\documentclass[12
pt,a4paper]{book}

\usepackage{amsmath}
\usepackage{amsfonts}
\usepackage{amsthm}
\usepackage{amssymb}
\usepackage{booktabs,amsmath}
\usepackage{fancyhdr}

\usepackage[T1]{fontenc}
\usepackage[utf8]{inputenc}
\usepackage[english]{babel}
\usepackage{emptypage}
\usepackage{imakeidx}
\makeindex[intoc]
\indexsetup{firstpagestyle=empty}
\usepackage[english]{varioref}
\usepackage{hyperref}
\usepackage{bookmark}

\def\Xint#1{\mathchoice
	{\XXint\displaystyle\textstyle{#1}}%
	{\XXint\textstyle\scriptstyle{#1}}%
	{\XXint\scriptstyle\scriptscriptstyle{#1}}%
	{\XXint\scriptscriptstyle\scriptscriptstyle{#1}}%
	\!\int}
\def\XXint#1#2#3{{\setbox0=\hbox{$#1{#2#3}{\int}$}
		\vcenter{\hbox{$#2#3$}}\kern-.5\wd0}}

\def\dashint{\Xint-}
\usepackage[english]{babel}
\usepackage{xcolor}
\usepackage{enumerate}
\usepackage{graphicx}
\usepackage{collectbox}

\newcommand{\mybox}{%
	\collectbox{%
		\setlength{\fboxsep}{4pt}%
		\fbox{\BOXCONTENT}%
	}%
}

\theoremstyle{plain}

\newtheorem{theo}{Theorem}[section]
\newtheorem{lem}[theo]{Lemma}
\newtheorem{prop}[theo]{Proposition}
\newtheorem{cor}[theo]{Corollary}

\theoremstyle{definition}
\newtheorem{definition}[theo]{Definition}
\theoremstyle{remark}

\newtheorem{rem}[theo]{Remark}

\numberwithin{equation}{section}

\numberwithin{figure}{chapter}

\begin{document}
		
\title{Notes on unique continuation properties for Partial Differential Equations -- Introduction to the stability estimates for inverse problems}
\author{Sergio Vessella \thanks{Dipartimento di Matematica e Informatica "U. Dini", Universit\`{a} di Firenze,
		viale Morgagni 67/A, 50134 Firenze (Italy), email: sergio.vessella@unifi.it}}
\date{}
\maketitle
\cleardoublepage
\thispagestyle{empty}
\vspace*{\stretch{1}}
\begin{flushright}
	\itshape Dedicated to my whole family and especially \\ to my dear wife Luisella and my dear son Luigi
\end{flushright}
\vspace{\stretch{2}}
\cleardoublepage

\tableofcontents 
\include{introduzione} 

\newpage
\textbf{Abstract.} These Notes are intended for graduate or undergraduate students who have familiarity with Lebesgue measure theory, partial differential equations, and functional analysis. The main topics covered in this work are the study of the Cauchy problem and unique continuation properties associated with partial differential equations. The primary objective is to familiarize students with stability estimates in inverse problems and quantitative estimates of unique continuation. The treatment is presented in a self-contained manner.

\bigskip 

\textbf{Mathematical Subject Classifications (2020) 35-01, 35-B60, 35-R25, 35-R30.} 

\medskip

\chapter*{Introduction}

 The main purpose of these Notes is to introduce the study of unique continuation properties and stability estimates for inverse problems for partial differential equations (PDEs). The topics covered in these Notes are all chosen with proximity to inverse problems in mind, but we believe that none of the topics should be neglected in the training of those interested in PDEs, especially with regard to the study of unique continuation properties. Despite the existence of excellent review articles and books on the subject, there is a lack of a truly introductory book starting from minimum basics for a graduate or undergraduate student, who should have some familiarity with Lebesgue measure theory, basic elements of functional analysis, and a first introductory course on PDEs.
 
 To facilitate the achievement of this purpose, we have covered other basic topics in the theory of PDEs, such as the theory of existence and regularity for second-order elliptic equations with real coefficients in Part I, and the classical theory of the Cauchy problem for equations with analytic coefficients in Part II. Part III focuses on the study of unique continuation properties for equations with non-analytic coefficients.
 
 We provide a brief description of the present Notes in the remainder of this Introduction, with more detailed descriptions of the topics covered in the individual chapters provided in their introductions.

 In Part I, in addition to Chapter \ref{prerequisiti} on Sobolev spaces, we have included Chapter \ref{preliminari-spazi-funz}, which serves as a connection and complement to elementary analysis topics. In Chapter \ref{preliminari-spazi-funz}, we recall the main definitions and theorems (without proof) of measure theory and prove some important theorems of real analysis, including the extension theorem in $C^{0,\alpha}$, the Lebesgue differentiation theorem, the Rademacher theorem, and the divergence theorem over open sets with Lipschitz boundary. Additionally, we study the distance function and the Hausdorff distance between compact sets, which is useful for studying the stability issue of inverse problems with unknown boundaries. In Chapter \ref{Lax-Milgram}, we provide the definition and first properties of the Dirichlet-to-Neumann map, in addition to the existence and regularity $L^2$ theory for second-order elliptic equations. We also introduce the inverse problem of inclusion detection and, in particular, size estimates, which have an interesting connection with the quantitative estimates of unique continuation developed in Part III. The books that inspired us the most in writing Part I are
 \cite{K-J-F}, \cite{Pu} (Chapter \ref{preliminari-spazi-funz}) and \cite{Br}, \cite{EV}, \cite{EV-GA} (Chapter \ref{prerequisiti}, Chapter \ref{Lax-Milgram}).    
 
 In Part II, we provide, in Chapter \ref{first-order-cauchy}, a concise discussion of the Cauchy problem for first-order PDEs. In Chapter \ref{funz-analitiche}, we have given the basic properties of real analytic functions which we need in Chapter \ref{Cauchy-chap}, where we give the formulation of the Cauchy problem for PDEs and prove the classical Cauchy-Kovalevskaya, Holmgren, and John theorems for (linear) PDEs with analytic coefficients. In Chapter \ref{PbI-froniera:22-11-22}, we apply the Holmgren and John theorems and the $L^2$ regularity theory to prove a uniqueness theorem for an inverse problem for the Laplace equation with unknown boundary. In Chapter \ref{LM}, we introduce the concept of a well-posed problem in the sense of Hadamard, and by means of the Lax-Mizohata Theorem, we highlight the important connection between uniqueness, solvability, and continuous dependence on the data in a Cauchy problem for equations with $C^{\infty}$ coefficients. In Chapter \ref{Stab-condizionata}, we give the definition of conditional stability (or ``well-posed problem in the sense of Tikhonov'') and some basic examples of conditional stability estimates for the calculus of derivatives and for the analytic continuation problem (in this area, the most famous theorem is the Hadamard three-circle inequality). Chapter \ref{Stab-condizionata} is a kind of ``laboratory'' in which we build some tools that should be kept in the toolbox of anyone who wants to study the conditional stability of not-well-posed problems in the sense of Hadamard. We conclude Part II with Chapter \ref{Holmgren-John}, in which we prove the John stability Theorem for the Cauchy problem for PDEs with analytic coefficients and discuss some of its consequences. The books that inspired us the most in writing Part II are 
\cite{Co-Hi}, \cite{EV}, \cite{Fo}, \cite{Joh}, \cite{Petr} (Chapter \ref{first-order-cauchy}, Chapter \ref{funz-analitiche}, Chapter \ref{Cauchy-chap}),  \cite[Ch. V]{HO63}, \cite[Vol. II]{HOII}, \cite{Mizo} (Chapter \ref{LM}), \cite{L-R-S} \cite{Tik} (Chapter \ref{Stab-condizionata}).
 
 In Part III, as we have already mentioned, we provide a gradual study of Carleman estimates and the main problems of unique continuation for PDEs. In Chapter \ref{Nirenberg}, we extensively explain the Nirenberg Theorem \cite{Ni} concerning the Cauchy problem for constant operators in the principal part. From an educational point of view, one of the merits of this theorem consists of its simple proof and, conversely, in the powerful consequences that allow us to solve the question of the uniqueness of solutions to the Cauchy problem for the equation
 
 \begin{equation*}
 	\Delta u=b(x)\cdot\nabla u+c(x)u,
 \end{equation*}
 where $b=(b_1,\cdots,b_n)\in L^{\infty}(\mathbb{R}^n;\mathbb{R}^n)$, $c\in
 L^{\infty}(\mathbb{R}^n)$.
 Furthermore, the Nirenberg Theorem allows for addressing some standard aspects involving Carleman estimates quite easily, especially with regard to how the aforementioned estimates are used to infer the unique continuation property for PDEs. The actual presentation of the Carleman estimates is carried out in chapters \ref{Carleman}, \ref{operatori-2ord}, and \ref{tre sfere-ellittiche}.   In Chapter \ref{Carleman}, we follow, with slight simplifications, the general and now classic approach developed by H\"{o}rmander \cite{HO63}, as it allows for a broad and general view of the issues concerning Carleman estimates.  The main theorem of Chapter \ref{Carleman} is the Carleman estimate for general elliptic operators, which corresponds to theorem 8.3.1 of \cite[Ch. VIII]{HO63}. Chapter \ref{Carleman} has its natural continuation in Chapter \ref{operatori-2ord}, in which we initially review the proofs of Chapter \ref{Carleman} in the simple case of the Laplace operator. From there, we move on to dealing with second-order operators that are not necessarily elliptic. Unlike Chapter \ref{Carleman}, where the integration by parts used to arrive at a Carleman estimate is based on a careful study of the quadratic differential forms, in Chapter \ref{operatori-2ord}, we adopt the Rellich identity and its natural generalization. This approach makes it easy to handle the case of second-order operators whose principal part has real Lipschitz continuous coefficients by providing a ``miniaturized'' proof of the uniqueness Calder'{o}n Theorem (here, Theorem \ref{oper-2ord-341}) for operators with simple characteristics. In Chapter \ref{operatori-2ord}, we provide a hint to the notion of pseudoconvex functions, which is particularly simplified in the case of second-order operators with real coefficients. Finally, in Chapter \ref{tre sfere-ellittiche}, we prove some Carleman estimates with a singular weight for the second-order elliptic operator 
 
 \begin{equation*}  \label{14-equ-fin}
 	Lu=\sum_{i,j=1}^n\partial_{x_i}\left(a^{ij}(x)\partial_{x^j}u\right),
 \end{equation*}
where  $\left\{a^{ij}(x)\right\}^n_{i,j=1}$ is a symmetric matrix whose entries are
real-valued Lipschitz-continuous functions. We use the Carleman estimates to deduce the optimal three sphere inequality, the doubling inequality, and the strong unique continuation property (corollaries \ref{SUCP} and \ref{SUCP-fin}) for the equation
$$Lu=b(x)\cdot\nabla u+c(x)u,$$ where $b\in L^{\infty}(\mathbb{R}^n;\mathbb{R}^n)$, $c\in
L^{\infty}(\mathbb{R}^n)$.
Of the various proofs in the literature for such Carleman estimates, we present the proof given in \cite{Ar}, \cite{AKS} which is based on transforming the elliptic operator into polar coordinates (Euclidean or Riemannian). We consider this elegant proof useful because it allows us to discuss the transformation into polar coordinates with respect to a Riemannian metric, which can be useful in other contexts of PDEs. In Chapter \ref{Misc:27-11-22}, we provide some brief and simple comments on the methods of log-convexity and the frequency function for studying the unique continuation property. In this chapter, we also mention some simple applications of $A_p$ weights in the stability and size estimates, and we conclude with the Runge property for the Laplace operator. The books that inspired us the most in writing Part III are \cite{HO63} and \cite{Lern}. 

I would like to conclude these Notes by thanking all those who provided me with useful advice on how to carry on this work, especially my friends Lorenzo Baldassari and Elisa Francini.

\part{THE SOBOLEV SPACES AND THE BOUNDARY VALUE PROBLEMS}\label{parte1}
\chapter{Main notation and basic formulas}\label{Notazioni richiami}
\section[Notation]{Notation} \label{Notazioni}

Let us denote  by $\mathbb{N}_0=\mathbb{N}\cup\{0\}$. We call multi-index \index{multi-index} 
any $n$-uple of elements of $\mathbb{N}_0$
$$\alpha=\left(\alpha_1,\alpha_2, \cdots,\alpha_n\right), \quad \alpha_j\in\mathbb{N}_0,\quad j=1,2,\cdots,n.$$
For any $\alpha \in \mathbb{N}^n_0$ we denote by 
$$\left\vert\alpha\right\vert=\left\vert\alpha_1\right\vert+\left\vert\alpha_2\right\vert \cdots,\left\vert\alpha_n\right\vert,
\quad\mbox{and}\quad \alpha!=\alpha_1!\alpha_2! \cdots\alpha_n!\mbox{ },$$ the lenght (modulus) and the factorial of $\alpha$, respectively. For any
$x\in\mathbb{R}^n$, $x=\left(x_1,x_2, \cdots,x_n\right)$, we set
$$x^{\alpha}=x_1^{\alpha_1}x_2^{\alpha_2}\cdots x_n^{\alpha_n}.$$ Let
$\alpha,\beta\in\mathbb{N}^n_0$ we write $\alpha\leq\beta$ provided
$\alpha_j\leq\beta_j$ for $j=1,2,\cdots,n$ and we write
$\alpha<\beta$ provided $\alpha\leq\beta$ and there exists
$j_0\in\{1,2,\cdots,n\}$ such that $\alpha_{j_0}<\beta_{j_0}$. For any
$x=(x_1,\cdots,x_n)\in \mathbb{R}^n$ we denote, unless otherwise stated, by $x'=(x_1,\cdots,x_{n-1})\in \mathbb{R}^{n-1}$
and we write $x=(x', x_n)$. Similar convention will be used for the multi-indices.

For any $\alpha,\beta\in\mathbb{N}_0$ and $\alpha\leq\beta$, let us denote by (the binomial "$\beta$ over $\alpha$") \index{binomial of multi--indices}
$$\binom{{\beta}}{{\alpha}}=\frac{\beta!}{\left(\beta-\alpha\right)!\alpha!}
\mbox{ }..$$ Let us denote by $\partial_k$ \index{$\partial_k$}the operator 
$\frac{\partial}{\partial x_k}$, $k=1,2,\cdots,n$ and by

$$\partial=\left(\partial_1,\partial_2, \cdots,\partial_n\right)\quad \mbox{(the gradient operator)}.$$
Hence, we set 

$$\partial^{\alpha}=\partial_1^{\alpha_1}\partial_2^{\alpha_2}\cdots\partial_n^{\alpha_n}=
\frac{\partial^{|\alpha|}}{\partial
x_1^{\alpha_1}\partial x_2^{\alpha_2}\cdots\partial x_n^{\alpha_n}}.$$
To denote the gradient operator we also use the notation $\nabla$,
however to denote $\partial^{\alpha}$ we \textbf{will not write} $\nabla^{\alpha}$. Of course, we will continue to denote by $u_{x_k}$ (or by other standard symbols) the partial derivative of $u$ with respect to $x_k$, $k=1,2,\cdots,n$. The Hessian matrix of a smooth function $u$ is denoted by \index{$\partial^2u$}
$$\partial^2u=\left\{\partial^2_{jk}u\right\}_{j,k=1}^n$$

We point out that some authors (and also in these notes in some
context) reserve the notation $D_k$ \index{$D_k$}to denote the operator
$\frac{1}{i}\frac{\partial}{\partial x_k}$, where $i=\sqrt{-1}$, consequently
$$D^{\alpha}=\left(\frac{1}{i}\right)^{|\alpha|}\partial_1^{\alpha_1}\partial_2^{\alpha_2}\cdots\partial_n^{\alpha_n}.$$
The latter notation is useful especially when an 
extensive use of the Fourier transform \index{Fourier transform}is done
$$\widehat{u}\left(\xi\right)=\int_{\mathbb{R}^n}u(x)e^{-ix\cdot\xi}dx.$$
Actually, we have
$$\widehat{D^{\alpha}u}\left(\xi\right)=\xi^{\alpha}\widehat{u}\left(\xi\right).$$
while, using the former notation, we have
$$\widehat{\partial^{\alpha}u}\left(\xi\right)=(i\xi)^{\alpha}\widehat{u}\left(\xi\right).$$

Let $\ell\in\mathbb{R}^n\setminus \{0\}$, for any
$j\in\mathbb{N}_0$ we set \index{$\frac{\partial^j}{\partial \ell^j}$}
$$\frac{\partial^j}{\partial \ell^j}=\sum_{|\alpha|=j}\ell^{\alpha}\partial^{\alpha}, $$
(we mean $\frac{\partial^0u}{\partial \ell^0}=u$). In particular
$$\frac{\partial}{\partial \ell}=\ell\cdot \partial= \ell\cdot\nabla.$$
As a consequence of the notations introduced above, we denote a polynomial $P$ of degree $m$  in the variables $\xi_1,\xi_2\cdots,\xi_n$

\begin{equation}\label{corr-10-3-23-1}
P(\xi)=\sum_{|\alpha|\leq m}a_{\alpha}\xi^{\alpha},
\end{equation}
where $a_{\alpha}\in \mathbb{R}$ (or $a_{\alpha}\in \mathbb{C}$) for any
$|\alpha|\leq m$. We say that the homogeneous polynomial  $$P_m(\xi)=\sum_{|\alpha|= m}a_{\alpha}\xi^{\alpha},$$ is the \textbf{principal part of a polynomial} \index{principal part:@{principal part:}!- of a polynomial@{- of a polynomial}} $P$  provided that there exists $\alpha_0\in\mathbb{N}_0^n$
such that $|\alpha_0|=m$ and $a_{\alpha_0}\neq 0$.

\section{Some useful formulas}\label{formule:1-11-22}
In this Section we recall some basic and useful formulas.

\begin{equation}\label{1-2N}
(x+y)^{\alpha}=\sum_{\beta\leq
\alpha}\binom{{\alpha}}{{\beta}}x^{\beta}y^{\alpha-\beta}, \quad
\forall x,y\in \mathbb{R}^n,\quad \forall\alpha\in \mathbb{N}_0^n.
\end{equation}
\begin{equation}\label{5-3N}
\partial^{\beta}x^{\alpha}=
\begin{cases}
\frac{\alpha!}{(\alpha-\beta)!}x^{\alpha-\beta}, \quad \mbox{ for } \alpha\geq\beta,\\
\\
0, \quad \quad\mbox{otherwise}. %
\end{cases}%
\end{equation}
\begin{equation}\label{3-3N}
(x_1+x_2+\cdots+x_n)^m=\sum_{|\alpha|=m}\frac{m!}{\alpha!}x^{\alpha},\quad
\forall m\in\mathbb{N}_0.
\end{equation}
\begin{equation}\label{4-3N}
\alpha!\leq|\alpha|!\leq n^{|\alpha|}\alpha!,\quad \forall\alpha\in
\mathbb{N}_0^n.
\end{equation}
Let us recall the following \textbf{Stirling formula} \index{Stirling formula}
\begin{equation}\label{Stirling-3N}
\lim_{n\rightarrow\infty}\frac{n!}{n^ne^{-n}\sqrt{n}}=\sqrt{2\pi}\ \ .
\end{equation}
Let $f$ be a smooth function and $m\in
\mathbb{N}_0$, we have
\begin{equation}\label{6-3N}
\frac{d^m}{dt^m}f(x+ty)=\left((\sum_{j=1}^n y_j\partial_j)^m
f\right)(x+ty)=
\sum_{|\alpha|=m}\frac{m!}{\alpha!}y^{\alpha}\left(\partial^{\alpha}
f\right)(x+ty).
\end{equation}
We recall the \textbf{Taylor formula}, \index{Taylor formula} centered at
$x_0\in\mathbb{R}^n$, of a polynomial $P$ of degree $m$   
\begin{equation}\label{2-2N}
P(x)=\sum_{|\alpha|\leq
m}\frac{1}{\alpha!}\partial^{\alpha}P(x_0)(x-x_0)^{\alpha}.
\end{equation}
Let $f$ and $g$ be two smooth functions and 
$\alpha\in \mathbb{N}_0^n$, we have the \textbf{Leibniz formula} for the $\alpha$-th derivative of the product $fg$
\begin{equation}\label{7-2N}
\partial^{\alpha}(fg)=\sum_{\beta\leq\alpha}\binom{{\alpha}}{{\beta}}\partial^{\beta}f\partial^{\alpha-\beta}g.
\end{equation}

\bigskip

Now we \textbf{check some formula.}

\medskip

Formula \eqref{1-2N} easily follows by the Newton binomial formula. Actually, we have 

\begin{equation*}
\begin{aligned}
(x+y)^{\alpha}&=(x_1+y_1)^{\alpha_1}\cdots (x_n+y_n)^{\alpha_n}=\\&
=\sum_{\beta_1\leq
\alpha_1}\binom{{\alpha_1}}{{\beta_1}}x_1^{\beta_1}y_1^{\alpha_1-\beta_1}
\cdots\sum_{\beta_n\leq
\alpha_n}\binom{{\alpha_n}}{{\beta_n}}x_n^{\beta_n}y^{\alpha_n-\beta_n}=\\&
=\sum_{\beta\leq
\alpha}\binom{{\alpha}}{{\beta}}x^{\beta}y^{\alpha-\beta}.
\end{aligned}
\end{equation*}

The proof of \eqref{5-3N} is immediate. Before checking 
\eqref{3-3N} let us notice that if $\alpha$ and $\beta$ are 
multi-indices such that $\alpha\leq\beta$ and $|\alpha|=|\beta|$, then
$\alpha=\beta$. Let us denote
$$S(x)=\sum_{j=1}^nx_j$$ and set
$$P(x)=\left(S(x)\right)^m.$$ Since $P$ is a homogeneous polynomial of degree $m$ we get
$$P(x)=\sum_{|\alpha|=m}c_{\alpha}x^{\alpha}.$$ Let us show that 
$c_{\alpha}=\frac{m!}{\alpha!}$ for every multi-indices $\alpha$
such that $|\alpha|=m$. Let $\beta$ be a multi-index  satisfying $|\beta|=m$. By what we notice above and by \eqref{5-3N}, we get

\begin{equation}\label{p1-N}
\partial^{\beta}P(x)=\sum_{|\alpha|=m}c_{\alpha}\partial^{\beta}x^{\alpha}=c_{\beta}\beta! \ .
\end{equation}
On the other hand,

\begin{equation}\label{p2-N}
	\begin{aligned}
		    &\partial^{\beta}P(x)=\partial^{\beta}S^m= \partial_1^{\beta_1}\cdots\partial_n^{\beta_n}S^m=\\&
    =m\cdots (m-\beta_n+1)\partial_1^{\beta_1}\cdots\partial_{n-1}^{\beta_{n-1}}  S^{m-\beta_n}=\\&
    =m\cdots  (m-\beta_n-\beta_{n-1}+1)\partial_1^{\beta_1}\cdots\partial_{n-2}^{\beta_{n-2}}  S^{m-\beta_n-\beta_{n-1}}=
    \cdots=m! \ .
\end{aligned}
\end{equation}
By \eqref{p1-N} and \eqref{p2-N} we get $$c_{\beta}\beta!=m!
\quad\mbox{ for every } \beta\in\mathbb{N}_0^n \mbox{ such that }
|\beta|=m,$$ from which we obtain \eqref{3-3N}. Of course, formula
\eqref{3-3N} can be proved more elementarly.
For instance, it can be proved by induction starting from the Newton binomial formula.

Concerning the inequality $\alpha!\leq|\alpha|!$ in 
\eqref{4-3N}, recalling that $h!k!\leq(h+k)!$ for every
$h,k\in\mathbb{N}_0$ we get
$$\alpha!=\alpha_1!\alpha_2!\cdots\alpha_n!\leq (\alpha_1+\alpha_2)!\alpha_3!\cdots\alpha_n!\leq\cdots \leq(\alpha_1+\alpha_2+\cdots+\alpha_n)!=|\alpha|! \ .$$
Regarding the inequality $|\alpha|!\leq
n^{|\alpha|}\alpha!$ it suffices to use formula \eqref{3-3N} and we have
$$n^{|\alpha|}=\underset{n }{(\underbrace{1+1+\cdots+1})}^{|\alpha|}=\sum_{|\beta|=|\alpha|}\frac{|\alpha|!}{\beta!}\geq
\frac{|\alpha|!}{\alpha!}. $$

The first equality in \eqref{6-3N} can be obtained by iterating the formula
$$\frac{d}{dt}f(x+ty)=\left(\sum_{j=1}^n y_j\partial_j f\right)(x+ty).$$ 
The second equality in \eqref{6-3N} can be obtained by a formal development of $$(\sum_{j=1}^n
y_j\partial_j)^m$$ through \eqref{3-3N}. 

Leibniz formula \eqref{7-2N} can be easily obtained by the namesake formula
for the one variable functions
$$\frac{d^k}{dt^k}(fg)=\sum_{h=0}^k \binom{{k}}{{h}}\frac{d^hf}{dt^h}\frac{d^{k-h}g}{dt^{k-h}},$$where $f$ and $g$ are two smooth functions in the variable $t$.

\bigskip

Let $m\in\mathbb{N}_0$ and let $a_{\alpha}$, $|\alpha|\leq m$, be some functions
defined in an open set $\Omega\subset\mathbb{R}^n$ with values
in $\mathbb{R}$ or in $\mathbb{C}$. We say that the operator

\begin{equation}\label{1-4N}
P(x,\partial)=\sum_{|\alpha|\leq m}a_{\alpha}(x)\partial^{\alpha}
\end{equation}
is a \textbf{linear differential operator of order} $m$ \index{order of a linear differential operator}
in $\Omega$, provided that there exists $\alpha\in \mathbb{N}_0^n$, $|\alpha|=m$ such that
$a_{\alpha}$ does not vanish identically in $\Omega$. We say that the
functions $a_{\alpha}$, $|\alpha|\leq m$, are the
\textbf{coefficients} of the differential operator \eqref{1-4N}.

We define the \textbf{symbol of operator} \eqref{1-4N} \index{symbol of an operator} as the following polynomial in the variable $\xi$
\begin{equation}\label{simbolo-4N}
P(x,\xi)=\sum_{|\alpha|\leq m}a_{\alpha}(x)(i\xi)^{\alpha}.
\end{equation}
Notice that if we write $P(x,\partial)$ as
$$\widetilde{P}(x,D)=\sum_{|\alpha|\leq m}a_{\alpha}(x)(i D)^{\alpha},$$ then the symbol \eqref{simbolo-4N} can be obtained by formally substituting $D$ to $\xi$.

We have

\begin{equation}\label{3ex-4N}
e^{-ix\cdot\xi}\widetilde{P}(x,D)e^{ix\cdot\xi}=P(x,\xi).
\end{equation}

It might seem more natural to define the symbol of
\eqref{1-4N} as simply $\sum_{|\alpha|\leq m}a_{\alpha}(x)\xi^{\alpha}$, actually the contexts in which it is
mostly used the definition of the symbol of a differential operator, are
often the same ones in which it is convenient to use $D_j=\frac{1}{i}\partial_j$ as derivative operator.
It is therefore advisable to stick to the standard definition of symbol for
do not stray from the current literature. For instance, the symbol of the 
\textbf{Laplace operator} \index{Laplace operator}
$$\Delta=\sum_{j=1}^n\partial_j^2=-\sum_{j=1}^n D_j^2$$
is given by
$$-\sum_{j=1}^n \xi_j^2,$$ the symbol of the \textbf{heat operator} \index{heat operator}

$$\sum_{j=1}^n\partial_j^2-\partial_t=-\sum_{j=1}^n D_j^2-iD_t,$$ is equal to
$$-\sum_{j=1}^n \xi_j^2-i\xi_{n+1}$$
and the symbol of the \textbf{wave operator} or \textbf{d'Alembertian operator} \index{d'Alembertian operator (wave operator)}

$$\Box = \Delta-\partial^2_t=\sum_{j=1}^n\partial_j^2-\partial^2_t=-\sum_{j=1}^n D_j^2+D^2_t$$
is equal to
$$-\sum_{j=1}^n \xi_j^2+\xi^2_0.$$ 

We will call the \textbf{principal part of operator} \index{principal part:@{principal part:}!- of a linear differential operator@{- of a linear operator}}  \eqref{1-4N},
the differential operator

\begin{equation}\label{2-4N}
P_m(x,\partial)=\sum_{|\alpha|=m}a_{\alpha}(x)\partial^{\alpha}.
\end{equation}

In the sequel, to simplify the notations, we will concentrate
our attention to the case in which the coefficients $a_{\alpha}$ are
real--valued functions. However, we warn that what we will establish, in
many cases, can easily be extended to the case where the
coefficients $a_{\alpha}$ is a complex--valued function.

If all the coefficients of the operator $P(x,\partial)$ are constants,
we will say that $P(x,\partial)$ is an \textbf{operator with 
constant coefficients}. In these cases, to
denote the operator  $P(x,\partial)$, we will just write $P(\partial)$.

We notice that, by the above definition,  the
symbol of the principal part of operator \eqref{1-4N} is the homogeneous polynomial

\begin{equation}\label{5-4N}
P_m(x,\xi)=i^m\sum_{|\alpha|= m}a_{\alpha}(x)\xi^{\alpha}.
\end{equation}

\bigskip

\underline{\textbf{Conventions on the constants.}} In the sequel, to denote a positive constant we will use very often
the letter $C$. We notice right now that the value of the constants will may
change from line to line, but we will generally indicate the dependence of the constants by the various
parameters. However, sometimes to be able to
better follow the various steps, we will put an index or a sign to $C$
and we will write $C_0, C_1, \overline{C}, \widetilde{C} \ldots$.
We will generally omit the dependence of the various constants on the
dimension of the space.

\newpage

\chapter{Review of some function spaces and measure theory}\label{preliminari-spazi-funz}

\section{The space $C^k$}
\label{Funz:Cont-diff} 

Let $X$ be a subset of $\mathbb{R}^n$. We will denote by $C^0(X)$ the vector space of continuous functions defined in $X$ with values in $\mathbb{R }$. If $u\in C^0(X)$, \index{$ C^0(X)$}we denote the \textbf{support of} \index{support of a continuous function} $u$ by 
$$\mbox{supp }u:=\overline{\left\{x\in X:\mbox{ } u(x)\neq 0\right\}} \quad \mbox{(closure in } \mathbb{R}^n).$$ 
 We will denote by $C_0^0(X)$ the space of continuous functions whose support is a compact set of $\mathbb{R}^n$ contained in $X$.

\begin{prop}\label{Contin:9-1}
	
	Let $u\in C^0(X)$. Let $K\subset X$ be a compact set of $\mathbb{R}^n$; for any $r>0$ we set
	$$K_r=\left\{x\in \mathbb{R}^n :\mbox{ } \mbox{dist} (x,K)\leq r \right\},$$ where dist($x,K)$ denotes the distance of $x$ from $K$. Let us suppose that there exists $r_0>0$ such that $K_{r_0}\subset X$. Then
	\begin{equation}\label{Contin:9-2}
	\lim_{r\rightarrow 0} \max_{K_r} u= \max_{K} u \ .
	\end{equation}
	\end{prop}
\textbf{Proof.} Since $K_{r_0}$ is a compact subset of $X$ and $u\in C^0(X)$, $u$ is uniformly continuous on $K_{r_0}$. Let $\varepsilon>0$ and $\delta>0$ such that

$$|u(x)-u(y)|<\varepsilon,\quad \mbox{for all } x,y\in K_{r_0}\mbox{ such that } |x-y|\leq\delta;$$ we may assume that $\delta<r_0$. Let $r\in (0,\delta]$ and let $x$ be any point of $K_r$. Hence there exists $y\in K$ such that $|x-y|\leq r$. Therefore
$$u(x)<u(y)+\varepsilon\leq \max_{K} u+\varepsilon.$$ By the arbitrariness of $x$ in $K_r$ we have
$$ \max_{K} u\leq \max_{K_r} u<\max_{K} u+\varepsilon,$$ which concludes the proof of \eqref{Contin:9-2}. $\blacksquare$ 

\bigskip 

Let $u\in C^0(X)$. We define the \textbf{modulus of continuity} \index{modulus of continuity} of $u$ in $X$

\begin{equation}\label{Contin:9-3}
	\omega(\delta)=\sup \left\{|u(x)-u(y)| :\mbox{ } x,y\in X,\mbox{ } |x-y|\leq \delta \right\},\quad \mbox{for }\delta>0.
\end{equation}
$\omega$ is an increasing function, defined on $[0,+\infty)$ and satisfies $\omega(0)=0$. Of course, $\omega$ may not be finite. It is easy to check that $u$ is uniformly continuous if and only if  

\begin{equation}\label{Contin:9-5-0}
	\lim_{\delta\rightarrow 0}\omega(\delta)=0.
\end{equation}

 If the function $\omega$ is bounded, it may be convenient to use the \textbf{concave modulus of continuity} which is defined as 

\begin{equation}\label{Contin:9-4}
\tilde{\omega}(\delta)=\inf \left\{f(\delta) :\mbox{ } f \mbox{ concave, } f\geq \omega, \mbox{ in } [0,+\infty) \right\},\quad \mbox{for }\delta>0.
\end{equation}

\medskip

 Now we check that  
 
 \begin{equation}\label{Contin:9-5}
 \lim_{\delta\rightarrow 0}\tilde{\omega}(\delta)=0.
 \end{equation}
 Let us denote
$$M=\sup_{\delta\in [0,+\infty)}\omega(\delta)<+\infty.$$ 
If $M=0$, \eqref{Contin:9-5} is trivial. Let us suppose therefore $M>0$. Let $0<\varepsilon<M$, from \eqref{Contin:9-5-0} it follows that there exists $\delta_0>0$ such that

\begin{equation*}
	0\leq \omega(\delta)<\frac{\varepsilon}{2},\quad \forall\delta\in [0,\delta_0].
\end{equation*}
Set 
\begin{equation*}
	g_{\varepsilon}(\delta)=
	\begin{cases}
		\frac{\varepsilon}{2}+\frac{2M-\varepsilon}{2\delta_0}\delta, \quad \mbox{ for } \delta \in [0,\delta_0],\\
		\\
		M, \quad \quad\mbox{for } \delta\in (\delta_0,+\infty). 
	\end{cases}%
\end{equation*}
It easy to check that $g_{\varepsilon}$ is concave and that $g_{\varepsilon}\geq \omega$ in $[0,+\infty)$. Furthermore, set
$$\delta_1=\frac{\varepsilon\delta_0}{2M-\varepsilon},$$ it turns out

$$g_{\varepsilon}(\delta)<\varepsilon, \quad\forall \delta \in [0,\delta_1).$$ Therefore
   
 $$\tilde{\omega}(\delta)<\varepsilon, \quad\forall \delta \in [0,\delta_1),$$ which gives  \eqref{Contin:9-5}. 

\bigskip

\textbf{Remark 1.} Let $\tilde{\omega}$ be a concave modulus of continuity, then

\begin{equation}\label{Contin:9-6}
	0<\eta_1<\eta_2 \ \Longrightarrow  \ \eta_1 \tilde{\omega}\left(\frac{1}{\eta_1}\right)\leq \eta_2 \tilde{\omega}\left(\frac{1}{\eta_2}\right).
\end{equation}
Let us check \eqref{Contin:9-6}. From the concavity of $\tilde{\omega}$ and recalling that $\tilde{\omega}(0)=0$ we have, for $0<\eta_1<\eta_2$,

\begin{equation*}
	  \eta_1 \tilde{\omega}\left(\frac{1}{\eta_1}\right)=\frac{\tilde{\omega}\left(\frac{1}{\eta_1}\right)-\tilde{\omega}(0)}{\frac{1}{\eta_1}-0}\leq \frac{\tilde{\omega}\left(\frac{1}{\eta_2}\right)-\tilde{\omega}(0)}{\frac{1}{\eta_2}-0}=\eta_2 \tilde{\omega}\left(\frac{1}{\eta_2}\right).
\end{equation*}
$\blacklozenge$

\bigskip

Let us denote by $C_{*}^0(X)$ the space of of the bounded functions of $C^0(X)$, let us define the norm

\begin{equation}\label{Contin:1}
	\left\Vert u\right\Vert_{C_{*}^0(X)}=\sup_{X}|u(x)|,\quad \forall u\in C_{*}^0(X).
\end{equation}
As it is well--known, the space $C_{*}^0(X)$ equipped with the norm \eqref{Contin:1}, is a \textbf{Banach space}.

\bigskip

The following Theorem holds true (see \cite[Corollary 1.3, Ch. 3]{Lang} for a proof)

\begin{theo}[\textbf{Weierstrass approximation}]\label{teor-approssimaz-weierst}
	\index{Theorem:@{Theorem:}!- Weierstrass approximation@{- Weierstrass approximation}} Let $X$ be a compact subset of $\mathbb{R}^n$. For every $u\in C^0(X)$ and for every $\varepsilon>0$ there exists a polynomial $P$ such that
	\begin{equation*}
		\left\Vert u-P\right\Vert_{C^0\left(X\right)}<\varepsilon.
		\end{equation*}
	\end{theo}
We notice that by approximating the polynomial $P$ given in the previous Theorem, by a polynomial with rational coefficients, we derive that $C^0(X)$, with compact $X$, is a \textbf{separable space} . Let us recall that a topological space $\mathcal{S}$ is said to be separable \index{separable topolical space} if there exists a countable set $D\subset \mathcal{S}$ such that $\overline{D}=\mathcal{S}$.

\medskip

The following Proposition holds true

\begin{prop}\label{Contin:11-10-1}
Let $\mathcal{S}$ be a metric space with distance $d$. If there exists $\mathcal{Y}$, uncountable subset of $\mathcal{S}$, and $\delta>0$ such that
\begin{equation}\label{Contin:11-10-2}
	d(x,y)>\delta,\quad \ \ \forall x,y\in \mathcal{Y},\ \ x\neq y, 
\end{equation}
then $\mathcal{S}$ is not a separable space.
\end{prop}
\textbf{Proof.} We argue by contradiction and we assume that $\mathcal{S}$ is separable.
Hence there exists $D=\left\{u_n\right\}_{n\in \mathbb{N}}$ such that $\overline{D}=\mathcal{S}$. Consequently, for every $x\in \mathcal{Y}$ there exists $u_{n_x}$ such that
\begin{equation*}
d\left(x,u_{n_x}\right)<\frac{\delta}{3} \ .
\end{equation*}
Therefore, if $x,y\in \mathcal{Y}$, $x\neq y$, the triangle inequality gives  

\begin{equation*}
	d\left(u_{n_x},u_{n_y}\right)\geq d(x,y)- d\left(x,u_{n_x}\right)-d\left(y,u_{n_y}\right)>\frac{\delta}{3} \ .
\end{equation*}
In particular, if $x,y\in \mathcal{Y}$, $x\neq y$, then $u_{n_x}\neq u_{n_y}$. Consequently, the map
$$\mathcal{Y}\ni x\rightarrow u_{n_x}\in D,$$ is injective, but this fact contradicts that $\mathcal{Y}$ is an uncountable set. Therefore  $\mathcal{S}$ is not separable. $\blacksquare$

\medskip

\textbf{Remark 2.} Let us note that the compactness assumption of $X$ cannot be dropped for $C_{*}^0(X)$ to be separable. We show, for instance, that $C_{*}^0(\mathbb{R})$ is not separable. 

For any $A\in \mathcal{P}(\mathbb{Z})\setminus\left\{\emptyset\right\}$ (where $\mathcal{P}(\mathbb{Z})$ denotes the power set of $\mathbb{Z}$) and any $\varepsilon\in\left(0,\frac{1}{2}\right)$, define

\begin{equation*}
	u_A=\sum_{g\in A}u_g,
\end{equation*}
where
\begin{equation*}
	u_g(t)=
	\begin{cases}
		1-\varepsilon^{-1}|t-g|, \quad \mbox{ for } t\in [g-\varepsilon,g+\varepsilon],\\
		\\
		0, \quad \quad\mbox{for } t\in \mathbb{R}\setminus [g-\varepsilon,g+\varepsilon]. %
	\end{cases}%
\end{equation*}
We have 
\begin{equation*}
	\left\Vert u_A-u_B\right\Vert_{C_{*}^0\left( \mathbb{R}\right)}=
1,\quad \forall A, B\in \mathcal{P}(\mathbb{Z})\setminus\left\{\emptyset\right\},\ \ A\neq B.
\end{equation*}
Since $\mathcal{P}(\mathbb{Z})\setminus\left\{\emptyset\right\}$  is uncountable, Proposition \ref{Contin:11-10-1} implies that $C_{*}^0\left( \mathbb{R}\right)$ is not separable. $\spadesuit$

\medskip

Generally we will be interested in the case when $X$ is an open or the closure of an open set $\Omega$ of $\mathbb{R}^n$. If $\Omega$ is a bounded open set then we may consider $$C^0\left(\overline{\Omega}\right)$$ as a subspace of $C_{*}^0(\Omega)$ and we will denote the norm of $C^0\left(\overline{\Omega}\right)$ by 

\begin{equation}\label{Contin:2}
	\left\Vert u\right\Vert_{C^0\left(\overline{\Omega}\right)}=
	\sup_{\overline{\Omega}}|u|,\quad \forall u\in C^0\left(\overline{\Omega}\right).
\end{equation}
$\blacklozenge$

\bigskip

In the sequel we will use the following classical theorems on relatively compact sets.

\begin{theo}\label{Sob:teopag.92}
	Let $(X,d)$ be a complete metric space and let $Y\subset X$. Then $Y$ is a
	relatively compact set (i.e., $\overline{Y}=X$) if and only if it is
	\textbf{totally bounded} \index{totally bounded set} that is, for every $\delta>0$
	there exists a finite set  $\left\{x_1,\cdots, x_N\right\}\subset X$ such that
	$$d\left(y,x_j\right)<\delta, \quad \forall y\in Y, \ \ j=1,\cdots N$$ or, equivalently,
	$$Y\subset \bigcup_{j=1}^N B_{\delta}(x_j).$$
\end{theo}

\bigskip

\begin{theo}[\textbf{Arzel\`{a}--Ascoli}]\label{Sob:teoAscoli-Arzel}
	\index{Theorem:@{Theorem:}!- Arzel\`{a}--Ascoli@{- Arzel\`{a}--Ascoli}}Let $\Omega$ a bounded open set of $\mathbb{R}^n$ and let
	$\left\{u_k\right\}$ be a a sequence of functions belonging \index{equibounded sequence of functions} to
	$C^0\left(\overline{\Omega}\right)$ such that:
	
	\smallskip
	
	(i) $\left\{u_k\right\}$ is equibounded, i.e., there exists $M>0$
	such that
	$$\left\Vert u_k\right\Vert_{C^0\left(\overline{\Omega}\right)}\leq M, \quad\mbox{ } \forall k\in \mathbb{N};$$
	
	\smallskip
	
	(ii) $\left\{u_k\right\}$ is equicontinuous, i.e., for every \index{equicontinuous of functions}
	$\eta>0$ there exists $\delta>0$ such that if $|x-y|<\delta$  and $x,y\in
	\overline{\Omega}$, then
	
	$$\left |u_k(x)-u_k(y)\right|\leq \eta, \quad \forall k\in \mathbb{N}.$$
	
	\smallskip
	
	\noindent Then there exists a subsequence $\left\{u_{k_j}\right\}$ of $\left\{u_k\right\}$ and a function $u\in C^0\left(\overline{\Omega}\right)$ such that
 $$ \lim_{j\rightarrow\infty}\left\Vert u_{k_j}-u\right\Vert_{C^0\left(\overline{\Omega}\right)}=0.$$
\end{theo}

\bigskip

\bigskip

Let $k\in \mathbb{N}$ and let $\Omega$ be an open set of $\mathbb{R}^n$, we will denote by $C^k(\Omega)$ \index{$C^k\left(\Omega\right)$, $C^k\left(\overline{\Omega}\right)$} the space of functions which satisfy $\partial^{\alpha}u\in C^0(\Omega)$ for every $\alpha\in \mathbb{N}^n_0$,  $|\alpha|\leq k$. Further, we will denote by $C^k\left(\overline{\Omega}\right)$ the space of the functions $u\in C^k(\Omega)$ such that, for every $\alpha\in \mathbb{N}^n_0$, $|\alpha|\leq k$, $\partial^{\alpha}u$ is extensible to a function $U_{\alpha} \in C^0\left(\overline{\Omega}\right)$. Of course, if such an extension exists it is unique and we will write $\partial^{\alpha}u$ instead of $U_{\alpha}$.  If $\Omega$ is a bounded \textbf{open set} of $\mathbb{R}^n$, we define the norm on $C^k\left(\overline{\Omega}\right)$ as follows
\begin{equation}\label{Contin:3}
	\left\Vert u\right\Vert_{C^k\left(\overline{\Omega}\right)}=\sum_{|\alpha|\leq k}
	\sup_{\overline{\Omega}}|\partial^{\alpha}u|,\quad \forall u\in C^k\left(\overline{\Omega}\right).
\end{equation}
As it is well--known, the space $C^k\left(\overline{\Omega}\right)$, equipped with the norm \eqref{Contin:3} is a \textbf{Banach space}.  In some contexts it turns out to be convenient to consider, instead of norm \eqref{Contin:3}, an equivalent "dimensionless" norm, e.g.
\begin{equation}\label{Contin:4}
	\sum_{|\alpha|\leq k}
	d_0^{|\alpha|}\sup_{\overline{\Omega}}|\partial^{\alpha}u|,\quad \forall u\in C^k\left(\overline{\Omega}\right),
\end{equation}
where $d_0$ is the diameter of $\Omega$.

\begin{prop}\label{Contin:10-1}
Let $\Omega$ be a bounded open set of $\mathbb{R}^n$ and $k\in \mathbb{N}_0$, then the space $C^{k}\left(\overline{\Omega}\right)$, with norm \eqref{Contin:3} is a separable space.	
\end{prop}
\textbf{Proof.} Recall that if a topological space is separable, then every subset of it is a separable space \cite[Ch. 3, Sec. 6]{CTV}.

We have already noticed (after Theorem \ref{teor-approssimaz-weierst}) that $C^{0}\left(\overline{\Omega}\right)$ is separable. We consider the case $k=1$ (the case $k>1$ can be
treated similarly). Let $\Psi$ be the map

$$\Psi:C^{1}\left(\overline{\Omega}\right) \ \rightarrow \ \mathcal{X} ,$$
where 

$$\mathcal{X}=\underset{(n+1) \mbox{ - times} }{\underbrace{C^{0}\left(\overline{\Omega}\right)\times\cdots \times C^{0}\left(\overline{\Omega}\right)}}, $$

\smallskip

$$\Psi(u)=(u,\partial_1 u,\cdots \partial_n u),\quad\forall u \in C^{1}\left(\overline{\Omega}\right).$$

\smallskip

if we equip $\mathcal{X}$ with the norm  
$$\left\Vert \mathbf{v}\right\Vert_{\mathcal{X}}=\sum_{j=0}^n \left\Vert v_j\right\Vert_{C^{0}\left(\overline{\Omega}\right)},\quad \forall \mathbf{v}=\left(v_0,v_1,\cdots,v_n\right)\in \mathcal{X},$$
$\Psi$ is an isometry. 

On the other hand, $\mathcal{X}$ is a separable space as a cartesian product of separable spaces. Thus
$\Psi\left(C^{1}(\overline{\Omega})\right)$ is separable
as a subspace of $\mathcal{X}$ and, since $\Psi$ is an isometry, also $C^1\left(\overline{\Omega}\right)$ is separable. $\blacksquare$

\medskip

It is evident that if $k,m\in \mathbb{N}_0$ and $k<m$ then $C^m\left(\overline{\Omega}\right)\subset C^k\left(\overline{\Omega}\right)$. We set

$$C^k_0(\Omega)=\left\{u\in C^k(\Omega):\mbox{ } \mbox{supp } u \mbox{ is a compact set contained in } \Omega\right\},$$

$$C^{\infty}(\Omega)=\bigcap_{k=0}^{\infty}C^k(\Omega),\quad C^{\infty}\left(\overline{\Omega}\right)=\bigcap_{k=0}^{\infty}C^k\left(\overline{\Omega}\right),
\quad C^{\infty}_0(\Omega)=\bigcap_{k=0}^{\infty}C^k_0(\Omega).$$

Let $k\in \mathbb{N}$ or $k=\infty$ and $\widetilde{\Omega}\supset \Omega$, we will often adopt the convention of identifying $C^{k}_0(\Omega)$ with the space of functions $u$ belonging to $C^{k}_0\left(\widetilde{\Omega}\right)$ and such that supp $u$ is a compact set contained in $\Omega$.

\section{The space $C^{k,\alpha}$ }\label{Contin:C-k-alpha}
\label{Funz:holder} Let $X$ be a subset of $\mathbb{R}^n$ and $\alpha\in (0,1]$, we will denote by $C^{0,\alpha}\left(X\right)$ \index{$C^{0,\alpha}\left(X\right)$}the space of the functions $u\in C^{0}\left(X\right)$ which satisfy
$$\left[u\right]_{C^{0,\alpha}\left(X\right)}=\sup \left\{\frac{\left\vert u(x)-u(y)
	\right\vert}{\left\vert x-y\right\vert^{\alpha}}:\mbox{ }x,y\in
X,\mbox{ }x\neq y \right\}<+\infty.$$

If $u$ is a function of $C^{0,\alpha}\left(X\right)$ it is also said that $u$ is a \textbf{H\"{o}lder function of order} $\alpha$\index{H\"{o}lder space}. The number $\alpha$ is said  the \textbf{H\"{o}lder exponent} of the space $C^{0,\alpha}\left(X\right)$. For any $u\in C^{0,\alpha}\left(X\right)$, the number $\left[u\right]_{C^{0,\alpha}\left(X\right)}$ is called the \textbf{H\"{o}lder constant} of $u$. If $\alpha=1$ we will also say that $u$ is a \textbf{Lipschitz function} \index{Lipschitz function} in $X$ and we call \textbf{Lipschitz constant} the number $\left[u\right]_{C^{0,1}\left(X\right)}$. We observe that if $\alpha>1$ and $X$ is a connected open set, then the space $C^{0,\alpha}\left(X\right)$ consists of only the constant functions (as a matter of fact, if $\alpha>1$ then any function of $C^{0,\alpha}\left(X\right)$ is differentiable with zero gradient in $X$). It can be easily checked that if $X$ is bounded, and the space $C^{0,\alpha}\left(X\right)$ is equipped with the norm  
\begin{equation}\label{Contin:5}
	\left\Vert u\right\Vert_{C^{0,\alpha}\left(X\right)}=\left\Vert u\right\Vert_{C^0\left(X\right)} + \left[u\right]_{C^{0,\alpha}\left(X\right)},
\end{equation}
then $C^{0,\alpha}\left(X\right)$ is a \textbf{Banach space}. Sometimes it is convenient to consider, instead of the norm \eqref{Contin:5}, an equivalent "dimensionless" norm, e.g.

\begin{equation*}
	\left\Vert u\right\Vert_{C^0\left(X\right)} + d_0^{\alpha}\left[u\right]_{C^{0,\alpha}\left(X\right)},
\end{equation*}
where $d_0$ is the diameter of $X$.

Let $m\in \mathbb{N}$, we denote by $C^{0,\alpha}\left(X;\mathbb{R}^m\right)$ the space of the functions \\ $u \in C^{0}\left(X;\mathbb{R}^m\right)$ satisfying $u_j\in C^{0,\alpha}\left(X\right)$, $j=1,\cdots,m$. We set

\begin{equation*}
	\left\Vert u\right\Vert_{C^{0,\alpha}\left(X; \mathbb{R}^m\right)}=\left\Vert u\right\Vert_{C^0\left(X;\mathbb{R}^n\right)} + \left[u\right]_{C^{0,\alpha}\left(X; \mathbb{R}^m\right)},
\end{equation*}
where $$\left[u\right]_{C^{0,\alpha}\left(X; \mathbb{R}^m\right)}=\sup \left\{\frac{\left\vert u(x)-u(y)\right\vert_{\mathbb{R}^m}}{\left\vert x-y\right\vert_{\mathbb{R}^n}^{\alpha}}:\mbox{ }x,y\in
X,\mbox{ }x\neq y \right\}<+\infty,$$
where $|\cdot|_{\mathbb{R}^m}$ is the Euclidean norm in $\mathbb{R}^m$ (in the sequel we will often omit the subscript $\mathbb{R}^m$ from this norm).

\bigskip

The following Proposition holds true.

\begin{prop}\label{Contin:prop-holder1}
	Let $X$ be a bounded set of  $\mathbb{R}^n$ and $0<\beta<\alpha\leq 1$ then
	
	\begin{equation}\label{Contin:6}
		\left[u\right]_{C^{0,\beta}\left(X\right)}\leq d_0^{\alpha-\beta} \left[u\right]_{C^{0,\alpha}\left(X\right)},
	\end{equation}
	where $d_0$ is the diameter of $X$;
	\begin{equation}\label{Contin:7}
		\left[u\right]_{C^{0,\beta}\left(X\right)}\leq \left(2\left\Vert u\right\Vert_{C^0\left(X\right)}\right)^{1-\frac{\beta}{\alpha}}\left(\left[u\right]_{C^{0,\alpha}\left(X\right)}\right)^{\frac{\beta}{\alpha}}.
	\end{equation}
	In particular we have 
	
	\begin{equation}\label{Contin:8}
		C^{0,\alpha}\left(X\right)\subset C^{0,\beta}\left(X\right).
	\end{equation}
	\end{prop}
\textbf{Proof.} It suffices to observe that for $x,y\in X$, $x\neq y$, by \eqref{Contin:6}, we have

$$\frac{\left\vert u(x)-u(y)
	\right\vert}{\left\vert x-y\right\vert^{\beta}}=\frac{\left\vert u(x)-u(y)
	\right\vert}{\left\vert x-y\right\vert^{\alpha}}\left\vert x-y\right\vert^{\alpha-\beta}\leq \left[u\right]_{C^{0,\alpha}\left(X\right)} d_0^{\alpha-\beta}.$$
Regarding \eqref{Contin:7}, we first note that the case 
 $$\left[u\right]_{C^{0,\alpha}\left(X\right)}=0$$ is trivial. Let us assume, then, 
$$\left[u\right]_{C^{0,\alpha}\left(X\right)}\neq 0$$
and let $r>0$ be chosen later. Let $x,y\in X$, $x\neq y$. If
$$\left\vert x-y\right\vert\leq r,$$ then

$$\frac{\left\vert u(x)-u(y)
	\right\vert}{\left\vert x-y\right\vert^{\beta}}\leq \left[u\right]_{C^{0,\alpha}\left(X\right)} r^{\alpha-\beta}.$$
If 

$$\left\vert x-y\right\vert> r,$$ then

$$\frac{\left\vert u(x)-u(y)
	\right\vert}{\left\vert x-y\right\vert^{\beta}}\leq 2r^{-\beta}\left\Vert u\right\Vert_{C^0\left(X\right)}.$$
In any case, we have 

$$\left[u\right]_{C^{0,\beta}\left(X\right)}\leq r^{-\beta}\max\left\{r^{\alpha}\left[u\right]_{C^{0,\alpha}\left(X\right)}, 2\left\Vert u\right\Vert_{C^0\left(X\right)}\right\}$$ and, choosing

$$r=\left(\frac{2\left\Vert u\right\Vert_{C^0\left(X\right)}}{\left[u\right]_{C^{0,\alpha}\left(X\right)}}\right)^{1/\alpha},$$ we obtain \eqref{Contin:7}. $\blacksquare$

\bigskip

\textbf{Remark 1.}  By using the Mean Value Theorem, it is easily shown that if $\Omega$ is a bounded, convex open set of $\mathbb{R}^n$ then
\begin{equation}\label{Contin:9}
	C^{1}\left(\overline{\Omega}\right)\subset C^{0,\beta}\left(\overline{\Omega}\right)
\end{equation}
and

\begin{equation*}
	\left[u\right]_{0,1, \Omega}\leq \left\Vert \nabla u\right\Vert_{C^0\left(\overline{\Omega}\right)},
\end{equation*}
where
$$\left\Vert \nabla u\right\Vert_{C^0\left(\overline{\Omega}\right)}=\left\Vert \left|\nabla u \right|\right\Vert_{C^0\left(\overline{\Omega}\right)}.$$

Nevertheless, for a bounded open set $\Omega$ it is not necessarily the case that the inclusion  \eqref{Contin:9} holds. Let us consider, for instance, the following example. Let

$$\Omega=\left\{(x,y)\in \mathbb{R}^2:\ x\leq \sqrt{|y|}, \ \  x^2+y^2<1\right\},$$   

$$1<\beta<2,$$

\begin{equation*}
	u(x,y)=
	\begin{cases}
		x^{\beta}\mbox{sgn}(y), \quad \mbox{ if } x\in \Omega,\ x>0,\\
		\\
		0, \quad \quad\mbox{if } x\in \Omega,\ x\leq 0. %
	\end{cases}%
\end{equation*}
We have $u\in C^{1}\left(\overline{\Omega}\right)$, however if $\alpha$ ssatisfies  $\frac{\beta}{2}<\alpha\leq 1$ then $$u\notin C^{0,\alpha}\left(\overline{\Omega}\right).$$ As a matter of fact, if $x=\sqrt{|y|}$, we get 

$$\frac{\left|u\left(x,y\right)-u\left(x,-y\right)\right|}{(2|y|)^{\alpha}}=2^{1-\alpha}|y|^{\frac{\beta}{2}-\alpha}\rightarrow +\infty,\quad \mbox{as } y\rightarrow 0.$$ $\blacklozenge$

\bigskip

\textbf{Remark 2.} The space $C^{0,\alpha}\left(\overline{\Omega}\right)$, where $\alpha\in (0,1]$ and $\Omega$ is a bounded open set, is not separable. Let us consider the case $n=1$ and  $\Omega=(0,1)$ and for any $a\in (0,1)$,  let us define

\begin{equation*}
	u_a(t)=
	\begin{cases}
		0, \quad \mbox{ for } t\in [0,a),\\
		\\
		(t-a)^{\alpha}, \quad \quad\mbox{for } t\in [a,1]. %
	\end{cases}%
\end{equation*}
We have
\begin{equation}\label{Contin:11-10-3}
	\left\Vert u_a-u_b\right\Vert_{C^{0,\alpha}\left([0,1]\right)}\geq 1,
	\quad \forall a,b \in [0,1] \ \ a\neq b.
\end{equation}

We check \eqref{Contin:11-10-3}. Let $a,b \in [0,1]$, $a<b$ and let us denote
$$v_{a,b}=u_a-u_b.$$
We have 
$$v_{a,b}(b)=u_a(b)-u_b(b)=u_a(b)=(b-a)^{\alpha}$$ and
$$v_{a,b}(a)=u_a(a)-u_b(a)=0.$$
Therefore
\begin{equation*}
	\left[u_a-u_b\right]_{C^{0,\alpha}([0,1])}=\left[v_{a,b}\right]_{C^{0,\alpha}([0,1]}\geq \frac{\left|v_{a,b}(b)-v_{a,b}(a)\right|}{\left|b-a\right|^{\alpha}}=1,
\end{equation*}
which implies \eqref{Contin:11-10-3}. Finally, from the latter and from Proposition \ref{Contin:11-10-1} it follows that $C^{0,\alpha}\left(\left[0,1\right]\right)$ is not separable.
$\blacklozenge$

 \bigskip

 \begin{theo}[\textbf{extension in $C^{0,\alpha}$}]\label{Contin:teo-est}
 	\index{Theorem:@{Theorem:}!- extension in $C^{0,\alpha}$@{- extension in $C^{0,\alpha}$}}
 	Let $X$ be a bounded set of $\mathbb{R}^n$ and $u\in C^{0,\alpha}(X)$, $\alpha\in (0,1]$, then there exists $U\in C^{0,\alpha}\left(\mathbb{R}^n\right)$ such that
 	
 	\begin{equation}\label{Contin:12-10-0}
 		U(x)=u(x),
 		\quad \forall x\in X,
 	\end{equation}
 	
 	\begin{equation}\label{Contin:12-10-00}
 		\left\Vert U\right\Vert_{C^{0}\left(\mathbb{R}^n\right)}=\left\Vert u\right\Vert_{C^{0}\left(X\right)},
 	\end{equation}
 	
 	\begin{equation}\label{Contin:12-10-000}
 		\left[U\right]_{C^{0,\alpha} \left(\mathbb{R}^n\right)}=\left[u\right]_{C^{0,\alpha} (X)}.
 	\end{equation}
 \end{theo}
 \textbf{Proof.} Let us denote
 
 $$M=\left\Vert u\right\Vert_{C^{0}\left(X\right)},\quad \quad m=\left[u\right]_{C^{0,\alpha} (X)}$$ and let us define the function
 
 \begin{equation*}
 	v(x)=\sup_{y\in X} \left\{u(y)-m|x-y|^{\alpha}\right\},\quad\mbox{ for } x\in \mathbb{R}^n.
 \end{equation*}
 We have
 
 \begin{equation}\label{Contin:12-10-1}
 	v(x)=u(x),
 	\quad \forall x\in X.
 \end{equation}

We check \eqref{Contin:12-10-1}. First note that we have trivially 
 \begin{equation}\label{Contin:12-10-2}
 	u(x)\leq v(x),
 	\quad \forall x\in X.
 \end{equation}
  On the other hand we have
 \begin{equation*}
 	u(y)-u(x)\leq m|x-y|^{\alpha}, \quad \ \ \forall x,y\in X
 	\quad \forall x\in X,
 \end{equation*}
 hence
 \begin{equation*}
 	u(y)-m|x-y|^{\alpha}\leq u(x), \quad \ \ \forall x,y\in X
 	\quad \forall x\in X.
 \end{equation*}
 Consequently 
 $$v(x)\leq u(x),\quad \ \ \forall x\in X.$$
 By the latter and by \eqref{Contin:12-10-2} we get \eqref{Contin:12-10-1}. 
 
 We also notice that
 \begin{equation}\label{Contin:12-10-3}
 	v(x)\leq M,
 	\quad \forall x\in \mathbb{R}^n.
 \end{equation}
 Now, for any $x\in \mathbb{R}^n$ let us define
 
 \begin{equation*}
 	U(x)=
 	\begin{cases}
 		v(x), \quad \mbox{for } |v(x)|\leq M ,\\
 		\\
 		-M, \quad \quad\mbox{for } v(x)<-M. %
 	\end{cases}%
 \end{equation*}
Let us note that, by \eqref{Contin:12-10-3}, $U$ is defined throughout $\mathbb{R}^n$. Let us note also that if $x \in X$, then \eqref{Contin:12-10-1} gives \eqref{Contin:12-10-0} and 
 
 $$\sup_{X} |U(x)|=M.$$
We also have 
$$\sup_{\mathbb{R}^n} |U(x)|=M,$$
 Concerning the latter, notice that, if $|v(x)|\leq M$ then $U(x)=v(x)$, hence $|U(x)|\leq M$ and if $v(x)<-M$ then $|U(x)|= M$.

 It only remains to prove that $U\in C^{0,\alpha}\left(\mathbb{R}^n\right)$ and that \eqref{Contin:12-10-000} holds.  Let, then, $x,y\in \mathbb{R}^n$ be such that $x\neq y$ and otherwise arbitrary. Let us suppose that $U(x)\neq U(y)$. For instance, let us assume

 \begin{equation}\label{Contin:12-10-4}
 	U(x)> U(y).
 \end{equation}
 Let us check that  
 \begin{equation}\label{Contin:12-10-5}
 	0<U(x)-U(y)\leq v(x)-v(y).
 \end{equation}
 The following cases occur.
 
 \smallskip
 
 \noindent (a) $v(x)< -M$ and $|v(y)|\leq M$, \\ (b) $v(x)< -M$ and $v(y)< -M$, \\ (c) $|v(x)|\leq M$ and $|v(y)|\leq M$, \\ (d)  $|v(x)|\leq M$ and $v(y)< -M$.
 
 \smallskip
 
Cases (a) and (b) cannot occur. As a matter of fact, in case (a) we would have
 $$U(x)-U(y)=-M-v(y)\leq -M+M=0,$$ that contradicts \eqref{Contin:12-10-4}. In case (b) we would have
 $$U(x)-U(y)=-M-(-M)=0,$$
 that contradicts \eqref{Contin:12-10-4}.
 
 In case (c) we have
 
 $$U(x)-U(y)=v(x)-v(y).$$
 Finally, in case (d) we have
 $$U(x)-U(y)=v(x)-(-M)=v(x)+M<v(x)-v(y).$$
 Therefore \eqref{Contin:12-10-5} holds true. 
 
 We have
 \begin{equation}\label{Contin:12-10-6}
 	\begin{aligned}
 		v(x)-v(y)&=\sup_{z\in X} \left\{u(z)-m|x-z|^{\alpha}\right\}-\sup_{z\in X} \left\{u(z)-m|y-z|^{\alpha}\right\}\leq\\&\leq m\sup_{z\in X} \left\{|y-z|^{\alpha}-|x-z|^{\alpha}\right\}\leq\\&\leq 
 		m\sup_{\zeta\in X} \left\{\left(|y-x|+|\zeta|\right)^{\alpha}-|\zeta|^{\alpha}\right\}.
 	\end{aligned}
 \end{equation}

Now, let us denote by $$\omega(t)=t^{\alpha}, \quad \mbox{if } t\in [0,+\infty).$$ Since $\omega$ is concave, we have
 $$ \omega(t+h)-\omega(t)\leq \omega(h)-\omega(0)=\omega(h),\quad \forall t,h\in   [0,+\infty),$$
 by the just obtained inequality, by \eqref{Contin:12-10-5} and \eqref{Contin:12-10-6} we get
 \begin{equation}
 	|U(x)-U(y)|\leq m|y-x|^{\alpha}
 \end{equation}
 and since we have proved \eqref{Contin:12-10-0}, we get \eqref{Contin:12-10-000}. 
 
 $\blacksquare$
 
 \bigskip
 
 A more general version of Theorem \ref{Contin:teo-est}, valid for uniformly continuous functions, can be found in \cite[Chapter 4]{Pu}. 
 
 \bigskip 
 
 Let us notice that if $$u:X\rightarrow\mathbb{R}^m$$ is a Lipschitz continuous function, then, Theorem \ref{Contin:teo-est} implies that there exists an extension 
 
 $$U:\mathbb{R}^n\rightarrow\mathbb{R}^m,$$ such that  
 
 $$[U]_{C^{0,1}\left(\mathbb{R}^n;\mathbb{R}^m\right)}\leq \sqrt{m} [u]_{C^{0,1}\left(X;\mathbb{R}^m\right)}.$$ Actually, this result can be improved. As a matter of fact, the following Theorem of \textbf{Kirszbraun} holds true, for the proof of which we refer to \cite[cap. 7]{Ma}.
 
 \begin{theo}[\textbf{Kirszbraun}]\label{Kirszb}
 	\index{Theorem:@{Theorem:}!- Kirszbraun@{- Kirszbraun}}
 	 Let $u:X\rightarrow \mathbb{R}^m$ be
 	a Lipschitz continuous function, where $X\subset \mathbb{R}^n$, then there
 	exists $U\in C^{0,1}\left(\mathbb{R}^n;\mathbb{R}^m\right)$ such that
 	$$U=u, \quad\mbox{in } X$$
 	and
 	$$[U]_{C^{0,1}\left(\mathbb{R}^n;\mathbb{R}^m\right)}=[u]_{C^{0,1}\left(X;\mathbb{R}^m\right)}.$$
 \end{theo}
 
 \bigskip
 
 \bigskip
 
 Let $\Omega$ be a bounded open set of $\mathbb{R}^n$, $k\in \mathbb{N}_0$ and $\alpha\in (0,1]$, we denote by  $C^{k,\alpha}\left(\overline{\Omega}\right)$ the space of the functions $u\in C^{k}\left(\overline{\Omega}\right)$, satisfying

 $$\left[\partial^{\beta}u\right]_{C^{0,\alpha} \left(\overline{\Omega}\right)}<+\infty,\quad \forall \beta\in \mathbb{N}_0^n, \ \ |\beta|=k.$$
 It is easily proven that $C^{k,\alpha}\left(\overline{\Omega}\right)$, equipped with the norm 

\begin{equation}\label{Contin:10}
	\left\Vert u\right\Vert_{C^{k,\alpha}\left(\overline{\Omega}\right)}=\left\Vert u\right\Vert_{C^k\left(\overline{\Omega}\right)}+ \left[u\right]_{C^{k,\alpha}\left(\overline{\Omega}\right)}
\quad \forall u\in C^k\left(\overline{\Omega}\right),
\end{equation}
is a \textbf{Banach space}, where

$$\left[u\right]_{C^{k,\alpha}\left(\overline{\Omega}\right)}=\sum_{|\beta|=k}\left[\partial^{\beta}u\right]_{C^{0,\alpha}\left(\overline{\Omega}\right)}.$$
Sometimes, instead of the norm \eqref{Contin:10} we will consider the dimensionless norm
$$	\left\Vert u\right\Vert_{C^{k,\alpha}\left(\overline{\Omega}\right)}=
\sum_{|\beta|\leq k} d_0^{|\beta|}
\sup_{\overline{\Omega}}|\partial^{\beta}u|+d_0^{k+\alpha}\left[u\right]_{C^{k,\alpha}\left(\overline{\Omega}\right)}.$$ 
where $d_0$ is the diameter of $\Omega$.

We define the space $C_{loc}^{k,\alpha}\left(\Omega\right)$, $k\in \mathbb{N}_0$, $0<\alpha\leq 1$, \index{$C_{loc}^{k,\alpha}\left(\Omega\right)$}as the space of functions $u\in C^0(\Omega)$ such that for every bounded open set $\omega \Subset \Omega$ (i.e. $\overline{\omega}$ compact and $\overline{\omega}\subset \Omega$ \index{$\Subset$}) we have 
$$u_{|\omega}\in C^{k,\alpha}\left(\overline{\omega}\right).$$

\section{Review of measure theory and $L^p$ spaces}\label{Sob:pag-r}

In this Section, we give, for the convenience of the reader, the main definitions and statements of the main theorems of the Measure Theory and of $L^p$ spaces.
Some reference texts are \cite{Pu}, \cite{Roy} (see also lecture notes \cite{Magn1} and \cite{Magn2}).  

\subsection{Measurable sets, measurable functions, positive measures}
	\label{Sob:RICHIAMI}
\begin{definition}\label{Sob:def1R}
	Let $X$ be a set and $\mathcal{M}$ be a family of subsets of
	$X$ with the following properties:
	
	\smallskip
		
	(i)$X\in\mathcal{M}$,
	
	(ii) $E\in\mathcal{M}$ $\Longrightarrow$ $\mathcal{C}E:=X\setminus
	E\in\mathcal{M}$,
	
	(iii) $E_j\in\mathcal{M}$, $j\in\mathbb{N}$ $\Longrightarrow$
	$\bigcup_{j\in\mathbb{N}}E_j\in \mathcal{M}$.
	
	\smallskip
	
	\noindent $\mathcal{M}$ is called a  \textbf{$\sigma$--algebra} \index{$\sigma$-algebra}and
	the couple $(X,\mathcal{M})$ is called a \textbf{measurable space}. \index{measurable space}
\end{definition}

We will be interested almost exclusively in the case where
$X=\mathbb{R}^n$ and $\mathcal{M}$ consists of the Lebesgue measurable subsets
of $\mathbb{R}^n$.

\begin{definition}\label{Sob:def2R}
	Let $(X,\mathcal{M})$ be a measurable spoce and $Y$ be a topologic space. We say that the function
	$$f:X\rightarrow Y,$$
	is a \textbf{measurable function}\index{measurable function}, provided we have
	$$f^{-1}(A)\in \mathcal{M},\quad\quad \mbox{for every open subset } A \mbox{ of }
	Y.$$
\end{definition}
Let us recall that if $(X,\mathcal{M})$ is a measurable space,
$Y$, $Z$ two topological spaces, $f:X\rightarrow Y$ is a measurable function and 
$g:Y\rightarrow Z$ is a continuous function, then $g\circ f:X\rightarrow Z$
is a measurable function.

\medskip

We denote by
$\overline{\mathbb{R}}=\mathbb{R}\cup\{-\infty,+\infty\}$, \index{$\overline{\mathbb{R}}$} the extended real line equipped with its usual topology.

\medskip

The following theorems hold true.
\begin{theo}\label{Sob:teo1R}
	Let $(X,\mathcal{M})$ be a measurable space and $f:X\rightarrow
	\overline{\mathbb{R}}$. Then $f$ is a measurable function if and only if for any 
	$t\in \mathbb{R}$ one of the following level sets	
	$$f^{-1}((t,+\infty]),\quad f^{-1}([t,+\infty]),\quad
	f^{-1}([-\infty,t)),\quad f^{-1}([-\infty,])),$$ is a measurable set.
\end{theo}

\begin{theo}\label{Sob:teo2R}
	Let $(X,\mathcal{M})$ be a measurable space. Then we have
	
	(i) If $f,g:X\rightarrow \overline{\mathbb{R}}$ are two measurable functions
	and $f+g$ is well defined, then  $f+g$ e $\lambda f$ are measurable functions
	(we use the convention that $0\cdot(\pm \infty)=0$).
	
	(ii) Let $\left\{f_k\right\}$ be a sequence of measurable functions, then	
	$$\sup_{k\in \mathbb{N}}f_k,\quad \inf_{k\in \mathbb{N}}f_k,\quad \liminf_{k\rightarrow\infty}f_k,\quad \limsup_{k\rightarrow\infty}f_k,$$
	are measurable functions.
\end{theo}

\bigskip

We define as \textbf{simple function} \index{simple function} on the measurable space $X$ a
function $$s:X\rightarrow \mathbb{R}$$ that assumes a
finite set of values.

\medskip

\begin{theo}\label{Sob:teo3R}
	Let $(X,\mathcal{M})$ be a measurable space and $f:X\rightarrow
	\overline{\mathbb{R}}$ be a measurable function, then there exists a sequence $\left\{s_k\right\}$ of simple functions such that
		$$\lim_{k\rightarrow\infty}s_k(x)=f(x),\quad \forall x\in X.$$ If
	$f$ is bounded then $\left\{s_k\right\}$ uniformly converges to $f$.
\end{theo}

\medskip

\begin{definition}\label{Sob:def3R} Let $(X,\mathcal{M})$ a measurable space. We say that
	
	$$\mu:\mathcal{M}\rightarrow[0,\infty],$$ is a \textbf{positive measure} \index{measure} provided that we have
	
	\smallskip
	
	(i) $\mu(\emptyset)=0$;
	
	(ii) if $\left\{E_j\right\}_{j\in \mathbb{N}}$ 
	is a countable family of measurable sets such that $$E_i\cap
	E_j=\emptyset,\quad \mbox{ for } i\neq j,$$ then we have 
	$$\mu\left(\bigcup_{j=1}^{\infty}E_j\right)=\sum_{j=1}^{\infty}\mu\left(E_j\right).$$
	The tern $(X,\mathcal{M},\mu)$ is called a \textbf{measure space} \index{measure space}.
\end{definition}

\bigskip

If $E$ is a Lebesgue measurable set of $\mathbb{R}^n$,
we will denote by $|E|$ its measure.

\medskip

The following theorems hold true.

\begin{theo}\label{Sob:teo4R}
	Let $(X,\mathcal{M},\mu)$ be a measure space. The following  properties hold true. 
	
	\smallskip

	(i) if $\left\{E_j\right\}_{1\leq j\leq N}$ is a finite family of measurable sets such that  $E_i\cap E_j=\emptyset$, for
	$i\neq j$, then 
	
	$$\mu\left(\bigcup_{j=1}^NE_j\right)=\sum_{j=1}^{N}\mu\left(E_j\right);$$
	
	(ii) if $E\subset F$ and $E,F\in \mathcal{M}$, then $\mu(E)\leq
	\mu(F)$;
	
	(iii) if $\left\{E_j\right\}_{j\in \mathbb{N}}$ is a countable  family of measurable sets such that $E_j \subset E_{j+1}$, for every $j\in \mathbb{N}$, then
	
	$$\lim_{j\rightarrow
		\infty}\mu\left(E_j\right)=\mu\left(\bigcup_{j=1}^{\infty}E_j\right);$$

	(iv) if $\left\{E_j\right\}_{j\in \mathbb{N}}$ is a countable  family of measurable sets such that $$E_{j+1} \subset E_{j},$$
	for every $j\in \mathbb{N}$, then
	
	$$\lim_{j\rightarrow
		\infty}\mu\left(E_j\right)=\mu\left(\bigcap_{j=1}^{\infty}E_j\right).$$
	
\end{theo}

\bigskip

\begin{theo}\label{Sob:teo5R}
	Let $E$ be a Lebesgue measurable subset of $\mathbb{R}^n$ whose measure be finite. Let  $\left\{f_j\right\}$ be a sequence of measurable functions such that there exists the limit     $$\lim_{j\rightarrow
		\infty}f_j(x)$$ and it is finite almost everywhere.
	Then, for every $\varepsilon>0$ there exists a compact set $K\subset E$ which satisfies
	$|E\setminus K|<\varepsilon$ and
	$$f_j\rightarrow f,
	\mbox{ as } j\rightarrow\infty,\quad\mbox{ \textbf{uniformly on } } K.$$
\end{theo}

\bigskip

\begin{theo}[\textbf{Lusin}]\label{Sob:teo6R}
	\index{Theorem:@{Theorem:}!- Lusin@{- Lusin}}
	Let $E$ be a Lebesgue mesurable subset of $\mathbb{R}^n$ 
	which has finite measure, and let $f:E\rightarrow
	\overline{\mathbb{R}}$ such that $$|f(x)|<+\infty,\quad \mbox{ a.e. }
	x\in E.$$ Then $f$ is a measurable function in $E$ if and only if for each
	$\varepsilon>0$ there exists $K\subset E$, $K$ closed, such that
	$|E\setminus K|<\varepsilon$ and $f_{|K}$ is a continuous function.
\end{theo}

\bigskip

Now, let us define \textbf{the integral over the measure space} \index{integral over a measure space}
$(X,\mathcal{M},\mu)$. Let $s$ be a nonnegative simple function

$$s(x)=\sum_{i=1}^Nc_j\chi_{E_j},$$
where $\left\{E_j\right\}_{1\leq j\leq N}$ is a finite family of measurable set
pairwise disjonts and $c_j\geq 0$,
$j=1,\cdots,N$. If $E\in \mathcal{M}$, we set by definition 

$$\int_Es(x)d\mu=\sum_{i=j}^Nc_j\mu\left(E\cap E_j\right),$$
in which the convention $0\cdot \infty=0$ occurs. We call 
$\int_E s(x)d\mu$ "the integral of $s$ over $E$".

\medskip

Let us consider the measurable function
$$f:X\rightarrow [0,+\infty].$$ We call
\textbf{the Lebesgue integral of $f$ with respect to the measure $\mu$} the following element
of $\overline{\mathbb{R}}$
$$\int_E f(x)d\mu:=\sup\left\{\int_E s(x)d\mu: s\mbox{ simple function, } 0\leq s\leq f \mbox{ in }
E\right\}.$$ We say that $f$ a summable function over $E$ if
$$\int_E f(x)d\mu<+\infty.$$

\bigskip

\begin{theo}\label{Sob:teo7R}
	Let $(X,\mathcal{M},\mu)$ be a measure space and let
	$f,g:X\rightarrow[0,+\infty]$ $E,F\in \mathcal{M}$. The following properties hold true:
	
	\smallskip
	
	(i) $\int_E f d\mu=\int_X f\chi_E d\mu$;
	
	\smallskip
	
	(ii) if $f\leq g$ in $E$ then $\int_E f d\mu\leq \int_E g d\mu$;
	
	\smallskip
	
	(iii) if $E\subset F$ then $\int_E f d\mu\leq \int_F f d\mu$;
	
	\smallskip
	
	(iv) if $f=0$ in $E$ then $\int_E f d\mu=0$;
	
	\smallskip
	
	(v) if $\mu(E)=0$ then $\int_E f d\mu=0$.
\end{theo}

\bigskip

\begin{theo}[\textbf{Monotone Convergence}]\label{Sob:teo8R}
	\index{Theorem:@{Theorem:}!- Monotone Convergence@{- Monotone Convergence}}
	Let $(X,\mathcal{M},\mu)$ be a measure space. Let
	$\left\{f_j\right\}$ be a sequence of nonnegative measurable functions which satisfy
	
	$$f_j(x)\leq f_{j+1}(x),\quad \forall x\in X,\quad \forall j\in
	\mathbb{N}.$$
	Then
	
	$$\int_X\lim_{j\rightarrow\infty}f_j(x)d\mu=\lim_{j\rightarrow\infty}\int_Xf_j(x)d\mu.$$
	
\end{theo}

\bigskip

\begin{theo}[\textbf{Fatou}]\label{Sob:teo9R}
	\index{Theorem:@{Theorem:}!- Fatou@{- Fatou}} Let $(X,\mathcal{M},\mu)$ be a measure space. Let
	$\left\{f_j\right\}$ be a sequence of nonnegative measurable functions, then we have
	
	$$\int_X\liminf_{j\rightarrow\infty}f_j(x)d\mu\leq\liminf_{j\rightarrow\infty}\int_Xf_j(x)d\mu.$$
	
\end{theo}

\medskip

\begin{theo}\label{Sob:teo10R}
	Let $(X,\mathcal{M},\mu)$ be a measure space and let
	$$f,g:X\rightarrow[0,+\infty]$$ be two measurable functions, then 
	$$\int_X(f+g)d\mu=\int_Xfd\mu+\int_Xgd\mu;$$
	$$\int_X\lambda fd\mu=\lambda\int_Xfd\mu,\quad\forall \lambda \in
	\mathbb{R},$$ with the convention $0\cdot\int_Xfd\mu=0$.
\end{theo}

\medskip

\begin{theo}\label{Sob:teo11R}
	Let $(X,\mathcal{M},\mu)$ be a measure space and let
	$\left\{f_j\right\}$ be a sequence of nonnegative measurable functions, then we have
	
	$$\int_X\sum_{j=1}^{\infty}f_jd\mu=\sum_{j=1}^{\infty}\int_Xf_jd\mu.$$
\end{theo}

\bigskip

Furthermore, recall that if
$$f:X\rightarrow[0,+\infty],$$
is a measurable function and defining 

$$\nu(E)=\int_Efd\mu,\quad \forall E\in \mathcal{M},$$
$\nu$ turns out to be a measure on $X$.

\medskip

\begin{definition}\label{Sob:def4R} Let $(X,\mathcal{M},\mu)$ be a measure space and let $$f:X\rightarrow \mathbb{R}.$$ We say that $f$ is
	summable over $X$ provided
	$$\int_X|f|d\mu<+\infty.$$
	In such a case we set
	$$\int_Xfd\mu=\int_Xf_+d\mu-\int_Xf_-d\mu,$$
	where $f_+=\max\{f,0\}$, $f_+=-\min\{f,0\}$. We denote by $\mathcal{L}^1(X)$ \index{$\mathcal{L}^1(X)$} the class of summable functions over $X$.
\end{definition}

\bigskip

\begin{theo}\label{Sob:teo12R}
	Let $(X,\mathcal{M},\mu)$ be a measure space and let
	$f\in\mathcal{L}^1(X)$, then 
	$$\mu\left(\left\{x\in X:\quad |f(x)|=+\infty\right\}\right)=0.$$
\end{theo}

\medskip

\begin{theo}\label{Sob:teo13R}
	$\mathcal{L}^1(X)$ is a vector space and
	$$\mathcal{L}^1(X)\ni f\rightarrow \int_Xfd\mu\in \mathbb{R},$$
	is a linear map. Furthermore, if $f,g\in \mathcal{L}^1(X)$ then 
	$$\max\{f,g\}\in \mathcal{L}^1(X)$$ and
	
	$$\left|\int_Xfd\mu\right|\leq \int_X|f|d\mu, \quad\quad \forall
	f\in \mathcal{L}^1(X).$$
\end{theo}

\medskip

\begin{theo}\label{Sob:teo14R}
	If $f\in\mathcal{L}^1(X)$ then
	$$\forall \varepsilon>0 \mbox{ } \exists \delta>0\mbox{ such that }\forall E\in \mathcal{M} \mbox{ and } \mu(E)<\delta\mbox{ we have } \int_E|f|d\mu<\varepsilon.$$
\end{theo}

\medskip

\begin{theo} [\textbf{Dominated Convergence}]\label{Sob:teo15R}
	\index{Theorem:@{Theorem:}!- Dominated Convergence@{- Dominated Convergence}}
	Let $\left\{f_j\right\}$ be a sequence of measurable functions in $\mathcal{L}^1(X)$. Let us assume
	
	\smallskip
	
	(i) $$\lim_{j\rightarrow\infty}f_j(x)=f(x),\quad\mbox{ a.e. } x\in
	X,$$
	
	\smallskip
	
	(ii) there exists $g\in \mathcal{L}^1(X)$ such that
	$$\left|f_j(x)\right|\leq g(x),\quad\mbox{ a.e. } x\in X,\mbox{
	}j\in \mathbb{N}.$$ 

Then $f\in \mathcal{L}^1(X)$ and

	$$\lim_{j\rightarrow\infty}\int_X\left|f_j-f\right|d\mu=0,$$
	
	$$\lim_{j\rightarrow\infty}\int_Xf_jd\mu=\int_Xfd\mu.$$
\end{theo}

\bigskip

The Monotone Convergence Theorem and the Dominated Convergence Theorem  give that if $\left\{f_j\right\}$ is a sequence of measurable functions in $\mathcal{L}^1(X)$ satisfying

$$\sum_{j=1}^{\infty}\int_X\left|f_j\right|d\mu<+\infty,$$
then $\sum_{j=1}^{\infty}f_j$ converges almost  everywhere to
a function of $\mathcal{L}^1(X)$ and 

$$\sum_{j=1}^{\infty}\int_Xf_jd\mu=\int_X\sum_{j=1}^{\infty}f_jd\mu.$$

\bigskip

\begin{theo}[\textbf{derivation under the integral sign}]\label{Sob:teo16R}
	\index{Theorem:@{Theorem:}!- derivation under the integral sign@{- derivation under the integral sign}}
Let $(X,\mathcal{M},\mu)$ be a measure space and let $A$ be an open set of $\mathbb{R}^n$ . Let
	$$F:A\times X\rightarrow \mathbb{R}$$
	satisfy
	
	\smallskip
	
	(i) $F(x,\cdot)\in \mathcal{L}^1(X)$ for every $x\in A$,
	
	(ii) $F(\cdot,y)\in C^1(A)$ for almost every $y\in X$.
	
	\smallskip
	
	If, for any $k=1,\cdots,n$, there exist $g_k\in \mathcal{L}^1(X)$, $g_k\geq
	0$ such that
	
	$$\left|\partial_{x_k}F(x,y)\right|\leq g_k(y),\quad \forall y\in
	X,\quad \forall x\in A,$$ then the function
	
	$$G(x):=\int_XF(x,y)d\mu(y),\quad x\in A$$
	is of  $C^1(A)$ class  and we have
	
	$$\partial_{x_k}G(x)=\int_X\partial_{x_k}F(x,y)d\mu(y),\quad \forall
	x\in A,\mbox{ } k=1, \cdots, n.$$
\end{theo}

\bigskip

\begin{theo}[\textbf{Fubini--Tonelli}]\label{Sob:teo17R}
	\index{Theorem:@{Theorem:}!- Fubini--Tonelli@{- Fubini--Tonelli}} Let
	$$f:\mathbb{R}^n\times \mathbb{R}^m\rightarrow
	\overline{\mathbb{R}}$$ be a measurable function. We have
	
	(i) If $f\geq 0$ then: $f(x,\cdot)$ is measurable for almost every $x\in \mathbb{R}^n$, in addition the function
	$$\mathbb{R}^n\ni x\rightarrow \int_{\mathbb{R}^m}f(x,y)dy\in[0,+\infty],$$
	is measurable over $\mathbb{R}^n$ and we have

	\begin{equation}\label{Sob:stella11R}
		\int_{\mathbb{R}^n}\left(\int_{\mathbb{R}^m}f(x,y)dy\right)dx=\int_{\mathbb{R}^{n+m}}f(x,y)dxdy.
	\end{equation}
	
	\smallskip
	
	(ii) If $f\in \mathcal{L}^1(\mathbb{R}^{n+m})$ then $f(x,\cdot)\in
	\mathcal{L}^1(\mathbb{R}^{m})$ for almost every $x\in \mathbb{R}^n$,
	furthermore
	
	$$\int_{\mathbb{R}^m}f(\cdot,y)dy\in \mathcal{L}^1(\mathbb{R}^{n})$$
	and \eqref{Sob:stella11R} holds true.
\end{theo}

\bigskip

\subsection{The $L^p$ spaces}\label{12-10-spaziLp}
Let $p\in[1,+\infty)$ and $(X, \mathcal{M},\mu)$ be a measurable space.
We say that $f\in \mathcal{L}^p(X)$ \index{$\mathcal{L}^p(X)$}if $f$ is measurable and
$|f|^p\in \mathcal{L}^1(X)$. $\mathcal{L}^p(X)$ is a
vector space. We define $L^p(X)$ as the quotient space $\left(\mathcal{L}^p(X)/\sim\right)$ where  "$\sim $" is the equivalence relation on $\mathcal{L}^p(X)$ defined as follows: $f\sim g$ if and only if $f=g$ almost everywhere. We  equip $L^p(X)$ with the norm \index{$\left\Vert \cdot \right\Vert_{L^p(X)}$}

$$\left\Vert f\right\Vert_{L^p(X)}=\left(\int_X |f|^p d\mu\right)^{1/p}.$$

We say that $f\in \mathcal{L}^{\infty}(X)$ provided

$$\mbox{ess}\sup|f|=\inf\left\{t\in \mathbb{R}: \quad
\mu(\left\{|f(x)|>t\right\})=0\right\}<+\infty.$$ 
We define $L^{\infty}(X)$ \index{$L^p(X)$} similarly as we have previously defined $L^{p}(X)$, $p<+\infty$. We equip $L^{\infty}(X)$ with the norm 
$$\left\Vert f\right\Vert_{L^{\infty}(X)}=\mbox{ess}\sup|f|.$$

\bigskip

\noindent\underline{\textbf{Minkowski inequality.}} \index{Minkowki inequality}If
$p\in[1,+\infty]$, $f,g\in L^p(X)$ then

$$\left\Vert f+g\right\Vert_{L^p(X)}\leq \left\Vert f\right\Vert_{L^p(X)}+\left\Vert g\right\Vert_{L^p(X)}.$$

\bigskip

\noindent\underline{\textbf{H\"{o}lder inequality.}} \index{H\"{o}lder inequality}Let
$p\in[1,+\infty]$, let us denote by $p'$ (the conjugate of $p$) the element of $[1,+\infty]$
 satisfying

$$\frac{1}{p}+\frac{1}{p'}=1.$$ If $f\in L^p(X)$ and $g\in L^{p'}(X)$
then $fg\in L^1(X)$ and
$$\left\Vert fg\right\Vert_{L^1(X)}\leq \left\Vert f\right\Vert_{L^p(X)}\left\Vert
g\right\Vert_{L^{p'}(X)}.$$

\bigskip

If $\mu(X)<+\infty$, we have

$$p_2\geq p_1\Longrightarrow L^{p_2}(X)\subset L^{p_1}(X)$$ and the function

$$p\rightarrow \left(\frac{1}{\mu(X)}\int_X |f|^p d\mu\right)^{1/p},$$
turns out to be an increasing function (just apply H\"{o}lder inequality).

Moreover

$$\lim_{p\rightarrow\infty}\left\Vert f\right\Vert_{L^p(X)}=\left\Vert
f\right\Vert_{L^{\infty}(X)}.$$

\bigskip

\begin{theo}\label{Sob:teo18R}
	Let $(X,\mathcal{M},\mu)$ be a measure space and let
	$p\in[1,+\infty]$. Then $L^p(X)$ is a \textbf{Banach space}. If $p=2$, $L^2(X)$ is a \textbf{Hilbert space} equipped with the scalar product	
	$$(f,g)_{L^2(X)}=\int_Xfgd\mu,\quad \forall f,g\in L^2(X).$$
\end{theo}

\bigskip

We say that the measure space $(X,\mathcal{M},\mu)$ is
$\sigma$--finite, if there exists a countable family 
$\left\{X_j\right\}_{j\in \mathbb{N}}\subset \mathcal{M}$ such that
$$X=\bigcup_{j\in \mathbb{N}}X_j,\quad\mbox{ and } \quad \mu(X_j)<+\infty.$$

\medskip

\begin{theo}\label{Sob:teo19R}
	Let $(X,\mathcal{M},\mu)$ be a $\sigma$--finite measure space and
	let \\ $p\in[1,+\infty)$, then $F$ is a bounded linear functional 
	from $L^p(X)$ to $\mathbb{R}$ if and only if there exists  $g\in
	L^{p'}(X)$ which satisfies
	
	$$F(f)=\int_Xgfd\mu,\quad \forall f\in L^p(X).$$
\end{theo}

\bigskip

\noindent\underline{ \textbf{Density and separability  in $L^p$.}} Let
$E\subset \mathbb{R}^n$ be a Lebesgue measurable set.

\bigskip

The following theorems hold true.

\begin{theo}[\textbf{density of simple functions in $L^p(E)$, $1\leq p\leq \infty$}]\label{Sob:teo20R}
	If $f\in L^p(E)$ and  $p\in[1,+\infty]$, then for every
	$\varepsilon>0$ there exists a simple function $s$ such that
	$$\left\Vert f-s\right\Vert_{L^p(E)}<\varepsilon.$$
\end{theo}

\medskip

\begin{theo}[\textbf{density of $C^0_0(E)$ in $L^p(X)$, $1\leq p< \infty$}]\label{Sob:teo21R}
	If $f\in L^p(E)$ and  $p\in[1,+\infty)$, then for every
	$\varepsilon>0$ there exists $g\in C^0_0\left(E\right)$, such that
	$$\left\Vert f-g\right\Vert_{L^p(E)}<\varepsilon.$$
\end{theo}

\medskip

\begin{theo}\label{Sob:teo22R}
	If $p\in[1,+\infty)$, then $L^p(X)$ is a separable space.
	$L^{\infty}(X)$ is not a separable space.
\end{theo}

\medskip

\begin{theo}\label{Sob:teo24R}
	If $p\in[1,+\infty)$ and $f\in L^p(\mathbb{R}^n)$ then
	$$\lim_{\delta\rightarrow 0} \left(\sup_{|h|<\delta}\int_{\mathbb{R}^n}
	|f(x-h)-f(x)|^pdx\right)=0.$$
\end{theo}

\bigskip

\noindent \underline{\textbf{Reflexivity.}}
\index{reflexive space}

Let $X$ be a normed space, let us denote by $X'$ the space of the
bounded linear functionals from $X$ to $\mathbb{R}$. As it is well--known, $X'$ is called the \textbf{dual} of $X$ and $X'$ turns out a Banach space (whether or not $X$ is a
Banach space) equipped with the norm

$$\left\Vert f\right\Vert_{X'}=\sup\left\{\langle f, u\rangle:\mbox{ }\left\Vert u\right\Vert_{X}\leq 1
\right\},\quad \forall f\in X',$$ where

$$\langle f, u\rangle=f(u),\quad \forall u\in X,$$
$\langle \cdot, \cdot\rangle$ is said the "duality bracket of 
of $X'$ and $X$".

\medskip

\begin{definition}\label{Sob:def5R} Let $X$ be a Banach space. We say that $X$ is a \textbf{riflexive space} \index{reflexive space} provided
	$$\forall v\in (X')'\quad \exists u\in X,\mbox{ such that } \langle v,
	f\rangle=\langle f, u\rangle\quad \forall f\in X'.$$
\end{definition}

\medskip

Keep in mind that if $u\in X$ then the map
\begin{equation*}
	X'\ni f \rightarrow j_u(f):=\langle f,
	u\rangle,\end{equation*}
is a bounded linear functional, hence
$j_u\in (X')'$ 
$j_u\in (X')'$ and it can prove (by the Hanh--Banach Theorem \cite{Br}) that

\begin{equation}\label{Sob:1-15R}
	\left\Vert j_u\right\Vert_{(X')'}=\left\Vert u\right\Vert_{X},\quad \forall u\in X.\end{equation}
Therefore, it is defined the map $j$
$$X\ni u\rightarrow j_u\in (X')'.$$
The map $j$ is injective (applying, again, Hanh--Banach Theorem  ).
Consequently, $j$ is an embedding and, by
\eqref{Sob:1-15R} it is an isometry. 
If $X$ is a reflexive space
 then $j$ is also suriective and, by the Open Map Theorem, $j^{-1}$ is continuous too. Ultimately, if $X$ is a reflexive space, we can identify $(X')'$ with $X$ by means of $j$. Recall that the
Hilbert spaces are reflexive.

\medskip

\begin{definition}\label{Sob:def6R} Let $X$ be a Banach space and let $\left\{u_k\right\}$ be a sequence of 
	$X$. We say tat $\left\{u_k\right\}$ weakly converges \index{weak convergence} to $u\in X$
	and we write
	
	$$u_k\rightharpoonup u,\quad \mbox{as }
	k\rightarrow\infty, \  \mbox{(or } \left\{u_k\right\}\rightharpoonup u,\mbox{)} ,$$ provided
	$$\langle f,
	u_k\rangle\rightarrow\langle f, u\rangle,\quad \mbox{as }
	k\rightarrow\infty,\mbox{ }\forall f\in X'.$$  If a weak limit exists, then it is unique
	(it can be again proved by Hanh--Banach Theorem).
\end{definition}

\medskip

\begin{prop}\label{Sob:pag16R}
	Let $X$ be a Banach space and $\left\{u_k\right\}$ be a
	sequence of $X$. 
	
	We have
	
	\smallskip
	
	(i) if $\left\{u_k\right\}$ weakly converges then it is bounded;
	
	\smallskip
	
	(ii) if $\left\{u_k\right\}$ weakly converges to $u$ then
	$$\left\Vert u\right\Vert_{X}\leq \liminf_{k\rightarrow \infty}\left\Vert
	u_k\right\Vert_{X}.$$
\end{prop}
\textbf{Proof.} Let us prove (i), since for every 
$f\in X'$, $\left\{\langle f, u_k\rangle\right\}$ is a converging
sequence of $\mathbb{R}$, it is bounded, that is

$$\sup_{k\in \mathbb{N}}\left|\langle f, u_k\rangle\right|\leq
C(f)<+\infty,\quad \forall f\in X'.$$ Now, applying the Banach--Steinhauss Theorem (see \cite{Magn2}) to the map
$$X'\ni f\rightarrow T_k(f):=\langle f, u_k\rangle\in \mathbb{R},\\ k\in \mathbb{N}$$
we have
$$\sup_{k\in \mathbb{N}}\left\Vert u_k\right\Vert_X=\sup_{k\in \mathbb{N}}\left\Vert
T_k\right\Vert_{(X')'}<+\infty.$$

Let us prove  (ii). Since $\left\{u_k\right\}\rightharpoonup
u$, $k\rightarrow\infty,$, by (i) we get 
$$\sup_{k\in \mathbb{N}}\left\Vert u_k\right\Vert_X<+\infty.$$
Moreover, for every $f\in X'$, we have

\begin{equation*}
	\begin{aligned}
		\langle f, u\rangle&=\lim_{k\rightarrow\infty}\langle f,
		u_k\rangle=\liminf_{k\rightarrow\infty}\langle f, u_k\rangle\leq
		\\&\leq \liminf_{k\rightarrow\infty}\left\Vert
		f\right\Vert_{X'}\left\Vert u_k\right\Vert_X=\\&=\left\Vert
		f\right\Vert_{X'}\liminf_{k\rightarrow\infty}\left\Vert
		u_k\right\Vert_X.
	\end{aligned}
\end{equation*}
Therefore

$$\left\Vert u\right\Vert_{X}\leq \liminf_{k\rightarrow\infty}\left\Vert
u_k\right\Vert_X.$$ $\blacksquare$

\medskip

We recall
\begin{theo}[\textbf{Banach--Alaoglu}]\label{Sob:teo23R}
	\index{Theorem:@{Theorem:}!- Banach--Alaoglu@{- Banach--Alaoglu}} Let $X$ be a reflexive Banach and let $\left\{u_k\right\}$ be a bounded sequence of $X$. Then there exists a
	subsequence $\left\{u_{k_j}\right\}$ which weakly converges.
\end{theo}

\medskip

\begin{theo}\label{Sob:teo23Rbis}
	Let $E$ be a measurable subset of $\mathbb{R}^n$. Then, 
	$L^p(E)$ is a reflexive space, for every $p\in (1,+\infty)$.
\end{theo}

\bigskip

\bigskip

\underline{\textbf{Convolution.}}\index{convolution product}

Let $f,g$ be two measurable functions defined in $\mathbb{R}^n$ with
values in $\mathbb{R}$. Let $x\in \mathbb{R}^n$, if the function of the
	variable $y$, $f(x-y)g(y)$ is summable, we set \index{$f\star g$}
	$$(f\star g)(x)=\int_{\mathbb{R}^n}f(x-y)g(y)dy.$$
	If the function $(f\star g)(x)$ is defined for almost every $x\in \mathbb{R}^n$, we will call it the \textbf{convolution product} (or, simply, the
	\textbf{convolution}) of $f$ and $g$.
	
\medskip

\begin{theo}[\textbf{the Young inequality}]\label{Sob:teoYoung}
	\index{Theorem:@{Theorem:}!- Young inequality@{- Young inequality}}
	Let $f\in L^p(\mathbb{R}^n)$ and $g\in L^q(\mathbb{R}^n)$, where
	$$\frac{1}{p}+\frac{1}{q}\geq 1,$$ then $f\star g \in
	L^r(\mathbb{R}^n)$ where
	$$r=\frac{1}{p}+\frac{1}{q}-1$$
	and we have
	
	$$\left\Vert f\star g\right\Vert_{L^r(\mathbb{R}^n)}\leq \left\Vert f\right\Vert_{L^p(\mathbb{R}^n)}\left\Vert
	g\right\Vert_{L^q(\mathbb{R}^n)}.$$
\end{theo}

\bigskip

\bigskip

\noindent Let $\eta\in C^{\infty}_0\left(\mathbb{R}^n\right)$ satisfy

\smallskip

(i) supp $\eta\subset B_1$,

\smallskip

(ii) $\eta\geq 0$

\smallskip

(iii) $\int_{\mathbb{R}^n} \eta(x)dx=1.$

\medskip

\noindent $\eta$ is named a \textbf{mollifier} \index{mollifier}. For instance, let 

\begin{equation*}
	\tilde{\eta}(s)=
	\begin{cases}
	c_n\exp \left\{-\frac{1}{1-4s^2} \right\}, \quad \mbox{ for } s\in [0,1/2),\\
		\\
		0, \quad \quad\mbox{otherwise}, %
	\end{cases}%
	\end{equation*}
where 

$$c_n=\left(\int^{1/2}_0s^{n-1}\exp \left\{-\frac{1}{1-4s^2}\right\}ds\right)^{-1},$$ then 
$$\eta(x)=\tilde{\eta}(|x|),$$
is a mollifier.

\medskip

\noindent Here and in the sequel we set, for $\varepsilon>0$,

$$\eta_{\varepsilon}(x)=\varepsilon^{-n}\eta\left(\varepsilon^{-1}x\right).$$

\begin{theo}\label{Sob:teo25R} (i) if $f\in L^p(\mathbb{R}^n)$, $p\in [1,+\infty)$, then
	$$\eta_{\varepsilon}\star f \rightarrow f,\quad\mbox{ as } \varepsilon\rightarrow 0, \quad\mbox{in }
	L^p(\mathbb{R}^n). $$ (ii) If
	$f$ is uniformly continuous and bounde in $\mathbb{R}^n$, then
	
	$$\eta_{\varepsilon}\star f \rightarrow f,\quad \mbox{ as } \varepsilon\rightarrow 0 \quad\mbox{uniformly in } \mathbb{R}^n.$$
\end{theo}

\medskip

\textbf{Remark.} Let $E$ be a measurable set of
$\mathbb{R}^n$, $f\in L^p(E)$, where $p\in [1,+\infty)$, then 

$$\left\Vert \eta_{\varepsilon}\star \left(f\chi_E\right)-f\right\Vert_{L^p(E)}\rightarrow 0,\quad\mbox{as } \varepsilon\rightarrow 0,$$
where

$$\eta_{\varepsilon}\star f=\eta_{\varepsilon}\star
\left(f\chi_E\right)=\int_E\eta_{\varepsilon}(x-y)f(y)dy.$$
$\blacklozenge$

\bigskip

\begin{theo}\label{Sob:teo26R} Let $\Omega$ be an open set of $\mathbb{R}^n$, $f\in L^p(\Omega)$, where $p\in [1,+\infty)$, then
	$$\eta_{\varepsilon}\star f \in C^{\infty}(\Omega),$$
	
	$$\partial^{\alpha}\left(\eta_{\varepsilon}\star f \right)=\left(\partial^{\alpha}\eta_{\varepsilon}\right)\star f$$
	and
	$$\eta_{\varepsilon}\star f \rightarrow f,\quad\mbox{in }
	L^p(\Omega),\mbox{ as } \varepsilon\rightarrow 0. $$ Hence
	$C^{\infty}(\Omega)$ is dense in $L^p(\Omega)$.
\end{theo}

We also have

\begin{theo}[\textbf{density of $C^{\infty}_0(\Omega)$ in $L^p(\Omega)$, $1\leq p<+\infty$}]\label{Sob:teo27R} 
	\index{Theorem:@{Theorem:}!- density of $C^{\infty}_0(\Omega)$ in $L^p(\Omega)$@{- density of $C^{\infty}_0(\Omega)$ in $L^p(\Omega)$}}
	Let $\Omega$ be an open set of
	$\mathbb{R}^n$, $p\in [1,+\infty)$, then $C^{\infty}_0(\Omega)$
	is dense in $L^p(\Omega)$.
\end{theo}

\bigskip

Let $\Omega$ be an open set of $\mathbb{R}^n$ and $p\in [1,+\infty]$.
We denote by $L^p_{loc}(\Omega)$ the space of measurable functions
defined on $\Omega$ such that for every compact set $K$ we have $f_{|K}\in L^p(K)$. Let
$\left\{u_k\right\}$ be a sequence of $L^p_{loc}(\Omega)$, we write

$$u_k\rightarrow u,\quad \mbox{as }
k\rightarrow\infty,\mbox{ in } L^p_{loc}(\Omega),$$ provided $u\in
L^p_{loc}(\Omega)$ and for every compact $K\subset \Omega$, we have
$$(u_k)_{|K}\rightarrow u_K,\quad \mbox{as }
k\rightarrow\infty,\mbox{ in } L^p(K).$$

Let us define, for any $\varepsilon>0$,

$$\Omega_{\varepsilon}=\left\{x\in \Omega:\mbox{ } \mbox{dist}(x,\partial\Omega)>\varepsilon\right\}.$$
If $f\in L^1_{loc}(\Omega)$, then
$\left(\eta_{\varepsilon}\star f\right)(x)$ is defined for every $x\in
\Omega_{\varepsilon}$, and we may rephrase theorems
\ref{Sob:teo25R} and \ref{Sob:teo26R} as follows.

\medskip

\begin{theo}\label{Sob:teo28R} Let $\Omega$ be an open set of $\mathbb{R}^n$ and $f\in
	L^1_{loc}(\Omega)$. Then
	
	\smallskip
	
	(i) $\eta_{\varepsilon}\star f\in C^{\infty}(\Omega_{\varepsilon})$;
	
	\smallskip
	
	(ii) if $f\in C^0(\Omega)$ then $\eta_{\varepsilon}\star
	f\rightarrow f$, uniformly on the compact sets of $\Omega$ as
	$\varepsilon\rightarrow 0$;
	
	\smallskip
	
	(iii) if $p\in [1,+\infty)$ and $f\in L^p_{loc}(\Omega)$ then
	$$\eta_{\varepsilon}\star
	f\rightarrow f,\quad\mbox{as } \varepsilon\rightarrow 0,\mbox{ in }
	L^p_{loc}(\Omega).$$
\end{theo}

\medskip

Let $f$ be a measurable function defined on $\mathbb{R}^n$ and let

$$\mathcal{O}=\left\{A\subset \mathbb{R}^n:\mbox{ } A \mbox{ open and }f=0 \mbox{ in } A \mbox{ a.e.
}\right\},$$ the set 
$$\mbox{supp}f=\mathbb{R}^n\setminus \bigcup_{A\in \mathcal{O}} A$$
is named the \textbf{essential support} \index{essential support} of $f$.
Hereafter, if there is no ambiguity, instead of "essential support of $f$" we will simply say "support of $f$". We recall that  then the essential support of a function if $u\in C^0\left(\mathbb{R}^n\right)$ is equal to the support defined in Section \ref{Funz:Cont-diff}.  

\section{Partition of unity \index{partition of unity}} \label{Partiz}
Let us start by the following

\begin{lem}\label{Partiz:22-10-22-1}
	Let $\Omega$ be an open set of $\mathbb{R}^n$ and $K$ be a compact set contained in $\Omega$, then there exists $\varphi\in C^{\infty}_0\left(\mathbb{R}^n\right)$ such that supp $\varphi\subset \Omega$, $0\leq \varphi\leq 1$ and $\varphi=1$ in a neighborhood of $K$.
\end{lem} 
\textbf{Proof.} For any $\varepsilon>0$, let us denote by 
$$K^{(\varepsilon)}=\left\{x\in \mathbb{R}^n:\mbox{ } \mbox{dist}(x,K)\leq \varepsilon\right\}.$$
Let $\varepsilon_0$ and $\varepsilon_1$ satisfy 
$$0<\varepsilon_0<\varepsilon_1<\varepsilon_0+\varepsilon_1<\mbox{dist}\left(K,\mathbb{R}^n\setminus \Omega\right).$$
Let us define

\begin{equation*}
	\varphi(x)=\int_{K^{(\varepsilon_1)}}\eta_{\varepsilon_0}(x-y)dy.
\end{equation*} 
It can be easily checked that $\varphi\in C^{\infty}_0\left(\mathbb{R}^n\right)$,
$$\mbox{supp }\varphi\subset K^{(\varepsilon_0+\varepsilon_1)}\subset \Omega$$
and
$$\varphi(x)=1,\quad \forall x\in K^{(\varepsilon_0)}.$$ $\blacksquare$ 

\medskip

\begin{lem}\label{Partiz:23-10-22-1}
	Let $K$ be a compact set of $\mathbb{R}^n$ and let $V_1, V_2,\cdots,
	V_l$ be some open sets of $\mathbb{R}^n$ satisfying
	$$K\subset \bigcup_{j=1}^lV_j.$$
	Then there exist the functions $\zeta_1, \cdots, \zeta_l\in
	C^{\infty}_0(\mathbb{R}^n)$ which satisfy 
	$$\mbox{supp }\zeta_j\subset V_j,\quad j=1,\cdots, l,$$
	$$ 0\leq \zeta_j,\quad j=1\cdots, l;\quad \sum_{j=1}^l
	\zeta_j\leq 1,\quad \mbox{on } \mathbb{R}^n$$
	$$ \sum_{j=1}^l
	\zeta_j=1,\quad \mbox{in a neighborhood of } K.$$ 
\end{lem}
\textbf{Proof.} Let us denote, for any $\varepsilon>0$ and $V\subset \mathbb{R}^n$, 

$$V_{\varepsilon}=\left\{x\in V:\mbox{ } \mbox{dist}\left(x,\partial V\right)>\varepsilon\right\}.$$
We have 
$$K\subset \bigcup_{\varepsilon>0}\bigcup_{j=1}^lV_{\varepsilon, j}$$ and by the compactness of $K$ it follows that there exists $\varepsilon_0>0$ such that  
$$K\subset \bigcup_{j=1}^lV_{\varepsilon_0,j}\subset \bigcup_{j=1}^l\overline{V_{\varepsilon_0,j}}\subset \bigcup_{j=1}^lV_j,$$
(because $\overline{V}_{\varepsilon_0, j}\subset V_j$). Hence, denoting

$$K_j=K\cap \overline{V_{\varepsilon_0.j}}\quad, j=1,\cdots, l$$ we have immediately that $K_j$ is a compact set, $K_j\subset V_j$, for any $j=1,\cdots, l$ and 
$$K\subset \bigcup_{j=1}^lK_j.$$
By Lemma \ref{Partiz:22-10-22-1}, we derive that for every $j\in \{1,\cdots,l\}$  there exist \\ $\varphi_j \in C^{\infty}_0\left(V_j\right)$ satisfying
$$0\leq \varphi_j\leq 1,\quad \varphi_j=1, \mbox{ in a neighborhood, } W_{j}, \mbox{ of } K_j.$$

Now, definining

$$\zeta_1=\varphi_1,\quad \zeta_2=\varphi_2\left(1-\varphi_1\right),\quad \cdots,\quad \zeta_j=\varphi_j\left(1-\varphi_1\right)\cdots \left(1-\varphi_{j-1}\right),$$ we get  

\begin{equation*}
	\begin{aligned}
		&\sum_{j=0}^l
		\zeta_j=\varphi_1+\varphi_2\left(1-\varphi_1\right)+\cdots+ \varphi_l\left(1-\varphi_1\right)\cdots \left(1-\varphi_{l-1}\right)=\\&=
		1- \left(1-\varphi_1\right)+\varphi_2\left(1-\varphi_1\right)+\cdots+ \varphi_l\left(1-\varphi_1\right)\cdots \left(1-\varphi_{l-1}\right)=\\&=
		1-\left(1-\varphi_1\right)\left(1-\varphi_2\right)+\varphi_3\left(1-\varphi_1\right)\left(1-\varphi_2\right)+\cdots +\varphi_l\left(1-\varphi_1\right)\cdots \left(1-\varphi_{l-1}\right)=\\&=
		1-\left(1-\varphi_1\right)\left(1-\varphi_2\right)\cdots \left(1-\varphi_l\right).
	\end{aligned}
\end{equation*}
Therefore, if $$x\in \bigcup_{j=1}^lW_{j},$$ there exists $\overline{j}\in \{1,\cdots, l\}$ such that $x\in W_{\overline{j}}$, hence $\varphi_{\overline{j}}(x)=1$ and
$$\sum_{j=0}^l
\zeta_j(x)=1,\quad \forall x\in \bigcup_{j=1}^lW_{j}.$$ Since $\bigcup_{j=1}^lW_{j}$ is a neighborhood of $K$, the Lemma is proved. $\blacksquare$

\bigskip

In what follows, we will say that the set of functions $\varphi_1,\cdots, \varphi_l$ is a \textbf{partition of the unity subordinate to the covering} $\left\{V_j\right\}_{1\leq j\leq l}$.

\bigskip

\begin{theo}[\textbf{partition of unity}]\label{Sob:lem3.3}
	\index{Theorem:@{Theorem:}!- partition of unity@{- partition of unity}}
	Let $\Omega$ be an open set of $\mathbb{R}^n$ and let $V_1, V_2,\cdots,
	V_l$ be open sets of $\mathbb{R}^n$ satisfying
	$$\partial \Omega\subset \bigcup_{j=1}^lV_j.$$ Then there exist the functions $\zeta_0, \zeta_1, \cdots, \zeta_l\in
	C^{\infty}(\mathbb{R}^n)$ such that 
	$$(\zeta_0)_{|\Omega}\in
	C^{\infty}_0(\Omega),$$
	$$\mbox{supp }\zeta_0\subset \mathbb{R}^n\setminus \partial\Omega;\quad \mbox{supp
	}\zeta_j\subset V_j,\quad j=1,\cdots, l,$$
	
	$$ \sum_{j=0}^l
	\zeta_j=1,\quad \mbox{on } \mathbb{R}^n;\quad 0\leq \zeta_j\leq
	1,\quad j=0,1\cdots, l.$$ 
\end{theo}
\textbf{Proof.} Let us consider a partition of unity subordinate to the covering, $\left \{V_j\right\}_{1\leq j\leq l}$, of $\partial \Omega$. Let us denote by $$\mathcal{O}_j= \mbox{supp} \varphi_j, \quad j=1,\cdots, l$$and set 
$$V_0=\mathbb{R}^n\setminus \bigcup_{j=1}^l \mathcal{O}_j.$$ Let
$$\zeta_0=1-\sum_{j=1}^l\zeta_j.$$
We get trivially $\zeta_0\in C^{\infty}_0\left(\mathbb{R}^n\right)$ and
\begin{equation*}
	\sum_{j=0}^l\zeta_j(x)=1,\quad\forall x\in \mathbb{R}^n,
\end{equation*}

\begin{equation*}
	\mbox{supp }\zeta_0\subset \mathbb{R}^n\setminus \partial\Omega.
\end{equation*}
Moreover, by  Lemma \ref{Partiz:23-10-22-1} there exists an open neighborhood, $\mathcal{U}$, of $\partial\Omega$ such that
$$\sum_{j=1}^l\zeta_j(x)=1,\quad \forall x\in \mathcal{U}.$$ Hence
$$\zeta_0(x)=0,\quad \forall x\in \mathcal{U}.$$ Therefore 
$$\mbox{supp} (\zeta_0)_{|\Omega}\subset \overline{\Omega}\setminus \mathcal{U}\subset \Omega,$$ this implies  $$(\zeta_0)_{|\Omega}\in
C^{\infty}_0(\Omega)$$ concluding the proof. $\blacksquare$ 

\bigskip

\textbf{Remark.} By Theorem \ref{Sob:lem3.3} it is evident that $V_0\cap \Omega, V_1,\cdots, V_l$  is a covering of $\overline{\Omega}$ and $(\zeta_0)_{|\Omega},\zeta_1,\cdots,\zeta_l$ is a partition of the unity subordinate to that covering. $ \blacklozenge$

\section{The Lebesgue differentiation Theorem}\label{Diff-Leb}
In this Section we prove the following

\begin{theo}[\textbf{Lebesgue differentiation}]\label{Diff-Leb:13-10-22-1} \index{Theorem:@{Theorem:}!- Lebesgue differentiation@{- Lebesgue differentiation}} If $f\in L_{loc}^1\left(\mathbb{R}^n\right)$ then 
	
	\begin{equation}\label{Diff-Leb:13-10-22-2}
	\lim_{r\rightarrow 0} \ \dashint_{B_r(x)}f(y)dy=f(x),\quad\mbox{a.e. } x\in \mathbb{R}^n, 
		\end{equation}
	where \index{$\dashint_{B_r(x)}f(y)dy$}
	\begin{equation*}
		\dashint_{B_r(x)}f(y)dy=\frac{1}{\left|B_r(x)\right|}\int_{B_r(x)}f(y)dy.
	\end{equation*}
	\end{theo}

\bigskip 
In order to prove Theorem \ref{Diff-Leb:13-10-22-1} we need some preliminary lemmas and propositions. We start by 

\begin{lem}[\textbf{Covering}]\label{Diff-Leb:13-10-22-3} 
	\index{Lemma:@{Lemma:}!- covering@{- covering}}Let $E$ be a Lebesgue measurable subset of $\mathbb{R}^n$ and let $\mathcal{B}$ be a family of balls of  $\mathbb{R}^n$ satisfying
\begin{equation*}
	E\subset\bigcup_{B\in \mathcal{B}}B.
\end{equation*}
and, 
	\begin{equation*}
		\sup_{B\in \mathcal{B}} d(B) <+\infty,
	\end{equation*}
where $d(B)$ is the diameter of $B$.
Then there exists a countable (or finite) family, $\left\{B_k\right\}_{k\in \Lambda}\subset \mathcal{B} $ which satisfies
\begin{equation*}
	B_j\cap B_k=\emptyset, \quad \mbox{for }j\neq k, \ j,k\in \Lambda
\end{equation*}
and
\begin{equation*}
	\sum_{k\in \Lambda}\left|B_k\right|\geq 5^{-n}|E|.
\end{equation*}

\end{lem}
\textbf{Proof of the Lemma.} Firstly we construct the family $\mathcal{B}_0:=\left\{B_k\right\}_{k\in \Lambda} $. Let $B_1\in \mathcal{B}$ satisfy
$$d\left(B_1\right)\geq\frac{1}{2}\sup\left\{d(B):\mbox{ }B\in \mathcal{B}\right\}$$
and set

$$d_1=\sup\left\{d(B):\mbox{ } B\in \mathcal{B}, \ B\cap B_1=\emptyset\right\}<+\infty.$$
Let us consider the set

$$\mathcal{F}_1=\left\{B\in \mathcal{B}: \mbox{ } B\cap B_1=\emptyset,\mbox{ }d(B)\geq \frac{1}{2}d_1\right\}
.$$
If $\mathcal{F}_1= \emptyset$, then we choose $\mathcal{B}_0:=\left\{B_1\right\}$; otherwise, if  $\mathcal{F}_1\neq \emptyset$, we choose $B_2\in \mathcal{F}_1$ and we continue the process. Let us suppose we have chosen $B_1,\cdots, B_i$, let us consider the family
$$\mathcal{F}_i=\left\{B\in \mathcal{B}: \mbox{ } B\cap \bigcup_{j=1}^i B_j=\emptyset,\mbox{ } d(B)\geq \frac{1}{2} d_i
\right\},$$
where
$$ d_i=\sup\left\{d(B):\mbox{ }B\in \mathcal{B}, \ B\cap \bigcup_{j=1}^i B_j=\emptyset\right\}.$$
If $\mathcal{F}_i= \emptyset$, then we choose $\mathcal{B}_0:=\left\{B_1,\cdots, B_i\right\}$; otherwise, if  $\mathcal{F}_i\neq \emptyset$, then we choose $B_{i+1}\in \mathcal{F}_i$ as above and continue the process. All in all, we construct a finite or infinite, countable family, $\left\{B_i\right\}_{i\in \Lambda}\subset \mathcal{B} $ such that
\begin{equation*}
	B_j\cap B_k=\emptyset, \quad \mbox{for }j\neq k, \ j,k\in \Lambda
\end{equation*}
and
\begin{equation*}
d\left(B_{i}\right)\geq \frac{1}{2}
\sup\left\{d(B):\mbox{ } B\in \mathcal{B}, \ B\cap \bigcup_{j=1}^{i-1}B_j=\emptyset\right\}, \quad \mbox{for } i\geq 2.
\end{equation*}

Let us consider the case in which $\mathcal{B}_0$ is finite. Let

\begin{equation}\label{Leb:14-10-22-1}
\mathcal{B}_0=\left\{B_1,\cdots, B_k\right\}.
\end{equation}
Let us denote by 
 $$ \mathcal{A}=\left\{B\in  \mathcal{B}: \mbox{ } B\cap \bigcup_{j=1}^k B_j =\emptyset\right\},$$ 
  $$ \mathcal{C}=\mathcal{B}\setminus \mathcal{A}=\left\{B\in  \mathcal{B}: \mbox{ } B\cap \bigcup_{j=1}^k B_j \neq \emptyset\right\},$$
 $$A=\bigcup_{B\in \mathcal{A} } B,$$

$$C=\bigcup_{B\in \mathcal{C} } B.$$

We have
\begin{equation}\label{Leb:14-10-22-2}
	E\subset \bigcup_{B\in \mathcal{B} } B=A\cup C.
\end{equation}

\medskip

\textbf{Claim.}
We have
$$ \mathcal{A}=\emptyset.$$

\medskip

\textbf{Proof of Claim.} We argue by contradiction. Let us assume that $\mathcal{A}\neq \emptyset$. Let us first observe that it cannot occur that  

$$d(B)< \frac{1}{2}
\sup\left\{d(B):\mbox{ } B\in \mathcal{B}, \ B\cap \bigcup_{j=1}^{k}B_j=\emptyset\right\}, \quad \forall B\in \mathcal{A}$$
otherwise we would have
\begin{equation}
	\begin{aligned}
		0<&\sup\left\{d(B):\mbox{ } B\in \mathcal{B}, \ B\cap \bigcup_{j=1}^{k}B_j=\emptyset\right\}\leq\\&\leq  \frac{1}{2}
		\sup\left\{d(B):\mbox{ } B\in \mathcal{B}, \ B\cap \bigcup_{j=1}^{k}B_j=\emptyset\right\},
	\end{aligned}
\end{equation}
Which is evidently absurd.

Therefore, there exists $\tilde{B}\in \mathcal{A}$ such that

$$d\left(\tilde{B}\right)\geq \frac{1}{2}
\sup\left\{d(B):\mbox{ } B\in \mathcal{B}, \ B\cap \bigcup_{j=1}^{k}B_j=\emptyset\right\}.$$
In particular, we have  $\tilde{B}\notin \left\{B_1,\cdots, B_k\right\}$ and, consequently, the process of construction of $\mathcal{B}_0$ does not stop, but we had assumed the opposite (i.e., $\mathcal{B}_0$ finite family consisting of $k$ elements) and thus we have a contradiction.

\medskip

By what is proved in the Claim and by \eqref{Leb:14-10-22-2}, we have

\begin{equation}\label{Leb:14-10-22-3}
	E\subset  C.
\end{equation}
Now, let us prove
\begin{equation}\label{Leb:14-10-22-4}
E\subset \bigcup_{j=1}^k B^*_j,
\end{equation}
where $B^*_j$ denotes the ball having the same center as $B_j$ with radius equal to $5$ times the radius of
$B_j$. Let $x\in E$, then \eqref{Leb:14-10-22-3} implies that there exists $\hat{B}
\in  \mathcal{B}$ satisfy $$\hat{B}\cap \bigcup_{j=1}^k B_j \neq \emptyset.$$ Let $j_0\in \{1,\cdots, k\}$ such that 

$$x\in \hat{B}, $$
$$\hat{B}\cap B_{j_0}\neq \emptyset$$ and $$\hat{B}\cap \bigcup_{j=1}^{j_0-1}B_j=\emptyset,$$ (if $j_0=1$, $\bigcup_{j=1}^{j_0-1}B_j$ is the empty set). By the third relationship we have

$$d\left(B_{j_0}\right)\geq \frac{1}{2}d\left(\hat{B}\right).$$ By the latter and by $\hat{B}\cap B_{j_0}\neq \emptyset$ we easily obtain
$$\hat{B}\subset B^*_{j_0}.$$ Hence 

\begin{equation*}
	x\in \bigcup_{j=1}^k B^*_j,
\end{equation*}
and \eqref{Leb:14-10-22-4} is proved. Moreover we have

$$|E|\leq \sum_{j=1}^k\left|B^*_j\right|=5^n\sum_{j=1}^k\left|B_j\right|.$$

\medskip 

Now let us consider the case where $\mathcal{B}_0$ is infinite. Hence, in such a case we have $\mathcal{B}_0=\left\{B_k\right\}_{k\in \mathbb{N}}$. If 

  $$\sum_{k=1}^{\infty}\left|B_k\right|=+\infty,$$
there is nothing to prove. Let us assume that
\begin{equation*}
	\sum_{k=1}^{\infty}\left|B_k\right|<+\infty,
\end{equation*}
which implies
\begin{equation*}
	\lim_{k\rightarrow\infty}d\left(B_k\right)=0.
\end{equation*}
If $\tilde{B}\in \mathcal{B}$ and 

\begin{equation}\label{Leb:14-10-22-5}
	k_0=\min\left\{j\in\mathbb{N}: \mbox{ } d\left(B_{j+1}\right) <\frac{1}{2}d\left(\tilde{B}\right)\right\},
\end{equation}
then

\begin{equation}\label{Leb:14-10-22-6}
\tilde{B}\cap\bigcup_{i=1}^{k_0}B_i\neq \emptyset.
\end{equation}
To prove \eqref{Leb:14-10-22-6} it suffices to notice that if it were
\begin{equation*}
	\tilde{B}\cap\bigcup_{i=1}^{k_0}B_i= \emptyset,
\end{equation*}
we would have
\begin{equation*}
	d\left(B_{k_0+1}\right)\geq \frac{1}{2}
	\sup\left\{d(B):\mbox{ } B\in \mathcal{B}, \ B\cap \bigcup_{j=1}^{k_0}B_j=\emptyset\right\}\geq \frac{1}{2}d\left(\tilde{B}\right),  
\end{equation*}
which contradicts \eqref{Leb:14-10-22-5}. 

Since \eqref{Leb:14-10-22-6} holds true, we set 
\begin{equation}\label{Leb:14-10-22-7}
	j_0=\min\left\{j\in \{1,\cdots, k_0\}: \mbox{ } B_{j}\cap \tilde{B}\neq \emptyset\right\}
\end{equation}
obtaining
$$\tilde{B}\cap \bigcup_{j=1}^{j_0-1}B_j=\emptyset.$$ Hence (by \eqref{Leb:14-10-22-7} and by the definition of $\mathcal{B}_0$)
\begin{equation*}
	\begin{cases}
	\tilde{B}\cap B_{j_0}\neq \emptyset ,\\
		\\
	d\left(B_{j_0}\right)\geq \frac{1}{2}d\left(\tilde{B}\right). %
	\end{cases}%
\end{equation*}
From which it follows that for every $B\in \mathcal{B}$ there exists $B_{j_0}\in \mathcal{B}_0$ such that 
$B\subset B^*_{j_0}$. Therefore, arguing as in the finite case, we have 

\begin{equation*}
	E\subset \bigcup_{B\in \mathcal{B} } B=\bigcup_{j=1}^{\infty}B^*_j
\end{equation*}
and by the latter the thesis follows. $\blacksquare$

\bigskip

Now we introduce the notion of \textbf{maximal function}\index{maximal function}. Let $f\in L^1\left(\mathbb{R}^n\right)$, the following function is called the maximal function associated to $f$

\begin{equation}\label{Diff-Leb:15-10-22-1}
	M(f)(x)=\sup_{r>0}\frac{1}{\left|B_r(x)\right|}\int_{B_r(x)}|f(y)|dy,\quad \forall x\in \mathbb{R}^n. 
\end{equation}
Let us observe that $M(f)$ is a measurable function. More precisely, the following Proposition holds true

\begin{prop}\label{Diff-Leb:prop-15-10-22-1}
If $f\in L^1\left(\mathbb{R}^n\right)$ then $M(f)$ is a lower semicontinuous function. 
\end{prop}
\textbf{Proof.} If $f$ is identically $0$, we have $M(f)\equiv 0$. Let us assume that $f$ is not identically $0$. Hence 

\begin{equation}\label{Diff-Leb:15-10-22-2}
M(f)(x)>0,\quad \forall x\in \mathbb{R}^n.
\end{equation}

Fix $t\geq 0$ and let us prove that 

\begin{equation*}
	A=\left\{x\in \mathbb{R}^n: \mbox{ } M(f)(x)>t \right\}
\end{equation*}
is an open set. 

In the case where $t=0$, we have $A= \mathbb{R}^n$. In the case where $t>0$, let $x_0\in A$ and $0<\varepsilon < M(f)(x_0)-t$. By the definition of $M(f)$, there exists $r_{\varepsilon}>0$ such that 

\begin{equation*}
\frac{1}{\left|B_{r_{\varepsilon}}(x_0)\right|}\int_{B_{r_{\varepsilon}}(x_0)}|f(y)|dy> M(f)(x_0)-\varepsilon> t.
\end{equation*}
Now, let $0<\eta<M(f)(x_0)-t-\varepsilon$. Since $f\in L^1\left(\mathbb{R}^n\right)$, there exists $\delta>0$ such that if $|x_0-x|<\delta$ then  

\begin{equation*}
	\left|\int_{B_{r_{\varepsilon}}(x_0)}|f(y)|dy-\int_{B_{r_{\varepsilon}}(x)}|f(y)|dy\right|<\eta\left|B_{r_{\varepsilon}}\right| .
\end{equation*}
Hence 

\begin{equation*}
\begin{aligned}
M(f)(x)&\geq \frac{1}{\left|B_{r_{\varepsilon}}(x)\right|}\int_{B_{r_{\varepsilon}}(x)}|f(y)|dy>\\&
>
\frac{1}{\left|B_{r_{\varepsilon}}(x_0)\right|}\int_{B_{r_{\varepsilon}}(x_0)}|f(y)|dy-\eta>\\&
>M(f)(x_0)-\varepsilon-\eta> t,
\end{aligned}
\end{equation*}
which implies
$$B_{\delta}(x_0)\subset A.$$
Therefore $A$ is open. $\blacksquare$

\bigskip

\begin{lem}\label{Leb:16-10-22-1-0}
	Let $f\in L^1\left(\mathbb{R}^n\right)$ and $M(f)$ its maximal function, then
	
	\begin{equation}\label{Diff-Leb:16-10-22-1}
	\left|\left\{x\in \mathbb{R}^n: \mbox{ } M(f)(x)>t \right\}\right|\leq \frac{5^n}{t}\int_{\mathbb{R}^n}|f(y)|dy,
	\end{equation}
	
\end{lem}
\textbf{Proof.} Set
$$E_t=\left\{x\in \mathbb{R}^n: \mbox{ } M(f)(x)>t \right\}.$$ If $x\in E_t$, then $M(f)(x)>t$. Hence there exists $r_x>0$ such that 
$$\int_{B_{r_{x}}(x)}|f(y)|dy>t\left|B_{r_{x}}(x)\right|.$$ For the sake of brevity, set $B_x=B_{r_{x}}(x)$ so that we have
\begin{equation}\label{Diff-Leb:16-10-22-2}
	\frac{1}{t}\int_{B_x}|f(y)|dy>\left|B_{x}\right|.
\end{equation}
Now, we have trivially that $\left\{B_x\right\}_{x\in E_t}$ is a covering of $E_t$; moreover \eqref{Diff-Leb:16-10-22-2} and $f\in L^1\left(\mathbb{R}^n\right)$ give 
$$\sup_{x\in E_t}\left|B_x\right|<+\infty.$$ Thus, the assumptions of Lemma \ref{Diff-Leb:13-10-22-3} are satisfied and therefore there exists a finite or countable, pairwise disjoint family of balls, $\left\{B_k\right\}_{k\in \Lambda}$ , such that 

\begin{equation*}
	\sum_{k\in \Lambda}\left|B_k\right|\geq 5^{-n}|E_t|.
\end{equation*}
Now, recalling that $B_k\cap  B_j=\emptyset$, for any $j\neq k$, we obtain

\smallskip

\begin{equation*}
	\begin{aligned}
		\int_{\mathbb{R}^n}|f(y)|dy&\geq \int_{\bigcup_{k\in \Lambda}B_k}|f(y)|dy=\\&
		=\sum_{k\in \Lambda}\int_{B_k}|f(y)|dy\geq\\&
		\geq t\sum_{k\in \Lambda}\left|B_k\right|\geq\\&
		\geq 
		5^{-n}t|E_t|.
	\end{aligned}
\end{equation*}
Therefore \eqref{Diff-Leb:16-10-22-1} is proved. $\blacksquare$
 
\bigskip

\textbf{Remark 1.} The function $M(f)$ may take the value $+\infty$, however it is almost everywhere finite. As a matter of fact, by Lemma \ref{Leb:16-10-22-1-0} we get

\smallskip

\begin{equation*}
	\begin{aligned}
	\left|\left\{x\in \mathbb{R}^n: \mbox{ } M(f)(x)=+\infty\right\}\right|&\leq \left|\left\{x\in \mathbb{R}^n: \mbox{ } M(f)(x)>t\right\}\right|\leq \\&
	\\&
	 \leq 
	\frac{5^n}{t}\int_{\mathbb{R}^n}|f(y)|dy	,\quad\forall t>0	
\end{aligned}
\end{equation*}
hence, passing to the limit as $t$ that goes to $+\infty$, we have
\begin{equation*}
	\left|\left\{x\in \mathbb{R}^n: \mbox{ } M(f)(x)=+\infty\right\}\right|=0.
\end{equation*}
$\blacklozenge$

\medskip

\textbf{Remark 2.} Inequality \eqref{Diff-Leb:16-10-22-1}, apart from the value of the constant $5^n$, cannot be improved. To show this it suffices to consider functions \\ $f\in L^1\left(\mathbb{R}^n\right)$ which approximate the Dirac measure concentrated at $0$. For instance 
$$f_{\varepsilon}=\frac{\chi_{B_{\varepsilon}}}{\left|B_{\varepsilon}\right|}.$$ 
Proceeding formally (the reader takes care of the details), we consider $f=\delta(x)$ (the Dirac delta). For this choice, we have
 
 \begin{equation*}
  M(f)(x)=\frac{1}{c_n|x|^n},
  \end{equation*}
 where $c_n$ is the measure of unit ball of $\mathbb{R}^n$. Therefore
 
 $$\left|\left\{x\in \mathbb{R}^n: \mbox{ } M(f)(x)>t \right\}\right|=\frac{1}{t}= \frac{1}{t}\int_{\mathbb{R}^n}|f(y)|dy.$$ 
 $\blacklozenge$
 
 \medskip
 
 \textbf{Remark 3.} Let us observe that, unless in the trivial case where $f$ is identically equal to $0$, we have  $$M(f)\notin  L^1\left(\mathbb{R}^n\right).$$ In this respect, we prove 
 
 \begin{equation}\label{Diff-Leb:16-10-22-4}
 	M(f)(x)\geq \frac{C}{|x|^n},\quad\mbox{for } |x|\geq 1. 
 \end{equation}
 Indeed, since $f$  does not vanish identically, there exists $t_0>0$ such that
 
 $$0<\left|E\right|<+\infty;$$ where
 $$E=\left\{x\in \mathbb{R}^n: \mbox{ } |f(x)|>t_0 \right\}.$$ Let $r_0>0$ satisfy
 
 $$\left|E\cap B_{r_0}\right|\geq \frac{1}{2}|E|.$$ For any $x\in \mathbb{R}^n$, we have (since $B_{r_0}\subset B_{r_{0}+|x|}(x)$)
 
 \begin{equation*}
 	\begin{aligned}
 		M(f)(x)&\geq \frac{1}{\left|B_{r_{0}+|x|}(x)\right|}\int_{B_{r_{0}+|x|}(x)}|f(y)|dy\geq\\&
 		\geq
 			\frac{1}{c_n\left(|x|+r_0\right)^n}\int_{B_{r_{0}}}|f(y)|dy\geq \\&
 			\geq\frac{t_0|E|}{2c_n\left(|x|+r_0\right)^n},
 	\end{aligned}
 \end{equation*}

\smallskip

\noindent from which we have \eqref{Diff-Leb:16-10-22-4} with $C=\frac{t_0|E|}{2c_n\left(1+r_0\right)^n}$.
 
 It can be proved that if $f\in L^p\left(\mathbb{R}^n\right)$, where $1<p\leq +\infty$, then $M(f)\in L^p\left(\mathbb{R}^n\right)$ and 
 
 \begin{equation}\label{7-12-1922}
 \left\Vert M(f)\right\Vert_{L^p\left(\mathbb{R}^n\right)}\leq C \left\Vert f\right\Vert_{L^p\left(\mathbb{R}^n\right)},\quad\forall f\in L^p\left(\mathbb{R}^n\right).
 \end{equation}
For more insights into the maximal function, we refer to \cite[Ch. 1]{ST}. $\blacklozenge$ 
 
\bigskip 

\textbf{Proof of Theorem \ref{Diff-Leb:13-10-22-1}.} Provided that $f$ is replaced by $f\chi_{B_R}$ with arbitrary $R$, we may assume $f\in L^1\left(\mathbb{R}^n\right)$. Let us denote

\begin{equation*}
	f_r(x)=\frac{1}{\left|B_r(x)\right|}\int_{B_r(x)}f(y)dy,\quad\mbox{ } x\in \mathbb{R}^n 
\end{equation*}
and notice that 

\begin{equation*}
	f_r=\varphi_r\star f, 
\end{equation*}
where
$$\varphi_r(x)=r^{-n}\varphi_1\left(r^{-1}x\right),$$

$$\varphi_1=\frac{1}{\left|B_1\right|}\chi_{B_1}.$$
Hence

\begin{equation*}
	\lim_{r\rightarrow 0}\left\Vert f_r-f\right\Vert_{L^1\left(\mathbb{R}^n\right)}=0.
\end{equation*}
Consequently, there exists a sequence $\left\{r_k\right\}$ such that
$$\left\{r_k\right\}\rightarrow 0^+$$ and

\begin{equation}\label{Leb:19-10-22-2}
	\lim_{k\rightarrow \infty}f_{r_k}(x)=f(x),\quad\mbox{ a.e. } x\in \mathbb{R}^n. 
\end{equation}

Now, let us denote

$$\Omega f(x)=\limsup_{r\rightarrow 0}f_r(x)-\liminf_{r\rightarrow 0}f_r(x).$$ 
By Remark 1 we get that $\limsup_{r\rightarrow 0^+}f_r(x)$, $\liminf_{r\rightarrow 0^+}f_r(x)$ are finite almost everywhere. As a matter of fact, we have

\begin{equation}\label{Leb:19-10-22-1}
	\left|\limsup_{r\rightarrow 0}f_r(x)\right|, \ \left|\liminf_{r\rightarrow 0}f_r(x)\right|\leq M(f)(x)<+\infty, \quad \mbox{ a.e. } x\in  \mathbb{R}^n.  
\end{equation} 
 
 Let us now prove that  
\begin{equation}\label{Leb:17-10-22-1}
	\Omega f(x)=0,\quad\mbox{ a.e. } x\in \mathbb{R}^n. 
\end{equation}

\medskip

\textbf{Claim.} If $g\in C^0_0\left(\mathbb{R}^n\right)$, then

\begin{equation*}
g_r\rightarrow g,\quad \mbox{ uniformly as } r\rightarrow 0, 
\end{equation*}
hence

\begin{equation}\label{Leb:18-10-22-1}
	\Omega g(x)=0,\quad\mbox{ } \forall x\in \mathbb{R}^n. 
\end{equation}

\medskip

\textbf{Proof of Claim.} Since $g$ is an uniformly continuous function, for any  $\varepsilon>0$ there exists $\delta>0$ such that if $|x-y|<\delta$ then 

$$|g(x)-g(y)|<\varepsilon.$$
Now

\begin{equation*}
	g_r(x)-g(x)=\int_{\mathbb{R}^n}\varphi_1(z)\left(g(x-rz)-g(x)\right)dz. 
\end{equation*}
Hence, if $0<r<\delta$, we get
\begin{equation*}
	\begin{aligned}
	\left|g_r(x)-g(x)\right|\leq \int_{\mathbb{R}^n}\varphi_1(z)\left|g(x-rz)-g(x)\right|dz\leq \varepsilon,\quad \forall x\in \mathbb{R}^n.
\end{aligned}
\end{equation*}
Claim is proved.

\medskip

Now, let $f\in L^1\left(\mathbb{R}^n\right)$.  Since $C^0_0\left(\mathbb{R}^n\right)$ is dense in $L^1\left(\mathbb{R}^n\right)$ (Theorem \ref{Sob:teo21R}) it follows that for any $\sigma>0$ there exists  $g\in C^0_0\left(\mathbb{R}^n\right)$ such that

\begin{equation*}
	\left\Vert f-g\right\Vert_{L^1\left(\mathbb{R}^n\right)}<\sigma.
\end{equation*}
Let $h=f-g$, we have trivially $f=g+h$, $\left\Vert h\right\Vert_{L^1\left(\mathbb{R}^n\right)}<\sigma$
 and by \eqref{Leb:18-10-22-1} we get
 
 \begin{equation*}
 \Omega f(x)\leq \Omega g(x)+ \Omega h(x)=\Omega h(x),\quad\mbox{ } \forall x\in \mathbb{R}^n.
 \end{equation*}
 Hence, for any $\eta>0$, we have 
 
 \begin{equation}\label{Leb:18-10-22-2}
 	\left|\left\{x\in \mathbb{R}^n: \mbox{ } \Omega f(x)>\eta\right\}\right|\leq \left|\left\{x\in \mathbb{R}^n: \mbox{ } \Omega h(x)>\eta\right\}\right|.
 \end{equation}
 On the other hand, we have trivially 
 \begin{equation*}
 	\Omega h(x)\leq 2M(h)(x),\quad\mbox{ } \forall x\in \mathbb{R}^n,
 \end{equation*}
 this inequality and \eqref{Leb:18-10-22-2} imply
 
 \begin{equation*}
 	\left|\left\{x\in \mathbb{R}^n: \mbox{ } \Omega f(x)>\eta \right\}\right|\leq \left|\left\{x\in \mathbb{R}^n: \mbox{ } M(h)(x)>\frac{\eta}{2}\right\}\right|.
 \end{equation*}
By the latter, by \eqref{Leb:18-10-22-2} and by Lemma \ref{Leb:16-10-22-1-0} we have

 \begin{equation*}
 	\left|\left\{x\in \mathbb{R}^n: \mbox{ } \Omega f(x)>\eta\right\}\right|\leq \frac{2C}{\eta}\left\Vert h\right\Vert_{L^1\left(\mathbb{R}^n\right)}\leq \frac{2C\sigma}{\eta},
 \end{equation*}
 where $C=5^n$. Hence, by choosing
 $$\sigma=\eta^2,$$ we obtain
 
 \begin{equation}\label{Leb:18-10-22-3}
 	\left|\left\{x\in \mathbb{R}^n: \mbox{ } \Omega f(x)>\eta\right\}\right|\leq 2C\eta,
 \end{equation}
 which yields
 \begin{equation*}
 	\left|\left\{x\in \mathbb{R}^n: \mbox{ } \Omega f(x)>0\right\}\right|=\lim_{j\rightarrow \infty}\left|\left\{x\in \mathbb{R}^n: \mbox{ } \Omega f(x)>\frac{1}{j}\right\}\right|=0
 \end{equation*}
 hence, \eqref{Leb:17-10-22-1} follows. Taking into account \eqref{Leb:19-10-22-1}, we have that the limit 
 $$\lim_{r\rightarrow 0}f_r(x)$$
there exists almost everywhere. Therefore, \eqref{Leb:19-10-22-2} implies

$$\lim_{r\rightarrow 0}f_r(x)=\lim_{k\rightarrow \infty}f_{r_k}(x)=f(x),\quad\mbox{ a.e. } x\in \mathbb{R}^n.$$ $\blacksquare$

\bigskip

\begin{cor}\label{Diff-Leb:29-10-22-1} Let $p\in [1,+\infty)$. If $f\in L_{loc}^p\left(\mathbb{R}^n\right)$,  then  
	\begin{equation}\label{Diff-Leb:29-10-22-1-00}
		\lim_{r\rightarrow 0}\ \dashint_{B_r(x)}|f(y)-f(x)|^pdy=0,\quad\mbox{a.e. } x\in \mathbb{R}^n. 
	\end{equation}
\end{cor}
\textbf{Proof.} For any $c\in \mathbb{R}$ let us denote by $D_c$ the subset of $\mathbb{R}^n$ of the points $x$ satisfying

\begin{equation*}
	\lim_{r\rightarrow 0}\left(\dashint_{B_r(x)}|f(y)-c|^pdy\right)^{1/p}=|f(x)-c|. 
\end{equation*} 
Set

\begin{equation*}
	E_c=\mathbb{R}^n\setminus D_c. 
\end{equation*} 
Theorem \ref{Diff-Leb:13-10-22-1} implies 
\begin{equation*}
\left|E_c\right|=0. 
\end{equation*} 
Consequently, setting

$$E=\bigcup_{q\in \mathbb{Q}}E_q,$$ we obtain
\begin{equation*}
	\left|E\right|=0. 
\end{equation*} 
Now, let us prove that 
\begin{equation}\label{Diff-Leb:29-10-22-3}
	\lim_{r\rightarrow 0}\left(\dashint_{B_r(x)}|f(y)-c|^pdy\right)^{1/p}=|f(x)-c|,\quad\forall x\in \mathbb{R}^n\setminus E, \ \forall c\in \mathbb{R}. 
\end{equation}
Let  $c\in \mathbb{R}$ and $x\in \mathbb{R}^n\setminus E$ and let $\delta>0$ and $q\in \mathbb{Q}$ satisfy

\begin{equation}\label{Diff-Leb:29-10-22-4}
|q-c|<\delta. 
\end{equation} 
The triangle inequality gives 

\begin{equation}\label{Diff-Leb:29-10-22-5}
	\begin{aligned}
	\left(\dashint_{B_r(x)}|f(y)-q|^pdy\right)^{1/p}-\delta& <\left(\dashint_{B_r(x)}|f(y)-c|^pdy\right)^{1/p}<\\&<	\left(\dashint_{B_r(x)}|f(y)-q|^pdy\right)^{1/p}+\delta.
	\end{aligned}
	\end{equation} 
Now, let us denote 

\begin{equation*}
	\Lambda'(x)= \liminf_{r\rightarrow 0}  \left(\dashint_{B_r(x)}|f(y)-c|^pdy\right)^{1/p},  
\end{equation*} 

\begin{equation*}
 \Lambda''(x)= \limsup_{r\rightarrow 0}  \left(\dashint_{B_r(x)}|f(y)-c|^pdy\right)^{1/p}.
\end{equation*} 
Passing to the limit as $r\rightarrow 0$, by \eqref{Diff-Leb:29-10-22-5} we obtain 
\begin{equation}\label{Diff-Leb:29-10-22-6}
|f(x)-q|-\delta \leq	\Lambda'(x) \leq \Lambda''(x)\leq |f(x)-q|+\delta. 
\end{equation}
Passing again to the limit as $\delta \rightarrow 0$ in \eqref{Diff-Leb:29-10-22-6} and taking into account \eqref{Diff-Leb:29-10-22-4}, we obtain  
$$\Lambda'(x)=\Lambda''(x)=|f(x)-c|,\quad\forall x\in \mathbb{R}^n\setminus E.$$ Therefore \eqref{Diff-Leb:29-10-22-3} holds true. Set therein $c=f(x)$ and we obtain \eqref{Diff-Leb:29-10-22-1-00}. $\blacksquare$

\section{The Rademacher Theorem} \label{funz-Lips}

Let us recall the definition of \textbf{absolutely continuous function} over the interval  $[a,b]$, where $a,b\in \mathbb{R}$ and $a<b$.  

We say that the function 

$$f: [a,b]\ \rightarrow \ \mathbb{R},$$
is absolutely continuous provided that for every $\varepsilon>0$ there exists $\delta>0$ such that, chosen anyway a finite family of pairwise disjoint intervals $(a_j,b_j)\subset [a,b]$, $j=1,\cdots, N$ satisfying 
$$\sum_{j=1}^{N}(b_j-a_j)<\delta,$$ we have

$$\sum_{j=1}^{N}\left|f\left(b_j\right)-f\left(a_j\right)\right|<\varepsilon.$$
We will denote by $AC\left([a,b]\right)$ the class of the
absolutely continuous functions \index{absolutely continuous functions} on $[a,b]$ and we will denote by
$AC_{loc}\left(\mathbb{R}\right)$ the class of functions
$f:\mathbb{R}\rightarrow\mathbb{R}$ such that for every interval
$[a,b]$ we have $f_{|[a,b]}\in
AC\left([a,b]\right)$. Let us recall that if $f\in AC_{loc}\left(\mathbb{R}\right)$, then $f$ is a differentiable function almost everywhere in $\mathbb{R}$. 

If $f$ is a Lipschitz continuous function in $\mathbb{R}$, then  $f\in AC_{loc}\left(\mathbb{R}\right)$, hence $f$ is a differentiable almost every everywhere. The main purpose of the present Section is to extend this result to several variable Lipschitz continuous functions. Precisely we want to prove
\begin{theo}[\textbf{Rademacher}]\label{Rade}
	\index{Theorem:@{Theorem:}!- Rademacher@{- Rademacher}} If $f\in C^{0,1}\left(\mathbb{R}^n\right)$, then $f$ is a differentiable function almost everywhere.
\end{theo}
\textbf{Proof.} Let $v\in \mathbb{R}^n$ satisfy $|v|=1$. Set 
\begin{equation*}
	\partial_vf(x)=\lim_{t\rightarrow 0} \frac{f(x+tv)-f(x)}{t},\quad \mbox{ provided the limit exists; }
\end{equation*}
notice that, since $f$ is Lipschitz continuous, if the limit above exists it is finite.

\medskip

\textbf{Claim 1.}   
 \begin{equation*}
 	\partial_vf(x),\quad \mbox{ exists a.e. } x\in \mathbb{R}^n.
 \end{equation*}

\bigskip

 \textbf{Proof of Claim 1.} Let us denote 
 
\begin{equation*}
	\begin{aligned}
\overline{\partial}_vf(x):=\limsup_{t\rightarrow 0} \frac{f(x+tv)-f(x)}{t}
\end{aligned}
\end{equation*} 
\\ 
\begin{equation*}
	\begin{aligned}
		\underline{\partial}_vf(x)&:=\liminf_{t\rightarrow 0} \frac{f(x+tv)-f(x)}{t}.
	\end{aligned}
\end{equation*} 
Since $f$ is a continuous function, we have
$$\overline{\partial}_vf(x)=\lim_{k\rightarrow \infty} \ \sup_{0<|t|<1/k, \ t\in \mathbb{Q}}\frac{f(x+tv)-f(x)}{t}$$
and 
$$\underline{\partial}_vf(x)=	\lim_{k\rightarrow \infty} \ \inf_{0<|t|<1/k, \ t\in \mathbb{Q}}\frac{f(x+tv)-f(x)}{t}.$$
hence $\overline{\partial}_vf$ and $\underline{\partial}_vf$ are measurable functions. Consequently

 \begin{equation*}
 	A_v:=\left\{x\in \mathbb{R}^n:\mbox{ } \underline{\partial}_vf(x)<\overline{\partial}_vf(x) \right\},
 \end{equation*}

  $$\ell_v:=\left\{tv:\mbox{ } t \in \mathbb{R} \right\}$$ are Lebesgue measurable sets. 
 
  Let us first consider the case $v=e_n$. The Fubini--Tonelli Theorem gives 
 
 \begin{equation}\label{Diff:19-10-22-3}
 	\begin{aligned}
 		\left|A_{e_n}\right|&=\int_{\mathbb{R}^n}\chi_{A_{e_n}}(x)dx=\\&=
 		\int_{\mathbb{R}^{n-1}}dx'\int_{R}\chi_{A_{e_n}}\left(x',x_n\right)dx_n=\\&=
 		\int_{\mathbb{R}^{n-1}}\left|\left(x'+\ell_{e_n}\right)\cap A_{e_n}\right|_1dx',
 	\end{aligned}
 \end{equation}

\smallskip

\noindent where $\left|\left(x'+\ell_{e_n}\right)\cap A_{e_n}\right|_1$ is the Lebesgue measure on $\mathbb{R}$ of $\left(x'+\ell_{e_n}\right)\cap A_{e_n}$. For any fixed $x'\in \mathbb{R}^{n-1}$, we have
 \begin{equation*}
 	\left(x'+\ell_{e_n}\right)\cap A_{e_n}=\{x'\}\times \left\{t\in \mathbb{R}:\mbox{ } \underline{\partial}_nf(x',t)<\overline{\partial}_nf(x',t) \right\}.
 \end{equation*}
 Now, by denoting
 $$\varphi(t)=f\left(x',t\right),$$ since $f$ is Lipschitz continuous, we have $\varphi\in AC_{loc}\left(\mathbb{R}\right)$. In particular, the function $\varphi$ is almost everywhere differentiable, hence 
 $$\left|\left(x'+\ell_{e_n}\right)\cap A_{e_n}\right|_1=0.$$ Therefore, by \eqref{Diff:19-10-22-3}, we get 
 
 $$\left|A_{e_n}\right|=0.$$
 Whenever $v\neq e_n$, let us consider a rotation $\mathcal{R}$  of $\mathbb{R}^n$ such that $$v=\mathcal{R}\left(e_n\right).$$ By setting 
 
 $$\tilde{f}=f\circ \mathcal{R},$$ we get (due to the invariance of the Lebesgue measure with respect to rotations)
 
 $$\left|A_v\right|=\left|\mathcal{R}^{-1}\left(A_v\right)\right|.$$ Moreover, it is easily checked that 
 $$ \mathcal{R}^{-1}\left(A_v\right)=\left\{y\in \mathbb{R}^n:\mbox{ } \underline{\partial}_n \tilde{f}(y)<\overline{\partial}_n\tilde{f}(y) \right\}.$$ Consequently, since $\tilde{f}$ is Lipschitz continuous, from what we have previously proved, we derive
 
$$\left|A_v\right|=\left|\mathcal{R}^{-1}\left(A_v\right)\right|=\left|\left\{y\in \mathbb{R}^n:\mbox{ } \underline{\partial}_n \tilde{f}(y)<\overline{\partial}_n\tilde{f}(y) \right\}\right|=0$$ that concludes the proof of Claim 1.
 
 \bigskip
 
 Set
 
 $$\nabla f(x)=\left(\partial_1 f(x),\cdots, \partial_n f(x) \right),\quad \mbox{ a.e. } x\in \mathbb{R}^n.$$
 
 \medskip
 
\textbf{Claim 2.} We have

\begin{equation}\label{Diff:21-10-22-1}
\partial_vf(x)=v\cdot \nabla f(x),\quad \mbox{ a.e. } x\in \mathbb{R}^n
\end{equation}
and
\begin{equation}\label{Diff:21-10-22-2}
|\nabla f(x)|\leq L,\quad\mbox{ a.e. } x\in \mathbb{R}^n,
\end{equation}
where
$$L=[f]_{0,1,\mathbb{R}^n},$$

\medskip
 
 \textbf{Proof of Claim 2.} Let $\zeta \in C^{\infty}_0\left(\mathbb{R}^n\right)$ be arbitrary. We easily get
 \begin{equation}\label{Diff:20-10-22-1}
 	\int_{\mathbb{R}^n}\frac{f(x+tv)-f(x)}{t}\zeta(x)dx=-\int_{\mathbb{R}^n}f(x)\frac{\zeta(x)-\zeta(x-tv)}{t}dx.
 \end{equation}
 Now, since 
 
 $$\lim_{t\rightarrow 0}\frac{f(x+tv)-f(x)}{t}=\partial_vf(x),\quad \mbox{ a.e. } x\in \mathbb{R}^n$$ and
 
 $$\left|\frac{f(x+tv)-f(x)}{t}\zeta(x)\right|\leq L|\zeta(x)|,\quad \forall x\in \mathbb{R}^n, \ \ \forall t\in \mathbb{R}\setminus\{0\},$$ we have by Dominated Convergence Theorem and by \eqref{Diff:20-10-22-1}, 
 
\begin{equation}\label{Diff:20-10-22-2}
	\begin{aligned}
	\int_{\mathbb{R}^n}\partial_vf(x)\zeta(x)dx&=\lim_{t\rightarrow 0}\int_{\mathbb{R}^n}\frac{f(x+tv)-f(x)}{t}\zeta(x)dx=\\&=-\lim_{t\rightarrow 0}\int_{\mathbb{R}^n}f(x)\frac{\zeta(x)-\zeta(x-tv)}{t}dx=\\&=
	-\sum_{j=1}^nv_j\int_{\mathbb{R}^n}f(x)\partial_j\zeta(x)dx=\\&=\sum_{j=1}^nv_j\int_{\mathbb{R}^n}\partial_jf(x)\zeta(x)dx.
	\end{aligned}
\end{equation}

\smallskip

\noindent Let us justify the last equality. For any $j\in\{1,\cdots,n\}$ we have that the function

$$x_j \ \rightarrow \ f(x_1,\cdots, x_{j-1},x_j, x_{j+1},\cdots,x_n),$$
belongs to $AC_{loc}(\mathbb{R})$. Therefore, by considering, for instance, the case $j=n$ (the other cases are similar), we have by Fubini--Tonelli Theorem

\begin{equation*}
	\begin{aligned}
		&-\int_{\mathbb{R}^n}f(x)\partial_n\zeta(x)dx=-\int_{\mathbb{R}^{n-1}}dx'\int_{\mathbb{R}}f(x',x_n)\partial_n\zeta(x',x_n)dx_n=\\&=-\int_{\mathbb{R}^{n-1}}dx'\int_{\mathbb{R}}\left(\partial_n\left(f(x',x_n)z(x',x_n)\right)-\partial_n f(x',x_n)\zeta(x',x_n)\right)dx_n=\\&=
		\int_{\mathbb{R}^{n-1}}dx'\int_{\mathbb{R}}\partial_n f(x',x_n)\zeta(x',x_n)dx_n=\\&=
		\int_{\mathbb{R}^n}\partial_n f(x)\zeta(x)dx.
	\end{aligned}
\end{equation*}

\smallskip

\noindent Consequently, \eqref{Diff:20-10-22-2} gives

\begin{equation*}
	\int_{\mathbb{R}^n}\partial_vf(x)\zeta(x)dx=\int_{\mathbb{R}^n}\left(v\cdot \nabla f(x)\right)\zeta(x)dx,\quad\forall \zeta \in C^{\infty}_0\left(\mathbb{R}^n\right)
\end{equation*}
which yields \eqref{Diff:21-10-22-1}. Concerning \eqref{Diff:21-10-22-2}, it is an immediate consequence of \eqref{Diff:21-10-22-1} and of the Cauchy--Schwarz inequality. The proof of Claim 2 is concluded.

\bigskip

Let $D=\left\{v_k\right\}_{k\in \mathbb{N}}$ be such that $$\left|v_k\right|=1,\quad \forall k\in \mathbb{N}$$ and
$$\overline{D}=\partial B_1.$$  
Moreover, let 

$$x\in \mathbb{R}^n\setminus \bigcup_{k=1}^{\infty}A_{v_{k}}.$$ Let us prove that $f$ is  differentiable in $x$ whereby we will conclude the proof of Thheorem, because

$$\left|\bigcup_{k=1}^{\infty}A_{v_{k}}\right|\leq \sum_{k=1}^{\infty}\left|A_{v_{k}}\right|=0.$$

Let $y\in \mathbb{R}^n\setminus \{x\}$ and set 
$$w=\frac{y-x}{|y-x|},\quad\quad  t=|y-x|,$$
we have trivially
$$y=x+tw.$$
Moreover, for any $k\in \mathbb{N}$ we have

\begin{equation}\label{Diff:21-10-22-3}
	\begin{aligned}
&\left|f(y)-f(x)-\nabla f(x)\cdot (y-x)\right|=\\&=\left|f(x+tw)-f(x)-t\nabla f(x)\cdot w\right|\leq\\&\leq\left|f(x+tv_k)-f(x)-t\nabla f(x)\cdot w \right|+\\&+\left|f(x+tw)-f(x+tv_k)\right|\leq \\&\leq 
\left|f(x+tv_k)-f(x)-t\nabla f(x)\cdot v_k \right|+\\&+
t\left|\nabla f(x)\right|\left|v_k-w\right|+\\&+\left|f(x+tw)-f(x+tv_k)\right|\leq \\&\leq 
t\left|\frac{f(x+tv_k)-f(x)-t\nabla f(x)\cdot v_k}{t}\right|+\\&+2L\left|w-v_k\right|t.
	\end{aligned}
\end{equation}

\bigskip

Let $\varepsilon>0$, since $D$ is dense in $\partial B_1$, there exists $k_{\varepsilon}$ such that

$$\left|w-v_{k_{\varepsilon}}\right|<\frac{\varepsilon}{2\left(2L+1\right)}.$$ Since 

$$x\in \mathbb{R}^n\setminus \bigcup_{k=1}^{\infty}A_{v_{k}},$$ we have

$$\lim_{\tau\rightarrow 0}\frac{f(x+\tau v_{k_{\varepsilon}})-f(x)}{\tau}=\nabla f(x)\cdot v_{k_{\varepsilon}}.$$ Therefore, there exists $\delta>0$ such that, if $0<|\tau|<\delta$ then

$$\left|\frac{f(x+\tau v_{k_{\varepsilon}})-f(x)-\tau \nabla f(x)\cdot v_{k_{\varepsilon}}}{\tau}\right|<\frac{\varepsilon}{2}. $$ Now taking into account  \eqref{Diff:21-10-22-3} (and $t=|x-y|$), we have
$$\left|f(y)-f(x)-\nabla f(x)\cdot (y-x)\right|<\varepsilon |x-y|,\quad \forall y\in B_{\delta}(x)\setminus\{x\},$$ which gives the differentiability of $f$ in $x$. $\blacksquare$

 \section{Description of the boundary of an open set of $\mathbb{R}^n$}\label{DescrBordo}
 Let $r>0$, $x\in\mathbb{R}^{n}$ and $x'\in\mathbb{R}^{n-1}$,
 we denote by $B_r(x)$ and $B'_r(x')$, the open ball of $\mathbb{R}^{n}$ centered in $x$ with radius $r$ and open ball of $\mathbb{R}^{n-1}$ centered in $x'$ with radius $r>0$ respectively.
 We will also write $B_r$ ($B'_r$) instead of $B_r(0)$ ($B'_r(0)$).
 For any $r, M>0$ and any $x\in \mathbb{R}^n$, here and in the sequel, we denote by 
 $$Q_{r,M}(x)=B_r'(x')\times (-Mr+x_n,Mr+x_n).$$ We will write also \index{$Q_{r,M}$, $Q_{r,M}(P)$}$Q_{r,M}$ instead of $Q_{r,M}(0)$.  
 
 \begin{definition}\label{Sob:def3.1}
 	\index{Definition:@{Definition:}!- boundary of class $C^m$, $C^{m,\alpha}$@{- boundary of class $C^m$, $C^{m,\alpha}$}}
 	
 	Let $\Omega$ be an open set of $\mathbb{R}^n$. Let $r_0, M_0$ be positive numbers and $m\in \mathbb{N}_0$.
 	
 	\smallskip
 	
 	(a) We say that $\Omega$ has the \textbf{boundary of $C^m$ class with constants}
 	$r_0, M_0$ (or, briefly, $\Omega$ is of class $C^m$ with constants $r_0, M_0$), if for every $P\in \partial \Omega$ there exists
 	an isometry $$ \Phi_P:\mathbb{R}^n\rightarrow \mathbb{R}^n,$$ such that $$\Phi_P(0)=P$$ and
 	\begin{equation}\label{Sob:for1.20r}
 	\Phi_P^{-1}	\left(\Omega\right)\cap Q_{r_0,2M_0}=\left\{x\in Q_{r_0,2M_0}:\mbox{ }
 		x_n>g_P(x') \right\},
 	\end{equation}
 	where $g_P\in C^m\left(\overline{B'_{r_0}}\right)$,

$$g_P(0)=0,\quad |\nabla g_P(0)|=0, \mbox{ for } m\geq 1$$
 	and
 	$$\left\Vert g_P\right\Vert_{C^m\left(\overline{B'_{r_0}}\right)}\leq
 	M_0r_0,$$ where
 	$$\left\Vert g_P\right\Vert_{C^m\left(\overline{B'_{r_0}}\right)}=\sum_{|\gamma|\leq m}
 	r_0^{|\gamma|}\left\Vert\partial^{\gamma}g_P\right\Vert_{L^{\infty}\left(B'_{r_0}\right)}.$$

 	\smallskip
 	
 	(b) Let $ \alpha\in (0,1]$. We say that $\Omega$ has the \textbf{boundary of $C^{m,\alpha}$ class with constants  $r_0, M_0$} (or, briefly,
 	$\Omega$ is of class $C^{m,\alpha}$ with constants $r_0, M_0$) if
 	for every $P\in
 	\partial \Omega$ there exists
 	an isometry $$ \Phi_P:\mathbb{R}^n\rightarrow \mathbb{R}^n,$$ such that $$\Phi_P(0)=P,$$ and
 	\begin{equation}\label{Sob:for1.20}
 	\Phi_P^{-1}\left(\Omega\right)\cap Q_{r_0,2M_0}	=\left\{x\in Q_{r_0,2M_0}:\mbox{ }
 		x_n>g_P(x') \right\},
 	\end{equation}
 	where $g_P\in C^m\left(\overline{B'_{r_0}}\right)$,
 	$$g_P(0)=0, \quad |\nabla g_P(0)|=0\quad\mbox{ (for } m\geq 1)$$
 	and
 	$$\left\Vert g_P\right\Vert_{C^{m,\alpha}\left(\overline{B'_{r_0}}\right)}\leq
 	M_0r_0$$ where
 	$$\left\Vert g_P\right\Vert_{C^{m,\alpha}\left(\overline{B'_{r_0}}\right)}=\left\Vert g_P\right\Vert_{C^m\left(\overline{B'_{r_0}}\right)}+r_0^{m}\sum_{|\gamma|= m}
 	r_0^{|\alpha|}\left[\partial^{\gamma}u\right]_{C^{0,\alpha}\left(\overline{B'_{r_0}}\right)}.$$ 	
 	(c) If there exists $r_0>0$ e $M_0>0$ such that $\partial \Omega$ has the boundary of class $C^{m}$ ($C^{m,\alpha}$) with constants $r_0>0$ e $M_0>0$, then we say that $\partial \Omega$ has the boundary of class $C^{m}$ ($C^{m,\alpha}$) .  
 	\end{definition}

 \underline{\textbf{Exercise.}} Prove that if $\Omega$ is an open set of class $C^0$ then $\overset{\circ }{\overline{\Omega }}=\Omega$. $\clubsuit$
 
 \bigskip
 
 Let us note that the graph of the function $g_P$ that occurs in the definition above is contained in $Q_{r_0,M_0}$ (see Figure 2.1).

 \begin{figure}\label{figura-p21sob}
 	\centering
 	\includegraphics[trim={0 0 0 0},clip, width=15cm]{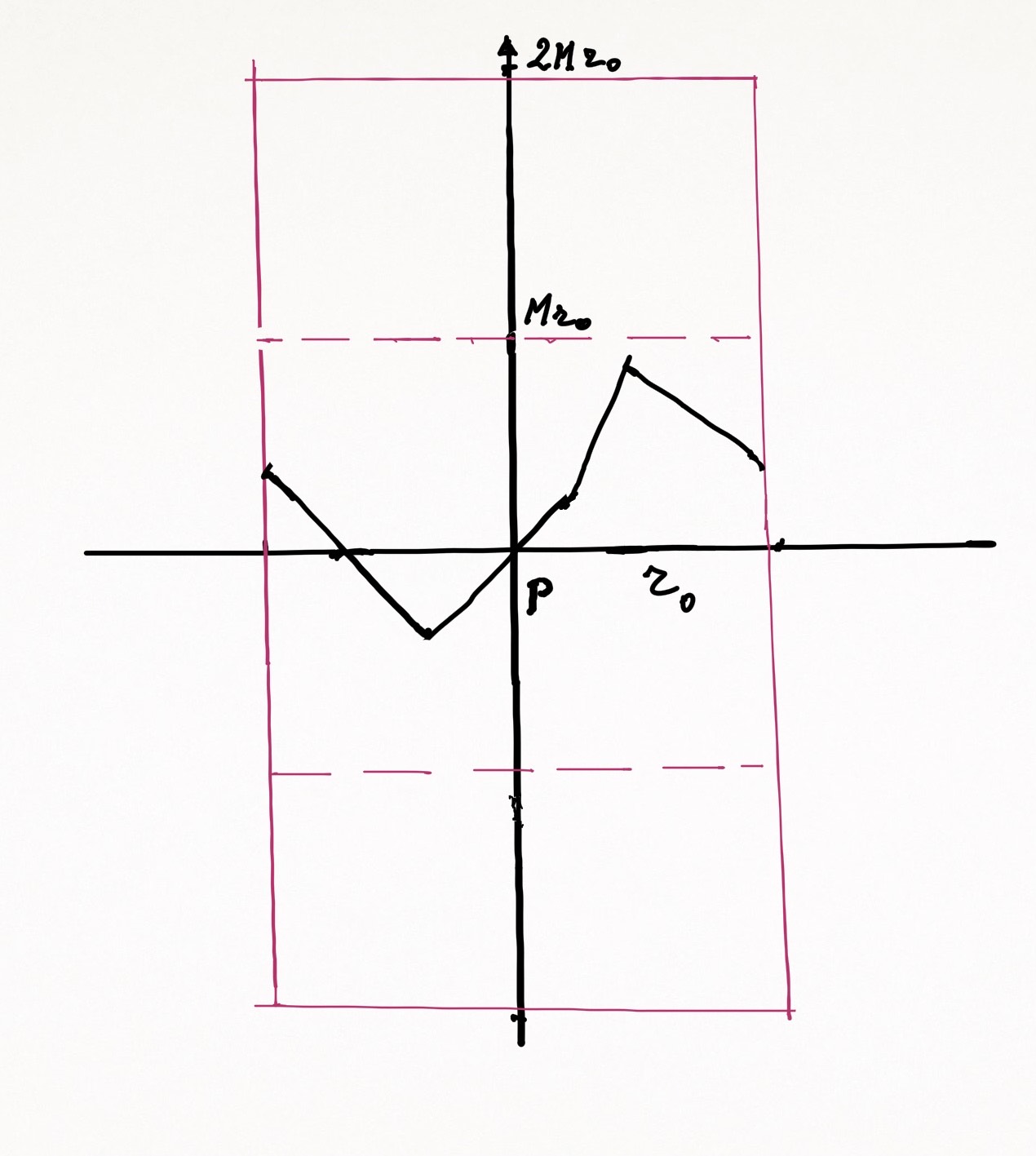}
 	\caption{}
 \end{figure}

\newpage

\section{The spaces $L^p(\partial \Omega)$}\label{correc-11-3-23-1}
We now provide a brief review of the definition of $L^p(\partial
\Omega)$ \index{$L^p(\partial \Omega)$} spaces, where $\Omega$ is a bounded open set of $\mathbb{R}^n$ of class $C^{0,1}$. For any $P\in \partial \Omega$ let us denote by $\Phi_P$ an isometry which satisfies \eqref{Sob:for1.20}. Let 
$$f:\partial\Omega\rightarrow \mathbb{R}.$$ We have that 
$\mathcal{U}=\left\{\Phi_{P_j}\left(Q_{r_0,2M_0}\right)\right\}_{1\leq j\leq l}$ is a finite open covering of $\partial \Omega$. Let us denote
$$g_j=g_{P_j}, \ \ j=1,\cdots, l,$$ and let $\zeta_1,\cdots,\zeta_l$ be a partition of unity subordinate to the covering $\mathcal{U}$. we say that $f\in L^1(\partial \Omega)$ provided that we have, for any $1\leq j\leq l$,  $$f_j\in L^1\left(B'_{r_0}\right),$$
where   
$$f_j(x')=\left(f\circ \Phi_{P_j}\right)\left(x',g_j(x')\right), \ \ \forall x'\in B'_{r_0}.$$

Let us denote 

\begin{equation}\label{int:1-11-22}
	\int_{\partial \Omega}fdS=\sum_{j=1}^l\int_{\partial \Omega}f\zeta_j dS,
	\end{equation}
where

\begin{equation}\label{Bordo:25-10-22-1}
\int_{\partial \Omega}f\zeta_j dS=\int_{B'_{r_0}}f_j\left(x'\right)\tilde{\zeta}_j(x')\sqrt{1+|\nabla g_j(x')|^2}dx',
\end{equation} 

$$\tilde{\zeta}_j(x')=\zeta_j\left(\Phi_{P_j}\left(x',g_j(x')\right)\right).$$

Concerning the last integral in \eqref{Bordo:25-10-22-1}, take into account that it is well defined,  because $g\in C^{0,1}\left(\overline{B'_{r_0}}\right)$ is differentiable almost everywhere and its gradient belongs to $L^{\infty}\left(B'_{r_0}; \mathbb{R}^n\right)$. Let us observe that  integral in \eqref{int:1-11-22} does not depend on the particular partition of unity that we use. As a matter of fact, if $\eta_1,\cdots,\eta_m$ is another partition of unity, then we have

\begin{equation*}
	\begin{aligned}
	\sum_{j=1}^l\int_{\partial \Omega}f\zeta_j dS&=	\sum_{j=1}^l\int_{\partial \Omega}\left(\sum_{k=1}^m\eta_k f\right)\zeta_j dS=\\&=\sum_{j=1}^l\sum_{k=1}^m\int_{\partial \Omega}\zeta_j\eta_k f dS=\\&=\sum_{k=1}^m\int_{\partial \Omega}\left(\sum_{j=1}^l \zeta_j f\right)\eta_kdS=\\&=\sum_{k=1}^m\int_{\partial \Omega} f\eta_k dS.
\end{aligned}
\end{equation*}

\medskip

Let us denote

\begin{equation*}
	\left\Vert
	u\right\Vert_{L^1(\partial\Omega)}=\int_{\partial \Omega}|f|
	dS.
\end{equation*}
Likewise, we define $L^p(\partial\Omega)$ for $1\leq p< +\infty$ and we set \index{$\left\Vert
	\cdot \right\Vert_{L^p(\partial\Omega)}$}

\begin{equation}\label{norma L-2 bordo}
	\left\Vert
	u\right\Vert_{L^p(\partial\Omega)}=\left(\int_{\partial \Omega}|f|^p
	dS\right)^{1/p}.
\end{equation}
The space $L^p(\partial\Omega)$ is
a separable Banach space, and if $p=2$, $L^2(\partial\Omega)$ is
a Hilbert space. 

Let $g\in C^{0,1}\left(\overline{B'_{r_0}}\right)$ such that

$$\left\Vert g\right\Vert_{C^{0,1}\left(\overline{B'_{r_0}}\right)}\leq
M_0r_0.$$ Let us consider the set

$$W=\left\{y\in Q_{r_0,2M_0}:\mbox{ }
y_n>g(y') \right\};$$
We define the field of unit outward normal \index{outward normal}to the graph of $g$ as follows

\begin{equation}\label{normal:26-10-22}
\nu^{(g)}(y',g(y'))=\frac{\left(\nabla_{y'}g(y'),-1\right)}{\sqrt{1+\left|\nabla_{y'}g(y')\right|^2}},\quad \mbox{a.e. }y'\in B'_{r_0}.
\end{equation}
Let $\Omega$ be a bounded open set of $\mathbb{R}^n$ of $C^{0,1}$ class; let $\mathcal{U}=\left\{\Phi_{P_j}\left(Q_{r_0,2M_0}\right)\right\}_{1\leq j\leq l}$ the open covering defined above. We define the field of unit outward normal on $\partial\Omega\cap \Phi_{P_j}\left(Q_{r_0,2M_0}\right)$, $j=1,\cdots, l$, as

\begin{equation}\label{Bordo:25-10-22-2}
	\nu(x)= \left(\Phi_{P_j}-P_j\right)\left(\nu^{(g_j)}(y',g_j(y'))\right),\quad x=\Phi_j\left(y',g_j(y')\right).
\end{equation}
If $f\in L^{1}(\partial\Omega)$, we define, by \eqref{Bordo:25-10-22-1},  

\begin{equation}\label{Bordo:25-10-22-3}
	\int_{\partial \Omega}f\nu dS=\sum_{j=1}^l\int_{\partial \Omega}f \zeta_j \nu dS.
\end{equation}

\section{The divergence Theorem} \label{TeoDiv:25-10-22}
The purpose of this Section is to prove the following 

\begin{theo}[\textbf{divergence}]\label{TeoDiv:25-10-22-1}
	\index{Theorem:@{Theorem:}!- divergence@{- divergence}}
	Let $\Omega$ be a bounded open set of $\mathbb{R}^n$ of class $C^{0,1}$ and let $f\in C^{0,1}\left(\overline{\Omega}\right)$, then we have
	
\begin{equation}\label{TeoDiv:25-10-22-2}
\int_{\Omega}\nabla f dx=\int_{\partial\Omega} f \nu dS.
\end{equation}
\end{theo}

\medskip

In order to prove Theorem \ref{TeoDiv:25-10-22-1} we need some lemmas.

\bigskip

\begin{lem}\label{TeoDiv:25-10-22-3}
	Let $\Omega$ be a bounded open set of $\mathbb{R}^n$ and let $f\in C^{0,1}\left(\overline{\Omega}\right)$. Then there exists a sequence $\left\{f_m\right\}$ in $C^{\infty}\left(\mathbb{R}^n\right)$ satisfying 
	\begin{equation}\label{TeoDiv:25-10-22-4}
	\lim_{m\rightarrow \infty}\left\Vert f_m-f\right\Vert_{C^{0}\left(\overline{\Omega}\right)}=0,
	\end{equation}

\begin{equation}\label{TeoDiv:25-10-22-5}
	\lim_{m\rightarrow \infty}\left\Vert \nabla f_m-\nabla f\right\Vert_{L^{2}\left(\Omega\right)}=0.
\end{equation}
\end{lem}
\textbf{Proof.} Set

$$L=\left[f\right]_{0,1, \overline{\Omega}}.$$
Let $\tilde{f}\in C^{0,1}\left(\mathbb{R}^n\right)$ an extension of $f$ which satisfies  (see Theorem \ref{Contin:teo-est})

$$\left[\tilde{f}\right]_{0,1, \mathbb{R}^n}=L.$$
For any $\varepsilon>0$, let us denote  
$$f_{\varepsilon}(x)=\int_{\mathbb{R}^n}\eta_{\varepsilon}(x-y)\tilde{f}(y)dy,$$
where $\eta$ is a mollifier. We have that $f_{\varepsilon}\in C^{\infty}\left(\mathbb{R}^n\right)$. Moreover, if $x\in \Omega$, we get  
\begin{equation*}
	\begin{aligned}
		\left|f_{\varepsilon}(x)-f(x)\right|&=\left|\int_{\mathbb{R}^n}\eta_{\varepsilon}(y) \left(\tilde{f}(x-y)-\tilde{f}(x)\right)dy\right|\leq \\&\leq \int_{\mathbb{R}^n}\eta_{\varepsilon}(y)\left|\tilde{f}(x-y)-\tilde{f}(x)\right|dy\leq\\&\leq L\int_{\mathbb{R}^n}\eta_{\varepsilon}(y)|y|dy\leq\\&\leq
		L\varepsilon.
	\end{aligned}
\end{equation*}

\smallskip

\noindent Therefore 

$$\lim_{\varepsilon\rightarrow 0}\left\Vert f_{\varepsilon}-f\right\Vert_{C^{0}\left(\overline{\Omega}\right)}=0,$$ which gives \eqref{TeoDiv:25-10-22-4}.

\medskip

Let us prove \eqref{TeoDiv:25-10-22-5}. Let $1\leq k\leq n$. Theorem \ref{Rade} gives 
\begin{equation*}
	\begin{aligned}
		\partial_{x_k}f_{\varepsilon}(x)&=\int_{\mathbb{R}^n}\partial_{x_k}\left(\eta_{\varepsilon}(x-y)\right)\tilde{f}(y)dy=\\&=-\int_{\mathbb{R}^n}\partial_{y_k}\left(\eta_{\varepsilon}(x-y)\right)\tilde{f}(y)dy=\\&=\int_{\mathbb{R}^n}\eta_{\varepsilon}(x-y)\partial_{y_k}\tilde{f}(y)dy.
	\end{aligned}
\end{equation*}
Now, for any $x\in \Omega$ we get 

\begin{equation*}
	\begin{aligned}
		\left|\partial_{x_k}f_{\varepsilon}(x)-\partial_{x_k}f(x)\right|&\leq \int_{\mathbb{R}^n}\eta_{\varepsilon}(y)\left|\left(\partial_{x_k}\tilde{f}\right)(x-y)-\partial_{x_k}\tilde{f}(x)\right|dy\leq \\&\\&\leq \left( \int_{\mathbb{R}^n}\eta_{\varepsilon}(y)\left|\left(\partial_{x_k}\tilde{f}\right)(x-y)-\partial_{x_k}\tilde{f}(x)\right|^2dy\right)^{1/2}.
	\end{aligned}
\end{equation*}
From which we derive
\begin{equation*}
	\begin{aligned}
		\int_{\Omega}\left|\partial_{x_k}f_{\varepsilon}(x)-\partial_{x_k}f(x)\right|^2dx&\leq \int_{\Omega}dx\int_{B_{\varepsilon}}\eta_{\varepsilon}(y)\left|\left(\partial_{x_k}\tilde{f}\right)(x-y)-\partial_{x_k}\tilde{f}(x)\right|^2dy= \\&=
		\int_{B_{\varepsilon}}\eta_{\varepsilon}(y)\int_{\Omega}\left|\left(\partial_{x_k}\tilde{f}\right)(x-y)-\partial_{x_k}\tilde{f}(x)\right|^2dx\leq \\&\leq 
		\sup_{|y|<\varepsilon}\int_{\Omega}\left|\left(\partial_{x_k}\tilde{f}\right)(x-y)-\partial_{x_k}\tilde{f}(x)\right|^2dx.
	\end{aligned}
\end{equation*}

\smallskip

\noindent By the just proved inequality and by Theorem \ref{Sob:teo24R} we have
$$\lim_{\varepsilon\rightarrow 0}\left\Vert \nabla f_{\varepsilon}-\nabla f\right\Vert_{L^{2}\left(\Omega\right)}=0,$$ that gives \eqref{TeoDiv:25-10-22-4}. $\blacksquare$

\medskip 

For any $f\in C^{0,1}\left(\overline{\Omega}\right)$ we say that a sequence $\left\{f_m\right\}$ which satisfies \eqref{TeoDiv:25-10-22-4} and \eqref{TeoDiv:25-10-22-5} is a \textbf{smooth approximating sequence} of $f$ \index{smooth approximating sequence of a Lipschitz function}.

\bigskip

\begin{lem}\label{Lemma:26-10-22}
Let $r, h$ be positive numbers. Let $g\in C^{\infty}\left(\overline{B'_r}\right)$ satisfy 
$$-h<g(x')<h, \quad \forall x'\in B'_{r}.$$ Let $f\in C^{\infty}_0\left(Q_{r,h}\right)$, where $Q_{r,h}=B'_r\times (-h,h)$. Denoting

$$W=\left\{x\in Q_{r,h}:\mbox{ }
x_n>g(x') \right\},$$ we have
\begin{subequations}
	\label{Lemma:26-10-22-1}
	\begin{equation}
		\label{Lemma:26-10-22-1-a}
		\int_{W} \partial_n fdx=-\int_{B'_r}f\left(x',g(x')\right)dx',
	\end{equation}
	\begin{equation}
		\label{Lemma:26-10-22-1-b}
		\int_{W} \partial_k fdx=\int_{B'_r}f\left(x',g(x')\right)\partial_kg(x')dx', \quad k=1,\cdots,n-1. 
	\end{equation}
	\end{subequations}
\end{lem}

\textbf{Proof.} Let $1\leq k\leq n$. If $k=n$, then we have

\begin{equation*}
	\int_{g(x')}^h \partial_nf\left(x',x_n\right)dx_n=-f\left(x',g(x')\right),\quad x'\in B'_r.
\end{equation*}
Hence 
\begin{equation}\label{Lemma:26-10-22-2}
	\begin{aligned}
	\int_{W} \partial_nfdx&=\int_{B'_r}dx'\int_{g(x')}^h\partial_nf(x',x_n)dx_n=\\&=
	-\int_{B'_r}f\left(x',g(x')\right)dx'.
	\end{aligned}
\end{equation}
If $1\leq k\leq n-1$, then we have 

\begin{equation*}
		\partial_k\int_{g(x')}^hf\left(x',x_n\right)dx_n=\int_{g(x')}^h \partial_kf\left(x',x_n\right)dx_n-f\left(x',g(x')\right)\partial_kg(x'),
\end{equation*}
from which we get

\begin{equation}\label{Lemma:26-10-22-3}
	\begin{aligned}
		\int_{W} \partial_kfdx&=\int_{B'_r}dx'\int_{g(x')}^h\partial_kf(x',x_n)dx_n=\\&=
		\int_{B'_r}dx'\left(\partial_k\int_{g(x')}^hf\left(x',x_n\right)dx_n\right)+\\&+
		\int_{B'_r}f\left(x',g(x')\right)\partial_kg(x')dx'=\\&=\int_{B'_r}f\left(x',g(x')\right)\partial_kg(x')dx,
	\end{aligned}
\end{equation}

\medskip

\noindent where in the fourth step we used the fact that the function
$$x'\rightarrow \  \int_{g(x')}^hf\left(x',x_n\right)dx_n,$$ has the support contained in $B'_r$. From what was obtained in \eqref{Lemma:26-10-22-2} and by \eqref{Lemma:26-10-22-3} we derive \eqref{Lemma:26-10-22-1}. $\blacksquare$

\bigskip

\textbf{Remark.} Under the same assumptions of Lemma \ref{Lemma:26-10-22}, taking into account \eqref{normal:26-10-22}, we have

\begin{equation}\label{2-11-22-1}
	\int_{W} \partial_s fdx=\int_{\Gamma(g)}  f \nu^{(g)}_sdS,\quad s=1,\cdots, n.
\end{equation} 
$\blacklozenge$

\bigskip

\textbf{Proof Theorem \ref{TeoDiv:25-10-22-1}.} Let us begin by considering the case where $f\in C^{\infty}\left(\mathbb{R}^n\right)$. 

Since $\partial \Omega$ is of class $C^{0,1}$, we may assume that there exist positive numbers $r_0, M_0$, such that for every $P\in 
\partial \Omega$ there is an isometry $$ \Phi_P:\mathbb{R}^n\rightarrow \mathbb{R}^n,$$ satisfying $$\Phi_P(0)=P$$ and
\begin{equation*}
	\Phi_P^{-1}\left(\Omega\right)\cap Q_{r_0,2M_0}	=W_P,
\end{equation*}
where
$$W_P=\left\{x\in Q_{r_0,2M_0}:\mbox{ }
x_n>g_P(x') \right\},$$
 $$g_P\in C^{0,1}\left(\overline{B'_{r_0}}\right),$$
$$g_P(0)=0, $$
and
$$\left\Vert g_P\right\Vert_{C^{0,1}\left(\overline{B'_{r_0}}\right)}\leq
M_0r_0.$$ 	

Since $\partial\Omega$ is compact, there exist $P_1,\cdots, P_l\in \partial \Omega$ such that the family of sets
 $\left\{\Phi_{P_j}\left(Q_{r_0,2M_0}\right)\right\}_{1\leq j\leq l}$ is a finite covering of  $\partial\Omega$.  
 Set
 $$V_j=\Phi_{P_j}\left(Q_{r_0,2M_0}\right),\quad\mbox{ } j=1,\cdots, l.$$
 By Theorem \ref{Sob:lem3.3} there exist $\zeta_0, \zeta_1, \cdots, \zeta_l\in
 C^{\infty}_0(\mathbb{R}^n)$ satisfying 
 $$(\zeta_0)_{|\Omega}\in
 C^{\infty}_0(\Omega),$$
 $$\mbox{supp }\zeta_0\subset \mathbb{R}^n\setminus \partial\Omega;\quad \mbox{supp
 }\zeta_j\subset V_j,\quad j=1,\cdots, l,$$
 
 $$ \sum_{j=0}^l
 \zeta_j=1,\quad \mbox{su } \mathbb{R}^n;\quad 0\leq \zeta_j\leq
 1,\quad j=0,1\cdots, l.$$ 
 We have 
 
 $$\int_{\Omega}\nabla\left(f\zeta_0\right)dx=0.$$ Hence
 
 \begin{equation}\label{TeoDiv:26-10-22-1}
 \int_{\Omega}\nabla fdx=\sum_{j=0}^l\int_{\Omega}\nabla\left(f\zeta_j\right)dx=\sum_{j=1}^l\int_{\Omega}\nabla\left(f\zeta_j\right)dx.
 \end{equation}
 Now, let us fix $j\in \left\{1\,\cdots,l\right\}$ and let us denote by $W:=W_{P_j}$, $g:=g_{P_j}$ and
 \begin{equation}\label{TeoDiv:26-10-22-1-0}
 	F:=f\zeta_j.
 \end{equation}
 Let $\left\{g_m\right\}$ be a sequence of smooth approximating of $g$.
 Set 
 $$\Phi_j=\Phi_{P_j}, \ \ j=1, \cdots l$$
 and denoting  
 $$U=\Phi_j\left(W\right), $$ 
 $$W_m=\left\{x\in Q_{r_0,2M_0}:\mbox{ }
 x_n>g_m(x') \right\},$$
 $$U_m=\Phi_j\left(W_m\right),$$
 we have
 
 \begin{equation}\label{TeoDiv:26-10-22-2} 
 	\begin{aligned}
 	\int_{\Omega}\nabla\left(f\zeta_j\right)dx=\int_{U}\nabla Fdx=
\lim_{m\rightarrow \infty} \int_{U_m}\nabla Fdx
	\end{aligned}
\end{equation}
(as a matter of fact, we have $|U_m\setminus U|\rightarrow 0$ as $m\rightarrow \infty$). 

Now let us deal with the last integral in \eqref{TeoDiv:26-10-22-2}. Since $\Phi_j$ is an isometry, there exists a matrix $A=\left\{a_{qs}\right\}_{1\leq q,s\leq n}$ such that  

\begin{equation*} 
	\Phi(y)=Ay+P_j,\quad A^TA=I_n,
\end{equation*} 
where $I_n$ is the the identity matrix $n\times n$. For the sake of brevity let us set $P:=P_j$, $\Phi:=\Phi_j$. We have
 \begin{equation}\label{TeoDiv:26-10-22-3}
\int_{U_m}\nabla Fdx=\int_{W_m}\left(\nabla_xF\right)\left(\Phi(y)\right)dy.
\end{equation}
By setting
\begin{equation}\label{TeoDiv:26-10-22-3-0}
	\overline{F}=F\circ \Phi,
\end{equation} we have 
$$\left(\partial_y\Phi(y)\right)^T \left(\nabla_x F\right)\left(\Phi(y)\right)=\nabla_y \overline{F}(y)$$ and, recalling that $A=\left(\partial_y\Phi(y)\right)$, $A^TA=I_n$, we get

\begin{equation}\label{TeoDiv:26-10-22-4}
\left(\nabla_x F\right)\left(\Phi(y)\right)=A\nabla_y \overline{F}(y),\quad \forall y\in Q_{r_0,2M_0}.
\end{equation}
 Now, since $\overline{F}\in C^{\infty}_0\left(Q_{r_0,2M_0}\right)$ and, for $m$ large enough, the graph of $g_m$ is contained in $Q_{r_0,2M_0}$ (recall that $\left\{g_m\right\}$ is a smooth approximating sequence of $g$), we obtain, by Lemma \ref{Lemma:26-10-22} and by \eqref{TeoDiv:26-10-22-4}
 \begin{equation}\label{TeoDiv:26-10-22-5}
 	\begin{aligned}
 \left(\int_{W_m}\left(\nabla_xF\right)\left(\Phi(y)\right)dy\right)_{q}&=
 \left(\int_{W_m}A\nabla_y \overline{F}(y)dy
 \right)_q=\\&=
 \sum_{s=1}^n\int_{W_m}a_{qs}\partial_s \overline{F}(y)dy=\\&=
 -\int_{B'_r}a_{qn}\overline{F}\left(y',g_m(y')\right)dy'+\\&+
 \sum_{s=1}^{n-1}\int_{B'_r}a_{qs}\overline{F}\left(y',g_m(y')\right)\partial_sg_m(y')dy'.
 \end{aligned}
 \end{equation}
Since $\left\{g_m\right\}$ is a smooth approximating sequence of $g$, we have

$$\overline{F}\left(y',g_m(y')\right) \rightarrow \overline{F}\left(y',g(y')\right),\quad \mbox{as } m\rightarrow\infty \mbox{ in } L^{\infty}\left(B'_r\right)$$
and, for any $s=1,\cdots, n-1$,
$$\overline{F}\left(y',g_m(y')\right)\partial_sg_m(y') \rightarrow \overline{F}\left(y',g(y')\right)\partial_sg(y'), \mbox{ as } m\rightarrow\infty, \mbox{ in } L^{2}\left(B'_r\right).$$
Hence
\begin{equation}\label{TeoDiv:26-10-22-6}
	\lim_{m\rightarrow\infty} \int_{B'_r}\overline{F}\left(y',g_m(y')\right)dy'=\int_{B'_r}\overline{F}\left(y',g(y')\right)dy'
\end{equation}
 and, for any $s=1,\cdots,n-1$,
 \begin{equation}\label{TeoDiv:26-10-22-7}
 	\lim_{m\rightarrow\infty} \int_{B'_r}\overline{F}\left(y',g_m(y')\right)\partial_sg_m(y')dy'=\int_{B'_r}\overline{F}\left(y',g(y')\right)\partial_sg(y')dy'.
 \end{equation}
By  \eqref{TeoDiv:26-10-22-2}, \eqref{TeoDiv:26-10-22-5}, \eqref{TeoDiv:26-10-22-6} and \eqref{TeoDiv:26-10-22-7}, denoting by $\nu^{(g)}$ the unit outward normal to $\Gamma(g)$ (recall \eqref{TeoDiv:26-10-22-1-0} and \eqref{TeoDiv:26-10-22-3-0}) we have 
\begin{equation*}
	\begin{aligned}
\int_{\Omega}\nabla\left(f\zeta_j\right)dx&=\lim_{m\rightarrow\infty}\int_{U_m}\nabla Fdx=\\&=\lim_{m\rightarrow\infty}\int_{W_m}\nabla_xF\left(\Phi(y)\right)dy=\\&=	\int_{\Gamma(g)} \overline{F}A\nu^{(g)}dS=\\&=\int_{\partial\Omega}f\zeta_j\nu dS.
	\end{aligned}
\end{equation*}

\smallskip

\noindent By what has just been obtained and by \eqref{TeoDiv:26-10-22-1} we derive that if $f\in C^{\infty}\left(\mathbb{R}^n\right)$, then

\begin{equation*}
	\int_{\Omega}\nabla f dx=\int_{\partial\Omega} f \nu dS.
\end{equation*}

Finally, let us consider the case where $f\in C^{0,1}\left(\overline{\Omega}\right)$. If $\left\{f_m\right\}$ is a smooth approximating sequence of $f$ given by Lemma \ref{TeoDiv:25-10-22-3}, then we have

\begin{equation*}
	\begin{aligned}
		\int_{\Omega}\nabla f dx&=\lim_{m\rightarrow\infty}\int_{\Omega}\nabla f_mdx=\\&=\lim_{m\rightarrow\infty}\int_{\partial \Omega}f_m\nu dS=\\&=	\int_{\partial\Omega}f\nu dS.
	\end{aligned}
\end{equation*}
$\blacksquare$

\section{The Hausdorff distance \index{Hausdorf distance}} \label{Hausd}

Let $(X,d)$ be a metric space. Let us recall that the distance of a point $x\in X$ from a subset $A$ of $X$, $A\neq \emptyset$, is given by 
\begin{equation}\label{Hausd:10-11-22-1}
	d(x,A)=\inf\left\{d(x,y): \mbox{ } y\in A\right\}.
\end{equation}

\begin{prop}\label{Hausd:10-11-22-2}
	If $A\subset X$, $A\neq\emptyset$, then we have  
	\begin{equation}\label{Hausd:10-11-22-3}
		\left|d(x,A)-d(y,A)\right|\leq d(x,y), \quad \forall x,y\in A.
	\end{equation}
In particular, the map
$$X\ni x\rightarrow d(x,A)\in \mathbb{R},$$
is Lipschitz continuous.	
\end{prop}
\textbf{Proof.} Let $x,y\in X$. By triangle inequality we have 
$$d(x,A)\leq d(x,z)\leq d(x,y)+d(y,z),\quad \forall z\in A, $$ from which we derive
$$d(x,A)\leq d(x,y)+d(y,A). $$
Hence 
$$d(x,A)-d(y,A)\leq d(x,y)$$ and interchanging $x$ and $y$, we obtain \eqref{Hausd:10-11-22-3}. $\blacksquare$

\bigskip

We denote by $\mathbf{K}(X)$ \index{$\mathbf{K}(X)$}the family of nonempty compact sets of $X$. \\  If $K\in \mathbf{K}(X)$ and $x\in X$ we have
$$d(x,K)=\min\left\{d(x,y): \mbox{ } y\in K\right\}.$$

Let $K\in \mathbf{K}(X)$. Let us denote by  $S(K)$ the set of the points $x\in X$ such that  
\begin{equation*}
	\left\{y\in K:\mbox{ } d(x,y)=d(x,K) \right\}, \mbox{ has only one point}.
\end{equation*}
Trivially, we have $K\subset S(K)$ and it is well--defined map 
\begin{equation}\label{Hausd:10-11-18-1}
	p_K:S(K)\rightarrow K,\quad \mbox{ such that }  d(x,p_K(x))=d(x,K).
\end{equation}
 If $x\in S(K)$, we call $p_K(x)$ the \textbf{point of minimum distance} of $x$ from $K$ or also \textbf{the projection} \index{projection of a point on a compact set} of $x$ on $K$.

\begin{prop}\label{funzdist:19-11-17-8}
	Let $K\in \mathbf{K}(X)$. Then $p_K$, defined by \eqref{Hausd:10-11-18-1}, is a continuous map.
\end{prop}
	
	\textbf{Proof.} 
	Let us argue by contradiction. Let us assume that $p_K$ is not continuous. Consequently,  there exists a point $x_0\in S(K)$ and a sequence $\left\{x_n\right\}$ of $S(K)$ satisfying
	\begin{equation}\label{funzdist:19-11-17-11}
		\left\{x_n\right\}\rightarrow x_0
	\end{equation}
	and
\begin{equation}\label{funzdist:19-11-17-12}
	\left\{p_K\left(x_n\right)\right\}	\nrightarrow p_K\left(x_0\right).
	\end{equation}
	The latter implies that there exists $\varepsilon>0$ and a subsequence $\left\{x^*_n\right\}$ of $\left\{x_n\right\}$  satisfying 
	
	\begin{equation}\label{funzdist:19-11-17-13}
		d\left( p_K\left(x^*_n\right),p_K(x_0)\right)\geq \varepsilon, \quad \forall n\in \mathbb{N}.
	\end{equation}
	
	Since for every $n\in \mathbb{N}$ we have $p_K\left(x^*_n\right)\in K$,  and $K$ is compact, there exists a subsequence $\left\{x^{**}_n\right\}$ of $\left\{x^{*}_n\right\}$ such that
	$\left\{p_K\left(x^{**}_n\right)\right\}$ converges to a point $z\in K$. On the other hand, by \eqref{funzdist:19-11-17-11} we get
	
	\begin{equation*}
		\left\{x^{**}_n\right\}\rightarrow x_0.
	\end{equation*}
	Hence
	
	\begin{equation*}
		d(x_0,K)=\lim_{n\rightarrow \infty}d\left(x^{**}_n, K\right)=\lim_{n\rightarrow \infty}d\left(x^{**}_n,p_K\left(x^{**}_n\right)\right) =d(x_0,z). 
	\end{equation*}
	Consequently 
	$$d(x_0,K)=d(x_0,z)$$
	and by the definition of $S(K)$ we have $$z=p_K(x_0).$$ On the other hand,  \eqref{funzdist:19-11-17-13} gives

	\begin{equation*}
		d(x_0,z)=\lim_{n\rightarrow}	d\left(p_K(x_0),p_k\left(x^{**}_n\right)\right) \geq \varepsilon.
	\end{equation*}
	We have actually reached a contradiction. Therefore the map $x \rightarrow p_K(x)$ is continuous. $\blacksquare$
	
\bigskip
	
\begin{definition}\label{Hausd:11-11-22-def1}
For any $K_1, K_2\in \mathbf{K}(X)$ we denote
\begin{equation}\label{Hausd:10-11-22-4}
	\delta\left(K_1,K_2\right)=\max\left\{d(x,K_2):\mbox{ } x\in K_1\right\}.
\end{equation}
\end{definition}

\bigskip

\begin{prop}\label{Hausd:10-11-22-5}
	If $K_1,K_2\in \mathbf{K}(X)$, then we have 
	\begin{equation}\label{Hausd:10-11-22-6}
	K_1\subset K_2 \ \ \Leftrightarrow \ \ \delta\left(K_1,K_2\right)=0.
	\end{equation}
\end{prop}
\textbf{Proof.} If $K_1\subset K_2$, then
$$d(x,K_2)=0,\quad \forall x\in K_1,$$
hence 
$$\delta\left(K_1,K_2\right)=0.$$ Conversely, if $\delta\left(K_1,K_2\right)=0$ then 
$$d(x,K_2)=0,\quad \forall x\in K_1.$$
Since $K_1$ is a closed set of $X$, we obtain 
$$x\in K_2,\quad\forall x\in K_1$$ that is $K_1\subset K_2$. $\blacksquare$

\bigskip 

Let us notice that $\delta(\cdot,\cdot)$ \textbf{does not} define a distance on $\mathbf{K}(X)$. Actually, by \eqref{Hausd:10-11-22-6} we have that $\delta(\cdot,\cdot)$ is not symmetric and
  
$$\delta\left(K_1,K_2\right)=0 \ \nRightarrow \ K_1= K_2.$$
 
\bigskip

For any $K_1,K_2\in \mathbf{K}(X)$, let us denote by \index{$d_{\mathcal{H}}$}

\begin{equation}\label{Hausd:10-11-22-7}
	 d_{\mathcal{H}}\left(K_1,K_2\right)=\max\left\{\delta\left(K_1,K_2\right), \delta\left(K_2,K_1\right)\right\}.
\end{equation}

\bigskip

\begin{prop}\label{Hausd:10-11-22-8}
	$d_{\mathcal{H}}\left(\cdot,\cdot\right)$, defined by \eqref{Hausd:10-11-22-7}, is a distance on $\mathbf{K}(X)$.
\end{prop}
\textbf{Proof.} It is obvious that $d_{\mathcal{H}}\left(K_1,K_2\right)=d_{\mathcal{H}}\left(K_2,K_1\right)$ and that \\ $d_{\mathcal{H}}\left(K_1,K_2\right)\geq 0$ for every $K_1,K_2\in \mathbf{K}(X)$. Furthermore, if 

$$d_{\mathcal{H}}\left(K_1,K_2\right)=0,$$
then $\delta \left(K_1,K_2\right)=0$ and $\delta \left(K_2,K_1\right)=0$ which imply, respectively, $K_1\subset K_2$ e $K_2\subset K_1$, hence $K_1=K_2$.

It only remains to prove the triangular inequality. We begin by proving that if $K_1,K_2\in \mathbf{K}(X)$ then for any $L\in \mathbf{K}(X)$ we have

\begin{equation}\label{Hausd:11-11-22-1}
	\delta \left(K_1,K_2\right)\leq \delta \left(K_1,L\right)+\delta \left(L,K_2\right).
\end{equation}   
Let $x\in K_1$. For any $y\in K_2$ and for any $z\in L$ we have
\begin{equation*}
d\left(x,K_2\right)\leq d\left(x,y\right)\leq d\left(x,z\right)+d\left(z,y\right),
\end{equation*}
from which we have

\begin{equation*}
	\begin{aligned}
		d\left(x,K_2\right)&\leq d\left(x,z\right)+d\left(z,K_2\right)\leq \\& \leq 
		d\left(x,z\right)+\delta \left(L,K_2\right).
	\end{aligned}
\end{equation*}

\smallskip

\noindent Therefore

\begin{equation*}
	\begin{aligned}
		d\left(x,K_2\right)&\leq  d\left(x,L\right)+\delta \left(L,K_2\right)\leq \\&\leq \delta \left(K_1,L\right)+\delta \left(L,K_2\right).
	\end{aligned}
\end{equation*}
By the latter we obtain \eqref{Hausd:11-11-22-1} and similarly

\begin{equation}\label{Hausd:11-11-22-2}
	\delta \left(K_2,K_1\right)\leq \delta \left(K_2,L\right)+\delta \left(L,K_1\right).
\end{equation}
Now, let us assume, for instance, that $d_{\mathcal{H}}\left(K_1,K_2\right)=\delta\left(K_1,K_2\right)$. We have, by \eqref{Hausd:11-11-22-1} and \eqref{Hausd:11-11-22-2},

\begin{equation*}
	\begin{aligned}
		d_{\mathcal{H}}\left(K_1,K_2\right)&= \delta\left(K_1,K_2\right)\leq \\&\leq \delta\left(K_1,L\right)+\delta \left(L,K_2\right)\leq \\&\leq
		d_{\mathcal{H}}\left(K_1,L\right)+d_{\mathcal{H}}\left(L,K_2\right).
		\end{aligned}
\end{equation*}
$\blacksquare$
\begin{definition}\label{Hausd:11-11-22-def2}
	Let $K_1, K_2\in \mathbf{K}(X)$, we call $d_{\mathcal{H}}\left(K_1,K_2\right)$, defined by \eqref{Hausd:10-11-22-7}, the \textbf{Hausdorff distance} between $K_1$ and $K_2$. 
\end{definition}

In Proposition \ref{Hausd:11-11-22-3} below we give an useful characterization of the Hausdorff distance. For this purpose we introduce the following notation. Let $r\geq 0$ and $K\in \mathbf{K}(X)$, set
$$[K]_r=\left\{u\in K: \mbox{ } d(x,K)\leq r\right\},$$
\index{$[K]_r$}
Let us call $[K]_r$ the $r$--\textbf{dilation} of $r$ of the set $K$ \index{dilation of a compact set}. 

\bigskip

\begin{prop}\label{Hausd:12-11-22-1}
If $K\in\mathbf{K}(X)$ and $r\geq 0$ then $[K]_r$ is a closed set of $X$.
\end{prop}
\textbf{Proof.} Let $\left\{x_n\right\}$ be any sequence in $[K]_r$ which satisfies 
$$\lim_{n\rightarrow\infty}x_n=x_0.$$
Since

$$d\left(x_n,K\right)\leq r, \quad \forall n\in \mathbb{N},$$ passing to the limit as $n\rightarrow \infty$ and taking into account that $x\rightarrow d\left(x,K\right)$ is continuous, we have

$$d\left(x_0,K\right)=\lim_{n\rightarrow\infty}d\left(x_n,K\right)\leq r.$$
Hence $x_0\in [K]_r$. Therefore $[K]_r$ is closed set of $X$. $\blacksquare$

\bigskip

\textbf{Remark.} It generally does not occur that $[K]_r$ is compact. For instance, let $X$ be a infinite dimensional Hilbert space, and let $K=\{a\}$, where $a\in X$, then we have $[K]_r=\overline{B_r(a)}$ and from Functional Analysis we know that $\overline{B_r(a)}$ is not compact. $\blacklozenge$

\begin{prop}\label{Hausd:11-11-22-3}
	Let $K_1, K_2\in \mathbf{K}(X)$ and $r\geq 0$. Then we have 
	\begin{equation}\label{Hausd:11-11-22-4}
		d_{\mathcal{H}}\left(K_1,K_2\right)\leq r \ \Longleftrightarrow \ \begin{cases}
K_1\subset \left[K_2\right]_r,\\
			\\
			K_2\subset \left[K_1\right]_r. 
		\end{cases}
	\end{equation}

\smallskip

	\begin{equation}\label{Hausd:11-11-22-5}
		d_{\mathcal{H}}\left(K_1,K_2\right)=\min \left\{r\geq 0: \mbox{ } K_1\subset \left[K_2\right]_r, \mbox{ } K_2\subset \left[K_1\right]_r\right\}.
	\end{equation}
\end{prop}

\textbf{Proof.}  The condition $$d_{\mathcal{H}}\left(K_1,K_2\right)\leq r$$ is equivalent to $\delta\left(K_1,K_2\right)\leq r$ and $\delta\left(K_2,K_1\right)\leq r$. On the other hand   
\begin{equation*}
	\delta\left(K_1,K_2\right)\leq r\Longleftrightarrow \left(d\left(x,K_2\right)\leq r, \ \forall x\in K_1\right) \Longleftrightarrow K_1\subset \left[K_2\right]_r.
\end{equation*}
Similarly, 

\begin{equation*}
	\delta\left(K_2,K_1\right)\leq r\Longleftrightarrow  K_2\subset \left[K_1\right]_r.
\end{equation*}
From what was obtained \eqref{Hausd:11-11-22-4} follows. 

Now, we prove \eqref{Hausd:11-11-22-5}. Let us denote
$$d=d_{\mathcal{H}}\left(K_1,K_2\right)$$ and
$$\rho=\inf \left\{r\geq 0: \mbox{ } K_1\subset \left[K_2\right]_r, \mbox{ } K_2\subset \left[K_1\right]_r\right\}.$$
\eqref{Hausd:11-11-22-4} implies ("$ \Longrightarrow$")
\begin{equation}\label{Hausd:11-11-22-6}
	\rho\leq d.
\end{equation}
On the other hand, for any $\varepsilon>0$, we have 
$$K_1\subset \left[K_2\right]_{\rho+\varepsilon},\quad \mbox{and}\quad K_2\subset \left[K_1\right]_{\rho+\varepsilon}$$ and  \eqref{Hausd:11-11-22-4} implies ("$\Longleftarrow$") 
$$d=d_{\mathcal{H}}\left(K_1,K_2\right)\leq \rho+\varepsilon.$$ Therefore, since  $\varepsilon$ is arbitrary, we obtain 
$$d\leq \rho.$$
By the latter and by \eqref{Hausd:11-11-22-6} we get
$$ d= \rho,$$
from which \eqref{Hausd:11-11-22-5} follows. $\blacksquare$

\bigskip

\textbf{Example.} Let us consider $$K_1=\overline{B_1}\setminus B_{\varepsilon},\quad K_2=\overline{B_1},$$ in $\mathbb{R}^n$, where $\varepsilon\in (0,1)$. 

We have 
$$K_1\subset K_2\subset \left[K_2\right]_{\varepsilon} $$ and

$$K_2\subset \left[K_1\right]_{r},\quad \forall r\geq \varepsilon.$$
Hence 
$$d_{\mathcal{H}}\left(K_1,K_2\right)=\varepsilon.$$
Let us see what happens with regard to 
$$d_{\mathcal{H}}\left(\partial K_1,\partial K_2\right).$$
We have
$$\partial K_2\subset \partial K_1$$ and
$$\partial K_1\subset \left[\partial K_2\right]_{r},\quad \forall r\geq 1-\varepsilon.$$
Hence 
$$d_{\mathcal{H}}\left(\partial K_1,\partial K_2\right)=1-\varepsilon.$$ 
Therefore, neither of the two relationships holds true
\begin{equation}\label{Hausd:11-11-22-7}
d_{\mathcal{H}}\left(K_1,K_2\right)\leq d_{\mathcal{H}}\left(\partial K_1,\partial K_2\right), 
\end{equation}
\begin{equation}\label{Hausd:11-11-22-8}
	d_{\mathcal{H}}\left(\partial K_1,\partial K_2\right)\leq d_{\mathcal{H}}\left(K_1,K_2\right). 
\end{equation}
As a matter of fact \eqref{Hausd:11-11-22-7} is false for $\frac{1}{2}\leq \varepsilon<1$, and \eqref{Hausd:11-11-22-8} is false for $ 0<\varepsilon<\frac{1}{2}$. $\spadesuit$

\bigskip

\begin{prop}\label{Hausd:11-11-22-9}
	$$\left[K_1\right]_{r}\cup \left[K_2\right]_{r}=\left[K_1\cup K_2\right]_{r}, \quad \forall K_1, K_2\in \mathbf{K}(X), \ \forall r\geq 0.$$
\end{prop}
\textbf{Proof.} Let $K_1, K_2\in \mathbf{K}(X)$. Let us begin to prove 
\begin{equation}\label{Hausd:12-11-22-2}
	\left[K_1\right]_{r}\cup \left[K_2\right]_{r}\subset\left[K_1\cup K_2\right]_{r}. 
\end{equation}
Let $x\in \left[K_1\right]_{r}\cup \left[K_2\right]_{r}$ and, for instance, let $x\in \left[K_1\right]_{r}$, then 
$$d\left(x, K_1 \cup K_2\right)\leq d\left(x, K_1 \right)\leq r,$$ which implies $x\in \left[K_1\cup K_2\right]_{r}$.  \eqref{Hausd:12-11-22-2} is proved. 

Now, let us prove
\begin{equation}\label{Hausd:12-11-22-3}
	\left[K_1\cup K_2\right]_{r}\subset \left[K_1\right]_{r}\cup \left[K_2\right]_{r}. 
\end{equation}
Let $x\in \left[K_1\cup K_2\right]_{r}$. We have 
$$d\left(x, K_1 \cup K_2\right)\leq r.$$ Let $y\in K_1 \cup K_2$ satisfy $d(x,y)=d\left(x, K_1 \cup K_2\right)$. Now, if $y\in K_1$, we have
$$d\left(x, K_1 \right)\leq d(x,y)\leq r,$$ 
hence $x\in \left[K_1\right]_{r}$. Similarly, if $y\in K_2$ then $x\in \left[K_2\right]_{r}$. In any case \\ $x\in \left[K_1\right]_{r}\cup \left[K_2\right]_{r}$.  Hence \eqref{Hausd:12-11-22-3} is proved. $\blacksquare$

\bigskip 

The following Theorem has been proved in  \textbf{Kuratowski}, \cite[ \S 15, VIII]{Kur}.

\begin{theo}\label{Hausd:Kurato}
	Let $x_0\in X$. Let us define, for any 
	$K\in \mathbf{K}(X)$, the function 
	\begin{equation}\label{Hausd:12-11-22-4}
		f_K(x)= d(x,K)-d(x,x_0). 
	\end{equation}
We have:

\smallskip

(i) 	$f_K$ is bounded,

\smallskip

\smallskip

(ii)
$$d_{\mathcal{H}}\left(K_1,K_2\right)=\sup_{x\in X}\left|f_{K_2}(x)-f_{K_1}(x)\right|,\quad \forall K_1, K_2\in \mathbf{K}(X).$$
\end{theo}
\textbf{Proof.} 

\noindent Let us prove (i). Let $x\in X$.  We have, by the triangle inequality,
$$d(x,K)\leq d(x,y)\leq d(y,x_0)+d(x_0,x), \quad \forall y\in K.$$
Hence
\begin{equation}\label{Hausd:12-11-22-5}
	d(x,K)\leq d(x_0,K)+d(x,x_0).
\end{equation}
By using again the triangle inequality we have, for any $y,z\in K$

\begin{equation*}
		d(x,x_0)\leq d(x,y)+d(y,z)+d(z,x_0)\leq d(x,y)+d(z,x_0)+d(K),
\end{equation*} 
where $d(K)$ is the diameter of $K$, by the last inequality we have
\begin{equation}\label{Hausd:12-11-22-6}
	d(x,x_0)\leq d(x,K)+d(x_0,K)+d(K).
\end{equation} 
Now, \eqref{Hausd:12-11-22-5} gives
$$f_K(x)= d(x,K)-d(x,x_0)\leq d(x_0,K)$$
and \eqref{Hausd:12-11-22-6} gives
 $$f_K(x)\geq d(x,K)-\left(d(x,K)+d(x_0,K)+d(K)\right)=-d(x_0,K)-d(K).$$ Therefore we have
 
 $$\left|f_K(x)\right|\leq d(x_0,K)+d(K), \quad\forall x\in X.$$
 
 Let us now prove (ii). Let $K_1, K_2\in \mathbf{K}(X)$. It is not restrictive to assume
 
 $$d_{\mathcal{H}}\left(K_1,K_2\right)=\delta \left(K_1,K_2\right)=\max_{x\in K_1}d\left(x,K_2\right).$$
 Let $\overline{x}\in K_1$ satisfy  
$$d\left(\overline{x},K_2\right)=d_{\mathcal{H}}\left(K_1,K_2\right).$$ Since we have trivially  
$d\left(\overline{x},K_1\right)=0$, we get 
\begin{equation}\label{Hausd:12-11-22-7}
	\begin{aligned}
	d_{\mathcal{H}}\left(K_1,K_2\right)&=d\left(\overline{x},K_2\right)- d\left(\overline{x},K_1\right)=\\&\\&= d\left(\overline{x},K_2\right)-d\left(\overline{x},x_0\right)-\left(d\left(\overline{x},K_1\right)-d\left(\overline{x},x_0\right)\right)=\\&\\&= f_{K_2}(\overline{x})-f_{K_1}(\overline{x})\leq \sup_{x\in X}\left|f_{K_1}(x)-f_{K_1}(x)\right|.
	\end{aligned}
\end{equation} 
Now, for any $x\in X$ let $y\in K_1$ satisfy 
$$ d(x,y)=d\left(x,K_1\right).$$ We have

$$d\left(x,K_2\right)\leq d(x,y)+d\left(y,K_2\right)=d\left(x,K_1\right)+d\left(y,K_2\right).$$
Hence 

\begin{equation*}
	d\left(x,K_2\right)-d\left(x,K_1\right)\leq d\left(y,K_2\right)\leq d_{\mathcal{H}}\left(K_1,K_2\right)
\end{equation*} 
and, by interchanging $K_1$ with $K_2$, we have. 
$$d\left(x,K_1\right)-d\left(x,K_2\right)\leq d_{\mathcal{H}}\left(K_1,K_2\right).$$
Hence
\begin{equation*}
	\begin{aligned}
		\left|f_{K_1}(x)-f_{K_2}(x)\right|&=\left|d\left(x,K_1\right)-d\left(x,K_2\right)\right|\leq d_{\mathcal{H}}\left(K_1,K_2\right), \quad \forall x\in X,
	\end{aligned}
\end{equation*} 
which implies

\begin{equation}\label{Hausd:12-11-22-8}
\sup_{x\in X}\left|f_{K_1}(x)-f_{K_1}(x)\right|\leq d_{\mathcal{H}}\left(K_1,K_2\right). 
\end{equation} 
Finally, \eqref{Hausd:12-11-22-7} and \eqref{Hausd:12-11-22-8} imply (ii). $\blacksquare$

\subsection{Completeness and compactness of $\left(\mathbf{K}(X), d_{\mathcal{H}}\right)$}\label{12-11-completezza}

The Main Theorem that we prove in the present Section is the following one.

\begin{theo}[\textbf{completeness}]\label{Hausd:12-11-22-9}
If $(X,d)$ is a complete metric space, then $\left(\mathbf{K}(X), d_{\mathcal{H}}\right)$ is complete.
\end{theo} 

\medskip

In order to prove Theorem \ref{Hausd:12-11-22-9} we need the following Lemma. 

\begin{lem}\label{Hausd:12-11-22-10}
	Let $(X,d)$ be a metric space and let $\left\{K_n\right\}$ be a Cauchy sequence in $\left(\mathbf{K}(X), d_{\mathcal{H}}\right)$. Let $\left\{n_j\right\}$ be a strictly increasing sequence in  $\mathbb{N}$ and let $\left\{x_{n_j}\right\}$ be a Cauchy sequence in $(X,d)$ satisfying
	\begin{equation*}
	x_{n_j}\in K_{n_j},\quad \forall j\in \mathbb{N}.
	\end{equation*}
Then there exists a Cauchy sequence  $\left\{\overline{x}_n\right\}$ in $(X,d)$ such that

\begin{equation}\label{Hausd:13-11-22-1}
\overline{x}_{n_j}=x_{n_j},\quad \forall j\in \mathbb{N},\quad \mbox{and}\quad  \overline{x}_n\in K_n, \quad \forall n\in \mathbb{N}.
\end{equation}
\end{lem} 
\textbf{Proof.} Let us define $\left\{\overline{x}_n\right\}$ as follows: if $1\leq n\leq n_1-1$, then we choose $\overline{x}_n$ satisfying 
$$d\left(x_{n_1},\overline{x}_n\right)=d\left(x_{n_1},K_n\right),$$
if $n_j+1\leq n\leq n_{j+1}-1$, $j\in \mathbb{N}$, then we choose $\overline{x}_n$ satisfying
$$d\left(x_{n_j},\overline{x}_n\right)=d\left(x_{n_j},K_n\right),$$ finally, if $n=n_j$, $j\in \mathbb{N}$, then we choose 
$$\overline{x}_n=x_{n_j}.$$
Notice that \eqref{Hausd:13-11-22-1} is satisfied by construction, so we are left to prove that $\left\{\overline{x}_n\right\}$ is a Cauchy sequence. 

Let us fix any $\varepsilon>0$ and let $\nu\in \mathbb{N}$ satisfy 
 
 \begin{equation}\label{Hausd:13-11-22-2}
 	d\left(x_{n_j},x_{n_h}\right)<\frac{\varepsilon}{3},\quad \forall j,h\geq \nu
 \end{equation}
and
 
\begin{equation}\label{Hausd:13-11-22-3}
	d_{\mathcal{H}}\left(K_n,K_m\right)<\frac{\varepsilon}{3},\quad \forall n,m\geq n_{\nu}.
\end{equation}
Let $n,m\geq n_{\nu}$ and let $j$ and $h$ be such that 
$$n_j\leq n\leq n_{j+1},\quad\quad n_h\leq m\leq n_{h+1}.$$
By the triangle inequality we have
\begin{equation}\label{Hausd:13-11-22-4}
d\left(\overline{x}_n,\overline{x}_m\right)\leq d\left(\overline{x}_n,x_{n_j}\right)+d\left(x_{n_j},x_{n_h}\right)+d\left(\overline{x}_m,x_{n_h}\right).
\end{equation}
Now, \eqref{Hausd:13-11-22-3} implies
$$d\left(x_{n_j},\overline{x}_n\right)=d\left(x_{n_j},K_n\right)\leq d_{\mathcal{H}}\left(K_{n_j},K_n\right)<\frac{\varepsilon}{3}.$$
Similarly, we have
$$d\left(x_{n_h},\overline{x}_m\right)<\frac{\varepsilon}{3}.$$ By these latter inequalities and by \eqref{Hausd:13-11-22-2}, \eqref{Hausd:13-11-22-4} we get 
$$d\left(\overline{x}_n,\overline{x}_m\right)<\varepsilon$$ and thereby we have also proved that $\left\{\overline{x}_n\right\}$ is a Cauchy sequence. $\blacksquare$

\bigskip

In what follows, for any sequence  $\left\{K_n\right\}$ in $\mathbf{K}(X)$ and any sequence  $\left\{x_n\right\}$ in $X$, we will write simply $\left\{x_n\in K_n\right\}$ to denote that
$$x_n\in K_n,\quad\forall n\in \mathbb{N}.$$

\bigskip

\textbf{Proof of Theorem \ref{Hausd:12-11-22-9}.} Let $\left\{K_n\right\}$ be a Cauchy sequence in $\mathbf{K}(X)$. Let us denote 

\begin{equation*}
	K=\left\{x\in X:\mbox{ there exists a sequence }  \left\{x_n\in K_n\right\} \mbox { such that } \lim_{n\rightarrow \infty}x_n=x\right\}.
\end{equation*}
Let us prove the following:

\smallskip

(\textbf{a}) $K\neq \emptyset$;

\smallskip

(\textbf{b}) $K$ is closed;

\smallskip

(\textbf{c}) $$\forall \varepsilon>0 \ \exists n_{\varepsilon}\in \mathbb{N} \mbox{ such that } \forall n\geq n_{\varepsilon} \ \ K\subset \left[K_n\right]_{\varepsilon};$$

\smallskip

(\textbf{d}) $K$ is compact;

\smallskip

\smallskip

\textbf{and}
$$\left\{K_n\right\}\rightarrow K,\quad \mbox{ in } \left(\mathbf{K}(X),d_{\mathcal{H}}\right).$$

\medskip

\textbf{Proof of (a).} Since $\left\{K_n\right\}$ is a  Cauchy sequence in 
 $\left(\mathbf{K}(X),d_{\mathcal{H}}\right)$, we have that, for any $\varepsilon>0$,  there exists $n_{\varepsilon}\in \mathbb{N}$, such that 
 
 \begin{equation}\label{Hausd:13-11-22-5}
 d_{\mathcal{H}}\left(K_n,K_m\right)<\varepsilon,\quad \forall n,m\geq n_{\varepsilon}.
 \end{equation}
We may assume $n_{\varepsilon}\in \mathbb{N}$ be strictly increasing w.r.t. $\varepsilon$.
For any $j\in \mathbb{N}$ let $$\varepsilon_j=\frac{1}{2^j}$$ and set $n_j=n_{\varepsilon_j}$. 
 
 Let $x_{n_1}\in K_{n_1}$ be chosen arbitrarily.  Since 
 
 $$d\left(x_{n_1},K_{n_2}\right)\leq d_{\mathcal{H}}\left(K_{n_1},K_{n_2}\right),$$ we may choose  $x_{n_2}\in K_{n_2}$ such that
 $$d\left(x_{n_1},x_{n_2}\right)=d\left(x_{n_1},K_{n_2}\right).$$ Similarly, after choosing $x_{n_1},\cdots, x_{n_{j-1}}$, we choose \\ $x_{n_j}\in K_{n_j}$. More precisely, let us suppose to have already chosen $x_{n_1},\cdots, x_{n_{j-1}}$, then we choose $x_{n_j}\in K_{n_j}$ so that 
 $$d\left(x_{n_{j-1}},x_{n_j}\right)=d\left(x_{n_{j-1}},K_{n_j}\right)\leq d_{\mathcal{H}}\left(K_{n_{j-1}},K_{n_j}\right)<\frac{1}{2^{j-1}}.$$
 
  Now, let us prove that $\left\{x_{n_j}\right\}$ is a  Cauchy sequence.  For any $h>j$ the triangle inequality gives
 
 $$d\left(x_{n_{j}},x_{n_h}\right)\leq \sum_{l=j}^{h-1}d\left(x_{n_{l}},x_{n_{l+1}}\right)<\sum_{l=j}^{h-1}\frac{1}{2^{l}}<\frac{1}{2^{j-1}}.$$
 On the other hand 
 $$x_{n_{j}}\in K_{n_{j}},\quad\forall j\in \mathbb{N}.$$ Hence, Lemma \ref{Hausd:12-11-22-10} implies that there exists a Cauchy sequence in $(X,d)$, $\left\{\overline{x}_n\in K_n\right\}$, which satisfies 
 $$\overline{x}_{n_j}=x_{n_{j}},\quad\forall j\in \mathbb{N}.$$
 Since $(X,d)$ is a complete metric space, there exists $x\in X$ such that
 
 $$\lim_{n\rightarrow \infty} \overline{x}_n=x.$$
Hence, $x$ belongs to $K$ (as we have defined $K$). Therefore $K\neq \emptyset$.

\medskip

\textbf{Proof of (b).} We prove that if $\left\{x_n\right\}$ is a sequence in $K$, which converges to $x_0$, then $x_0\in K$. By the definition of $K$, we have that for every $n\in \mathbb{N}$ there exists a sequence $\left\{y^{(n)}_j\right\}_{j\in \mathbb{N}}$ which satisfies

\begin{equation}\label{Hausd:14-11-22-0}
y^{(n)}_j\in K_j, \quad \forall  j\in \mathbb{N}
\end{equation}
and
$$y^{(n)}_j\rightarrow x_n ,\mbox{ as } j\rightarrow \infty, \quad \forall n\in \mathbb{N}.$$
 Since $\left\{x_n\right\}$  converges to $x_0$, there exists a stricly increasing sequence in $\mathbb{N}$, $\left\{n_h\right\}_{h\in \mathbb{N}}$, which satisfies
\begin{equation}\label{Hausd:14-11-22-1}
d\left(x_{n_{h}},x_0\right)< \frac{1}{h},\quad \forall h\in \mathbb{N}. 
\end{equation} 
In addition, since 
$$y^{(n_h)}_j\rightarrow x_{n_h} ,\mbox{ as } j\rightarrow \infty, \quad \forall h\in \mathbb{N},$$
there exists a stricly increasing sequence in $\mathbb{N}$, $\left\{m_h\right\}_{h\in \mathbb{N}}$, which satisfies

\begin{equation}\label{Hausd:14-11-22-2}
	d\left(y^{(n_{h})}_{m_h},x_{n_h}\right)< \frac{2}{h},\quad \forall h\in \mathbb{N}.
\end{equation} 
By \eqref{Hausd:14-11-22-1} and \eqref{Hausd:14-11-22-2} we have

\begin{equation}\label{Hausd:14-11-22-3}
	d\left(y^{(n_{h})}_{m_h},x_{0}\right)< \frac{1}{h}\quad \forall h\in \mathbb{N}.
\end{equation} 
Now, let us consider the sequence $\left\{y^{(n_{h})}_{m_h}\right\}_{h\in \mathbb{N}}$. Since it is convergent, it is a Cauchy sequence and, by \eqref{Hausd:14-11-22-0}, we have
$$y^{(n_{h})}_{m_h}\in K_{m_h}, \quad \forall h\in \mathbb{N}.$$ Therefore by Lemma \ref{Hausd:12-11-22-10}, there exists a Cauchy sequence, $\left\{\overline{y}_{n}\right\}$ which satisfies 
$$\overline{y}_{m_h}=y^{(n_{h})}_{m_h}, \quad \forall h\in \mathbb{N}$$ and
$$\overline{y}_{n}\in K_n, \quad\forall n\in \mathbb{N}.$$ 
Since $X$ is a complete space and $\left\{\overline{y}_{n}\right\}$ is a Cauchy sequence, it converges. On the other hand, \eqref{Hausd:14-11-22-3} implies that the subsequence $\left\{y^{(n_{h})}_{m_h}\right\}_{h\in \mathbb{N}}$ converges to $x_0$. Therefore the whole sequence $\left\{\overline{y}_{n}\right\}$ converges to $x_0$ and by the definition of $K$ we have $x_0\in K$.

\medskip

\textbf{Proof of (c).} Let $\varepsilon >0$ and let $n_{\varepsilon}$ satisfy

\begin{equation*}
	d_{\mathcal{H}}\left(K_n,K_m\right)<\varepsilon,\quad \forall n,m\geq n_{\varepsilon}.
\end{equation*}
Proposition \ref{Hausd:11-11-22-3} gives

\begin{equation}\label{Hausd:14-11-22-4}
	K_m\subset \left[K_n\right]_{\varepsilon},\quad \forall n,m\geq n_{\varepsilon}.
\end{equation}
Let $x\in K$. Let us prove that $x\in \left[K_n\right]_{\varepsilon}$ for every $n\geq n_{\varepsilon}$. Fix $n\geq n_{\varepsilon}$. By the definition of $K$, there exists $\left\{x_m\in K_m\right\}$ such that 
$$x_m\rightarrow x, \quad\mbox{ as } \  m \rightarrow \infty.$$  Now, since  $x_m\in K_m\subset \left[K_n\right]_{\varepsilon}$ (by \eqref{Hausd:14-11-22-4}), for every $m\geq n_{\varepsilon}$ and taking into account that $\left[K_n\right]_{\varepsilon}$ is a closed set (Proposition
\ref{Hausd:12-11-22-1}), we have $x\in \left[K_n\right]_{\varepsilon}$.

\medskip

\textbf{Proof of (d).} Since $X$ is a complete metric space and $K$ is a closed set (by (b)), by Theorem \ref{Sob:teopag.92}, it suffices to prove that $K$ is totally bounded. Let us argue by contradiction. Let us assume that $K$ is not totally bounded. Hence, let us assume that there exists $\delta>0$ and there exists a sequence $\left\{x_n\right\}$ in $K$ so that 
$$d\left(x_n,x_m\right)\geq \delta.$$ Now, by (c), there exists $\nu\in \mathbb{N}$ such that, if $n\neq m$,
$$K\subset \left[K_{\nu}\right]_{\frac{\delta}{4}}.$$ From which we have that for every $n\in \mathbb{N}$ there exists $y_n\in K_{\nu}$ such that
$$d\left(y_n,x_n\right)<\frac{\delta}{4}.$$ On the other hand, since $K_{\nu}$ is compact, there exists a subsequence of $\left\{y_{n}\right\}$, $\left\{y_{n_j}\right\}$, which converges, consequently there exists $\nu'\geq \nu$ such that  
$$d\left(y_{n_j},y_{n_h}\right)<\frac{\delta}{4}, \quad \forall j,h \geq \nu'.$$ Therefore, by the triangle inequality we have, if $j\neq h$,
$$\delta \leq d\left(x_{n_j},x_{n_h}\right)\leq d\left(x_{n_j},y_{n_j}\right)+d\left(y_{n_j},y_{n_h}\right)+d\left(y_{n_h},x_{n_h}\right)< \frac{3\delta}{4},$$
for  $j,h \geq \nu'$, $j\neq h$.  This is clearly a contradiction.

\medskip

\textbf{Proof of (e).} Since we have proved (c), it suffices to prove 

\begin{equation}\label{Hausd:15-11-22-1}
	\forall \varepsilon>0 \ \exists \nu_{\varepsilon}\in \mathbb{N} \mbox{ such that }\forall n\geq \nu_{\varepsilon} \ \  K_n\subset \left[K\right]_{\varepsilon}.
\end{equation}
Since $\left\{K_n\right\}$ is a Cauchy sequence we have that, for any $\varepsilon>0$ there exists $\nu_{\varepsilon} \in \mathbb{N}$  such that 

\begin{equation}\label{Hausd:15-11-22-2}
	d_{\mathcal{H}}\left(K_n,K_m\right)<\frac{\varepsilon}{2},\quad \forall n,m\geq \nu_{\varepsilon}.
\end{equation}
Hence 

\begin{equation*}
	K_m\subset \left[K_n\right]_{\frac{\varepsilon}{2}}, \quad \forall n,m\geq \nu_{\varepsilon}.
\end{equation*}
Let us fix $\overline{n}\geq \nu_{\varepsilon}$. Inequality  \eqref{Hausd:15-11-22-2} implies that there exists a strictly increasing sequence $\left\{n_j\right\}$ in $\mathbb{N}$ such that $n_j\geq \nu_{\varepsilon}$, for every $j\in \mathbb{N}$, and 
\begin{equation*}
	d_{\mathcal{H}}\left(K_{n_{j-1}},K_{n_j}\right)<\frac{\varepsilon}{2^j}.
\end{equation*}
Since $n_1, \overline{n}\geq\nu_{\varepsilon}$, we get by \eqref{Hausd:15-11-22-2} 

\begin{equation*}
	d_{\mathcal{H}}\left(K_{\overline{n}},K_{n_1}\right)<\frac{\varepsilon}{2}.
\end{equation*}
Hence
\begin{equation}\label{Hausd:15-11-22-3}
	K_{\overline{n}}\subset \left[K_{n_1}\right]_{\frac{\varepsilon}{2}}.
\end{equation}
Now, let us fix $y\in K_{\overline{n}}$ and let us prove that $y\in \left[K\right]_{\varepsilon}$. By \eqref{Hausd:15-11-22-3} we have 
$$y\in \left[K_{n_1}\right]_{\frac{\varepsilon}{2}},$$ hence there exists $x_{n_1}\in K_{n_1}$ such that

\begin{equation}\label{Hausd:15-11-22-4}
d\left(x_{n_1},y\right)<\frac{\varepsilon}{2}.
\end{equation}
Generally speaking, since

\begin{equation*}
	K_{n_{j-1}}\subset \left[K_{n_j}\right]_{\frac{\varepsilon}{2^j}},\quad \forall j\geq 2,
\end{equation*}
there exists a sequence $\left\{x_{n_j}\right\}$ which satisfies $x_{n_j}\in K_{n_j}$, for every $j\in \mathbb{N}$ and 
$$d\left(x_{n_{j-1}},x_{n_{j}}\right)<\frac{\varepsilon}{2^j},\quad\forall j\in \mathbb{N}.$$ By the latter and by \eqref{Hausd:15-11-22-4} we get

\begin{equation}\label{Hausd:15-11-22-5}
	d\left(y,x_{n_{j}}\right)\leq d\left(y,x_{n_{1}}\right)+d\left(x_{n_{1}}, x_{n_{2}}\right)+\cdots + d\left(x_{n_{j-1}}, x_{n_{j}}\right)<\varepsilon
\end{equation}
and 
\begin{equation*}
	d\left(x_{n_{j}}, x_{n_{h}}\right)\leq \sum_{l=j}^{h-1}d\left(x_{n_{l}}, x_{n_{l+1}}\right)<\frac{\varepsilon}{2^j},\quad\forall h>j\geq \nu_{\varepsilon}\in \mathbb{N}.
\end{equation*}
In particular, the just obtained inequality implies that for every $\delta>0$ there exists $n_{\delta}$ such that
$$d\left(x_{n_{j}}, x_{n_{h}}\right)<\delta, \quad \forall j,h\geq n_{\delta}.$$
Hence $\left\{x_{n_j}\right\}$ is a Cauchy sequance and it satisfies 
$$ x_{n_j}\in K_{n_j},\quad \forall j\in \mathbb{N}.$$ Now, Lemma \ref{Hausd:12-11-22-10} implies that there exists a Cauchy sequence $\left\{\overline{x}_n\in K_n\right\}$ which satisfies
$$ \overline{x}_{n_j}=x_{n_j}, \quad \forall j\in \mathbb{N}.$$ Consequently $\left\{\overline{x}_n\in K_n\right\}$ converges to a point $x$ and such a point $x$, by the definition of $K$, belongs to $K$. In addition, since  $\left\{x_{n_j}\right\}$ is a subsequence of $\left\{\overline{x}_n\right\}$, we have 
$$\left\{x_{n_j}\right\}\rightarrow x.$$
Hence, by \eqref{Hausd:15-11-22-5}, we have
$$d(x,y)=\lim_{j\rightarrow \infty}d\left(y,x_{n_{j}}\right)\leq \varepsilon.$$ Therefore
$$y\in [K]_{\varepsilon}.$$
Hence \eqref{Hausd:15-11-22-1} is proved. $\blacksquare$

\bigskip

\begin{theo}[\textbf{compactness}]\label{Hausd:15-11-22-6}
	If $(X,d)$ is a compact metric space, then $\left(\mathbf{K}(X), d_{\mathcal{H}}\right)$ is a compact metric space.
\end{theo} 
\textbf{Proof.} Since $(X,d)$ is a compact space it is complete and totally bounded. On the other hand, by Theorem \ref{Hausd:12-11-22-9}, $\left(\mathbf{K}(X),d_{\mathcal{H}}\right)$ is complete, hence for proving that it is compact, it suffices to prove that $\left(\mathbf{K}(X),d_{\mathcal{H}}\right)$ is totally bounded. 

Let $\varepsilon$ be any positive number, since $X$ is totally bounded, there exists a finite set
$F_{\varepsilon}$ which satisfies

\begin{equation}\label{Hausd:15-11-22-7}
	d\left(x,F_{\varepsilon}\right)<\varepsilon,\quad \forall x\in X.
\end{equation}
Let $\mathcal{G}_{\varepsilon}$ be the family of all subsets of $F_{\varepsilon}$ 
($\mathcal{G}_{\varepsilon}$ is finite because $F_{\varepsilon}$ is finite). Let $K\in \mathbf{K}(X)$. Let us consider the set 
$$G=\left\{p\in F_{\varepsilon}:\mbox{ } d(p,K)<\varepsilon\right\}.$$ 
We have $G\neq \emptyset$. As a matter of fact, \eqref{Hausd:15-11-22-7} implies that if $x\in K$ then there exists $y\in  F_{\varepsilon}$ such that 
$$d(x,y)<\varepsilon,$$
hence
$$d(y,K)\leq d(x,y)<\varepsilon,$$
therefore $y\in G$. 
Notice, that we have trivially $G\in\mathcal{G}_{\varepsilon}$ and 
\begin{equation}\label{Hausd:15-11-22-8}
	\delta(G,K)=\max_{p\in G}d(p, K)<\varepsilon.
\end{equation}
Now, let us prove
\begin{equation}\label{Hausd:15-11-22-9}
	\delta(K,G)=\max_{x\in K}d(x, G)<\varepsilon.
\end{equation}
Let $x\in K$. Relationship \eqref{Hausd:15-11-22-7} implies that there exists $y\in F_{\varepsilon}$ such that $$d(x,y)<\varepsilon,$$ consequently
$$d(y,K)\leq d(y,x)<\varepsilon.$$
Therefore $y\in G$ which yields
$$d(x,G)\leq d(x,y)<\varepsilon,$$
and \eqref{Hausd:15-11-22-9} follows. Hence
$$d_{\mathcal{H}}(K,G)=\max\left\{\delta(K,G),\delta(G,K)\right\}<\varepsilon.$$
All in all, we have proved 
$$\forall K\in \mathbf{K}(X) \ \ \exists G\in \mathcal{G}_{\varepsilon} \mbox{ such that } d_{\mathcal{H}}(K,G)<\varepsilon,$$
which, since $\varepsilon$ is arbitrary, implies that  $\mathbf{K}(X)$ is a compact metric space. $\blacksquare$

\section{The distance function \index{distance function}} \label{funzdist:17-11-22}
\index{distance function} \index{$d_{\partial \Omega}(x)$}

In this Section we will give some properties of the function 
$$\mathbb{R}^n\ni x\rightarrow d_{\partial \Omega}(x):=d(x,\partial\Omega),$$ 
where $\Omega$ is a bounded open set of $\mathbb{R}^n$ whose boundary is of class $C^{1,1}$. When there is no risk of ambiguity, we simply write $d(x)$. In Proposition \ref{Hausd:10-11-22-2} we have proved that $d_{\partial \Omega}(x)$ is a Lipschitz continuous function, because it satisfies the inequality

\begin{equation}\label{funzdist:17-11-22-1}
	\left\vert d_{\partial \Omega}(x)-d_{\partial \Omega}(y)\right\vert\leq \left\vert x-y\right\vert, \quad \forall x,y\in \mathbb{R}^n.
\end{equation}

We say that an open set $A$ of $\mathbb{R}^n$ enjoys the \textbf{interior ball property} \index{interior ball property}, if for every point $P\in \partial A$ there exists $P'\in A$ and $r>0$ such that
\begin{equation*}
\overline{B_r(P')}\cap \overline{\Omega}=\{P\}.
\end{equation*}
We say that $A$  enjoys the \textbf{exterior ball property} \index{exterior ball property} if  $\mathbb{R}^n\setminus\overline{\Omega}$ enjoys the property of the interior ball.

\medskip

The following Proposition holds true.

\begin{prop}\label{funzdist:17-11-22-2}
Let $\Omega$ be a bounded open set of $\mathbb{R}^n$ whose boundary is of class $C^{1,1}$ whith constants $r_0, M_0$. Then $\Omega$ enjoys the interior ball property and the exterior ball property. More precisely we have what follows. Denoting by
\begin{equation}\label{funzdist:17-11-22-2-0}
	\mu_0=\min\left\{\frac{1}{M_0}, M_0\right\},
\end{equation}
for any $P\in \partial \Omega$ and for any  $r\in \left(0,\mu_0r_0\right)$, we have 
\begin{equation}\label{funzdist:17-11-22-3}
	\overline{B_r(P-r\nu(P))}\cap \overline{\Omega}=\{P\}
\end{equation}
and 
\begin{equation}\label{funzdist:19-11-22-1}
	\overline{B_r(P+r\nu(P))}\cap \overline{\mathbb{R}^n\setminus\Omega}=\{P\},
\end{equation}
where $\nu(P)$ is the unit outward normal to $\partial\Omega$ in $P$
 \end{prop}
\textbf{Proof.} 
Let $P\in \partial \Omega$. Let us consider a local representation of $\partial \Omega$. Hence, let us assume $P=0$ and let us assume, up to a isometry, 

\begin{equation*}
	\Omega\cap Q_{r_0, 2M_0}=\left\{x\in Q_{r_0,2M_0}:\mbox{ }
	x_n>g(x') \right\},
\end{equation*}  
where $g\in C^{1,1}\left(\overline{B'}_{r_0}\right)$ satisfies
$$g(0)=|\nabla g(0)|=0$$
and
$$\left\Vert g\right\Vert_{C^{1,1}\left(\overline{B'_{r_0}}\right)}\leq
M_0r_0.$$
We have $$\nu(0)=-e_n.$$ Now, notice that
$$g\left(x'\right)=\int^1_0 (1-s)\partial^2 g\left(sx'\right)x'\cdot x'ds.$$
Hence 
$$g\left(x'\right)\leq \frac{M_0\left|x'\right|^2}{2r_0}.$$
Therefore, in order to satisfy \eqref{funzdist:17-11-22-3} it suffices that, besides the condition $r\leq r_0$, the following conditions are satisfied

$$\frac{M_0\left|x'\right|^2}{2r_0}<r-\sqrt{r^2-\left|x'\right|^2}, \quad \forall x'\in \overline{B'}_r\setminus \{0\}$$
and 
$$r+\sqrt{r^2-\left|x'\right|^2}< 2M_0r_0,\quad \forall x'\in \overline{B'}_r.$$ It is easy to check that if $r\in \left(0,\mu_0r_0\right)$, the above conditions are satisfied. In a similar way we proceed for the property of the exterior ball. $\blacksquare$

\bigskip

For any $\rho\in (0,\mu_0r_0)$, set

\begin{equation}\label{funzdist:19-11-22-2}
	S_{\rho}=\left\{x\in \Omega:\mbox{ } d_{\partial \Omega}(x)<\rho\right\}.
\end{equation} 
We observe  that for every $x\in S_{\rho}$ there exists a unique point $p(x)\in \partial\Omega$ such that 
\begin{equation*}
	|x-p(x)|=d_{\partial \Omega}(x).
\end{equation*} 
As a matter of fact, let $\overline{x}\in S_{\rho}$ and let $p\in \partial\Omega$ a point which satisfies 
$$d_{\partial \Omega}\left(\overline{x}\right)=\left|\overline{x}-p\right|.$$ We may assume  that $y$ belongs to the graph, $\Gamma^{(g)}$, of a function $g\in C^{1,1}\left(\overline{B'}_{r_0}\right)$ such that
\begin{equation*}
	\Omega\cap Q_{r_0, 2M_0}=\left\{x\in Q_{r_0,2M_0}:\mbox{ }
	x_n>g(x') \right\},
\end{equation*}  
and $g(0)=|\nabla g(0)|=0$.  Since $p$ is a minimum point on $\Gamma^{(g)}$  of the function 
$$y\rightarrow \frac{1}{2} \left|\overline{x}-y\right|^2.$$
By The Lagrange Multiplier Theorem, the following conditions need to be fulfilled

\begin{equation*}
	\begin{cases}
		p_j-\overline{x}_j+\lambda \partial_jg\left(p'\right)=0, \quad 1\leq j\leq n-1,\\
		\\
		p_n-\overline{x}_n-\lambda=0, \\
		\\
		g\left(p'\right)-p_n=0.
	\end{cases}%
\end{equation*}
Hence 
$$ p-\overline{x}=\left(p'-\overline{x}',p_n-\overline{x}_n\right)=\lambda \left(-\nabla_{p'}g\left(p'\right),1\right),$$
which implies
\begin{equation}\label{funzdist:19-11-22-4}
	p-\overline{x}=\left|p-\overline{x}\right|\nu(y).
\end{equation}
This relation, in turn, implies that $p$ is the unique point of $\partial\Omega$
that achieves the minimum distance. Indeed, let $B_{\rho}(z)$ be the interior ball tangent in $y$ to $\partial \Omega$ (such a ball exists by Proposition \ref{funzdist:19-11-22-1}), the equality  \eqref{funzdist:19-11-22-4} ensures us that  $\overline{x}$ lies on the segment of extremes $z$ and $p$. Consequently, by setting $\rho_1=\left|\overline{x}-p\right|$ we have $B_{\rho_1}\left(\overline{x}\right)\subset B_{\rho}\left(z\right)$. It is, therefore, evident that the distance of $\overline{x}$ from $\partial B_{\rho}\left(z\right)$ is greater than or equal to $\rho_1$ and, recalling that $B_{\rho}\left(\overline{y}\right)$ is an interior ball to $\Omega$, tangent to $\partial \Omega$ at the unique point $y$, we obtain that
$$ \left|\overline{x}-p\right|=\rho_1\leq \rho<\left|\overline{x}-\xi\right|,\quad \forall \xi\in \partial \Omega\setminus \{y\}.$$ 
Therefore we have proved

\begin{prop}\label{funzdist:19-11-22-5}
If $\rho\in (0,\mu_0r_0)$, then for any $x\in S_{\rho}$ there exists an unique point $p(x)\in \partial\Omega$ which attains the minimum of distance from $x$ to $\partial\Omega$. Moreover we have \begin{equation}\label{funzdist:19-11-22-7}
	x=p(x)-d_{\partial \Omega}(x)\nu(p(x)).
\end{equation}
\end{prop}

\bigskip

The following Proposition holds true 

\begin{prop}\label{funzdist:20-11-22-8}
	If $\rho\in (0,\mu_0r_0)$, then 
	\begin{equation}\label{funzdist:19-11-22-9}
		S_{\rho}=\left\{y-t\nu(y):\mbox{ } y\in \partial \Omega, \ \ 0\leq t<\rho\right\}.
	\end{equation}
\end{prop}
\textbf{Proof.} Proposition \ref{funzdist:19-11-22-5} implies 

\begin{equation*}
	S_{\rho}\subset \left\{y-t\nu(y):\mbox{ } y\in \partial \Omega, \ \ 0\leq t<\rho\right\}.
\end{equation*}

Now, let  $x=y-t\nu(y)$, where $y\in \partial \Omega$ e $0\leq t<\rho$. We have $$d_{\partial\Omega}(x)\leq |x-y|=t<\rho.$$
Hence $x\in S_{\rho}$. Therefore 
\begin{equation*}
	\left\{y-t\nu(y):\mbox{ } y\in \partial \Omega, \ \ 0\leq t<\rho\right\}\subset S_{\rho}.
\end{equation*}
Therefore \eqref{funzdist:19-11-22-9} is proved. $\blacksquare$ 

\bigskip

 Let us prove the following
\begin{lem}\label{funzdist:19-11-22-8}
	There exists $\mu_1\leq \mu_0$ such that if $\rho \in (0,\mu_1 r_0)$ then the maps
	\begin{equation}\label{funzdist:19-11-22-10}
		S_{\rho}\ni x \rightarrow p(x)\in\partial\Omega, \quad \mbox{ and } \quad S_{\rho}\ni x \rightarrow \nu(p(x))\in \mathbb{S}^{n-1}
	\end{equation}
	are Lipschitz continuous, where $p(x)$ is the point that realizes the minimum distance of $x\in S_{\rho}$ from $\partial \Omega$.
\end{lem}
\textbf{Proof.} We begin by proving that the map

\begin{equation}\label{funzdist:19-11-22-10-0}
	\partial\Omega\ni y \rightarrow \nu(y)\in \mathbb{S}^{n-1},
\end{equation}
is Lipschitz continuous. To prove this, let $y_1,y_2\in \partial\Omega$ and distinguish two cases

\smallskip

(a) $\left|y_1-y_2\right|\geq r_0$,

\smallskip

(b) $\left|y_1-y_2\right|< r_0$.

\smallskip

\smallskip

In case (a), we have trivially 
\begin{equation}\label{funzdist:19-11-22-10-00}
	\left|\nu\left(y_1\right)-\nu\left(y_2\right) \right|\leq 2\leq 2\frac{\left|y_1-y_2\right|}{r_0}.
\end{equation} 

In case (b), we may employ a local representation of $\partial \Omega$ assuming that $y_2=0$ and $y_1=g\left(x'\right)$ where $g\in C^{1,1}\left(\overline{B'}_{r_0}\right)$ and $g(0)=|\nabla g(0)|=0$. Hence
$$\nu\left(y_2\right)=-e_n$$
and
$$\nu\left(y_1\right)=\left(\frac{\nabla_{x'}g(x')}{\sqrt{1+\left|\nabla_{x'}g(x')\right|^2}},
\frac{-1}{\sqrt{1+\left|\nabla_{x'}g(x')\right|^2}}\right).$$
Now it is easy to check that 

\begin{equation}\label{funzdist:19-11-22-10-000}
	\left|\nu\left(y_1\right)-\nu\left(y_2\right) \right|\leq \frac{\sqrt{2} M_0}{r_0}\left|x'\right|=\frac{\sqrt{2}M_0}{r_0}\left|y_1-y_2\right|.
\end{equation} 
Therefore, by \eqref{funzdist:19-11-22-10-00} and \eqref{funzdist:19-11-22-10-000} we get
\begin{equation}\label{funzdist:19-11-22-10a}
	\left|\nu\left(y_1\right)-\nu\left(y_2\right) \right|\leq \frac{M_1}{r_0}\left|y_1-y_2\right|,
\end{equation}
where

$$M_1=\max\left\{2, \sqrt{2} M_0\right\}.$$

\medskip 

Now, let us prove that $x\rightarrow p(x)$ is Lipschitz continuous. By \eqref{funzdist:19-11-22-7} we have (we omit subscript in $d_{\partial \Omega}$)
$$p(x)=x+d(x)\nu(p(x)).$$
Therefore, recalling  \eqref{funzdist:17-11-22-1} and \eqref{funzdist:19-11-22-10a}, we get
\begin{equation*}
	\begin{aligned}
	|p(x)-p(y)|&\leq |x-y|+d(x)|\nu(p(x))-\nu(p(y))|+|d(x)-d(y)||\nu(p(y))|\leq\\&\leq 
	2|x-y|+\rho \frac{M_1}{r_0} |p(x)-p(y)|.
	\end{aligned}
\end{equation*}
Hence
\begin{equation*}
	\left(1-\rho \frac{M_1}{r_0}\right)|p(x)-p(y)|\leq 2|x-y|. 
\end{equation*}
Moreover for any 
$$\rho<\min\left\{\frac{1}{2M_1}, M_0, \frac{1}{M_0}\right\}$$
we have 
\begin{equation}\label{funzdist:19-11-22-13}
	|p(x)-p(y)|\leq 4|x-y|, \quad \forall x,y\in S_{\rho}. 
\end{equation}
The above inequality proves that the map $x\rightarrow p(x)$ is Lipschitz continuous provided $$\mu_1=\min\left\{\frac{1}{2M_1}, M_0, \frac{1}{M_0}\right\}.$$
Therefore  \eqref{funzdist:19-11-22-10a} and \eqref{funzdist:19-11-22-13} imply that $\nu(p(x))$ is Lipschitz continuous. $\blacksquare$

\bigskip

\begin{lem}\label{funzdist:19-11-22-14}
	Let $A$ be on open set of $\mathbb{R}^n$ and let $f\in C_{loc}^{0,1}(A)$. If there exists a function $g\in C^0(A;\mathbb{R}^n)$ which satisfies
	\begin{equation}\label{funzdist:19-11-22-15}
		\nabla f(x)= g(x),\quad \mbox{ a.e. } x\in A,
	\end{equation}
	then $f\in C^1(A)$ and 
	\begin{equation}\label{funzdist:20-11-22-3}
\nabla f(x)= g(x),\quad  \forall x\in A.
	\end{equation}
\end{lem}
\textbf{Proof.} Let $x_0\in A$ and $\delta=\frac{1}{4}d(x,\partial A)$. For any $\varepsilon\in (0,\delta)$ let us consider the function
\begin{equation*}
	f_{\varepsilon}(x)=\int_A f(y)\eta_{\varepsilon}(x-y)dy,\quad \forall x\in B_{\delta}(x_0),
\end{equation*}
where $\eta$ is a mollifier. It turns out that $f_{\varepsilon}\in C^{\infty}\left(\overline{B_{\delta}(x_0)}\right)$ and, moreover, the divergence Theorem gives
\begin{equation*}
	\begin{aligned}
	\partial_jf_{\varepsilon}(x)&=-\int_A f(y)\partial_{y_j}\eta_{\varepsilon}(x-y)dy=\\&=-\int_{B_{\varepsilon}(x)} \left[\partial_{y_j}\left(f(y)\eta_{\varepsilon}(x-y)\right)-\partial_{y_j}f(y)\eta_{\varepsilon}(x-y)\right]dy=\\&=\int_{B_{\varepsilon}(x)}\partial_{y_j}f(y)\eta_{\varepsilon}(x-y)dy=\\&=\int_A\partial_{y_j}f(y)\eta_{\varepsilon}(x-y)dy, \quad j=1,\cdots n, \quad \forall x\in B_{\delta}(x_0) .
	\end{aligned}
\end{equation*}
By what has just been obtained and by \eqref{funzdist:19-11-22-15} we have

\begin{equation}\label{funzdist:20-11-22-4}
	\nabla f_{\varepsilon}(x)=g_{\varepsilon}(x):=\int_Ag(y)\eta_{\varepsilon}(x-y)dy, \quad \forall x\in B_{\delta}(x_0).
\end{equation}
Let $v$ be a versor of $\mathbb{R}^n$. By \eqref{funzdist:20-11-22-4} we get 
\begin{equation}\label{funzdist:20-11-22-5}
	f_{\varepsilon}(x_0+tv)-f_{\varepsilon}(x_0)=\int^t_0g_{\varepsilon}(x_0+sv)\cdot vds, \quad \forall t\in [-\delta,\delta].
\end{equation}
Now, Theorem \ref{Sob:teo25R} implies that $f_{\varepsilon}$ and $g_{\varepsilon}$ uniformly converge  in $B_{\delta}(x_0)$. Therefore passing to the limit in \eqref{funzdist:20-11-22-5} as $\varepsilon\rightarrow 0$, we obtain 
\begin{equation*}
	f(x_0+tv)-f(x_0)=\int^t_0g(x_0+sv)\cdot vds, \quad \forall t\in [-\delta,\delta].
\end{equation*}
On the other hand, since $g$ is continuous, we have

$$\lim_{t\rightarrow 0}\frac{f(x_0+tv)-f(x_0)}{t}=\lim_{t\rightarrow 0}\frac{1}{t} \int^t_0 g(x_0+sv)\cdot vds=g(x_0)\cdot v.$$
Therefore
$$\frac{\partial f(x_0)}{\partial v}=g(x_0)\cdot v, \ \ \forall v\in \mathbb{R}^n, \ \ |v|=1.$$
which implies
$$\nabla f(x_0)= g(x_0)$$ 
so that, since $x_0$ is arbitrary in $A$, we have $\nabla f= g$ in $A$. That, in turn, by the continuity of  $g$ implies $f\in C^1(A)$. $\blacksquare$

\bigskip

\begin{theo}\label{funzdist:19-11-22-16}
Let $\mu_1$ be the same of Proposition \ref{funzdist:19-11-22-8}, then we have 	$$d_{\partial \Omega}\in C^{1,1}\left(\overline{S_{\mu_1r_0}}\right)$$ and 
	\begin{equation}\label{funzdist:19-11-22-17}
		\nabla d_{\partial \Omega}(x)=-\nu(p(x)), \quad \forall x\in S_{\mu_1r_0}.
	\end{equation} 
\end{theo}
\textbf{Proof.} For the sake of brevity, we omit the subscript in $d_{\partial \Omega}$. Inequality \eqref{funzdist:17-11-22-1} implies that $d$ is almost everywhere  differentiable, in addition in the points where it is differentiable we have
\begin{equation}\label{funzdist:20-11-22-1}
	|\nabla d|\leq 1. 
\end{equation}
Let $x\in S_{\mu_1r_0}\setminus\partial\Omega$ be a point in which $d$ is differentiable. Let $t\in (0, d(x))$. Since $x+t\nu(p(x))$ lies on the segment of endpoints $x$ and $p(x)$, we have
$$d(x+t\nu(p(x)))=d(x)-t.$$
Hence
\begin{equation*}
	\nabla d(x)\cdot \nu(p(x))=\lim_{t\rightarrow 0^+}\frac{d(x+t\nu(p(x)))-d(x)}{t}=-1.
\end{equation*}
Therefore
\begin{equation}\label{funzdist:20-11-22-2}
	\nabla d(x)\cdot \nu(p(x))=-1.
\end{equation}
Consequently we get
\begin{equation*}
1=|\nabla d(x)\cdot \nu(p(x))|\leq |\nabla d(x)| |\nu(p(x))|\leq 1.
\end{equation*}
Hence, there exits $\lambda\in \mathbb{R}$ such that $\nabla d(x)=\lambda  \nu(p(x))$ and by \eqref{funzdist:20-11-22-2} we have $\lambda=-1$. This implies

\begin{equation}\label{funzdist:20-11-22-6}
	\nabla d(x)=-\nu(p(x)), \quad \mbox{ a.e. } x\in S_{\mu_1r_0}\setminus\partial\Omega.
\end{equation}
Now, by Lemma \ref{funzdist:19-11-22-8} we know that $\nu(p(x))$ is Lipschitz continuous, therefore by Lemma \ref{funzdist:19-11-22-14} we obtain $d\in C^1\left(S_{\mu_1r_0}\setminus\partial\Omega\right))$. Finally, exploit again \eqref{funzdist:20-11-22-6}, we have $d\in C^{1,1}\left(\overline{S_{\mu_1r_0}}\right)$.  $\blacksquare$

\bigskip

\begin{cor}\label{funzdist:20-11-22-10}
For any $\rho\in (0,\mu_1r_0)$ the boundary of the open set 
\begin{equation}\label{funzdist:20-11-22-11}
	(\Omega)_{\rho}=\left\{x\in \Omega: \mbox{ } d_{\partial \Omega}(x)>\rho \right\}
\end{equation}
is of class $C^{1,1}$ and we have

\begin{equation}\label{funzdist:20-11-22-12}
	\partial(\Omega)_{\rho}=\Gamma_{\rho},
\end{equation}
where
\begin{equation}\label{funzdist:20-11-22-11-01}
	\Gamma_{\rho}=\left\{x\in \Omega: \mbox{ } d_{\partial \Omega}(x)=\rho \right\}.
\end{equation}
 Moreover
 \begin{equation}\label{funzdist:20-11-22-11-02}
 	\Gamma_{\rho}=\left\{y-\rho\nu(\rho): \mbox{ } y\in \partial\Omega \right\}.
 \end{equation}

\end{cor}
\textbf{Proof.} 
We first prove \eqref{funzdist:20-11-22-12}. To prove that $\Gamma_{\rho}\subset\partial(\Omega)_{\rho}$ we argue by contradiction. Let us assume that there exists $x\in \partial(\Omega)_{\rho}$ sucht that $d_{\partial \Omega}(x)>\rho$. Consequently, $x$ would be an interior point of  $(\Omega)_{\rho}$. If $d_{\partial \Omega}(x)<\rho$ then $x$ would be exterior to $(\Omega)_{\rho}$. Therefore, $x\in  \Gamma_{\rho}$ and we have $\partial(\Omega)_{\rho}\subset\Gamma_{\rho}$. Now, let $x\in \Gamma_{\rho}$, since $\rho\in (0,\mu_1r_0)$ (and $\mu_1\leq \mu_0$), Proposition \ref{funzdist:19-11-22-5} implies that there is an unique point $p(x)\in \partial \Omega$ which attains the minimum of distance of $x$ from $\partial\Omega$, moreover 
$$x=p(x)-\rho\nu(p(x)).$$
For any $\epsilon>0$ small enough, we have
$$x-\varepsilon \nu(p(x)=p(x)-(\rho+\varepsilon)\nu(p(x))\in (\Omega)_{\rho}$$
and 
$$x+\varepsilon \nu(p(x)=p(x)-(\rho-\varepsilon)\nu(p(x))\notin (\Omega)_{\rho}.$$
Hence $x\in \partial(\Omega)_{\rho}$. Therefore  $\Gamma_{\rho}\subset\partial(\Omega)_{\rho}$.

In order to prove that $\partial(\Omega)_{\rho}$ is of class $C^{1,1}$, we exploit  \eqref{funzdist:20-11-22-12}. By Theorem \ref{funzdist:19-11-22-16} we derive $\left|\nabla d_{\partial\Omega}(x)\right|=1$, for every $x\in \Gamma_{\rho}$, and by applying Implicit Function Theorem  we easily reach the assertion. 

Concerning \eqref{funzdist:20-11-22-11-02}, let us note that if $x=y-\rho\nu(y)$, taking into account $\rho<\mu_1r_0\leq \mu_0r_0$, then $d_{\partial\Omega}(x)=\rho$. Conversely, if $x\in \Gamma_{\rho}$  Proposition \ref{funzdist:19-11-22-5} gives
$$x=p(x)-d_{\partial\Omega}(x)\nu(p(x))=(x)-\rho\nu(p(x))$$ 
so that, since $p(x)\in \partial\Omega$ we have $x\in \left\{y-\rho\nu(\rho): \mbox{ } y\in \partial\Omega \right\}$. $\blacksquare$ 

\bigskip

We now provide a few words about the map
\begin{equation*}
\Phi: \partial\Omega\times \left(0,\mu_1r_0\right)\rightarrow \ \mathbb{R}^n,\quad 
\end{equation*}
such that
\begin{equation}\label{funzdist:20-11-22-13}
	\Phi(y,t)=y-t\nu(y), \ \ \forall (y,t)\in \partial\Omega\times \left(0,\mu_1r_0\right).
\end{equation}

\medskip

The following Proposition holds true

\begin{prop}\label{funzdist:20-11-22-14}
If $\Omega$ is a bounded open set of $\mathbb{R}^n$ of class $C^{1,1}$ then we have: 

\smallskip

(a) $\Phi\left(\partial\Omega\times \left(0,\mu_1r_0\right)\right)=S_{\mu_1r_0}$,

\smallskip

(b) $\Phi\in C^{0,1}\left(\partial\Omega\times \left[0,\mu_1r_0\right]\right)$,

\smallskip

(c) $\Phi$ is injective on $\partial\Omega\times \left(0,\mu_1r_0\right)$ and it inverse is Lipschitz continuous map.
\end{prop}
\textbf{Proof.}
(a) is a consequence of Proposition \ref{funzdist:20-11-22-8}. 
 (b) is a consequence of Lemma \ref{funzdist:19-11-22-8}. Now let us prove (c). Let $x\in S_{\mu_1r_0}$ satisfy
 $$x=\Phi(y,t)=y-t\nu(y), \quad (y,t)\in \partial\Omega\times \left(0,\mu_1r_0\right).$$
 By Proposition \ref{funzdist:19-11-22-5} and by the interior ball property we get
 $$y=p(x),\quad\quad t=d_{\partial \Omega}(x).$$
 Hence 
 $$\Phi^{-1}(x)=p(x)-d_{\partial \Omega}(x)\nu(p(x)).$$
 By the latter and by  Lemma \ref{funzdist:19-11-22-8} it follows that $\Phi^{-1}$ is Lipschitz continuous. $\blacksquare$
 
\bigskip

If $\Omega$ is of class $C^k$, $k\geq 2$, other properties of the distance function and the map $\Phi$ can be proved. For instance, one can prove that $d_{\partial\Omega}\in C^k$ and $\Phi\in C^{k-1}$. For further details, we refer to \cite[Ch. 14, Sect. 6]{G-T}.

\bigskip

We say that a continuous map 
$$\gamma:[0,1]\rightarrow A,$$
is a \textbf{continuous path in a set $A\subset \mathbb{R}^n$} \textbf{continuous path}.
Let $B\subset A$ and $x,y\in B$, if $\gamma([0,1])\subset B$ and $\gamma(0)=x$, $\gamma(1)=y$, we say that the path $\gamma$ joins  $x$ and $y$ in $B$. If $\gamma_1$ e $\gamma_2$ are two continuous paths in $A$ which satisfy $\gamma_1(1)=\gamma_2(0)$, we denote by $\gamma_1\vee \gamma_2$ the following continuous path

\begin{equation}\label{funzdist:26-11-22-1-0}
	\left(\gamma_1\vee \gamma_2\right)(t)=\begin{cases}
	\gamma_1(2t)	, \quad\mbox{for } t\in \left[0,\frac{1}{2}\right), \\
		\\
		\gamma_2(2t-1), \quad\mbox{for } t\in \left[\frac{1}{2},1\right].
	\end{cases}
\end{equation}
If $\gamma_1, \cdots \gamma_k$ are $k\geq 2$ continuous paths in $A$ such that $\gamma_{j-1}(1)=\gamma_j(0)$, $j=2,\cdots k$, we set
$$\gamma_1\vee \cdots \vee\gamma_k:=\left(\gamma_1\vee \cdots \vee\gamma_{k-1}\right)\vee\gamma_k.$$ 
We say that $\gamma_1\vee \cdots \vee\gamma_k$ is the  \index{union of continuous paths}  $\gamma_1, \cdots \gamma_k$.

\bigskip

\begin{prop}\label{funzdist:26-11-22-1}
Let us assume that $\Omega$ and $\partial \Omega$ are connected. If $\rho\in (0,\mu_1r_0)$, then the $(\Omega)_{\rho}$, defined by \eqref{funzdist:20-11-22-11}, is connected.
\end{prop}
\textbf{Proof.} Let $z,w\in (\Omega)_{\rho}$ and let $\varepsilon>0$ such that $$\rho+\varepsilon<\min\left\{\mu_1r_0, d_{\partial\Omega}(z), d_{\partial\Omega}(w) \right\}.$$ Then $$z,w\in (\Omega)_{\rho+\varepsilon},$$
\begin{equation}\label{funzdist:26-11-22-2}
\Gamma_{\rho+\varepsilon}\subset (\Omega)_{\rho},
\end{equation}
and $\Gamma_{\rho+\varepsilon}$ is connected, as it is the image by $\Phi$, defined in Proposition \ref{funzdist:20-11-22-14}, of the connected set $\partial\Omega\times\{\rho+\varepsilon\}$. 

Now, since $\Omega$ is connected, $\overline{\Omega}$ is also connected (path connected, because $\partial \Omega$ \`{e} of class $C^{1,1}$). Be, therefore, $x\in \partial \Omega$ and be \\ $$\gamma_1:[0,1]\rightarrow \overline{\Omega}\quad \mbox{and}\quad \gamma_2:[0,1]\rightarrow \overline{\Omega}$$ two continuous paths such that 
$$\gamma_1(0)=z,\ \ \gamma_1(1)=x,\quad \gamma_2(0)=x,\ \ \gamma_2(1)=w.$$
Let 
$$t_1=\inf\left\{t\in [0,1]:\mbox{ } d\left(\gamma_1(t),\partial\Omega\right)<\rho+\varepsilon \right\},$$
we have (because $d\left(\gamma_1(\cdot),\partial\Omega\right)$ is continuous) 
$$y':=\gamma_1(t_1)\in \Gamma_{\rho+\varepsilon}.$$ 
Similarly, let
$$t_2=\sup\left\{t\in [0,1]:\mbox{ } d\left(\gamma_2(t),\partial\Omega\right)<\rho+\varepsilon \right\},$$
we have
$$y'':=\gamma_2(t_2)\in \Gamma_{\rho+\varepsilon}.$$
Since $\Gamma_{\rho+\varepsilon}$ is connected, there exists a continuous path $\widetilde{\gamma}:[0,1]\rightarrow \Gamma_{\rho+\varepsilon}$, such that
$$\widetilde{\gamma}(0)=y',\ \ \widetilde{\gamma}(1)=y''.$$ 
It is now evident that the path 
$$\gamma:=\gamma_1\vee \widetilde{\gamma}\vee \gamma_2,$$ 
is continuous and it joins $z$ e $w$ in $(\Omega)_{\rho}$. $\blacksquare$ 

\bigskip

\textbf{Remark.} In Proposition \ref{funzdist:26-11-22-1}, the assumption that $\partial \Omega$ is connected is not necessary.  The proof of this assertion may follow arguing likewise  the  proof of Proposition \ref{funzdist:26-11-22-1}, taking into account that  due to the boundedness of $\Omega$ and the $C^{1,1}$ character of $\partial \Omega$, the connected components of $\partial \Omega$ are finite in number. We invite the reader to develop the details. 
$\blacklozenge$

\chapter{The Sobolev spaces}\label{prerequisiti}

\section{Weak derivatives} \label{definizione H-k}

Let us give the definition of weak derivative.
\begin{definition}\label{weak-der}
	Let $\Omega$ be an open set of $\mathbb{R}^n$ and $\alpha\in
	\mathbb{N}_0^n$. Let $u, v\in L_{loc}^1(\Omega)$. We say that
	$v$ is the \textbf{$\alpha$--th weak derivative} \index{weak derivative} of $u$ and we write
	
	$$\partial^{\alpha}u=v,$$
	if
	\begin{equation}\label{derivata debole}
		\int_{\Omega}u\partial^{\alpha}\phi dx=(-1)^{\alpha}
		\int_{\Omega}v\phi dx, \quad\quad \forall\phi\in
		C^{\infty}_0(\Omega).
	\end{equation}
\end{definition}

\medskip

Definition \ref{weak-der} is justified by the integration by parts formula  that, in the case of $u\in
C^{|\alpha|}(\Omega)$, gives precisely the derivative
$\partial^{\alpha}u$ in the classical sense. For instance, if $u\in
C^{1}(\Omega)$,  we have

\begin{equation*}
	\int_{\Omega}u\partial_j\phi dx=
	\int_{\Omega}\left[\partial_j\left(u\phi \right)-\phi\partial_j
	u\right]dx=-\int_{\Omega}\phi\partial_ju  dx, \quad\quad
	\forall\phi\in C^{\infty}_0(\Omega).
\end{equation*}

\bigskip

\begin{prop}\label{Sob:1.1}
	If $u\in L_{loc}^1(\Omega)$ admits the $\alpha$--th weak derivative, it is unique (up to a set of measure zero).
	\end{prop}
\textbf{Proof.} Let us assume that $v_1,v_2\in
L_{loc}^1(\Omega)$ are two $\alpha$--th weak derivative of $u$,
then
\begin{equation*}
	(-1)^{\alpha} \int_{\Omega}v_1\phi
	dx=\int_{\Omega}u\partial^{\alpha}\phi dx=(-1)^{\alpha}
	\int_{\Omega}v_2\phi dx, \quad\quad \forall\phi\in
	C^{\infty}_0(\Omega),
\end{equation*}
which implies

\begin{equation*}
	\int_{\Omega}\left(v_1-v_2\right)\phi dx=0, \quad\quad
	\forall\phi\in C^{\infty}_0(\Omega),
\end{equation*}
however, $v_1-v_2 \in L_{loc}^1(\Omega)$, so
$$v_1=v_2,\quad \mbox{ a.e. in } \Omega.$$ $\blacksquare$

\bigskip

\textbf{Example 1.} Let $\Omega=(-1,1)$, $u(x)=|x|$; let us show that
$$u'=\mbox{sgn}(x):=\begin{cases}
	1, \quad \mbox{ for } x>0,\\
	\\
	0, \quad \mbox{ for } x=0,\\
	\\
	-1, \quad \quad\mbox{for }  x<0,
\end{cases}
,\quad\mbox{ in the weak sense}.$$

As a matter of fact we have sgn$(\cdot)$$\in L^1(-1,1)$ and  

\begin{equation*}
	\begin{aligned}
		\int^1_{-1}|x|\phi'(x)dx&=\int^1_{0}x\phi'(x)dx-\int^{0}_{-1}x\phi'(x)dx=\\&
		=\left[x\phi(x)\right]^1_0-\int^1_{0}\phi(x)dx-\left[x\phi(x)\right]^0_{-1}+\int^0_{-1}\phi(x)dx=\\&
		=-\int^1_{-1}\mbox{sgn}(x)\phi(x)dx,\quad\quad \forall \phi\in
		C_0^{\infty}(-1,1).
	\end{aligned}
\end{equation*}

$\spadesuit$

\bigskip

\textbf{Example 2.} Let $\Omega=(-1,1)$, $u(x)=$sgn$(x)$. Let us prove that
 $u$ has \textbf{not} the weak derivative. Let us assume the contrary and be $v\in
L_{loc}^1(-1,1)$ such that

\begin{equation}\label{Sob:1.2}
	\int^1_{-1}u(x)\phi'(x)dx=-\int^1_{-1}v(x)\phi(x)dx,\quad\quad
	\forall \phi\in C_0^{\infty}(-1,1).
\end{equation}
Let $\phi\in C_0^{\infty}(-1,1)$ arbitrary. We have 

\begin{equation*}
	\begin{aligned}
		\int^1_{-1}u(x)\phi'(x)dx=\int^1_{0}\phi'(x)dx-\int^0_{-1}\phi(x)dx=-2\phi(0).
	\end{aligned}
\end{equation*}
Taking into account \eqref{Sob:1.2}, we get 

\begin{equation}\label{Sob:1.2bis}
	\int^1_{-1}v(x)\phi(x)dx=2\phi(0),\quad\quad \forall \phi\in
	C_0^{\infty}(-1,1).
\end{equation}
Let now $\left\{\phi_k\right\}_{k\geq 2}$ be the following sequence of functions 

\begin{equation*}
	\phi_k(x)=
	\begin{cases}
		e^{k^2-\frac{k^2}{1-k^2x^2}}, \quad \mbox{ for } |x|<\frac{1}{k},\\
		\\
		0, \quad \quad\mbox{for } \frac{1}{k}\leq |x|<1.
	\end{cases}%
\end{equation*}
we have $\phi_k\in C_0^{\infty}(-1,1)$,
supp $\phi_k\subset\left[-\frac{1}{2}, \frac{1}{2}\right]$ for every
$k\geq 2$. Moreover

$$\phi_k(0)=1,\quad \mbox{and}\quad \lim_{k\rightarrow
	\infty}\phi_k(x)=0,\quad \mbox{ for } x\neq 0.$$ On the other hand, by
\eqref{Sob:1.2bis} we have

\begin{equation}\label{Sob:1.3}
	2=2\phi_k(0)=\int^1_{-1}v(x)\phi_k(x)dx,\quad\quad \forall k\geq 2,
\end{equation}
but $v\in L_{loc}^1(-1,1)$, hence the Dominated Convergence Theorem implies

$$\lim_{k\rightarrow
	\infty}\int^1_{-1}v(x)\phi_k(x)dx=0.$$ By the latter and by \eqref{Sob:1.3} we reach a contradiction. $\spadesuit$

\section{Definition of the Sobolev spaces} \label{Sob:paragrafo2}
Let us give the following
\begin{definition}\label{Sob:def1.2}
	Let $1\leq p\leq \infty$, $k\in \mathbb{N}_0$ and let $\Omega$ be an open set of $\mathbb{R}^n$, $n\geq 1$. If $k=0$, set \index{$W^{k,p}(\Omega)$}
	$$W^{0,p}(\Omega)=L^{p}(\Omega).$$
	If $k\geq 1$, $W^{k,p}(\Omega)$ is the set of functions $u
	\in L_{loc}^1(\Omega)$ satisfying
	\begin{equation}\label{Sob:2.3}
		\partial^{\alpha}u\in L^{p}(\Omega), \quad\quad\mbox{for }
		|\alpha|\leq k,
	\end{equation}
	where $\partial^{\alpha}u$  is the $\alpha$--th weak derivative of 
	$u$.
\end{definition}

\medskip

It is easy to check that $W^{k,p}(\Omega)$ is a vector space. Furthermore we define the following norms. If $1\leq
p<+\infty$, we set
\index{$\left\Vert \cdot\right\Vert_{W^{k,p}(\Omega)}$}
\begin{equation}\label{Sob:1.4}
	\left\Vert u\right\Vert_{W^{k,p}(\Omega)}=\left(\sum_{|\alpha|\leq
		k}\int_{\Omega}\left|\partial^{\alpha}u\right|^{p} dx\right)^{1/p}.
\end{equation}
If $p=+\infty$, we set
\begin{equation}\label{Sob:2.4}
	\left\Vert u\right\Vert_{W^{\infty,p}(\Omega)}=\sum_{|\alpha|\leq
		k}\left\Vert \partial^{\alpha}u\right\Vert_{L^{\infty}(\Omega)}.
\end{equation}

\medskip

\noindent If $p=2$, we also set \index{$H^k(\Omega)$, ($k$ positive integer number)}
$$H^k(\Omega)=W^{k,2}(\Omega).$$
Let us observe that $H^k(\Omega)$ is a pre--Hilbertian space equipped with the
scalar product \index{$(\cdot,\cdot)_{H^k(\Omega)}$, ($k$ positive integer number)}

\begin{equation}\label{Sob:3.4}
	(u,v)_{H^k(\Omega)}=
	\int_{\Omega}\sum_{|\alpha|\leq k} \partial^{\alpha}u
	\partial^{\alpha}v dx, \quad \forall u,v\in H^k(\Omega).
\end{equation}

\medskip

Here and in the sequel, for any $k\in \mathbb{N}$, $p\in [1,\infty]$,
we denote by $W^{k,p}_{loc}(\Omega)$, ($H^{k}_{loc}(\Omega)$) \index{$W^{k,p}_{loc}(\Omega)$, $H^{k}_{loc}(\Omega)$ }the subspace of 
$L^{p}_{loc}(\Omega)$ ($L^{2}_{loc}(\Omega)$) of the functions $u$ such that for every open set
$\omega\Subset \Omega$ (i.e. $\overline{\omega}\subset \Omega$) we have $u_{|\omega}\in
W^{k,p}(\omega)$ ($u_{|\omega}\in
H^{k}(\omega)$). Let $\left\{u_m\right\}$ be a sequence in $W_{loc}^{k,p}(\Omega)$ and $u\in W_{loc}^{k,p}(\Omega)$, we say that \index{convergence of a sequence in $W_{loc}^{k,p}(\Omega)$}

$$u_m\rightarrow u,\quad \mbox{ as } m\rightarrow\infty,\mbox{ in }
W_{loc}^{k,p}(\Omega),$$ if

$$(u_m)_{|\omega}\rightarrow u_{|\omega},
\quad \mbox{ as } m\rightarrow\infty,\mbox{ in }
W^{k,p}(\omega),\quad\forall \omega\Subset \Omega.$$

\bigskip

\underline{\textbf{Exercise 1.}} Check that, if $1\leq p\leq \infty$,
$k\in \mathbb{N}_0$, then $W^{k,p}(\Omega)$ is a vector subspace of $L^p(\Omega)$ and $\left\Vert
\cdot\right\Vert_{W^{k,p}(\Omega)}$ defines a norm on
$W^{k,p}(\Omega)$. $\clubsuit$

\bigskip

\begin{prop}\label{Sob:prop1.2.6}
	If $u\in W^{k,p}(\Omega)$, then we have
	
	\noindent (i) $\partial^{\alpha}u\in W^{k-|\alpha|,p}(\Omega) $ for
	$|\alpha|\leq k$ and
	$\partial^{\beta}\partial^{\alpha}u=\partial^{\alpha}\partial^{\beta}u=\partial^{\alpha+\beta}u$
	for $|\alpha|+|\beta|\leq k$,
	
	\noindent (ii) for any $\zeta \in
	C^{\infty}\left(\overline{\Omega}\right)$ we have $\zeta u \in
	W^{k,p}(\Omega)$ and
	
	$$\partial^{\alpha}(\zeta u)=\sum_{\beta\leq\alpha}\binom{{\alpha}}{{\beta}}\partial^{\beta}\zeta\partial^{\alpha-\beta}u.$$
\end{prop}
\textbf{Proof.} (i) Let $u\in W^{k,p}(\Omega)$, $|\alpha|\leq k$, and let
$\beta$ satisfy $|\beta|\leq k-|\alpha|$. For any $\phi\in
C^{\infty}_0(\Omega)$, we have

\begin{equation*}
	\begin{aligned}
		\int_{\Omega}\partial^{\alpha}u\partial^{\beta}\phi
		dx&=(-1)^{|\alpha|}\int_{\Omega}u\partial^{\alpha+\beta}\phi dx=\\&
		=(-1)^{|\alpha|}(-1)^{|\alpha|+|\beta|}\int_{\Omega}\partial^{\alpha+\beta}u\phi
		dx=\\&=(-1)^{|\beta|}\int_{\Omega}\partial^{\alpha+\beta}u\phi dx.
	\end{aligned}
\end{equation*}
Hence
$$\int_{\Omega}\partial^{\alpha}u\partial^{\beta}\phi
dx=(-1)^{|\beta|}\int_{\Omega}\partial^{\alpha+\beta}u\phi dx,\quad
\forall \phi\in C^{\infty}_0(\Omega),$$ consequently
$$\partial^{\beta}\partial^{\alpha}u=\partial^{\alpha+\beta}u.$$ The latter implies $$\partial^{\alpha}u\in
W^{k-|\alpha|,p}(\Omega),$$ for any $|\alpha|\leq k$.

\medskip

\noindent (ii) Let us consider the case $|\alpha|=1$.  Let  $\alpha=e_j$, for
$j=1,\cdots, n$. We have, for any $\phi\in C^{\infty}_0(\Omega)$,

\begin{equation*}
	\begin{aligned}
		\int_{\Omega}\zeta u\partial_{j}\phi
		dx&=\int_{\Omega}u\left[\partial_{j}(\zeta\phi)-(\partial_{j}\zeta)\phi\right]
		dx=\\& =-\int_{\Omega}(\partial_{j}u)\zeta\phi
		dx-\int_{\Omega}u(\partial_{j}\zeta)\phi dx =\\&=
		-\int_{\Omega}\left[(\partial_{j}u)\zeta+u\partial_{j}\zeta\right]\phi
		dx.
	\end{aligned}
\end{equation*}
Now, let us notice that
$$(\partial_{j}u)\zeta+u\partial_{j}\zeta\in L^p(\Omega),$$ hence
$$\partial_j(\zeta u)=(\partial_{j}u)\zeta+u\partial_{j}\zeta,\quad \mbox{ in the weak sense}.$$
If $|\alpha|>1$, one proceeds by induction, and we leave the details to the
reader. $\blacksquare$

\bigskip

\begin{theo}[\textbf{completeness of $W^{k,p}(\Omega)$}]\label{Sob:2.2.8}
\index{Theorem:@{Theorem:}!- completeness of $W^{k,p}(\Omega)$@{- completeness of $W^{k,p}(\Omega)$}}	
	The space $W^{k,p}(\Omega)$, $k\in \mathbb{N}_0$, $1\leq p\leq
	\infty$, equipped with the norm \eqref{Sob:1.4}, \eqref{Sob:2.4}, is a Banach space. If $p=2$, $H^{k}(\Omega)$ is a Hilbert space.
	\end{theo}
\textbf{Proof.} We limit ourselves to the case $k=1$. Similarly it can be handle the
case $k>1$.
Let $\left\{u_m\right\}$ be a Cauchy sequence in
$W^{1,p}(\Omega)$. From the definition of norm of
$W^{1,p}(\Omega)$ we have that
$$\left\{u_m\right\}\quad\mbox{ and }\quad
\left\{\partial_ju_m\right\},\quad j=1,\cdots,n,$$ are Cauchy sequences in $L^p(\Omega)$. On the other hand, $L^p(\Omega)$ is
complete; hence there exist $u$, $v_1,\cdots,v_n \in L^p(\Omega)$ such that

\begin{subequations}
	\label{Sob:1.8}
	\begin{equation}
		\label{Sob:1.8a} u_m\rightarrow u,\quad\mbox{ as }
		m\rightarrow\infty,\quad\mbox{in } L^p(\Omega),
	\end{equation}
	\begin{equation}
		\label{Sob:1.8b} \partial_ju_m\rightarrow
		v_j,\quad\mbox{ as } m\rightarrow\infty, \quad\mbox{in }
		L^p(\Omega),\mbox{ } j=1,\cdots,n.
	\end{equation}
\end{subequations}
Now, \eqref{Sob:1.8a} and \eqref{Sob:1.8b} imply that, for any $\phi\in
C^{\infty}_0(\Omega)$, we have
\begin{equation*}
	\begin{aligned}
		\int_{\Omega} u\partial_{j}\phi dx&=\lim_{m\rightarrow
			\infty}\int_{\Omega}u_m \partial_{j}\phi dx=\\& =-\lim_{m\rightarrow
			\infty}\int_{\Omega}\partial_{j}u_m\phi dx=\\&=
		-\int_{\Omega}v_j\phi dx.
	\end{aligned}
\end{equation*}
Hence
$$\partial_ju=v_j,\quad\mbox{ for } j=1,\cdots,n. $$
Therefore by \eqref{Sob:1.8a} e \eqref{Sob:1.8b} we have
$$u_m\rightarrow u,\quad\mbox{ as }
m\rightarrow\infty, \quad \mbox{in }W^{1,p}(\Omega).$$
$\blacksquare$

\bigskip

\begin{prop}\label{Sob:pag113}
	The space $W^{k,p}(\Omega)$, $k\in \mathbb{N}_0$, $1\leq p< \infty$,
	equipped with the norm \eqref{Sob:1.4} is a	separable space.
\end{prop}
\textbf{Proof.} The proof is similar to the one of Proposition \ref{Contin:10-1}.  Let us consider the case $k=1$. Let 

$$\Phi: W^{1,p}(\Omega)\rightarrow L^{p}(\Omega)\times
L^{p}(\Omega;\mathbb{R}^n),$$
$$\Phi(u)=(u,\nabla u),\quad\forall u \in W^{1,p}(\Omega).$$
$\Phi$ turns out to be an isometry, provided that we equip  $L^{p}(\Omega)\times
L^{p}(\Omega;\mathbb{R}^n)$ by the norm

$$\left(\int_{\Omega}|v_0|^pdx+\sum_{j=1}^n\int_{\Omega}\left|v_j\right|^pdx\right)^{1/p},$$
for every $v=\left(v_0,v_1,\cdots,v_n\right)\in L^{p}(\Omega)\times
L^{p}(\Omega;\mathbb{R}^n)$. Now, since $p<+\infty$, $L^{p}(\Omega)\times
L^{p}(\Omega;\mathbb{R}^n)$ is separable because it is the cartesian product of separable spaces. Hence
$\Phi\left(W^{1,p}(\Omega)\right)$ is separable as a subspace of $L^{p}(\Omega)\times
L^{p}(\Omega;\mathbb{R}^n)$ and, since $\Phi$ is an isometry,
 $W^{1,p}(\Omega)$ is separable too. $\blacksquare$

\medskip

It can be proved that if $1<p<\infty$, then $W^{1,p}(\Omega)$
is a \textit{reflexive space}. In the sequel we will not make
explicitly use this property, however for a
proof we refer to  \cite[Proposizione IX.1]{Br}.

\bigskip

\textbf{Example 1.} Let $\alpha>0$.  Let us consider
$$u(x)=\frac{1}{|x|^{\alpha}}.$$
We prove that $u\in W^{1,p}(B_1)$ if and only if $p<n$ and
$\alpha<\frac{n}{p}-1$.

We begin by assuming that $u\in W^{1,p}(B_1)$. Then $u\in L^{p}(B_1)$ and consequently $\alpha p<n$, therefore
$p<\infty$. Moreover, for any $j=1,\cdots,n$ there exists $v_j\in
L^{p}(B_1)$ such that

$$\int_{B_1}\frac{1}{|x|^{\alpha}}\partial_j \phi dx=-\int_{B_1}v_j \phi
dx,\quad \forall \phi\in C^{\infty}_0(B_1).$$ In particular
we have, for any $\phi\in C^{\infty}_0(B_1\setminus \{0\})$,

\begin{equation*}
	\begin{aligned}
		-\int_{B_1} v_j \phi
		dx=\int_{B_1}\frac{1}{|x|^{\alpha}}\partial_j \phi dx=
		\int_{B_1} \frac{\alpha x_j}{|x|^{\alpha+2}}\phi dx.
	\end{aligned}
\end{equation*}
Hence, for any $j=1,\cdots, n$,

$$v_j(x)=-\frac{\alpha x_j}{|x|^{\alpha+2}}, \quad\mbox{ a.e. in }
B_1.$$ Now, $v_j\in L^{p}(B_1)$, therefore

\begin{equation}\label{Sob:1.10}
	\sum_{j=1}^n\int_{B_1} \left|\frac{\alpha
		x_j}{|x|^{\alpha+2}}\right|^pdx<\infty.
\end{equation}
Let us observe that if $a=\left(a_1, \cdots, a_n\right)\in
\mathbb{R}^n$, then 

\begin{equation}\label{Sob:2.10}
	\frac{1}{n^{p-1}}\left|a\right|^p\leq\sum_{j=1}^n \left|a_j\right|^p\leq n
	\left|a\right|^p,
\end{equation}
(the first inequality is just a consequence of H\"{o}lder inequality, yje second is trivial). Hence \eqref{Sob:1.10} is satisfied if and only if

\begin{equation*}
	\int_{B_1} \frac{dx}{|x|^{(\alpha+1)p}}=\int_{B_1}
	\left|\frac{ x}{|x|^{\alpha+2}}\right|^pdx<\infty,
\end{equation*}
from which we derive

\begin{equation}\label{Sob:0.11}
	\alpha<\frac{n}{p}-1.
\end{equation}

Conversely, let us assume that \eqref{Sob:0.11} holds true and that $p<n$. Let us show that $u\in W^{1,p}(B_1)$.  Inequality \eqref{Sob:0.11} implies
$\alpha<\frac{n}{p}$ that, in turn implies $u\in L^{p}(B_1)$. Now, let
$\phi\in C^{\infty}_0(B_1)$ be arbitrary. We have, for any $j=1,\cdots,
n$,

\begin{equation}\label{Sob:1.11}
	\begin{aligned}
		\int_{B_1}u\partial_j\phi dx&=\lim_{\varepsilon\rightarrow
			0}\int_{B_1\setminus B_{\varepsilon} }u\partial_j\phi dx=\\&=
		\lim_{\varepsilon\rightarrow 0}\int_{B_1\setminus
			B_{\varepsilon} }\left(\partial_j(u\phi)-
		\phi\partial_ju\right)dx=\\&= \lim_{\varepsilon\rightarrow
			0}\left\{\int_{\partial B_{\varepsilon}
		}u\phi\nu_jdS-\int_{B_1\setminus B_{\varepsilon}
		}\phi\partial_j udx\right\}.
	\end{aligned}
\end{equation}
On the other hand we have
\begin{equation*}
	\left |\int_{\partial B_{\varepsilon}
	}u\phi\nu_jdS \right|\leq
	\omega_n\left\Vert\phi\right\Vert_{L^{\infty}(B_1)}
	\varepsilon^{-\alpha+n-1},
\end{equation*}
where $\omega_n$ is the measure of $\left|\partial B_{1}\right|$.
Now, by \eqref{Sob:0.11} and $p\geq 1$ we have $\alpha<n-1$. Hence
\begin{equation*}
	\lim_{\varepsilon\rightarrow 0}\left |\int_{\partial
		B_{\varepsilon} }u\phi\nu_jdS\right|=0;
\end{equation*}
coming back to \eqref{Sob:1.11} and keeping in mind that (by
\eqref{Sob:0.11}) 
$$\frac{\alpha x_j}{|x|^{\alpha+2}}\in L^p(B_1)\subset
L^1(B_1),$$ we get

\begin{equation*}
	\begin{aligned}
		\int_{B_1}u\partial_j\phi dx=\lim_{\varepsilon\rightarrow
			0}\int_{B_1\setminus B_{\varepsilon} }\frac{\alpha
			x_j}{|x|^{\alpha+2}}\phi dx=\int_{B_1}\frac{\alpha
			x_j}{|x|^{\alpha+2}}\phi dx.
	\end{aligned}
\end{equation*}
All in all we have

$$\partial_j u=-\frac{\alpha x_j}{|x|^{\alpha+2}}\in L^p(B_1),$$
therefore $u\in W^{1,p}(B_1)$. $\spadesuit$

\bigskip

\textbf{Remark.} Similarly to
in Example 1, it can be proved that if

$$u\in C^0\left(\overline{B_1}\right)\cap C^1\left(\overline{B_1}\setminus
\{0\}\right)$$ then we have
$$u\in W^{1,p}(B_1)\Longleftrightarrow \nabla u \in L^p(B_1),$$
here  $\nabla u$ is the gradient of $u$ in $B_1\setminus \{0\}$
in the classic sense,
Let us consider the case $n=1$ only, because the case $n>1$ can be treated
in precisely the same way as Example 1 (the proof is left to the reader). Let $\phi\in
C_0^{\infty}(-1,1)$, we have

\begin{equation*}
	\begin{aligned}
		\int_{-1}^1u\phi' dx&=\lim_{\varepsilon\rightarrow
			0}\int_{(-1,1)\setminus [-\varepsilon,\varepsilon] }u\phi'
		dx=\\&=
		\lim_{\varepsilon\rightarrow
			0}\left\{-u(\varepsilon)\phi(\varepsilon)+u(-\varepsilon)\phi(-\varepsilon)-\int_{(-1,1)\setminus
			[-\varepsilon,\varepsilon] }u'\phi dx\right\},
	\end{aligned}
\end{equation*}
but $u\in C^0([-1,1])$, hence
$$\lim_{\varepsilon\rightarrow
	0}\left(-u(\varepsilon)\phi(\varepsilon)+u(-\varepsilon)\phi(-\varepsilon)\right)=0$$
and $u'\in L^p(-1,1)$, implies
$$\int_{-1}^1u\phi' dx=-\int_{-1}^1u'\phi dx.$$
Therefore, $u'$ is the weak derivative of $u$ and $u\in W^{1,p}(-1,1)$.
$\blacklozenge$

\subsection{The spaces $W_0^{k,p}(\Omega)$ }
\label{Sob:kp0} We give the following
\begin{definition}\label{Sob:49-0}
	Let $1\leq p\leq \infty$, $k\in \mathbb{N}_0$ and $\Omega$ be an open set of $\mathbb{R}^n$, $n\geq 1$. Let us denote by \index{$W_0^{k,p}(\Omega)$}
	$$W_0^{k,p}(\Omega),$$
	the closure of $C_0^{\infty}(\Omega)$ in $W^{k,p}(\Omega)$.
	We write \index{$H_0^{k}(\Omega)$}
	$$H_0^{k}(\Omega)=W_0^{k,2}(\Omega).$$
\end{definition}

\bigskip

Let us notice that $W_0^{0,p}(\Omega)=L^{p}(\Omega)$. Moreover
$W_0^{k,p}(\Omega)$, equipped with the norm
$\left\Vert\cdot\right\Vert_{W^{k,p}(\Omega)}$, as it is a closed subspace of
 $W^{k,p}(\Omega)$, is a Banach space.

\section{Approximation and density theorems} \label{teore-dens}
Let $\eta$ be a mollifier, namely $\eta\in
C^{\infty}_0\left(\mathbb{R}^n\right)$ satisfies (i) supp
$\eta\subset B_1$, \\ (ii) $\eta\geq 0$, (iii) $\int_{\mathbb{R}^n}
\eta(x)dx=1$. Set, for any $\varepsilon>0$,

$$\eta_{\varepsilon}(x)=\varepsilon^{-n}\eta\left(\varepsilon^{-1}x\right).$$
Let $\Omega$ be an open set of $\mathbb{R}^n$. Set

$$\Omega_{\varepsilon}=\left\{x\in \Omega:\mbox{ } \mbox{dist}(x,\partial\Omega)>\varepsilon\right\}.$$

\medskip

The following Theorem holds true

\begin{theo}[\textbf{local approximation by $C^{\infty}$ functions}]\label{Sob:teo1.3}
	 Let $k\in \mathbb{N}_0$,
	$p\in[1,+\infty)$. Let us assume that $u\in W^{k,p}\left(\Omega\right) $. Let us denote by
	$$u^{\varepsilon}=\eta_{\varepsilon}\star u,\quad\mbox{ in }
	\Omega_{\varepsilon}.$$ Then we have $$u^{\varepsilon}\in
	C^{\infty}\left(\Omega_{\varepsilon}\right)\cap
	W^{k,p}_{loc}\left(\Omega_{\varepsilon}\right),\quad
	\forall\varepsilon>0$$ and
	$$u^{\varepsilon}\rightarrow u,\quad\mbox{ as } \varepsilon\rightarrow 0,\mbox{ in } W^{k,p}_{loc}\left(\Omega\right).$$
\end{theo}
\textbf{Proof.} The fact that $u^{\varepsilon}\in
C^{\infty}\left(\Omega_{\varepsilon}\right)$ is an immediate consequence of Theorem \ref{Sob:teo28R}. Concernig
$u^{\varepsilon}\in W^{k,p}_{loc}\left(\Omega_{\varepsilon}\right)$,
we have by Theorem \ref{Sob:teo26R}

$$\partial^{\alpha}u^{\varepsilon}=\left(\partial^{\alpha}\eta_{\varepsilon}\right)\star
u,\quad\forall \alpha\in \mathbb{N}_0^n.$$ Now, if $|\alpha|\leq k$,
we have, for any $x\in \Omega_{\varepsilon}$,

\begin{equation*}
	\begin{aligned}
		\left(\partial^{\alpha}\eta_{\varepsilon}\right)\star
		u&=\int_{\Omega}\left(\partial_x^{\alpha}\eta_{\varepsilon}\right)(x-y)u(y)dy=\\&=
		(-1)^{|\alpha|}\int_{\Omega}\partial_y^{\alpha}\left(\eta_{\varepsilon}(x-y)\right)u(y)dy=\\&=
		\int_{\Omega}\eta_{\varepsilon}(x-y)\partial_y^{\alpha}u(y)dy=\\&=
		\left(\eta_{\varepsilon}\star \partial^{\alpha}u\right)(x).
	\end{aligned}
\end{equation*}

Let $\omega \Subset \Omega$,  since
$\partial^{\alpha}u\in L^p(\omega)$ we have, for any $|\alpha|\leq k$,

$$\eta_{\varepsilon}\star
\partial^{\alpha}u\rightarrow\partial^{\alpha}u,\quad\mbox{ as }
\varepsilon\rightarrow 0,\mbox{ in } L^p(\omega).$$ Therefore

$$u\rightarrow u^{\varepsilon},\quad\mbox{ as }
\varepsilon\rightarrow 0,\mbox{ in } W^{k,p}(\omega).$$
$\blacksquare$

\bigskip

\begin{theo}[\textbf{Meyers -- Serrin}]\label{densit 1}
	\index{Theorem:@{Theorem:}!- Meyers -- Serrin@{- Meyers -- Serrin}}
	 Let $\Omega$ be an bounded open set of $\mathbb{R}^n$. Let $k\in \mathbb{N}_0$,
	$p\in[1,+\infty)$. If $u\in W^{k,p}\left(\Omega\right)$ then
	there exists a sequence $\left\{u_m\right\}$  in $W^{k,p}\left(\Omega\right)\cap C^{\infty}\left(\Omega\right)$ which satisfies
	$$u_{m}\rightarrow u,\quad\mbox{ as } m\rightarrow \infty,\mbox{ in } W^{k,p}\left(\Omega\right).$$
\end{theo}
\textbf{Proof.} Let

$$\Omega_{j}=\left\{x\in \Omega:\mbox{ }
\mbox{dist}(x,\partial\Omega)>\frac{1}{j}\right\}, \quad j\in \mathbb{N},$$
It is not restrictive to assume $\Omega_{j}\neq \emptyset$, for every $j\in \mathbb{N}$. We have

\begin{equation}\label{Sob:0.17}
	\overline{\Omega}_{j}\subset \Omega_{j+1},\quad \forall j\in
	\mathbb{N},\quad \bigcup_{j=1}^{\infty}\Omega_{j}=\Omega.
\end{equation} Let
$\phi_j\in C^{\infty}\left(\mathbb{R}^n\right)$, $j\in \mathbb{N}$,
satisfy supp $\phi_j\subset \Omega_{j+1}$; $\phi_j(x)=1$ for every
$x\in \Omega_{j}$; $0\leq \phi_j(x)\leq 1$, for every $x\in
\mathbb{R}^n$. 

For any $j\in \mathbb{N}$ we get
$$x\in \mathbb{R}^n\setminus \Omega_{j+1} \mbox{ }\Longrightarrow\mbox{
} \phi_j(x)=0\leq \phi_{j+1}(x)$$ and

$$x\in \Omega_{j+1} \mbox{ }\Longrightarrow\mbox{
} \phi_j(x)\leq 1=\phi_{j+1}(x).$$ Hence

$$\phi_j\leq \phi_{j+1}, \quad \forall j\in \mathbb{N},\quad \mbox{ in } \mathbb{R}^n.$$
Set
$$\zeta_0=\phi_2,\quad \zeta_{j}=\phi_{j+1}-\phi_{j},\quad \forall j\in
\mathbb{N}$$ and
$$V_0=\Omega_2,\quad V_{j}=\Omega_{j+3}\setminus \overline{\Omega}_{j},\quad \forall j\in
\mathbb{N}.$$ We have

$$\zeta_{j}\in C^{\infty}_0\left(\mathbb{R}^n\right),\quad\mbox{supp }\zeta_{j}\subset V_j, \quad \forall j\in
\mathbb{N}_0.$$ Moreover
$$0\leq \zeta_{j}\leq 1,\quad\forall x\in \mathbb{R}^n,\quad \forall j\in
\mathbb{N}_0$$ and

\begin{equation}\label{Sob:for1.17}
	\sum_{j=0}^{\infty}\zeta_{j}(x)=1, \quad\forall x\in \Omega.
\end{equation}
Let us check \eqref{Sob:for1.17}. Let $x\in \Omega$, by
\eqref{Sob:0.17} we have that there exists $\overline{m}\in \mathbb{N}$ such that $x\in \Omega_j$, for every $j\geq \overline{m}$. Let $m\geq
\overline{m}$, we have $x\in \Omega_m$, hence $\phi_{m+1}(x)=1$.
Consequently, we have

\begin{equation*}
	\begin{aligned}
		\sum_{j=0}^m
		\zeta_{j}(x)&=\zeta_{0}(x)+\zeta_{1}(x)+\cdots+\zeta_{m}(x)=\\&=
		\phi_{2}(x)+(\phi_{3}(x)-\phi_{2}(x))+\cdots+(\phi_{m+1}(x)-\phi_{m}(x))=\\&
		=\phi_{m+1}(x)=1.
	\end{aligned}
\end{equation*}
Therefore we have checked \eqref{Sob:for1.17}.

Now, let $u\in W^{k,p}(\Omega)$ and let us consider the functions
$\zeta_ju$, $j\in \mathbb{N}_0$. Proposition gives
\ref{Sob:prop1.2.6} we get  $\zeta_ju\in W^{k,p}(\Omega)$, in addition

$$\zeta_ju=0,\quad\mbox{ in } \Omega\setminus \overline{V}_j.$$
Let us denote by $W_0=\Omega_4$, $W_1=\Omega_5$, $W_j=\Omega_{j+4}\setminus\overline{\Omega}_{j-1}$, $j\geq 2$. Let $\delta>0$ be fixed and let $0<\varepsilon_j<\frac{1}{j+4}-\frac{1}{j+3}$ satisfy

$$u^j=\eta_{\varepsilon_j}\star (\zeta_ju)\in C^{\infty}(\Omega)\cap W^{k,p}(\Omega),$$ we have

$$u^j=0,\quad\mbox{ in } \Omega\setminus \overline{W}_j.$$ Theorem \ref{Sob:teo1.3} implies that for every $j\in \mathbb{N}_0$ there exists $\varepsilon_j>0$ such that

\begin{equation}\label{Sob:for1.19}
	\left\Vert u^j-\zeta_ju\right\Vert_{W^{k,p}(\Omega)}=\left\Vert
	u^j-\zeta_ju\right\Vert_{W^{k,p}(W_j)}\leq \frac{\delta}{2^{j+1}},
	\quad j\in \mathbb{N}_0.
\end{equation}
We now set

\begin{equation}\label{Sob:for2.19}
	v(x)=\sum_{j=0}^{\infty}u^j(x).
\end{equation}
Notice that, for any $x\in \Omega$, only a finite number of
terms of series \eqref{Sob:for2.19} is different from $0$. Moreover, as
 $u^j \in C^{\infty}(\Omega)\cap W^{k,p}(\Omega)$, for every $j\in \mathbb{N}_0$, we have $v\in C^{\infty}(\Omega)\cap
W^{k,p}(\Omega)$.

Now, taking into account that
$$u=\sum_{j=0}^{\infty}\zeta_ju,$$
for any $h\in \mathbb{N}$,  
\eqref{Sob:for1.19} and \eqref{Sob:for2.19} give

\begin{equation*}
	\begin{aligned}
		\left\Vert v-u\right\Vert_{W^{k,p}(\Omega_h)}&=\left\Vert
		\sum_{j=0}^{\infty}u^j-\sum_{j=0}^{\infty}\zeta_ju\right\Vert_{W^{k,p}(\Omega_h)}\leq\\&\leq
		\sum_{j=0}^{\infty}\left\Vert
		u^j-\zeta_ju\right\Vert_{W^{k,p}(\Omega_h)}\leq\\&\leq
		\sum_{j=0}^{\infty}\left\Vert
		u^j-\zeta_ju\right\Vert_{W^{k,p}(\Omega)}\leq\\&\leq 
		\sum_{j=0}^{\infty}\frac{\delta}{2^{j+1}}=\delta.
	\end{aligned}
\end{equation*}

All in all, we have

$$\left\Vert v-u\right\Vert_{W^{k,p}(\Omega_h)}\leq
\delta,\quad\forall h\in \mathbb{N}.$$ Hence

$$\left\Vert v-u\right\Vert_{W^{k,p}(\Omega)}=\lim_{h\rightarrow
	\infty} \left\Vert v-u\right\Vert_{W^{k,p}(\Omega_h)}\leq \delta.$$
Therefore, the sequence $$u_m=\sum_{j=0}^{m}u^j(x),\quad
m\in\mathbb{N},$$ satisfies the thesis. $\blacksquare$

\bigskip

\underline{\textbf{Exercise.}} Prove Theorem  \ref{densit 1} without the assumption that $\Omega$ is bounded. [Hint: consider $\Omega_j\cap B_j(x_0)$, $x_0$ fixed point, instead of $\Omega_j$.] $\clubsuit$

\bigskip

The following Theorem holds true

\begin{theo}[\textbf{$C^{\infty}$ approximation to the boundary}]\label{densit 2}
	\index{Theorem:@{Theorem:}!- approximation $C^{\infty}\left(\overline{\Omega}\right)$ of $W^{k,p}(\Omega)$ function@{- approximation $C^{\infty}\left(\overline{\Omega}\right)$ of $W^{k,p}(\Omega)$ function}}
	Let $\Omega$ be a bounded open set of $\mathbb{R}^{n}$ whose boundary is of
	class $C^{0,1}$ with constants $r_0,M_0$. Let $u\in W^{k,p}(\Omega)$,
	$1\leq p<+\infty$. Then there exists a sequence of functions
	$\left\{u_j\right\}\subset C^{\infty}\left(\overline{\Omega}\right)$
	such that
	
	$$u_j\rightarrow u,\quad\mbox{ as }
	j\rightarrow\infty,\mbox{ in }W^{k,p}(\Omega).$$
\end{theo}

\medskip

To preparare the proof of Theorem \ref{densit 2}, we introduce some
notations and we prove a Proposition,

Let $\Omega$ be a bounded open set of $\mathbb{R}^{n}$ whose boundary is of
class $C^{0,1}$ with constants $r_0,M_0$. Let $x_0\in \partial\Omega$. We may assume (up to isometry) that
$x_0=0$ and
\begin{equation*}
	\Omega\cap Q_{r_0,2M_0}=\left\{x\in Q_{r_0,2M_0}:\mbox{ }
	x_n>\varphi(x') \right\},
\end{equation*}
where $\varphi\in C^{0,1}\left(B'_{r_0}\right)$ satisfies
$$\varphi(0)=0$$
and
$$\left\Vert\varphi\right\Vert_{C^{0,1}\left(\overline{B'_{r_0}}\right)}\leq
M_0r_0.$$
Set $$V=\Omega\cap Q_{\frac{r_0}{2},\frac{M_0}{2}}.$$ 

\begin{figure}\label{figura-p23sob}
	\centering
	\includegraphics[trim={0 0 0 0},clip, width=15cm]{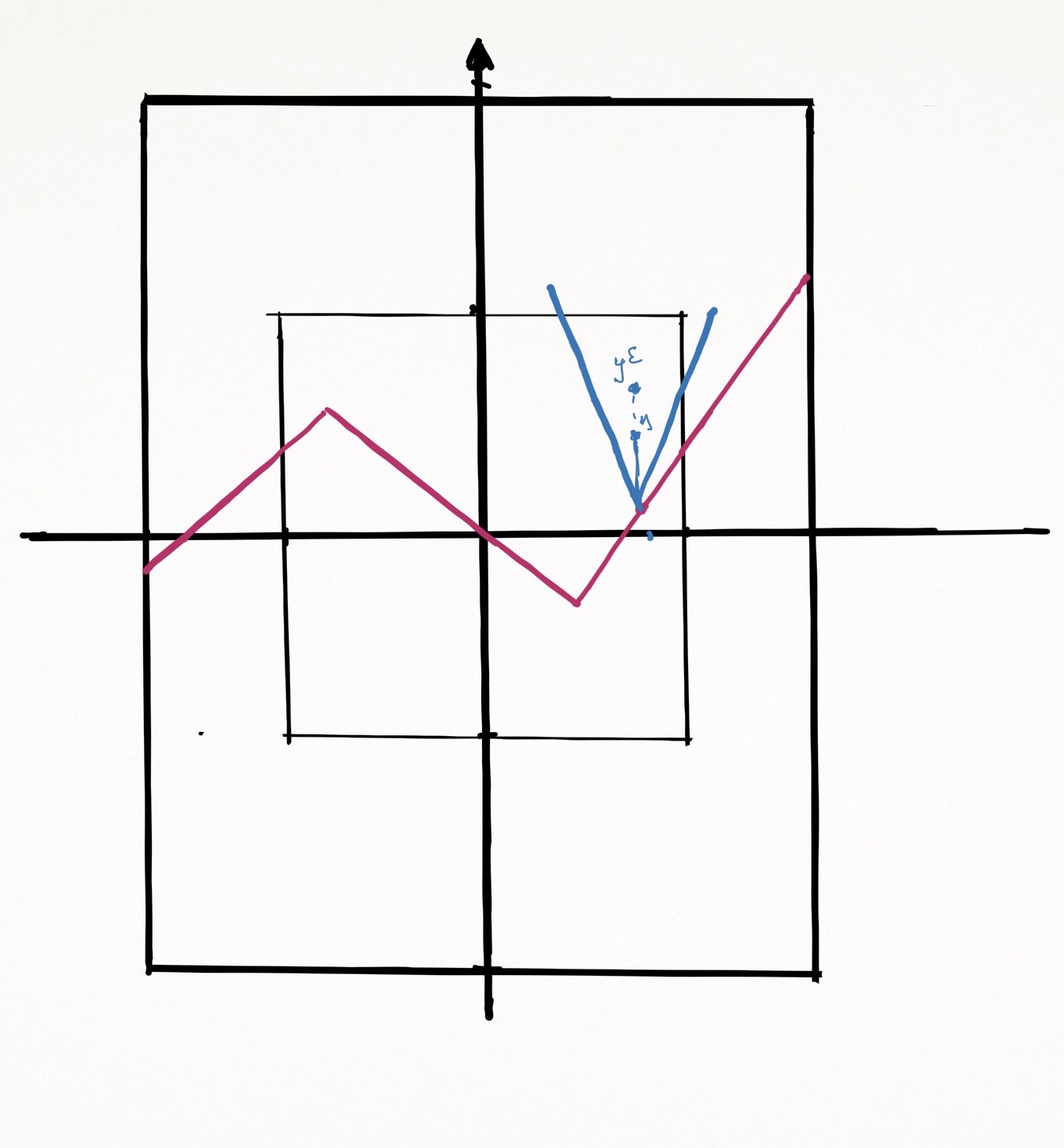}
	\caption{}
\end{figure}
Let $y$ be any point of $V$, we look for what $\lambda>0$ and $\varepsilon>0$ we have (see Figure 3.2)
\begin{equation}\label{Sob:inclusione24}
	B_{\varepsilon}\left(y^{\varepsilon}\right)\subset \Omega\cap
	Q_{r_0,M_0},
\end{equation}
where
$$y^{\varepsilon}=y+\varepsilon\lambda e_n,$$

\smallskip

\noindent \textbf{1.} Let us check that if $\varepsilon$ and $\lambda$
satisfy

\begin{equation}\label{Sob:stella24}
	\varepsilon<\frac{r_0}{2},\quad \varepsilon
	(1+\lambda)<\frac{M_0r_0}{2},
\end{equation}
then we have

\begin{equation}\label{Sob:stella24bis}
	B_{\varepsilon}\left(y^{\varepsilon}\right)\subset Q_{r_0,M_0}.
\end{equation}
Since $B_{\varepsilon}\left(y^{\varepsilon}\right)\subset
B'_{\varepsilon}\left((y^{\varepsilon})'\right)\times\left[y^{\varepsilon}_n-\varepsilon,y^{\varepsilon}_n+\varepsilon\right]$,
we have that the first condition of \eqref{Sob:stella24} implies
\begin{equation}\label{Sob:stella24ter}
	B'_{\varepsilon}\left((y^{\varepsilon})'\right)\subset
	B'_{r_0}\end{equation} and second condition of 
\eqref{Sob:stella24} implies
$$y^{\varepsilon}_n+\varepsilon\leq
\frac{M_0r_0}{2}+\varepsilon\lambda+\varepsilon<M_0r_0$$ and,
similarly,

$$y^{\varepsilon}_n-\varepsilon\geq
-\frac{M_0r_0}{2}+\varepsilon\lambda-\varepsilon>-M_0r_0.$$ Hence
$$\left[y^{\varepsilon}_n-\varepsilon,y^{\varepsilon}_n+\varepsilon\right]\subset \left[-M_0r_0,
M_0r_0\right],$$ which gives \eqref{Sob:stella24bis}.

\smallskip

\noindent \textbf{2.} In order that
$B_{\varepsilon}\left(y^{\varepsilon}\right)\subset \Omega\cap
Q_{r_0,M_0}$, it is suffices that, besides conditions \eqref{Sob:stella24},
$y^{\varepsilon}$ have a distance greater or equal to $\varepsilon$ from the cone
$$x_n=M_0\left|x'-y'\right|+\varphi(y').$$
Now, denoting by $d_{\varepsilon}$ this distance, we have

\begin{equation*}
	\begin{aligned}
		d_{\varepsilon}&=\frac{\left|y_n+\varepsilon\lambda-\varphi(y')\right|}{\sqrt{1+M_0^2}}=\\&=
		\frac{y_n+\varepsilon\lambda-\varphi(y')}{\sqrt{1+M_0^2}}\geq\\&\geq\frac{\varepsilon\lambda}{\sqrt{1+M_0^2}}.
	\end{aligned}
\end{equation*}
Hence, in order that $d_{\varepsilon}>\varepsilon$ it suffices that 
$\lambda>\sqrt{1+M_0^2}$. Therefore, by choosing
$$\lambda=\lambda_0:=2\sqrt{1+M_0^2}$$ and by requiring that
$$\varepsilon<\varepsilon_0:=\min\left\{\frac{r_0}{2},\frac{M_0r_0}{2},
\frac{M_0r_0}{2\sqrt{1+M_0^2}}\right\}$$ we obtain
\eqref{Sob:inclusione24}.

\bigskip

For any $u\in L^p(\Omega)$ and $\varepsilon<\varepsilon_0$ we set

\begin{equation}\label{Sob:for2.25}
	u_{\varepsilon}(x)=u\left(x^{\varepsilon}\right)=u(x+\lambda_0\varepsilon
	e_n),\quad \forall x\in V
\end{equation}
and

\begin{equation}\label{Sob:for3.25}
	v^{\varepsilon}(x)=\int_{B_{\varepsilon}\left(x^{\varepsilon}\right)\cap
		\Omega}\eta_{\varepsilon}(x+\lambda_0\varepsilon e_n-y)u(y)dy,\quad
	\forall x\in V,
\end{equation}
where, we recall,
$\eta_{\varepsilon}(x)=\varepsilon^{-n}\eta\left(\varepsilon^{-1}x\right)$.
Now, since
$B_{\varepsilon}\left(x^{\varepsilon}\right)\subset \Omega\cap
Q_{r_0,M_0}$, we have

\smallskip

\begin{equation}\label{Sob:for4.25}
	\begin{aligned}
		v^{\varepsilon}(x)&=\int_{B_{\varepsilon}\left(x^{\varepsilon}\right)}\eta_{\varepsilon}(x+\lambda_0\varepsilon
		e_n-y)u(y)dy=\\&\\&=
		\int_{B_{\varepsilon}}\eta_{\varepsilon}(y)u(x+\lambda_0\varepsilon
		e_n-y)dy ,\quad \forall x\in V.
	\end{aligned}
\end{equation}

\smallskip

\noindent The first equality in \eqref{Sob:for4.25} gives
$$v^{\varepsilon}(x)=\int_{\Omega}\eta_{\varepsilon}(x+\lambda_0\varepsilon
e_n-y)u(y)dy.$$ Hence

\begin{equation}\label{Sob:for5.25}
	\partial^{\alpha}v^{\varepsilon}(x)=\int_{\Omega}\partial_x^{\alpha}\eta_{\varepsilon}(x+\lambda_0\varepsilon e_n-y)u(y)dy,\quad
	\forall x\in V
\end{equation}
so that

$$\partial^{\alpha}v^{\varepsilon}\in C^{\infty}\left(\overline{V}\right).$$
Moreover, for any $u\in W^{k,p}(\Omega)$, we have
\begin{equation}\label{Sob:for6.25}
	\begin{aligned}
		\partial^{\alpha}v^{\varepsilon}(x)&=\int_{\Omega}
		\eta_{\varepsilon}(x+\lambda_0\varepsilon
		e_n-y)\partial^{\alpha}u(y)dy=\\&\\&=
		\int_{B_{\varepsilon}\left(x^{\varepsilon}\right)}\eta_{\varepsilon}(x+\lambda_0\varepsilon
		e_n-y)\partial^{\alpha}u(y)dy ,\quad \forall x\in V,
	\end{aligned}
\end{equation}
for  every $|\alpha|\leq k$.

\medskip

We have the following

\begin{prop}\label{Sob:prop4.3}
	If $u\in W^{k,p}(\Omega)$ and $p\in [1,+\infty)$ then
	$$v^{\varepsilon}\rightarrow u,\quad\mbox{ as }
	\varepsilon\rightarrow 0,\mbox{ in } W^{k,p}(V).$$
\end{prop}
\textbf{Proof.} First of all we prove
$$v^{\varepsilon}\rightarrow u,\quad\mbox{ as }
\varepsilon\rightarrow 0,\mbox{ in } L^{p}(V).$$ The triangle inequality gives 

\begin{equation}\label{Sob:for0.26}\left\Vert v^{\varepsilon}-u\right\Vert_{L^{p}(V)}\leq \left\Vert
	u_{\varepsilon}-u\right\Vert_{L^{p}(V)}+\left\Vert
	v^{\varepsilon}-u_{\varepsilon}\right\Vert_{L^{p}(V)}.\end{equation}
Now

$$\left\Vert
u_{\varepsilon}-u\right\Vert^p_{L^{p}(V)}=\int_V\left|u(x+\lambda_0\varepsilon
e_n)-u(x)\right|^pdx$$ and, by Theorem \ref{Sob:teo24R}, we have 

\begin{equation}\label{Sob:for00.26}\lim_{\varepsilon \rightarrow0}\left\Vert
	u_{\varepsilon}-u\right\Vert_{L^{p}(V)}=0.\end{equation} Moreover,
by the second equality in \eqref{Sob:for4.25} we have, for any $x\in
V$,

$$v^{\varepsilon}(x)-u_{\varepsilon}(x)=\int_{B_{\varepsilon}}\eta_{\varepsilon}(y)\left(u(x+\lambda_0\varepsilon
e_n-y)-u(x+\lambda_0\varepsilon e_n)\right)dy.$$ 
In order to prove that
the second term on the right hand side in \eqref{Sob:for0.26} goes to $0$ it suffices to
repeat the same steps which provide the proof of Theorem
\ref{Sob:teo25R}. For completeness, let us repeat these
steps.

\begin{equation*}
	\begin{aligned}
		&\int_{V}\left|v^{\varepsilon}-u_{\varepsilon}\right|^pdx\leq
		\int_{V}\left(\int_{B_{\varepsilon}}\eta_{\varepsilon}(y)\left|u(x+\lambda_0\varepsilon
		e_n-y)-u(x)\right|dy\right)^pdx=\\&=
		\int_{V}\left(\int_{B_{\varepsilon}}\eta^{1/p'}_{\varepsilon}(y)\eta^{1/p}_{\varepsilon}(y)\left|u(x+\lambda_0\varepsilon
		e_n-y)-u(x)\right|dy\right)^pdx\leq\\&\leq \int_{V}dx
		\left(\int_{B_{\varepsilon}}\eta_{\varepsilon}(y)dy\right)^{p/p'}\int_{B_{\varepsilon}}\eta_{\varepsilon}(y)\left|u(x+\lambda_0\varepsilon
		e_n-y)-u(x)\right|^pdy=\\&=
		\int_{V}dx\int_{B_{\varepsilon}}\eta_{\varepsilon}(y)\left|u(x+\lambda_0\varepsilon
		e_n-y)-u(x)\right|^p=\\&=
		\int_{B_{\varepsilon}}\left(\eta_{\varepsilon}(y)\int_{V}\left|u(x+\lambda_0\varepsilon
		e_n-y)-u(x)\right|^pdx\right)dy\leq\\&\leq \sup_{|y|\leq
			\varepsilon}\int_{V}\left|u(x+\lambda_0\varepsilon
		e_n-y)-u(x)\right|^pdx.
	\end{aligned}
\end{equation*}
All in all, we have
$$\int_{V}\left|v^{\varepsilon}-u_{\varepsilon}\right|^pdx\leq \sup_{|y|\leq
	\varepsilon}\int_{V}\left|u(x+\lambda_0\varepsilon
e_n-y)-u(x)\right|^pdx.$$  Theorem \ref{Sob:teo24R} now yields

$$\lim_{\varepsilon\rightarrow 0}\left\Vert
v^{\varepsilon}-u_{\varepsilon}\right\Vert_{L^{p}(V)}=0.$$ By the latter, by \eqref{Sob:for0.26} and by \eqref{Sob:for00.26} we have
$$\lim_{\varepsilon\rightarrow 0}\left\Vert
v^{\varepsilon}-u\right\Vert_{L^{p}(V)}=0.$$

If $u\in W^{k,p}(\Omega)$, we obtain from what has been proven above and from
\eqref{Sob:for6.25}
$$\partial^{\alpha}v^{\varepsilon}\rightarrow \partial^{\alpha} u,\quad\mbox{ per }
\varepsilon\rightarrow 0,\mbox{ in } L^{p}(V),\mbox{ for
}|\alpha|\leq k.$$ Hence

$$v^{\varepsilon}\rightarrow u,\quad\mbox{ as }
\varepsilon\rightarrow 0,\mbox{ in } W^{k,p}(V).$$ $\blacksquare$

\bigskip

\textbf{Proof of Theorem \ref{densit 2}.} Let
$x_0\in\partial \Omega$, let us denote by 
$\widetilde{Q}_{r_0,2M_0}(x_0)$ the  cylinder isometric to
$Q_{r_0,2M_0}$ such that  \eqref{Sob:for1.20} holds. As a consequence,
$\left\{\widetilde{Q}_{\frac{r_0}{2},\frac{M_0}{2}}(x_0)\right\}_{x_0\in
	\partial\Omega}$ is an open covering of the compact set
$\partial\Omega$. Let
$$\left\{\widetilde{Q}_{\frac{r_0}{2},\frac{M_0}{2}}(x_i)\right\}_{1\leq i\leq N}$$ be a finite subcovering
of $\partial\Omega$. For any $1\leq i\leq N$ and
let us denote $$V_i=\Omega\cap \widetilde{Q}_{\frac{r_0}{2},\frac{M_0}{2}}(x_i).$$ For any fixed $\delta>0$ let
$v_i\in C^{\infty}\left(\overline{V_i}\right)$ be the function  constructed
in \eqref{Sob:for3.25} which satisfies

\begin{equation}\label{Sob:for0.28}
	\left\Vert v_i-u\right\Vert_{W^{k,p}\left(V_i\right)}\leq\delta.\end{equation}
Moreover, let  $V_0\subset \Omega$ be such that 
$$\Omega\subset \bigcup_{i=0}^N V_i$$ and let $\left\{\zeta_i\right\}_{0\leq i\leq
	N}$ be a partition of unity  (compare Theorem \ref{Sob:lem3.3}) which satisfies
 $\zeta_{i}\in C^{\infty}_0\left(\mathbb{R}^n\right)$,
$\mbox{supp }\zeta_{i}\subset V_i$ per $1\leq i\leq N$ and
\begin{equation*}
	\sum_{i=0}^{N}\zeta_{i}(x)=1, \quad\forall x\in \Omega.
\end{equation*}
Let us denote by $v_0=\zeta_0 u$ and 
$$v=\sum_{i=0}^{N}\zeta_{i}v_i.$$
We have $v\in C^{\infty}(\overline{\Omega})$ and, taking into account  \eqref{Sob:for0.28}, 
\begin{equation*}
	\begin{aligned}
		\left\Vert
		\partial^{\alpha}v-\partial^{\alpha}u\right\Vert_{L^{p}(\Omega)}&=\left\Vert
		\sum_{i=0}^{N}\partial^{\alpha}(\zeta_{i}v_i)-\sum_{i=0}^{N}\partial^{\alpha}(\zeta_iu)\right\Vert_{L^{p}(\Omega)}\leq\\&\leq
		\sum_{i=0}^{N}\left\Vert
		\partial^{\alpha}(\zeta_{i}v_i)-\partial^{\alpha}(\zeta_iu)\right\Vert_{L^{p}(V_i)}\leq\\&\leq
		C\sum_{i=0}^{N}\left\Vert v_i-u\right\Vert_{W^{k,p}(V_i)}\leq\\&\leq
		CN\delta,
	\end{aligned}
\end{equation*}
for every $|\alpha|\leq k$.

Therefore, if $u\in W^{k,p}(\Omega)$ then for every $\eta>0$ there exists
$v\in C^{\infty}\left(\overline{\Omega}\right)$ such that
$$\left\Vert v-u\right\Vert_{W^{k,p}(\Omega)}<\eta.$$
The Theorem is proved. $\blacksquare$

\bigskip

We conclude this Section with some propositions and exercises.

\bigskip

\begin{prop}\label{Sob:Es5.39}
	Let $\Omega$ be a connected open set of $\mathbb{R}^n$ and let $u\in
	W_{loc}^{1,1}(\Omega)$ satisfy
	$$\nabla u=0,\quad \mbox{in } \Omega,$$ then $u$ is almost everywhere equal to a constant.
\end{prop}
\textbf{Proof.} Let us first consider the case in which
$\Omega=B_r$, $r>0$, and $u\in W^{1,1}(B_r)$. Let
$\delta$ be any number in $(0,r)$ and  let $\varepsilon\in (0,\delta)$.
Set

$$u_{\varepsilon}(x)=\int_{B_r}\eta_{\varepsilon}(x-y)u(y)dy.$$
By Theorem \ref{Sob:teo1.3} we derive that $u_{\varepsilon}\in
C^{\infty}\left(B_{r-\delta}\right)$ and that, for any $x\in
B_{r-\delta}$,

\begin{equation*}
	\begin{aligned}
		\nabla
		u_{\varepsilon}(x)=-\int_{B_r}\nabla\left(\eta_{\varepsilon}(x-y)\right)u(y)dy=\int_{B_r}\eta_{\varepsilon}(x-y)\nabla_yu(y)dy=0.
	\end{aligned}
\end{equation*}
Hence
$$u_{\varepsilon}(x)=C_{\varepsilon},\quad\mbox{in }
B_{r-\delta},$$ where $C_{\varepsilon}$ is a constant
which depends on $\varepsilon$. On the other hand
$$u_{\varepsilon}\rightarrow u,\quad\mbox{as
}\varepsilon\rightarrow 0, \mbox{ in }
L^1\left(B_{r-\delta}\right).$$  Since the limit, in
$L^1\left(B_{r-\delta}\right)$, of a sequence of constant functions is a
constant function, we have $u=\widetilde{C}_{\delta}$, almost everywhere in
$B_{r-\delta}$, where $\widetilde{C}_{\delta}$ is a constant. Triavially, $\widetilde{C}_{\delta}$ does not depend on  $\delta$ and, as  $\delta$ is arbitrary in $(0,r)$, we have that $u$ is constant almost everywhere in $B_r$.

Now, let us consider the general case and let us assume that $u\in
W_{loc}^{1,1}(\Omega)$. Let $\overline{x}\in \Omega$ and
$B_r\left(\overline{x}\right)\Subset \Omega$, for what
proved before we have that there is $C\in \mathbb{R}$ such that

\begin{equation}\label{Sob:for0.39ter}
	u=C,\quad\mbox{ a.e. in }B_{r}\left(\overline{x}\right).
\end{equation}
Let $y$ be any point of $\Omega$, $y\neq \overline{x}$. We prove that there exists $\rho>0$ such that

\begin{equation}\label{Sob:for1.39ter}
	u=C,\quad\mbox{ a.e. in }B_{\rho}\left(y\right).
\end{equation}
Since $\Omega$ is a connected open set, there exists a continuous path
$\gamma:[0,1]\rightarrow\Omega$, $\gamma$ such that
$\gamma(0)=\overline{x}$, $\gamma(1)=y$. Since
$\gamma\left([0,1]\right)$ is a compact, we have
$$r_0:=\mbox{dist }
\left(\gamma\left([0,1]\right),\partial\Omega\right)>0.$$ Moreover,
let $\rho=\min\{r_0,r\}$, we can estract a finite subcovering by the open covering
$\left\{B_{\rho}(x)\right\}_{x\in \gamma\left([0,1]\right)}$, of $\gamma\left([0,1]\right)$. Let
$\left\{B_{\rho}(x_j)\right\}_{1\leq j\leq N}$ be such a finite subcovering of $\gamma\left([0,1]\right)$, where $x_j\in\gamma\left([0,1]\right)$. It is not restrictive to assume
$x_1=\overline{x}$, $x_N=y$. For this purpose it suffices, eventually,
to add to the family $\left\{B_{\rho}(x_j)\right\}_{1\leq j\leq N}$, the balls $B_{\rho}\left(\overline{x}\right)$ and $B_{\rho}\left(y\right)$) and, rearranging the remaining points $x_2,\cdots, x_{N-1}$, we may assume that (as $\gamma\left([0,1]\right)$ is connected)

\begin{equation}\label{Sob:for-stella.39ter}
	B_{\rho}\left(x_j\right)\cap B_{\rho}\left(x_{j+1}\right)\neq
	\emptyset,\quad j=1,\cdots, N-1.
\end{equation}
In each ball $B_{\rho}\left(x_j\right)$, $u$ is constant almost everywhere and, since $B_{\rho}\left(x_j\right)\cap
B_{\rho}\left(x_{j+1}\right)$ has positive measure, for $j=1,\cdots,
N-1$, we have by \eqref{Sob:for0.39ter} that $u=C$ almost everywhere in $B_{\rho}\left(x_j\right)$, $j=1,\cdots, N$.  Therefore we obtain
\eqref{Sob:for1.39ter}. $\blacksquare$

\bigskip

\begin{prop}\label{Sob:Es1.29}
	Let $F\in C^1(\mathbb{R})$ be such that $F'$ is bounded.  Let
	$\Omega$ be a bounded open set of $\mathbb{R}^n$ and let $u\in
	W^{1,p}(\Omega)$, $p\in[1,+\infty)$. Let us denote
	$$v:=F(u),$$
	we have $v\in W^{1,p}(\Omega)$ and
	$$\partial_jv=F'(u)\partial_ju,\quad j=1,\cdots,n.$$
\end{prop}
\textbf{Proof.} Since $F'$ is a bounded function  and $\Omega$
is a bounded set, we have $v\in L^{p}(\Omega)$. As a matter of fact
$$|v|\leq |F(u)-F(0)|+|F(0)|\leq \left\Vert
F'\right\Vert_{L^{\infty}(\mathbb{R})}|u|+|F(0)|\in L^{p}(\Omega).$$
Now, we apply  Theorem \ref{densit 2} and let
$\left\{u_m\right\}\subset C^{\infty}(\Omega)\cap
W^{1,p}(\Omega)$ be a sequence such that 
\begin{equation}\label{Sob:for1.30}
	u_m\rightarrow u,\quad\mbox{as } m\rightarrow \infty,\mbox{ in }
	W^{1,p}(\Omega).
\end{equation}
We have
\begin{equation}\label{Sob:for2.30}
	F(u_m)\rightarrow F(u),\quad\mbox{as } m\rightarrow \infty,\mbox{
		in } L^{p}(\Omega).
\end{equation}
Concerning the latter we have

\begin{equation*}
	\lim_{m\rightarrow\infty}\int_{\Omega}\left|
	F(u_m)-F(u)\right|^pdx\leq \left\Vert
	F'\right\Vert^p_{L^{\infty}(\mathbb{R})}\lim_{m\rightarrow\infty}\int_{\Omega}\left|
	u_m-u\right|^pdx=0.
\end{equation*}

Now let us check that

\begin{equation}\label{Sob:for3.30}
	F'(u_m)\partial_ju_m\rightarrow F'(u)\partial_ju,\quad\mbox{as }
	m\rightarrow \infty,\mbox{ in } L^{p}(\Omega),
\end{equation}
for $j=1,\cdots,n$.

We have
\begin{equation*}
	\begin{aligned}
		\left\Vert
		F'(u_m)\partial_ju_m-F'(u)\partial_ju\right\Vert_{L^{p}(\Omega)}&\leq
		\left\Vert
		F'(u_m)\left(\partial_ju_m-\partial_ju\right)\right\Vert_{L^{p}(\Omega)}+\\&+\left\Vert
		\left(F'(u_m)-F'(u)\right)\partial_ju\right\Vert_{L^{p}(\Omega)}\leq\\&
		\leq \left\Vert F'\right\Vert_{L^{\infty}(\mathbb{R})}\left\Vert
		\partial_ju_m-\partial_ju\right\Vert_{L^{p}(\Omega)}+\\&+\left\Vert
		\left(F'(u_m)-F'(u)\right)\partial_ju\right\Vert_{L^{p}(\Omega)}.
	\end{aligned}
\end{equation*}
Since \eqref{Sob:for1.30} holds, the second-to-last term on the right goes to
zero as $m\rightarrow \infty$, concerning the last term, it goes to zero by the Dominated Convergence Theorem. Thus, we have checked \eqref{Sob:for3.30}.

Now, by \eqref{Sob:for2.30} and \eqref{Sob:for3.30} we have 
\begin{equation*}
	\begin{aligned}
		\int_{\Omega}v\partial_j\phi
		dx&=\lim_{m\rightarrow\infty}\int_{\Omega}F(u_m)\partial_j\phi
		dx=\\&=-\lim_{m\rightarrow\infty}\int_{\Omega}\partial_j\left(F(u_m)\right)\phi
		dx=\\&=
		-\lim_{m\rightarrow\infty}\int_{\Omega}F'(u_m)\partial_ju_m\phi
		dx=\\&=-\int_{\Omega}F'(u)\partial_ju\phi dx,
	\end{aligned}
\end{equation*}
for every $\phi\in C_0^{\infty}(\Omega)$ and every $j=1,\cdots,n$. Hence

$$\partial_jv=F'(u)\partial_ju\in L^p(\Omega),\quad j=1,\cdots,n$$
so that, taking into account that $v\in L^p(\Omega)$, we have $v\in W^{1,p}(\Omega)$.
$\blacksquare$

\bigskip

\begin{prop}\label{Sob:Es1.31}
	Let $\Omega$ be a bounded open set of  $\mathbb{R}^n$, let $u\in
	W^{1,p}(\Omega)$ and $p\in [1,+\infty)$. Let us denote by
	$u_+=\max\{u,0\}$, $u_-=\min\{u,0\}$. We have $u_+,u_-\in
	W^{1,p}(\Omega)$ and
	
	\begin{equation}\label{Sob:for1.32}
		\nabla u_+=
		\begin{cases}
			\nabla u, \mbox{ for } u>0, \\
			\\
			0,\quad \mbox{ for } u\leq 0,
		\end{cases}
	\end{equation}
	
	\begin{equation}\label{Sob:for2.32}
		\nabla u_-=
		\begin{cases}
			0,\quad \mbox{ for } u\geq0, \\
			\\
			-\nabla u,\quad \mbox{ for } u< 0,
		\end{cases}
	\end{equation}
	
	\begin{equation}\label{Sob:for3.32}
		\nabla |u|=
		\begin{cases}
			\nabla u, \mbox{ for } u>0, \\
			0, \quad\mbox{ for } u=0,\\
			-\nabla u, \mbox{ for } u< 0.
		\end{cases}
	\end{equation}
\end{prop}
\textbf{Proof}. Let us prove \eqref{Sob:for1.32}. For any
$\varepsilon>0$ let us define

\begin{equation*}
	f_{\varepsilon}(t)=
	\begin{cases}
		\sqrt{t^2+\varepsilon^2}-\varepsilon, \mbox{ for } t >0, \\
		\\
		0,\quad\quad\quad \quad\quad\mbox{ for } t\geq 0.
	\end{cases}
\end{equation*}
Recalling that $u\in W^{1,p}(\Omega)$, by the Dominated Convergence Theorem we get

$$\lim_{\varepsilon\rightarrow
	0}\int_{\Omega}\left|f_{\varepsilon}(u)-u_+\right|^pdx=0,$$ As a matter of fact we have
$$\lim_{\varepsilon\rightarrow
	0}\left|f_{\varepsilon}(u)-u_+\right|^p=0, \quad \mbox{in } \Omega$$ and
$$\left|f_{\varepsilon}(u)-u_+\right|^p\leq
2^p\left(\left|f_{\varepsilon}(u)\right|^p+\left|u_+\right|^p\right)\leq
2^{p+1}|u|^p\in L^1(\Omega).$$ Now, we have

\begin{equation*}
	f'_{\varepsilon}(t)=
	\begin{cases}
		\frac{t}{\sqrt{t^2+\varepsilon^2}}, \mbox{ for } t >0, \\
		\\
		0,\quad\quad \mbox{ for } t\geq 0
	\end{cases}
\end{equation*}
and
$$\left|f'_{\varepsilon}(t)\right|\leq 1.$$ Hence, Proposition \ref{Sob:Es1.29}
implies
$$f_{\varepsilon}(u)\in W^{1,p}(\Omega).$$
Now

\begin{equation*}
	\partial_jf_{\varepsilon}(u)=
	\begin{cases}
		\frac{u\partial_ju}{\sqrt{u^2+\varepsilon^2}}, \mbox{ for } u >0, \\
		\\
		0,\quad\quad \mbox{ for } u\geq 0,
	\end{cases}
\end{equation*}
for $j=1,\cdots, n$. Hence, for any $\phi\in
C^{\infty}_0(\Omega)$,

\begin{equation}\label{Sob:for1.33}
	\begin{aligned}
		\int_{\Omega}u_+\partial_j\phi dx&=\lim_{\varepsilon\rightarrow
			0}\int_{\Omega}f_{\varepsilon}(u)\partial_j\phi
		dx=\\&=-\lim_{\varepsilon\rightarrow
			0}\int_{u>0}\frac{u\partial_ju}{\sqrt{u^2+\varepsilon^2}}\phi
		dx=\\&= -\int_{u>0}\partial_ju\phi dx
	\end{aligned}
\end{equation}
in the last step we have applied the Dominated Convergence Theorem. Therefore

$$\int_{\Omega}u_+\partial_j\phi dx=-\int_{\Omega}\partial_ju\chi_{u>0}\phi
dx,\quad \forall\phi \in C^{\infty}_0(\Omega),$$ from which we get 
\eqref{Sob:for1.32}. Concerning \eqref{Sob:for2.32},
it suffices to notice that $u_-=(-u)_+$. All in all, \eqref{Sob:for3.32}
follows by \eqref{Sob:for1.32} and \eqref{Sob:for2.32}  (recall that $|u|=u_++u_-$). $\blacksquare$

\bigskip

\underline{\textbf{Exercise 1.}} We say that $f:\mathbb{R}\rightarrow\mathbb{R}$ is a piecewise $C^1$ function, provided that $f$ satisfies what follows: $f$ is a continuous function, it has a continuous derivative  in $\mathbb{R}\setminus\left\{a_1,\cdots, a_l\right\}$, where $a_j\in \mathbb{R}$ and $f$ has the right and the left derivatives in $a_j$, for $j=1,\cdots,l$ and such derivatives are finite.
	
	Prove that if $f$ is a piecewise $C^1$ function, $f'\in
	L^{\infty}(\mathbb{R})$, $\Omega$ is a bounded open set of
	$\mathbb{R}^n$ and $u\in W^{1,p}(\Omega)$, $p\in [1,+\infty)$, then we have
	$f(u)\in W^{1,p}(\Omega)$ and

	\begin{equation*}
		\nabla(f(u))=
		\begin{cases}
			f'(u)\nabla u, \mbox{ for } u \notin \left\{a_1,\cdots, a_l\right\} , \\
			\\
			0,\quad \mbox{ for } u\in  \left\{a_1,\cdots, a_l\right\}.
		\end{cases}
	\end{equation*}
	[Hint: consider preliminarly the case $l=1$ and, in doing so, first address to the
	case in which $f(0)=0$; observe that
	
	\begin{equation*}
		f(t)=
		\begin{cases}
			f_1(t), \mbox{ for } t>0 , \\
			f_2(t), \mbox{ for } t\leq 0,
		\end{cases}
	\end{equation*}
	where $f_1,f_2\in C^1(\mathbb{R})$ and $f'_1,f'_2\in
	L^{\infty}(\mathbb{R})$. Let us note that $f(t)=f_1(t_+)+f_2(-t_-)$ and use Proposition \ref{Sob:Es1.29} ...].

\bigskip

\underline{\textbf{Exercise 2.}} Let $\Omega$ be an open set of $\mathbb{R}^n$ and let $u,v\in W^{1,p}(\Omega)\cap L^{\infty}(\Omega)$. Prove that
	$$uv\in W^{1,p}(\Omega)\cap L^{\infty}(\Omega)$$ and
	\begin{equation}\label{Sob:for1.37}
		\nabla(uv)=v\nabla u+u\nabla v.
	\end{equation}
\textbf{Solving.} Let $\phi\in C_0^{\infty}(\Omega)$ and let $V$ be
an open set such that
$$\mbox{supp } \phi\subset V\Subset \Omega.$$ Let

$$u_{\varepsilon}(x)=\int_{\Omega}\eta_{\varepsilon}(x-y)u(y)dy,\quad
v_{\varepsilon}(x)=\int_{\Omega}\eta_{\varepsilon}(x-y)v(y)dy.$$
We have $(u_{\varepsilon})_{|_V},(v_{\varepsilon})_{|_V}\in
C^{\infty}\left(\overline{V}\right)$ and

\begin{equation*}
	\begin{aligned}
		\int_{\Omega}uv\partial_j\phi dx&=\int_{V}uv\partial_j\phi
		dx=\\&=\lim_{\varepsilon\rightarrow
			0}\int_{V}u_{\varepsilon}v_{\varepsilon}\partial_j\phi dx=\\&=
		-\lim_{\varepsilon\rightarrow 0}\int_{V}\left[\left(\partial_j
		u_{\varepsilon}\right)v_{\varepsilon}+u_{\varepsilon}\partial_jv_{\varepsilon}\right]\phi
		dx=\\&=-\int_{V}\left[\left(\partial_j u\right)v+u\partial_jv\right]\phi dx,
	\end{aligned}
\end{equation*}
first limit is justified by the Dominated Convergence Theorem, the second limit is justified as follows

\begin{equation*}
	\begin{aligned}
		\left\Vert \left(\partial_j u_{\varepsilon}\right)v_{\varepsilon}\phi-\left(\partial_j
		u \right)v\phi \right\Vert_{L^1(V)}&\leq \left\Vert
		\phi\right\Vert_{L^{\infty}(V)} \left\Vert
		v\right\Vert_{L^{\infty}(V)}|V|^{1/p'}\left\Vert
		\partial_j u_{\varepsilon}-\partial_j u
		\right\Vert_{L^p(V)}+\\&+\left\Vert \phi\right\Vert_{L^{\infty}(V)}
		\left\Vert
		\partial_ju\right\Vert_{L^{p}(V)}\left\Vert
		v_{\varepsilon}-v \right\Vert_{L^{p'}(V)}\rightarrow 0,\mbox{ as }
		\varepsilon\rightarrow 0,
	\end{aligned}
\end{equation*}
similarly we argue for $\left\Vert\left( \partial_j
v_{\varepsilon}\right)u_{\varepsilon}\phi-\left(\partial_j v\right) u\phi
\right\Vert_{L^1(V)}$. Hence, we get

\begin{equation}
	\partial_j(uv)=\left(\partial_jv\right)u+\left(\partial_ju\right) v.
\end{equation}
Since (by H\"{o}lder inequality) $\left(\partial_jv\right) u+\left(\partial_ju\right)
v\in L^p(\Omega)$, we get \eqref{Sob:for1.37}. $\clubsuit$

\section{The extension theorems} \label{Sob:sec4.40}

Let us start by some propositions about the space $W_0^{k,p}(\Omega)$.
\begin{prop}\label{Sob:prop2.4}
	If $k\in \mathbb{N}$ and $p\in [1,+\infty)$ then
	$$W^{k,p}(\mathbb{R}^{n})=W_0^{k,p}(\mathbb{R}^{n}).$$
\end{prop}
\textbf{Proof.} We limit ourselves to the case $k=1$, the case $k>1$ can be
proved in a similar way and is left to the reader. 

Let $R>1$ and let $\zeta_R\in C_0^{\infty}(\mathbb{R}^n)$ satisfy 
$$0\leq \zeta_R\leq 1,\quad\mbox{in } \mathbb{R}^n,$$
$$\zeta_R(x)=1,\quad \forall x\in B_R;\quad \zeta_R(x)=0,\quad \forall x\in \mathbb{R}^n\setminus
B_{2R},$$
$$\left|\nabla \zeta_R\right|\leq C,\quad\mbox{in } \mathbb{R}^n,$$
where $C$ is independent of $R$.

If $u\in W^{1,p}(\mathbb{R}^{n})$, we have

\begin{equation}\label{Sob:for1.51}
	\left\Vert u-\zeta_Ru
	\right\Vert_{W^{1,p}(\mathbb{R}^{n})}\rightarrow 0,\quad\mbox{as }
	R\rightarrow\infty.
\end{equation}
Let us check \eqref{Sob:for1.51}.
$$\left\Vert u-\zeta_Ru
\right\Vert_{L^{p}(\mathbb{R}^{n})}\leq \left\Vert u
\right\Vert_{L^{p}(\mathbb{R}^{n}\setminus B_R)}\rightarrow
0,\quad\mbox{as } R\rightarrow\infty$$ and, for any $j=1,\cdots,n$,
\begin{equation*}
	\begin{aligned}
		\left\Vert \partial_ju-\partial_j(\zeta_Ru)
		\right\Vert_{L^{p}(\mathbb{R}^{n})}&=\left\Vert
		(1-\zeta_R)\partial_ju-u\partial_j\zeta_R
		\right\Vert_{L^{p}(\mathbb{R}^{n})}\leq\\&\leq \left\Vert
		\partial_ju
		\right\Vert_{L^{p}(\mathbb{R}^{n}\setminus B_R)}+C\left\Vert u
		\right\Vert_{L^{p}(\mathbb{R}^{n}\setminus B_R)}\rightarrow 0,\quad
		\mbox{as } R\rightarrow \infty.
	\end{aligned}
\end{equation*}

In order to complete the proof, firstly we observe (by 
Theorem \ref{Sob:teo26R})

\begin{equation}\label{Sob:for2.51}
	\left(\zeta_Ru\right)\star\eta_{\varepsilon}\rightarrow
	\zeta_R u,\quad\mbox{as }\varepsilon\rightarrow 0,\mbox{in }
	W^{1,p}(\mathbb{R}^{n}),
\end{equation}
where $\eta_{\varepsilon}$ a mollifier. Moreover, let
$\delta$ be any positive number, let $R_0>1$ be such that

$$\left\Vert u-\zeta_{R_0}u
\right\Vert_{W^{1,p}(\mathbb{R}^{n})}<\frac{\delta}{2}$$ and let 
$\varepsilon_0>0$ be such that

\begin{equation*}
	\left\Vert\left(\zeta_{R_0} u\right)\star\eta_{\varepsilon_0}-
	\zeta_{R_0}u \right\Vert_{W^{1,p}(\mathbb{R}^{n})}<\frac{\delta}{2}.
\end{equation*}
From the last two inequalities and the triangle inequality
we get
\begin{equation*}
	\left\Vert\left(\zeta_{R_0}u\right)\star\eta_{\varepsilon_0}- u
	\right\Vert_{W^{1,p}(\mathbb{R}^{n})}<\delta.
\end{equation*}
Since
$\left(\zeta_{R_0}u\right)\star\eta_{\varepsilon_0}\in
C^{\infty}_0(\mathbb{R}^n)$, the Proposition is proved.
$\blacksquare$

\bigskip

\begin{theo}[\textbf{The first Poincar\'{e} inequality}]\label{Poincar}
	\index{Theorem:@{Theorem:}!- first Poincar\'{e} inequality@{- first Poincar\'{e} inequality}}
	Let $\Omega$ be a bounded open set of $\mathbb{R}^n$. Let $p\in
	[1,+\infty]$, $k\in \mathbb{N}$. The following inequality holds true, for any $\alpha\in
	\mathbb{N}_0^n$, $|\alpha|\le k-1$, 
	\begin{equation}\label{Poinc}
		\left\Vert \partial^{\alpha}u\right\Vert_{L^p(\Omega)}\leq C
		d^{k-|\alpha|}\sum_{|\beta|= k}\left\Vert
		\partial^{\beta}u \right\Vert_{L^p(\Omega)}, \quad\forall u\in
		W^{k,p}_0(\Omega),
	\end{equation}
	where $d$ is the diameter of $\Omega$ and $C$ depends on  $n$ and
	 $k$ only.
\end{theo}
\textbf{Proof.} We restrict ourselves to the case $k=1$, actually starting from this case \eqref{Poinc} can easily be deduced by induction. It is not restrictive to assume $0\in \Omega$ and

$$\Omega \subset [-d,d]^n.$$ Let $u\in W^{1,p}_0(\Omega)$. Let
$\left\{u_j\right\}$ be a sequence in
$C^{\infty}_0(\Omega)$ such that

$$\left\{u_j\right\}\rightarrow u,\quad \mbox{ in } W^{1,p}(\Omega).$$ For any $p\in
[1,+\infty)$ and any $j\in \mathbb{N}$, we have

\begin{equation*}
	\begin{aligned}
		\left|u_j(x)\right|&=\left|u_j(x)-u_j(x',-d)\right|=\\&=
		\left|\int_{-d}^{x_n}\partial_{y}u_j(x',y)dy\right|\leq\\&\leq
		\int_{-d}^{d}\left|\partial_{y}u_j(x',y)\right|dy\leq\\&\leq
		(2d)^{1/p'}\left(\int_{-d}^{d}\left|\partial_{y}u_j(x',y)\right|^pdy\right)^{1/p}.
	\end{aligned}
\end{equation*}
Hence
\begin{equation}\label{Sob:for1.53}
	\left|u_j(x)\right|^p\leq
	(2d)^{p-1}\int_{-d}^{d}\left|\partial_{y}u_j(x',y)\right|^pdy.
\end{equation}
Let us integrate both the sides of \eqref{Sob:for1.53} over $[-d,d]$
w.r.t. $x_n$. We get

\begin{equation*}
	\int_{-d}^{d}\left|u_j(x',x_n)\right|^pdx_n\leq
	(2d)^{p}\int_{-d}^{d}\left|\partial_{y}u_j(x',y)\right|^pdy.
\end{equation*}
Now, let us integrate both the sides of the last inequality over $[-d,d]^{n-1}$. We get

\begin{equation*}
	\left(\int_{\Omega}\left|u_j(x)\right|^pdx\right)^{1/p}\leq
	2d\left(\int_{-d}^{d}\left|\partial_{x_n}u_j(x)\right|^pdx\right)^{1/p}.
\end{equation*}
Passing to the limit as $j\rightarrow\infty$, we obtain

\begin{equation*}
	\left\Vert u\right\Vert_{L^p(\Omega)}\leq C d \left\Vert \nabla u
	\right\Vert_{L^p(\Omega)}.
\end{equation*}
If $p=+\infty$, then we have
$$\left|u_j(x)\right|=\left|u_j(x)-u_j(x',-d)\right|\leq 2d \left\Vert
\partial_{x_n}u_j\right\Vert_{L^{\infty}(\Omega)},$$ from which,
passing to the limit as $j\rightarrow\infty$, we have

\begin{equation*}
	\left\Vert u\right\Vert_{L^{\infty}(\Omega)}\leq C d \left\Vert
	\nabla u \right\Vert_{L^{\infty}(\Omega)}.
\end{equation*}
$\blacksquare$

\bigskip

\textbf{Remarks.}

\smallskip

\noindent\textbf{1.} From the proof of Proposition
\ref{Poincar} it is evident that inequality \eqref{Poinc}
also holds if $\Omega$ is contained in a strip of
$\mathbb{R}^n$ of the type $\mathbb{R}^{n-1}\times [-d,d]$ or isometric to it.

\smallskip

\noindent\textbf{2.} Proposition \ref{Poincar} implies that, the following seminorms are actually norms on
$W_0^{k,p}(\Omega)$  
$$\sum_{|\beta|= k}\left\Vert
\partial^{\beta}u \right\Vert_{L^p(\Omega)},\quad \left(\sum_{|\beta|= k}\left\Vert
\partial^{\beta}u \right\Vert^p_{L^p(\Omega)}\right)^{1/p}.$$
Moreover such norms are equivalent to the norm
$$\left\Vert u\right\Vert_{W^{k,p}(\Omega)}.$$
$\blacklozenge$

\bigskip

\begin{prop}\label{Sob:propWextra}
	Let $p\in [1,+\infty]$, $k\in \mathbb{N}$ and let $\Omega$ and
	$\widetilde{\Omega}$ be open sets of $\mathbb{R}^n$ such that
	$\Omega\subset\widetilde{\Omega}$. Let $u\in W_0^{k,p}(\Omega)$.
	
	Denoting
	
	\begin{equation}\label{estensione-00}
		\widetilde{u}=
		\begin{cases}
			u, \mbox{ in } \Omega, \\
			\\
			0, \mbox{ in } \widetilde{\Omega}\setminus \Omega,%
		\end{cases}
	\end{equation}
	we have $\widetilde{u}\in W_0^{k,p}\left(\widetilde{\Omega}\right)$.
\end{prop}
\textbf{Proof.} Let $u\in W_0^{k,p}(\Omega)$ and let
$\left\{u_j\right\}\subset C^{\infty}_0(\Omega)$ be a sequence such that 
$$u_j\rightarrow u, \quad \mbox{as } j\rightarrow\infty,\mbox{in }
W^{k,p}(\Omega).$$ Hence,  by denoting

\begin{equation}
	\widetilde{u}_j=
	\begin{cases}
		u_j, \mbox{ in } \Omega, \\
		\\
		0, \mbox{ in } \widetilde{\Omega}\setminus \Omega,%
	\end{cases}
\end{equation}
we have $\left\{\widetilde{u}_j\right\}\subset
C^{\infty}_0\left(\widetilde{\Omega}\right)$ and
$$\widetilde{u}_j\rightarrow \widetilde{u}, \quad \mbox{as } j\rightarrow\infty,\mbox{in }
W^{k,p}\left(\widetilde{\Omega}\right).$$ Therefore $\widetilde{u}\in
W_0^{k,p}\left(\widetilde{\Omega}\right)$. $\blacksquare$

\bigskip

The Main Theorem of the present Section is the following one.

\begin{theo}[\textbf{extension in $W^{1,p}$}]\label{Sob:teo1.4}
	\index{Theorem:@{Theorem:}!- extension in $W^{1,p}$@{- extension in $W^{1,p}$}}
	Let $\Omega$ be a bounded open set of $\mathbb{R}^{n}$ whose boundary is of
	class $C^{0,1}$ with cosntants $r_0,M_0$. Let $d_0$ be the diameter of
	$\Omega$. Let $\widetilde{\Omega}$ be an open set of $\mathbb{R}^{n}$
	such that $\Omega\Subset \widetilde{\Omega}$ and let $p\in[1,+\infty)$.
	
	Then there exists a linear bounded operator 
	
	\begin{equation}\label{Sob:for1.40}
		E:W^{1,p}(\Omega)\rightarrow W^{1,p}(\mathbb{R}^{n}),
	\end{equation}
	which satisfies, for any $u\in W^{1,p}(\Omega)$,
	
	\begin{equation}\label{Sob:for2.40}
		Eu=u,\quad\mbox{in } \Omega,
	\end{equation}
	
	\begin{equation}\label{Sob:for3.40}
		\mbox{supp }(Eu)\subset \widetilde{\Omega}.
	\end{equation}
	Moreover, there exists a constant $C$ depending on $r_0$, $M_0$, $d_0$, $n$ and
	$p$ only such that
	
	\begin{equation}\label{Sob:for4.40r}
		\left\Vert Eu \right\Vert_{W^{1,p}(\mathbb{R}^{n})}\leq C\left\Vert
		u \right\Vert_{W^{1,p}(\Omega)},\quad\forall u\in W^{1,p}(\Omega).
	\end{equation}
\end{theo}
\textbf{Proof.} Let $x_0\in \partial \Omega$. We may assume (up tp an isometry) that $x_0=0$ and

$$\Omega\cap Q_{r_0,2M_0}=\left\{x\in Q_{r_0,2M_0}: \mbox{ }
x_n>\varphi(x')\right\},$$ where $\varphi\in
C^1\left(\overline{B'}_{r_0}\right)$ satisfies
$$\varphi (0)=0$$ 
and $$\left\Vert \varphi
\right\Vert_{C^1\left(\overline{B'}_{r_0}\right)}=\left\Vert \varphi
\right\Vert_{C^0\left(\overline{B'}_{r_0}\right)}+r_0[\varphi]_{0,1,B'_{r_0}}\leq M_0r_0.$$

Set

$$V^+=Q_{\frac{r_0}{4},\frac{M_0}{4}}\cap \Omega,\quad
V^-=Q_{\frac{r_0}{4},\frac{M_0}{4}}\setminus \overline{V^+}.$$ 
Notice that, for every $x'\in \overline{B'}_{r_0/4}$ we have
\begin{equation}\label{correct:2-4-23-1}
\left|\varphi(x')\right|=\left|\varphi(x')-\varphi(0)\right|\leq [\varphi]_{0,1,B'_{r_0}}\left|x'\right|\leq \frac{M_0r_0}{4}.
\end{equation}

\begin{figure}\label{figura-p41sob}
	\centering
	\includegraphics[trim={0 0 0 0},clip, width=12cm]{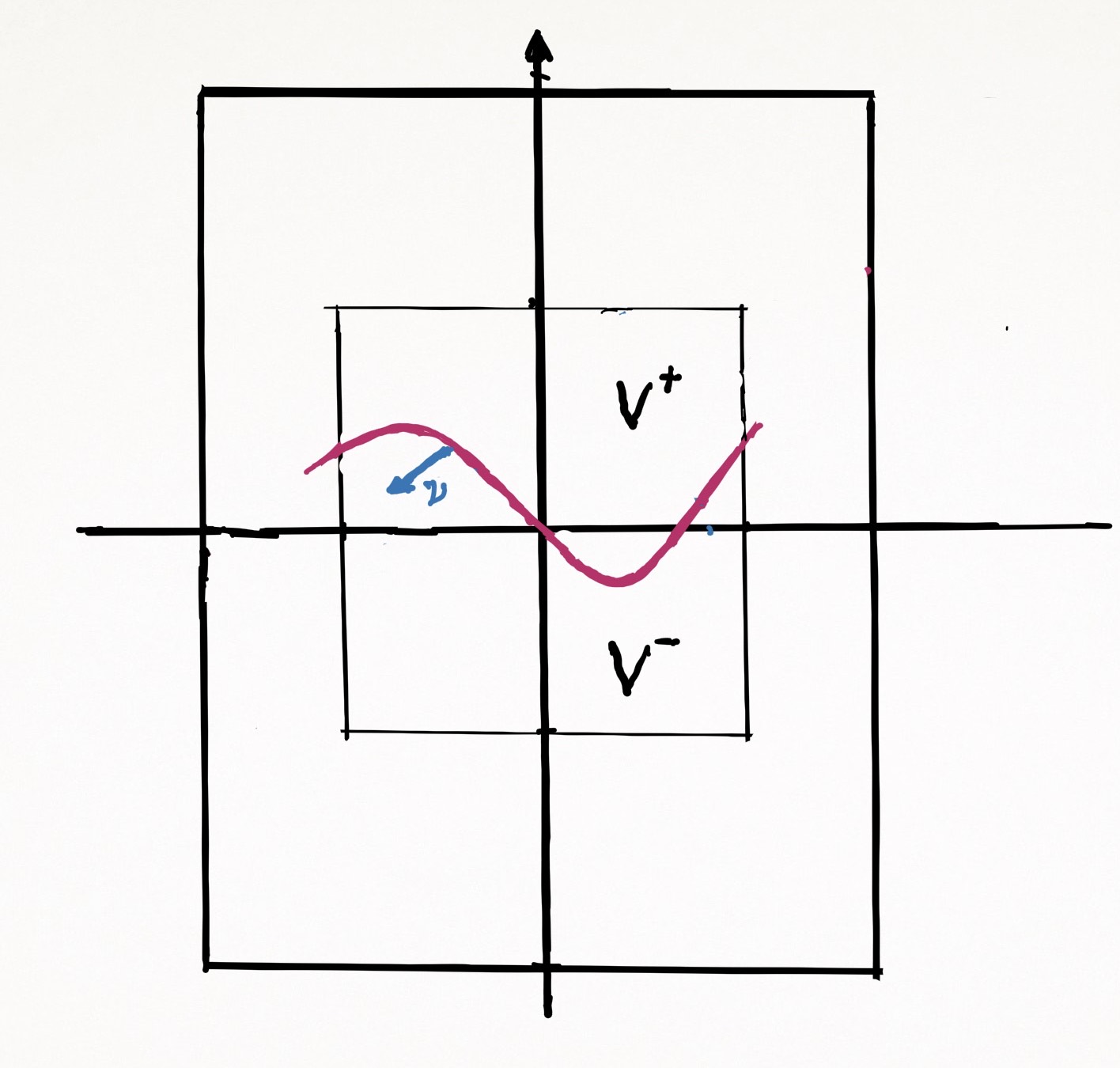}
	\caption{}
\end{figure}

First, we assume that $u\in
C^{\infty}\left(\overline{\Omega}\right)$ and we define
$$\nu(x')=\frac{\left(\nabla
	\varphi(x'),-1\right)}{\sqrt{\left|\nabla_{x'}\varphi\right|^2+1}},
\quad x'\in B'_{r_0}, $$

\begin{equation*}
	\overline{u}(x)=
	\begin{cases}
		u(x), \mbox{ in } V^+, \\
		\\
		v(x), \mbox{ in } V^-,%
	\end{cases}
\end{equation*}
where

$$v(x)=u(x',2\varphi(x')-x_n).$$

\medskip

\textbf{Claim.} $\overline{u}\in
W^{1,p}\left(Q_{\frac{r_0}{4},\frac{M_0}{4}}\right)$ and

\begin{equation}\label{Sob:for1.42}
	\left\Vert \overline{u}
	\right\Vert_{W{1,p}\left(Q_{\frac{r_0}{4},\frac{M_0}{4}}\right)}\leq C
	\left\Vert u \right\Vert_{W{1,p}\left(\Omega\right)},
\end{equation}
where $C$ depends on $M_0$ only.

\medskip

\textbf{Proof of Claim.} Let $\Phi\in
C^{\infty}_0\left(Q_{\frac{r_0}{4},\frac{M_0}{4}}\right)$ and $1\leq i\leq n$.
Denoting by $\Gamma$ the graph of $\varphi_{|B'_{r_0/4}}$, by the divergence Theorem we get
\begin{equation*}
	\begin{aligned}
		\int_{Q_{\frac{r_0}{4},\frac{M_0}{4}}}\overline{u}\partial_{i}\Phi
		dx&=\int_{V^+}u\partial_{i}\Phi dx+\int_{V^-}v\partial_{i}\Phi
		dx=\\&=-\int_{V^+}\partial_{i}u\Phi dx+\int_{\Gamma} u\Phi
		(\nu\cdot e_i) dS-\\&-\int_{V^-}\partial_{i}v\Phi dx-\int_{\Gamma}
		v\Phi (\nu\cdot e_i) dS=\\&=-\int_{Q_{\frac{r_0}{4},\frac{M_0}{4}}}w_i\Phi
		dx+\int_{\Gamma} (u-v)\Phi (\nu\cdot e_i) dS,
	\end{aligned}
\end{equation*}
where

\begin{equation}\label{Sob:for-stella-43}
	w_i(x)=
	\begin{cases}
		\partial_{i}u(x), \mbox{ in } V^+, \\
		\\
		\partial_{i} v(x), \mbox{ in } V^-.
	\end{cases}
\end{equation}
On the other hand, 

$$(u-v)(x',\varphi(x'))=0,\quad\forall x\in B'_{r_0/4},$$ hence

\begin{equation}\label{Sob:for1.43}
	\int_{Q_{\frac{r_0}{4},\frac{M_0}{4}}}\overline{u}\partial_{i}\Phi
	dx=-\int_{Q_{\frac{r_0}{4},\frac{M_0}{4}}}w_i\Phi dx,\quad\forall \Phi \in
	C^{\infty}_0\left(Q_{\frac{r_0}{2},M_0}\right).
\end{equation}
Therefore
\begin{equation}\label{Sob:forExtra.44}
	\partial_{i}\overline{u}=w_i(x),\quad\forall x\in
	Q_{\frac{r_0}{4},\frac{M_0}{4}}.
\end{equation}
Now, let us notice that

\begin{equation}\label{Sob:for2.43}
	\int_{Q_{\frac{r_0}{4},\frac{M_0}{4}}}\left|\overline{u}(x)\right|^pdx=\int_{V^+}\left|u(x)\right|^pdx+\int_{V^-}\left|v(x)\right|^pdx
\end{equation}
and
\begin{equation}\label{Sob:for3.43}
	\begin{aligned}
		\int_{V^-}\left|v(x)\right|^pdx&=\int_{B'_{r_0/4}}dx'\int^{\varphi(x')}_{-\frac{M_0r_0}{4}}\left|u(x',2\varphi(x')-x_n)\right|^pdx_n
		=\\&=\int_{B'_{r_0/4}}dx'\int^{2\varphi(x')+\frac{M_0r_0}{4}}_{\varphi(x')}\left|u(x',\xi_n)\right|^pd\xi_n\leq\\&\leq
		\int_{\Omega}\left|u(x)\right|^pdx.
	\end{aligned}
\end{equation}
In the last inequality we have used \eqref{correct:2-4-23-1}.  By
\eqref{Sob:for2.43} and \eqref{Sob:for3.43}, we have,
\begin{equation}\label{Sob:for1.44}
	\int_{Q_{\frac{r_0}{4},\frac{M_0}{4}}}\left|\overline{u}(x)\right|^pdx\leq
	2\int_{\Omega}\left|u(x)\right|^pdx.
\end{equation}
Now,  \eqref{Sob:for-stella-43} gives 
$$\left|\nabla \overline{u}(x)\right|\leq C\left|(\nabla u)
(x',2\varphi(x')-x_n)\right|,\quad \forall x\in V^-,$$ where $C$
depends on $M_0$ only. Hence

\begin{equation*}
	\begin{aligned}
		\int_{Q_{\frac{r_0}{4},\frac{M_0}{4}}}\left|\nabla\overline{u}(x)\right|^pdx&\leq\int_{V^+}\left|\nabla
		u(x)\right|^pdx+C\int_{V^-}\left|(\nabla u)
		(x',2\varphi(x')-x_n)\right|^pdx\leq\\&\leq
		C\int_{\Omega}\left|\nabla u(x)\right|^pdx.
	\end{aligned}
\end{equation*}
From the just obtained inequality e from \eqref{Sob:for1.44} we obtain
\eqref{Sob:for1.42}. Claim is proved.

\medskip

Since $\partial \Omega$ is a compact set, there exist $x_{0,
	1},\cdots, x_{0, N}\in \partial \Omega$ such that
$$\partial \Omega\subset \bigcup_{j=1}^N \widetilde{Q}_{\frac{r_0}{4},\frac{M_0}{4}}(x_{0,j})$$
where, for any $j=1,\cdots,N$,
$\widetilde{Q}_{\frac{r_0}{4},\frac{M_0}{4}}(x_{0,j})$ are suitable cylinders which are
isometric to $Q_{\frac{r_0}{4},\frac{M_0}{4}}$. Moreover, let us denote 
$\overline{u}_j$ the extensions of $u$ on
$\widetilde{Q}_{\frac{r_0}{4},\frac{M_0}{4}}(x_{0,j})$. Let us employ the partition of unity  (Lemma \ref{Sob:lem3.3}). Set
$V_j=\widetilde{Q}_{\frac{r_0}{4},\frac{M_0}{4}}(x_{0,j})$, $j=1,\cdots,N$,
we have that there exist $\zeta_0,\zeta_1,\cdots,\zeta_N\in
C_0^{\infty}(\mathbb{R}^n)$ such that
$$0\leq \zeta_j(x)\leq 1,\quad j=1,\cdots,N,\ \ \forall x\in \mathbb{R}^n$$
$$\mbox{supp } \zeta_j\subset V_j, \quad j=1,\cdots,N, \quad \mbox{supp }
\zeta_0\subset \mathbb{R}^n\setminus \partial\Omega,$$ $$\sum_{j=0}^N
\zeta_j(x)=1,\quad\forall x\in \mathbb{R}^n,$$

and
$$\sum_{j=1}^N
\zeta_j(x)=1,\quad\mbox{for every }\  x  \ \mbox{ in a neighborhood of } \ \partial \Omega.$$
 Be, also,
$\eta\in C_0^{\infty}\left(\widetilde{\Omega}\right)$, such that
$0\leq \eta\leq 1$, $\eta(x)=1$, for $x\in \Omega$. Set
$$ \widetilde{u}=\eta\left(\zeta_0u+\sum_{j=1}^N\zeta_j \overline{u}_j\right),$$
By \eqref{Sob:for1.42} and by the triangle inequality, we have

\begin{equation}\label{Sob:for1.45}
	\left\Vert \widetilde{u}
	\right\Vert_{W^{1,p}\left(\mathbb{R}^n\right)}\leq C \left\Vert u
	\right\Vert_{W^{1,p}\left(\Omega\right)},
\end{equation}
where $C$ depends on $r_0$, $n$ and $M_0$ only. Moreover, we have

\begin{equation}\label{Sob:for2.45}
	\widetilde{u}(x)=u(x),\quad\forall x\in \Omega,
\end{equation}

\begin{equation}\label{Sob:for3.45}
	\mbox{supp } \widetilde{u} \subset\widetilde{\Omega}.
\end{equation}
Now, let us denote

$$Eu:=\widetilde{u},\quad\forall u\in
C^{\infty}\left(\overline{\Omega}\right).$$ $E$ is a linear operator and
satisfies the inequality

\begin{equation}\label{Sob:for1.46}
	\left\Vert Eu \right\Vert_{W^{1,p}\left(\mathbb{R}^n\right)}\leq C
	\left\Vert u \right\Vert_{W^{1,p}\left(\Omega\right)}.
\end{equation}

Now, let $u\in W^{1,p}(\Omega)$ and apply Theorem \ref{densit
	2}. Let therefore be \\ $\left\{u_m\right\}\subset
C^{\infty}\left(\overline{\Omega}\right)$ such that

$$\left\{u_m\right\}\rightarrow u, \quad \mbox{in
}W^{1,p}\left(\Omega\right).$$ We have, by \eqref{Sob:for1.46}, 

\begin{equation*}
	\left\Vert Eu_m-Eu_{m'}
	\right\Vert_{W^{1,p}\left(\mathbb{R}^n\right)}\leq C \left\Vert
	u_m-u_{m'} \right\Vert_{W^{1,p}\left(\Omega\right)}.
\end{equation*}
Hence  $\left\{Eu_m\right\}$ is a Cauchy sequence in
$W^{1,p}\left(\mathbb{R}^n\right)$ consequently it converges to a function which we continue to denote by $\widetilde{u}\in
W^{1,p}\left(\mathbb{R}^n\right)$ which satisfies trivially  
\eqref{Sob:for1.45}--\eqref{Sob:for3.45}. $\blacksquare$

\bigskip

We merely state, with some comments, the Theorem of
extension for $W^{k,p}(\Omega)$, where $k\geq 1$. We refer to \cite{G-T} for a proof.

\begin{theo}[\textbf{extension in $W^{k,p}$}]\label{Sob:teo2.4}
	\index{Theorem:@{Theorem:}!- extension in $W^{k,p}$@{- extension in $W^{k,p}$}}
	Let $k\geq 1$ and $p\in[1,+\infty)$. Let $\Omega$ be a bounded open set
	of $\mathbb{R}^{n}$ whose boundary is of class $C^{k-1,1}$ with costants
	$r_0,M_0$. Let $d_0$ be the diameter of $\Omega$. Let
	$\widetilde{\Omega}$ an open set of $\mathbb{R}^{n}$ such that
	$\Omega\Subset \widetilde{\Omega}$. 
	
	Then there exists a bounded linear operator
	
	\begin{equation*}
		E:W^{k,p}(\Omega)\rightarrow W^{k,p}(\mathbb{R}^{n}),
	\end{equation*}
	such that, for any $u\in W^{k,p}(\Omega)$ we have
	
	\begin{equation*}
		Eu=u,\quad\mbox{on } \Omega,
	\end{equation*}
	
	\begin{equation*}
		\mbox{supp }(Eu)\subset \widetilde{\Omega}.
	\end{equation*}
	Moreover, there exists a constant $C$ depending on $r_0$, $M_0$, $d_0$, $n$
	$k$ and $p$ only, such that
	
	\begin{equation}\label{Sob:for4.40}
		\left\Vert Eu \right\Vert_{W^{k,p}(\mathbb{R}^{n})}\leq C\left\Vert
		u \right\Vert_{W^{k,p}(\Omega)},\quad\forall u\in W^{k,p}(\Omega).
	\end{equation}
\end{theo}

\bigskip

\underline{\textbf{Exercise 1.}}
	\textbf{(i)} Let $u\in C^{\infty}\left(\overline{B^+_r}\right)$,
	where $B^{\pm}_r=\left\{x\in \mathbb{R}^n: \mbox{ } |x|<r\mbox{,
	}x_n\gtrless 0\right\}$. Let us define
	
	\begin{equation*}
		\overline{u}(x)=
		\begin{cases}
			u(x), \mbox{ in } B^+_r, \\
			\\
			v(x), \mbox{ in } B^-_r,
		\end{cases}
	\end{equation*}
	where
	$$v(x)=-3u(x',-x_n)+4u\left(x',-\frac{x_n}{2}\right).$$ Prove that, if $p\in
	[1,+\infty)$ then $\overline{u}(x)\in W^{2,p}(B_r)$ and the following inequality holds true
	
	\begin{equation*}
		\left\Vert \overline{u} \right\Vert_{W^{2,p}(B_r)}\leq C\left\Vert u
		\right\Vert_{W^{2,p}(\Omega)},\quad\forall u\in W^{2,p}(B^+_r).
	\end{equation*}
	\textbf{(ii)} Let us define
	$$C^{\infty}\left(\overline{B^+_r}\right)\ni u\rightarrow Eu=\overline{u}\in
	W^{2,p}(B_r).$$ Prove that the operator $E$ can be extended
	to $W^{2,p}\left(B^+_r\right)$ and that it satisfies  $Eu=u$ in
	$B^+_r$.
	
	\noindent \textbf{(iii)} Let $k\in \mathbb{N}$ and let
	$c_1,\cdots,c_k$ be such that 
	$$\sum_{j=1}^kc_j\left(-\frac{1}{j}\right)^m=1,\quad, m=0,1,\cdots, k-1,$$
	(check that such $c_1,\cdots,c_k$ exist); let us define for any
	$u\in C^{\infty}\left(\overline{B^+_r}\right)$
	\begin{equation*}
		\overline{u}(x)=
		\begin{cases}
			u(x), \mbox{ in } B^+_r, \\
			\\
			w(x), \mbox{ in } B^-_r,
		\end{cases}
	\end{equation*}
	where
	$$w(x)=\sum_{j=1}^kc_ju\left(x',-\frac{x_n}{j}\right).$$
	Prove that if $p\in [1,+\infty)$ then $\overline{u}(x)\in
	W^{k,p}(B_r)$ and the following inequality holds true
	
	\begin{equation*}
		\left\Vert \overline{u} \right\Vert_{W^{k,p}(B_r)}\leq C\left\Vert u
		\right\Vert_{W^{k,p}(B^+_r)},\quad\forall u\in W^{k,p}(B^+_r).
	\end{equation*}
	Moreover, deduce that the following operator
	$$C^{\infty}\left(\overline{B^+_r}\right)\ni u\rightarrow E_ku=\overline{u}\in
	W^{k,p}(B_r)$$ can be extended to
	$W^{k,p}\left(B_r\right)$and $E_ku=u$ in $B^+_r$.

\section{Traces in $W^{1,p}(\Omega)$} \label{tracce}

 It is well--known that if
$u\in C^0\left(\overline{\Omega}\right)$, then we can define its
trace on $\partial \Omega$, namely $u_{|\partial\Omega}$. If, on the other hand, $u\in
L^p\left(\Omega\right)$, generally, it does not make sense to consider 
its trace on $\partial\Omega$. In the present Section we will see that we can define a notion of
trace that extends the known one for the functions of
$W^{1,p}\left(\Omega\right)\cap C^0\left(\overline{\Omega}\right)$.

\bigskip

More precisely we have

\begin{theo} [\textbf{trace Theorem}]\label{traccia}
	\index{Theorem:@{Theorem:}!- trace@{- trace}}
	Let $\Omega$ be a bounded open set of class $C^{0,1}$ with constants $r_0$ and
	$M_0$. Let $p\in[1,+\infty)$. Let $d_0$ be the diameter of $\Omega$.
	Then there exists an unique bounded linear operator 
	$$T:W^{1,p}(\Omega)\rightarrow L^p(\partial\Omega),$$
	which satisfies:
	
	\medskip
	
	\noindent (i) $T(u)=u_{|\partial \Omega}$ for every $u\in
	C^0\left(\overline{\Omega}\right)\cap W^{1,p}(\Omega)$;
	
	\smallskip
	
	\noindent(ii) $$\left\Vert T(u)\right\Vert_{L^p(\partial\Omega)}\leq
	\left\Vert u\right\Vert_{W^{1,p}(\Omega)}, \quad \forall u\in
	W^{1,p}(\Omega),$$ where $C$ depends by $r_0$, $M_0$, $d_0$ and $p$;
	
	\smallskip
	
	\noindent (iii) 
	
	$$\int_{\Omega}u \mbox{div } \Phi dx=-\int_{\Omega}\nabla u\cdot\Phi
	dx+\int_{\partial\Omega}(\Phi\cdot\nu)Tu dS,$$ for every $u\in
	W^{1,p}(\Omega)$ and for every $\Phi\in C^1\left(\overline{\Omega},
	\mathbb{R}^n\right)$.
	
	The function $Tu$ is called the \textbf{trace} of $u$ on $\partial \Omega$. \index{trace of function of $W^{1,p}(\Omega)$} 
\end{theo}
\textbf{Proof.} We first notice that, since $C^{\infty}\left(\overline{\Omega}\right)\subset W^{1,p}\left(\Omega\right)$ (Theorem  \ref{densit 1}) and $\Phi$ is arbitrary in $C^1\left(\overline{\Omega},
\mathbb{R}^n\right)$, if $T$ there exists, then
it is unique.

Let us prove the existence of $T$. First, let us consider the case where $ u\in C^{\infty}\left(\overline{\Omega}\right)$. Since $\partial \Omega$ is
compact, we may consider a partition of unity
subordinate to a finite covering $\left\{V_j\right\}_{1\leq
	j\leq N}$, where
$V_j=\widetilde{Q}_{\frac{r_0}{2},M_0}\left(x_{0,j}\right)$,
$j=1,\cdots, N$
($\widetilde{Q}_{\frac{r_0}{2},M_0}\left(x_{0,j}\right)$ is a cylinder isometric to $Q_{\frac{r_0}{2},M_0}$) $$\zeta_j\in
C^{\infty}_0(\mathbb{R}^n), \quad 0\leq \zeta_j\leq 1,\quad
\mbox{supp }\zeta_j\subset V_j, \quad \sum_{j=1}^n\zeta_j=1 \mbox{
	on }
\partial\Omega.$$
Let $j\in \left\{1,\cdots, N\right\}$ be fixed. Up to isometries
we may assume
$V_j=Q_{\frac{r_0}{2},M_0}$ and

$$Q_{r_0,M_0}\cap \Omega=\left\{(x',x_n)\in Q_{r_0,M_0}:\mbox{ }
x_n>\varphi(x')\right\},$$ where
$\varphi\in C^1\left(\overline{B'}_{r_0}\right)$ satisfies
$\varphi(0)=0$ and
$$\left\Vert \varphi\right\Vert_{L^{\infty}(B'_{r_0})}+r_0\left\Vert
\nabla\varphi\right\Vert_{L^{\infty}(B'_{r_0})}\leq M_0r_0.$$

\bigskip

Let $v=\zeta_ju$. For any $t\in \left[0,\frac{M_0r_0}{2}\right]$ and any $x'\in B'_{r_0/2}$ we have

$$v\left(x',\varphi(x')\right)=v\left(x',\varphi(x')+t\right)-\int_{\varphi(x')}^{\varphi(x')+t}\partial_{x_n}v(x',x_n)dx_n.$$
H\"{o}lder inequality gives 

\begin{equation*}
	\left\vert v\left(x',\varphi(x')\right)\right\vert^p\leq
	2^{p-1}\left\vert
	v\left(x',\varphi(x')+t\right)\right\vert^p+2^{p-1}t^{p-1}\int_{\varphi(x')}^{M_0r_0}\left\vert
	\partial_{x_n}v\left(x',x_n\right)\right\vert^pdx_n
\end{equation*}
and if we integrate with respect to $x_n$ over $\left[0, \frac{M_0r_0}{2}\right]$ both the sides of the last inequality we get

\begin{equation*}
	\begin{aligned}
		\frac{M_0r_0}{2}\left\vert
		v\left(x',\varphi(x')\right)\right\vert^p&\leq
		2^{p-1}\int_{\varphi(x')}^{\varphi(x')+\frac{M_0r_0}{2}}\left\vert
		v\left(x',x_n\right)\right\vert^pdx_n+\\&+\frac{2^{p-1}}{p}\left(\frac{M_0r_0}{2}\right)^{p}\int_{\varphi(x')}^{M_0r_0}\left\vert
		u\partial_{x_n}\zeta_j+\zeta_j\partial_{x_n}u\right\vert^pdx_n.
	\end{aligned}
\end{equation*}
Now, we multiply both the sides of the last
inequality by $\sqrt{1+\left|\nabla_{x'}\varphi(x')\right|^2}$,
and we integrate over $B'_{r_0/2}$ obtaining
$$\int_{\partial \Omega}\left\vert
u\zeta_j\right\vert^pdS\leq C\int_{\Omega}\left\vert
u\right\vert^pdx+C\int_{\Omega}\left\vert \nabla u\right\vert^pdx,$$
where $C$ depends by $M_0$ e $r_0$. Therefore
\begin{equation}\label{Sob:for1.59}
	\begin{aligned}
		\left\Vert u\right\Vert_{L^{p}(\partial\Omega)}&=\left\Vert
		\sum_{j=1}^N\zeta_j u\right\Vert_{L^{p}(\partial\Omega)}\leq\\&\leq
		\sum_{j=1}^N\left\Vert \zeta_j
		u\right\Vert_{L^{p}(\partial\Omega)}\leq \\&\leq C\left\Vert
		u\right\Vert_{W^{1,p}(\Omega)},\quad\forall u\in
		C^{\infty}\left(\overline{\Omega}\right).
	\end{aligned}
\end{equation}
Set
$$Tu=u_{|\partial\Omega},\quad\forall u\in
C^{\infty}\left(\overline{\Omega}\right).$$  Inequality \eqref{Sob:for1.59}
implies
\begin{equation}\label{Sob:for2.59}
	\left\Vert Tu\right\Vert_{L^{p}(\partial\Omega)}\leq C \left\Vert
	u\right\Vert_{W^{1,p}(\Omega)}, \ \ \forall u\in
	C^{\infty}\left(\overline{\Omega}\right).
\end{equation}
Let now  $u \in W^{1,p}(\Omega)$. From Theorem
\ref{densit 2} we have that there exists a sequence $\left\{u_m\right\}$ in
$C^{\infty}\left(\overline{\Omega}\right)$ such that

$$\left\{u_m\right\}\rightarrow u,\quad\mbox{in
}W^{1,p}(\Omega).$$ In particular, \eqref{Sob:for2.59} implies that
$\left\{Tu_m\right\}$ is a Cauchy sequence in
$L^{p}(\partial\Omega)$. Set
$$Tu=\lim_{m\rightarrow\infty} Tu_m.$$ Now, we observe
that if $u\in C^{0}\left(\overline{\Omega}\right)\cap
W^{1,p}(\Omega)$, then the sequence $\left\{u_m\right\}$ constructed in the proof of Proposition  \ref{Sob:prop4.3} (with $\varepsilon=1/m$)
uniformly converges to $u$. Hence, for any $u\in
C^{0}\left(\overline{\Omega}\right)\cap W^{1,p}(\Omega)$, we have 

$$Tu=\lim_{m\rightarrow \infty}Tu_m=\lim_{m\rightarrow
	\infty}u_m=u,\quad \mbox{in } L^p(\partial\Omega).$$ Therefore, we have proved (i) and (ii), now
let us prove (iii). Let $u\in W^{1,p}(\Omega)$ and let
$\left\{u_m\right\}$ be a sequence in
$C^{\infty}\left(\overline{\Omega}\right)$ which converges to $u$ in
$W^{1,p}(\Omega)$. We have, for any $\Phi\in
C^1\left(\overline{\Omega}, \mathbb{R}^n\right)$,

\begin{equation*}
	\begin{aligned}
		\int_{\Omega}u \mbox{div } \Phi
		dx&=\lim_{m\rightarrow\infty}\int_{\Omega}u_m \mbox{div } \Phi
		dx=\\&= \lim_{m\rightarrow\infty}\left(-\int_{\Omega}\nabla
		u_m\cdot\Phi dx+\int_{\partial\Omega}(\Phi\cdot\nu)u_m
		dS\right)=\\&=-\int_{\Omega}\nabla u\cdot\Phi
		dx+\int_{\partial\Omega}(\Phi\cdot\nu)Tu dS.
	\end{aligned}
\end{equation*}
$\blacksquare$

\bigskip

\begin{rem}\label{Sob:oss-traccia} 
	If $u\in W_0^{1,p}(\Omega)$ then $Tu=0$. Actually,
	under the same assumption of Theorem \ref{traccia} the conversely is also valid, but here we omit the proof and refer to Theorem 2 of Ch. 5 of \cite{EV}. $\blacklozenge$
\end{rem}

\bigskip

The issue of traces will be taken up in Section
\ref{Sobolev-Fourier-Tracce}.

\section{The Sobolev spaces of function of one variable}
\label{Sob:sec6.61} In the present Section we will dwell briefly on the
Sobolev spaces in the case where the space dimension is equal to $1$. Let us observe that if  $I\subset
\mathbb{R}$ is a bounded open interval, then the theorems
proved in the previous sections remain valid: it is
certainly a useful exercise (left to the reader) to adapt the
proofs of these theorems and observe that they turn out to be
simplified with respect to the general case. In particular, by Theorem
\ref{densit 2} we have that $C^{\infty}\left(\overline{I}\right)$
is dense in $W^{k,p}(I)$, for every $p\in [1,+\infty)$ and, by the extension Theorem
\ref{Sob:teo1.4} it turns out that if $\widetilde{I}\Supset
I$, where $\widetilde{I}$ is an open interval of $\mathbb{R}$, then there 
exists a bounded linear operator
$$E:W^{1,p}(I)\rightarrow W^{k,p}\left(\widetilde{I}\right)$$ which satisfies

$$\left\Vert Eu\right\Vert_{W^{1,p}\left(\widetilde{I}\right)}\leq
\left\Vert u\right\Vert_{W^{1,p}(I)},\quad \forall u\in
W^{1,p}(I),$$
$$Eu_{|I}=u,\quad \mbox{supp } (Eu)\subset I.$$

Now we investigate the basic relations between the absolutely continuous functions \index{absolutely continuous functions}
 and the functions of $W^{1,p}(I)$. By the above mentioned extension Theorem, we may consider the space $W^{1,p}(\mathbb{R})$, instead of $W^{1,p}(I)$.

Let us recall that if $u\in L_{loc}^1(\mathbb{R})$ then

\begin{equation}\label{Sob:for1.62}
	\lim_{r\rightarrow 0}\frac{1}{2r}\int^{x+r}_{x-r}\left\vert
	u(t)-u(x)\right\vert dt=0,\quad\mbox{ a.e. } x\in \mathbb{R},
\end{equation}
and
\begin{equation}\label{Sob:for2.62}
	\lim_{r\rightarrow 0}\frac{1}{2r}\int^{x+r}_{x-r}u(t)
	dt=u(x),\quad\mbox{ a.e. } x\in \mathbb{R}.
\end{equation}
When \eqref{Sob:for2.62} holds true in $x$, we say that $x$ is a \textbf{Lebesgue point} of $u$ \index{Lebesgue point}.\\ For any $c\in \mathbb{R}$ we have

\begin{equation}\label{Sob:for3.62}
	\mathbb{R}\ni x\rightarrow\int^{x}_{c}u(t) dt\in
	AC_{loc}(\mathbb{R})
\end{equation}
and
\begin{equation}\label{Sob:for4.62}
	\left(\int^{x}_{c}u(t) dt\right)'=u(x),\quad\mbox{ a.e. }x \in
	\mathbb{R}.
\end{equation}
When $u\in L_{loc}^1(\mathbb{R})$, we set
\begin{equation*}
	u^{*}(x)=
	\begin{cases}
		\lim_{r\rightarrow 0}\frac{1}{2r}\int^{x+r}_{x-r}u(t)
		dt, \mbox{ provided the limit exists, } \\
		\\
		0, \quad\mbox{ otherwise}.
	\end{cases}
\end{equation*}
By \eqref{Sob:for2.62} we have
$$u^{*}(x)=u(x),\quad\mbox{ a.e. } x\in \mathbb{R}.$$
The function $u^{*}$ is called the \textbf{precise representative} of $u$ \index{precise representative}. In the sequel to this Section, if $f\in
AC_{loc}(\mathbb{R})$ we will denote by $f'$ its derivative, and if
$g\in W_{loc}^{1,p}(\mathbb{R})$, we will denote by  $\frac{d}{dx}g$ its weak derivative.

\bigskip

\begin{theo}\label{Sob:teo1.6} Let $p\in [1,+\infty)$. We have
	
	\smallskip
	
	\noindent (i) if $u\in W_{loc}^{1,p}(\mathbb{R})$, then $u^{*}\in
	AC_{loc}(\mathbb{R})$; moreover
	$$\left(u^{*}\right)'\in L^p_{loc}(\mathbb{R}),\quad \mbox{and }\quad
	\frac{d}{dx}u=\left(u^{*}\right)';$$
	
	\smallskip
	
	\noindent (ii) let $u\in L_{loc}^{p}(\mathbb{R})$. If $v\in
	AC_{loc}(\mathbb{R})$ satisfies
	
	$$u=v,\quad\mbox{ a.e. in } \mathbb{R}$$ and $v'\in L_{loc}^{p}(\mathbb{R})$
	then $u\in W_{loc}^{1,p}(\mathbb{R})$ and $\frac{d}{dx}u=v$.
\end{theo}

\bigskip

In order to prove the Theorem above we need the following Lemma

\begin{lem}\label{Sob:lem2.6} Let $u\in L^1_{loc}(\mathbb{R})$ and let $\eta_{\varepsilon}$ be
	a mollifier. If $x$ is a Lebesgue point of $u$ we have
	
	$$u^{\varepsilon}(x):=\left(\eta_{\varepsilon}\star
	u\right)(x)\rightarrow u^*(x),\quad \mbox{as }
	\varepsilon\rightarrow 0^+.$$ Hence
	$$u^{\varepsilon}(x)\rightarrow u(x),\quad \mbox{as }
	\varepsilon\rightarrow 0^+,\quad\mbox{a.e. } x\in \mathbb{R}.$$
\end{lem}

\textbf{Proof of Lemma.} Let us recall
$$\eta_{\varepsilon}=\varepsilon^{-1}\eta\left(\varepsilon^{-1}x\right),$$
where supp $\eta\subset (-1,1)$, $\eta\in C^{\infty}_0(\mathbb{R})$,
$\eta\geq 0$ and

$$\int_{\mathbb{R}}\eta(x)dx=1.$$

If $x\in \mathbb{R}$ is a Lebesgue point of $u$, then we have 
\begin{equation*}
	\begin{aligned}
		\left|u^{\varepsilon}(x)-u(x)\right|&=\left|\frac{1}{\varepsilon}\int_{\mathbb{R}}\eta\left(\frac{x-t}{\varepsilon}\right)\left(u(t)-u(x)\right)dt\right|\leq\\&\leq
		\left\Vert
		\eta\right\Vert_{L^{\infty}(\mathbb{R})}\frac{1}{\varepsilon}\int_{x-\varepsilon}^{x+\varepsilon}\left|u(t)-u(x)\right|dt\rightarrow
		0,\quad \mbox{ as } \varepsilon\rightarrow 0^+.
	\end{aligned}
\end{equation*}
$\blacksquare$

\bigskip

\textbf{Proof of Theorem \ref{Sob:teo1.6}.}

\textbf{1.} Let $u\in W_{loc}^{1,p}(\mathbb{R})$. We have
$u^{\varepsilon}\in C^{\infty}(\mathbb{R})$ and

\begin{equation}\label{Sob:for1.64}
	u^{\varepsilon}(y)=u^{\varepsilon}(x)+\int_{x}^y\left(u^{\varepsilon}\right)'(t)dt,\quad
	\forall x,y\in \mathbb{R}.
\end{equation}
Let $x_0\in \mathbb{R}$ be a Lebesgue point of $u$ (hence
$u^*(x_0)=u(x_0)$). By \eqref{Sob:for1.64} we have
\begin{equation}\label{Sob:for1.65}
	u^{\varepsilon}(x)=u^{\varepsilon}(x_0)+\int_{x}^y\left(u^{\varepsilon}\right)'(t)dt
\end{equation}
and, for any $\varepsilon,\delta>0$,
\begin{equation}\label{Sob:for2.65}
	\left|u^{\varepsilon}(x)-u^{\delta}(x)\right|\leq\left|u^{\varepsilon}(x_0)-u^{\delta}(x_0)\right|+\left|\int_{x_0}^x\left|\left(u^{\varepsilon}\right)'
	-\left(u^{\delta}\right)'\right|dt\right|.
\end{equation}
Now, Lemma \ref{Sob:lem2.6} yields
\begin{equation}\label{Sob:for3.65}
	u^{\varepsilon}(x_0)\rightarrow u(x_0),\quad \mbox{ as }
	\varepsilon\rightarrow 0.
\end{equation}
Moreover

\begin{equation*}
	\begin{aligned}
		\left(u^{\varepsilon}\right)'(t)=\int_{\mathbb{R}}-\partial_y\left(\eta_{\varepsilon}(t-y)\right)u(y)dy=\int_{\mathbb{R}}\eta_{\varepsilon}(t-y)\frac{du(y)}{dy}dy,
	\end{aligned}
\end{equation*}
which in turn implies
\begin{equation}\label{Sob:for4.65}
	\left(u^{\varepsilon}\right)'\rightarrow\frac{du}{dt} ,\quad \mbox{
		as } \varepsilon\rightarrow 0,\quad\mbox{ in }
	L^p_{loc}(\mathbb{R}).
\end{equation}
Now, by \eqref{Sob:for2.65}--\eqref{Sob:for4.65} we have that
$\left\{u^{\varepsilon}\right\}$ satifies the Cauchy property on every compact set of
$\mathbb{R}$. Therefore $\left\{u^{\varepsilon}\right\}$ uniformly converges to a contiuous function $v$ 
on every compact set of $\mathbb{R}$. Since we have
\begin{equation*}
	u^{\varepsilon}\rightarrow u ,\quad \mbox{ as }
	\varepsilon\rightarrow 0,\quad\mbox{ in } L^p_{loc}(\mathbb{R}),
\end{equation*}
we get
$$u=v,\quad\mbox{a.e. in } \mathbb{R}.$$
On the other hand, by \eqref{Sob:for1.65} and by \eqref{Sob:for4.65} (taking into account that $\left\{u^{\varepsilon}\right\}\rightarrow v$
on any compact), we have
$$v(x)=v(x_0)+\int_{x_0}^x\frac{du(t)}{dt}dt,\quad\forall x\in
\mathbb{R}$$ which implies that $v\in AC_{loc}(\mathbb{R})$ e
\begin{equation}\label{Sob:for1.66}
	v'(x)=\frac{du(x)}{dx},\quad\mbox{ a.e. in } \mathbb{R}.
\end{equation}
Moreover,

$$\frac{1}{2r}\int_{x-r}^{x+r}u(t)dt=\frac{1}{2r}\int_{x-r}^{x+r}v(t)dt,\quad\forall
x\in \mathbb{R}.$$

Hence, passing to the limit as $r\rightarrow 0$ we have

$$u^*(x)=v(x), \quad\forall x\in \mathbb{R}.$$
By the latter and by \eqref{Sob:for1.66} we get

$$\left(u^*\right)'(x)=\frac{du(x)}{dx}, \quad\mbox{ a.e. in } \mathbb{R}.$$

\textbf{2.} Let $u\in L^p_{loc}(\mathbb{R})$ satisfy
$$u=v,\quad\mbox{ a.e. in } \mathbb{R},$$ where  $v\in
AC_{loc}(\mathbb{R})$ and $v'\in L^p_{loc}(\mathbb{R})$. We have

$$\int_{\mathbb{R}}u\Phi'dx=\int_{\mathbb{R}}v\Phi'dx=-\int_{\mathbb{R}}v'\Phi
dx,\quad \forall \Phi\in C^{\infty}_0(\mathbb{R}).$$ Hence, as $v'\in L^p_{loc}(\mathbb{R})$, we get

$$\frac{du}{dx}=v',\quad\mbox{ and }\quad u\in
W^{1,p}_{loc}(\mathbb{R}).$$ $\blacksquare$

\bigskip

\textbf{Remark.} The Extension Theorem implies that if
$1\leq p<+\infty$,  $a,b\in \mathbb{R}$, $a<b$, 
denoting by $E$ the extension operator, then we have

\smallskip

(i') if $u\in W^{1,p}(a,b)$, $(Eu)^*\in AC([a,b])$; in particular
$u$ is almost everywhere equal to an absolutely continuous function in
$[a,b]$ and the weak derivative of $u$ is equal to the classic derivative
of $(Eu)^*$ in $(a,b)$;

\smallskip

(ii') if $u\in L^{p}(a,b)$ and $v\in AC([a,b])$ satisfies
$$u=v,\quad\mbox{ a.e. in } [a,b]$$
and $v'\in L^{p}(a,b)$ then $u\in W^{1,p}(a,b)$ and the weak derivative of $u$ is equal to the classic derivative of $v$. $\blacklozenge$

\bigskip

Theorem \ref{Sob:teo1.6} can be accomplished by the following

\begin{prop}\label{Sob:prop3.6}
	Let $p>1$ and $u\in W_{loc}^{1,p}(a,b)$. We have that $u$ is almost everywhere equal to a function $C^{0,\alpha}_{loc}(\mathbb{R})$, where
	$\alpha=1-1/p$. Here $C^{0,\alpha}_{loc}(\mathbb{R})$ denotes the 
	space of the functions $u$ satisfying $u_{|I}\in
	C^{0,\alpha}(I)$ for every $I$ compact  interval of $\mathbb{R}$.
\end{prop}
\textbf{Proof.} Let $I$ be a bounded interval. By
\eqref{Sob:for1.64}, (by using H\"{o}lder inequality),
we have, for any $x,y\in I$

\begin{equation*}
	\begin{aligned}
		\left|u^{\varepsilon}(x)-u^{\varepsilon}(y)\right|&=\left|\int_{x}^y\left(u^{\varepsilon}\right)'(t)dt\right|\leq\\&\leq
		|x-y|^{1-1/p}\left(\int_I\left|\left(u^{\varepsilon}\right)'(t)\right|^pdt\right)^{1/p}.
	\end{aligned}
\end{equation*}
Hence, passing to the limit as $\varepsilon\rightarrow 0$, taking into account
that $u^{\varepsilon}\rightarrow u^*$ and $u^*(x)=u(x)$ almost everywhere, we obtain

$$\left|u(x)-u(y)\right|\leq
|x-y|^{1-1/p}\left(\int_I\left|\frac{du}{dt}\right|^pdt\right)^{1/p},\quad
\mbox{ a.e. } x,y\in I.$$ $\blacksquare$

\section{The embedding theorems} \label{Sobolev-emb}

In this Section we are interested in proving non
trivial embedding theorems of $W^{k,p}(\Omega)$ in other function spaces. For instance we will be interested in establishing when it happens that
$W^{k,p}(\Omega)\subset L^q(\Omega)$, for $q\neq p$ as well
$W^{k,p}(\Omega)\subset C^{m,\gamma}(\Omega)$ for appropriate $m\in
\mathbb{N}_0$, $0<\gamma\leq 1$.

First we consider the space $W^{1,p}(\Omega)$ and
we distinguish the following three cases

$$ \mbox{(a)}\quad 1\leq p<n,\quad \mbox{(b)} \quad n<p\leq +\infty,
\quad \mbox{(c)}\quad  p=n.$$ About case (c) we will just give brief hints.

\subsection{Case $1\leq p<n$. The Gagliardo -- Nirenberg inequality} \label{Sob:sec7.1}
Let us assume

\begin{equation}\label{Sob:for1.68}
	1\leq p<n
\end{equation}
and let us ask ourselves for what $q\in [1,+\infty]$ can be true an estimate like

\begin{equation}\label{Sob:for2.68}
	\left\Vert u\right\Vert_{L^q(\mathbb{R}^n)}\leq C \left\Vert \nabla
	u\right\Vert_{L^p(\mathbb{R}^n)}, \quad\forall u\in C^{\infty}_0(\mathbb{R}^n),
\end{equation}
where $C$ and $q$ do not depend on $u$.

Let us assume that \eqref{Sob:for2.68} is true  and let us prove that,
\textit{necessarily}

$$q=\frac{np}{n-p}.$$
First, we examine the case $q\in [1,+\infty)$. Let $u\in
C^{\infty}_0(\mathbb{R}^n)$ be \textbf{not identically equal to $0$}, and, for any
$\lambda>0$, let
$$u_{\lambda}(x)=u(\lambda x),\quad \forall x\in \mathbb{R}^n.$$
Of course if \eqref{Sob:for2.68} holds true, then

\begin{equation}\label{Sob:for1.69}
	\left\Vert u_{\lambda}\right\Vert_{L^q(\mathbb{R}^n)}\leq C \left\Vert
	\nabla u_{\lambda}\right\Vert_{L^p(\mathbb{R}^n)}, \quad\forall u\in
	C^{\infty}_0(\mathbb{R}^n) \ \forall \lambda >0.
\end{equation}
Now

\begin{equation*}
	\left\Vert
	u_{\lambda}\right\Vert^q_{L^q(\mathbb{R}^n)}=\int_{\mathbb{R}^n}\left|u(\lambda
	x)\right|^qdx=\lambda^{-n}\int_{\mathbb{R}^n}\left|u(x)\right|^qdx;
\end{equation*}
hence

\begin{equation}\label{Sob:for2.69}
	\left\Vert
	u_{\lambda}\right\Vert_{L^q(\mathbb{R}^n)}=\lambda^{-\frac{n}{q}}\left\Vert
	u\right\Vert_{L^q(\mathbb{R}^n)}
\end{equation}
and
\begin{equation*}
	\left\Vert \nabla
	u_{\lambda}\right\Vert^p_{L^p(\mathbb{R}^n)}=\int_{\mathbb{R}^n}\lambda^p\left|(\nabla
	u)(\lambda
	x)\right|^pdx=\lambda^{p-n}\int_{\mathbb{R}^n}\left|\nabla
	u(x)\right|^pdx;
\end{equation*}
therefore

\begin{equation}\label{Sob:for3.69}
	\left\Vert \nabla
	u_{\lambda}\right\Vert_{L^p(\mathbb{R}^n)}=\lambda^{1-\frac{n}{p}}\left\Vert
	\nabla u\right\Vert_{L^p(\mathbb{R}^n)}.
\end{equation}
Since \eqref{Sob:for2.69} and \eqref{Sob:for3.69} hold true, we may write
\eqref{Sob:for1.69} as follows
\begin{equation}\label{Sob:for4.69}\left\Vert
	u\right\Vert_{L^q(\mathbb{R}^n)}\leq
	C\lambda^{1+\frac{n}{q}-\frac{n}{p}}\left\Vert \nabla
	u\right\Vert_{L^p(\mathbb{R}^n)},\quad\forall
	\lambda>0.\end{equation} Now, if
$$q>\frac{np}{n-p},$$ we get
$$1+\frac{n}{q}-\frac{n}{p}<1+n\left(\frac{n-p}{np}\right)-\frac{n}{p}=0.$$
Consequently

$$\left\Vert
u\right\Vert_{L^q(\mathbb{R}^n)}\leq \lim_{\lambda\rightarrow
	+\infty} C\lambda^{1+\frac{n}{q}-\frac{n}{p}}\left\Vert \nabla
u\right\Vert_{L^p(\mathbb{R}^n)}=0,$$ this is a contradiction because
$u$ does not vanish identically.

On the other hand, if
$$q<\frac{np}{n-p},$$ we have
$$1+\frac{n}{q}-\frac{n}{p}>0,$$ hence
$$\left\Vert u\right\Vert_{L^q(\mathbb{R}^n)}\leq
\lim_{\lambda\rightarrow 0^+}
C\lambda^{1+\frac{n}{q}-\frac{n}{p}}\left\Vert \nabla
u\right\Vert_{L^p(\mathbb{R}^n)}=0,$$ this is again a contradiction.

Finally, if $q=+\infty$, instead of \eqref{Sob:for4.69} we have

\begin{equation*}
	\left\Vert u\right\Vert_{L^{\infty}(\mathbb{R}^n)}\leq
	C\lambda^{1-\frac{n}{p}}\left\Vert \nabla
	u\right\Vert_{L^p(\mathbb{R}^n)},\quad\forall
	\lambda>0.\end{equation*} and, by \eqref{Sob:for1.68}, passing to the limit as $\lambda\rightarrow +\infty$ we have a contradiction.

\bigskip

Here and in the sequel, if $1\leq p< n$, we denote by $p^{\star}$ the number

$$\frac{1}{p^{\star}}=\frac{1}{p}-\frac{1}{n}$$ and we call
$p^{\star}$ the \textbf{Sobolev exponent} or the \textbf{Sobolev conjugate} of $p$ \index{$p^{\star}$}. Let us notice
$$p^{\star}=\frac{pn}{n-p}>p.$$

\bigskip

The Main Theorem of the present Subsection is the following one

\begin{theo}[\textbf{The Gagliardo -- Nirenberg inequality}]\label{Sob:teo1.7}
	\index{Theorem:@{Theorem:}!- Gagliardo -- Nirenberg inequality@{- Gagliardo -- Nirenberg inequality}}
	Let $$1\leq p<n.$$ Then there exists $C$ depending on $p$ and $n$ only
	such that
	
	\begin{equation}\label{Sob:for1.70}\left\Vert
		u\right\Vert_{L^{p^{\star}}(\mathbb{R}^n)}\leq C\left\Vert \nabla
		u\right\Vert_{L^p(\mathbb{R}^n)},\quad\forall u\in
		C^1_0\left(\mathbb{R}^n\right).\end{equation}
\end{theo}

The most challenging part of the proof of Theorem
\ref{Sob:for1.70} concerns the case $p=1$ and this, in turn, is
based on the following

\begin{lem}\label{Sob:lem2.7}
	Let $n\geq 2$ and
	$$g_j:\mathbb{R}^{n-1}\rightarrow [0,+\infty),\quad j=1, \cdots,
	n,$$ be  measurable functions. Then
	
	\begin{equation}\label{Sob:for2.70}
		\begin{aligned}
			\int_{\mathbb{R}^n}\prod_{j=1}^n
			g_j(x_1,\cdots,&x_{j-1},x_{j+1},\cdots, x_n)dx_1\cdots
			dx_n\leq\\&\leq \prod_{j=1}^n
			\left(\int_{\mathbb{R}^{n-1}}g_j^{n-1}(y)dy\right)^{\frac{1}{n-1}}.
		\end{aligned}
	\end{equation}
\end{lem}
\textbf{Proof of Lemma \ref{Sob:lem2.7}.} Let us proceed by
induction on $n$. If $n=2$, we have 
$$\int_{\mathbb{R}^2}g_1(x_2)g_2(x_1)dx_1dx_2=\int_{\mathbb{R}}g_1(x_2)dx_2
\int_{\mathbb{R}}g_2(x_1)dx_1.$$ Therefore, if $n=2$,
\eqref{Sob:for2.70} holds true. Now, let us assume that
\eqref{Sob:for2.70} holds for $n$ and let us prove it for $n+1$. Hence, let us assume
that \textbf{for any} nonnegative measurable functions $g_1,g_2 \cdots$, we have

\begin{equation*}
	\begin{aligned}
		\int_{\mathbb{R}^n}\prod_{j=1}^n g_j dx_1\cdots dx_n\leq
		\prod_{j=1}^n
		\left(\int_{\mathbb{R}^{n-1}}g_j^{n-1}dy\right)^{\frac{1}{n-1}},
	\end{aligned}
\end{equation*}
let us notice that, to shorten the formula, we have omitted the
variables. However, it is important to recall that $g_j$ does not depend
on $x_j$.

By the H\"{o}lder inequality we get

\begin{equation}\label{Sob:for1.71}
	\begin{aligned}
		&\int_{\mathbb{R}^{n+1}}\prod_{j=1}^{n+1} g_jdx_1\cdots
		dx_{n+1}=\\&=\int_{\mathbb{R}}dx_{n+1}\int_{\mathbb{R}^{n}}g_{n+1}\prod_{j=1}^{n}
		g_jdx_1\cdots dx_{n}\leq\\&\leq
		\int_{\mathbb{R}}dx_{n+1}\left(\int_{\mathbb{R}^{n}}g^n_{n+1}dx_1\cdots
		dx_{n}\right)^{\frac{1}{n}}\times \\&\times\left(\int_{\mathbb{R}^{n}}\prod_{j=1}^{n}
		g^{\frac{n}{n-1}}_jdx_1\cdots dx_{n}\right)^{\frac{n-1}{n}}=\\&=
		\left(\int_{\mathbb{R}^{n}}g^n_{n+1}dx_1\cdots
		dx_{n}\right)^{\frac{1}{n}}\times \\&\times\int_{\mathbb{R}}dx_{n+1}\left(\int_{\mathbb{R}^{n}}\prod_{j=1}^{n}
		g^{\frac{n}{n-1}}_jdx_1\cdots dx_{n}\right)^{\frac{n-1}{n}}.
	\end{aligned}
\end{equation}

\medskip

Now let us apply the inductive assumption to the functions
$g_j^{\frac{n}{n-1}}(\cdot,x_{n+1})$, $j=1, \cdots, n$. We get

$$\int_{\mathbb{R}^{n}}\prod_{j=1}^{n}
g^{\frac{n}{n-1}}_jdx_1\cdots dx_{n}\leq
\prod_{j=1}^{n}\left(\int_{\mathbb{R}^{n-1}}g_j^n(y,x_{n+1})dy\right)^{\frac{1}{n-1}}.$$
Hence, we have trivially

$$\left(\int_{\mathbb{R}^{n}}\prod_{j=1}^{n}
g^{\frac{n}{n-1}}_jdx_1\cdots dx_{n}\right)^{\frac{n-1}{n}}\leq
\prod_{j=1}^{n}\left(\int_{\mathbb{R}^{n-1}}g_j^n(y,x_{n+1})dy\right)^{\frac{1}{n}}.$$
The last inequality and \eqref{Sob:for1.71} imply

\begin{equation}\label{Sob:for0.72}
	\begin{aligned}
		\int_{\mathbb{R}^{n+1}}\prod_{j=1}^{n+1} g_jdx_1\cdots dx_{n+1}&\leq
		\left(\int_{\mathbb{R}^{n}}g^n_{n+1}dx_1\cdots
		dx_{n}\right)^{\frac{1}{n}}\times\\&\times\int_{\mathbb{R}}dx_{n+1}\prod_{j=1}^{n}\left(\int_{\mathbb{R}^{n-1}}g_j^n(y,x_{n+1})dy\right)^{\frac{1}{n}}.
	\end{aligned}
\end{equation}
Now, we set
$$h_j\left(x_{n+1}\right)=\left(\int_{\mathbb{R}^{n-1}}g_j^n(y,x_{n+1})dy\right)^{\frac{1}{n}}$$ and we use the extended H\"{o}lder inequality:
\begin{equation*}
	\begin{aligned}
		\int_{\mathbb{R}}\prod_{j=1}^{n}\left(\int_{\mathbb{R}^{n-1}}g_j^n(y,x_{n+1})dy\right)^{\frac{1}{n}}&dx_{n+1}=
		\int_{\mathbb{R}}\prod_{j=1}^{n}
		h_j\left(x_{n+1}\right)dx_{n+1}\leq\\&\leq
		\prod_{j=1}^{n}\left(\int_{\mathbb{R}}h_j^n\left(x_{n+1}\right)dx_{n+1}\right)^{\frac{1}{n}}=\\&=
		\prod_{j=1}^{n}\left(\int_{\mathbb{R}}dx_{n+1}\int_{\mathbb{R}^{n-1}}g_j^n(y,x_{n+1})dy\right)^{\frac{1}{n}}=\\&=
		\prod_{j=1}^{n}\left(\int_{\mathbb{R}^{n}}g_j^n(y)dy\right)^{\frac{1}{n}}.
	\end{aligned}
\end{equation*}
The just obtained inequality and \eqref{Sob:for0.72} yield

\begin{equation*}
	\begin{aligned}
		\int_{\mathbb{R}^{n+1}}\prod_{j=1}^{n+1} g_j dx_1\cdots dx_{n+1}\leq
		\prod_{j=1}^{n+1}
		\left(\int_{\mathbb{R}^n}g_j^{n}dy\right)^{\frac{1}{n}}.
	\end{aligned}
\end{equation*}
Proof of Lemma is concluded. $\blacksquare$

\bigskip

\textbf{Proof of Theorem \ref{Sob:teo1.7}.} Let $u\in
C^1_0\left(\mathbb{R}^n\right)$. For any $j=1,\cdots,n$, we set

\begin{equation}\label{Sob:for1.73}
	\begin{aligned}
	&f_j\left(x_1,\cdots,x_{j-1},x_{j+1},\cdots,x_n\right)=\\&=\int_{\mathbb{R}}\left|\nabla
	u\left(x_1,\cdots,x_{j-1},y_j,x_{j+1},\cdots,x_n\right)\right|dy_j.
	\end{aligned}
	\end{equation}

\smallskip

\noindent We obtain

\begin{equation*}
	|u(x)|\leq f_j\left(x_1,\cdots,x_{j-1},x_{j+1},\cdots,x_n\right),
	\quad j=1,\cdots,n.
\end{equation*}
As a matter of fact, we have

\begin{equation*}
	u(x)=\int_{-\infty}^{x_j}
	\partial_{x_j}u\left(x_1,\cdots,x_{j-1},y_j,x_{j+1},\cdots,x_n\right)dy_j,
	\quad j=1,\cdots,n,
\end{equation*}
from which, for any $j=1,\cdots,n$ we have

\begin{equation}\label{Sob:for2.73}
	\begin{aligned}
		|u(x)|&\leq \int_{\mathbb{R}}\left|\nabla
		u\left(x_1,\cdots,x_{j-1},y_j,x_{j+1},\cdots,x_n\right)\right|dy_j=\\&=f_j\left(x_1,\cdots,x_{j-1},x_{j+1},\cdots,x_n\right).
	\end{aligned}
\end{equation}
Now by multiplying all \eqref{Sob:for2.73} we get

$$|u(x)|^n\leq\prod_{j=1}^nf_j,$$
hence

$$|u(x)|^{\frac{n}{n-1}}\leq\prod_{j=1}^nf^{\frac{1}{n-1}}_j$$ and, by integrating over $\mathbb{R}^n$ we have
$$\int_{\mathbb{R}^n}|u(x)|^{\frac{n}{n-1}}dx\leq
\int_{\mathbb{R}^n}\prod_{j=1}^nf^{\frac{1}{n-1}}_jdx.$$ At this stage
let us exploit Lemma \ref{Sob:lem2.7}. Set 
$$g_j=f^{\frac{1}{n-1}}_j,\quad j=1,\cdots, n$$ and we obtain
\begin{equation}\label{Sob:for1.74}
	\begin{aligned}
		\int_{\mathbb{R}^n}|u(x)|^{\frac{n}{n-1}}dx\leq \prod_{j=1}^n
		\left(\int_{\mathbb{R}^{n-1}}f_j(\eta)d\eta\right)^{\frac{1}{n-1}}.
	\end{aligned}
\end{equation}
Now, we notice that

\begin{equation*}
	\begin{aligned}
		&\int_{\mathbb{R}^{n-1}}f_j(\eta)d\eta=\\&=\int_{\mathbb{R}^{n-1}}d\eta\int_{\mathbb{R}}\left|\nabla
		u\left(\eta_1,\cdots,\eta_{j-1},y_j,\eta_{j+1},\cdots,\eta_n\right)\right|dy_j=\\&=\int_{\mathbb{R}^{n}}
		\left|\nabla u(x)\right|dx,\quad j=1,\cdots,n.
	\end{aligned}
\end{equation*}
By the just obtained equality and by \eqref{Sob:for1.74} we get

$$\int_{\mathbb{R}^n}|u(x)|^{\frac{n}{n-1}}dx\leq \left(\int_{\mathbb{R}^{n}}\left|\nabla
u(x)\right|dx\right)^{\frac{n}{n-1}},$$ which implies

\begin{equation}\label{Sob:for2.74}
	\begin{aligned}
		\left\Vert u \right\Vert_{L^{1^{\star}}(\mathbb{R}^n)}\leq
		\left\Vert \nabla u\right\Vert_{L^{1}(\mathbb{R}^n)}.
	\end{aligned}
\end{equation}
Therefore \eqref{Sob:for1.70} is proved for $p=1$. Now, let
$1<p<n$ and $\alpha>1$ to be chosen. By applying 
\eqref{Sob:for2.74} to $|u|^{\alpha}$ we have

\begin{equation*}
	\begin{aligned}
		\left(\int_{\mathbb{R}^{n}}\left| u(x)\right|^{\frac{\alpha
				n}{n-1}}dx\right)^{\frac{n-1}{n}}&\leq
		\int_{\mathbb{R}^{n}}\left|\nabla\left(\left|
		u(x)\right|^{\alpha}\right)\right|dx=\\&=\int_{\mathbb{R}^{n}}\alpha
		\left| u(x)\right|^{\alpha-1}|\nabla u|dx\leq \\&\leq
		\alpha\left(\int_{\mathbb{R}^{n}}\left|
		u(x)\right|^{\frac{(\alpha-1) p}{p-1}}dx\right)^{\frac{p}{p-1}}
		\left(\int_{\mathbb{R}^{n}}\left| \nabla
		u(x)\right|^{p}dx\right)^{\frac{1}{p}}.
	\end{aligned}
\end{equation*}
Now, let us choose $\alpha$ satisfying

$$\frac{\alpha n}{n-1}=(\alpha-1)\frac{p}{p-1}$$
that is

$$\alpha=\frac{p(n-1)}{n-p},$$
notice that $\alpha>1$, as $n>p>1$. The above choice of
$\alpha$ gives

$$\frac{\alpha n}{n-1}=(\alpha-1)\frac{p}{p-1}=\left(\frac{p(n-1)}{n-p}-1\right)\frac{p}{p-1}=\frac{pn}{n-p}.$$
Hence

$$\left(\int_{\mathbb{R}^{n}}\left| u(x)\right|^{\frac{p
		n}{n-p}}dx\right)^{\frac{n-1}{n}}\leq
\frac{p(n-1)}{n-p}\left(\int_{\mathbb{R}^{n}}\left|
u(x)\right|^{\frac{p
		n}{n-p}}dx\right)^{\frac{p-1}{p}}\left(\int_{\mathbb{R}^{n}}\left|
\nabla u(x)\right|^{p}dx\right)^{\frac{1}{p}}$$ that is

$$\left(\int_{\mathbb{R}^{n}}\left| u(x)\right|^{p^{\star}}dx\right)^{\frac{n-1}{n}-\frac{p-1}{p}}\leq
\frac{p(n-1)}{n-p}\left(\int_{\mathbb{R}^{n}}\left| \nabla
u(x)\right|^{p}dx\right)^{\frac{1}{p}}.$$ On the other hand
$$\frac{n-1}{n}-\frac{p-1}{p}=\frac{1}{p}-\frac{1}{n}=\frac{1}{p^{\star}}.$$
Therefore
$$\left(\int_{\mathbb{R}^{n}}\left| u(x)\right|^{p^{\star}}dx\right)^{\frac{1}{p^{\star}}}\leq
\frac{p(n-1)}{n-p}\left(\int_{\mathbb{R}^{n}}\left| \nabla
u(x)\right|^{p}dx\right)^{\frac{1}{p}}.$$ $\blacksquare$

\bigskip

\begin{theo}[\textbf{The Sobolev inequality}]\label{Sob:teo3.7}
	\index{Theorem:@{Theorem:}!- Sobolev inequality@{- Sobolev inequality}}
	Let $\Omega$ be a bounded open set of $\mathbb{R}^n$ whose boundary is of class $C^{0,1}$ with constants
	$M_0$ and $r_0$. Let $d_0$ be the diameter of $\Omega$.
	Let us assume $1\leq p<n$. 
	
	Then there exists $C$ depending on
	$M_0$, $r_0$, $d_0$, $p$ and $n$ only such that 
	
	\begin{equation}\label{Sob:forSTELLa.75}\left\Vert
		u\right\Vert_{L^{p^{\star}}(\Omega)}\leq C\left\Vert
		u\right\Vert_{W^{1,p}(\Omega)},\quad\forall u\in
		W^{1,p}(\Omega).
	\end{equation}
\end{theo}
\textbf{Proof.} Since $\Omega$ is a bounded open set of
$\mathbb{R}^n$ of class $C^{0,1}$, we can apply the extension Theorem \ref{Sob:teo1.4}. Hence, there exists $\widetilde{u}\in
W^{1,p}(\mathbb{R}^n)$ such that
\begin{equation}\label{Sob:for1.75}
	\widetilde{u}=u, \mbox{ in } \Omega,\quad\quad \mbox{supp }
	\widetilde{u} \mbox{ }\mbox{ compact of } \mathbb{R}^n.
\end{equation}
Moreover, there exists $C$ depending on $M_0$, $r_0$, $d_0$, $p$ and
$n$ only such that

\begin{equation}\label{Sob:for2.75}\left\Vert
	\widetilde{u}\right\Vert_{W^{1,p}(\mathbb{R}^n)}\leq C\left\Vert
	u\right\Vert_{W^{1,p}(\Omega)},\quad\forall u\in W^{1,p}(\Omega).
\end{equation}
Proposition \ref{Sob:prop2.4} implies that there exists a sequence $\left\{v_j\right\}\subset C_0^{\infty}(\mathbb{R}^n)$ such that

$$\left\{v_j\right\}\rightarrow \widetilde{u},\quad\mbox{in }
W^{1,p}(\mathbb{R}^n).$$ Now, by Theorem \ref{Sob:teo1.7} we have

\begin{equation*}
	\left\Vert v_j-v_m\right\Vert_{L^{p^{\star}}(\mathbb{R}^n)}\leq
	C\left\Vert \nabla v_j-\nabla
	v_m\right\Vert_{W^{1,p}(\mathbb{R}^n)},\quad\forall j,m\in
	\mathbb{N},
\end{equation*}
from which it follows that $\left\{v_j\right\}$ is a
Cauchy sequence in $L^{p^{\star}}(\mathbb{R}^n)$, hence there exists  $v\in
L^{p^{\star}}(\Omega)$ such that

$$\left\{v_j\right\}\rightarrow v,\quad\mbox{in }
L^{p^{\star}}(\mathbb{R}^n).$$ Hence $\widetilde{u}=v$; as a matter of fact,
for any $R>0$ we have $L^{p^{\star}}(B_R)\subset L^{p}(B_R)$, consequently
\begin{equation*}
	\left\Vert v-\widetilde{u}\right\Vert_{L^{p}(B_R)}\leq \left\Vert
	v-v_j\right\Vert_{L^{p}(B_R)}+\left\Vert
	v_j-\widetilde{u}\right\Vert_{L^{p}(B_R)}\rightarrow 0,\quad \mbox{as }
	j\rightarrow\infty
\end{equation*}
and, as $R$ is arbitrary, we get $\widetilde{u}=v$.

Therefore $\widetilde{u}\in L^{p^{\star}}(\mathbb{R}^n)$ and passing to the limit in the following inequality 
\begin{equation*}
	\left\Vert v_j\right\Vert_{L^{p^{\star}}(\mathbb{R}^n)}\leq
	C\left\Vert \nabla
	v_j\right\Vert_{W^{1,p}(\mathbb{R}^n)},\quad\forall j\in \mathbb{N},
\end{equation*}
we obtain

\begin{equation}\label{Sob:for1.76}
	\left\Vert
	\widetilde{u}\right\Vert_{L^{p^{\star}}(\mathbb{R}^n)}\leq
	C\left\Vert \nabla \widetilde{u}\right\Vert_{W^{1,p}(\mathbb{R}^n)}.
\end{equation}
On the other hand \eqref{Sob:for1.75} implies

\begin{equation}\label{Sob:for2.76}
	\left\Vert u\right\Vert_{L^{p^{\star}}(\Omega)}\leq \left\Vert
	\widetilde{u}\right\Vert_{L^{p^{\star}}(\mathbb{R}^n)}
\end{equation}
and \eqref{Sob:for2.75} yields

\begin{equation*}\left\Vert
	\nabla\widetilde{u}\right\Vert_{L^{p}(\mathbb{R}^n)}\leq C\left\Vert
	u\right\Vert_{W^{1,p}(\Omega)},\quad\forall u\in W^{1,p}(\Omega).
\end{equation*}
so that, by the latter, by \eqref{Sob:for1.76} and by
\eqref{Sob:for2.76} we get

\begin{equation*}
	\left\Vert u\right\Vert_{L^{p^{\star}}(\Omega)}\leq C\left\Vert
	u\right\Vert_{W^{1,p}(\Omega)}.
\end{equation*}
$\blacksquare$

\bigskip

\begin{cor}\label{Sob:cor4.7}
	Let $\Omega$ be a bounded open set of $\mathbb{R}^n$ and let $d_0$ be its 
	diameter. Let us assume $1\leq p<n$. We have that, if $u\in
	W^{1,p}_0(\Omega)$, then  $u\in L^{p^{\star}}_0(\Omega)$ and there exists
	$C$ depending on $p$, $n$ and $d_0$ only, such that 
		\begin{equation*}
		\left\Vert u\right\Vert_{L^{p^{\star}}(\Omega)}\leq C\left\Vert
		u\right\Vert_{W^{1,p}(\Omega)}.
	\end{equation*}
\end{cor}
\textbf{Proof.} Let $x_0\in \Omega$ and $R=2d_0$.
Since $u\in W^{1,p}_0(\Omega)$ we have that the function

\begin{equation*}
	\widetilde{u}(x)=
	\begin{cases}
		u(x), \mbox{ for } x\in \Omega \\
		\\
		0, \quad\mbox{ for } x\in B_R\setminus \Omega,
	\end{cases}
\end{equation*}
belongs to $W^{1,p}_0(B_R(x_0))$. By the first
Poincar\'{e} inequality (Theorem \ref{Poincar}) and by Theorem \ref{Sob:teo3.7}
we obtain

\begin{equation*}
	\begin{aligned}
		\left\Vert u\right\Vert_{L^{p^{\star}}(\Omega)}&=\left\Vert
		\widetilde{u}\right\Vert_{L^{p^{\star}}(B_R(x_0))}\leq \\&\leq C
		\left\Vert \widetilde{u}\right\Vert_{W^{1,p}(B_R(x_0))}\leq \\&\leq
		C\widetilde{C}\left\Vert
		\nabla\widetilde{u}\right\Vert_{L^{p}(B_R(x_0))}=\\&=C\widetilde{C}\left\Vert
		\nabla u\right\Vert_{L^{p}(\Omega)},
	\end{aligned}
\end{equation*}
where $C$ is the same constant that occurs in the inequality
\eqref{Sob:forSTELLa.75}, hence, it depends on $p$ e da $n$ only,  and
$\widetilde{C}$ is the same constant that occurs in the Poincar\'{e} inequality, hence, it depends on $R$, that is, on $d_0$, only.
$\blacksquare$

\bigskip

\textbf{Counterexample.} If $n>1$ it \textbf{does not} happen that
$W^{1,n}(\Omega)\subset L^{\infty}(\Omega)$. Let us prove what is claimed.

Let
$$u(x)=\log\log\left(1+\frac{1}{|x|}\right),\quad\mbox{in } B_1.$$
We have $u\notin L^{\infty}(B_1)$ and

$$\nabla u(x)=\frac{x}{\left(|x|^3+|x|^2\right)\log
	\left(1+\frac{1}{|x|}\right)}.$$

Let us check that $u\in W^{1,n}(B_1)$. We have

\begin{equation*}
	\begin{aligned}
		\int_{B_1}|u|^ndx&=\int_{B_1}\left|\log\log\left(1+\frac{1}{|x|}\right)\right|^n
		dx=\\&=\omega_n\int^1_0\rho^{n-1}\left|\log\log\left(1+\frac{1}{\rho}\right)\right|^nd\rho<+\infty
	\end{aligned}
\end{equation*}
and (as $n>1$)
\begin{equation*}
	\begin{aligned}
		\int_{B_1}|\nabla
		u|^ndx&=\omega_n\int^1_0\frac{\rho^{n-1}}{\left(\rho^2+\rho\right)^n\left(\log\left(1+\frac{1}{\rho}\right)\right)^n}
		d\rho=\\&=\omega_n\int^1_0\frac{d\rho}{\left(\rho+1\right)^n\rho\left(\log\left(1+\frac{1}{\rho}\right)\right)^n}<+\infty.
	\end{aligned}
\end{equation*}
Therefore $u\in W^{1,n}(B_1)$, but $u\notin L^{\infty}(B_1)$.

\medskip

Let us observe that if $n=1$, for what was proved in Section \ref{Sob:sec6.61}, we have $W^{1,1}(I)\subset L^{\infty}(I)$,
where $I$ is a bounded open interval of $\mathbb{R}$. $\spadesuit$

\subsection{The Sobolev--Poincar\'{e} inequality} \label{Sob:sec7.2}

We begin by the following

\begin{lem}\label{Sob:lem5.7}
	If $p\in[1,+\infty)$, then there exists a constant $C$ depending on  $p$ e da $n$ only, such that 
	\begin{equation}\label{Sob:for1.79}
		\int_{B_r(x)}\left\vert u(y)-u(z)\right\vert^pdy\leq Cr^{n+p-1}
		\int_{B_r(x)}\left\vert \nabla u(y)\right\vert^p||y-z|^{1-n}dy,
	\end{equation}
	for every $u\in C^1\left(\overline{B_r(x)}\right)$, $r>0$, $x\in
	\mathbb{R}^n$ and for every $z\in \overline{B_r(x)}$.
\end{lem}
\textbf{Proof.} It is not restrictive to assume $x=0$. By the Fundamental Theorem of Calculus we have, for any $y,z\in B_r$,
\begin{equation*}
	\begin{aligned}
		u(y)-u(z)=\int^1_0\frac{d}{dt}u(z+t(y-z))dt=\int^1_0\nabla
		u(z+t(y-z))dt\cdot (y-z);
	\end{aligned}
\end{equation*}
which implies
\begin{equation*}
	\begin{aligned}
		\left|u(y)-u(z)\right|^p\leq \left|y-z\right|^p \int^1_0\left|\nabla
		u(z+t(y-z))\right|^pdt.
	\end{aligned}
\end{equation*}
Let $\rho>0$ and let us integrate over
$B_r\cap\partial B_{\rho}(z)$ both the sides  of the last inequality. We get

\begin{equation*}
	\begin{aligned}
		\int_{B_r\cap \partial B_{\rho}(z)}\left|u(y)-u(z)\right|^p&
		dS_y\leq \int_{B_r\cap \partial B_{\rho}(z)}dS_y \left|y-z\right|^p
		\int^1_0\left|\nabla u(z+t(y-z))\right|^pdt=\\&= \rho^p \int^1_0 dt
		\int_{B_r\cap \partial B_{\rho}(z)}\left|\nabla
		u(z+t(y-z))\right|^pdS_y=\\&=\rho^p \int^1_0 dt \int_{B_r\cap
			\partial B_{t\rho}(z)}\left|\nabla u(\xi)\right|^pt^{1-n}dS_{\xi}:=(\bigstar),
	\end{aligned}
\end{equation*}

\smallskip

\noindent where in the last inequality we set $\xi=z+t(y-z)$, so that
$t^{n-1}dS_y=dS_{\xi}$. Now, if $\xi\in\partial B_{t\rho}(z)$, we have
$\frac{|\xi-z|}{\rho}=t$; hence

\begin{equation*}
	\begin{aligned}
		(\bigstar)&=\rho^p \int^1_0 dt \int_{B_r\cap
			\partial B_{t\rho}(z)}\left|\nabla
		u(\xi)\right|^p\frac{|\xi-z|^{1-n}}{\rho^{1-n}}dS_{\xi}=\\&=
		\rho^{n+p-1} \int^1_0 dt \int_{B_r\cap
			\partial B_{t\rho}(z)}\left|\nabla
		u(\xi)\right|^p|\xi-z|^{1-n}dS_{\xi}=\\&= \rho^{n+p-2} \int^{\rho}_0
		d\tau \int_{B_r\cap
			\partial B_{\tau}(z)}\left|\nabla
		u(\xi)\right|^p|\xi-z|^{1-n}dS_{\xi}=\\&= \rho^{n+p-2} \int_{B_r\cap
			B_{\rho}(z)}\left|\nabla u(\xi)\right|^p|\xi-z|^{1-n}d\xi\leq\\&\leq
		\rho^{n+p-2} \int_{B_r}\left|\nabla
		u(\xi)\right|^p|\xi-z|^{1-n}d\xi.
	\end{aligned}
\end{equation*}
Therefore

\begin{equation*}
	\int_{B_r\cap \partial B_{\rho}(z)}\left|u(y)-u(z)\right|^p dS_y\leq
	\rho^{n+p-2} \int_{B_r}\left|\nabla
	u(\xi)\right|^p|\xi-z|^{1-n}d\xi.
\end{equation*}
Now, let us integrate w.r.t. $\rho$ over $[0,r]$ both the sides of the previous inequality. We have

\begin{equation*}
	\begin{aligned}
		\int^r_0d\rho\int_{B_r\cap \partial
			B_{\rho}(z)}\left|u(y)-u(z)\right|^p dS_y&\leq
		\int^r_0\rho^{n+p-2}d\rho \int_{B_r}\left|\nabla
		u(\xi)\right|^p|\xi-z|^{1-n}d\xi=\\&=
		\frac{r^{n+p-1}}{n+p-1}\int_{B_r}\left|\nabla
		u(\xi)\right|^p|\xi-z|^{1-n}d\xi,
	\end{aligned}
\end{equation*}
on the other hand, by using the polar coordinates, we have

\begin{equation*}
	\int^r_0d\rho\int_{B_r\cap \partial
		B_{\rho}(z)}\left|u(y)-u(z)\right|^p
	dS_y=\int_{B_r}\left|u(y)-u(z)\right|^p dy.
\end{equation*}
Therefore, for any $z\in B_r$, we have
\begin{equation*}
	\int_{B_r}\left|u(y)-u(z)\right|^p dy\leq
	\frac{r^{n+p-1}}{n+p-1}\int_{B_r}\left|\nabla
	u(\xi)\right|^p|\xi-z|^{1-n}d\xi
\end{equation*}
from which we obtain \eqref{Sob:for1.79}. $\blacksquare$

\bigskip

In what follows, for any $g\in L^1\left(B_r(x)\right)$, we set

$$(g)_{x,r}=\dashint_{B_r(x)}g(y)dy=\frac{1}{\left|B_r(x)\right|}\int_{B_r(x)}g(y)dy. $$

\bigskip

\begin{theo}[\textbf{The Sobolev--Poincar\'{e} inequality}]\label{Sob:teo6.7}
	\index{Theorem:@{Theorem:}!- Sobolev--Poincar\'{e} inequality@{- Sobolev--Poincar\'{e} inequality}}
	 Let $1\leq p<n$. Then there exists a constant $C$ depending on $p$ and $n$ only such that
	\begin{equation}\label{Sob:for1.82}
		\left(\dashint_{B_r(x)}\left|u(y)-(u)_{x,r}\right|^{p\star}
		dy\right)^{\frac{1}{p^{\star}}}\leq C
		r\left(\dashint_{B_r(x)}\left|\nabla
		u(y)\right|^pdy\right)^{\frac{1}{p}},
	\end{equation}
	for every $u\in W^{1,p}\left(B_r(x)\right)$.
\end{theo}
\textbf{Proof.} First, we assume $x=0$ and $r=1$. Let $u\in
C^1\left(\overline{B_1}\right)$. Lemma \ref{Sob:lem5.7} and H\"{o}lder inequality give

\begin{equation*}
	\begin{aligned}
		\dashint_{B_1}\left|u(y)-\dashint_{B_1}u(z)dz\right|^{p}
		dy&=\dashint_{B_1}\left|\dashint_{B_1}(u(y)-u(z))dz\right|^{p}
		dy\leq\\&\leq
		\dashint_{B_1}dy\dashint_{B_1}\left|u(y)-u(z)\right|^{p}dz\leq\\&\leq
		C\dashint_{B_1}dy\dashint_{B_1}\left|\nabla
		u(z)\right|^{p}\left|y-z\right|^{1-n}dz=\\&=
		C\dashint_{B_1}dz\left|\nabla
		u(z)\right|^{p}\dashint_{B_1}\left|y-z\right|^{1-n}dy\leq\\&\leq
		C_1\dashint_{B_1}dz\left|\nabla
		u(z)\right|^{p}\int_{B_2(z)}\left|y-z\right|^{1-n}dy=\\&=
		C_2\dashint_{B_1}\left|\nabla u(z)\right|^{p}dz,
	\end{aligned}
\end{equation*}
where $C_1=\frac{\omega_n}{n}C$, $C_2=2\omega_n C_1$. Hence, we have obtained the inequality 

\begin{equation*}
	\dashint_{B_1}\left|u(y)-(u)_{0,1}\right|^{p} dy\leq
	C_2\dashint_{B_1}\left|\nabla u(z)\right|^{p}dz.
\end{equation*}
Now, let us apply the Sobolev inequality (Theorem
\ref{Sob:teo3.7}) to $u-(u)_{0,1}$. We get

\begin{equation*}
	\begin{aligned}
		\left(\dashint_{B_1}\left|u(y)-(u)_{0,1}\right|^{p\star}
		dy\right)^{\frac{1}{p\star}}&\leq
		C_3\left[\dashint_{B_1}\left|u(y)-(u)_{0,1}\right|^{p\star}
		dy+\dashint_{B_1}\left|\nabla
		u(y)\right|^{p}dy\right]^{\frac{1}{p}}\leq\\&\leq
		C_4\left(\dashint_{B_1}\left|\nabla
		u(y)\right|^{p}dy\right)^{\frac{1}{p}}.
	\end{aligned}
\end{equation*}
Thus, we have obtained
\begin{equation}\label{Sob:for1.83}
	\left(\dashint_{B_1}\left|u(y)-(u)_{0,1}\right|^{p\star}
	dy\right)^{\frac{1}{p\star}}\leq C_4
	\left(\dashint_{B_1}\left|\nabla
	u(y)\right|^{p}dy\right)^{\frac{1}{p}},
\end{equation}
where $C_4$ depend on $n$ and $p$ only.

Let now $u\in C^1\left(\overline{B_r}(0)\right)$. Set

$$v(x)=u(rx),\quad x\in B_1$$ and apply \eqref{Sob:for1.83} to function
$v$. We have
\begin{equation}\label{Sob:for1.84}
	(v)_{0,1}=(u)_{0,r},
\end{equation}

\begin{equation}\label{Sob:for2.84}
	\dashint_{B_1}\left|\nabla
	v(y)\right|^{p}dy=r^p\dashint_{B_r}\left|\nabla u(x)\right|^{p}dx
\end{equation}
and
\begin{equation}\label{Sob:for3.84}
	\dashint_{B_1}\left|v(y)-(v)_{0,1}\right|^{p\star}
	dy=\dashint_{B_r}\left|u(x)-(u)_{0,r}\right|^{p\star} dx.
\end{equation}
By applying inequality \eqref{Sob:for1.83} to $v$ and taking into account
\eqref{Sob:for1.84}--\eqref{Sob:for3.84} we obtain
\eqref{Sob:for1.82} by density. $\blacksquare$

\subsection{The Morrey inequality} \label{Sob:sec7.3}

Let $E$ be any measurable of $\mathbb{R}^n$, here and in the sequel we say that $u^*:E \rightarrow\mathbb{R}$ is a \textit{version}
of a given function $u:E \rightarrow\mathbb{R}$ if

$$u=u^*,\quad \mbox{ a.e. in } E.$$

\begin{lem}\label{Sob:lem7.7}
	If $n<p\leq +\infty$, then there exists a constant $C$ depending on $p$ and $n$ only, such that we have
	\begin{equation}\label{Sob:for1.85}
		\left\vert u(y)-u(z)\right\vert\leq Cr^{1-\frac{n}{p}}
		\int_{B_r(x)}\left\vert \nabla u(\xi)\right\vert^p|d\xi,
	\end{equation}
	for every $u\in C^1\left(\overline{B_r(x)}\right)$, $r>0$,  $x\in
	\mathbb{R}^n$ and for every $y,z\in
	\overline{B_r(x)}$.
\end{lem}
\textbf{Proof.} Let $u\in
C^1\left(\overline{B_r(x)}\right)$. Let us apply Lemma
\ref{Sob:lem5.7} for $p=1$. We have, for any $y,z\in \overline{B_r(x)}$,

\begin{equation}\label{Sob:for2.85}
	\begin{aligned}
		\left\vert
		u(y)-u(z)\right\vert&=\dashint_{B_r}\left|u(y)-u(z)\right|
		d\xi\leq\\&\leq \dashint_{B_r}\left|u(y)-u(\xi)\right| d\xi+
		\dashint_{B_r}\left|u(z)-u(\xi)\right| d\xi\leq\\&\leq
		C\int_{B_r(x)}\left(|y-\xi|^{1-n}+|z-\xi|^{1-n}\right)\left\vert
		\nabla u(\xi)\right\vert d\xi,
	\end{aligned}
\end{equation}
where $C$ depends on $n$ only. Let $n<p<+\infty$, by the H\"{o}lder inequality we have

\begin{equation}\label{Sob:forstella.85}
	\begin{aligned}
		&\int_{B_r(x)}\left(|y-\xi|^{1-n}+|z-\xi|^{1-n}\right)\left\vert
		\nabla u(\xi)\right\vert d\xi\leq\\&\leq
		\left[\int_{B_r(x)}\left(|y-\xi|^{1-n}+|z-\xi|^{1-n}\right)d\xi\right]^{\frac{p-1}{p}}\left(\int_{B_r(x)}\left\vert
		\nabla u(\xi)\right\vert^p d\xi\right)^{\frac{1}{p}}\leq\\&\leq
		2^{\frac{1}{p-1}} \underset{I
		}{\left[\underbrace{\int_{B_r(x)}\left(|y-\xi|^{-\frac{(n-1)p}{p-1}}+|z-\xi|^{-\frac{(n-1)p}{p-1}}\right)d\xi}\right]}^{\frac{p-1}{p}}\left(\int_{B_r(x)}\left\vert
		\nabla u(\xi)\right\vert^p d\xi\right)^{\frac{1}{p}}.
	\end{aligned}
\end{equation}

\smallskip

Now, let us check that

\begin{equation}\label{Sob:for1.86}
	I\leq
	C2^{\frac{p-n}{p-1}}\left(\frac{p-1}{p-n}\right)r^{\frac{p-n}{p-1}},
\end{equation}
where $C$ depends on $n$ only. We have $B_r(x)\subset B_{2r}(y)$, for any $y\in B_r(x)$. Hence (taking into account that $p>n$
implies $\frac{(n-1)p}{p-1}<n$)

\begin{equation*}
	\begin{aligned}
		\int_{B_r(x)}|y-\xi|^{-\frac{(n-1)p}{p-1}}d\xi&\leq
		\int_{B_{2r}(y)}|y-\xi|^{-\frac{(n-1)p}{p-1}}d\xi=\\&=
		\int_{\partial B_1}d S
		\int^{2r}_0\rho^{-\frac{(n-1)p}{p-1}}\rho^{n-1}d\rho=\\&=
		\omega_n2^{\frac{p-n}{p-1}}\left(\frac{p-1}{p-n}\right)r^{\frac{p-n}{p-1}}.
	\end{aligned}
\end{equation*}
Since a similar  estimate holds true for the integral
$$\int_{B_r(x)}|z-\xi|^{-\frac{(n-1)p}{p-1}}d\xi$$
provided $z\in B_r(x)$, we get \eqref{Sob:for1.86}. By estimate \eqref{Sob:for1.86} and by \eqref{Sob:for2.85} we obtain

\begin{equation}\label{Sob:for1.87}
	\left\vert u(y)-u(z)\right\vert\leq
	C_n2^{\frac{1}{p-1}}2^{1-\frac{n}{p}}r^{1-\frac{n}{p}}\left(\int_{B_r(x)}\left\vert
	\nabla u(\xi)\right\vert^p d\xi\right)^{\frac{1}{p}}.
\end{equation}
The last obtained estimate implies \eqref{Sob:for1.85} for every $n<p<+\infty$.
If $p=+\infty$, we may pass to the limit as $p\rightarrow +\infty$ in \eqref{Sob:for1.87}. $\blacksquare$

\bigskip

\begin{lem}\label{Sob:lem8.7}
	If $n<p\leq +\infty$, then there exists a constant $C$ depending on $p$ and $n$ such that
	\begin{equation}\label{Sob:for2.87}
		\left\vert u(x)\right\vert\leq C\left\Vert \nabla
		u\right\Vert_{W^{1,p}(B_1(x)},
	\end{equation}
	for every $u\in C^1\left(\overline{B_1(x)}\right)$, $r>0$ and for every $x\in
	\mathbb{R}^n$.
\end{lem}
\textbf{Proof.} By \eqref{Sob:for1.85} ($y=x$, $r=1$ and
$z\in B_1(x)$) we get
$$\left\vert u(x)\right\vert\leq C\left\Vert
\nabla u\right\Vert_{L^p(B_1(x))}+\left\vert u(z)\right\vert.$$ Now, by
integrating both the side w.r.t. $z$ on $B_1(x)$ and by applying the H\"{o}lder inequality, we have

\begin{equation*}
	\begin{aligned}
		|B_1(x)|\left\vert u(x)\right\vert&\leq C\left\Vert
		\nabla u\right\Vert_{L^p(B_1(x))}+\int_{B_1(x)}|u(z)|dz\leq \\&\leq
		C'\left(\left\Vert \nabla u\right\Vert_{L^p(B_1(x))}+\left\Vert
		u\right\Vert_{L^p(B_1(x))}\right)\leq\\&\leq C''\left\Vert
		u\right\Vert_{W^{1,p}(B_1(x))}.
	\end{aligned}
\end{equation*}
Which gives \eqref{Sob:for2.87}.
$\blacksquare$

\bigskip

\begin{theo}[\textbf{the Morrey inequality}]\label{Sob:teo9.7} 
	\index{Theorem:@{Theorem:}!- Morrey inequality@{- Morrey inequality}}
    Let $n<p<+\infty$ and let $\Omega$ be a bounded open set of $\mathbb{R}^n$
	whose baundary is of class $C^{0,1}$ with constants $M_0,r_0$. Then there exists a constant $C$, depending on $p$, $n$, $M_0$ and $r_0$ only, such that
	for every $u\in W^{1,p}(\Omega)$ there exists a version of $u$, $u^*\in
	C^{0,\gamma}\left(\overline{\Omega}\right)$, where
	$$\gamma=1-\frac{n}{p}.$$ Moreover
	\begin{equation}\label{Sob:forStella.88}
		\left\Vert
		u^*\right\Vert_{C^{0,\gamma}\left(\overline{\Omega}\right)}\leq C
		\left\Vert u\right\Vert_{W^{1,p}\left(\Omega\right)}.
	\end{equation}
\end{theo}
\textbf{Proof.} Let us begin by proving that

\begin{equation}\label{Sob:for1.88}
	\left\Vert v\right\Vert_{C^{0,\gamma}\left(\mathbb{R}^n\right)}\leq
	C \left\Vert v\right\Vert_{W^{1,p}\left(\mathbb{R}^n\right)},\quad
	\forall v\in C_0^{\infty}(\mathbb{R}^n),
\end{equation}
where $C$ depends on $p$ and $n$ only. Indeed, by
\eqref{Sob:for2.87} we get trivially
\begin{equation}\label{Sob:for2.88}
	\left\Vert v\right\Vert_{L^{\infty}\left(\mathbb{R}^n\right)}\leq C
	\left\Vert v\right\Vert_{W^{1,p}\left(\mathbb{R}^n\right)},\quad
	\forall v\in C_0^{\infty}(\mathbb{R}^n)
\end{equation}
Now, we set  $$x=\frac{y+z}{2},\quad \mbox{and}\quad r=|y-z|,$$
and by \eqref{Sob:for1.85} we get, for any that $y,z\in \mathbb{R}^n$, 

\begin{equation}\label{Sob:for3.88}
	\begin{aligned}
		\left\vert v(y)-v(z)\right\vert\leq C \left\vert
		y-z\right\vert^{\gamma}\left(\int_{B_r(x)}|\nabla
		v(\xi)|^pd\xi\right)^{\frac{1}{p}}\leq C\left\vert
		y-z\right\vert^{\gamma}\left\Vert
		v\right\Vert_{W^{1,p}\left(\mathbb{R}^n\right)}.
	\end{aligned}
\end{equation}

\smallskip

\noindent Hence \eqref{Sob:for2.88} and \eqref{Sob:for3.88} give
\eqref{Sob:for1.88}.

Now, as $\partial \Omega$ is of class $C^{0,1}$ with constants
$M_0,r_0$, by extension Theorem  \ref{Sob:teo1.4} we have that if
$u\in W^{1,p}(\Omega)$ there exists $\widetilde{u}\in
W^{1,p}(\mathbb{R}^n)$ such that

\begin{equation}\label{Sob:for1.89}
	\begin{cases}
		\widetilde{u}(x)=
		u(x), \mbox{ for } x\in \Omega, \\
		\\
		\mbox{supp }\widetilde{u}, \mbox{ is a compact of } \mathbb{R}^n,\\
		\\
		\left\Vert
		\widetilde{u}\right\Vert_{W^{1,p}\left(\mathbb{R}^n\right)}\leq
		C\left\Vert u\right\Vert_{W^{1,p}\left(\Omega\right)}.
	\end{cases}
\end{equation}

\smallskip

\noindent By Proposition \ref{Sob:prop2.4} we derive that there exists a sequence 
$\left\{v_j\right\}$ in $C^{\infty}_0(\mathbb{R}^n)$ such that
\begin{equation}\label{Sob:for2.89}
	\left\{v_j\right\}\rightarrow \widetilde{u},\quad\mbox{in }
	W^{1,p}\left(\mathbb{R}^n\right)
\end{equation}
and \eqref{Sob:for1.88} implies

$$\left\Vert
v_j-v_k\right\Vert_{C^{0,\gamma}\left(\mathbb{R}^n\right)}\leq
C \left\Vert v_j-v_k\right\Vert_{W^{1,p}\left(\mathbb{R}^n\right)}$$
for every $j,k\in \mathbb{N}$.
Hence $\left\{v_j\right\}$ is a Cauchy sequence in
$C^{0,\gamma}\left(\mathbb{R}^n\right)$ and therefore there exists $u^*\in
C^{0,\gamma}\left(\mathbb{R}^n\right)$ such that

$$v_j\rightarrow u^*,\quad\mbox{ as } j\rightarrow\infty, \mbox{ in
} C^{0,\gamma}\left(\mathbb{R}^n\right).$$ By the latter and by
\eqref{Sob:for1.89}, \eqref{Sob:for2.89} we obtain

$$u^*_{|\Omega}=\widetilde{u}_{|\Omega}=u,\quad\mbox{ a.e. in }
\Omega.$$ Since \eqref{Sob:for1.88} yields

$$\left\Vert
v_j\right\Vert_{C^{0,\gamma}\left(\mathbb{R}^n\right)}\leq
C \left\Vert
v_j\right\Vert_{W^{1,p}\left(\mathbb{R}^n\right)},\quad\forall j\in
\mathbb{N},$$ passing to the limit, we have

$$\left\Vert
u^*\right\Vert_{C^{0,\gamma}\left(\overline{\Omega}\right)}\leq
C\left\Vert u^*\right\Vert_{W^{1,p}\left(\mathbb{R}^n\right)}\leq C
\left\Vert
\widetilde{u}\right\Vert_{W^{1,p}\left(\mathbb{R}^n\right)}\leq C' C
\left\Vert u\right\Vert_{W^{1,p}\left(\Omega\right)},$$ which concludes the proof. $\blacksquare$

\subsection{The General Sobolev inequalities} \label{Sob:sec7.4}

By Theorems \ref{Sob:teo3.7} \ref{Sob:teo9.7}, proceeding by
iteration we obtain the following general theorem, whose
proof we leave to the reader.

\begin{theo} [\textbf{Sobolev embedding}]\label{Sobolev-inequ}
	\index{Theorem:@{Theorem:}!- Sobolev embedding@{- Sobolev embedding}}
	Let $\Omega$ be a bounded open set of class $C^{0,1}$ with constants $M_0, r_0$ and let
	$u\in W^{k,p}(\Omega)$.
	
	(i) If
	\begin{equation}\label{Sob-1}
		k<\frac{n}{p},
	\end{equation}
	then $u\in L^{q}(\Omega)$, where
	\begin{equation}\label{Sob-2}
		\frac{1}{q}=\frac{1}{p}-\frac{k}{n}.
	\end{equation}
	Moreover
	\begin{equation}\label{stimaSobolev-1}
		\left\Vert u\right\Vert_{L^{q}(\Omega)}\leq C\left\Vert u
		\right\Vert_{W^{k,p}(\Omega)},
	\end{equation}
	where $C$ depends on $M_0,r_0$, $k$ and $n$ only.
	
	(ii) If
	\begin{equation}\label{Sob-11}
		k>\frac{n}{p},
	\end{equation}
	then $u\in C^{m,\alpha}\left(\overline{\Omega}\right)$, where
	$m=k-[\frac{n}{p}]-1$ and
	\begin{equation}\label{Sob-20}
		\alpha=
		\begin{cases}
			[\frac{n}{p}]+1-\frac{n}{p}, \  \mbox{ if } \frac{n}{p} \mbox{ is not an integer number} , \\
			\\
			\mbox{any positive number, if } \alpha<1 \mbox{ and } \frac{n}{p} \mbox{ is an integer number}%
		\end{cases}%
	\end{equation}
	and
	\begin{equation}\label{stimaSobolev-12}
		\left\Vert u\right\Vert_{C^{m,\alpha}(\overline{\Omega})}\leq
		C\left\Vert u \right\Vert_{W^{k.p}(\Omega)},
	\end{equation}
	where $C$ depends on $M_0,r_0$, $k$ and $n$ only.
\end{theo}

\textbf{Examples.}

If $n=1$ and $u\in H^1(0,1)$, then $u\in C^{0,1/2}([0,1])$. If $n=2$
and $u\in H^1(\Omega)$, then  $u\in L^{q}(\Omega)$ for every $1\leq
q<\infty$ and,  if $u\in H^2(\Omega)$ then $u\in
C^{0,\alpha}(\overline{\Omega})$ for every $\alpha<1$. Finally, if
$u\in H^k(\Omega)$ for every $k\in \mathbb{N}$, then $u\in
C^{\infty}(\overline{\Omega})$. $\spadesuit$

\section{The compactness theorems} \label{teoremi-compatt}
In the previous Section we have proved that if $\Omega$ is a bounded open set of class $C^{0,1}$ and  $1\leq p<n$ then
$W^{1,p}(\Omega)\subset L^{p^{\star}}(\Omega)$. Moreover the embedding

$$W^{1,p}(\Omega) \hookrightarrow L^{p^{\star}}(\Omega),$$
is continuous, as inequality
\eqref{Sob:forSTELLa.75} holds true. Similarly, (Theorem \ref{Sob:teo9.7}),
for $n<p<+\infty$, the embedding

$$W^{1,p}(\Omega) \hookrightarrow C^{0,1-\frac{n}{p}}\left(\overline{\Omega}\right),$$
is continuous.
In this Section we will prove compact embedding theorems,
in particular, the Rellich -- Kondrachov Theorem, which gives
the compactness of the embedding

$$W^{1,p}(\Omega) \hookrightarrow L^{q}(\Omega),$$
for $1\leq p<n$ and $q<p^{\star}$. This means that any
bounded subset $Y$ di $W^{1,p}(\Omega)$ is relatively
compact in $L^{q}(\Omega)$ (namely, $\overline{Y}$ is compact
in $L^{q}(\Omega)$).

\begin{theo}[\textbf{Rellich -- Kondrachov}]\label{Sob:teo2.8}
	\index{Theorem:@{Theorem:}!- Rellich -- Kondrachov@{- Rellich -- Kondrachov}}
	Let $\Omega$ be a bounded open set of $\mathbb{R}^n$ with boundary of class $C^{0,1}$ and let $1\leq p<n$, $1\leq
	q<p^{\star}=\frac{np}{n-p}$. Then the embedding of
	$W^{1,p}(\Omega)$ in $L^{q}(\Omega)$ is compact.
\end{theo}

\bigskip

In order to prove Theorem \ref{Sob:teo2.8} we need the following.

\begin{lem}\label{Sob:lem1.8}
	Let $1\leq q<+\infty$, let $\Omega$ be a bounded open set of
	$\mathbb{R}^n$ and let $\Lambda$ be the subset of $L^q(\Omega)$
	defined as follows
	$$\Lambda=\left\{u\in L^q(\Omega):\mbox{ } \left\Vert u\right\Vert_{L^q(\Omega)}\leq 1  \right\}.$$
	Let us assume 
	
	\begin{equation}\label{Sob:for1.92}
		\lim_{\varepsilon\rightarrow 0}\left(\sup_{u\in \Lambda}\left\Vert
		u_{\varepsilon}-u\right\Vert_{L^{q}(\Omega)}\right)=0,
	\end{equation}
	where
	$$u_{\varepsilon}(x)=\int_{\Omega}\eta_{\varepsilon}(x-y)u(y)dy$$
	and
	$\eta_{\varepsilon}=\varepsilon^{-n}\eta\left(\varepsilon^{-1}x\right)$
	where $\eta\in C^{\infty}_0(\mathbb{R}^n)$, supp $\eta \subset
	B_1$, $\int_{\mathbb{R}^n}\eta(x)dx=1$.
	
	Then $\Lambda$ is relatively compact in $L^{q}(\Omega)$.
\end{lem}

\textbf{Proof of Lemma \ref{Sob:lem1.8}.} We prove that 
$\Lambda$ is a totally bounded set in $L^{q}(\Omega)$.

Let $\delta>0$. By \eqref{Sob:for1.92} we have that there exists
$\varepsilon_0>0$ so that
\begin{equation}\label{Sob:for1.93}
	\left\Vert
	u_{\varepsilon_0}-u\right\Vert_{L^{q}(\Omega)}<\frac{\delta}{2},\quad\forall
	u\in \Lambda.
\end{equation}
Let

$$\Lambda_0=\left\{u_{\varepsilon_0}:\mbox{ } u\in \Lambda
\right\}$$ Now we prove that $\Lambda_0$ is relatively compact
in $C^0\left(\overline{\Omega}\right)$.

Let us denote

$$M_0=\sup_{\mathbb{R}^n} \left|\eta_{\varepsilon_0}\right|,\quad M_1=\sup_{\mathbb{R}^n}
\left|\nabla\eta_{\varepsilon_0}\right|.$$ We have, for any $u\in
\Lambda$,

\begin{equation*}
	\begin{aligned}
		\left\vert u_{\varepsilon_0}(x)\right\vert&=\left\vert
		\int_{\Omega}\eta_{\varepsilon_0}(x-y)u(y)dy\right\vert\leq \\&\leq
		M_0|\Omega|^{1-\frac{1}{q}}\left\Vert u\right\Vert_{L^q(\Omega)}\leq
		\\&\leq
		M_0|\Omega|^{1-\frac{1}{q}}
	\end{aligned}
\end{equation*}
and, similarly,
\begin{equation*}
	\begin{aligned}
		\left\vert \nabla u_{\varepsilon_0}(x)\right\vert&=\left\vert
		\int_{\Omega}\nabla\eta_{\varepsilon_0}(x-y)u(y)dy\right\vert\leq
		M_1|\Omega|^{1-\frac{1}{q}}.
	\end{aligned}
\end{equation*}
Therefore $\Lambda_0$ is equibounded and equicontinuous. Hence,  the Arzel\`{a}--Ascoli Theorem implies that
$\Lambda_0$ is relatively compact in
$C^0\left(\overline{\Omega}\right)$. 

Now we prove that
$\Lambda_0$ is relatively compact in $L^q(\Omega)$. The inequality

$$\left\Vert w\right\Vert_{L^q(\Omega)}\leq |\Omega|^{\frac{1}{q}}  \left\Vert
w\right\Vert_{C^0\left(\overline{\Omega}\right)},\quad \forall w\in
C^0\left(\overline{\Omega}\right),$$ implies that, for any $w\in
C^0\left(\overline{\Omega}\right)\subset L^q(\Omega)$ and for any
$r>0$

$$B^{C^0}_{r'}(w)\subset B^{L^q}_{r}(w),$$ where $r'=r|\Omega|^{-\frac{1}{q}}$,
$B^{C^0}_{r'}(w)$ is the open ball of
$C^0\left(\overline{\Omega}\right)$ centered at $w$ with radius $r'$ and
$B^{L^q}_{r}(w)$ is the open ball of $L^q(\Omega)$ centered at $w$ with radius $r$.
Since $\Lambda_0$ is relatively compact in 
$C^0\left(\overline{\Omega}\right)$, there exist $w_1,\cdots,
w_{N_r}\in C^0\left(\overline{\Omega}\right)$ such that

$$\Lambda_0\subset \bigcup_{j=1}^{N_r}
B^{C^0}_{r'}\left(w_j\right)\subset \bigcup_{j=1}^{N_r}
B^{L^q}_{r}\left(w_j\right).$$ All in all, $\Lambda_0$ is totally bounded
set of $L^q(\Omega)$. Hence

$$\Lambda_0\subset \bigcup_{j=1}^{N}
B^{L^q}_{\delta/2}\left(w_j\right),$$ where $N$ depends by
$\delta>0$. Consequently, if $u\in \Lambda$, then there exists $j_u\in
\left\{1,\cdots, N\right\}$ so that

$$\left\Vert
u_{\varepsilon_0}-w_{j_u}\right\Vert_{L^q(\Omega)}<\frac{\delta}{2}.$$
By this inequality and by  \eqref{Sob:for1.93} we derive that, if $u\in
\Lambda$ then
$$\left\Vert
u-w_{j_u}\right\Vert_{L^q(\Omega)}<\delta.$$ Hence
$$\Lambda\subset \bigcup_{j=1}^{N}
B^{L^q}_{\delta}\left(w_j\right),$$ which implies compactness of
$\Lambda$. $\blacksquare$

\bigskip

\textbf{Proof of Theorem \ref{Sob:teo2.8}.} Let us apply
Lemma \ref{Sob:lem1.8}. Set

\begin{equation}\label{Sob:for1.95}
	\Lambda=\left\{u\in W^{1,p}(\Omega):\mbox{ } \left\Vert
	u\right\Vert_{W^{1,p}(\Omega)}\leq 1  \right\}.
\end{equation}
We begin by proving the Theorem for $q=1$. Hence, let us prove
that

\begin{equation}\label{Sob:for2.95}
	\lim_{\varepsilon\rightarrow 0}\left(\sup_{u\in \Lambda}\left\Vert
	u_{\varepsilon}-u\right\Vert_{L^{1}(\Omega)}\right)=0,
\end{equation}
where

$$u_{\varepsilon}=\left(\overline{u}\star
\eta_{\varepsilon}\right),$$ being $\overline{u}$ the extension of
$u$ to $0$ in $\mathbb{R}^n\setminus \Omega$.

Let $\delta>0$ and let $\widetilde{\Omega}\Subset \Omega$ satisfy

\begin{equation}\label{Sob:for3.95}
	\left|\Omega\setminus
	\widetilde{\Omega}\right|<\delta^{\frac{p^{\star}}{p^{\star}-1}}.
\end{equation}
We have
\begin{equation}\label{Sob:for4.95}
	\int_{\Omega}\left|u_{\varepsilon}(x)-u(x)\right|dx=\int_{\Omega\setminus
		\widetilde{\Omega}}\left|u_{\varepsilon}(x)-u(x)\right|dx+\int_{
		\widetilde{\Omega}}\left|u_{\varepsilon}(x)-u(x)\right|dx.
\end{equation}
Now, by the H\"{o}lder inequality we derive 

\begin{equation}\label{Sob:for1.96}
	\begin{aligned}
		\int_{\Omega\setminus
			\widetilde{\Omega}}\left|u_{\varepsilon}(x)-u(x)\right|dx&\leq
		\left|\Omega\setminus
		\widetilde{\Omega}\right|^{1-\frac{1}{p^{\star}}}\left\Vert
		u_{\varepsilon}-u\right\Vert_{L^{p^{\star}}(\Omega)}\leq\\&\leq
		\delta\left(\left\Vert
		u_{\varepsilon}\right\Vert_{L^{p^{\star}}(\Omega)}+\left\Vert
		u\right\Vert_{L^{p^{\star}}(\Omega)}\right),
	\end{aligned}
\end{equation}
On the other hand by the Young inequality for convolutions, we have
\begin{equation}\label{Sob:for2.96}
	\left\Vert u_{\varepsilon}\right\Vert_{L^{p^{\star}}(\Omega)}\leq
	\left\Vert u\right\Vert_{L^{p^{\star}}(\Omega)}.
\end{equation}
Now, Theorem \ref{Sob:teo3.7} gives

\begin{equation}\label{Sob:for3.96}
	\left\Vert u\right\Vert_{L^{p^{\star}}(\Omega)}\leq C_1\left\Vert
	u\right\Vert_{W^{1,p}(\Omega)}\leq C_1, \quad\forall u\in \Lambda.
\end{equation}
Therefore by \eqref{Sob:for1.96} -- \eqref{Sob:for3.96} we obtain

\begin{equation}\label{Sob:for4.96}
	\int_{ \Omega\setminus\widetilde{\Omega}}\left|u_{\varepsilon}(x)-u(x)\right|dx\leq
	2C_1\delta,\quad\forall u\in \Lambda.
\end{equation}
Now, let us consider  second addend on the right--hand side of \eqref{Sob:for4.95}. Let
$\overline{\varepsilon}=$dist$\left(\widetilde{\Omega},\partial
\Omega\right)$. For any $x\in \widetilde{\Omega}$ and for any
$\varepsilon\in \left(0,\overline{\varepsilon}\right)$, we have

$$u_{\varepsilon}(x)=\int_{\Omega}\eta_{\varepsilon}(x-y)u(y)dy=\int_{\Omega}\eta(\xi)u(x-\varepsilon
\xi)d\xi.$$ Hence

\begin{equation}\label{Sob:for5.96}
	\begin{aligned}
		\int_{\widetilde{\Omega}}\left|u(x-\varepsilon\zeta)-u(x)\right|dx&=\int_{\widetilde{\Omega}}dx\left|\int_{\Omega}\eta(\xi)(u(x-\varepsilon
		\xi)-u(x))d\xi\right|\leq\\&\leq \int_{\Omega}d\xi\int_{\widetilde{\Omega}}\eta(\xi)\left|u(x-\varepsilon
		\xi)-u(x)\right|dx=\\&=\int_{\Omega}\eta(\xi)d\xi\int_{\widetilde{\Omega}}\left|u(x-\varepsilon \xi)-u(x)\right|dx.
	\end{aligned}
\end{equation}
Let now $y\in B_1$ and $\varepsilon<\overline{\varepsilon}$,
by applying Theorem \ref{densit 2} we have, for almost every $x\in
\widetilde{\Omega}$ and for every $\xi\in B_1$,

\begin{equation*}
	\begin{aligned}
		\left|u(x-\varepsilon \xi)-u(x)\right|=\left|\int^1_0\nabla
		u(x-t\varepsilon \xi)\cdot \varepsilon \xi dt\right|\leq
		\varepsilon\int^1_0 \left|\nabla u(x-t\varepsilon \xi)\right|dt.
	\end{aligned}
\end{equation*}
Hence, for any $u\in \Lambda$, $\xi\in B_1$ and $\varepsilon\in
\left(0,\overline{\varepsilon}\right)$ we have

\begin{equation}\label{Sob:for1.97}
	\begin{aligned}
		\int_{\widetilde{\Omega}}\left|u_{\varepsilon}(x)-u(x)\right|dx&\leq
		\varepsilon\int^1_0 dt \int_{\widetilde{\Omega}}\left|\nabla
		u(x-t\varepsilon \xi)\right|dx=\\&= \varepsilon\int^1_0 dt
		\int_{\widetilde{\Omega}-t\varepsilon \xi}\left|\nabla
		u(z)\right|dz\leq\\&\leq \varepsilon \int_{\Omega}\left|\nabla
		u(z)\right|dz\leq \varepsilon
		|\Omega|^{1-\frac{1}{p}}\left\Vert
		u\right\Vert_{W^{1,p}(\Omega)}\leq
		\\&\leq \varepsilon
		|\Omega|^{1-\frac{1}{p}}.
	\end{aligned}
\end{equation}
From what we obtained in \eqref{Sob:for5.96} and by \eqref{Sob:for1.97}
we get

$$\int_{\widetilde{\Omega}}\left|u_{\varepsilon}(x)-u(x)\right|dx\leq \varepsilon
|\Omega|^{1-\frac{1}{p}}.$$ By the latter, by \eqref{Sob:for4.95}
and by \eqref{Sob:for4.96} we have

\begin{equation}\label{Sob:for1.98}
	\int_{\Omega}\left|u_{\varepsilon}(x)-u(x)\right|dx\leq
	2C_1\delta+\varepsilon |\Omega|^{1-\frac{1}{p}}, \quad\forall u\in
	\Lambda,\mbox{ } \forall \varepsilon\in
	\left(0,\overline{\varepsilon}\right),
\end{equation}
which implies

\begin{equation*}
	\limsup_{\varepsilon\rightarrow 0}\left(\sup_{u\in
		\Lambda}\left\Vert
	u_{\varepsilon}-u\right\Vert_{L^{1}(\Omega)}\right)\leq 2C_1\delta
\end{equation*}
and, as $\delta$ is arbitrary, we get

\begin{equation*}
	\lim_{\varepsilon\rightarrow 0}\left(\sup_{u\in \Lambda}\left\Vert
	u_{\varepsilon}-u\right\Vert_{L^{1}(\Omega)}\right)=0.
\end{equation*}
Hence, by Lemma \ref{Sob:lem1.8}, $\Lambda$ is relatively 
compact in $L^1(\Omega)$.

Now we consider the case $1<q<p^{\star}$. Denoting
$\theta= \frac{q-1}{p^{\star}-1}$,  we have $0<\theta<1$,
$q=1-\theta+\theta p^{\star}$, and

\begin{equation*}
	\begin{aligned}
		\int_{\Omega}\left|u_{\varepsilon}(x)-u(x)\right|^qdx&=\int_{\Omega}\left|u_{\varepsilon}(x)-u(x)\right|^{1-\theta+\theta
			p^{\star}}dx\leq\\&\leq \left\Vert
		u_{\varepsilon}-u\right\Vert_{L^1(\Omega)}^{1-\theta} \left\Vert
		u_{\varepsilon}-u\right\Vert_{L^{p^{\star}}(\Omega)}^{\theta}.
	\end{aligned}
\end{equation*}
Hence

\begin{equation}\label{Sob:for2.98}
	\left\Vert u_{\varepsilon}-u\right\Vert_{L^q(\Omega)}\leq \left\Vert
	u_{\varepsilon}-u\right\Vert_{L^1(\Omega)}^{\frac{1-\theta}{q}}
	\left\Vert
	u_{\varepsilon}-u\right\Vert_{L^{p^{\star}}(\Omega)}^{\frac{\theta}{q}}.
\end{equation}
On the other hand, by \eqref{Sob:for3.96} we have
$$\left\Vert
u\right\Vert_{L^{p^{\star}}(\Omega)}\leq C_1,\quad\forall u\in
\Lambda$$ and by the Young inequality we get

$$\left\Vert
u_{\varepsilon}\right\Vert_{L^{p^{\star}}(\Omega)}\leq \left\Vert
u\right\Vert_{L^{p^{\star}}(\Omega)}\leq C_1,\quad\forall u\in
\Lambda.$$ Hence, \eqref{Sob:for1.98} and \eqref{Sob:for2.98}
give
$$\left\Vert u_{\varepsilon}-u\right\Vert_{L^q(\Omega)}\leq
(2C_1)^{\frac{\theta}{q}}\left(2 C_1\delta+\varepsilon
|\Omega|^{1-\frac{1}{p}}\right)^{\frac{1-\theta}{q}},\quad\forall
u\in \Lambda.$$ Consequently

\begin{equation*}
	\lim_{\varepsilon\rightarrow 0}\left(\sup_{u\in \Lambda}\left\Vert
	u_{\varepsilon}-u\right\Vert_{L^{q}(\Omega)}\right)=0
\end{equation*}
and by Lemma \ref{Sob:lem1.8} we have that $\Lambda$ is relatively compact in $L^q(\Omega)$. $\blacksquare$

\bigskip

Now we state and prove a compactness theorem in the case $p>n$.

\begin{theo}\label{Sob:teo3.8}
	Let $\Omega$ be a bounded open set of $\mathbb{R}^n$ with boundary of 
	class $C^{0,1}$ and let $p>n$, $\alpha\in (0,\gamma)$, where
	$\gamma=1-\frac{n}{p}$. Then the embedding $$W^{1,p}(\Omega)
	\hookrightarrow C^{0,\alpha}\left(\overline{\Omega}\right),$$  is
	compact.
\end{theo}
\textbf{Proof.} Let $\left\{u_j\right\}$ be a sequence in  $W^{1,p}(\Omega)$ satisfying

$$\left\Vert
u_j\right\Vert_{W^{1,p}(\Omega)}\leq 1,\quad\forall j\in
\mathbb{N}.$$ By Theorem \ref{Sob:teo9.7} we have

\begin{equation}\label{Sob:for1.100}
	\left\Vert
	u_j\right\Vert_{C^{0,\alpha}\left(\overline{\Omega}\right)}\leq
	C_1,\quad\forall j\in \mathbb{N},
\end{equation}
where $C_1$ depends on $p,n$ and $\Omega$. The Arzel\`{a}--Ascoli Theorem
yields that there exists a subsequence
$\left\{u_{k_j}\right\}$ and $u\in
C^{0}\left(\overline{\Omega}\right)$ which satisfy

\begin{equation}\label{Sob:for2.100}
	\left\{u_{k_j}\right\}\rightarrow u,\quad\mbox{ uniformly}.
\end{equation}
By \eqref{Sob:for1.100} and \eqref{Sob:for2.100} we have, for any
$x,y\in \overline{\Omega}$, $x\neq y$,

\begin{equation}\label{Sob:for3.100}
	\frac{|u(x)-u(y)|}{|x-y|^{\gamma}}=\lim_{j\rightarrow\infty}
	\frac{\left|u_{k_j}(x)-u_{k_j}(y)\right|}{|x-y|^{\gamma}}\leq C_1.
\end{equation}
Hence $u\in C^{0,\gamma}\left(\overline{\Omega}\right)$. Therefore
$u\in C^{0,\alpha}\left(\overline{\Omega}\right)$ for
$0<\alpha<\gamma$. Now, let us recall the following inequality (see
Proposition \ref{Contin:prop-holder1}):

$$\left\Vert
f\right\Vert_{C^{0,\alpha}\left(\overline{\Omega}\right)}\leq C
\left\Vert
f\right\Vert^{\frac{\alpha}{\gamma}}_{C^{0,\gamma}\left(\overline{\Omega}\right)}\left\Vert
f\right\Vert^{1-\frac{\alpha}{\gamma}}_{C^{0}\left(\overline{\Omega}\right)},$$
For every $f\in C^{0,\gamma}\left(\overline{\Omega}\right)$, where $C$
depends by $\alpha$, $\gamma$ and $\Omega$ only. By applying such an inequality to $u_{k_j}-u$, taking into account \eqref{Sob:for1.100}--\eqref{Sob:for3.100}, we easily obtain

\begin{equation*}
	\left\{u_{k_j}\right\}\rightarrow u,\quad\mbox{ in }
	C^{0,\alpha}\left(\overline{\Omega}\right).
\end{equation*} $\blacksquare$

\subsection{Counterexamples} \label{controesempi-compatt}

\textbf{1.} The Rellich--Kondrachov Theorem does not hold for
$q=p^{\star}$. Indeed, we have the following counterexample. Let $u\in
C^{\infty}_0\left(B_1\setminus \{0\}\right)$ be a \textbf{not identically
	vanishing} function and let
$$u_j(x)=j^{\frac{n}{p^{\star}}}u(jx),\quad\forall j\in \mathbb{N}\mbox{, }\forall x\in B_1.$$
We have (see beginning of Section \ref{Sob:sec7.1})

\begin{equation}\label{Sob:for1.101}
	\left\Vert u_j\right\Vert_{L^{p^{\star}}(B_1)}=\left\Vert
	u\right\Vert_{L^{p^{\star}}(B_1)},\quad\forall j\in \mathbb{N},
\end{equation}

\begin{equation}\label{Sob:for3.101}
	\left\Vert u_j\right\Vert_{L^{p}(B_1)}=j^{-1}\left\Vert
	u\right\Vert_{L^{p}(B_1)},\quad\forall j\in \mathbb{N},
\end{equation}

\begin{equation}\label{Sob:for2.101}
	\left\Vert\nabla u_j\right\Vert_{L^{p}(B_1)}=\left\Vert \nabla
	u\right\Vert_{L^{p}(B_1)},\quad\forall j\in \mathbb{N}.
\end{equation}
Hence, by \eqref{Sob:for3.101} and \eqref{Sob:for2.101} we have

\begin{equation}\label{Sob:for4.101}
	\left\Vert u_j\right\Vert_{W^{1,p}(B_1)}=C\left\Vert
	u\right\Vert_{W^{1,p}(B_1)}<+\infty,\quad\forall j\in \mathbb{N}.
\end{equation}
Moreover
\begin{equation}\label{Sob:for5.101}
	\lim_{j\rightarrow\infty}u_j(x)=0,\quad\forall x\in
	B_1.
\end{equation}
Now, if the embedding

$$W^{1,p}(B_1)
\hookrightarrow L^{p^{\star}}(B_1)$$ were compact, by \eqref{Sob:for4.101} there should exist a subsequence
$\left\{u_{k_j}\right\}$ and \\ $v\in L^{p^{\star}}(B_1)$ such that

\begin{equation}\label{Sob:for5.101.0}
	\left\{u_{k_j}\right\}\rightarrow v,\quad \mbox{ in }
	L^{p^{\star}}(B_1).\end{equation}
Hence, by \eqref{Sob:for1.101} we should have
\begin{equation}\label{Sob:for1.102}
	\left\Vert v\right\Vert_{L^{p^{\star}}(B_1)}=\left\Vert
	u\right\Vert_{L^{p^{\star}}(B_1)}.
\end{equation}
On the other hand, passing eventually to another
subsequence, by \eqref{Sob:for5.101.0} we should have  

$$v(x)=\lim_{j\rightarrow\infty} u_{k_j}(x)\quad\mbox{ a.e. }
x\in B_1,$$ from the latter and from \eqref{Sob:for5.101} we should have

$$v(x)=0,\quad\mbox{ a.e. }
x\in B_1,$$ that would contradict \eqref{Sob:for1.102}.

\bigskip

\textbf{2.} Now, let us consider the case where $\Omega=\mathbb{R}^n$ and let us show that if $1\leq p<n$ and $q\leq p^{\star}$, then the embedding 

$$W^{1,p}\left(\mathbb{R}^n\right)
\hookrightarrow L^{q}\left(\mathbb{R}^n\right),$$ \textbf{is not}
compact. 

Let  $u\in C^{\infty}_0\left(\mathbb{R}^n\right)$ be a not identically vanishing function 
such that supp $u\subset B_1$ and let

$$u_j(x)=u(x-2je_1),\quad\forall j\in \mathbb{N}.$$
We obtain

$$\left\Vert u_j\right\Vert_{W^{1,p}\left(\mathbb{R}^n\right)}=\left\Vert
u\right\Vert_{W^{1,p}\left(\mathbb{R}^n\right)} ,\quad\forall j\in
\mathbb{N}$$ and

\begin{equation}\label{Sob:for2.102}
	\left\Vert
	u_j-u_k\right\Vert_{L^{q}\left(\mathbb{R}^n\right)}=2\left\Vert
	u\right\Vert_{L^{q}\left(\mathbb{R}^n\right)}>0,\quad\forall j, k\in
	\mathbb{N}\mbox{, } j\neq k.
\end{equation}
Hence $\left\{u_{j}\right\}$ is a bounded sequence in
$W^{1,p}\left(\mathbb{R}^n\right)$ but, as \eqref{Sob:for2.102} holds true, we cannot extract any subsequence that converges in
$L^{q}\left(\mathbb{R}^n\right)$.

\bigskip

\textbf{3.} Let us prove that if $\Omega$ is a bounded open set and $p>n$,
$\gamma=1-\frac{n}{p}$ then the embedding

$$W^{1,p}\left(\Omega\right)
\hookrightarrow C^{0,\gamma}\left(\overline{\Omega}\right)$$
\textbf{is not} compact.

Let $u\in C^{\infty}_0\left(B_1\right)$, not identically equal to $0$. Denote by $\overline{u}$ the extension of $u$ to $0$ in $\mathbb{R}^n\setminus B_1$. Let us denote
$$u_j(x)=\frac{1}{j^{\gamma}}\overline{u}(jx),\quad\forall j\in
\mathbb{N}\mbox{, } \forall x\in \overline{B_1}.$$ Now, let us notice (the reader check as an exercise)
$$\left[\overline{u}\right]_{0,\gamma,\mathbb{R}^n }=\left[u\right]_{0,\gamma,B_1 } $$
and
\begin{equation*}
	\begin{aligned}
		\left[u_j\right]_{0,\gamma,B_1 }&=\sup_{x,y\in B_1, x\neq y}\mbox{
		}\frac{\left|u_j(x)-u_j(y)\right|}{|x-y|^{\gamma}}=\\&=\sup_{x,y\in
			B_1, x\neq y}\mbox{
		}\frac{\left|\overline{u}(jx)-\overline{u}(jy)\right|}{|jx-jy|^{\gamma}}=\\&
		=\left[\overline{u}\right]_{0,\gamma, \mathbb{R}^n }.
	\end{aligned}
\end{equation*}

In addition we have
\begin{equation*}
	\begin{aligned}
		\left\Vert
		u_j\right\Vert_{L^{p}(B_1)}&=\frac{1}{j^{\gamma}}\left(\int_{B_1}|u(jx)|^pdx\right)^{\frac{1}{p}}=\\&=
		\frac{1}{j^{\gamma+\frac{n}{p}}}\left(\int_{B_{1/j}}|u(x)|^pdx\right)^{\frac{1}{p}}=\\&=
		\frac{1}{j}\left(\int_{B_{1/j}}|u(x)|^pdx\right)^{\frac{1}{p}}
	\end{aligned}
\end{equation*}
and
\begin{equation*}
	\begin{aligned}
		\left\Vert \nabla
		u_j\right\Vert_{L^{p}(B_1)}&=\frac{1}{j^{\gamma}}\left(\int_{B_1}|(\nabla
		u)(jx)|^pj^pdx\right)^{\frac{1}{p}}=\\&=
		\frac{j^{1-\frac{n}{p}}}{j^{\gamma}}\left(\int_{B_{1/j}}|\nabla
		u(x)|^pdx\right)^{\frac{1}{p}}=\\&= \left(\int_{B_{1/j}}|\nabla
		u(x)|^pdx\right)^{\frac{1}{p}}.
	\end{aligned}
\end{equation*}
Hence
$$\left\{u_j\right\}\rightarrow 0, \quad\mbox{ in } W^{1,p}(B_1),$$
in particular, $\left\{u_j\right\}$ is a bounded sequence in
$W^{1,p}(B_1)$. On the other hand, if there was a subsequence
$\left\{u_{k_j}\right\}$ of $\left\{u_j\right\}$ and $v\in
C^{0,\gamma}\left(\overline{B_1}\right)$ such that

$$\left\{u_{k_j}\right\}\rightarrow v, \quad\mbox{ in }
C^{0,\gamma}\left(\overline{B_1}\right),$$ we should necessarily
have $v\equiv 0$ and

$$\left[u_j-v\right]_{0,\gamma, B_1}=\left[u_j\right]_{0,\gamma,
	B_1}=\left[\overline{u}\right]_{0,\gamma, \mathbb{R}^n}>0$$ Which is a contradiction.

\bigskip

\textbf{4.} The case $p>n$, $\Omega=\mathbb{R}^n$, can be handle similarly to the case $p<n$. Let $u\in
C^{\infty}_0\left(\mathbb{R}^n\right)$, supp $u\subset B_1$, $u$ not identically vanishing function; let $u_j(x)=u(x-2je_1)$. We have
$$\left\Vert
u_j\right\Vert_{W^{1,p}\left(\mathbb{R}^n\right)}=\left\Vert
u\right\Vert_{W^{1,p}\left(\mathbb{R}^n\right)},\quad\forall j\in
\mathbb{N}$$ and
$$\left[u_j-u_k\right]_{0,\alpha, \mathbb{R}^n}\geq \left[u_j\right]_{0,\alpha,
	\mathbb{R}^n}=\left[u\right]_{0,\alpha, \mathbb{R}^n}>0,\quad \mbox{
	for } j\neq k,$$ where $\alpha\leq 1-\frac{n}{p}$. Hence, no extracted sequence of $ \left\{u_j\right\}$ can be a Cauchy sequence in
$C^{0,\alpha}\left(\mathbb{R}^n\right)$.

\section{The second Poincar\'{e} inequality}\label{Sob:poincar}

In Theorem \ref{Sob:teo6.7} we have proved the Sobolev--Poincar\'{e} inequality that, in particular, holds in the following form (see the proof
of the above mentioned Theorem)

\begin{equation}\label{Sob:for1.106}
	\left(\dashint_{B_r(x)}\left|u(y)-(u)_{x,r}\right|^{p\star}
	dy\right)^{\frac{1}{p^{\star}}}\leq C
	r\left(\dashint_{B_r(x)}\left|\nabla
	u(y)\right|^pdy\right)^{\frac{1}{p}},
\end{equation}
for every $u\in W^{1,p}\left(B_r(x)\right)$, where $p\in [1,+\infty)$
(actually it holds true for $p=+\infty$). We will now prove a
more general version of  \eqref{Sob:for1.106}.

\begin{theo}[\textbf{The second Poincar\'{e} inequality}]\label{Poincar1}
	\index{Theorem:@{Theorem:}!- second Poincar\'{e} inequality@{- second Poincar\'{e} inequality}}
	Let $\Omega$ be a bounded connected open set of $\mathbb{R}^n$ with $\partial
	\Omega$ of class $C^{0,1}$. Let $p\in [1,+\infty)$ and
	$$u_{\Omega}=\frac{1}{|\Omega|}\int_{\Omega} udx.$$
	Then there exists a constant $C$ depending on $p, n$ and $\Omega$ only, such that
	\begin{equation}\label{Poinc1}
		\left\Vert u-u_{\Omega}\right\Vert_{L^p(\Omega)}\leq C \left\Vert
		\nabla u\right\Vert_{L^p(\Omega)},\quad\forall u\in W^{1,p}(\Omega).
	\end{equation}
\end{theo}
\textbf{Proof.} We argue by contradiction. Let us assume that
\eqref{Poinc1} does not hold. Consequently for any $k\in \mathbb{N}$ there exists
$u_k\in W^{1,p}(\Omega)$ such that

$$\left\Vert u_k-(u_k)_{\Omega}\right\Vert_{L^p(\Omega)}> k \left\Vert
\nabla u_k\right\Vert_{L^p(\Omega)}.$$ Let us denote

$$v_k=\frac{u_k-(u_k)_{\Omega}}{\left\Vert
	u_k-(u_k)_{\Omega}\right\Vert_{L^p(\Omega)}},\quad \forall k\in
\mathbb{N}.$$ We have
$$(v_k)_{\Omega}=0,$$
$$\left\Vert v_k\right\Vert_{L^p(\Omega)}=1$$ and

\begin{equation}\label{Sob:for1.107}
	k\left\Vert \nabla v_k\right\Vert_{L^p(\Omega)}<1.
\end{equation}
Hence, there exists $M<+\infty$ such that 

$$\left\Vert v_k\right\Vert_{W^{1,p}(\Omega)}\leq M.$$
Therefore, by Rellich--Kondrachov Theorem we have that there exists a subsequence $\left\{v_{k_j}\right\}$ of $\left\{v_{k}\right\}$, and  $v\in L^p(\Omega)$ which satisfy
$$\left\{v_{k_j}\right\}\rightarrow v,\quad\mbox{ in } L^p(\Omega).$$
Hence
\begin{equation}\label{Sob:for2.107}
	v_{\Omega}=0
\end{equation}
and
\begin{equation}\label{Sob:for3.107}
	\left\Vert v\right\Vert_{L^p(\Omega)}=\lim_{j\rightarrow\infty}
	\left\Vert v_{k_j}\right\Vert_{L^p(\Omega)}=1.
\end{equation}
On the other hand by \eqref{Sob:for1.107} we have

$$\int_{\Omega}v\partial_l\phi dx=\lim_{j\rightarrow\infty}\int_{\Omega}v_{k_j}\partial_l\phi
dx=-\lim_{j\rightarrow\infty}\int_{\Omega}\partial_lv_{k_j}\phi
dx=0,\quad\forall \phi\in C^{\infty}_0\left(\Omega\right),$$ for
$l=1,\cdots, n$. Consequently
$$\int_{\Omega}v\partial_l\phi dx=0,\quad\forall \phi\in C^{\infty}_0\left(\Omega\right)\mbox{, }l=1,\cdots, n.$$
Therefore $\nabla v=0$ in $\Omega$ (and, trivially, $v\in
W^{1,p}(\Omega)$). Since $\Omega$ is a connected open set, Proposition \ref{Sob:Es5.39} yields that there exists a constant $c_0\in \mathbb{R}$ such that
$$v\equiv c_0,$$
and by \eqref{Sob:for2.107} we have $c_0=0$ that contradicts
\eqref{Sob:for3.107}. Therefore \eqref{Poinc1} holds true.
$\blacksquare$

\bigskip

\textbf{Remark.} 
The proof of \eqref{Poinc1} that we have given before is not constructive and this does not allow us to
further specify the dependence of the constant $C$ on
$\Omega$. To fill this gap we refer the reader to  \cite{A-R-M-08}. $\blacklozenge$

\section{The difference quotients}\label{quozienti di differenze} \index{difference quotients} In this Section we provide the definition and the main properties of the difference quotients. These topics
will turn out to be useful in the study of the regularity of the
solutions of second order elliptic equations.

\begin{definition}\label{Sob:def108}
	Let $V$ and $\Omega$ be open sets $\mathbb{R}^n$ such that $V\Subset
	\Omega$. Let $j\in\left\{1,\cdots,n \right\}$ and let $u\in
	L^1_{loc}(\Omega)$. The following function \index{$\delta^h_ju(x)$, $\delta^hu(x)$}
	
	\begin{equation}\label{Sob:for1.108}
		\delta^h_ju(x)=\frac{u(x+he_j)-u(x)}{h},\quad \forall x\in V.
	\end{equation}
is called $j$-th partial quotient  of $u$
with increment $h\in \mathbb{R}\setminus\{0\}$, \\ $|h|<$ dist $(V,\partial
\Omega)$. 
	We denote
	\begin{equation}\label{Sob:for2.108}
		\delta^hu(x)=\left(\delta^h_1u(x),\cdots,
		\delta^h_nu(x)\right),\quad \forall x\in V.
	\end{equation}
\end{definition}

\bigskip

We have the following

\begin{theo}\label{Sob:teo1.9}
	Let $\Omega$ be an open set of $\mathbb{R}^n$.
	
	\smallskip
	
	(i) If  $p\in [1,+\infty)$, $u\in W^{1,p}(\Omega)$,  then
	
	\begin{equation}\label{Sob:for3.108}
		\left\Vert \delta^h u\right\Vert_{L^p(V)}\leq C\left\Vert \nabla
		u\right\Vert_{L^p(\Omega)},\quad\mbox{ for }
		|h|<\frac{1}{2}\mbox{dist }(V,\partial \Omega), \ h\neq 0
	\end{equation}
	where $C$ depends on $n$ only.
	
	\smallskip
	
	(ii) Let us assume $p\in (1,+\infty)$, $u\in L^{p}(\Omega)$ and let us assume that
	there exists $C>0$ satisying
	\begin{equation}\label{Sob:for1.109}
		\left\Vert \delta^h u\right\Vert_{L^p(V)}\leq C,\quad\mbox{ for }
		|h|<\frac{1}{2}\mbox{dist }(V,\partial \Omega), \ h\neq 0
	\end{equation}
	then
	$$u\in W^{1,p}(V)\quad\mbox{ and }\quad \left\Vert \nabla u\right\Vert_{L^p(V)}\leq
	C.$$
\end{theo}
\textbf{Proof.} In order to prove (i) it suffices to assume $u\in
C^{\infty}(\Omega)\cap W^{1,p}(\Omega)$ and to apply Theorem
\ref{densit 1}.

If $0<|h|<$ dist $(V,\partial \Omega)$, $j=1,\cdots, n$ and $x\in V$,
we have

$$u(x+he_j)-u(x)=\int_0^1\nabla u (x+the_j)\cdot(he_j)dt,$$
from which we have

$$\left\vert u(x+he_j)-u(x)\right\vert\leq |h|\int_0^1\left\vert \nabla u
(x+the_j)\right\vert dt.$$ By using H\"{o}lder inequality and by integrating both the sides of the last inequality over  $V$, we get

\begin{equation*}
	\begin{aligned}
		\int_V \left\vert \delta^h_ju\right\vert^pdx&\leq \int_V dx
		\int^1_0\left\vert \nabla u (x+the_j)\right\vert^p
		dt=\\&=\int_0^1dt\int_V\left\vert \nabla u
		(x+the_j)\right\vert^pdx\leq\\&\leq \int_{\Omega}\left\vert \nabla u
		\right\vert^pdx.
	\end{aligned}
\end{equation*}

\medskip

Now, let us prove (ii). Let us assume that for some $C>0$ we have

\begin{equation}\label{correct:4-5-23}
	\left\Vert \delta^h u\right\Vert_{L^p(V)}\leq C,\quad\mbox{ for }
	0<|h|<\frac{1}{2}\mbox{dist }(V,\partial \Omega).
\end{equation}

\medskip

\textbf{Claim}

If $\phi\in
C^{\infty}_0(V)$ and let us denote $K=$ supp $\phi$, then for any $j\in \left\{1,\cdots,n\right\}$  we have

\begin{equation}\label{Sob:for1.110}
	\int_Vu\delta^h_j\phi dx=-\int_V \delta^{-h}_ju\phi dx,\quad\mbox{
		for } 0<|h|<\mbox{dist }(K,\partial V).
\end{equation}

\textbf{Proof of the Claim}. Let us notice that $$K-he_j\subset V\quad\mbox{
	for } 0<|h|<\mbox{dist }(K,\partial V),$$ for $j=1,\cdots,n$. Hence we have

\begin{equation*}
	\begin{aligned}
		\int_Vu\delta^h_j\phi
		dx&=\frac{1}{h}\left\{\int_Vu(x)\phi(x+he_j)dx-\int_Vu(x)\phi(x)dx\right\}=\\&=
		\frac{1}{h}\left\{\int_{K-he_j}u(x)\phi(x+he_j)dx-\int_Vu(x)\phi(x)dx\right\}=\\&=
		\frac{1}{h}\left\{\int_{K}u(x-he_j)\phi(x)dx-\int_Vu(x)\phi(x)dx\right\}=\\&=
		\frac{1}{h}\left\{\int_Vu(x-he_j)\phi(x)dx-\int_Vu(x)\phi(x)dx\right\}=\\&=
		-\int_V \delta^{-h}_ju\phi dx.
	\end{aligned}
\end{equation*}
Claim is proved.

\medskip

Let us fix $j\in\left\{1,\cdots,n \right\}$. Since $L^p(V)$ is a
a reflexive Banach space for $1<p<+\infty$, by

$$\sup\left\Vert \delta^{-h}_ju\right\Vert_{L^p(V)}\leq C$$
(recalling Theorems \ref{Sob:teo23R} and \ref{Sob:teo23Rbis}) there exists a sequence $\left\{h_k\right\}$ which goes to $0$ and $v_j\in L^p(V)$, such that

\begin{equation}\label{Sob:for2.110}
	\left\{\delta^{-h_k}_ju\right\}\rightharpoonup v_j,\quad\mbox{
		weakly in } L^p(V).
\end{equation}
On the other hand, by the Dominated Convergence Theorem we have, for any
$\phi\in C^{\infty}_0(\Omega)$ such that supp $\phi \subset V$,

$$\int_{\Omega}u\partial_j\phi dx=\lim_{k\rightarrow
	\infty}\int_{\Omega}u \delta^{h_k}_j\phi dx.$$ As a matter of fact

$$u(x)\delta^{h_k}_j\phi(x)\rightarrow
u(x)\partial_j\phi(x),\quad \forall x\in \Omega\mbox{  as } k\rightarrow\infty$$
and
$$\left|u\delta^{h_k}_j\phi\right|\leq |u|\left\Vert \nabla
\phi\right\Vert_{L^{\infty}(\Omega)}
\chi_{\widetilde{V}},\quad\forall k\in \mathbb{N},$$ where

$$\widetilde{V}=\left\{x\in \Omega:\mbox{ } \mbox{dist }(x,V)\leq \frac{1}{2} \mbox{dist }(V,\partial
\Omega)\right\}.$$ Therefore
\begin{equation*}
	\begin{aligned}
		\int_Vu\partial_j\phi dx&=\int_{\Omega}u\partial_j\phi dx=\\&=
		\lim_{k\rightarrow\infty}\int_{\Omega}u\left(\delta^{h_k}_j\phi\right) dx=\\&=
		-\lim_{k\rightarrow\infty}\int_{\Omega}\left(\delta^{-h_k}_ju\right)\phi dx=\\&=
		-\lim_{k\rightarrow\infty}\int_{V}\left(\delta^{-h_k}_ju\right)\phi dx=\\&=
		-\int_{V}v_j\phi dx.
	\end{aligned}
\end{equation*}
Consequently $$\partial_ju=v_j,\quad \mbox{in the weak sense for }
j=1,\cdots,n.$$ Hence $\nabla u\in L^p\left(V,\mathbb{R}^n\right)$,
but $u\in L^p\left(V\right)$. Therefore $u\in
W^{1,p}\left(V\right)$. 

Finally, by \eqref{Sob:for2.110} we have   
\begin{equation*}
\left\Vert \nabla u\right\Vert_{L^p(V)}\leq \liminf_{k\rightarrow\infty} \left\Vert \delta^{-h_k}u\right\Vert_{L^p(V)}\leq C,
\end{equation*}
($C$ is the same constant that occurs in \eqref{correct:4-5-23}).
$\blacksquare$

\bigskip

\textbf{Remark.} If $p=1$, then (ii) of Theorem \ref{Sob:teo1.9} does not hold. As a matter of fact, let $\Omega=(-2,2)$ and
$$u(t)=\chi_{(-1,1)}.$$
We have $u\in L^1(-2,2)$. Let
$V=\left(-\frac{3}{2},\frac{3}{2}\right)$. Now, dist $(V,\partial
\Omega)=\frac{1}{2}$ and for $0<|h|<\frac{1}{4}$ we have (for $h>0$)

\begin{equation*}
	\begin{aligned}
		\delta^hu(t)&=\frac{\chi_{(-1,1)}(t+h)-\chi_{(-1,1)}(t)}{h}=\\&=
		\frac{\chi_{(-1-h,1-h)}(t)-\chi_{(-1,1)}(t)}{h}=\\&=
		\frac{1}{h}\chi_{(-1,-1-h)\cup (1-h,1)}.
	\end{aligned}
\end{equation*}
Hence
$$\int_V\left|\delta^hu(t)\right|dt=2,\quad\mbox{ for } 0<|h|<\frac{1}{4},$$
but (see Example 2 of Section \ref{definizione H-k}) $$u\notin W^{1,1}(V).$$
$\blacklozenge$

\bigskip

In the sequel we will use the following variant of Theorem \ref{Sob:teo1.9}.
\begin{theo}\label{Sob:teo1.9bis}
	Let $r>0$.  We have
	
	(i) If  $p\in [1,+\infty)$ and $u\in W^{1,p}\left(B^+_r\right)$,  then for any $k\in\{1,\cdots,n-1\}$ we have, for $0<|h|<\frac{r}{2}$, 
	
	\begin{equation}\label{Sob:for3.108bis}
		\left\Vert \delta^h_k u\right\Vert_{L^p\left(B^+_{r/2}\right)}\leq C\left\Vert \partial_k u
		\right\Vert_{L^p\left(B^+_r\right)},
	\end{equation}
	where $C$ depends on $n$ only.
	
	\smallskip
	
	(ii) Let $k\in\{1,\cdots,n-1\}$. Let $p\in (1,+\infty)$, $u\in L^{p}\left(B^+_r\right)$ and let us suppose that there exists $C>0$ such that
	\begin{equation}\label{Sob:for1.109bis}
		\left\Vert \delta^h u\right\Vert_{L^p\left(B^+_{r/2}\right)}\leq C,\quad \mbox{ for } \ 0<|h|<\frac{r}{2},
	\end{equation}
	then
	$$\partial_ku\in L^p\left(B^+_r\right)\quad\mbox{ and }\quad \left\Vert \partial_ku\right\Vert_{L^p\left(B^+_r\right)}\leq
	C.$$
\end{theo}
The proof of the above Theorem is completely analogous to the one of Theorem \ref{Sob:teo1.9} and it is left to the reader as an exercise.

\section{The dual space of $H_0^{1}(\Omega)$}

Let $\Omega$ be an open set of $\mathbb{R}^n$. We denote by $H^{-1}(\Omega)$  the dual space of $H^1_0(\Omega)$ (i.e. the space of the linear bounded form from $H^1_0(\Omega)$ to $\mathbb{R}$).
If $F\in H^{-1}(\Omega)$, we write \index{$H^{-1}(\Omega)$}

\begin{equation*}\langle F,v\rangle:=F(v),\quad v\in
	H^1_0(\Omega)
\end{equation*}
and
\begin{equation*}
	\left\Vert F\right\Vert_{H^{-1}(\Omega)}=\sup\left\{\langle
	F,v\rangle:\mbox{ } v\in H^1_0(\Omega)\mbox{, } \left\Vert
	u\right\Vert_{H^{1}_0(\Omega)}\leq 1 \right\}.
\end{equation*}

\bigskip

The following Theorem holds true.

\begin{theo}[\textbf{characterization of $H^{-1}(\Omega)$}]\label{Sob:teo2.10}
	\index{Theorem:@{Theorem:}!- characterization of $H^{-1}(\Omega)$@{- characterization of $H^{-1}(\Omega)$}}
	Let $\Omega$ be an open set of $\mathbb{R}^n$.
	
	\smallskip
	
	\noindent (i) $F\in H^{-1}(\Omega)$  if and only if there exist $f_0,
	f_1, \cdots, f_n\in L^2(\Omega)$ satisfying
	
	\begin{equation}\label{duale}
		\langle F,v\rangle=\int_{\Omega} f_0vdx+\sum_{j=1}^n\int_{\Omega}
		f_j v_{x_j}dx,\quad \forall v\in H^1_0(\Omega).
	\end{equation}
	
	\smallskip
	
	\noindent(ii)
	\begin{equation*}
		\begin{aligned}
			&\left\Vert F\right\Vert_{H^{-1}(\Omega)}=\\&
			=\inf\left\{\left(\sum_{j=0}^n\int_{\Omega}|f_j|^2dx\right)^{1/2}: F
			\mbox{ satisfies } \eqref{duale} \mbox{ for }  f_0, f_1, \cdots,
			f_n\in L^2(\Omega)\right\}.
		\end{aligned}
	\end{equation*}
\end{theo}

\bigskip

We also write  $$F=f_0-\sum_{j=1}^n\partial_j f_j.$$ If
$$f_1, \cdots, f_n=0,$$ we will identify
the functional

\begin{equation*}
	\langle F,v\rangle=\int_{\Omega} f_0vdx,\quad\quad \forall v\in
	H^1_0(\Omega)
\end{equation*}
with $f_0$ and we will write $F\in L^2(\Omega)$. Similarly, if $f_0\in
H^k(\Omega)$, we will write $F\in H^k(\Omega)$.

Let us note that $f_0, f_1, \cdots, f_n$ are \emph{not} uniquely
determined. For instance, if $\Omega$ is bounded, then the functional
$$\langle F,v\rangle=\int_{\Omega} f_0vdx,$$
where $f_0\in L^2(\Omega)$, can also be represented by
$$\langle F,v\rangle=\int_{\Omega} (f_0+2x_1)vdx+\int_{\Omega} x^2_1\partial_1vdx, \quad \forall v\in H^1_0(\Omega)$$
and in infinite other ways.

\bigskip

\textbf{Proof of Theorem \ref{Sob:teo2.10}.} Let us equip
$H^1_0(\Omega)$ with the scalar product

$$(u,v)_{H^1_0(\Omega)}=\int_{\Omega}(uv+\nabla u\cdot\nabla v)dx,\quad\forall u,v\in H^1_0(\Omega).$$
It is clear that if $F$ is like \eqref{duale}, then $F\in
H^{-1}(\Omega)$, as a matter of fact, by applying 
the Cauchy--Schwarz inequality we get 
\begin{equation}\label{Sob:for1.115}
	\begin{aligned}
		\left\vert\langle F,v\rangle\right\vert&= \left\vert\int_{\Omega}
		f_0vdx+\sum_{j=1}^n\int_{\Omega} f_j
		\partial_jvdx\right\vert\leq\\&\leq
		\left(\sum_{j=0}^n\int_{\Omega}|f_j|^2dx\right)^{1/2} \left\Vert
		v\right\Vert_{H^1_0(\Omega)},\quad\forall v\in H^1_0(\Omega).
	\end{aligned}
\end{equation}

Conversely, let us assume $F\in H^{-1}(\Omega)$. By the Riesz representation Theorem we have that there exists a unique $u\in
H^1_0(\Omega)$ such that

\begin{equation}\label{Sob:for2.115}
	\langle F,v\rangle=(u,v)_{H^1_0(\Omega)}=\int_{\Omega}(uv+\nabla
	u\cdot\nabla v)dx ,\quad\forall v\in H^1_0(\Omega).
\end{equation}
Hence, denoting

\begin{equation}\label{Sob:for3.115}
	f_0=u,\quad f_j=\partial_ju,\quad j=1,\cdots,n,
\end{equation}
we have
\begin{equation}\label{Sob:for4.115}
	\langle F,v\rangle=\int_{\Omega} f_0vdx+\sum_{j=1}^n\int_{\Omega}
	f_j\partial_jvdx,\quad\forall v\in H^1_0(\Omega).
\end{equation}
The proof of (i) is concluded.

\medskip

Now, let us prove (ii). Let $u\in H^1_0(\Omega)$ and $f_j\in
L^2(\Omega)$, $j=0,1,\cdots,n$, be like in \eqref{Sob:for3.115}. Let
$g_j\in L^2(\Omega)$, $j=0,1,\cdots,n$ satisfy 
\begin{equation*}
	\langle F,v\rangle=\int_{\Omega} \left( g_0v+\sum_{j=1}^n
	g_j\partial_jv\right)dx,\quad\forall v\in H^1_0(\Omega).
\end{equation*}
Let us check that

\begin{equation}\label{Sob:for1.116}
	\int_{\Omega}\sum_{j=0}^n \left|f_j\right|^2dx\leq
	\int_{\Omega}\sum_{j=0}^n \left|g_j\right|^2dx.
\end{equation}
We have

\begin{equation*}
	\begin{aligned}
		\int_{\Omega}\left(\left|u\right|^2 +\left|\nabla u\right|^2\right)
		dx&=\langle F, u\rangle=\int_{\Omega} \left( g_0u+\sum_{j=1}^n
		g_j\partial_ju\right)dx\leq\\&\leq \left(\int_{\Omega}\sum_{j=0}^n
		\left|g_j\right|^2
		dx\right)^{1/2}\left(\int_{\Omega}\left(\left|u\right|^2
		+\left|\nabla u\right|^2\right) dx\right)^{1/2}.
	\end{aligned}
\end{equation*}
Hence
\begin{equation}\label{Sob:for2.116}
	\begin{aligned}
		\int_{\Omega}\sum_{j=0}^n \left|f_j\right|^2
		dx&=\int_{\Omega}\left(\left|u\right|^2 +\left|\nabla
		u\right|^2\right) dx\leq\\&\leq \int_{\Omega}\sum_{j=0}^n
		\left|g_j\right|^2 dx,
	\end{aligned}
\end{equation}
which proves \eqref{Sob:for1.116}.

In order to complete the proof, let us notice that by 
\eqref{Sob:for1.115} we get

\begin{equation}\label{Sob:for3.116}
	\begin{aligned}
		\left\Vert F\right\Vert_{H^{-1}(\Omega)}\leq
		\left(\sum_{j=0}^n\int_{\Omega}|f_j|^2dx\right)^{1/2}.
	\end{aligned}
\end{equation}
On the other hand, setting
$$\widetilde{v}=\frac{u}{\left\Vert
	u\right\Vert_{H^{1}_0(\Omega)}},$$ we obtain, by \eqref{Sob:for4.115}
(recall that $u$ satisfies \eqref{Sob:for3.115}),

\begin{equation}\label{Sob:for1.117}
	\begin{aligned}
		\langle F,\widetilde{v}\rangle&=\int_{\Omega}(u\widetilde{v}+\nabla
		u\cdot\nabla \widetilde{v})dx=\\&=\left\Vert
		u\right\Vert_{H^{1}_0(\Omega)}=\\&=\left(\sum_{j=0}^n\int_{\Omega}|f_j|^2dx\right)^{1/2}.
	\end{aligned}
\end{equation}
Hence, by \eqref{Sob:for3.116} and \eqref{Sob:for1.117} we have

\begin{equation*}
	\begin{aligned}
		\left\Vert F\right\Vert_{H^{-1}(\Omega)}=
		\left(\sum_{j=0}^n\int_{\Omega}|f_j|^2dx\right)^{1/2}.
	\end{aligned}
\end{equation*}
By the just obtained equality and by \eqref{Sob:for2.116} we obtain (ii).
$\blacksquare$

\bigskip

\textbf{Remark.} Let us note that the  greatest lower bound in (ii) is actually
the minimum. $\blacklozenge$

\bigskip

\underline{\textbf{Exercise.}} Let us denote by $H^{-m}(\Omega)$, $m\in \mathbb{N}$, the  dual space of
$H^{m}_0(\Omega)$. Prove that $F\in H^{-m}(\Omega)$ if and only if there exist $f_{\alpha}\in L^2(\Omega)$, $|\alpha|\leq m$ such that

$$\langle F,\widetilde{v}\rangle=\int_{\Omega}\sum_{|\alpha|\leq
	m}f_{\alpha}\partial^{\alpha}vdx,\quad \forall v\in
H^{m}_0(\Omega)$$ and prove the analogue of the part (ii) of
Theorem \ref{Sob:teo2.10}. $\clubsuit$

\section{The Sobolev spaces with noninteger exponents and traces}
\label{Sobolev-Fourier-Tracce} 
In the present Section we will provide a
characterization of the traces of $H^k(\Omega)$ function,
$k\in \mathbb{N}$. For this purpose we need to extend the notion of the
Sobolev space that we have studied so far to the spaces with non
integer exponents. First of all, we provide brief reminders of the
Fourier transform.

\subsection{Review of the Fourier transform} \label{Trasformata
	di Fourier} \index{Fourier transform} Let us denote by $D_j=\frac{1}{i}\partial_j$, 
$j=1,\cdots, n$.

\medskip

\begin{definition}\label{Fourier:def1.7.1} We denote by $\mathcal{S}$ the space of functions $f\in C^{\infty}(\mathbb{R}^n)$ which satisfy
	\begin{equation}\label{Four:for1.7.3}
		p_{\alpha,\beta}(f):=\sup_{x\in \mathbb{R}^n }\left\vert
		x^{\alpha}D^{\beta}f\right\vert<\infty, \quad \forall \alpha,\beta
		\in \mathbb{N}_0^n.
	\end{equation}
	The topology on $\mathcal{S}$ is defined by the seminorms
	$p_{\alpha,\beta}(f)$.
\end{definition}
According to Definition \ref{Fourier:def1.7.1}, a sequence
$\left\{f_k\right\}\subset \mathcal{S}$ converges to $f\in
\mathcal{S}$ if and only if
$$\lim_{k\rightarrow\infty}p_{\alpha,\beta}(f_k-f)=0, \quad \forall \alpha,\beta
\in \mathbb{N}_0^n.$$ The space $\mathcal{S}$ is known as
\textbf{the Schwartz space} \index{Schwartz space} or, also, the space of
rapidly decreasing functions; equipped with the family of seminorms $\left\{p_{\alpha,\beta}\right\}$, $\mathcal{S}$ is a Frech\'{e}t space \index{Fr\'{e}chet space}
(for the definition of Fr\'{e}chet space we refer to \cite{Ru} and, in the present Notes, Section \ref{LM-24}). We have

$$C_0^{\infty}\left(\mathbb{R}^n\right)\subset \mathcal{S}\subset L^p\left(\mathbb{R}^n\right), \ \ \forall p\in [1,+\infty].$$
It is simple to prove that  $C_0^{\infty}\left(\mathbb{R}^n\right)$
is dense in $\mathcal{S}$. The function $f(x)=e^{-|x|^2}$ is an example of function that does not belong to $C_0^{\infty}\left(\mathbb{R}^n\right)$, but 
belongs to $\mathcal{S}$.

\bigskip

\begin{definition}\label{Fourier:deftrasformata} Let $f\in L^1\left(\mathbb{R}^n\right)$,
	we define its  Fourier transform\index{Definition:@{Definition:}!- Fourier transform@{- Fourier transform}} by \index{$\widehat{f}$} 
	$$\widehat{f}(\xi):=\mathcal{F}(f)(\xi):=\int_{\mathbb{R}^n}f(x)e^{-ix\cdot\xi}dx, \ \ \forall \xi\in\mathbb{R}^n.$$
	\end{definition}
\bigskip

We have 
\begin{equation*}
	\left\Vert\widehat{f}\right\Vert_{L^{\infty}\left(\mathbb{R}^n\right)}\leq \left\Vert f\right\Vert_{L^{\infty}\left(\mathbb{R}^n\right)}, \ \ \forall f\in L^1\left(\mathbb{R}^n\right).
\end{equation*}
Actually, we have $\widehat{f}\in C^0\left(\mathbb{R}^n\right)$ and
\begin{equation}\label{correct:15-4-23-RiemannLebesgue}
	\widehat{f}(\xi) \rightarrow 0, \ \ \mbox{as } |\xi|\rightarrow +\infty.
\end{equation} 
Property \eqref{correct:15-4-23-RiemannLebesgue} is known as \textbf{Riemann--Lebesgue Lemma}. \index{Riemann--Lebesgue Lemma}

\bigskip

Let us recall that, if $f(x)=e^{-\frac{|x|^2}{2}}$ then
$$\widehat{f}(\xi)=(2\pi)^{n/2}e^{-\frac{|\xi|^2}{2}}.$$

\bigskip

\begin{theo}\label{Fourier:propI} Let $f\in\mathcal{S}$. Then we have $\widehat{f}\in\mathcal{S}$. Moreover the following properties hold.
	
	\smallskip
	
	a) The map $$\mathcal{S}\ni f\rightarrow \widehat{f}\in
	\mathcal{S}$$ is one--to--one, continuous, with continuous inverse and 
	\begin{equation}\label{correct:15-4-23-inversionformula}
	f(x)=\frac{1}{(2\pi)^n}
	\int_{\mathbb{R}^n}\widehat{f}(\xi)e^{ix\cdot\xi}d\xi, \ \ \forall x\in \mathbb{R}^n;	
\end{equation}

	\smallskip
	
	b) $$\widehat{D_x^{\alpha}f}(\xi)=\xi^{\alpha}\widehat{f}(\xi),
	\quad \forall f\in \mathcal{S};$$
	
	\smallskip
	
	c)
	$$\widehat{(x^{\alpha}f)}(\xi)=(-1)^{|\alpha|}D_{\xi}^{\alpha}\widehat{f}(\xi)\quad \forall f\in
	\mathcal{S}.$$
\end{theo}

Formula \eqref{correct:15-4-23-inversionformula} is known as the \textbf{inversion formula for the Fourier transform}\index{inversion formula for the Fourier transform}. If $f\in L^1\left(\mathbb{R}^n\right)$ and $\widehat{f}\in L^1\left(\mathbb{R}^n\right)\cap L^{\infty}\left(\mathbb{R}^n\right)$, we have
\begin{equation}\label{correct:15-4-23-inversionformula-l1}
	f(x)=\frac{1}{(2\pi)^n}
	\int_{\mathbb{R}^n}\widehat{f}(\xi)e^{ix\cdot\xi}d\xi, \ \ \mbox{a.e. in } \mathbb{R}^n.	
\end{equation}  

\bigskip

\begin{theo}\label{Fourier:propII} Let $f,g\in\mathcal{S}$.
	We have

	\begin{equation}\label{Four:for1.7.6}
		\int_{\mathbb{R}^n}f(x)\widehat{g}(x)dx=\int_{\mathbb{R}^n}\widehat{f}(x)g(x)dx,
	\end{equation}
	\begin{equation}\label{Four:for1.7.7}
		\int_{\mathbb{R}^n}f(x)\overline{g(x)}dx=\frac{1}{(2\pi)^n}\int_{\mathbb{R}^n}\widehat{f}(\xi)\overline{\widehat{g}(\xi)}d\xi,
	\end{equation}
	
	\begin{equation}\label{Four:for1.7.8}
		\widehat{(f\star g)}(\xi)=\widehat{f}(\xi)\widehat{ g}(\xi),
	\end{equation}
	
	\begin{equation}\label{Four:for1.7.9}
		\widehat{(fg)}(\xi)=\frac{1}{(2\pi)^n}(\widehat{f}\star\widehat{
			g})(\xi).
	\end{equation}
\end{theo}

\bigskip

Formula \eqref{Four:for1.7.7} is known as \textbf{Parseval formula}\index{Parseval formula} and it is equivalent to the following one
\begin{equation}\label{Four:for1.7.7-pars}
	\left\Vert
	f\right\Vert_{L^2(\mathbb{R}^n)}=\frac{1}{(2\pi)^{n/2}}\left\Vert
	\widehat{f}\right\Vert_{L^2(\mathbb{R}^n)}, \ \ \forall f\in\mathcal{S}.
\end{equation}

Let us notice that the restriction of the linear operator $\mathcal{F}$ over $\mathcal{S}$ acts as follows  
$$ \mathcal{S}\ni f\rightarrow \mathcal{F}(f):=\widehat{f}\in \mathcal{S}.$$
Moreover $\mathcal{F}$ is bijective and by \eqref{Four:for1.7.7-pars} we have
\begin{equation}\label{Four:for1.7.7-pars-15-4}
	\left\Vert
	\mathcal{F}(f)\right\Vert_{L^2(\mathbb{R}^n)}=(2\pi)^{n/2}\left\Vert
	f\right\Vert_{L^2(\mathbb{R}^n)}, \ \ \forall f\in\mathcal{S}.
\end{equation}
Let us observe that, since $C_0^{\infty}\left(\mathbb{R}^n\right)$ is dense in $L^{2}\left(\mathbb{R}^n\right)$ and $$C_0^{\infty}\left(\mathbb{R}^n\right)\subset \mathcal{S} \subset L^{2}\left(\mathbb{R}^n\right),$$ then $\mathcal{S}$ is dense in $L^{2}\left(\mathbb{R}^n\right)$. Hence \eqref{Four:for1.7.7-pars-15-4} implies that the linear operator $\mathcal{F}$  
 can be extended to a bounded linear operator from $L^{2}\left(\mathbb{R}^n\right)$ to $L^{2}\left(\mathbb{R}^n\right)$. We continue to denote by $\mathcal{F}$ such an extension. Hence it is defined 
 $$\widehat{f}:=\mathcal{F}(f), \ \ \forall f\in L^{2}\left(\mathbb{R}^n\right).$$ 
 It can be proved that the operator $$\mathcal{F}:L^2(\mathbb{R}^n)\rightarrow L^2(\mathbb{R}^n)$$ is bijective and Theorem \ref{Fourier:propII} continue to holds. Moreover, denoting by
 $$\mathcal{C}: L^{2}\left(\mathbb{R}^n\right)\rightarrow L^{2}\left(\mathbb{R}^n\right),$$
 $$ \left(\mathcal{C}(f)\right)(x)=\frac{1}{(2\pi)^n}f(-x), \  \ \forall f\in L^{2}\left(\mathbb{R}^n\right), \ \forall x\in \mathbb{R}^n,$$
 we have 
 \begin{equation}\label{correct:16-4-23-1}
 	f=\mathcal{C}\mathcal{F}(f), \ \ \forall f\in L^{2}\left(\mathbb{R}^n\right).
 \end{equation}
If $f\in \mathcal{S}$ or $f\in L^{2}\left(\mathbb{R}^n\right)\cap L^{1}\left(\mathbb{R}^n\right)$ then formula \eqref{correct:16-4-23-1} is nothing but inversion formula \eqref{correct:15-4-23-inversionformula} or \eqref{correct:15-4-23-inversionformula-l1} respectively. For the proofs and much more details we refer the reader to \cite[Vol. I]{HOII}, \cite{EV}, \cite{Ru} 
\subsection{Fourier transform and $H^k\left(\mathbb{R}^n\right)$ spaces, $k\in \mathbb{N}_0$}
\label{SobFourInt} Let us state and prove 

\begin{theo}\label{Sob:teo1.11}
Let $k\in \mathbb{N}_0$. The following properties hold
	
	\smallskip
	
	\noindent (i) Let $u\in L^2\left(\mathbb{R}^n\right)$. We have that $u\in
	H^k\left(\mathbb{R}^n\right)$ if and only if
	
	\begin{equation}\label{Sob:for1.118}
		\left(1+|\xi|^2\right)^{k/2}\widehat{u}(\xi)\in
		L^2\left(\mathbb{R}^n\right).
	\end{equation}
	\noindent (ii) There exists a constant $C\geq 1$ depending on $k$ and
	$n$ only, such that 
	\begin{equation}\label{Sob:for2.118}
		C^{-1}\left\Vert u\right\Vert_{H^k\left(\mathbb{R}^n\right)}\leq
		\left\Vert
		\left(1+|\xi|^2\right)^{k/2}\widehat{u}(\xi)\right\Vert_{L^2\left(\mathbb{R}^n\right)}\le
		C\left\Vert u\right\Vert_{H^k\left(\mathbb{R}^n\right)},
	\end{equation}
	for every $u\in H^k\left(\mathbb{R}^n\right)$.
\end{theo}
\textbf{Proof.} If $k=0$, then (i) and (ii) are obvious. Let us assume $k\geq 1$ and let us begin to prove  (i). Since $\mathcal{S}\subset H^k\left(\mathbb{R}^n\right)$, we have  that if $u\in H^k\left(\mathbb{R}^n\right)$ then

\begin{equation}\label{corr:20-4-23-1}
	\int_{\mathbb{R}^n}\partial^{\alpha}u\varphi dx=(-1)^{|\alpha|}\int_{\mathbb{R}^n}u\partial^{\alpha} \varphi dx, \ \ \mbox{for } |\alpha|\leq k , \ \forall \varphi \in \mathcal{S}.
\end{equation}

\bigskip

\textbf{Claim I} 

\begin{equation}\label{corr:20-4-23-2}
u \in H^k\left(\mathbb{R}^n\right) \ \Longleftrightarrow \ (i\xi)^{\alpha} \widehat{u}(\xi)\in L^2\left(\mathbb{R}^n\right), \ \mbox{for } |\alpha|\leq k.
\end{equation}

\medskip

\textbf{Proof of Claim I.} Let us prove "$\Longrightarrow$". Let $u\in H^k\left(\mathbb{R}^n\right)$, $R$ be an arbitrary positive number and $\psi\in C^{\infty}_0\left(B_R\right)$. Let us define

\begin{equation*}
\varphi(x)=\frac{1}{(2\pi)^n}
\int_{\mathbb{R}^n}\overline{\psi (\xi)}e^{ix\cdot\xi}d\xi, \ \ \forall x\in \mathbb{R}^n.
\end{equation*}
 We have $$\varphi \in \mathcal{S}, \ \ \mbox{and} \ \ \overline{\widehat{\varphi}(\xi)}=\psi(\xi), \ \forall \xi \in \mathbb{R}^n.$$
 In addition, as $u\in L^2\left(\mathbb{R}^n\right)$, 
 $$(i\xi)^{\alpha} \widehat{u}(\xi)_{|B_R}\in L^2\left(B_R\right), \ \mbox{for } |\alpha|\leq k,$$
hence $(i\xi)^{\alpha} \widehat{u}(\xi)\psi(\xi)\in L^2\left(\mathbb{R}^n\right)$ and by the Parseval identity we have, for $|\alpha|\leq k$,

\begin{equation*}
	\begin{aligned}
	\int_{\mathbb{R}^n} (i\xi)^{\alpha} \widehat{u}(\xi) \psi(\xi)d\xi &=\int_{\mathbb{R}^n} (i\xi)^{\alpha} \widehat{u}(\xi) \overline{\widehat{\varphi}(\xi)} d\xi=\\&=
	(-1)^{|\alpha|} \int_{\mathbb{R}^n}  \widehat{u}(\xi) \overline{(i\xi)^{\alpha}\widehat{\varphi}(\xi)} d\xi=\\&=(-1)^{|\alpha|}(2\pi)^n \int_{\mathbb{R}^n} u(x)\overline{\partial^{\alpha}\varphi(x)}dx=\\&=(2\pi)^n \int_{\mathbb{R}^n} \partial^{\alpha}u(x)\overline{\varphi(x)}dx=\\&=\int_{\mathbb{R}^n}\widehat{\partial^{\alpha}u}(\xi)\overline{\widehat{\varphi}(\xi)} d\xi=\\&=\int_{\mathbb{R}^n}\widehat{\partial^{\alpha}u}(\xi)\psi(\xi)d\xi
	\end{aligned}
\end{equation*}
Hence
\begin{equation*}
	\int_{\mathbb{R}^n}\left((i\xi)^{\alpha} \widehat{u}-\widehat{\partial^{\alpha}u}\right)\psi d\xi=0, \ \forall \psi \in C^{\infty}_0\left(B_R\right),
\end{equation*}
 but $\widehat{\partial^{\alpha}u}_{|B_R}\in L^2\left(B_R\right)$ for $|\alpha|\leq k$, (because $u\in H^k\left(\mathbb{R}^n\right)$), therefore 
 
 \begin{equation*}
 	(i\xi)^{\alpha} \widehat{u}=\widehat{\partial^{\alpha}u}, \ \mbox{in } B_R
  \end{equation*}
and since $R$ is arbitrary,
 $$(i\xi)^{\alpha} \widehat{u}=\widehat{\partial^{\alpha}u}, \ \mbox{in } \mathbb{R}^n.$$
In particular 
$$(i\xi)^{\alpha} \widehat{u}\in L^2\left(\mathbb{R}^n\right),$$
hence "$\Longrightarrow$" is proved.

\medskip

Now we prove"$\Longleftarrow$". Let us assume that

\begin{equation*}
	(i\xi)^{\alpha} \widehat{u}(\xi)\in L^2\left(\mathbb{R}^n\right), \ \mbox{for } |\alpha|\leq k.
\end{equation*}
Let $u_{\alpha}\in L^2\left(\mathbb{R}^n\right)$, $|\alpha|\leq
k$, be defined by

$$u_{\alpha}(x)=\frac{1}{(2\pi)^n}\int_{\mathbb{R}^n}(i\xi)^{\alpha}
\widehat{u}(\xi)e^{ix\cdot\xi}d\xi.$$ For any $\phi\in
C_0^{\infty}\left(\mathbb{R}^n\right)$ and any $|\alpha|\leq k$, we have
(recall \eqref{Four:for1.7.7})

\begin{equation*}
	\begin{aligned}
		\int_{\mathbb{R}^n}\partial^{\alpha}\phi(x)\overline{u}(x)dx&=\frac{1}{(2\pi)^n}\int_{\mathbb{R}^n}\widehat{\partial^{\alpha}\phi}(\xi)\overline{\widehat{u}(\xi)}d\xi=\\&=
		\frac{1}{(2\pi)^n}\int_{\mathbb{R}^n}(i\xi)^{\alpha}\widehat{\phi}(\xi)\overline{\widehat{u}(\xi)}d\xi=\\&=\frac{(-1)^{|\alpha|}}{(2\pi)^n}\int_{\mathbb{R}^n}\widehat{\phi}(\xi)\overline{(i\xi)^{\alpha}\widehat{u}(\xi)}d\xi=\\&=
		\frac{(-1)^{|\alpha|}}{(2\pi)^n}\int_{\mathbb{R}^n}\widehat{\phi}(\xi)\overline{\widehat{u_{\alpha}}(\xi)}d\xi=\\&=
		(-1)^{|\alpha|}\int_{\mathbb{R}^n}\phi(x) \overline{u_{\alpha}(x)}dx.
	\end{aligned}
\end{equation*}
Hence
$$\partial^{\alpha}u=u_{\alpha}\in L^2\left(\mathbb{R}^n\right),\quad\mbox{ for
} |\alpha|\leq k.$$ Therefore $u\in H^k\left(\mathbb{R}^n\right)$. Claim I is proved.

\bigskip

\textbf{Claim II.}

The following conditions are equivalent

\medskip

\noindent (a) $(i\xi)^{\alpha} \widehat{u}\in
L^2\left(\mathbb{R}^n\right)$ for $|\alpha|\leq k$,

\medskip

\noindent (b) $\left(1+|\xi|^2\right)^{k/2}\widehat{u}(\xi)\in
L^2\left(\mathbb{R}^n\right)$.

\medskip

\textbf{Proof of Claim II.} First, let us note that (a) is equivalent to

$$\int_{\mathbb{R}^n}\sum_{|\alpha|\leq k}\left|\xi^{\alpha}\right|^2
\left|\widehat{u}(\xi)\right|^2d\xi<+\infty,$$ hence, in order to prove
that (a) and (b) are equivalent it suffices to prove that there exists $C\geq 1$
such that

\begin{equation}\label{Sob:for1.120}
	C^{-1}\leq \frac{\sum_{|\alpha|\leq
			k}\left|\xi^{\alpha}\right|^2}{\left(1+|\xi|^2\right)^k}\leq
	C,\quad\forall\xi\in \mathbb{R}^n.
\end{equation}
To this purpose we notice that the function

$$g(\xi,\tau)=\frac{\sum_{|\alpha|\leq
		k}\tau^{2(k-|\alpha|)}\left|\xi^{\alpha}\right|^2}{\left(\tau^2+|\xi|^2\right)^k},$$
is homogeneous of degree $0$, it is continuous in
$\mathbb{R}^{n+1}\setminus\{(0,0)\}$, and
$$g(\xi,\tau)>0,\quad \forall
(\xi,\tau)\in \mathbb{R}^{n+1}\setminus\{(0,0)\}.$$ Hence there exists $C\geq
1$ such that
$$C^{-1}\leq g(\xi,\tau)\leq C$$
so that, if $\tau=1$ we get \eqref{Sob:for1.120}, which, in turn implies the equivalence of (a) an (b). Claim II is proved. 

\medskip

By Claim I and Claim II we obtain (i).

\bigskip

Concerning (ii), it is enough to observe that by
\eqref{correct:16-4-23-1} we have

$$\partial^{\alpha}u(x)=\frac{1}{(2\pi)^n}\int_{\mathbb{R}^n}(i\xi)^{\alpha}
\widehat{u}(\xi)e^{ix\cdot\xi}d\xi, \ \mbox{for } |\alpha|\leq k$$ and by the Parseval identity we have

$$\left\Vert
u\right\Vert_{H^k\left(\mathbb{R}^n\right)}=\frac{1}{(2\pi)^n}\int_{\mathbb{R}^n}\sum_{|\alpha|\leq
	k}\left|\xi^{\alpha}\right|^2 \left|\widehat{u}(\xi)\right|^2d\xi.$$
Hence by \eqref{Sob:for1.120} we derive  \eqref{Sob:for2.118}.
$\blacksquare$

\subsection{The Sobolev spaces with noninteger exponents}
\label{Sob:frazionari} \index{Sobolev spaces with noninteger exponent}Theorem \ref{Sob:teo1.11} justifies the
the following definition (instead of $\mathbb{R}^n$ we will consider
$\mathbb{R}^m$, with $m\in \mathbb{N}$ to avoid confusion
later on, when we will need to set $m=n-1$)

\begin{definition}\label{Sob:def1.11} Let $m\in \mathbb{N}$ and let $s$ be a
	real positive number, we say that $u\in H^s\left(\mathbb{R}^m\right)$ \index{$H^s\left(\mathbb{R}^m\right)$, ($s>0$ real number)}if
	
	\begin{equation*}
		\left(1+|\xi|^2\right)^{s/2}\widehat{u}(\xi)\in
		L^2\left(\mathbb{R}^m\right),
	\end{equation*}
	in this case we denote \index{$\left\Vert
		\cdot\right\Vert_{H^s\left(\mathbb{R}^m\right)}$, ($s>0$ real number)}
	\begin{equation}\label{Sob:for1.122}
		\left\Vert
		u\right\Vert_{H^s\left(\mathbb{R}^m\right)}=\left(\frac{1}{(2\pi)^m}\int_{\mathbb{R}^m}\left(1+|\xi|^2\right)^{s}\left|\widehat{u}(\xi)\right|^2d\xi\right)^{1/2}.
	\end{equation}
\end{definition}
It is evident that if $s\in \mathbb{N}$ we again obtain the
Sobolev spaces with integer exponents that we have studied so far,
nevertheless if $s\notin \mathbb{N}$ we obtain some new spaces, namely the \textbf{Sobolev spaces with noninteger exponents}, also known as "the Sobolev spaces with
fractional exponent". It is simple to check that the norm 
$\left\Vert \cdot\right\Vert_{H^s\left(\mathbb{R}^m\right)}$ is
induced by the scalar product \index{$(\cdot,\cdot)_{H^s\left(\mathbb{R}^m\right)}$}

$$(u,v)_{H^s\left(\mathbb{R}^m\right)}=\frac{1}{(2\pi)^m}\int_{\mathbb{R}^m}\left(1+|\xi|^2\right)^{s}\widehat{u}(\xi)\overline{\widehat{v}(\xi)}d\xi.$$
We leave the reader to verify that $H^s\left(\mathbb{R}^m\right)$
is a  Hilbert space.

\bigskip

\begin{theo}\label{Sob:teo2.11}
	If $0<s<1$ then the norm $\left\Vert
	u\right\Vert_{H^s\left(\mathbb{R}^m\right)}$ is equivalent to the 
	norm
	\begin{equation}\label{Sob:for1.123}
		\left\Vert u\right\Vert=\left(\left\Vert
		u\right\Vert^2_{L^2\left(\mathbb{R}^m\right)}+\left\vert
		u\right\vert^2_{s,\mathbb{R}^m}\right)^{1/2},
	\end{equation}
	where \index{$\left\vert \cdot\right\vert_{s,\mathbb{R}^m}$}
	\begin{equation}\label{Sob:for2.123}
		\left\vert u\right\vert^2_{s,\mathbb{R}^m}=\int_{\mathbb{R}^m}dx
		\int_{\mathbb{R}^m}\frac{|u(x)-u(y)|^2}{|x-y|^{m+2s}}dy.
	\end{equation}
\end{theo}
\textbf{Proof.} We need to prove that there exists $C\geq 1$ such that, for every $u\in H^s\left(\mathbb{R}^m\right)$ we have

\begin{equation}\label{Sob:for3.123}
	C^{-1}\int_{\mathbb{R}^m}|\xi|^{2s}\left|\widehat{u}(\xi)\right|^2d\xi\leq\left\vert
	u\right\vert^2_{s,\mathbb{R}^m}\leq
	C\int_{\mathbb{R}^m}|\xi|^{2s}\left|\widehat{u}(\xi)\right|^2d\xi.
\end{equation}
Let us begin by observing that
\begin{equation}\label{Sob:for4.123}
	\begin{aligned}
		\left\vert u\right\vert^2_{s,\mathbb{R}^m}&=\int_{\mathbb{R}^m}dx
		\int_{\mathbb{R}^m}\frac{|u(x+z)-u(x)|^2}{|z|^{m+2s}}dz=\\&=
		\frac{1}{(2\pi)^m}\int_{\mathbb{R}^m}\frac{1}{|z|^{m+2s}}dz\int_{\mathbb{R}^m}\left|\widehat{u(\cdot+z)}-\widehat{u(\cdot)}\right|^2d\xi=\\&=
		\frac{1}{(2\pi)^m}\int_{\mathbb{R}^m}\frac{1}{|z|^{m+2s}}dz\int_{\mathbb{R}^m}\left|e^{iz\cdot\xi}-1\right|^2
		\left|\widehat{u}(\xi)\right|^2d\xi=\\&=\frac{1}{(2\pi)^m}\int_{\mathbb{R}^m}\phi(\xi)
		\left|\widehat{u}(\xi)\right|^2d\xi,
	\end{aligned}
\end{equation}
where

$$\phi(\xi)=\int_{\mathbb{R}^m}\frac{\left|e^{iz\cdot\xi}-1\right|^2}{|z|^{m+2s}}dz.$$
Notice that $\phi$ is a homogeneous with degree $2s$.
As a matter of fact, for any $t>0$, we have

\begin{equation*}
	\begin{aligned}
		\phi(t\xi)&=\int_{\mathbb{R}^m}\frac{\left|e^{itz\cdot\xi}-1\right|^2}{|z|^{m+2s}}dz=\\&=
		\int_{\mathbb{R}^m}\frac{\left|e^{iy\cdot\xi}-1\right|^2}{|t^{-1}y|^{m+2s}}\frac{dy}{t^m}=
		t^{2s}\int_{\mathbb{R}^m}\frac{\left|e^{iy\cdot\xi}-1\right|^2}{|y|^{m+2s}}dy=\\&=
		t^{2s}\phi(t\xi), \quad \forall \xi\in \mathbb{R}^m.
	\end{aligned}
\end{equation*}
Moreover, $\phi$ is a continuous function in $\mathbb{R}^m$. In order to prove this, let $\xi_0\in \mathbb{R}^m$ and let us check that

\begin{equation}\label{Sob:correz.1}
	\lim_{\xi\rightarrow\xi_0}\int_{\mathbb{R}^m}\frac{\left|e^{iz\cdot\xi}-1\right|^2}{|z|^{m+2s}}dz=\int_{\mathbb{R}^m}\frac{\left|e^{iz\cdot\xi_0}-1\right|^2}{|z|^{m+2s}}dz.
\end{equation}
We have

$$\lim_{\xi\rightarrow\xi_0}\frac{\left|e^{iz\cdot\xi}-1\right|^2}{|z|^{m+2s}}=\frac{\left|e^{iz\cdot\xi_0}-1\right|^2}{|z|^{m+2s}},\quad \forall \xi\in \mathbb{R}^m$$
and, if $|\xi-\xi_0|<1$, we have
\begin{equation*}
	\begin{aligned}
		\frac{\left|e^{itz\cdot\xi}-1\right|^2}{|z|^{m+2s}}&=\frac{\left|e^{itz\cdot\xi}-1\right|^2}{|z|^{m+2s}}\chi_{B_1}(z)+
		\frac{\left|e^{itz\cdot\xi}-1\right|^2}{|z|^{m+2s}}\chi_{\mathbb{R}^m\setminus
			B_1}(z)\leq \\&\leq
		\frac{C(1+\left|\xi_0\right|)^2}{|z|^{m-2(1-s)}}\chi_{B_1}(z)+\frac{4}{|z|^{m+2s}}\chi_{\mathbb{R}^m\setminus
			B_1}(z)\in L^1\left(\mathbb{R}^m\right).
	\end{aligned}
\end{equation*}
Therefore by the Dominated Convergence Theorem we get \eqref{Sob:correz.1}.

Now, since $\phi$ is continuous in $\mathbb{R}^m$ and $\phi(\xi)>0$, for every $|\xi|=1$, we have that there exists
$C\geq 1$ such that
$$C^{-1}|\xi|^{2s}\leq \phi(\xi)\leq C|\xi|^{2s},\quad\forall\xi\in
\mathbb{R}^m.$$ By the last inequality and by \eqref{Sob:for4.123} we obtain
\eqref{Sob:for3.123}. $\blacksquare$

\bigskip

Similarly to the previous Theorem the following one can be proved
\begin{theo}\label{Sob:teo3.11}
	If $s>0$, $s\notin \mathbb{N}$, then  the norm
	$\left\Vert u\right\Vert_{H^s\left(\mathbb{R}^m\right)}$ is
	equivalent to the norm
	\begin{equation}\label{Sob:for1.123r}
		\left\Vert u\right\Vert=\left(\left\Vert
		u\right\Vert^2_{H^{[s]}\left(\mathbb{R}^m\right)}+\left\vert
		u\right\vert^2_{s,\mathbb{R}^m}\right)^{1/2},
	\end{equation}
	where
	\begin{equation}\label{Sob:for2.123r}
		\left\vert u\right\vert^2_{s,\mathbb{R}^m}=\sum_{|\alpha|=[s]}\int_{\mathbb{R}^m}dx
		\int_{\mathbb{R}^m}\frac{|\partial^{\alpha}u(x)-\partial^{\alpha}u(y)|^2}{|x-y|^{m+2(s-[s])}}dy.
	\end{equation}
\end{theo}

\bigskip

Theorems \ref{Sob:teo2.11}, \ref{Sob:teo3.11} justify the following definition.

\begin{definition}\label{Sob:def2.11} Let $\Theta$ be a bounded open set of $\mathbb{R}^m$ of class $C^{0,1}$. Let $s\notin \mathbb{N}$ be a positive real number. We define $H^s(\Theta)$ as the space of functions
	$u\in H^{[s]}(\Theta)$ such that
	\begin{equation*}
		\left\vert u\right\vert^2_{s,\Theta}=\sum_{|\alpha|=[s]}\int_{\Theta}dx
		\int_{\Theta}\frac{|\partial^{\alpha}u(x)-\partial^{\alpha}u(y)|^2}{|x-y|^{m+2(s-[s])}}dy<+\infty,
	\end{equation*}
	equipped with the norm
	
	$$\left\Vert u\right\Vert_{H^s(\Theta)}=\left(\left\Vert u\right\Vert_{H^{[s]}(\Theta)}+\left\vert
	u\right\vert^2_{s,\Theta}\right)^{1/2}.$$
\end{definition}

\medskip

It is not difficult to prove that the space $H^s(\Theta)$ is
complete.

\medskip

Now let us define $H^s(\partial\Omega)$, \index{$H^s(\partial\Omega)$} where $\Omega$ is
a bounded open set. If $s\in \mathbb{N}$, we assume that $\partial\Omega$ is of class $C^{s}$. If $s\notin \mathbb{N}$ we assume that  $\partial \Omega$ is of class $C^{[s],1}$.

We proceed basically as we did in Section \ref{DescrBordo} to define $L^p(\partial \Omega)$ (we will use the same notations as Section \ref{DescrBordo}).

Let us begin by the case \index{$\left\Vert \cdot \right\Vert_{H^s(\partial\Omega)}$ ($s>0$ real number)}  $s:=k$, positive integer number. Thus, let us assume $\partial\Omega$ of class $C^k$ with constants
$r_0,M_0$ and let us cover $\partial\Omega$ by a finite number, $N$,
of cylinders $\widetilde{Q}_{r_0,2M_0}(X_i)$, $i=1\cdots, N$,  where
$X_i\in\partial\Omega$ isometric to $Q_{r_0,2M_0}$.  Moreover let us assume that:
$\Sigma_i:=\widetilde{Q}_{r_0,2M_0}(X_i)\cap \partial\Omega$, for any
$i=1\cdots, N$, up to isometry for which $X_i$ is mapped in $0$,  is the graph of a function $\varphi_i\in
C^k(\overline{B'}_{r_0})$ likewise Definition \ref{Sob:def3.1}. We say that $f\in H^k(\partial \Omega)$ provided that the functions
$f(x',\varphi_i(x'))$, $i=1,\cdots, N$, belong to $H^k(B'_{r_0})$ and we denote

\begin{equation*}
	\left\Vert f\right\Vert_{H^k(\partial\Omega)}=\left(\sum_{i=1}^N
	\left\Vert f\right\Vert^2_{H^k(\Sigma_i)}\right)^{1/2},
\end{equation*}
where
$$\left\Vert f\right\Vert_{H^k(\Sigma_i)}=\left\Vert f(\cdot,\varphi_i(\cdot))\right\Vert_{H^k(B_{r_0})}.$$ Of course, the norm $\left\Vert \cdot\right\Vert_{H^k(\partial\Omega)}$ depends on the particular family of cylinders that we use as a covering, but they are all equivalent norms. Moreover,
$H^k(\partial\Omega)$ is a separable Hilbert space.

If $s\notin \mathbb{N}$, we say that $f\in H^s(\partial\Omega)$ provided that the
functions $f(x',\varphi_i(x'))$ belong to $H^s(B'_r)$, for any
$i=1,\cdots, N$ and we denote

\begin{equation*}
	\left\Vert f\right\Vert_{H^s(\partial\Omega)}=\left(\sum_{i=1}^N
	\left\Vert f\right\Vert^2_{H^s(\Sigma_i)}\right)^{1/2}.
\end{equation*}
The space $H^s(\partial\Omega)$ is complete and can be equipped of a Hilbert structure.

\bigskip

For an extended discussion of Sobolev spaces with noninteger exponent 
we refer the reader to \cite[Cap. 6]{Fo}, \cite[Cap. 6]{K-J-F},
\cite[Cap. 2]{Ne}. Here we limit ourselves to prove a Theorem that will be
useful in the next Section.

\bigskip

\begin{theo}[\textbf{density of $C^{\infty}(\mathbb{R}^m)$ in $H^s(\mathbb{R}^m)$}]\label{Sob:densita-Hs}
\index{Theorem:@{Theorem:}!- density of $C^{\infty}(\mathbb{R}^m)$ in $H^s(\mathbb{R}^m)$@{- density of $C^{\infty}(\mathbb{R}^m)$ in $H^s(\mathbb{R}^m)$}}	If $s$ is a positive real number then
	$C^{\infty}_0\left(\mathbb{R}^m\right)$ is dense in
	$H^{s}\left(\mathbb{R}^m\right)$.
\end{theo}

\bigskip

We premise the following

\begin{lem}\label{Sob:Lemma139}
	Let $\eta\in C^{\infty}_0\left(\mathbb{R}^m\right)$ satisfy
	
	\smallskip
	
	\noindent(i) supp $\eta\subset B_1$,
	
	\smallskip
	
	\noindent(ii) $\eta\geq 0$,
	
	\smallskip
	
	\noindent(iii) $\int_{\mathbb{R}^n} \eta(x)dx=1$.
	
	\smallskip
	
	Let us denote, for any $\varepsilon>0$, $v\in H^s\left(\mathbb{R}^m\right)$,
	$s>0$,
	
	$$\eta_{\varepsilon}(x)=\varepsilon^{-m}\eta\left(\varepsilon^{-1}x\right)$$
	and
	$$v^{\varepsilon}=\eta_{\varepsilon}\star v,\quad\mbox{ in }
	\mathbb{R}^m.$$ Then we have
	\begin{equation}\label{Sob:for139.7}
		\lim_{\varepsilon\rightarrow 0}\left\Vert \eta_{\varepsilon}\star
		v-v \right\Vert_{H^s\left(\mathbb{R}^m\right)}=0.
	\end{equation}
	\end{lem}

\textbf{Proof of the Lemma.} We have
$$\widehat{\eta_{\varepsilon}}(\xi)=\int_{\mathbb{R}^m}\varepsilon^{-m}\eta\left(\varepsilon^{-1}x\right)e^{-ix\cdot\xi}dx=\int_{\mathbb{R}^m}\eta\left(y\right)e^{-i\varepsilon y\cdot\xi}dy=\widehat{\eta}(\varepsilon\xi).$$
Moreover
\begin{equation}\label{Sob:for139.8}
	\lim_{\varepsilon\rightarrow
		0}\widehat{\eta}(\varepsilon\xi)=\widehat{\eta}(0)=\int_{\mathbb{R}^m}\eta(x)dx=1
\end{equation} and
\begin{equation}\label{Sob:for139.9}\left\vert\widehat{\eta}(\varepsilon\xi)\right\vert\leq
	\int_{\mathbb{R}^m}\eta(x)dx=1,\quad\forall \xi\in
	\mathbb{R}^m\mbox{, } \forall \varepsilon>0.
\end{equation}
Hence, for any $v\in H^s\left(\mathbb{R}^m\right)$, we have

\begin{equation*}
	\left\Vert \eta_{\varepsilon}\star v-v
	\right\Vert^2_{H^s\left(\mathbb{R}^m\right)}=\frac{1}{(2\pi)^m}\int_{\mathbb{R}^m}\left(1+|\xi|^2\right)^{s}\left\vert\widehat{\eta}(\varepsilon\xi)-1\right\vert^2
	\left|\widehat{v}(\xi)\right|^2d\xi.
\end{equation*}
By the last equality, by the Dominated Convergence Theorem
(take into account \eqref{Sob:for139.8} and \eqref{Sob:for139.9}) we obtain
\eqref{Sob:for139.7}. $\blacksquare$

\bigskip

\textbf{Proof of Theorem \ref{Sob:densita-Hs}.} The case where
$s\in \mathbb{N}$ has been proved in Proposition
\ref{Sob:prop2.4}. Let us consider the case $0<s<1$ (if $s>1$, the proof proceeds in a similar way, and we leave the details to the reader).

\bigskip

\textbf{Claim.} Let us denote by $\mathcal{H}_s$ the subspace of the
functions of $H^{s}\left(\mathbb{R}^m\right)$ with compact support. Then $\mathcal{H}_s$ is dense in
$H^{s}\left(\mathbb{R}^m\right)$.

\bigskip

\textbf{Proof of the Claim.}

Let $\zeta\in C_0^{\infty}(\mathbb{R}^m)$ satisfy
$$0\leq \zeta\leq 1,\quad\mbox{in } \mathbb{R}^m,$$
$$\zeta(x)=1,\quad \forall x\in B_1,\quad \zeta(x)=0,\quad \forall x\in \mathbb{R}^m\setminus
B_{2}$$ and
$$\left|\nabla \zeta\right|\leq C_0,\quad\mbox{in } \mathbb{R}^m,$$
where $C_0$ is a constant.  Let $R>1$ and 
$$\zeta_R(x)=\zeta\left(R^{-1}x\right).$$
We have
$$\left\Vert u-\zeta_Ru
\right\Vert_{L^{2}(\mathbb{R}^{m})}\leq \left\Vert u
\right\Vert_{L^{2}(\mathbb{R}^{m}\setminus B_R)}\rightarrow
0,\quad\mbox{as } R\rightarrow\infty.$$ Now, we prove 
\begin{equation}\label{Sob:for139.1}
	\lim_{R\rightarrow\infty} \left\vert
	u-u\zeta_R\right\vert_{s,\mathbb{R}^m}=0.
\end{equation}
In proving the latter, we will obtain as a by-product of the
performed calculations that $u\zeta_R$ belongs to $H^s\left(\mathbb{R}^m\right)$.

We apply Theorem \ref{Sob:teo2.11} and we write

\begin{equation*}
	\left\vert
	u-u\zeta_R\right\vert^2_{s,\mathbb{R}^m}=\int_{\mathbb{R}^m}dx
	\int_{\mathbb{R}^m}\frac{\Phi^2_R(x,y)}{|x-y|^{m+2s}}dy,
\end{equation*}
where
$$\Phi_R(x,y)=\left|\left(1-\zeta_R(x)\right)u(x)-\left(1-\zeta_R(y)\right)u(y)\right|.$$
We have
\begin{equation*}
	\begin{aligned}
		\Phi_R(x,y)&\leq
		\left|1-\zeta_R(y)\right|\left|u(x)-u(y)\right|+\left|\zeta_R(x)-\zeta_R(y)\right|
		\left|u(x)\right|\leq\\&\leq \chi_{\mathbb{R}^m\setminus
			B_R}(y)\left|u(x)-u(y)\right|+\left|\zeta_R(x)-\zeta_R(y)\right|
		\left|u(x)\right|.
	\end{aligned}
\end{equation*}
Hence
\begin{equation}\label{Sob:for139.2}
	\begin{aligned}
		\left\vert u-u\zeta_R\right\vert^2_{s,\mathbb{R}^m}&\leq
		2\int_{\mathbb{R}^m}dx\int_{\mathbb{R}^m\setminus
			B_R}\frac{\left|u(x)-u(y)\right|^2}{|x-y|^{m+2s}}dy+\\&+
		2\int_{\mathbb{R}^m}dx\int_{\mathbb{R}^m}\frac{\left|\zeta_R(x)-\zeta_R(y)\right|^2\left|u(x)\right|^2}{|x-y|^{m+2s}}dy.
	\end{aligned}
\end{equation}
Set
$$I:=\int_{\mathbb{R}^m}dx\int_{\mathbb{R}^m}\frac{\left|\zeta_R(x)-\zeta_R(y)\right|^2\left|u(x)\right|^2}{|x-y|^{m+2s}}dy.$$
Let us notice that
$$\left|\zeta_R(x)-\zeta_R(y)\right|=0,\quad \mbox{ for } |x|\geq
2R\mbox{ and } |y|\geq 2R$$ moreover, 

$$\left|\zeta_R(x)-\zeta_R(y)\right|\leq \min\left\{2,\frac{C_0|x-y|}{R}\right\},\quad \forall x\in \mathbb{R}^m, \forall y\in \mathbb{R}^m.$$
Hence
\begin{equation}\label{Sob:for139.3}
	\begin{aligned}
		I&\leq \int_{|x|\leq
			2R}dx\int_{\mathbb{R}^m}\frac{\left|\zeta_R(x)-\zeta_R(y)\right|^2\left|u(x)\right|^2}{|x-y|^{m+2s}}dy+\\&+\int_{\mathbb{R}^m}dx\int_{|y|\leq
			2R}\frac{\left|\zeta_R(x)-\zeta_R(y)\right|^2\left|u(x)\right|^2}{|x-y|^{m+2s}}dy:=I_1+I_2.
	\end{aligned}
\end{equation}
We have
\begin{equation}\label{Sob:for139.4}
	\begin{aligned}
		I_1&\leq \frac{C_0^2}{R^2}\int_{|x|\leq
			2R}\left|u(x)\right|^2dx\int_{|x-y|\leq
			\frac{2R}{C_0}}\frac{dy}{|x-y|^{m-2+2s}}+\\&+ 4\int_{|x|\leq
			2R}\left|u(x)\right|^2dx\int_{|x-y|>
			\frac{2R}{C_0}}\frac{dy}{|x-y|^{m+2s}}\leq\\&\leq
		\frac{\omega_m}{2-2s}\frac{C_0^2}{R^2}\left(\frac{2R}{C_0}\right)^{2-2s}\int_{|x|\leq
			2R}\left|u(x)\right|^2dx+\\&+\frac{\omega_m}{2s}\left(\frac{2R}{C_0}\right)^{-2s}\int_{|x|\leq
			2R}\left|u(x)\right|^2dx\leq \\&\leq CR^{-2s}\left\Vert
		u\right\Vert^2_{L^2\left(\mathbb{R}^m\right)},
	\end{aligned}
\end{equation}
where $C$ depends by $m$ and $s$ only.

\smallskip

Concerning $I_2$, we have

\begin{equation*}
	\begin{aligned}
		I_2&=\int_{|y|\leq
			2R}dy\int_{\mathbb{R}^m}\frac{\left|\zeta_R(x)-\zeta_R(y)\right|^2\left|u(x)\right|^2}{|x-y|^{m+2s}}dx=\\&=
		\int_{|y|\leq 2R}dy\int_{|x-y|\leq
			\frac{2R}{C_0}}\frac{\left|\zeta_R(x)-\zeta_R(y)\right|^2\left|u(x)\right|^2}{|x-y|^{m+2s}}dx+\\&+
		\int_{|y|\leq 2R}dy\int_{|x-y|>
			\frac{2R}{C_0}}\frac{\left|\zeta_R(x)-\zeta_R(y)\right|^2\left|u(x)\right|^2}{|x-y|^{m+2s}}dx\leq\\&\leq
		\frac{C_0^2}{R^2}\int_{|y|\leq 2R}dy\int_{|x|\leq
			2R(1+1/C_0)}\frac{\left|u(x)\right|^2}{|x-y|^{m-2+2s}}dx+\\&+
		4\left(\frac{2R}{C_0}\right)^{-m-2s}\int_{|y|\leq
			2R}dy\int_{|x-y|>
			\frac{2R}{C_0}}\left|u(x)\right|^2dx.
	\end{aligned}
\end{equation*}
By interchanging the order of integration in the second--to--last integral and trivially estimating from above the last integral, we obtain

\begin{equation}\label{Sob:for139.5}
	\begin{aligned}
		I_2&\leq \frac{C_0^2}{R^2}\int_{|x|\leq
			2R(1+1/C_0)}\left|u(x)\right|^2dx\int_{|y|\leq
			2R}\frac{dy}{|x-y|^{m-2+2s}}+\\&+
		4\frac{\omega_m}{m}\left(\frac{2R}{C_0}\right)^{-m-2s}(2R)^{m}\left\Vert
		u\right\Vert^2_{L^2\left(\mathbb{R}^m\right)}\leq\\&\leq
		\frac{C_0^2}{R^2}\int_{|x|\leq
			2R(1+1/C_0)}\left|u(x)\right|^2dx\int_{|y-x|\leq
			2R(2+1/C_0)}\frac{dy}{|x-y|^{m-2+2s}}+\\&+CR^{-2s}\left\Vert
		u\right\Vert^2_{L^2\left(\mathbb{R}^m\right)}\leq\\&\leq
		C'R^{-2s}\left\Vert u\right\Vert^2_{L^2\left(\mathbb{R}^m\right)},
	\end{aligned}
\end{equation}
where $C$ e $C'$ depend by $m$ and $s$ only.

By \eqref{Sob:for139.3}, \eqref{Sob:for139.4} and \eqref{Sob:for139.5}
we get

\begin{equation*}
	I=I_1+I_2\leq CR^{-2s}\left\Vert
	u\right\Vert^2_{L^2\left(\mathbb{R}^m\right)},
\end{equation*}
where $C$ depends by $s$ and $m$ only. 

By the just obtained inequality and by \eqref{Sob:for139.2} we have

\begin{equation*}
	\begin{aligned}
		\left\vert u-u\zeta_R\right\vert^2_{s,\mathbb{R}^m}&\leq
		2\int_{\mathbb{R}^m}dx\int_{\mathbb{R}^m\setminus
			B_R}\frac{\left|u(x)-u(y)\right|^2}{|x-y|^{m+2s}}dy+CR^{-2s}\left\Vert
		u\right\Vert^2_{L^2\left(\mathbb{R}^m\right)}.
	\end{aligned}
\end{equation*}
Since the following function belongs to $L^2\left(\mathbb{R}^m\times \mathbb{R}^m\right)$ (as $u\in H^s\left(\mathbb{R}^m\right)$) 
$$\mathbb{R}^m\times \mathbb{R}^m\ni(x,y)\rightarrow \frac{\left|u(x)-u(y)\right|^2}{|x-y|^{m+2s}},$$
we obtain \eqref{Sob:for139.1}. The \textbf{Claim} is proved.

\medskip

Now, let $\delta>0$ and let 
$R_0$ be (recall \eqref{Sob:for139.1}) such that

\begin{equation*}
	\left\Vert
	u-u\zeta_{R_0}\right\Vert_{s,\mathbb{R}^m}<\frac{\delta}{2}.
\end{equation*}
Lemma \ref{Sob:Lemma139} implies that there exists $\varepsilon_0>0$
such that
\begin{equation*}
	\left\Vert
	u\zeta_{R_0}-\left(u\zeta_{R_0}\right)\star\eta_{\varepsilon_0}\right\Vert_{s,\mathbb{R}^m}<\frac{\delta}{2}.
\end{equation*}
Hence
\begin{equation*}
	\left\Vert
	u-\left(u\zeta_{R_0}\right)\star\eta_{\varepsilon_0}\right\Vert_{s,\mathbb{R}^m}<\delta
\end{equation*}
Since $\left(u\zeta_{R_0}\right)\star\eta_{\varepsilon_0}\in
C^{\infty}_0\left(\mathbb{R}^m\right)$, the last inequality concludes the proof. $\blacksquare$

\subsection{The Theorem of characterization of the traces.} \label{Sob:tracce}
We preliminarily examine the extension of the notion of trace of a function belonging to $H^1\left(\mathbb{R}^n_+\right)$, where $\mathbb{R}^n_+=\left\{x=(x',x_n)\in
\mathbb{R}^n:\mbox{ } x_n>0 \right\}$.

If $u\in
H^1\left(\mathbb{R}^n_+\right)$, then the function

$$\widetilde{u}\left(x',x_n\right):=u \left(x',|x_n|\right),$$
belongs to $H^1\left(\mathbb{R}^n\right)$ and, as
$H^1\left(\mathbb{R}^n\right)=H_0^1\left(\mathbb{R}^n\right)$,
there exists a sequence $\left\{v_j\right\}$ in
$C^{\infty}_0\left(\mathbb{R}^n\right)$ such that
$$\left\{v_j\right\}\rightarrow \widetilde{u},\quad\mbox{ in } H^1\left(\mathbb{R}^n\right).$$
Hence

$$\left\{(v_j)_{|\mathbb{R}^n_+}\right\}\rightarrow u,\quad\mbox{ in }
H^1\left(\mathbb{R}^n_+\right).$$ Denoting
$$w=v_j-v_k,\quad j,k\in \mathbb{N},$$
we have

$$w\left(x',0\right)=w\left(x',x_n\right)-\int^{x_n}_0\partial_{x_n}w\left(x',y\right)dy,\quad
\forall x_n>0,$$ from which we obtain

$$\left|w\left(x',0\right)\right|^2\leq
2\left|w\left(x',x_n\right)\right|^2+2x_n
\int^{x_n}_0\left|\partial_{x_n}w\left(x',y\right)\right|^2dy.$$ Now, integrating both the sides of the last inequality over $(0,\delta)$, $\delta>0$,  \\ w.r.t. $x_n$, we have

\begin{equation*}
	\begin{aligned}
		\delta\left|w\left(x',0\right)\right|^2&\leq
		2\int^{\delta}_0\left|w\left(x',x_n\right)\right|^2dx_n+2\delta
		\int^{\delta}_0dx_n\int^{x_n}_0\left|\partial_{x_n}w\left(x',y\right)\right|^2dy\leq\\&\leq
		2\int^{\delta}_0\left|w\left(x',x_n\right)\right|^2dx_n+2\delta^2
		\int^{\delta}_0\left|\partial_{x_n}w\left(x',y\right)\right|^2dy.
	\end{aligned}
\end{equation*}
Integrating both the sides of the last inequality over $\mathbb{R}^{n-1}$  w.r.t. $x'$,
we have

\begin{equation*}
	\int_{\mathbb{R}^{n-1}}\left|w\left(x',0\right)\right|^2dx'\leq
	\frac{2}{\delta}\int_{\mathbb{R}^{n}}\left|w\left(x\right)\right|^2dx+2\delta
	\int_{\mathbb{R}^{n}}\left|\nabla w\left(x\right)\right|^2dx.
\end{equation*}
Starting from this inequality we proceed as in the
proof of Theorem \ref{traccia} (inequality
\eqref{Sob:for2.59}) and we obtain the extension of the trace operator  from the space

$$C^{\infty}_{\ast}\left(\mathbb{R}^n_+\right):=\left\{ \widetilde{u}_{|\mathbb{R}^n_+}
:\mbox{ } \widetilde{u}\in C_0^{\infty}\left(\mathbb{R}^n\right)
\right\}$$ to the space $H^1\left(\mathbb{R}^n_+\right)$. In particular, we have  $$Tu\in
L^2\left(\mathbb{R}^{n-1}\right),\quad\forall u\in
H^1\left(\mathbb{R}^n_+\right)$$ and

\begin{equation}\label{Sob:correz.2}
	\left\Vert Tu\right\Vert_{L^2\left(\mathbb{R}^{n-1}\right)}\leq C\left\Vert u\right\Vert_{H^1\left(\mathbb{R}^{n}\right)},\quad \forall u\in
	H^1\left(\mathbb{R}^n_+\right),
\end{equation}
where $C$ depends on $n$ only.

\bigskip

The following Theorem provides a
characterization of the image of $H^1(\Omega)$ by means the trace operator. 
\begin{theo}[\textbf{characterization of the trace}]\label{Sob:teo4.11}
	\index{Theorem:@{Theorem:}!- characterization of the trace@{- characterization of the trace}}
	Let $\Omega$ be either a bounded open set of class $C^{0,1}$ with constants
	$M_0,r_0$, or $\Omega=\mathbb{R}^n_+$. Let
	
	$$T:H^1(\Omega)\rightarrow L^2(\partial\Omega)$$
	be the trace operator defined in Theorem \ref{traccia} (in the case	where $\Omega=\mathbb{R}^n_+$ the definition is given at the beginning of this Section).
	
	Then we have
	
	\begin{equation*}
		T\left(H^1(\Omega)\right)=H^{1/2}(\partial\Omega).
	\end{equation*}
	Moreover
	
	\noindent(i)
	
	\begin{equation*}
		\left\Vert T(u)\right\Vert_{H^{1/2}(\partial\Omega)}\leq C
		\left\Vert u\right\Vert_{H^{1}(\Omega)},\quad\forall u\in
		H^{1}(\Omega),
	\end{equation*}
	where $C$ depends on $M_0$, $r_0$ and $n$ only.
	
	\noindent(ii) There exists a bounded, linear map
	
	\begin{equation*}
		\mathcal{T}:H^{1/2}(\partial\Omega)\rightarrow H^{1}(\Omega)
	\end{equation*}
	such that
	
	\begin{equation*}
		T(\mathcal{T}(h))=h,\quad\forall h\in H^{1/2}(\partial\Omega).
	\end{equation*}
	In particular, denoting by $u=\mathcal{T}(h)$, we have
	\begin{equation*}
		\left\Vert u\right\Vert_{H^{1}(\Omega)}\leq C \left\Vert
		h\right\Vert_{H^{1/2}(\partial\Omega)},\quad\forall h\in
		H^{1/2}(\Omega),
	\end{equation*}
	where $C$ depends on $M_0$, $r_0$ and $n$ only.
	
\end{theo}

\medskip

In the general case the proof of Theorem \ref{Sob:teo4.11}, is quite technical. 
Here we limit ourselves to the case
$\Omega=\mathbb{R}^n_+$. A
complete treatment (including the traces of the functions belonging to
$W^{k,p}(\Omega)$, $k\in \mathbb{N}$) can be founded in \cite[Ch. 6]{K-J-F} and in
\cite[Ch. 2, Secs. 2.3 -- 2.5]{Ne}.

\bigskip

\textbf{Proof of Theorem \ref{Sob:teo4.11} in the case
	$\Omega=\mathbb{R}^n_+$.}

We have proved, in \eqref{Sob:correz.2}, that $T\in L^{2}\left(\mathbb{R}^{n-1}\right)$. Now we prove $$Tu \in H^{1/2}\left(\mathbb{R}^{n-1}\right).$$

\smallskip

\textbf{Claim}

Let $v\in C^{\infty}_{\ast}\left(\mathbb{R}^n_+\right)$. Let us denote

$$h\left(x'\right)=v\left(x',0\right),\quad \forall x'\in \mathbb{R}^{n-1}.$$
We have

\begin{equation}\label{Sob:for1.130}
	\widehat{h}\left(\xi'\right)=\frac{1}{2\pi}\int_{\mathbb{R}}\widehat{v}\left(\xi',\xi_n\right)d\xi_n.
\end{equation}

\medskip

\textbf{Proof of the Claim.}

Since

$$v\left(x',x_n\right)=\frac{1}{(2\pi)^n}\int_{\mathbb{R}^n}\widehat{v}\left(\xi\right)e^{ix\cdot\xi}d\xi,$$
we have

\begin{equation*}
	\begin{aligned}
		h\left(x'\right)&=v\left(x',0\right)=\\&=\frac{1}{(2\pi)^n}\int_{\mathbb{R}^n}\widehat{v}\left(\xi\right)e^{ix'\cdot\xi'}d\xi=\\&=
		\frac{1}{(2\pi)^{n-1}}\int_{\mathbb{R}^{n-1}}\left(\frac{1}{2\pi}\int_{\mathbb{R}}\widehat{v}\left(\xi',\xi_n\right)d\xi_n\right)e^{ix'\cdot\xi'}d\xi'
	\end{aligned}
\end{equation*}
and by proposition (a) of Theorem \ref{Fourier:propI} we obtain
\eqref{Sob:for1.130}. Claim is proved.

\medskip

Now, by \eqref{Sob:for1.130} we have

\begin{equation*}
	\begin{aligned}
		&\int_{\mathbb{R}^{n-1}}\left\vert
		\widehat{h}(\xi')\right\vert^2\left(1+|\xi'|^2\right)^{1/2}d\xi'\leq\frac{1}{4\pi^2}\int_{\mathbb{R}^{n-1}}
		\left(\int_{\mathbb{R}}\left\vert\widehat{v}(\xi',\xi_n)\right\vert d\xi_n \right)^2\left(1+|\xi'|^2\right)^{1/2}d\xi' =\\&=
		\frac{1}{4\pi^2}\int_{\mathbb{R}^{n-1}}\left(1+|\xi'|^2\right)^{1/2}d\xi'\left(\int_{\mathbb{R}}\left(1+|\xi|^2\right)^{1/2}\left\vert\widehat{v}(\xi',\xi_n)\right\vert
		\left(1+|\xi|^2\right)^{-1/2} d\xi_n\right)^2\leq\\&\leq
		\frac{1}{4\pi^2}\int_{\mathbb{R}^{n-1}}\left(1+|\xi'|^2\right)^{1/2}d\xi'\left(\int_{\mathbb{R}}\left(1+|\xi|^2\right)\left\vert\widehat{v}(\xi',\xi_n)\right\vert^2
		d\xi_n \int_{\mathbb{R}}\left(1+|\xi|^2\right)^{-1}d\xi_n\right).
	\end{aligned}
\end{equation*}
Moreover
$$\int_{\mathbb{R}}\left(1+|\xi|^2\right)^{-1}d\xi_n=\int_{\mathbb{R}}\frac{d\xi_n}{1+|\xi'|^2+\xi_n^2}
=\frac{\pi}{\left(1+|\xi'|^2\right)^{1/2}}.$$ Therefore

$$\int_{\mathbb{R}^{n-1}}\left\vert
\widehat{h}(\xi')\right\vert^2\left(1+|\xi'|^2\right)^{1/2}d\xi'\leq
\frac{1}{4\pi}\int_{\mathbb{R}^n}\left(1+|\xi|^2\right)\left\vert\widehat{v}(\xi)\right\vert^2
d\xi$$ that is

\begin{equation*}
	\left\Vert T(v)\right\Vert_{H^{1/2}(\mathbb{R}^{n-1})}\leq
	\frac{1}{2\sqrt{\pi}} \left\Vert
	v\right\Vert_{H^{1}(\mathbb{R}^{n}_+)},\quad\forall v\in
	C^{\infty}_{\ast}\left(\mathbb{R}^n_+\right),
\end{equation*}
from which, by density we have
\begin{equation}\label{SOBBB}
	\left\Vert
	T(u)\right\Vert_{H^{1/2}\left(\mathbb{R}^{n-1}\right)}\leq
	\frac{1}{2\sqrt{\pi}} \left\Vert
	u\right\Vert_{H^{1}\left(\mathbb{R}^{n}\right)},\quad\forall u\in
	H^{1}\left(\mathbb{R}^{n}_+\right).
\end{equation}

\medskip

Now we prove (ii).

Let $h\in C^{\infty}_0\left(\mathbb{R}^{n-1}\right)$ and, for any
$\varepsilon>0$, let

\begin{equation*}
	u_{\varepsilon}\left(x',x_n\right)=\frac{1}{(2\pi)^{n-1}}\int_{\mathbb{R}^{n-1}}e^{-\left(1+\left|\xi'\right|\right)\left(x_n+\varepsilon\right)}
	\widehat{h}\left(\xi'\right)e^{ix'\cdot\xi'}d\xi',\quad \forall x\in
	\overline{\mathbb{R}^n_+}.
\end{equation*}
By applying  Theorem \ref{Fourier:propI} and by performing the derivative under the integral sign, it can be easily checked that

$$u_{\varepsilon}\in
C^{\infty}\left(\overline{\mathbb{R}^n_+}\right)$$ and

\begin{equation*}
	u_{\varepsilon}\left(x',0\right)=\frac{1}{(2\pi)^{n-1}}\int_{\mathbb{R}^{n-1}}e^{-\left(1+\left|\xi'\right|\right)\varepsilon}
	\widehat{h}\left(\xi'\right)e^{ix'\cdot\xi'}d\xi',\quad \forall
	x'\in \mathbb{R}^{n-1}.
\end{equation*}

\medskip

Now we prove what follows

\smallskip

\noindent \textbf{(a)} $u_{\varepsilon}\in H^{1}\left(\mathbb{R}^n_+\right)$
and, denoting

\begin{equation*}
	u\left(x',x_n\right)=\frac{1}{(2\pi)^{n-1}}\int_{\mathbb{R}^{n-1}}e^{-\left(1+\left|\xi'\right|\right)x_n}
	\widehat{h}\left(\xi'\right)e^{ix'\cdot\xi'}d\xi',\quad \forall x\in
	\mathbb{R}^{n},
\end{equation*}
(let us notice that $u(\cdot,0)=h$)  we have, 

$$u_{\varepsilon}\in H^{1}\left(\mathbb{R}^n_+\right)\quad\mbox{ and }\quad u\in H^{1}\left(\mathbb{R}^n_+\right).$$
Moreover

\begin{equation}\label{Sob:for1.137}
	\begin{aligned}
		\left\Vert
		u_{\varepsilon}\right\Vert_{H^1\left(\mathbb{R}^{n}\right)}\leq
		C\left\Vert h\right\Vert_{H^{1/2}\left(\mathbb{R}^{n-1}\right)},
	\end{aligned}
\end{equation}
and
\begin{equation}\label{Sob:for1bis.137}
	\begin{aligned}
		\left\Vert u\right\Vert_{H^1\left(\mathbb{R}^{n}\right)}\leq
		C\left\Vert h\right\Vert_{H^{1/2}\left(\mathbb{R}^{n-1}\right)},
	\end{aligned}
\end{equation}
where $C$ depends on $n$ only.

\smallskip

\noindent \textbf{(b)} $$u_{\varepsilon}\rightarrow u,\quad\mbox{ as }
\varepsilon\rightarrow 0\mbox{ in }
H^{1}\left(\mathbb{R}^n_+\right).$$

\medskip

\textbf{Proof of (a).}

The Parseval identity implies

\begin{equation*}
	\int_{\mathbb{R}^{n-1}}\left|u_{\varepsilon}\left(x',x_n\right)\right|^2dx'=c_n\int_{\mathbb{R}^{n-1}}e^{-2\left(1+\left|\xi'\right|\right)\left(x_n+\varepsilon\right)}
	\left|\widehat{h}\left(\xi'\right)\right|^2d\xi',
\end{equation*}
where $c_n$ depends on $n$ only. Hence

\begin{equation*}
	\begin{aligned}
		\int_{\mathbb{R}^{n}_+}\left|u_{\varepsilon}\left(x',x_n\right)\right|^2dx'dx_n&=c_n\int_{\mathbb{R}^{n-1}}\left|\widehat{h}\left(\xi'\right)\right|^2d\xi'
		\left(\int^{+\infty}_0e^{-2\left(1+\left|\xi'\right|\right)\left(x_n+\varepsilon\right)}dx_n\right)=\\&=
		c_n\int_{\mathbb{R}^{n-1}}\left|\widehat{h}\left(\xi'\right)\right|^2\frac{e^{-2\varepsilon\left(1+\left|\xi'\right|\right)}}{2\left(1+\left|\xi'\right|\right)}d\xi'\leq\\&\leq
		\frac{(2\pi)^{n-1}c_n}{2}\left\Vert
		h\right\Vert^2_{L^2\left(\mathbb{R}^{n-1}\right)}.
	\end{aligned}
\end{equation*}
Therefore

\begin{equation}\label{Sob:for1.135}
	\left\Vert
	u_{\varepsilon}\right\Vert_{L^2\left(\mathbb{R}^{n}\right)}\leq
	C\left\Vert h\right\Vert_{L^2\left(\mathbb{R}^{n-1}\right)},
\end{equation}
where $C$ depends on $n$ only.

\medskip

Now we estimate from above $\left\Vert \nabla
u_{\varepsilon}\right\Vert_{L^2\left(\mathbb{R}^{n}\right)}$.
We have

$$\partial_{x_n} u_{\varepsilon}\left(x',x_n\right)=-\frac{1}{(2\pi)^n}\int_{\mathbb{R}^{n-1}}\left(1+\left|\xi'\right|\right)
e^{-\left(1+\left|\xi'\right|\right)\left(x_n+\varepsilon\right)}\widehat{h}\left(\xi'\right)e^{ix'\cdot\xi'}d\xi'.$$
Hence, arguing as above, we get

\begin{equation}\label{Sob:for1.136}
	\begin{aligned}
		&\int_{\mathbb{R}^{n}_+}\left|\partial_{x_n}u_{\varepsilon}\left(x',x_n\right)\right|^2dx'dx_n=\\&=c_n\int_{\mathbb{R}^{n-1}}\left|\widehat{h}\left(\xi'\right)\right|^2d\xi'
		\left(\int^{+\infty}_0\left(1+\left|\xi'\right|\right)^2
		e^{-2\left(1+\left|\xi'\right|\right)\left(x_n+\varepsilon\right)}dx_n\right)=\\&=
		\frac{c_n}{2}\int_{\mathbb{R}^{n-1}}\left(1+\left|\xi'\right|\right)\left|\widehat{h}\left(\xi'\right)\right|^2e^{-2\varepsilon\left(1+\left|\xi'\right|\right)}d\xi'\leq\\&\leq
		\frac{(2\pi)^{n-1}c_n}{2}\left\Vert
		h\right\Vert^2_{H^{1/2}\left(\mathbb{R}^{n-1}\right)}.
	\end{aligned}
\end{equation}

\smallskip

Now, if $1\leq j\leq n-1$, we have
$$\partial_{x_j}
u_{\varepsilon}\left(x',x_n\right)=\frac{1}{(2\pi)^n}\int_{\mathbb{R}^{n-1}}i\xi_j
e^{-2\left(1+\left|\xi'\right|\right)\left(x_n+\varepsilon\right)}\widehat{h}\left(\xi'\right)e^{ix'\cdot\xi'}d\xi'$$
and arguing as in \eqref{Sob:for1.136}, we get
\begin{equation}\label{Sob:for2.136}
	\begin{aligned}
		\int_{\mathbb{R}^{n}_+}\left|\partial_{x_j}u_{\varepsilon}\left(x',x_n\right)\right|^2dx'dx_n\leq
		\frac{(2\pi)^{n-1}c_n}{2}\left\Vert
		h\right\Vert^2_{H^{1/2}\left(\mathbb{R}^{n-1}\right)}.
	\end{aligned}
\end{equation}
By \eqref{Sob:for1.136} and \eqref{Sob:for2.136} we get
$$\left\Vert \nabla
u_{\varepsilon}\right\Vert^2_{L^2\left(\mathbb{R}^{n}\right)}\leq
C\left\Vert h\right\Vert^2_{H^{1/2}\left(\mathbb{R}^{n-1}\right)},$$
where $C$ depends on $n$ only. By the just obtained inequality and by
\eqref{Sob:for1.135} we derive

\begin{equation*}
	\begin{aligned}
		\left\Vert
		u_{\varepsilon}\right\Vert_{H^1\left(\mathbb{R}^{n}\right)}\leq
		C\left\Vert h\right\Vert_{H^{1/2}\left(\mathbb{R}^{n-1}\right)}
	\end{aligned}
\end{equation*}
where $C$ depends on $n$ only.

As can be easily observed, the calculations performed above also apply to $\varepsilon=0$ and
similarly yield
\begin{equation*}
	\begin{aligned}
		\left\Vert u\right\Vert_{H^1\left(\mathbb{R}^{n}\right)}\leq
		C\left\Vert h\right\Vert_{H^{1/2}\left(\mathbb{R}^{n-1}\right)},
	\end{aligned}
\end{equation*}
where $C$ depends on $n$ only. Hence we have proved
\eqref{Sob:for1.137} and \eqref{Sob:for1bis.137}. Proof of (a)
is concluded.

\medskip

\textbf{Proof of (b).}

Since
\begin{equation*}
	\begin{aligned}
	\left(u-u_{\varepsilon}\right)\left(x',x_n\right)=\frac{1}{(2\pi)^{n-1}}\int_{\mathbb{R}^{n-1}}e^{-\left(1+\left|\xi'\right|\right)x_n}\left(1-
	e^{-\varepsilon\left(1+\left|\xi'\right|\right)}\right)\widehat{h}\left(\xi'\right)e^{ix'\cdot\xi'}d\xi,
	\end{aligned}
\end{equation*}
we easily obtain

\begin{equation*}
	\begin{aligned}
		\left\Vert
		u-u_{\varepsilon}\right\Vert_{L^2\left(\mathbb{R}^{n}\right)}\leq
		C\int_{\mathbb{R}^{n-1}}\frac{\left(1-
			e^{-\varepsilon\left(1+\left|\xi'\right|\right)}\right)^2}{1+\left|\xi'\right|}\left|\widehat{h}\left(\xi'\right)\right|^2d\xi'
	\end{aligned}
\end{equation*}
and
\begin{equation*}
	\begin{aligned}
		\left\Vert
		\nabla\left(u-u_{\varepsilon}\right)\right\Vert_{L^2\left(\mathbb{R}^{n}\right)}\leq
		C\int_{\mathbb{R}^{n-1}}\left(1+\left|\xi'\right|\right)\left(1-
		e^{-\varepsilon\left(1+\left|\xi'\right|\right)}\right)^2\left|\widehat{h}\left(\xi'\right)\right|^2d\xi'.
	\end{aligned}
\end{equation*}
Therefore, by the Dominated Convergence Theorem we get

$$\lim_{\varepsilon\rightarrow 0}\left\Vert
u-u_{\varepsilon}\right\Vert^2_{H^1\left(\mathbb{R}^{n}\right)}=0.$$
By the previous limit and by the trace Theorem  (that is, by (i)) we have

\begin{equation}\label{Sob:for1.138}
	T(u)=\lim_{\varepsilon\rightarrow
		0}T\left(u_{\varepsilon}\right),\quad\mbox{ in }
	L^2\left(\mathbb{R}^{n-1}\right). \end{equation} On the other hand
\begin{equation}\label{Sob:for2.138}
	\begin{aligned}
		&\lim_{\varepsilon\rightarrow 0}\left\Vert
		T\left(u_{\varepsilon}\right)-h
		\right\Vert_{L^2\left(\mathbb{R}^{n-1}\right)}=\\&=\lim_{\varepsilon\rightarrow
			0}c_n\int_{\mathbb{R}^{n-1}}\left(1-
		e^{-\varepsilon\left(1+\left|\xi'\right|\right)}\right)^2\left|\widehat{h}\left(\xi'\right)\right|^2d\xi'=0
	\end{aligned}
\end{equation}
Therefore, by \eqref{Sob:for1.138} and \eqref{Sob:for2.138} we get

\begin{equation}\label{Sob:for3.138}
	T(u)=h,\quad\forall h\in C^{\infty}_0\left(\mathbb{R}^{n-1}\right).
\end{equation}

Now, we set

$$\mathcal{T}(h)=u.$$
By \eqref{Sob:for1bis.137} we have

\begin{equation}\label{SOBBB1}
	\begin{aligned}
		\left\Vert
		\mathcal{T}(h)\right\Vert_{H^1\left(\mathbb{R}^{n}\right)}\leq
		C\left\Vert
		h\right\Vert_{H^{1/2}\left(\mathbb{R}^{n-1}\right)},\quad \forall
		h\in C_0^{\infty}\left(\mathbb{R}^{n-1}\right)
	\end{aligned}
\end{equation}
and by \eqref{Sob:for3.138} we have

\begin{equation}\label{SOBBB2}
	T(\mathcal{T}(h))=h,\quad \forall h\in
	C_0^{\infty}\left(\mathbb{R}^{n-1}\right).
\end{equation}
Now, \eqref{SOBBB} and \eqref{SOBBB1} give

\begin{equation}\label{SOBBB3}
	\left\Vert
	T(\mathcal{T}(h))\right\Vert_{L^2\left(\mathbb{R}^{n}\right)}\leq
	C\left\Vert
	h\right\Vert_{H^{1/2}\left(\mathbb{R}^{n-1}\right)},\quad \forall
	h\in C^{\infty}\left(\mathbb{R}^{n-1}\right).\end{equation}
Finally, by density Theorem \ref{Sob:densita-Hs}, by \eqref{SOBBB1},
\eqref{SOBBB2} and by \eqref{SOBBB3} the thesis follows. $\blacksquare$

\section{Final comments and supplements}\label{complementi-Sob}

In this Section we will state, without proof, some theorems
concerning the traces and Lipschitz continuous functions.

Concerning the traces, if $k\geq 1$ and if $\Omega$ is a bounded open set of $\mathbb{R}^n$ of
class $C^{k-1,1}$, it can be proved that (see \cite{K-J-F},
\cite{Ne})
$$H^{k-1/2}(\partial\Omega)=T\left(H^k(\Omega)\right).$$ In the 
sequel we will be mainly interested in the cases  $k=1,2$. In
particular, if $u\in H^2(\Omega)$, then $\partial_ku\in
H^1(\Omega)$ for $k=1,\cdots,n$, hence
$\partial_ku_{|\partial\Omega}\in H^{1/2}(\partial\Omega)$. Moreover,
we have
$$H^{2}(\partial\Omega)\subset H^{3/2}(\partial\Omega) \subset
H^{1}(\partial\Omega).$$ If $u\in H^2(\Omega)$, we can define $\frac{\partial u}{\partial \nu}=\nabla u\cdot\nu$ on
$\partial\Omega$ ($\nu$ unit outward normal vector) and we have

\begin{equation}\label{traccia-1}
	\left\Vert \frac{\partial u}{\partial
		\nu}\right\Vert_{H^{1/2}(\partial\Omega)}\leq C\left\Vert
	u\right\Vert_{H^{2}(\Omega)}, \quad \forall u\in H^{2}(\Omega),
\end{equation}
where $C$ depends by $\Omega$ only.

\medskip

Similarly to what we saw in Section \ref{tracce}, the following Theorem can be
proved

\begin{theo} \label{traccia-10}
	Let $\Omega$ be a bounded open set of $\mathbb{R}^n$ of class $C^{1,1}$. If $u\in H^2(\Omega)$ then
	\begin{equation}\label{traccia nulla-1}
		u\in H_0^2(\Omega) \quad\mbox{ if and only if }\quad u_{|\partial\Omega}=0 \mbox{ and } \frac{\partial u}{\partial \nu}=0 \mbox{ on } \partial\Omega.
	\end{equation}
\end{theo}

\bigskip

Theorem \ref{Sob:teo4.11} can be generalized as follows 
\begin{theo}\label{traccia-inv-k}
	Let $\Omega$ be a bounded open set of $\mathbb{R}^n$ of class $C^{1,1}$. Then, for every
	$$(\psi_0,\psi_1)\in H^{3/2}(\partial\Omega)\times
	H^{1/2}(\partial\Omega)$$ there exists $u\in H^2(\Omega)$ such that
	$$u_{|\partial\Omega}=\psi_0 \quad\mbox{ and }\quad \frac{\partial u}{\partial \nu}=\psi\quad \mbox{ on }\quad \partial\Omega$$
	and
	$$\left\Vert u\right\Vert_{H^{2}(\Omega)} \leq C\left(\left\Vert \psi_0\right\Vert_{H^{3/2}(\partial\Omega)}+\left\Vert \psi_1\right\Vert_{H^{1/2}(\partial\Omega)}\right),$$
	where $C$ depends on $\Omega$ only.
\end{theo}

\bigskip






\subsection{The space $H^{-1/2}(\partial \Omega)$} \label{H--1/2}
\index{$H^{-1/2}(\partial \Omega)$}
Let $\Omega$ be a bounded open set of $\mathbb{R}^n$ of class $C^{0,1}$. We denote by $H^{-1/2}(\partial
\Omega)$ the dual space of $H^{1/2}(\partial \Omega)$. Thus,
$H^{-1/2}(\partial \Omega)$ is the space of the linear  functionals
$$\Phi:H^{1/2}(\partial \Omega)\rightarrow \mathbb{R}$$ such that for a constant $C$ we have
\begin{equation}\label{duale di H1/2}
	|\Phi(\varphi)|\leq C\left\Vert \varphi\right\Vert_{H^{1/2}(\partial\Omega)} \quad\forall  \varphi\in H^{1/2}(\partial\Omega).
\end{equation}

We define the norm, $\left\Vert
\Phi\right\Vert_{H^{-1/2}(\partial\Omega)}$, of $\Phi$ in
$H^{-1/2}(\partial \Omega)$ as the greatest lower bound of $C$  satisfying \eqref{duale di H1/2}.

\subsection{The space $W_{loc}^{1,\infty}\left(\mathbb{R}^n\right)$ and $C_{loc}^{0,1}\left(\mathbb{R}^n\right)$}  \label{sob:29-10-22}
 We say that $u\in C_{loc}^{0,1}\left(\mathbb{R}^n\right)$ provided that for any $x_0\in \mathbb{R}^n$ there exists $r>0$ such that 
$$u_{|\overline{B_r(x_0)}} \in C^{0,1}\left(\overline{B_r(x_0)}\right),$$

\bigskip

\begin{theo}\label{W-lip:29-10-22-1} Let $u:\mathbb{R}^n\rightarrow
	\mathbb{R}$. We have that $u\in W_{loc}^{1,\infty}\left(\mathbb{R}^n\right)$  if and only if
	$u\in C_{loc}^{0,1}\left(\mathbb{R}^n\right)$.
\end{theo}
\textbf{Proof.} Let  $u\in C_{loc}^{0,1}\left(\mathbb{R}^n\right)$. Theorem \ref{Rade} implies that $u$ is differentiable almost everywhere and $\nabla u=\left(\partial_1 u,\cdots \partial_nu\right)\in L_{loc}^{\infty}\left(\mathbb{R}^n\right)$. Hence, for any $1\leq k\leq n$, we have

$$\int_{\mathbb{R}^n}u\partial_k\varphi dx=\int_{\mathbb{R}^n}\left(\partial_k\left(u\varphi\right)-\partial_ku\varphi\right)dx=-\int_{\mathbb{R}^n}\partial_ku\varphi dx,\quad \forall \varphi\in C^{\infty}_0\left(\mathbb{R}^n\right).$$ Therefore $\partial_ku$ is the weak derivative of $u$, hence $u\in W_{loc}^{1,\infty}\left(\mathbb{R}^n\right)$. 

Conversely, let $u\in W_{loc}^{1,\infty} \left(\mathbb{R}^n\right)$ and let $B_{r}(x_0)$ be a ball of $\mathbb{R}^n$. Let us denote

$$ u_{\varepsilon}(x)= \int_{\mathbb{R}^n} \eta_{\varepsilon} (x-y) u(y) dy.$$

We have  $u_{\varepsilon}\in  C^{\infty}\left(\mathbb{R}^n\right)$ and 
$$\nabla u_{\varepsilon}(x)=\int_{\mathbb{R}^n} \eta_{\varepsilon} (x-y) \nabla u(y)dy.$$ Consequently, for every $\varepsilon \in (0,r]$, we get

$$\left|\nabla u_{\varepsilon}(x)\right|\leq \left\Vert\nabla u \right\Vert_{L^{\infty}\left(B_{2r}(x_0)\right)}<+\infty,\quad \forall x\in B_r(x_0).$$ 
Hence

\begin{equation}\label{L-W:29-10-22-3}
\left|u_{\varepsilon}(x)-u_{\varepsilon}(y)\right|\leq \left\Vert\nabla u \right\Vert_{L^{\infty}\left(B_{2r}(x_0)\right)} |x-y|,\quad\forall x,y\in B_r(x_0).
\end{equation}
Moreover, Theorem \ref{Sob:teo28R} implies that 
$$u_{\varepsilon}\rightarrow u,\quad  \mbox{as } \varepsilon\rightarrow 0 \ \ \mbox{ (uniformly). }$$ By this and by \eqref{L-W:29-10-22-3}, we have

\begin{equation}\label{L-W:29-10-22-4}
	\left|u(x)-u(y)\right|\leq \left\Vert\nabla u \right\Vert_{L^{\infty}\left(B_{2r}(x_0)\right)} |x-y|,\quad\forall x,y\in B_r(x_0).
\end{equation}
Therefore $u_{|\overline{B_r(x_0)}}\in C^{0,1}\left(\overline{B_r(x_0)}\right)$. Since $B_r(x_0)$ is arbitrary, the proof is complete. $\blacksquare$

\subsection{Almost everywhere differenziability of function belonging to  $W_{loc}^{1,p}\left(\mathbb{R}^n\right)$ with $p>n$.}  \label{sob:29-10-22-4}

\begin{theo}\label{sob:29-10-22-5}
	Let $n<p\leq +\infty$. If $u\in W_{loc}^{1,p}\left(\mathbb{R}^n\right)$, then  $u$ is almost everywhere differentiable in $\mathbb{R}^n$.
\end{theo}
\textbf{Proof.} Since $W_{loc}^{1,\infty}\left(\mathbb{R}^n\right) \subset W_{loc}^{1,p}\left(\mathbb{R}^n\right)$ for $p<+\infty$, we may assume that $n<p<+\infty$. Let $x\in \mathbb{R}^n$ be such that (Corollary \ref{Diff-Leb:29-10-22-1} )
\begin{equation}\label{Diff:30-10-22-1}
		\lim_{r\rightarrow 0}\ \dashint_{B_r(x)}|\nabla u(\xi)-\nabla u(x)|^pd\xi=0. 
\end{equation}
Let us denote
\begin{equation*}
	v(\xi)=u(\xi)-u(x)-\nabla u(x)\cdot (\xi-x),\quad \xi \in  \mathbb{R}^n,
\end{equation*}
let $y\in \mathbb{R}^n$, $y\neq x$ and $r=|x-y|$. Theorem \ref{densit 2} and Lemma \ref{Sob:lem7.7} give

\begin{equation*}
	\begin{aligned}
&|u(y)-u(x)-\nabla u(x)\cdot (y-x)|=|v(y)-v(x)|\leq\\&\leq  Cr^{1-\frac{n}{p}}\left(\int_{B_r(x)}|\nabla v(\xi)|^pd\xi\right)^{1/p}=\\&=
C'|x-y|\left(\dashint_{B_r(x)}|\nabla u(\xi)-\nabla u(x)|^pd\xi\right)^{1/p}.
\end{aligned}
\end{equation*}
Hence
$$u(y)-u(x)-\nabla u(x)\cdot (y-x)=o(|y-x|),\quad \mbox{ as } y\rightarrow x.$$
Therefore $u$ is differentiable in every point $x$ which satisfies \eqref{Diff:30-10-22-1} and by Corollary \ref{Diff-Leb:29-10-22-1} we conclude the proof. $\blacksquare$

\medskip 

\textbf{Remark.} From Theorems \ref{W-lip:29-10-22-1} and \ref{sob:29-10-22-5} one immediately obtains the Rademacher Theorem by a proof different from the one followed in Section \ref{funz-Lips}. The reader is invited to make sure that following this new proof does not lead to "vicious circles". $\blacklozenge$ 

\chapter{The boundary value problems for second order elliptic equations and the Dirichlet to Neumann map }\markboth{Chapter 4. The boundary value problems}{}
\label{Lax-Milgram}
\section{Introduction}\label{introd-Dirich}

Let $\Omega$ be a bounded open set of $\mathbb{R}^n$. Let us denote by
$\mathbb{M}(n)$ \index{$\mathbb{M}(n)$}the vector space of the matrices $n\times n$ whose entries are real numbers and let
$A\in L^{\infty}(\Omega;\mathbb{M}(n))$, \index{$L^{\infty}(\Omega;\mathbb{M}(n))$} i.e.
$A=\left\{a^{jk}\right\}_{j,k=1}^n$ is a matrix whose entries  $a^{jk}$ belong to $L^{\infty}(\Omega)$, for $j,k=1,\cdots,n$.

Throughout this Chapter we will assume that $A$ satisfies the following
condition of \textbf{uniform ellipticity} \index{uniform ellipticity}

\begin{equation}\label{gamma-eq}
	\lambda^{-1}|\xi|^2\leq A(x)\xi\cdot \xi \leq\lambda |\xi|^2
	\quad\mbox{ a.e. in } \Omega,\mbox{ } \forall \xi\in \mathbb{R}^n,
\end{equation}
where $\lambda\geq 1$ is a given number. We define 
\begin{equation*}
	|A(x)|_{\mathbb{M}(n)}=\left(\sum_{j,k=1}^n\left(a^{jk}(x)\right)^2\right)^{\frac{1}{2}}, \ \ \mbox{a.e. } x\in \Omega
\end{equation*}
and
\begin{equation*}
\left\Vert
A\right\Vert_{L^{\infty}(\Omega;\mathbb{M}(n))}=\left\Vert
|A(\cdot)|_{\mathbb{M}(n)}\right\Vert_{L^{\infty}(\Omega)}.
\end{equation*}

Let us notice that if $A$ is symmetric, then the second inequality of \eqref{gamma-eq} implies also
\begin{equation}\label{corr:25-4-23-1}
	\sup_{|\xi|=1,|\eta|=1}\left\vert
	A(x)\xi\cdot\eta \right\vert\leq\lambda, \ \ \forall x\in \Omega
\end{equation}
and
\begin{equation}\label{corr:25-4-23-2}
	\left\Vert
	a^{jk}\right\Vert_{L^{\infty}(\Omega)}\leq\lambda, \ \ \mbox{for } j,k=1,\cdots, n.
\end{equation}

Let us check \eqref{corr:25-4-23-1} and \eqref{corr:25-4-23-1}. Let $\xi,\eta\in \mathbb{R}^n$ such that $|\xi|=1$ and $|\eta|=1$. By the symmetry of $A(x)$ we have

\begin{equation*}
	\begin{aligned}
		|A(x)\xi\cdot \eta|&=\frac{1}{4}\left|A(x)(\xi+\eta)\cdot (\xi+\eta)-A(x)(\xi-\eta)\cdot (\xi\eta)\right|\leq \\&\leq \frac{1}{4}\left(\lambda |\xi+\eta|^2+\lambda |\xi-\eta|^2\right)=\\&=
		\frac{\lambda}{4}\left(2|\xi|^2+2|\eta|^2\right)=\\&=
		\lambda.
	\end{aligned}
\end{equation*}
Hence, \eqref{corr:25-4-23-1} follows. Concerning \eqref{corr:25-4-23-2}, we have, for $ j,k=1,\cdots, n$ 
\begin{equation*}
	\left\Vert
	a^{jk}(x)\right\Vert_{L^{\infty}(\Omega)}=|A(x)e_j\cdot e_k|\leq \lambda, \ \ \mbox{a.e. } x\in \Omega.
\end{equation*}

\smallskip

In this Chapter we will tackle the Dirichlet problem \index{Dirichlet problem} for the operator $-\mbox{div}(A \nabla u)$, which formally consists in determining $u\in H^{1}(\Omega)$ such that
\begin{equation}\label{Ell:Dirichlet}
	\begin{cases}
		-\mbox{div}(A\nabla u)=F, \quad\mbox{ in } \Omega, \\
		\\
		u=\varphi, \quad\mbox{ on } \partial\Omega,
	\end{cases}
\end{equation}
where $\varphi\in H^{1/2}(\partial\Omega) $ and $F\in H^{-1}(\Omega)$.
We will deal with the variational formulation of problem
\eqref{Ell:Dirichlet}. Next we will deal with the existence and the uniqueness of
the solutions in $H^{1}(\Omega)$ and subsequently we prove some
regularity results for the same problem, i.e., in coarse terms,
we will prove that if $\Omega$, $A$, and $F$ have greater
regularity, then $u$ also acquires more regularity.

The investigation on  problem \eqref{Ell:Dirichlet} will guide us
to deal with the more general case in which, instead of $-\mbox{div} (\nabla u)$, we will have the operator

\begin{equation}\label{Ell:intro1}
	-\sum_{j,k=1}^n\partial_j\left( a^{jk}\partial_k u
	+d^ju\right)+\sum_{j=1}^nb^{j}\partial_j u+cu,\end{equation} where $b^j,d^j, c\in
L^{\infty}(\Omega)$, for $j=1,\cdots, n$.

\section{The Lax--Milgram Theorem and the Fredholm Theorem}\label{Prob-Dirichlet}

Let $H$ be a real Hilbert space, let us denote by $\left\Vert
\cdot\right\Vert$ and $(\cdot,\cdot)$ the scalar product on $H$ and the induced norm induced respectively. As usual we denote by $H'$ the dual space of $H$.

If $A$ is a linear operator we denote by $\mathcal{R}(A)$ \index{$\mathcal{R}(A)$}
the range of $A$, that is
$$\mathcal{R}(A):=\left\{Au:\mbox{ } u\in H\right\}$$
and by $\mathcal{N}(A)$ \index{$\mathcal{N}(A)$}the kernel of $A$
$$\mathcal{N}(A):=\left\{u\in H:\mbox{ } Au=0\in H\right\}.$$

Let

\begin{equation}\label{forma biline}
	a:H\times H\rightarrow \mathbb{R},
\end{equation}
a bilinear form. We say that $a$ is \textbf{continuous} if there exists
$C>0$ such that

\begin{equation}\label{forma biline cont}
	|a(u,v)|\leq C \left\Vert u\right\Vert \left\Vert
	v\right\Vert,\quad\forall u,v\in H.
\end{equation}

We say that the bilinear form $a$ is \textbf{coercive} \index{coercive bilinear form}if there exists $\alpha>0$ such that

\begin{equation}\label{forma biline cont-1}
	\alpha\left\Vert u\right\Vert^2\leq a(u,u),\quad\forall
	u\in H.
\end{equation}

\medskip

The following Theorem holds true.

\begin{theo}[\textbf{Lax--Milgram}]\label{laxmilgram}
	\index{Theorem:@{Theorem:}!- Lax -- Milgram@{- Lax -- Milgram}}
	Let $a$ be a coercive bilinear form and let $F\in H'$.
	Then there exists a unique $u\in H$ such that
	\begin{equation}\label{lm}
		a(u,v)=F(v), \quad\forall v\in H.
	\end{equation}
	Moreover
	\begin{equation}\label{lm-cont}
		\left\Vert u\right\Vert\leq \frac{1}{\alpha}\left\Vert
		F\right\Vert_{H'},
	\end{equation}
	where $\left\Vert \cdot\right\Vert_{H'}$ is the norm of $H'$.
\end{theo}
\textbf{Proof.} 
The Riesz Representation Theorem implies that there exists a unique $f\in H$
such that
\begin{equation}\label{ELLlaxmilgram.2}
	F(v)=(f,v), \quad \forall v\in H.
\end{equation}
Moreover
\begin{equation*}
	\left\Vert F\right\Vert_{H'}=\left\Vert f\right\Vert_{H}.
\end{equation*}
Let $u\in H$ be fixed and observe that, as \eqref{forma biline cont} holds, the map
$$H\ni v\rightarrow a(u,v)\in \mathbb{R},$$
is linear and bounded.
The Riesz Representation Theorem implies that there exists a unique
$Au\in H$ such that $$a(u,v)=(Au,v),\quad\forall v\in H.$$
Hence, we have defined the map
$$A:H\rightarrow H,$$
such that
\begin{equation}\label{ELLlaxmilgram.1}
	a(u,v)=(Au,v), \quad \forall u,v\in H.
\end{equation}
By \eqref{ELLlaxmilgram.2} and \eqref{ELLlaxmilgram.1} we have that
\eqref{lm} is equivalent to
\begin{equation}\label{ELLlaxmilgram.3}
	Au=f,\quad u\in H.
\end{equation}

Now, we prove that the map $A$ is 
\textbf{(i) linear, (ii) bounded} and \\ \textbf{(iii) bijective}.

\medskip

\textbf{(i)} Let $u_1,u_2\in H$, $\lambda_1, \lambda_2\in
\mathbb{R}$. We have
\begin{equation*}
	\begin{aligned}
		\left(A\left(\lambda_1u_1+\lambda_2u_2\right),v\right)&=a\left(\lambda_1u_1+\lambda_2u_2,v\right)\\&=\lambda_1a\left(u_1,v\right)+\lambda_2a\left(u_2,v\right)=\\&=
		\lambda_1\left(Au_1,v\right)+\lambda_2\left(Au_2,v\right)=\\&=\left(\lambda_1Au_1+\lambda_2
		Au_2,v\right),\quad\forall v\in H.
	\end{aligned}
\end{equation*}
Hence
\begin{equation*}
	A\left(\lambda_1u_1+\lambda_2u_2\right)=\lambda_1Au_1+\lambda_2
	Au_2,\quad \forall u_1,u_2\in H.
\end{equation*}

\medskip

\textbf{(ii)} By \eqref{forma biline cont} we get

\begin{equation*}
	\left|\left(Au,v\right)\right|=|a(u,v)|\leq C \left\Vert
	u\right\Vert \left\Vert v\right\Vert,\quad\forall u,v\in
	H.
\end{equation*}
By the Cauchy--Schwarz inequality we get

\begin{equation*}
	\left\Vert Au\right\Vert\leq C \left\Vert u\right\Vert
	,\quad\forall u\in H.
\end{equation*}
Therefore $A$ is bounded and

\begin{equation*}
	\left\Vert A\right\Vert_{\mathcal{L}(H)}\leq C,
\end{equation*}
where $\mathcal{L}(H)$ \index{$\mathcal{L}(H)$}is the space of bounded linear map  from $H$ in itself.
\medskip

\textbf{(iii)} Condition \eqref{forma biline cont-1} implies

\begin{equation*}
	\alpha\left\Vert u\right\Vert^2\leq a(u,u)=(Au,u)\leq \left\Vert A
	u\right\Vert\left\Vert u\right\Vert ,\quad\forall u\in H.
\end{equation*}
Hence

\begin{equation}\label{ELLlaxmilgram.4}
	\alpha\left\Vert u\right\Vert\leq \left\Vert A
	u\right\Vert,\quad\forall u\in H.
\end{equation}
Since $\alpha>0$, $A$ is inijective.

In order to prove that $A$ is onto, we first prove
that $\mathcal{R}(A)$ is closed. Let  $\left\{w_k\right\}$ be a sequence in $\mathcal{R}(A)$ such that
\begin{equation}\label{ELLlaxmilgram.5}
	\left\{w_k\right\}\rightarrow w, \end{equation} and let us check that
$w\in \mathcal{R}(A)$. Let $u_k\in H$, $k\in \mathbb{N}$, satisfy
$$Au_k=w_k,\quad\forall k\in \mathbb{N}.$$
By \eqref{ELLlaxmilgram.4} we have
\begin{equation*}
	\left\Vert u_k-u_j\right\Vert\leq \frac{1}{\alpha}\left\Vert A
	u_k-Au_j\right\Vert=\frac{1}{\alpha}\left\Vert
	w_k-w_j\right\Vert,\quad\forall k,j\in \mathbb{N}.
\end{equation*}
Since $\left\{w_k\right\}$ converges, it is a Cauchy sequence. Consequently, $\left\{u_k\right\}$ is  a Cauchy sequence too. Therefore there exists  $u\in H$ such that
$$\left\{u_k\right\}\rightarrow u.$$ Now, by \eqref{ELLlaxmilgram.5} and, as  $A$ is continuous, we obtain

$$w=\lim_{k\rightarrow\infty}w_k=\lim_{k\rightarrow\infty}Au_k=Au.$$
Therefore $w\in \mathcal{R}(A)$.

\smallskip

Now we prove 

\begin{equation}\label{ELLlaxmilgram.6}
	\mathcal{R}(A)=H. \end{equation} We argue by contradiction.  Let us assume that $\mathcal{R}(A)\subsetneqq H$. Since $\mathcal{R}(A)$ is closed, there exists $w\in H\setminus\{0\}$  such that $w\perp
\mathcal{R}(A)$, (by this we mean $(w,h)=0$ for every $h\in \mathcal{R}(A)$). Hence 

\begin{equation*}
	\alpha\left\Vert w\right\Vert^2\leq a(w,w)=(Aw,w)=0.
\end{equation*}
Consequently, we should have  $w=0$ that contradicts $w\neq 0$. Thus
\eqref{ELLlaxmilgram.6} is proved.

\medskip

Now, since \eqref{lm} and
\eqref{ELLlaxmilgram.3} are equivalent, there exists one and only one solution $u\in
H$ of the problem 

\begin{equation}\label{lm:1}
	a(u,v)=F(v), \quad\forall v\in H.
\end{equation}

Concerning estimate \eqref{lm-cont}, it follows immediately by
\begin{equation}\label{lm-cont:1}
	\alpha\left\Vert u\right\Vert^2\leq a(u,u)=F(u)\leq \left\Vert
	F\right\Vert_{H'}\left\Vert u\right\Vert.
\end{equation}
$\blacksquare$

\bigskip

\textbf{Remark.} We observe that, as $a$ is a
bilinear form, by \eqref{lm-cont} we obtain that if $F_j\in
H'$, $j=1,2$ and $u_j\in H$ are solutions to

\begin{equation*}
	a(u_j,v)=F_j(v) \quad\forall v\in H,
\end{equation*}
then

\begin{equation}\label{lm-cont:2}
	\left\Vert u_1-u_2\right\Vert\leq \frac{1}{\alpha}\left\Vert
	F_1-F_2\right\Vert_{H'}
\end{equation}
from which, in particular, we get again the uniqueness.
$\blacklozenge$

\bigskip

We now recall the Fredholm Alternative Theorem \cite{Magn2}. Meanwhile, we recall that a linear operator  $\mathcal{K}$ from $H$ in itself is said compact \index{compact linear operator}, provided that for every bounded set $M\subset H$, $\mathcal{K}(M)$
is relatively compact in $H$. We recall that a compact operator is necessarily bounded, neverthless if $H$ does not have finite dimension, the identity on $H$ is a bounded operator, but it is not compact.

\begin{theo}[\textbf{Fredholm Alternative}]\label{Fredholm}
	\index{Theorem:@{Theorem:}!- Fredholm Alternative@{- Fredholm Alternative}}
	Let
	$$\mathcal{K}:H\rightarrow H,$$
	a linear compact operator. Then we have:
	
	\smallskip
	
	(i) the dimension of $\mathcal{N}(I-\mathcal{K})$ is finite;
	
	\smallskip
	
	(ii) $\mathcal{R}(I-\mathcal{K})$ is a closed subspace;
	
	\smallskip
	
	(iii)
	$\mathcal{R}(I-\mathcal{K})=\mathcal{N}(I-\mathcal{K}^{\star})^{\perp},$
	
	\smallskip
	
	(iv) $\mathcal{N}(I-\mathcal{K})=\{0\}$ if and only if
	$\mathcal{R}(I-\mathcal{K})=H$;
	
	\smallskip
	
	(v) dim $\mathcal{N}(I-\mathcal{K})$ = dim
	$\mathcal{N}(I-\mathcal{K}^{\star})$.
\end{theo}
Recall that $\mathcal{K}^{\star}$ is the adjoint operator
$\mathcal{K}$ defined by

\begin{equation*}
	\left(\mathcal{K}u,v\right)=\left(u,\mathcal{K}^{\star}v\right),\quad\forall
	u,v\in H.
\end{equation*}

\section[The variational formulation of the Dirichlet problem]{The variational formulation of the Dirichlet problem. Existence theorems}\label{Ell:for-variazionale}
Let us begin by clarifying what we mean by \textbf{the
	variational formulation} \index{variational formulation of Dirichlet problem}of Dirichlet problem \eqref{Ell:Dirichlet}.
We begin by the case in which the \textbf{condition at the boundary is
	homogeneous}. In this case \eqref{Ell:Dirichlet} can be written
(formally):
\begin{equation}\label{Dirichlet}
	\begin{cases}
		-\mbox{div}(A\nabla u)=F, \quad\mbox{ in } \Omega, \\
		\\
		u=0, \quad\mbox{ on } \partial\Omega
	\end{cases}
\end{equation}
and the variational formulation of problem above is the following:

\medskip

\noindent Determine $u$ such that

\begin{equation}\label{variazDirichlet}
	\begin{cases}
		\int_{\Omega} A\nabla u\cdot \nabla v  dx=F(v), \quad\forall v \in H_0^{1}(\Omega), \\
		\\
		u\in H_0^{1}(\Omega).
	\end{cases}
\end{equation}

\medskip 

Let us observe that if $\Omega$ is an open set of class $C^{0,1}$,
$a^{jk}\in C^{1}(\overline{\Omega})$, $j,k=1,\cdots,n$ and $u\in
C^{2}\left(\overline{\Omega}\right)$ with $u_{|\partial\Omega}=0$ and
$F\in C^{0}(\overline{\Omega})$, then \eqref{Dirichlet} is equivalent to \eqref{variazDirichlet}. Indeed, under such assumptions, the divergence Theorem implies

\begin{equation}\label{div}
	\int_{\Omega}\mbox{div}(A\nabla u) v dx=-\int_{\Omega}A\nabla
	u\cdot\nabla v dx,\quad\quad \forall v \in C_0^{\infty}(\Omega).
\end{equation}
Therefore, \eqref{Dirichlet} implies \eqref{variazDirichlet}
(taking into account Remark \ref{Sob:oss-traccia} and that $C_0^{\infty}(\Omega)$ is dense in
$H_0^{1}(\Omega)$). Conversely, if \eqref{variazDirichlet} holds true, then
Theorem \ref{traccia} implies $u_{|\partial\Omega}=0$ and
\eqref{variazDirichlet} gives

\begin{equation}\label{div-1}
	\int_{\Omega}\left(\mbox{div}(A\nabla u)+F\right) v dx=0,
	\quad\forall v \in C_0^{\infty}(\Omega),
\end{equation}
that gives  \eqref{Dirichlet}. Of course, under the general assumptions on $\Omega$, $A$ and $F$, the formulation \eqref{Dirichlet} makes no sense, while the formulation \eqref{variazDirichlet} makes
perfectly sense, and the Lax--Milgram Theorem will tell us
easily that it is a well--posed problem. Actually, let $H=H_0^{1}(\Omega)$ and

\begin{equation}\label{a-dirich}
	a(u,v)=\int_{\Omega}A\nabla u\cdot\nabla v dx, \quad\forall u,v \in
	H_0^{1}(\Omega).
\end{equation}
The form \eqref{a-dirich} is bilinear. Moreover by the Cauchy--Schwarz inequality we have, for any $u, v \in
H_0^{1}(\Omega)$,
\begin{equation*}\label{a-dirich-cont}
	\begin{aligned}
		|a(u,v)|&\leq \left\Vert
	A\right\Vert_{L^{\infty}(\Omega;\mathbb{M}(n))}\int_{\Omega}|\nabla
	u||\nabla v| dx\leq \\&\leq \left\Vert
	A\right\Vert_{L^{\infty}(\Omega;\mathbb{M}(n)}\left\Vert \nabla
	u\right\Vert_{L^{2}(\Omega)}\left\Vert \nabla
	v\right\Vert_{L^{2}(\Omega)}
	\end{aligned}
\end{equation*}
and by \eqref{gamma-eq} we have
\begin{equation*}\label{a-dirich-coerc}
	\lambda^{-1}\left\Vert \nabla u\right\Vert^2_{L^{2}(\Omega)}\leq
	\int_{\Omega}A\nabla u\cdot\nabla u dx=a(u,u), \quad \forall u \in
	H_0^{1}(\Omega).
\end{equation*}
Now, recalling that $\left\Vert \nabla u\right\Vert_{L^{2}(\Omega)}$
and $\left\Vert u\right\Vert_{H^{1}(\Omega)}$ are two equivalent norms of $H_0^{1}(\Omega)$, we have only to apply the Lax--Milgram Theorem to conclude that problem \eqref{variazDirichlet} has an unique solution
in  $H_0^{1}(\Omega)$. Moreover the following inequality holds true 
\begin{equation}\label{stab}
	\left\Vert \nabla u\right\Vert_{L^{2}(\Omega)}\leq \lambda\left\Vert
	F\right\Vert_{H^{-1}(\Omega)}.
\end{equation}

Inequality \eqref{stab}, together with the already existence and uniqueness results, implies that problem \eqref{variazDirichlet} is well--posed in $H_0^{1}(\Omega)$.

\bigskip

Now we consider the case where the boundary condition is \textbf{not
	homogeneous}, but the equation is still homogeneous.

Let $\Omega$ be an open set of $\mathbb{R}^n$ of class $C^{0,1}$, let $\varphi\in
H^{1/2}(\partial\Omega)$ and let us assume 
that \eqref{gamma-eq} is satisfied. Formally, the Dirichlet problem can be written as

\begin{equation}\label{Dirichlet-non om}
	\begin{cases}
		-\mbox{div}(A\nabla u)=0, \quad\mbox{ in } \Omega, \\
		\\
		u=\varphi, \quad\mbox{ on } \partial\Omega.
	\end{cases}
\end{equation}
We wish to give the variational formulation of problem \eqref{Dirichlet-non om} and to prove
the existence of the solutions in $H^1(\Omega)$ to this problem.

The variational formulation of \eqref{Dirichlet-non om} is

\begin{equation}\label{variazDirichlet-100}
	\begin{cases}
		\int_{\Omega} A\nabla u\cdot \nabla v  dx=0, \quad\forall v \in H_0^{1}(\Omega), \\
		\\
		u=\varphi, \quad\mbox{ on } \partial\Omega\quad\mbox{ (in the sense of the traces)}.
	\end{cases}
\end{equation}
Notice that, in the case where $A$, $u$ and $\varphi$ are
sufficiently regular, the first equation in
\eqref{variazDirichlet-100} is equivalent to the first equation
of \eqref{Dirichlet-non om}.

\medskip

In order to solve \eqref{Dirichlet-non om} (\eqref{variazDirichlet-100}),
we proceed in the following way. Recalling Theorem \ref{Sob:teo4.11}, there exists $\Phi\in H^{1}(\Omega)$ such that
$$\Phi_{|\partial\Omega}=\varphi,\quad\mbox{(in the sense of the traces)}$$ 
which in turn (by (ii) of Theorem \ref{Sob:teo4.11}) implies

\begin{equation}\label{Dirichlet-non om-1}
	\left\Vert \Phi\right\Vert_{H^{1}(\Omega)}\leq C\left\Vert
	\varphi\right\Vert_{H^{1/2}(\partial \Omega)},
\end{equation}
where $C$ is a constant depending on $\Omega$ only. Set
$w=u-\Phi$. Since $u$ and $\Phi$ have the same trace on
$\partial \Omega$ we have $w_{|\partial \Omega}=0$, so that problem
\eqref{Dirichlet-non om} can be written (formally),
\begin{equation}\label{Dirichlet-non om-10}
	\begin{cases}
		-\mbox{div}(A\nabla w)=\mbox{div}(A\nabla \Phi), \quad\mbox{ in } \Omega, \\
		\\
		w=0, \quad\mbox{ on } \partial\Omega,
	\end{cases}
\end{equation}
whose variational formulation is

\begin{equation}\label{variazDirichlet-1}
	\begin{cases}
		\int_{\Omega} A\nabla w\cdot \nabla v  dx=\int_{\Omega} A\nabla \Phi\cdot \nabla v  dx, \quad\forall v\in H_0^{1}(\Omega),   \\
		\\
		w\in H_0^{1}(\Omega).
	\end{cases}
\end{equation}

Let us note that the bilinear form is still given by
\eqref{a-dirich}. The solution of problem \eqref{Dirichlet-non om}
(\eqref{variazDirichlet-100}) is given by

\begin{equation}\label{variazDirichlet-25}
	u=w+\Phi\in H^1(\Omega).
\end{equation}
Moreover, denoting
\begin{equation*}\label{variazDirichlet-10}
	F(v)=\int_{\Omega} A\nabla \Phi\cdot \nabla v  dx, \quad\forall
	v\in H_0^{1}(\Omega),
\end{equation*}
it turns out that $F\in H^{-1}(\Omega)$. As a matter of fact by the Cauchy--Schwarz inequality we have

\begin{equation*}
	|F(v)|\leq \lambda \left\Vert
	\nabla\Phi\right\Vert_{L^{2}(\Omega)}\left\Vert \nabla
	v\right\Vert_{L^{2}(\Omega)} , \quad\forall v\in H_0^{1}(\Omega).
\end{equation*}
By \eqref{Dirichlet-non om-1} and \eqref{Dirichlet-non om-10} we have
\begin{equation*}\label{variazDirichlet-20}
	|F(v)|\leq  C\lambda\left\Vert
	\varphi\right\Vert_{H^{1/2}(\partial\Omega)}\left\Vert
	v\right\Vert_{H^{1}(\Omega)}, \quad\forall v\in H_0^{1}(\Omega).
\end{equation*}
Therefore $F\in H^{-1}(\Omega)$ and the following inequality holds
\begin{equation}\label{variazDirichlet-21}
	\left\Vert F\right\Vert_{H^{-1}(\Omega)}\leq C'\left\Vert
	\varphi\right\Vert_{H^{1/2}(\partial\Omega)},
\end{equation}
where $C'$ depends on $\Omega$ and $\lambda$ only .

\medskip

The Lax--Milgram Theorem implies that  problem
\eqref{variazDirichlet-1} has a unique solution $w\in
H_0^{1}(\Omega)$, moreover by \eqref{lm-cont} we have
\begin{equation}\label{variazDirichlet-29}
	\left\Vert \nabla w\right\Vert_{L^{2}(\Omega)}\leq \lambda\left\Vert
	F\right\Vert_{H^{-1}(\Omega)}\leq \lambda C'\left\Vert
	\varphi\right\Vert_{H^{1/2}(\partial\Omega)}.
\end{equation}
Now, by using \eqref{Dirichlet-non om-1},
\eqref{variazDirichlet-25} and \eqref{variazDirichlet-29} we get
\begin{equation}\label{variazDirichlet-30}
	\left\Vert  u \right\Vert_{H^{1}(\Omega)}\leq \left\Vert w
	\right\Vert_{H^{1}(\Omega)}+\left\Vert \Phi
	\right\Vert_{H^{1}(\Omega)}\leq C''\left\Vert
	\varphi\right\Vert_{H^{1/2}(\partial\Omega)},
\end{equation}
where $C''$ depends by $\lambda$ and $\Omega$ only. Notice that
\eqref{variazDirichlet-30} implies, in particular, the uniqueness of solution of problem \eqref{Dirichlet-non om}.

From what we have proved so far, we have the following
\begin{theo}\label{teorema-es-unic}
	Let $\Omega$ be an open set of class $C^{0,1}$. Let us assume that $A\in
	L^{\infty}(\Omega;\mathbb{M}(n))$ and $A$ satisfies \eqref{gamma-eq}. Let
	$F\in H^{-1}(\Omega)$ and $\varphi\in H^{1/2}(\partial\Omega)$.
	
	Then the following problem
	\begin{equation}\label{Dirichlet-completo}
		\begin{cases}
			-\mbox{div}(A\nabla u)=F, \quad\mbox{ in } \Omega, \\
			\\
			u=\varphi, \quad\mbox{ on } \partial\Omega,
		\end{cases}
	\end{equation}
	whose variational formulation is
	
	\begin{equation}\label{variazDirichlet-100F}
		\begin{cases}
			\int_{\Omega} A\nabla u\cdot \nabla v  dx=F(v), \quad\forall v \in H_0^{1}(\Omega), \\
			\\
			u=\varphi, \quad\mbox{ su } \partial\Omega\quad\mbox{ (in the sense of
				traces)},
		\end{cases}
	\end{equation}
	has a unique solution $u\in H^{1}(\Omega)$ and we have
	
	\begin{equation}\label{stima-Dirichlet-completo}
		\left\Vert u\right\Vert_{H^{1}(\Omega)}\leq C\left(\left\Vert
		F\right\Vert_{H^{-1}(\Omega)}+\left\Vert
		\varphi\right\Vert_{H^{1/2}(\partial\Omega)}\right),
	\end{equation}
	where $C$ depends on $\lambda$ and $\Omega$ only.
\end{theo}

\bigskip

Now let $L$ be the following operator
\begin{equation}\label{Ell:for-variazionale.1}
	Lu=-\sum_{j,k=1}^n\partial_j\left( a^{jk}\partial_k u
	+d^ju\right)+\sum_{j=1}^nb^{j}\partial_j u+cu,\end{equation} where  $A\in
L^{\infty}(\Omega;\mathbb{M}(n))$,
$A=\left\{a^{jk}\right\}_{j,k=1}^n$, satisfies \eqref{gamma-eq} and
$$b^j, \ d^j, \ c\in L^{\infty}(\Omega),$$ for $j=1,\cdots, n$. 

\smallskip

Let us consider the Dirichlet problem 

\begin{equation}\label{Ell:for-variazionale.2}
	\begin{cases}
		Lu=f, \quad\mbox{ in } \Omega, \\
		\\
		u=0, \quad\mbox{ su } \partial\Omega,
	\end{cases}
\end{equation}
where $f\in L^{2}(\Omega)$. The variational formulation of above problem is:

\medskip

\noindent Determine $u$ such that

\begin{equation}\label{Ell:for-variazionale.3}
	\begin{cases}
		a(u,v)=(f,v), \quad\forall v \in H_0^{1}(\Omega), \\
		\\
		u\in H_0^{1}(\Omega),
	\end{cases}
\end{equation}
where

\begin{equation}\label{Ell:for-variazionale.4}
	a(u,v)=\int_{\Omega}\left(A\nabla u\cdot \nabla v+ud\cdot v-b\cdot
	\nabla u v-cuv\right)dx,
\end{equation}
 $d=\left(d^1,\cdots,d^n\right)$ and
$b=\left(b^1,\cdots,b^n\right)$. 

\medskip

Problem \eqref{Ell:for-variazionale.2} does not always have existence and uniqueness. To show this fact, let us consider the following simple example 

\begin{equation}\label{Ell:for-variazionale.5}
	\begin{cases}
		-u''-u=f, \quad\mbox{ in } (0,\pi), \\
		\\
		u(0)=u(\pi)=0,
	\end{cases}
\end{equation}
where $f\in L^2(0,\pi)$. The solutions to
\eqref{Ell:for-variazionale.5} have to be found among the functions of the
type

$$C_1\sin x +C_2 \cos x-\int^x_0\sin(x-t)f(t)dt.$$
By the boundary conditions we have $C_2=0$
and
\begin{equation}\label{Ell:for-variazionale.6}\int^{\pi}_0f(t)\sin
	tdt=0.\end{equation} Therefore, if \eqref{Ell:for-variazionale.6}
is satisfied, then  \eqref{Ell:for-variazionale.5} has infinite solutions, given by

$$C\sin x -\int^x_0\sin(x-t)f(t)dt, \ \ C\in\mathbb{R}$$
Whereas if \eqref{Ell:for-variazionale.6} is not satisfied, then  \eqref{Ell:for-variazionale.5} has no  solutions. So we cannot  expect that the bilinear form \eqref{Ell:for-variazionale.4} is always coercive.
We can, however, prove the following Theorem that will be useful for establish some conditions of existence and uniqueness to problem \eqref{Ell:for-variazionale.3}

\begin{theo}\label{Ell:stimaenergia}
	Let $\Omega$ be a bounded open set of $\mathbb{R}^n$. Let us assume
	that $A\in L^{\infty}(\Omega;\mathbb{M}(n))$ and that $A$ satisfies condition
 \eqref{gamma-eq} and $b,d \in
	L^{\infty}\left(\Omega;\mathbb{R}^{n}\right)$, $c\in L^{\infty}(\Omega)$, for
	$j=1,\cdots, n$. Moreover, let
	\begin{equation*}
		a(u,v)=\int_{\Omega}\left(A\nabla u\cdot \nabla v+ud\cdot \nabla
		v-b\cdot \nabla u v-cuv\right)dx.
	\end{equation*}
	Then $a$ is a continuous bilinear form and there exist
	$\gamma_0\geq 0$ and $\alpha_0>0$ such that for any $\gamma\geq
	\gamma_0$ we have
	
	\begin{equation}\label{Ell:stimaenergia.1}
		\alpha_0\left\Vert u\right\Vert^2_{H^{1}(\Omega)}\leq
		a(u,u)+\gamma\left\Vert u\right\Vert^2_{L^{2}(\Omega)},\quad \forall
		u\in H_0^{1}(\Omega).
	\end{equation}
	In particular,
	\begin{equation}\label{Ell:stimaenergia.1bis}
		a_{\gamma}(u,v):=a(u,u)+\gamma\left(u,v\right)_{L^{2}(\Omega)},
	\end{equation}
	is a continuous and coercive bilinear form for every $\gamma\geq \gamma_0$.
\end{theo}
\textbf{Proof.} By the Cauchy--Schwarz inequality we have, for any $u,v\in H^1_0(\Omega)$

\begin{equation*}
	\begin{aligned}
		&|a(u,v)|\leq \int_{\Omega}\left\Vert
		A\right\Vert_{L^{\infty}(\Omega;\mathbb{M}(n))}\left\vert \nabla
		u\right\vert\left\vert \nabla v\right\vert
		dx+\\&+\int_{\Omega}\left(\left\Vert
		d\right\Vert_{L^{\infty}(\Omega;\mathbb{R}^n)}\left\vert u\right\vert\left\vert
		\nabla v\right\vert +\left\Vert
		b\right\Vert_{L^{\infty}(\Omega;\mathbb{R}^n)}\left\vert \nabla
		u\right\vert\left\vert v\right\vert+\left\Vert
		c\right\Vert_{L^{\infty}(\Omega)}\left\vert u\right\vert\left\vert
		v\right\vert\right)dx\leq\\&\leq \left\Vert
		A\right\Vert_{L^{\infty}(\Omega;\mathbb{M}(n))}\left\Vert \nabla
		u\right\Vert_{L^{2}(\Omega)}\left\Vert \nabla
		v\right\Vert_{L^{2}(\Omega)}+ \left\Vert
		d\right\Vert_{L^{\infty}(\Omega;\mathbb{R}^n)}\left\Vert
		u\right\Vert_{L^{2}(\Omega)}\left\Vert \nabla
		v\right\Vert_{L^{2}(\Omega)}+\\&+\left\Vert
		b\right\Vert_{L^{\infty}(\Omega;\mathbb{R}^n)}\left\Vert \nabla
		u\right\Vert_{L^{2}(\Omega)}\left\Vert
		v\right\Vert_{L^{2}(\Omega)}+\left\Vert
		c\right\Vert_{L^{\infty}(\Omega)}\left\Vert
		u\right\Vert_{L^{2}(\Omega)}\left\Vert
		v\right\Vert_{L^{2}(\Omega)}\leq\\&\leq C\left\Vert
		u\right\Vert_{H^{1}(\Omega)}\left\Vert v\right\Vert_{H^{1}(\Omega)},
	\end{aligned}
\end{equation*}
where
$$C=\left(\left\Vert
A\right\Vert_{L^{\infty}(\Omega;\mathbb{M}(n))}+\left\Vert
d\right\Vert_{L^{\infty}(\Omega;\mathbb{R}^n)}+\left\Vert
b\right\Vert_{L^{\infty}(\Omega;\mathbb{R}^n)}+\left\Vert
c\right\Vert_{L^{\infty}(\Omega)}\right).$$ Therefore

\begin{equation*}
	|a_{\gamma}(u,v)|\leq (C+|\gamma|)\left\Vert
	v\right\Vert_{H^{1}(\Omega)}\left\Vert u\right\Vert_{H^{1}(\Omega)},
	\quad\forall u,v\in H_0^{1}(\Omega)
\end{equation*}
which implies the continuity of $a_{\gamma}$ for any $\gamma\in \mathbb{R}$.

Concerning coercivity of $a_{\gamma}$, we notice
firstly that \eqref{gamma-eq} gives

\begin{equation}\label{Ell:stimaenergia.2}
	\int_{\Omega}A\nabla u\cdot \nabla u dx\geq \lambda^{-1}\left\Vert
	\nabla u\right\Vert^2_{L^{2}(\Omega)}, \quad \forall u\in
	H_0^{1}(\Omega).
\end{equation}
Moreover, let $\varepsilon$ be a positive number which we will choose
later on. We get

\begin{equation}\label{Ell:stimaenergia.3}
	\begin{aligned}
		\left|\int_{\Omega}ud\cdot \nabla u dx\right|&\leq \left\Vert
		d\right\Vert_{L^{\infty}(\Omega;\mathbb{R}^n)}\left\Vert
		u\right\Vert_{L^{2}(\Omega)}\left\Vert \nabla
		u\right\Vert_{L^{2}(\Omega)}\leq\\&\leq
		\frac{\varepsilon}{2}\left\Vert \nabla
		u\right\Vert^2_{L^{2}(\Omega)}+\frac{1}{2\varepsilon}\left\Vert
		d\right\Vert^2_{L^{\infty}(\Omega;\mathbb{R}^n)}\left\Vert
		u\right\Vert^2_{L^{2}(\Omega)} \end{aligned}
\end{equation} 
similarly,
\begin{equation}\label{Ell:stimaenergia.4}
	\left|\int_{\Omega}b \cdot\nabla u u dx\right|\leq
	\frac{\varepsilon}{2}\left\Vert \nabla
	u\right\Vert^2_{L^{2}(\Omega)}+\frac{1}{2\varepsilon}\left\Vert
	b\right\Vert^2_{L^{\infty}(\Omega;\mathbb{R}^n)}\left\Vert
	u\right\Vert^2_{L^{2}(\Omega)}.
\end{equation}
Furthermore
\begin{equation}\label{Ell:stimaenergia.5}
	\left|\int_{\Omega}c  u^2 dx\right|\leq \left\Vert
	c\right\Vert_{L^{\infty}(\Omega)}\left\Vert
	u\right\Vert^2_{L^{2}(\Omega)}.
\end{equation}
Hence, by \eqref{Ell:stimaenergia.2}--\eqref{Ell:stimaenergia.5}
we obtain

\begin{equation*}
	\begin{aligned}
		a_{\gamma}(u,u)&\geq \left(\lambda^{-1}-\varepsilon\right)\left\Vert
		\nabla
		u\right\Vert^2_{L^{2}(\Omega)}+\\&+\left[\gamma-\left(\frac{1}{2\varepsilon}\left\Vert
		d\right\Vert^2_{L^{\infty}(\Omega;\mathbb{R}^n)}+\frac{1}{2\varepsilon}\left\Vert
		b\right\Vert^2_{L^{\infty}(\Omega;\mathbb{R}^n)}+\left\Vert
		c\right\Vert_{L^{\infty}(\Omega)}\right) \right]\left\Vert
		u\right\Vert^2_{L^{2}(\Omega)}.
	\end{aligned}
\end{equation*}
Now, by choosing
$$\varepsilon=\frac{\lambda^{-1}}{2}$$ and denoting

$$\gamma_0=\lambda\left\Vert
d\right\Vert^2_{L^{\infty}(\Omega;\mathbb{R}^n)}+\lambda\left\Vert
b\right\Vert^2_{L^{\infty}(\Omega;\mathbb{R}^n)}+\left\Vert
c\right\Vert_{L^{\infty}(\Omega)},$$ we have, for any $\gamma\geq
\gamma_0$
\begin{equation*}
	\begin{aligned}
		a_{\gamma}(u,u)\geq \frac{\lambda^{-1}}{2}\left\Vert \nabla
		u\right\Vert^2_{L^{2}(\Omega)},\quad \forall u\in H_0^{1}(\Omega).
	\end{aligned}
\end{equation*}
Finally, the first Poincar\'{e} inequality (Theorem
\ref{Poincar})  implies that there exists $\alpha_0$, depending on the
diameter of $\Omega$, such that \eqref{Ell:stimaenergia.1} is satisfied.
$\blacksquare$

\bigskip

\underline{\textbf{Exercise.}} Prove there that there exists $\delta>0$ such that if
diam($\Omega)\leq \delta$ (diam($\Omega)$ is the diameter of) then problem
\eqref{Ell:for-variazionale.2}, where $f\in H^{-1}(\Omega)$, has a unique solution in $H_0^{1}(\Omega)$. [Hint: use the first Poincar\'{e} inequality].

\bigskip

The following operator is called \textbf{the (formal) adjoint} \index{adjoint (formal)}of the operator
$L$ 

\begin{equation}\label{Ell:for-variazionale.7}
	L^{\star}v=-\sum_{j,k=1}^n\partial_k\left( a^{jk}\partial_j v
	+b^kv\right)+\sum_{j=1}^nd^{j}\partial_j v+cv.\end{equation} To the operator
$L^{\star}$ corresponds the bilinear form

\begin{equation*}
	a^{\star}(v,u)=\int_{\Omega}\left(A^T\nabla v\cdot \nabla u-vb\cdot
	\nabla u+d\cdot \nabla v u-cvu\right)dx, \quad \forall u,v\in
	H^1_0(\Omega),\end{equation*} where $A^T$ is the transposed of
	the matrix  $A$. Let us notice that

\begin{equation*}
	a^{\star}(v,u)=a(u,v), \quad \forall u,v\in
	H^1_0(\Omega).\end{equation*} Finally, let $f\in L^{2}(\Omega)$,
we say that $v\in H^1_0(\Omega)$ is a weak solution of the adjoint problem

\begin{equation}\label{Ell:for-variazionale.8}
	\begin{cases}
		L^{\star} v=f, \quad\mbox{ in } \Omega, \\
		\\
		u=0, \quad\mbox{ on } \partial\Omega,
	\end{cases}
\end{equation}
provided that we have

\begin{equation}\label{Ell:for-variazionale.3-0}
	\begin{cases}
		a^{\star}(v,u)=(f,u)_{L^{2}(\Omega)}, \quad\forall u \in H_0^{1}(\Omega), \\
		\\
		v\in H_0^{1}(\Omega).
	\end{cases}
\end{equation}

\medskip 

The following Theorem holds true

\begin{theo}\label{Ell:esistenza}
	Let $L$ be operator  \eqref{Ell:for-variazionale.1} and $L^{\star}$
	its formal adjoint.
	
	\noindent(i) The following alternative holds true.
	
	\smallskip
	
	either
	
	(a) for any $f\in L^{2}(\Omega)$ there exists a unique $u\in
	H_0^{1}(\Omega)$ such that
	
	\begin{equation}\label{Ell:esistenza1}
		\begin{cases}
			Lu=f, \quad\mbox{ in } \Omega, \\
			\\
			u=0, \quad\mbox{ on } \partial\Omega
		\end{cases}
	\end{equation}
	
	or
	
	(b) there exists at least one \textbf{not identically vanishing solution } $u\in
	H_0^{1}(\Omega)$ to the homogeneous problem
	\begin{equation}\label{Ell:esistenza2}
		\begin{cases}
			Lu=0, \quad\mbox{ in } \Omega, \\
			\\
			u=0, \quad\mbox{ on } \partial\Omega.
		\end{cases}
	\end{equation}
	
	\noindent(ii) If (b) holds true, then, denoting by $N$
	the subspace $H_0^{1}(\Omega)$ of the solutions to
	\eqref{Ell:esistenza2} and by $N^{\star}$ the subspace of
	$H_0^{1}(\Omega)$ of the solutions to
	\begin{equation}\label{Ell:esistenza3}
		\begin{cases}
			L^{\star}v=0, \quad\mbox{ in } \Omega, \\
			\\
			v=0, \quad\mbox{ on } \partial\Omega,
		\end{cases}
	\end{equation}
we have that $N$ and $N^{\star}$ have finite \textbf{dimension}, moreover 
	
	$$\mbox{dimension of }  N=  \mbox{dimension of }N^{\star}.$$
	
	\noindent (iii) Finally, problem \eqref{Ell:esistenza1} admits
	a solution in $H_0^{1}(\Omega)$ if and only if
	
	\begin{equation}\label{Ell:esistenza4}
		(f,v)_{L^{2}(\Omega)}=0,\quad\forall v\in N^{\star}.
	\end{equation}
\end{theo}
\textbf{Proof.} Let $\gamma_0$ be the same of  Theorem
\ref{Ell:stimaenergia}. Let us fix $\gamma\geq \gamma_0$. Since
 $a_{\gamma}$, defined by
\eqref{Ell:stimaenergia.1bis}, is a continuous and coercive bilinear form, we have that for any
 $g\in L^{2}(\Omega)$ there exists a unique $u\in
H_0^{1}(\Omega)$ such that 

$$L_{\gamma}u=g.$$
Set

$$L^{-1}_{\gamma}g=u.$$

Now, we notice that $u\in H_0^{1}(\Omega)$ solves the
boundary value problem
\begin{equation}\label{Ell:esistenza-p5}
	\begin{cases}
		Lu=f, \quad\mbox{ in } \Omega, \\
		\\
		u=0, \quad\mbox{ su } \partial\Omega
	\end{cases}
\end{equation} if and only if

\begin{equation*}
	a_{\gamma}(u,v)=(\gamma u+f,v)_{L^{2}(\Omega)}, \quad \forall v\in
	H_0^{1}(\Omega)
\end{equation*}
which, in turn, is equivalent to
\begin{equation}\label{Ell:esistenza5}
	u=L^{-1}_{\gamma}(\gamma u+f).
\end{equation}
Let us denote
\begin{equation}\label{Ell:esistenza6}
	\mathcal{K}u=\gamma L^{-1}_{\gamma} u
\end{equation}
and
\begin{equation}\label{Ell:esistenza6bis}
	h= L^{-1}_{\gamma} f.
\end{equation}
Let us notice that $\mathcal{K}$ is linear and it satisfies

\begin{equation*}
	a_{\gamma}(\mathcal{K}g,v)=\gamma(g,v)_{L^{2}(\Omega)}, \quad
	\forall g\in L^{2}(\Omega), \mbox{ } \forall v\in H_0^{1}(\Omega).
\end{equation*}
Moreover, as \eqref{Ell:esistenza5} and \eqref{Ell:esistenza6} hold,
\eqref{Ell:esistenza-p5} is equivalent to
\begin{equation}\label{Ell:esistenza7}
	u-\mathcal{K}u=h.
\end{equation}
Also, we notice that we have $h\in
H_0^{1}(\Omega)$ and by the definition of $\mathcal{K}$, we have $\mathcal{K}g \in
H_0^{1}(\Omega)$, for any $g\in L^{2}(\Omega)$, . Therefore, every function of $L^{2}(\Omega)$  which is a solution to
\eqref{Ell:esistenza7} belongs to $H_0^{1}(\Omega)$. 

\smallskip

Now we examine the solvability of \eqref{Ell:esistenza7} in $L^{2}(\Omega)$.

\medskip

Let us begin to check that
$$\mathcal{K}:L^{2}(\Omega)\rightarrow L^{2}(\Omega),$$ is a well--defined compact operator. Let $g\in L^{2}(\Omega)$, set
$$w=\mathcal{K}g.$$
By \eqref{Ell:stimaenergia.1} we have
\begin{equation*}
	\alpha_0\left\Vert w\right\Vert^2_{H^{1}(\Omega)}\leq
	a_{\gamma}(w,w)=\gamma(g,w)_{L^{2}(\Omega)}\leq \gamma\left\Vert
	g\right\Vert_{L^{2}(\Omega)}\left\Vert w\right\Vert_{L^{2}(\Omega)}.
\end{equation*}
Hence \begin{equation}\label{Ell:esistenza8}\left\Vert
	\mathcal{K}g\right\Vert_{H^{1}(\Omega)}=\left\Vert
	w\right\Vert_{H^{1}(\Omega)}\leq \frac{\gamma}{\alpha_0}\left\Vert
	g\right\Vert_{L^{2}(\Omega)}.\end{equation} Therefore the operator
$\mathcal{K}$ is well--defined from $L^{2}(\Omega)$ in itself. In addition $\mathcal{K}$ is compact. As a matter of fact, let $M$ be a bounded set of
$L^{2}(\Omega)$, by \eqref{Ell:esistenza8}, we have that
$\mathcal{K}(M)$ is bounded in $H^{1}(\Omega)$ so that the Rellich--Kondrachov Theorem
implies that $\mathcal{K}(M)$ is relatively compact in $L^{2}(\Omega)$.

Let us apply Theorem \ref{Fredholm}. By proposition (iv) of such a  Theorem
we have tthe following alternative:

\medskip

\noindent either

\smallskip

(j)  the equation

\begin{equation}\label{Ell:esistenza9}
	u-\mathcal{K}u=\widetilde{h}.
\end{equation}
has a unique solution in $L^{2}(\Omega)$, for every $\widetilde{h}\in L^{2}(\Omega)$

\medskip

\noindent or

\smallskip

(jj) there exists at least a not identically vanishing solution (in $L^{2}(\Omega)$) to the equation
\begin{equation}\label{Ell:esistenza10}
	u-\mathcal{K}u=0.
\end{equation}

\medskip

If  proposition (j) holds true, then, as  
\eqref{Ell:esistenza-p5} and \eqref{Ell:esistenza7} are equivalent, we have that there exists a unique solution to problem \eqref{Ell:esistenza1}.
Whereas, if  proposition (jj) holds true, then, by (i) and (v) of
Theorem \ref{Fredholm}, we have that the subspace
$N\neq\{0\}$ of solutions of \eqref{Ell:esistenza10} (hence,
of the solutions to \eqref{Ell:esistenza2}) has finite dimension.
Let us observe that, in the latter case we have
\begin{equation}\label{Ell:esistenza10bis}
	\gamma\neq 0.
\end{equation}
otherwise we should have $\mathcal{K}=0$ and, by \eqref{Ell:esistenza10},
$N=\{0\}$. Moreover the dimension of $N$ is equal to the dimension
of $N^{\star}$, where $N^{\star}$ is the subspace of the solutions to
\begin{equation}\label{Ell:esistenza11}
	u-\mathcal{K}^{\star}u=0.
\end{equation}

At this point, to conclude (b), let us examine what relationship
holds true between $\mathcal{K}^{\star}$ and $L^{\star}$.

\medskip

\textbf{Claim.} Let us denote

\begin{equation}\label{Ell:esistenza11-0}
	L_{\gamma}^{\star}=L^{\star}+\gamma,
\end{equation}
we have

\begin{equation}\label{Ell:esistenza11bis}
	\mathcal{K}^{\star}=\gamma \left( L_{\gamma}^{\star}\right)^{-1}.
\end{equation}

\textbf{Proof of the Claim.} Firstly, recall that 
\begin{equation*}
	\begin{aligned}
		&\forall g\in L^{2}(\Omega)\quad \exists u\in H^{1}_0(\Omega)\mbox{
		(unique) such that } \\&
	a_{\gamma}(u,v)=(g,v)_{L^{2}(\Omega)},\quad\forall v\in
	H^{1}_0(\Omega)
	\end{aligned}
\end{equation*}
that is
\begin{equation}\label{Ell:esistenza12}
	u=L^{-1}_{\gamma} g.
\end{equation}
By the Lax--Milgram Theorem we have
\begin{equation*}
	\begin{aligned}
			&\forall \widetilde{g}\in L^{2}(\Omega)\quad \exists \widetilde{u}\in
	H^{1}_0(\Omega)\mbox{ (unique) such that } \\&
	a^{\star}_{\gamma}\left(\widetilde{u},\widetilde{v}\right)=\left(\widetilde{g},\widetilde{v}\right)_{L^{2}(\Omega)},\quad\forall
	\widetilde{v}\in H^{1}_0(\Omega)
	\end{aligned}
\end{equation*}
that is
\begin{equation}\label{Ell:esistenza13}
	\widetilde{u}=\left(L^{\star}_{\gamma}\right)^{-1} \widetilde{g}.
\end{equation}
Now, let us recall that

\begin{equation*}
	a^{\star}_{\gamma}\left(\widetilde{u},\widetilde{v}\right)=a_{\gamma}\left(\widetilde{v},\widetilde{u}\right),\quad\forall
	\widetilde{v}\in H^{1}_0(\Omega),
\end{equation*}
and let us choose  $\widetilde{v}=u$, where $u$ is given by \eqref{Ell:esistenza12}. We get

\begin{equation*}
	\left(\widetilde{g},u\right)_{L^{2}(\Omega)}=a^{\star}_{\gamma}\left(\widetilde{u},u\right)=a_{\gamma}\left(u,\widetilde{u}\right)=\left(g,\widetilde{u}\right)_{L^{2}(\Omega)}.
\end{equation*}
Hence,  \eqref{Ell:esistenza12} and \eqref{Ell:esistenza13} give

\begin{equation*}
	\left(\widetilde{g}, L^{-1}_{\gamma}
	g\right)_{L^{2}(\Omega)}=\left(\left(L^{\star}_{\gamma}\right)^{-1}
	\widetilde{g},g\right)_{L^{2}(\Omega)}, 
\end{equation*}
for every $\widetilde{g}\in L^{2}(\Omega)$ and for every $g \in L^{2}(\Omega)$. Now, as $\mathcal{K} g=\gamma L^{-1}_{\gamma} g$ we obtain
\eqref{Ell:esistenza11bis}. The Claim is proved.

\medskip

By \eqref{Ell:esistenza11-0} and \eqref{Ell:esistenza11bis} (taking into account that $\mathcal{K}^{\star}$ assumes its values in
$H^{1}_0(\Omega)$) we have the equivalences 

\begin{equation*}
	v-\mathcal{K}^{\star}v=0\Longleftrightarrow
	L_{\gamma}^{\star}v-\gamma v=0\Longleftrightarrow L^{\star}v=0.
\end{equation*}
Therefore, the solutions \eqref{Ell:esistenza11} to are all and only
the solutions to \eqref{Ell:esistenza3} and by that also (b) is
proved

\medskip

Now, we prove (iii). Proposition (iii) of Theorem \ref{Fredholm}
implies that the boundary value problem  \eqref{Ell:esistenza-p5} (which,
we recall, is equivalent to equation \eqref{Ell:esistenza7})
admits a solution if and only if
\begin{equation*}
	\left(h, v\right)_{L^{2}(\Omega)}=0,\quad\forall v\in N^{\star} .
\end{equation*}
On the other hand, by $v=\mathcal{K}^{\star}v$,
\eqref{Ell:esistenza6} and by \eqref{Ell:esistenza6bis}, we have

\begin{equation*}
	\left(f, v\right)_{L^{2}(\Omega)}=\left(f,
	\mathcal{K}^{\star}v\right)_{L^{2}(\Omega)}=\left(\mathcal{K}f,
	v\right)_{L^{2}(\Omega)}=\gamma\left(h, v\right)_{L^{2}(\Omega)}
\end{equation*}
and, taking into account \eqref{Ell:esistenza10bis}, we have that problem
\eqref{Ell:esistenza1} has a solution in $H_0^{1}(\Omega)$ if and only if
\begin{equation*}
	(f,v)_{L^{2}(\Omega)}=0,\quad\forall v\in N^{\star}.
\end{equation*}

$\blacksquare$
\section{The Neumann problem}\label{prob-Neumann}
Let $\Omega$ be a connected open set of $\mathbb{R}^n$  whose boundary is of class $C^{0,1}$ and let $A$ be a matrix whose entries are functions of $L^{\infty}(\Omega)$. Let us assume that $A$ satisfies
\eqref{gamma-eq}. Formally the Neumann problem \index{variational formulation of Neumann problem} for the equation
$$-\mbox{div}(A\nabla u)=F,$$ may be written as follows
\begin{equation}\label{Neumann}
	\begin{cases}
		-\mbox{div}(A\nabla u)=F, \quad\mbox{ in } \Omega, \\
		\\
		A\nabla u\cdot \nu=g, \quad\mbox{ on }
		\partial\Omega.
	\end{cases}
\end{equation}
Concerning the variational formulation we first have to specify that
$$F\in \left(H^{1}(\Omega)\right)'\quad\mbox{ and }\quad g\in
H^{-1/2}(\partial\Omega), $$ where $\left(H^{1}(\Omega)\right)'$ is the dual space of $H^{1}(\Omega)$. Having done this, arguing similarly to the Dirichlet problem, we 
formulate the Neumann problem as follows:

\medskip

Determine $u\in H^{1}(\Omega)$ such that

\begin{equation}\label{variazNeumann}
	\begin{cases}
		\int_{\Omega} A\nabla u\cdot \nabla v  dx=F(v)+\langle g, v \rangle_{H^{-1/2}, H^{1/2}}, \quad\forall v \in H^{1}(\Omega), \\
		\\
		u\in H^{1}(\Omega),
	\end{cases}
\end{equation}
where we mean 
$$\langle g, v \rangle_{H^{-1/2}, H^{1/2}}=\langle g, T(v) \rangle_{H^{-1/2},
	H^{1/2}},$$ here $T(v)$ is the trace of $v$ on $\partial
\Omega$.

\medskip
Let us notice at once that, by setting $v=1$, in \eqref{variazNeumann} we have
that a necessary condition (and, as we will see in Theorem
\ref{esistenza-unic-Neum}, also sufficient) to ensure that the
problem \eqref{variazNeumann} admits solutions, is

\begin{equation}\label{cond-Neumann}
	F(1)+\langle g, 1 \rangle_{H^{-1/2}, H^{1/2}}=0.
\end{equation}
Also, we notice that if $u_0 \in H^{1}(\Omega)$ is a solution to 
problem \eqref{variazNeumann}, then all the solutions to
\eqref{variazNeumann} are given by
$$u_0+C,$$
where $C$ is any constant. Indeed, it is immediately checked that $u_0+C$ is a solution to \eqref{variazNeumann}.
Conversely, if $u\in H^{1}(\Omega)$ is a solution to
\eqref{variazNeumann}, then
$$\int_{\Omega} A\nabla (u-u_0)\cdot \nabla v  dx=0, \quad\forall v \in
H^{1}(\Omega).$$ Now, we choose $v=u-u_0$ and 
\eqref{gamma-eq} gives
$$\lambda^{-1}\int_{\Omega} |\nabla (u-u_0)|^2\leq \int_{\Omega} A\nabla (u-u_0)\cdot \nabla(u-u_0)
dx=0.$$ Since $\Omega$ is connected,  we obtain that $u-u_0$ is a constant in $\Omega$. 

Thus, to ensure the uniqueness to the Neumann problem we may formulate it as follows

\begin{equation}\label{variazNeumann-1}
	\begin{cases}
		\int_{\Omega} A\nabla u\cdot \nabla v  dx=F(v)+\langle g, v \rangle_{H^{-1/2}, H^{1/2}}, \quad\forall v \in H^{1}(\Omega), \\
		\\
		u\in \left\{w\in H^{1}(\Omega):\quad \int_{\Omega}wdx=0\right\}.
	\end{cases}
\end{equation}

\medskip

Concerning the existence,  we have

\begin{theo}\label{esistenza-unic-Neum}
Let $\Omega$ be a connected open set of $\mathbb{R}^n$  whose boundary is of class $C^{0,1}$. Let us assume  that $A\in L^{\infty}(\Omega;\mathbb{M}(n))$ and that $A$ satisfies \eqref{gamma-eq}. Let us assume that $F\in \left(H^{1}(\Omega)\right)'$ and
	$g\in H^{-1/2}(\partial\Omega)$ satisfy \eqref{cond-Neumann}.
	
	Then problem \eqref{variazNeumann-1} has a unique solution and the following inequality holds true
	
	\begin{equation}\label{stima-Neumann-completo}
		\left\Vert u\right\Vert_{H^{1}(\Omega)}\leq C\left(\left\Vert
		F\right\Vert_{\left(H^{1}(\Omega)\right)'}+\left\Vert
		g\right\Vert_{H^{-1/2}(\partial\Omega)}\right),
	\end{equation}
	where $C$ depends on $\lambda$ and $\Omega$ only.
\end{theo}
\textbf{Proof.} Set
$$\widetilde{H}:=\left\{w\in H^{1}(\Omega):\quad
\int_{\Omega}wdx=0\right\},$$ Theorem
\ref{Poincar1} implies that $\widetilde{H}$ is a Hilbert space equipped with the 
norm
$$\left\Vert w\right\Vert_{\widetilde{H}}=\left(\int_{\Omega}|\nabla
w|^2dx\right)^{1/2}.$$ Moreover, the bilinear form on
$\widetilde{H}$
\begin{equation}\label{stima-Neumann-completo-1}
	a(u,v)=\int_{\Omega} A\nabla u\cdot \nabla v  dx,
\end{equation}
is coercive and continuous. Now we check that the linear functional
$$\widetilde{H}\ni v\rightarrow \widetilde{F}(v):= F(v)+\langle g, v \rangle_{H^{-1/2},
	H^{1/2}}\in \mathbb{R},$$ is well--defined and continuous on
$\widetilde{H}$. As a matter of fact, by Theorem \ref{Poincar1} we have

\begin{equation}\label{variazNeumann-2}
	|F(v)|\leq \left\Vert F\right\Vert_{\left(H^{1}(\Omega)\right)'}\left\Vert
	v\right\Vert_{H^{1}(\Omega)}\leq C\left\Vert
	F\right\Vert_{H^{-1}(\Omega)}\left\Vert
	v\right\Vert_{\widetilde{H}},\mbox{ }\forall v\in \widetilde{H},
\end{equation}
where $C$ depends on $\Omega$ only. Moreover, recalling that
$H^{-1/2}(\Omega)$ is the dual space of $H^{1/2}(\Omega)$ (compare
Section \ref{H--1/2}), inequality (ii) of
Theorem \ref{traccia} gives

\begin{equation}\label{variazNeumann-3}
	\begin{aligned}
		\left|\langle g, v \rangle_{H^{-1/2}, H^{1/2}}\right|&\leq
		\left\Vert g\right\Vert_{H^{-1/2}(\partial\Omega)}\left\Vert
		v\right\Vert_{H^{1/2}(\partial\Omega)}\leq \\&\leq C\left\Vert
		g\right\Vert_{H^{-1/2}(\partial\Omega)}\left\Vert
		v\right\Vert_{H^{1}(\Omega)}\leq\\&\leq C\left\Vert
		g\right\Vert_{H^{-1/2}(\partial\Omega)}\left\Vert
		v\right\Vert_{\widetilde{H}},
	\end{aligned}
\end{equation}
where $C$ depends by $\Omega$ only. Therefore, \eqref{variazNeumann-2} and
\eqref{variazNeumann-3} give

\begin{equation}\label{variazNeumann-4}
	|\widetilde{F}(v)|\leq C\left(\left\Vert
	F\right\Vert_{H^{-1}(\Omega)}+\left\Vert
	g\right\Vert_{H^{-1/2}(\partial\Omega)}\right)\left\Vert
	v\right\Vert_{\widetilde{H}},\mbox{ }\forall v\in \widetilde{H}.
\end{equation}
Now, since bilinear form  \eqref{stima-Neumann-completo-1}
is continuous and coercive on $\widetilde{H}$ and since 
\eqref{variazNeumann-4} holds, $\widetilde{F}$ is a bounded linear functional on $\widetilde{H}$. Therefore by the Lax--Milgram Theorem we have that there exists a unique $u\in \widetilde{H}$ which satisfies

\begin{equation}\label{variazNeumann-5}
	a(u,v)=\widetilde{F}(v),\quad\forall v\in \widetilde{H}.
\end{equation}
Moreover

\begin{equation}\label{stima-Neumann-completo-2}
	\left\Vert u\right\Vert_{H^{1}(\Omega)}\leq C \left\Vert \nabla
	u\right\Vert_{L^{2}(\Omega)}\leq C\left(\left\Vert
	F\right\Vert_{\left(H^{1}(\Omega)\right)'}+\left\Vert
	g\right\Vert_{H^{-1/2}(\partial\Omega)}\right),
\end{equation}
where $C$ depends on $\lambda$ and $\Omega$ only.

Now, let $v$ be any function of $H^1(\Omega)$ and let us denote

$$v_{\Omega}=\frac{1}{|\Omega|}\int_{\Omega}v dx,\quad \widetilde{v}=v-v_{\Omega}.$$
Since $\widetilde{v}\in \widetilde{H}$, by \eqref{cond-Neumann} and
\eqref{variazNeumann-5}  we obtain
\begin{equation}\label{variazNeumann-6}
	\begin{aligned}
		&\int_{\Omega} A\nabla u\cdot \nabla v  dx=\int_{\Omega} A\nabla
		u\cdot \nabla\widetilde{v}
		dx=\\&=a\left(u,\widetilde{v}\right)=F\left(\widetilde{v}\right)+\langle
		g, \widetilde{v} \rangle_{H^{-1/2}, H^{1/2}}=\\&= F(v)+\langle g, v
		\rangle_{H^{-1/2}, H^{1/2}}-v_{\Omega}\left(F(1)+\langle g, 1
		\rangle_{H^{-1/2}, H^{1/2}}\right)=\\&= F(v)+\langle g, v
		\rangle_{H^{-1/2}, H^{1/2}}.
	\end{aligned}
\end{equation}
Therefore $u$ is a solution to problem
\eqref{variazNeumann-1}. Estimate \eqref{variazNeumann-1} follows by
\eqref{stima-Neumann-completo-2}. $\blacksquare$

\section{The Caccioppoli inequality}\label{Ell:caccioppoli}
\begin{theo}[\textbf{the Caccioppoli inequality}]\label{dis-Caccioppoli}
	\index{Theorem:@{Theorem:}!- Caccioppoli inequality@{- Caccioppoli inequality}}
	Let $x_0\in \mathbb{R}^n$ and $R>0$. Let  $A$ be a symmetric matrix whose entries
	are  measurable functions on $B_R(x_0)$. Let us assume $A$ satisfies 
	\eqref{gamma-eq} (with $\Omega=B_R(x_0)$). \\ Let  $b\in
	L^{\infty}(B_R(x_0);\mathbb{R}^n)$ and $c\in L^{\infty}(B_R(x_0))$.
	Let $u\in H_{loc}^1 (B_R(x_0))$ satisfy
	\begin{equation}\label{Cacciopp-1}
		\int_{B_R(x_0)}A\nabla u\cdot\nabla v dx=\int_{B_R(x_0)}\left(b\nabla
		u+cu\right)v dx,\quad \forall v\in H_{0}^1 (B_R(x_0)).
	\end{equation}
	If $0<r<\rho<R$, then  we have
	\begin{equation}\label{Cacciopp-2}
		\int_{B_{r}(x_0)}|\nabla u|^2 dx\leq \frac{C}{(\rho-r)^2}
		\int_{B_{\rho}(x_0)}u^2 dx,
	\end{equation}
	where $C$ depends on $\lambda$, $R\left\Vert
	b\right\Vert_{L^{\infty}(B_R(x_0);\mathbb{R}^n)}$ and $R^2\left\Vert
	c\right\Vert_{L^{\infty}(B_R(x_0);\mathbb{R}^n)}$ only.
\end{theo}

\textbf{Proof.} It is not restrictive to assume $x_0=0$. Let $\eta\in
C^{\infty}_0(B_{\rho}$) satisfy
\begin{equation}\label{Cacciopp-3}
	0\leq \eta\leq 1;\qquad \eta=1,\quad \mbox{ in } B_{r}
\end{equation}
and
\begin{equation}\label{Cacciopp-4}
	|\nabla \eta|\leq \frac{K}{\rho-r},
\end{equation}
where $K$ is a positive constant. We choose in \eqref{Cacciopp-1}
$$v=\eta^2u$$
and we have
\begin{equation}\label{Cacciopp-5}
	\int_{B_R}A\nabla u\cdot\nabla\left(\eta^2u\right)
	dx=\int_{B_R}\left(b\nabla u+cu\right)\eta^2u dx.
\end{equation}
Hence

\begin{equation*}
	\begin{aligned}
		&\int_{B_R}(A\nabla u\cdot\nabla
		u)\eta^2dx=\int_{B_R}\left(b\nabla u+cu\right)\eta^2u
		dx-\\&-2\int_{B_R}(A\nabla u\cdot\nabla\eta) \eta udx\leq\\&\leq
		\int_{B_R}\left(|b||\nabla
		u||u|\eta^2+|c|u^2\eta^2\right)dx+\\&+2\int_{B_R}(A \nabla
		u\cdot\nabla u)^{1/2}(A \nabla \eta\cdot\nabla \eta)^{1/2} |u|\eta
		dx\leq \\& \leq \int_{B_R}\left(|b||\nabla
		u||u|\eta^2+|c|u^2\eta^2\right)dx+\\&+\frac{1}{2}\int_{B_R}(A\nabla
		u\cdot\nabla u)\eta^2dx+2\int_{B_R}(A\nabla \eta\cdot\nabla
		\eta)u^2dx.
	\end{aligned}
\end{equation*}
By moving to the left--hand side the second-to-last integral  and, by estimating from above the last
integral, we have
\begin{equation}\label{Cacciopp-6}
	\begin{aligned}
		&\frac{1}{2}\int_{B_R}(A\nabla u\cdot\nabla u)\eta^2dx\leq\\&
		\leq \int_{B_R}\left(|b||\nabla
		u||u|\eta^2+|c|u^2\eta^2\right)dx+\frac{2K^2\lambda}{(\rho-r)^2}\int_{B_{\rho}}
		u^2dx.
	\end{aligned}
\end{equation}
Now let us estimate from above the first integral on the right hand side of \eqref{Cacciopp-6}. We obtain, for $\varepsilon >0$ to be choosen,

\begin{equation}\label{correct-13-3-23-1}
	\begin{aligned}
		\int_{B_R}\left(|b||\nabla u||u|\eta^2+|c|u^2\eta^2\right)dx&\leq
		\frac{\varepsilon}{2}\int_{B_R}|\nabla
		u|^2\eta^2dx+\\&+\frac{1}{2\varepsilon}\left\Vert
		b\right\Vert^2_{L^{\infty}(B_R;\mathbb{R}^n)}\int_{B_R}u^2\eta^2dx+\\&+\left\Vert
		c\right\Vert_{L^{\infty}(B_R)}\int_{B_R}u^2\eta^2dx\leq\\&\leq
		\frac{\varepsilon\lambda}{2}\int_{B_R}(A\nabla u\cdot\nabla
		u)\eta^2dx+\\&+C_{\varepsilon}\int_{B_R}u^2\eta^2dx,
	\end{aligned}
\end{equation}
where

$$C_{\varepsilon}=\frac{1}{2\varepsilon}\left\Vert
b\right\Vert^2_{L^{\infty}(B_R,\mathbb{R}^n)}+\left\Vert
c\right\Vert_{L^{\infty}(B_R)}.$$ Using inequality \eqref{correct-13-3-23-1} in \eqref{Cacciopp-6}, after a few easy calculations, we have 

\begin{equation}\label{Cacciopp-7}
	\begin{aligned}
		&\frac{1}{2}\left(1-\varepsilon\lambda\right)\int_{B_R}(A\nabla
		u\cdot\nabla u)\eta dx\leq\\& \leq
		C_{\varepsilon}\int_{B_R}u^2\eta^2dx+\frac{2K^2\lambda}{(\rho-r)^2}\int_{B_{\rho}}
		u^2dx.
	\end{aligned}
\end{equation}
Now we choose
$$\varepsilon=\varepsilon_0:=\frac{1}{2\lambda},$$ and we get

$$C_{\varepsilon_0}=\lambda\left\Vert
b\right\Vert^2_{L^{\infty}(B_R;\mathbb{R}^n)}+\left\Vert
c\right\Vert_{L^{\infty}(B_R)}$$ and by \eqref{Cacciopp-7} we have

\begin{equation}\label{Cacciopp-8}
	\begin{aligned}
		\frac{\lambda^{-1}}{4}\int_{B_{r}}|\nabla
		u|^2dx&\leq\frac{1}{4}\int_{B_R}(A\nabla u\cdot\nabla
		u)\eta^2dx\leq\\& \leq\frac{K_1}{(\rho-r)^2}\int_{B_{\rho}} u^2,
	\end{aligned}
\end{equation}
where
$$K_1=2K^2\lambda+R^2C_{\varepsilon_0}.$$ By \eqref{Cacciopp-8} we obtain
immediately

\begin{equation*}
	\begin{aligned}
		\int_{B_{r}}|\nabla u|^2dx
		\leq\frac{4K_1\lambda}{(\rho-r)^2}\int_{B_{\rho}} u^2dx,
	\end{aligned}
\end{equation*}
so that \eqref{Cacciopp-2} follows. $\blacksquare$

\bigskip

\underline{\textbf{Exercise 1.}} Let $x_0\in \mathbb{R}^n$ and $R>0$. Let $L$ be the
operator

\begin{equation}\label{Ell:EserCacc1}
	Lu=-\sum_{j,k=1}^n\partial_j\left( a^{jk}\partial_k u
	+d^ju\right)+\sum_{j=1}^nb^{j}\partial_j u+cu,\end{equation} where  $A\in
L^{\infty}(B_{R}(x_0);\mathbb{M}(n))$,
$A=\left\{a^{jk}\right\}_{j,k=1}^n$, satisfies \eqref{gamma-eq} and \\
$b^j,d^j, c\in L^{\infty}(B_{R}(x_0))$, for $j=1,\cdots, n$. Let
$f\in L^{2}(B_{R}(x_0))$ and let us assume that $u\in H^1(B_{R}(x_0))$ is a weak solution to

$$Lu=f,\quad \mbox{ in } B_{R}(x_0).$$
Prove that, if $0<r<\rho<R$ then we have
\begin{equation}\label{Cacciopp-9}
	\int_{B_{r}(x_0)}|\nabla u|^2 dx\leq \frac{C}{(\rho-r)^2}
	\int_{B_{\rho}(x_0)}u^2 dx+C\rho^2\int_{B_{\rho}(x_0)}f^2dx,
\end{equation}
where $C$ depends on $\lambda$, $\left\Vert
A\right\Vert_{L^{\infty}(B_R(x_0);\mathbb{M}(n))}$,  $R\left\Vert
d\right\Vert_{L^{\infty}(B_R(x_0);\mathbb{R}^n)}$, $R\left\Vert
b\right\Vert_{L^{\infty}(B_R(x_0;\mathbb{R}^n)}$ and $R^2\left\Vert
c\right\Vert_{L^{\infty}(B_R(x_0)}$ only. $\clubsuit$

\bigskip

\underline{\textbf{Exercise 2.}} Let $R>0$ and $x_0\in \{(x',0):\mbox{ } x'\in
\mathbb{R}^{n-1}\}$. Let $L$ be operator \eqref{Ell:EserCacc1} and let 
$u\in H^1(B^+_{R}(x_0))$ satisfy

\begin{equation*}
	\begin{cases}
		Lu=f, \quad\mbox{ in } B^+_{R}(x_0),\mbox{ in weak sense,} \\
		\\
		u(x',0)=0,\quad x'\in B'_R(x_0) \mbox{ in the traces sense,}
	\end{cases}
\end{equation*}
then
\begin{equation}\label{Cacciopp-10}
	\int_{B^+_{r}(x_0)}|\nabla u|^2 dx\leq \frac{C}{(\rho-r)^2}
	\int_{B^+_{\rho}(x_0)}u^2 dx+C\rho^2\int_{B^+_{\rho}(x_0)}f^2dx,
\end{equation}
where $C$ depends on $\lambda$, $\left\Vert
A\right\Vert_{L^{\infty}(B^+_R(x_0);\mathbb{M}(n))}$,  $R\left\Vert
d\right\Vert_{L^{\infty}(B^+_R(x_0);\mathbb{R}^n)}$, $R\left\Vert
b\right\Vert_{L^{\infty}(B^+_R(x_0);\mathbb{R}^n)}$ and $R^2\left\Vert
c\right\Vert_{L^{\infty}(B^+_R(x_0))}$ only. $\clubsuit$

\bigskip

\underline{\textbf{Exercise 3.}} Let $x_0$, $R$ and $A$ be like Exercise 2.

(a) Give the variational formulation of problem

\begin{equation}\label{Neuma000}
	\begin{cases}
		\mbox{div}(A\nabla u)=0, \quad\mbox{ in } B^+_R(x_0), \\
		\\
		(A\nabla u)(x',0)\cdot e_n=0, \quad\mbox{ for } x'\in B'_R(x_0).
	\end{cases}
\end{equation}

(b) Prove that

\begin{equation}\label{Cacciopp-11}
	\int_{B^+_{r}(x_0)}|\nabla u|^2 dx\leq \frac{C}{(\rho-r)^2}
	\int_{B^+_{\rho}(x_0)}u^2 dx,
\end{equation}
where $C$ depends on $\lambda$ and $\left\Vert
A\right\Vert_{L^{\infty}(B^+_R(x_0);\mathbb{M}(n))}$ only. $\clubsuit$

\section{The regularity theorems} \label{regolarit} In
Section \ref{Lax-Milgram} (Theorem \ref{teorema-es-unic}) we have proved that if $\Omega$ is a bounded open set of
$\mathbb{R}^n$, $A\in L^{\infty}(\Omega;\mathbb{M}(n))$ satisfies
\eqref{gamma-eq} and $F\in H^{-1}(\Omega)$, then the
Dirichlet problem

\begin{equation}\label{Dirich-reg}
	\begin{cases}
		-\mbox{div}(A\nabla u)=F, \quad\mbox{ in } \Omega, \\
		\\
		u=0, \quad\mbox{ on } \partial\Omega,
	\end{cases}
\end{equation}
is well--posed in $H_0^{1}( \Omega)$.

It is natural to ask whether with more restrictive assumptions on the
data $\Omega, A$ and $F$, $u$ is more regular. More precisely, we ask whether there exists $k>1$ such that $u\in H^{k}( \Omega)$. In carrying out this
investigation it is convenient to distinguish between the \textbf{regularity
in the interior} and \textbf{regularity at the boundary} \index{regularity in the interior, at the boundary}.
In the investigation of the regularity in the interior we are interested in whether for
some $k>1$ we have $u\in H_{loc}^{k}( \Omega)$, while in the investigation of regularity at the boundary we are interested in knowing whether
for some $k>1$ it happens that for every $x_0\in \partial\Omega$ there exists a neighborhood of $x_0$,
 $\mathcal{U}$,  such that $u_{|\mathcal{U}} \in
H^{k}(\Omega\cap \mathcal{U})$.
As we should expect, in the study of
regularity in the interior, the regularity of
$\partial\Omega$ plays no role. In contrast, in the study of regularity at the
boundary, the regularity of $\partial\Omega$ plays a crucial role.

Before going on to the rigorous treatment, let us illustrate in a
rough manner the main idea that drives the study of the
regularity in the interior; similar arguments can be made for the
regularity at the boundary.

Let us consider the equation
\begin{equation}\label{Ell:reg1}
	-\Delta u=f,\quad \mbox{ in } \Omega
\end{equation}
where $f\in L^2(\Omega)$. Let $u\in H^1(\Omega)$ be a solution to 
\eqref{Ell:reg1}; that is.
\begin{equation}\label{Ell:reg2}
	\int_{\Omega}\nabla u\cdot \nabla vdx=\int_{\Omega}fvdx, \quad
	\forall v\in H_0^1(\Omega).
\end{equation}

\underline{\textbf{Let us suppose}} that we know $u$ be sufficiently regular
(say $u\in H_{loc}^3(\Omega)$) so that the operations that we will make are allowed. 

Let $x_0\in \Omega$ and $R>0$ satisfy 
$B_{2R}(x_0)\subset \Omega$. Let $\eta\in C^{\infty}_0(B_{2R}(x_0)$
such that
\begin{equation}\label{Ell:reg2bis}
	0\leq \eta\leq 1;\qquad \eta=1,\quad \mbox{ in } B_{R}(x_0)
\end{equation}
and
\begin{equation}\label{Ell:reg2ter}
	|\nabla \eta|\leq \frac{K}{R},
\end{equation}
where $K$ is a positive constant. Let $k\in \{1,\cdots,n\}$.
Multiply both the sides of \eqref{Ell:reg1} by
$\partial_k\left(\eta^2\partial_k u\right)$ and integrate over $\Omega$
or, equivalently, choose in \eqref{Ell:reg2})
\begin{equation}\label{Ell:reg3bis} v=\partial_k\left(\eta^2\partial_ku\right)\end{equation}
 obtaining

\begin{equation}\label{Ell:reg3}
	\int_{\Omega} \Delta u
	\partial_k\left(\eta^2\partial_ku\right)dx=\int_{\Omega}f\partial_k\left(\eta^2\partial_ku\right)dx.
\end{equation}
Let us consider the left--hand side of \eqref{Ell:reg3}, integration by per parts yields
\begin{equation*}
	\begin{aligned}
		\int_{\Omega} \Delta u
		\partial_k\left(\eta^2\partial_ku\right)dx&=\int_{\Omega} \sum_{j=1}^n\partial^2_j u
		\partial_k\left(\eta^2\partial_ku\right)dx=\\&=-\int_{\Omega} \sum_{j=1}^n\partial_j\left(\partial^2_{jk}
		u\right) \eta^2\partial_kudx=\\&=\int_{\Omega}
		\sum_{j=1}^n\partial^2_{jk} u
		\partial_{j}\left(\eta^2\partial_ku\right)dx=\\&=\int_{\Omega}
		\sum_{j=1}^n\left|\partial^2_{jk} u\right|^2\eta^2dx+\\&+
		2\int_{\Omega} \sum_{j=1}^n\left(\partial^2_{jk}
		u\right)\eta\partial_j\eta
		\partial_k u dx.
	\end{aligned}
\end{equation*}
Concerning the right--hand side of \eqref{Ell:reg3}, we get
\begin{equation*}
	\begin{aligned}
		\int_{\Omega}f\partial_k\left(\eta^2\partial_ku\right)dx&=\int_{\Omega}f\eta^2\partial^2_kudx+2\int_{\Omega}f\eta\partial_k\eta\partial_kudx.
	\end{aligned}
\end{equation*}
Using in \eqref{Ell:reg3} the last two obtained equalities and summing up over $k$, we have

\begin{equation}\label{Ell:reg4}
	\begin{aligned}
		\int_{\Omega} \sum_{j,k=1}^n\left|\partial_{jk}^2
		u\right|^2&\eta^2dx=-2\int_{\Omega}
		\sum_{j,k=1}^n\left(\partial^2_{jk} u\right)\eta\partial_j\eta
		\partial_k u dx+\\&+
		\int_{\Omega}f\eta^2\sum_{k=1}^n\partial^2_kudx+2\int_{\Omega}f\eta\sum_{k=1}^n\partial_k\eta\partial_kudx:=I.
	\end{aligned}
\end{equation}
Let $\varepsilon>0$ to be choosen later on, let us denote by $\partial^2u$ the Hessian matrix $\left\{\partial_{jk}^2 u\right\}_{j,k=1}^n$. We have
\begin{equation*}
	\begin{aligned}
		I&\leq \varepsilon\int_{\Omega}\left|\partial^2
		u\right|^2\eta^2dx+\frac{1}{\varepsilon}\int_{\Omega}\left|\nabla
		u\right|^2\left|\nabla \eta\right|^2dx+\\&+
		\frac{\varepsilon}{2}\int_{\Omega}\left|\partial^2
		u\right|^2\eta^2dx+\frac{1}{2\varepsilon}\int_{\Omega}\left|f\right|^2\eta^2dx+\\&+
		\int_{\Omega}\left|f\right|^2\eta^2dx+\int_{\Omega}\left|\nabla
		\eta\right|^2 \left|\nabla u\right|^2dx.
	\end{aligned}
\end{equation*}
Now, in \eqref{Ell:reg4}, we move to the left--hand side the terms
that contain the second derivatives, and we get

\begin{equation*}
	\begin{aligned}
		\left(1-\frac{3\varepsilon}{2}\right) \int_{\Omega}\left|\partial^2
		u\right|^2\eta^2dx&\leq
		\left(1+\frac{1}{2\varepsilon}\right)\int_{\Omega}\left|\nabla
		u\right|^2\left|\nabla \eta\right|^2dx+\\&+
		\left(1+\frac{1}{2\varepsilon}\right)\int_{\Omega}\left|f\right|^2\eta^2dx.
	\end{aligned}
\end{equation*}
At this point, we choose $\varepsilon=\frac{1}{3}$ and by
\eqref{Ell:reg2bis}, \eqref{Ell:reg2ter} we get

\begin{equation}\label{Ell:reg5}
	\int_{B_{R}(x_0)}\left|\partial^2 u\right|^2dx\leq
	\frac{5K^2}{R^2}\int_{B_{2R}(x_0)}\left|\nabla u\right|^2dx+
	5\int_{B_{2R}(x_0)}\left|f\right|^2dx.
\end{equation}
We observe that \eqref{Ell:reg5} allows us to estimate the
second derivatives of $u$ in $L^2\left(B_{R}(x_0)\right)$ by means of
the finite quantity that occurs on the right. \textbf{It is evident that
this estimate by itself do not provide  a proof}
that $u\in H^2\left(B_{R}(x_0)\right)$, $j,k=1,\cdots, n$, since
to obtain the estimate we exploited a regularity of
$u$ even greater than was proved (!). However, in the
rigorous proofs that we will present soon in this
Chapter, we will "retrace", in a sense, the previous steps
by considering as test function, instead of the \eqref{Ell:reg3bis},
the function

\begin{equation*}
	v=-\delta^{-h}_k\left(\eta^2 \delta^{h}_ku\right)
\end{equation*}
where $\delta^{-h}_k$ and $\delta^{h}_k$, are the difference quotients
studied in Section \ref{quozienti di differenze}.

\subsection{The regularity theorems in the interior} \label{regolarit-interno}

The Main Theorem of the present Subsection is the following.

\begin{theo} [\textbf{regularity in the interior}]\label{reg-interno}
	\index{Theorem:@{Theorem:}!- regularity in the interior@{- regularity in the interior}}
	Let $\Omega$ be a bounded open set of $\mathbb{R}^n$. Let
	\begin{equation}\label{Ell:E1.3}
		f\in L^2(\Omega).
	\end{equation}
	Let $A$ be a symmetric matrix. Let us assume that $A$ satisfies \eqref{gamma-eq},  \\ $A\in
	C^{0,1}(\Omega;\mathbb{M}(n))$ and it satisfies
	
	\begin{equation}\label{Ell:E1.2}
		\left|A(x)-A(y)\right|\leq E|x-y|,\quad \forall x,y\in
		\Omega,
	\end{equation}
	where $E$ is a positive number. Let us assume $u\in H^1(\Omega)$ is a
	solution to
	
	\begin{equation}\label{Ell:E1.4}
		-\mbox{div}(A\nabla u)=f, \quad\mbox{ in } \Omega.
	\end{equation}
	Then we have
	\begin{equation}\label{Ell:E1.5}
		u\in H^2_{loc}(\Omega),
	\end{equation}
	and, for any $B_{2R}(x_0)\subset \Omega$, the following estimate holds true
	\begin{equation}\label{Ell:E2.1}
		\begin{aligned}
			\sum_{|\alpha|\leq
				2}R^{2|\alpha|}\int_{B_R(x_0)}\left|\partial^{\alpha}u
			\right|^2dx&\leq C\left(1+E^2R^2\right)\int_{B_{2R}(x_0)}u
			^2dx+\\&+CR^4\int_{B_{2R}(x_0)}f ^2dx,
		\end{aligned}
	\end{equation}
	where $C$ depends on $\lambda$ only.
\end{theo}

\textbf{Proof.} It is not restrictive to assume  $0\in \Omega$ and $R<
\frac{1}{2}$dist ($0,\partial\Omega)$. Let $\eta\in
C^{\infty}_0(B_{3R/2})$ satisfy 
\begin{equation}\label{Ell:E2.1-0}
	0\leq \eta\leq 1;\qquad \eta=1,\quad \mbox{ in } B_{R}
\end{equation}
and
\begin{equation}\label{Ell:E2.1-00}
	|\nabla \eta|\leq \frac{K}{R},
\end{equation}
where $K$ is a positive constant.

Since $u\in H^1(\Omega)$ satisfies \eqref{Ell:E1.4}, we have

\begin{equation}\label{Ell:E2.2}
	\int_{\Omega} A\nabla u\cdot \nabla v  dx=\int_{\Omega}fvdx,
	\quad\forall v \in H_0^{1}(\Omega).
\end{equation}

Let $h\in \left(-\frac{R}{8},\frac{R}{8}\right)\setminus\{0\}$. Let us note that, if $w_1,w_2\in
H^1(\Omega)$ and supp $w_1\subset B_{3R/2}$ (or supp $w_2\subset
B_{3R(2}$), then, for any $k\in \left\{1,\cdots,n\right\}$, we have

\begin{equation}\label{Ell:E3.2}
	\int_{\Omega}w_1\delta^{-h}_kw_2dx=
	-\int_{\Omega}w_2\delta^{h}_kw_1dx
\end{equation}
and
\begin{equation}\label{Ell:E3.3}
	\delta^{h}_k\left(w_1w_2\right)=w_1^h\delta^{h}_kw_2+w_2\delta^{h}_kw_1,
\end{equation}
where $w_1^h(x)=w_1(x+he_k)$. Concerning \eqref{Ell:E3.2},
just argue like in the the Claim of the
proof of Theorem \ref{Sob:teo1.9}. While equalty \eqref{Ell:E3.2} follows
easily by 
\begin{equation*}
	\begin{aligned}
		h\delta^{h}_k\left(w_1w_2\right)&=w_1(x+he_k)w_2(x+he_k)-w_1(x)w_2(x)=\\&=
		w_1(x+he_k)w_2(x+he_k)-w_1(x+he_k)w_2(x)+\\&+w_1(x+he_k)w_2(x)-w_1(x)w_2(x)=\\&=
		h\left(w^h_1\delta^{h}_kw_2+w_2\delta^{h}_kw_1\right).
	\end{aligned}
\end{equation*}

Now, let us choose as test function in \eqref{Ell:E2.2}

\begin{equation}\label{Ell:E3.1}
	v=-\delta^{-h}_k\left(\eta^2 \delta^{h}_ku\right).
\end{equation}
We get
\begin{equation}\label{Ell:E4.1}
	\begin{aligned}
		\int_{\Omega} A\nabla u\cdot \nabla v  dx&=-\int_{\Omega} A\nabla
		u\cdot\left[\delta^{-h}_k\nabla\left(\eta^2
		\delta^{h}_ku\right)\right]dx=\\&=
		\int_{\Omega}\delta^{h}_k\left(A\nabla u\right)\cdot
		\nabla\left(\eta^2 \delta^{h}_ku\right)dx=\\&=\int_{\Omega}A^h
		\left(\delta^{h}_k\nabla u\right)\cdot\nabla\left(\eta^2
		\delta^{h}_ku\right)dx+\\&+ \int_{\Omega}
		\left(\delta^{h}_kA\right)\nabla u\cdot\nabla\left(\eta^2
		\delta^{h}_ku\right)dx=\\&=\int_{\Omega}\left[A^h
		\left(\delta^{h}_k\nabla u\right)\cdot \left(\delta^{h}_k\nabla
		u\right)\right]\eta^2dx+\mathcal{R},
	\end{aligned}
\end{equation}
where
\begin{equation*}
	\begin{aligned}
		\mathcal{R}&=\int_{\Omega}A^h \left(\delta^{h}_k\nabla
		u\right)\cdot\left(2\eta\nabla\eta\delta^{h}_ku\right)dx+\\&+\int_{\Omega}\left(\delta^{h}_k
		A\right)\nabla u\cdot \left(\delta^{h}_k\nabla
		u\right)\eta^2dx+\\&+\int_{\Omega}\left(\delta^{h}_k A\right)\nabla
		u\cdot\left(2\eta\nabla\eta\delta^{h}_ku\right)dx.
	\end{aligned}
\end{equation*}
Now, by \eqref{gamma-eq} we get

\begin{equation}\label{Ell:E6.1}
	\begin{aligned}
		\int_{\Omega}\left[A^h \left(\delta^{h}_k\nabla u\right)\cdot
		\left(\delta^{h}_k\nabla u\right)\right]\eta^2dx\geq
		\lambda^{-1}\int_{\Omega}\left|\delta^{h}_k\nabla
		u\right|^2\eta^2dx.
	\end{aligned}
\end{equation}
Concerning $\mathcal{R}$, let $\varepsilon$ be a positive number which will choose later on.
By \eqref{Ell:E1.2} we have

\begin{equation}\label{Ell:E6.2}
	\begin{aligned}
		|\mathcal{R}|&\leq
		\frac{c_n\lambda}{R}\int_{\Omega}\left|\delta^{h}_k\nabla
		u\right|\left|\delta^{h}_k u\right|\eta dx+ E
		\int_{\Omega}\left|\delta^{h}_k\nabla u\right|\left|\nabla
		u\right|\eta^2 dx+\\&+\frac{E}{R}\int_{\Omega}\left|\nabla
		u\right|\left|\delta^{h}_k u\right|\eta^2 dx\leq\\&\leq
		\varepsilon\int_{\Omega}\left|\delta^{h}_k\nabla u\right|^2\eta^2
		dx+\\&+\frac{C}{\varepsilon}\left(R^{-2}+E^2\right)\int_{B_{3R/2}}\left(\left|\delta^{h}_k
		u\right|^2+\left|\nabla u\right|^2\right)dx,
	\end{aligned}
\end{equation}
where $c_n$ depends on $n$ only and $C$ depends on $\lambda$ and $n$ only. Now let us choose

$$\varepsilon=\frac{\lambda^{-1}}{2}$$
so that, by \eqref{Ell:E6.1} and \eqref{Ell:E6.2}, we get
\begin{equation}\label{Ell:E7.1}
	\begin{aligned}
		\int_{\Omega}&\left[A^h \left(\delta^{h}_k\nabla u\right)\cdot
		\left(\delta^{h}_k\nabla u\right)\right]\eta^2dx+\mathcal{R}\geq
		\frac{\lambda^{-1}}{2}\int_{\Omega}\left|\delta^{h}_k\nabla
		u\right|^2\eta^2dx-\\&-2C\lambda\left(R^{-2}+E^2\right)\int_{B_{3R/2}}\left(\left|\delta^{h}_k
		u\right|^2+\left|\nabla u\right|^2\right)dx.
	\end{aligned}
\end{equation}

\smallskip

Hence, \eqref{Ell:E2.2}, \eqref{Ell:E3.1}, \eqref{Ell:E4.1} and
\eqref{Ell:E7.1} yield

\begin{equation}\label{Ell:E7.2}
	\begin{aligned}
		\frac{\lambda^{-1}}{2}\int_{\Omega}\left|\delta^{h}_k\nabla
		u\right|^2\eta^2dx&\leq - \int_{\Omega}f\delta^{-h}_k\left(\eta^2
		\delta^{h}_ku\right)dx+\\&+C\left(R^{-2}+E^2\right)\int_{B_{3R/2}}\left(\left|\delta^{h}_k
		u\right|^2+\left|\nabla u\right|^2\right)dx,
	\end{aligned}
\end{equation}
where $C$ depends on $\lambda$ and $n$ only.

Now, by Theorem \ref{Sob:teo1.9}--(i) (with $V=B_{3R/2}$ and
$\Omega=B_{7R/4}$) we have
\begin{equation}\label{Ell:E7.3}
	\begin{aligned}
		\int_{B_{3R/2}}\left|\delta^{h}_k u\right|^2dx\leq
		\int_{B_{7R/4}}|\nabla u|^2dx.
	\end{aligned}
\end{equation}
Moreover

\begin{equation}\label{Ell:E8.1}
	\begin{aligned}
		&\int_{\Omega}\left|\delta^{-h}_k\left(\eta^2
		\delta^{h}_ku\right)\right|^2dx\leq
		C\int_{\Omega}\left|\nabla\left(\eta^2
		\delta^{h}_ku\right)\right|^2dx\leq\\&\leq
		C'\int_{\Omega}\left|\eta\nabla\eta\left(\eta^2
		\delta^{h}_ku\right)\right|^2dx+C'\int_{\Omega}\eta^2\left|\nabla
		\delta^{h}_ku \right|^2dx\leq\\&\leq
		\frac{C''}{R^2}\int_{B_{3R/2}}\left|\delta^{h}_k
		u\right|^2dx+C'\int_{\Omega}\eta^2\left|\delta^{h}_k\nabla u
		\right|^2dx.
	\end{aligned}
\end{equation}

\smallskip

Let now $\sigma>0$ to be choosen, \eqref{Ell:E8.1} implies

\begin{equation}\label{Ell:E9.1}
	\begin{aligned}
		&\left|\int_{\Omega}f\delta^{-h}_k\left(\eta^2
		\delta^{h}_ku\right)dx\right|\leq
		\frac{1}{2\sigma}\int_{B_{2R}}f^2dx+\frac{\sigma}{2}\int_{\Omega}\left|\delta^{-h}_k\left(\eta^2
		\delta^{h}_ku\right)\right|^2dx\leq\\&\leq
		\frac{1}{2\sigma}\int_{B_{2R}}f^2dx+\\&+
		C\sigma\left(R^{-2}\int_{B_{3R/2}}\left|\delta^{h}_k
		u\right|^2dx+\int_{\Omega}\eta^2\left|\delta^{h}_k\nabla u
		\right|^2dx\right).
	\end{aligned}
\end{equation}
Now we apply Theorem \ref{Sob:teo1.9}--(i) and inequality, \eqref{Cacciopp-9}, so that we have

\begin{equation*}
	\begin{aligned}
		R^{-2}\int_{B_{3R/2}}\left|\delta^{h}_k u\right|^2dx\leq
		R^{-2}\int_{B_{7R/4}}\left|\nabla u\right|^2dx\leq
		CR^{-4}\int_{B_{2R}}u^2dx+C\int_{B_{2R}}f^2dx.
	\end{aligned}
\end{equation*}
By the last obtained estimate and by \eqref{Ell:E9.1} we have

\begin{equation*}
	\begin{aligned}
		\left|\int_{\Omega}f\delta^{-h}_k\left(\eta^2
		\delta^{h}_ku\right)dx\right|&\leq
		\frac{1}{2\sigma}\int_{B_{2R}}f^2dx+\\&+
		C_{\ast}\sigma\left(R^{-4}\int_{B_{2R}}
		u^2dx+\int_{\Omega}\eta^2\left|\delta^{h}_k\nabla u
		\right|^2dx\right),
	\end{aligned}
\end{equation*}
where $C_{\ast}$ is a constant depending on $\lambda$ and $n$ only. Inserting what we have just obtained into \eqref{Ell:E7.2}
we get 
\begin{equation*}
	\begin{aligned}
		\left(\frac{\lambda^{-1}}{2}-C_{\ast}\sigma\right)\int_{\Omega}\left|\delta^{h}_k\nabla
		u\right|^2\eta^2dx&\leq \frac{1}{2\sigma}\int_{B_{2R}}f^2dx+\\&+
		\frac{C\left(1+\sigma+E^2R^2\right)}{R^4}\int_{B_{2R}} u^2dx.
	\end{aligned}
\end{equation*}
Now, choosing

$$\sigma=\frac{\lambda^{-1}}{4C_{\ast}}$$ and we have, for $k=1,\cdots,n$,

\begin{equation*}
	\begin{aligned}
		\int_{B_{R}}\left|\delta^{h}_k\nabla u\right|^2\eta^2dx\leq
		C\int_{B_{2R}}f^2dx +\frac{C\left(1+E^2R^2\right)}{R^4}\int_{B_{2R}}
		u^2dx.
	\end{aligned}
\end{equation*}
By the last obtained inequality and by Theorem \ref{Sob:teo1.9}--(ii) we obtain
$$\partial_k u\in H^1\left(B_R\right),$$ for $k=1,\cdot,n$. Hence
$u\in H^2\left(B_R\right)$ and

\begin{equation*}
	\sum_{|\alpha|=2 }R^{2|\alpha|}\int_{B_R}\left|\partial^{\alpha}u
	\right|^2dx\leq C\left(1+E^2R^2\right)\int_{B_{2R}}u
	^2dx+CR^4\int_{B_{2R}}f ^2dx.
\end{equation*}
Finally, by the latter and by \eqref{Cacciopp-9} we get
\eqref{Ell:E2.1}. $\blacksquare$

\bigskip

\underline{\textbf{Exercise 1.}} Under the same assumptions of Theorem \ref{reg-interno}, prove that, if $0<r<\rho$ and $B_{\rho}(x_0)\subset \Omega$ then
\begin{equation}\label{Ell:reg-0}
	\begin{aligned}
		\sum_{|\alpha|\leq
			2}(\rho-r)^{2|\alpha|}\int_{B_r(x_0)}\left|\partial^{\alpha}u
		\right|^2dx&\leq C\left(1+E^2\rho^2\right)\int_{B_{\rho}(x_0)}u
		^2dx+\\&+C(\rho-r)^4\int_{B_{\rho}(x_0)}f ^2dx,
	\end{aligned}
\end{equation}
where $C$ depends on $\lambda$ only. [Hint: consider a finite covering $\overline{B_r(x_0)}$ consisting
of balls of the type $B_{\frac{\rho-r}{2}}(x)$, $x\in B_r(x_0)$, and apply \eqref{Ell:E2.1}].

\bigskip

\underline{\textbf{Exercise 2.}} Under the same assumption of Theorem \ref{reg-interno}, prove that, if $$\Omega'\Subset \Omega$$ then
\begin{equation}\label{Ell:reg-00}
	\begin{aligned}
		\sum_{|\alpha|\leq
			2}\delta_0^{2|\alpha|}\int_{\Omega'}\left|\partial^{\alpha}u
		\right|^2dx&\leq C\left(1+E^2d_0^2\right)\int_{\Omega}u
		^2dx+\\&+Cd_0^4\int_{\Omega}f ^2dx,
	\end{aligned}
\end{equation}
where $d_0$ is the diameter of $\Omega$, $\delta_0=$dist$\left(\Omega', \partial \Omega\right)$ and $C$ depends on $\lambda$ and $d_0\delta_0^{-1}$ only.
[Hint: use Exercise 1 and a partition of unity].

\bigskip

\underline{\textbf{Exercise 3.}} (a) Generalize Theorem \ref{reg-interno}
to the equation
\begin{equation*}
	-\sum_{j,k=1}^n\partial_j\left( a^{jk}\partial_k u
	+d^ju\right)+\sum_{j=1}^nb^{j}\partial_j u+cu =f
\end{equation*}
where $A=\left\{a^{jk}\right\}_{j,k=1}^n$ and $f$ satisfy the same assumptions of Theorem \ref{reg-interno},
$d, b\in L^{\infty}(\Omega;\mathbb{R}^n)$, $c\in L^{\infty}(\Omega)$.

(b) Generalize (a) to the case where $A$ is a nonsymmetric matrix. [Hint to (b): write the operator div $(A\nabla u)$ like $$ \mbox{div}(A^s\nabla u)+\mbox{terms of order less than } 2,$$
where $A^s$ is symmetric part of $A$]. $\clubsuit$

\bigskip

\begin{theo} [\textbf{improved regularity in the interior}]\label{Ell:reg-m}
	\index{Theorem:@{Theorem:}!- improved regularity in the interior@{- improved regularity in the interior}}
Let $\Omega$ be a bounded open set of
	$\mathbb{R}^n$ with diameter $d_0$. Let
	\begin{equation}\label{Ell:reg-m1}
		f\in H^m(\Omega).
	\end{equation}
	Let  $A$ be a symmetric matrix.
	Let us assume that $A$ satisfies \eqref{gamma-eq}, let us assume that $A\in
	C^{m,1}(\overline{\Omega};\mathbb{M}(n))$ and it satisfies
	
	\begin{equation}\label{Ell:reg-m2}
		\left\Vert
		A\right\Vert_{C^{m,1}\left(\overline{\Omega};\mathbb{M}(n)\right)}\leq
		E_m,
	\end{equation}
	where, $E_m$ is a positive number and, recall,
	$$\left\Vert
	A\right\Vert_{C^{m,1}\left(\overline{\Omega};\mathbb{M}(n)\right))}=\sum_{|\alpha|\leq
		m}d_0^{|\alpha|}\left\Vert
	\partial^{\alpha}A\right\Vert_{L^{\infty}\left(\Omega;\mathbb{M}(n)\right))}+d_0^{m+1}\sum_{|\alpha|=
		m}[\partial^{\alpha}A]_{1, \Omega}.$$ 
	Let us assume that $u\in H^1(\Omega)$ is a solution to
	
	\begin{equation}\label{Ell:reg-m3}
		-\mbox{div}(A\nabla u)=f, \quad\mbox{ in } \Omega.
	\end{equation}
	Then we have
	\begin{equation}\label{Ell:reg-m4}
		u\in H^{m+2}_{loc}(\Omega),
	\end{equation}
	moreover, if $B_{2R}(x_0)\subset \Omega$, then the following inequality holds true
	\begin{equation}\label{Ell:reg-m5}
		\begin{aligned}
			\sum_{|\alpha|\leq
				m+2}R^{2|\alpha|}\int_{B_R(x_0)}\left|\partial^{\alpha}u
			\right|^2dx&\leq C\left(1+E_m^2\right)\int_{B_{2R}(x_0)}u
			^2dx+\\&+C\sum_{|\alpha|\leq
				m}R^{2(|\alpha|+4)}\int_{B_{2R}(x_0)}\left|\partial^{\alpha}f
			\right|^2dx,
		\end{aligned}
	\end{equation}
	where $C$ depends on $\lambda$ only.
\end{theo}
\textbf{Proof.} We simply consider the case $m=1$,
leaving the reader to complete the proof by induction.
Let $l\in\{1,\cdots,n\}$. Let $\widetilde{v}$ be any function belonging to $C_0^{\infty}\left(B_{3R/2}(x_0)\right)$. Choose, as a
test function,
$$v=-\partial_l \widetilde{v}.$$ We get

\begin{equation}\label{Ell:reg-m6}
	-\int_{\Omega} A\nabla u\cdot \nabla\partial_l \widetilde{v}
	dx=-\int_{\Omega}f\partial_l \widetilde{v}dx.
\end{equation}
Since $f\in H^1(\Omega)$, we have
\begin{equation}\label{Ell:reg-m7}
	-\int_{\Omega}f\partial_l \widetilde{v}dx=\int_{\Omega}\partial_lf
	\widetilde{v}dx.
\end{equation}
Moreover, since $u\in
H_{loc}^2(\Omega)$ (by Theorem \ref{reg-interno}), we have

\begin{equation*}
	\begin{aligned}
		-\int_{\Omega} A\nabla u\cdot \nabla\partial_l
		\widetilde{v}dx&=-\int_{\Omega}\sum_{j,k=1}^na^{jk}\partial_ku\partial_l\left(\partial_j\widetilde{v}\right)dx=\\&=
		\int_{\Omega}\sum_{j,k=1}^n\partial_l\left(a^{jk}\partial_ku\right)\partial_j\widetilde{v}dx=\\&=
		\int_{\Omega}\sum_{j,k=1}^na^{jk}\partial_k\left(\partial_lu\right)\partial_j\widetilde{v}dx+\\&+
		\int_{\Omega}\sum_{j,k=1}^n\partial_la^{jk}\partial_ku\partial_j\widetilde{v}dx=\\&=
		\int_{\Omega}\sum_{j,k=1}^na^{jk}\partial_k\left(\partial_lu\right)\partial_j\widetilde{v}dx-\\&-
		\int_{\Omega}\sum_{j,k=1}^n\partial_j\left[\left(\partial_la^{jk}\right)\partial_ku\right]\widetilde{v}dx.
	\end{aligned}
\end{equation*}
By the equality obtained above, by \eqref{Ell:reg-m6} and by \eqref{Ell:reg-m7} we have
(recall that  $C_0^{\infty}\left(B_{3R/2}(x_0)\right)$ is dense in $H_0^{1}\left(B_{3R/2}(x_0)\right)$

\begin{equation}\label{Ell:reg-m8}
	\int_{B_{3R/2}(x_0)}A\nabla\left(\partial_lu\right)\cdot\nabla wdx=\int_{B_R(x_0)}
	\widetilde{f}wdx,\ \ \forall w\in H_0^{1}\left(B_{3R/2}(x_0)\right),
\end{equation}
where
\begin{equation*}
	\widetilde{f}=\partial_lf+\sum_{j,k=1}^n\partial_j\left[\left(\partial_la^{jk}\right)\partial_ku\right].
\end{equation*}
Now $\widetilde{f}\in L^{2}(B_{3R/2}(x_0))$. As a matter of fact

\begin{equation}\label{Ell:reg-m9}
	\begin{aligned}
		&\int_{B_{3R/2}(x_0)}\left|\widetilde{f}\right|^2dx\leq\\&\leq 2\int_{B_{3R/2}(x_0)}\left|\partial_lf\right|^2dx+cE^2_1d_0^{-4}\int_{B_{3R/2}(x_0)}\left|\partial^2u\right|^2dx+
		\\&+cE^2_1d_0^{-2}\int_{B_{3R/2}(x_0)}\left|\nabla u\right|^2dx.
	\end{aligned}
\end{equation}

\medskip

By the latter, by Theorem \ref{reg-interno} (more precisely, by \eqref{Ell:reg-0}) and by \eqref{Ell:reg-m8}, we obtain

\begin{equation}
	\begin{aligned}
		R^4\sum_{|\alpha|=
			3}\int_{B_R(x_0)}\left|\partial^{\alpha}u
		\right|^2dx&\leq C\left(1+E_1^2\right)\sum_{|\alpha|=
			2}\int_{B_{3R/2}(x_0)}\left|\partial^{\alpha}u
		\right|^2dx+\\&+CR^4\int_{B_{3R/2}(x_0)}\left|\widetilde{f}\right|^2dx.
	\end{aligned}
\end{equation}
At this point, we again apply \eqref{Ell:reg-0} to estimate from above the derivatives of order less than or equal to $2$ we obtain \eqref{Ell:reg-m5}. $\blacksquare$

\bigskip

\underline{\textbf{Exercise 4.}} Generalize Theorem \ref{Ell:reg-m}
to the equation
\begin{equation*}
	-\sum_{j,k=1}^n\partial_j\left( a^{jk}\partial_k u
	+d^ju\right)+\sum_{j,k=1}^n b^{j}\partial_j u+cu =f
\end{equation*}
where $A$ and $f$ satisfy the same assumption of Theorem \ref{Ell:reg-m},
$d, b\in C^{m-1,1}\left(\overline{\Omega};\mathbb{R}^n\right)$, $c\in C^{m-1,1}\left(\overline{\Omega}\right)$, $m\geq 1$. $\clubsuit$

\bigskip

\begin{cor} [\textbf{$C^{\infty}$ regularity in the interior}]\label{Ell:reg-infin}
	\index{Corollary:@{Corollary:}!- $C^{\infty}$ regularity in the interior@{- $C^{\infty}$ regularity in the interior}}
	Let $\Omega$ an open set of
	$\mathbb{R}^n$. Let
	\begin{equation*}
		f\in C^{\infty}(\Omega).
	\end{equation*}
	Let $A\in
	C^{\infty}(\Omega;\mathbb{M}(n))$. Let us assume that $A$ satisfies \eqref{gamma-eq}.
	
	Let $u\in H^1(\Omega)$ be a solution to 
	
	\begin{equation*}
		-\mbox{div}(A\nabla u)=f, \quad\mbox{ in } \Omega.
	\end{equation*}
	Then we have
	\begin{equation*}
		u\in C^{\infty}(\Omega).
	\end{equation*}
\end{cor}
\textbf{Proof.} Let $B_R(x_0)\Subset \Omega$. For any  $m\geq 0$, we have $f\in H^m(B_R(x_0))$ and $A\in
C^{m,1}(\overline{B_R(x_0)};\mathbb{M}(n))$. Therefore Theorem \ref{Ell:reg-m} implies

$$u\in \bigcap_{m=0}^{\infty}H^m(B_R(x_0))=C^{\infty}(B_R(x_0)),$$
Where the last equality is due to Theorem \ref{Sobolev-inequ}. Since $B_R(x_0)$ is arbitrary, the thesis follows. $\blacksquare$

\bigskip

\underline{\textbf{Exercise 5.}} Prove Corollary \ref{Ell:reg-infin}
for the equation
\begin{equation*}
	-\sum_{j,k=1}^n\partial_j\left( a^{jk}\partial_k u
	+d^ju\right)+\sum_{j=1}^nb^{j}\partial_j u+cu =f
\end{equation*}
where $A$ satisfies the same assumption of Corollary \ref{Ell:reg-infin} and
$d, b\in C^{\infty}\left(\Omega;\mathbb{R}^n\right)$, $c\in C^{\infty}\left(\Omega\right)$. $\clubsuit$

\subsection{Regularity teorems at the boundary--global regularity} \label{regolarit-bordo}

In this Section we will study the regularity at the boundary. The following Lemma is a crucial step in the proof of the forthcoming theorems

\begin{lem}[\textbf{local regularity at the boundary}]\label{Ell:reg-bordo-loc}
	\index{Lemma:@{Lemma:}!- local regularity at the boundary@{- local regularity at the boundary}}
	Let $R>0$ and
	\begin{equation}\label{Ell:E12.1L}
		f\in L^2\left(B^+_{2R}\right)
	\end{equation}
	and let $A$ be a symmetric matrix. Let us assume that $A$ satisfies \eqref{gamma-eq}, $A\in
	C^{0,1}(\Omega;\mathbb{M}(n))$ and it satisfies
	
	\begin{equation}\label{Ell:E12.2L}
		\left|A(x)-A(\overline{x})\right|\leq E|x-\overline{x}|,\quad \forall x,\overline{x}\in
		B^+_{2R},
	\end{equation}
	where $E$ is a positive number. Let us assume that $u\in H^1\left(B^+_{2R}\right)$ satisfies
	
	\begin{equation}\label{Ell:E12.3L}
		-\mbox{div}(A\nabla u)=f, \quad\mbox{ in } B^+_{2R}
	\end{equation}
	and 
	\begin{equation}\label{Ell:E12.4L}
		u(\cdot,0)=0, \quad\mbox{ in the sense of the traces in }  B'_{2R}.
	\end{equation}
	Then we have
	\begin{equation}\label{Ell:E12.5L}
		u\in H^2\left(B^+_{R}\right)
	\end{equation}
	and the following estimate holds true
	\begin{equation}\label{Ell:E12.6L}
		\begin{aligned}
			\sum_{|\alpha|\leq
				2}R^{2|\alpha|}\int_{B^+_R}\left|\partial^{\alpha}u
			\right|^2dx&\leq C\left(1+E^2R^2\right)\int_{B^+_{2R}}u
			^2dx+\\&+CR^4\int_{B^+_{2R}}f ^2dx,
		\end{aligned}
	\end{equation}
	where $C$ depends on $\lambda$ only.
\end{lem}
\textbf{Proof.} Let $\eta\in
C^{\infty}_0(B_{3R/2})$ satisfy
\begin{equation}\label{Ell:E2.12-0}
	0\leq \eta\leq 1;\qquad \eta=1,\quad \mbox{ in } B_{R}
\end{equation}
and
\begin{equation}\label{Ell:E2.12-00}
	|\nabla \eta|\leq \frac{K}{R},
\end{equation}
where $K$ is a positive constant. Let $h\in \left(-\frac{R}{8},\frac{R}{8}\right)\setminus\{0\}$ and let $$k\in\{\ 1,\cdots,n-1\}.$$ Let us denote

\begin{equation}\label{Ell:E2.12-000}
	v=-\delta^{-h}_k\left(\eta^2 \delta^{h}_ku\right).
\end{equation}
By \eqref{Ell:E12.4L}, taking into account that $\eta\in
C^{\infty}_0\left(B_{3R/2}\right)$, we have 
$$v\in H^{1}_0\left(B^+_{2R}\right).$$
Therefore, by $u\in H^{1}\left(B^+_{2R}\right)$, by \eqref{Ell:E12.1L}, \eqref{Ell:E12.3L} and by \eqref{Ell:E12.4L}, we get

\begin{equation}\label{Ell:E13.1}
	\begin{aligned}
		\int_{B^+_{2R}} A\nabla u\cdot \nabla v  dx=\int_{B^+_{2R}} f v  dx.
	\end{aligned}
\end{equation}
At this point we may argue likewise in the proof of Theorem \ref{reg-interno}. Actually, by using  Theorem \ref{Sob:teo1.9bis} instead of \ref{Sob:teo1.9} we have
\begin{equation}\label{Ell:E14.1}
	\begin{aligned}
		\frac{\lambda^{-1}}{2}\int_{B^+_{2R}}\left|\delta^{h}_k\nabla
		u\right|^2&\eta^2dx\leq - \int_{B^+_{2R}}f\delta^{-h}_k\left(\eta^2
		\delta^{h}_ku\right)dx+\\&+C\left(R^{-2}+E^2\right)\int_{B^+_{3R/2}}\left(\left|\delta^{h}_k
		u\right|^2+\left|\nabla u\right|^2\right)dx,
	\end{aligned}
\end{equation}
for any $k\in\{1,\cdots,n-1\}$, where $C$ depends on $\lambda$ and $n$ only.
Let $\sigma>0$ to be choosen,  we get (compare with \eqref{Ell:E9.1})

\begin{equation}\label{Ell:E15.1}
	\begin{aligned}
		\left|\int_{B^+_{2R}}f\delta^{-h}_k\left(\eta^2
		\delta^{h}_ku\right)dx\right|&\leq
		\frac{1}{2\sigma}\int_{B^+_{2R}}f^2dx+\\&+
		C\sigma\left(R^{-2}\int_{B^+_{3R/2}}\left|\delta^{h}_k
		u\right|^2dx+\right. \\&+
		\left.\int_{B^+_{2R}}\eta^2\left|\nabla \delta^{h}_ku
		\right|^2dx\right).
	\end{aligned}
\end{equation}
By applying Theorem \ref{Sob:teo1.9bis} and inequality
\eqref{Cacciopp-10}, we get

\begin{equation*}
	\begin{aligned}
		R^{-2}\int_{B^+_{3R/2}}\left|\delta^{h}_k u\right|^2dx&\leq
		R^{-2}\int_{B^+_{7R/4}}\left|\nabla u\right|^2dx\leq\\&\leq
		CR^{-4}\int_{B^+_{2R}}u^2dx+C\int_{B^+_{2R}}f^2dx.
	\end{aligned}
\end{equation*}
By the latter and by \eqref{Ell:E15.1} we have

\begin{equation*}
	\begin{aligned}
		\left|\int_{\Omega}f\delta^{-h}_k\left(\eta^2
		\delta^{h}_ku\right)dx\right|&\leq
		\frac{1}{2\sigma}\int_{B_{2R}}f^2dx+\\&+
		C_{\star}\sigma\left(R^{-4}\int_{B_{2R}}
		u^2dx+\int_{\Omega}\eta^2\left|\nabla \delta^{h}_ku,
		\right|^2dx\right),
	\end{aligned}
\end{equation*}
where $C_{\star}$  depends on  $\lambda$ and
and $n$ only. By using in \eqref{Ell:E14.1} the just obtained inequality, 
we obtain
\begin{equation*}
	\begin{aligned}
		\left(\frac{\lambda^{-1}}{2}-C_{\star}\sigma\right)\int_{B^+_{2R}}\left|\delta^{h}_k\nabla
		u\right|^2\eta^2dx&\leq \frac{1}{2\sigma}\int_{B^+_{2R}}f^2dx+\\&+
		\frac{C\left(1+\sigma+E^2R^2\right)}{R^4}\int_{B^+_{2R}} u^2dx.
	\end{aligned}
\end{equation*}
Now let us choose

$$\sigma=\frac{\lambda^{-1}}{4C_{\ast}}$$ and we get, for $k=1,\cdots,n-1$,

\begin{equation*}
	\begin{aligned}
		\int_{B^+_{R}}\left|\delta^{h}_k\nabla u\right|^2\eta^2dx\leq
		C\int_{B^+_{2R}}f^2dx +\frac{C\left(1+E^2R^2\right)}{R^4}\int_{B^+_{2R}}
		u^2dx.
	\end{aligned}
\end{equation*}
By this inequality and by  Theorem \ref{Sob:teo1.9bis}--(ii) we get
$$\partial_k u\in H^1\left(B^+_R\right),\quad \mbox{for } k=1,\cdots,n-1$$ and applying again  \eqref{Cacciopp-10} we have

\begin{equation}\label{Ell:E16.1}
	\sum_{\underset{k+j<2n}{k,j=1}}^n\int_{B^+_R}\left|\partial^{2}_{jk}u
	\right|^2dx\leq \frac{C\left(1+E^2R^2\right)}{R^4}\int_{B^+_{2R}}u
	^2dx+CR^4\int_{B^+_{2R}}f ^2dx,
\end{equation}
where $C$ depends on $\lambda$  $n$ only.

Now, Theorem \ref{reg-interno} implies that $u\in H_{loc}^2\left(B^+_R\right)$, this allows us to write the equation \eqref{Ell:E12.2L} in the form

$$\sum_{k,j=1}^na^{jk}\partial^2_{jk}u+\sum_{k,j=1}^n\partial_ja^{jk}\partial_{k}u=-f, \quad \mbox{a.e. in } B^+_R.$$
By the last equality we can find (taking into account that $a^{nn}\geq\lambda^{-1}>0$)
\begin{equation}\label{Ell:E16.1bis}
	\partial^2_nu=-\frac{1}{a^{nn}}\left(\sum_{\underset{k+j<2n}{k,j=1}}^n a^{jk}\partial^2_{jk}u+\sum_{k,j=1}^n\partial_ja^{jk}\partial_{k}u+f\right).\end{equation}
By \eqref{Ell:E16.1bis}, \eqref{Ell:E16.1}, \eqref{Ell:E12.2L} we have $\partial^2_nu\in H^2\left(B^+_R\right)$. Finally, by \eqref{Cacciopp-10} we obtain
\eqref{Ell:E12.6L}. $\blacksquare$

\begin{theo} [\textbf{global regularity}]\label{reg-bordo}
	\index{Theorem:@{Theorem:}!- global regularity@{- global regularity}}
	Let $\Omega$ be a bounded open set of $\mathbb{R}^n$ whose boundary is of class $C^{1,1}$ with constants $M_0,r_0$. Let
	\begin{equation}\label{Ell:E11.2}
		f\in L^2(\Omega)
	\end{equation}
	and let $A$ be a symmetric matrix. Let us assume that $A$ satisfies \eqref{gamma-eq},  $A\in
	C^{0,1}(\Omega;\mathbb{M}(n))$ and it satisfies
	
	\begin{equation}\label{Ell:E11.3}
		\left|A(x)-A(\overline{x})\right|\leq E|x-\overline{x}|,\quad \forall x,\overline{x}\in
		\Omega,
	\end{equation}
	where $E$ is a positive number. Let us suppose that $u\in H_0^1(\Omega)$ is the solution
	to
	
	\begin{equation}\label{Ell:E11.3bis}
		-\mbox{div}(A\nabla u)=f, \quad\mbox{ in } \Omega.
	\end{equation}
	Then we have
	\begin{equation}\label{Ell:E11.3ter}
		u\in H^2(\Omega),
	\end{equation}
	and the following estimate holds true
	\begin{equation}\label{Ell:E11.4}
		\begin{aligned}
			\sum_{|\alpha|\leq
				2}r_0^{2|\alpha|}\int_{\Omega}\left|\partial^{\alpha}u
			\right|^2dx&\leq Cd_0^4\int_{\Omega}f ^2dx,
		\end{aligned}
	\end{equation}
	where $d_0$ is the diameter of $\Omega$ and $C$ depends by $\lambda, E, M_0$ and $\frac{d_0}{r_0}$ only.
\end{theo}

\textbf{Proof.} Let $P\in
\partial \Omega$. There exists a rigid transformation of coordinates under which we have
 $P=0$ and
\begin{equation*}
	\Omega\cap Q_{r_0,2M_0}=\left\{x\in Q_{r_0,2M_0}:\mbox{ }
	x_n>g(x') \right\},
\end{equation*}
where $g\in C^{1,1}\left(\overline{B'_{r_0}}\right)$,
$$g(0)=0, \quad |\nabla_{x'} g(0)|=0$$
and
$$\left\Vert g\right\Vert_{C^{1,1}\left(\overline{B'_{r_0}}\right)}\leq
M_0r_0.$$ 
Let us consider the change of coordinates
$$\Phi: Q_{r_0,2M_0}\rightarrow\mathbb{R}^n,\quad \Phi(x)=\left(x',x_n-g(x')\right).$$
Let us note that $\Phi$ "flattens the boundary" i.e.

$$Q_{r_0,2M_0}\cap\Phi\left(\partial \Omega\right)  =\left\{(y',0): \ y'\in B'_{r_0}\right\}.$$
$\Phi\in C^{1,1}\left(Q_{r_0,2M_0}\cap \Omega\right)$ is injiective and it is a local diffeomorphism. Let
$$W=Q_{r_0,2M_0}\cap\Phi\left(\Omega\right)$$
and

$$\Psi: W\rightarrow Q_{r_0,2M_0}\cap\Omega,\quad \Psi=\Phi^{-1}.$$
Let $J$ be the jacobian matrix  of $\Psi$. Then

$$\det J(y)=1,\quad \forall y\in W.$$

Let $u\in H_0^1(\Omega)$ be the weak solutions to \eqref{Ell:E11.3bis}. Let $\widetilde{v}$ be any function belonging to $H_0^1(W)$ and set $$v(x)=\widetilde{v}\left(\Psi(x)\right).$$ We have

\begin{equation}\label{Ell:E18.0}
	\int_{Q_{r_0,2M_0}\cap\Omega} A\nabla u\cdot \nabla v  dx=\int_{Q_{r_0,2M_0}\cap\Omega}fvdx.
\end{equation}
By the change of variables $x=\Psi(y)$ equation \eqref{Ell:E18.0} becomes

\begin{equation}
	\int_{W} \widetilde{A}(y)\nabla w(y)\cdot \nabla \widetilde{v}(y)  dy=\int_{W}\widetilde{f}(y)\widetilde{v}(y)dy,
\end{equation}
where

$$w(y)=u(\Psi(y)), \quad \forall y\in W,$$

$$\widetilde{A}(y)=(J(y))^{-1}A(\Psi(y))\left((J(y))^{-1}\right)^t,\quad \forall y\in W$$
and
$$\widetilde{f}(y)=f(\Psi(y)),\quad \forall y\in W.$$
It is easy to check that

$$\widetilde{\lambda}^{-1}|\xi|^2\leq \widetilde{A}(y)\xi\cdot\xi\leq \widetilde{\lambda}|\xi|^2,\quad \forall y\in W,\ \forall \xi\in \mathbb{R}^n,$$
where $\widetilde{\lambda}\geq 1$ depends on $\lambda$ and $M_0$ only. Moreover

$$w\left(y',0\right)=0,\quad y'\in B'_{r_0}\mbox{ in the sense of the traces. }$$

Also we have $\widetilde{A}\in C^{0,1}\left(W\right)$ and

$$\left|\widetilde{A}(y)-\widetilde{A}\left(\overline{y}\right)\right|\leq \widetilde{E}\left|y-\overline{y}\right|, \quad \forall y,\overline{y}\in W,$$
where
$$\widetilde{E}=C\left(E+M_0r_0^{-1}\right),$$
 and $C$ depends on $\lambda$ and $M_0$ only.

At this stage (compare Fig. 4.1) we introduce the quantity

\begin{figure}\label{figura-Ell-pagE19}
	\centering
	\includegraphics[trim={0 0 0 0},clip, width=15cm]{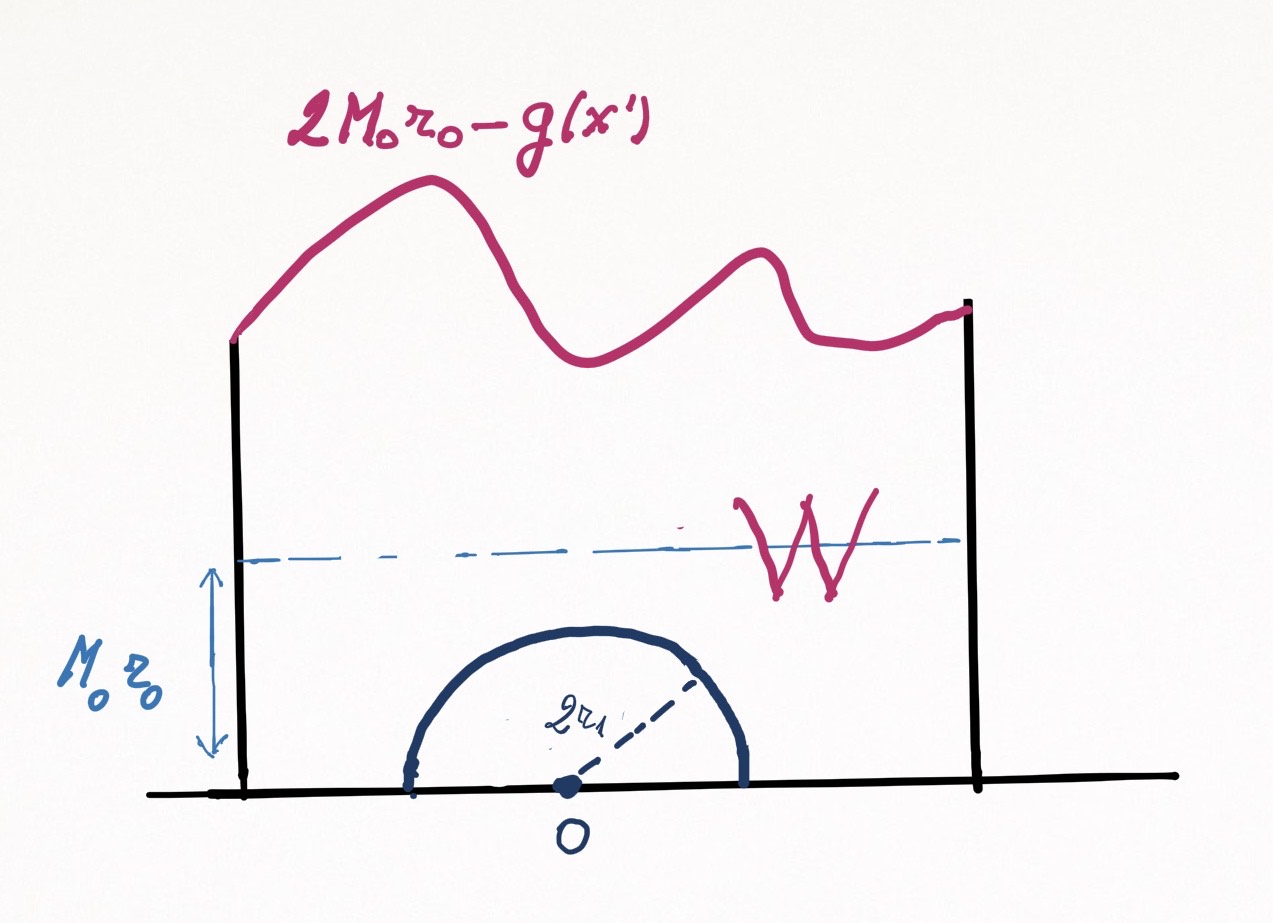}
	\caption{}
\end{figure}

$$r_1=\frac{r_0}{2}\min\{1, M_0\}$$ in such a way that we have
$$B^+_{2r_1}\subset W.$$ Applying Lemma \ref{Ell:reg-bordo-loc} we get

\begin{equation}\label{Ell:E20.1}
	\begin{aligned}
		\sum_{|\alpha|\leq
			2}r_1^{2|\alpha|}\int_{B^+_{r_1}}\left|\partial^{\alpha}w
		\right|^2dy&\leq C\left(1+\widetilde{E}^2r_1^2\right)\int_{B_{2r_1}}w
		^2dy+\\&+Cr_1^4\int_{B^+_{2r_1}}\widetilde{f} ^2dy,
	\end{aligned}
\end{equation}
where $C$ depends on $\widetilde{\lambda}$ only. Coming back to the original variables, after some calculation and simple estimates, we have

\begin{equation}\label{Ell:E20.1-0}
	\begin{aligned}
		\sum_{|\alpha|\leq
			2}r_0^{2|\alpha|}\int_{\Psi\left(B^+_{r_1}\right)}\left|\partial^{\alpha}u
		\right|^2dx&\leq C\left(\int_{\Omega}u
		^2dx+r_0^4\int_{\Omega}f ^2dx\right),
	\end{aligned}
\end{equation}
where $C$ depends on $\lambda, M_0$ and $Er_0$ only. On the other hand, as is easily checked, there is $\overline{C}\geq1$, $C$ depending on $\lambda$ and $M_0$ only, such that, if
$r_2=\frac{r_0}{C}$, we have

$$\Omega\cap B_{r_2}(P)\subset \Psi\left(B^+_{r_1}\right).$$ Therefore by \eqref{Ell:E20.1-0} we get trivially

\begin{equation}\label{Ell:E20.2}
	\begin{aligned}
		\sum_{|\alpha|\leq
			2}r_0^{2|\alpha|}\int_{\Omega\cap B_{r_2}(P)}\left|\partial^{\alpha}u
		\right|^2dx&\leq C\left(\int_{\Omega}u
		^2dx+r_0^4\int_{\Omega}f ^2dx\right).
	\end{aligned}
\end{equation}
Now, by the compactness of  $\partial\Omega$, we can extract a finite covering by $\left\{B_{r_2}(P)\right\}_{P\in \partial\Omega}$.  Let
$\left\{B_{r_2}(P_j)\right\}_{1\leq j\leq N}$,  such a finite covering, where $P_j\in \partial \Omega$, and let
$$\Lambda=\Omega\cap\bigcup_{j=1}^N B_{r_2}(P_j)\quad \mbox{e}\quad \Omega'=\Omega \setminus \Lambda.$$
We can make dist($\Omega',\partial \Omega)\geq r_0/C$, where $C\geq 1$ depends on $M_0$ only, furthermore $N$ depends on $M_0$ and $\frac{d_0}{r_0}$ only.
Inequality \eqref{Ell:E20.2} implies that there is a constant $C$ depending on $\lambda, E, M_0$ and $\frac{d_0}{r_0}$ so that

\begin{equation*}
	\begin{aligned}
		\sum_{|\alpha|\leq
			2}r_0^{2|\alpha|}\int_{\Lambda}\left|\partial^{\alpha}u
		\right|^2dx&\leq C\left(\int_{\Omega}u
		^2dx+r_0^4\int_{\Omega}f ^2dx\right).
	\end{aligned}
\end{equation*}
By the just obtained inequality and by \eqref{Ell:reg-00} we have

\begin{equation}\label{Ell:E20.3}
	\begin{aligned}
		\sum_{|\alpha|\leq
			2}r_0^{2|\alpha|}\int_{\Omega}\left|\partial^{\alpha}u
		\right|^2dx&\leq C\left(\int_{\Omega}u
		^2dx+d_0^4\int_{\Omega}f ^2dx\right),
	\end{aligned}
\end{equation}
By the first Poincar\'{e} inequality (Proposition \eqref{Poincar}) and by inequality \eqref{stab}, we find

\begin{equation*}
	\int_{\Omega}u
	^2dx\leq Cd_0^4\int_{\Omega}f ^2dx,
\end{equation*}
where $C$ depends on $\lambda$ only. By the last obtained inequality and by \eqref{Ell:E20.3} we get \eqref{Ell:E11.3ter}. $\blacksquare$

\bigskip

\underline{\textbf{Exercise 1.}} (a) Generalize Theorem \ref{reg-bordo} to the boundary value problem
\begin{equation*}
	\begin{cases}
		-\sum_{j,k=1}^n\partial_j\left( a^{jk}\partial_k u
		+d^ju\right)+\sum_{j=1}^nb^{j}\partial_j u+cu =f, \quad\mbox{ in } \Omega\\
		\\
		u=0, \quad\mbox{ su } \partial\Omega.
	\end{cases}
\end{equation*}
where $A=\left\{a^{jk}\right\}_{j,k=1}^n$ and $f$ satisfy the same assumptions of  Theorem \ref{reg-interno}
and $d, b\in L^{\infty}(\Omega;\mathbb{R}^n)$, $c\in L^{\infty}(\Omega)$.

(b) Generalize the result obtained in (a) to the case where $A$ is not a symmetric matrix.  $\clubsuit$

\bigskip

\begin{theo} [\textbf{improved global regularity}]\label{Ell:reg-bordo-m}
	\index{Theorem:@{Theorem:}!- improved global regularity@{- improved global regularity}}
	Let $\Omega$ be a bounded open set of	$\mathbb{R}^n$ of class $C^{m+1,1}$ with constants $M_0,r_0$, let $d_0$ be the diameter of $\Omega$. Let
	\begin{equation}\label{Ell:reg-bordo-m1}
		f\in H^m(\Omega)
	\end{equation}
	and let $A$ be a a symmetric matrix. Let us assume that $A$ satisfies \eqref{gamma-eq},  $A\in
	C^{m,1}(\overline{\Omega};\mathbb{M}(n))$ and it satisfies
	
	\begin{equation*}
		\left\Vert
		A\right\Vert_{C^{m,1}\left(\overline{\Omega};\mathbb{M}(n)\right)}\leq
		E_m,
	\end{equation*}
	where $E_m$ is a positive number.	
	
	Let us suppose that $u\in H_0^1(\Omega)$ is the solution to
	
	\begin{equation}\label{Ell:reg-bordo-m2}
		-\mbox{div}(A\nabla u)=f, \quad\mbox{ in } \Omega.
	\end{equation}
	Then we have
	\begin{equation}\label{Ell:reg-bordo-m3}
		u\in H^{m+2}(\Omega),
	\end{equation}
	and the following estimate holds true
	\begin{equation}\label{Ell:reg-bordo-m4}
		\begin{aligned}
			\sum_{|\alpha|\leq
				m+2}r_0^{2|\alpha|}\int_{\Omega}\left|\partial^{\alpha}u
			\right|^2dx&\leq C\sum_{|\alpha|\leq
				m}d_0^{2(|\alpha|+4)}\int_{\Omega}\left|\partial^{\alpha}f
			\right|^2dx,
		\end{aligned}
	\end{equation}
	where $C$ depends on $\lambda, E_m, M_0$ and $\frac{d_0}{r_0}$ only.
\end{theo}
\textbf{Proof.} The proof is mostly similar to that of Theorem \ref{Ell:reg-m}, so we focus here by considering, in the case $m=1$, the steps in which the new proof differs from the proof of Theorem \ref{Ell:reg-m}, leaving the details to the care of the reader. First of all, in analogy to Lemma \ref{Ell:reg-bordo-loc}, let us consider
the following special situation. Let $R>0$ and
\begin{equation*}
	f\in H^1\left(B^+_{2R}\right).
\end{equation*}
Let $A$ be a symmetric matrix whose entries are measurable functions in $B^+_{2R}$.
Let us suppose that $A$ satisfies \eqref{gamma-eq} and that $A\in
C^{1,1}\left(B_{2R},\mathbb{M}(n)\right)$. Let us suppose that $u\in H^1\left(B^+_{2R}\right)$
is a solution to

\begin{equation*}
	-\mbox{div}(A\nabla u)=f, \quad\mbox{ in } B^+_{2R},
\end{equation*}
and
\begin{equation}\label{Ell:E23.1}
	u(\cdot,0)=0, \quad\mbox{ in } B'_{2R}, \mbox{ (in the sense of the traces) }.
\end{equation}
Let us prove that
\begin{equation}\label{Ell:E23.2}
	u\in H^3\left(B^+_{r}\right),\quad\forall r<R.
\end{equation}
To this purpose, let us prove 

\medskip

\textbf{Claim.} Let $l=1,\cdots,n-1$. We have, for every $r\in (0,2R)$,
\begin{equation}\label{Ell:E23.3}
	\partial_lu(\cdot,0)=0, \quad\mbox{ in } B'_{r}, \mbox{ (in the sense of the traces)}.
\end{equation}

\medskip

\textbf{Proof of the Claim.} First let us note that, since $u\in H^2\left(B^+_{r}\right)$, for every $r<R$,  $\partial_lu(x',0)$ is well--defined in the sense of traces. Now, let
$$v=\partial_lu$$ and let us denote by $T(v)$ the trace of $v$ on $\{x_n=0\}$. As a consequence of Theorem \ref{traccia}, $T(v)$ is  characterised by the identity.

\begin{equation}\label{Ell:E25.0}
	-\int_{B'_R}\Phi_n (x',0) T(v)dx'=\int_{B^+_R}v \mbox{div}\Phi dx+\int_{B^+_R}\nabla v \cdot\Phi dx,
\end{equation}
for every $\Phi\in C^{\infty}_0\left(B_R; \mathbb{R}^n\right)$.

\smallskip

Let now $\Phi$ be any function belonging to $C^{\infty}_0\left(B_R; \mathbb{R}^n\right)$. By \eqref{Ell:E23.2} and by Theorem \ref{traccia} (applied to $u$), we get

\begin{equation*}
	\begin{aligned}
		\int_{B^+_R}v (\mbox{div}\Phi) dx&=\int_{B^+_R}\partial_l u (\mbox{div}\Phi) dx=-\int_{B^+_R}u(\partial_l\mbox{div}\Phi) dx=\\&=-\int_{B^+_R}u(\mbox{div}\partial_l\Phi) dx=\\&
		=\int_{B^+_R}\nabla u\cdot \partial_l\Phi dx+\int_{B'_R}\partial_l\Phi_n (x',0) T(u)dx'=\\&=\int_{B^+_R}\nabla u\cdot \partial_l\Phi dx=
		-\int_{B^+_R}\partial_l\nabla u\cdot \Phi dx=\\&=-\int_{B^+_R}\nabla v\cdot \Phi.
	\end{aligned}
\end{equation*}
Hence, \eqref{Ell:E25.0} implies

\begin{equation*}
	T(v)=0,\quad \mbox{ in } B'_r, \quad\forall r\in (0,2R).
\end{equation*}
Claim is proved.

\medskip
We now briefly and only formally show the most significant steps to complete the proof; we encourage the reader to treat the steps in a rigorous manner using appropriately the variational formulation in a similar way as we did in the proof of Theorem \eqref{Ell:reg-m}.

By calculating the derivative w.r.t. $x_l$, $l\in \{1,\cdots, n-1\}$, of both the sides of the equation

\begin{equation}\label{Ell:E25.00}
	-\sum_{j,k=1}^n\partial_j\left(a^{jk}\partial_ku\right)=f,\end{equation}
we obtain

\begin{equation}\label{Ell:E25.1}
	-\sum_{j,k=1}^n\partial_j\left(a^{jk}\partial_kv\right)=\widetilde{f},\end{equation}
where $v=\partial_l u$ and

$$\widetilde{f}=\partial_lf+\sum_{j,k=1}^n\partial_j\left((\partial_la^{jk})\partial_ku\right).$$
Now, as $f\in H^1\left(B^+_{2R}\right)$, $a^{jk}\in C^{1,1}\left(B^+_{2R}\right)$, $u\in H^2\left(B^+_r\right)$ for every $r\in (0,2R)$, we get

$$\widetilde{f}\in L^2\left(B^+_r\right),\quad\forall r\in (0,2R).$$
on the other hand, for every $r\in (0,2R)$, we have $v\in H^1\left(B^+_r\right)$ and
\begin{equation*}
	v(\cdot,0)=0, \quad\mbox{ in } x'\in B'_{r}, \mbox{ (in the sense of the traces) }.
\end{equation*}
Therefore by Lemma \ref{Ell:reg-bordo-loc} we have
\begin{equation}\label{Ell:E25.2}\partial_l u=v\in H^2\left(B^+_r\right),\quad\forall r\in (0,2R),\end{equation}
for $l\in \{1,\cdots, n-1\}$.
Moreover by \eqref{Ell:E25.1} we get (likewise to \eqref{Ell:E16.1bis})
\begin{equation}\label{Ell:E25.3}
	\partial^3_nu=\partial_n\left[-\frac{1}{a^{nn}}\left(\sum_{\underset{k+j<2n}{k,j=1}}^n a^{jk}\partial^2_{jk}u+
	\sum_{k,j}^n\partial_ja^{jk}\partial_{k}u+f\right)\right],\end{equation}
from which, taking into account \eqref{Ell:E25.2}, we get 
$$u\in H^3\left(B^+_r\right),\quad\forall r\in (0,2R).$$
Finally, inequality \eqref{Ell:reg-bordo-m4} (for $m=1$) is obtained applying inequality \eqref{Ell:E12.6L} to \eqref{Ell:E25.1}. In addition, inequality \eqref{Ell:E12.6L} gives the estimates

\begin{equation*}
	\begin{aligned}
		R^6\sum_{\underset{k+j+l<3n}{k,j,l=1}}^n\left\Vert\partial^3_{jkl}u\right\Vert^2_{L^2\left(B^+_R\right)}&\leq C\left(1+E_1^2\right)\int_{B_{2R}}u
		^2dx+\\&+C\sum_{|\alpha|\leq 1}R^{2|\alpha|+4}\int_{B^+_{2R}}\left|\partial^{\alpha}\right| ^2dx.
	\end{aligned}
\end{equation*}By the last inequality, by means of \eqref{Ell:E25.3} (applying
again inequality \eqref{Ell:E12.6L} to equation \eqref{Ell:E25.00}), we obtain \eqref{Ell:reg-bordo-m4} (for $m=1$).

In order to complete the proof, simply follow the proof of Theorem \ref{reg-bordo} taking into account that diffeomorphisms $\Phi$ e $\Psi$ are, in this case, of class $C^{2,1}$ ($C^{m+1,1}$ in the general case). $\blacksquare$

\bigskip

\begin{cor} [\textbf{$C^{\infty}$ global regularity}]\label{Ell:reg-infin-bordo}
	\index{Corollary:@{Corollary:}!- $C^{\infty}$ global regularity@{- $C^{\infty}$ global regularity}}
	Let $\Omega$ be a bounded open set of
	$\mathbb{R}^n$ with boundary of class $C^{\infty}$. Let
	\begin{equation*}
		f\in C^{\infty}(\overline{\Omega}).
	\end{equation*}
	Let $A\in
	C^{\infty}(\overline{\Omega};\mathbb{M}(n))$. Let us assume that $A$ satisfies \eqref{gamma-eq}.
	
	Let us assume that $u\in H_0^1(\Omega)$ is a solution to
	
	\begin{equation*}
		-\mbox{div}(A\nabla u)=f, \quad\mbox{ in } \Omega
	\end{equation*}
	Then we have
	\begin{equation*}
		u\in C^{\infty}(\overline{\Omega}).
	\end{equation*}
\end{cor}

\textbf{Proof.} By Theorem \ref{Ell:reg-bordo-m} and by the Embedding Theorem \ref{Sobolev-inequ} we have

$$u\in \bigcap_{m=0}^{\infty}H^m(\Omega)=C^{\infty}(\overline{\Omega}).$$ $\blacksquare$

\bigskip

\section{The Dirichlet to Neumann Map} \label{Mappa Dir-Neum} Denote by $\mathbb{M}^S(n)$ the vector space of symmetric matrix $n\times$ with real entries\index{$\mathbb{M}^S(n)$}. Let
$\Omega$ be a bounded open set of $\mathbb{R}^n$ of class $C^{0,1}$. Let $A\in
L^{\infty}(\Omega;\mathbb{M}^S(n))$ and let us suppose that \eqref{gamma-eq} holds. Formally, the Dirichlet to Neumann Map \index{Dirichlet to Neumann Map}can be constructed in the following way: let $\varphi\in
H^{1/2}(\partial\Omega)$ and let $u\in H^1(\Omega)$ be the solution to the
problem
\begin{equation}\label{D-N-0}
	\begin{cases}
		\mbox{div}(A\nabla u)=0, \quad\mbox{ in } \Omega, \\
		\\
		u=\varphi, \quad\mbox{ on } \partial\Omega.
	\end{cases}
\end{equation}
We have seen that $u$  is uniquely determined by
$\varphi$ hence, \textit{if it would make sense}, we could define the
map 
\begin{equation}\label{D-N-1}
	\varphi\rightarrow A \nabla u\cdot \nu,\quad \mbox{
		(conormal derivative of $u$ on $\partial \Omega$). }
\end{equation}
Let us observe that if $u\in C^{2}(\overline{\Omega})$ and $A\in
C^{1}\left(\overline{\Omega};\mathbb{M}^S(n)\right)$, by
\eqref{D-N-0} we have

\begin{equation}\label{D-N-1-1}
	\int_{\partial\Omega}(A\nabla u\cdot \nu) v
	dS=\int_{\Omega}A\nabla u\cdot\nabla v dx, \mbox{ }
	\forall v \in C^{1}(\overline{\Omega}).
\end{equation}
As a matter of fact, if $v \in C^{1}(\overline{\Omega})$, by the divergence Theorem and by \eqref{D-N-0} we get

\begin{equation*}
	\begin{aligned}
		&\int_{\Omega}A\nabla u\cdot\nabla v
		dx=\int_{\Omega}\left(\mbox{div}(v A\nabla u)-
		\mbox{div}(A\nabla u)v\right)dx= \int_{\partial\Omega}(A\nabla
		u\cdot \nu)v dS.
	\end{aligned}
\end{equation*}
Equality \eqref{D-N-1-1} allows us to "read" $A \frac{\partial
	u}{\partial \nu}$ by means integral on the right--hand side. Now,
if $A\in L^{\infty}(\Omega;\mathbb{M}^S(n))$ and $u\in
H^1(\Omega)$ then the integral on the right--hand side of \eqref{D-N-1-1} makes perfectly sense. Based on these insights we will define
below $A \nabla u\cdot\nu$ as an element of
$H^{-1/2}(\partial\Omega)$ (the dual space of $H^{1/2}(\partial\Omega)$).

\medskip

First, notice that 
\begin{equation}\label{D-N-11}
	\int_{\Omega}A\nabla u\cdot\nabla v dx=0, \quad \forall v
	\in H_0^1(\Omega).
\end{equation}
As a matter of fact, recalling \eqref{variazDirichlet-1} and \eqref{variazDirichlet-25},
we have, for each \\ $v \in H_0^1(\Omega)$,

\begin{equation*}\label{D-N-20}
	\int_{\Omega}A\nabla u\cdot\nabla v dx=\int_{\Omega}\{A\nabla
	\Phi\cdot\nabla v +A\nabla w\cdot\nabla v\}dx=0.
\end{equation*}
By \eqref{D-N-11}, recalling Theorem \ref{traccia}, we have that if
$\varphi \in H^1(\Omega)$ then the integral

\begin{equation}\label{D-N-12}
	\int_{\Omega}A\nabla u\cdot\nabla v dx
\end{equation}
depends only by the trace of $v$ on $\partial\Omega$. As a matter of fact,
if $v_1,v_2 \in H^1(\Omega)$  have the same trace on
$\partial\Omega$ then $v_1-v_2\in H_0^1(\Omega)$.
Hence \eqref{D-N-11} implies $$\int_{\Omega}A\nabla
u\cdot\nabla v_1 dx=\int_{\Omega}A\nabla u\cdot\nabla v_2
dx.$$ Therefore, for every $\varphi\in H^{1/2}(\partial\Omega)$ it turns out to be  well--defined the functional  $$L_{\varphi}:
H^{1/2}(\partial\Omega)\rightarrow \mathbb{R},$$ which maps $\phi\in
H^{1/2}(\partial\Omega)$ in the real number
$$L_{\varphi}(\phi)=\int_{\Omega}A\nabla u\cdot\nabla v dx,$$
where $v_{|\partial\Omega}=\phi$ (in the sense of traces).

\medskip

\textbf{We prove that  $L_{\varphi}$ is a linear and bounded functional.}

\medskip

The linearity of $L_{\varphi}$ is trivial. Concerning the boundedness, let \\ $$\phi\in H^{1/2}(\partial\Omega),$$ by Theorem \ref{Sob:teo4.11}
we know that there exists $v\in H^1(\Omega)$ so that \\
$$v_{|\partial\Omega}=\phi, \ \  \mbox{(in the sense of traces)}$$ which satisfies

\begin{equation}\label{D-N-2}
	\left\Vert v\right\Vert_{H^{1}(\Omega)}\leq C\left\Vert
	\phi\right\Vert_{H^{1/2}(\partial \Omega)},
\end{equation}
where $C$ is a constant  depending on $\Omega$ only. Proceeding in a
similar way to what we did to obtain
\eqref{variazDirichlet-21} and taking into account
\eqref{variazDirichlet-30}, we get

\begin{equation}\label{D-N-10}
	\begin{aligned}
		|L_{\varphi}(\phi)|&=\left\vert\int_{\Omega}A\nabla u\cdot\nabla v
		dx\right\vert \leq \lambda \left\Vert\nabla
		u\right\Vert_{L^{2}(\Omega)}\left\Vert
		\nabla v\right\Vert_{L^{2}(\Omega)}\leq\\&\leq
		\overline{C}\left\Vert \varphi\right\Vert_{H^{1/2}(\partial
			\Omega)}\left\Vert \phi\right\Vert_{H^{1/2}(\partial \Omega)},
	\end{aligned}
\end{equation}
where $\overline{C}$ is a constant depending on $\Omega$
and $\lambda$ only. Therefore, the functional $L_{\varphi}$ is bounded and it satisfies

\begin{equation}\label{Funz-D-N-2}
	\left\Vert L_{\varphi}\right\Vert_{H^{-1/2}(\partial\Omega)}\leq
	\overline{C}\left\Vert \varphi\right\Vert_{H^{1/2}(\partial \Omega)},\quad \forall \varphi\in H^{1/2}(\partial \Omega).
\end{equation}
Moreover, notice that \eqref{Funz-D-N-2} implies that the linear operator $$H^{1/2}(\partial \Omega)\ni\varphi\rightarrow
L_{\varphi}\in H^{-1/2}(\partial \Omega),$$ is bounded.

Now we set 
$$A\nabla u\cdot \nu:=L_{\varphi}$$
consequently, we write, for any $v \in H^1(\Omega)$ such that
$v_{|\partial\Omega}=\phi$,
\begin{equation}\label{Funz-D-N-3}
	\langle {A\nabla u\cdot \nu,\phi} \rangle_{H^{-1/2},
		H^{1/2}}=L_{\varphi}(\phi)=\int_{\Omega}A\nabla u\cdot\nabla v dx,
\end{equation}
where, by $\langle {\cdot,\cdot} \rangle_{H^{-1/2}, H^{1/2}}$ we denote
the scalar product in the duality. With the notations used so far
we should have written $(A\nabla u\cdot \nu)(\phi)$ instead of
$\langle {A\nabla u\cdot \nu,\phi} \rangle_{H^{-1/2}, H^{1/2}}$, but the latter
notation it is certainly more handleable in the present context.
.

\medskip

Finally, we define the \textbf{Dirichlet to Neumann Map} as

\begin{equation}\label{Mappa-D-N}
	\Lambda_{A}: H^{1/2}(\partial \Omega)\rightarrow H^{-1/2}(\partial
	\Omega),\quad \Lambda_{A}(\varphi)=A\nabla u\cdot \nu.
\end{equation}
By \eqref{Funz-D-N-2} it follows that $\Lambda_{A}\in \mathcal{L}\left(
H^{1/2}(\partial \Omega), H^{-1/2}(\partial\Omega)\right)$, where, we recall, $\mathcal{L}\left(H^{1/2}(\partial
\Omega),H^{-1/2}(\partial \Omega)\right)$ denotes the space of the linear and bounded
operators from $H^{1/2}(\partial \Omega)$ to
$H^{-1/2}(\partial\Omega)$. From the that construction we have performed so far we get

\begin{equation}\label{Mappa-D-N-bilin}
	\langle {\Lambda_{A}(\varphi),\phi} \rangle_{H^{-1/2},
		H^{1/2}}=\int_{\Omega}A\nabla u\cdot\nabla v dx, \quad \forall \varphi, \phi\in H^{1/2}(\partial \Omega),
\end{equation}
where $u\in H^1(\Omega)$ is the solution to \eqref{D-N-0} and
$v$ is any function of $H^1(\Omega)$ which satisfies
$v_{|\partial\Omega}=\phi$. It is simple to check that

\begin{equation*}
	\begin{aligned}
		&\left\Vert
	\Lambda_{A}\right\Vert_{\mathcal{L}\left(H^{1/2},H^{-1/2}\right)}=\\&=\sup\left\{\langle
	{\Lambda_{A}(\varphi),\phi} \rangle_{H^{-1/2}, H^{1/2}}: \left\Vert
	\varphi\right\Vert_{H^{1/2}(\partial \Omega)}\leq 1, \left\Vert
	\phi\right\Vert_{H^{1/2}(\partial \Omega)}\leq 1\right\}.
\end{aligned}
\end{equation*}

\medskip

In what follows we prove other simple but important properties of
$\Lambda_{A}$.

\medskip

We first observe that since the right-hand integral
in \eqref{Mappa-D-N-bilin} is independent of the choice of
$v$ (as long as it is a trace of $\phi$) we can choose
$v=w$ where $w\in H^1(\Omega)$ is the solution of the Dirichlet problem

\begin{equation}\label{D-N-0-v}
	\begin{cases}
		\mbox{div}(A\nabla w)=0, \quad\mbox{ in } \Omega, \\
		\\
		w=\phi, \quad\mbox{ on } \partial\Omega.
	\end{cases}
\end{equation}
From this it follows that the bilinear form
$$H^{1/2}(\partial \Omega)\times H^{1/2}(\partial\Omega)\ni(\varphi,\phi)\rightarrow\langle {\Lambda_{A}(\varphi),\phi} \rangle_{H^{-1/2}, H^{1/2}}\in\mathbb{R}$$
is  \textbf{symmetric} that is

\begin{equation}\label{Simmetria}
	\langle {\Lambda_{A}(\varphi),\phi} \rangle_{H^{-1/2}, H^{1/2}}=\langle
	{\Lambda_{A}(\phi),\varphi} \rangle_{H^{-1/2}, H^{1/2}}, \quad \forall \varphi, \phi\in H^{1/2}(\partial \Omega).
\end{equation}
As a matter of fact
\begin{equation*}
	\begin{aligned}
	\langle {\Lambda_{A}(\varphi),\phi} \rangle_{H^{-1/2},
		H^{1/2}}&=\int_{\Omega}A\nabla u\cdot\nabla w dx=\\&=\int_{\Omega}A\nabla
	w\cdot\nabla u dx=\\&=\langle {\Lambda_{A}(\phi),\varphi}
	\rangle_{H^{-1/2}, H^{1/2}}.	
\end{aligned}
	\end{equation*}
An important consequence of \eqref{Simmetria} is the following identity

\begin{theo}[\textbf{the Alessandrini identity}]\label{id-Alessandrini}
	\index{Theorem:@{Theorem:}!- Alessandrini identity@{- Alessandrini identity}}
	Let $A_1,A_2\in L^{\infty}(\Omega;\mathbb{M}^S(n))$. Let us assume that $A_1,A_2$
	satisfy \eqref{gamma-eq}. Let $\varphi,\phi\in
	H^{1/2}(\partial\Omega)$ and let $u_1,u_2\in H^{1}(\Omega)$ be the
	solutions to 
	\begin{equation}\label{D-N-0-u-1}
		\begin{cases}
			\mbox{div}(A\nabla u_1)=0, \quad\mbox{ in } \Omega, \\
			\\
			u_1=\varphi, \quad\mbox{ on } \partial\Omega
		\end{cases}
	\end{equation}
	and
	\begin{equation}\label{D-N-0-u-2}
		\begin{cases}
			\mbox{div}(A\nabla u_2)=0, \quad\mbox{ in } \Omega, \\
			\\
			u_2=\phi, \quad\mbox{ on } \partial\Omega,
		\end{cases}
	\end{equation}
	then
	
	\begin{equation}\label{id-A}
		\langle {(\Lambda_{A_1}-\Lambda_{A_2})(\varphi),\phi}
		\rangle_{H^{-1/2}, H^{1/2}}=\int_{\Omega}(A_1-A_2)\nabla
		u_1\cdot\nabla u_2 dx.
	\end{equation}
\end{theo}

\medskip

\textbf{Proof.} By \eqref{Simmetria} we have (by setting for sake of brevity, $\langle
{\cdot,\cdot} \rangle=\langle {\cdot,\cdot} \rangle_{H^{-1/2},
	H^{1/2}}$)

\begin{equation}\label{id-A-1}
	\begin{aligned}
		\langle {(\Lambda_{A_1}-\Lambda_{A_2})(\varphi),\phi} \rangle&=\langle
		{\Lambda_{A_1}(\varphi),\phi} \rangle-\langle
		{\Lambda_{A_2}(\varphi),\phi} \rangle=\\&=\langle
		{\Lambda_{A_1}(\varphi),\phi} \rangle-\langle
		{\Lambda_{A_2}(\phi),\varphi} \rangle.
	\end{aligned}
\end{equation}
Now, by \eqref{D-N-0-u-1} and \eqref{D-N-0-u-2} we have,
respectively,

$$\langle {\Lambda_{A_1}(\varphi),\phi}\rangle=\int_{\Omega}A_1\nabla u_1\cdot\nabla u_2 dx$$
and

$$\langle {\Lambda_{A_2}(\phi),\varphi}\rangle=\int_{\Omega}A_2\nabla u_2\cdot\nabla u_1 dx= \int_{\Omega}A_2\nabla u_1\cdot\nabla u_2$$
by inserting the latter in \eqref{id-A-1} we get

\begin{equation*}
	\langle {(\Lambda_{A_1}-\Lambda_{A_2})(\varphi),\phi}
	\rangle=\int_{\Omega}(A_1-A_2)\nabla u_1\cdot\nabla u_2 dx.
\end{equation*}
$\blacksquare$

\bigskip

A simple consequence of the Alessandrini identity is the continuity of the map $A\rightarrow\Lambda_{A}$. Precisely we have the following
\begin{prop}\label{continuita mappa}
	There exists a constant $C$  dpending on $\lambda$ and $\Omega$ only
	so that, if $A_1,A_2\in L^{\infty}(\Omega;\mathbb{M}^S(n))$
	satisfy \eqref{gamma-eq}, then 
	\begin{equation}\label{continuita mappa-1}
		\left\Vert
		\Lambda_{A_1}-\Lambda_{A_2}\right\Vert_{\mathcal{L}\left(H^{1/2},H^{-1/2}\right)}\leq
		C\left\Vert A_1-A_2\right\Vert_{L^{\infty}(\Omega,\mathbb{M}^S(n))}.
	\end{equation}
\end{prop}

\medskip

\textbf{Proof.} By the Alessandrini identity and by the Cauchy--Schwarz inequality we have

\begin{equation}\label{continuita mappa-2}
	\begin{aligned}
		&\left\vert\langle {(\Lambda_{A_1}-\Lambda_{A_2})(\varphi),\phi}
		\rangle_{H^{-1/2}, H^{1/2}}\right\vert=
		\left\vert\int_{\Omega}(A_1-A_2)\nabla u_1\cdot\nabla u_2 dx\right\vert\leq \\
		& \leq\left\Vert
		A_1-A_2\right\Vert_{L^{\infty}(\Omega;\mathbb{M}^S(n)))}\left\Vert\nabla
		u_1\right\Vert_{L^{2}(\Omega)}\left\Vert\nabla
		u_2\right\Vert_{L^{2}(\Omega)}.
	\end{aligned}
\end{equation}
Now, \eqref{variazDirichlet-30} gives

$$\left\Vert\nabla u_1\right\Vert_{L^{2}(\Omega)}\leq C\left\Vert \varphi\right\Vert_{H^{1/2}(\partial \Omega)}
\mbox{ and }\left\Vert\nabla u_2\right\Vert_{L^{2}(\Omega)}\leq C\left\Vert \phi\right\Vert_{H^{1/2}(\partial \Omega)},$$
where $C$ depends on $\lambda$ and $\Omega$ only. Inserting the last obtained  inequalities in \eqref{continuita mappa-2} we have
\begin{equation*}
	\begin{aligned}
	&\left\vert\langle {(\Lambda_{A_1}-\Lambda_{A_2})(\varphi),\phi}
	\rangle_{H^{-1/2}, H^{1/2}}\right\vert\leq \\& \leq C\left\Vert
	A_1-A_2\right\Vert_{L^{\infty}(\Omega;\mathbb{M}^S(n))}\left\Vert
	\varphi\right\Vert_{H^{1/2}(\partial \Omega)}\left\Vert
	\phi\right\Vert_{H^{1/2}(\partial \Omega)}
	\end{aligned}
\end{equation*}
From which \eqref{continuita mappa-1} follows. $\blacksquare$

\bigskip

Similarly, one can also define the \textbf{Neumann to
	Dirichlet Map}. Let $g\in H^{-1/2}(\partial \Omega)$ satisfy
$$\langle g,1\rangle=0.$$
Let us consider the solution $u\in H^1(\Omega)$ to the Neumann problem

\begin{equation}\label{variazNeumann-1-10122}
	\begin{cases}
		\int_{\Omega} A\nabla u\cdot \nabla v  dx=\langle g, \varphi \rangle_{H^{-1/2}, H^{1/2}}, \quad\forall v \in H^{1}(\Omega), \\
		\\
		u\in \left\{w\in H^{1}(\Omega):\quad \int_{\Omega}wdx=0\right\},
	\end{cases}
\end{equation}
we define the Neumann to Dirichlet Map as follows
\begin{equation}\label{Mappa-D-N-10122}
	\begin{aligned}
	&\mathcal{N}_{A}: H^{-1/2}(\partial \Omega)\rightarrow
	H^{1/2}(\partial \Omega),\quad \\& \mathcal{N}_{A}(g)=u_{|\partial
		\Omega}, \quad \mbox{(in the sense of the traces)}.
	\end{aligned}
\end{equation}

\bigskip

\section{The inclusion inverse problem}\label{size-est}
Let $n=2$ or $n=3$ and let us assume that $\Omega$ represents an electrically
conductor of constant conductivity, say, $1$ and
let us suppose that $\Omega$ contains an \textit{inclusion} $D$ of
different conductivity, say $k$, with $k>0$ and $k\neq 1$.
We consider the problem of determining $D$ \index{inclusion inverse problem}from the knowledge of a
density of prescribed current on $\partial \Omega$ and of the
corresponding voltage $u$ measured on $\partial \Omega$.

We provide a mathematical formulation of the problem. Let us assume that
$\Omega$ is a bounded open set of $\mathbb{R}^n$
($n=2,3$) whose boundary is of class $C^{0,1}$, let $\phi\in
H^{-1/2}(\partial \Omega)$ satisfy

\begin{equation}\label{size-est-1}
	\int_{\partial \Omega}\phi dS=0,
\end{equation}
Where we have denoted
$$\int_{\partial \Omega}\phi d S=\langle
\phi,1 \rangle_{H^{-1/2}, H^{1/2}}.$$ $\phi$ represents the 
density of prescribed current on $\partial \Omega$. If
the inclusion is present, the electrostatic potenzial  $u$ is determined, up to an additive constant, as a solution of Neumann problem

\begin{equation}\label{size-est-2}
	\begin{cases}
		\mbox{div}\left(\left(1+(k-1)\chi_D\right)\nabla u\right)=0, \quad\mbox{ in } \Omega, \\
		\\
		\frac{\partial u}{\partial\nu}=\phi, \quad\mbox{ on }
		\partial\Omega,
	\end{cases}
\end{equation}
where $D$ is a measurable subset of $\Omega$. In what follows we assume
\begin{equation}\label{size-est-3}
	\int_{\Omega}u(x) dx=0,
\end{equation}
which yields with the uniqueness of the boundary value problem  \eqref{size-est-2}.

The \textbf{inverse problem} consists in determining $D$,
 by assigning a nontrivial input $\phi$ and
measuring the corresponding trace $u_{|\partial\Omega}$. The uniqueness of $D$ is still an open question.
We point out that if we dispose of the entire Neumann to Dirichlet Map
	(or the Dirichlet to Neumann Map), and $\partial D$ is enough regular, the uniqueness can be proved (see \cite{Is2}). Keep in mind that having the entire Dirichlet to Neumann Map is equivalent to being able to make
	infinite measurements on $\partial \Omega$. However, the ideas (developed
	detailed in \cite{A-R-M}) that we will present here allow us to find
	 \textbf{size estimates (of volume or area)} \index{size estimates}of the inclusion from
	certain integrals of the data, $\phi$ and $u_{|\partial\Omega}$ as we show below. 
	
	 Let us consider the quantity 

\begin{equation}\label{size-est-4}
	W=\int_{\partial\Omega}\phi u,
\end{equation}
and compare it with 
\begin{equation}\label{size-est-5}
	W_0=\int_{\partial\Omega}\phi u_0,
\end{equation}
where $u_0$ represents the solution to the Neumann problem
\begin{equation}\label{size-est-6}
	\begin{cases}
		\Delta u_0=0, \quad\mbox{ in } \Omega, \\
		\\
		\frac{\partial u_0}{\partial\nu}=\phi, \quad\mbox{ on }
		\partial\Omega,
	\end{cases}
\end{equation}
with $$\int_{\Omega}u_0=0.$$

\smallskip

$W$ and $W_0$ represent the power required to maintain the current
$\phi$, when the inclusion $D$ is present and it is not present, respectively. 
Partially anticipating the results that
\textbf{we will prove later in this Section, we have that if
	$k\neq 1$, $k>0$, then the following inequality holds true}

\begin{equation}\label{size-est-7}
	C_1\int_{D}\left|\nabla u_0\right|^2\leq |W_0-W|\leq C_2
	\int_{D}\left|\nabla u_0\right|^2,
\end{equation}
where $C_1$ and $C_2$ are positive constants  depending on $k$ only.
Let us assume, for instance+, that $D$ is a connected open set.
 Inequalities \eqref{size-est-7} implies that if $\phi\neq 0$
 then $|D|$  is $0$  if and only if $W_0-W=0$. Let us prove this claim. If $|D|=0$ then by the second
 inequality we immediately have $W_0-W=0$. Conversely, if
$W_0-W=0$, then  by the first inequality we have that if $D\neq
\emptyset$ then 

$$u_0=\mbox{constant in } D.$$
Since $u_0$ is an analytic function, we have
$$u_0=\mbox{constant in } \Omega,$$
consequently $\phi=0$, but we have assumed $\phi\neq 0$,
therefore $D=\emptyset$.

\bigskip

Inequalities  \eqref{size-est-7} can be proved as a consequence of
general properties of the continuous symmetric coercive bilinear forms on a Hilbert space.

\medskip

Introduce some notation. Let $H$ be a real Hilbert space
and $H'$ its dual space.

Let $\lambda_0,\lambda_1\in [1,+\infty)$ and let $a_0(\cdot,\cdot)$ and
$a_1(\cdot,\cdot)$ two symmetric bilinear forms on $H$ which satisfy
the conditions

\begin{subequations}
	\label{size-est-8}
	\begin{equation}
		\label{size-est-8a} \lambda^{-1}_0\|u\|^2\leq
		a_0(u,u)\leq\lambda_0\|u\|^2,\qquad \forall\,u\in H,
	\end{equation}
	\begin{equation}
		\label{size-est-8b} \lambda^{-1}_1\|u\|^2\leq a_1(u,u)\leq
		\lambda_1\|u\|^2,\qquad \forall\,u\in H.
	\end{equation}
\end{subequations}
Let us note that \eqref{size-est-8a} and \eqref{size-est-8b}
imply, respectively, the continuity of $a_0$ and of $a_1$. Just
check this for $a_0(\cdot,\cdot)$. We have that \eqref{size-est-8a} implies
$$|a_0(u,v)|\leq \sqrt{a_0(u,u)}\sqrt{a_0(v,v)}\leq \lambda_0\|u\|\|v\|,\quad
\forall u,v\in H.$$ Moreover, let
\begin{equation}\label{size-est-9}
	\alpha(u,v)=a_1(u,v)-a_0(u,v),\qquad u,v\in H.
\end{equation}

Let $F\in H'$. By the Lax--Milgram Theorem, there exist $u_1,u_0\in H$ such that
\begin{equation}\label{size-est-10}a_j(u_j,v)=<F,v>\qquad \forall\,v\in H,\quad j=0,1.
\end{equation} Define
\begin{equation}\label{size-est-10-1}
	W_0=<F,u_0>,\quad W_1=<F,u_1>,\quad \delta W=W_0-W_1.
\end{equation}

Now we prove two simple lemmas.
\begin{lem}\label{correct:8-4-23-1}
	The following equalities hold true.
	\begin{subequations}
		\label{size-est-11}
		\begin{equation}
			\label{size-est-11a} a_0(u_1-u_0,u_1-u_0)-\alpha(u_0,u_0)=-\delta W,
		\end{equation}
		\begin{equation}
			\label{size-est-11b} a_0(u_0-u_1,u_0-u_1)+\alpha(u_1,u_1)=\delta W,
		\end{equation}
		\begin{equation}
			\label{size-est-11c} \alpha(u_1,u_0)=-\delta W.
		\end{equation}
	\end{subequations}
\end{lem}
\textbf{Proof.}
Let us check \eqref{size-est-11a}.
\begin{equation*}
	\begin{aligned}
		&a_0(u_1-u_0,u_1-u_0)-\alpha(u_0,u_0)=\\&=a_0(u_1-u_0,u_1-u_0)-[a_1(u_0,u_0)-a_0(u_0,u_0)]=\\&
		=a_1(u_1,u_1)-2a_1(u_1,u_0)+a_1(u_0,u_0)-a_1(u_0,u_0)+a_0(u_0,u_0)=\\
		&=a_1(u_1,u_1)-2a_1(u_1,u_0)+a_0(u_0,u_0)=\\
		&=<F,u_1>-2<F,u_0>+<F,u_0>=\\&=<F,u_1-u_0>=-\delta W.
	\end{aligned}
\end{equation*}
Equality \eqref{size-est-11b} can be obtained similarly and
\eqref{size-est-11c} is an immediate consequence of
\eqref{size-est-9}. $\blacksquare$

\begin{lem}
	\label{size-160122-33} If one of the following conditions is
	satisfied
	$$\alpha(u,u)\geq 0, \quad \forall u\in H, $$ or $$\alpha(u,u)\leq 0, \quad \forall u\in H, $$
	then we have
	\begin{equation}
		\label{size-160122-eq4}
		|\alpha(u,v)|\leq|\alpha(u,u)|^{1/2}|\alpha(v,v)|^{1/2},\quad\forall
		u,v\in H.
	\end{equation}
\end{lem}

\medskip

\textbf{Proof.} Let $u,v\in H$. If $\alpha(u,u)=0$ and
$\alpha(v,v)=0$, then, assuming $\alpha(w,w)\geq 0$, for every
$w\in H$, we have
$$0\leq\alpha(u+tv,u+tv)=2t\alpha(u,v),\qquad \forall\,t\in \mathbb R,$$
which implies $\alpha(u,v)=0$ and \eqref{size-160122-eq4} is
proved.\\
If either $\alpha(u,u)\neq0$ or $\alpha(v,v)\neq0$, then, assuming,
for instance, $\alpha(v,v)>0$, we have
$$0\leq\alpha(u+tv,u+tv)=t^2\alpha(v,v)+2t\alpha(u,v)+\alpha(u,u), \qquad \forall\,t\in\mathbb R,$$
hence
$$\left(\alpha(u,v)\right)^2-\alpha(u,u)\alpha(v,v)\leq0$$
which gives \eqref{size-160122-eq4}.\\
If $\alpha(w,w)\leq0$, for every $w\in H$, the thesis follows easily
by applying the previous procedure to $-\alpha(\cdot,\cdot)$.
$\blacksquare$

\bigskip

Now we prove 

\begin{theo}
	\label{teo-size} Let $F\in H'$. Let us assume that the bilinear forms
	$a_0(\cdot,\cdot)$ and $a_1(\cdot,\cdot)$ satisfy conditions
	\eqref{size-est-8} and let us assume that $u_0,u_1$ satisfy
	\eqref{size-est-10}. If $\alpha(\cdot,\cdot)$ (defined by
	\eqref{size-est-9}), satisfies 
	\begin{equation}
		\label{size-160122-4a} 0\leq\alpha(u,u)\leq C_0
		a_0(u,u),\quad\forall\,u\in H,
	\end{equation}
	where $C_0$ is a positive constant, then
	$$\delta W\geq 0$$ and
	\begin{equation}
		\label{size-160422-1-6} \delta W\leq\alpha(u_0,u_0)\leq(1+C_0)\delta
		W.
	\end{equation}
	If $\alpha(\cdot,\cdot)$ satisfies the condition
	\begin{equation}
		\label{size-160122-4b} \alpha(u,u)\leq 0,\quad\forall\,u\in H,
	\end{equation}
	then $$\delta W\leq 0$$ and
	\begin{equation}
		\label{size-160122-2-7} -C\delta W\leq -\alpha(u_0,u_0)\leq-\delta
		W,
	\end{equation}
	where $C$ is a positive constant depending on $\lambda_0$
	and $\lambda_1$ only.
\end{theo}
\textbf{Proof.} First we consider the case in which \eqref{size-160122-4a} holds. By \eqref{size-est-11b} we have $\delta W\geq 0$ and by \eqref{size-est-11a}
we have $-\alpha(u_0,u_0)\leq-\delta W$. Therefore
\begin{equation}
	\label{size-160122-5-5} \delta W\leq\alpha(u_0,u_0).
\end{equation}
Now, let us estimate $\alpha(u_0,u_0)$ from above. Lemma \ref{correct:8-4-23-1} -- b and Lemma
\ref{size-160122-33} give
\begin{equation*}
	\begin{aligned}
		&\alpha(u_0,u_0) =
		\alpha(u_0-u_1,u_0-u_1)+\alpha(u_1,u_1)+2\alpha(u_0-u_1,u_1)\leq\\[2mm]&\leq\alpha(u_0-u_1,u_0-u_1)+\alpha(u_1,u_1)
		+2|\alpha(u_0-u_1,u_0-u_1)|^{1/2}|\alpha(u_1,u_1)|^{1/2}\leq\\[2mm]
		&\leq \alpha(u_0-u_1,u_0-u_1)+\alpha(u_1,u_1)
		+\frac{1}{C_0}\alpha(u_0-u_1,u_0-u_1)+C_0\alpha(u_1,u_1)=\\&
		=(1+C_0)[\frac{1}{C_0}\alpha(u_0-u_1,u_0-u_1)+\alpha(u_1,u_1)]\leq\\[2mm]&
		\leq(1+C_0)\max\left\{1,\frac{1}{C_0}\right\}[a_0(u_0-u_1,u_0-u_1)+\alpha(u_1,u_1)]\leq\\[2mm]
		&\leq C_1\delta W,
	\end{aligned}
\end{equation*}
where, in the last inequality we have set
$$C_1=(1+C_0)\max\left\{1,\frac{1}{C_0}\right\}.$$
Hence we have
$$\alpha(u_0,u_0)\leq C_1\delta W.$$
By the just obtained inequality and by \eqref{size-160122-5-5} we get
\eqref{size-160422-1-6}.

Now we consider the case in which \eqref{size-160122-4b} holds. By
\eqref{size-est-11a} we get $\delta W\leq0$ and also
\begin{equation}
	\label{size-160122-2-6} |\alpha(u_0,u_0)|\le-\delta W.
\end{equation}
Now we estimate $|\alpha(u_0,u_0)|$ from below. By
\eqref{size-est-11c} we obtain, for $\varepsilon>0$ to be choosen,
\begin{eqnarray}
	\label{size-160122-3-6}
	-\delta W&=&\alpha(u_0,u_1)\leq\left(-\alpha(u_0,u_0)\right)^{1/2}\left(-\alpha(u_1,u_1)\right)^{1/2}\leq\nonumber\\[2mm]
	&\leq&\frac{\varepsilon}{2}(-\alpha(u_1,u_1))+\frac{1}{2\varepsilon}(-\alpha(u_0,u_0)).
\end{eqnarray}
By \eqref{size-est-11b} we have
\begin{equation}
	\label{size-160122-4-6} -\alpha(u_1,u_1)=a_0(u_1-u_0,u_1-u_0)-\delta
	W.
\end{equation}
Moreover, \eqref{size-est-8} gives

\begin{equation*}
	\begin{aligned}
		a_0(u_1-u_0,u_1-u_0)&\leq \lambda_0\left\Vert
		u_1-u_0\right\Vert^2\leq \lambda_0\lambda_1 a_1(u_1-u_0,u_1-u_0).
	\end{aligned}
\end{equation*}
By the just obtained inequality and by \eqref{size-160122-4-6}
we get
$$-\alpha(u_1,u_1)\leq \lambda_0\lambda_1a_1(u_1-u_0,u_1-u_0)-\delta W.$$
The last inequality togheter with \eqref{size-160122-3-6} and
\eqref{size-est-11a},  give (we denote $A=\lambda_0\lambda_1$)
\begin{equation*}
	\begin{aligned}
		-\delta W &\leq
		\frac{\varepsilon}{2}\left[Aa_1(u_1-u_0,u_1-u_0)-\delta W\right]
		+\frac{1}{2\varepsilon}(-\alpha(u_0,u_0))=\\&=
		\frac{\varepsilon}{2}\left[A\left(a_1(u_1-u_0,u_1-u_0)-\alpha(u_0,u_0)\right)+A\alpha(u_0,u_0)-\delta W\right]+\\&
		+\frac{1}{2\varepsilon}(-\alpha(u_0,u_0))=\\&=
		-\frac{\varepsilon}{2}(1+A)\delta
		W+\left(\frac{1}{2\varepsilon}-A\frac{\varepsilon}{2}\right)(-\alpha(u_0,u_0)).
	\end{aligned}
\end{equation*}

Therefore
$$\left(1-\frac{\varepsilon}{2}(1+A)\right)|\delta W|
\leq\frac{1-A\varepsilon^2}{2\varepsilon}|\alpha(u_0,u_0)|.$$ If
$$\varepsilon=\min\left\{\frac{1}{\sqrt{2A}},\frac{1}{1+A}\right\}$$ we have
$$\frac{2\varepsilon\left(1-\frac{\varepsilon}{2}(1+A)\right)}{1-A\varepsilon^2}|\delta
W| \leq|\alpha(u_0,u_0)|.$$  Ultimately, we have
$$C|\delta W|\leq|\alpha(u_0,u_0)|,$$
where $C$ depends on $\lambda_0,\lambda_1$ only. $\blacksquare$

\bigskip

\textbf{Remark 1.} If \eqref{size-160122-4a} holds, condition \eqref{size-est-8} can be weakened by assuming
that $a_0(\cdot,\cdot)$, $a_1(\cdot,\cdot)$ are semidefinite positive. In turn, if case \eqref{size-160122-4b} occurs, it suffices to assume that $a_0(\cdot,\cdot)$, $a_1(\cdot,\cdot)$
are semidefinite positive and satisfy
$$a_0(u,u)\leq C_1 a_1(u,u),\qquad\forall\,u\in H,$$
where $C_1$ is a positive constant. $\blacklozenge$

\bigskip

Now we apply Theorem \ref{teo-size} to inclusion inverse problem. Let
$$H=\left\{v\in H^{1}(\Omega):\quad \int_{\Omega}vdx=0\right\},$$

$$a_1(u,v)=\int_{\Omega}\left(1+(k-1)\chi_D\right)\nabla
u\cdot\nabla v,\quad u,v\in H,$$

$$a_0(u,v)=\int_{\Omega}\nabla
u\cdot\nabla v,\quad u,v\in H.$$

Let us assume that $\phi\in H^{-1/2}(\partial \Omega)$ satisfies
\eqref{size-est-1}. Then Neumann problems \eqref{size-est-2},
\eqref{size-est-6} can be formulated (see Section
\ref{prob-Neumann} and the proof of Theorem
\ref{esistenza-unic-Neum}) as follows
\begin{equation}\label{size-est-12}
	a_1(u,v)=\langle g, v \rangle_{H^{-1/2}, H^{1/2}},\quad \forall v\in
	H,\end{equation}

\begin{equation}\label{size-est-13}
	a_0(u,v)=\langle g, v \rangle_{H^{-1/2}, H^{1/2}},\quad \forall v\in
	H.
\end{equation}
Hence, in our case we have

$$H\ni v\rightarrow \langle F, v \rangle=\langle g, v \rangle_{H^{-1/2},
	H^{1/2}}\in\mathbb{R},$$

$$\alpha(u,v)=(1-k)\int_{D}\nabla u\cdot \nabla v, \quad u,v\in H,$$

$$\delta W=\int_{\partial\Omega}\phi (u_0-u_1),$$
where $u_0$ e $u_1$ are solutions to 
\eqref{size-est-13} and \eqref{size-est-12} respectively.

\medskip

\noindent \textbf{Case $k<1$.} In this case, if $C_0=1$ then
\eqref{size-160122-4a} is satisfied. Hence by
\eqref{size-160422-1-6} we have
\begin{equation}\label{size-est-14}
	\frac{\delta W}{1-k}\leq\int_{D}|\nabla u_0|^2\leq \frac{2\delta
		W}{1-k}.
\end{equation}

\medskip

\noindent \textbf{Case $k>1$.} In this case
\eqref{size-160122-4b} holds true. Hence, by 
\eqref{size-160122-2-7} we get
\begin{equation}\label{size-est-15}
	-\frac{C\delta W}{k-1}\leq\int_{D}|\nabla u_0|^2\leq -\frac{\delta
		W}{k-1},
\end{equation}
where $C$ depends on $k$ only. Inequalities \eqref{size-est-14} and
\eqref{size-est-15} imply \eqref{size-est-7}.

\bigskip

Let us observe that estimates \eqref{size-est-7} can be used to
easily find a size estimate of $D$. For instance, if
$u_0=x_1+c$ in $\Omega$ we have $\nabla u_0=e_1$ and then by
\eqref{size-est-7} we get.

\begin{equation}\label{size-est-16}
	C_2^{-1}|W_0-W|\leq |D|\leq C_1^{-1}|W_0-W|.
\end{equation}
Obviously,
to assign the value of $u_0$ on $\Omega$ is equivalent to assign some
stringent conditions on the "input" current density. In the case
in question, $\phi=e_1\cdot \nu$ and it is not certain that, in practice, one
one can make such a choice. For this reason it is useful to
find the estimates (from above and below) of the measure of $D$ for a generic
nontrivial $\phi$.

In order to examine this issue a little more deeply, let us begin by
observing that  to formulate Neumann problem
\eqref{size-est-2} it is not necessary that $D$ be an open, but it suffices that
\textbf{$D$ be a Lebesgue measurable set} of $\mathbb{R}^n$. If, for instance we know that $$\mbox{dist}(D,\partial \Omega)\geq d>0$$ it is not
difficult to find an estimate from below of $|D|$ by exploiting the
second inequality in \eqref{size-est-7}, i.e. the inequality

\begin{equation}\label{size-est-7-n}
	|\delta W|\leq C_2 \int_{D}\left|\nabla u_0\right|^2.
\end{equation}

Let us examine in which manner we can find an estimate from below of $|D|$.

We have
\begin{equation}\label{size-est-17-n}\int_D|\nabla u_0|^2\leq
	|D|\max_{\overline{D}}|\nabla u_0|^2.\end{equation} Now, let $x_0\in
\overline{D}$ satisfy

\begin{equation}\label{size-est-17}
	|\nabla u_0(x_0)|=\max_{\overline{D}}|\nabla u_0|.
\end{equation}
By the Mean Property for harmonic functions, we have

\begin{equation*}
	\begin{aligned}
		\nabla u_0(x_0)&=\frac{1}{|B_{d/2}(x_0)|}\int_{B_{d/2}(x_0)}\nabla
		u_0(x)dx=\\&=\frac{1}{|B_{d/2}(x_0)|}\int_{\partial
			B_{d/2}(x_0)}u_0(x)\nu dS,
	\end{aligned}
\end{equation*} which implies

\begin{equation}\label{size-est-19}|\nabla u_0(x_0)|\leq \frac{2}{d}\max_{\partial
		B_{d/2}(x_0)}|u_0|.\end{equation} Now, let us estimate  $\max_{\partial
	B_{d/2}(x_0)}|u_0|$ from above. Let $\overline{x}\in \partial B_{d/2}(x_0)$
fulfill

$$|u_0\left(\overline{x}\right)|=\max_{\partial
	B_{d/2}(x_0)}|u_0|,$$ by using again the Mean Property and the Cauchy--Schwarz inequality we get

\begin{equation}\label{size-est-20}
	\begin{aligned}
		|u_0(\overline{x})|&=\left|\frac{1}{|B_{d/4}(\overline{x})|}\int_{B_{d/4}(\overline{x})}u(y)dy\right|\leq
		\frac{|B_{d/4}(\overline{x})|^{1/2}}{|B_{d/4}(\overline{x})|}\left(\int_{B_{d/4}(\overline{x})}|u(y)|^2dy\right)^{1/2}\leq\\&\leq
		\frac{1}{|B_{d/4}(\overline{x})|^{1/2}|}\left\Vert
		u_0\right\Vert_{L^2(\Omega)}\leq C \left\Vert \nabla
		u_0\right\Vert_{L^2(\Omega)},
	\end{aligned}
\end{equation}

\smallskip

\noindent where $C$ depends  by $\Omega$ and $d$ only; in the last inequality of
\eqref{size-est-20} we have applied Theorem \ref{Poincar1}. On the other hand, by \eqref{stima-Neumann-completo} we have
$$\left\Vert \nabla
u_0\right\Vert_{L^2(\Omega)}\leq C\left\Vert
\phi\right\Vert_{H^{-1/2}(\partial\Omega)}.$$ The just obtained inequality and
\eqref{size-est-20} yield

$$\max_{\partial
	B_{d/2}(x_0)}|u_0|=|u_0\left(\overline{x}\right)|\leq C\left\Vert
\phi\right\Vert_{H^{-1/2}(\partial\Omega)}.$$ Now, by this inequality and by
\eqref{size-est-17-n}-- \eqref{size-est-19} we get
\begin{equation}
	\label{size-est-21}\int_D|\nabla u_0|^2\leq C_{*}|D|\left\Vert
	\phi\right\Vert^2_{H^{-1/2}(\partial\Omega)},\end{equation} where
$C_{*}$ is a constant depending on $\Omega$ and $d$ only.
Finally, by \eqref{size-est-7-n} and \eqref{size-est-21} we have

\begin{equation*}
	\frac{|\delta W|}{C_2C_{*}\left\Vert
		\phi\right\Vert^2_{H^{-1/2}(\partial\Omega)}}\leq
	|D|.\end{equation*}
To find an estimate from above of $|D|$  (of course, in terms of
$\delta W$) is definitely more challenging and, as we have
already mentioned in the case where $D$ is an open set, such estimate from above has
inevitably to do with the \textbf{unique continuation property} \index{unique continuation property} of solution to the Laplace equation.

When $D$ is only a Lebesgue measurable set, even
prove that

\begin{equation}
	\label{size-est-22}\delta W=0 \quad \Longrightarrow\quad |D|=0,
\end{equation}
is not trivial. To present here a proof of \eqref{size-est-22} we need the differentiation Lebesgue Theorem \ref{Diff-Leb:13-10-22-1}. In particular such a Theorem
\ref{Diff-Leb:13-10-22-1} implies that if $D$ is a Lebesgue  measurable set, then

\begin{equation}
	\label{size-est-24} \lim_{r\rightarrow 0}\frac{|D\cap B_r(x)|
	}{|B_r(x)|}=1, \quad \mbox{ a.e. } x\in D.
\end{equation}
Let us set \begin{equation}
	\label{size-est-24n}\widetilde{D}=\left\{x\in D:\quad
	\lim_{r\rightarrow 0}\frac{|D\cap B_r(x)|
	}{|B_r(x)|}=1\right\}.\end{equation} The following Proposition holds true
(\cite{DeFi-Go}).
\begin{prop}\label{prop-Lebesgue}
	Let $\Omega$ be an open set of $\mathbb{R}^n$ and let $D$ be a Lebesgue measurable set such that $\overline{D}\subset \Omega$ and
	$|D|>0$. Let $u\in H_{loc}^1(\Omega)$. Let us assume that $u$ satisfies the condition
	\begin{equation}
		\label{size-est-26} u(x)=0,\quad \forall x\in D.
	\end{equation}
	Moreover, let us assume that, in a given point $x_0\in \widetilde{D}$ we have
	
	\begin{equation}
		\label{size-est-25} \int_{B_r(x_0)}|\nabla u|^2dx\leq
		\frac{C}{r^2}\int_{B_{2r}(x_0)}u^2dx,
	\end{equation}
	for every $r>0$ such that
	$\overline{B_{2r}(x_0)}\subset \Omega$, where $C$ is independent of $r$.
	
	Then we have
	\begin{equation}
		\label{size-est-27} \int_{B_r(x_0)}u^2dx\leq
		\mathcal{O}\left(r^k\right), \quad \mbox{ as } r\rightarrow 0,\quad
		\forall k\in \mathbb{N}.
	\end{equation}
\end{prop}

To prove Proposition \ref{prop-Lebesgue} we need the following

\begin{lem}\label{imm-Sobol}
	If $R>0$ and $u\in H^1\left(B_R\right)$ then
	\begin{equation}\label{size-est-28}
		\left(\int_{B_R}|u|^qdx\right)^{\frac{1}{q}}\leq
		\frac{C_{n,q}}{|B_R|^{\frac{1}{2}-\frac{1}{q}}}
		\left(\int_{B_R}\left[R^2|\nabla u|^2+
		u^2\right]dx\right)^{\frac{1}{2}},
	\end{equation}
	where $q$ is an arbitrary number of $(1,+\infty)$ for $n=2$, and it is
equal to $\frac{2n}{n-2}$ for $n\geq 3$. Moreover, $C_{n,q}$ depends
	on $q$ and $n$ only.
\end{lem}
\textbf{Proof of Lemma \ref{imm-Sobol}.} Set
$$v(y)=u(Ry),\quad \forall y\in B_1,$$
it turns out that $v\in H^1\left(B_1\right)$. Now, by the Embedding Sobolev Theorem (Theorem \ref{Sobolev-inequ}) and performing  the
change of variables $y=Rx$, we get
\begin{equation}\label{size-est-29}
	\begin{aligned}
		\left(\int_{B_R}|u|^qdx\right)^{\frac{1}{q}}&=\left(R^n\int_{B_1}|u(Ry)|^qdy\right)^{\frac{1}{q}}=\\&=
		R^{\frac{n}{q}}\left(\int_{B_1}|v(y)|^qdy\right)^{\frac{1}{q}}\leq\\&\leq
		C R^{\frac{n}{q}}\left(\int_{B_1}\left[|\nabla
		v|^2+|v|^2\right]dy\right)^{\frac{1}{2}}=\\&= C
		R^{\frac{n}{q}}\left(\int_{B_1}\left[R^2|(\nabla
		u)(Ry)|^2+|u(Ry)|^2\right]dy\right)^{\frac{1}{2}}=\\&= C
		R^{\frac{n}{q}}\left(R^{-n}\int_{B_R}\left[R^2|\nabla u|^2+|
		u|^2\right]dx\right)^{\frac{1}{2}}=\\&=
		\frac{C\omega_n^{\frac{1}{2}-\frac{1}{q}}}{|B_R|^{\frac{1}{2}-\frac{1}{q}}}\left(\int_{B_R}\left[R^2|\nabla
		u|^2+| u|^2\right]dx\right)^{\frac{1}{2}}.
	\end{aligned}
\end{equation}
$\blacksquare$

\bigskip

\textbf{Proof of Proposition \ref{prop-Lebesgue}.}
Let us denote by $\rho=$ dist$(x_0,\partial \Omega)$. Since
$x_0\in\widetilde{D}$, we have  that for any $\varepsilon>0$ there exists
$r_{\varepsilon}<2\rho$ such that
\begin{equation*} \frac{| B_r(x_0)\setminus D|
	}{|B_r(x_0)|}<\varepsilon, \quad \forall r\in
	(0,r_{\varepsilon}].
\end{equation*}
By the above inequality and by Lemma \ref{imm-Sobol} we have

\begin{equation*}
	\begin{aligned}
		\int_{B_r(x_0)}u^2dx&=\int_{B_r(x_0)\setminus D}u^2dx\leq
		|B_r(x_0)\setminus D|^{1-\frac{2}{q}}\left(\int_{B_r(x_0)\setminus
			D}|u|^qdx\right)^{\frac{2}{q}}\leq
		\\& \leq |B_r(x_0)\setminus
		D|^{1-\frac{2}{q}}\left(\int_{B_r(x_0)}|u|^qdx\right)^{\frac{2}{q}}\leq\\&\leq
		C\left(\frac{|B_r(x_0)\setminus
			D|}{|B_r(x_0)|}\right)^{1-\frac{2}{q}}\int_{B_r(x_0)}\left[r^2|\nabla
		u|^2+u^2\right]dx\leq\\&\leq
		C\varepsilon^{1-\frac{2}{q}}\int_{B_r(x_0)}\left[r^2|\nabla u|^2+
		u^2\right]dx,
	\end{aligned}
\end{equation*}
where $q$ is an arbitrary number of $(1,+\infty)$ for $n=2$, and it is
equal to $\frac{2n}{n-2}$ for $n\geq 3$. In addition, $C$ depends on $q$ and $n$ only.

By the just obtained inequality and by \eqref{size-est-25} we get

\begin{equation}\label{size-est-30}
	\int_{B_r(x_0)}u^2dx\leq C
	\varepsilon^{1-\frac{2}{q}}\int_{B_{2r}(x_0)} u^2dx,\quad
	r\in(0,r_{\varepsilon}],
\end{equation}
where $C$ is  independent on $r$ and $\varepsilon$.

Let $k\in \mathbb{N}$ be arbitrary and let $\varepsilon>0$ satisfy $$C\varepsilon^{1-\frac{2}{q}}=2^{-k}.$$ 
We further let us denote by $r_k$ the
value of $r_{\varepsilon}$ that corresponds to this choice of $\varepsilon$. Let us intruduce the following function
$$f(r)=\int_{B_r(x_0)}u^2dx,\quad r\in(0,2r_k].$$ Then
\eqref{size-est-30} can be written as
\begin{equation}\label{size-est-31}
	f(r)\leq 2^{-k}f(2r),\quad r\in(0,r_k].
\end{equation}

Now, for any $0<r<r_k$, let $m\in \mathbb{N}$ satisfy
\begin{equation}\label{size-est-32}
	2^{-m}r_k\leq r<2^{1-m}r_k.
\end{equation}
Iteration of \eqref{size-est-31} gives

$$f(r)\leq 2^{-k}f(2r)\leq \cdots \leq 2^{-km}f\left(2^mr\right)\leq 2^{-km}f\left(2r_k\right).$$
On the other hand \eqref{size-est-32} implies
$$2^{-m}\leq \frac{r}{r_k},$$ hence
$$f(r)\leq \left(\frac{r}{r_k}\right)^{k}f\left(2r_k\right),$$
which gives \eqref{size-est-27}. $\blacksquare$

\bigskip

\textbf{Remark 2.} By the proof of Proposition \ref{prop-Lebesgue} it is clear that assumption \eqref{size-est-25} can be replaced by the assumption that there is $p>2$ such that

\begin{equation}\label{7-12-1921}
\left(\dashint_{B_r}|u|^pdx\right)^{1/p} \leq C \left(\dashint_{B_{2r}}|u|^2dx\right)^{1/2}.
\end{equation}
See also Section \ref{Ap:7-12-22}. $\blacklozenge$

\begin{theo}\label{harmonic-UCP}
	Let $\Omega$ be a connected bounded open set of $\mathbb{R}^n$ an let
	$$D\subset \Omega$$ be a Lebesgue measurable set of positive measure.
	Let $u$ be a harmonic function in $\Omega$ which vanishes on $D$.  Then
	
	\begin{equation}\label{size-est-33}
		u\equiv 0.
	\end{equation}
\end{theo}
\textbf{Remark 3.} As will become clear from the
proof, the boundedness assumption of $\Omega$ is not essential: we have introduced it solely for the purpose of easing the
proof. We leave the simple extension to the reader. 

$\spadesuit$

\medskip

\textbf{Proof of Theorem \ref{harmonic-UCP}.} Since 
\begin{equation}\label{size-est-34}
	D=\bigcup_{j=1}^{\infty} \left(D\cap \Omega_j\right),\end{equation}
where
$$\Omega_j=\left\{x\in \Omega:\quad \mbox{dist}(x,\partial
\Omega)>1/j\right\},$$ we have $$0<|D|=\lim_{j\rightarrow
	\infty}|D\cap \Omega_j|,$$ hence $|D\cap \Omega_j|>0$ for
$j$ large enough, in addition we have $\overline{D\cap \Omega_j}\subset
\Omega$. Hence, provided to replace $D\cap \Omega_j,$ to $D$,
we may always assume  $\overline{D}\subset\Omega$.

Let us apply Proposition \ref{prop-Lebesgue}. Inequality
\eqref{size-est-25} is nothing more than the Caccioppoli inequality proved in Theorem \ref{dis-Caccioppoli}. Hence
we have, for $x_0\in \widetilde{D}$ ($\widetilde{D}$ is defined
by \eqref{size-est-24n}),

\begin{equation}
	\label{size-est-35} \int_{B_r(x_0)}u^2dx=
	\mathcal{O}\left(r^k\right), \quad \mbox{ as } r\rightarrow 0,\quad
	\forall k\in \mathbb{N}.
\end{equation}
Let now $x\in\Omega$ satisfy
$|x-x_0|<\frac{1}{2}$dist$(x_0,\partial\Omega)$. Set
$r=|x-x_0|$, by the Mean Property and by the Cauchy--Schwarz inequality we have

\begin{equation*}
	\begin{aligned}
		|u(x)|&=\left|\frac{1}{|B_r(x)|}\int_{B_r(x)}u(y)dy\right|\leq
		\frac{1}{|B_r(x)|}\int_{B_r(x)}|u(y)|dy\leq\\&\leq
		\frac{1}{|B_r(x)|}\int_{B_{2r}(x_0)}|u(y)|dy\leq
		\frac{|B_{2r}(x_0)|}{|B_r(x)|}\left(
		\int_{B_{2r}(x_0)}|u(y)|^2dy\right)^{1/2}\leq\\&\leq
		c_nr^{n/2}\left( \int_{B_{2r}(x_0)}|u(y)|^2dy\right)^{1/2},
	\end{aligned}
\end{equation*}
where $c_n$ depends on $n$ only. From what has just been obtained and from
\eqref{size-est-35}, recalling that $r=|x-x_0|$, we have

$$u(x)=\mathcal{O}\left(|x-x_0|^k\right), \quad \mbox{ as } x\rightarrow x_0,\quad
\forall k\in \mathbb{N}.$$ Therefore, as $u$ is an analytic function
the thesis follows. $\blacksquare$

\part{CAUCHY PROBLEM FOR PDEs AND STABILITY ESTIMATES}\label{parte2}

\chapter{The Cauchy problem for the first order PDEs}\label{first-order-cauchy}
\section{Review of ordinary differential equations}
\label{richiami-ode} In this Section we give, without proof, some results
on ordinary differential equations that we will need
later on. For further discussion we refer to \cite{Pi-St-Vi}.

Let $t_0\in \mathbb{R}$, $x_0\in \mathbb{R}^n$,
$x_0=\left(x_{0,1},\cdots, x_{0,n}\right)$, $a>0$, $b>0$. Set
$$Q=\left\{(t,x)\in
\mathbb{R}^{n+1}:\quad |t-t_0|\leq a,\quad \left\vert
x_i-x_{0,i}\right\vert\leq b\right\}$$ and let
$$f:Q\rightarrow \mathbb{R}^n,$$ a \textbf{continuous function}  in $Q$ which is 
\textbf{Lipschitz continuous with respect to the variable $x$}, that is 

\begin{equation}\label{1-R1}
\left\vert f(t,x)-f(t,y)\right\vert\leq L \left\vert x-y\right\vert,
\quad\quad\forall (t,x), (t,y)\in Q.
\end{equation}
Let us consider the following Cauchy problem: determine the function
$x(t)$ differentiable in a neighborhood of $t_0$ and satisfying

\begin{equation}\label{2-R1}
\begin{cases}
\overset{\cdot }{x}(t)=f\left(t,x(t)\right),\\
\\
x(t_0)=x_0,
\end{cases}
\end{equation}
here $\overset{\cdot }{x}$ is the derivative of $x$ w.r.t. $t$ \index{$\overset{\cdot }{x}$}
The following Theorem holds true

\begin{theo}\label{Teor-R1}
Let $f\in C^0(Q)$ satisfy \eqref{1-R1}. Then there exists 
$\delta>0$ and there exists a unique solution $x\in
C^1\left([t_0-\delta,t_0+\delta];\mathbb{R}^n\right)$ to problem
\eqref{2-R1}. Moreover, setting
$$M_i=\max_{Q} \left\vert f_i\right\vert,\quad \quad M=\max_{1\leq i\leq
n}M_i,$$ we can choose $\delta=\min\left\{a,
\frac{b}{M}\right\}$.
\end{theo}

The following Lemma will be very useful
\begin{lem}[\textbf{Gronwall}]\label{Lemma-R2}
	\index{Lemma:@{Lemma:}!- Gronwall@{- Gronwall}}
Let $I$ be an interval of $\mathbb{R}$, $\alpha\in I$ and $c\geq 0$.
Moreover, let $u,v\in C^0(I,\mathbb{R})$ where $v(t)\geq 0$ and $u(t)\geq 0$, for every
$t\in I$.

What follows holds true.

\noindent(i) If

$$v(t)\leq c+\int^{t}_{\alpha} u(s)v(s)ds,\quad\forall t\geq \alpha,$$
then
$$v(t)\leq ce^{\int^{t}_{\alpha} u(s)ds},\quad\forall t\geq \alpha.$$

\medskip

\noindent(ii) If

$$v(t)\leq c+\int_{t}^{\alpha} u(s)v(s)ds,\quad\forall t\leq \alpha$$
then
$$v(t)\leq ce^{\int_{t}^{\alpha} u(s)ds},\quad\forall t\leq \alpha.$$
\end{lem}

\bigskip

The Gronwall Lemma makes it simple to
prove the continuous dependence result of the solution to
\eqref{2-R1} by the data $t_0, x_0, f$. More precisely, we have

\begin{theo}[\textbf{Continuous dependence by the data}]\label{Teor-R2}
	\index{Theorem:@{Theorem:}!- continuous dependence by the data (for ODEs)@{- continuous dependence by the data (for ODEs)}}
Let $f, \widetilde{f}\in C^0(Q)$ satisfy \eqref{1-R1}.
Let $\sigma_1,\sigma_2, \varepsilon$ be positive numbers. Let us suppose that

$$\left\vert t_0-\widetilde{t}_0\right\vert\leq \sigma_1,\quad \left\vert x_0-\widetilde{x}_0\right\vert\leq
\sigma_2,\quad \max_{Q} \left\vert f-\widetilde{f}\right\vert\leq
\varepsilon.$$ Set
$$M=\max_{Q} \left\vert f\right\vert,\quad \widetilde{M}=\max_{Q} \left\vert
\widetilde{f}\right\vert.$$

Then the following fact occurs:

There exists $\sigma_0>0$ depending on $a,b, M, \widetilde{M}$
such that if $\sigma_1, \sigma_2<\sigma_0$, then there is $\delta>0$ and
 $x,\widetilde{x}\in
C^1\left([t_0-\delta,t_0+\delta],\mathbb{R}^n\right)$ that
satisfy what follows:

\begin{equation*}
\begin{cases}
\overset{\cdot }{x}(t)=f\left(t,x(t)\right),\\
\\
x(t_0)=x_0, %
\end{cases}
\end{equation*}
\begin{equation*}
\begin{cases}
\overset{\cdot }{\widetilde{x}}(t)=f\left(t,\widetilde{x}(t)\right),\\
\\
\widetilde{x}\left(\widetilde{t}_0\right)=\widetilde{x}_0 %
\end{cases}
\end{equation*}
and
$$\left\vert x(t)-\widetilde{x}(t)\right\vert\leq
C\left(\sigma_1+\sigma_2+\varepsilon\right),\quad\quad \forall t\in
\left[t_0-\delta,t_0+\delta\right],$$ where $C$ is a constant that depends
on $a,b,L, M, \widetilde{M}$ only.
\end{theo}

\begin{theo}[\textbf{regularity}]\label{Teor-R3}
	\index{Theorem:@{Theorem:}!- regularity (for ODEs)@{- regularity (for ODEs)}}
	Let $k\in\mathbb{N}$ and $f\in C^k(Q)$. Then the solution to Cauchy 
problem \eqref{2-R1} belongs to
$C^{k+1}\left([t_0-\delta,t_0+\delta];\mathbb{R}^n\right)$, where
$\delta=\min\left\{a, \frac{b}{M}\right\}$.
\end{theo}

\bigskip

We now describe the Theorem of \textbf{differentiability of the
	solution} of the Cauchy problem with respect to a parameter and with respect
to the inital values. Let $\lambda_1,\lambda_2\in \mathbb{R}$ be such that  $\lambda_1<\lambda_2$ and let

$$\widetilde{Q}=\left\{(t,x;\lambda)\in
\mathbb{R}^{n+2}:\quad |t-t_0|\leq a,\quad \left\vert
x_i-x_{0,i}\right\vert\leq b,\quad
\lambda_1\leq\lambda\leq\lambda_2\right\}.$$ Moreover, let  $f\in
C^1\left(\widetilde{Q},\mathbb{R}^{n} \right)$. If
$\left(t_0,x_0,\overline{\lambda}\right)\in \widetilde{Q}$ then there exists
a unique solution to the Cauchy problem

\begin{equation}\label{1-R4}
\begin{cases}
\overset{\cdot }{x}(t)=f\left(t,x(t);\overline{\lambda}\right),\\
\\
x(t_0)=x_0.
\end{cases}
\end{equation}
Let us denote by
$$\overline{x}\left(t,t_0,x_0;\overline{\lambda}\right)$$ this solution. It can be proved  (and for this we refer to
\cite[Ch. 1]{Pi-St-Vi}) that $\overline{x}$ is differentiable with respect to all the variables and the derivatives are continuous. In order to calculate the derivatives
$$\frac{\partial \overline{x}}{\partial t_0}, \quad\quad \frac{\partial \overline{x}}{\partial
x_{0,j}}, \quad\quad \frac{\partial \overline{x}}{\partial
\overline{\lambda}},\quad j=1,\cdots,n,$$ we proceed in the
following way: we write the system \eqref{1-R4} in the form
\begin{equation}\label{2extr-R4}
\begin{cases}
\frac{\partial}{\partial t}\overline{x}\left(t,t_0,x_0;\overline{\lambda}\right)=f\left(t,\overline{x}\left(t,t_0,x_0;\overline{\lambda}\right);\overline{\lambda}\right)\\
\\
\overline{x}\left(t_0,t_0,x_0;\overline{\lambda}\right)=x_0,
\end{cases}
\end{equation}
next, we make the derivatives  of both the sides of \eqref{2extr-R4}
obtaining a linear first-order system with Cauchy conditions  in the "new unknowns"

 $$\frac{\partial
\overline{x}}{\partial t_0}, \ \ \frac{\partial \overline{x}}{\partial
x_{0,k}}, \  \ \frac{\partial \overline{x}}{\partial \lambda}.$$

For instance, in case $n=1$, whether we are interested in calutating $\frac{\partial
\overline{x}}{\partial t_0}$, we make the derivatives of both the sides of equation
\eqref{2extr-R4} with respect to $t_0$ and we get

$$\frac{\partial}{\partial t_0}\frac{\partial}{\partial t}\overline{x}\left(t,t_0,x_0;\overline{\lambda}\right)=
\frac{\partial f}{\partial
x}\left(t,\overline{x}\left(t,t_0,x_0;\overline{\lambda}\right);\overline{\lambda}\right)\frac{\partial}{\partial
t_0}\overline{x}\left(t,t_0,x_0;\overline{\lambda}\right)$$ and by the initial datum, making the derivative with respect to  $t_0$, we have

$$\left(\frac{\partial}{\partial
t_0}\overline{x}\left(t,t_0,x_0;\overline{\lambda}\right)+\frac{\partial}{\partial
t}\overline{x}\left(t,t_0,x_0;\overline{\lambda}\right)\right)_{|t=t_0}=0.$$
From which, taking into account \eqref{2extr-R4}, we have
$$\frac{\partial}{\partial
t_0}\overline{x}\left(t_0,t_0,x_0;\overline{\lambda}\right)=-f\left(t_0,x_0;\overline{\lambda}\right).$$
Now, set

\begin{equation}\label{R5}
U\left(t,t_0,x_0;\overline{\lambda}\right)=\frac{\partial}{\partial
t_0}\overline{x}\left(t,t_0,x_0;\overline{\lambda}\right),
\end{equation}

$$A\left(t,t_0,x_0;\overline{\lambda}\right)=\frac{\partial f}{\partial
x}\left(t,\overline{x}\left(t,t_0,x_0;\overline{\lambda}\right);\overline{\lambda}\right)$$
and interchanging the order of derivatives $\frac{\partial}{\partial
t_0}\frac{\partial}{\partial t}\overline{x}$ (in the rigorous proof
it is proved that this step is admissible, compare \cite[Cap.
1]{Pi-St-Vi}) we have
\begin{equation*}
\begin{cases}
\frac{\partial U}{\partial t}=A\left(t,t_0,x_0;\overline{\lambda}\right)U,\\
\\
U\left(t_0,t_0,x_0;\overline{\lambda}\right)=-f\left(t_0,x_0;\overline{\lambda}\right),
\end{cases}
\end{equation*}
that is a Cauchy problem for an ordinary differential equation in the new unkwnon $U$.
When $U$ is determined, also   $\frac{\partial}{\partial
t_0}\overline{x}$ turns out determined by \eqref{R5}. Similarly  we proceed when $n>1$ and for the others derivatives.

Likewise, if $f\in C^k\left(\widetilde{Q},\mathbb{R}^{n} \right)$,
it can be proved that
$\overline{x}\left(\cdot,t_0,x_0;\overline{\lambda}\right)$ has continuous 
derivatives w.r.t. $t_0$, $x_{0,j}$ and
$\overline{\lambda}$ up to order $k$ .

\section{First order linear PDEs}\label{linear-quasilinear}

We begin by giving some definitions. Let $\Omega$ be a connected open set of $\mathbb{R}^n$ and let
$$a:\Omega\rightarrow \mathbb{R}^n,\quad\quad
a(x)=\left(a_1(x),\cdots,a_n(x)\right),$$ be a vector field on
$\Omega$. Let us denote by $P(x,\partial)$ the following linear differential operator
\begin{equation}\label{1-1}
P(x,\partial)=a\cdot \nabla =\sum_{j=1}^n a_j(x)\partial_j,\quad
x\in \Omega.
\end{equation}

Throughout this Chapter we will call the \textbf{symbol of the operator}
$P(x,\partial)$ the following homogeneous polynomial of first degree w.r.t. the variables $\xi_1,\cdots,\xi_n$

\begin{equation}\label{1-2}
P(x,\xi)=a(x)\cdot \xi =\sum_{j=1}^n a_j(x)\xi_j,\quad x\in
\Omega,\quad \xi\in \mathbb{R}^n.
\end{equation}
We say that $\xi\in \mathbb{R}^n\setminus\{0\}$ is a
\textbf{characeristic direction  for $P(x,\partial)$ at the point} \index{characteristic:@{characteristic:}!- direction@{- direction}} $x_0\in \Omega$ if
\begin{equation}\label{1-3-24}
P(x_0,\xi)=0.
\end{equation}
A surface
$$\Gamma=\left\{x\in \Omega: \quad\phi(x)=\phi(x_0)\right\},$$
where $\phi\in C^1(\Omega)$ and
$$\nabla \phi(x_0)\neq 0$$ is said a
\textbf{characteristic surface at the point $x_0\in \Omega$ for
$P(x,\partial)$ } \index{characteristic:@{characteristic:}!- surface@{- surface}}if $\nabla \phi(x_0)$ is a characteristic direction
for $P(x,\partial)$ at the point $x_0$. That is

\begin{equation}\label{1-4-24}
P(x_0,\nabla \phi(x_0))=a\cdot \nabla \phi(x_0) =\sum_{j=1}^n
a_j(x_0)\partial_j\phi(x_0)=0,
\end{equation}
let us note $$P(x,\nabla \phi)=P(x,\partial)\phi.$$

We say that $\Gamma$ is a characteristic surface for $P(x,\partial)$ if

\begin{equation}\label{2-1}
P(x,\nabla \phi(x))=P(x,\partial)\phi(x)=a(x)\cdot \nabla \phi(x)=0,
\quad \forall x\in \Gamma.
\end{equation}
We say that the vectors 
$$ \nu(x_0)=-\frac{\nabla \phi(x_0)}{|\nabla \phi(x_0)|} \ \ \mbox{and} \  -\nu(x_0)=\frac{\nabla \phi(x_0)}{|\nabla \phi(x_0)|}$$
Let us note that the versor $\nu(x_0)$ is directed toward the region 
$\left\{x\in \Omega: \phi(x)<0 \right\}$. We agree to say that $\nu(x_0)$ is the
\textit{unit outward normal} to $\Gamma$ in $x_0$ and $-\nu(x_0)$ is the
\textit{inner outward normal} to $\Gamma$ in $x_0$, respectively.

Notice that $\Gamma$ is a characteristic surface at $x_0$ for  $P(x,\partial)$ if and only if
\begin{equation}\label{2-1-24}
	P(x_0,\nu(x_0))=a(x_0)\cdot \nabla \nu(x_0)=0,
\end{equation}
In other words,
$\Gamma$ is a characteristic surface at $x_0$ if and only if $a(x_0)$
is a tangent vector to $\Gamma$ at $x_0$.

\begin{figure}
	\centering
	\includegraphics[trim={1cm 5cm 0 0},clip, width=8cm]{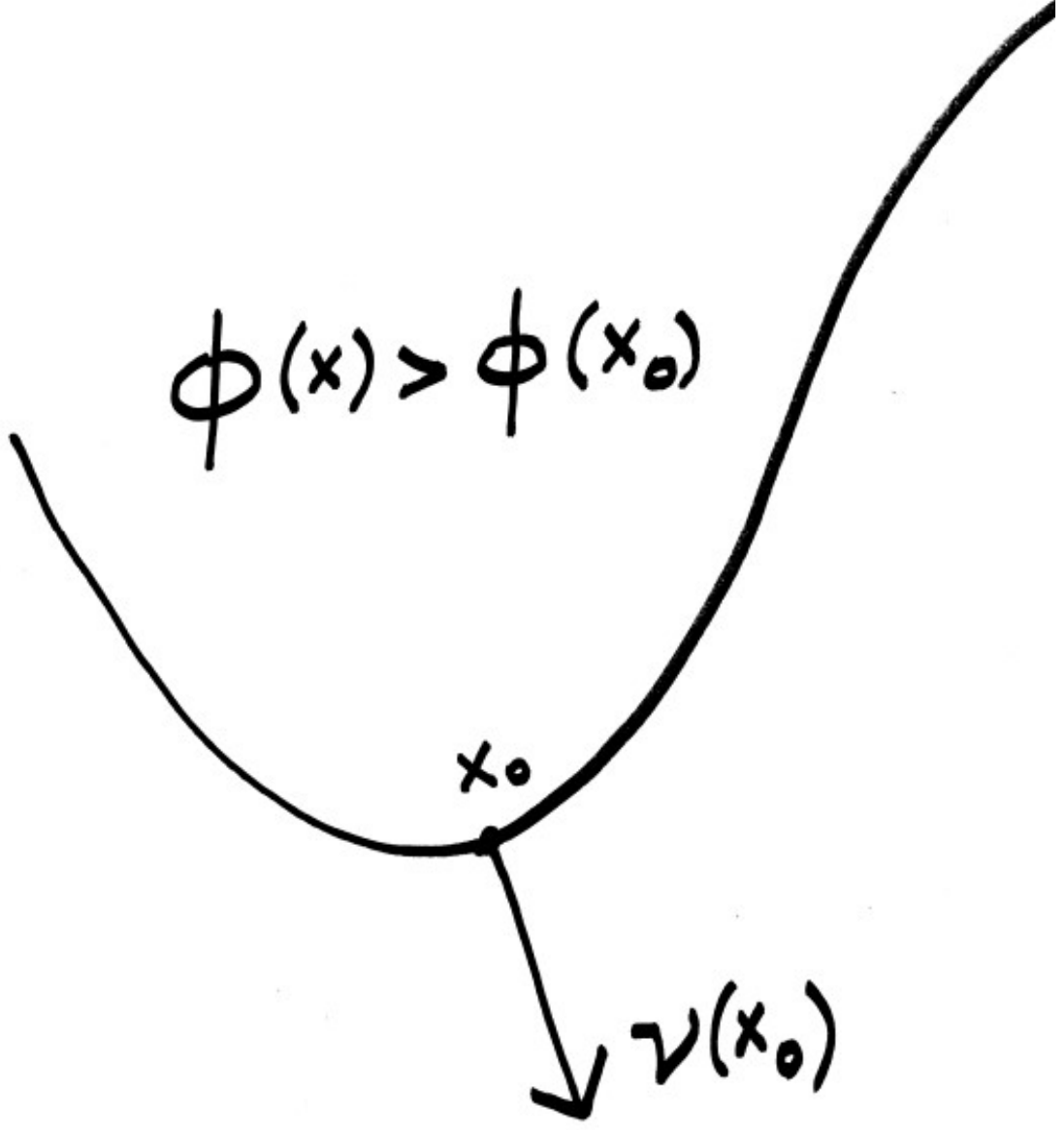}
	\label{50122}
	\caption{}
\end{figure}

Of course, when $n=2$, instead of the characteristic surfaces
we will simply speak of the characteristic lines (or curves).

\bigskip

In this Section we will study the following \textbf{Cauchy problem} for the \textbf{first order linear differential equation}  \index{Cauchy problem:@{Cauchy problem:}!- first order linear differential equations@{- first order linear differential equations}} \index{first order linear differential equation}.  Given the vector field $a(\cdot)$,  and the function $h:\Gamma\rightarrow \mathbb{R}$ determine $u$ such that

\begin{equation}\label{3-2}
\begin{cases}
P(x,\partial)u=c(x)u+f(x),\\
\\
u_{|\Gamma}=h,
\end{cases}
\end{equation}
We will specify the assumptions on $a, c, f, h$ in a while. Let us now
premise some simple example to the investigation of problem \eqref{3-2}. We call $\Gamma$ the \textbf{initial surface} and $h$ the \textbf{initial datum} of Cauchy problem \eqref{3-2}.\index{intital surface of Cauchy problem} \index{intial datum of Cauchy problem}

\bigskip

\textbf{Example 1.}

Let $\Omega=\mathbb{R}^2$ and let

$$P(x,y,\partial)=\partial_x+\partial_y.$$
Let us consider the equation

\begin{equation}\label{1-3}
\partial_xu+\partial_yu=0,\quad\quad\mbox{in }\mbox{ }\mathbb{R}^2.
\end{equation}
We can easily determine all the $C^1\left(\mathbb{R}^2\right)$ solutions
 to equation \eqref{1-3}. Actually,
setting
$$\mu=(1,1),$$  we can write \eqref{1-3}
\begin{equation}\label{2-3}
\frac{\partial u}{\partial\mu}=0,\quad\quad\mbox{in  }\mbox{
}\mathbb{R}^2,
\end{equation}
where $\frac{\partial u}{\partial\mu}$ denotes the derivatives of $u$
w.r.t. direction $\mu$ (defined in Section \ref{Notazioni}).

It is clear that the functions $u \in C^1\left(\mathbb{R}^2\right)$ which satisfy \eqref{2-3} are all and only the functions
constant on the lines parallel to the vector $\mu$. Hence, for any fixed
$x_0$, we have

$$u(x_0+t,t)=F(x_0),\quad t\in\mathbb{R}$$ from which by eliminating $t$,
we have
\begin{equation}\label{3-3}
u(x,y)=F(x-y).
\end{equation}
Therefore it suffices to assume that $F\in C^1(\mathbb{R})$ for obtaining, by
\eqref{3-3}, all the solutions to \eqref{1-3}.

Having \eqref{3-3} available, the study of the Cauchy problem for
equation \eqref{1-3} is quite simple, and here
we take the opportunity to highlight some important facts.

Let us consider the following Cauchy problem 

\begin{equation}\label{1-4}
\begin{cases}
u_x+u_y=0,\\
\\
u(x,0)=h(x),
\end{cases}
\end{equation}
where $h\in C^1(\mathbb{R})$. By \eqref{3-3}, taking into account of the initial datum in \eqref{1-4}, we have

$$h(x)=u(x,0)=F(x),\quad x\in \mathbb{R}.$$
Therefore, the unique solution of problem \eqref{1-4} is given by
\begin{equation}\label{2-4}
u(x,y)=h(x-y),\quad (x,y)\in \mathbb{R}^2.
\end{equation}

Let us now consider a somewhat more general situation and let us assume
that the Cauchy datum is assigned on a regular curve $\Gamma$
of parametric equations

\begin{equation}\label{3-4}
x=\overline{x}(\tau),\quad\quad y=\overline{y}(\tau), \quad \tau\in
I,
\end{equation}
where $I$ is an interval. Let us examine what happens when $\Gamma$
is a \textbf{characteristic line}. We therefore consider the
problem

\begin{equation}\label{4-4}
\begin{cases}
u_x+u_y=0,\\
\\
u(\tau,\tau)=h(\tau), \quad \tau\in \mathbb{R}.
\end{cases}
\end{equation}
We immediately realize that if $h$ is not constant, the
problem \eqref{4-4} has no solutions: we had, indeed, already
observed that every  $C^1\left(\mathbb{R}^2\right)$ which is a solution
of the equation $u_x+u_y=0$ must be constant on the lines
parallel to the vector $\mu=(1,1)$ and therefore, in particular, they must
be constant on the line
$$\{(\tau,\tau):\tau\in \mathbb{R}\}.$$
Moreover, if $h$ is constant then Cauchy problem
\eqref{4-4} has \textbf{infinite solutions} as we easily, if $h\equiv 0$
then every solutions to \eqref{4-4} is given by 
$$u(x,y)=F(x-y),\quad\mbox{with } F(0)=0.$$

\begin{figure}
	\centering
	\includegraphics[trim={0 5cm 5cm 5cm},clip, width=8cm]{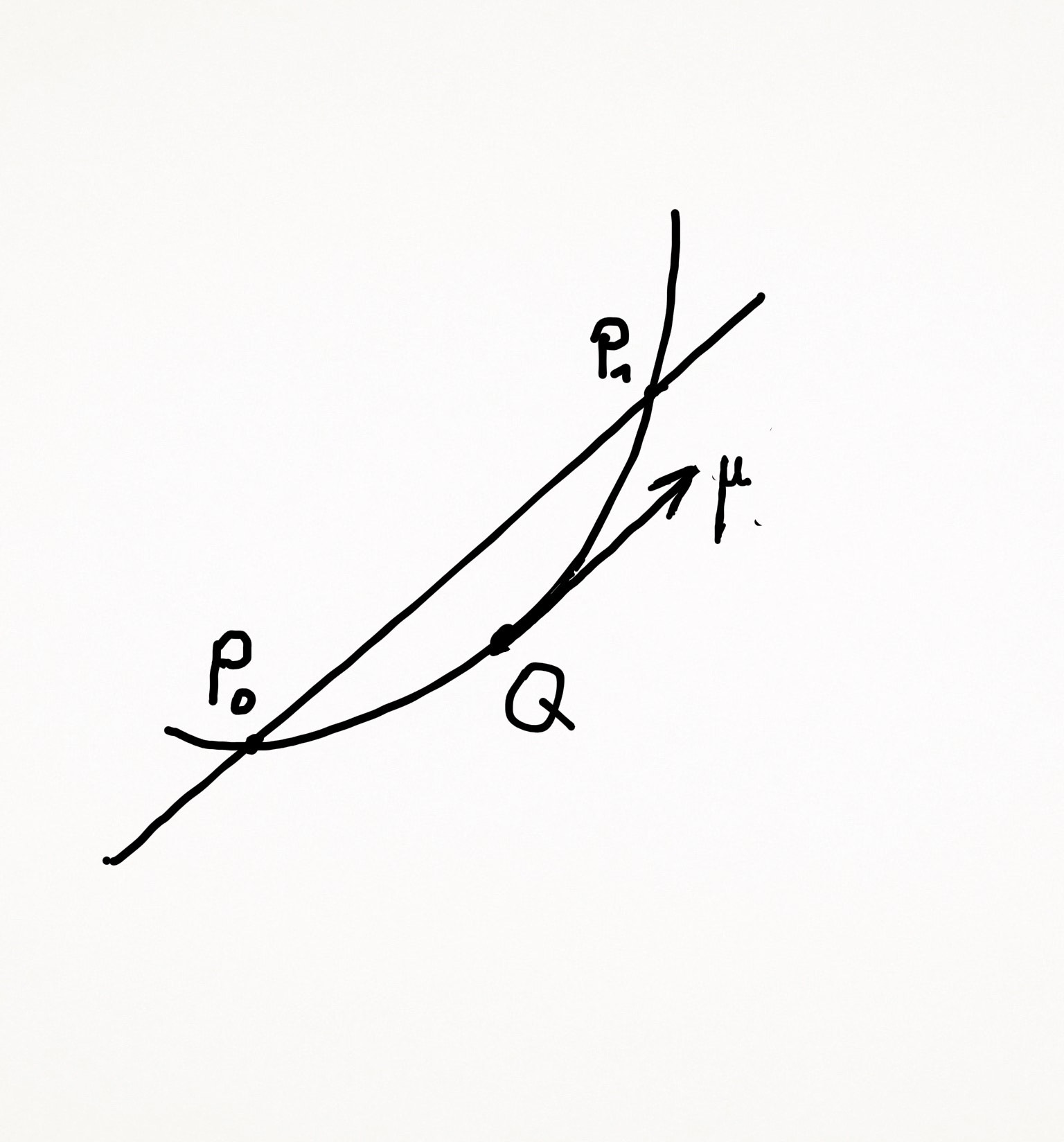}
	\label{22122-1}
	\caption{}
\end{figure}

More generally, if we have to face the Cauchy problem

\begin{equation}\label{1-5}
\begin{cases}
u_x+u_y=0,\\
\\
u\left(\overline{x}(\tau),\overline{y}(\tau)\right)=h(\tau), \quad
\tau\in \mathbb{R}
\end{cases}
\end{equation}
and if a characteristic line intersects $\Gamma$ at two
distinct, say $P_0=\left(\overline{x}(\tau_0),\overline{y}(\tau_0)\right)$ and
$P_1=\left(\overline{x}(\tau_1),\overline{y}(\tau_1)\right)$ where
$\tau_0\neq \tau_1$, then in order that problem \eqref{1-5}
has solution, it is necessary that $h$ satisfies the condition
$$h(\tau_0)=h(\tau_1)$$ and this imposes, in turn, some restrictions on the datum
$h$ itself. In other words, \textbf{not every initial value is admissible
for Cauchy problem \eqref{1-5}}. This is because the values of
$u$ on $\Gamma$ are determined by the values of $u$ on a portion
smaller of $\Gamma$ itself. Keep in mind that in the
situation we have just considered, between the points $P_0$ and $P_1$
there must be a point $Q\in \Gamma$ at which the direction
characteristic is tangent to $\Gamma$. That is, $\Gamma$ 
characteristic line w.r.t. the operator $\partial_x+\partial_y$ at the point
$Q$ (Figure \ref{22122-1}). 

Let us further illustrate what has just been said. Let it be, then, (Figure $5.3$)

$$\Gamma=\left\{\left(\tau,\frac{\tau^2}{2}\right):\quad
0<\tau<2\right\}$$ and let us consider the following Cauchy problem

\begin{equation}\label{1-6}
\begin{cases}
u_x+u_y=0,\\
\\
u\left(\tau,\frac{\tau^2}{2}\right)=h(\tau), \quad \tau\in (0,2).
\end{cases}
\end{equation}

\bigskip

\begin{figure}
	\centering
	\includegraphics[trim={0 0 0 0},clip, width=8cm]{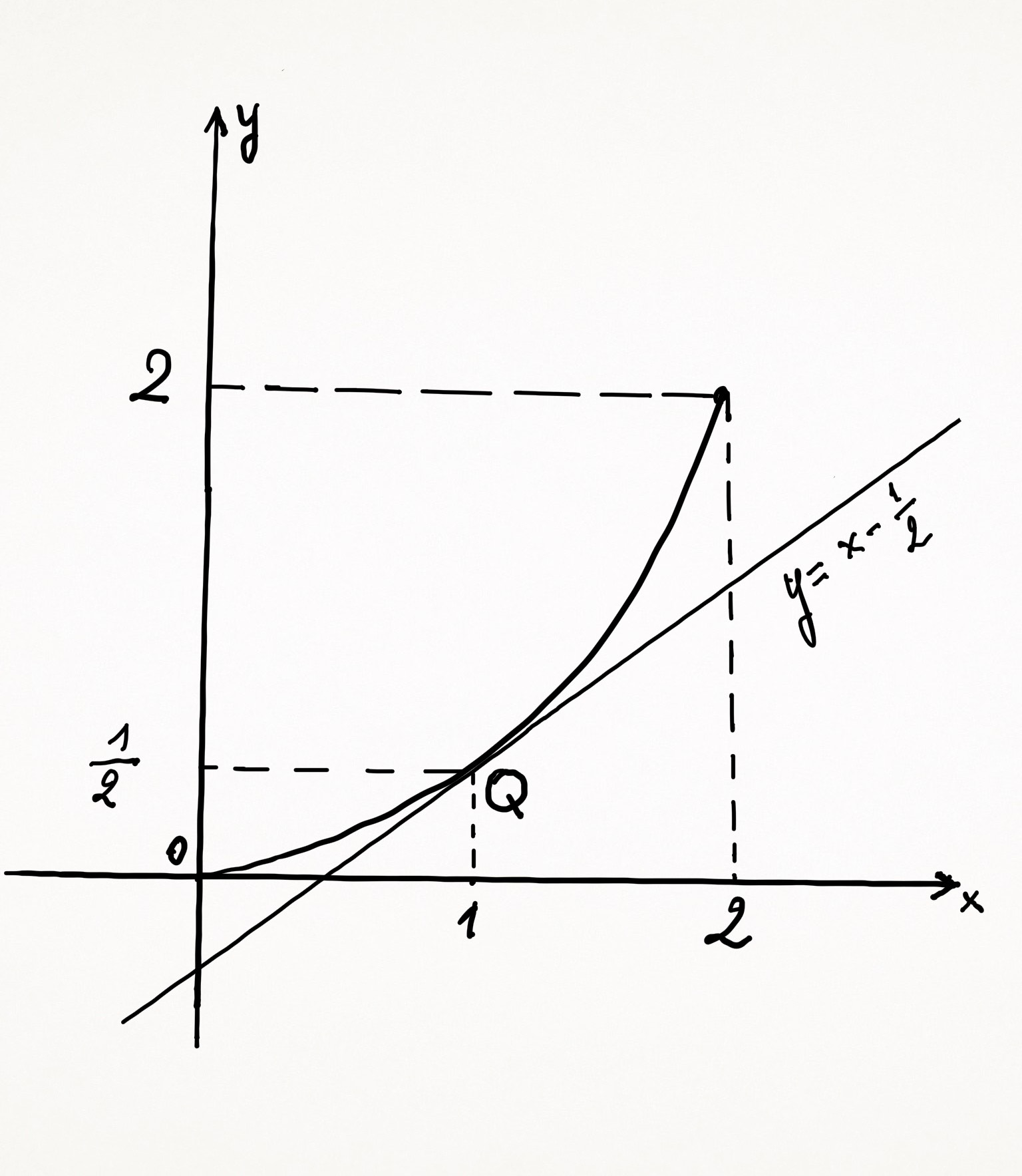}
	\label{figurap6}
	\caption{}
\end{figure}

\bigskip

Notice that $\Gamma$ is a characteristic line w.r.t. the operator
$\partial_x+\partial_y$ at the point $Q=\left(1,\frac{1}{2}\right)$.
Moreover, we check what follows.

\medskip

\noindent \textbf{(i)} There exists a solution $u\in
C^1\left(\mathbb{R}^2\right)$ to problem \eqref{1-6} if and only if
$h\in C^1((0,2))$ and
\begin{equation}\label{2-6}
h(x)=h(2-x), \quad,\quad\forall x\in (0,2).
\end{equation}

\medskip

\noindent \textbf{(ii)} If $h\in C^1(0,2)$ satisfies \eqref{2-6},
then for any $0<r<1$, there exist infinite solutions to the Cauchy problem

\begin{equation}\label{extra0-7}
\begin{cases}
u_x+u_y=0,\quad \mbox{in } B_r(Q),\\
\\
u_{|\Gamma\cap B_r(Q)}=h.
\end{cases}
\end{equation}

\bigskip

\textbf{Let us check (i).} By \eqref{3-3} we have $u(x,y)=F(x-y)$. Consequently, in order to
$$u\left(x,\frac{x^2}{2}\right)=h(x),$$ we need to have 
$$ h(x)=F\left(x-\frac{x^2}{2}\right)=F\left(\frac{1}{2}-
\frac{1}{2}(1-x)^2\right), \quad\quad\forall x\in (0,2),$$ from which \eqref{2-6} follows.

\bigskip

\textbf{Let us check (ii).} By the linearity of  problem
\eqref{extra0-7},  we can choose $h\equiv 0$. Moreover, let $g\in
C^1\left( (0,2)\right)$ be an arbitrary function which satisfies
3
$$g\left(\frac{1}{2}\right)=g'\left(\frac{1}{2}\right)=0.$$ Then it is easily checked that the functions

\begin{equation}\label{extra0-8}
u(x,y)=
\begin{cases}
0,\quad\quad \mbox{in } B_r(Q)\cap\left\{y\geq x-\frac{1}{2}\right\},\\
\\
g(x-y),\quad \mbox{in } B_r(Q)\cap\left\{y< x-\frac{1}{2}\right\},
\end{cases}
\end{equation}
are all solutions of  Cauchy problem (Figure $5.4$)

\begin{equation}\label{extra0-7-24}
\begin{cases}
u_x+u_y=0,\quad \mbox{in } B_r(Q),\\
\\
u_{|\Gamma\cap B_r(Q)}=0.
\end{cases}
\end{equation}

\begin{figure}
	\centering
	\includegraphics[trim={0 0 0 0},clip, width=8cm]{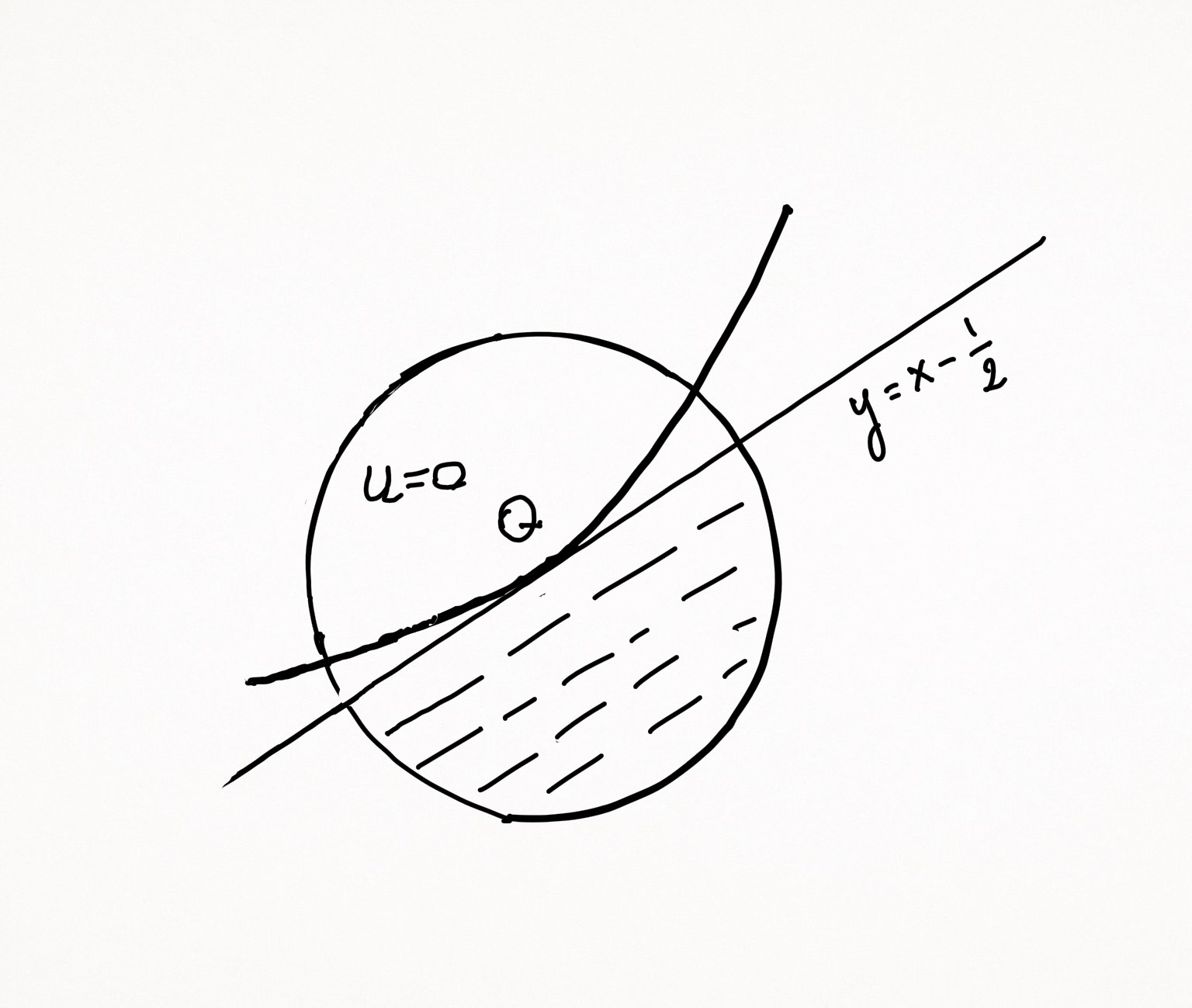}
	\label{figura-p7}
	\caption{}
\end{figure}

\section{The method of characteristics -- the linear case}\label{metodo-caratt}

Let $\Omega$ be a connected open set of $\mathbb{R}^n$ and let $a\in
C^1\left(\Omega, \mathbb{R}^n\right),$
$$a(x)=\left(a_1(x),\cdots,a_n(x)\right).$$ Let us consider the operator
\begin{equation}\label{1-8}
P(x,\partial)=\sum_{j=1}^n a_j(x)\partial_j,\quad x\in \Omega.
\end{equation}
We will call \textbf{characteristic line} \index{characteristic:@{characteristic:}!- equations@{- equations}}of $P(x,\partial)$
any solution to the system of ordinary differential equations --
\textbf{characteristic equations}

\begin{equation}\label{2-8}
\frac{d X(t)}{dt}=a\left(X(t)\right).
\end{equation}
Let $u\in C^1(\Omega)$. Obxerve that, if we set
$$z(t)=u\left(X(t)\right),$$
then we have

\begin{equation}\label{1-9}
\frac{d z(t)}{dt}=\frac{d X(t)}{dt}\cdot \left(\nabla
u\right)(X(t))=a\left(X(t)\right)\cdot \left(\nabla u\right)(X(t)).
\end{equation}
This simple relationship is the starting point of the
method of characteristics by which we will tackle and
we will solve the Cauchy problem
\begin{equation}\label{2-9}
\begin{cases}
P(x,\partial)u=c(x)u+f(x),\\
\\
u_{|\Gamma}=h,
\end{cases}
\end{equation}
where $\Gamma$ is a portion of regular surface of parametric equations

\begin{equation}\label{3-9}
x=\overline{x}(y), \quad\mbox{ } y\in B'_1,
\end{equation}
$\overline{x}\in C^1\left(B'_1\right)$ and

\begin{equation}\label{4-9}
\mbox{Rank}\left(\frac{\partial \overline{x}}{\partial
y}(y)\right)=\mbox{Rank}\left(
\begin{array}{ccc}
\partial_{y_1}\overline{x}_1&\cdots  & \partial_{y_{n-1}}\overline{x}_1 \\
\vdots& \cdots & \vdots \\
\partial_{y_1}\overline{x}_n&  \cdots&
\partial_{y_{n-1}}\overline{x}_n
\end{array}%
\right)=n-1.
\end{equation}
Moreover, we assume
\begin{equation}\label{5-9}
c\in C^1(\Omega),\quad\mbox{ }\quad f\in C^1(\Omega).
\end{equation}
The initial condition $u_{\Gamma}=h$ is expressed by
\begin{equation}\label{6-9}
u\left(\overline{x}(y)\right)=h\left(\overline{x}(y)\right):=\overline{h}(y),\quad\mbox{
} \forall y\in B'_1,
\end{equation}
where

\begin{equation}\label{7-9}
\overline{h}\in C^1\left(B'_1\right).
\end{equation}

The \textbf{method of characteristics} \index{method of characteristics}consists of constructing a
local change of coordinates of class $C^1$,

\begin{equation}\label{1-10} (-\delta,\delta)\times B'_{r}\ni
(t,y)\rightarrow X(t,y)\in\mathbb{R}^n,
\end{equation}
for suitable $\delta>0$ and $r\in (0,1)$. Where $X$ has the following properties 

\begin{equation}\label{2-10}
X(0,y)=\overline{x}(y),\quad\forall  y\in B'_{r}
\end{equation}
and
\begin{equation}\label{3-10}
\partial_t X(t,y)=a(X(t,y)),\quad\forall (t,y)\in (-\delta,\delta)\times
B'_{r}.
\end{equation}
In this way, setting

$$z(t,y)=u(X(t,y)),\quad C(t,y)=c(X(t,y)),\quad F(t,y)=f(X(t,y))$$
and taking into account \eqref{1-9}, Cauchy problem \eqref{2-9}
becomes

\begin{equation}\label{4-X}
\begin{cases}
\partial_tz(t,y)=C(t,y)z(t,y)+F(t,y),\\
\\
z(0,y)=\overline{h}(y).
\end{cases}
\end{equation}
By the assumptions made on $c$, $f$, and $h$, see Section
\ref{richiami-ode}, it turns out that $z\in C^1\left((-\delta,\delta)\times
B'_{r}\right)$ and, as we will see,  that the function (under appropriate
conditions)

\begin{equation}\label{5-10}
u(x)=z\left(X^{-1}(x)\right)
\end{equation}
is a solution to \eqref{2-9}.

Let us begin to see under what conditions, the transformation
defined by \eqref{1-10}--\eqref{3-10} is a diffeomorphism in a
neighborhood of $0\in \mathbb{R}^n$.
Since, $a\in C^1(\Omega)$ and 
$\overline{x}\in C^1\left(B'_1\right)$ we have (for the
differentiability w. r, t.  the parameters, see Section \ref{richiami-ode})
$$X\in C^1\left(J\times B'_1\right),$$
where $J$ is a suitable neighborhood of $0$.

Moreover

\begin{equation*}
\begin{aligned}
\frac{\partial X}{\partial (t,y)}(0,0)&=\underset{\mbox{column vectors} }{\left(\underbrace{\partial_t
X(0,0),\partial_{y_1}X(0,0),\cdots,\partial_{y_{n-1}}X(0,0)}
\right)}=\\&\\& =\left(a(x_0),\partial_{y_1}\overline{x}(0),\cdots,
\partial_{y_{n-1}}\overline{x}(0) \right)
\end{aligned}
\end{equation*}
and, since

$$\partial_t
X(0,0)=a(X(0,0))=a(x_0),$$

$$\mbox{Rank}\left(a(x_0),\partial_{y_1}\overline{x}(0),\cdots,
\partial_{y_{n-1}}\overline{x}(0) \right)=n-1,$$ we have that the following conditions are equivalent 
\begin{equation}\label{etra-11}
\mbox{Rank}\frac{\partial X}{\partial (t,y)}(0,0)=n\end{equation} and

\begin{equation}\label{etra0-11}
a(x_0)\notin\left\langle\partial_{y_1}\overline{x}(0),\cdots,
\partial_{y_{n-1}}\overline{x}(0)\right\rangle,\end{equation} where $\langle
v_1,v_2,\cdots,v_{n-1}\rangle$ is the vector space  generated
by $v_1,v_2,\cdots,v_{n-1}$. Condition \eqref{etra0-11} is, in turn, equivalent to the condition
 that $a(x_0)$ \textbf{is not tangent to
$\Gamma$ in} $x_0$.

All in all, \textbf{if $\Gamma$ is noncharacteristic in $x_0$ for operator \eqref{1-8}} then there exists $\delta>0$ and
$r\in(0,1)$ such that $X$ is a diffeomorphism in
$(-\delta,\delta)\times B'_{r}$.

Now, we denote by $\mathcal{U}_{x_0}=X\left((-\delta,\delta)\times
B'_{r}\right)$ and by

\begin{equation}\label{1-12}
\Psi(x)=X^{-1}(x),\quad x\in \mathcal{U}_{x_0}.
\end{equation}
Let us check that  $$u(x)=z(\Psi(x))$$ solves Cauchy problem \eqref{2-9}.

Regarding the \textbf{initial condition}, by \eqref{2-10} we have
immediately

$$\overline{x}(y)=X(0,y),\quad\quad \forall y\in B'_{r}.$$ Hence, by
\eqref{1-12} and recalling \eqref{4-X}, we get

\begin{equation}\label{2-12}
u\left(\overline{x}(y)\right)=z\left(\Psi(X(0,y))\right)=z(0,y)=\overline{h}(y),\quad\quad
\forall y\in B'_{r}.
\end{equation}
Concerning \textbf{the equation}

$$\sum_{j=1}^na_j(x)\partial_ju=c(x)u+f(x),$$ recall that by
\eqref{1-12} we have

\begin{equation}\label{3-12}
\left(\frac{\partial\Psi(x)}{\partial x}\right)\left(\frac{\partial
X\left(\Psi(x)\right)}{\partial (t,y)}\right)=I_n,
\end{equation}
where $I_n$ is the identity matrix $n\times n$. 
In particular, considering the first column on the right--hand side and the first
column on the left --hand side of \eqref{3-12}, we have

\begin{equation}\label{1-13}
\begin{cases}
\sum_{j=1}^n\partial_{x_j}\Psi_1\partial_t X_j=1,\\
\\
\sum_{j=1}^n\partial_{x_j}\Psi_k\partial_t X_j=0, \quad
k=2,\cdots,n.
\end{cases}
\end{equation}
Now

\begin{equation*}
\begin{aligned}
\partial_{x_j}u(x)&=\partial_tz(\Psi(x))\partial_{x_j}\Psi_1(x)+\partial_{y_1}z(\Psi(x))\partial_{x_j}\Psi_2(x)+\\&\cdots+\partial_{y_{n-1}}z(\Psi(x))\partial_{x_j}\Psi_{n-1}(x).
\end{aligned}
\end{equation*}
Hence (multiplying by $a_j(x)$ and summing up on $j$)

\begin{equation}\label{2-13}
\begin{aligned}
\sum_{j=1}^na_j(x)\partial_{x_j}u&=\partial_tz(\Psi(x))\sum_{j=1}^n\partial_{x_j}\Psi_1a_j(x)+
\\&+\partial_{y_1}z(\Psi(x))\sum_{j=1}^n\partial_{x_j}\Psi_2a_j(x)+\cdots\\&+\partial_{y_{n-1}}z(\Psi(x))\sum_{j=1}^n\partial_{x_j}\Psi_{n-1}(x)a_j(x).
\end{aligned}
\end{equation}
On the other hand by \eqref{3-10} we know

$$a_j(x)=\left(\partial_t X_j\right)(\Psi(x)).$$
By this equality, by \eqref{1-13} and by the equation in
\eqref{4-X} we have

\begin{equation*}\label{2extra-13}
\sum_{j=1}^na_j(x)\partial_{x_j}u(x)=\partial_tz(\Psi(x))=c(x)u(x)+f(x).
\end{equation*}
Hence $u$ is solution
to Cauchy problem \eqref{2-9}.

Finally, we observe that, by the hypotheses \eqref{5-9} (actually,
it suffices $c\in C^0(\Omega)$), $u$ \textbf{is the unique
	solution} of class $C^1$ to problem \eqref{2-9} in the neighborhood
$\mathcal{U}_{x_0}$. Indeed, if $u_1,u_2\in
C^1\left(\mathcal{U}_{x_0}\right)$ are two solutions then, setting
$$w=u_1-u_2,$$ we have

\begin{equation*}
\begin{cases}
P(x,\partial)w=c(x)w,\\
\\
w_{|\Gamma}=0,
\end{cases}
\end{equation*}
and setting
$$\widetilde{z}(t,y)=w(X(t,y)),$$ by \eqref{4-X} we have
\begin{equation*}\label{4-10}
\begin{cases}
\partial_t\widetilde{z}(t,y)=C(t,y)\widetilde{z}(t,y),\\
\\
\widetilde{z}(0,y)=0.
\end{cases}
\end{equation*}
From which we have $\widetilde{z}=0$ in $(-\delta,\delta)\times B'_{r}$,
therefore $w=0$ in $\mathcal{U}_{x_0}$.

The construction we have illustrated and the local uniqueness hold
for any point of $\Gamma$.

\medskip

Hence we have proved
\begin{theo}\label{esist-unic-24}
Let $a\in C^1\left(\Omega, \mathbb{R}^n\right)$, $c\in
C^1\left(\Omega\right)$ and $f\in C^1\left(\Omega\right)$. Let
$\Gamma$ be a non characteristic  surface of parametric equations $x=\overline{x}(y)$, 
where $\overline{x}\in
C^1\left( B_1'\right)$ and satisfying \eqref{4-9}. Let $h$ be a
function $C^1$ on $\Gamma$ (i.e.   $h\circ
\overline{x}\in C^1\left(B'_1\right)$.

Then there exists a neighborhood  $\mathcal{U}$ of $\Gamma$ such that
there exists unique solution $u$ in $C^1(\mathcal{U})$ to the  Cauchy problem

\begin{equation}\label{ext-unic}
\begin{cases}
\sum_{j=1}^n a_j(x)\partial_ju=c(x)u+f(x),\quad x\in \mathcal{U},\\
\\
u_{|\Gamma}=h.
\end{cases}
\end{equation}
\end{theo}

\smallskip

Given a surface $\Gamma$ in $\mathbb{R}^n$ we call \textbf{domain of
dependence} \index{domain of dependence}of $\Gamma$ with respect to the equation
$$\sum_{j=1}^n a_j(x)\partial_ju=c(x)u,$$ the largest closed set
$D_{\Gamma}$ for which we have

\begin{equation*}
\begin{cases}
\sum_{j=1}^n a_j(x)\partial_ju=c(x)u,\quad x\in D_{\Gamma},\\
\\
u_{|\Gamma}=0,%
\end{cases}\quad\Longrightarrow \quad u=0\quad\mbox{in } D_{\Gamma}.
\end{equation*}

\bigskip

\underline{\textbf{Exercise 1.}} Let $0<r<1$ and
$$\Gamma=\left\{\left(x',-\sqrt{1-|x'|^2}\right):\quad
|x'|<r\right\}.$$  Construct a vector field $a\in
C^1\left(\overline{B_1'}\right)$ such that the domain of
dependence of $\Gamma$ with respect to the equation
$$a(x)\cdot \nabla u=0,$$
contains $\overline{B_1}$. $\clubsuit$

\bigskip

\underline{\textbf{Exercise 2.}} Apply the characteristic method to prove
that the functions $u\in
C^1\left(\mathbb{R}^n\setminus\{0\}\right)$ which satisfy

$$\sum_{j=1}^n x_j\partial_ju=\alpha u,$$ are the homogeneous function
of degree $\alpha$. $\clubsuit$

\bigskip

\underline{\textbf{Exercise 3.}} Let $b$ be a vector of $\mathbb{R}^n$ and let $f\in
C^0\left(\mathbb{R}^{n+1}\right)$. Apply the characteristic method to solve the following Cauchy problem

\begin{equation*}
\begin{cases}
\partial_tu+b\cdot\nabla u=f(x,t),\\
\\
u(x,0)=0.
\end{cases}
\end{equation*}
The equation $\partial_tu+b\cdot\nabla u=f(x,t)$ is known as the
\textbf{transport equation} \index{equation:@{equation:}!- transport@{- transport}}. $\clubsuit$

\section{The method of characteristics -- quasilinear case}\label{metodo-caratt-ql} Let
$J$ be an open interval of $\mathbb{R}$, let $\Omega$ be a connected open set
of $\mathbb{R}^n$. Let $a\in
C^1\left(J\times\Omega,\mathbb{R}^{n}\right)$ and $c\in
C^1\left(J\times\Omega\right)$. The following equation

\begin{equation}\label{1-17}
a(x,u)\cdot\nabla u=c(x,u),\end{equation} is called
a \textbf{first-order quasilinear equation}\index{first order quasilinear equation}. Of course, a linear equations are special case of the
quasilinear equations.

With minor modifications, the characteristics method studied in the
Section \ref{metodo-caratt} can be adapted to handle
equation \eqref{1-17} and the related Cauchy problem. In the case
of equation \eqref{1-17}, the characteristic equation \eqref{1-17}
is the following one

\begin{equation}\label{2-17}
\begin{cases}
\frac{d X}{dt}(t)=a(X(t),z(t)),\\
\\
\frac{d z}{dt}(t)=c(X(t),z(t)).
\end{cases}
\end{equation}
Let us note that in the linear case, the equation $$\frac{d
	z}{dt}(t)=c(X(t),z(t)),$$ is precisely the one satisfied by
$z(t)=u(X(t))$ when $u$ is a solution of the linear equation
$$a(x)\cdot \nabla u=c(x)u+f(x).$$
We continue to call \textbf{characteristic line}, the curve of
parametric equations

\begin{equation}\label{2extra-17}
(X,z)=(X(t),z(t))
\end{equation}
where $X(t),z(t)$ is a solution of the system \eqref{2-17}.
When there is no risk of ambiguity, we will call "characteristic line"\index{characteristic:@{characteristic:}!- line@{- line}} also the projection on $\mathbb{R}^n$ of the line
\eqref{2extra-17}.

Let us consider the Cauchy problem \index{Cauchy problem:@{Cauchy problem:}!- first order quasilinear differential equations@{- first order quasilinear differential equations}}
\begin{equation}\label{1-18}
\begin{cases}
a(x,u)\cdot\nabla u=c(x,u),\\
\\
u_{|\Gamma}=h,
\end{cases}
\end{equation}
where $\Gamma$ is a portion of regular surface of parametric equations

\begin{equation*}
x=\overline{x}(y), \quad\mbox{ } \forall y\in B'_1.
\end{equation*}
To solve \eqref{1-18}, we proceed similarly to what we did in the linear case. Namely, we consider $X(t,y)$ and $z(t,y)$ such that

\begin{equation}\label{2-18}
\begin{cases}
\partial_tX(t,y)=a(X(t,y),z(t,y)),\\
\\
\partial_tz(t,y)=c(X(t,y),z(t,y)),\\
\\
X(0,y)=\overline{x}(y),\\
\\
z(0,y)=h(y)
\end{cases}
\end{equation}
and it can be checked, exactly as in Section \ref{metodo-caratt}
that if
\begin{equation*}
(-\delta,\delta)\times B'_{r}\ni (t,y)\rightarrow
X(t,y)\in\mathbb{R}^n,
\end{equation*}
for some $\delta>0$ and $r\in (0,1)$, is local change of coordinates of $\mathbb{R}^n$, then the function

\begin{equation}\label{1-19}
u(x):=z\left(X^{-1}(x) \right),\end{equation} is a solution
to Cauchy problem \eqref{1-18}. More precisely: setting
$x_0=\overline{x}(0)$ there exists a neighborhood $\mathcal{U}_{x_0}$ such that the function
$u$ defined by \eqref{1-19} satisfies

\begin{equation}\label{2-19}
\begin{cases}
a(x,u)\cdot\nabla u=c(x,u),\quad\mbox{in } \mathcal{U}_{x_0},\\
\\
u_{|\Gamma\cap \mathcal{U}_{x_0}}=h.
\end{cases}
\end{equation}

In order to the map $X$ be a local diffeomorphism, it suffices
to have

\begin{equation}\label{3-19}
\mbox{Rank}\frac{\partial X}{\partial (t,y)}(0,0)=n
\end{equation}
and since

\begin{equation*}
\begin{aligned}
\frac{\partial X}{\partial (t,y)}(0,0)&=\left(\partial_t
X(0,0),\partial_{y_1}X(0,0),\cdots,\partial_{y_{n-1}}X(0,0)
\right)=\\&= \left(a(x_0,
h(x_0)),\partial_{y_1}\overline{x}(0),\cdots,
\partial_{y_{n-1}}\overline{x}(0) \right)
\end{aligned}
\end{equation*}
and
\begin{equation*}
\mbox{Rank}\left(\frac{\partial \overline{x}}{\partial
y}(0)\right)=n-1,
\end{equation*}
condition \eqref{3-19} is equivalent to

\begin{equation*}
a(x_0,h(x_0))\notin\left\langle\partial_{y_1}\overline{x}(0),\cdots,
\partial_{y_{n-1}}\overline{x}(0)\right\rangle,
\end{equation*}
(compare this condition with \eqref{etra0-11}).

\bigskip

Now we briefly consider the issue of continuous dependence by initial datum in
problem \eqref{1-18}.

Let $u_k$, $k=1,2$, satisfy

\begin{equation}\label{1-20}
\begin{cases}
a(x,u_k)\cdot\nabla u_k=c(x,u_k),\\
\\
u_{k_{|\Gamma}}=h_k.
\end{cases}
\end{equation}
We have

$$c(x,u_1)-c(x,u_2)=\left(u_1-u_2\right)\int^1_0\partial_uc\left(x,
u_2(x)+t\left(u_1(x)-u_2(x)\right)\right)dt$$ and, similarly

$$a(x,u_1)\cdot\nabla u_1-a(x,u_2)\cdot\nabla
u_2=\overline{a}(x)\cdot\nabla\left(u_1-u_2\right)+b(x),$$ where
$$\overline{a}(x)= a(x,u_1(x))$$
and
$$b(x)=-\nabla u_2(x)\cdot \int^1_0\partial_ua\left(x,
u_2(x)+t\left(u_1(x)-u_2(x)\right)\right)dt.$$ Set
$$\overline{c}(x)=b(x)+\int^1_0\partial_uc\left(x,
u_2(x)+t\left(u_1(x)-u_2(x)\right)\right)dt,$$
$$w=u_1-u_2,$$
$$\overline{h}=h_1-h_2.$$ By \eqref{1-20} we have

\begin{equation}\label{1-21}
\begin{cases}
\overline{a}(x)\cdot\nabla w=\overline{c}(x)w,\\
\\
w_{|\Gamma}=\overline{h}.
\end{cases}
\end{equation}

Now, if $\overline{h}\equiv 0$ and if at a point $x_0\in \Gamma$, then $a(x_0,h_1(x_0))$ is not tangent to $\Gamma$, then
there exists a neighborhood $\mathcal{U}_{x_0}$ of $x_0$ such that $w\equiv 0$
in $\mathcal{U}_{x_0}$, that is

$$u_1\equiv u_2,\quad\mbox{in } \mathcal{U}_{x_0}.$$
As a matter of fact, the equation
$$\overline{a}(x)\cdot\nabla w=\overline{c}(x)w$$
is linear and Theorem \ref{esist-unic-24} applies.

We leave as an exercise to the reader to prove that if $a(x_0,h_1(x_0))$
is not tangent to $\Gamma$ in $x_0$ then there exists a neighborhood
$\mathcal{V}_{x_0}$ such that

\begin{equation}\label{2-21}
\left\Vert
w\right\Vert_{L^{\infty}\left(\mathcal{V}_{x_0}\right)}\leq K
\left\Vert
\overline{h}\right\Vert_{L^{\infty}\left(\Gamma\cap\mathcal{V}_{x_0}\right)},
\end{equation}
that is

\begin{equation*}
\left\Vert
u_1-u_2\right\Vert_{L^{\infty}\left(\mathcal{V}_{x_0}\right)}\leq K
\left\Vert
h_1-h_2\right\Vert_{L^{\infty}\left(\Gamma\cap\mathcal{V}_{x_0}\right)},
\end{equation*}
where $K$ is a constant which depends on $C^1$ norm of $a$ and 
$c$ and on the (convex) angle between the vector $a(x_0,h_1(x_0))$ and the  unit outward normal to $\Gamma$ in $x_0$.

\bigskip

We conclude this Section by the following

\medskip

\textbf{Example.} Let us consider the following Cauchy problem

\begin{equation}\label{1-22}
\begin{cases}
u_y+uu_x=0,\\
\\
u(x,0)=h(x),\quad x\in \mathbb{R},
\end{cases}
\end{equation}
where $h\in C^1(\mathbb{R})$.

The characteristic equations are given by

\begin{equation}\label{2-22}
\left \{
\begin{array}{c}
\frac{\partial x(t,s)}{\partial t}=z,\\
\\
\frac{\partial y (t,s)}{\partial t}=1,\\
\\
\frac{\partial z(t,s)}{\partial t}=0,
\end{array}%
\right.
\end{equation}
the initial conditions are

\begin{equation}\label{3-22}
x(0,s)=s,\quad\quad y(0,s)=0\quad\quad z(0,s)=h(s).
\end{equation}
By \eqref{2-22} and \eqref{3-22} we have easily

\begin{equation}\label{4-22}
\begin{cases}
x(t,s)=s+th(s),\\
\\
y (t,s)=t,\\
\\
z(t,s)=h(s).
\end{cases}
\end{equation}
By the method shown in this Section, the solution to
\eqref{1-22} is given by a function $u$ such that

\begin{equation}\label{5-22}
u\left(x(t,s),y(t,s)\right)=z(t,s)=h(s).
\end{equation}
To express $u$ in the variables $x$ and $y$ we eliminate $s$ and $t$
from the first two equations of \eqref{4-22}. We have
\begin{figure}\label{figura-p24}
	\centering
	\includegraphics[trim={0 0 0 0},clip, width=12cm]{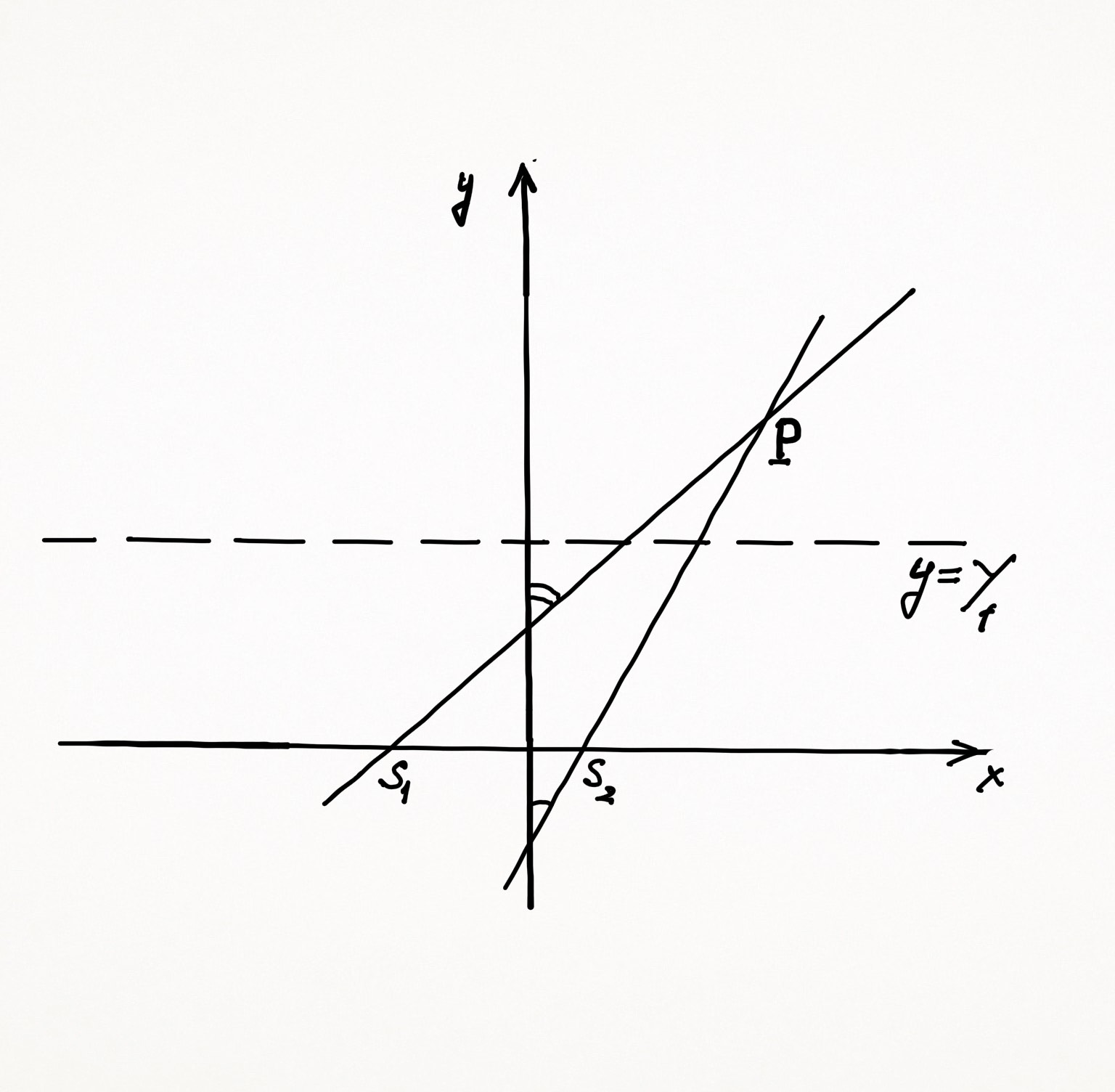}
	\caption{}
\end{figure}

\begin{equation*}
\begin{cases}
x=s+yh(s),\\
\\
t=y.
\end{cases}
\end{equation*}
To obtain $s$ from the first equation it is necessary that $s\rightarrow
s+yh(s)$ be injective, that is, it is necessary that
$$0\neq \frac{d}{ds}\left(s+yh(s)\right)=1+yh'(s).$$
For instance, if $h'>0$,  let us assume 
$$Y_0:=\sup_{s\in\mathbb{R}}-\frac{1}{h'(s)} <0,$$ then
we have that the solution to \eqref{1-22} is defined for all
$y$ such that
$$y>Y_0.$$
If $h'(s)<0$, we assume

$$Y_1:=\inf_{s\in\mathbb{R}}-\frac{1}{h'(s)} >0,$$ then we have the solution to
\eqref{1-22} is defined for all $y$ such that $$y<Y_1.$$

Let us dwell for a while on the latter case and examine what happens above the line $y=Y_1$. Let us come back to \eqref{5-22};
this relation tells us that $u$ is constant on the projection
of the characteristic passing through the point $(s,0)$ and there it is equal to $h(s)$. Let now $s_1, s_2\in \mathbb{R}$ satisfy $s_1<s_2$,
then the straight lines whose equations
are given by
$$x=s_1+yh(s_1)$$ and

$$x=s_2+yh(s_2)$$ intersect at the point

$$P=\left(-\frac{s_2h(s_1)}{h(s_2)-h(s_1)},
-\frac{s_2-s_1}{h(s_2)-h(s_1)}\right),$$ that implies that the
function $u$ \textit{cannot be continuous} in $P$.  Let us observe
that the point $P$ is situated either on the line $y=Y_1$ or above
it since, for an appropriate $\overline{s}\in(s_1,s_2)$, we have (Figure 5.5)

$$-\frac{s_2-s_1}{h(s_2)-h(s_1)}=-\frac{1}{h'(\overline{s})}\geq
Y_1.$$

It is, actually, of some physical interest to include (in an appropriate sense) the discontinuous solution among  the solutions to problem \eqref{1-22} since they correspond to "shock waves" \index{shock waves}. We  refer for
insights to \cite[Ch. 3, Sec. 4]{EV}.

\section{Brief review on the fully nonlinear case}\label{NonlinEq}
In this Section we wish briefly consider to the method of the
characteristics to solve the Cauchy problem for the fully nonlinear equation
$$	F\left(x,u(x),\nabla u(x)\right)=0.$$
\index{Cauchy problem:@{Cauchy problem:}!- first order fully nonlinear equationS@{- first order fully nonlinear equations}}
Namely, the Cauchy problem
\begin{equation}\label{NonlinEq-0}
	\begin{cases}
		F\left(x,u(x),\nabla u(x)\right)=0, \quad \mbox{in } \Omega,\\
		\\
		u_{|\Gamma}=g, 
	\end{cases}
\end{equation}
where $\Omega$ is an open set $\mathbb{R}^n$, $\Gamma$ is a regular
surface of $\mathbb{R}^n$ contained $\Omega$,
\begin{equation}\label{NonlinEq-5}
	F:\Omega\times \mathbb{R}\times \mathbb{R}^n\rightarrow
	\mathbb{R},\end{equation} is a function of class 
$C^{2}(\Omega\times \mathbb{R}\times \mathbb{R}^n)$. The variables of
$F$ are $x\in \Omega$, $z\in \mathbb{R}$ and $p\in \mathbb{R}^n$.
Moreover
$g:\Gamma\rightarrow \mathbb{R}$ is a function defined on $\Gamma$. The assumptions on $\Gamma$, $F$ and
$g$ will be specified in more detail later on. We notice that in the
linear and quasilinear cases investigated, respectively, in Sections
\ref{linear-quasilinear} and \ref{metodo-caratt-ql}, we have
$F(x,z,p)=a(x)\cdot p-c(x)z-f(x)$, $F(x,z,p)=a(x,z)\cdot p-c(x,z)$.
Notice that in the nonlinear case, generally, we cannot expect the uniqueness of the solutions to
Cauchy \eqref{NonlinEq-0}. The following simple example will
help us to understand this fact. Let us consider the Cauchy problem
\begin{equation}\label{NonlinEq-6}
	\begin{cases}
		u^2_{x_1}+u^2_{x_2}=1, \quad \mbox{in } \mathbb{R}^2,\\
		\\
		u(x_1,0)=g(x_1), \quad \mbox{for } x_1\in \mathbb{R}.
	\end{cases}
\end{equation}
We note that all we can say about
$u_{x_2}\left(x_1,0\right)$ is that it satisfies the condition

\begin{equation*}
	g'^2(x_1)+u^2_{x_2}\left(x_1,0\right)=1,
\end{equation*}
which leaves undetermined the sign of $u_{x_2}\left(x_1,0\right)$.
If, for instance, $g=0$, then $u=x_2$ and $u=-x_2$ are both
solutions of Cauchy problem \eqref{NonlinEq-6}.

\medskip

We have already studied the method of characteristics for the linear and the quasilinear
linear equations. In the nonlinear case we follow a
procedure similar to the previous two cases, but in the
nonlinear case it is less obvious which are the
characteristic equations. For this purpose some geometrical considerations may be
useful, which, neverthless, we do not take up here,
referring the interested reader to \cite[Ch. 1]{Joh}. We start
by the following

\begin{definition}\label{NonlinEq-10} Let $F$ be the function
	\eqref{NonlinEq-5}. Let us assume that $F$ is of class
	$C^{2}\left(\Omega\times \mathbb{R}\times \mathbb{R}^n\right)$.
	We call the characteristic equations related to the partial differential equation
	
	\begin{equation}\label{NonlinEq-15}
		F\left(x,u(x),\nabla u(x)\right)=0,\mbox{ in } \Omega,\end{equation}
	the following system of ordinary differential equations
	
	\begin{equation}\label{NonlinEq-20}
		\begin{cases}
			\frac{d X(t)}{dt}=\nabla_pF(X(t),z(t),p(t)),\\
			\\
			\frac{d z(t)}{dt}=\nabla_pF(X(t),z(t),p(t))\cdot p(t),\\
			\\
			\frac{d p(t)}{dt}=-\partial_zF(X(t),z(t),p(t))p(t)-\nabla_x
			F(X(t),z(t),p(t)).
		\end{cases}
	\end{equation}
	The function $X(\cdot), z(\cdot), p(\cdot)$ are called the
	\textbf{characteristic lines of the equation} \index{characteristic:@{characteristic:}!- line@{- line}}\eqref{NonlinEq-15}.
	$X(\cdot)$ is called the \textbf{ray} \index{characteristic:@{characteristic:}!- projected line (ray)@{- projected line (ray)}}or the
	\textbf{projected characteristic lines} on $\mathbb{R}^n$.
\end{definition}
If $F$ does not depend on $z$, system \eqref{NonlinEq-20} is
decoupled in $z$, while the first and the third equations
constitute the \textbf{Hamilton--Jacobi system}:\index{Hamilton--Jacobi system}

\begin{equation}\label{NonlinEq-21}
	\begin{cases}
		\frac{d X(t)}{dt}=\nabla_pF(X(t),p(t)),\\
		\\
		\frac{d p(t)}{dt}=-\nabla_x F(X(t),p(t)).
	\end{cases}
\end{equation}

\bigskip

\bigskip

\textbf{Remarks.}

\noindent\textbf{1.} Let us note that if $(X(\cdot),z(\cdot),
p(\cdot))$ is a characteristic line for equation
\eqref{NonlinEq-15} then

\begin{equation}\label{NonlinEq-25}
	F(X(t),z(t),p(t))=\mbox{constant}.
\end{equation}
As a matter of fact, exploiting \eqref{NonlinEq-20}, we have

\begin{equation*}
	\begin{aligned}
		&\frac{d}{dt}F(X(t),z(t),p(t))=\nabla_xF(X(t),z(t),p(t))\cdot\frac{dX}{dt}+\\&
		+F_z(X(t),z(t),p(t))\frac{dz}{dt}+\nabla_pF(X(t),z(t),p(t))\cdot\frac{dp}{dt}=\\&
		=\nabla_xF\cdot \nabla_pF+F_z\nabla_pF\cdot p(t)+\\&
		+\nabla_pF\cdot\left(-F_zp(t)-\nabla_xF\right)=0,
	\end{aligned}
\end{equation*}
In the last step, for the sake of brevity, we have omitted the arguments
$X(t),z(t),p(t)$ in $F$.

\noindent\textbf{2.} Let us suppose that $u$ is a solution of class
$C^2$ to the equation \begin{equation} \label{NonlinEq-30}
	F\left(x,u(x),\nabla u(x)\right)=0,
\end{equation}
we wish to look for $X(t)$ (or, more precisely,
for an equation for $X(t)$) such that, setting

\begin{equation*}
	z(t)=u(X(t)),\quad\quad p(t)=\nabla u(X(t)),
\end{equation*}
it happens that $(X(t),z(t),p(t))$ solves system
\eqref{NonlinEq-20}. 

We have

\begin{equation}\label{NonlinEq-35}
	\frac{d z(t)}{dt}=p(t)\cdot \frac{dX(t)}{dt} ,
\end{equation}

\begin{equation}\label{NonlinEq-41}
	\frac{d
		p_i(t)}{dt}=\sum_{j=1}^n\partial_{ij}^2u(X(t))\frac{dX_j(t)}{dt},\quad
	i=1,\cdots, n.
\end{equation}
Now, calculating the derivatives of both the sides of equation \eqref{NonlinEq-30} w.r.t. $x_i$,
$i=1,\cdots n$, we have

\begin{equation}\label{NonlinEq-42}
	\sum_{j=1}^n\partial_{p_j}F\partial_{ij}^2u(x)=-\partial_{x_i}F-\partial_zF\partial_{x_i}u(x),
\end{equation}
where the argument of $F$ in \eqref{NonlinEq-42} is
$\left(x,u(x),\nabla u(x)\right)$.  Now, let us observe what follows: if

\begin{equation}\label{NonlinEq-45}
	\frac{d X_j(t)}{dt}=\partial_{p_j}F(X(t),z(t),p(t)),\quad
	j=1,\cdots, n,
\end{equation}
then by \eqref{NonlinEq-41} and by \eqref{NonlinEq-42},
calculated for $x=X(t)$, we have, for $i=1,\cdots, n$

\begin{equation}\label{NonlinEq-50}
	\frac{d
		p_i(t)}{dt}=-\partial_{x_i}F(X(t),z(t),p(t))-\partial_zF(X(t),z(t),p(t))p_i(t)
\end{equation}
and by \eqref{NonlinEq-35} we have

\begin{equation}\label{NonlinEq-55}
	\frac{d z(t)}{dt}=\sum_{j=1}^n \partial_{p_j}F
	(X(t),z(t),p(t))p_j(t).
\end{equation}
Equations \eqref{NonlinEq-45}, \eqref{NonlinEq-50} and \eqref{NonlinEq-55}
are just the equations of the system \eqref{NonlinEq-20}.
$\blacklozenge$

\bigskip

In order to solve Cauchy problem \eqref{NonlinEq-0} we will follow an
approach similar to that followed in the linear (and quasilinear) case
by letting the projected characteristic lines,
$X(t)$, play a similar role to that played, in the linear case, by the
characteristic lines.

In what follows we will consider the case

\begin{equation}\label{NonlinEq-60}
	\Gamma=\left\{x\in \Omega: x_n=0 \right\}.
\end{equation}
We observe that we can always lead back to this situation,
at least locally, even if $\Gamma$ is given by

\begin{equation*}
	\Gamma=\left\{x\in \Omega: \phi(x)=0 \right\},
\end{equation*}
where $\phi\in C^{3}(\Omega)$ and $\phi(x_0)=0$ for a given $x_0\in
\Omega$ and

\begin{equation}\label{NonlinEq-65}
	\nabla\phi(x_0)\neq 0.
\end{equation}
Indeed, thanks to \eqref{NonlinEq-65}, there exists a
 neighborhood, $\mathcal{U}$, of $x_0$ such that $\Gamma\cap \mathcal{U}$ is
is a graph of a function of $n-1$ variables. If, for instance,
let us suppose that $\phi_{x_n}(x_0)\neq 0$ then, up to a
translation that moves $x_0$ to $0$, we may assume that for an
appropriate $\delta>0$, we have

\begin{equation}\label{NonlinEq-70}
	\Gamma\cap \mathcal{U}=\left\{(x',\varphi(x')): x'\in
	B'_{\delta}\right\},
\end{equation}
where $\varphi\in C^3(B'_{\delta})$, $\varphi(0)=\left|\nabla_{x'}\varphi(0)\right|=0$. Now, let
																											$$\Lambda: B_{\delta}\subset\mathbb{R}^n_x\rightarrow\mathbb{R}^n_y,\quad \Lambda(x)=\left(x',x_n-\varphi(x')\right),$$

$$\Lambda(\Gamma)=\left\{(y',0): y'\in B'_{\delta}\right\}=\left\{y\in B_{\delta}:-y_n=0\right\}$$
and, setting
$$v(y)=u\left(\Lambda^{-1}(y)\right),$$ we easily obtain that the
problem \eqref{NonlinEq-0} takes the form

\begin{equation}\label{NonlinEq-75}
	\begin{cases}
		\widetilde{F}\left(y,v(y),\nabla_y v(y)\right)=0, \quad \mbox{in } \mathcal{V},\\
		\\
		v(y)=\widetilde{g}(y), \quad \mbox{for } y\in \Lambda(\Gamma)\cap
		\mathcal{V},
	\end{cases}
\end{equation}
where $\mathcal{V}$ is a neighborhood of $0$ and
$\widetilde{F}:\mathcal{V}\times \mathbb{R}\times
\mathbb{R}^n\rightarrow \mathbb{R}$ is a function of class
$C^{2}(\mathcal{V}\times \mathbb{R}\times \mathbb{R}^n)$.

\begin{theo}\label{teor-NonL}
	Let $R>0$ and $F\in C^{2}(B_R\times \mathbb{R}\times
	\mathbb{R}^n)$. Let $g\in C^{2}(B'_R)$.  Let $\eta\in \mathbb{R}$
	satisfy
	
	\begin{equation}\label{NonlinEq-80}
		F\left(0,g(0),\nabla_{x'}g(0), \eta\right)=0
	\end{equation}
	and
	\begin{equation}\label{NonlinEq-85}
		F_{p_n}\left(0,g(0),\nabla_{x'}g(0), \eta\right)\neq 0,
	\end{equation}
	then for some $r\in (0,R)$ there exists a unique solution 
	$u\in C^{2}(B_r)$ to the initial--value problems 
	
	\begin{equation}\label{NonlinEq-90}
		\begin{cases}
			F\left(x,u(x),\nabla u(x)\right)= 0, \quad \mbox{in } B_r,\\
			\\
			u(x',0)=g(x'), \quad \mbox{for } x'\in B_r',\\
			\\
			u_{x_n}(0,0)=\eta.
		\end{cases}
	\end{equation}
	
\end{theo}

\bigskip

\textbf{Proof.} 
Let us begin by proving \textbf{the uniqueness}. It will suffice to prove that if $u$ is a
solution to

\begin{equation}\label{NonlinEq-95}
	\begin{cases}
		F\left(x,u(x),\nabla u(x)\right)= 0,\quad \mbox{in } B_R,\\
		\\
		u(x',0)=g(x'),\quad \mbox{for } x'\in B_R', \\
		\\
		u_{x_n}(0,0)=\eta,
	\end{cases}
\end{equation}
then there exists a neighborhood of $0$ in which $u$ is uniquely
determined.

Set

\begin{equation*}
	z_0=g(0),\quad\quad p'_0=\nabla_{x'}g(0).
\end{equation*}
By \eqref{NonlinEq-80} and \eqref{NonlinEq-85} we have
\begin{equation*}
	F\left(0,z_0,p'_0, \eta\right)=0,\quad\quad F_{p_n}\left(0,z_0,p'_0,
	\eta\right)\neq 0.
\end{equation*}
By applying the Implicit Function Theorem, we have that there exists
$\delta\in (0,R]$ such that, setting
$$\mathcal{U}=B_{\delta}\times
\left(z_0-\delta,z_0+\delta\right)\times B'_{\delta}(p_0')\times
\left(\eta-\delta,\eta+\delta\right),$$ we have that the set

$$\left\{(x,z,p)\in \mathcal{U}:\mbox{ } F(x,z,p)=0 \right\},$$ is equal to the graph of the function

\begin{equation}\label{NonlinEq-96}
	\psi:B_{\delta}\times
	\left(z_0-\delta,z_0+\delta\right)\times
	B'_{\delta}(p_0')\rightarrow
	\left(\eta-\delta,\eta+\delta\right),\end{equation} where $\psi$ is of
class $C^{2}$ and $$\psi(0,z_0,p'_0)=\eta.$$ Now, since $u$
satisfies \eqref{NonlinEq-95}, we have
$$F\left(x',0,g(x'),\nabla_{x'}g(x'),u_{x_n}(x',0)\right)=0,$$

$$u_{x_n}(0,0)=\eta,$$
and we have

$$u_{x_n}(x',0)=\psi\left(x',0,g(x'),\nabla_{x'}g(x')\right),\quad
\forall x'\in B'_{\delta}.$$ Set
$$p_n(x')=\psi\left(x',0,g(x'),\nabla_{x'}g(x')\right)$$
and 
\begin{equation}\label{NonlinEq-104}
	p^{(0)}(x')=\left(\nabla_{x'}g(x'),p_n(x') \right),
\end{equation}
we have

\begin{equation}\label{NonlinEq-105}
	\begin{cases}
		F\left(x,u(x),\nabla u(x)\right)= 0, \quad \mbox{in } B_{\delta},\\
		\\
		u(x',0)=g(x'), \quad \mbox{for } x'\in B_{\delta}',\\
		\\
		\nabla u (x',0)=p^{(0)}(x'),\quad \mbox{for } x'\in B_{\delta}'.
	\end{cases}
\end{equation}
We denote by $y$ an arbitrary point of $B_{\delta}'$ and
recalling that by Remark 2 of the present Section, the function $$t\rightarrow
(X(t),z(t),p(t)):= (X(t,y),u(X(t,y)), \nabla u(X(t,y))$$ is a solution to the characteristic equations (for each $y\in B_{\delta}'$) 

\begin{equation}\label{NonlinEq-110}
	\begin{cases}
		\frac{d X}{dt}=\nabla_pF(X,z,p),\\
		\\
		\frac{d z}{dt}=\nabla_pF(X,z,p)\cdot p,\\
		\\
		\frac{d p}{dt}=-\partial_zF(X,z,p)p-\nabla_x F(X,z,p)
	\end{cases}
\end{equation}
and

\begin{equation}\label{NonlinEq-115}
	\begin{cases}
		X(0,y)=(y,0),\\
		\\
		z(0,y)=u(X(0,y))=g(y),\\
		\\
		p(0,y)=\nabla u(X(0,y))=p^{(0)}(y).
	\end{cases}
\end{equation}

Therefore, due to the uniqueness of the solution to Cauchy problem
\eqref{NonlinEq-110}--\eqref{NonlinEq-115} it turns out
that $u(X(t,y))$ is determined for every $y\in B_{\delta}'$ and
for every $t$ in a neighborhood of $0$ (this neighborhood depends on $y$).
To conclude the proof, it suffices, therefore, to prove that
the map

\begin{equation}\label{NonlinEq-120}
	(t,y)\rightarrow X(t,y),
\end{equation} 
is a diffeomorphism in a neighborhood of
$0$. From what we said in Section \ref{richiami-ode} (final part), map \eqref{NonlinEq-120} is of class $C^{2}$. To
establish that it is a local diffeomorphism, it suffices to check that the
Jacobian matrix of $(t,y)\rightarrow X(t,y)$ is nonsingular
in $0$. Now from \eqref{NonlinEq-115} we have

\begin{equation}\label{NonlinEq-125}
	\frac{\partial X(0,0)}{\partial(t,y)}=\left(
	\begin{array}{ccc}
		\partial_{t}X_1(0,0) \mbox{ } \partial_{y_1}X_1(0,0)&\cdots  & \partial_{y_{n-1}}X_1(0,0) \\
		\vdots\mbox{ }\vdots &\cdots  &\vdots \\
		\partial_{t}X_n(0,0) \mbox{ } \partial_{y_1}X_n(0,0)&  \cdots&
		\partial_{y_{n-1}}X_n(0,0)
	\end{array}%
	\right).
\end{equation}
On the other hand
$$\partial_{y_i}X_j(0,0)=\delta_{ij}, \mbox{for  } 1\leq i\leq n-1,\mbox{ } 1\leq j\leq n-1,$$
$$\partial_{y_i}X_n(0,0)=0, \mbox{ for } 1\leq i\leq n-1$$
and, for $1\leq j\leq n-1$,
$$\partial_{t}X_j(0,0)=\partial_{p_j}F(X(0,0),z(0,0),p(0,0))=\partial_{p_j}F\left(0,g(0),\nabla_{x'}g(0), \eta\right).$$
Hence

\begin{equation}\label{NonlinEq-130}
	\det\left(\frac{\partial
		X(0,0)}{\partial(t,y)}\right)=(-1)^n\partial_{p_n}F\left(0,g(0),\nabla_{x'}g(0)
	,\eta\right)\neq 0,
\end{equation}
from which it follows that map \eqref{NonlinEq-120} is a
local diffeomorphism. The proof of uniqueness is complete.

\medskip

Now, let us prove \textbf{the esistence} of the solution to problem
\eqref{NonlinEq-90}. Let $p^{(0)}(x')$ be defined by
\eqref{NonlinEq-104} and let $(X(t,y),z(t,y),p(t,y))$ be the solution 
of the Cauchy problem comprising the system \eqref{NonlinEq-110} and
the initial conditions

\begin{equation}\label{NonlinEq-135}
	\begin{cases}
		X(0,y)=(y,0),\\
		\\
		z(0,y)=g(y),\\
		\\
		p(0,y)=p^{(0)}(y).
	\end{cases}
\end{equation}
Set
\begin{equation*}
	f(t,y)=F\left(X(t,y),z(t,y),p(t,y)\right),
\end{equation*}
we have

\begin{equation*}
	\begin{aligned}
		f(0,y)&=F\left((y,0),g(y),p^{(0)}(y)\right)=\\&=F\left((y,0),g(y),\nabla_yg(y),\psi\left((y,0),g(y),\nabla_{y}g(y)\right)\right)=0,
	\end{aligned}
\end{equation*}
where $\psi$ is given by \eqref{NonlinEq-96}. Hence, by
\eqref{NonlinEq-25}, we have

\begin{equation}\label{NonlinEq-140}
	f(t,y)=F\left(X(t,y),z(t,y),p(t,y)\right)=0.
\end{equation}
Moreover, in a completely similar way to what has been done above for the uniqueness
we have that there exists $\delta_1>0$ and a neighbourhood of $0$,
$\mathcal{U}_0$, such that
$$B_{\delta_1}\ni(t,y)\rightarrow X(t,y)\in \mathcal{U}_0,$$
is a diffeomorphism of class $C^2\left(B_{\delta_1}(0)\right)$.
We denote the inverse of $X(\cdot,\cdot)$ by
$$X^{-1}(x)=(t(x),y(x))$$ and set

\begin{equation*}
	u(x)=z(t(x),y(x)), \quad p(x)=p(t(x),y(x)).
\end{equation*}
The remaining part of the proof consists of proving that $u$ satisfies \eqref{NonlinEq-90}. First of all, we check that

\begin{equation}\label{NonlinEq-145}
	u(x',0)=g(x').
\end{equation}
To this purpose we note that 

\begin{equation*}
	t(x',0)=0,\quad y(x',0)=x'.
\end{equation*}
Hence
\begin{equation*}
	u(x',0)=z\left(t(x',0),y(x',0)\right)=z\left(0,x'\right)=g\left(x'\right).
\end{equation*}
Therefore, we have \eqref{NonlinEq-145}. We will check the condition
$\partial_{x_n}u(0,0)=\eta$ later on,  now
we check that $u$ satisfies the equation

\begin{equation}\label{NonlinEq-150}
	F\left(x,u(x),\nabla u(x)\right)= 0.
\end{equation}
Firstly, we observe that from \eqref{NonlinEq-140} we have 
\begin{equation}\label{NonlinEq-155}
	F\left(x,u(x),p(x)\right)= f(t(x),y(x))=0, \quad \forall x\in
	\mathcal{U}_0.
\end{equation}
Therefore, to prove \eqref{NonlinEq-150} it suffices to prove that
\begin{equation}\label{NonlinEq-156}
	p(x)=\nabla u(x), \quad \forall x\in \mathcal{U}_0.
\end{equation}
To this purpose we prove the following claims:

\medskip

\textbf{Claim I}
\begin{equation}\label{NonlinEq-157}
	\partial_tz(t,y)=\sum_{j=1}^np_j(t,y)\partial_tX_j(t,y),\quad \forall (t,y)\in B_{\delta_1}.
\end{equation}

\medskip

\textbf{Claim II}
\begin{equation}\label{NonlinEq-160}
	\partial_{y_i}z(t,y)=\sum_{j=1}^np_j(t,y)\partial_{y_i}X_j(t,y),\quad \forall (t,y)\in B_{\delta_1}.
\end{equation}

\bigskip

\textbf{Claim I} follows by the first and the second equation of \eqref{NonlinEq-110}. As a matter of fact, we have
$$\partial_tz(t,y)=\nabla_p F\left(X(t,y),z(t,y),p(t,y)\right)\cdot
p=\partial_tX(t,y)\cdot p(t,y).$$

\bigskip

The proof of \textbf{Claim II} is less immediate than Claim I. Set

\begin{equation}\label{NonlinEq-165}
	h_i(t,y)=\partial_{y_i}z(t,y)-\sum_{j=1}^np_j(t,y)\partial_{y_i}X_j(t,y).
\end{equation}
By \eqref{NonlinEq-104} and recalling that

\begin{equation}
	\begin{aligned}
	&\partial_{y_i}X_j(0,y)=\delta_{ij}, \mbox{ for } 1\leq i,j\leq
	n-1,\\&  \partial_{y_i}X_n(0,y)=0 , \mbox{ for } 1\leq i\leq n-1,	
\end{aligned}
\end{equation}

we have, for $i=1,\cdots,n-1$,

\begin{equation}\label{NonlinEq-170}
	h_i(0,y)=\partial_{y_i}z(0,y)-\sum_{j=1}^np_j(0,y)\partial_{y_i}X_j(0,y)=\partial_{y_i}g(y)-p^{(0)}_i(y)=0.
\end{equation}
Now, we prove that $h_i(\cdot,y)$ satisfies

\begin{equation}\label{NonlinEq-175}
	\partial_th_i(t,y)=-\partial_zF\left(X(t,y),z(t,y),p(t,y)\right)h_i(t,y).
\end{equation}
By \eqref{NonlinEq-157} we have
\begin{equation}\label{NonlinEq-180}
	\partial^2_{ty_i}z=\sum_{j=1}^n\left(\partial_{y_i}p_j\partial_tX_j+p_j\partial^2_{ty_i}X_j\right).
\end{equation}
Now making the derivative w.r.t. $t$ of both the sides of \eqref{NonlinEq-165} we have

\begin{equation*}
	\begin{aligned}
		\partial_{t}h_i=\partial^2_{ty_i}z-\sum_{j=1}^n\left(\partial_{t}p_j\partial_{y_i}X_j+p_j\partial^2_{ty_i}X_j\right).
	\end{aligned}
\end{equation*}
By this equality, by \eqref{NonlinEq-180} and by \eqref{NonlinEq-110}
we have

\begin{equation}\label{NonlinEq-185}
	\begin{aligned}
		\partial_{t}h_i&=\sum_{j=1}^n\left(\partial_{y_i}p_j\partial_tX_j-\partial_{t}p_j\partial_{y_i}X_j\right)=\\&=
		\sum_{j=1}^n\left(\partial_{y_i}p_j\partial_{p_j}F-\left(-\partial_{x_j}F-\partial_zFp_j\right)\partial_{y_i}X_j\right)=\\&=
		\sum_{j=1}^n\left(\partial_{y_i}p_j\partial_{p_j}F+\partial_{x_j}F\partial_{y_i}X_j+\partial_zFp_j\partial_{y_i}X_j\right).
	\end{aligned}
\end{equation}
Now by \eqref{NonlinEq-140}, making the derivative w.r.t.  $y_i$ of both the sides, we get 

\begin{equation*}
	\begin{aligned}
		\sum_{j=1}^n\left(\partial_{y_i}p_j\partial_{p_j}F+\partial_{x_j}F\partial_{y_i}X_j\right)=-\partial_zF\partial_{y_i}z
	\end{aligned}
\end{equation*}
and inserting the latter into \eqref{NonlinEq-185} we get

\begin{equation*}
	\begin{aligned}
		\partial_{t}h_i&=-\partial_zF\partial_{y_j}z+\sum_{j=1}^n\partial_zFp_j\partial_{y_i}X_j
		=\\&=-\partial_zF\left(\partial_{y_i}z-\sum_{j=1}^np_j\partial_{y_i}X_j\right)=\\&=
		-\partial_zFh_i.
	\end{aligned}
\end{equation*}
All in all, by the latter and by \eqref{NonlinEq-170} we get,
for $i=1,\cdot,n-1$,

\begin{equation*}
	\begin{cases}
		\partial_{t}h_i=-\partial_zFh_i,\\
		\\
		h_i(0,y)=0,
	\end{cases}
\end{equation*}
from which we have
$$h_i(t,y)=0,\quad \mbox{for }i=1,\cdot,n-1.$$
Claim II is proved.

\medskip

Now, let us prove \eqref{NonlinEq-156}. First, we recall that
$u(x)=z\left(X^{-1}(x)\right)=z(t(x),y(x))$. We have by 
\eqref{NonlinEq-157} and \eqref{NonlinEq-160},

\begin{equation*}
	\begin{aligned}
		\partial_{x_i}u&=\partial_tz\partial_{x_i}t+\sum_{j=1}^{n-1}\partial_{y_j}z\partial_{x_i}y_j
		=\\&=\left(\sum_{k=1}^{n}p_k\partial_t
		X_k\right)\partial_{x_i}t+\sum_{j=1}^{n-1}\sum_{k=1}^{n}p_k\partial_{y_j}X_k\partial_{x_i}y_j=\\&=
		\sum_{k=1}^{n}p_k\left(\partial_t
		X_k\partial_{x_i}t+\sum_{j=1}^{n-1}\partial_{y_j}X_k\partial_{x_i}y_j\right)=\\&=
		\sum_{k=1}^{n}p_k\partial_{x_i}\left(X_k\left(X^{-1}(x)\right)\right)=\\&=
		\sum_{k=1}^{n}p_k\delta_{ik}=p_i
	\end{aligned}
\end{equation*}
for $i=1,\cdots, n$. From which we have \eqref{NonlinEq-156} and, in
particular,
$$\partial_{x_n}u(0)=p_n(0)=\eta.$$ Which concludes the proof. $\blacksquare$

\section{Appendix: geodesics and Hamilton--Jacobi equations}\label{appendice}
We warn that throughout this Appendix we will adopt the convention
of repeated indices. In addition, we will strictly adhere to the
notation on indices (upper or lower) for the components of a
tensor.

In the first part of this Appendix we will present the rudiments of the
theory of Hamilton-Jacobi equations, these topics can
be carried out in a more general way, for more details we refer
to \cite{EV}.

Let $\Omega$ be an open set of $\mathbb{R}^n$, we say that a real--valued function, $$L\in
C^{\infty}\left(\Omega\times \mathbb{R}^n\right)$$ is a 
\textbf{Lagrangian} \index{Lagrangian} on $\Omega$. An example of
Lagrangian that we are interested in is given by

\begin{equation}  \label{NonlinEq-285}
	L(x,q)=\frac{1}{2}g_{ij}\left(x\right)q^iq^j,\quad \forall x\in
	\Omega, \mbox{ } \forall q\in \mathbb{R}^n,
\end{equation}
where  $\left\{g_{ij}\left(x\right) \right\}_{i,j=1}^{n}$ is a real symmetric nonsingular
matrix, $n\times n$, whose entries belong to $C^{\infty}\left(\Omega\right)$.

Given a Lagrangian $L$ we will call \textbf{equation of
	Euler--Lagrange} \index{equation:@{equation:}!- Euler--Lagrange@{- Euler--Lagrange}} 
 the differential equation in the unknown $x=x(t)$

\begin{equation}  \label{NonlinEq-300n}
	\frac{d}{dt}\left(\nabla_qL\left(x(t),\overset{\cdot
	}{x}(t)\right)\right)-\nabla_xL\left(x(t),\overset{\cdot
	}{x}(t)\right)=0.
\end{equation}
Here and in the sequel we will indistinctly let us denote by $\frac{df}{dt}$
or by $\overset{\cdot }{f}(t)$ the derivative with respect to $t$ of a
differentiable function $f$ .
The solutions $x:[t_0,t_1]\rightarrow
\mathbb{R}^n$ of \eqref{NonlinEq-300n} are also called the \textbf{extremal} 
of the functional \index{extremal of a functional}

\begin{equation}  \label{NonlinEq-220n}
	\int^{t_1}_{t_0}L\left(x(t),\overset{\cdot }{x}(t)\right)dt.
\end{equation}

\noindent \textbf{Assumption I.} In what follows we suppose that, for every
$p\in \mathbb{R}^n$, the equation
\begin{equation}  \label{NonlinEq-290}
	\nabla_qL(x,q)=p,
\end{equation}
has a unique solution of class
$C^{\infty}\left(\Omega\times \mathbb{R}^n\right)$. We denote such a solution by $q(x,p)$.

\medskip

The function

\begin{equation}  \label{NonlinEq-295}
	H(x,p)=p\cdot q(x,p)-L\left(x,q(x,p)\right).
\end{equation}
is called the \textbf{Hamiltonian} \index{Hamiltonian}associated to $L$
\medskip

We have

\begin{theo}\label{prop51-geod} Let $x=x(t)$ be a solution 
	to  Euler--Lagrange equation
	
	\begin{equation}  \label{NonlinEq-300}
		\frac{d}{dt}\left(\nabla_qL\left(x(t),\overset{\cdot
		}{x}(t)\right)\right)-\nabla_xL\left(x(t),\overset{\cdot
		}{x}(t)\right)=0.
	\end{equation}
	Then, setting
	$$p(t)=\nabla_qL\left(x(t),\overset{\cdot
	}{x}(t)\right),$$ it turns out that $(x(t),p(t))$ is a solution of the 
	Hamilton--Jacobi system
	
	\begin{equation}\label{NonlinEq-305}
		\begin{cases}
			\frac{dx}{dt}=\nabla_pH(x(t),p(t)),\\
			\\
			\frac{dp}{dt}=-\nabla_xH(x(t),p(t)).
		\end{cases}
	\end{equation}
	Moreover
	\begin{equation}\label{NonlinEq-310}
		H(x(t),p(t))=\mbox{constant}.
	\end{equation}
\end{theo}
\textbf{Proof.} Let $x(t)$ be a solution to equation 
\eqref{NonlinEq-300}. Then, since equation
\eqref{NonlinEq-290} has a unique solution, $q(x,p)$,  and since 
$$p(t)=\nabla_qL\left(x(t),\overset{\cdot
}{x}(t)\right),$$ we have
\begin{equation}\label{NonlinEq-315}
	\overset{\cdot }{x}(t)=q(x(t),p(t)).
\end{equation}
Now we have, for $i=1,\cdots,n$,

\begin{equation*}
	\begin{aligned}
		&\partial_{p_i}H(x,p)=\partial_{p_i}\left(p\cdot
		q(x,p)-L\left(x,q(x,p)\right)\right)=\\&=p_k\partial_{p_i}q^k(x,p)+q^i(x,p)-\partial_{q^k}L(x,q(x,p))\partial_{p_i}q^k(x,p)=\\&=
		\partial_{p_i}q^k(x,p)\left(p_k-\partial_{q^k}L(x, q(x,p))\right)+q^i(x,p)=\\&=
		q^i(x,p).
	\end{aligned}
\end{equation*}
Hence, recalling \eqref{NonlinEq-315}, we have, for
$i=1,\cdots,n$,
$$\frac{dx^i(t)}{dt}=\partial_{p_i}H(x(t),p(t)),$$
which is the first equation of system \eqref{NonlinEq-305}.
Concerning the second equation, for $i=1,\cdots,n$,
we have
\begin{equation}\label{NonlinEq-316}
	\begin{aligned}
		&\partial_{x^i}H(x,p)=\partial_{x^i}\left(p\cdot
		q(x,p)-L\left(x,q(x,p)\right)\right)=\\&=
		p_k\partial_{x^i}q^k(x,p)-\partial_{x^i}L(x,q(x,p))-\partial_{q^k}L(x,q(x,p))\partial_{x^i}q^k(x,p)=\\&=
		\partial_{x^i}q^k(x,p)\left(p_k-\partial_{q^k}L(x,q(x,p))\right)-\partial_{x^i}L(x,q(x,p))=\\&=
		-\partial_{x^i}L\left(x,q(x,p)\right).
	\end{aligned}
\end{equation}
On the other hand, by \eqref{NonlinEq-300} and \eqref{NonlinEq-315},
we have, for $i=1,\cdots,n$,

\begin{equation*}
	\begin{aligned}
		\partial_{x^i}L(x(t),q(x(t),p(t))&=\partial_{x^i}L(x(t),\overset{\cdot
		}{x}(t))=\\&= \frac{d}{dt}\partial_{q^i}L(x(t),\overset{\cdot
		}{x}(t))=\\&= \frac{dp_i(t)}{dt}.
	\end{aligned}
\end{equation*}
By the just obtained equality and by \eqref{NonlinEq-316} we get
\begin{equation}\label{NonlinEq-320}
	\frac{dp_i(t)}{dt}=-\partial_{x^i}H(x(t),p(t)),\quad\mbox{ for }
	i=1,\cdots,n.
\end{equation}
 which is the second equation of system \eqref{NonlinEq-305}.

Finally, \eqref{NonlinEq-310} follows by  \eqref{NonlinEq-305} and by 
\begin{equation*}
	\begin{aligned}
		&\frac{d}{dt}H(x(t),p(t))=\partial_{x^i}H(x(t),p(t))\frac{dx^i(t)}{dt}+\partial_{p_i}H(x(t),p(t))\frac{dp_i(t)}{dt}=\\&=\partial_{x^i}H(x(t),p(t))\partial_{p_i}H(x(t),p(t))-\partial_{p_i}H(x(t),p(t))
		\partial_{x^i}H(x(t),p(t))=0.
	\end{aligned}
\end{equation*}
$\blacksquare$

\bigskip

\noindent \textbf{Assumption II.} Let $H$ be the Hamiltonian associated to
the Lagrangian $L$ which satisfies Assumption I. Let us suppose that, for every $q\in \mathbb{R}^n$, the equation
\begin{equation*}
	\nabla_pH(x,p)=q
\end{equation*}
has a unique solution of class
$C^{\infty}\left(\Omega\times \mathbb{R}^n\right)$. We denote by $p(x,q)$ such solution. Let us notice that \eqref{NonlinEq-295} implies trivially
\begin{equation}\label{NonlinEq-399}
	L(x,q)=p(x,q)\cdot q-H(x,p(x,q)).
\end{equation}

We now prove the converse of Theorem \ref{prop51-geod}.

\begin{theo}\label{prop52-geod} Let $H$ the Hamiltonian associated to
	$L$ and let us suppose that Assumption I and II hold true.
	Moreover, let us suppose  $(x(t),p(t))$ that is a solution to the Hamilton--Jacobi equation
	
	\begin{equation}\label{NonlinEq-400}
		\begin{cases}
			\frac{dx(t)}{dt}=\nabla_pH(x(t),p(t)),\\
			\\
			\frac{dp(t)}{dt}=-\nabla_xH(x(t),p(t)).
		\end{cases}
	\end{equation}
	Then $x(t)$ is a solution to Euler--Lagrange equation
	\begin{equation}  \label{NonlinEq-405}
		\frac{d}{dt}\left(\nabla_qL\left(x(t),\overset{\cdot
		}{x}(t)\right)\right)-\nabla_xL\left(x(t),\overset{\cdot
		}{x}(t)\right)=0.
	\end{equation}
	
\end{theo}

\textbf{Proof.} By Assumpyion II and by \eqref{NonlinEq-400},
in particular by

$$\nabla_pH(x(t),p(t))=\overset{\cdot
}{x}(t),$$ we have
$$p(t)=p\left(x(t),\overset{\cdot
}{x}(t)\right).$$ Now, by \eqref{NonlinEq-399} we have
\begin{equation}\label{NonlinEq-401}
	L\left(x,\overset{\cdot }{x}(t)\right)=p(t)\cdot \overset{\cdot
	}{x}(t)-H\left(x(t),p(t)\right).
\end{equation}
On the other hand

$$\nabla_q L(x,q)=p(x,q),$$
hence
$$\nabla_q L\left(x(t),\overset{\cdot
}{x}(t)\right)=p\left(x(t),\overset{\cdot }{x}(t)\right)=p(t).$$
Therefore, taking into account \eqref{NonlinEq-400}, we have

\begin{equation}\label{NonlinEq-410}
	\frac{d}{dt}\nabla_qL\left(x(t),\overset{\cdot
	}{x}(t)\right)=\overset{\cdot }{p}(t)=-\nabla_xH(x(t),p(t)).
\end{equation}
Now, let us make the derivatives w.r.t. $x^i$ of both the sides of \eqref{NonlinEq-399}

\begin{equation*}
	\begin{aligned}
		&\partial_{x^i}L(x,q)=q^k\partial_{x^i}p_k(x,q)-\partial_{x^i}H(x,p(x,q))-\\&-\partial_{p_k}H(x,p(x,q))\partial_{x^i}p_k(x,q)=\\&=
		\left(q^k-\partial_{p_k}H(x,p(x,q))\right)\partial_{x^i}p_k(x,q)-\partial_{x^i}H(x,p(x,q)),
	\end{aligned}
\end{equation*}
from which we have

\begin{equation*}
	\begin{aligned}
		&\partial_{x^i}L\left(x(t),\overset{\cdot
		}{x}(t)\right)=\\&=\left(\frac{dx^k}{dt}-\partial_{p_k}H(x(t),p(t))\right)\partial_{x^i}p_k(x(t),p(t))-\\&-\partial_{x^i}H(x(t),p(t))=\\&=
		-\partial_{x^i}H(x(t),p(t)).
	\end{aligned}
\end{equation*}
By the just obtained inequality and by \eqref{NonlinEq-410} we have

\begin{equation*}
	\frac{d}{dt}\nabla_qL\left(x(t),\overset{\cdot
	}{x}(t)\right)=\nabla_xL\left(x(t),\overset{\cdot }{x}(t)\right).
\end{equation*}
$\blacksquare$

\bigskip

Now let us consider the \textbf{geodesic lines} with respect to the Riemannian metric \index{Riemannian}
$$g_{ij}( x) dx^{i}\otimes dx^{j},$$ where $\left\{g_{ij}\left(x\right)
\right\}_{i,j=1}^{n}$ is a symmetric real matrix $n\times n$
whose entries belong to $C^{\infty}\left(\Omega\right)$. Let us suppose

\begin{equation}  \label{NonlinEq-190}
	\lambda^{-1}\left|\xi\right|^{2}\leq
	g_{ij}\left(x\right)\xi^i\xi^j\leq\lambda\left|\xi\right|^{2},\quad
	\forall\xi\in\mathbb{R}^{n} \mbox{, }\forall x\in \Omega,
\end{equation}
where $\lambda\geq 1$. Let us denote by $\left\{ g^{ij}\left(
x\right) \right\}_{i,j=1}^{n}$ the inverse matrix of
$\left\{g_{ij}\left(x\right) \right\}_{i,j=1}^{n}$.
\begin{definition}
	\index{Definition:@{Definition:}!- geodesic line@{- geodesic line}}
	We say that the path 
	$$\gamma:[t_0,t_1]\rightarrow \Omega,$$
	Is a \textbf{geodesic line} with respect to the Riemannian metric $g_{ij}( x)
	dx^{i}\otimes dx^{j}$, if $\gamma\in
	C^{\infty}\left([t_0,t_1],\Omega\right)$ and it solves the equations \index{equation:@{equation:}!- of geodesic lines@{- geodesic lines}}
	
	\begin{equation}  \label{NonlinEq-195}
		\frac{d^2\gamma^h(t)}{dt^2}+\Gamma^h_{ij}(\gamma(t))\frac{d
			\gamma^i(t)}{dt}\frac{d \gamma^j(t)}{dt}=0,\quad h=1,\cdots,n,
	\end{equation}
	where (Christoffel symbols) \index{Christoffel symbols}, for $i,j,h=1,\cdots,n$,
	
	\begin{equation}  \label{NonlinEq-200}
		\Gamma^h_{ij}(x)=\frac{1}{2}g^{hk}(x)\left[\partial_{i}g_{kj}(x)+\partial_{j}g_{ki}(x)-\partial_{k}g_{ij}(x)\right].
	\end{equation}
\end{definition}

\bigskip

The following Proposition holds true.

\begin{prop}\label{prop1-geod}
	The path $\gamma:[t_0,t_1]\rightarrow \Omega$ is a geodesic line w.r.t.
	the Riemannian metric $g_{ij}( x) dx^{i}\otimes dx^{j}$ if and only if
	$\gamma$ is an extremal of the Lagrangian
	\begin{equation}  \label{NonlinEq-205}
		L(x,q)=g_{ij}(x)q^iq^j,\quad x\in \Omega,\mbox{ } q\in \mathbb{R}.
	\end{equation}
\end{prop}
\textbf{Proof.} Let us write the Euler--Lagrange equation

\begin{equation}  \label{NonlinEq-215}
	\frac{d}{dt}\left(\nabla_qL\left(x(t),\overset{\cdot
	}{x}(t)\right)\right)-\nabla_xL\left(x(t),\overset{\cdot
	}{x}(t)\right)=0.
\end{equation}
We have

\begin{equation*}
		\partial_{q_k}L\left(x(t),\overset{\cdot
}{x}(t)\right)=2g_{kj}(x(t))\frac{dx^j(t)}{dt},
\end{equation*}

\medskip

\begin{equation*}
	\begin{aligned}
		\frac{d}{dt}\left(\partial_{q_k}L\left(x(t),\overset{\cdot
		}{x}(t)\right)\right)&=2g_{kj}(x(t))\frac{d^2x^j(t)}{dt^2}+\\&+2\partial_{x_i}g_{kj}(x(t))\frac{dx^i(t)}{dt}\frac{dx^j(t)}{dt}
	\end{aligned}
\end{equation*}

\medskip

\noindent and
$$\partial_{x_k}L\left(x(t),\overset{\cdot
}{x}(t)\right)=\partial_{x_k}g_{ij}(x(t))\frac{dx^i(t)}{dt}\frac{dx^j(t)}{dt}.$$
Therefore the Euler--Lagrange equation can be written (we omit the variables, for the sake of brevity)

\begin{equation*}
	\begin{aligned}
		g_{kj}\frac{d^2x^j}{dt^2}=\left(\frac{1}{2}\partial_{x_k}g_{ij}-\partial_{i}g_{kj}\right)\frac{dx^i}{dt}\frac{dx^j}{dt}.
	\end{aligned}
\end{equation*}
From which we have

\begin{equation*}
	\begin{aligned}
		\frac{d^2x^h}{dt^2}&=g^{hk}\left(\frac{1}{2}\partial_{x_k}g_{ij}-\partial_{i}g_{kj}\right)\frac{dx^i}{dt}\frac{dx^j}{dt}=\\&=
		\frac{1}{2}g^{hk}\left(\partial_{x_k}g_{ij}-\partial_{x_i}g_{kj}-\partial_{x_j}g_{ki}\right)\frac{dx^i}{dt}\frac{dx^j}{dt}=\\&=-\Gamma^h_{ij}\frac{dx^i}{dt}\frac{dx^j}{dt}.
	\end{aligned}
\end{equation*}
Hence the equation

\begin{equation*}
	\begin{aligned}
		\frac{d^2x^h}{dt^2}+\Gamma^h_{ij}\frac{dx^i}{dt}\frac{dx^j}{dt}=0,
	\end{aligned}
\end{equation*}
is equivalent to the Euler--Lagrange equation   related to $L$
and by \eqref{NonlinEq-195} the thesis follows. $\blacksquare$

\bigskip

\textbf{Remark.} Let $\left\{g_{ij}\left(x\right)
\right\}_{i,j=1}^{n}$ be a matrix like in Proposition \ref{prop1-geod}, set
$$L(x,q)=\frac{1}{2}g_{ij}\left(x\right)q^iq^j,\quad \forall x\in \Omega,
\mbox{ } \forall q\in \mathbb{R}^n.$$ Notice that, by
\eqref{NonlinEq-285}, we can write \eqref{NonlinEq-290} as

\begin{equation*}
	g_{ij}q^j=p_i,\quad \mbox{ for } i=1,\cdots,n.
\end{equation*}
Hence
\begin{equation*}
	H(x,p)=g^{ij}p_jp_i-\frac{1}{2}g_{ij}g^{ih}p_hg^{ik}p_k=\frac{1}{2}g^{ij}p_jp_i.
\end{equation*}

$\blacklozenge$

\bigskip

The following Theorem holds true

\bigskip

\begin{theo}\label{prop5-geod}
	Let $\left\{g_{ij}\left(x\right) \right\}_{i,j=1}^{n}$ be a real symmetric matrix $n\times n$ whose entries belong to
	$C^{\infty}\left(\Omega\right)$ and let us assume that it satisfies
	\eqref{NonlinEq-190}. Let $u\in C^{\infty}(\Omega)$ be a solution
	to the eikonal equation \index{equation:@{equation:}!- eikonal@{- eikonal}}
	\begin{equation}  \label{NonlinEq-280}
		g^{ij}\left(x\right)\partial_{x^i}u\partial_{x^j}u=1
	\end{equation}
	and let $x=\gamma(t)$ be a solution to the system
	\begin{equation*}
		\frac{dx^i(t)}{dt}=g^{ij}\left(x(t)\right)\partial_{x^j}u(x(t)),\quad
		i=1,\cdots,n.
	\end{equation*} Then $x=\gamma(t)$ is a geodesic line w.r.t. 
	the Riemannian metric $$g_{ij}( x) dx^{i}\otimes dx^{j}.$$
\end{theo}
\textbf{Proof.} Let
$$L(x,q)=\frac{1}{2}g_{ij}\left(x\right)q^iq^j,\quad \forall x\in \Omega,
\mbox{ } \forall q\in \mathbb{R}^n.$$ Let $H$ be the Hamiltonian of $L$,
that is
\begin{equation*}
	H(x,p)=\frac{1}{2}g^{ij}(x)p_jp_i\quad \forall x\in \Omega, \mbox{ }
	\forall p\in \mathbb{R}^n.
\end{equation*}

Set

\begin{equation}\label{NonlinEq-1a}
	p(t)=\nabla u(\gamma(t)),
\end{equation}
where $u$  is a solution to equation \eqref{NonlinEq-280} and
$\gamma$ is a  solution to the equations

\begin{equation}  \label{NonlinEq-281a}
	\frac{d\gamma^i(t)}{dt}=g^{ij}\left(\gamma(t)\right)\partial_{x^j}u(\gamma(t))\left(=\partial_{p_i}H(\gamma(t),p(t))\right),
\end{equation}
for $i=1, \cdots, n$.

\medskip

Now, we make the derivative w.r.t. $x^k$ of both the sides of equation \eqref{NonlinEq-280}
and we get

\begin{equation}  \label{NonlinEq-2a}
	2\left(g^{ij}\left(x\right)\partial^2_{x^ix^k}u\right)\partial_{x^j}u+\left(\partial_{x^k}g^{ij}\left(x\right)\right)\partial_{x^i}u\partial_{x^j}u=0.
\end{equation}
By \eqref{NonlinEq-1a}, \eqref{NonlinEq-281a} and \eqref{NonlinEq-2a}
we have, for $i=1,\cdots,n$,
\begin{equation*}
	\begin{aligned}
		\frac{dp_i(t)}{dt}&=\partial^2_{x^ix^k}u(\gamma(t))\frac{d\gamma^k(t)}{dt}=\\&=
		\partial^2_{x^ix^k}u(\gamma(t))g^{kj}\left(\gamma(t)\right)\partial_{x_j}u(\gamma(t))=\\&=
		\left(g^{kj}\left(\gamma(t)\right)\partial^2_{x^ix^k}u(\gamma(t))\right)\partial_{x^j}u(\gamma(t))=\\&=
		-\frac{1}{2}\partial_{x^i}g^{jk}\left(\gamma(t)\right)\partial_{x^j}u(\gamma(t))\partial_{x^k}u(\gamma(t))=\\&=
		-\partial_{x^i} H(\gamma(t),p(t)).
	\end{aligned}
\end{equation*}
The just obtained equality and \eqref{NonlinEq-281a} implies that
$(\gamma(t),p(t))$ is a solution to the system

\begin{equation}\label{NonlinEq-400a}
	\begin{cases}
		\frac{d\gamma(t)}{dt}=\nabla_pH(x(t),\gamma(t)),\\
		\\
		\frac{dp(t)}{dt}=-\nabla_xH(\gamma(t),p(t)).
	\end{cases}
\end{equation}

\medskip

Therefore, by Theorem \ref{prop52-geod} and by Proposition
\eqref{prop1-geod} the thesis follows. $\blacksquare$

\bigskip

\textbf{Remark.} Let us observe that if $u$ is a solution to the equation

\begin{equation*}
	g^{ij}\left(x\right)\partial_{x^i}u\partial_{x^j}u=1
\end{equation*}
and $\gamma(t)$ is a solution to the system
\begin{equation*}
	\frac{d\gamma^i(t)}{dt}=g^{ij}\left(\gamma(t)\right)\partial_{x^j}u(\gamma(t)),\quad
	i=1,\cdots,n,
\end{equation*}
then $t$ is the natural parameter (in the Riemannian metric)
of the path $x=\gamma(t)$. As a matter of fact we have

\begin{equation*}
	\begin{aligned}
		g_{ij}\left(\gamma(t)\right)\frac{d\gamma^i(t)}{dt}\frac{d\gamma^j(t)}{dt}&=g_{ij}\left(\gamma(t)\right)\left(g^{ik}\left(\gamma(t)\right)\partial_{x^k}u(\gamma(t))\right)
		\left(g^{jl}\left(\gamma(t)\right)\partial_{x^l}u(\gamma(t))\right)=\\&=
		\delta^k_j\partial_{x^k}u(\gamma(t))g^{jl}\left(\gamma(t)\right)\partial_{x^l}u(\gamma(t))=\\&=
		g^{jl}\left(\gamma(t)\right)\partial_{x^l}u(\gamma(t))\partial_{x^j}u(\gamma(t))=1.
	\end{aligned}
\end{equation*}
Moreover, for fixed $\overline{x}\in \mathbb{R}^n$, set
$$u\left(\overline{x}\right)=R_0.$$
If $\gamma$ is the solution to the Cauchy
problem

\begin{equation*}
	\begin{cases}
		\frac{d\gamma^i(t)}{dt}=g^{ij}\left(\gamma(t)\right)\partial_{x^j}u(\gamma(t)),\quad
		i=1,\cdots,n,\\
		\\
		\gamma(R_0)=\overline{x},
	\end{cases}
\end{equation*}
then

\begin{equation}\label{ossGeod}
	u(\gamma(t))=t,\quad \forall t\in I
\end{equation} ($I$ is the maximal interval
of the solution $\gamma$). As a matter of fact we have 
\begin{equation*}
	\begin{aligned}
		\frac{d}{dt}u(\gamma(t))&=\partial_{x^i}u(\gamma(t))\frac{d\gamma^i(t)}{dt}=\\&=
		\partial_{x^i}u(\gamma(t))
		g^{ij}\left(\gamma(t)\right)\partial_{x^j}u(\gamma(t))=1,
	\end{aligned}
\end{equation*}
hence

$$u(\gamma(t))=t+C,$$ where $C$ is a constant which can be
determined easily in the following way

$$R_0=u\left(\overline{x}\right)=u\left(\gamma(R_0)\right)=R_0+C,$$
hence $C=0$ and \eqref{ossGeod} is proved. These
comments will be used in Ch. \ref{tre
	sfere-ellittiche}. $\blacklozenge$

\bigskip

Let us conclude this Appendix by some propositions on the extremal
and other comments on the geodetics lines.

\begin{prop}\label{prop2-geod}
	Let $L\in C^{\infty}\left(\Omega\times\mathbb{R}^n\right)$. We have what follows.
	
	\smallskip
	
	\noindent (i) if $\overline{x}\in
	C^{\infty}\left([t_0,t_1],\mathbb{R}^n\right)$ is an  extremal
	of the functional
	\begin{equation}  \label{NonlinEq-220}
		\int^{t_1}_{t_0}L\left(x(t),\frac{dx(t)}{dt}\right)dt,
	\end{equation}
	we have
	\begin{equation}  \label{NonlinEq-225}
		\nabla_qL\left(\overline{x}(t),\frac{d\overline{x}(t)}{dt}\right)\cdot
		\frac{d\overline{x}(t)}{dt}-L\left(\overline{x}(t),\frac{d\overline{x}(t)}{dt}\right)=\mbox{
			constant}.
	\end{equation}
	\noindent (ii) If $L(x,q)$ is an homogeneous function w.r.t.  $q$ of degree $\alpha\neq 1$ and $\overline{x}(t)$ is an extremal of functional \eqref{NonlinEq-220} then
	
	\begin{equation}  \label{NonlinEq-230}
		L\left(\overline{x}(t),\frac{d\overline{x}(t)}{dt}\right)=\mbox{
			constant}.
	\end{equation}
\end{prop}
\textbf{Proof.}

(i) Set
$$F(t)=\nabla_qL\left(\overline{x}(t),\frac{d\overline{x}(t)}{dt}\right)\cdot
\frac{d\overline{x}(t)}{dt}-L\left(\overline{x}(t),\frac{d\overline{x}(t)}{dt}\right).$$
We have
\begin{equation*}
	\begin{aligned}
		&\frac{d
			F(t)}{dt}=\frac{d}{dt}\left(\nabla_qL\left(\overline{x}(t),\frac{d\overline{x}(t)}{dt}\right)\right)\cdot
		\frac{d\overline{x}(t)}{dt}
		+\nabla_qL\left(\overline{x}(t),\frac{d\overline{x}(t)}{dt}\right)\cdot
		\frac{d^2\overline{x}(t)}{dt^2}-\\&-\nabla_xL\left(\overline{x}(t),\frac{d\overline{x}(t)}{dt}\right)\cdot
		\frac{d\overline{x}(t)}{dt}-\nabla_qL\left(\overline{x}(t),\frac{d\overline{x}(t)}{dt}\right)\cdot
		\frac{d^2\overline{x}(t)}{dt^2}=\\&=
		\left[\frac{d}{dt}\left(\nabla_qL\left(\overline{x}(t),\frac{d\overline{x}(t)}{dt}\right)\right)-\nabla_xL\left(\overline{x}(t),\frac{d\overline{x}(t)}{dt}\right)\right]\cdot
		\frac{d\overline{x}(t)}{dt}=0.
	\end{aligned}
\end{equation*}
From which the thesis follows.

\medskip

(ii) By point (i) we have

$$\nabla_qL\left(\overline{x}(t),\frac{d\overline{x}(t)}{dt}\right)\cdot
\frac{d\overline{x}(t)}{dt}-L\left(\overline{x}(t),\frac{d\overline{x}(t)}{dt}\right)=\mbox{constant}.$$
On the other hand by the homogeneity of $L(x,\cdot)$ we get

$$\nabla_qL\left(\overline{x}(t),\frac{d\overline{x}(t)}{dt}\right)\cdot
\frac{d\overline{x}(t)}{dt}=\alpha
L\left(\overline{x}(t),\frac{d\overline{x}(t)}{dt}\right).$$
Therefore
$$(\alpha-1)L\left(\overline{x}(t),\frac{d\overline{x}(t)}{dt}\right)
=\mbox{constant}$$ and recalling that $\alpha\neq 1$, the thesis follows.
$\blacksquare$

\bigskip

\begin{prop}\label{prop3-geod} Let $L\in
	C^{\infty}\left(\Omega\times\mathbb{R}^n\right)$ be an homogeneous function w.r.t. $q$ of degree
	$2$. If $\overline{x}\in
	C^{\infty}\left([t_0,t_1],\mathbb{R}^n\right)$ is an extremal
	of the functional
	\begin{equation}  \label{NonlinEq-235}
		\int^{t_1}_{t_0}L\left(x(t),\frac{dx(t)}{dt}\right)dt,
	\end{equation}
	and \begin{equation} \label{NonlinEq-236}
		L\left(\overline{x}(t),\frac{d\overline{x}(t)}{dt}\right)>0,\quad
		\forall t\in[t_0,t_1],\end{equation} then $\overline{x}$ is an
	extremal of the functional
	\begin{equation}  \label{NonlinEq-240}
		\int^{t_1}_{t_0}\sqrt{L\left(x(t),\frac{dx(t)}{dt}\right)}dt.
	\end{equation}
\end{prop}
\textbf{Proof.} By Proposition \eqref{prop2-geod} and by
\eqref{NonlinEq-236} we may set
$$c^2_0=L\left(\overline{x}(t),\frac{d\overline{x}(t)}{dt}\right)>0,$$
where $c_0$ is a positive constant. Since $\overline{x}$ is an
extremal of the functional \eqref{NonlinEq-235}, we have
\begin{equation*}
	\begin{aligned}
		\frac{d}{dt}\left(\nabla_q\sqrt{L\left(\overline{x}(t),\frac{d\overline{x}(t)}{dt}\right)}\right)&
		=\frac{1}{2c_0}\frac{d}{dt}\left(\nabla_q
		L\left(\overline{x}(t),\frac{d\overline{x}(t)}{dt}\right)\right)=\\&=\frac{1}{2c_0}\nabla_xL\left(\overline{x}(t),\frac{d\overline{x}(t)}{dt}\right)=\\&=
		\nabla_x\sqrt{L\left(\overline{x}(t),\frac{d\overline{x}(t)}{dt}\right)}.
	\end{aligned}
\end{equation*}
Hence $\overline{x}$ is an extremal of functional
\eqref{NonlinEq-240}. $\blacksquare$

\bigskip

\begin{prop}\label{prop4-geod} Let us suppose that $L$ satisfies the same assumptions of Proposition \ref{prop3-geod}. Let $\varphi$ an extremal
	of functional
	
	\begin{equation}  \label{NonlinEq-245}
		\int^{t_1}_{t_0}\sqrt{L\left(x(t),\frac{dx(t)}{dt}\right)}dt.
	\end{equation}
	Let us suppose $$L\left(\varphi(t),\frac{d\varphi(t)}{dt}\right)>0,\quad
	\forall t\in [t_0,t_1],$$ then there exists a unique
	parametrization $t(\tau)$, $t'(\tau)>0$ in $[\tau_0,\tau_1]$
	($t(\tau_0)=t_0$ and $t(\tau_1)=t_1$), such that, setting
	$\psi(\tau)=\varphi(t(\tau))$, we have
	\begin{equation}  \label{NonlinEq-246}
		L\left(\psi(\tau),\frac{d\psi(\tau)}{d\tau}\right)=\mbox{constant}.
	\end{equation}
	Moreover the path $x=\psi(\tau)$ is an extremal of the functional
	\begin{equation}  \label{NonlinEq-250}
		\int^{\tau_1}_{\tau_0}L\left(x(\tau),\frac{dx(\tau)}{d\tau}\right)d\tau.
	\end{equation}
\end{prop}
\textbf{Proof.} Let

$$f(t)=\sqrt{L\left(\varphi(t),\frac{d\varphi(t)}{dt}\right)},\quad \forall t\in
[t_0,t_1].$$ Let $c_0>0$ a be constant and let $t(\tau)$ satisfy 

$$\int^{t(\tau)}_{t_0}f(t)dt=c_0\tau, \ \ \forall \tau\in\left[\tau_0,\tau_1\right]. $$ Set
$\psi(\tau)=\varphi(t(\tau))$; we get, by the homogeneity of
$L(x,\cdot)$,

\begin{equation*}
	\begin{aligned}
		\sqrt{L\left(\psi(\tau),\frac{d\psi(\tau)}{d\tau}\right)}&=\sqrt{L\left(\varphi(t(\tau)),\frac{d\varphi}{dt}(t(\tau))\right)t'(\tau)}=\\&=
		t'(\tau)\sqrt{L\left(\varphi(t(\tau)),\frac{d\varphi}{dt}(t(\tau))\right)}=c_0.
	\end{aligned}
\end{equation*}
Hence \eqref{NonlinEq-246} is proved.

Now, let us prove that $x=\psi(\tau)$ is an  extremal of the functional
\eqref{NonlinEq-250}. Set

$$F(t)=\frac{1}{f(t)}\nabla_q
L\left(\varphi(t),\frac{d\varphi(t)}{dt}\right)$$ and

$$G(t)=\frac{1}{f(t)}\nabla_x
L\left(\varphi(t),\frac{d\varphi(t)}{dt}\right).$$ Since
$\varphi$ is an extremal of functional \eqref{NonlinEq-245}, we have

\begin{equation}  \label{NonlinEq-255}
	\frac{dF(t)}{dt}=G(t).
\end{equation}
Now, recalling
$$\sqrt{L\left(\psi(\tau),\frac{d\psi(\tau)}{d\tau}\right)}=c_0,$$
we have (by the homogeneity of $L$ w.r.t. $q$)
\begin{equation*}
	\begin{aligned}
		\nabla_q
		L\left(\psi(\tau),\frac{d\psi(\tau)}{d\tau}\right)&=t'(\tau)\nabla_q
		L\left(\varphi(t(\tau)),\frac{d\varphi}{dt}(t(\tau))\right)=\\&=
		t'(\tau)\sqrt{L\left(\varphi(t(\tau)),\frac{d\varphi}{dt}(t(\tau))\right)}
		F(t(\tau))=\\&=c_0F(t(\tau)).
	\end{aligned}
\end{equation*}
Hence, recalling \eqref{NonlinEq-255} (and the homogeneity of $L$
w.r.t. $q$), we have

\begin{equation*}
	\begin{aligned}
		\frac{d}{d\tau}\left(\nabla_qL\left(\psi(\tau),
		\frac{d\psi(\tau)}{d\tau}\right)\right)&=c_0\frac{d}{d\tau}(F(t(\tau)))=\\&=
		c_0t'(\tau)\frac{d F}{dt}(t(\tau))=c_0t'(\tau)G(t(\tau))=\\&=
		\frac{c_0t'(\tau)}{f(t(\tau))}\nabla_x
		L\left(\varphi(t(\tau)),\frac{d\varphi}{dt}(t(\tau))\right)=\\&=
		\frac{c_0}{t'(\tau)f(t(\tau))}\nabla_x
		L\left(\psi(\tau),\frac{d\psi(\tau)}{d\tau}\right)=\\&= \nabla_x
		L\left(\psi(\tau),\frac{d\psi(\tau)}{d\tau}\right).
	\end{aligned}
\end{equation*}
Hence $\psi$ is an extremal of functional
\eqref{NonlinEq-250}. $\blacksquare$

\bigskip

\textbf{Remark.} Let $\left\{g_{ij}\left(x\right)
\right\}_{i,j=1}^{n}$ be a real symmetric matrix  $n\times n$ whose entries belong to $C^{\infty}\left(\Omega\right)$ and let us assume that it satisfies \eqref{NonlinEq-190}. Let
$$L(x,q)=g_{ij}\left(x\right)q^iq^j,\quad \forall x\in \Omega,
\mbox{ } \forall q\in \mathbb{R}^n.$$ By Proposition
\ref{prop3-geod} we have that, if $x=\gamma (t)$ is a geodesic line,
i.e. it is an extremal of the functional

\begin{equation}  \label{NonlinEq-260}
	\int^{t_1}_{t_0}L\left(x(t),\frac{dx(t)}{dt}\right)dt=
	\int^{t_1}_{t_0}g_{ij}\left(x(t)\right)\frac{dx^i}{dt}\frac{dx^j}{dt}dt,
\end{equation}
then $x=\gamma (t)$ is also an extremal of the functional
\begin{equation}  \label{NonlinEq-265}
	\int^{t_1}_{t_0}L\left(x(t),\frac{dx(t)}{dt}\right)dt=
	\int^{t_1}_{t_0}\sqrt{g_{ij}\left(x(t)\right)\frac{dx^i}{dt}\frac{dx^j}{dt}}dt.
\end{equation}
On the other hand, by Proposition \ref{prop4-geod}, we have that if
$x=\gamma (t)$ is an extremal of functional
\eqref{NonlinEq-265} and if $t(\tau)$ is strictly increasing and it satisfies

\begin{equation}  \label{NonlinEq-270}
	t'(\tau)\sqrt{g_{ij}\left(\gamma(t(\tau))\right)\frac{d\gamma^i}{dt}(t(\tau)\frac{dx^j}{dt}(t(\tau)}=c,
\end{equation}
where $c>0$ is a positive constant, then $x=\gamma(t(\tau)$
is an extremal of

\begin{equation}  \label{NonlinEq-275}
	\int^{t_1}_{t_0}L\left(x(t),\frac{dx(t)}{dt}\right)dt=
	\int^{\tau_1}_{\tau_0}g_{ij}\left(x(\tau)\right)\frac{dx^i}{d\tau}\frac{dx^j}{d\tau}d\tau,
\end{equation}
where $\tau_0$ and $\tau_1$ satisfy $t(\tau_0)=t_0$ and
$t(\tau_1)=t_1$. Let us notice that if $c=1$, then condition
\eqref{NonlinEq-270} means that $\tau$ is the natural parameter
of the path $x=\gamma(t)$ (extremal of
\eqref{NonlinEq-265}) in the riemanniann metric $g_{ij}( x)
dx^{i}\otimes dx^{j}$. $\blacklozenge$

\newpage

\chapter{Real analytic functions}\label{funz-analitiche}
\section{Power series}\label{funz-analitiche-Intr}

In this chapter we will consider the multiple series \index{multiple series}

\begin{equation}\label{1-16C}
\sum_{\alpha\in \mathbb{N}^n_0}c_{\alpha},
\end{equation}
where $c_{\alpha}\in \mathbb{R}$ (or $c_{\alpha}\in \mathbb{C}$).

When we say that the series \eqref{1-16C} converges,
we will mean \textit{always} that it is \textbf{absolutely convergent} \index{absolutely convergent multiple series}. That is 
\begin{equation*}
\sum_{\alpha\in \mathbb{N}^n_0}|c_{\alpha}|<+\infty.
\end{equation*}
Therefore, if the series \eqref{1-16C} converges, the value
of the sum in \eqref{1-16C} does not depend on the order of the terms
$c_{\alpha}$. If $c_{\alpha}(x)$ are functions, we will naturally extend the
notions of uniform, total convergence, $C^k(\overline{\Omega})$ convergence and so on. For instance, we will say
that

\begin{equation}\label{0-16C}
\sum_{\alpha\in \mathbb{N}^n_0}c_{\alpha}(x),
\end{equation}
uniformly converges to a function $f$ in a set $K\subset
\mathbb{R}^n$ provided that: 

(i) for every $x\in K$,  $\sum_{\alpha\in
\mathbb{N}^n_0}|c_{\alpha}(x)|$ converges, 

(ii) we have
$$f(x)=\sum_{\alpha\in \mathbb{N}^n_0}c_{\alpha}(x),\quad\forall x\in K$$
and

(iii)
$$\lim_{N\rightarrow +\infty}\sup_{x\in K}\left\vert f(x)-\sum_{|\alpha|\leq N}c_{\alpha}(x)\right\vert=0.$$
 Let $c_{\alpha}\in \mathbb{R}$ (or $c_{\alpha}\in \mathbb{C}$) we call \textbf{power series} \index{power series}a series like  
\begin{equation}\label{2-16C}
\sum_{\alpha\in \mathbb{N}^n_0}c_{\alpha}x^{\alpha}.
\end{equation}
For any $y\in \mathbb{R}^n$ set
$$Q_y=\{x\in \mathbb{R}^n: |x_j|\leq |y_j|, j=1,\cdots,n\}.$$
We have 

\begin{prop}\label{prop0-16C}
If series \eqref{2-16C} converges at a point $y\in \mathbb{R}^n$
then the series uniformly converges in $Q_y$.
\end{prop}

\textbf{Proof.} The convergence at $y$ of \eqref{2-16C}
is equivalent to

$$\sum_{\alpha\in \mathbb{N}^n_0}|c_{\alpha}||y^{\alpha}|<+\infty.$$
Hence, we have

$$\sum_{\alpha\in \mathbb{N}^n_0}\sup_{Q_y}|c_{\alpha}x^{\alpha}|\leq\sum_{\alpha\in \mathbb{N}^n_0}|c_{\alpha}||y^{\alpha}|<+\infty.$$
From which we get the total convergence and, consequently, the uniform convergence, of series
\eqref{2-16C}.$\blacksquare$

\bigskip

Proposition \ref{prop0-16C} implies that the sum of series
\eqref{2-16C} is continuous in $Q_y$.

The differentiability will be proved in Proposition
\ref{prop1-17C} to prove such a Proposition we need  

\begin{lem}\label{Es1-19C}
Let us denote by $\upsilon=(1,1,\cdots,1)$. If $|x_j|<1$ 
$j=1,\cdots,n$, we have  
\begin{equation}\label{1-19C}
\sum_{\alpha\in
\mathbb{N}^n_0}x^{\alpha}=\frac{1}{(\upsilon-x)^{\upsilon}},
\end{equation}
and
\begin{equation}\label{2-19C}
\sum_{\alpha\geq \beta}\frac{\alpha!}{(\alpha-\beta)!}
x^{\alpha-\beta}=\partial^{\beta}\left(\frac{1}{(\upsilon-x)^{\upsilon}}\right)=
\frac{\beta!}{(\upsilon-x)^{\upsilon+\beta}}.
\end{equation}
\end{lem}
\textbf{Proof.} Concerning the convergence of
series \eqref{1-19C}, we have, for $|x_j|<1$, $j=1,\cdots,n$,

\begin{equation*}
\begin{aligned}
\sum_{\alpha\leq \upsilon N}\left|x^{\alpha}\right|&=
\sum_{\alpha\leq \upsilon N}\left|x_j^{\alpha_j}\right|=\\&
=\prod_{j=1}^n\sum_{\alpha_j\leq N}\left|x_j^{\alpha_j}\right|=\\&
=\prod_{j=1}^n\frac{1-|x_j|^{N+1}}{1-|x_j|}\rightarrow
\prod_{j=1}^n\frac{1}{1-|x_j|},\quad \mbox{as } N\rightarrow\infty.
\end{aligned}
\end{equation*}
Concerning the sum of the series we have, similarly,

\begin{equation*}
\sum_{\alpha\leq \upsilon N}x^{\alpha}=
\lim_{N\rightarrow\infty}\prod_{j=1}^n\sum_{\alpha_j\leq
N}x_j^{\alpha_j}= \frac{1}{(1-x_1)\cdots
(1-x_n)}=\frac{1}{(\upsilon-x)^{\upsilon}}.
\end{equation*}

Now, let us prove \eqref{2-19C}. Recalling \eqref{5-3N} 
\begin{equation*}
\partial^{\beta}x^{\alpha}=
\begin{cases}
	\frac{\alpha!}{(\alpha-\beta)!}x^{\alpha-\beta}, \quad \mbox{ for } \alpha\geq\beta,\\
	\\
	0, \quad \quad\mbox{otherwise}, %
\end{cases}%
\end{equation*}

we have

\begin{equation*}
\begin{aligned}
\sum_{\alpha\geq \beta}\frac{\alpha!}{(\alpha-\beta)!}
x^{\alpha-\beta}&=
\sum_{\alpha\in\mathbb{N}^n_0}\partial^{\beta}x^{\alpha}=\\&
=\prod_{j=1}^n\sum_{\alpha_j\in\mathbb{N}_0}\partial^{\beta_j}x_j^{\alpha_j}=\\&
=\prod_{j=1}^n\partial_j^{\beta_j}\frac{1}{1-x_j}=
\partial^{\beta}\left(\frac{1}{(\upsilon-x)^{\upsilon}}\right)=\\&=
\frac{\beta!}{(\upsilon-x)^{\upsilon+\beta}}.
\end{aligned}
\end{equation*}
$\blacksquare$

\begin{prop}\label{prop1-17C}
If series \eqref{2-16C} converges at the point $y\in \mathbb{R}^n$ and
$y_j>0$ for every $j=1,\cdots,n$ then, denoted by $f$ the sum of such a series, we have $f\in C^{\infty}\left(Int(Q_y)\right)$, where $Int(Q_y)$ is the interior part
of $Q_y$. 

Moreover

\begin{equation}\label{1-17C}
\partial^{\alpha}f(0)=\frac{1}{\alpha!}c_{\alpha}.
\end{equation}
\end{prop}

\textbf{Proof.} In order to prove that $f\in
C^{\infty}\left(Int(Q_y)\right)$ it suffices to prove that for every
$q\in(0,1)$ and for every $\beta\in \mathbb{N}_0^n$ we have
$$\sum_{\alpha\in \mathbb{N}^n_0}\sup_{Q_{qy}}\left|\partial^{\beta}\left(c_{\alpha}x^{\alpha}\right)\right|<+\infty.$$
By \eqref{5-3N} we have

\begin{equation*}
\partial^{\beta}\left(c_{\alpha}x^{\alpha}\right)=\left \{
\begin{array}{c}
\frac{\alpha!}{(\alpha-\beta)!}c_{\alpha}x^{\alpha-\beta}, \quad \mbox{ for } \alpha\geq\beta,\\
\\
0, \quad \mbox{otherwise}. %
\end{array}%
\right.
\end{equation*}
Hence

\begin{equation}\label{1-18C}
\begin{aligned}
\sum_{\alpha\in
\mathbb{N}^n_0}\sup_{Q_{qy}}\left|\partial^{\beta}\left(c_{\alpha}x^{\alpha}\right)\right|&\leq
\sum_{\alpha\geq \beta}
\frac{\alpha!}{(\alpha-\beta)!}|c_{\alpha}|\left|\left(q
y\right)^{\alpha-\beta}\right|=\\&
=\frac{1}{\left|y^{\beta}\right|}\sum_{\alpha\geq
\beta}\frac{\alpha!|c_{\alpha}y^{\alpha}|}{(\alpha-\beta)!}q^{|\alpha-\beta|}.
\end{aligned}
\end{equation}
Now, since \eqref{2-16C} converges in $y\in \mathbb{R}^n$,
we get

\begin{equation}\label{2-18C}
\left|c_{\alpha}y^{\alpha}\right|\leq \mu_y:=\sum_{\alpha\in
\mathbb{N}^n_0}\left|c_{\alpha}y^{\alpha}\right|<+\infty.
\end{equation}
By the above obtained inequality and by \eqref{1-18C} we have (for $y\neq 0$)

\begin{equation}\label{3-18C}
\sum_{\alpha\in
\mathbb{N}^n_0}\sup_{Q_{qy}}\left|\partial^{\beta}\left(c_{\alpha}x^{\alpha}\right)\right|
\leq \frac{\mu_y}{\left|y^{\beta}\right|}\sum_{\alpha\geq
\beta}\frac{\alpha!}{(\alpha-\beta)!}q^{|\alpha-\beta|}.
\end{equation}
Applying \eqref{2-19C} with $x=q\upsilon=q(1,1\cdots,1)$,
we have

\begin{equation}\label{4-18C}
	\begin{aligned}
	\sum_{\alpha\geq
\beta}\frac{\alpha!}{(\alpha-\beta)!}q^{|\alpha-\beta|}&=
\sum_{\alpha\geq
\beta}\frac{\alpha!}{(\alpha-\beta)!}(q\upsilon)^{\alpha-\beta}=\\&=
\frac{\beta!}{(\upsilon-\upsilon
q)^{\upsilon+\beta}}=\\&=\frac{\beta!}{(1-q)^{n+|\beta|}}.
\end{aligned}
\end{equation}

\smallskip

By the just obtained equality and by \eqref{3-18C} we obtain
\begin{equation}\label{Asterisco-19C}
\sum_{\alpha\in
\mathbb{N}^n_0}\sup_{Q_{qy}}\left|\partial^{\beta}\left(c_{\alpha}x^{\alpha}\right)\right|\leq
\frac{\mu_y}{\left|y^{\beta}\right|}\frac{\beta!}{(1-q)^{n+|\beta|}}<+\infty.
\end{equation}
We have proved so far that for each $q\in(0,1)$ we have $f\in
C^{\infty}\left(Q_{yq}\right)$. Therefore $f\in
C^{\infty}\left(\mbox{Int}(Q_y)\right).$

Concerning \eqref{1-17C} we have

$$\partial^{\beta}f(x)=\sum_{\alpha\in \mathbb{N}^n_0}\partial^{\beta}\left(c_{\alpha}x^{\alpha}\right)=
\sum_{\alpha\geq
\beta}\frac{\alpha!}{(\alpha-\beta)!}c_{\alpha}x^{\alpha-\beta}.$$
Therefore
$$\partial^{\beta}f(0)=\frac{1}{\beta!}c_{\beta}.$$$\blacksquare$

\bigskip

\underline{\textbf{Exercise 1.}} Let $f(x)$ be the sum of the series 
$$\sum_{\alpha\in\mathbb{N}^n_0}c_{\alpha}x^{\alpha},$$
in $Q_{\varrho}$ where $\varrho=(\varrho_1,\cdots,\varrho_n)$, with
$\varrho_j>0$, $j=1,\cdots,n$. Let $q\in (0,1)$ and set
\begin{equation}\label{1-20C}
r=(1-q)\min_{1\leq j\leq n}\varrho_j,\quad\mbox{ and }\quad
\mu_{\varrho}=\sum_{\alpha\in \mathbb{N}^n_0}
\left|c_{\alpha}\varrho^{\alpha}\right|.
\end{equation}
Then we have

\begin{equation}\label{2-20C}
\sup_{Q_{q\varrho}}\left|\partial^{\beta}f\right|\leq
(1-q)^{-n}\mu_{\varrho}r^{-|\beta|}\beta!.
\end{equation}

\textbf{Solving Exercise 1.}

Is an immediate consequence of \eqref{Asterisco-19C}. $\clubsuit$

\bigskip

\underline{\textbf{Exercise 2.}} Prove that

\begin{equation}\label{1-21C}
\sum_{\alpha\in\mathbb{N}^n_0}\frac{|\alpha|!}{\alpha!}x^{\alpha}=\frac{1}{1-\sum_{j=1}^nx_j},
\quad\mbox{ for } \ \  \sum_{j=1}^n|x_j|<1,
\end{equation}

\begin{equation}\label{2-21C}
\sum_{\alpha\geq
\beta}\frac{|\alpha|!}{(\alpha-\beta)!}x^{\alpha-\beta}=\frac{|\beta|!}{\left(1-\sum_{j=1}^nx_j\right)^{1+|\beta|}},
\quad\mbox{ for }  \ \  \sum_{j=1}^n|x_j|<1.
\end{equation}

\textbf{Solving Exercise 2.}

By \eqref{3-3N} we have
\begin{equation*}
(x_1+x_2+\cdots+x_n)^m=\sum_{|\alpha|=m}\frac{m!}{\alpha!}x^{\alpha},\quad
m\in\mathbb{N}_0.
\end{equation*}
Hence, for $\sum_{j=1}^n|x_j|<1$, we have
\begin{equation*}
\begin{aligned}
\sum_{\alpha\in
\mathbb{N}^n_0}\frac{|\alpha|!}{\alpha!}x^{\alpha}&=\sum_{m=0}^{\infty}
\sum_{|\alpha|=m} \frac{|\alpha|!}{\alpha!}x^{\alpha}=\\&
=\sum_{m=0}^{\infty}(x_1+x_2+\cdots+x_n)^m=\\&=\frac{1}{1-(x_1+x_2+\cdots+x_n)},
\end{aligned}
\end{equation*}
from which we get \eqref{1-21C}.

Concerning \eqref{2-21C}, it suffices to note that by
\eqref{1-21C},  we have, for $\sum_{j=1}^n|x_j|<1$,

\begin{equation*}
\begin{aligned}
\sum_{\alpha\geq
\beta}\frac{|\alpha|!}{(\alpha-\beta)!}x^{\alpha-\beta}&=\sum_{\alpha\in
\mathbb{N}^n_0}
\frac{|\alpha|!}{\alpha!}\partial^{\beta}(x^{\alpha})=\\&
=\partial^{\beta}(\frac{1}{1-(x_1+x_2+\cdots+x_n)})=\\&=\frac{|\beta|!}{\left(1-(x_1+x_2+\cdots+x_n)\right)^{1+|\beta|}}.
\end{aligned}
\end{equation*}

\bigskip

\section{Analytic functions in an open set of $\mathbb{R}^n$}\label{funz-analitiche-in aperto}
\begin{definition}\label{def1-22C}
	\index{Definition:@{Definition:}!- analytic function of real variables@{- analytic function of real variables}}
Let $\Omega$ be an open set of  $\mathbb{R}^n$ and $f:\Omega\rightarrow
\mathbb{R}$ (or $\mathbb{C}$) a function. We say that $f$ is
an \textbf{analytic function of real variables} in $x_0\in \Omega$ if there exist a neighborhood $\mathcal{U}_{x_0}$ of $x_0$, $c_{\alpha}\in \mathbb{R}$
($\mathbb{C}$),$\alpha\in \mathbb{N}_{0}^n$, such that
$$f(x)=\sum_{\alpha\in \mathbb{N}^n_0}c_{\alpha}(x-x_0)^{\alpha},\quad\forall x\in \mathcal{U}_{x_0}.$$
We say that $f$ is a \textbf{analytic function of real variables} in
$\Omega$ if it is  analytic function of real variables in every
$x_0\in\Omega$.

We say that $f:\Omega\rightarrow \mathbb{R}^m$ (or
$\mathbb{C}^m$) where $m\in\mathbb{N}$, $f=(f_1,\cdots,f_m)$ is
an \textbf{analytic function of real variables} in $x_0\in \Omega$ (or in
$\Omega$) provided $f_j$ are analytic functions of real variables in $x_0\in
\Omega$ (or in $\Omega$).
\end{definition}

\bigskip

In what follows, if there is no ambiguity, we will simply say
\textbf{"analytic functions"}, omitting the expression "of real variables".

We will denote by  $C^{\omega}(\Omega; \mathbb{R}^m)$
($C^{\omega}(\Omega; \mathbb{C}^m)$) \index{$C^{\omega}(\Omega; \mathbb{R}^m)$, ($C^{\omega}(\Omega; \mathbb{C}^m)$)} the class of analytic function defined on $\Omega$ with values in  $\mathbb{R}^m$ ($\mathbb{C}^m$), if $m=1$ we write, if there is no ambiguity,
simply $C^{\omega}(\Omega)$ to denote
$C^{\omega}(\Omega; \mathbb{R})$ ($C^{\omega}(\Omega; \mathbb{C})$).

 Proposition \ref{prop1-17C} immediately gives:

\begin{equation}\label{1-22C}
C^{\omega}(\Omega)\subset C^{\infty}(\Omega).
\end{equation}
Moreover, for every $x_0\in \Omega$ there exists a neighborhood 
$\mathcal{U}_{x_0}$ of $x_0$ such that
\begin{equation}\label{2-22C}
f(x)=\sum_{\alpha\in
\mathbb{N}^n_0}\frac{1}{\alpha!}\partial^{\alpha}f(x_0)(x-x_0)^{\alpha},\quad\forall
x\in \mathcal{U}_{x_0}
\end{equation}
and there exist $M>0$, $r>0$ (depending on $x_0$) and
$\widetilde{\mathcal{U}}_{x_0}$, neighborhood of $x_0$, with
$\widetilde{\mathcal{U}}_{x_0}\Subset \mathcal{U}_{x_0}$, such that
\begin{equation}\label{3-22C}
\left|\partial^{\alpha}f(x)\right|\leq
M|\alpha|!r^{-|\alpha|},\quad\forall\alpha\in \mathbb{N}_0^n,\quad
\forall x\in \mathcal{U}_{x_0},
\end{equation}
(\eqref{3-22C} follows by \eqref{2-20C}).

We recall that the inclusion \eqref{1-22C} is proper. As a matter of fact, the
function

\begin{equation}\label{4-22C}
f(t)=
\begin{cases}
e^{-1/t^2}, \quad \mbox{for } t\neq 0,\\
\\
0, \quad\quad\quad \mbox{for } t=0, %
\end{cases}
\end{equation}
belongs to $C^{\infty}(\mathbb{R})$, but it is not analytic. As a matter of fact, since we have
 $f^{(k)}(0)=0$, for every $k\in \mathbb{N}_0$, we have

$$\sum_{k=0}^{\infty}\frac{1}{k!}f^{(k)}(0)t^k=0\neq f(t),\quad\forall
t\in\mathbb{R}\setminus\{0\}.$$ Therefore, for every $t \neq 0$, $f(t)$
is different from the sum of its Taylor series.

\bigskip

The analytic functions enjoy the \textbf{unique continuation property.}\index{unique continuation property} Indeed, the following holds true.

\begin{theo}\label{teo1-23C}
Let $\Omega$ be a connected open set of $\mathbb{R}^n$ and $f\in
C^{\omega}(\Omega)$. Let $x_0\in \Omega$. Then we have

$$\partial^{\alpha}f(x_0)=0,\quad \forall \alpha\in \mathbb{N}_0^n\quad\Longrightarrow \quad f\equiv 0, \mbox{ in } \Omega.$$
In particular, if $f$ vanishes identically in an open set (not empty) of
$\Omega$ then $f$ vanishes identically in $\Omega$.
\end{theo}

\textbf{Proof.} Let
$$\widetilde{\Omega}=\left\{x\in \Omega:\mbox{ } \partial^{\alpha}f(x)=0,\quad \forall \alpha\in \mathbb{N}_0^n \right\}.$$
It is clear that $\widetilde{\Omega}\neq\emptyset$ because
$x_0\in\widetilde{\Omega}$. Therefore, whether
 we prove that $\widetilde{\Omega}$ is at the same time an open and a closed set in $\Omega$ (in the topology induced by  $\mathbb{R}^n$),
 as $\Omega$ is connected, we have $\Omega=\widetilde{\Omega}$ and the thesis follows.

Since $f$ is continuous, the set $\widetilde{\Omega}$ is closed in)
$\Omega$. Now we prove that $\widetilde{\Omega}$ is an open set of
$\Omega$.

Let $\widetilde{x}\in \widetilde{\Omega}$. By the analyticity of $f$ we have that there exists a neighborhood $\mathcal{U}_{\widetilde{x}}$ of $\widetilde{x}$ such that

$$f(x)=\sum_{\alpha\in \mathbb{N}^n_0}\frac{1}{\alpha!}\partial^{\alpha}f(\widetilde{x})(x-\widetilde{x})^{\alpha},
\quad\forall x\in \mathcal{U}_{\widetilde{x}}.$$ On the other hand,
$\partial^{\alpha}f(\widetilde{x})=0$ for every $\alpha\in
\mathbb{N}_0^n$, hence $f(x)=0$ for every $x\in
\mathcal{U}_{\widetilde{x}}$. Hence
$$\partial^{\alpha}f(x)=0,\quad \forall \alpha\in \mathbb{N}_0^n,\quad \forall x\in \mathcal{U}_{\widetilde{x}},$$
from which $\mathcal{U}_{\widetilde{x}}\subset \widetilde{\Omega}$.
Therefore $\Omega=\widetilde{\Omega}$. $\blacksquare$

\bigskip

In Theorem \ref{teo2-24C} we will prove that the analytic functions
can be characterized by the growth of their derivatives.
We premise the following

\begin{definition}\label{def2-24C}
Let $m\in \mathbb{N}$ and let $\Omega$ be an open set of $\mathbb{R}^n$. Let
$f:\Omega\rightarrow \mathbb{R}^m$ (or $\mathbb{C}^m$),
$f=(f_1,\cdots,f_m)$. Let $x_0\in \Omega$ and $M, r>0$,
$j=1,\cdots,m$. We write $$f\in\mathcal{C}_{M,r}(x_0),$$ \index{$\mathcal{C}_{M,r}(x_0)$}provided that $f\in
C^{\infty}$ in a neighborhood of $x_0$ and we have
$$\left|\partial^{\alpha}f_j(x_0)\right|\leq M|\alpha|!r^{-|\alpha|},\quad\forall\alpha\in \mathbb{N}_0^n,\quad j=1,\cdots,m.$$
\end{definition}

\bigskip

\begin{theo}\label{teo2-24C}
Let $\Omega$ be an open set of $\mathbb{R}^n$ and $f:\Omega\rightarrow
\mathbb{R}^m$ (or $\mathbb{C}^m$). The following conditions are equivalent:

(i) $f$ is analytic in $\Omega$,

(ii) for any compact $K$, $K\subset\Omega$, there exist  $M,r>0$
(depending on $K$) such that
$$f\in\mathcal{C}_{M,r}(x_0),\quad\forall x_0\in K.$$
\end{theo}

\textbf{Proof.} It suffices to consider the case $m=1$.

We prove that (i)$\Rightarrow$(ii).

By \eqref{3-22C} we know that if $y\in \Omega$, there exist
$M_y,r_y>0$ and a neighborhood $\mathcal{U}_y$ of $y$ such that

\begin{equation}\label{1-25C}
\sup_{\mathcal{U}_y}\left|\partial^{\alpha}f\right|\leq
M_y|\alpha|!r_y^{-|\alpha|},\quad\forall\alpha\in \mathbb{N}_0^n.
\end{equation}
Let $K\subset\Omega$ be a compact, then
$\left\{\mathcal{U}_y\right\}_{y\in K}$ is an open covering of $K$,
hence there exist $\mathcal{U}_{y_1},\cdots, \mathcal{U}_{y_N}$ such that $$K\subset \bigcup_{j=1}^N\mathcal{U}_{y_j}.$$ Set
$$M=\max_{1\leq j\leq N}M_{y_j}, \quad r=\min_{1\leq j\leq N}r_{y_j}.$$
By \eqref{1-25C} we have
\begin{equation}\label{2-25C}
\sup_{K}\left|\partial^{\alpha}f\right|\leq\sup_{\bigcup_{j=1}^N\mathcal{U}_{y_j}}\left|\partial^{\alpha}f\right|\leq
M|\alpha|!r^{-|\alpha|}, \quad\forall\alpha\in \mathbb{N}_0^n.
\end{equation}

\medskip

We prove that (ii)$\Rightarrow$(i).

Let us suppose that (ii) holds. Let $x_0\in \Omega$ and $\rho>0$ satisfy
$\overline{B_{\rho}(x_0)}\subset\Omega$. It is not restrictive to assume that $x_0=0\in \Omega$. Let us choose 
$K=\overline{B_{\rho}}$ and let $M,r>0$ satisfy

\begin{equation}\label{3-25C}
\left|\partial^{\alpha}f(x)\right|\leq
M|\alpha|!r^{-|\alpha|},\quad\forall x\in
\overline{B_{\rho}},\quad \forall\alpha\in\mathbb{N}_0^n.
\end{equation}
Set $d$, $0<d<\min\{r,\rho\}$, we now prove 

$$f(x)=\sum_{\alpha\in \mathbb{N}^n_0}\frac{1}{\alpha!}\partial^{\alpha}f(0)x^{\alpha},\quad\mbox{ for every } x \mbox{ such that }
\sum_{j=1}^n\left|x_j\right|\leq d.$$

Let $x$ satisfy $\sum_{j=1}^n\left|x_j\right|\leq d$ and set
$$\phi(t)=f(tx), \quad t\in [0,1].$$
We have, for any $m\in \mathbb{N}$,
\begin{equation}\label{1-26C}
f(x)=\phi(1)=\sum_{k=0}^{m-1}\frac{1}{k!}\phi^{(k)}(0)+\mathcal{R}_m=
\sum_{|\alpha|\leq
m-1}\frac{1}{\alpha!}\partial^{\alpha}f(0)x^{\alpha}+\mathcal{R}_m.
\end{equation}
where

$$\mathcal{R}_m=\frac{1}{(m-1)!}\int^1_0(1-t)^{(m-1)}\phi^{(m)}(t)dt.$$
By \eqref{3-3N}, \eqref{3-25C} and by
$\sum_{j=1}^n\left|x_j\right|\leq d$ we have

\begin{equation*}
\begin{aligned}
\left|\frac{1}{m!}\phi^{(m)}(t)\right|&=\left|\sum_{|\alpha|=m}\frac{1}{\alpha!}\partial^{\alpha}f(0)x^{\alpha}\right|\leq\\&
\leq\sum_{|\alpha|=m}\frac{|\alpha|!}{\alpha!}M
r^{-|\alpha|}\left|x^{\alpha}\right|=\\&
=Mr^{-m}\sum_{|\alpha|=m}\frac{|\alpha|!}{\alpha!}\left|x_1\right|^{\alpha_1}\cdots\left|x_n\right|^{\alpha_n}=\\&
=Mr^{-m}\left(\sum_{j=1}^n\left|x_j\right|\right)^m\leq \\&\leq
M\left(\frac{d}{r}\right)^m.
\end{aligned}
\end{equation*}
Therefore
\begin{equation}\label{1-27C}
\begin{aligned}
\left|\mathcal{R}_m\right|&\leq\frac{1}{(m-1)!}\int^1_0(1-t)^{(m-1)}\left|\phi^{(m)}(t)\right|dt\leq\\&
\leq\frac{1}{(m-1)!}\int^1_0(1-t)^{(m-1)}m!M
\left(\frac{d}{r}\right)^mdt=\\&=M\left(\frac{d}{r}\right)^m.
\end{aligned}
\end{equation}

All in all, by \eqref{1-26C} and \eqref{1-27C} we get, if
$\sum_{j=1}^n\left|x_j\right|\leq d$ (recall $d<r$),

$$\left\vert f(x)-\sum_{|\alpha|\leq m}\frac{1}{\alpha!}\partial^{\alpha}f(0)x^{\alpha}\right\vert
\leq M\left(\frac{d}{r}\right)^m\rightarrow 0, \quad\mbox{as }
m\rightarrow \infty.$$ $\blacksquare$

\bigskip

\section{Majorant functions}\label{funz-analitiche-maggioranti}
In the proof  of the Cauchy--Kowalevkaya Theorem we will make use the \textbf{method of majorant series}. We give the following

\begin{definition}\label{def3-24C}
	\index{Definition:@{Definition:}!- majorant function@{- majorant function}}
Let $m\in \mathbb{N}$. Let $\mathcal{U}_{x_0}$ be a neighborhood of
$x_0\in\mathbb{R}^n$ and let $f:\mathcal{U}_{x_0}\rightarrow
\mathbb{R}^m$ (or $\mathbb{C}^m$), $F:\mathcal{U}_{x_0}\rightarrow
\mathbb{R}^m$, $f=(f_1,\cdots,f_m)$,
$F=(F_1,\cdots,F_m)$. Let us suppose that $f_j,F_j\in
C^{\infty}(\mathcal{U}_{x_0})$, $j=1,\cdots,m$. We say that $F$ is \textbf{a majorant} of $f$ or $F$ 
\textbf{majorazes} $f$  \textbf{in} $x_0$ and we write \index{$\preccurlyeq$}
$$f\preccurlyeq F, \quad \mbox{ in } x_0,$$
provided
$$\left|\partial^{\alpha}f_j(x_0)\right|\leq \partial^{\alpha}F_j(x_0),\quad\forall\alpha\in \mathbb{N}_0^n,\quad j=1,\cdots,m.$$
\end{definition}

\medskip

\textbf{Remark.} Notice that if $f\preccurlyeq F$ in $x_0$, then,  we have $f\preccurlyeq \widetilde{F}$ in
$x_0$, where
$$\widetilde{F}(x)=F(a^{-1}_1x_1,\cdots, a^{-1}_nx_n)$$ for every $a_j\in
(0,1]$, $j=1,\cdots,n$. $\blacklozenge$

\medskip

In what follows we will assume, without any restriction that  $x_0=0$
and we will write simply (if there is no ambiguity)
$f\preccurlyeq F$ instead of $f\preccurlyeq F$ in $0$.

\bigskip

\begin{prop}\label{prop1-28C}
Let $f:\mathcal{U}_{0}\rightarrow \mathbb{C}^m$, where
$\mathcal{U}_{0}$ is a neighborhood of $0\in\mathbb{R}^n$ and, for given $M,
r>0$, let
$$\phi_{M,r}(x)=\frac{Mr}{r-(x_1+\cdots+x_n)}.$$ Then we have

\smallskip

(i) $f\in \mathcal{C}_{M,r}(0)$ if and only if $f\preccurlyeq
\upsilon_m \phi_{M,r}$ (here $\upsilon_m=\underset{m
}{(\underbrace{1\cdots,1})}$),

\smallskip

(ii) $f\in \mathcal{C}_{M,r}(0)$ and $f(0)=0$ if and only if
$f\preccurlyeq \upsilon_m\left( \phi_{M,r}-M\right)$.
\end{prop}
\textbf{Proof}

 (i) is an immediate consequence of \eqref{1-21C}. We have,
indeed,

$$\phi_{M,r}(x)=\frac{M}{1-(\frac{x_1+\cdots+x_n}{r})}=\sum_{\alpha\in \mathbb{N}^n_0}\frac{M|\alpha|!r^{-|\alpha|}}{\alpha!}x^{\alpha}.$$
Hence $\partial^{\alpha}\phi_{M,r}(0)=M|\alpha|!r^{-|\alpha|}$, from which we get
 (i).  (ii) is a trivial consequence  of (i).
$\blacksquare$

\bigskip

From the derivation rules we have

\begin{prop}\label{prop0-29C}
Let $f,g:\mathcal{U}_{0}\rightarrow \mathbb{C}$; $F,G:\mathcal{U}_{0}\rightarrow \mathbb{R}$ where
$\mathcal{U}_{0}$ is a neighborhood of $0\in\mathbb{R}^n$. If,
$$f\preccurlyeq F, \quad\mbox{and}\quad g\preccurlyeq G$$
then

$$f+g\preccurlyeq F+G,\quad fg\preccurlyeq FG.$$
\end{prop}
\textbf{Proof.} The proof is left as an exercise to the reader. $\blacksquare$
\bigskip

\begin{prop}\label{prop1-29C}
Let $n,m,p\in\mathbb{N}$ and 

$$f,F:\mathcal{U}_0\subset\mathbb{R}^n\rightarrow \mathbb{R}^m,$$
where $\mathcal{U}_0$ is a neighborhood  of $0$. Let us suppose 
$$f(0)=F(0)=0,$$ and let $\mathcal{V}$ be an open set of $\mathbb{R}^m$ such that $\mathcal{U}_0\subset\mathcal{V}$. Let

$$g,G:\mathcal{V}\rightarrow \mathbb{R}^p.$$
Let us assume that $f,F,g,G$ be of class $C^{\infty}$ and 
$$f\preccurlyeq F, \quad g\preccurlyeq G.$$
Then we have
$$g\circ f \preccurlyeq G\circ F.$$
\end{prop}
\textbf{Proof.} Set $h=g\circ f$ e $H=G\circ F$. By
the chain rule we have that for every $\alpha\in
\mathbb{N}_0^n$ there exists a polynomial $P_{\alpha}$ with positive coefficients \textit{indipendent} of $f,g,F,G$ such that, for
$j=1,\cdots,p,$ we have

$$\partial^{\alpha}h_j(0)=P_{\alpha}\left(\partial^{\beta}g_l(0),\cdots, \partial^{\gamma}f_k(0)\right),$$

$$\partial^{\alpha}H_j(0)=P_{\alpha}\left(\partial^{\beta}G_l(0),\cdots, \partial^{\gamma}F_k(0)\right).$$

Since

$$\left|\partial^{\gamma}f_k(0)\right|\leq \partial^{\gamma}F_k(0), \quad \forall \beta\in \mathbb{N}_0^n,\quad k=1,\cdots,m,$$

$$\left|\partial^{\beta}g_l(0)\right|\leq \partial^{\beta}G_l(0), \quad \forall \beta\in \mathbb{N}_0^m,\quad l=1,\cdots,p,$$
and taking into account that the coefficients of $P_{\alpha}$ are positive, we get

\begin{equation*}
\begin{aligned}
\left|\partial^{\alpha}h_j(0)\right|&\leq
P_{\alpha}\left(\left|\partial^{\beta}g_l(0)\right|,\cdots,
\left|\partial^{\gamma}f_k(0)\right|\right)\leq\\& \leq
P_{\alpha}\left(\left|\partial^{\beta}G_l(0)\right|,\cdots,
\left|\partial^{\gamma}F_k(0)\right|\right)=\\&
=\partial^{\alpha}H_j(0).
\end{aligned}
\end{equation*}
$\blacksquare$

\bigskip

\begin{prop}\label{prop1-30C}
Let $n,m,p\in\mathbb{N}$. Let $M,r, \mu,\rho$ be  positive numbers and $\widetilde{x}\in \mathbb{R}^n$, $\widetilde{y}\in
\mathbb{R}^m$. Let $f\in \mathcal{C}_{M,r}(\widetilde{x})$ be a function with values in
$\mathbb{R}^m$, $g\in \mathcal{C}_{\mu,\rho}(\widetilde{y})$ be a function with
values in $\mathbb{R}^p$, $\widetilde{y}=f(\widetilde{x})$. Then
\begin{equation}\label{1-30C}
h:=g\circ f\in \mathcal{C}_{\mu,\frac{\rho
r}{Mm+\rho}}(\widetilde{x}).
\end{equation}
\end{prop}
\textbf{Proof.} Set

$$g^*(y)=g(y+\widetilde{y}),\quad f^*(x)=f(x+\widetilde{x})-f(\widetilde{x}).$$
We have
$$h(x+\widetilde{x})=g\left(\widetilde{y}+f(x+\widetilde{x})-f(\widetilde{x})\right)=g^*\left(f^*(x)\right)$$
and

$$f^*\in \mathcal{C}_{M,r}(0),\quad g^*\in \mathcal{C}_{\mu,\rho}(0).$$
Proposition \ref{prop1-28C} implies

\begin{equation}\label{1-31C}
f^*\preccurlyeq \upsilon_m\left(\phi_{M,r}-M\right).
\end{equation}

\begin{equation}\label{2-31C}
g^*\preccurlyeq \upsilon_p\phi_{\mu,\rho}.
\end{equation}
Set
$$ \chi(x)=\phi_{\mu,\rho}\left(\left(\phi_{M,r}(x)-M\right)\right).$$
By Proposition \ref{prop1-29C}, by \eqref{1-31C} and
\eqref{2-31C} we have

\begin{equation}\label{3-31C}
h(x+\widetilde{x})=\left(g^*\circ f^*\right)(0)\preccurlyeq
\upsilon_p\chi(x).
\end{equation}
Moreover it is simple to obtain
\begin{equation*}
\begin{aligned}
\chi(x)&=\phi_{\mu,\rho}\left(\phi_{M,r}(x)-M\right)=\\&
=\frac{\mu\rho}{\rho-m\left(\phi_{M,r}(x)-M\right)}=\\&
=\frac{\mu\rho\left(r-(x_1+\cdots+x_n)\right)}{\rho
r-(\rho+mM)(x_1+\cdots+x_n)}.
\end{aligned}
\end{equation*}
Now we check that

\begin{equation}\label{1-32C}
\chi(x)\preccurlyeq \frac{\mu\rho r}{\rho
r-(\rho+mM)(x_1+\cdots+x_n)}.
\end{equation}
To this purpose, set $A=\rho r$, $B=\rho+mM$, $t=x_1+\cdots+x_n$,
$$\widetilde{\chi}(t)=\mu\rho\frac{r-t}{A-Bt}.$$ We have

\begin{equation}\label{0-32C}
\chi(x)= \widetilde{\chi}(x_1+\cdots+x_n).
\end{equation}
and

\begin{equation}\label{00-32C}
\begin{aligned}
\frac{r-t}{A-Bt}&=\frac{r/A}{1-Bt/A}-\frac{t}{A}\frac{1}{\left(1-\frac{B}{A}t\right)}=\\&
=\frac{r}{A}\sum_{k=0}^{\infty}\left(\frac{Bt}{A}\right)^k-\frac{t}{A}\sum_{k=0}^{\infty}\left(\frac{Bt}{A}\right)^k=\\&
=\frac{r}{A}+\frac{1}{A}\sum_{k=1}^{\infty}\left[\left(\frac{B}{A}\right)^kr-\left(\frac{B}{A}\right)^{k-1}\right]t^k\preccurlyeq\\&
\preccurlyeq
\frac{r}{A}+\frac{r}{A}\sum_{k=1}^{\infty}\left(\frac{B}{A}\right)^kt^k=\frac{r}{A-Bt}.
\end{aligned}
\end{equation}
By \eqref{0-32C} and \eqref{00-32C}
we obtain  \eqref{1-32C}. Finally, by \eqref{3-31C}, \eqref{1-32C} and by Proposition
\ref{prop1-28C} we get \eqref{1-30C}. $\blacksquare$

\bigskip

\begin{theo}
Let $n,m,p\in\mathbb{N}$. Let $\Omega_1$ be an open set of 
$\mathbb{R}^n$, and $\Omega_2$ be an open set of $\mathbb{R}^m$. Let $f\in
C^{\omega}(\Omega_1,\mathbb{R}^m)$ satisfy $f(\Omega_1)\subset
\Omega_2$. Let $g\in C^{\omega}(\Omega_2,\mathbb{R}^p)$. Then 
$g\circ f\in C^{\omega}(\Omega_1,\mathbb{R}^p)$.
\end{theo}

\textbf{Proof.} Is an immediate consequence of Proposition \ref{prop1-30C} and of
Theorem \ref{teo2-24C}. $\blacksquare$

\bigskip

In what follows we will use the \textbf{Inverse Function Theorem} 
and the \textbf{Implicit Function Theorem for analytic functions} \index{Implicit Function Theorem for analytic functions}.

For instance, we will exploit the following fact. If
$\phi:\Omega\rightarrow\mathbb{R}$ is an analytic function in $\Omega$ and
$$\nabla\phi(x_0)\neq 0,$$  $x_0\in \Omega$, then there exist
$r,\delta>0$  and an isometry $$ \Psi:\mathbb{R}^n\rightarrow \mathbb{R}^n,$$ such that 

$$\Psi(0)=x_0,$$ 

\begin{equation*}
	\Psi^{-1}\left(\left\{x\in \Omega: \phi(x)=\phi(x_0)\right\}\right)\cap Q_{r, 2M}=\left\{(x',\varphi(x')):x'\in B'_r\right\} 
\end{equation*}
where $\varphi\in C^{\omega}\left(B'_r; \mathbb{R}\right)$ and it satisfies

$$\varphi(0)=0,\quad |\nabla \varphi(0)|=0, $$
and
$$\left\Vert \varphi\right\Vert_{C^1\left(\overline{B'_{r_0}}\right)}\leq
Mr.$$ 
We will not prove this Theorem which can be proved by the method of the majorant functions that
we will learn to use in the next Chapter.



\chapter[The Cauchy problem for PDEs with
analytic coefficients]{The Cauchy problem for PDEs with
analytic coefficients}\label{Cauchy-chap} 
\section{Formulation of the Cauchy problem} \label{formulazione-cauchy}
In this Section we will give a fairly general formulation of the
Cauchy problem. Although we are mainly interested in the
linear operators, the formulation that we will give also applies to the fully
nonlinear operators.

Let $\Omega$ be a connected open set of  $\mathbb{R}^n$. Let $x_0\in \Omega$,
$\phi\in C^m(\Omega;\mathbb{R})$, where $m\in \mathbb{N}$. Let us suppose that
\begin{equation}\label{0-5C}
\phi(x_0)=0.
\end{equation}
Set

\begin{equation}\label{1-5C}
\Gamma=\left\{x\in \Omega: \phi(x)=0 \right\}
\end{equation}
and let us assume

\begin{equation}\label{2-5C}
\nabla\phi(x)\neq 0,\quad \forall x\in\Gamma.
\end{equation}
Let us denote by
\begin{equation}\label{3-5C}
\nu(x)=-\frac{\nabla\phi(x)}{\left\vert
\nabla\phi(x)\right\vert},\quad \forall x\in\Gamma.
\end{equation}

Let be given the function  $g_0,g_1,\cdots,g_{m-1}$, defined on
$\Gamma$, and let $F\left(x,(p_{\alpha})_{|\alpha|\leq
m}\right)$ be a function defined on $\Omega\times \mathbb{R}^{N_m}$, where $N_m\in
\mathbb{N}$ depends on $m$ only. The \textbf{Cauchy problem} is formulated as follows. 
\index{Cauchy problem:@{Cauchy problem:}!- general equations@{- general equations}}

\smallskip

Determine $u$ of class $C^m$ in a neighborhood $\mathcal{U}$ of $x_0$
such that

\begin{equation}\label{4-1C}
\begin{cases}
F\left(x,(\partial^{\alpha} u)_{|\alpha|\leq m}\right)=0, \quad \forall x\in \mathcal{U}, \\
\\
\frac{\partial^ju(x)}{\partial \nu^j}=g_j(x), \quad
j=0,1,\cdots,m-1, \mbox{ } \forall x\in \Gamma\cap\mathcal{U}.
\end{cases}
\end{equation}

\smallskip

The functions $g_0,g_1,\cdots,g_{m-1}$ and $\Gamma$ are called, respectively, the \textbf{initial data} \index{initial data (values) of the Cauchy problem}or the \textbf{initial values} and the \textbf{initial surface} \index{initial surface of the Cauchy problem} of Cauchy problem \eqref{4-1C}. The equations
\begin{equation}\label{4a-1C}
\frac{\partial^ju(x)}{\partial \nu^j}=g_j(x), \quad
j=0,1,\cdots,m-1, \mbox{ } \forall x\in \Gamma\cap\mathcal{U},
\end{equation}
are called the \textbf{initial conditions of the Cauchy problem} \index{initial conditions of the Cauchy problem}. Of course, it makes sense and interest to set more general initial conditions. For instance, instead of the vector field
$\nu(x)$, we may consider a vector field  $\ell(x)$ in a neighborhood
of $\Gamma$, of class $C^{m-1}$, requiring
that
\begin{equation}\label{4b-1C}
\frac{\partial^ju(x)}{\partial \ell^j}=g_j(x), \quad
j=0,1,\cdots,m-1, \mbox{ }\forall x\in \Gamma\cap\mathcal{U}.
\end{equation}
We can easly check that if the vector field $\ell(x)$ and the functions
$g_j$ are smooth enough and if $\ell(x)\cdot\nu(x)\neq 0$,
for any $x\in \Gamma$, then conditions \eqref{4a-1C} and \eqref{4b-1C} are
equivalent. To realize this, let us consider the simple case where $m=2$, $\Omega=\mathbb{R}^n$, $\phi(x)=x_n$, hence
$\nu(x)=-e_n$ for every $x\in \Gamma$ and $\ell=(\ell',\ell_n)$, where
$\ell'=(\ell_1,\cdots,\ell_{n-1})$, is a vector field such that
$\ell_n(x')\neq 0$ for every $x'\in \mathbb{R}^{n-1}$. Let us assume
$u\in C^2(\mathbb{R}^n)$ and

\begin{equation}\label{4c-1C}
u(x',0)=g_0(x'), \quad \partial_nu(x',0)=-g_1(x'),\quad \forall
x'\in \mathbb{R}^{n-1},
\end{equation}
where $g_0\in C^1(\mathbb{R}^{n-1})$ and $g_1\in
C^0(\mathbb{R}^{n-1})$. By the first equation in \eqref{4c-1C} we have
$$\nabla_{x'} u(x',0)=\nabla_{x'} g_0(x'), \quad \forall x'\in \mathbb{R}^{n-1} $$ that,  together with $\partial_n u(x',0)=-g_1(x')$, gives 
\begin{equation}\label{4ca-1C}
\frac{\partial u}{\partial \ell}(x',0)=\ell'\cdot \nabla_{x'}
g_0(x')-\ell_ng_1(x'), \quad \mbox{ } \forall x'\in \mathbb{R}^{n-1}.
\end{equation}
Conversely, let us suppose that 
\begin{equation}\label{4d-1C}
u(x',0)=\widetilde{g}_0(x'), \quad
\frac{\partial}{\partial\ell}u(x',0)=\widetilde{g}_1(x'),\quad
\forall x'\in \mathbb{R}^{n-1},
\end{equation}
where $\widetilde{g}_0\in C^1(\mathbb{R}^{n-1})$ and
$\widetilde{g}_1\in C^0(\mathbb{R}^{n-1})$. We have, by the first equation in
\eqref{4d-1C} and taking into account that $\ell_n\neq 0$,
$$\partial_nu(x',0)=\frac{1}{\ell_n(x')}\left(-\ell'(x')\cdot\nabla_{x'} \widetilde{g}_0(x')+\widetilde{g}_1(x')\right). $$
We also notice that if $\ell_n=0$ at some point $x'_0\in
\mathbb{R}^{n-1}$ then between \eqref{4c-1C} and \eqref{4d-1C} there is no
equivalence. Actually, if \eqref{4c-1C} holds, then
we can equally get \eqref{4ca-1C}, but by
\eqref{4d-1C} we see that between $\widetilde{g}_0$ and $\widetilde{g}_1$
the following condition of compatibility needs to be fulfilled

$$-\ell'(x'_0)\cdot\nabla_{x'} \widetilde{g}_0(x'_0)+\widetilde{g}_1(x'_0)=0.$$
On the other hand, if this condition is satisfied, it is undeterminate
the value of $\partial_nu(x'_0,0)$.

\section{The characteristic surfaces} \label{superf-caratteristiche}
The notion of the \textbf{characteristic surface} \index{characteristic:@{characteristic:}!- surface@{- surface}}has a fundamental importance in the investigation of the Cauchy problem. Roughly speaking, we say that the \textbf{surface} $\Gamma$
is noncharacteristic if, assuming $\phi, u, F, g_j\in
C^{\infty}$, \\ for $j=0,1,\cdots,m-1$,  all the
derivatives of $u$ on $\Gamma$ can be determined from  by \eqref{4-1C}.

Of course, we need to specify this notion and arrive at a
formal definition, however in what follows we will not tackle
problem \eqref{4-1C} in its full generality, but we will
limit ourselves to the linear case, i.e.

$$F\left(x,(\partial^{\alpha} u)_{|\alpha|\leq m}\right)=P(x,\partial)u-f(x)=\sum_{|\alpha|\leq m}a_{\alpha}(x)\partial^{\alpha}u-f(x),$$
Where we will assume, unless explicitly otherwise stated, 
\begin{equation}\label{0-2C}
\phi,a_{\alpha}, f\in C^{\infty}(\Omega), \quad |\alpha|\leq m
\end{equation}
and that \eqref{2-5C} holds.

To motivate the definition of a characteristic surface that we will give,
we begin by considering the following Cauchy problem

\begin{equation}\label{00-2C}
\begin{cases}
P(x,\partial)u=f(x), \quad \forall x\in B_R, \\
\\
\partial_n^ju(x',0)=g_j(x'), \quad j=0,1,\cdots,m-1,\mbox{ } \forall x'\in
B'_R.
\end{cases}
\end{equation}
Here $\phi(x)=-x_n$, $\Gamma=\{x\in
B_R: x_n=0\}$.

Now, if $g_j$, $j=0,1,\cdots,m-1$, are of class  $C^{\infty}$,
and if there exists a solution $u\in C^{\infty}(B_R)$ to problem
\eqref{00-2C}, we have, by the initial condizions, 

\begin{equation}\label{1-3C}
\partial^{\alpha'}\partial_n^ju(x',0)=\partial^{\alpha'}g_j(x'),\quad \forall \alpha'\in \mathbb{N}_0^{n-1}, \quad j=0,1,\cdots,m-1.
\end{equation}
Let us notice that by  \eqref{1-3C} we \textit{cannot} determine the
derivatives

\begin{equation}\label{2-3C}
\partial_n^ju(x',0),  \quad j\geq m
\end{equation}
and a fortiori, we cannot determine the derivatives
$\partial^{\alpha'}\partial_n^ju(x',0)$ for $j\geq m$, $\alpha'\in
\mathbb{N}_0^{n-1}$. To gain such derivatives we should exploit
the equation $$P(x,\partial)u=f(x)$$ which we write in the form

\begin{equation}\label{3-3C}
a_{(0',m)}(x',0)\partial_n^mu(x',0)=-\sum_{|\alpha ^{\prime }|\leq
m,\phantom{n} \alpha
_{n}<m}a_{\alpha}(x',0)\partial^{\alpha}u(x'0)+f(x',0).
\end{equation}
Let us note that to the right--hand side of \eqref{3-3C}, by
\eqref{1-3C}, all the derivatives that appear can be expressed
in terms of the initial data and their derivatives. Therefore, if

\begin{equation}\label{4-3C}
a_{(0',m)}(x',0)\neq 0,\quad \forall x'\in B'_R,
\end{equation}
by \eqref{3-3C} we determine $\partial_n^mu(x',0)$, hence we can calculate
\begin{equation}\label{5-3C}
\partial^{\alpha'}\partial_n^mu(x',0),\quad \forall \alpha'\in \mathbb{N}_0^{n-1}, \forall x'\in B'_R.
\end{equation}
Actually, condition \eqref{4-3C} allows us to deduce all the
derivatives of $u(x',0)$. As a matter of fact, by calculating the derivatives w.r.t.  $x_n$ of both the
sides of \eqref{3-3C}, we have

\begin{equation}\label{1-4C}
\begin{aligned}
a_{(0',m)}(x)\partial_n^{m+1}u(x)&=\underset{\mbox{known for } x_n=0
\mbox{ (by \eqref{5-3C}) }
}{\underbrace{-\partial_na_{(0',m)}(x)\partial_n^{m}u(x)}}-\\&
\\&
-\underset{\mbox{known for } x_n=0 \mbox{ (by \eqref{1-3C}) }
}{\underbrace{\partial_n(\sum_{|\alpha ^{\prime }|\leq
m,\phantom{n}) \alpha
_{n}<m}a_{\alpha}(x)\partial^{\alpha}u(x))}}+\\& \\& +\underset{\mbox{known for } x_n=0
  }{\underbrace{\partial_nf(x)}}.
\end{aligned}
\end{equation}
Again, condition \eqref{4-3C} allows us to derive from
\eqref{1-4C} the derivative
$$\partial_n^{m+1}u(x',0),\quad \forall x'\in B'_R.$$
Of course, we can further make the derivatives of both the sides of
\eqref{1-4C} and we determine $\partial_n^{m+2}u(x',0)$. Iterating the
procedure we obtain the derivatives

$$\partial_n^{k}u(x',0),\quad \forall k\in\mathbb{N}_0,\mbox{ }\forall x'\in B'_R $$
and then calculating the derivatives w.r.t.  $x_1,x_2\cdots,x_{n-1}$ we can
calculate all the derivatives of $u$ at the points $(x',0)\in B'_R$.

At this point we note that condition \eqref{4-3C} can be written
\begin{equation}\label{2-4C}
P_m((x',0),\nu)=i^m\sum_{|\alpha|= m}a_{\alpha}(x',0)\nu^{\alpha}\neq
0,\quad \forall \alpha'\in \mathbb{N}_0^{n-1},\mbox{ }\forall x'\in
B'_R
\end{equation}
being, in this case, $\nu=e_n$ and recalling that
$P_m(x,\xi)=i^m\sum_{|\alpha|= m}a_{\alpha}(x)\xi^{\alpha}$.

\bigskip

Now we give the following general definition
\begin{definition}\label{def-5C}
	\index{Definition:@{Definition:}!- characteristic surface for linear operator of order $m$@{- characteristic surface of linear operator of order $m$}}
Let $\Omega$ be an open set of $\mathbb{R}^n$ and $x_0\in \Omega$. Let
$a_{\alpha}\in C^0(\Omega)$, for any $|\alpha|\leq m$. Let
$$P(x,\partial)=\sum_{|\alpha|\leq m}a_{\alpha}(x)\partial^{\alpha}$$ be
a linear differential operator of order $m$.

We say that $\ell\in \mathbb{R}^n$ is a \textbf{characteristic direction} \index{characteristic:@{characteristic:}!- direction@{- direction}} for the operator $P(x,\partial)$ in $x_0$ if
\begin{equation}\label{24-7-21-4C}
P_m(x_0,\ell)=0.
\end{equation}
Let $\phi\in C^1(\Omega)$ satisfy
$$\nabla \phi(x_0)\neq 0.$$
We say that $\Gamma=\left\{x\in \Omega: \phi(x)=\phi(x_0) \right\}$
is a \textbf{characteristic surface in} $x_0$ for the 
operator $P(x,\partial)$ provided
$\nu(x_0)=-\frac{\nabla\phi(x_0)}{\left\vert
\nabla\phi(x_0)\right\vert}$ is a characteristic direction for
$P(x,\partial)$ in $x_0$ that is, if $$P_m(x_0,\nu(x_0))=i^m\sum_{|\alpha|= m}a_{\alpha}(x_0)\nu^{\alpha}(x_0)=0,$$
or, equivalently, 
$$P_m(x_0,\nabla\phi(x_0)= 0.$$
\end{definition}

\bigskip

In the sequel, we will say that $\Gamma$ is a \textbf{noncharacteristic surface in} $x_0$ for the operator $P(x,\partial)$, provided that
$$P_m(x_0,\nabla\phi(x_0)\neq 0.$$

We will say that $\Gamma$ is a \textbf{characteristic surface} for the operator $P(x,\partial)$ if it is characteristic \textbf{at each point of} $\Gamma$.
Finally, we say that $\Gamma$ is a  \textbf{noncharacteristic surface} \index{noncharacteristic surface}for
the operator $P(x,\partial)$ as long as it is noncharacteristic \textbf{at every point of $\Gamma$}.
Let us notice that being a "non characteristic" is more
restrictive than the negation of "characteristic surface." This little abuse will simplify the form of expression later on.

\bigskip

We notice that the previous definition \textbf{involves only
the principal part of the operator} $P(x,\partial)$.

\bigskip

For completeness, we also give a definition of a non
characteristic surface for the quasilinear operator of order $m$

\begin{equation}\label{quasilinear}
\mathcal{P}(u)=\sum_{|\alpha|=m}a_{\alpha}\left(x,(\partial^{\beta}u)_{||\beta|\le
m-1}\right)\partial^{\alpha}u+a_{0}
\left(x,(\partial^{\beta}u)_{||\beta|\le m-1}\right),
\end{equation}
where $a_{\alpha}\left(x,(p^{\beta})_{||\beta|\le m-1}\right), a_0$
are given functions.

\begin{definition}\label{def-quasilin-5C}
	\index{characteristic:@{characteristic:}!- surface for qualinear operator@{- surface for quasilinear operator}}
Let $\Omega$, $x_0$, $\phi$, $\Gamma$ be like in Definition
\ref{def-5C} and let $\mathcal{P}$ be like  in \eqref{quasilinear}. 
We say
that $\Gamma$ is a \textbf{noncharacteristic surface in}
$x_0$ for the operator $\mathcal{P}$ if $\Gamma$ is
a noncharacteristic surface in $x_0$ for
to the operator
$$\sum_{|\alpha|=m}a_{\alpha}\left(x,(p^{\beta})_{||\beta|\le m-1}\right)\partial^{\alpha}$$
for each value of $p^{\beta}$, for $|\beta|\le m-1$. We say that
$\Gamma$ is a  \textbf{noncharacteristic surface} for the operator $\mathcal{P}$ if it is a noncharacteristic in each point of $\Gamma$.
\end{definition}

\bigskip

\textbf{Remark.} Definition \ref{def-quasilin-5C} is actually more restrictive than that would be needed to determine $\partial^{\gamma}u(x',0)$ for all $\gamma\in \mathbb{N}_0^n$. Let us consider, for instance, the Cauchy problem
\begin{equation}\label{1-6C}
\begin{cases}
\mathcal{P}(u)=0, \quad \forall x\in B_R, \\
\\
\partial_n^ju(x',0)=g_j(x'), \quad j=0,1,\cdots,m-1, \mbox{  } \forall x'\in
B'_R,
\end{cases}.
\end{equation}
Let us suppose that $a_{\alpha}, a_0, g_j$ are functions $C^{\infty}$. We have proved before that the 
	derivatives $\partial^{\gamma}u(x',0)$, for $|\gamma|\leq m-1$, depend on the initial data only. Hence the value of
$a_{\alpha}\left((x',0),(\partial^{\beta}u(x',0))_{||\beta|\le m-1}\right)$ are determined by the Cauchy data only. Therefore, it should be more natural to say
that $\Gamma$ is a noncharacteristic surface for the operator $\mathcal{P}$ provided that $\Gamma$ is
a noncharacteristic surface in $x_0$ for
to the operator

$$\sum_{|\alpha|=m}a_{\alpha}\left((x',0),(\partial^{\beta}u (x',0) )_{||\beta|\le m-1}\right)\partial^{\alpha}.$$
 $\blacklozenge$

\bigskip

\section{Transformation of a linear differential operator.} \label{trasformazione}
We wish to examine the transformation of the principal part
of the linear differential operator
$$P(x,\partial)=\sum_{|\alpha|\leq m}a_{\alpha}(x)\partial^{\alpha},$$
under the action of $$\Lambda\in C^m(\Omega,\mathbb{R}^n),$$ where
$\Omega$ is an open set of $\mathbb{R}^n$ and
$\Lambda=(\Lambda_1,\cdots,\Lambda_n)$ is a diffeomorphism of
class $C^m$. By this we mean that $\Lambda$ is injective and it
satisfies

$$\det\left(\partial_x \Lambda(x)\right)\neq 0, \quad\forall x\in \Omega,$$
where $\partial_x \Lambda(x)$ is the jacobian matrix of
$\Lambda$. Let $u\in C^m(\Omega)$ and set
$$v(y)=u\left(\Lambda^{-1}(y)\right), \quad\forall y\in\widetilde{\Omega}:=\Lambda(\Omega),$$
we have $v\in C^m(\widetilde{\Omega})$. 

Moreover by
$$u(x)=v\left(\Lambda(x)\right), \quad\forall x\in\Omega,$$
we have
$$\partial_{x_j}u(x)=\sum_{k=1}^n(\partial_{y_k}v)\left(\Lambda(x)\right)\partial_{x_j}\Lambda_k(x), \quad j=1,\cdots,n,\quad \forall x\in\Omega$$
$$\partial_{x_jx_i}^2u(x)=\sum_{h,k=1}^n(\partial_{y_hy_k}v)\left(\Lambda(x)\right)\partial_{x_i}\Lambda_h(x)\partial_{x_j}\Lambda_k(x)+$$

\begin{equation*}
\begin{aligned}
\partial_{x_jx_i}^2u(x)&=\sum_{h,k=1}^n(\partial_{y_hy_k}v)\left(\Lambda(x)\right)\partial_{x_i}\Lambda_h(x)\partial_{x_j}\Lambda_k(x)+\\&
+\sum_{k=1}^n(\partial_{y_k}v)\left(\Lambda(x)\right)\partial_{x_jx_i}\Lambda_k(x)=\\&
=\left(\left((\partial_{x}\Lambda(x))^{t}\partial_{y}\right)_i\left((\partial_{x}\Lambda(x))^{t}\partial_{y}\right)_j\right)v(\Lambda(x))+\\&
+\mbox{(first order terms)}.
\end{aligned}
\end{equation*}

\medskip

\noindent In general we have
\begin{equation}\label{1-11C}
\partial^{\alpha}_{x}u(x)=\left(\left((\partial_{x}\Lambda(x))^{t}\partial_{y}\right)^{\alpha}\right)v(\Lambda(x))
+\mbox{(terms of order less than } |\alpha| \mbox{)}.
\end{equation}
Now, let us denote by $\widetilde{P}(y,\partial_y)$ the transformed operator \index{transformed operator}
of $P(x,\partial_x)$ through $\Lambda$, that is
the operator satisfying

\begin{equation}\label{2-11C}
\left(\widetilde{P}(y,\partial_y)v(y)\right)_{|y=\Lambda(x)}=P(x,\partial_x)u(x).
\end{equation}
By \eqref{2-11C} we have that the principal part of $\widetilde{P}(y,\partial_y)$,
$\widetilde{P}_m(y,\partial_y)$, is given by
\begin{equation}
\begin{aligned}\label{3-11C}
\widetilde{P}_m(y,\partial_y)&=\sum_{|\alpha|=
m}a_{\alpha}(\Lambda^{-1}(y))
\left(\left(\partial_x\Lambda(x)\right)^t\partial_y\right)_{|x=\Lambda^{-1}(y)}^{\alpha}=\\&=
P(x,\partial_x(\Lambda (x))^t\partial_y)_{|x=\Lambda^{-1}(y)}
\end{aligned}
\end{equation}
and its symbol $\widetilde{P}_m(y,\eta)$ is given (up to the multiplicative constant $i^m$) by
\begin{equation}\label{4-11C}
\widetilde{P}_m(y,\eta)=\sum_{|\alpha|=
m}a_{\alpha}(\Lambda^{-1}(y))\left(\left(\partial_x\Lambda(x)\right)^t\eta\right)_{|x=\Lambda^{-1}(y)}^{\alpha}.
\end{equation}

From what we have so far established, we easily obtain

\begin{theo}[\textbf{invariant property of the characteristic surfaces}]\label{inv-caratt}
	\index{Theorem:@{Theorem:}!- invariant property of characteristic surfaces@{- invariant property of the characteristic surfaces}}
Let $\Omega$ be an open set of $\mathbb{R}^n$, and let $x_0\in \Omega$,
$\phi\in C^m(\Omega)$,  $$\Gamma=\left\{x\in \Omega:
\phi(x)=\phi(x_0) \right\}.$$ Let us suppose that
$$\nabla \phi(x)\neq 0, \quad \forall x\in \Gamma.$$

Moreover, let  $P(x,\partial)$ be a linear differential operator of order $m$ and $\Lambda\in C^m(\Omega,\mathbb{R}^n)$ be a
diffeomorphism of class $C^m$.

Then $\Gamma$ is a noncharacteristic surface  for
$P(x,\partial)$ if and only if $\Lambda(\Gamma)$ is a noncharacteristic surface for the operator $\widetilde{P}(y,\partial_y)$.
\end{theo}
\textbf{Proof.} It is not restrictive to assume  $x_0=0$ and
$\Lambda(x_0)=0$. Now, since  $\Gamma=\left\{x\in \Omega:
\phi(x)=0 \right\}$, we have
$$\Lambda(\Gamma)=\left\{y\in \Lambda(\Omega): \widetilde{\phi}(y)=0 \right\},$$
where $\widetilde{\phi}=\phi\circ\Lambda$. On the other hand, $\Gamma$
is a noncharacteristic surface for $P(x,\partial)$ if and only if
$$P_m(x,\nabla\phi(x))\neq 0,\quad\forall x\in \Gamma,$$
but we have
\begin{equation*}
\nabla_y\widetilde{\phi}(y)=\left(\partial_x\Lambda(x)\right)^t_{|x=\Lambda^{-1}(y)}(\nabla_x\phi)\left(\Lambda^{-1}(y)\right).
\end{equation*}
Therefore by \eqref{4-11C} we have
\begin{equation*}
\widetilde{P}_m(y,\nabla_y\widetilde{\phi}(y))=i^m\sum_{|\alpha|=
m}a_{\alpha}(x)\left(\nabla_x\phi\right)^{\alpha}=P_m(x,\nabla_x\phi(x))
\end{equation*}
from which the thesis follows. $\blacksquare$

\bigskip

The following Proposition holds true.
\begin{prop}\label{prop-6C}
Let $\Omega$ be an open set of $\mathbb{R}^n$ and let  $x_0\in \Omega$ and
$\phi\in C^{\infty}(\Omega)$ satisfy
$$\nabla \phi(x_0)\neq 0.$$
Let $f, a_{\alpha}\in C^{\infty}(\Omega)$ for $|\alpha|\leq m$.
Let us assume that $$\Gamma=\left\{x\in \Omega: \phi(x)=\phi(x_0)
\right\}$$ is a \textbf{noncharacteristic surface in} $x_0$
for the operator $P(x,\partial)$. Let  \\ $g_j\in
C^{\infty}(\Omega) $, for $j=0,1,\cdots,m-1$.

If $u$  is a  $C^{\infty}$ solution in a neighborhood $\mathcal{U}$
of $x_0$ of the Cauchy problem
\begin{equation}\label{0-6C}
\begin{cases}
P(x,\partial)u=f(x), \quad \forall x\in \mathcal{U}, \\
\\
\frac{\partial^ju(x)}{\partial \nu^j}=g_j(x), \quad
j=0,1,\cdots,m-1, \forall x\in \Gamma\cap\mathcal{U},
\end{cases}.
\end{equation}
then the derivatives $\partial^{\alpha}u(x_0)$  are uniquely determined for every $\alpha \in\mathbb{N}_0^{n}$.
\end{prop}
\textbf{Proof.} Although the proof follows the line
indicated in the particular case $\phi(x)=-x_n$ we want to
dwell on some details that may be useful later.

We split the proof into two steps. In \textbf{Step
	I} we determine the derivatives of $u$ of order less than $m$ on
$\Gamma$ in, \textbf{Step II} we determine the higher order 
derivatives of $u$.

\textbf{Step I.} We begin by proving that if $\phi\in C^k(\Omega)$,
$k\geq 1$, $\nabla\phi(x_0)\neq 0$ and $v$ is any
$C^k$ function in a neighborhood of $x_0$ then there exists a neighborhood
$\mathcal{V}$ of $x_0$ such that the derivatives
$\frac{\partial^jv(x)}{\partial \nu^j}$ on $\Gamma\cap \mathcal{V}$,
for $j=0,1,\cdots,k$ determine the derivatives $\partial^{\alpha}v$ on
$\Gamma\cap \mathcal{V}$ for each $|\alpha|\leq k$.

Since $\nabla\phi(x_0)\neq 0$, we may limit ourselves, up to
isometries, to consider the case where, we have, for a suitable $\delta>0$, 

\begin{equation}\label{0-7C}
\Gamma=\left\{(x',\varphi(x')): x'\in B'_{\delta}\right\},
\end{equation}
where $\varphi\in C^k(\overline{B'_{\delta}})$, $$\varphi(0)=|\nabla
\varphi(0)|=0.$$

We have
\begin{equation}\label{1-7C}
\nu((x',\varphi(x')))=\frac{(\nabla_{x'}\varphi(x'),-1)}{\sqrt{1+|\nabla_{x'}\varphi|^2}}.
\end{equation}
Also set 
\begin{equation}\label{2-7C}
\mu(x'):=\nu((x',\varphi(x')))\sqrt{1+|\nabla_{x'}\varphi|^2}=(\nabla_{x'}\varphi(x'),-1).
\end{equation}

We proceed by induction on the order $s$ of
derivatives. If $s=0$, then  $v$ is known on $\Gamma$, but to better understand the procedure, we also consider the case $k=1$. In
such a case we have, for $x'\in B'_{\delta}$,
\begin{equation}\label{2-9C}
v(x',\varphi(x'))=g_0(x'),
\end{equation}
\begin{equation}\label{3-9C}
	\begin{aligned}
&\sum_{j=1}^{n-1}(\partial_jv)(x',\varphi(x'))\partial_j\varphi(x')-(\partial_nv)(x',\varphi(x'))=\\&=g_1(x')\sqrt{1+|\nabla_{x'}\varphi(x')|^2}.
\end{aligned}
\end{equation}

\medskip

Making the derivatives w.r.t. $x_i$ of both the sides of \eqref{2-9C}, for
$i=1,\cdots,n-1$, we get, (for the sake of brevity, omit the variables)
$$\partial_iv+\partial_nv\partial_i\varphi=\partial_ig_0.$$
Hence
$$\partial_iv=\partial_ig_0-\partial_nv\partial_i\varphi$$
and, inserting these derivatives in \eqref{3-9C}, we have
$$\sum_{j=1}^{n-1}\partial_jg_0\partial_j\varphi-\left(1+|\nabla_{x'}\varphi|^2\right)\partial_nv=g_1\sqrt{1+|\nabla_{x'}\varphi|^2}$$
from which we get

$$\partial_nv=\frac{1}{1+|\nabla_{x'}\varphi|^2}\left(\nabla_{x'}g_0\cdot\nabla_{x'}\varphi-g_1\sqrt{1+|\nabla_{x'}\varphi|^2}\right).$$

Now let us prove that if $\partial^{\alpha}v$ are determined by
$|\alpha|\leq s$ on $\Gamma$ (with $s\leq k-1$) then $\partial^{\alpha}v$ are 
determined  on $\Gamma$ for $|\alpha|\leq s+1$.

Let
\begin{equation}\label{3-7C}
\partial^{\alpha}v(x',\varphi(x'))=h_{\alpha}(x'),\quad \mbox{for } |\alpha|=s, \ \mbox{ on } \Gamma
\end{equation}
and set, for $j=0,1,\cdots,k$,  
$$\widetilde{g}_j(x')=\left(1+|\nabla_{x'}\varphi|^2\right)^{j/2}\frac{\partial^j v}{\partial\nu^j}=
\sum_{|\alpha|=
j}\mu^{\alpha}(x')\partial^{\alpha}v(x',\varphi(x')).$$ Then, besides \eqref{3-7C}, we know that
\begin{equation}\label{4-7C}
\sum_{|\alpha|=
s+1}\mu^{\alpha}\partial^{\alpha}v=\widetilde{g}_{s+1}\quad\mbox{ on
} \Gamma,
\end{equation}
which we write 

\begin{equation}\label{1-8C}
\sum_{j=0}^{s+1}\sum_{|\alpha'|= s+1-j}\mu'^{\alpha'}
\mu^j\partial_n^{j}\partial^{\alpha'}v=\widetilde{g}_{s+1}\quad\mbox{
on } \Gamma,
\end{equation}
where $\alpha'=(\alpha_1,\cdots,\alpha_{n-1})$ and
$\mu'=(\mu_1,\cdots,\mu_{n-1})$.

At this point we express $\partial_n^{j}\partial^{\alpha'}v$, for
$|\alpha'|= s+1-j$, through the functions $\partial_n^{s+1}v$ and
$h_{\gamma}$ for $|\gamma|\leq s$.

If $j<s+1$, there exists $i\in\{1,\cdots,n-1\}$ such that
$\alpha'-e_i\geq 0$. Set $\beta'=\alpha'-e_i$ and let us recall that by \eqref{3-7C} we have

$$\partial_n^{j}\partial^{\beta'}v=h_{(\beta',j)},$$
from which, making the derivative w.r.t.  $x_i$ we have

$$\partial_n^{j}\underset{\partial^{\alpha'}v }{\underbrace{\partial_i \partial^{\beta'}v}}+
\mu_i\partial_n^{j+1}\partial^{\beta'}v=\partial_i h_{(\beta',j)},$$
hence

$$\partial_n^{j}\partial^{\alpha'}v=\partial_i h_{(\beta',j)}-\mu_i\partial_n^{j+1}\partial^{\beta'}v.$$
Now, if $\beta'\neq 0$, we proceed in a similar manner for
$\partial_n^{j}\partial^{\alpha'}v$ and then we iterate. Denoting by
$H_{(\alpha',j)}$ the functions of the type
$$\sum_{k=0}^{j}\sum_{|\beta'|= s-j}c_{(\beta',k)}h_{(\beta',k)},$$
where $c_{(\beta',k)}$ are known functions expressable by means of 
$\mu_i$, $i=1,\cdots,n-1$, we obtain

\begin{equation}\label{2-8C}
\partial_n^{j}\partial^{\alpha'}v=H_{(\alpha',j)}+(-1)^{s+1-j}\mu'^{\alpha'}\partial_n^{s+1}v\quad\mbox{ for }
|\alpha|=s+1-j.
\end{equation}
Inserting \eqref{2-8C} in \eqref{1-8C} we get (recall $\mu_n=-1$)

\begin{equation*}
\begin{aligned}
&\widetilde{g}_{s+1}=\sum_{j=0}^{s+1}\sum_{|\alpha'|=
s+1-j}(-1)^{j}\mu'^{\alpha'}
\left(H_{(\alpha',j)}+(-1)^{s+1-j}\mu'^{\alpha'}\partial_n^{s+1}v\right)=\\&
=\sum_{j=0}^{s+1}\sum_{|\alpha'|=
s+1-j}(-1)^{j}\mu'^{\alpha'}H_{(\alpha',j)}+\left(\sum_{j=0}^{s+1}\sum_{|\alpha'|=
s+1-j}(-1)^{j}\mu'^{2\alpha'}\right)\partial_n^{s+1}v.
\end{aligned}
\end{equation*}
From which we have
\begin{equation}\label{1-9C}
	\begin{aligned}
&(-1)^{s+1}\left(1+\sum_{1\leq|\alpha'|\leq
s+1}\mu'^{2\alpha'}\right)\partial_n^{s+1}v=\\&=\widetilde{g}_{s+1}-
\sum_{j=0}^{s+1}\sum_{|\alpha'|=
s+1-j}(-1)^{j}\mu'^{\alpha'}H_{(\alpha',j)}.
\end{aligned}
\end{equation}
Since
$$\sum_{1\leq|\alpha'|\leq s+1}\mu'^{2\alpha'}\geq 0,$$
by \eqref{1-9C} we determine $\partial_n^{s+1}v$.

\bigskip

\textbf{Step II.} Let us consider the Cauchy problem

\begin{equation}\label{3-12C}
\begin{cases}
P(x,\partial)u=f(x), \quad \forall x\in \mathcal{U}, \\
\\
\frac{\partial^ju(x)}{\partial \nu^j}=g_j(x), \quad
j=0,1,\cdots,m-1, \ \ \forall x\in \Gamma\cap\mathcal{U}.
\end{cases}.
\end{equation}
Like Step I we assume that $\Gamma$ is the graph
 \eqref{0-7C}, where $\varphi\in
C^{\infty}(\overline{B'_{\delta}})$, $\varphi(0)=|\nabla_{x'}\varphi(0)|=0$.
By what was proved in Step I we can determine, from
the functions $g_0,g_1,\cdots,g_{m-1}$ only, the derivatives

\begin{equation}\label{4-12C}
\partial^{\alpha}u,\quad \mbox{ for } |\alpha|\leq m-1, \ \mbox{ on } \Gamma.
\end{equation}

Let now us show in which a way  we determine the other derivatives of $u$ on
$\Gamma$. Let
$$\Lambda: B_{\delta}\subset\mathbb{R}^n_x\rightarrow\mathbb{R}^n_y,\quad \Lambda(x)=\left(x',x_n-\varphi(x')\right),$$

\begin{figure}\label{figura-Gamma}
	\centering
	\includegraphics[trim={0 0 0 0},clip, width=10cm]{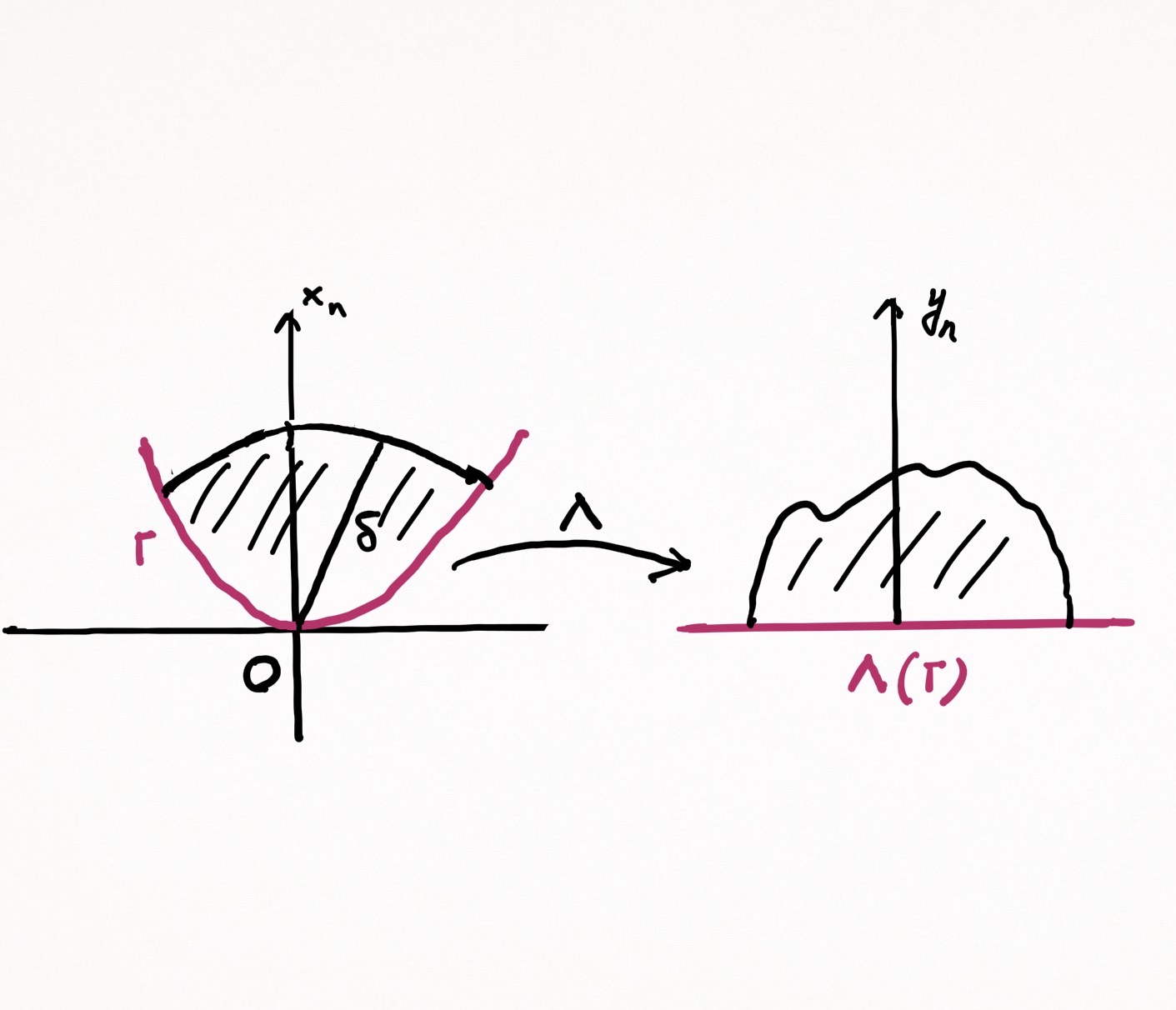}
	\caption{}
\end{figure}

\noindent we have
 
$$\Lambda(\Gamma)=\left\{(y',0): y'\in B'_{\delta}\right\}=\left\{y\in B_{\delta}:-y_n=0\right\}$$
and, set
$$v(y)=u\left(\Lambda^{-1}(y)\right).$$
Equality \eqref{1-11C} implies that the derivatives

\begin{equation*}
\partial^{\alpha}v(y',0), \  \mbox{ for } |\alpha|\leq m-1, \quad\forall  y'\in B'_{\delta},
\end{equation*}
are all uniquely determined. In particular, the following derivatives are
determined 
\begin{equation}\label{1-13C}
\partial_n^{j}v(y',0), \ \mbox{ for } j=0,1,\cdots,m-1,\quad\forall y'\in B'_{\delta},.
\end{equation}
On the other hand, by \eqref{2-11C} we have that $v$, 
solves the following equation in a neighborhood of  $0$

$$\widetilde{P}(y,\partial_y)v(y):=\sum_{|\alpha|\leq m}b_{\alpha}\partial_y^{\alpha}v(y)=\widetilde{f}(y),$$
where $\widetilde{f}(y)=f\left(\Lambda^{-1}(y)\right)$. Now, by
Proposition \ref{inv-caratt} we have that $\Lambda(\Gamma)$ is a noncharacteristic surface for $\widetilde{P}(y,\partial_y)$.
Since $\Lambda(\Gamma)=\left\{y\in
B_{\delta}:-y_n=0\right\}$ we have
$$b_{(0,m)}(y',0)=\sum_{|\alpha|= m}b_{\alpha}(y',0)e_n^{\alpha}\neq 0 \quad\forall y'\in B'_{\delta}.$$
We are therefore reduced to the same situation examined
at the beginning of this Section \ref{superf-caratteristiche} and we 
then calculate the derivatives $\partial_y^{\alpha}v(y',0)$ for $|\alpha|\geq m$ from
the derivatives of $g_j$, $j=0,1,\cdots,m-1$ and the
coefficients (and their derivatives) of
$\widetilde{P}(y,\partial_y)$. Finally, by exploiting formula
\eqref{1-11C} we obtain the derivatives $\partial_x^{\alpha}u$ on
$\Gamma$ for $|\alpha|\geq m$. $\blacksquare$

\bigskip

We conclude this Section with some examples and remarks.

\medskip

\textbf{Remarks.} Let $P(x,\partial)$, $x\in\Omega$, be a linear differential operator
whose principal part is $P_m(x,\partial)$.

\medskip

\noindent\textbf{1.} We call  \textbf{characteristic equation} \index{equation:@{equation:}!- characteristic@{- characteristic}}the non linear first order equation 
in the unknown $\phi$
\begin{equation}\label{0-14C}
P_m(x,\nabla\phi(x))=0.
\end{equation}

\medskip

\noindent\textbf{2.} We say that $P(x,\partial)$  is
\textbf{elliptic in the point} $x_0\in\Omega$ \index{elliptic operator}if
$$P_m(x_0,\xi)\neq 0, \quad\forall\xi\in\mathbb{R}^n\setminus \{0\}.$$
It is evident that the elliptic operators have not characteric surface. If $P(x,\partial)$  is elliptic in each point of $\Omega$ we say that $P(x,\partial)$ is 
\textbf{elliptic in} $\Omega$. 
Each linear differential operators of one variable are elliptic
$$P(t,\frac{d}{dt})=a_m(t)\frac{d^m}{dt^m}+\cdots+a_0(t),\quad t\in I,$$
where $a_m(t)\neq 0$ for $t\in I$, where $I$ is an interval of
$\mathbb{R}$. As a matter of fact
$$P(t,\xi)=a_m(t)\xi^m\neq 0,\quad t\in I,\quad \forall \xi\in\mathbb{R}\setminus \{0\}.$$
A remarkable example of an elliptic operator is the operator of
Cauchy-Riemann
$$P((x,y),\partial_x,\partial_y)=\partial_x+i\partial_y.$$
If $m=2$ and $a_{jk}$, where $a_{jk}(x)=a_{kj}(x)$ for $x\in\Omega$,
$j,k=1,\cdots,n$ are real--valued functions in $\mathbb{R}$, we define
a \textbf{uniformly elliptic} operator with bounded coefficients an operator of the type
$$P(x,\partial)=\sum_{j,k=1}^na_{jk}(x)\partial^2_{jk}+\sum_{j=1}^nb_j(x)\partial_{j}+c(x)$$ such that there exists a constant \index{uniformly elliptic}
$\lambda\geq 1$ satisfying
$$\lambda^{-1}|\xi|^2\leq \sum_{j,k=1}^na_{jk}(x)\xi_j\xi_k\leq \lambda|\xi|^2, \quad \forall x\in \Omega, \quad \forall\xi\in \mathbb{R}^n.$$



$\blacklozenge$

\bigskip

\underline{\textbf{Exercise.}} Let us consider the following operator with constant coefficients

\begin{equation}\label{4-15C}
P(\partial)=\sum_{j,k=1}^na_{jk}\partial^2_{jk},
\end{equation}
where $\{a_{jk}\}$ is a real symmetric matrix. If $\det\{a_{jk}\}\neq 0$
then there exists a nonsingular matrix $C$  such that, setting $y=Cx$, operator \eqref{4-15C}, is transformed in 

\begin{equation*}
P(\partial_y)=\sum_{j=1}^n\varkappa_j \partial^2_{y_j},
\end{equation*}
where $\varkappa_j$, $j=1,\cdots,n$, is equal either to $1$ or to
$-1$. $\clubsuit$

\section{The Cauchy-Kovalevskaya Theorem} \label{Cauchy-Kov}
In what follows we denote by $\Omega$ an open set of $\mathbb{R}^n$, $x_0\in \Omega$ and by \\ $\phi\in C^{\omega}(\Omega,\mathbb{R})$ a function such that
\begin{equation}\label{1-34C}
	\nabla \phi(x)\neq 0, \quad \forall x\in \Gamma :=\left\{x\in \Omega:
	\phi(x)=\phi(x_0) \right\}.
\end{equation}
By \eqref{1-34C} we have that for every $\widehat{x}\in\Gamma$ there exist $r, M>0$ and an isometry $\Psi$ under which we have $\Psi(0)=\widehat{x}$, and 

\begin{equation*}
	\Psi^{-1}(\Gamma)\cap Q_{r, 2M}=\left\{(x',\varphi(x')):x'\in B'_r\right\} 
\end{equation*}
where $\varphi\in C^{\omega}\left(B'_r; \mathbb{R}\right)$ and it satisfies

$$\varphi(0)=0,\quad |\nabla \varphi(0)|=0, $$

$$\left\Vert \varphi\right\Vert_{C^1\left(\overline{B'_{r}}\right)}\leq
Mr.$$ 

We say that a function 
$$h:\Gamma \rightarrow \mathbb{C}$$ is analitic on $\Gamma$, provided that

$$\left(h\circ \Psi\right)\left(\cdot,\varphi(\cdot)\right)\in C^{\omega}\left(B'_r\right).$$ 
 
\bigskip  

In this Section we will prove
\begin{theo}[\textbf{Cauchy--Kovalevskaya}]\label{teor-CK}
	\index{Theorem:@{Theorem:}!- Cauchy--Kovalevskaya@{- Cauchy--Kovalevskaya}}
Let $m\in\mathbb{N}$ and let  $\Omega$ be an open set of $\mathbb{R}^n$, $x_0\in \Omega$,
$\phi\in C^{\omega}(\Omega,\mathbb{R})$ which satisfies \eqref{1-34C}.
Moreover, let  $P(x,\partial)$ be the linear differential operator
\begin{equation}\label{2-34C}
P(x,\partial)=\sum_{|\alpha|\leq m}a_{\alpha}(x)\partial^{\alpha},
\end{equation}
where $a_{\alpha}\in C^{\omega}(\Omega)$, for $|\alpha|\leq m$. Let
$g_0,g_1,\cdots, g_{m-1}$ be \textbf{analytic functions on} $\Gamma$. Let us assume that $\Gamma$ is a noncharacteristic surface
for the operator $P(x,\partial)$. Let
$f \in C^{\omega}(\Omega)$.

Then for every $\widetilde{x}\in \Gamma$ there exists a neighborhood 
$\mathcal{U}_{\widetilde{x}}$ such that the Cauchy problem
\begin{equation}\label{3-34C}
\begin{cases}
P(x,\partial)u(x)=f(x), \quad \forall  x\in \mathcal{U}_{\widetilde{x}}, \\
\\
\frac{\partial^ju(x)}{\partial \nu^j}=g_j(x), \quad
j=0,1,\cdots,m-1,\quad \forall x\in
\Gamma\cap\mathcal{U}_{\widetilde{x}}
\end{cases}
\end{equation}
has a unique analytic solution in $\mathcal{U}_{\widetilde{x}}$.
\end{theo}

\bigskip

In order to prove Cauchy-Kovalevskaya Theorem we need two preliminary steps

\smallskip

(i) local flatness of  initial surface;

(ii) transformation of problem \eqref{3-34C} to a Cauchy problem  for a first order system.

\bigskip

We have already considered \textbf{point (i)} in the context of
of the proof of Step II of Proposition \ref{prop-6C}.
Here it suffices to add that, referring to the notations used
in the above proof, the function $\varphi$ introduced there
is, not only $C^{\infty}$, but also analytic in $B'_{\delta}$ and, consequently, the map
\begin{equation}\label{2-24-5}
\Lambda:
B_{\delta}\subset\mathbb{R}^n_x\rightarrow\mathbb{R}^n_y,\quad
\Lambda(x)=\left(x',x_n-\varphi(x')\right),
\end{equation} which, we recall, flatten $\Gamma$ in the sense that

$$\Lambda(\Gamma)=\left\{(y',0): y'\in B'_{\delta}\right\}=\left\{y\in B_{\delta}:-y_n=0\right\}.$$
Moreover, setting
\begin{equation}\label{3-36C}
v(y)=u\left(\Lambda^{-1}(y)\right),
\end{equation}
the operator $\widetilde{P}(y,\partial_y)$, defined by
\begin{equation}\label{4-36C}
\left(\widetilde{P}(y,\partial_y)v(y)\right)_{|y=\Lambda(x)}=P(x,\partial_x)u(x),
\end{equation}
has its principal part  $\widetilde{P}_m(y,\partial_y)$, given by

\begin{equation}\label{5-36C}
\widetilde{P}_m(y,\partial_y)=\sum_{|\alpha|=
m}b_{\alpha}(y)\partial^{\alpha},
\end{equation}
where
\begin{equation}\label{5n-36C}
b_{\alpha}(y)=a_{\alpha}(\Lambda^{-1}(y)), \quad \mbox{ for }
|\alpha|\leq m,
\end{equation}
are analytic functions on $B_{\delta}$.

Again by Proposition \ref{prop-6C} we know that the new initial data \\
$\widetilde{g}_j(y')=\partial^j_n v(y',0)$,
$j=0,1,\cdots,m-1$, are determined by the initial data $g_j$ and that 
$\widetilde{g}_j$ are analytic in a neighborhood of $0\in
\mathbb{R}^{n-1}$. Moreover, we know that the surface $$\left\{y\in
B_{\delta}:-y_n=0\right\}$$ is noncharacteristic. Therefore

$$b_{(0,m)}(y',0)\neq 0, \quad\forall y'\in B'_{\delta}(0).$$
Hence, setting

$$\widetilde{b}_{\alpha}(y)=-\frac{b_{\alpha}(y)}{b_{(0,m)}(y)},\quad \widetilde{f}(y)=\frac{f\left(\Lambda^{-1}(y)\right)}{b_{(0,m)}(y)},$$
we write problem \eqref{3-34C} as

\begin{equation}\label{1-38C}
\begin{cases}
\partial_n^mv=\sum_{|\alpha|\leq m, \alpha_n\leq m-1}\widetilde{b}_{\alpha}(y)\partial^{\alpha}v+\widetilde{f}(y),  \\
\\
\partial^j_n v(y',0)=\widetilde{g}_j(y'), \quad j=0,1,\cdots,m-1.
\end{cases}
\end{equation}

We may easily transform problem \eqref{1-38C} into a Cauchy problem 
with homogeneous initial conditions. To this purpose it suffices to define the
function

$$H(y)=\sum_{j=0}^{m-1}\frac{y_n^j}{j!}\widetilde{g}_j(y')$$
and set
$$w=v-H,$$
obtaining

\begin{equation}\label{2-38C}
\begin{cases}
\partial_n^mw=\sum_{|\alpha|\leq m, \alpha_n\leq m-1}\widetilde{b}_{\alpha}(y)\partial^{\alpha}w+F(y),  \\
\\
\partial^j_n w(y',0)=0, \quad j=0,1,\cdots,m-1.
\end{cases}
\end{equation}
where

\begin{equation}\label{2n-38C}
F(y)=\widetilde{f}(y)-\partial_n^mH(y)+\sum_{|\alpha|\leq m,
\alpha_n\leq m-1}\widetilde{b}_{\alpha}(y)\partial^{\alpha}H(y).
\end{equation}

\bigskip

\noindent \textbf{(ii)}. 
The idea of the tranformation is simple and it partly replicates the one
usually followed to reduce a Cauchy problem for ordinary differential equations
of order $m$ to a Cauchy problem for a first-order system. However, in our case the unknown
depends on $n>1$ variables, and this requires further arrangements.

In order to highlight the main steps we illustrate 
the procedure in the case where $n=2$ and
the operator is equal to its principal part only. Next
we will outline how to proceed in the general case.

\medskip

Let us consider the Cauchy problem
\begin{equation}\label{1-39C}
\begin{cases}
\partial_t^mu=\sum_{j\leq m-1}\sum_{i+j=m}a_{i,j}(x,t)\partial^i_x\partial^j_t u+ f(x,t),  \\
\\
\partial_t^j u(x,0)=0, \quad j=0,1,\cdots,m-1, \ \ \forall x\in \mathbb{R}.
\end{cases}
\end{equation}

Let us assume that $m\geq2$ and let $u(x,t)$ be a $C^{\infty}$ solution
of \eqref{1-39C}. We set

$$V_{i,j}=\partial^i_x\partial^j_t u,\quad\mbox{ for } i+j\leq m-1.$$
It is simple to check what follows

\begin{subequations}
\label{2-39C}
\begin{equation}
\label{2a-39C}
\partial_tV_{i,j}=V_{i,j+1},\quad\mbox{ for } i+j< m-1,
\end{equation}
\begin{equation}
\label{2b-39C}
\partial_tV_{i,j}=\partial_xV_{i-1,j+1},\quad\mbox{ for } i+j=m-1,\quad
i>0,
\end{equation}
\begin{equation}
\label{2c-39C}
\partial_tV_{0,m-1}=\sum_{i+j=m-1}a_{i+1,j}(x,t)\partial_x V_{i,j} +
f(x,t),
\end{equation}
\begin{equation}
\label{3-39C} V_{i,j}(x,0)=0 \quad\mbox{ for } i+j\leq m-1.
\end{equation}
\end{subequations}

Hence, if $V_{0,0}$ is a solution to Cauchy problem
\eqref{1-39C}, then  it is a solution to Cauchy problem
\eqref{2a-39C}--\eqref{3-39C}. Thus, in order to prove 
the equivalence of problem \eqref{1-39C} and problem
\eqref{2a-39C}--\eqref{3-39C} it suffices to prove the converse. Let us suppose, therefore, that $V=\left(V_{i,j}\right)_{i+j\leq m-1}$ is a
$C^{\infty}$ a solution to problem
\eqref{2a-39C}--\eqref{3-39C} and let us prove that $V_{0,0}$ is a
solution to problem \eqref{1-39C}.

\bigskip

\noindent\textbf{1. We prove that}

\begin{equation}\label{1-40C}
\mbox{If }\quad i+j=m-1, i>0\quad\mbox{ then }\quad
V_{i,j}=\partial_xV_{i-1,j}.
\end{equation}

\medskip

\noindent\textbf{Proof.} Let $i+j\leq m-1$ and $i>0$. By \eqref{2b-39C}
we have

\begin{equation}
\label{2-40C}
\partial_tV_{i,j}=\partial_xV_{i-1,j+1},\quad\mbox{ for } i+j=m-1,\quad i>0.
\end{equation}
Now, since $(i-1)+j=m-2$, by \eqref{2a-39C} we have

\begin{equation}
\label{3-40C} V_{i-1,j+1}=\partial_tV_{i-1,j}.
\end{equation}
Hence, by \eqref{2-40C} e \eqref{3-40C} we have

\begin{equation}
\label{4-40C}
\partial_t\left(V_{i,j}-\partial_xV_{i-1,j}\right)=0.
\end{equation}
On the other hand, by \eqref{3-39C} we have

$$\left(V_{i,j}-\partial_xV_{i-1,j}\right)(x,0)=0$$ this equality and \eqref{4-40C} implies

$$V_{i,j}=\partial_xV_{i-1,j}.$$

\bigskip

\noindent\textbf{2. We prove that}

\begin{equation}\label{1-41C}
\mbox{If }\quad i+j\leq m-1, i>0\quad\mbox{ then }\quad
V_{i,j}=\partial_xV_{i-1,j}.
\end{equation}

\medskip

\noindent\textbf{Proof.} Set $l=(m-1)-(i+j)$ and let us proceed by
induction on $l$. If $l=0$, then \eqref{1-41C} holds true, because it is nothing but  \eqref{1-40C}. Let now let us suppose that \eqref{1-41C} holds true
for $l$ and prove it for $l+1$. Hence,  let us suppose that

\begin{equation}\label{2-41C}
\mbox{if }\quad i+j=(m-1)-l, i>0\quad\mbox{ then }\quad
V_{i,j}=\partial_xV_{i-1,j}.
\end{equation}
Let $i,j$ satisfy $i+j=(m-1)-(l+1)$ and $i>0$. Since $i+j<m-1$
by \eqref{2a-39C} we have

\begin{equation}
\label{3-41C}
\partial_tV_{i,j}=V_{i,j+1}.
\end{equation}
Since we have $i+(j+1)=(m-1)-l$,  \eqref{2-41C} gives

\begin{equation}\label{4-41C}
V_{i,j+1}=\partial_xV_{i-1,j+1}.
\end{equation}
Since $(i-1)+j=(m-1)-l-2<m-1$,  \eqref{2a-39C} gives

\begin{equation}
\label{5-41C} V_{i-1,j+1}=\partial_tV_{i-1,j},\quad\mbox{ for } i+j<
m-1.
\end{equation}
Hence, by \eqref{4-41C} e \eqref{5-41C} we have
$$V_{i,j+1}=\partial_xV_{i-1,j+1}=\partial_t\partial_xV_{i-1,j},$$
by the latter and by \eqref{3-41C} we have

\begin{equation}
\label{6-41C}
\partial_t\left(V_{i,j}-\partial_xV_{i-1,j}\right)=0
\end{equation}
so, by \eqref{3-39C} we have that, if $i+j=(m-1)-(l+1)$, $i>0$, then

 $$V_{i,j}=\partial_xV_{i-1,j}.$$
 \eqref{1-41C} is proved.

\medskip

\noindent\textbf{Conclusions.} Iteration of  \eqref{1-41C} implies what follows:
\begin{equation}\label{1-42C}
\mbox{if }\quad i+j\leq m-1, i>0\quad\mbox{ then }\quad
V_{i,j}=\partial^i_xV_{0,j}.
\end{equation}
On the other hand, \eqref{2a-39C} gives

\begin{equation}\label{2-42C}
 V_{0,j}=\partial_tV_{0,j-1}=\cdots=\partial^j_tV_{0,0}.
\end{equation}
All in all, by \eqref{1-42C} and \eqref{2-42C} we get
$$V_{i,j}=\partial^i_x\partial^j_tV_{0,0}$$
and using this equality into \eqref{2c-39C}--\eqref{3-39C} we have that
$V_{0,0}$ solves \eqref{1-39C}.

We may rewrite problem \eqref{2a-39C}--\eqref{3-39C} in a
more concentrated form as follows

\begin{equation}\label{0-42C}
\begin{cases}
\partial_tV(x,t)=B(x,t)\partial_x V+ F(x,t),  \\
\\
V(x,0)=0.
\end{cases}
\end{equation}
where, for an appropriate $N\in\mathbb{N}$, $V$ is a function with
values in $\mathbb{R}^N$, $B$ is an $N\times N$ matrix whose
entries are analytic and $F=fe_N$.

In the case $n>2$ one may similarly reduce Cauchy problem \eqref{3-34C} to a Cauchy problem for a first order system. We outiline the procedure (the details of which we leave to the reader). First, it is convenient to introduce the following notation. If $\alpha\in \mathbb{N}_0^{n-1}\setminus\{0\}$ is a  multi-index, we set
$$i(\alpha)=\min\{i:\alpha_i>0\}.$$ In addition, we set  $t=x_n$
and
$$V_{\alpha,j}=\partial_t^j\partial_{x'}^{\alpha}w,\quad\mbox{ for } |\alpha|+j\leq m-1.$$
By \eqref{2-38C} we have

\begin{subequations}
\label{a-42C}
\begin{equation}
\label{aa-42C}
\partial_tV_{\alpha,j}=V_{\alpha,j+1}, \quad\mbox{ for } |\alpha|+j<
m-1,
\end{equation}
\begin{equation}
\label{ab-42C}
\partial_tV_{\alpha,j}=\partial_{x_{i(\alpha)}}V_{\alpha-e_{i(\alpha)},j+1}, \quad\mbox{ for } |\alpha|+j= m-1,
|\alpha|>0,
\end{equation}
\begin{equation}
\label{ac-43C}
\partial_tV_{0,m-1}=\sum_{|\alpha|+j=m, j<m}c_{\alpha,j}\partial_{x_{i(\alpha)}} V_{\alpha-e_{i(\alpha)},j} +
\sum_{|\alpha|+j\leq m}d_{\alpha,j}V_{\beta,j}+ F,
\end{equation}
\begin{equation}
\label{ad-43C} V_{\alpha,j}(x',0)=0 \quad\mbox{ for } |\alpha|+j\leq
m-1,
\end{equation}
\end{subequations}
where $c_{\alpha,j}$, $d_{\alpha,j}$, $F$ are analytic functions in the variables
 $x'$ and $t$.

\bigskip

\noindent\textbf{Proof of the Cauchy--Kovalevskaya Theorem.}

Taking into account what has been done in (i) and (ii) and changing the
notations a little, we may reformulate Cauchy problem
\eqref{2-38C} as follows.

\begin{equation}\label{1-43C}
\begin{cases}
\partial_tU_k(x,t)=\sum_{j=1}^{n-1} \sum_{l=1}^{N}B_j^{lk}\partial_{x_j} U_l+ \sum_{l=1}^{N}C^{lk} U_l+F_k, \quad k=1,\cdots, N, \\
\\
U_k(x,0)=0, \quad k=1,\cdots, N,
\end{cases}
\end{equation}

\medskip

\noindent where $B_j^{lk}$, $C^{lk}$, $F_k$, $j=1,\cdots, n-1$, $l,k=1,\cdots,N$ 
are analytic functions in a neighborhood of $0$.

\medskip

In order to prove the existence and the uniqueness for Cauchy problem
\eqref{1-43C} we proceed as follows:

\medskip

\noindent\textbf {Step I.} For any function $U\in C^{\infty}$ that 
\textit{we suppose} to satisfy \eqref{1-43C}, we will calculate the
derivatives
$$\partial^{\alpha}U_k(0,0):=U_{k,\alpha},\quad \forall \alpha\in \mathbb{N}_0^n,\quad k=1,\cdots,N.$$
Setting, $\partial^{\alpha}=\partial_x^{\alpha'}\partial_t^{\alpha_n}$, for $\alpha=(\alpha',\alpha_n)$,
we will have

\begin{equation}\label{0-44C}
U_{k,(\alpha',0)}=0,\quad \forall \alpha'\in
\mathbb{N}_0^{n-1},\quad k=1,\cdots,N.
\end{equation}

\medskip

\noindent\textbf{Step II.} We will show what follows. Let us assume that the functions
$\widetilde{B}_{j}^{lk}$, $\widetilde{C}^{lk}$, $\widetilde{F}$ and
$\widetilde{\varphi}$ (the latter is independent of $t$) satisfy
the following conditions

\smallskip

\textbf{(a)} $B_j^{lk}\preccurlyeq \widetilde{B}_{j}^{lk}$,
$C^{lk}\preccurlyeq \widetilde{C}^{lk}$, for $j=1,\cdots, n-1$,
$l,k=1,\cdots,N$, $F\preccurlyeq \widetilde{F}$, \indent $0\preccurlyeq
\widetilde{\varphi}$

\noindent and let us assume that

\textbf{(b)} it occurs that for any $ C^{\infty}$ solution $\widetilde{U}$ to the Cauchy problem

\begin{equation}\label{1-44C}
\begin{cases}
\partial_t\widetilde{U}_k=\sum_{j=1}^{n-1} \sum_{l=1}^{N}\widetilde{B}_j^{lk}\partial_{x_j} \widetilde{U}_l+
\sum_{l=1}^{N}\widetilde{C}^{lk} \widetilde{U}_l+\widetilde{F}_k,  \\
\\
\widetilde{U}_k(x,0)=\widetilde{\varphi}_k(x),
\end{cases}
\end{equation}
 we will have

\begin{equation}\label{2-44C}
\left|U_{k,\alpha}\right|\leq \partial^{\alpha}
\widetilde{U}_k(0,0),\quad \forall \alpha\in \mathbb{N}_0^n,\quad
k=1,\cdots,N.
\end{equation}

\medskip

\noindent\textbf{Step III.} We will construct some majorants
$\widetilde{B}_{j}^{lk}$, $\widetilde{C}^{lk}$, $\widetilde{F}$ and
$\widetilde{\varphi}$ for which Cauchy problem
\eqref{1-44C} \textbf{does indeed have an analitic solution}. 
Let us denote again by $\widetilde{U}$ such a solution. By Step II and, in particular, by \eqref{2-44C} it will follow that the power series
$$\sum_{\alpha\in \mathbb{N}_0^n}\frac{1}{\alpha !}U_{\alpha}x^{\alpha'}t^{\alpha_n},$$
will converge in a neighborhood of $0$ and its sum, which we denote
by $V$, will satisfy, the system

\begin{equation}\label{3-44C}
\partial_tV_k(x,t)=\sum_{j=1}^{n-1} \sum_{l=1}^{N}B_j^{lk}(x,t)\partial_{x_j} V_l+ \sum_{l=1}^{N}C^{lk}(x,t) V_l+F_k(x,t).  \\
\end{equation}
Indeed, for $k=1,\cdots, N$, the  analytic functions

$$\partial_t V_k$$ and
$$\sum_{l=1}^{N}B_j^{lk}(x,t)\partial_{x_j} V_l+ \sum_{l=1}^{N}C^{lk}(x,t) V_l+F_k(x,t),$$ will have (by construction) all the derivatives equal at $0$  .
Furthermore, by \eqref{0-44C}, we will have

\begin{equation}\label{4-44C}
V_k(x,0)=0, \quad k=1,\cdots, N
\end{equation}
and we will then have proved the existence of a solution to
Cauchy problem \eqref{1-43C}.

\medskip

\noindent\textbf{Step IV.} The uniqueness in the class of analytic functions in a connected neighborhood of $0$ for Cauchy problem
\eqref{1-43C} will be a consequence of \textbf{Step I} and of the unique continuation property for the analytic functions (Theorem \ref{teo1-23C}).

\bigskip

\bigskip

\noindent\textbf{Step I.} By the initial conditions $U(x,0)=0$ we have

\begin{equation}\label{1-45C}
\partial_x^{\alpha'}U_k(x,0)=0, \quad k=1,\cdots, N,
\end{equation}
which implies \eqref{0-44C}. Now, for every $\alpha \in
\mathbb{N}_0^n$, where $\alpha_n>0$, we have, for $k=1,\cdots, N$,

\begin{equation}\label{2-45C}
	\begin{aligned}
	\partial^{\alpha}U_k(0,0)&=P_{k,\alpha}\left(\partial^{\gamma}B_{j}^{kl},\cdots,\partial^{\delta}C^{kl},\cdots
\partial^{\beta}U_h\right)_{x=0,t=0}+\\&+\partial_x^{\alpha'}\partial_t^{\alpha_n-1}F_k(0,0),
	\end{aligned}
\end{equation}

\smallskip

\noindent where $P_{k,\alpha}$ \textbf{is a polynomial with positive coefficients} and the multi-indices $\beta$ in the derivatives
$\partial^{\beta}U_h$ satisfy
$$|\beta|\leq |\alpha|,\quad \mbox{ e }\quad \beta_n\leq \alpha_n-1.$$
In particular, \eqref{2-45C} is a recursive relation on the
derivatives of $U$. To prove \eqref{2-45C}
it suffices to make the derivatives of both the sides of the equations in \eqref{1-43C}. Concerning the positivity of the coefficients of $P_{k,\alpha}$ it suffices to keep in mind that we only use the rules of derivation of a product and a sum of functions.

For instance, if $i=1,\cdots,n-1$, we have

\begin{equation}\label{0-46C}
\begin{aligned}
\partial_t\partial_{x_i}U_k&=\sum_{j=1}^{n-1} \sum_{l=1}^{N}\left(B_j^{lk}\partial^2_{x_ix_j} U_l
+\partial_{x_i}B_j^{lk}\partial_{x_j} U_l\right)+\\& +\sum_{l=1}^{N}
\left(C^{lk} \partial_{x_i}U_l+\partial_{x_i}C^{lk}
U_l\right)+\partial_{x_i}F_k
\end{aligned}
\end{equation}
and, taking into account \eqref{1-45C}, we have

\begin{equation*}
\partial_t\partial_{x_i}U_k(0,0)=\partial_{x_i}F_k(0,0).
\end{equation*}
Analogously,
$$\partial_t\partial_{x}^{\alpha'}U_k(0,0)=\partial_{x}^{\alpha'}F_k(0,0)$$
and

\begin{equation}\label{00-46C}
\begin{aligned}
\partial^2_tU_k&=\sum_{j=1}^{n-1} \sum_{l=1}^{N}\left(B_j^{lk}\partial_t\partial_{x_j} U_l+
\partial_{t}B_j^{lk}\partial_{x_j} U_l\right)+\\&
+\sum_{l=1}^{N} \left(C^{lk} \partial_{t}U_l+\partial_{t}C^{lk}
U_l\right)+\partial_{t}F_k.
\end{aligned}
\end{equation}
Let us observe that all the derivatives of $U$ in $(0,0)$ that occur in
\eqref{00-46C} can be obtained by \eqref{0-46C} and by \eqref{1-45C}. A similar argument applies to
$\partial_t^2\partial^{\alpha'}U(0,0)$, $\cdots$
$\partial_t^j\partial^{\alpha'}U(0,0)$, $j=1,\cdots$, $\alpha'\in
\mathbb{N}_0^{n-1}$.

\bigskip

\noindent\textbf{Step II.} Let
$\widetilde{U}$ be a solution to problem \eqref{1-44C}. In a similar way
to what we have done in \textbf{ Step I} we obtain, for $\alpha \in
\mathbb{N}_0^n$ with $\alpha_n>0$, that for each $k=1,\cdots, N$,

\begin{equation}\label{3-45C}
	\begin{aligned}
		\partial^{\alpha}\widetilde{U}_k(0,0)&=P_{k,\alpha}
\left(\partial^{\gamma}\widetilde{B}_{j}^{kl},\cdots,\partial^{\delta}\widetilde{C}^{kl},\cdots
\partial^{\beta}\widetilde{U}_h\right)_{x=0,t=0}+\\&+\partial_x^{\alpha'}\partial_t^{\alpha_n-1}\widetilde{F}_k(0,0)+,
\end{aligned}
\end{equation}

\medskip

\noindent where $P_{k,\alpha}$ is the same polynomial with
positive coefficients that occurs in \eqref{2-45C} and (as in
\eqref{2-45C}) the multi-indices $\beta$ in the derivatives
$\partial^{\beta}U_h$ satisfy $|\beta|\leq |\alpha|$ and
$\beta_n\leq \alpha_n-1.$ Furthermore, we have that

\begin{equation}\label{4-47C}
\partial_x^{\alpha'}\widetilde{U}_k(0,0)=\partial_x^{\alpha'}\widetilde{\varphi}_k(0)
\end{equation}
and by  $0\preccurlyeq \widetilde{\varphi}$ we have

\begin{equation}\label{4n-47C-24}
0\leq \partial_x^{\alpha'}\widetilde{\varphi}_k(0).
\end{equation}
In order to prove \eqref{2-44C} one can proceed by
induction on the order $\alpha_n$ of the derivative with respect to $t$. If
$\alpha_n=0$ then we have

\begin{equation}\label{5-47C-24-1}
\left|U_{k,(\alpha',0)}\right|\leq \partial_x^{\alpha'}
\widetilde{\varphi}_k(0)=
\partial_x^{\alpha'} \widetilde{U}_k(0,0),\quad \forall \alpha\in \mathbb{N}_0^{n-1},\quad k=1,\cdots,N.
\end{equation}
Now, let us suppose that for a given $\alpha_n$ we have

\begin{equation}\label{5-47C}
\left|U_{k,(\alpha',\alpha_n)}\right|\leq
\partial_t^{\alpha_n}\partial_x^{\alpha'} \widetilde{U}_k(0,0),
\quad \forall \alpha\in \mathbb{N}_0^{n-1}\quad k=1,\cdots,N.
\end{equation}
We have, for $k=1,\cdots,N$,

\begin{equation*}
\begin{aligned}
U_{k,(\alpha',\alpha_n+1)}&=\partial_t^{\alpha_n+1}\partial_x^{\alpha'}
U_k(0,0)=\partial_x^{\alpha'}\partial_t^{\alpha_n}F_k(0,0)+\\&
+P_{k,(\alpha',\alpha_n+1)}\left(\partial^{\gamma}B_{j}^{kl},\cdots,\partial^{\delta}C^{kl},\cdots
\partial^{\beta}U_h\right)_{x=0,t=0},
\end{aligned}
\end{equation*}
where
$$|\beta|\leq |\alpha|,\quad \mbox{ and } \quad\beta_n\leq \alpha_n.$$
Then, since the coefficients of
$P_{k,(\alpha_n+1)}$ are positive, using 
\textbf{(a)} of \textbf{Step II} and the inductive hypothesis \eqref{5-47C},
we have

\begin{equation}\label{5n-47C}
\begin{aligned}
&\left\vert
U_{k,(\alpha',\alpha_n+1)}\right\vert\leq\left\vert\partial_x^{\alpha'}\partial_t^{\alpha_n}F_k(0,0)\right\vert+\\&
+P_{k,(\alpha',\alpha_n+1)}\left(\left\vert\partial^{\gamma}B_{j}^{kl}\right\vert,\cdots,
\left\vert\partial^{\delta}C^{kl}\right\vert,\cdots \left\vert
U_{h,\beta}\right\vert\right)_{x=0,t=0}\leq\\& \leq
\partial_x^{\alpha'}\partial_t^{\alpha_n}\widetilde{F}_k(0,0)+\\&
+P_{k,(\alpha',\alpha_n+1)}\left(\partial^{\gamma}\widetilde{B}_{j}^{kl},\cdots,
\partial^{\delta}\widetilde{C}^{kl},\cdots \partial^{\beta}\widetilde{U}_h\right)_{x=0,t=0}=\\&
=\partial_t^{\alpha_n+1}\partial_x^{\alpha'} \widetilde{U}_k(0,0).
\end{aligned}
\end{equation}

\bigskip

\noindent\textbf{Step III.} We may assume that for appropriate positive numbers $M_1,M_2$ and $\rho_1,\rho_2$ with $M_1\geq 1$, we have

\begin{equation}\label{Stella-48C}
B_j^{lk}, C^{lk}\in \mathcal{C}_{M_1,\rho_1}(0), \quad j=1,\cdots,
n-1,\quad l,k=1,\cdots,N,
\end{equation}

\begin{equation}\label{Stella1-48C}
F_k\in \mathcal{C}_{M_2,\rho_2}\quad k=1,\cdots,N.
\end{equation}
We set

$$M=\frac{M_2}{M_1},\quad \rho=\min\{\rho_1,\rho_2\}.$$

By Proposition \ref{prop1-28C}, we may choose, for $j=1,\cdots, n-1$, $l,k=1,\cdots,N$,

$$\widetilde{B}_j^{lk}=\widetilde{C}^{lk}=\frac{M_1\rho}{\rho-(\sigma^{-1}t+x_1+\cdots+x_{n-1})}$$
and
$$\widetilde{F}_k=\frac{M_2\rho}{\rho-(\sigma^{-1}t+x_1+\cdots+x_{n-1})},$$
where $\sigma\in (0,1]$ is to be chosen. System 
\eqref{1-44C} becomes, for $k=1,\cdots, N$,
\begin{equation}\label{1-48C}
	\begin{aligned}
	&\partial_t\widetilde{U}_k(x,t)=\\&=\frac{M_1\rho}{\rho-(\sigma^{-1}t+x_1+\cdots+x_{n-1})}
\left(\sum_{j=1}^{n-1} \sum_{l=1}^{N}\partial_{x_j} \widetilde{U}_l+
\sum_{l=1}^{N} \widetilde{U}_l+M\right).
\end{aligned}
\end{equation}
At this point we search for a solution to equation \eqref{1-48C}
of the form

\begin{equation}\label{2-48C}
\widetilde{U}_k(x,t)=w\left(\sigma^{-1}t+x_1+\cdots+x_{n-1}\right).
\end{equation}
Set $s=\sigma^{-1}t+x_1+\cdots+x_{n-1}$ and

\begin{equation}\label{phi-48C}
\phi(s)=\frac{M_1}{1-\rho^{-1}s}
\end{equation}
we have

\begin{equation}\label{3-48C}
\left(\sigma^{-1}-\phi(s)N(n-1)\right)\frac{dw}{ds}=N\phi(s)w+M\phi(s).
\end{equation}
Now we choose $\sigma>0$ so that $$\sigma^{-1}-\phi(s)N(n-1)>0$$
in a neighborhood of $0$. For instance, we choose 
\begin{equation}\label{sigma0-48C}
\sigma=\sigma_0:=\frac{1}{2NM_1(n-1)}.
\end{equation}
We get

$$\sigma_0^{-1}-\phi(s)N(n-1)=NM_1(n-1)\left(1-2s/\rho\right)>0,\quad\mbox{ for } |s|<\frac{\rho}{2}.$$

Setting
$$h(s)=\frac{\phi}{\sigma_0^{-1}-N(n-1)\phi}=\frac{1}{N(n-1)}\frac{\rho}{\rho-2s},$$
equation \eqref{3-48C} becomes

\begin{equation}\label{3n-48C}
\frac{dw}{ds}=Nh(s)w+Mh(s).
\end{equation}
Let $w_0$ the solution to \eqref{3n-48C} such that

\begin{equation}\label{3nn-48C}
w_0(0)=0.
\end{equation}
We have

\begin{equation}\label{3nnn-48C}
w_0=\frac{M}{N}\left[\exp\left(N\int^s_0h(\eta)d\eta\right)
-1\right]=
\frac{M}{N}\left[\left(\frac{\rho}{\rho-2s}\right)^{\frac{\rho}{2(n-1)}}-1\right].
\end{equation}
In particular, $w_0$ is analytic in
$\left(-\frac{\rho}{2},\frac{\rho}{2}\right)$ and

\begin{equation}\label{w-0-48C}
0\preccurlyeq w_0.
\end{equation}
The latter relationship can be easily checked by using formula \eqref{3nnn-48C} or can be also easily derived from \eqref{3n-48C} and \eqref{3nnn-48C}, by expressing the derivatives of
$w_0$ in $0$ by means of those of lower order and noticing
that they are all nonnegative.
 Now, for $k=1,\cdots,N$, let us consider the following functions
\begin{equation}\label{u-tilde-k-48C}
\widetilde{U}_k(x,t)=w_0\left(\sigma_0^{-1}t+x_1+\cdots+x_{n-1}\right),
\end{equation}
we have that $\widetilde{U}_k$  are solutions to equations
\eqref{1-48C}  (when $\sigma=\sigma_0$) and by Proposition
\ref{prop1-29C}, they are analytic. Moreover by \eqref{3-3N} we have

\begin{equation}\label{1ex-48C}
\begin{aligned}
\widetilde{\varphi}_k(x)&=\widetilde{U}_k(x,0)=w_0\left(x_1+\cdots+x_{n-1}\right)=\\&
=\sum_{m=0}^{\infty}w_0^{(m)}(0)\sum_{|\alpha'|=m}\frac{1}{\alpha'!}x^{\alpha'}.
\end{aligned}
\end{equation}
From which, taking into account \eqref{w-0-48C}, it is obvious that

\begin{equation}\label{100-48C}
0\preccurlyeq \widetilde{\varphi}_k,\quad\mbox{ per } k=1,\cdots,N.
\end{equation}
All in all, $\widetilde{U}$ is an analytic solution in a
neighborhood of $(0,0)$, of the Cauchy problem

\begin{equation}\label{1-44C-corr}
\begin{cases}
\partial_t\widetilde{U}_k=H(x,t)\left(\sum_{j=1}^{n-1} \sum_{l=1}^{N}\partial_{x_j} \widetilde{U}_l+
\sum_{l=1}^{N} \widetilde{U}_l+M\right),
  \\
\\
\widetilde{U}_k(x,0)=\widetilde{\varphi}_k(x),
\end{cases}
\end{equation}
where
$$ H(x,t)=\frac{M_1\rho}{\rho-(\sigma_0^{-1}t+x_1+\cdots+x_{n-1})}.$$
Since $\widetilde{U}$ is analytic, \eqref{5n-47C} implies that the following power series converges in a neighborhood $\mathcal{U}$ of $(0,0)$ 
$$\sum_{\alpha\in \mathbb{N}_0^n}\frac{1}{\alpha !}U_{\alpha}x^{\alpha'}t^{\alpha_n},$$
in addition its sum, $U$, solves Cauchy problem \eqref{1-43C} in
$\mathcal{U}$.

\bigskip

\noindent\textbf{Step IV.} The uniqueness to Cauchy problem  \eqref{1-43C} in the class of
analytic functions in a connected neighborhood of $0$  is a consequence of \textbf{Step I} and of the unique continuation property for analytic functions.
As a matter of fact, if $V'$
$V''$ are analytic solutions of \eqref{1-43C},+ then
$$\partial^{\alpha}V'_k(0,0)=U_{k,\alpha}=\partial^{\alpha} V''_k(0,0),\quad k=1,\cdots,N$$
so that, by Theorem \ref{teo1-23C}, we have $V'=V''$ in a neighborhood of $0$.

$\blacksquare$

\bigskip

\textbf{Remark on the neighborhood in which there exist
	solutions of the Cauchy problem}

In what follows we will be interested in having some detailed information about the dependence of the neighborhood $\mathcal{U}$ by the known term $f$ and, consequently, by the initial data of Cauchy problem \eqref{3-34C}. From  Step III of the previous proof we can say that
the neighborhood $\mathcal{U}$ does not depend on the constant $M_2$. To
clarify what we have just claimed, it suffices to apply Proposition
\ref{prop1-30C} to the composite function

$$w_0\left(\sigma_0^{-1}t+x_1+\cdots+x_{n-1}\right).$$ 
Let us observe that by \eqref{3nnn-48C}, setting
$\kappa=\frac{\rho}{2(n-1)}$, we have there exists a constant
$c_{\kappa}\geq 1$ such that
$$0\leq w_0^{(m)}(0)\leq  \frac{c_{\kappa}M}{N} \left(c_{\kappa}\rho^{-1}\right)^m m!,\quad \forall m\in \mathbb{N}_0.$$
Hence $$w_0\in\mathcal{C}_{\frac{c_{\kappa}M}{N},
\frac{\rho}{c_{\kappa}}}.$$ On the other hand we get trivially

$$\left(\sigma_0^{-1}t+x_1+\cdots+x_{n-1}\right)\in \mathcal{C}_{(\sigma_0^{-1}+n-1),1}.$$
Hence $$w_0\left(\sigma_0^{-1}t+x_1+\cdots+x_{n-1}\right)\in
\mathcal{C}_{\frac{c_{\kappa}M}{N},R},$$ where
$$R=\frac{\rho}{c_{\kappa}\sigma_0^{-1}+n-1}$$
and by \eqref{sigma0-48C} it turns out that $\sigma_0$ does not depend by
$M_2$. Therefore $R$ depends on $M_1$, $\rho$,  $n$ and $N$ only. Moreover, we can choose  $\mathcal{U}=\{(x,t)\in
\mathbb{R}^n:|t|+|x_1|+\cdots+|x_{n-1}|<R\}$ as the neighborhood in which the
Cauchy problem \eqref{1-43C} admits a solution.

\bigskip
In preparation for what we will do later, let us go back to consider the
case of a linear differential operator of order $m$ given by
\begin{equation}\label{0-51C}
P(x,\partial)=\sum_{|\alpha|\leq m}a_{\alpha}(x)\partial^{\alpha}
\end{equation}
and let us consider the following Cauchy problem

\begin{equation}\label{1-51C}
\begin{cases}
P(x,\partial)u(x)=f(x), \\
\\
\partial^j_n u(x',0)=0, \quad j=0,1,\cdots,m-1, \forall x'\in
B'_1.
\end{cases}
\end{equation}

Let us suppose that, for given $M_0,\rho_0$, we have
\begin{equation}\label{2a-51C}
P_m((x',0),e_n)\neq 0,\quad\forall x'\in \overline{B'_1},
\end{equation}

\begin{equation}\label{2b-51C}
a_{\alpha}\in \mathcal{C}_{M_0,\rho_0}(z),\quad |\alpha|\leq m, \quad\forall z\in
\overline{B'_1}\times[-\delta_0,\delta_0],
\end{equation}

\begin{equation}\label{2c-51C}
f\mbox{ be a polynomial}.
\end{equation}
Then the solution to problem \eqref{1-51C} there
exists in $\overline{B'_1}\times[-\delta,\delta]$ (actually,
in a neighborhood of  $\overline{B'_1}\times\{0\}$) where $\delta>0$
depends on $M_0,\rho_0,\delta_0$ and $\min\{|P_m((x',0),e_n)|:
x'\in \overline{B'_1}\}$, \textbf{but does not  depend by the polynomial}
$f$. In order to check this assertion, let $h$ be the degree of the polynomial $f$. Set

$$K=1+\sum_{|\beta|\leq h}\max_{\overline{B'_1}\times[-\delta_0,\delta_0]}\left|\partial^{\beta}f\right|.$$
It is evident that, setting.
$$\widetilde{u}=\frac{u}{K}, \quad\quad \widetilde{f}=\frac{f}{K},$$
$u$ solves Cauchy problem \eqref{1-51C} if and only if $v$ is solves the following Cauchy problem

\begin{equation}\label{1-51bisC}
\begin{cases}
P(x,\partial)\widetilde{u}(x)=\widetilde{f}(x), \\
\\
\partial^j_n \widetilde{u}(x',0)=0, \quad j=0,1,\cdots,m-1, \mbox{ }\forall x'\in
B'_1.
\end{cases}
\end{equation}
On the other hand, because of the way we defined $\widetilde{f}$
we can certainly state that

\begin{equation}\label{2-51bisC}
\widetilde{f}\in \mathcal{C}_{1,1}(z), \quad\forall z\in
\overline{B'_1}\times[-\delta_0,\delta_0].
\end{equation}
We can then return to problem \eqref{1-43C}. Hence by
applying the Cauchy--Kovalevskaya Theorem and taking into account that
$a_{\alpha}$, $|\alpha|\leq m$, are analytic functions in a neighborhood of
$\overline{B'}_1(0)\times\{0\}$, we conclude that
the solution of Cauchy problem \eqref{1-51C} exists and it is
analytic in $\overline{B'_1}\times[-\delta,\delta]$ where
$\delta>0$ depends on $M_0,\rho_0,\delta_0$ and on
$$\min\{|P_m((x',0),e_n)|: x'\in \overline{B'_1}\}$$ (but does not depend on
$f$).

Obviously, if \eqref{2a-51C}, \eqref{2b-51C},
\eqref{2c-51C} hold,
similar conclusions are valid to the Cauchy problem

\begin{equation}\label{1-52C}
\begin{cases}
P(x,\partial)u(x)=f(x), \quad \mbox{ se } x\in \mathcal{U}_{\widetilde{x}}, \\
\\
\partial^j_n u(x',0)=g_j, \quad j=0,1,\cdots,m-1, \mbox{ }\forall x'\in
B'_1.
\end{cases}
\end{equation}
provided that $g_j$ are polynomials for $j=0,1,\cdots,m-1$.
$\blacklozenge$

\bigskip

\section[Further comments on the Cauchy--Kovalevskaya Theorem]{Further comments on the Cauchy--Kovalevskaya Theorem. Examples} \label{ulteriori}
\subsection{A few brief note on the qualisinear and the nonlinear case} \label{nonlinear}

It is not difficult to adapt the proof of the Cauchy--Kovalevskaya Theorem 
 to the case of a quasilinear operator

$$\sum_{|\alpha|=m}a_{\alpha}\left(x,(\partial^{\beta}u)_{||\beta|\le m-1}\right)\partial^{\alpha}u+
a_{0}\left(x,(\partial^{\beta}u)_{||\beta|\le m-1}\right).$$ In
this case, we recall, the Cauchy problem is
\begin{equation}\label{1-24-5}
\begin{cases}
\mathcal{P}(u)=0, \\
\\
\partial_n^ju(x)=g_j(x), \quad j=0,1,\cdots,m-1, \mbox{ }\forall x\in\Gamma.
\end{cases}
\end{equation}
One can proves that if $a_{\alpha}, a_0, g_j$ are 
analytic functions, $\Gamma$ is analytic and noncharacteristic
 (Definition \ref{def-quasilin-5C}) then for every  $\widetilde{x}\in \Gamma$ there exists a neighborhood
$\mathcal{U}_{\widetilde{x}}$ in which Cauchy problem
\eqref{1-24-5} has analytic solution. For the proof  we refer to \cite{EV}. \index{Cauchy problem:@{Cauchy problem:}!- general quasilinear equations@{- general quasilinear equations}}

In \eqref{4-1C} we have formulated the Cauchy problem for the fully nonlinear equation
\index{Cauchy problem:@{Cauchy problem:}!- general fully nonlinear equations@{- general fully nonlinear equations}} 
\begin{equation}\label{1-55C}
\begin{cases}
F\left(x,(\partial^{\alpha} u)_{|\alpha|\leq m}\right)=0,  \\
\\
\frac{\partial^ju(x)}{\partial \nu^j}=g_j(x), \quad
j=0,1,\cdots,m-1,\quad \forall x\in \Gamma.
\end{cases}
\end{equation}
Here we only outline the proof of the existence of the solutions to problem \eqref{1-55C}
referring for more details to \cite[Chapter 1]{Co-Hi}, \cite[Chapter 1]{Fo}.

Let us consider the case where $\Gamma=\{x_n=0\}$. We know that
we may always reduce to this case by "flattening" $\Gamma$
(by map \eqref{2-24-5}). So, instead of the
 conditions $\frac{\partial^ju(x)}{\partial \nu^j}=g_j(x)$,
$j=0,1,\cdots,m-1$, for $x\in \Gamma$, we may consider

\begin{equation}\label{2-55C}
\partial_n^ju(x',0)=g_j(x') \quad j=0,1,\cdots,m-1,\quad \forall x'\in
B'_r,
\end{equation}
where $r>0$. 
In the first part of the proof of Proposition
\ref{prop-6C} we have already seen that conditions \eqref{2-55C}
allow us to determine the derivatives

\begin{equation}\label{3-55C}
	\partial_{x'}^{\alpha'}\partial_n^ju(x',0)=\partial_{x'}^{\alpha'}g_j(x')
	\quad j=0,1,\cdots,m-1,\ \alpha'\in \mathbb{N}_0^{n-1}\
	x'\in B'_r.
\end{equation}
\smallskip

\noindent \textit{without involving the equation} $F\left(x,(\partial^{\alpha} u)_{|\alpha|\leq m}\right)=0$. Let us recall that in order to calculate the derivative $\partial_n^m u(x',0)$ we need to use
the equation. More precisely, we have

\begin{equation}\label{4-55C}
F\left(x',0,(\partial_{x'}^{\alpha'} g_j)_{|\alpha'|+j\leq m, j\leq
m-1}(x'),\partial_n^m u(x',0)\right)=0.
\end{equation}
To find $z=\partial_n^m u(x',0)$ from equation \eqref{4-55C} we need that the equation

\begin{equation}\label{1-56C}
F\left(x',0,(\partial_{x'}^{\alpha'} g_j)_{|\alpha'|+j\leq m, j\leq
m-1}(x'),z\right)=0,
\end{equation}
admits a solution.
If, for instance, we require 
\begin{equation}\label{2-56C}
(\partial_zF)\left(x',0,(\partial_{x'}^{\alpha'}
g_j)_{|\alpha'|+j\leq m, j\leq m-1}(x'),z\right)\neq 0
\end{equation}
then we may express $z$ as a
a function of variable $x'$ of class $C^1$ (provided $F$ is of class $C^1$) or
as an analytic function, provided  $F$ is analytic. Let us observe that 
condition \eqref{2-56C} makes it possible to write, in a neighborhood of
$\Gamma$, the equation

$$F\left(x,(\partial^{\alpha} u)_{|\alpha|\leq m}\right)=0$$
like

\begin{equation}\label{3-56C}
\partial_n^m u=G\left(x,(\partial^{\alpha} u)_{|\alpha|\leq m,
\alpha_n<m}\right),
\end{equation}
where $G$ is analytic (provided $F$ is analytic). 
Hence, in the fully nonlinear case, condition
\eqref{2-56C} may replace  the condition that $\Gamma$ is
noncharacteristic surface for a linear (or quasilinear) operator.
Furthermore, since we can find $\partial_n^m u(x',0)$
from \eqref{4-55C}, we set $g_m(x')=\partial_n^m u(x',0)$,
by making the derivatives of both the sides of  \eqref{3-56C} w.r.t. $x_n$ we have

\begin{equation*}
\begin{cases}
\partial_n^{m+1}u=\sum_{|\alpha|\leq m, j<m}a_{\alpha',j}\partial_{x'}^{\alpha'}\partial_n^{j+1}u+
(\partial_n G)\left(x,(\partial^{\alpha}u)_{||\alpha|\le m-1, \alpha_n}\right),  \\
\\
\frac{\partial^ju(x',0)}{\partial \nu^j}=g_j(x'), \quad
j=0,1,\cdots,m,\quad \forall x\in \Gamma.
\end{cases}
\end{equation*}
where
$$a_{\alpha',j}=\left(\partial_{p_{\alpha',j}} G \right)
\left(x,(\partial_{x'}^{\alpha'}\partial_n^ju)_{||\alpha'|+j \leq m,
j<m}\right),$$ for $|\alpha'|+j \leq m, j<m$. 

In other words, if  conditions \eqref{1-56C} and \eqref{2-56C} hold we may
reformulate Cauchy problem \eqref{1-55C} as a Cauchy problem for a quasilinear equation of order $m+1$ and we may use the existence results that obtained in the quasilinear case.

\subsection{Comments about the existence and the uniqueness of solutions. Examples and counterexamples}\label{esist-unic}
In this Subsection we return to the \textbf{linear case.} In general, if $\Gamma$ is a characteristic surface for the operator $P(x,\partial)$, neither existence nor
uniqueness for Cauchy problem \eqref{3-34C} can be expected. Let us look at some example.

\noindent\textbf{(a)} Let $P_m(\xi)$ be a homogeneous polynomial of degree
$m$ and let us assume
$$P_m(N)=0,$$
where $N\in \mathbb{R}^n\setminus\{0\}$. Then the hyperplane
$$\pi=\{x\cdot N=0\}$$ is a characteristic surface for the operator $P_m(\partial)$. Now, let us consider the functions

$$u^{(t)}(x)=e^{tx\cdot N}-\sum_{k=0}^{m-1}\frac{t^k(x\cdot N)^k}{k!},\quad t\in \mathbb{R}$$
where $t\neq 0$. It easy to check that, for every  $t\in \mathbb{R}$, $u^{(t)}$ solves the  Cauchy problem

\begin{equation}\label{1-53C}
\begin{cases}
P_m(\partial)u(x)=0, \\
\\
\frac{\partial^ju}{\partial N^j}=0, \quad j=0,1,\cdots,m-1, \mbox{
on }\pi,
\end{cases}
\end{equation}
Hence, does not hold the uniqueness for problem \eqref{1-53C}.

\medskip

\noindent\textbf{(b)} The Cauchy--Kovalevskaya Theorem gives
the existence and uniqueness of solutions to the Cauchy problem
having initial data on a noncharacteristic surface \textit{under the assumption
of analyticity of all the data of the problem}. Regarding the
existence of the solutions, if we desire to preserve the same
generality of the Theorem, the  assumptions of analyticty cannot
be reduced. To prove this, it suffices to consider the following
Cauchy problem 

\begin{equation}\label{1-57C}
\begin{cases}
\partial^2_yu+\partial^2_xu=0, \\
\\
u(x,0)=g_0(x), \mbox{ for } x\in (-r,r),\\
\\
\partial_yu(x,0)=g_1(x), \mbox{ for } x\in(-r,r).
\end{cases}
\end{equation}
Since $u$ satisfies the Laplace equation 
$\partial^2_yu+\partial^2_xu=0$ in a neighborhood of $(0,0)$ it is
analytic in such a neighborhood, therefore the initial data,  $u(x,0)=g_0(x)$
and $\partial_yu(x,0)=g_1(x)$ also need to be analytic.

Incidentally, even though we consider the "one--sided" Cauchy problem \index{Cauchy problem:@{Cauchy problem:}!- one--sided@{- one--sided}}the situation
do not change in a significant way. Let us consider, indeed, the problem
\begin{equation*}
\begin{cases}
\partial^2_yu+\partial^2_xu=0, \mbox{ for } (x,y)\in B_r^+,\\
\\
u(x,0)=g_0(x), \mbox{ for } x\in (-r,r),\\
\\
\partial_yu(x,0)=g_1(x), \mbox{ for } x\in(-r,r),
\end{cases}
\end{equation*}
where $u\in C^2(B_r^+)\cap C^0(\overline{B_r^+})$. Let us suppose for simplicity that
$g_1\equiv0$ and let us consider the even reflection w.r.t. $x$--axis of
$u$ 
$$v(x,y)=u(x,|y|)$$ then, by the Schwarz reflection principle, we have
$$\Delta v=0, \mbox{ in } B_r$$
and again we have that $v$ is an analytic function and therefore
$g_0$ is an analytic function.

\bigskip

\underline{\textbf{Exercise 1.}} Let us suppose $g_0,g_1\in C^{\omega}(-r,r)$ in \eqref{1-57C}.
Prove that the solution to Cauchy problem \eqref{1-57C} is given by

\begin{equation}\label{0-58C}
u(x,y)=\sum_{n=0}^{\infty}(-1)^n\left(\frac{g_0^{(2n)}(x)y^{2n}}{(2n)!}+\frac{g_1^{(2n+1)}(x)y^{2n+1}}{(2n+1)!}\right).
\end{equation}
$\clubsuit$

\bigskip

\bigskip

\noindent\textbf{(c)} Let us consider the Cauchy problem

\begin{equation}\label{2-53C}
\begin{cases}
\partial_tu-\partial^2_xu=0, \\
\\
u(x,0)=g_0(x), \mbox{ for } x\in\mathbb{R},\\
\\
\partial_tu(x,0)=g_1(x), \mbox{ for } x\in\mathbb{R}.
\end{cases}
\end{equation}
In this case the straight line $\{t=0\}$ is a characteristic line for
the operator $\partial_t-\partial^2_x$. It is evident that if we
do not require

\begin{equation}\label{3-53C}
g_1(x)=g''_0(x), \mbox{ }\forall x\in\mathbb{R},
\end{equation}

\smallskip

\noindent then Cauchy problem \eqref{2-53C} has no solutions. As a matter of fact, if problem
 \eqref{2-53C} admits solutions, even just of class
$C^2$, then
$$g''_0(x)=\partial^2_xu(x,0)=\partial_tu(x,0)= g_1(x), \mbox{ for } |x|<1$$
and therefore \eqref{3-53C} hold.

We now check that even condition \eqref{3-53C} is satisfied
we can exhibit a $g_0\in C^{\omega}$ such that problem
\eqref{2-53C} has no solutions. First, it is evident that if
\eqref{3-53C} is satisfied, then problem \eqref{2-53C} can be
be formulated as

\begin{equation}\label{2n-53C}
\begin{cases}
\partial_tu-\partial^2_xu=0, \\
\\
u(x,0)=g_0(x), \mbox{ for } x\in\mathbb{R}.\\
\end{cases}
\end{equation}
Let
$$g_0(x)=\frac{1}{1+x^2}.$$
We have that $g_0\in C^{\omega}(\mathbb{R})$ and

$$g^{(2k)}_0(0)=(-1)^{k}(2k)!\ , \quad\mbox{ for } k\in \mathbb{N}_0.$$
Now, if there exists a solution $u$ \textbf{analytic} in a neighborhood
of $(0,0)$ of problem \ref{2n-53C}, then $u(0,t)$ should
be expanded in Taylor series in $t=0$

$$\partial^j_tu(0,0)=\partial^{2j}_xu(0,0)=g^{(2j)}_0(0)=(-1)^{j}(2j)!,$$
on the other side, the power series

$$\sum_{j=0}^{\infty}\frac{(-1)^{j}(2j)!}{j!}t^j$$
has the radius of convergence equal to zero. Hence $u(0,t)$ is not analytic.

Also, we note that the Cauchy problem

\begin{equation}\label{0-54C}
\begin{cases}
\partial_tu-\partial^2_xu=0, \\
\\
u(x,0)=0, \mbox{ for } x\in\mathbb{R}.\\
\end{cases}
\end{equation}
admits only one analytic solution, namely the null solution.
As a matter of fact, the null function is trivially solution of problem
\eqref{0-54C} and if $u$ is an analytic solution of
\eqref{0-54C} then we have, for each $i,j\in \mathbb{N}_0$

$$\partial^{i}_x\partial^j_tu(0,0)=\partial_x^{i+2j}u(0,0)=0$$
hence $u\equiv 0$.

Of course the one-sided Cauchy problems are also of interest
\begin{equation}\label{00-54C}
\begin{cases}
\partial_tu-\partial^2_xu=0,\mbox{ for } x\in\mathbb{R}, \quad t>0 \\
\\
u(x,0)=g_0(x), \mbox{ for } x\in\mathbb{R}.\\
\end{cases}
\end{equation}
For this problem it turns out that (if $g_0$ is regular enough)
$$u(x,t)=\frac{1}{\sqrt{4\pi}}\int^{+\infty}_{-\infty}g_0(\xi)e^{-\frac{(x-\xi)^2}{4t}}d\xi$$
is a solution to \eqref{00-54C}. In particular if
$g_0=\frac{1}{1+x^2}$, then Cauchy problem \eqref{00-54C} admits
solutions (of course,  nonanalytic w.r.t. $t$).

Keep in mind that the problems that we have considered in the
point (c) are not written in the form \eqref{1-39C}. As a matter of fact,  the term on the right--hand side of equation
$\partial_tu=\partial^2_xu$ has order $2$ greater than the order of
derivative  $\partial_tu$, on the left--hand side.

\bigskip

\noindent\textbf{(d)} In a strong contrast with the example considered in (b)
we present now the following example for the vibrating string equation. Let us consider the Cauchy problem

\begin{equation}\label{1-25-5}
\begin{cases}
\partial^2_tu-\partial^2_xu=0, \\
\\
u(x,0)=g_0(x), \mbox{ for } x\in (-1,1),\\
\\
\partial_tu(x,0)=g_1(x), \mbox{ for } x\in(-1,1).
\end{cases}
\end{equation}

Let us first prove problem \eqref{1-25-5} has at most
one solution $u\in C^2(\overline{Q})$ where  $Q=\{(x,t)\in
\mathbb{R}^2: |x|+|t|\leq 1\}$.

\begin{figure}\label{figura-Q}
	\centering
	\includegraphics[trim={0 0 0 0},clip, width=12cm]{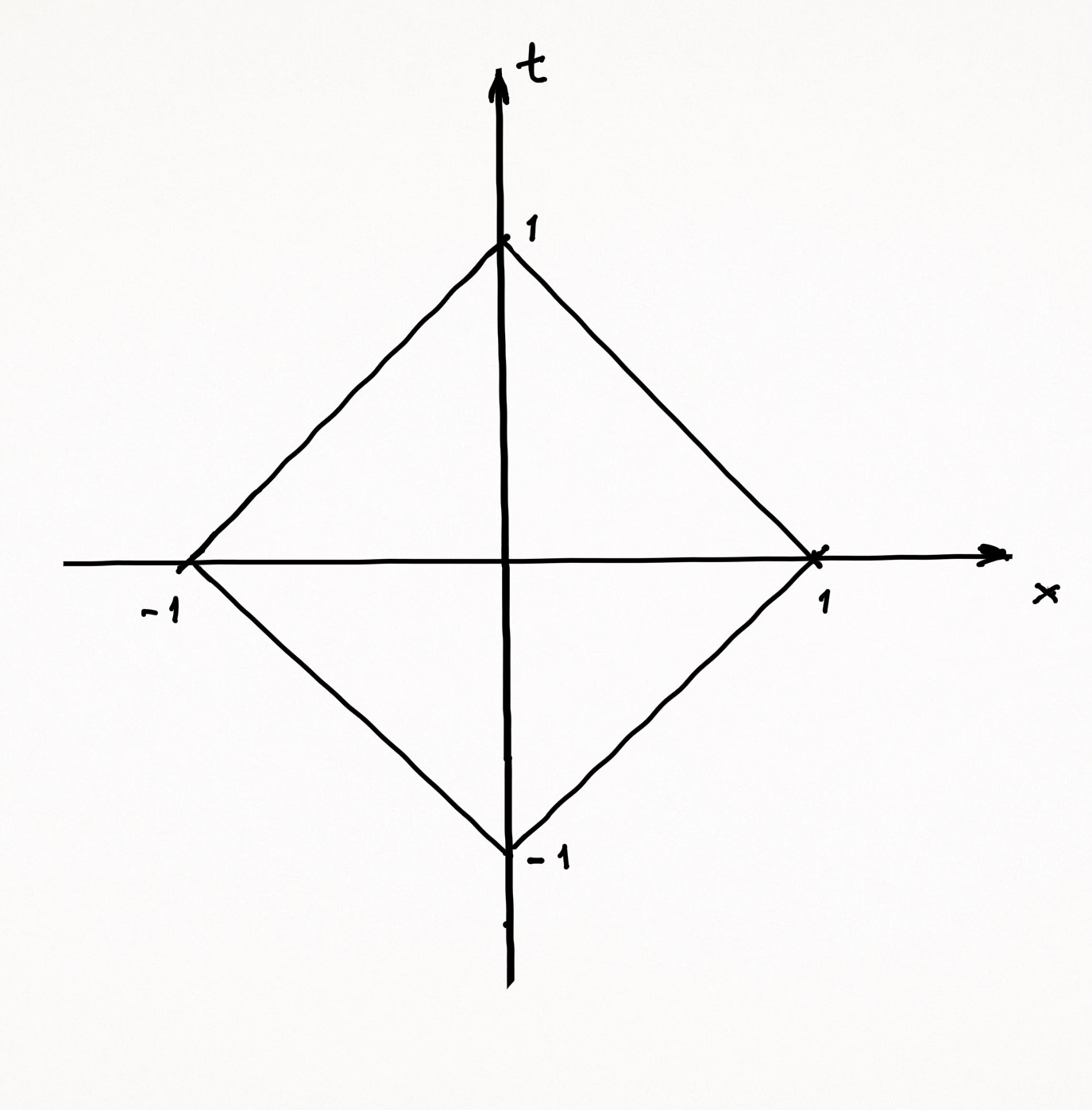}
	\caption{$Q=\{(x,t)\in
		\mathbb{R}^2: |x|+|t|\leq 1\}$}
\end{figure}

To this purpose we observe that

\begin{equation}\label{2-25-5}
0=\left(\partial^2_tu-\partial^2_xu\right)\partial_tu=\frac{1}{2}\left(\partial_t(\partial_tu)^2-2\partial_x(\partial_xu\partial_tu)+
\partial_t(\partial_xu)^2\right).
\end{equation}
Moreover, let us suppose that $g_0=0$ and $g_1=0$ in \eqref{1-25-5}. We integrate both the sides of 
\eqref{2-25-5} over

$$\mathcal{T}_{\delta}=Q\cap\{(x,t)\in \mathbb{R}^2:
0<t<1-\delta\},$$ where $\delta\in (0,1)$ is arbitrary. We obtain, by the divergence Theorem,

\begin{equation}\label{3-25-5}
\begin{aligned}
&0=\iint_{\mathcal{T}_{\delta}}\left(\partial^2_tu-\partial^2_xu\right)\partial_tudxdt=\\&=\frac{1}{2}\int_{\partial\mathcal{T}_{\delta}}
\left((\partial_tu)^2\nu_t-2(\partial_xu\partial_tu)\nu_x+
(\partial_xu)^2\nu_x\right)dS=\\& = -\frac{1}{2}\int_{-1}^1
\left((\partial_tu)^2(x,0)+ (\partial_xu)^2(x,0)\right)dx+\\&
+\frac{1}{2\sqrt{2}}\int_{1-\delta}^{1}
\left(\partial_tu(x,1-x)-\partial_xu(x,1-x)\right)^2dx+\\&
+\frac{1}{2}\int_{-1+\delta}^{1-\delta}
\left((\partial_tu)^2(x,\delta)+
(\partial_xu)^2(x,\delta)\right)dx+\\&
+\frac{1}{2\sqrt{2}}\int_{-1}^{-1+\delta}
\left(\partial_tu(x,1+x)+\partial_xu(x,1+x)\right)^2dx\geq\\& \geq
\frac{1}{2}\int_{-1+\delta}^{1-\delta}
\left((\partial_tu)^2(x,\delta)+ (\partial_xu)^2(x,\delta)\right)dx.
\end{aligned}
\end{equation}

\smallskip

\noindent Hence
$$\int_{-1+\delta}^{1-\delta}
\left((\partial_tu)^2(x,\delta)+
(\partial_xu)^2(x,\delta)\right)dx=0$$ from which we have
$(\partial_tu)^2(x,\delta)+ (\partial_xu)^2(x,\delta)=0$ and, since  $\delta$ is arbitrary, we have $\partial_xu=\partial_tu=0$
in $Q\cap \{(x,t)\in \mathbb{R}^2: 0\leq t\}$. Finally, since
$u(x,0)=\partial_xu(x,0)=0$, we have $u=0$ in $Q\cap \{(x,t)\in
\mathbb{R}^2: 0\leq t\}$. Similarly, we obtain $u=0$ in
$Q\cap \{(x,t)\in \mathbb{R}^2: 0\geq t\}$. Therefore $u=0$ in $Q$.

The existence of solutions also does not require that the Cauchy data
$g_0$ and $g_1$ are analytic.  As a matter of fact it is checked straightforwardly that,
if $g_0\in C^2([-1,1])$ and $g_1\in C^1([-1,1])$, then the solution
to Cauchy problem \eqref{1-25-5} is given by

\begin{equation}\label{4-25-5}
u(x,t)=\frac{g_0(x+t)+g_0(x-t)}{2}+\frac{1}{2}\int^{x+t}_{x-t}g_1(\eta)d\eta.
\end{equation}

\bigskip

\textbf{Conclusion.} From the short discussion of this Section
we can say that, in the context of partial differential equations, the Cauchy--Kovalevskaya Theorem represents for us more a starting point than an ending point.
Starting with the next Section we will focus more on the issue
of the uniqueness, under assumptions which are weaker than the analyticity of all the data.

\section{The Holmgren Theorem} \label{Holmgren}
Let us start by recalling the \textbf{divergence Theorem}. Let
$D$ a bounded open $\mathbb{R}^n$ such that its boundary
$\partial D$ is of class $C^{0,1}$. Then we have

\begin{equation}\label{1-60C}
\int_{D}\partial_j udx=\int_{\partial D}u\nu_j dS, \quad
j=1,\cdots,n\quad \forall u\in C^1\left(\overline{D}\right),
\end{equation}
where $\nu=(\nu_1,\cdots,\nu_n)$ is the unit outward normal to
$\partial D$ and $dS$ is the $(n-1)$--element of surface.

Let $m\in \mathbb{N}$, $a_{\alpha}\in
C^m\left(\overline{D}\right)$, $|\alpha|\leq m$ and

\begin{equation}\label{2-60C}
P(x,\partial)=\sum_{|\alpha|\leq m}a_{\alpha}(x)\partial^{\alpha}.
\end{equation}
We call the \textbf{(formal) adjoint operator} \index{adjoint (formal)} of $P(x,\partial)$ the following operator

\begin{equation}\label{3-60C}
C^m\left(\overline{D}\right)\ni u \rightarrow
P^{\ast}(x,\partial)u=\sum_{|\alpha|\leq
m}(-1)^{|\alpha|}\partial^{\alpha}\left(a_{\alpha}(x)u\right).
\end{equation}
Let us note that, up to the sign, the principal part of
$P^{\ast}(x,\partial)$ is equal to the principal part of
$P(x,\partial)$.

The following \textbf{Green identity} holds true, for any $u,v\in
C^m\left(\overline{D}\right)$,
\begin{equation}\label{1-61C}
\int_{D}\left(vP(x,\partial)u-uP^{\ast}(x,\partial)v\right)dx=\int_{\partial
D} \mathcal{M}(u,v;\nu) dS,
\end{equation}

\smallskip

\noindent where $\mathcal{M}(u,v;\nu)$ is linear w.r.t. $u$, $v$ and $\nu$. Moreover
\begin{equation}\label{M-61C}
\mathcal{M}(u,v;\nu)=\sum_{|\beta|+|\gamma|\leq
m-1}c_{\beta\gamma}(x)\partial^{\beta}u\partial^{\gamma}v,
\end{equation}
where $c_{\beta\gamma}$, $|\beta|+|\gamma|\leq m-1$, belong
(at least) to $C^0\left(\partial D\right)$. 


Identity \eqref{1-61C} can be  obtained by applying repeatedly the following
simple identity

$$v(x)a(x)\partial_j u(x)=\partial_j\left(v(x)a(x) u(x)\right)-\partial_j\left(v(x)a(x)\right)  u(x)$$
obtaining
\begin{equation}\label{2-61C}
\begin{aligned}
&v(x)a_{\alpha}(x)\partial^{\alpha}
u(x)=v(x)a_{\alpha}(x)\partial_j\partial^{\alpha-e_j} u(x)=\\&
=\partial_j\left(v(x)a_{\alpha}(x) \partial^{\alpha-e_j}
u(x)\right)-\partial_j\left(v(x)a_{\alpha}(x)\right)
\partial^{\alpha-e_j} u(x)= \\& =\cdots=\\&
=\mbox{div}\left(F_{\alpha}\right)+(-1)^{|\alpha|}\partial^{\alpha}\left(a_{\alpha}(x)v(x)\right)u(x),
\end{aligned}
\end{equation}
where $F_{\alpha}$ is a suitable vector field. Next we add up the identities obtained in
\eqref{2-61C}, we integrate over $D$ the obtained new identity, and by the divergence Theorem we get \eqref{1-61C}.

\medskip

Before stating the Holmgren Theorem, let us introduce some
notation.

Let $\Omega$ be an open set of $\mathbb{R}^n$, $x_0\in\Omega$
and let $\phi\in C^2(\Omega)$. Let us denote by
$$\Gamma=\{x\in \Omega: \phi(x)=\phi(x_0)\}.$$
We will assume
\begin{equation}\label{1-62C}
\nabla\phi(x)\neq 0,\quad \forall x\in \Gamma.
\end{equation}
If $\mathcal{U}$ is a neighborhood of $x_0$ we denote by 
$\mathcal{U}_+$ the set

\begin{equation}\label{2-62C}
\mathcal{U}_+=\mathcal{U}\cap\{x\in \Omega: \phi(x)\geq\phi(x_0)\}.
\end{equation}

\begin{theo}[\textbf{Holmgren}]\label{teor-Holm}
Let $a_{\alpha}\in C^{\omega}\left(\Omega\right)$, $|\alpha|\leq m$. Let us suppose that $\Gamma$ is a noncharacteristic surface in $x_0$ for the operator
\begin{equation}\label{3-62C}
P(x,\partial)=\sum_{|\alpha|\leq m}a_{\alpha}(x)\partial^{\alpha}.
\end{equation}
Then there exists a neighborhood  $\mathcal{U}$ of $x_0$ such that we have:

\smallskip
 
 if $u\in C^m\left(\overline{\mathcal{U}_+}\right)$ satisfies
\begin{equation}\label{4-62C}
\begin{cases}
P(x,\partial)u=0, & \mbox{in}\quad \mathcal{U}_+ \\
\\
\partial^{\alpha} u=0, & \mbox{for } |\alpha|\leq m-1,\quad x\in
\Gamma\cap\mathcal{U}_+.
\end{cases}
\end{equation}
Then we have
\begin{equation}\label{4n-62C}
u\equiv 0\quad\mbox{in}\quad \mathcal{U}_+.
\end{equation}
\end{theo}

\textbf{Remarks.} Before starting with the proof of
Theorem \ref{teor-Holm} we observe what follows.

(i) In \eqref{4-62C}, $u$ is required to be a solution
to the equation $P(x,\partial)u=0$ in $\mathcal{U}_+$, unlike the Cauchy-Kovalevskaya Theorem in which $u$ is required to be a solution in a full neighborhood of $x_0$.
Furthermore, it is only required that $u\in
C^m\left(\overline{\mathcal{U}_+}\right)$.

(ii) The initial surface $\Gamma$ is assumed to be of class $C^2$, thus,
not analytic like in the Cauchy-Kovalevskaya Theorem. 
Also, let us note that in \eqref{4-62C} we require $\partial^{\alpha} u=0$ for $|\alpha|\leq m-1$, on $\Gamma\cap\mathcal{U}_+$ and not just that
$$\frac{\partial^ju}{\partial\nu^j}=0, \quad\mbox{for }j=0,1,\cdots, m-1.$$
Of course, if we want to assume the latter conditions we should
require $\phi\in C^{m-1}(\Omega)$ (compare with the proof of the first part of
 Proposition \ref{prop-6C}). $\blacklozenge$

\begin{figure}\label{figura-p63}
	\centering
	\includegraphics[trim={0 0 0 0},clip, width=10cm]{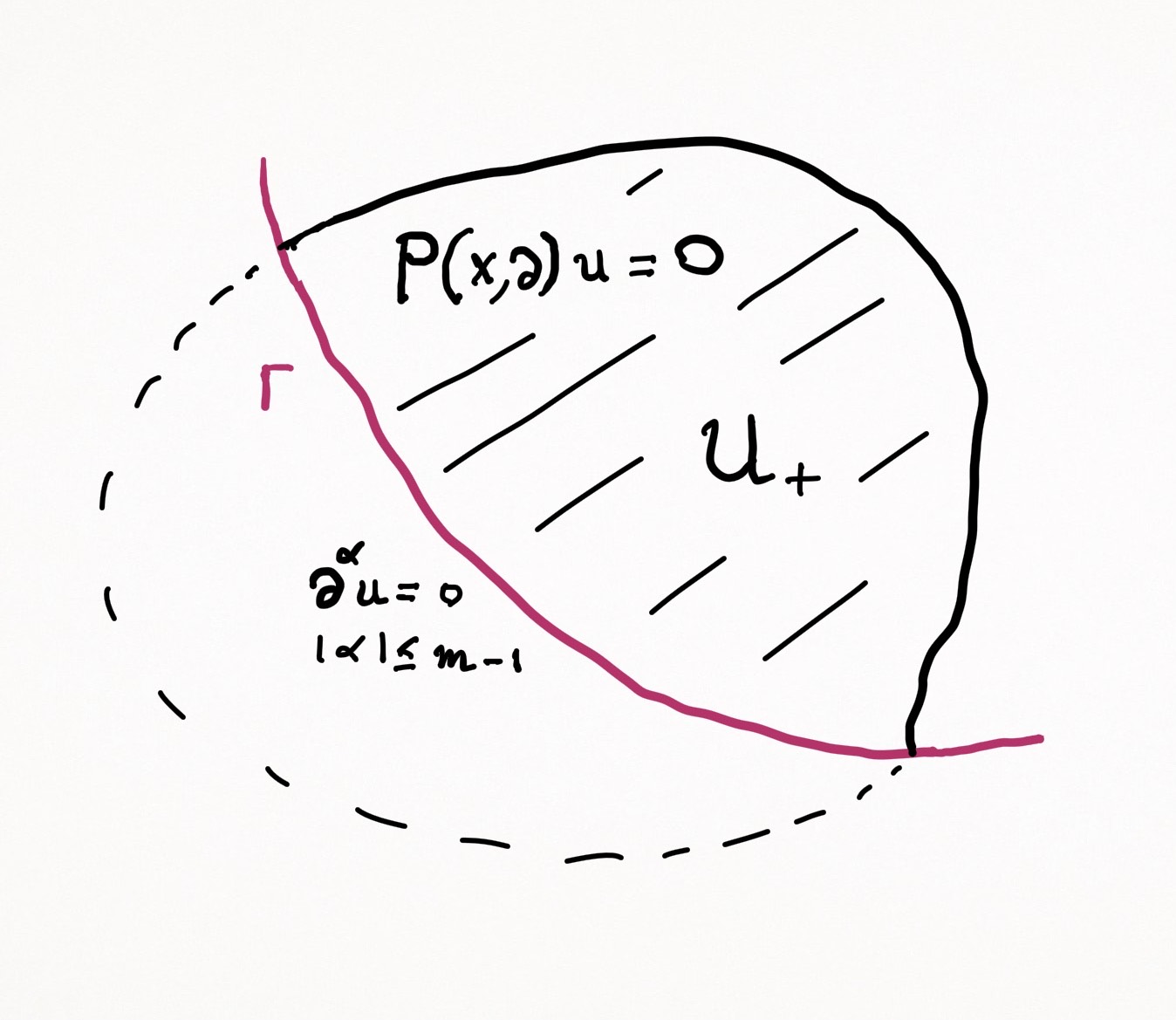}
	\caption{}
\end{figure}

\bigskip

\textbf{Proof of the Holmgren Theorem.}

We may assume $\Gamma$ be a graph of a function. More precisely we may assume that.

\begin{equation}\label{1-63C}
\Gamma=\{(x',\psi(x')): x'\in B'_{r_0}\},
\end{equation}
where $\psi\in C^{2}\left(\overline{B'_{r_0}}\right)$ satisfies

\begin{equation}\label{2-63C}
\psi(0)=\left|\nabla_{x'}\psi(0)\right|=0.
\end{equation}

We will divide the proof of the Theorem into two steps. In the first
step we will assume that $\psi$ is a \textbf{strictly convex function}.
In the second step we will reduce to the first part by means of the so--called Holmgren transformation.

\medskip

\textbf{Step I.}  Let $\psi$ strictly convex, let $R_0$ satisfy

\begin{equation}\label{1-64C}
\Gamma\subset Q_{R_0}:=B'_{R_0}\times(-R_0, R_0)
\end{equation}
and, by assumptions, 
\begin{equation}\label{2-64C}
a_{\alpha}\in C^{\omega}\left(Q_{2R_0}\right).
\end{equation}

\begin{figure}\label{figura-p64}
	\centering
	\includegraphics[trim={0 0 0 0},clip, width=10cm]{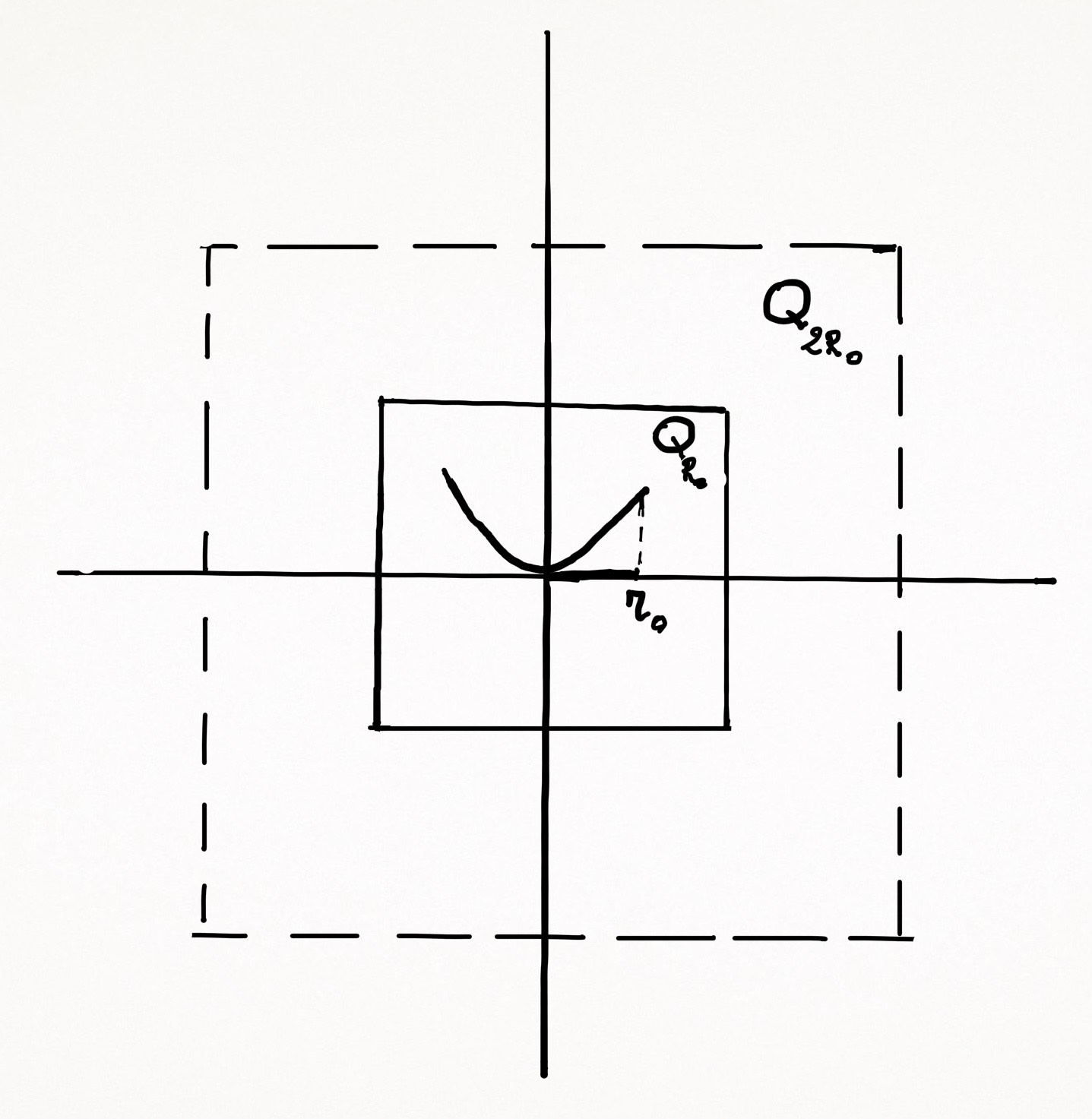}
	\caption{}
\end{figure}

Since $\Gamma$ is a  noncharacteristic surface in $0$ we may assume, recalling \eqref{2-63C},

\begin{equation}\label{3-64C}
\left|P_m(0,e_n)\right|=c_0>0.
\end{equation}
By the continuity of the coefficients of $P_m(x,\partial)$ and by 
\eqref{3-64C} we have that there exists $\rho_1>0$ such that 

\begin{equation}\label{1-65C}
\left|P_m\left((x',h),e_n\right)\right|\geq \frac{c_0}{2}, \quad
\forall x'\in B'_{\rho_1}(0), |h|\leq \rho_1.
\end{equation}
This implies that for every $h\in\left[-\rho_1,\rho_1\right]$ the
flat surface
$$\{(x',h): x'\in B'_{r_0}\}$$
is noncharacteristic.

Let now $f$ be a polynomial. By the Cauchy--Kovalevskaya Theorem and by the Remark which follows such a Theorem, there exists $\rho_2$,
$0<\rho_2<\rho_1$, $\rho_2$ \textit{independent of}  $f$ such that, there exists an analytic solution in
$\overline{B'_{\rho_2}(0)}\times [h-\rho_2,h+\rho_2]$ to the following Cauchy problem

\begin{equation}\label{2-65C}
\begin{cases}
P^{\ast}(x,\partial)w=f, & \\
\\
\partial^{\alpha} w(x',h)=0, & |\alpha|\leq m-1.
\end{cases}
\end{equation}
Let

$$h_0=\min_{\partial B'_{\rho_2}(0)}\psi,$$

$$h_1=\min\left\{h_0, \frac{\rho_2}{2}\right\}.$$
Let us notice that by the strict convexity of $\psi$, $h_1$ is
positive. Let us choose in \eqref{2-65C}

\begin{equation}\label{h-65C}
h=h_1.
\end{equation}
Let us consider the set

$$D=\left\{(x',x_n):x'\in B'_{\rho_2}(0), \quad \psi(x')<x_n<h_1\right\}.$$
We have that $D$ has a "lens" shape in particular on the boundary of $D$ there are no vertical segments.
Let us note that, because of the way we choose $h_1$,
we have $w\in C^{\omega}(\overline{D})$.
\begin{figure}\label{figura-p66}
	\centering
	\includegraphics[trim={0 0 0 0},clip, width=9cm]{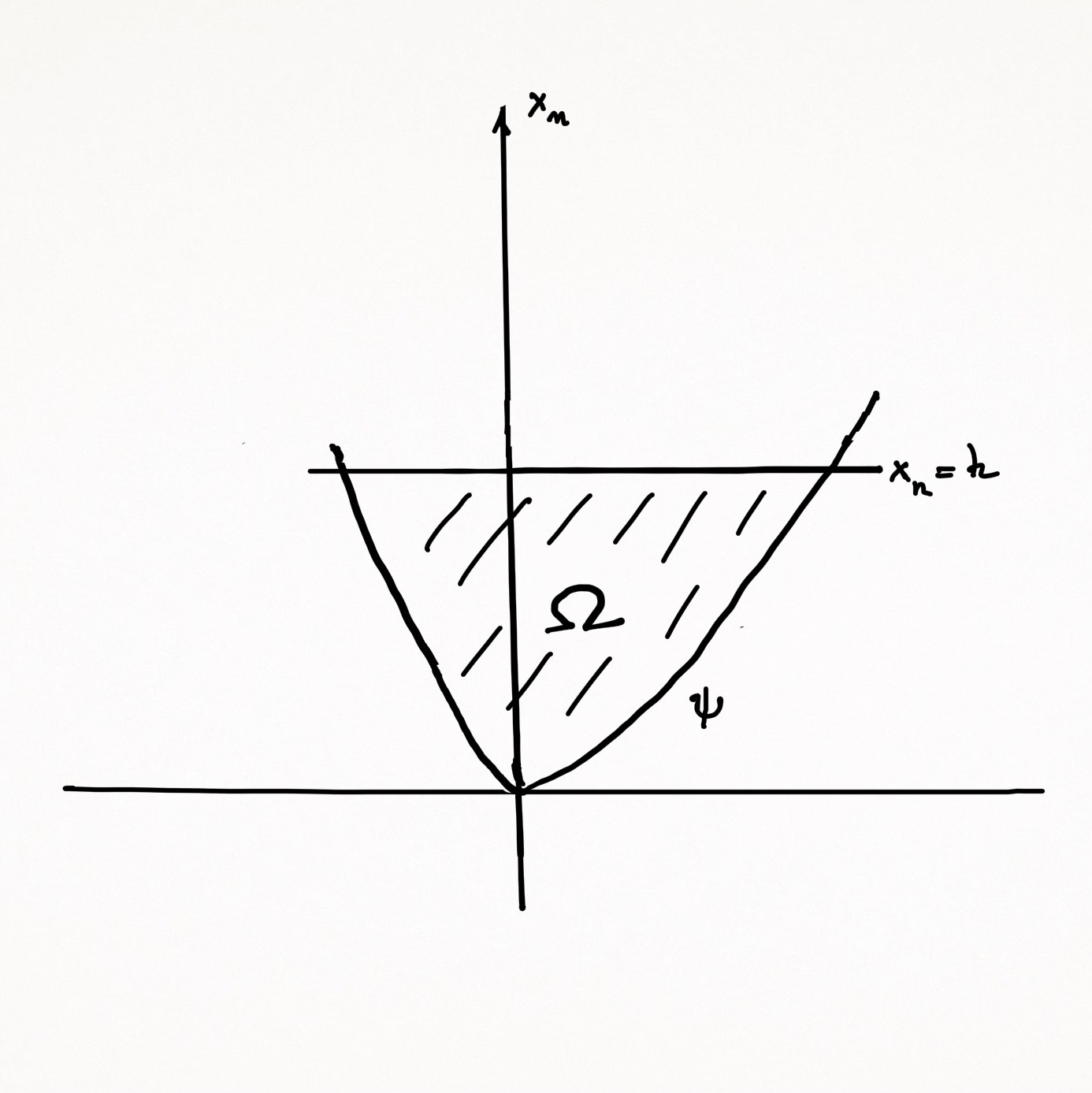}
	\caption{}
\end{figure}

Now, by the
 assumptions $u$ belongs to $C^m\left(\overline{D}\right)$ and it is
a solution to the Cauchy problem
\begin{equation}\label{Cauchy-66C}
\begin{cases}
P(x,\partial)u=0, & \mbox{in}\quad D, \\
\\
\partial^{\alpha} u=0, & \mbox{for } |\alpha|\leq m-1,\quad x\in \Gamma\cap
D.
\end{cases}
\end{equation}
By the Green identity  \eqref{1-61C}, we have

\begin{equation}\label{1-66C}
\begin{aligned}
\int_{D}fudx&=\int_{D}uP^{\ast}(x,\partial)w dx=\\&
=\int_{D}\left(uP^{\ast}(x,\partial)w-wP(x,\partial)u\right)dx=\\&
=\int_{\partial D} \mathcal{M}(u,w;\nu) dS=0.
\end{aligned}
\end{equation}
To prove that

$$\int_{\partial D} \mathcal{M}(u,w;\nu) dS=0,$$
it suffices to write the integral over $\partial D$ as a sum of two
integrals, say $I_1$ and $I_2$, with the same integrand
$\mathcal{M}(u,w;\nu)$, where $I_1$ is the integral over a portion of
the graph of $\psi$, on which $\partial^{\alpha}u=0$, for
$|\alpha|\leq m-1$, and $I_2$ is
 the integral over a portion of hyperplane $\{x_n=h_1\}$ on which, by \eqref{2-65C} and \eqref{h-65C}, we have $\partial^{\alpha}w=0$,
 for $|\alpha|\leq m-1$. Hence, by \eqref{M-61C}, both $I_1$ and $I_2$ is equal to zero.
Therefore, by \eqref{1-66C} we get

$$\int_{D}fudx=0,\quad\mbox{ for every polynomial } f$$
and since the set of polynomials is dense in
$C^0\left(\overline{D}\right)$ (Theorem \ref{teor-approssimaz-weierst}) we have
$$u\equiv 0, \quad\mbox{in } D.$$
The first part of proof is concluded.

\medskip

\textbf{Step II.} Now, we suppose that $\Gamma$ satisfies
\eqref{2-63C} e \eqref{1-63C}, but we \textit{do not suppose} that $\psi$
is strictly convex.  We may reduce to the case discussed
in Step t I using the following \textbf{Holmgren transformation} \index{Holmgren transformation}
\begin{equation}\label{3-4-67C}
\Lambda:\mathbb{R}^n_x\rightarrow \mathbb{R}^n_y,\quad x\rightarrow
y=\Lambda(x',x_n)=\left(x',x_n+\frac{A}{2}|x'|^2\right),
\end{equation}
where $A>0$ is to be chosen. Let us note that $\Lambda$ is a
diffeomorphism. Let it be further

\begin{equation}\label{5-67C}
\widetilde{\Gamma}:=\Lambda(\Gamma)=\left\{(x',\widetilde{\psi}(x')):x'\in
B'_{r_0}\right\},
\end{equation}
where
\begin{equation}\label{6-67C}
\widetilde{\psi}:=\psi(x')+\frac{A}{2}|x'|^2.
\end{equation}
Let us choose $A$ in such a way that $\widetilde{\psi}$ is strictly
convex. For this purpose it suffices to have
$$A>\left\Vert\partial^2\psi\right\Vert_{L^{\infty}(B'_{r_0})},$$
where $\partial^2\psi$ is the Hessian matrix of $\psi$. Fix such a number $A$. Let us denote by $\widetilde{P}(y,\partial_y)$
the transformed operator by means of $\Lambda$

$$\widetilde{P}(y,\partial_y)v(y)_{|y=\Lambda(x)}=P(x,\partial_x)u(x),$$
where $v$ is defined by
$$v(\Lambda(x))=u(x).$$
Let us notice that the coefficients of  $\widetilde{P}$ are analytic functions. Setting
$\widetilde{\mathcal{U}}_+=\Lambda(\mathcal{U}_+)$, we have that
$v\in C^m\left(\overline{\widetilde{\mathcal{U}}_+}\right)$ is a
solution to the Cauchy problem

\begin{equation*}
\begin{cases}
\widetilde{P}(y,\partial_y)v=0, & \mbox{in}\quad\mathcal{U}_+, \\
\\
\partial^{\alpha} v(y)=0, & \mbox{for } |\alpha|\leq m-1,\quad y\in \widetilde{\Gamma}\cap \mathcal{U}_+.%
\end{cases}
\end{equation*}
Let us recall (compare with \eqref{4-11C}) that the symbol of the principal part
of $\widetilde{P}_m(y,\eta)$ is given by

\begin{equation*}
\widetilde{P}_m(y,\eta)=P(x,\partial_x(\Lambda
(x))^t\eta)_{|x=\Lambda^{-1}(y)}.
\end{equation*}
Since
$$(\Lambda (0))^te_n=e_n$$
we have
$$\widetilde{P}_m(0,e_n)=P_m(0,e_n)\neq 0.$$
In short, we come back to the case already treated in Step I.
Therefore, for a suitable neighborhood $\mathcal{U}$ of $0$, we have $v\equiv 0$ in
$\mathcal{U}_+$ which implies $u\equiv 0$ in $U_+$. $\blacksquare$

\bigskip

\textbf{Remarks about the Holmgren Theorem.}

\textbf{1.} If $\Gamma$ is a noncharacteristic surface,
Theorem \ref{teor-Holm} allows us to say that there exists an open set $S$ such that $\Gamma\subset S$ and such that, denoting by
$S_+=S\cap\{x\in \Omega: \phi(x)\geq\phi(x_0)\}$, it occurs that if
$v\in C^m\left(\overline{S_+}\right)$ is a solution of the Cauchy  problem

\begin{equation}\label{0-68C}
\begin{cases}
P(x,\partial)u=0, & \mbox{in}\quad S_+, \\
\\
\partial^{\alpha} u(x)=0, & \mbox{for } |\alpha|\leq m-1,\quad x\in \Gamma,%
\end{cases}
\end{equation}
then $u\equiv 0$ in $S_+$.

Nevertheless, the statement of Holmgren Theorem does not clarify sufficiently how large
the set $S$ (or $S_+$) can be. Actually, one would expect that $\partial S\setminus \Gamma$
should consist of characteristic surfaces or, in other words, that the uniqueness for the Cauchy problem
would hold until a characteristic surface is encountered.

If $P(x,\partial)$ is an
elliptic operator with analytic coefficients in an open connected set $\Omega$ of $\mathbb{R}^n$ (Section \ref{trasformazione})
and if $\Gamma$ is a portion of a regular  surface, then, since the ellipticity of $P(x,\partial)$ guarantees us
that $\Gamma$ is not characteristic, we would expect the same ellipticity of $P(x,\partial)$ guarantees
that a solution of $P(x,\partial)u=0$ in $\Omega$, with null Cauchy data on $\Gamma$, is identically null on
$\Omega$. For instance, if $P(x,\partial)=\Delta$, the above occurs. As a matter of fact, it is enough to keep in mind that the solutions of
$\Delta u=0$ are analytic in $\Omega$ to obtain that $u\equiv 0$ in $\Omega$.

If we have the vibrating string operator
$\partial_t^2-\partial_x^2$ we know that if $u$ is a solution of
$\partial_t^2u-\partial_x^2u=0$ with zero initial conditions on $\Gamma= (-R, R)$ then $u$ vanishes in the 
square $\{|x|+|t|<R\}$ that is, $u$ vanishes in a region
bounded by characteristic lines parallel to  $\{x+t=0\}$, $\{x-t=0\}$.

Neither the situation described for the Laplace equation nor
the one described for the vibrating string operator are a
direct consequence of the statement of Theorem \ref{teor-Holm}. A
general answer to the problems is given  by the
Global Uniqueness Theorem proved by F. John, of which we will here provide
the statement and examine some of its consequences.

\medskip

\textbf{2.} The assumptions of Theorem \ref{teor-Holm} can be
weakened. Here we merely give a few hints and refer to
\cite[Theorem 5.3.1]{HO63} the interested reading in learning more about
the topic. We point out, in particular, that:

(i) one may assume $\psi\in C^1(\Omega)$.

(ii) one may give a distributional formulation of Cauchy problem \eqref{4-62C} and in this framework prove the uniqueness
of the solution.   

\medskip

\textbf{3.} It is worth mentioning that the Holmgren uniqueness Theorem cannot be extended to the nonlinear case. For further information we refer the reader to \cite{Garding}, \cite{Met}. $\blacklozenge$

\subsection{Statement of the Holmgren--John Theorem. Examples} \label{holm-globale}

Here we we only state the Holmgren--John global uniqueness Theorem, for the proof we refer to \cite{Joh}) or, in these notes, to the Chapter \ref{Holmgren-John} in which we will prove the Stability Theorem due to F. John himself and from which we can trivially deduce the uniqueness.

In order to state the Global Uniqueness Theorem, we need the
following definition (in it we follow the terminology introduced in \cite{Joh}).

\begin{definition}
	\index{Definition:@{Definition:}!- analytic field@{- analytic field}}
Let
\begin{equation}\label{0-71C}
F:B'_1\times(0,1)\rightarrow \mathbb{R}^n,
\end{equation}
a function satisfying the following properties (let us denote by
$y'\in B'_1$ and $\lambda\in (0,1)$ the independent variables):

(i) $F$ is injective,

(ii) $F$ is analytic in $B'_1\times(0,1)$,

(iii) for every $(y',\lambda)\in B'_1\times(0,1)$ the jacobian matrix $\partial_{y',\lambda}F(y',\lambda)$ is nonsingular.

We call \textbf{analytic field} in $\mathbb{R}^n$ the family of sets
$\left\{S_{\lambda}\right\}_{\lambda\in (0,1)}$, where

\begin{equation}\label{0-72C}
S_{\lambda}=\left\{F(y',\lambda): y'\in (0,1)\right\}, \quad\mbox{
for } \lambda\in (0,1).
\end{equation}
We call  \textbf{support of the analytic field} the open set (see Figure 7.6)

$$\Sigma:=\bigcup_{\lambda\in (0,1)} S_{\lambda}.$$
\end{definition}

\begin{figure}\label{figura-p72}
	\centering
	\includegraphics[trim={0 0 0 0},clip, width=12cm]{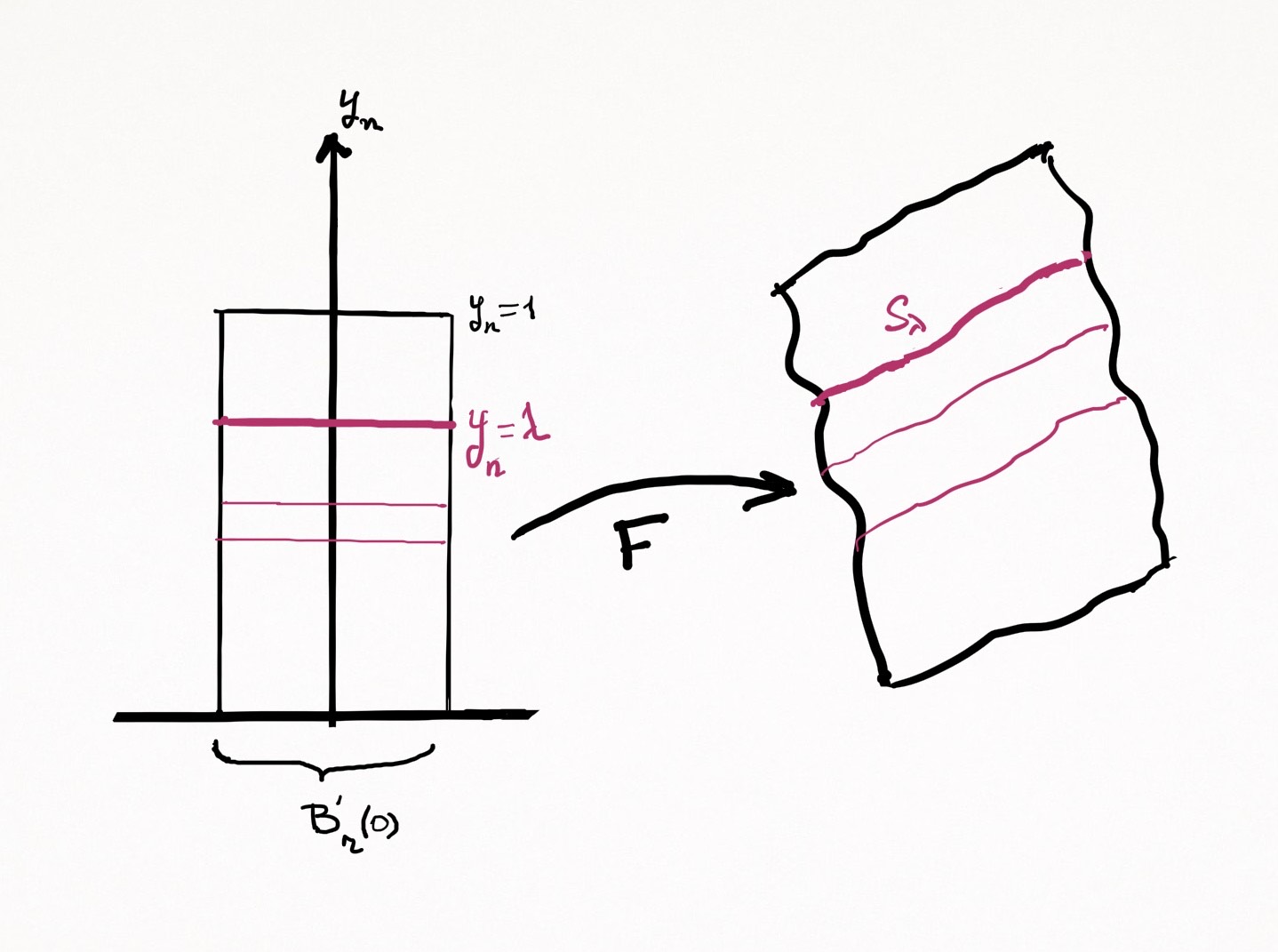}
	\caption{}
\end{figure}
\noindent We denote, for any $\mu \in (0,1)$,

$$\Sigma_{\mu}:=\bigcup_{\lambda\in (0,\mu]} S_{\lambda}.$$
\begin{figure}\label{figura-p73}
	\centering
	\includegraphics[trim={0 0 0 0},clip, width=9cm]{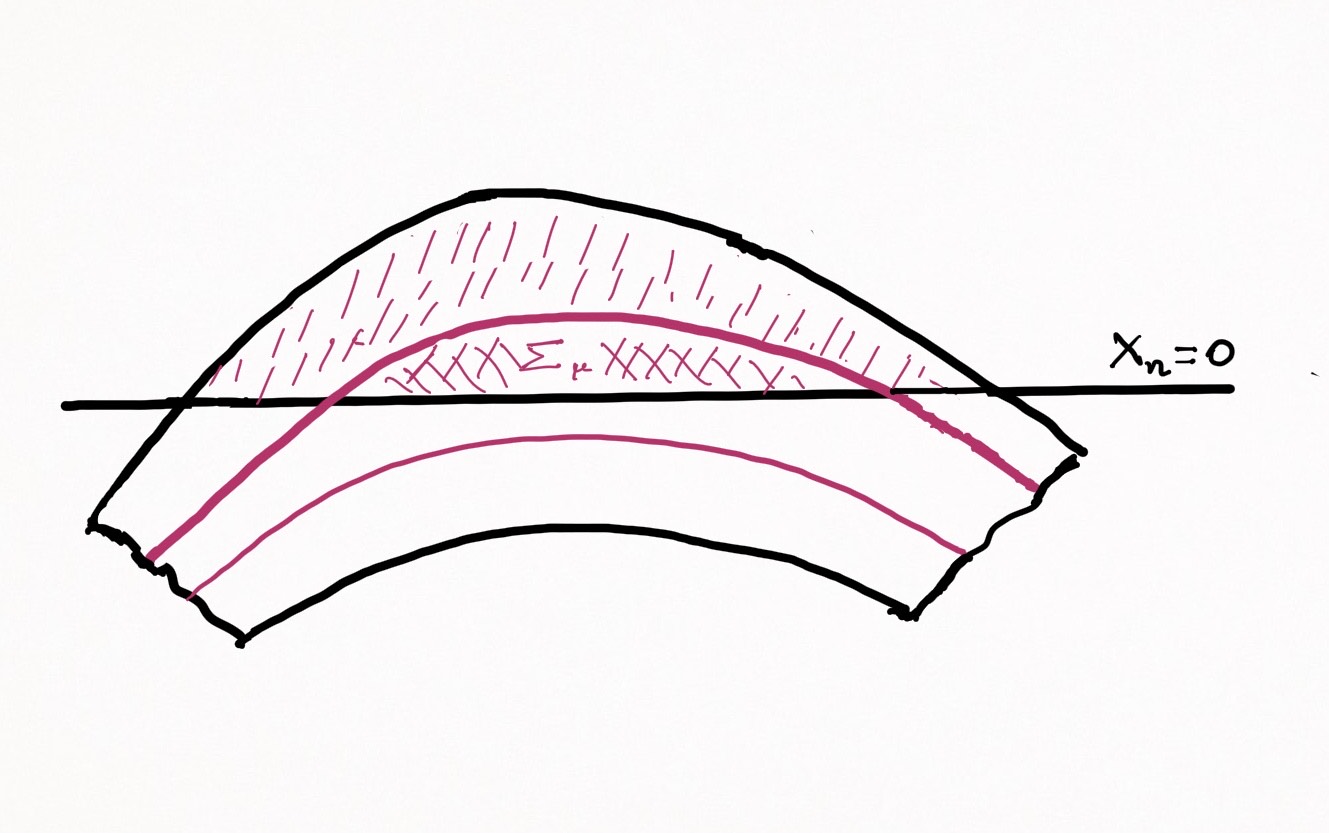}
	\caption{}
\end{figure}

\begin{theo}[\textbf{Holmgren--John}]\label{teor-Jo-Holm}
	\index{Theorem:@{Theorem:}!- Holmgren--John@{- Holmgren--John}}
Let $S_{\lambda}$ be an analytic field  in $\mathbb{R}^n$ and let
$\Sigma$ be its support. Let
$$P(x,\partial)=\sum_{|\alpha|\leq m}a_{\alpha}(x)\partial^{\alpha}$$
be a linear differential operator of order $m$, where $a_{\alpha}\in
C^{\omega}(\Sigma)$. Let us define the sets  (see Figure 7.7)
\begin{equation}\label{1-73C}
\mathcal{R}=\left\{x\in \mathbb{R}^n: x\in \Sigma, \quad x_n\geq
0\right\},
\end{equation}

\begin{equation}\label{2-73C}
\mathcal{Z}=\left\{(x',0): (x',0)\in \Sigma\right\}.
\end{equation}
Let us suppose

\medskip

(a) $\mathcal{Z}$ and $S_{\lambda}$, $\lambda\in (0,1)$, are noncarhacteristic for $P(x,\partial)$,

(b) for every $\mu\in (0,1)$, $\mathcal{R}\cap \Sigma_{\mu}$, is a
\textbf{closed set} of $\mathbb{R}^n$.

\medskip

Then we have that if $u\in C^m(\mathcal{R})$ is a solution to the Cauchy problem

\begin{equation}\label{3-73C}
\begin{cases}
P(x,\partial)u=0, & \mbox{in}\quad \mathcal{R}, \\
\\
\partial^{\alpha} u=0, & \mbox{for } |\alpha|\leq m-1,\quad x\in \mathcal{Z},%
\end{cases}
\end{equation}
we have $$u\equiv 0.$$
\end{theo}

Let us illustrate somewhat the assumptions of Theorem \ref{teor-Jo-Holm}.
We observe that $\mathcal{Z}$ is a portion of the hyperplane
$\left\{x_n=0\right\}$: this is exclusively a expository choice, actually $\mathcal{Z}$ can be any
$C^2$ noncharacteristic surface for $P(x,\partial)$. Furthermore,
hypothesis (b) assures us that in the boundary of $\mathcal{R}\cap
\Sigma_{\mu}$, for $\mu\in (0,1)$, there are no vertical segments,
in other words, $\mathcal{R}\cap \Sigma_{\mu}$ has a
"lens" shape which we have already encountered in the proof of the
Theorem \ref{teor-Holm}.

\bigskip

\textbf{Example 1 -- Wave equation.}

Let us denote by $x\in \mathbb{R}^n$ and $t\in \mathbb{R}$ the independent variables, let

\begin{equation}\label{1-77C}
K=\left\{(x,t)\in \mathbb{R}^{n+1}: |x|<|1-t|\right\}.
\end{equation}
Let us prove that if $u\in C^2(K)$ is a solution to the Cauchy problem

\begin{equation}\label{1n-77C}
\begin{cases}
\partial^2_tu-\Delta_x u=0, & \mbox{in}\quad K, \\
\\
u(x,0)=0, & \mbox{for }  |x|<1,\\
\\
\partial_tu(x,0)=0, & \mbox{for }  |x|<1,
\end{cases}
\end{equation}
then $u=0$ in $K$.

We apply Theorem \ref{teor-Jo-Holm}

In this case we have
$$P(\partial_t,\partial_x)=\partial^2_tu-\Delta_x.$$
Let $\varepsilon \in (0,1)$ be fixed and let, for $\lambda \in (0,1)$ (Figure 7.8)

$$S^{\varepsilon}_{\lambda}=\left\{(x,t)\in \mathbb{R}^{n+1}:
t=1-\sqrt{(1-\lambda+\varepsilon)^2+|x|^2},\quad|x|<\sqrt{1-\varepsilon^2}\right\}.$$

\begin{figure}\label{figura-p77}
	\centering
	\includegraphics[trim={0 0 0 0},clip, width=9cm]{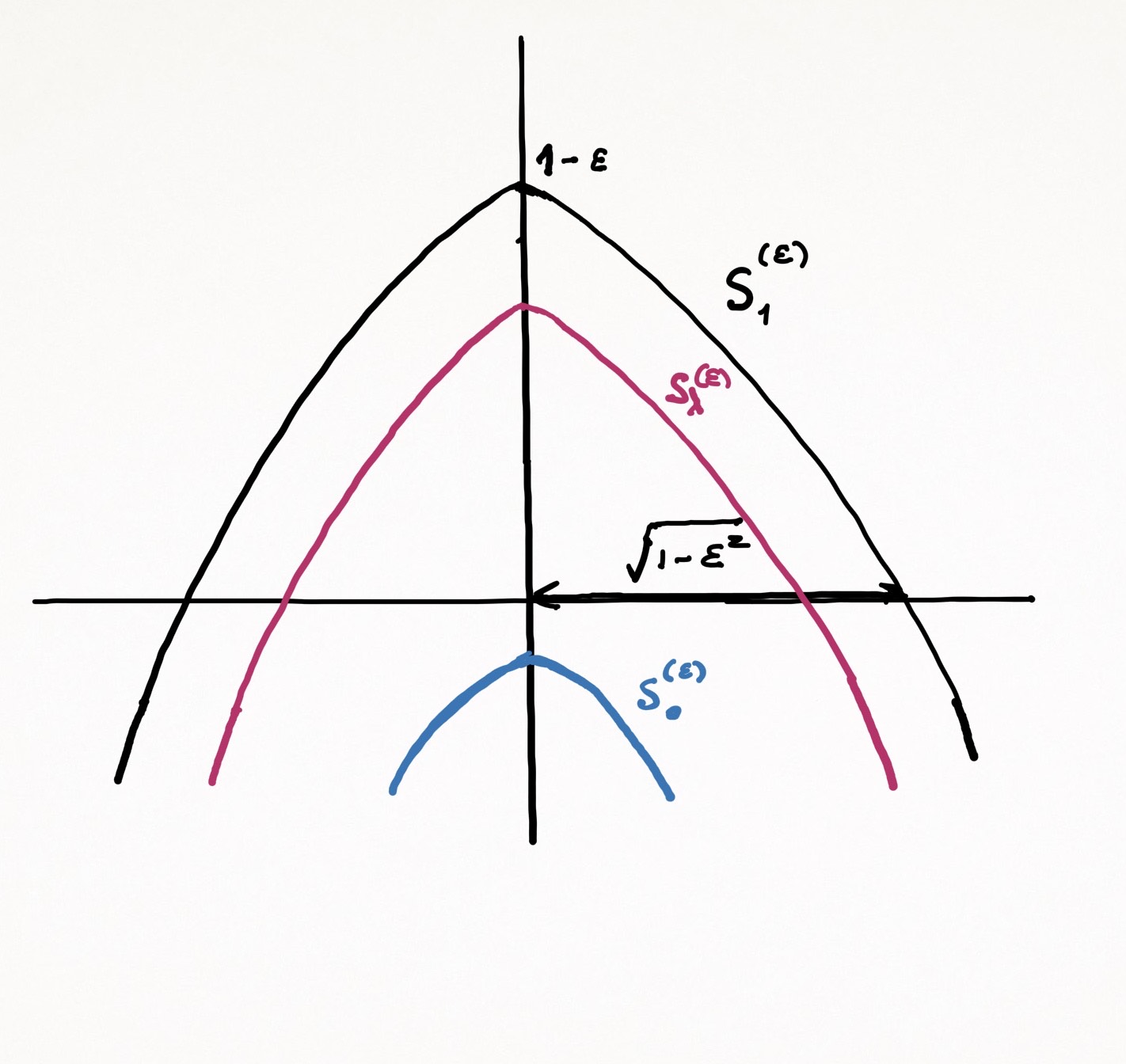}
	\caption{}
\end{figure}

\noindent Let

\begin{equation*}
F:B_{\sqrt{1-\varepsilon^2}}(0)\times(0,1)\rightarrow
\mathbb{R}^{n+1},
\end{equation*}

\begin{equation*}
F(y,\lambda)=\left(y, 1-\sqrt{(1-\lambda+\varepsilon)^2+|y|^2}\right).
\end{equation*}

Moreover, set

\begin{equation*}
\begin{aligned}
&\Sigma^{\varepsilon}  =\bigcup_{\lambda\in (0,1)}
S^{\varepsilon}_{\lambda}=\\& =\left\{
1-\sqrt{(1-\lambda+\varepsilon)^2+|x|^2}<t<1-\sqrt{\varepsilon^2+|x|^2},\quad|x|<\sqrt{1-\varepsilon^2}\right\},
\end{aligned}
\end{equation*}

\begin{equation*}
\mathcal{R}^{\varepsilon}=\left\{(x,t)\in \mathbb{R}^{n+1}: 0\leq
t<1-\sqrt{\varepsilon^2+|x|^2}\right\},
\end{equation*}

\begin{equation*}
\mathcal{Z}^{\varepsilon}=\left\{(x,0):|x|<\sqrt{1-\varepsilon^2}\right\}
\end{equation*}
and
\begin{equation*}
\Sigma^{\varepsilon}_{\mu}=\left\{(x,t)\in \mathbb{R}^{n+1}: 0\leq
t\leq1-\sqrt{(1-\mu+\varepsilon)^2+|x|^2}\right\}.
\end{equation*}
For each $0<\mu<\varepsilon$ we have
$\mathcal{R}^{\varepsilon}\cap\Sigma^{\varepsilon}_{\mu}=\emptyset$
and, for each $\varepsilon\leq \mu<1$, we have
\begin{equation*}
\Sigma^{\varepsilon}_{\mu}=\mathcal{R}^{\varepsilon}\cap\Sigma^{\varepsilon}_{\mu}.
\end{equation*}

Finally, let us check that $\mathcal{Z}$ and $S^{\varepsilon}_{\lambda}$
are noncharacteristic. We have trivially that $\mathcal{Z}$ is a noncharacteristic surface.  Concerning
$S^{\varepsilon}_{\lambda}$, we notice that

\begin{equation*}
S^{\varepsilon}_{\lambda}=\left\{(x,t)\in \mathbb{R}^{n+1}:
\phi(x,t)=(1-\lambda+\varepsilon)^2, \quad t<1\right\},
\end{equation*}
where
\begin{equation*}
\phi(x,t)=(1-t)^2-|x|^2.
\end{equation*}
Therefore we have

$$\nabla_{x,t}\phi(x,t)=(-2x,2(t-1)),\quad P(\nabla_{x,t}\phi(x,t))=4\left((1-t)^2-|x|^2\right).$$
Hence, if $(x,t)\in S^{\varepsilon}_{\lambda}$, then

$$P(\nabla_{x,t}\phi(x,t))=4(1-\lambda+\varepsilon)^2>0.$$
Therefore $S^{\varepsilon}_{\lambda}$ is noncharacteristic.

By Theorem \ref{teor-Jo-Holm} we get

$$u=0,\quad\mbox{ in } \mathcal{R}^{\varepsilon}$$
and, since $\varepsilon$ is arbitrary, we have

$$u=0,\quad\mbox{ for } 0\leq t<1-|x|.$$
Similarly we can check that $u=0$ for  $0\leq |x|-1<t\leq 0$.
Therefore $u=0$ in $K$. $\spadesuit$

\bigskip

\textbf{Example 2 -- Elliptic equations with analytic coefficients.}
\begin{figure}\label{figura-p83}
	\centering
	\includegraphics[trim={0 0 0 0},clip, width=11cm]{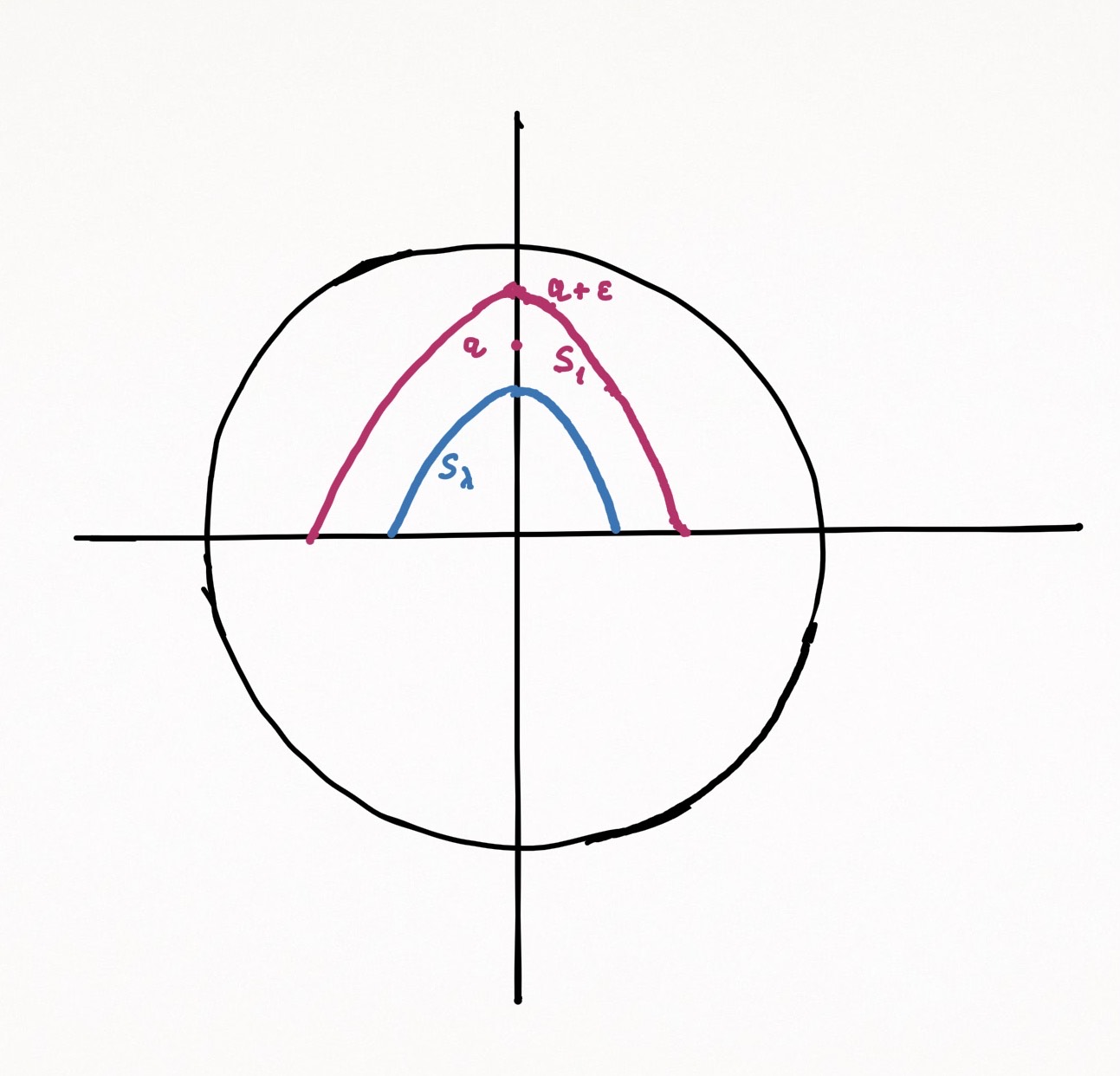}
	\caption{}
\end{figure}

Let $$P(x,\partial)$$ be a
\textbf{linear elliptic operator} of order $m$ whose coefficients are analytic functions.
Let us begin by proving the following unique continuation property \index{unique continuation property} .
\begin{prop}\label{prop-UC-82C}
Let $\rho, R$ be such that  $0<\rho<R$. Let $u\in
C^m\left(B_R\right)$ and let us suppose that 
\begin{equation}\label{3-82C}
\begin{cases}
P(x,\partial)u=0, & \mbox{in}\quad B_R, \\
\\
u=0, & \mbox{in } B_{\rho},%
\end{cases}
\end{equation}
then
$$u\equiv 0,\quad\mbox{in } B_R.$$
\end{prop}

\textbf{Proof.} Let $x_0\in B_R\setminus \overline{B_{\rho}}$ and let us prove that $u(x_0)=0$. 
We may assume that $x_0$ lies on the $x_n$--axis, because, by means of a rotation of $\mathbb{R}^n$, we may always reduce to this case (Figure 7.9). Hence let us suppose that
$$x_0=ae_n,$$ where $\rho<a<R$.

Trivially we have

\begin{equation}\label{1-83C-ccc}
\partial^{\alpha}u(x',0)=0, \quad |\alpha|\leq m,\quad x'\in B'_{\rho}.
\end{equation}
Let $\varepsilon$ be such that $0<\varepsilon<\min\{R-a,\rho\}$ and let
$S_{\lambda}$ be defined by
$$x_n=(a+\varepsilon)^2\left[\lambda-\frac{|x'|^2}{(\rho-\varepsilon)^2}\right],\quad x'\in B'_{\rho-\varepsilon},\quad 0<\lambda<1.$$
Recalling that an elliptic operator has no characteristics, by
Theorem \ref{teor-Jo-Holm} we have $u(ae_n)=0$. $\blacksquare$

\bigskip

It is evident from the proof of the above Proposition that if $u\in
C^m\left(B_R\right)$, is a solution to $P(x,\partial)u=0$ in
$B_{R}$ and it is null in $B_{\rho}(\widetilde{x})\subset B_{R}$
we have that $u\equiv 0$ in $B_{R}$.

\begin{theo}\label{UCP-ell-anal}
Let $\Omega$ be a connected open set of $\mathbb{R}^n$ and let
$$P(x,\partial)=\sum_{|\alpha|\leq m}a_{\alpha}(x)\partial^{\alpha},$$
be a linear differential operator of order $m$, elliptic in
$\Omega$, where $$a_{\alpha}\in C^{\omega}(\Omega\cup \Gamma), $$ and
$\Gamma\subset\partial\Omega$. Let us assume that $\Gamma$ is a local
graph of class $C^2$. \\ If $u\in C^m(\Omega\cup\Gamma)$
is a solution to the Cauchy problem
\begin{equation}\label{2-82C}
\begin{cases}
P(x,\partial)u=0, & \mbox{in }\quad \Omega, \\
\\
\partial^{\alpha}u=0, & |\alpha|\leq m-1 \mbox{ on } \Gamma,%
\end{cases}
\end{equation}
then
$$u=0,\quad\mbox{in } \Omega.$$
\end{theo}
\textbf{Proof.} By Theorem \ref{teor-Holm} (or \ref{teor-Jo-Holm}) we that have
there exists an open set (in the topology induced on $\Omega\cup\Gamma$) $\mathcal{U}_+\subset \Omega$ such that $\partial\mathcal{U}_+\cap \Gamma \neq \emptyset$ and such that 

\begin{equation}\label{1-84C}
u=0,\quad\mbox{in } \mathcal{U}_+.
\end{equation}
Set

\begin{equation}\label{1n-84C}
A=\left\{x\in \Omega: \exists \rho_x>0 \mbox{ such that } u=0, \mbox{
in }B_{\rho_x}(x)\right\}.
\end{equation}
By \eqref{1-84C} we have $A\neq\emptyset$ and, trivially,
we have that $A$ is an open set in $\Omega$. In order to prove the assertion, it suffices to prove
that $A$ is also closed in  $\Omega$  and,
since $\Omega$ is connected, we have $A=\Omega$ from which we will obtain the thesis. 
In order to prove that $A$ is closed in $\Omega$
it suffices to prove that if $\left\{x_j\right\}_{j\in \mathbb{N}}$ is
a sequence of $A$ such that

\begin{equation}\label{1nn-84C}
\lim_{j\rightarrow\infty}x_j=x_0,
\end{equation}
where $x_0\in \Omega$, then $x_0\in A$. 

Let $\varepsilon>0$ satisfy 
$B_{\varepsilon}(x_0)\subset \Omega$. By \eqref{1nn-84C}, 
there exists $j_0$ such that
$$\left|x_0-x_{j_0}\right|<\frac{\varepsilon}{4}.$$
Since  $x_{j_0}\in A$ there exists $\rho_{x_{j_0}}>0$ such that

\begin{equation}\label{1nnn-84C}
u=0, \quad\mbox{in } B_{\rho_{x_{j_0}}}(x_{j_0})\subset \Omega.
\end{equation}
Set
$\overline{\rho}=\min\left\{\frac{\varepsilon}{4},\rho_{x_{j_0}}\right\}$,
we have $B_{\overline{\rho}}(x_{j_0})\subset B_{\varepsilon}(x_0)$ and, by
\eqref{1nnn-84C}, we have $u=0$ in $B_{\overline{\rho}}(x_{j_0})$,  then, by
Proposition \ref{prop-UC-82C} and subsequent remarks, we have
$u=0$ in $B_{\varepsilon}(x_0)$. Therefore $x_0\in A$ and the thesis follows. $\blacksquare$

\bigskip

\textbf{Remarks on Theorem \ref{UCP-ell-anal}.}

\textbf{1.} At this point we report, for information only, that if
$P(x,\partial)$ is an eliptic operator in $\Omega$, $f\in C^{\omega}(\Omega)$ and $u\in C^m(\Omega)$ is a solution
of the equation $P(x,\partial)u=f$ in $\Omega$, then $u\in
C^{\omega}(\Omega)$, see \cite{Mo-Ni} for a proof. It is evident that if we had
used this regularity property of the solutions,
Theorem \ref{UCP-ell-anal} would be a consequence of the
unique continuation property for analytic functions (and the Holmgren Theorem).
However, to prove Theorem \ref{UCP-ell-anal} we \textit{ did not} need the above regularity result.

\bigskip

\textbf{2.} Let us consider some interesting consequences of Theorem
\ref{UCP-ell-anal}. Let $\Omega$, $\Gamma$ and $P(x,\partial)$ be as
in Theorem \ref{UCP-ell-anal}, further let us suppose that $\Gamma$ is of class $C^m$.
We denote by $\mathcal{X}_{\Gamma}$ the class of functions
$$g:\Gamma\rightarrow \mathbb{R}^m, \ \ g(x)=(g_0(x),g_1(x),\cdots,
g_{m-1}(x)), \ \ \forall x\in \Gamma$$ such that there exists $u\in
C^m(\Omega\cup\Gamma)$ solution to the Cauchy problem

\begin{equation}\label{2-86C}
\begin{cases}
P(x,\partial)u=0, & \mbox{in}\quad \Omega, \\
\\
\frac{\partial^{j}u}{\partial\nu^j}=g_j, & j=0,1,\cdots,m-1 \mbox{
on } \Gamma.
\end{cases}
\end{equation}
The class $\mathcal{X}_{\Gamma}$ enjoys the following
property.

For every $\Gamma_0\subset \Gamma$, $\Gamma_0$ open in $\Gamma$
in the induced topology, it occurs that

\begin{equation}\label{3-86C}
g\in \mathcal{X}_{\Gamma},\ \mbox{and} \ g=0 \ \mbox{on} \
\Gamma_0\quad\Longrightarrow\quad g=0 \ \mbox{on}\ \Gamma
\end{equation}
and, therefore, by linearity, if $g, \widetilde{g}\in
\mathcal{X}_{\Gamma}$ and $g=\widetilde{g}$ on $\Gamma_0$ then
$g=\widetilde{g}$ on $\Gamma$. In other words, the class
$\mathcal{X}_{\Gamma}$ must enjoy the unique continuation property \eqref{3-86C}. By this same fact we
deduce that if the initial data of a Cauchy problem for $P(x,\partial)$ belong to the class of functions $C^k(\Gamma,\mathbb{R}^m)$, for
any $k\geq m$, such a Cauchy problem cannot, in general, admits solutions. As a matter of fact, the class $C^k(\Gamma,\mathbb{R}^m)$  does not  enjoy
property \eqref{3-86C}.

The proof of the assertion above is very simple. As a matter of fact, let
$g\in \mathcal{X}_{\Gamma}$ and let $u\in C^k(\Omega\cup\Gamma)$, $k\geq m$ be a solution
of problem \eqref{3-86C}. Let us suppose
$$g=0,\quad \mbox{ on }\quad \Gamma_0.$$
Then applying Theorem \ref{UCP-ell-anal} (with $\Gamma_0$ in the
place of $\Gamma$) we have that $u=0$ in $\Omega$ and being $u\in
C^k(\Omega\cup\Gamma)$, we have $g=u_{|\Gamma}=0$.
 $\blacklozenge$

\bigskip

\textbf{Example 3 -- One--dimensional heat equation.}

We consider the following Cauchy problem

\begin{equation}\label{1-88C}
\begin{cases}
\partial_x^2u-\partial_tu=0, & \mbox{in}\quad D:=(0,1)\times(0,1), \\
\\
u(0,t)=0,
 &  \mbox{ for } t\in (a,b),\\
 \\
 \partial_xu(0,t)=0, &  \mbox{ for } t\in (a,b),\\
\end{cases}
\end{equation}
where $a,b$ are given numbers and such that $0<a<b<1$. Let us notice that problem \eqref{1-88C} is a \textbf{noncharacteristic Cauchy problem}, since the initial line is $\{x=0\}$.

By using the same arguments exploited in \textbf{Example 1} and
in \textbf{Example 2} (Proposition \ref{prop-UC-82C}) it can be proved
easily that if $u\in C^2(\overline{D})$ then $u=0$ in
$[0,1]\times [a,b]$. The details are left to the reader (it is
useful to keep in mind that the only characteristics of the operator
$\partial_x^2-\partial_t$ are the straight lines $t=t_0$, for any $t_0\in \mathbb{R}$).

It is quite natural to wonder whether a solution
of \eqref{1-88C} is null in $D$. The answer to this
question is negative as proved in an example due to
\textbf{Tychonoff}, \index{Tychonoff example}\cite{Tik1}, which we will discuss below.

\medskip

First, let us consider the following Cauchy problem.

\begin{equation}\label{1-89C}
\left \{
\begin{array}{cc}
\partial_x^2u-\partial_tu=0, &  \\
\\
u(0,t)=\varphi(t),
 &  \mbox{ for } t\in \mathbb{R},\\
 \\
 \partial_xu(0,t)=0, &  \mbox{ for } t\in \mathbb{R},\\
\end{array}%
\right.
\end{equation}
For the time being, let us just assume that $\varphi\in
C^{\infty}(\mathbb{R})$ and let us search, at first just formally, a
solution of the type

\begin{equation}\label{1n-89C}
u(x,t)=\sum_{j=0}^\infty a_j(t)x^j,
\end{equation}
where $a_j$ are functions to be found.

Obviously we have to require that

$$a_0(t)=\varphi(t),\quad\mbox{ and }\quad a_1(t)=0$$
and, by requiring that \eqref{1n-89C} is a solution to the equation
$\partial_x^2u-\partial_tu=0$, we need to require 

\begin{equation*}
\sum_{j=0}^\infty a'_j(t)x^j-\sum_{j=2}^\infty
j(j-1)a_j(t)x^{j-2}=0,
\end{equation*}
from which we have

\begin{equation*}
a_{j+2}(t)=\frac{1}{(j+2)(j+1)}a'_j(t), \quad \forall j\in \mathbb{N}_0.
\end{equation*}
Hence

\begin{equation*}
a_{2k}(t)=\frac{\varphi^{(k)}(t)}{(2k)!}, \quad a_{2k+1}(t)=0\quad
\forall k\in \mathbb{N}_0.
\end{equation*}
Therefore

\begin{equation}\label{2-89C}
u(x,t)=\sum_{k=0}^\infty\frac{\varphi^{(k)}(t)}{(2k)!}x^{2k}.
\end{equation}
In order that \eqref{2-89C} is actually a solution to 
Cauchy problem \eqref{1-89C} (in a neighborhood of $\{0\}\times
\mathbb{R}$) it suffices to require that there exist two positive numbers $c$ and
$M$ such that

\begin{equation}\label{1-90C}
\left\vert\varphi^{(k)}(t)\right\vert\leq cM^k(2k)!,\quad \forall
k\in\mathbb{N}_0.
\end{equation}
Let now

\begin{equation}\label{1n-90C}
\varphi(t)=
\begin{cases}
e^{-\frac{1}{t^2}}, & \mbox{ for } t> 0, \\
\\
0,
 &  \mbox{ for } t\leq 0.\\
\end{cases}
\end{equation}
We show now that \eqref{1-90C} is satisfied. We will prove,
indeed,

\begin{equation}\label{2-90C}
\left\vert\varphi^{(k)}(t)\right\vert\leq
\left(\frac{9k}{2e}\right)^{\frac{k}{2}}k!,\quad \forall
k\in\mathbb{N}_0.
\end{equation}
which (by the Stirling  formula, \eqref{Stirling-3N}) implies \eqref{1-90C}.

To prove \eqref{2-90C} we use the Cauchy formula for the
holomorphic functions. Let therefore $t>0$ and let $S$ be the circumference centered at $t+i0$ and with radius $\frac{t}{2}$ in the complex plane, i.e.

$$S=\left\{t\left(1+\frac{1}{2}e^{i\vartheta}\right): \vartheta\in [0,2\pi)\right\}.$$
We have, for any $k\in\mathbb{N}_0$,

\begin{equation*}
\varphi^{(k)}(t)= \frac{k!}{2\pi
i}\int_S\frac{e^{-\frac{1}{z^2}}}{(z-t)^{k+1}}dz.
\end{equation*}
Now, when $z=t\left(1+\frac{1}{2}e^{i\vartheta}\right)\in S$, it is easily checked that
$$ \frac{1}{z}=\frac{4}{3t}+\frac{2}{3t}e^{i\vartheta}.$$
Therefore
\begin{equation*}
\begin{aligned}
\Re\left(\frac{1}{z^2}\right)  & =\left(\frac{4}{3t}\right)^2
\left[\left(1+\frac{1}{2}\cos\vartheta\right)^2-\left(\frac{1}{2}\sin\vartheta\right)^2\right]=\\&
=\left(\frac{4}{3t}\right)^2\left[\frac{1}{4}+\frac{1}{2}\left(1+\cos\vartheta\right)^2\right]\geq
\frac{4}{9t^2}.
\end{aligned}
\end{equation*}
Therefore we have

\begin{equation}\label{0-91C}
\begin{aligned}
\left\vert\varphi^{(k)}(t)\right\vert  & \leq
\frac{k!}{2\pi}\int_{|z-t|=
\frac{t}{2}}\left\vert\frac{e^{-\frac{1}{z^2}}}{(z-t)^{k+1}}\right\vert
ds\leq\\& \leq
\frac{k!}{2\pi}\int_{|z-t|=\frac{t}{2}}\frac{e^{-\frac{4}{9t^2}}}{|z-t|^{k+1}}ds=\\&
=
\frac{k!}{2\pi}\left(2\pi\frac{t}{2}\right)\frac{1}{(t/2)^{k+1}}e^{-\frac{4}{9t^2}=}\\&
=\frac{2^kk!}{t^k}e^{-\frac{4}{9t^2}}.
\end{aligned}
\end{equation}
Since
$$\sup \frac{1}{t^k} e^{-\frac{4}{9t^2}}=\left(\frac{9k}{8e}\right)^{k/2},$$
by \eqref{0-91C} we get \eqref{2-90C} which in turn implies that
 the series in \eqref{2-89C} converges for every $x\in\mathbb{R}$ and  its
sum, $u$, is indeed the solution to Cauchy problem \eqref{1-89C}.

From what we have just established, it turns out that, denoting
by $\psi$ the following function

\begin{equation*}
\psi(t)=
\begin{cases}
e^{-\frac{1}{(t-b)^2}}, & \mbox{for } t>b,  \\
\\
0,
 &  \mbox{ per } a\leq t\leq b,\\
 \\
 e^{-\frac{1}{(a-t)^2}}, &  \mbox{ for } t<a,\\
\end{cases}
\end{equation*}
we have that

\begin{equation*}
\widetilde{u}=\sum_{k=0}^\infty\frac{\psi^{(k)}(t)}{(2k)!}x^{2k},
\end{equation*}
is the solution to the Cauchy problem

\begin{equation*}
\left \{
\begin{array}{cc}
\partial_x^2\widetilde{u}-\partial_t\widetilde{u}=0, & \mbox{in } \mathbb{R}^2, \\
\\
\widetilde{u}(0,t)=\psi(t),
 &  \mbox{ for } t\in \mathbb{R},\\
 \\
 \partial_x\widetilde{u}(0,t)=0, &  \mbox{ for } t\in \mathbb{R}\\
\end{array}%
\right.
\end{equation*}
and $\widetilde{u}=0$ in $\mathbb{R}\times[a,b]$, but if $t_0\notin
[a,b]$ then $\widetilde{u}(\cdot,t_0)$ does not  identically
vanish, more precisely $\widetilde{u}(\cdot,t_0)$ does not
vanish in any open set of $\mathbb{R}$ because, as
$\widetilde{u}(\cdot,t_0)$ is analytic (as it is the sum of a
series of powers in the variable $x$) one would have that
$\psi^{(k)}(t_0)=0$ for every $k\in \mathbb{N}_0$ that is false.

We conclude this discussion about the heat equation
by observing that does not hold the uniqueness in $C^2(\mathbb{R}\times[0,+\infty))$ to the following Cauchy problem \textbf{characteristic}

\begin{equation}\label{Cauchy-caratt}
\begin{cases}
\partial_x^2u-\partial_tu=0, & \mbox{in } \mathbb{R}^2, \\
\\
u(x,0)=0,
 &  \mbox{ for } x\in \mathbb{R}.\\
\end{cases}
\end{equation}
It suffices to consider the function $u$ defined by \eqref{2-89C} where
$\varphi$ defined by \eqref{1n-90C} and we have that $u$ is the solution
of problem \eqref{Cauchy-caratt}, but it does not  vanish identically. $\spadesuit$

\bigskip

\underline{\textbf{Exsercise}}. Let $T,r>0$ and denote
$$K=\left\{(x,t)\in \mathbb{R}^{n+1}: \ |x|+|t|<T+r\right\}.$$ Prove that if $u\in C^2(K)$ satisfies 
\begin{equation*}
	\begin{cases}
		\partial^2_tu-\Delta_xu=0, \ \ \mbox{in } K \\
		\\
		u=0, \ \ \mbox{in } B_r\times (-T,T),
	\end{cases}
\end{equation*}
then
\begin{equation*}
	u=0, \ \ \mbox{in } K.
\end{equation*}
$\clubsuit$

\chapter{Uniqueness for an inverse problem}\label{PbI-froniera:22-11-22}
\section{Introduction}\label{PbI-froniera:22-11-22-1}

In this short chapter we present an inverse problem \index{inverse problem} for the Laplace equation.
The direct problem is nothing but the Dirichlet problem

\begin{equation}\label{PbI-froniera:22-11-22-2}
	\begin{cases}
		\Delta u=0, \quad \mbox{ in } \Omega, \\
		\\
		u=\varphi, \quad\mbox{ on } \partial\Omega.\\
	\end{cases}
\end{equation}  
For simplicity, we assume that $\Omega$ is a bounded and connected open of $\mathbb{R}^n$ whose boundary is of class  $C^{\infty}$ and $\varphi\in C^{\infty}(\partial\Omega)$. In Chapter \ref{Lax-Milgram} we saw (see in particular, Corollary \ref{Ell:reg-infin-bordo}) that, under these assumptions, there exists a unique $u\in C^{\infty}\left(\overline{\Omega}\right)$, which is the solution to \eqref{PbI-froniera:22-11-22-2}.

Now, let us suppose that a portion of $\partial\Omega$, which we will call $\Gamma^{(i)}$, is unknown and that we have
\begin{equation}\label{PbI-froniera:22-11-22-3}
	\varphi(x)=0, \quad \forall x\in \Gamma^{(i)}
\end{equation}
and let us suppose that we know  

\begin{equation}\label{PbI-froniera:22-11-22-4}
	\frac{	\partial u}{\partial \nu}=\psi, \quad \mbox{su}\quad  \Sigma,
\end{equation} 
where $\Sigma\subset \partial\Omega\setminus \Gamma^{(i)}$. We are interested in determining $\Gamma^{(i)}$. This is our \textbf{inverse problem}. Let us note that if we consider all the data of the problem we should write 

\begin{equation}\label{PbI-froniera:22-11-22-5}
	\begin{cases}
		\Delta u=0, \quad \mbox{ in } \Omega, \\
		\\
		u=\varphi, \quad\mbox{ on } \partial\Omega,\\
		\\
		\frac{	\partial u}{\partial \nu}=\psi, \quad \mbox{on}\quad  \Sigma,\\
	\end{cases}
\end{equation}  
which evidently is an overdetermined problem from which there is to be expected a compatibility relation between $\varphi,\psi$  and $\partial\Omega$ so that, if we assume to known $\varphi,\psi$ and $\partial\Omega\setminus \Gamma^{(i)}$, we can reasonably hope to obtain some information about $\Gamma^{(i)}$ itself. It is evident that if $\varphi\equiv 0$ on $\partial\Omega$ then we have no information on $\Gamma^{(i)}$.  We will see that this trivial case is (under precise assumptions) the only case in which $\Gamma^{(i)}$ is not uniquely determined.

Instead of considering the inverse problem as an overdetermined problem, it turns out to be more efficient to consider the inverse problem from the point of view which we now illustrate. Since the direct problem has a unique
solution $u\in C^{\infty}\left(\overline{\Omega}\right)$ for every bounded open set of class $C^{\infty}$,
and for any  $\varphi\in C^{\infty}(\partial\Omega)$ satisfying \eqref{PbI-froniera:22-11-22-3}, it turns out that the derivative $$\frac{	\partial u}{\partial \nu}_{|\Sigma}$$ is a function of $\Gamma^{(i)}$. Let us let us denote such a function by $$\mathbf{U}\left(\Gamma^{(i)}\right),$$ the inverse problem we are intersted in may be formulated as follows

\begin{equation*}
	\mybox{Determine the solution to the equation $\mathbf{U}\left(\Gamma^{(i)}\right)=\psi$.}
\end{equation*}     

\medskip

In the next Section we will specify the assumptions and we will formulate the uniqueness theorem for the inverse problem above.

\section{Statement of the uniqueness theorem for the inverse problem}\label{PbI-froniera:22-11-22-6}

Let us suppose that $\Omega$ is as above and let us suppose that $\partial\Omega$ is the union of two internally disjoint portions, $\Gamma^{(a)}$ ("accessible" portion) and $\Gamma^{(i)}$ ("inaccessible" portion). More precisely, let us suppose that:

\begin{equation}\label{PbI-froniera:22-11-22-7}
	\partial\Omega=\Gamma^{(a)}\cup \Gamma^{(i)}
\end{equation}

\begin{equation}\label{PbI-froniera:22-11-22-8}
	\Gamma^{(a)} \mbox{ and } \Gamma^{(i)} \mbox{ closed in } \partial\Omega \ \ \mbox{ and } \ \ \overline{\overset{\circ }{\Gamma^{(a)}}}=\Gamma^{(a)}, \ \ \overline{\overset{\circ }{\Gamma^{(i)}}}=\Gamma^{(i)},
\end{equation}

\smallskip
\noindent we equip $\partial\Omega$ with the topology induced by the Euclidean topology of $\mathbb{R}^n$,
\begin{equation}\label{PbI-froniera:22-11-22-9}
	\overset{\circ }{\Gamma^{(a)}}\cap \overset{\circ }{\Gamma^{(i)}}=\emptyset.
\end{equation}
Hence 

\begin{equation}\label{PbI-froniera:22-11-22-9-0}
	\Gamma^{(i)}=\partial\Omega\setminus \overset{\circ }{\Gamma^{(a)}}\ \ \mbox{ and } \ \ \Gamma^{(a)}=\partial\Omega\setminus \overset{\circ }{\Gamma^{(i)}}
\end{equation}
and
\begin{equation}\label{PbI-froniera:22-11-22-10}
	\overset{\circ }{\Gamma^{(a)}} \ \ \mbox{ \textbf{is  connected,}} 
\end{equation}
therefore, \eqref{PbI-froniera:22-11-22-8} implies $\Gamma^{(a)}$ is connected. 

\smallskip

Let $\Sigma$ be a compact subset of $\partial\Omega$ such that
\begin{equation}\label{PbI-froniera:23-11-22-1}
	\Sigma\Subset	\overset{\circ }{\Gamma^{(a)}},\quad\quad \overset{\circ }{\Sigma}\neq \emptyset. 
\end{equation}

\smallskip

Let $\varphi\in C^{\infty}(\partial\Omega)$ satisfy

\begin{equation}\label{PbI-froniera:23-11-22-2}
	\mbox{supp } \varphi \subset\overset{\circ }{\Gamma^{(a)}}. 
\end{equation}

\bigskip

We now state the Theorem

\begin{theo}[\textbf{uniqueness}]\label{PbI-froniera:23-11-22-3}
	Let $\Omega_k$, $k=1,2$, be two bounded connected open sets of $\mathbb{R}^n$ whose boundary is of class $C^{\infty}$. Let us assume that
	\begin{equation}\label{PbI-froniera:23-11-22-4}
		\partial\Omega_1=\Gamma^{(a)}\cup \Gamma^{(i)}_1 \quad \mbox{and}\quad \partial\Omega_2=\Gamma^{(a)}\cup \Gamma^{(i)}_2, 
	\end{equation}
	where $\Gamma^{(a)},\Gamma^{(i)}_k$, $k=1,2$, satisfies \eqref{PbI-froniera:22-11-22-8},  \eqref{PbI-froniera:22-11-22-9} and \eqref{PbI-froniera:22-11-22-10}. Let us assume that $\varphi_k\in C^{\infty}\left(\partial\Omega_k\right)$ do not vanish identically and satisfying \eqref{PbI-froniera:23-11-22-2}.
	Moreover, let us assume 
	
	\begin{equation}\label{PbI-froniera:23-11-22-5}
		\varphi_1=\varphi_2,\quad \mbox{on}\quad \Gamma^{(a)}.
	\end{equation}
	Let $u_k\in C^{\infty}\left(\overline{\Omega_k}\right)$, $k=1,2$, be the solutions to
	
	\begin{equation}\label{PbI-froniera:23-11-22-6}
		\begin{cases}
			\Delta u_k=0, \quad \mbox{ in } \Omega_k, \\
			\\
			u_k=\varphi, \quad\mbox{ on } \partial\Omega_k.\\
		\end{cases}
	\end{equation}  
	Let us assume that
	
	\begin{equation}\label{PbI-froniera:23-11-22-7}
		\frac{\partial u_1}{\partial \nu}=\frac{\partial u_2}{\partial \nu}, \quad\mbox{on}  \ \ \Sigma,
	\end{equation}
	where $\Sigma$ is a compact subset of $\partial\Omega$ which satisfies \eqref{PbI-froniera:23-11-22-1}.
	
	Then we have
	\begin{equation}\label{PbI-froniera:23-11-22-8}
		\Omega_1=\Omega_2	,
	\end{equation}
	and, consequently,
	\begin{equation*}
		\Gamma^{(i)}_1=\Gamma^{(i)}_2.
	\end{equation*} 
\end{theo}

The proof will be given in the next Section.

\section{Proof of the uniqueness}\label{PbI-froniera:23-11-22-9}
The idea of the proof of uniqueness Theorem is quite simple, but it requires propositions of general topology that we will prove separately for the purpose of not breaking the main argument of the proof. 

We argue by contradiction.  We assume that one of the two sets $\Omega_1\setminus\overline{\Omega}_2$, $\Omega_2\setminus\overline{\Omega}_1$, is not empty. For instance, let us assume  

\begin{equation}\label{PbI-froniera:23-11-22-10}
	\Omega_1\setminus\overline{\Omega}_2\neq \emptyset.
\end{equation}
Define $G$ as 
\begin{equation}\label{PbI-froniera:23-11-22-11}
	G=\bigcup_{A\in \mathcal{A}}A,
\end{equation}
where 
\begin{equation}\label{PbI-froniera:23-11-22-12}
	\mathcal{A}=\left\{A \mbox{ open set of } \mathbb{R}^n:\mbox{ } A\subset \Omega_1\cap\Omega_2, \  \Gamma^{(a)}\subset \overline{A},  \ A \mbox{ connected }\right\}.
\end{equation}

\smallskip

As we will prove in Proposition \ref{PbI-froniera:23-11-22-16} (and as can be expected), $G$ is a connected open set, moreover $\Gamma^{(a)}\subset \overline{G}$ and $G\subset \Omega_1\cap \Omega_2$. From the latter and from \eqref{PbI-froniera:23-11-22-10} we have
\begin{equation}\label{PbI-froniera:23-11-22-12-0}
	\Omega_1\setminus\overline{G}\neq \emptyset.
\end{equation}
Set 
$$u:=u_1-u_2,\quad \mbox{in}\quad G.$$
We have

\begin{equation*}
	\begin{cases}
		\Delta u=0, \quad \mbox{ in } G, \\
		\\
		u=0, \quad\mbox{ on } \Sigma,\\
		\\
		\frac{	\partial u}{\partial \nu}=0, \quad \mbox{on}\quad  \Sigma.\\
	\end{cases}
\end{equation*}  
Theorem \ref{UCP-ell-anal} implies
\begin{equation}\label{PbI-froniera:23-11-22-13}
	u\equiv 0,\quad \mbox{in}\quad G
\end{equation}

\medskip

\textbf{Claim.}
\begin{equation}\label{PbI-froniera:23-11-22-14}
	u_1= 0,\quad \mbox{su}\quad \partial\left(\Omega_1\setminus \overline{G}\right).
\end{equation}

\medskip

\textbf{Proof of Claim.} We exploit the following relationship (proved in Proposition \ref{PbI-froniera:25-11-22-3})
\begin{equation}\label{PbI-froniera:23-11-22-15}
	\partial\left(\Omega_1\setminus \overline{G}\right)\subset \Gamma_1^{(i)}\cup \left(\Gamma_2^{(i)}\cap \partial G\right).
\end{equation}

Let   $x\in \partial\left(\Omega_1\setminus \overline{G}\right)$, then by \eqref{PbI-froniera:23-11-22-15} we distinguish two cases: 

\smallskip 

(a) $x\in \Gamma_1^{(i)}$,

\smallskip

(b) $x\in \Gamma_2^{(i)}\cap \partial G$.

\smallskip

In case (a) $u_1(x)=0$ because $u=\varphi_1$ on $\partial\Omega_1$ and $\varphi_1=0$ on $\Gamma_1^{(i)}$.

In case (b), we have $u_2(x)=0$ because $u=\varphi_2$ on $\partial\Omega_2$ and $\varphi_2=0$ on $\Gamma_2^{(i)}$. On the other hand, by \eqref{PbI-froniera:23-11-22-13} and taking into account that $u$ is continuous in $\overline{G}$, we have 

$$u_1(x)=u_1(x)-u_2(x)=u(x)=0.$$
The Claim is proved.

\medskip

Therefore $u_1$ solves the Dirichlet problem

\begin{equation*}
	\begin{cases}
		\Delta u_1=0, \quad \mbox{ in } \Omega_1\setminus \overline{G}, \\
		\\
		u_1=0, \quad\mbox{ on } \partial\left(\Omega_1\setminus \overline{G}\right),\\
	\end{cases}
\end{equation*}  
and the maximum principle implies

$$u_1=0, \quad\mbox{in}\quad \Omega_1\setminus \overline{G}.$$
Now, taking into account \eqref{PbI-froniera:23-11-22-12-0}, the unique continuation property gives $$u_1=0, \quad\mbox{in}\quad \Omega_1.$$
From which we have
$$\varphi_1=0, \quad\mbox{in}\quad \partial\Omega_1,$$
But this contradicts the assumption that $\varphi_1$ does not vanish identically. Therefore  $\Omega_1\setminus\overline{\Omega}_2\neq \emptyset$. Hence $\Omega_1\subset\overline{\Omega}_2$. Similarly we have  
$\Omega_2\subset\overline{\Omega}_1$. Therefore, $\overline{\Omega}_1=\overline{\Omega}_2$. On the other hand (see Exercise of Section \ref{DescrBordo}) $$\Omega_1=\overset{\circ }{\overline{\Omega }}_1=\overset{\circ }{\overline{\Omega }}_2=\Omega_2.$$ $\blacksquare$

\bigskip

\begin{prop}\label{PbI-froniera:23-11-22-16}
	$G$, defined by \eqref{PbI-froniera:23-11-22-11} is an open nonempty set and it enjoys the following properties 
	
	\smallskip
	(a) $G\subset\Omega_1\cap \Omega_2 $,
	
	\smallskip
	
	(b) $\Gamma^{(a)}\subset \overline{G}$,
	
	\smallskip
	
	(c) $G$ is connected. 
\end{prop}
\textbf{Proof.} Since $\Omega_k$, $k=1,2$, are of class $C^{\infty}$, there exist $r_0, M_0$, positive numbers, such that $\Omega_k$, $k=1,2$, are of class $C^{1,1}$ with constant  $r_0, M_0$. Proposition \ref{funzdist:20-11-22-14} implies that there exists  $\mu_1>0$ such that, the following map is continuous

\begin{equation*}
	\Phi: \partial\Omega\times \left(0,\mu_1r_0\right)\rightarrow \ \mathbb{R}^n, 
\end{equation*}

\begin{equation*}
	\Phi(y,t)=y-t\nu(y), \ \ \forall (y,t)\in \partial\Omega\times \left(0,\mu_1r_0\right).
\end{equation*}

For any $\delta\in \left(0,\mu_1r_0\right)$ we denote

\begin{equation*}
	\Lambda_{\delta}=\Phi\left(\overset{\circ }{\Gamma^{(a)}}\times \left(0,\delta\right) \right).
\end{equation*}
Let us check 
\begin{equation}\label{PbI-froniera:23-11-22-17}
	\Lambda_{\delta}\in \mathcal{A},
\end{equation}
where $\mathcal{A}$  is defined in \eqref{PbI-froniera:23-11-22-12}.
It evident that $\Lambda_{\delta}\neq \emptyset$. Moreover, by the continuity of $\Phi^{-1}$ (see (c) of Proposition \ref{funzdist:20-11-22-14}), $\Lambda_{\delta}$ is an open set (see (a) of Proposition \ref{funzdist:20-11-22-14})
$$\Lambda_{\delta}\subset \Omega_1,\quad\mbox{and}\quad \Lambda_{\delta}\subset \Omega_2.$$
Also,  $\Lambda_{\delta}$ is connected, as the image of the connected set $\Gamma^{(a)}\times \left(0,\delta\right)$ by means of the continuous map $\Phi$. To complete the proof of \eqref{PbI-froniera:23-11-22-17}, it suffices to check
\begin{equation}\label{PbI-froniera:23-11-22-18}
	\Gamma^{(a)}\subset \overline{\Lambda}_{\delta}.
\end{equation}
Let $x\in \Gamma^{(a)}$ and let $r>0$ arbitrary. Since, by \eqref{PbI-froniera:22-11-22-8},
 
\begin{equation*}
	\overline{\overset{\circ }{\Gamma^{(a)}}}=\Gamma^{(a)},
\end{equation*}
we have $B_r(x)\cap \overset{\circ }{\Gamma^{(a)}}\neq \emptyset$. Let $y\in B_r(x)\cap \overset{\circ }{\Gamma^{(a)}}$, if $t$  is a positive number small enough, we have 
$$y-t\nu(y)\in  B_r(x)\cap \Lambda_{\delta}.$$
Therefore
$$B_r(x)\cap \Lambda_{\delta}\neq \emptyset,\quad\forall r>0.$$
Hence $x\in \overline{\Lambda}_{\delta}$ and \eqref{PbI-froniera:23-11-22-18} is proved. Now, let us notice that

\begin{equation}\label{PbI-froniera:23-11-22-19}
	 \Lambda_{\delta}\cap A\neq \emptyset, \quad\forall A\in \mathcal{A}.
\end{equation}
 Let us fix $A\in \mathcal{A}$ and let $x\in \overset{\circ }{\Gamma^{(a)}}$. Since $x$ is an interior point of $\Gamma^{(a)}$,  $\Gamma^{(a)}\subset \partial\Omega_1\cap \partial\Omega_2$ and since $\Omega_1$ and $\Omega_2$ are of class $C^{1,1}$, there exists  $\overline{r}>0$ such that
$$B_{\overline{r}}(x)\cap \Omega_1=B_{\overline{r}}(x)\cap \Omega_2\subset \Lambda_{\delta}.$$ Hence
$$B_{\overline{r}}(x)\cap A  \subset \Lambda_{\delta}.$$ On the other hand, since $\overset{\circ }{\Gamma^{(a)}}\subset \overline{A}$, we have $x\in \overline{A}$. Hence
$$\emptyset\neq B_{\overline{r}}(x)\cap A\subset  \Lambda_{\delta}\cap A,$$
which implies \eqref{PbI-froniera:23-11-22-19}. 

\medskip

Now, since
\begin{equation}\label{PbI-froniera:23-11-22-20}
G=\bigcup_{A\in \mathcal{A}}A
\end{equation}
 and since $\Lambda_{\delta}\in \mathcal{A}$, we have that $G\neq \emptyset$ and (trivially) \\ $$G\subset \Omega_1\cap\Omega_2,\ \ \Gamma^{(a)}\subset\overline{\Lambda}_{\delta}\subset \overline{G}.$$ Hence (a) and (b) are proved.
 It remains to prove (c). Let $x,y\in G$ and let $A, B\in \mathcal{A}$ satisfy $x\in A$ e $y\in B$. By \eqref{PbI-froniera:23-11-22-19}, we have $\Lambda_{\delta}\cap A\neq\emptyset$ and $\Lambda_{\delta}\cap B\neq\emptyset$. Let 
 $$z\in  \Lambda_{\delta}\cap A\quad\mbox{and}\quad w\in \Lambda_{\delta}\cap B$$ and let 
  $\gamma_1$ be a continuous path that joins $x$ and $z$ in $A$, $\gamma_2$ be a continuous path that joins $z$ and $w$ in $\Lambda_{\delta}$ and  $\gamma_3$ be a continuous path that joins $w$ and $y$ in $B$.
 Set  $$\gamma=\gamma_1\vee\gamma_2\vee \gamma_3.$$ 
$\gamma$ is a continuous path that joins $x$ and $y$ in $A\cup \Lambda_{\delta}\cup B$. Hence  $A\cup \Lambda_{\delta}\cup B\in \mathcal{A}$ and we have $A\cup \Lambda_{\delta}\cup B\subset G$. Consequently  $\gamma$ joins $x$ and $y$ in $G$. All in all, $G$ is connected. $\blacksquare$  

\bigskip

We have

\begin{prop}\label{PbI-froniera:25-11-22-1}
Let $C$ a nonempty set of $\mathbb{R}^n$ and $A$ an open set of $\mathbb{R}^n$.
then
\begin{equation}\label{PbI-froniera:25-11-22-2}
	A\cap \overline{C}\neq \emptyset \Longleftrightarrow A\cap C\neq \emptyset
\end{equation} 
\end{prop}
\textbf{Proof.} The implication "$\Longleftarrow$" is trivial. Concerning the implication "$\Longrightarrow$", let $z\in A\cap \overline{C}$ and $B_r(z)\subset A$. Since $z\in\overline{C}$, we have $\emptyset\subsetneq B_r(z)\cap C\subset A\cap C$. Therefore $A\cap C\neq \emptyset$. $\blacksquare$  

\begin{prop}\label{PbI-froniera:25-11-22-3}
Let $G$ defined by \eqref{PbI-froniera:23-11-22-11}. Let us suppose 
\begin{equation}\label{PbI-froniera:25-11-22-3-0}
	\Omega_1\setminus\overline{\Omega}_2\neq \emptyset.
\end{equation}

Then we have

\begin{equation}\label{PbI-froniera:25-11-22-4}
	\partial\left(\Omega_1\setminus \overline{G}\right)\subset \Gamma_1^{(i)}\cup \left(\Gamma_2^{(i)}\cap \partial G\right).
\end{equation}
\end{prop}
\textbf{Proof.} 

\textbf{Step I.} First, let us notice that, since $\Omega_1\setminus G\supset  \Omega_1\setminus \left(\overline{\Omega_1\cap\Omega_2}\right)\varsupsetneq \emptyset$ and $\Omega_1\setminus \overline{G}\neq \mathbb{R}^n$, we have $\partial\left(\Omega_1\setminus \overline{G}\right)\neq \emptyset$. 

Now we prove that
\begin{equation}\label{PbI-froniera:25-11-22-5}
	\partial\left(\Omega_1\setminus \overline{G}\right)\subset \partial\Omega_1\cup \partial\Omega_2.
\end{equation}
 We argue by contradiction. Let $x_0\in \partial\left(\Omega_1\setminus \overline{G}\right)$ and let us suppose that  

\begin{equation}\label{PbI-froniera:25-11-22-6}
	x_0\notin \partial\Omega_1\cup \partial\Omega_2.
\end{equation}
It cannot be the case that $x_0\notin \overline{\Omega}_1$ because if it were, it would exist $\rho>0$ such that $B_{\rho}(x_0)\cap \left(\Omega_1\setminus \overline{G}\right)\subset B_{\rho}(x_0)\cap \Omega_1=\emptyset$, that contradicts $\partial\left(\Omega_1\setminus \overline{G}\right)\neq \emptyset$. Let now examine the following two cases:

\smallskip

(a) $x_0\in \Omega_1\cap \Omega_2$
 
\smallskip

(b) $x_0\in \Omega_1\cap \left(\mathbb{R}^n\setminus\overline{\Omega}_2\right)$.

\medskip

\textbf{Case (a)}. Let $\rho>0$ such that 

\begin{equation}\label{PbI-froniera:25-11-22-7}
	B_{\rho}(x_0)\subset \Omega_1\cap \Omega_2. 
\end{equation}
On the other hand, $x_0\in \partial\left(\Omega_1\setminus \overline{G}\right)$, hence

\begin{equation}\label{PbI-froniera:25-11-22-8}
	\begin{cases}
		B_{\rho}(x_0)\cap \left(\Omega_1\setminus \overline{G}\right)\neq \emptyset, \\
		\\
		B_{\rho}(x_0)\cap\left(\mathbb{R}^n\setminus\left(\Omega_1\setminus \overline{G}\right)\right)\neq \emptyset.\\
\end{cases} 
\end{equation}
By \eqref{PbI-froniera:25-11-22-7} and by the second of \eqref{PbI-froniera:25-11-22-8}, we have
\begin{equation*}
	\begin{aligned}
\emptyset\subsetneq B_{\rho}(x_0)\cap\left(\mathbb{R}^n\setminus\left(\Omega_1\setminus \overline{G}\right)\right)&=B_{\rho}(x_0)\cap\left[\left(\mathbb{R}^n\setminus\Omega_1\right)\cup\left(\overline{G}\cap\Omega_1\right)\right]=\\&=
B_{\rho}(x_0)\cap \left(\overline{G}\cap\Omega_1\right)\subset\\&\subset
B_{\rho}(x_0)\cap \overline{G}.
	\end{aligned}
\end{equation*} 
All in all, $B_{\rho}(x_0)\cap \overline{G}\neq \emptyset$ and, by Proposition \ref{PbI-froniera:25-11-22-1} we get
\begin{equation}\label{PbI-froniera:25-11-22-9}
B_{\rho}(x_0)\cap G\neq \emptyset.
\end{equation}
Moreover, by the first relationship of \eqref{PbI-froniera:25-11-22-8}, we obtain
\begin{equation}\label{PbI-froniera:25-11-22-9-0}
	B_{\rho}(x_0)\cap\left(\mathbb{R}^n\setminus \overline{G}\right)\supset B_{\rho}(x_0)\cap\left(\Omega_1\setminus\overline{G}\right)\supsetneq \emptyset.
\end{equation}

Therefore (recalling \eqref{PbI-froniera:25-11-22-7}) 
\begin{equation}\label{PbI-froniera:25-11-22-10}
	G\subsetneq B_{\rho}(x_0)\cup G\subset \Omega_1\cap\Omega_2.
\end{equation}

\medskip

\textbf{Claim.} $B_{\rho}(x_0)\cup G$ is a connected set. 

\medskip

\textbf{Proof of the Claim.} Let $x,y\in B_{\rho}(x_0)\cup G$. Let us prove that there exists a continuous path  $\gamma$  that joins $x$ and $y$ in $B_{\rho}(x_0)\cup G$. If  $x,y\in G$ we have nothing to prove because (by Proposition \ref{PbI-froniera:23-11-22-16}) $G$ is connected.  If $x,y\in B_{\rho}(x_0)$, of course we have nothing to prove. Hence, let us suppose that $x\in G$ and $y\in B_{\rho}(x_0)$. By \eqref{PbI-froniera:25-11-22-9} there exists $z\in B_{\rho}(x_0)\cap G$. Let $\gamma_1$ be a continuous path that joins $x$ and $z$ in $G$ and $\gamma_2$ be a continuous path that joins $z$ and $y$ in $B_{\rho}(x_0)$. Let $$\gamma=\gamma_1\vee\gamma_2.$$ $\gamma$ is a continuous path that joins $x$ and $y$ in $B_{\rho}(x_0)\cup G$. Claim is proved.

\medskip

Now, from (b) of Proposition \ref{PbI-froniera:23-11-22-16} we have  $\overline{B_{\rho}(x_0)\cup G}\supset \Gamma^{(a)}$ and recalling the definition of $G$ we would have
$$G=B_{\rho}(x_0)\cup G,$$
which contradicts \eqref{PbI-froniera:25-11-22-10}. 
\medskip

\textbf{Case (b).} Let $\overline{\rho}>0$ satisfy

\begin{equation}\label{PbI-froniera:25-11-22-11}
B_{\overline{\rho}}(x_0)\subset \Omega_1\cap\left(\mathbb{R}^n\setminus\overline{\Omega}_2\right).
\end{equation}
Since $x_0\in \partial\left(\Omega_1\setminus G\right)$, we again obtain \eqref{PbI-froniera:25-11-22-9}. Therefore

$$\emptyset\subsetneq B_{\overline{\rho}}(x_0)\cap G\subset \Omega_1\cap\Omega_2\subset\Omega_2,$$ hence 
$B_{\rho}(x_0)\cap \Omega_2\neq \emptyset$, that contradicts \eqref{PbI-froniera:25-11-22-10}. Hence \eqref{PbI-froniera:25-11-22-11} is proved.

\bigskip

\textbf{Step II.} Now we prove 
\begin{equation}\label{PbI-froniera:25-11-22-12-0}
	\partial\left(\Omega_1\setminus G\right)\subset \Gamma^{(i)}_1\cup \Gamma^{(i)}_2.
\end{equation}
To this aim, let us prove
\begin{equation}\label{PbI-froniera:25-11-22-12}
\overset{\circ }{\Gamma^{(a)}}\cap \partial\left(\Omega_1\setminus\overline{G}\right)=\emptyset.
\end{equation}
Let us suppose that \eqref{PbI-froniera:25-11-22-12} does not hold. Hence, there exists  $\overline{x}$ such that
\begin{equation}\label{PbI-froniera:25-11-22-13}
	\overline{x}\in \overset{\circ }{\Gamma^{(a)}}\cap \partial\left(\Omega_1\setminus\overline{G}\right).
\end{equation}
Since $\overline{x}\in \overset{\circ }{\Gamma^{(a)}}$ and $\Gamma^{(a)}\subset \left(\partial\Omega_1\right)\cap\left(\partial\Omega_2\right)$, taking into account that $\partial\Omega_k$, $k=1,2$, are of class $C^{1,1}$, there exist $\overline{r}, \overline{M}$ positive numbers, such that 

\smallskip

(i) $Q_{\overline{r}, 2\overline{M}}\left(\overline{x}\right)\cap \partial\Omega_1 = Q_{\overline{r}, 2\overline{M}}\left(\overline{x}\right)\cap \partial\Omega_2\subset \overset{\circ }{\Gamma^{(a)}}$,

\smallskip

 (ii) $Q_{\overline{r}, 2\overline{M}}\left(\overline{x}\right)\cap \Omega_1 = Q_{\overline{r}, 2\overline{M}}\left(\overline{x}\right)\cap \Omega_2$ and they are connected.
 
\smallskip

Now we prove that 
\begin{equation}\label{PbI-froniera:25-11-22-14}
	G\subsetneq G\cup\left( Q_{\overline{r}, 2\overline{M}}\left(\overline{x}\right)\cap \Omega_1\right)\subset \Omega_1\cap\Omega_2
\end{equation}
and
\begin{equation}\label{PbI-froniera:25-11-22-14-0}
	G\cup\left( Q_{\overline{r}, 2\overline{M}}\left(\overline{x}\right)\cap \Omega_1\right)\quad \mbox{is connected. }
\end{equation}
 From these we will arrive to a contradiction.

First, by $G \subset \Omega_1\cap\Omega_2$ we have
\begin{equation}\label{PbI-froniera:25-11-22-15}
G\cap\left(Q_{\overline{r}, 2\overline{M}}\left(\overline{x}\right)\cap \Omega_1\right)\subset \Omega_1\cap\Omega_2.
\end{equation}
On the other hand, since, by \eqref{PbI-froniera:25-11-22-13} and $\overline{x}\in \partial\left(\Omega_1\setminus\overline{G}\right)$, we obtain 
\begin{equation*}
\emptyset\subsetneq Q_{\overline{r}, 2\overline{M}}\left(\overline{x}\right)\cap \left(\Omega_1\setminus\overline{G}\right)\subset \left(Q_{\overline{r}, 2\overline{M}}\left(\overline{x}\right)\cap \Omega_1\right)\setminus G.
\end{equation*}
Hence 
\begin{equation*}
G \subsetneq G\cup \left(Q_{\overline{r}, 2\overline{M}}\left(\overline{x}\right)\cap \Omega_1\right)
\end{equation*}
and \eqref{PbI-froniera:25-11-22-14} is proved. 

In order to prove that the set $G\cup \left(Q_{\overline{r}, 2\overline{M}}\left(\overline{x}\right)\cap \Omega_1\right)$ is connected,  we first of all check that 
\begin{equation}\label{PbI-froniera:25-11-22-16}
	\left(Q_{\overline{r}, 2\overline{M}}\left(\overline{x}\right)\cap \Omega_1\right) \cap G \neq \emptyset,
\end{equation}
As a matter of fact, by $\overline{x}\in \partial\left(\Omega_1\setminus\overline{G}\right)$ we have

\begin{equation*}
	\begin{aligned}
		\emptyset& \subsetneq \left(Q_{\overline{r}, 2\overline{M}}\left(\overline{x}\right)\cap \Omega_1\right) \cap \left[\mathbb{R}^n\setminus\left(\Omega_1\setminus\overline{G}\right)\right]=\\&=\left(Q_{\overline{r}, 2\overline{M}}\left(\overline{x}\right)\cap \Omega_1\right)\cap \left[\left(\mathbb{R}^n\setminus\Omega_1\right)\cup\left(\Omega_1\cap \overline{G}\right)\right]=\\&=\left(Q_{\overline{r}, 2\overline{M}}\left(\overline{x}\right)\cap \Omega_1\right)\cap \overline{G}.
	\end{aligned}
\end{equation*}
Hence
$$\left(Q_{\overline{r}, 2\overline{M}}\left(\overline{x}\right)\cap \Omega_1\right)\cap \overline{G}\neq \emptyset$$
and by Proposition \ref{PbI-froniera:25-11-22-1} we get \eqref{PbI-froniera:25-11-22-16}. 

At this point, in order to prove that $\left(Q_{\overline{r}, 2\overline{M}}\left(\overline{x}\right)\cap \Omega_1\right) \cup G$ is connected, we need only recall (ii) and to proceed as in the proof of the Claim in Step I. Therefore \eqref{PbI-froniera:25-11-22-14-0} is proved. Now, since $$\overline{G\cup\left(Q_{\overline{r}, 2\overline{M}}\left(\overline{x}\right)\cap \Omega_1\right)}\supset \Gamma^{(a)},$$ by \eqref{PbI-froniera:25-11-22-14}, \eqref{PbI-froniera:25-11-22-14-0} and by the definition of  $G$, we arrive to a contradiction. Hence \eqref{PbI-froniera:25-11-22-12} holds, which combined with \eqref{PbI-froniera:25-11-22-5} implies

\begin{equation*}
	\begin{aligned}
		\partial\left(\Omega_1\setminus G\right)&\subset \left(\partial\Omega_1\cup \partial\Omega_2\right)\setminus \overset{\circ }{\Gamma^{(a)}}=\\&=
		\left(\partial\Omega_1\setminus \overset{\circ }{\Gamma^{(a)}} \right)\cup \left(\partial\Omega_2\setminus \overset{\circ }{\Gamma^{(a)}} \right)=\\&=
		\Gamma^{(i)}_1\cup \Gamma^{(i)}_2,  
	\end{aligned}
\end{equation*}
which gives \eqref{PbI-froniera:25-11-22-12-0}. 

\bigskip

\textbf{Step III.} We conclude the proof of \eqref{PbI-froniera:25-11-22-4}.  Let $$x\in \partial\left(\Omega_1\setminus \overline{G}\right).$$
Let us distinguish two cases:

\smallskip

(j) $x\in \mathbb{R}^n\setminus\overline{\Omega}_1$,

\smallskip

(jj) $x\in \overline{\Omega}_1$.

\smallskip

Case (j) cannot occur, because if $x\in \mathbb{R}^n\setminus\overline{\Omega}_1$ then there exists $r>0$ such that $B_r(x)\subset  \mathbb{R}^n\setminus\overline{\Omega}_1$, hence 
$$B_r(x)\cap\left(\Omega_1\setminus \overline{G}\right)\subset B_r(x)\cap \Omega_1=\emptyset,$$ but this cannot hold because $x\in \partial\left(\Omega_1\setminus \overline{G}\right)$. 

\smallskip

Let us consider case (jj). If $x\in \Gamma^{(i)}_1$ then trivially $x\in  \Gamma^{(i)}_1\cup \left(\Gamma^{(i)}_2\cap \partial G\right)$. Instead, if  $x\notin \Gamma^{(i)}_1$ then, by \eqref{PbI-froniera:25-11-22-12-0}, we have

\begin{equation}\label{PbI-froniera:26-11-22-1}
	x\in \Gamma^{(i)}_2.
\end{equation}
Moreover, since, by \eqref{PbI-froniera:25-11-22-12}, we have $x\notin \overset{\circ }{\Gamma^{(a)}}$ and since, by \eqref{PbI-froniera:22-11-22-9-0}, we have $\partial\Omega_1=\Gamma^{(i)}_1\cup \overset{\circ }{\Gamma^{(a)}}$, we get $x\notin \partial \Omega_1$ and, taking into account that $x\in \overline{\Omega}_1$, we have $x\in \Omega_1$. Hence there exists $r>0$ such that

\begin{equation}\label{PbI-froniera:25-11-22-17}
B_r(x)\subset \Omega_1.
\end{equation}
Now, let $s \in (0,r]$ be arbitrary. Since $x\in \partial\left(\Omega_1\setminus \overline{G}\right)$ we have

\begin{equation*}
	\begin{aligned}
		\emptyset &\subsetneq B_{s}(x)\cap\left[\mathbb{R}^n\setminus \left(\Omega_1\setminus \overline{G}\right)\right]= \\&=
 B_{s}(x)\cap \left(\Omega_1\cap \overline{G}\right).
	\end{aligned}
\end{equation*}
Hence
\begin{equation}\label{PbI-froniera:25-11-22-18}
	B_{s}(x)\cap \overline{G}\neq \emptyset,\quad \forall s \in (0,r],
\end{equation}
therefore
\begin{equation}\label{PbI-froniera:25-11-22-19}
	x\in \overline{G}.
\end{equation}
On the other hand, since $x\in \partial\left(\Omega_1\setminus \overline{G}\right)$, we have for every $s>0$
$$\emptyset \subsetneq B_{s}(x)\cap \left(\Omega_1\setminus \overline{G}\right)\subset B_{s}(x)\cap \left(\mathbb{R}^n\setminus G\right).$$
Therefore
\begin{equation*}
x\in \overline{\mathbb{R}^n\setminus G},
\end{equation*}
that combined with \eqref{PbI-froniera:25-11-22-19} gives $x\in \partial G$, which in turn combined with \eqref{PbI-froniera:26-11-22-1} concludes proof. $\blacksquare$

\chapter{The Hadamard example. Solvability of the Cauchy problem and continuous dependence by the data}\markboth{Chapter 9. The Hadamard example. Solvability}{}
\label{LM}

\section{The Hadamard example}\label{esempio-Had} We present the
Hadamard example \index{Hadamard example}relating to the Cauchy problem. 

Let us consider the following Cauchy problem

\begin{equation}\label{Cauchy-Had}
\begin{cases}
\partial^2_xu+\partial^2_yu=0, &  \\
\\
u(x,0)=\varphi(x),
 &  \mbox{for } x\in (0,1),\\
 \\
 \partial_yu(x,0)=\psi(x), &  \mbox{ for } x\in (0,1).\\
\end{cases}
\end{equation}
If $\varphi$ and $\psi$ are
analytic, by the Cauchy--Kovalevskaya Theorem,  there exists a unique analytic solution to problem
\eqref{Cauchy-Had}. Such a solution, by the  Holmgren Theorem, is the unique
solution of class $C^2$ in a neighborhood of $(0,1)\times\{0\}$.
Let now
$$\varphi_{\nu}(x)=e^{-\sqrt{\nu}}\sin\nu x,\quad\mbox{and}\quad \psi_{\nu}(x)=0,\quad\nu\in \mathbb{N}.$$
It is immediately checked that
$$u_{\nu}(x,y)=e^{-\sqrt{\nu}}\sin\nu x\sinh \nu y,\quad\nu\in \mathbb{N},$$
is the solution to the Cauchy problem

\begin{equation}\label{Cauchy-Had-24}
\begin{cases}
\partial^2_xu_{\nu}+\partial^2_yu_{\nu}=0, &  \\
\\
u_{\nu}(x,0)=\varphi_{\nu}(x),
 &  \mbox{for } x\in (0,1),\\
 \\
 \partial_yu_{\nu}(x,0)=\psi_{\nu}(x), &  \mbox{ for } x\in (0,1).\\
\end{cases}
\end{equation}
Let us note that for every $k\in \mathbb{N}_0$

$$\sup_{(0,1)}\left|\frac{d^k\varphi_{\nu}}{dx^k}\right|\rightarrow 0,\quad\mbox{as } \nu\rightarrow\infty,$$
on the other hand, for every $a,b\in (0,1)$, $a<b$ and for every $\delta>0$ we have

$$\sup_{[a,b]\times [-\delta,\delta]}\left|u_{\nu}\right|\rightarrow +\infty,\quad\mbox{as } \nu\rightarrow\infty.$$
In other words, \textbf{"small errors" on the data of
Cauchy problem \eqref{Cauchy-Had} yield uncontrollable errors
on the solution}. This phenomenon makes problem \eqref{Cauchy-Had} essentially intractable
for the applications. Actually, in any problem of an applied nature, the data, in the present case the initial data, are derived through measurements and these are necessarily approximated with
some error, so in order to be able to practically use the
mathematical solution it is necessary that it depends continuously by the
data.

\medskip

In a broad way, we may present the notion of \textbf{well--posed problem
 in the sense of Hadamard} \index{well-posed problem in the sense of Hadamard} as follows. Let $X$ and $Y$ be two metric spaces and be
$$A:X\rightarrow Y$$
a map from $X$ to $Y$. Let us consideri the problem of
determining $x\in X$ such that
\begin{equation}\label{pb-Had}
A(x)=f,
\end{equation}
where $f\in Y$.

\medskip

We say that problem \eqref{pb-Had} is well--posed in the sense of Hadamard provided that we have

\begin{enumerate}
    \item \textbf{(Existence)} for any $f\in Y$ there exists at least one $x\in X$ such that $A(x)=f$.
    \item \textbf{(Uniqueness)} for any $x_1, x_2\in X$ which satisfy $A(x_1)=A(x_2)$ we have $x_1=x_2$;
    \item \textbf{(Continuous dependence by the data)} let us suppose that condition 2 is satisfied, then the map 
    $$A^{-1}: A(X)\rightarrow X$$ is continuous ($A(X)$ with the topology induced by $Y$);
\end{enumerate}

\medskip
In the next Section we will see that in the Cauchy problem there is an interesting relationship between the
first two points above (existence and uniqueness) and the third point (continuous dependence by the data).

\bigskip

\section[Solvability of the Cauchy problem]{Solvability of the Cauchy problem and its
	relations to the continuous dependence on the data} \label{LM-24}
In this Section we will use some theorems from Functional Analysis
on topological vector spaces of which, however, we will not
give the proof. As a reference book we will use
W. Rudin's book \cite{Ru} to which we refer for the above-mentioned
proofs and for further consideration. We will quote
detailed references from of \cite{Ru} as we go along. 

\bigskip

Let us recall the following Theorem of General Topology.
\begin{theo}[\textbf{Baire}]\label{Baire}
	\index{Theorem:@{Theorem:}!- Baire@{- Baire}}
Let $\mathcal{X}$ be a complete metric space. Then for any
countable family of closed subset, $\{F_n\}_{n\in
\mathbb{N}}$, which satisfies
\begin{equation}\label{1-LM}
\mbox{Int}(F_n)=\emptyset,\quad\forall n\in \mathbb{N},
\end{equation}
we have 
\begin{equation}\label{2-LM}
\mbox{Int}\left(\bigcup_{n\in \mathbb{N}}F_n\right)=\emptyset.
\end{equation}
\end{theo}

\bigskip

\begin{definition}[\textbf{topological vector space}]\label{svt}
	\index{Definition:@{Definition:}!- topological vector space@{- topological vector space}}
Let $\mathcal{X}$ be a vector space on $\mathbb{C}$ (or
$\mathbb{R}$) equipped with a topology $\tau$. We say that $\mathcal{X}$
is a topological vector space

(a) for every $x\in \mathcal{X}$, $\{x\}$  is closed w.r.t. $\tau$

(b) the maps $$\mathcal{X}\times\mathcal{X}\ni
(x,y)\rightarrow x+y\in \mathcal{X}$$ and
$$\mathcal{X}\times\mathbb{C}\ni (x,\lambda)\rightarrow \lambda x\in \mathcal{X}$$
(or,  $\mathcal{X}\times\mathbb{R}\ni (x,\lambda)\rightarrow
\lambda x\in \mathcal{X}$) are continuous.
\end{definition}

\medskip

Let $\mathcal{X}$ be a vector space, we say that
$p:\mathcal{X}\rightarrow [0,+\infty)$ is a \textbf{seminorm} \index{seminorm}
on $\mathcal{X}$ provided we have $$p(x+y)\leq p(x)+p(y),\quad \forall x,y\in
\mathcal{X}$$ and $$p(\lambda x)\leq |\lambda|p(x),\quad \forall \lambda
\in \mathbb{C},\mbox{ }\forall x\in \mathcal{X}.$$

\medskip

Let $\mathcal{X}$ be a topological vector space:

\noindent\textbf{(i)} We say that  $\mathcal{X}$ is \textbf{locally
convex}\index{locally convex topological vector space} if there exists a local base of neighborhoods of  $0$ whose members are convex. By the condition (b) of Definition \ref{svt}
it is clear that if  $\mathcal{U}$ is a neighborhood of $0$ then,
for every $x\in \mathcal{X}$, $x+\mathcal{U}$ is a neighborhood of $x$
and conversely. Hence, given a local base of neighborhoods of  $0$
it turns out defined trivially a local base of neighborhoods of each point of $\mathcal{X}$;

\noindent\textbf{(ii)} We say that $\mathcal{X}$ is a
\textbf{F--space} \index{F--space}if its topology is induced by a complete metric
$d$ which is invariant w.r.t. translation (i.e. $d(x+z,y+z)=d(x,y)$ for every $x,y,z\in \mathcal{X}$);

\noindent\textbf{(iii)} We say that $\mathcal{X}$ is a
\textbf{Fr\'{e}chet space} \index{Fr\'{e}chet space} if it is a locally convex F--space.

\bigskip

Let $\Omega$ be an open set of $\mathbb{R}^n$. If $\Omega$ is bounded, We set, as usual,  for any $u\in
C^k\left(\overline{\Omega}\right)$, $k\in \mathbb{N}_0$,

$$\left\Vert u\right\Vert_{C^k\left(\overline{\Omega}\right)}=\sum_{j=0}^k
\sum_{|\alpha|\leq k} \max_{\overline{\Omega}}\left\vert
\partial^{\alpha}u\right\vert.$$
As it is well known, $C^k\left(\overline{\Omega}\right)$, equipped with
$\left\Vert \cdot\right\Vert_{C^k\left(\overline{\Omega}\right)}$, is a Banach
space. Also, we recall that  $C^{k,\sigma}\left(\overline{\Omega}\right)$, $0<\sigma\leq 1$, is a Banach space equipped with the norm

 $$\left\Vert u\right\Vert_{C^{k,\sigma}\left(\overline{\Omega}\right)}=\left\Vert u\right\Vert_{C^k\left(\overline{\Omega}\right)}+[u]_{\Omega;k,\sigma}.$$ 

\medskip

Now we equip $C^k(\Omega)$, where $k\in
\mathbb{N}_0$, or $k=\infty$, with a topology that makes it a space of Fr\'{e}chet space.

We start by $C^0(\Omega)$ ($C^0(\Omega, \mathbb{C})$ or
$C^0(\Omega, \mathbb{R})$, \cite[Ch. 1, Sect. 1.44]{Ru}). Let
$\{K_j\}_{j\in\mathbb{N}}$ be a family of compact sets contained in $\Omega$  such that $K_j\neq \emptyset$,
$K_j\subset\overset{\circ }{K}_{j+1}$ for every $j\in\mathbb{N}$ and
$$\bigcup_{j=1}^{\infty}K_j=\Omega.$$ For any $f\in C^0(\Omega)$,
let us denote by
$$p_{0,j}(f)=\max_{K_j}|f|.$$
$\{p_{0,j}\}_{j\in\mathbb{N}}$ is a family of \textbf{separating seminorms}, \index{sepating seminorms} that is, for each $f\in C^0(\Omega)$, which does not vanish identically, there exists $j\in\mathbb{N}$ such that $p_j(f)\neq 0$.
The collection of sets \cite[Thm. 1.37]{Ru}

 $$\mathcal{V}_j=\left\{f\in C^0(\Omega): p_{0,j}(f)<\frac{1}{j} \right\},$$
 make up a local base in $C^0(\Omega)$ of convex neighborhoods of $0$
 which in turn defines a topology induced by the distance

\begin{equation}\label{metr1-LM}
d_0(f,g)=\sum_{j=1}^{\infty}\frac{2^{-j}p_{0,j}(f-g)}{1+p_{0,j}(f-g)},\quad
\forall f,g\in C^0(\Omega).
\end{equation}
It is proved that $C^0(\Omega)$ with the metric \eqref{metr1-LM} is a  Fr\'{e}chet space (it is quite simple and is left as an exercise).

In a similar way we proceed for $C^k(\Omega)$, $k$ finite. More
precisely, we start, rather than from the seminorms $p_{0,j}$, from the
seminorms

$$p_{k,j}(f)=\max\left\{\max_{K_j}|\partial^{\alpha}f|: |\alpha|\leq k\right\}$$
and we define the distance $d_k$ on $C^k(\Omega)$ by substituting 
$p_{0,j}$ by $p_{k,j}$ in \eqref{metr1-LM}, that is

\begin{equation}\label{metr2-LM}
d_k(f,g)=\sum_{j=1}^{\infty}\frac{2^{-j}p_{k,j}(f-g)}{1+p_{0,j}(f-g)},\quad\forall
f,g\in C^k(\Omega).
\end{equation}
Similarly, it is proved that $C^k(\Omega)$ with the metric
\eqref{metr2-LM} is a Fr\'{e}chet space.

Finally, concerning $C^{\infty}(\Omega)$, the following seminorms (with the corresponding metric),  are defined

$$q_{N}(f)=\max\left\{|\partial^{\alpha}f(x)|:x\in K_{N},\quad |\alpha|\leq N\right\}.$$

\begin{equation}\label{metr3-LM}
d_{\infty}(f,g)=\sum_{N=1}^{\infty}\frac{2^{-N}q_{N}(f-g)}{1+q_{N}(f-g)},\quad\forall
f,g\in C^{\infty}(\Omega).
\end{equation}
in the sequel, when we are dealing
with the convergence of sequences, it will be more convenient to use directly the
seminorms instead of the distances $d_k$, $0\leq k\leq \infty$. For
instance, the sequence $\left\{f_k\right\}$ in
$C^{\infty}(\Omega)$ converges to $f$ in the topology induced by the
norm $d_{\infty}$ if and only if
$$\lim_{k\rightarrow \infty}q_{N}(f_k-f)=0,\quad\forall N\in \mathbb{N}.$$

\medskip

\begin{theo}[\textbf{closed graph}]\label{graf-chiuso}
	\index{Theorem:@{Theorem:}!- closed graph@{- closed graph}}
Let $\mathcal{X}$, $\mathcal{Y}$ be two F-spaces and
\begin{equation}\label{10-LM}
\Lambda:\mathcal{X}\rightarrow \mathcal{Y}
\end{equation}
be a linear map. Then $\Lambda$ is continuous if and only if the graph of $\Lambda$
$$\mathcal{G}=\left\{(x,\Lambda x):x\in \mathcal{X}\right\}$$
is closed in $\mathcal{X}\times\mathcal{Y}$.
\end{theo}
We refer to \cite[Prop. 2.14, Thm 2.15]{Ru} for a proof.
Keep in mind that the most significant implication of
Theorem \ref{graf-chiuso} consists of $$\mathcal{G} \ \mbox{closed}
\Longrightarrow\Lambda \ \mbox{continuous}.$$ The reverse is true even if
$\Lambda$ is nonlinear, with $\mathcal{X}$ and $\mathcal{Y}$ topological spaces
topological and  $\mathcal{Y}$ is a  Hausdorff space

\bigskip

Let
\begin{equation}\label{11-LM}
P(x,\partial)=\sum_{|\alpha|\leq m}a_{\alpha}(x)\partial^{\alpha},
\end{equation}
a linear differential operator of order $m$ and $a_{\alpha}\in C^{\infty}(\mathbb{R}^n,\mathbb{C})$,
for $|\alpha|\leq m$, we say that the Cauchy problem with initial
surface $\{x_n=0\}$ for  $P(x,\partial)$ enjoys
the \textbf{local uniqueness property (in $0$)} \index{local uniqueness property} provided there exists
 $\delta>0$ such that we have: if $u\in
C^{m}\left(\overline{B_{\delta}}\right)$ satisfies

\begin{equation}\label{12-LM}
\begin{cases}
P(x,\partial)u=0, & \mbox{in}\quad B_{\delta}, \\
\\
\partial_n^j u(x',0)=0, & \mbox{for } j=0,1,\cdots, m-1,\quad \forall x'\in B'_{\delta},
\end{cases}
\end{equation}
 then $$u\equiv 0 \ \ \mbox{in } B_{\delta}.$$

 For instance, if the coefficients of $P(x,\partial)$
 are analytic  (in a neighborhood of $0$) and $P(0,e_n)\neq 0$, the Holmgren Theorem
 implies that the \textbf{local uniqueness property} is satisfied.

\medskip

We say that the following Cauchy problem

\begin{equation}\label{13-LM}
\begin{cases}
P(x,\partial)u=0, & \\
\\
\partial_n^j u(x',0)=g_j(x'), & \mbox{for } j=0,1,\cdots, m-1,
\end{cases}
\end{equation}
is \textbf{locally solvable (in the origin)} \index{local solvable Cauchy problem} in $C^{\infty}$
if the following occurs:

\smallskip

\noindent for every open neighborhood $\mathcal{U}$ of $0$  in $\mathbb{R}^{n-1}$
and for every $$g=(g_0,g_1,\cdots, g_{m-1})\in C^{\infty}(\mathcal{U},
\mathbb{C}^m)$$ there exists  $\mathcal{V}$, open neighborhood of $0$ in
$\mathbb{R}^{n}$ and there exists $u\in C^{\infty}(\mathcal{V})$ such that
\begin{equation}\label{13n-LM}
\begin{cases}
P(x,\partial)u=0, & \mbox{in}\quad \mathcal{V}, \\
\\
\partial_n^j u=g_j, & \mbox{for } j=0,1,\cdots, m-1,\quad \mbox{ in } \mathcal{V}\cap\{x_n=0\}.
\end{cases}
\end{equation}

\medskip

The definitions of the \textbf{local uniqueness property} and
of the \textbf{solvability in the origin} of the Cauchy problem with
initial surface the hyperplane $$\left\{x\in \mathbb{R}^n: N\cdot
x=0 \right\},$$ where $N$ is a versor of $\mathbb{R}^n$, is formulated in
obvious way or, simply, by reconducting them to the case $N=e_n$ by means of
an isometry of $\mathbb{R}^n$.

\medskip

The following two theorems and their immediate consequences are known in the literature as the Lax--Mizohata Theorem. Here we present them in a slightly modified form. In particular, Theorem \ref{thm-Lax} is due to Lax, and Theorem \ref{thm-preMizohata} is due to Mizohata. The above theorems are treated, for instance, in \cite{Mizo}.

\begin{theo}\label{thm-Lax}
Let $P(x,\partial)$ be operator \eqref{11-LM} with $C^{\infty}(\mathbb{R}^n,\mathbb{C})$ coefficients. Let us suppose that
$P(x,\partial)$ enjoys the \textbf{local uniqueness property} and that Cauchy problem \eqref{13-LM} is
\textbf{locally solvable} in $C^{\infty}$. Then
for every $\mathcal{U}$, neighborhood of $0$ in $\mathbb{R}^{n-1}$, there exists 
$r>0$ such that for every $$g=(g_0,g_1,\cdots, g_{m-1})\in
C^{\infty}\left(\mathcal{U}, \mathbb{C}^m\right),$$ there exists a unique
$u\in C^m\left(\overline{B_r}\right)$ such that
\begin{equation}\label{20-LM}
\begin{cases}
P(x,\partial)u=0, & \mbox{in } B_r, \\
\\
\partial_n^j u=g_j & \mbox{for } j=0,1,\cdots, m-1,\quad \mbox{ in } B'_r.
\end{cases}
\end{equation}
\end{theo}

\medskip

\noindent\textbf{Remark 1.} Let us observe that in Theorem
\ref{thm-Lax}, unlike in the definition of local solvability, $r$, hence the neighborhood of $0$, $B_r$, \textbf{does not} depend on $g$. Of course $r$
depends on $\mathcal{U}$. $\blacklozenge$

\bigskip

\textbf{Proof.} Let $\mathcal{U}\in \mathbb{R}^{n-1}$ be a neighborhood of $0$ and let
$\sigma\in (0,1)$ be fixed. Since Cauchy problem
\eqref{13-LM} is locally solvable in $0$, by local uniqueness property
we have that, for every $$g=(g_0,g_1,\cdots, g_{m-1})\in C^{\infty}\left(\mathcal{U},
\mathbb{C}^m\right),$$ there exists $\delta>0$ such that
$B'_{\delta}\subset \mathcal{U}$ and such that for every $\rho\in
(0,\delta]$ there exists a unique solution $u\in
C^{m,\sigma}\left(\overline{B_{\rho}}\right)$ of
problem $\mathcal{P}_{g,\rho}$:

\begin{equation*}
(\mathcal{P}_{g,\rho})\quad \quad\quad
\begin{cases}
P(x,\partial)u=0, & \mbox{in } B_{\rho}, \\
\\
\partial_n^j u=g_j, & \mbox{for } j=0,1,\cdots, m-1,\quad \mbox{ in } B'_{\rho}.
\end{cases}
\end{equation*}

 Let $\left\{\rho_k\right\}$ be a strictly decreasing sequence such that 
 $$\lim_{k\rightarrow \infty}\rho_k=0.$$ For any $k,M\in \mathbb{N}$, let us consider the sets
 \begin{equation*}
 	\begin{aligned}
 	&\mathcal{A}_{k,M}= \\&= \left\{g\in C^{\infty}(\mathcal{U}, \mathbb{C}^m):
 	\mbox{ exists } u \mbox{ solution to }\mathcal{P}_{g,\rho_k} \mbox{
 		and } \left\Vert
 	u\right\Vert_{C^{m,\sigma}\left(\overline{B_{\rho_k}}\right)}
 	\leq M\right\} 	\end{aligned}
 \end{equation*}
 
$$.$$

It is evident that $\mathcal{A}_{k,M}$ is symmetric, for any $k,M\in \mathbb{N}$
($g\in\mathcal{A}_{k,M}\Rightarrow -g\in\mathcal{A}_{k,M}$)
and convex for every $k,M\in \mathbb{N}$.

\medskip

\noindent\textbf{Step 1.} Let us check that

\begin{equation}\label{21-3LM}
C^{\infty}\left(\mathcal{U}, \mathbb{C}^m\right)=\bigcup_{k,M\in
\mathbb{N}}\mathcal{A}_{k,M}.
\end{equation}
Of course, it suffices to check "$\subset$". Let $g\in
C^{\infty}(\mathcal{U}, \mathbb{C}^m)$. By the local solvability there exists $\mathcal{V}\subset \mathbb{R}^n$,
neighborhood of $0$, such that Cauchy problem \eqref{13n-LM} admits
a solution $u\in C^{\infty}(\mathcal{V})$. It is enough then to choose
$k\in \mathbb{N}$ such that $\overline{B_{\rho_k}}\subset
\mathcal{V}$ and $M$ such that $M\geq \left\Vert
u\right\Vert_{C^{m,\sigma}\left(\overline{B_{\rho_k}}\right)}$ and we have $g\in \mathcal{A}_{k,M}$.

\medskip

\noindent\textbf{Step 2.} Now we prove that for every $k,M\in
\mathbb{N}$, $\mathcal{A}_{k,M}$ is closed in
$C^{\infty}\left(\mathcal{U}, \mathbb{C}^m\right)$ equipped with the topology
induced by the metric $d_{\infty}$.

We fix  $k,M\in \mathbb{N}$ and let $\{g_{\nu}\}_{\nu\in
\mathbb{N}}$ be a sequence in $\mathcal{A}_{k,M}$ such that
\begin{equation}\label{30-4LM}
\{g_{\nu}\}\rightarrow \widetilde{g}, \quad \mbox{ in } C^{\infty}(\mathcal{U},
\mathbb{C}^m).
\end{equation}
let us check that $\widetilde{g}\in \mathcal{A}_{k,M}$.

First of all, we have
\begin{equation}\label{gtilde-4LM}
\widetilde{g}\in C^{\infty}(\mathcal{U}, \mathbb{C}^m).
\end{equation}
Moreover, denoting by $u_{\nu}$ the solution of problem
$\mathcal{P}_{g_{\nu},\rho_k}$ we have, by the definition of $\mathcal{A}_{k,M}$,
\begin{equation*}\left\Vert u_{\nu}\right\Vert_{C^{m,\sigma}\left(\overline{B_{\rho_k}}\right)}
\leq M,\quad \forall \nu\in\mathbb{N}.
\end{equation*}
Hence, by the Arzel\`{a}--Ascoli Theorem there exists a subsequence  $\left\{u_{\nu_{q}}\right\}$ of
$\left\{u_{\nu}\right\}$ and $\widetilde{u}\in
C^{m,\sigma}\left(\overline{B_{\rho_k}}\right)$ such that

\begin{equation}\label{stella-4LM}
u_{\nu_{q}}\rightarrow \widetilde{u},\quad\mbox{as }
q\rightarrow\infty,\mbox{ in }
C^{m}\left(\overline{B_{\rho_k}}\right).
\end{equation}
Moreover we have
\begin{equation}\label{u-limitata-4LM}\left\Vert \widetilde{u}\right\Vert_{C^{m,\sigma}\left(\overline{B_{\rho_k}}\right)} \leq M,
\end{equation}

Concerning the justification of the latter inequality, we notice that the sequence $\left\{u_{\nu_{q}}\right\}$ converges to
 $\widetilde{u}$ in $C^{m}\left(\overline{B_{\rho_k}}\right)$, but not necessarily it converges in $C^{m,\sigma}\left(\overline{B_{\rho_k}}\right)$, neverthless  \eqref{u-limitata-4LM} holds, since for $x\neq y\in \overline{B_{\rho_k}}$ we have, for any $q\in \mathbb{N}$, 
$$\sum_{|\alpha|=m}\frac{\left\vert\partial^{\alpha}u_{\nu_{q}}(x)-\partial^{\alpha}u_{\nu_{q}}(y)\right\vert}{|x-y|^{\sigma}}\leq M-\left\Vert u_{\nu_q}\right\Vert_{C^{m}\left(\overline{B_{\rho_k}}\right)};$$
hence, from the punctual convergence of $\left\{\partial^{\alpha}u_{\nu_{q}}\right\}$ for $|\alpha|=m$ and by
$$\left\Vert u_{\nu_q}\right\Vert_{C^{m}\left(\overline{B_{\rho_k}}\right)}\rightarrow \left\Vert u\right\Vert_{C^{m}\left(\overline{B_{\rho_k}}\right)},\quad\mbox{as }
q\rightarrow\infty,$$
we obtain \eqref{u-limitata-4LM}.

By \eqref{stella-4LM} we have easily 
\begin{equation*}
0=P(x,\partial)u_{\nu_{q}}\rightarrow
P(x,\partial)\widetilde{u},\quad\mbox{as }
q\rightarrow\infty,\mbox{ in }
 C^{0}\left(\overline{B_{\rho_k}}\right).
\end{equation*}
Therefore
\begin{equation}\label{stella-5LM}
 P(x,\partial)\widetilde{u}=0,\quad\mbox{in } \overline{B_{\rho_k}}.
\end{equation}
By \eqref{30-4LM} and \eqref{stella-4LM} we have, for
$j=0,1,\cdots,m-1$,

\begin{equation}\label{1-6LM}
 \partial_n^ju(x',0)=\lim_{q\rightarrow \infty}\partial_n^ju_{\nu_{q}}u(x',0)=
 \lim_{q\rightarrow \infty}g_{j,\nu_{q}}(x')=\widetilde{g},\quad\forall x'\in\overline{B'_{\rho_k}}.
\end{equation}
By \eqref{gtilde-4LM}, \eqref{u-limitata-4LM}--\eqref{1-6LM} we get $\widetilde{g}\in \mathcal{A}_{k,M}$.

\medskip

\noindent\textbf{Step 3.} Now, recalling that
$C^{\infty}\left(\mathcal{U}, \mathbb{C}^m\right)$ is a Frech\'{e}t space,  by Theorem \ref{Baire} we have that there exist
$k_0,M_0\in \mathbb{N}$ such that
\begin{equation}\label{1-7LM}
\mbox{Int}\left(\mathcal{A}_{k_0,M_0}\right)\neq \emptyset.
\end{equation}
To prove \eqref{1-7LM} we argue by contradiction. Let us assume \eqref{1-7LM} does not hold. Consequently we have
$$ \mbox{Int}\left(\mathcal{A}_{k,M}\right)= \emptyset, \ \ \ \forall k,M\in \mathbb{N},$$
and recalling that $\mathcal{A}_{k,M}$ is closed for every $k,M\in
\mathbb{N}$,  \eqref{21-3LM} and Theorem \ref{Baire} imply

$$ C^{\infty}\left(\mathcal{U}, \mathbb{C}^m\right)=
\mbox{Int}\left(C^{\infty}\left(\mathcal{U},
\mathbb{C}^m\right)\right)=\mbox{Int}\left(\bigcup_{k,M\in\mathbb{N}}\mathcal{A}_{k,M}\right)=\emptyset.$$
This is clearly a contradiction, then  \eqref{1-7LM} needs to hold. 

On the other hand, \eqref{1-7LM} implies that there exists $\psi\in
\mbox{Int}\left(\mathcal{A}_{k_0,M_0}\right)$, and since
$\mathcal{A}_{k_0,M_0}$ is symmetric and convex, we have  $-\psi\in
\mbox{Int}\left(\mathcal{A}_{k_0,M_0}\right)$ therefore
$$0=\frac{1}{2}(\psi-\psi)\in
\mbox{Int}\left(\mathcal{A}_{k_0,M_0}\right).$$ 

All in all, we have 
$$0\in \mbox{Int}\left(\mathcal{A}_{k_0,M_0}\right).$$
Consequently there exists
$\mathcal{W}_0\subset\mbox{Int}\left(\mathcal{A}_{k_0,M_0}\right)$,
where $\mathcal{W}_0$ is an element of local base of neighborhood of 
$0$ in $C^{\infty}(\mathcal{U},\mathbb{C}^m)$, hence  it is of the type 
$$\mathcal{W}_0=\bigcap_{|\alpha|\leq h, 1\leq j\leq h'}\left\{f\in C^{\infty}(\mathcal{U},\mathbb{C}^m):\quad \max_{K_j}|\partial^{\alpha}f|<\varepsilon_{\alpha, j}\right\},$$
where $h, h', \varepsilon_{\alpha, j}$ are suitable positive numbers and $\{K_j\}_{j\in\mathbb{N}}$ is a family of compact subset of  $\mathcal{U}$ which satisfies  $K_j\neq \emptyset$,
$K_j\subset\overset{\circ }{K}_{j+1}$ for every $j\in\mathbb{N}$ and
$$\bigcup_{j=1}^{\infty}K_j=\mathcal{U}.$$

We recall that the local uniqueness property implies that there exists $r>0$, that we may assume less or equal to $\rho_{k_0}$, such that if $w\in C^{m}\left(\overline{B_r}\right)$ satisfies

\begin{equation*}
	\begin{cases}
		P(x,\partial)w=0, & \mbox{in } B_r, \\
		\\
		\partial_n^j w=0, & \mbox{for } j=0,1,\cdots, m-1,\quad \mbox{ in } B'_{r},
	\end{cases}
\end{equation*}
then we have
$$w\equiv 0 \ \ \mbox{in } B_r.$$

Now, let $g\in C^{\infty}\left(\mathcal{U}\right)$, then there exists trivially  $\lambda>0$ such that $$\lambda^{-1}g\in \mathcal{W}_0\subset
\mbox{Int}\left(\mathcal{A}_{k_0,M_0}\right)$$ and by the definition of
$\mathcal{A}_{k_0,M_0}$ we have there exists $v\in
C^m\left(\overline{B_{\rho_{k_0}}}\right)$ (actually, $v\in
C^{m,\sigma}\left(\overline{B_{\rho_{k_0}}}\right)$) solution to
$\mathcal{P}_{\lambda^{-1}g,\rho_{k_0}}$. Hence $v_{|B_r}$ is the unique solution in $C^{m}\left(\overline{B_r}\right)$ to the Cauchy problem
\begin{equation*}
	\begin{cases}
		P(x,\partial)v=0, & \mbox{in } B_r, \\
		\\
		\partial_n^j v=\lambda^{-1}g_j, & \mbox{for } j=0,1,\cdots, m-1,\quad \mbox{ in } B'_{r},
	\end{cases}
\end{equation*}

Therefore, $u=\lambda v$ is the unique solution in
$C^{m}\left(\overline{B_{r}}\right)$ to the Cauchy problem
\begin{equation*}
\begin{cases}
P(x,\partial)u=0, & \mbox{in } B_{r}, \\
\\
\partial_n^j u=g_j & \mbox{for } j=0,1,\cdots, m-1,\quad \mbox{ in } B'_{r}.
\end{cases}
\end{equation*}

$\blacksquare$

\medskip

\textbf{Remark 2.} 
It is evident from the proof of
Theorem \ref{thm-Lax} that if in the definition of local
solvability we require the existence of $u$ in $
C^{m'}(\mathcal{V})$ with $m'>m$, then we reach the same
conclusions. $\blacklozenge$

\medskip
\begin{theo}[\textbf{Lax--Mizohata}]\label{thm-preMizohata}
	\index{Theorem:@{Theorem:}!- Lax--Mizohata@{- Lax--Mizohata}}
Let us suppose that $P(x,\partial)$ satisfies the same assumption
of Theorem \ref{thm-Lax}. Let $\mathcal{U}\subset \mathbb{R}^{n-1}$ be
a neighborhood of $0$ and let $r>0$ be defined in Theorem \ref{thm-Lax}. Let
$\Lambda$ be the map 
$$\Lambda:C^{\infty}\left(\mathcal{U}, \mathbb{C}^m\right)\rightarrow C^m\left(\overline{B_r}\right),$$

$$C^{\infty}\left(\mathcal{U}, \mathbb{C}^m\right)\ni g\rightarrow \Lambda(g)=u
\mbox{ solution to the Cauchy problem}:$$

\begin{equation}\label{Cauchy-miz}
\begin{cases}
P(x,\partial)u=0, & \mbox{in } B_{r}, \\
\\
\partial_n^j u=g_j, & \mbox{for } j=0,1,\cdots, m-1,\quad \mbox{ in } B'_{r}.
\end{cases}
\end{equation}
Then $\Lambda$ is continuous.
\end{theo}

\textbf{Proof.} Let us prove that the graph of $\Lambda$ is closed and then we apply Theorem \ref{graf-chiuso}. Let $\left\{\left(g_{\nu},
u_{\nu}\right)\right\}_{\nu\in \mathbb{N}}$ be a sequence in 
$C^{\infty}\left(\mathcal{U}, \mathbb{C}^m\right)\times
C^m\left(\overline{B_r}\right)$ which satisfies

\begin{equation}\label{Cauchy-mizLM11}
\begin{cases}
P(x,\partial)u_{\nu}=0, & \mbox{in } B_{r}, \\
\\
\partial_n^j u_{\nu}=g_{\nu, j}, & \mbox{for } j=0,1,\cdots, m-1,\quad \mbox{ in }
B'_{r}
\end{cases}
\end{equation}
and

\begin{equation}\label{1-LM11}
\left\{\left(g_{\nu}, u_{\nu}\right)\right\}\rightarrow \left(g,
v\right), \quad\mbox{ in }
C^{\infty} \left(\mathcal{U}, \mathbb{C}^m\right)\times
C^m\left(\overline{B_r}\right).
\end{equation}
Then
$$P(x,\partial)v=\lim_{\nu\rightarrow \infty}
P(x,\partial)u_{\nu}=0$$ and, for $0\leq j\leq m-1$,
$$\partial_n^jv(x',0)=\lim_{\nu\rightarrow \infty}\partial_n^ju_{\nu}(x',0)=\lim_{\nu\rightarrow \infty}g_{\nu,j}(x')=g_j(x').$$
 Hence $v$ solves Cauchy problem
\eqref{Cauchy-miz} and by Theorem \ref{thm-Lax}, we have
$$v=\Lambda g.$$ Therefore the graph of $\Lambda$ is closed, hence
Theorem \ref{graf-chiuso} implies that $\Lambda$ is continuous.
$\blacksquare$

\bigskip

We observe that by Theorem \ref{thm-preMizohata} and the Hadamard example it follows (again) that the Cauchy problem for the Laplace equation,
\eqref{Cauchy-Had}, is not locally solvable in $C^{\infty}$. Indeed, by the Holmgren Theorem, such a Cauchy problem for the Laplace equation enjoys the
property of local uniqueness in $C^{2}$, but, as
shown in the Hadamard example,
the map $\Lambda$ defined for this Cauchy problem is not continuous.

Similarly, Theorem \ref{thm-preMizohata} may be applied to obtain some necessary condition for the local solvability of the Cauchy problem. Here we limit ourselves only to consider the case of the operators with constant coefficients  which are equal to their principal part. We refer to \cite[Ch.5, Sect.4]{HO63},
\cite[Ch.12, Sect.3]{HOII} for the general operators with constant coefficients, and to
\cite[Ch.23, Sect.3]{HOII} for the operators with $C^{\infty}$ coefficients, warning the reader, however, that (especially in \cite[Ch.23, Sect. 3]{HOII}) quite advanced tools are used.

\medskip

Let
\begin{equation}\label{1-Iper}
P_m(\xi)=\sum_{|\alpha|=m}a_{\alpha}\xi^{\alpha},
\end{equation}
be a \textbf{homogeneous polynomial} of degree $m$ with coefficients
$a_{\alpha}\in\mathbb{C}$, for $|\alpha|= m$ and let $N$ be a versor of
$\mathbb{R}^n$. We say that $P_m$ is hyperbolic
\index{hyperbolic:@{hyperbolic:}!- homogeneous polynomial w.r.t. a direction@{- homogeneous polynomial w.r.t. a direction}} with respect to the direction $N$ provided we have

\smallskip 

\noindent(a) $P_m(N)\neq 0$

\noindent(b) for every $\xi\in \mathbb{R}^n$ the algebraic equation in $z$
$$P_m(\xi+zN)=0,$$ has real roots only.

\medskip

\begin{theo}\label{thm-Iper-1}
Let  $N$ be a versor of
$\mathbb{R}^n$ and let $P_m(\partial)$ be the differential operator with constant coefficients
\begin{equation}\label{2-Iper}
P_m(\partial)=\sum_{|\alpha|=m}a_{\alpha}\partial^{\alpha}.
\end{equation}
Let us suppose that the Cauchy problem for operator \eqref{2-Iper}
with initial surface
\begin{equation}\label{3-Iper}
\left\{x\in \mathbb{R}^n: N\cdot x=0 \right\}
\end{equation}
enjoys the local sovability and the local uniqueness property.

Then the polynomial $P_m(\xi)$ is hyperbolic w.r.t. 
the direction $N$.
\end{theo}

\textbf{Proof.} We have already seen in Section \ref{esist-unic} that
by the local uniqueness for the Cauchy problem for operator \eqref{2-Iper} with initial surface
\eqref{3-Iper}, we have

\begin{equation}\label{4-Iper-24}
P_m(N)\neq 0.
\end{equation}
We argue by contradiction to prove that, for every $\xi\in \mathbb{R}^n$, there exist real roots only to  the equation
\begin{equation}\label{4-Iper}
P_m(\xi+zN)=0.
\end{equation}
Let us suppose that there exist
$\xi_0\in\mathbb{R}^n$ and $z_0=\Re z_0+i\Im z_0$, such that
$\Im z_0\neq 0$ and
\begin{equation}\label{5-Iper}
P_m(\xi_0+z_0N)=0.
\end{equation}
Let us denote by $$\eta=\xi_0+\Re z_0 N, \quad\mbox{ and }\quad \tau=\Im
z_0$$ and by

\begin{equation*}
u_{\nu}(x)=\exp\left\{-|\nu|^{1/2}+\nu\left(i\eta\cdot x-\tau N\cdot
x\right)\right\},\quad \nu\in\mathbb{Z}.
\end{equation*}
We have, by \eqref{5-Iper} (recalling that  $\xi_0+z_0N=\eta+i\tau N$ and
$P_m(\xi)$ is a homogeneous polynomial)

\begin{equation}\label{6-Iper}
	\begin{aligned}
	P_m(\partial)u_{\nu}(x)&=u_{\nu}(x)P_m\left(\nu\left(i\eta-\tau
N\right)\right)=\\&=u_{\nu}(x)(i\nu)^mP_m\left(\eta+i\tau N\right)=0.
	\end{aligned}
\end{equation}
Set
$$g_{\nu,j}=\frac{\partial^j}{\partial N^j}u_{\nu},\quad \mbox{for }\nu\in\mathbb{Z}\quad \mbox{and }j=0,1\cdots,m-1,$$
we have trivially that $u_{\nu}$ solves the Cauchy problem

\begin{equation*}
\begin{cases}
P(x,\partial)u_{\nu}=0, & \mbox{in }  \mathbb{R}^n,\\
\\
\frac{\partial^j}{\partial N^j}u_{\nu}=g_{\nu, j}, & \mbox{for }
j=0,1,\cdots, m-1,\quad N\cdot x=0.
\end{cases}
\end{equation*}
On the other hand it is easy to check that 
\begin{equation*}
g_{\nu}\rightarrow 0,\quad\mbox{as }
|\nu|\rightarrow +\infty,\mbox{ in } C^{\infty}.
\end{equation*}
Moreover, if $\tau N\cdot x>0$ then we have
$$\left|u_{\nu}(x)\right|\rightarrow +\infty ,\quad\mbox{as } \nu\rightarrow -\infty$$
and if $\tau N\cdot x<0$ then we have
$$\left|u_{\nu}(x)\right|\rightarrow +\infty ,\quad\mbox{as } \nu\rightarrow +\infty.$$
Hence, the map $\Lambda$ defined in Theorem \ref{thm-preMizohata}
is not continuous. Since $P_m(\partial)$ enjoys the local uniqueness property, $P_m(\partial)$ cannot enjoy at the same time the
local sovability. Thus we have a contradiction. Therefore equation \eqref{4-Iper} has real roots only. $\blacksquare$

\medskip

If the polynomial $P(\xi)$ is not homogeneous, one can prove (with more efforts) a theorem which is similar to Theorem
\ref{thm-Iper-1} that we merely state here (see \cite[Ch.5,
Sect. 4]{HO63} , \cite[Ch.12, Sect. 3]{HOII})

\begin{theo}\label{thm-Iper-2}
Let $P(\xi)$ be a polynomial of degree $m$, $N$ a versor of
$\mathbb{R}^n$ and $P(\partial)$ the following differential operator with 
constant coefficient 
\begin{equation}\label{20-Iper}
P(\partial)=\sum_{|\alpha|\leq m}a_{\alpha}\partial^{\alpha},
\end{equation}
where $a_{\alpha}\in \mathbb{C}$ for $|\alpha|\leq m$.  Let us suppose that
the Cauchy problem for operator \eqref{20-Iper} with initial surface
\begin{equation*}
\left\{x\in \mathbb{R}^n: N\cdot x=0 \right\},
\end{equation*}
enjoys local sovability  property and let us suppose
\begin{equation}\label{21-Iper}
P_m(N)\neq 0,
\end{equation}
where $P_m$ is the principal part of $P$.

Then
\begin{equation}
\label{22-Iper}  \exists\tau_0\in\mathbb{R} \mbox{ such that }
\forall\xi\in \mathbb{R}^n \mbox{ and } \forall\tau<\tau_0\quad
P(i(\xi+i\tau N))\neq 0.
\end{equation}
\end{theo}

A polynomial that enjoys the properties \eqref{21-Iper} and
\eqref{22-Iper} is called \textbf{hyperbolic polynomial w.r.t.
the direction} $N$ \index{hyperbolic:@{hyperbolic:}!- operator w.r.t. a direction@{- operator with respect to a direction}}. It is not difficult to prove (see the
literature quoted above) that in the case where $P(\xi)$ is a homogeneous polynomial
the two definitions of hyperbolicity coincide. Moreover, if $P(\xi)$ is a hyperbolic polynomial (not necessarily
homogeneous), denoted by $P_m(\xi)$ its principal part,  then we have

\smallskip

\noindent(a) $P_m(N)\neq 0$, 

\smallskip
\noindent(b) for every $\xi\in
\mathbb{R}^n$ the equation in $z$ $$P_m(\xi+zN)=0,$$ has real roots only. 

\smallskip

However, if the principal part of $P(\xi)$ satisfies (a)
and (b) it does not imply that $P(x,\partial)$ is hyperbolic w.r.t. $N$. It can be proved that if 
the equation $P_m(\xi+zN)=0$, has  \textbf{simple real roots} only, then the
 hyperbolicity condition with respect to a direction $N$ is
also sufficient for the local solvability of the
Cauchy problem with initial surface $\{N\cdot x=0\}$.

\section{Concluding Remarks}\label{concl-hadamard}

\noindent\textbf{1.} In this Chapter we have considered the
properties of the local uniqueness and the local solvability for the
Cauchy problem, but with simple and natural changes
one could consider the same properties for the one--sided Cauchy problem. Actually, it suffices to replace each neighborhood of $0\in \mathbb{R}^n$, say $\mathcal{V}$, by $\mathcal{V}_+=\mathcal{V}\cap \mathbb{R}^{n-1}\times[0,+\infty)$, and in defininig the seminorms of $C^{\infty}(\mathcal{V}_+)$ we
consider a family of compacts $K_j\subset\mbox{Int}_{\mathcal{V}_+}(K_{j+1})$ for every
$j\in\mathbb{N}$ e $\bigcup_{j=1}^{\infty}K_j=\mathcal{V}_+$. In
this way provided that the appropriate modifications are introduced, Theorems
\ref{thm-preMizohata} and \ref{thm-Lax} preserve their validity,
in particular, the ball $B_r$ should be replaced with the half-ball
$B_r\cap \left(\mathbb{R}^{n-1}\times[0,+\infty)\right)$.

\medskip

\noindent\textbf{2.} In Theorem \ref{thm-preMizohata} we have seen
that if the properties of the local uniqueness and of the
local solvability in $C^{infty}$ hold true then the map
$\Lambda$, defined in Theorem \ref{thm-preMizohata}, is continuous.
In what follows, we check that if the operator $P(x,\partial)$ has 
analytic coefficients, the converse also holds (in a sense).  

Let us suppose that
$$P_m(0,e_n)\neq 0,$$ and let us suppose that there exists $r>0$ such that only the null function solves (in $C^{\infty}\left(B_{r}\right)$)  the  Cauchy problem 

\begin{equation}\label{1-Oss}
\begin{cases}
P(x,\partial)u=0, & \mbox{in } B_{r}, \\
\\
\partial_n^j u=0, & \mbox{for } j=0,1,\cdots, m-1,\quad \mbox{ in }
B'_{r}.
\end{cases}
\end{equation}
Then by the Cauchy-Kovalevskaya Theorem, there exists $\widetilde{r}\leq r$ such that the
following map $\widetilde{\Lambda}$ is well--defined

$$C^{\omega}\left(B'_{r}, \mathbb{C}^m\right)\ni g\rightarrow \widetilde{\Lambda}(g)=u\in C^{\omega}\left(B_{\widetilde{r}}\right)
\mbox{ solution to the Cauchy problem:}$$

\begin{equation}\label{2-Oss}
\begin{cases}
P(x,\partial)u=0, & \mbox{in } B_{\widetilde{r}}, \\
\\
\partial_n^j u=g_j, & \mbox{for } j=0,1,\cdots, m-1,\quad \mbox{ in } B'_{r}.
\end{cases}
\end{equation}

We prove what follows: let us assume that $\widetilde{\Lambda}$ is a continuous map provided that we equip the
spaces $C^{\omega}\left(B'_{r}, \mathbb{C}^m\right)$ and $C^{\omega}\left(B_{\widetilde{r}}\right)$ with the metric $d_{\infty}$, we have that Cauchy problem \eqref{2-Oss} satisfies the
local solvability property in $C^{\infty}$ (i.e. if the initial data of the Cauchy belong to $C^{\infty}$
there exist solutions to the Cauchy problem).

 Let us suppose that $\widetilde{\Lambda}$ is continuous and let
$g\in C^{\infty}\left(B'_{r}, \mathbb{C}^m\right)$. Let  $\left\{g_{,\nu}\right\}$ be the sequence  in $C^{\omega}\left(B'_{r},
\mathbb{C}^m\right)$ defined as follows
$$g_{,\nu}(x')=\left(\frac{\nu}{2\pi}\right)^{(n-1)/2}\int_{B'_{r}}e^{-\frac{\nu|x'-y'|^2}{2}}g(x')dx'.$$
The sequence $\left\{g_{,\nu}\right\}$ converges
to $g$ in $C^{\infty}\left(B'_{r}, \mathbb{C}^m\right)$ and since $g_{,\nu}\in C^{\omega}\left(B'_{r},
\mathbb{C}^m\right)$, for every $\nu\in \mathbb{N}$, the Cauchy-Kovalevskaya Theorem yields the existence of a solution
$u_{\nu}$ which is unique in $C^{\omega}\left(B_{\widetilde{r}}\right)$.
By the continuity of $\widetilde{\Lambda}$ and since
$\left\{g_{\nu}\right\}_{\nu\in \mathbb{N}}$ is a Cauchy sequence in
$C^{\infty}\left(B'_{r}, \mathbb{C}^m\right)$ we derive that
$\left\{u_{\nu}\right\}$ is a Cauchy sequence in
$C^{\infty}\left(B_{\widetilde{r}}(0)\right)$. The completeness of
$C^{\infty}\left(B_{\widetilde{r}}\right)$ implies that $\left\{u_{\nu}\right\}_{\nu\in \mathbb{N}}$ converges to a function $u\in C^{\infty}\left(B_{\widetilde{r}}\right)$.

Moreover, by

$$P(x,\partial)u_{\nu}\rightarrow P(x,\partial)u, \quad\mbox{as }
\nu \rightarrow \infty$$ and

$$ \partial^j_nu_{\nu}(x',0)= g_{j,\nu}(x')\rightarrow
g_{j}(x'),\quad\quad \partial^j_nu_{\nu}(x',0)\rightarrow
\partial^j_nu(x',0),\mbox{ as } \nu \rightarrow \infty,$$ for
$j=0,1,\cdots,m-1$, we derive that $u\in
C^{\infty}\left(B_{\widetilde{r}}\right) $ solves the Cauchy
problem

\begin{equation*}
\begin{cases}
P(x,\partial)u=0, & \mbox{in } B_{\widetilde{r}}, \\
\\
\partial_n^j u=g_j, & \mbox{for } j=0,1,\cdots, m-1,\quad \mbox{ in } B'_{\widetilde{r}}.
\end{cases}
\end{equation*}

\medskip

\noindent\textbf{3.} Keep in mind that the definition of
hyperbolicity that we provided above, involves not only
the operator (or, more precisely, its symbol), but also the
direction $N$. Let us consider, for instance, the wave operator in space dimension
two 
$$P(\partial_t,\partial_{x_1}, \partial_{x_2})=\partial^2_t-\left(\partial^2_{x_1}+\partial^2_{x_2}\right).$$
The symbol of $P(\partial_t,\partial_{x_1}, \partial_{x_2})$ is the polynomial
$$P(\eta,\xi_1,\xi_2)=-\eta^2+\left(\xi_1^2+\xi^2_2\right),$$
which is hyperbolic w.r.t. the direction $(1,0,0)$, but \textbf{is not hyperbolic
w.r.t. the directions $(0,1,0)$ e $(0,0,1)$}.

Elliptic operators (with constant coefficients), as is checked
easily, are not hyperbolic with respect to any direction.
Hence, the Cauchy problem for the elliptic operators does not enjoy the
property of local solvability in $C^{\infty}$, nor,
at the light of what was shown in Section \ref{LM-24}, it may happen that there is a continuous dependence in $C^{\infty}$ (or from $C^{\infty}$
in $C^{m}$).

\newpage

\chapter{Well--posed problems. Conditional stability
}\label{Stab-condizionata}
\section{Introduction}\label{Stab-condizionata-intro}

In the previous Chapter, we introduced the notion of well-posed problem
in the sense of Hadamard and we observed that some Cauchy problems
are not well posed. In particular, we observed that in such problems may
fail some kind of continuous dependence of the solutions
with respect to the data. This phenomenon represents a serious
obstacle in the study of the problems originating from the
applications. Indeed, in these problems the measurements of the data
are, apart for trivial cases, affected by the errors of approximation
the effect of which must always be taken into account if the
theoretical results obtained have any reasonable application.

Moreover, one should not believe that the phenomena of noncontinuous dependence of the
solutions with respect to the data are present only in
particularly complicated situations as is the case of Cauchy  problems. Indeed, such phenomena are encountered even in the
approximate calculation of a derivative or, to put ourselves in an
"applicative" perspective, in the approximate calculation of the velocity from a given time law. Let us suppose that a certain object moves
with rectilinear motion with a time law
$$x=s(t)$$
and let us let us suppose that we are interested in determining its velocity $v(t)$.
As we know very well
$$v(t)=s'(t).$$
However, now,  let us suppose that we only have an approximation of the time law. Let $\varepsilon>0$ and let $s_{\varepsilon}$ be an
approximate measure of the time law of our object. For instance let us suppose that

\begin{equation}\label{stabilit-1}
\sup_{t\in [0,T]}\left\vert
s(t)-s_{\varepsilon}(t)\right\vert\leq\varepsilon,
\end{equation}
where $T>0$ is the time in which the motion is performed (initial time $0$). It would be desired that
$v_{\varepsilon}=s'_{\varepsilon}$ be itself an approximation of $v$.
 Nevertheless, this does not happen. As a matter of fact, let 
$$s_{\varepsilon}(t)=s(t)+\varepsilon\sin
\left(\varepsilon^{-2}t\right), \quad t\in [0,T].$$ Then
\eqref{stabilit-1} is satisfied, hence

\begin{equation}\label{stabilit-2}
\lim_{\varepsilon\rightarrow 0}\sup_{t\in [0,T]}\left\vert
s(t)-s_{\varepsilon}(t)\right\vert=0,
\end{equation}
on the other side

\begin{equation*}
v(t)-
v_{\varepsilon}(t)=s'(t)-s'_{\varepsilon}(t)=\varepsilon^{-1}\cos\left(\varepsilon^{-2}t\right)\nrightarrow
0, \quad\mbox{as } \varepsilon\rightarrow 0.
\end{equation*}
Thus, without additional information about the motion of the object, we cannot obtain an approximation of its velocity based only on an approximation of its time law.

The example we have just considered may be expressed more formally by saying that the operator

$$\frac{d}{dt}: X\rightarrow Y,$$ where
$$X=C^1\left([0,T]\right), \quad\mbox{equipped with norm } \left\Vert
\cdot\right\Vert_{C^0\left([0,T]\right)}$$ and
$$Y=C^0\left([0,T]\right), \quad\mbox{equipped with norm } \left\Vert
\cdot\right\Vert_{C^0\left([0,T]\right)}$$ is not continuous.

\textit{In what follows when we say that a problem is not well--posed,\index{not well--posed problem in the sense of Hadamard} or ill--posed,\index{ill--posed problem in the sense of Hadamard}
 in the sense of Hadarmard we will always mean (if we do not risk 
ambiguities) that although the existence and uniqueness of the
solutions occur, the problem does not enjoy the continuous dependence with respect
to the data. The continuous dependence should be considered with respect to
the topologies suggested by the same nature of the applied problem under investigation. 
}

The Functional Analysis is a rich repository of examples of
not well--posed problems in the sense specified above. For instance, it is known that.
if $X$ is a Banach space and $A\in \mathcal{L}(X)$ (where
$\mathcal{L}(X)$ denotes the space of linear and continuous operators
from $X$ into itself) is injective then

$$A^{-1}: \mathcal{R}(A)\rightarrow X,$$
where $$\mathcal{R}(A)=\left\{Au: u\in X \right\},$$ is continuous
if and only if $\mathcal{R}(A)$ is closed in $X$. In particular,
if $\mathcal{R}(A)$ is not closed, then the problem

\begin{equation}\label{stabilit-3}
A(u)=f,
\end{equation}
is not well--posed in the sense of Hadamard. This is the case of the integral operator

\begin{equation}\label{stabilit-4}
L^2(0,1)\ni u\rightarrow (Au)(t)=\int_0^tu(s)ds\in L^2(0,1).
\end{equation}
Similarly, if $X$ has not a finite dimension
 and $A$ is a compact and injective operator from
$X$ into itself, then again $A^{-1}$ is not
continuous and therefore \eqref{stabilit-3} is an ill-posed problem.

As it is well known, operator \eqref{stabilit-4} is compact.
More generally, if $k\in L^2([0,1]\times [0,1])$ then
the operator

\begin{equation}\label{stabilit-5}
L^2(0,1)\ni u\rightarrow (Ku)(t)=\int_0^1k(t,s)u(s)ds\in L^2(0,1),
\end{equation}
is a compact operator. Hence the integral equation
\begin{equation}\label{stabilit-6}
\int_0^1k(t,s)u(s)ds=f(t), \quad t\in [0,1],
\end{equation}
is certainly an ill-posed problem in the sense of Hadamard since,
even if it admits solutions they do not depend continuously (in
$L^2(0,1)$) by $f$.

\bigskip

As we have already observed, a problem for which there is no
continuous dependence with respect to the data, without further information,
cannot be treated practically. In order to treat it,
additional information is needed. This introduces the notion of
\textbf{conditionally well-posed problem} whose formal definition
formal is as follows.

\begin{definition}\label{stabilit-7}
	\index{Definition:@{Definition:}!- conditional well--posed problem@{- conditional well--posed problem}}
	\index{Definition:@{Definition:}!- well--posed problem problem in the sense of Tikhonov@{- well--posed problem in the sense of Tikhonov}}
Let $(X,d_1)$ e $(Y,d_2)$ be two metric spaces and let
$$A:X\rightarrow Y$$ be a map
from $X$ to $Y$. we say that the problem of
determining $u\in X$ such that
\begin{equation}\label{stabilit-8}
A(u)=f,
\end{equation}
where $f\in Y$, is a \textbf{condizional well--posed problem} or likewise \textbf{well--posed problem in the sense of Tikhonov} \cite{Tik} with respect to $\mathcal{K} \subset X$ provided that we have

\noindent (i) $A_{|\mathcal{K}}: \mathcal{K}\rightarrow Y$ is
injective,

\noindent (ii)
$\left(A_{|\mathcal{K}}\right)^{-1}:A(\mathcal{K})\rightarrow
\mathcal{K}$ is continuous.
\end{definition}

\medskip

\textbf{Remark.} Let us note that in the above Definition
\textit{there is no requirement for the existence of solutions} of the
problem $A(u)=f$. In essence, the problem that is considered is the following one.

Given $f\in A(\mathcal{K})$, determine $u$ such that

\begin{equation}\label{stabilit-9}
\begin{cases}
A(u)=f, \\
\\
u\in\mathcal{K}.\\
\end{cases}
\end{equation}
The introduction of the set $\mathcal{K}$ into the definition of a
well-posed problem in the sense of Tikhonov is equivalent to the introduction of a
\textbf{additional information} \index{additional information} or an \textbf{a priori information}\index{a priori information}
(as it is often said in the literature) to the problem under investigation. We will synthesize the above requirements by saying

\begin{equation}\label{stabilit-9n}
	\begin{aligned}
&\mbox{ Determine } u\in X\mbox{ such that } A(u)=f \\& \mbox{
\textit{with the a priori bound } } u\in \mathcal{K}.
\end{aligned}
\end{equation}

\medskip

\noindent The a priori information pertains to the specific character of the problem
under investigation and it is suggested by the applied nature of that
problem. To find \textbf{a stability estimate for problem} \eqref{stabilit-9} or \eqref{stabilit-9n} means to find an appropriate estimate of the modulus of continuity $\omega$ of
$\left(A_{|\mathcal{K}}\right)^{-1}$, which, we recall, is defined
by

$$\omega(\delta)=\sup\left\{d_1\left(A^{-1}(f_1),A^{-1}(f_2)\right): f_1, f_2\in
A(\mathcal{K}), \quad d_2\left(f_1,f_2\right)\leq \delta\right\}.$$
Obviously the best that one can do is to determine
exactly $\omega$, but very often to arrive at an accurate
asymptotic estimate of $\omega(\delta)$ as $\delta$ goes to $0$,
can be considered a satisfactory result for many
applications.

$\blacklozenge$

\bigskip

If the problem
$$A(u)=f,$$
is not well--posed in the sense of Hadamard and $A$ is a continuous and injective map , there always exist sets $\mathcal{K}$ for which problem \eqref{stabilit-8} is conditionally well--posed. The following Theorem holds true.

\begin{theo}[\textbf{Tikhonov}]\label{teor-Tichonov}
	\index{Theorem:@{Theorem:}!- Tikhonov@{- Tikhonov}}
Let $\mathcal{K}$ and $Y$ be two metric spaces. Let us assume
$\mathcal{K}$ is a compact. Moreover, let
$$F:\mathcal{K}\rightarrow Y$$
which satisfies

\smallskip

\noindent (i) $F$ is injective,

\noindent (ii) $F$ is continuous.

\noindent Then
$$F^{-1}:F(\mathcal{K})\rightarrow \mathcal{K},$$
is continuous.
\end{theo}
\textbf{Proof}. Let $C$ be a closed subset of
$\mathcal{K}$, let us prove that $$F(C)=\left(F^{-1}\right)^{-1}(C)$$
is closed in $Y$. We have

\begin{equation*}
\begin{aligned}
C\subset \mathcal{K}, C\mbox{ closed } \mathcal{K} &\Rightarrow
F(C) \mbox{ compact in } Y\Rightarrow\\& \Rightarrow F(C) \mbox{
closed in } Y.
\end{aligned}
\end{equation*}
Therefore $F$ is continuous. $\blacksquare$

\bigskip

\textbf{Example 1.} To illustrate what we have said so far, let us return
to the problem of calculating the derivative of a function and we consider a  simplified version of it, i.e. the following one: calculate the derivative of $f\in
C^1([0,1])$ which satisfies

\begin{equation}\label{stabilit-10}
f'(0)=f'(1)=0.
\end{equation}
Let us note that also this problem is ill--posed in $L^2(0,1)$.
As a matter of fact, let

$$f_n(t)=\frac{1}{\sqrt{n}}\cos \pi n t, \quad t\in [0,1],\quad n\in
\mathbb{N}.$$ Then \eqref{stabilit-10} is satisfied and
$$\left\Vert
f_n\right\Vert_{L^2(0,1)}=\frac{1}{\sqrt{2n}}\rightarrow 0,
\quad\mbox{as}\quad n\rightarrow \infty$$ and
$$\left\Vert
f'_n\right\Vert_{L^2(0,1)}=\frac{\pi\sqrt{n}}{\sqrt{2}}\rightarrow
+\infty, \quad\mbox{as } n\rightarrow \infty.$$

We reformulate the problem as a
functional equation. Let be

$$X=\left\{u\in C^0([0,1]): u(0)=u(1)=0 \right\},$$

$$Y=\left\{f\in C^1([0,1]): f'(0)=f'(1)=0 \right\}$$
and
$$T:X\rightarrow Y, \quad \left(T(u)\right):=\int_{0}^tu(s)ds.$$
Equip $X$ and $Y$ with the  $L^2(0,1)$ norm

 The above problem becomes: given $f\in Y$ determine $u\in X$ such that

\begin{equation}\label{stabilit-11}
T(u)=f,
\end{equation}
which has, trivially, a unique solution given by
$$u(t)=f'(t), \quad t\in [0,1].$$ We have seen above that
there is no continuous dependence of the solutions by the datum. On the other hand $T$ is continuous and, setting
$$\mathcal{K}=\left\{u\in X: u\in C^1([0,1]), \quad \left\Vert
u'\right\Vert_{L^2(0,1)}\leq E \right\}, $$ where $E$ is a positive number, since $\mathcal{K}$ is a compact of $X$, the problem
\begin{equation}\label{stabilit-12}
\begin{cases}
T(u)=f, \\
\\
u\in\mathcal{K},\\
\end{cases}
\end{equation}
is well--posed in the sense of Tikhonov. Hence, denoting by $\omega_E$ the
modulus of continuity of $\left(T_{|\mathcal{K}}\right)^{-1}$
we have

\begin{equation}\label{stabilit-13}
	\begin{aligned}
		\left\Vert T^{-1}(f_1)-T^{-1}(f_2)\right\Vert_{L^2(0,1)}&=\left\Vert
f'_1-f'_2\right\Vert_{L^2(0,1)}\leq \\&\leq \omega_E\left(\left\Vert
f_1-f_2\right\Vert_{L^2(0,1)}\right),
\end{aligned}
\end{equation}

\smallskip

\noindent for every $f_1,f_2 \in T(\mathcal{K})$.

An estimate of $\omega_E$ can be easily proved as follows.  If $f\in
C^2([0,1])$ and $f'(0)=f'(1)=0$ then integrating by parts and applying
the Cauchy--Schwarz inequality we get

\begin{equation*}
 \begin{aligned}
\int_0^1f'^2(t)dt&=\int_0^1f'(t)f'(t)dt=-\int_0^1f''(t)f(t)dt\leq
\\&\leq
\left(\int_0^1f''^2(t)dt\right)^{1/2}\left(\int_0^1f^2(t)dt\right)^{1/2}.
\end{aligned}
\end{equation*}
Hence
\begin{equation}\label{stabilit-14}
 \int_0^1f'^2(t)dt\leq
 \left(\int_0^1f''^2(t)dt\right)^{1/2}\left(\int_0^1f^2(t)dt\right)^{1/2}.
\end{equation}
Let now $f_1,f_2 \in T(\mathcal{K})$, satisfy
$$\left\Vert
f_1-f_2\right\Vert_{L^2(0,1)}\leq \varepsilon,$$ then by
\eqref{stabilit-14} we have

\begin{equation}\label{stabilit-15}
\left\Vert T^{-1}(f_1)-T^{-1}(f_2)\right\Vert_{L^2(0,1)}\leq
(2E\varepsilon )^{1/2}.
\end{equation}
Therefore
\begin{equation}\label{stabilit-16}
\omega_E(\varepsilon)\leq (2E\varepsilon )^{1/2}.
\end{equation}
In particular, for fixed $E$, we have

\begin{equation}\label{stabilit-17}
\omega_E(\varepsilon)=\mathcal{O}\left((\varepsilon
)^{1/2}\right)\quad \mbox{as } \varepsilon\rightarrow 0,
\end{equation}
the reader is invited to check that the exponent $1/2$ in the
estimate \eqref{stabilit-17} cannot be improved, however on
this kind of issue we will return to further in this
Chapter.

\section{Interpolation estimates for the derivatives of a function.}\label{stime-interp}

In this Section we will prove some estimates among functions and their
derivatives. These estimates can be considered as conditional stability estimates of some not well--posed problem. We start by the following.

\begin{prop}\label{prop1-s1}
If $f\in C^2([a,b])$, where $a,b\in \mathbb{R}$, $a<b$, then we have
\begin{equation}\label{s1-1}
	\begin{aligned}
	&\left\Vert f'\right\Vert_{L^{\infty}(a,b)}\leq\\& \leq
c_0\left(((b-a)^{-2}\left\Vert
f\right\Vert_{L^{\infty}(a,b)}+\left\Vert
f''\right\Vert_{L^{\infty}(a,b)}\right)^{1/2}\left\Vert
f\right\Vert_{L^{\infty}(a,b)}^{1/2},
\end{aligned}
\end{equation}
where $c_0\leq 8\sqrt{2}$ is a positive constant. 
\end{prop}
\textbf{Proof.} The proof of \eqref{s1-1} can
be reduced to the case where \\ $[a,b]=[0,1]$. To this aim it suffices to consider , instead of $f$, the function $$[0,1]\ni t\rightarrow
f\left(a+(b-a)t\right)\in \mathbb{R}.$$ Let us continue to denote by
$f$ this function. Let us fix $x\in \left[0,\frac{1}{2}\right]$ and
let $h\in \left(0,\frac{1}{2}\right]$. We have

\begin{equation}\label{s1-2}
f'(x)=\left(f'(x)-\frac{f(x+h)-f(x)}{h}\right)+\frac{f(x+h)-f(x)}{h}.
\end{equation}
The Lagrange Theorem implies that there exist  $\xi,\eta$ such that
$$x<\eta<\xi<x+h$$ and
\begin{equation*}\label{s1-3}
f'(x)-\frac{f(x+h)-f(x)}{h}=f'(x)-f'(\xi)=(x-\xi)f'(\eta).
\end{equation*}
Hence
\begin{equation}\label{s1-4}
\left\vert f'(x)-\frac{f(x+h)-f(x)}{h}\right\vert \leq h \left\Vert
f''\right\Vert_{L^{\infty}(0,1)}.
\end{equation}
By the just obtained inequality and by  \eqref{s1-2} we have, for every $x\in
\left[0,\frac{1}{2}\right]$ and for every $h\in
\left(0,\frac{1}{2}\right]$,

\begin{equation}\label{s2-1}
 \begin{aligned}
\left\vert f'(x)\right\vert&\leq h \left\Vert
f''\right\Vert_{L^{\infty}(0,1)}+h^{-1}\left(\left\vert
f(x+h)\right\vert+\left\vert f(x)\right\vert\right)\leq
\\&\leq
h \left\Vert f''\right\Vert_{L^{\infty}(0,1)}+2h^{-1}\left\Vert
f\right\Vert_{L^{\infty}(0,1)}.
\end{aligned}
\end{equation}
Instead, if $x\in \left[\frac{1}{2},1\right]$, then it suffices to replace 
\eqref{s1-2} by

\begin{equation*}\label{s1-2-24}
f'(x)=\left(f'(x)-\frac{f(x-h)-f(x)}{-h}\right)+\frac{f(x-h)-f(x)}{-h},
\end{equation*}
for every $h\in 
\left(0,\frac{1}{2}\right]$ and we obtain

\begin{equation}\label{s2-1n}
\left\vert f'(x)\right\vert\leq h \left\Vert
f''\right\Vert_{L^{\infty}(0,1)}+2h^{-1}\left\Vert
f\right\Vert_{L^{\infty}(0,1)},
\end{equation}
for every $h\in 
\left(0,\frac{1}{2}\right]$. By \eqref{s2-1} and
\eqref{s2-1n} we get
\begin{equation}\label{s2-2}
\left\Vert f'\right\Vert_{L^{\infty}(0,1)}\leq h \left\Vert
f''\right\Vert_{L^{\infty}(0,1)}+2h^{-1}\left\Vert
f\right\Vert_{L^{\infty}(0,1)},\quad \forall h\in
\left(0,\frac{1}{2}\right].
\end{equation}
Now, set
\begin{equation}\label{s2-3}
E=\left\Vert f''\right\Vert_{L^{\infty}(0,1)}\quad\mbox{e}\quad
\varepsilon=\left\Vert f\right\Vert_{L^{\infty}(0,1)}
\end{equation}
and let us determine the  minimum of the function (in $h$ variable) on the right--hand side of \eqref{s2-2}, i.e.

$$
\left(0,\frac{1}{2}\right]\ni h\rightarrow\Phi(h)=hE+2h^{-1}\varepsilon.$$
It turns out that if $\left(\frac{2\varepsilon}{E}\right)^{1/2}\leq
\frac{1}{2}$ then, for
$h=h_0:=\left(\frac{2\varepsilon}{E}\right)^{1/2}$,

\begin{equation}\label{s3-1}
\min_{\left[0,1/2\right]}\Phi=\Phi(h_0)=2\sqrt{2E\varepsilon},
\end{equation}
while, if $\left(\frac{2\varepsilon}{E}\right)^{1/2}\geq
\frac{1}{2}$ then
$$
\min_{\left[0,1/2\right]}\Phi=\Phi(1/2)=\frac{1}{2}E+4\varepsilon,$$
but, since in this case $8\varepsilon\geq E$, we get

\begin{equation}\label{s3-2}
\min_{\left[0,1/2\right]}\Phi\leq 8\varepsilon.
\end{equation}
By \eqref{s3-1}, \eqref{s3-2} and \eqref{s2-2} we have

\begin{equation*}
\left\Vert f'\right\Vert_{L^{\infty}(0,1)}\leq
2\sqrt{2E\varepsilon}+8\varepsilon\leq
8\sqrt{2}(E+\varepsilon)^{1/2}\varepsilon^{1/2}
\end{equation*}
and recalling \eqref{s2-3} we obtain \eqref{s1-1}.$\blacksquare$

\bigskip

\textbf{Remarks}

\noindent \textbf{1.} By Proposition \ref{prop1-s1} it follows
the equivalence  of the norms
\begin{equation*}
\left\Vert f\right\Vert_{C^2([a,b])}= \left\Vert
f\right\Vert_{L^{\infty}(a,b)}+(b-a) \left\Vert
f'\right\Vert_{L^{\infty}(a,b)}+(b-a)^{2}\left\Vert
f''\right\Vert_{L^{\infty}(a,b)},
\end{equation*}
and
\begin{equation*}
\left\Vert f\right\Vert= \left\Vert
f\right\Vert_{L^{\infty}(a,b)}+(b-a)^{2}\left\Vert
f''\right\Vert_{L^{\infty}(a,b)},
\end{equation*}
(reader chek: use the inequality
$2AB\leq A^2+B^2$).

\medskip

\noindent \textbf{2.} Inequality \eqref{s1-1} is a \textbf{stability estimate for the calculation of the first derivative} \index{stability estimate:@{stability estimate:}!- for the calculation of the first derivative@{- for the calculation of the first derivative}} provided we have the a priori information

\begin{equation}\label{s4-3}
(b-a)^{2}\left\Vert f''\right\Vert_{L^{\infty}(a,b)}\leq E.
\end{equation}
As a matter of fact, if

\begin{equation}\label{s4-4}
\left\Vert f\right\Vert_{L^{\infty}(a,b)}\leq \varepsilon,
\end{equation}
then

$$\left\Vert
f'\right\Vert_{L^{\infty}(a,b)}\leq\frac{c_0}{b-a}(E+\varepsilon)^{1/2}\varepsilon^{1/2}.$$
More precisely, setting

$$\mathcal{K}_E=\left\{f\in C^2([a,b]): (b-a)^{2}\left\Vert f''\right\Vert_{L^{\infty}(a,b)}\leq
E\right\}$$ and denoting by $\omega$ the modulus of continuity of
$$\mathcal{K}_E\ni f\rightarrow f'\in C^1([a,b]),$$ we have

\begin{equation}\label{s5-1}
\omega(\varepsilon)\leq
\frac{c_0}{b-a}(E+\varepsilon)^{1/2}\varepsilon^{1/2}, \quad \forall
\varepsilon>0.
\end{equation}
Also, we observe 

\begin{equation}\label{s5-2}
\omega(\varepsilon)\geq
\frac{1}{b-a}(E+\varepsilon)^{1/2}\varepsilon^{1/2}, \quad \forall
\varepsilon>0.
\end{equation}
In order to check \eqref{s5-2}, for the sake of brevity, let us consider the case $[a,b]=[0,1]$ and let us denote

\begin{equation}\label{s5-3}
f_{\varepsilon}(x)=
\frac{E\varepsilon}{E+\varepsilon}\sin\frac{x}{\sqrt{\varepsilon/(\varepsilon+E)}}.
\end{equation}
We have 
$$\left\Vert f_{\varepsilon}\right\Vert_{L^{\infty}(0,1)}\leq
\varepsilon,\quad \left\Vert
f''_{\varepsilon}\right\Vert_{L^{\infty}(0,1)}\leq E.$$ Hence
$f_{\varepsilon}\in \mathcal{K}_E$ and we have
$$\omega(\varepsilon)\geq \left\Vert
f'_{\varepsilon}\right\Vert_{L^{\infty}(0,1)}=(E+\varepsilon)^{1/2}\varepsilon^{1/2},$$
from which we get \eqref{s5-2}.

Inequality \eqref{s5-2} implies that the exponent  $\frac{1}{2}$ in stability estimate
\eqref{s1-1} is \textbf{optimal},\index{optimal stability estimate} that is it cannot be improved by a bigger exponent.  $\blacklozenge$

\bigskip

In what follows we extend Proposition
\ref{prop1-s1} to higher order derivatives. To this aim we need some
notations. Let $a\in \mathbb{R}$. We recall that the translation operator, $\tau_a$, is defined as
\begin{equation}\label{traslazione}
(\tau_a f)(x)=f(x-a)
\end{equation}
where $f$ is any one real variable function. Also we denote 

\begin{equation}\label{s6-1b}
\tau^{*}_a=\tau_{-a}.
\end{equation}
Let us denote by $I$ the identity operator. Moreover, if $h\in
\mathbb{R}$ let us denote by $\Delta_hf$ the difference operator 

\begin{equation}\label{s6-1c}
\left(\Delta_h f\right)(x)=\left(\left(\tau^{*}_h-I\right)
f\right)(x)=f(x+h)-f(x).
\end{equation}
Let  $\mathfrak{e}_0, \mathfrak{e}_1, \cdots \mathfrak{e}_j $,
$j\in \mathbb{N}_0$, be the polynomials 
$$\mathfrak{e}_0(x)=1,$$
$$\mathfrak{e}_1(x)=x,$$
$$\mathfrak{e}_2(x)=x(x-h),$$
$$\cdots$$
$$\mathfrak{e}_j(x)=x(x-h)\cdots\left(x-(j-1)h\right).$$
Let us note that $\mathfrak{e}_j,$ has degree $j$ for every $j\in
\mathbb{N}_0$. Let $h\neq 0$ and let us denote by $P_{n-1}$, $n\in
\mathbb{N}$, the \textbf{Newton interpolation polynomial centered at
$x_0$ and with degree $n-1$},\index{Newton interpolation polynomial} that is

\begin{equation}\label{s6-2}
P_{n-1}(x)=\sum_{j=0}^{n-1}\frac{\left(\Delta^j_h
f\right)(x_0)}{j!h^j}\mathfrak{e}_j(x-x_0).
\end{equation}
Let us notice that if $1\leq j\leq n-1$ and $s\in \left\{0,1, \cdots,
n-1\right\}$ then we have

\begin{equation}\label{s7-1}
\mathfrak{e}_j(sh)=
\begin{cases}
0, \quad \mbox{ for } 0\leq s\leq j-1, \\
\\
\binom{{s}}{{j}}j!h^j, \quad\mbox{ for } j\leq s\leq n-1.\\
\end{cases}
\end{equation}
We check \eqref{s7-1}. If $0\leq s\leq j-1$ then one of the factors of
$\mathfrak{e}_j(sh)$ vanishes. While, if $j\leq s\leq n-1$ (where
$1\leq j\leq n-1$) we have
\begin{equation*}
 \begin{aligned}
\mathfrak{e}_j(sh)&=(sh)(sh-h)\cdots\left(sh-(sh-(j-1)h)\right)=
\\&=
h^js(s-1)\cdots\left(h-(s-(j-1))\right)=\\&=h^j\binom{{s}}{{j}}j! \ .
\end{aligned}
\end{equation*}
Also, let us note that if $0\leq s\leq n-1$, we get
\begin{equation}\label{s7-2}
f(x_0+sh)=\left(\left(\tau^{*}_h\right)^{s}f\right)(x_0)
\end{equation}
and
\begin{equation}\label{s7-3}
\left(\tau^{*}_h\right)^{s}=\left(\Delta_h+I\right)^{s}=\sum_{j=0}^{n-1}\binom{{s}}{{j}}\Delta_h^j.
\end{equation}
Hence \eqref{s7-2} and \eqref{s7-3} yield

\begin{equation*}
 \begin{aligned}
f(x_0+sh)&=\sum_{j=0}^{n-1}\binom{{s}}{{j}}\left(\Delta_h^jf\right)(x_0)=
\\&=
\sum_{j=0}^{n-1}\frac{1}{j!h^j}\left(\Delta_h^jf\right)(x_0)\mathfrak{e}_j(x_0+sh-x_0)=\\&=P_{n-1}(x_0+sh).
\end{aligned}
\end{equation*}
All in all, we have
\begin{equation}\label{s8-1}
f(x_0+sh)=P_{n-1}(x_0+sh),\quad \mbox{per } s=0,1, \cdots, n-1.
\end{equation}

\bigskip

We can now state and prove the following
\begin{prop}\label{prop2-s8}
Let $n\geq 2$ and $f\in C^n([a,b])$, where $a,b\in \mathbb{R}$, $a<b$.
For $1\leq k\leq n-1$  we have
\begin{equation}\label{s8-2}
	\begin{aligned}
	&\left\Vert f^{(k)}\right\Vert_{L^{\infty}(a,b)}\leq\\&\leq
c_{k,n}\left(((b-a)^{-n}\left\Vert
f\right\Vert_{L^{\infty}(a,b)}+\left\Vert
f^{(n)}\right\Vert_{L^{\infty}(a,b)}\right)^{\frac{k}{n}}\left\Vert
f\right\Vert_{L^{\infty}(a,b)}^{1-\frac{k}{n}},
\end{aligned}
\end{equation}
where $c_{k,n}$ is a positive constant which depends on $k$ and
$n$ only.
\end{prop}
\textbf{Proof.} As in the proof of Proposition
\ref{prop1-s1} we may reduce to the case $[a,b]=[0,1]$. We begin to prove  \eqref{s8-2} when $k=n-1$. Let $f\in
C^n([0,1])$. Fix $x_0\in \left[0,\frac{1}{2}\right]$ and let
$h\in \left(0, \frac{1}{2(n-1)}\right]$ and set 
$$R(x)=f(x)-P_{n-1}(x), \quad x\in [0,1].$$ By \eqref{s8-1}
we have
$$R(x_0)=R(x_0+h)=\cdots=R\left(x_0+(n-1)h\right)=0$$ from which,
By repeatedly applying the Rolle Theorem, we have that there exists 
$\xi\in \left(x_0,x_0+(n-1)h\right)$ such that 

$$f^{(n-1)}(\xi)-P_{n-1}^{(n-1)}(\xi)=R^{(n-1)}(\xi)=0.$$ Therefore there exists $\xi\in \left(x_0,x_0+(n-1)h\right)$ such that

\begin{equation}\label{s9-1}
f^{(n-1)}(\xi)=P_{n-1}^{(n-1)}(\xi).
\end{equation}
On the other hand we have
$$P_{n-1}^{(n-1)}(x)=\frac{\mathfrak{e}_{n-1}^{(n-1)}(x-x_0)}{(n-1)!}\frac{\left(\Delta_h^{n-1}f\right)(x_0)}{h^{n-1}},$$
and $$\mathfrak{e}_{n-1}^{(n-1)}(x)=(n-1)! \ .$$ By the latter and by 
\eqref{s9-1} we have
\begin{equation}\label{s9-2}
f^{(n-1)}(\xi)=\frac{\left(\Delta_h^{n-1}f\right)(x_0)}{h^{n-1}}.
\end{equation}
Therefore

\begin{equation}\label{s9-3}
f^{(n-1)}(x_0)=\left(f^{(n-1)}(x_0)-f^{(n-1)}(\xi)\right)+\frac{\left(\Delta_h^{n-1}f\right)(x_0)}{h^{n-1}}
\end{equation}
and by the Lagrange Theorem, we get

\begin{equation}\label{s9-4}
	\begin{aligned}
\left|f^{(n-1)}(x_0)-f^{(n-1)}(\xi)\right|&\leq\left\Vert
f^{(n)}\right\Vert_{L^{\infty}(0,1)}\left|x_0-\xi\right|\leq\\&\leq
\left\Vert f^{(n)}\right\Vert_{L^{\infty}(0,1)}(n-1)h.
\end{aligned}
\end{equation}
Moreover
\begin{equation}\label{s10-1}
	\begin{aligned}
\left(\Delta_h^{n-1}f\right)(x_0)&=\left(\left(\tau^*_h-I\right)^{n-1}f\right)(x_0)=\\&=\sum_{j=0}^{n-1}\binom{{n-1}}{{j}}(-1)^{n-1-j}f(x_0+jh)\\
\end{aligned}
\end{equation}
therefore
\begin{equation}\label{s10-2}
\left|\left(\Delta_h^{n-1}f\right)(x_0)\right|\leq
\sum_{j=0}^{n-1}\binom{{n-1}}{{j}}\left|f(x_0+jh)\right|\leq
2^{n-1}\left\Vert f\right\Vert_{L^{\infty}(0,1)}.
\end{equation}
By \eqref{s9-3}, \eqref{s9-4} and \eqref{s10-2} we have

\begin{equation}\label{s10-3}
\left|f^{(n-1)}(x_0)\right|\leq \left\Vert
f^{(n)}\right\Vert_{L^{\infty}(0,1)}(n-1)h+\left(\frac{2}{h}\right)^{n-1}\left\Vert
f\right\Vert_{L^{\infty}(0,1)}.
\end{equation}
Applying \eqref{s10-3} to $f(1-x)$ we also obtain the estimate 
for $x_0\in \left[\frac{1}{2},1\right]$. All in all, we have
for $0<h\leq \frac{1}{2(n-1)}$,

\begin{equation}\label{s10-4}
\left\Vert f^{(n-1)}\right\Vert_{L^{\infty}(0,1)}\leq
E(n-1)h+\left(\frac{2}{h}\right)^{n-1}\varepsilon,
\end{equation}
here we set

$$E:=\left\Vert f^{(n)}\right\Vert_{L^{\infty}(0,1)},\quad
\varepsilon:=\left\Vert f\right\Vert_{L^{\infty}(0,1)}.$$ Now we find the minimum of the function
$$\left(0,\frac{1}{2(n-1)}\right]\ni\rightarrow\Phi(h)=E(n-1)h+\left(\frac{2}{h}\right)^{n-1}\varepsilon.$$
By elementary calculation we have

\begin{equation*}\label{s11-2}
\min_{\left(0,\frac{1}{2(n-1)}\right]}\Phi \leq c_n
(E+\varepsilon)^{1-\frac{1}{n}}\varepsilon^{\frac{1}{n}},
\end{equation*}
where $c_n$ depends on $n$ only. Hence

\begin{equation}\label{s12-1}
\left\Vert f^{(n-1)}\right\Vert_{L^{\infty}(0,1)} \leq c_n
\left(\left\Vert f^{(n)}\right\Vert_{L^{\infty}(0,1)}+\left\Vert
f\right\Vert_{L^{\infty}(0,1)}\right)^{1-\frac{1}{n}}\left\Vert
f\right\Vert_{L^{\infty}(0,1)}^{\frac{1}{n}}.
\end{equation}

Now, let $1\leq k\leq n-1$, by iteration of \eqref{s12-1} we obtain

\begin{equation*}
 \begin{aligned}
&\left\Vert f^{(k)}\right\Vert_{L^{\infty}(0,1)}\leq
c_{k+1}\left(\left\Vert
f^{(k+1)}\right\Vert_{L^{\infty}(0,1)}+\left\Vert
f\right\Vert_{L^{\infty}(0,1)}\right)^{\frac{k}{k+1}}\left\Vert
f\right\Vert_{L^{\infty}(0,1)}^{\frac{1}{k+1}} \leq
\\&\leq
c_{k+1}c_{k+2}^{\frac{1}{k+2}}\left(\left\Vert
f^{(k+2)}\right\Vert_{L^{\infty}(0,1)}+\left\Vert
f\right\Vert_{L^{\infty}(0,1)}\right)^{\frac{k}{k+2}}\left\Vert
f\right\Vert_{L^{\infty}(0,1)}^{\frac{2}{k+2}}\leq\\&\leq \cdots
\leq \\&\leq c_{k,n}\left(\left\Vert
f^{(n)}\right\Vert_{L^{\infty}(0,1)}+\left\Vert
f\right\Vert_{L^{\infty}(0,1)}\right)^{\frac{k}{n}}\left\Vert
f\right\Vert_{L^{\infty}(0,1)}^{1-\frac{k}{n}}.
\end{aligned}
\end{equation*}
where $c_{k,n}$ depends on $k$ and $n$ only. By the above obtained inequality, coming back $[a,b]$ we get \eqref{s8-2}. $\blacksquare$

\bigskip

Estimates like \eqref{s8-2} can be easily derived
for  $L^p$ norm, $1\leq p\leq \infty$. Let
\begin{equation}\label{s12-2}
\left\Vert
f\right\Vert_{L^{p}(a,b)}=\left(\int^b_a\left|f(x)\right|^p\right)^{1/p}.
\end{equation}
We have
\begin{prop}\label{prop3-s12}
Let $n\geq 2$ and $f\in C^n([a,b])$, where $a,b\in \mathbb{R}$, $a<b$.
For $1\leq k\leq n-1$ we have
\begin{equation}\label{s12-3}
	\begin{aligned}
&\left\Vert f^{(k)}\right\Vert_{L^{p}(a,b)}\leq\\&\leq
c_{k,n}\left(((b-a)^{-n}\left\Vert
f\right\Vert_{L^{p}(a,b)}+\left\Vert
f^{(n)}\right\Vert_{L^{p}(a,b)}\right)^{\frac{k}{n}}\left\Vert
f\right\Vert_{L^{p}(a,b)}^{1-\frac{k}{n}},
\end{aligned}
\end{equation}
where $c_{k,n}$ is a positive constant which depends on $k$ and
 $n$ only.
\end{prop}

\textbf{Proof.} Similarly to the proof of the previous Proposition, we may reduce to the case
$[a,b]=[0,1]$. Let us prove \eqref{s12-3} when $k=n-1$. Let $t\in
\left[0,\frac{1}{2}\right]$  $h\in \left(0,
\frac{1}{2(n-1)}\right]$, by \eqref{s9-3} (with $t$ replacing
$x_0$) we have

\begin{equation}\label{s13-1}
f^{(n-1)}(t)=\left(f^{(n-1)}(t)-f^{(n-1)}(\xi)\right)+\frac{\left(\Delta_h^{n-1}f\right)(t)}{h^{n-1}}.
\end{equation}
Set

\begin{equation*}
\widetilde{f^{(n)}}(\tau)=\left \{
\begin{array}{c}
f^{(n)}(\tau), \quad \mbox{ for } \tau\in [0,1], \\
\\
0, \quad\mbox{ for } \tau\notin [0,1].\\
\end{array}%
\right.
\end{equation*}
We have

\begin{equation*}
 \begin{aligned}
&\left|f^{(n-1)}(t)-f^{(n-1)}(\xi)\right|=
\left|\int^t_{\xi}f^{(n)}(\tau)d\tau\right| \leq
\\&\leq \int_t^{t+(n-1)h}
\left|f^{(n)}(\tau)\right|d\tau=\int_{\mathbb{R}}\left|\widetilde{f^{(n)}}(\tau)\right|\chi_{(0,(n-1)h)}(\tau-t)d\tau,
\end{aligned}
\end{equation*}
where $\chi_{(0,(n-1)h)}$ is the characteric function 
$(0,(n-1)h)$. Hence, by \eqref{s13-1}, for $t\in
\left[0,\frac{1}{2}\right]$ and $h\in \left(0,
\frac{1}{2(n-1)}\right]$, we have
\begin{equation}\label{s13-2}
\left|f^{(n-1)}(t)\right|\leq
\int_{\mathbb{R}}\left|\widetilde{f^{(n)}}(\tau)\right|\chi_{(0,(n-1)h)}(\tau-t)d\tau+\frac{\left(\Delta_h^{n-1}f\right)(t)}{h^{n-1}}.
\end{equation}
At this point we use the triangle inequality in  $L^p$ and the
Young inequality for convolutions:
$$\left\Vert
F\star G\right\Vert_{L^{p}(\mathbb{R})}\leq \left\Vert
F\right\Vert_{L^{p}(\mathbb{R})} \left\Vert
G\right\Vert_{L^{1}(\mathbb{R})}$$ where $F=\widetilde{f^{(n)}}$,
$G=\chi_{(0,(n-1)h)}$ and we get

\begin{equation*}
\left(\int^{1/2}_0\left|f^{(n-1)}(t)\right|^pdt\right)^{1/p}\leq(n-1)h\left\Vert
f^{(n)}\right\Vert_{L^{p}(0,1)}+\left(\frac{2}{h}\right)^{n-1}\left\Vert
f\right\Vert_{L^{p}(0,1)}.
\end{equation*}
Similarly, we have
\begin{equation*}
\left(\int^{1}_{1/2}\left|f^{(n-1)}(t)\right|^pdt\right)^{1/p}\leq(n-1)h\left\Vert
f^{(n)}\right\Vert_{L^{p}(0,1)}+\left(\frac{2}{h}\right)^{n-1}\left\Vert
f\right\Vert_{L^{p}(0,1)}.
\end{equation*}
Hence, by $h\in \left(0, \frac{1}{2(n-1)}\right]$, we have
\begin{equation*}
\left\Vert f^{(n-1)}\right\Vert_{L^{p}(0,1)}\leq(n-1)h\left\Vert
f^{(n)}\right\Vert_{L^{p}(0,1)}+\left(\frac{2}{h}\right)^{n-1}\left\Vert
f\right\Vert_{L^{p}(0,1)}.
\end{equation*}
From now on we proceed as in the proof of Proposition
\ref{prop2-s8}. $\blacksquare$

\bigskip

\textbf{Remark 4.} It can be proved (see exercise below) that the exponent $1-\frac{k}{n}$ of the estimate \eqref{s8-2} is optimal. Regarding the constant we report here, without a proof (we refer to \cite{Gorny}), the following sharp estimate ($k,m\in \mathbb{N}_0$, $m>0$) 
  
\begin{equation*}
\left\Vert f^{(k)}\right\Vert_{L^{\infty}(a,b)}\leq 4e^{2k}m^k \left\Vert f\right\Vert^{1-\frac{1}{m}}_{L^{\infty}(a,b)}M_{km}^{\frac{1}{m}},
\end{equation*}
where
\begin{equation*}
M_{km}=\max\left\{\frac{(km)!}{(b-a)^{nm}}\left\Vert f\right\Vert^{1-\frac{1}{m}}_{L^{\infty}(a,b)}, \left\Vert f^{(km)}\right\Vert_{L^{\infty}(a,b)} \right\}.
\end{equation*} $\blacklozenge$

\bigskip

\underline{\textbf{Exercise 1.}} Prove  the optimality of the exponent $1-\frac{k}{n}$
in inequalities \eqref{s8-2} and \eqref{s12-3}. (hint:
 note that inequalities \eqref{s8-2} and \eqref{s12-3} hold for complex--valued functions. After that, instead of trigonometric functions like \eqref{s5-3} use complex exponential). $\clubsuit$

\bigskip

\underline{\textbf{Exercise 2.}} 
\noindent (i) Let $0<\alpha<\beta\leq 1$. Prove that for every $f\in
C^{0,\beta}([a,b])$ we have
$$\left|f\right|_{\alpha,[a,b]}\leq C(b-a)^{-\alpha}\left[(b-a)^{\beta}\left|f\right|_{\beta,[a,b]}+\left\Vert
f\right\Vert_{L^{\infty}(a,b)}\right]^{\frac{\alpha}{\beta}}\left\Vert
f\right\Vert_{L^{\infty}(a,b)}^{1-\frac{\alpha}{\beta}},$$ where $C$
depends on $\alpha$ and $\beta$ only.

\noindent (ii) Let $0<\alpha\leq 1$. Prove that for every $f\in
C^{1,\alpha}([a,b])$ we have
$$\left\Vert f'\right\Vert_{L^{\infty}(a,b)}\leq C(b-a)^{-1}\left[(b-a)^{1+\alpha}\left|f'\right|_{\alpha,[a,b]}+\left\Vert
f\right\Vert_{L^{\infty}(a,b)}\right]^{\frac{1}{1+\alpha}}\left\Vert
f\right\Vert_{L^{\infty}(a,b)}^{\frac{\alpha}{1+\alpha}},$$ where  $C$
depends on $\alpha$ only.

\noindent  \textbf{hint to (ii):} note that instead of
\eqref{s1-4} we have
$$ \left\vert f'(x)-\frac{f(x+h)-f(x)}{h}\right\vert\leq
\left|f\right|_{\alpha,[a,b]}h^{\alpha}.$$ $\clubsuit$

\bigskip

We conclude this Section with two estimates for the derivatives of several variables
functions.

\begin{prop}\label{Stima1-c16}
Let $f\in C^2\left(\overline{B_1}\right)$. We have
\begin{equation}\label{c16-2}
\left\Vert \nabla f\right\Vert_{L^{\infty}(B_1)}\leq
c\left(\left\Vert \partial^2
f\right\Vert_{L^{\infty}(B_1)}+\left\Vert
f\right\Vert_{L^{\infty}(B_1)}\right)^{\frac{1}{2}}\left\Vert
f\right\Vert_{L^{\infty}(B_1)}^{\frac{1}{2}},
\end{equation}
where $c$ is a positive constant whhich depends on $n$ only.
\end{prop}
\textbf{Proof.} Let $h\in (0,1]$ be to choose later on and let $x\in
B_1$. Set

$$\Omega_h(x)=B_h(x)\cap B_1.$$
For $j=1,\cdots, n$ we have

\begin{equation}\label{c17-1}
 \begin{aligned}
f_{x_j}(x)&=
f_{x_j}(x)-\frac{1}{\left|\Omega_h(x)\right|}\int_{\Omega_h(x)}f_{x_j}(y)dy+\\&+\frac{1}{\left|\Omega_h(x)\right|}\int_{\Omega_h(x)}f_{x_j}(y)dy=
\\&=\frac{1}{\left|\Omega_h(x)\right|}\int_{\Omega_h(x)}\left(f_{x_j}(x)-f_{x_j}(y)\right)dy+\\&+\frac{1}{\left|\Omega_h(x)\right|}\int_{\partial\Omega_h(x)}f(y)\nu_jdS.
\end{aligned}
\end{equation}
We have, for a suitable  $\overline{x}$ on the segment of  extremes $x$
and $y$,

\begin{equation}\label{c17-2}
f_{x_j}(x)-f_{x_j}(y)=\nabla f_{x_j}\left(\overline{x}\right)\cdot (x-y),
\end{equation}
in addition we have
\begin{equation}\label{c17-3}
\left|\Omega_h(x)\right|\geq C_1 h^{n},\quad
\left|\partial\Omega_h(x)\right|\geq C_2 h^{n-1},
\end{equation}
where $C_1$ and $C_2$ depend on $n$ only. From what was obtained in
\eqref{c17-1}, \eqref{c17-2} and \eqref{c17-3} we get
$$\left\vert
\nabla f(x)\right\vert\leq C\left(h\left\Vert
\partial^2 f\right\Vert_{L^{\infty}(B_1)}+h^{-1}\left\Vert
f\right\Vert_{L^{\infty}(B_1)}\right),$$ where $C$ depends on
$n$ only. Now we minimize the function on the right--hand side of \eqref{c17-3} and
we obtain \eqref{c16-2}. $\blacksquare$

\bigskip

\begin{prop}\label{Stima2-c18}
Let $f\in C^1\left(\overline{B_1}\right)$. We have
\begin{equation}\label{c18-1}
\left\Vert f\right\Vert_{L^{\infty}(B_1)}\leq c\left(\left\Vert
\nabla f\right\Vert_{L^{\infty}(B_1)}+\left\Vert
f\right\Vert_{L^{2}(B_1)}\right)^{\frac{n}{n+2}}\left\Vert
f\right\Vert_{L^{2}(B_1)}^{\frac{n}{n+2}},
\end{equation}
where $c$ depends on $n$ only.
\end{prop}
\textbf{Proof.} Let $h\in(0,1]$ be to choose, $x\in B_1$
and $\Omega_h(x)$ like in the previous proof. We have
\begin{equation}\label{c18-2}
f(x)
=\frac{1}{\left|\Omega_h(x)\right|}\int_{\Omega_h(x)}\left(f(x)-f(y)\right)dy+\frac{1}{\left|\Omega_h(x)\right|}\int_{\Omega_h(x)}f(y)dy.
\end{equation}
By
\begin{equation*}
\left|f(x)-f(y)\right|\leq\left\Vert \nabla
f\right\Vert_{L^{\infty}(B_1)}|x-y|
\end{equation*}
we have
\begin{equation}\label{c18-3}
\frac{1}{\left|\Omega_h(x)\right|}\int_{\Omega_h(x)}\left|f(x)-f(y)\right|dy\leq
h\left\Vert \nabla f\right\Vert_{L^{\infty}(B_1)}.
\end{equation}
On the other hand, by the Cauchy--Schwarz inequality we have

\begin{equation}\label{c18-4}
	\begin{aligned}
		\left\vert\frac{1}{\left|\Omega_h(x)\right|}\int_{\Omega_h(x)}\left|f(y)\right|dy\right\vert&\leq
\frac{1}{\left|\Omega_h(x)\right|^{1/2}}\left\Vert
f\right\Vert_{L^{2}(B_1)}\leq \\&\leq \frac{1}{(c_1h^n)^{1/2}}\left\Vert
f\right\Vert_{L^{2}(B_1)}.
\end{aligned}
\end{equation}
Hence
\begin{equation}\label{c18-5}
\left\vert f(x)\right\vert\leq C\left(h^{-n/2}\left\Vert
f\right\Vert_{L^{2}(B_1)}+h\left\Vert \nabla
f\right\Vert_{L^{\infty}(B_1)}\right), \end{equation} where $C$
depends on $n$ only. Now we minimize the function on the right--hand side of
\eqref{c18-5} and we get \eqref{c18-1}. $\blacksquare$

\section{Stability estimates for the continuation of holomorphic functions}\label{stime-prolungamento-olo} 
\index{stability estimate:@{stability estimate:}!- for the continuation of holomorphic functions@{- for the continuation of holomorphic functions}}
In what follows we will identify
$\mathbb{C}$ with $\mathbb{R}^2$. Let us recall very quickly the
definition and some properties of the holomorphic functions.\index{holomorphic functions}

\bigskip

\noindent \textbf{1.} Let $\Omega$ be an open set of 
$\mathbb{C}$, $f:\Omega\rightarrow\mathbb{C}$ be a complex--valued function defined 
in $\Omega$ and $z_0\in \Omega$. We say that $f$
is \textbf{holomorphic in $z_0$} if the following 
limit there exists (in $\mathbb{C}$) 
\begin{equation}\label{s18-1}
\lim_{z\rightarrow z_0} \frac{f(z)-f(z_0)}{z-z_0}.
\end{equation}
In such a case we denote by $f'(z_0)$ the value of limit
\eqref{s18-1} and we say that it is the \textbf{derivative of $f$ in
$z_0$}. We say that $f$ is \textbf{holomorphic in $\Omega$} provided it is holomorphic
in each point of $\Omega$. For instance, $z, z^n$ are
holomorphic functions in $\mathbb{C}$ while $\Re z, \Im z,
\overline{z}$ are not.

\bigskip

\noindent \textbf{2.} From what we say in \textbf{1} it follows that
if $f$ is holomorphic in $z_0=x_0+iy_0$ then $f$, as a function
of the real variables $x$ and $y$, is differentiable in
$(x_0,y_0)$ and

\begin{equation}\label{s18-2}
\frac{\partial f}{\partial x}+i\frac{\partial f}{\partial
y}=0,\quad\mbox{in } (x_0,y_0).
\end{equation}
Denoting by $P=\Re f$, $Q=\Im f$,  \eqref{s18-2} we may write
\begin{equation}\label{s19-1}
\frac{\partial P}{\partial x}=\frac{\partial Q}{\partial
y},\quad\mbox{ } \frac{\partial P}{\partial y}=-\frac{\partial
Q}{\partial x}, \quad\mbox{in }(x_0,y_0).
\end{equation}
Equations \eqref{s18-2} and \eqref{s19-1} are known as the \textbf{Cauchy-Riemann equations (or condition)}.\index{equation:@{equation:}!- Cauchy--Riemann@{- Cauchy--Riemann}} By introducing the notations
\begin{equation}\label{s19-2}
\frac{\partial}{\partial
z}=\frac{1}{2}\left(\frac{\partial}{\partial
x}-i\frac{\partial}{\partial y}\right), \quad
\frac{\partial}{\partial
\overline{z}}=\frac{1}{2}\left(\frac{\partial}{\partial
x}+i\frac{\partial}{\partial y}\right),
\end{equation}
the Cauchy-Riemann conditions  can be written as

\begin{equation}\label{s19-3}
\frac{\partial f}{\partial \overline{z}}=0, \quad \mbox{in
}(x_0,y_0).
\end{equation}
Also we have, setting $dz=dx+idy$, $d\overline{z}=dx-idy$ and
by considering $f$ as a function of $z$ and $\overline{z}$
\begin{equation}\label{s19-3n}
df=\frac{\partial f}{\partial x}dx+ \frac{\partial f}{\partial
y}dy=\frac{\partial f}{\partial z}dz+ \frac{\partial f}{\partial
\overline{z}}d\overline{z}.
\end{equation}

\bigskip

\noindent \textbf{3.} \textbf{The Cauchy Theorem.} If
$f:\Omega\rightarrow\mathbb{C}$ is holomorphic in $\Omega$ then the differential form  $$f(z)dz=f(x+iy)dx+if(x+iy)dy$$ is locally exact. That is, for every $(x_0,y_0)\in \Omega$ there exist  $\delta>0$ and
 $$F:B_{\delta}(x_0,y_0)\rightarrow\mathbb{C},$$ $F$
differentiable in $B_{\delta}(x_0,y_0)$ such that 
$$\frac{\partial F}{\partial x}=f,\quad \frac{\partial F}{\partial
y}=if.$$

\bigskip

\noindent \textbf{4.} It can be proved that if $f\in
C^0(\Omega)$ is holomorphic in $\Omega\setminus L$, where $L$
is a straight line then $fdz$ is holomorphic in $\Omega$. In
particular, if $f\in C^0(\Omega)$ and $f$ is holomorphic in
$\Omega\setminus \{a\}$ where $a\in \Omega$, then $fdz$ is
locally exact. This in turn enables the proof of the

\noindent \textbf{Cauchy integral formula.} \index{Cauchy formula} Let $f$ be holomorphic in
$\Omega$. Let $a\in \Omega$ and $r>0$ satisfy
$\overline{B_r(a)}\subset \Omega$. Setting $\gamma(t)=a+re^{it}$,
$t\in [0,2\pi)$, we have

\begin{equation}\label{s21-1}
f(a)=\frac{1}{2\pi i}\int_{\gamma}\frac{f(z)}{z-a}dz.\end{equation}

\bigskip

\noindent \textbf{5.} The Cauchy formula implies that if $f$ is
holomorphic in $B_{\rho}$ then $f$ can be expanded in a power series in $B_{\rho}$, that is there exists $\left\{a_n\right\}_{n\geq
0}$, sequence of $\mathbb{C}$, such that

\begin{equation}\label{s21-2}
f(z)=\sum_{n=0}^{\infty}a_nz^n,\quad \forall z\in
B_{\rho}.\end{equation} Since a holomorphic function can be expanded in a power series in each point of an open set $\Omega$, we have
$f:\Omega\rightarrow\mathbb{C}$ is holomorphic in $\Omega$ if and only if
$f$ is a complex analitic function  in $\Omega$, i.e. if and only if
for every $a\in \Omega$ there exists $\delta$ such that
\begin{equation}\label{s21-3}
f(z)=\sum_{n=0}^{\infty}\frac{f^{(n)}(a)}{n!}(z-a)^n,\quad \forall
z\in B_{\delta}(a).\end{equation} The "if $\cdots$ then" part
of the equivalence follows immediately by the properties of
differentiability of power series and by \eqref{s18-1}.
Keep in mind, however, that the expression "analytic complex function" should not be confused with the expression "analytic complex--valued function " For instance $f(z,\overline{z})=z^2-\overline{z}^2$ is
complex-valued analytic function, but not  analytic complex function,
as it is not holomorphic.

\bigskip

\noindent \textbf{6.} From what we have said in \textbf{4}, we get
the converse of the Cauchy Theorem.  That is to say: if $f\in
C^0(\Omega)$ and $f(z)dz$ is locally exact in $\Omega$ then
$f$ is holomorphic in $\Omega$. As a matter of fact, if $f(z)dz$ is
locally exact in $\Omega$ then it has locally a
primitive hence,  there exists locally, $F\in C^1$ such that
$\frac{\partial F}{\partial x}=f,\quad \frac{\partial F}{\partial y}=if $ from which we have that $F$ satisfies the Cauchy-Riemann conditions, hence $F$ is holomorphic and $f=F'$, on the other hand since the
derivative of a complex analytic function is still a complex analytic function, hence holomorphic, $f=F'$ is holomorphic.

We also have that if $f\in C^0(\Omega)$ and $\overline{B_r(a)}\subset \Omega$ then

\begin{equation}\label{s22-1}
f^{(n)}(a)=\frac{n!}{2\pi
i}\int_{\gamma}\frac{f(z)}{(z-a)^{n+1}}dz,\quad\forall n\in
\mathbb{N}_0,\end{equation} where $\gamma(t)=a+re^{it}$, $t\in
[0,2\pi)$. Moreover, again by \eqref{s21-1}, we obtain
the \textbf{mean property} \index{mean property}
\begin{equation}\label{s22-1-45}
f(a)=\frac{n!}{2\pi }\int_{0}^{2\pi}f(a+re^{it})dt,\end{equation}
as soon as $\overline{B_r(a)}\subset \Omega$. From the mean property
it follows the

\medskip

\noindent \textbf{Maximum modulus principle.}\index{maximum modulus principle} Let $\Omega\subset
\mathbb{C}$ a bounded open set and \\ $f\in C^0(\overline{\Omega})$ be a holomorphic function in $\Omega$, then

$$\max_{\overline{\Omega}} |f|=\max_{\partial\Omega} |f|.$$ Moreover,
if $\Omega$ is connected and there exists $a\in \Omega$ such that $$|f(a)|=\max_{\overline{\Omega}} |f|$$ then $f$ is constant in $\Omega$.

\bigskip

\noindent \textbf{7.} Let us now return our attention
to the analyticity of holomorphic functions and let us recall what follows.
 If $$f:\Omega\rightarrow \mathbb{C},$$ is holomorphic in
$\Omega$, connected open set of $\mathbb{C}$, then if $a\in \Omega$
we have

\begin{equation}\label{s23-1}
f^{(n)}(a)=0,\quad \forall n\in \mathbb{N}_0\quad
\Longrightarrow\quad f\equiv 0\quad \mbox{ in }
\Omega.\end{equation} From which we have that, if  $D$ is nonempty open set contained in $\Omega$, then

\begin{equation}\label{s23-2}
f=0,\quad \quad \mbox{ in } D \quad \Longrightarrow\quad f\equiv
0\quad \mbox{ in } \Omega.\end{equation} Moreover, if $f$ does not vanish identically in $\Omega$ then \textbf{the set of zeros
of $f$ has no accumulation points in $\Omega$}. As a matter of fact, if $f$
does not vanish identically in $\Omega$ then \eqref{s23-1}
implies that for every $a\in \Omega$ there exists  $k\in \mathbb{N}_0$ such that
$$f^{(k)}(a)\neq 0.$$ Denoting by $k_0\in \mathbb{N}_0$ the minimum
of such $k$, \eqref{s21-3} gives
\begin{equation*}
f(z)=(z-a)^{k_0}\sum_{n=0}^{\infty}\frac{f^{(n)}(a)}{n!}(z-a)^{n-k_0}:=(z-a)^{k_0}\varphi(z),\quad
\forall z\in B_{\delta}(a),\end{equation*} where $\varphi(a)\neq 0$,
hence $f(z)\neq 0$ in $B_{\delta_1}(a)\setminus \{a\}$ for some $\delta_1>0$. The above point can also be
expressed in the following way: let $$\left\{z_n:n\in
\mathbb{N}\right\}$$ a \textbf{infinite set which has at least an accumulation point in $\Omega$}, then

\begin{equation}\label{s23-3}
f(z_n)=0,\quad \forall n\in \mathbb{N}\quad \Longrightarrow\quad
f\equiv 0\quad \mbox{ in } \Omega\end{equation} or also

\begin{equation}\label{s24-1}
f(z)=\mathcal{O}\left(|z-a|^N\right),\quad \forall N\in
\mathbb{N}\quad \Longrightarrow\quad f\equiv 0\quad \mbox{ in }
\Omega.\end{equation}

Let us notice that \eqref{s23-3} does not hold if $\left\{z_n:n\in
\mathbb{N}\right\}$ has  accumulations points on $\partial
\Omega$ only. Let us consider, for instance,  $$\Omega=\left\{x+iy:x>0, \quad
|y|<x\right\}$$ and let $$f(z)=e^{-\frac{1}{z}}.$$ We have
$f(z)=\mathcal{O}\left(|z|^N\right)$, for every $N\in \mathbb{N}$, but
$f\not\equiv 0$.

\bigskip

\noindent \textbf{8.} In this concluding part of this
summary, we prove

\begin{prop}\label{Prop-s24}
Let $\Omega$ be a connected open set of $\mathbb{C}$. Let us suppose that

\begin{equation}\label{s24-2}
\overset{o}{\overline{\Omega}}=\Omega.
\end{equation}
Let $z_0\in\partial\Omega$ and $\Gamma=\partial \Omega\cap
B_R(z_0)\neq \emptyset$, where $R>0$. In addition, Let $f\in
C^0(\Omega\cup\overline{\Gamma})$, $f$ be holomorphic in $\Omega$ which satisfies

\begin{equation}\label{s24-3}
f=0\quad\mbox{on } \Gamma,
\end{equation}
then $f\equiv 0$ in $\Omega$.
\end{prop}
\textbf{Proof.} In what follows we will need some
simple topological relationships that we will prove (for the convenience
of the reader) in the concluding part of the main
proof. By \eqref{s24-2} we have immediately

\begin{equation}\label{s25-1}
\partial \Omega=\overline{\Omega}\cap \overline{\left(\mathbb{C}\setminus
\overline{\Omega}\right)}.
\end{equation}
Now, let us fix $\delta\in \left(0,\frac{R}{4}\right)$. Since
$z_0\in
\partial\Omega$, we have by  \eqref{s25-1} 
\begin{equation}\label{s25-2}
B_{\delta}(z_0)\cap \Omega\neq\emptyset, \quad
B_{\delta}(z_0)\cap\left(\mathbb{C}\setminus
\overline{\Omega}\right)\neq \emptyset.
\end{equation}
Now, let $a$ and $b$ be such that 
$$b\in B_{\delta}(z_0)\cap \Omega, \quad a\in B_{\delta}(z_0)\cap\left(\mathbb{C}\setminus
\overline{\Omega}\right).$$ Since
$$|b-a|\leq |b-z_0|+|z_0-a|<2\delta<R-2\delta,$$ we get
$$b\in B_{R-2\delta}(a)\cap \Omega\subset  B_{R}(z_0)\cap
\Omega,$$ \textbf{in particular} $B_{R-2\delta}(a)\cap
\Omega\neq\emptyset$, (the second inclusion relationship follows
from the triangle inequality). Let us denote by
$$r=R-2\delta.$$

\begin{figure}\label{figura-ps25}
	\centering
	\includegraphics[trim={0 0 0 0},clip, width=11cm]{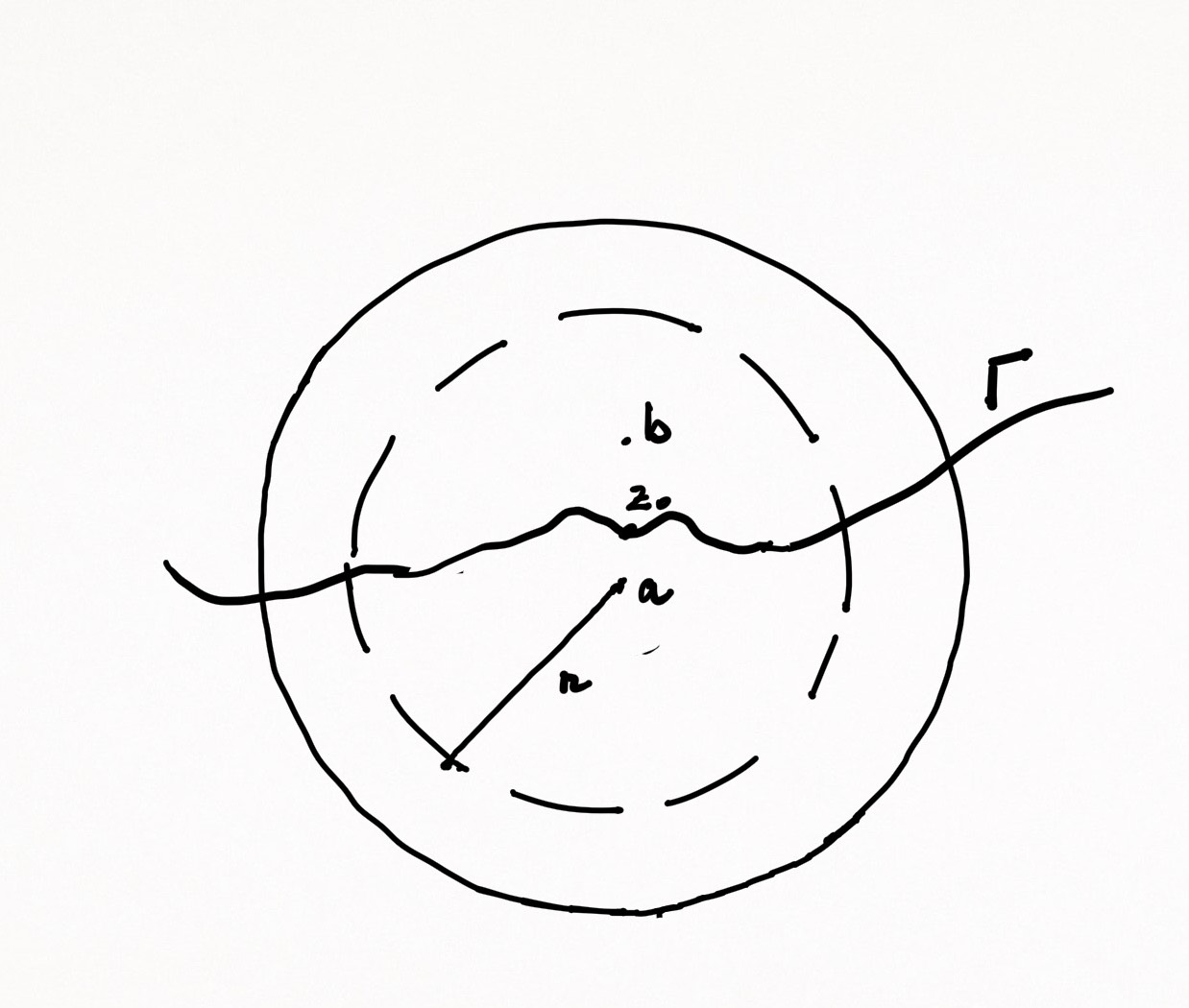}
	\caption{}
\end{figure}

We have also (see "concluding part")

\begin{equation}\label{s25-3}
\partial\left(B_r(a)\cap\Omega\right)\setminus \Gamma \subset
\left(\partial B_r(a)\right)\cap\overline{\Omega}.
\end{equation}
Set now

\begin{equation}\label{s26-1}
\varphi(z)=(z-a)^{-n}f(z),\quad n\in \mathbb{N}.
\end{equation}
It turns out that $\varphi$ is holomorphic in $B_r(a)\cap\Omega$ and continuous
in $\overline{B_r(a)\cap\Omega}$. By the maximum modulus principle we obtain

\begin{equation}\label{s26-2}
\left|\varphi(z)\right|\leq
\max_{\partial\left(B_r(a)\cap\Omega\right)}\left|\varphi\right|,\quad
\forall z\in B_r(a)\cap\Omega
\end{equation}
and by \eqref{s25-3} we get, recalling that $f=0$ on $\Gamma$,

\begin{equation}\label{s26-3}
 \begin{aligned}
\max_{\partial\left(B_r(a)\cap\Omega\right)}\left|\varphi\right|&\leq
\max_{\partial\left(B_r(a)\cap\Omega\right)\setminus
\Gamma}\left|\varphi\right|+\max_{\Gamma}\left|\varphi\right|
\leq\\&\leq
\max_{\partial\left(B_r(a)\right)\cap\overline{\Omega}}\left|\varphi\right|=\max_{\partial\left(B_r(a)\right)\cap\overline{\Omega}}\left(|z-a|^{-n}\left|f\right|\right)
=\\&=r^{-n}\max_{\partial\left(B_r(a)\right)\cap\overline{\Omega}}|f|.
\end{aligned}
\end{equation}
Now let $z\in B_r(a)\cap\Omega$, by \eqref{s26-1}, \eqref{s26-2} and
\eqref{s26-3} we derive
\begin{equation}\label{s26-4}
\left|f(z)\right|\leq \left(\frac{|z-a|}{r}\right)^n
\max_{\partial\left(B_r(a)\right)\cap\overline{\Omega}}\left|f\right|,
\quad  \forall n \in \mathbb{N}
\end{equation}
and passing to the limit as $n$ goes to infinity we deduce

$$f(z)=0, \quad \forall z\in B_r(a)\cap\Omega.$$
Since $B_r(a)\cap\Omega$ is a nonempty open set and $\Omega$
is a connected open set, by \eqref{s23-2} we have $f\equiv 0$ in
$\Omega$.

\medskip

\textbf{Concluding part of the proof.} We prove 
\eqref{s25-3}. Let us recall

\begin{equation}\label{s27-0}
\partial (A\cap B)\subset \partial A
\cup \partial B. \end{equation} 

Let now $$x\in \partial \left(B_r(a)\cap
\Omega\right)\setminus \Gamma,$$ we wish to prove that  
\begin{equation}\label{correct:15-4-23-2}
x\in
\left(\partial B_r(a)\right)\cap\overline{\Omega}.
\end{equation}
 First, we have trivially 
\begin{equation}\label{correct:15-4-23-1}
\partial \left(B_r(a)\cap
\Omega\right)\setminus \Gamma\subset \overline{\Omega}
\end{equation}

By \eqref{s27-0} we have

\begin{equation}\label{s27-1}
x\in \left(\partial B_r(a)\right) \cup \partial \Omega\quad \mbox{ and
}\quad x\notin \Gamma.
\end{equation}
Now, by \eqref{s27-1} we have that if $x\notin \partial B_r(a)$ then $x\in \partial
\Omega$. Moreover we have $x\in B_r(a)$. Because, if $x\notin B_r(a)$,
as $x\notin \partial B_r(a)$, we would have $x\notin
\overline{B_r(a)}$, hence it would exist  $\rho>0$ such that
$B_{\rho}(x)\cap B_r(a)=\emptyset$. Consequently, for such $\rho$ we would have $B_{\rho}(x)\cap \left(B_r(a)\cap \Omega\right)=\emptyset$
which contradicts \\ $x\in \partial\left(B_r(a)\cap \Omega\right)$.

All in all, if $x\notin \partial B_r(a)$ then
$$x\in \partial\Omega\cap B_r(a)\subset\partial\Omega\cap
B_R(z_0)=\Gamma,$$ But this cannot occur because, by \eqref{s27-1},
$x\notin \Gamma$. Therefore, \eqref{s27-1} implies that $$x\in
\partial B_r(a).$$ Finally, since \eqref{correct:15-4-23-1} holds we get \eqref{correct:15-4-23-2}, hence 
\eqref{s25-3} is proved. $\blacksquare$

\bigskip

\textbf{Remark.} Assumption
$\overset{o}{\overline{\Omega}}=\Omega$ excludes, for instance, that
$\Omega=B_1\setminus \left\{z_1, \cdots, z_n\right\}$ where $z_j\in
B_R$ for $j=1,\cdots,n$, where $R<1$. In this case, Proposition \ref{Prop-s24} does not hold for $R<1$ and
$\Gamma=\left\{z_1, \cdots, z_n\right\}$. $\blacklozenge$

\section[The Hadamard three circle inequality]{The Hadamard three circle inequality and other examples of stability estimates. }\label{esempi-stime-errore} In the previous Section we
focused on the unique continuation property for the \index{unique continuation property}
holomorphic functions. As we have seen it takes on several facets
corresponding to \eqref{s23-1}--\eqref{s24-1} and to Proposition \ref{Prop-s24}. In particular, we have that if a holomorphic function  $f:\Omega\rightarrow\mathbb{C}$
is known in a set $D\subset \Omega$ which  admits at least one accumulation point in $\Omega$ and if $\Omega$ is connected, then $f$
is uniquely determined in $\Omega$. The problem of determining
effectively the values of $f$ on $\Omega$ from $f_{|D}$ has an
interest in applications.
However this problem is not well posed
in the sense of Hadamard as can be inferred from the following simple
example.

\bigskip

\textbf{Example 1.}

\noindent Let $\Omega=B_1$ and $D=B_r$ where $r\in (0,1)$.
Then any holomorphic function on $B_1$ is uniquely determined by
$f_{|B_r}$. Nevertheless, small errors in the evaluation of
$f_{|B_r}$ may produce uncontrollable errors on $f$. Let
indeed

$$f_n(z)=\frac{1}{n}\left(\frac{z}{r}\right)^n, \quad n\in
\mathbb{N}.$$ We have

$$\max_{\overline{B_r}}\left\vert f_n\right\vert=\frac{1}{n}\rightarrow 0, \quad \mbox{as } n\rightarrow
\infty,$$ on the other hand, if $|z|>r$, then we have

$$\left\vert f_n(z)\right\vert=\frac{1}{n}\left(\frac{|z|}{r}\right)^n\rightarrow\infty, \quad \mbox{as } n\rightarrow
\infty.$$ $\spadesuit$

\bigskip

The conditional stability question for the analytic extension problem may be formulated as
follows

Let $\Omega$ be a connected open set of $\mathbb{C}$ and $D\subset
\Omega$ which has at least one accumulation point in $\Omega$. Let $f$
be any holomorphic function in $\Omega$, continuous on
$\overline{\Omega}$, which satisfies

\begin{equation}\label{s30-1}
\max_{\overline{\Omega}}\left\vert f\right\vert\leq E
\end{equation}
and
\begin{equation}\label{s30-2}
\max_{\overline{D}}\left\vert f\right\vert\leq \varepsilon.
\end{equation}
We are interested in finding a stability estimate like the following one

\begin{equation}\label{s30-3}
\left\vert f(z)\right\vert\leq E\eta\left(\frac{\varepsilon}{E};
z\right), \quad\forall z\in \Omega,
\end{equation}
where
$$\eta\left(s;
z\right)\rightarrow 0 \quad\mbox{as } s\rightarrow 0, \mbox{ 
 } \forall z\in \Omega.$$

\bigskip

There is a fairly general treatment of estimates of
stability \eqref{s30-3}, but here we will examine only a few examples
that are particularly significant.

\bigskip

\textbf{Example 2: The Hadamard three circle inequality.}

\noindent Let $0<r<\rho<R$. Let $f$ be a holomorphic function in
$B_R$ and continuous in $\overline{B_R}$. Let us denote by

\begin{equation}\label{s30-4}
M(s):=\max_{\overline{B_s}}\left\vert f\right\vert,
\quad\mbox{for } 0<s\leq R.
\end{equation}
We have

\begin{equation}\label{s30-5}
M(\rho)\leq (M(r))^{\theta_0}(M(R))^{1-\theta_0},
\end{equation}
where
\begin{equation}\label{s30n-5}
\theta_0=\frac{\log\frac{R}{\rho}}{\log\frac{R}{r}}.
\end{equation}

\textbf{Proof of \eqref{s30-5}.}

\noindent Let $n$ and $m$ be two integer numbers, $m>0$. Let us consider
the function

\begin{equation}\label{s31-1}
F(z)=z^{-n}\left(f(z)\right)^m, \quad\mbox{for } z\in
B_R\setminus \{0\}.
\end{equation}
$F$ is holomorphic in $B_R\setminus \{0\}$ and it is continuous in
$\overline{B_R}\setminus B_r$. We can  apply the maximum modulus principle. Set
\begin{equation*}
\widetilde{M}(s):=\max_{\partial B_s}\left\vert f\right\vert,
\quad\mbox{for } 0<s\leq R.
\end{equation*}
We have, for any $\rho\in (r,R)$,
\begin{equation*}
 \begin{aligned}
\rho^{-n}\left(\widetilde{M}(\rho)\right)^{m}&=\max_{\partial
B_{\rho}}\left|F\right|\leq \\&\leq \max\left\{\max_{\partial
B_{r}}\left|F\right|,\max_{\partial
B_{R}}\left|F\right|\right\}=\\&=\max\left\{r^{-n}\left(\widetilde{M}(r)\right)^{m},R^{-n}\left(\widetilde{M}(R)\right)^{m}\right\},
\end{aligned}
\end{equation*}
which gives

\begin{equation}\label{s31-3}
\widetilde{M}(\rho)\leq
\max\left\{\left(\frac{\rho}{r}\right)^{\frac{n}{m}}\widetilde{M}(r),\left(\frac{\rho}{R}\right)^{\frac{n}{m}}\widetilde{M}(R)\right\}.
\end{equation}
Since $\mathbb{Q}$ is dense in $\mathbb{R}$, by
\eqref{s31-3} we have

\begin{equation}\label{s31-4}
\widetilde{M}(\rho)\leq
\max\left\{\left(\frac{\rho}{r}\right)^{\alpha}\widetilde{M}(r),\left(\frac{\rho}{R}\right)^{\alpha}\widetilde{M}(R)\right\},\quad\mbox{ } \forall \alpha\in \mathbb{R}.
\end{equation}
Now, let us choose  $\alpha$ in such a way that
$$\left(\frac{\rho}{r}\right)^{\alpha}\widetilde{M}(r)=\left(\frac{\rho}{R}\right)^{\alpha}\widetilde{M}(R),$$
that is, let

$$\alpha=\frac{\log\left(\frac{\widetilde{M}(R)}{\widetilde{M}(r)}\right)}{\log\frac{R}{\rho}}$$
and  \eqref{s31-4} implies

\begin{equation}\label{s32-1}
\widetilde{M}(\rho)\leq
\left(\widetilde{M}(r)\right)^{\theta_0}\left(\widetilde{M}(R)\right)^{1-\theta_0},
\end{equation}
where $\theta_0$ is given by \eqref{s30n-5}. Finally, by the maximum modulus principle, we get \eqref{s30-5}. $\blacksquare$

\begin{figure}\label{trepalle}
	\centering
	\includegraphics[trim={0 0 0 0},clip, width=10cm]{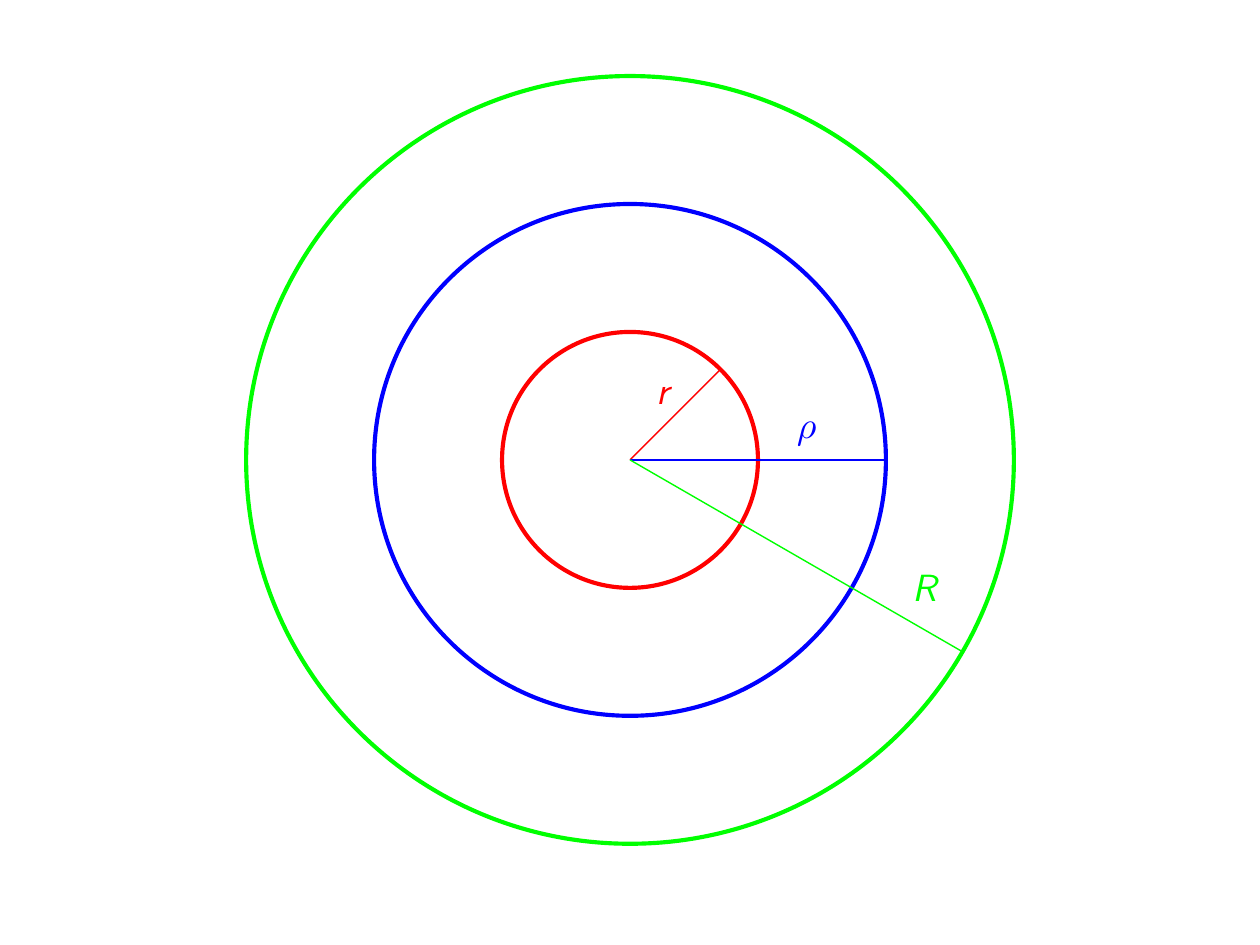}
	\caption{}
\end{figure}

\bigskip

\textbf{Remarks}

\noindent \textbf{1.} If $f$ is holomorphic in $B_R\setminus
\overline{B_r}$ and it is continuous in $\overline{B_R}\setminus
B_r$, inequality \eqref{s32-1} still applies.

\medskip

\noindent \textbf{2.} It is evident that  \eqref{s30-5} is a stability estimate for the problem:

Determine $f\in C^0\left(\overline{B_R}\right)$, $f$ holomorphic 
in $B_R$ which satisfies

$$\max_{\overline{B_r}}\left\vert f\right\vert\leq \varepsilon$$
and

$$\max_{\overline{B_R}}\left\vert f\right\vert\leq E.$$

\medskip

\noindent \textbf{3.} Let us notice that inequality \eqref{s30-5}
is equivalent to the convexity of the function

$$t\rightarrow \log M(e^t).$$

\medskip

\noindent \textbf{4.} The inequality \eqref{s30-5} cannot be
improved. More precisely, the following facts apply.

For every $C>0$ independent of $f$ we have
\begin{equation}\label{s32-a}
\theta_0=\sup\left\{\theta: M(\rho)\leq C\left(M(r)\right)^{\theta},
\quad M(R)=1\right\}.
\end{equation}
In other words, the exponent $\theta_0$ in \eqref{s30-5} is the best exponent.
Moreover
\begin{equation}\label{s32-b}
\inf\left\{C>0: M(\rho)\leq C\left(M(r)\right)^{\theta_0}, \quad
M(R)=1\right\}=1,
\end{equation}
that is the constant $1$ in \eqref{s30-5} is the best constant.

\bigskip

\textbf{Proof of \eqref{s32-a}}

\noindent It suffices to prove that if $\theta\in \mathbb{R}$ satisfies

\begin{equation}\label{s33-1}
M(\rho)\leq C\left(M(r)\right)^{\theta},
\end{equation}
for every $f\in C^0\left(\overline{B_R}\right)$, $f$ holomorphic in
$B_R$ and such that $M(R)=1$ then

\begin{equation}\label{s33-2}
\theta\leq \theta_0.
\end{equation}
Now, let

$$f_n(z)=\left(\frac{z}{R}\right)^n,\quad n\in \mathbb{N}.$$
We have

$$M(\rho)=\left(\frac{\rho}{R}\right)^n, \quad
M(r)=\left(\frac{r}{R}\right)^n$$ and by \eqref{s33-1} we have

$$n\log\frac{\rho}{R}\leq \log C+\theta n \log\frac{r}{R},\quad
\forall n\in \mathbb{N},$$ from which (recalling that $r<\rho<R$) we have

$$\frac{\log (\rho/R)}{\log(r/R)}\geq \theta + \frac{\log
C}{n\log(r/R)},\quad n\in \mathbb{N}$$ and passing to the limit as 
$n\rightarrow \infty$ we obtain \eqref{s33-2}.

\bigskip

\textbf{Proof of \eqref{s32-b}}

\noindent It suffices to prove that if $C>0$ satisfies

\begin{equation*}
M(\rho)\leq C\left(M(r)\right)^{\theta_9},
\end{equation*}
for every $f\in C^0\left(\overline{B_R}\right)$, $f$ holomorphic in
$B_R$ and $M(R)=1$ then

\begin{equation}\label{s33-cost}
C\geq 1.
\end{equation}
It suffices to choose $$f(z)=\frac{z}{R}$$ and we have trivially

$$\left(M(r)\right)^{\theta_0}=\left(\frac{r}{R}\right)^{\theta_0}=\frac{\rho}{R}=M(\rho),$$
from which \eqref{s33-cost} follows.

\medskip

\noindent \textbf{5.} It is interesting to note that the mere
inequality \eqref{s30-5} implies the following 
unique continuation property \index{unique continuation property}

\begin{equation}\label{s34-1}
f(z)=\mathcal{O}\left(|z|^N\right), \  \mbox{as }
z\rightarrow 0, \ \ \forall N\in \mathbb{N}\Longrightarrow f\equiv 0 \mbox{ in } B_R.
\end{equation}
indeed, let us assume that

\begin{equation}\label{s34-2}
f(z)=\mathcal{O}\left(|z|^N\right),   \mbox{ as }
z\rightarrow 0, \ \ \forall N\in \mathbb{N}
\end{equation}
and, arguing by contradiction, let us suppose that
\begin{equation}\label{s34n-2}
f\not\equiv 0\quad \mbox{in } B_R.
\end{equation}
Then there exists $\rho\in (0,R)$ such that
\begin{equation}\label{s34nn-2}
M(\rho)>0.\end{equation} On the other hand,  \eqref{s34-2} implies 

$$M(r)\leq C_N\left(\frac{r}{R}\right)^N, \ \ \forall N\in \mathbb{N}$$
for some constant $C_N$ (independent of $r$).
By this inequality and by \eqref{s30-5} we have

\begin{equation*}
 \begin{aligned}
&\frac{M(\rho)}{M(R)}\leq
\left(\frac{M(r)}{M(R)}\right)^{\theta_0}\leq\\&\leq
\left(C_N\left(\frac{r}{R}\right)^N\right)^{\theta_0}=\\&=\exp\left\{\frac{\log
R/\rho}{\log R/r}\left[-N\log\frac{R}{r}+\log C_N \right]\right\},
\quad \forall r\in (0,\rho), \forall N\in \mathbb{N}
\end{aligned}
\end{equation*}
and passing to the limit as $r\rightarrow 0$, we have

\begin{equation}\label{s35-1}
\frac{M(\rho)}{M(R)})\leq \exp \left[-N \log
\frac{R}{\rho}\right],\quad \forall N\in \mathbb{N}
\end{equation}
now, again passing to the limit as $N\rightarrow \infty$ we have
$$M(\rho)=0$$
which contradicts \eqref{s34nn-2}. Hence $f\equiv 0$ in $B_R$.

\medskip

\noindent \textbf{6.} We can easily prove an
inequality in $L^2$ similar to \eqref{s30-5}. More
precisely the following inequalities ($0<r<\rho<R$) hold true

\begin{equation}\label{s36-1}
	\begin{aligned}
&\int^{2\pi}_0\left|f\left(\rho e^{i\phi}\right)\right|^2d\phi\leq \\&\leq 
\left(\int^{2\pi}_0\left|f\left(r
e^{i\phi}\right)\right|^2d\phi\right)^{\theta_0}\left(\int^{2\pi}_0\left|f\left(R
e^{i\phi}\right)\right|^2d\phi\right)^{1-\theta_0},
\end{aligned}
\end{equation}

\begin{equation}\label{s36n-1}
\int_{B_{\rho}}\left|f\right|^2dxdy\leq
\left(\int_{B_{r}}\left|f\right|^2dxdy\right)^{\theta_0}\left(\int_{B_{R}}\left|f\right|^2dxdy\right)^{1-\theta_0}.
\end{equation}

\textbf{Proof of \eqref{s36-1}}.

First, let us observe that 

\begin{equation}\label{s36-2}
\rho=r^{\theta_0}R^{1-\theta_0}.
\end{equation}

By the assumption on $f$ we have

$$f(z)=\sum_{n=0}^{\infty} a_n z^n, \quad |z|\leq R.$$

Hence, by  \eqref{s36-2} and by the H\"{o}lder inequality
we have

\begin{equation*}
 \begin{aligned}
&\int^{2\pi}_0\left|f\left(\rho e^{i\phi}\right)\right|^2d\phi=
\sum_{n=0}^{\infty}\rho^{2n}|a_n|^2=\\&=
\sum_{n=0}^{\infty}\left(r^{2n}|a_n|^2\right)^{\theta_0}\left(R^{2n}|a_n|^2\right)^{1-\theta_0}\leq\\&\leq
\left(\sum_{n=0}^{\infty}r^{2n}|a_n|^2\right)^{\theta_0}\left(\sum_{n=0}^{\infty}R^{2n}|a_n|^2\right)^{1-\theta_0}=\\&=
\left(\int^{2\pi}_0\left|f\left(r
e^{i\phi}\right)\right|^2d\phi\right)^{\theta_0}\left(\int^{2\pi}_0\left|f\left(R
e^{i\phi}\right)\right|^2d\phi\right)^{1-\theta_0}.
\end{aligned}
\end{equation*}

Now, \textbf{let us prove \eqref{s36n-1}}. We have

\begin{equation*}
	\begin{aligned}
\int_{B_{\rho}}\left|f\right|^2dxdy&=\int_{0}^{\rho} s\left(\int_0^{2\pi}\left|f\left(s
e^{i\phi}\right)\right|^2d\phi\right)ds=\\&=\int_{0}^{1}
t\rho\left(\int_0^{2\pi}\left|f\left(t\rho
e^{i\phi}\right)\right|^2d\phi\right)dt.
\end{aligned}
\end{equation*} 
From which, by using 
\eqref{s36-1}, \eqref{s36-2} and by H\"{o}lder inequality,
we get
\begin{equation*}
	\begin{aligned}
&\int_{B_{\rho}}\left|f\right|^2dxdy\leq \\&\leq \rho
\left\{\int_0^1\left[\int^{2\pi}_0 t\left|f\left(tr
e^{i\phi}\right)\right|^2d\phi\right]^{\theta_0}\left[\int^{2\pi}_0t\left|f\left(t
R e^{i\phi}\right)\right|^2d\phi\right]^{1-\theta_0}dt\right\}\leq\\&
\leq\rho\left[\int_0^1\int^{2\pi}_0 t\left|f\left(tr
e^{i\phi}\right)\right|^2d\phi
dt\right]^{\theta_0}\left[\int_0^1\int^{2\pi}_0 t\left|f\left(tR
e^{i\phi}\right)\right|^2d\phi dt\right]^{1-\theta_0}=\\&
=\left[\int_0^1\int^{2\pi}_0 tr\left|f\left(tr
e^{i\phi}\right)\right|^2d\phi
dt\right]^{\theta_0}\left[\int_0^1\int^{2\pi}_0 tR\left|f\left(tR
e^{i\phi}\right)\right|^2d\phi dt\right]^{1-\theta_0}=\\&
=\left(\int_{B_{r}}\left|f\right|^2dxdy\right)^{\theta_0}\left(\int_{B_{R}}\left|f\right|^2dxdy\right)^{1-\theta_0}.
\end{aligned}
\end{equation*}
One can also prove $L^p$ versions of the inequalities
\eqref{s36-1} and \eqref{s36n-1}, for these we refer the interested reader
to \cite[Ch. 1]{Du}. $\spadesuit$

\bigskip

\underline{\textbf{Exercise.}} Let $u$ be a harmonic function in $B_R\subset \mathbb{R}^2$ such that \\ $u\in C^0\left(\overline{B}_R\right)$. Prove that if $0<r<\rho<R$ then the following inequality holds true 

\begin{equation}\label{3-12-22-s36n-1-0}
	\int_{\partial B_{\rho}}u^2dS\leq
	\left(\int_{\partial B_{r}}u^2dS\right)^{\theta_0}\left(\int_{\partial B_{R}}u^2dS\right)^{1-\theta_0},
\end{equation}

\begin{equation}\label{3-12-22-s36n-1}
	\int_{B_{\rho}}u^2dxdy\leq
	\left(\int_{B_{r}}u^2dxdy\right)^{\theta_0}\left(\int_{B_{R}}u^2dxdy\right)^{1-\theta_0},
\end{equation}
where $\theta_0$ is given by \eqref{s30n-5}.
[Hint: recall the solution formula for Dirichlet problem in polar coordinates
$$u(\varrho,\phi)=\frac{a_0}{2}+\sum_{n=1}^{\infty}\varrho^n\left(a_n\cos n\phi+b_n\sin n\phi\right)$$ and apply it to obtain \eqref{3-12-22-s36n-1}]. $\clubsuit$

\bigskip

\bigskip  

\textbf{Example 3: Stability estimate on the bisector of an angle (\cite{Car0}).}

Let $S$ be a bounded open set  of $\mathbb{C}$ whose
boundary is is made up of two segments, $OA$ e $OB$ such that
$\widehat{AOB}=\pi\alpha$, $0<\alpha<2$ and by a Jordan curve 
$\Gamma$ of extremes $A$ and $B$ 
Let $z_0\in S$ and let us assume that $z_0$ belongs to the  bisector
of the angle $\widehat{AOB}$. Let $f\in
C^0\left(\overline{S}\right)$ be holomorphic in $S$. Let us denote by

$$E=\max_{\overline{S}}|f|,\quad
\varepsilon=\max_{\overline{\Gamma}}|f|.$$ Then

\begin{equation}\label{s38-1}
|f(z_0)|\leq
E^{1-\left(\frac{|z_0|}{R}\right)^{1/\alpha}}\varepsilon^{\left(\frac{|z_0|}{R}\right)^{1/\alpha}},
\end{equation}
where $R$ denotes the diameter of $S$.

\medskip

\textbf{Proof of \eqref{s38-1}.}

 Let $\sigma>0$ be to choose and
$$F(z)=f(z)\exp\sigma\left(\frac{z}{z_0}\right)^{1/\alpha}.$$
Set $|z_0|=r$. We have, for $z=\rho e^{\pm\frac{i\alpha\pi}{2}}$,
\begin{figure}\label{figura-ps38}
	\centering
	\includegraphics[trim={0 0 0 0},clip, width=9cm]{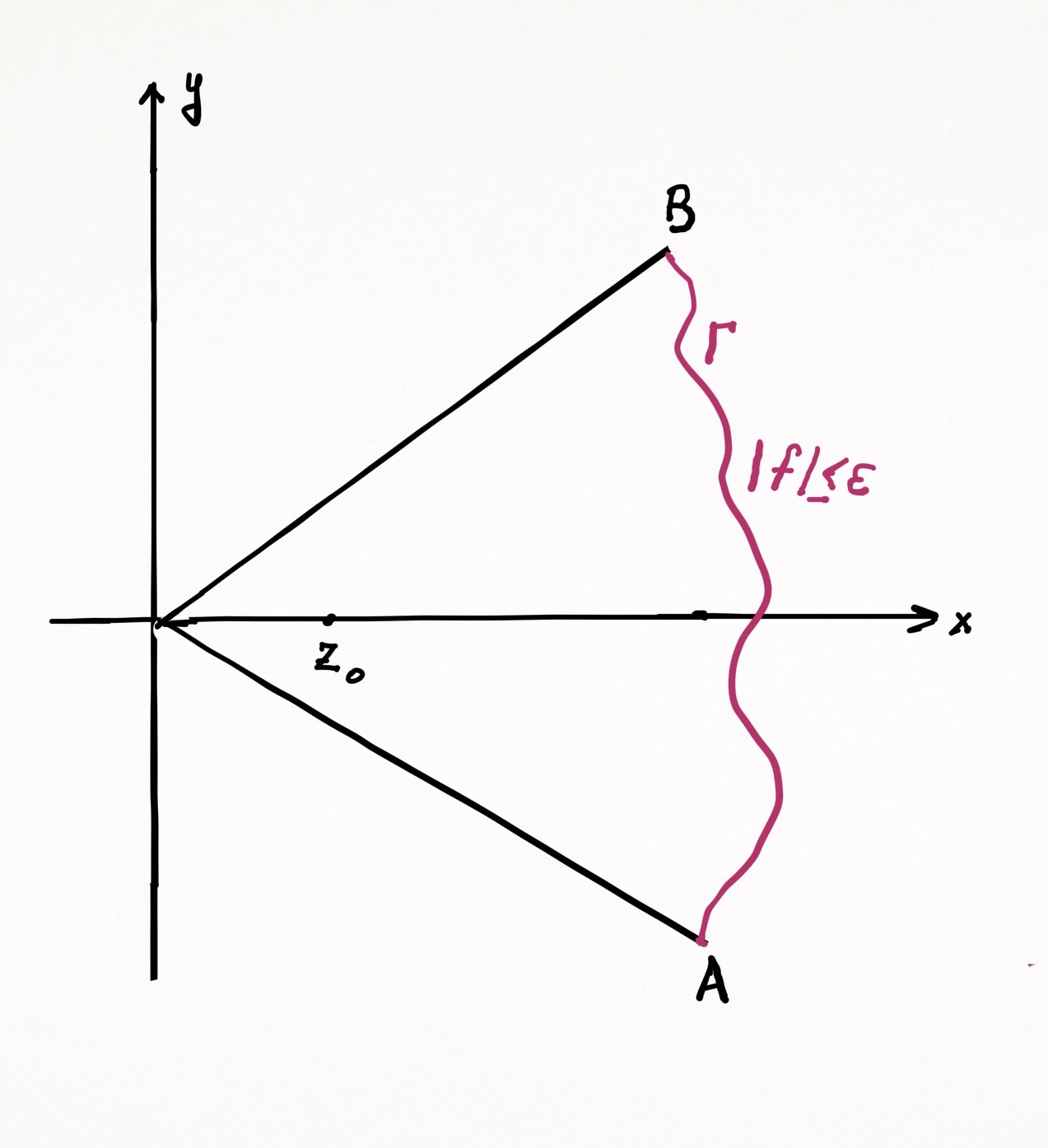}
	\caption{}
\end{figure}

$$\left|F\left(\rho e^{\pm\frac{i\alpha\pi}{2}}\right)\right|=\left|f\left(\rho e^{\pm\frac{i\alpha\pi}{2}}\right)\right|.$$ Hence
$$|F(z)| \leq E,\quad \mbox{on } OA \mbox{ and } OB.$$
Moreover
$$|F(z)| \leq \varepsilon \exp\sigma\left(\frac{R}{r}\right)^{1/\alpha},\quad \mbox{on }
\Gamma.$$ Hence, by the maximum modulus principle, we have

$$|f(z_0)|e^{\sigma}=|F(z_0)|\leq \max\left\{E,\varepsilon \exp\sigma\left(\frac{R}{r}\right)^{1/\alpha}
\right\},$$ from which we have

$$|f(z_0)|=e^{-\sigma}|F(z_0)|\leq \max\left\{Ee^{-\sigma},\varepsilon
\exp\sigma\left[\left(\frac{R}{r}\right)^{1/\alpha}-1\right]
\right\}.$$ Now, we choose  $\sigma$ such that
$$e^{\sigma}=\left(\frac{E}{\varepsilon}\right)^{\left(\frac{r}{R}\right)^{1/\alpha}}$$
and we obtain \eqref{s38-1}. $\blacksquare$

\medskip

Stability estimates for the analytic continuation can be
proved even for more general sets than those considered
in examples 1 and 2. We report, without proof, the
following result (see \cite[cap. III]{L-R-S}, \cite{Is1}):

Let $\Omega\subset \mathbb{C}$ be a bounded simply connected open set whose boundary is of $C^1$ class. Let
$\Gamma=\partial \Omega$ and let us assume that
$\Gamma=\Gamma_1\cup\Gamma_2$ where $\Gamma_1$ and $\Gamma_2$ is a regular path
such that $\Gamma_1\cap\Gamma_2=\emptyset$.
Let us assume that $f\in C^0\left(\overline{\Omega}\right)$, $f$ 
holomorphic in $\Omega$ satisfying

$$|f(z)|\leq \varepsilon, \quad \forall z\in \Gamma_1, \quad\quad |f(z)|\leq E, \quad \forall z\in
\Gamma_2$$ then

\begin{equation}\label{s40-1}
|f(z)|\leq E^{1-\omega(z)}\varepsilon^{\omega(z)},\quad \forall z\in
\Omega,
\end{equation}
where $\omega(z)$ is the  harmonic function in $\Omega$ such that
$$\omega(z)=1, \quad \forall z\in \Gamma_1, \quad\quad \omega(z)=0, \quad \forall z\in
\Gamma_2.$$ $\omega$ is called \textit{harmonic measure associated to
$\Gamma_1$ in $\Omega$}

\bigskip

\textbf{Example 4: Stability estimate for the continuation of real analytic 
functions.}

In Theorem \ref{teo1-23C} we have seen that if $f:\Omega\rightarrow
\mathbb{R}$ (or $\mathbb{C}$), where $\Omega$ is a connected open set of
$\mathbb{R}^n$, is an analytic function and  $D\subset \Omega$
is a (nonempty) open set then

\begin{equation}\label{s41-1}
f=0\quad\mbox{ in } D\quad \Longrightarrow \quad f=0 \quad\mbox{ in
} \Omega.
\end{equation}
That is, a real analytic function on a connected open $\Omega$
is determined by its values on any nonempty open set
$D\subset \Omega$. Nevertheless small errors on $f_{|D}$ can have
uncontrollable effects on $f(z)$ for $z\in \Omega\setminus
\overline{D}$. 

In the present Example 4, as application of Example 3,
we will find an error estimate for the analytic continuation  problem.

Let us consider the following particular situation: let $\Omega$ be a
star shaped open set of $\mathbb{R}^n$ w.r.t.  $x_0\in \Omega$ (i.e. for every $x\in\Omega$ we have $x_0+t\left(x-x_0\right)\in \Omega$ for every $x\in \Omega$).
Let us assume that $\overline{B_r(x_0)}\subset \Omega$ for some
$r>0$. Let $f:\Omega\rightarrow \mathbb{R}$ satisfy
($E,\varepsilon>0$)

\begin{equation}\label{s41-2}
f\in \mathcal{C}_{E, \rho}(x), \quad \forall x\in \Omega,
\end{equation}
that is (compare Definition \ref{def2-24C})

\begin{equation}\label{s41n-2}
\left|\partial^{\alpha}f(x)\right| \leq E\rho^{-|\alpha|}|\alpha|! \,\quad
\forall \alpha \in \mathbb{N}^n_0, \ \forall x\in \Omega,
\end{equation}
and 
\begin{equation}\label{s42-1}
|f(x)|\leq \varepsilon, \quad \forall x\in B_r.
\end{equation}
We want to prove the following stability estimate

\begin{equation}\label{s42-2}
\left|f(x)\right| \leq (2E)^{1-\theta}\varepsilon^{\theta}\quad
\forall x\in \Omega,
\end{equation}
where $\theta\in (0,1)$ and $\theta$ depends on $n$,
$\frac{\rho}{r}$ and $\frac{\rho}{d}$ only, where $d$ is the diameter of
$\Omega$.

\bigskip

\textbf{Proof of \eqref{s42-1}.}

It is not restrictive to assume $x_0=0$. The idea of the proof
is as follows: let us fix $x\in \Omega\setminus \overline{B_r}$
and let us consider the function
\begin{equation}\label{s42-3}
\varphi(t)=f(tx), \quad t\in[0,1];
\end{equation}
we extend such a function holomorphically  to a function $\varphi$ in a neighborhood (in $\mathbb{C}$) of $\left\{t+i0: t\in[0,1]\right\}$ and by \eqref{s41n-2},
\eqref{s42-1} and the result of Example 3, we reach
\eqref{s42-2}.

By formula \eqref{6-3N} we have, for $t_0\in[0,1]$ and $k\in
\mathbb{N}_0$
\begin{equation*}
\varphi^{(k)}(t_0)=
\sum_{|\alpha|=k}\frac{k!}{\alpha!}x^{\alpha}\left(\partial^{\alpha}
f\right)(t_0x).
\end{equation*}
By \eqref{s41n-2} we have

\begin{equation}\label{s43-1}
 \begin{aligned}
\left|\varphi^{(k)}(t_0)\right|&\leq
\sum_{|\alpha|=k}\frac{k!}{\alpha!}|x|^{|\alpha|}\left|\left(\partial^{\alpha}
f\right)(t_0x)\right|\leq\\&\leq
Ek!\left(\frac{|x|}{\rho}\right)^k\sum_{|\alpha|=k}\frac{k!}{\alpha!}
=\\& =Ek!\left(\frac{n|x|}{\rho}\right)^k.
\end{aligned}
\end{equation}
Hence the power series

\begin{equation}\label{s43-2}
\sum_{k=0}^{\infty}\frac{1}{k!}\varphi^{(k)}(t_0)\left(z-t_0\right)^k,
\quad z=t+i\tau\in \mathbb{C},
\end{equation}
converges in

\begin{equation*}
B_h(t_0)=\left\{z\in \mathbb{C}: \left|z-t_0\right|<h\right\},
\end{equation*}
where $h=\frac{\rho}{2n|x|}$. Moreover the sum of power series 
\eqref{s43-2} is holomorphic in $B_h(t_0)$. The above extension can be performed
for every $t_0\in [0,1]$ therefore the function $\varphi$ can be holomorphically extended in (Figure 10.4)

$$K=\left\{z\in \mathbb{C}:\quad \mbox{dist}(z, I)<h\right\},$$
where $I:=\left\{t+i0: t\in[0,1]\right\}$. The extension $\varphi$
to $K$ is formally written as $\varphi(t+i\tau)$ and by
\eqref{s43-1} we have

\begin{figure}\label{figura-ps44}
	\centering
	\includegraphics[trim={0 0 0 0},clip, width=9cm]{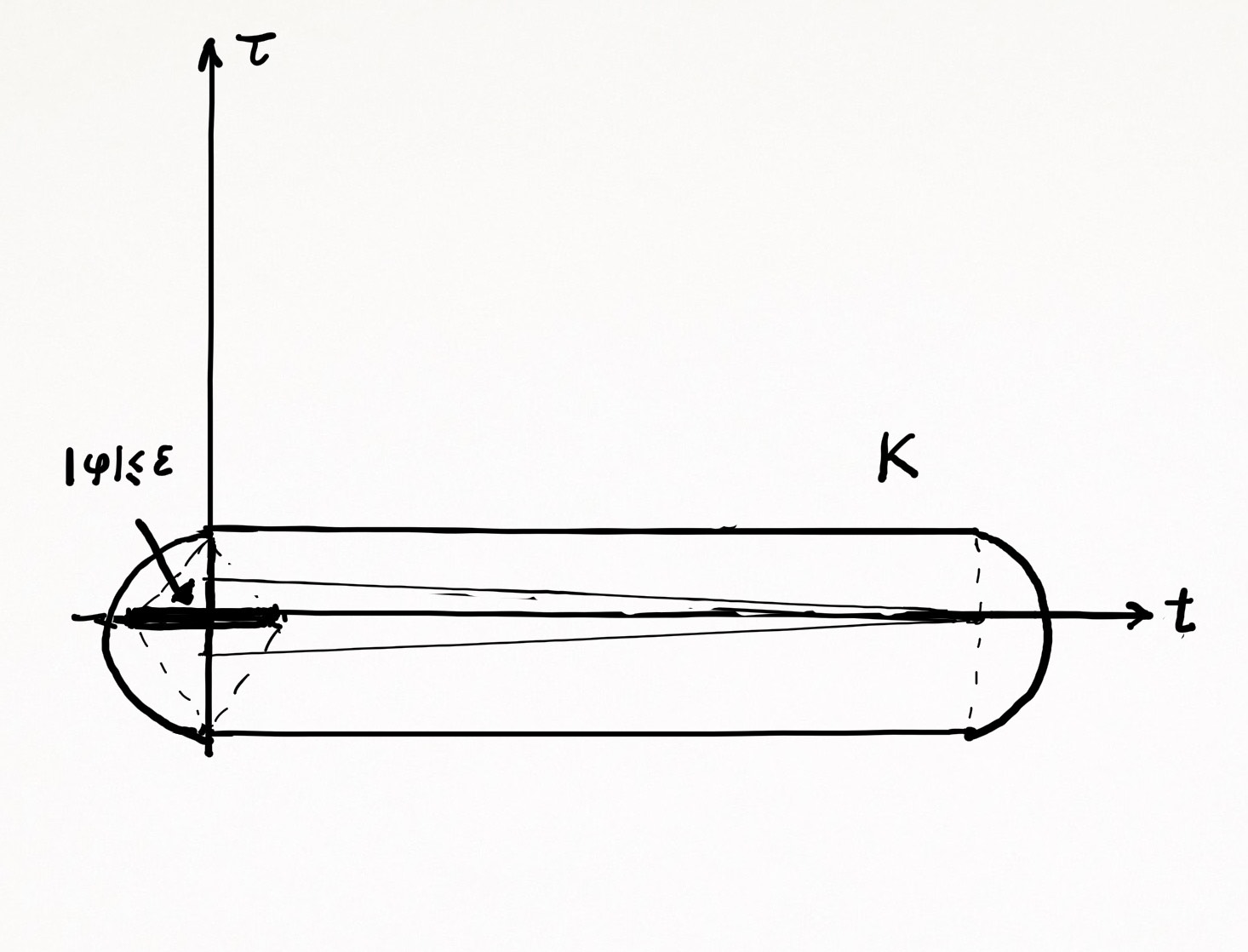}
	\caption{}
\end{figure}

\begin{equation}\label{s44-1}
 |\varphi(t+i\tau)|\leq 2E, \quad\mbox{for every } t+i\tau\in K.
\end{equation}
On the other hand by \eqref{s42-1} and by \eqref{s42-3} we have

\begin{equation}\label{s44-2}
 |\varphi(t+i0)|\leq \varepsilon, \quad\mbox{for}\quad |t|\leq \frac{r}{|x|}.
\end{equation}
At this point it suffices to prove an estimate from above of
$|\varphi(t+i0)|=|f(x)|$. 
This estimate can be obtained
by applying twice estimate \eqref{s38-1}. First we apply estimate 
\eqref{s38-1} in the triangles

$$S_+=\left\{t+i\tau : \quad |t|\leq s \mbox{, }
0\leq\tau \leq h\left(1-\frac{|t|}{s}\right)\right\},$$

$$S_-=\left\{t+i\tau : \quad |t|\leq s \mbox{, }
-h\left(1-\frac{|t|}{s}\right)\leq\tau\leq 0\right\},$$
 where
$$s:=\min\{\frac{r}{|x|}, h\}.$$ We get
\begin{equation}\label{s45-1}
 |\varphi(i\tau)|\leq (2E)^{1-\vartheta}\varepsilon^{\vartheta}, \quad\mbox{for } |\tau|\leq
 \frac{h}{2},
\end{equation}
where

\begin{equation}\label{s45-2}
\vartheta=\left(\frac{h}{2R_0}\right)^{1/\alpha_0}=\frac{1}{\left(2\sqrt{\frac{s^2}{h^2}+1}\right)^{1/\alpha_0}},
\end{equation}
$$\alpha_0=\frac{2}{\pi}\arctan \frac{s}{h},$$
$$R_0=\sqrt{s^2+h^2},$$
(hence, $\vartheta$ depends on $\frac{r}{\rho}$ and $n$ only).

Now, we use  \eqref{s45-1} and \eqref{s44-1} to apply
\eqref{s38-1} in the triangle

$$T=\left\{t+i\tau : \quad |\tau|\leq\frac{h}{2} \mbox{, }
0\leq t\leq (h+1)\left(1-\frac{2|\tau|}{h}\right)\right\}.$$ To this aim, set

$$\alpha_1=\frac{2}{\pi}\arctan \frac{h}{h+1},$$

$$R_1=\sqrt{\left(\frac{h}{2}\right)^2+(1+h)^2},$$ and we have
\begin{equation}\label{s46-1}
 |\varphi(1+i0)|\leq (2E)^{1-\vartheta\widetilde{\vartheta}}\varepsilon^{\vartheta\widetilde{\vartheta}},
\end{equation}
where
$$\widetilde{\vartheta}=\left(\frac{1}{R_1}\right)^{1/\alpha_1}.$$
Therefore we have proved  \eqref{s42-2} with
$\theta=\vartheta\widetilde{\vartheta}$. $\blacksquare$

\bigskip

\textbf{Concluding Remarks.}

Come back to the holomorphic functions. The Hadamard three circle inequality  allows us to estimate  $|f(z)|$ for $z\in B_R$ provided we know that 
$f$ is holomorphic in $B_R$ and, in addition, we know
\begin{equation}\label{l1-1} \sup_{|z|\leq r}|f(z)|\leq
\varepsilon \mbox{ (error) },
\end{equation}
and
\begin{equation}\label{l1-2} \sup_{|z|\leq R}|f(z)|\leq E \mbox{ (a priori information) }.
\end{equation}

As a matter of fact, we have

\begin{equation}\label{l1-3}
|f(z)|\leq \varepsilon^{\theta_{|z|}}E^{1-\theta_{|z|}},
\end{equation}
where
\begin{equation}\label{l1-4}\theta_{|z|}=\frac{\log R/|z|}{\log R/r}.
\end{equation}

It is immediately checked that the a priori information

$$\sup_{|z|\leq R}|f(z)|\leq E, $$ is not sufficient to
control the error on $\{|z|=R\}$. It is enough to consider
$f_n(z)=\left(\frac{z}{R}\right)^n$, obtaining 

$$\sup_{|z|\leq r}|f_n(z)|\rightarrow 0 \mbox{ as } n\rightarrow
\infty,\quad\mbox{ and }\quad |f_n(z)|=1, \mbox{ as } |z|=R.$$

We will now show that by "strengthening" the a priori information
\eqref{l1-2} we can find a stability estimate for $|f(z)|$
when $|z|=R$. Let us consider, for instance, the following a
priori information:
 $f\in C^{0,\alpha}\left(\overline{B_R}\right)$, $0<\alpha\leq 1$

\begin{equation}\label{l1-alpha} \sup_{|z|\leq R}|f(z)|+\left[f\right]_{0,\alpha}\leq
E_{\alpha},
\end{equation}
where
\begin{equation}\label{l1-alpha-1}
\left[f\right]_{0,\alpha}=\sup_{z,w\in B_R, z\neq w
}\frac{|f(z)-f(w)|}{|z-w|}.
\end{equation}

\medskip

Let us assume, for brevity, that $R=1$. Let $z_0\in \partial B_1$. Set
$z_t=z_0(1-t)$, with $t\in [0,1)$ to be chosen, we have

\begin{equation*}
 \begin{aligned}
|f(z_0)|&\leq
\left|f(z_t)-f(z_0)\right|+\left|f(z_t)\right|\leq\\&\leq
E_{\alpha}\left|z_t-z_0\right|^{\alpha}+\left|f(z_t)\right|\leq\\&
\leq E_{\alpha} t^{\alpha}+\left|f(z_t)\right|.
\end{aligned}
\end{equation*}
On the other hand by \eqref{l1-3} we have

\begin{equation*}
|f(z)|\leq
\varepsilon^{\widetilde{\theta}_{t}}E_{\alpha}^{1-\widetilde{\theta}_{t}},
\end{equation*}
where
\begin{equation*}
\widetilde{\theta}_{t}=\frac{\log \frac{1}{1-t}}{\log 1/r}.
\end{equation*}
Hence

\begin{equation*}
|f(z_0)|\leq
E_{\alpha}\left(t^{\alpha}+\varepsilon_1^{\log\frac{1}{1-t}}\right),\quad
\forall t\in [0,1),
\end{equation*}
where

\begin{equation}\label{l3-6}
\varepsilon_1=\left(\frac{\varepsilon}{E_{\alpha}}\right)^{\frac{1}{|\log
r|}}.
\end{equation}
Now we have

$$t^{\alpha}+\varepsilon_1^{\log\frac{1}{1-t}}=t^{\alpha}+\exp\left(\log(1-t)\left|\log\varepsilon_1\right|\right).$$
On the other hand, we have
$$\log(1-t)\leq -t.$$
Hence

\begin{equation}\label{l3-1}
|f(z_0)|\leq
E_{\alpha}\left(t^{\alpha}+\exp\left(-t\left|\log\varepsilon_1\right|\right)\right),\quad
\forall t\in [0,1).
\end{equation}
Now we note that, if $\varepsilon_1<1$ then
$$0<\left|\log\varepsilon_1\right|^{-1}\left(\log\left|\log\varepsilon_1\right|\right)\leq
e^{-1}$$ and we can choose 

\begin{equation*}
t=\left|\log\varepsilon_1\right|^{-1}\left(\log\left|\log\varepsilon_1\right|\right)\in[0,1).
\end{equation*}
We obtain

\begin{equation}\label{l3-3}
\begin{aligned}
|f(z_0)|&\leq
E_{\alpha}\left(\left|\log\varepsilon_1\right|^{-\alpha}\left(\log\left|\log\varepsilon_1\right|\right)+\left|\log\varepsilon_1\right|^{-1}\right)\leq
\\& \leq C
E_{\alpha}\left|\log\varepsilon_1\right|^{-\alpha}\left(\log\left|\log\varepsilon_1\right|\right),
\end{aligned}
\end{equation}
where $C$ depends on $\alpha$ only. If $\varepsilon_1\geq 1$ then we have trivially

\begin{equation}\label{l3-4}
|f(z_0)|\leq E_{\alpha}\leq E_{\alpha}\varepsilon_1.
\end{equation}
Thus, by \eqref{l3-3} and \eqref{l3-4} we have the following
\textbf{stability estimate}, for every $z_0\in \partial B_1$

\begin{equation}\label{l3-5}
|f(z_0)|\leq \widetilde{C}
E_{\alpha}\left|\log\varepsilon_1\right|^{-\alpha}\left(\log\left|\log\varepsilon_1\right|\right).
\end{equation}
where $\widetilde{C}$ depends on $\alpha$ and $\varepsilon_1$ is given by \eqref{l3-6}. $\blacklozenge$

\newpage

\chapter{The John stability Theorem for the Cauchy problem for PDEs with analytic coefficients} \markboth{Chapter 11. The John stability Theorem for the Cauchy problem}{}\label{Holmgren-John}
\section{Statement of the Theorem}\label{Holmgren-John-enunc} The stability estimate that
we present in this Chapter is due to F. John \cite{Joh60}. The basic elements
of the proof are as follows.

\begin{enumerate}
    \item The Green identity \index{Green identity} and the construction of an appropriate  solution of the adjoint operator.
    \item The stability estimates for the analytic continuation problem.
\end{enumerate}

In what follows we will consider the following linear system

\begin{equation}\label{Cs1-1}
u_t(x,t)=\sum_{j=1}^n A_j(x,t)u_{x_j}(x,t)+A_0(x,t)u(x,t),
\end{equation}
where $u:=\left(u^1, \cdots, u^N \right)^T$, $x\in \mathbb{R}^n$,
$t\in \mathbb{R}$, $A_j(x,t)$, $j=0,1,\cdots, n$  are
$N\times N$ matrices . Moreover, let us introduce the following notations 

\begin{equation}\label{Cs1-2}
\gamma(x)=\left(1-|x|^2 \right)^{n+1},
\end{equation}
For any $\lambda\in \mathbb{R}$ let us denote by $S_{\lambda}$ the surface
\begin{equation}\label{Cs2-1}
S_{\lambda}=\left\{(x,\lambda\gamma(x))|x\in B_1 \right\}
\end{equation}
and, for any $\lambda_1<\lambda_2$, let

\begin{equation}\label{Cs2-2}
\mathcal{R}_{\lambda_1,\lambda_2}=\left\{(x,t)\in \mathbb{R}^{n+1}
|x\in B_1 \quad \lambda_1\gamma(x)<t<\lambda_2\gamma(x) \right\}.
\end{equation}

\bigskip

\begin{theo}[\textbf{John stability estimate}]\label{Cs2-teorema-st}
	\index{Theorem:@{Theorem:}!- John stability estimate@{- John stability estimate}}
Let $c_0, L, M, \rho, E, \varepsilon$ positive numbers. Let
$A_j$, $j=1,\cdots, n$ and $B$ matrices $N \times N$ whose entries are analytic in $\overline{\mathcal{R}_{0,L}}$ and satisfy

\begin{equation}\label{Cs2-3}
A_j\in \mathcal{C}_{M,\rho}\left(\overline{x},\overline{t}\right),
\quad j=0,1,\cdots, n\quad
\forall \left(\overline{x},\overline{t}\right)\in \mathcal{R}_{0,L}.
\end{equation}
 Set
\begin{equation}\label{Cs2n-3}
A(x,\lambda)=I+\lambda\sum_{j=1}^n
A_j(x,\lambda\gamma(x))\gamma_{x_j}(x).
\end{equation}
 Let us assume that 

\begin{equation}\label{Cs2-4}
\left|\det A(x,\lambda)\right|\geq c_0,\quad \forall x\in
\overline{\mathcal{R}_{0,L}},\quad \forall \lambda \in [0,L],
\end{equation}
(that is $S_{\lambda}$ is a noncharacteristic surface for every 
$\lambda \in [0,L]$).

Let $u\in C^{n+1}\left(\overline{\mathcal{R}_{0,L}}\right)$ satisfy
\begin{subequations}
\label{Cs2-5}
\begin{equation}
\label{Cs2-5a} u_t(x,t)=\sum_{j=1}^n
A_j(x,t)u_{x_j}(x,t)+A-0(x,t)u(x,t),\quad \mbox{ } \forall (x,t)\in
\mathcal{R}_{0,L}
\end{equation}
\begin{equation}
\label{Cs2-5b} \left\Vert
u(\cdot,0)\right\Vert_{L^{\infty}(B_1)}\leq\varepsilon,
\end{equation}
\begin{equation}
\label{Cs2-5c} \left\Vert
u\right\Vert_{C^{n+1}\left(\overline{\mathcal{R}_{0,L}}\right)}\leq
E.
\end{equation}
\end{subequations}
Then, for every $r\in (0,1)$, we have

\begin{equation}\label{Cs3-1}
\left|u(x,t)\right|\leq
\frac{C(E+2\varepsilon)}{(1-r)^{n+1}}\left|\log\frac{\varepsilon}{E+2\varepsilon}\right|^{-1},\quad
\forall x\in \overline{\mathcal{R}_{0,L}}\cap
\left(\overline{B_{1-r}}\times \mathbb{R}\right),
\end{equation}
where $C$ depends on $M, L, \rho, c_0$ and $n$ only.
\end{theo}

\section{Proof of the Theorem}\label{Holmgren-John-dimost}

Let us premise the following 

\begin{lem}\label{stabilita-john}
Let $A_j$ be as in Theorem \ref{Cs2-teorema-st}. Let us assume that

\begin{equation}\label{51ter-1}
A_j,\in \mathcal{C}_{M_0,\rho_0}\left(\overline{x},
\overline{t}\right), \quad \forall \left(\overline{x},
\overline{t}\right)\in \overline{B_1}\times[-\delta_0,\delta_0].
\end{equation}
Let $W\in \mathbb{R}^N$ be such that $|W|=1$ and let $\xi\in
\mathbb{R}^n$ arbitrary.

Let $U$ be the solution to the Cauchy problem

\begin{equation}\label{51ter-1-28}
\begin{cases}
U_t=\sum_{j=1}^n A_j(x,t)U_{x_j}+A_0(x,t)U, \\
\\
U(x,0)=e^{-ix\cdot\xi}W,\quad \forall x\in B_1.
\end{cases}
\end{equation}
Then there exist $M,\rho,\delta$ positive numbers  depending by
$M_0,\rho_0,\delta_0$, but independent of  $\xi$ such that

\begin{equation}\label{51ter-2}
U\in
\mathcal{C}_{Me^{|\xi|},\rho}\left(\overline{x},\overline{t}\right),
\quad \forall \left(\overline{x},\overline{t}\right)\in
\overline{B_1}\times [-\delta,\delta],
\end{equation}
that is

\begin{equation}\label{51ter-3}
\left|\partial^{\alpha}U\right|\leq Me^{|\xi|}\rho^{-|\alpha|}
|\alpha|!, \quad \forall \left(\overline{x},\overline{t}\right)\in
\overline{B_1}\times [-\delta,\delta], \mbox{ }\forall \alpha\in
\mathbb{N}^{n+1}_0.
\end{equation}
\end{lem}

\textbf{Proof of Lemma \ref {stabilita-john}.} Let
$$\psi(x)=e^{-|\xi|}e^{-i\xi\cdot x}W.$$ Let us consider the Cauchy problem

\begin{equation}\label{51ter-4}
\begin{cases}
V_t=\sum_{j=1}^n A_j(x,t)V_{x_j}+A_0(x,t)V, \\
\\
V(x,0)=\psi(x),\quad \forall x\in B_1.
\end{cases}
\end{equation}
We have, trivially,

\begin{equation}\label{51ter-5}
U(x,t)=e^{|\xi|}V(x,t).
\end{equation}
On the other hand

\begin{equation*}
\begin{aligned}
\left|\partial^{\alpha'}\psi(x)
\right|&=\left|(i\xi)^{\alpha'}e^{-|\xi|}e^{-\xi\cdot
x}W\right|=|\xi|^{|\alpha'|}e^{-|\xi|}\leq\\&\leq
|\alpha'|^{|\alpha'|}e^{-|\alpha'|}\leq |\alpha'|!, \quad \forall
\alpha'\in \mathbb{N}_0^n, \mbox{ } \forall x\in \overline{B_1},
\end{aligned}
\end{equation*}
hence
$$\psi\in \mathcal{C}_{1,1}\left(\overline{x}\right),\quad\forall \overline{x}\in \overline{B_1}.$$
Therefore there exist $M,\rho,\delta$ which depend on
$M_0,\rho_0,\delta_0$, but independent of $\xi$ such that

$$V\in \mathcal{C}_{1,1}\left(\overline{x}, \overline{t}\right),
\quad \forall \left(\overline{x}, \overline{t}\right)\in
\overline{B_1}\times [-\delta,\delta]$$ and by \eqref{51ter-5}
we have \eqref{51ter-2}. $\blacksquare$

\bigskip

 \textbf{Proof of Theorem \ref{Cs2-teorema-st}.}

\noindent \textbf{Step 1. The Green identity.}

We have 
\begin{equation}\label{Cs3-2}
	\begin{aligned}
	&v^T\left(u_t-\sum_{j=1}^n A_ju_{x_j}-A_0u\right)=\\&
	=\partial_t(v^Tu)-
	\partial_{x_j}\left(v^T \sum_{j=1}^n A_ju\right)-\\&-\left(v^T_t-\sum_{j=1}^n
	\left(v^TA_j\right)_{x_j}+v^TA_0\right)u.	
\end{aligned}
\end{equation}
Let now  $v\in
C^1\left(\overline{\mathcal{R}_{\lambda_1,\lambda_2}}\right)$ be a
solution to the adjoint system
\begin{equation}\label{Cs3-3}
v_t-\sum_{j=1}^n \left(A^T_jv\right)_{x_j}+A_0^Tv=0, \quad \mbox{in }
\mathcal{R}_{\lambda_1,\lambda_2}.
\end{equation}
Since $u$ is a solution to system \eqref{Cs2-5a},
integrating both the sides of \eqref{Cs3-2} over
$\mathcal{R}_{\lambda_1,\lambda_2}$ we have

\begin{equation*}
	\begin{aligned}
	&0=\int_{\mathcal{R}_{\lambda_1,\lambda_2}} \left[\partial_t(v^Tu)-
	\partial_{x_j}\left(v^T \sum_{j=1}^n A_ju\right)\right]dxdt=\\&
	=\int_{\partial\mathcal{R}_{\lambda_1,\lambda_2}}\left[(v^Tu)(\nu\cdot
	e_{n+1})- \left(v^T \sum_{j=1}^n A_ju (\nu\cdot
	e_{j})\right)\right]dS=\\& 
	=\int_{B_1}\left(v^Tu+\lambda_2 v^T\sum_{j=1}^n
	\gamma_{x_j}A_ju\right)(x,\lambda_2\gamma(x))dx-\\& 
	-\int_{B_1}\left(v^Tu+\lambda_1 v^T\sum_{j=1}^n
	\gamma_{x_j}A_ju\right)(x,\lambda_1\gamma(x))dx.
	\end{aligned}
\end{equation*}

Hence, recalling \eqref{Cs2n-3} we get

\begin{equation}\label{Cs4-1}
	\begin{aligned}
	&\int_{B_1}
	v^T(x,\lambda_2\gamma(x))A(x,\lambda_2)u(x,\lambda_2\gamma(x))dx=\\& 
	=\int_{B_1}
	v^T(x,\lambda_1\gamma(x))A(x,\lambda_1)u(x,\lambda_1\gamma(x))dx.	
\end{aligned}
\end{equation}

\bigskip

\noindent \textbf{Step 2. Construction of an appropriate solution to
\eqref{Cs3-3}.} Let $W\in \mathbb{R}^N$ be such that $|W|=1$ and $\xi\in \mathbb{R}^n$ arbitrary. Let us denote by $w$ the function

\begin{equation}\label{Cs4-2}
w(x)=e^{-ix\cdot \xi}W.
\end{equation}
Let $\lambda\in [0,L]$. Let us consider the following Cauchy problem

\begin{equation}\label{Cs4-3}
\begin{cases}
v_t=\sum_{j=1}^n \left(A^T_j(x,t)v\right)_{x_j}-A_0^T(x,t)v, \\
\\
v(x,\lambda\gamma(x))=\gamma(x)\left(A^T(x,\lambda)\right)^{-1}w(x),\quad
x\in B_1.
\end{cases}
\end{equation}
Let us prove that there exists  $\delta>0$, depending on $M, \rho, c_0$ only,
such that there exists the solution $v(x,t;\lambda)$ of \eqref{Cs4-3},  and it is 
analytic in $\overline{\mathcal{R}_{\lambda-\delta,\lambda+\delta}}$. To this purpose we perform some change of variables. First, we set

\begin{equation*}
s=\frac{t}{\gamma(x)}, \quad
v(x,t;\lambda)=\gamma(x)V\left(x,\frac{t}{\gamma(x)};\lambda
\right)
\end{equation*}
and we have
\begin{equation*}
v_t(x,t;\lambda)=V_s\left(x,\frac{t}{\gamma(x)};\lambda \right).
\end{equation*}

\begin{equation*}
 \begin{aligned}
v_{x_j}&=\gamma
V_{x_j}+\gamma_{x_j}V-\frac{t\gamma_{x_j}}{\gamma}V_s=
\\&= \gamma
V_{x_j}+\gamma_{x_j}V-s\gamma_{x_j}V_s.
\end{aligned}
\end{equation*}
Inserting what obtained above in system \eqref{Cs4-3}, we get

\begin{equation*}
 \begin{aligned}
V_{s}&=v_t=\sum_{j=1}^n
A^T_jv_{x_j}+\left(\sum_{j=1}^nA^T_{j,x_j}-A_0^T\right)v=\\&=\sum_{j=1}^n
A^T_j\left(\gamma V_{x_j}-s\gamma_{x_j}V_s\right)+
\left[\sum_{j=1}^n\left(A^T_j\gamma_{x_j}+A^T_{j,x_j}\gamma\right)-A_0^T\gamma\right]V,
\end{aligned}
\end{equation*}
From which (recalling \eqref{Cs2n-3}), we get

\begin{equation*}
 \begin{aligned}
A^T(x,s)V_{s}&=\gamma\sum_{j=1}^n
A^T_jV_{x_j}+\left(\sum_{j=1}^nA^T_{j,x_j}-A_0^T\right)V=\\&=\sum_{j=1}^n
A^T_j\left(\gamma V_{x_j}-s\gamma_{x_j}V_s\right)+
\left[\sum_{j=1}^n\left(A^T_j\gamma_{x_j}+A^T_{j,x_j}\gamma\right)-A_0^T\gamma\right]V.
\end{aligned}
\end{equation*}
Now, set

\begin{equation*}
\overline{A}_j(x,s)=\gamma(x)\left(A^T(x,s)\right)^{-1}A^T_j(x,s\gamma(x)),
\quad j=1,\cdots,n,
\end{equation*}

\begin{equation*}
 \begin{aligned}
&\overline{A}_0(x,s)=\\&=\left(A^T(x,s)\right)^{-1}\left[\sum_{j=1}^nA^T_j(x,s\gamma(x))\gamma_{x_j}+\gamma
\left(\sum_{j=1}^n
A^T_{j,x_j}(x,s\gamma(x))-A_0^T(x,s\gamma(x))\right)\right].
\end{aligned}
\end{equation*}
Therefore problem \eqref{Cs4-3} can be written as

\begin{equation}\label{Cs6-1}
\begin{cases}
V_s(x,s;\lambda)=\sum_{j=1}^n \overline{A}_j(x,s)V_{x_j}(x,s;\lambda)+\overline{A}_0(x,s)V(x,s;\lambda), \\
\\
V(x,s;\lambda)_{|s=\lambda}=\left(A^T(x,\lambda)\right)^{-1}w(x),\quad
\forall x\in B_1.
\end{cases}
\end{equation}
Now we denote

\begin{equation}\label{Cs6-2}
Z(x,s;\lambda)=\left(A^T(x,\lambda)\right)^{-1}V(x,s+\lambda;\lambda)
\end{equation}
and by \eqref{Cs4-2}, \eqref{Cs6-2} we have

\begin{equation}\label{Cs6-3}
\begin{cases}
Z_s(x,s;\lambda)=\sum_{j=1}^n \widetilde{A}_j(x,s;\lambda)Z_{x_j}(x,s;\lambda)+\widetilde{A}_0(x,s;\lambda)Z(x,s;\lambda), \\
\\
Z(x,0;\lambda)=e^{-ix\cdot \xi} W,\quad x\in B_1,
\end{cases}
\end{equation}
where

\begin{equation*}
\widetilde{A}_j(x,s;\lambda)=A^T(x,\lambda)\overline{A}_j(x,s+\lambda)\left(A^T(x,\lambda)\right)^{-1},\quad
j=1,\cdots,n,
\end{equation*}
and
\begin{equation*}
	\begin{aligned}
&\widetilde{A}_0(x,s;\lambda)=\\&=A^T(x,\lambda)\left[\overline{B}(x,s+\lambda)\left(A^T(x,\lambda)\right)^{-1}+\sum_{j=1}^n
\overline{A}_j(x,s+\lambda)\partial_{x_j}\left(A^T(x,\lambda)\right)^{-1}\right].
\end{aligned}
\end{equation*}
Now $\widetilde{A}_j(x,s;\lambda)$, $j=1,\cdots,n$ and
$\widetilde{A}_0(x,s;\lambda)$ are analitic functions in $(x,s,\lambda)$. In addition, for every
$\left(\overline{x},\overline{s},\overline{\lambda}\right)\in
\overline{B_1}\times [0,L]\times[0,L]$, we have
\begin{equation*}
\widetilde{A}_j\in
\mathcal{C}_{M',\rho'}\left(\overline{x},\overline{s},\overline{\lambda}\right),\quad
j=0,1,\cdots,n,
\end{equation*}
where $M'$ e $\rho'$ depend on $M,\rho,c_0$ (and $n$, which we will omit in the sequel) only.

By Lemma \ref{stabilita-john}, there exist
$M'',\rho'',\delta$ depending on $M,\rho, c_0$ and $L$ only such that
there is  $Z$ which is the solution to \eqref{Cs6-3}, it is analitic in
$\overline{B_1}\times[-\delta,\delta]\times[0,L]$ and satisfies

\begin{equation}\label{Cs7-1}
Z\in
\mathcal{C}_{M''e^{|\xi|},\rho''}\left(\overline{x},\overline{s},\overline{\lambda}\right),
\quad \forall
\left(\overline{x},\overline{s},\overline{\lambda}\right)\in
\overline{B_1}\times [-\delta,\delta]\times[0,L].
\end{equation}
Coming back to problem \eqref{Cs4-3}, we have that there exists $v$,
solution to \eqref{Cs4-3} in
$\mathcal{R}_{\lambda-\delta,\lambda+\delta}$, such that

\begin{equation}\label{Cs7-2}
\left|\partial^{\alpha}_{x,t,\lambda}v(x,t;\lambda)\right|\leq M_0
e^{|\xi|}\rho^{-|\alpha|}|\alpha|!, \mbox{ } 
\forall \alpha \in \mathbb{N}_0^{n+2}, \end{equation}
for all
$\left(x,t;\lambda\right)\in
\overline{\mathcal{R}_{\lambda-\delta,\lambda+\delta}}\times[0,L],$ where $M_0,\delta_0$ and $\delta$ depend on $M,\rho, c_0$ and $L$
only. Let us note that to obtain \eqref{Cs7-2} for every
$\lambda\in [0,L]$ it suffices to consider  \eqref{Cs7-1} and the similar relationships on the coefficients corresponding to $\lambda=0$.

\bigskip

\noindent \textbf{Step 3. Planning the concluding part of the
	proof.}

We employ \eqref{Cs4-1}, where $v$ is the solution to  problem
\eqref{Cs4-3} for some $\lambda\in[0,L]$. Let $\lambda_0\in
[0,L]$ be fixed and let $\lambda$ satisfy
$|\lambda-\lambda_0|<\delta$. By \eqref{Cs4-1} we have

\begin{equation}\label{Cs8-1}
\begin{aligned}
g(\lambda)&:=\int_{B_1}
\gamma(x)w^T(x)u(x,\lambda\gamma(x))dx=\\&=\int_{B_1}
v^T(x,\lambda_0\gamma(x);\lambda)A(x,\lambda_0)u(x,\lambda_0\gamma(x))dx.
\end{aligned}
\end{equation}
The function $g$ is \textbf{analitic} because the integrand in the
second integral of \eqref{Cs8-1} depends analytically by
$\lambda$.

We are first interested in proving an estimates from above for
$g(\lambda)$ from which, subsequently, we will derive the estimates from above for
$u$. Setting $\lambda_0=0$ in \eqref{Cs8-1} we have

\begin{equation}\label{Cs8-2}
g(\lambda)=\int_{B_1} v^T(x,0;\lambda)u(x,0)dx, \quad \forall
\lambda\in[0,\delta).
\end{equation}
By \eqref{Cs2-5b} and by \eqref{Cs7-2}, we have, for $\alpha=0$,

\begin{equation}\label{Cs8-3}
|g(\lambda)|\leq\int_{B_1}
\left|v^T(x,0;\lambda)\right||u(x,0)|dx\leq
cM_0e^{|\xi|}\varepsilon,\ \ \forall \lambda\in[0,\delta),
\end{equation}
where $c\geq 1$ depends on $n$ only.

Let now $\lambda_0$ be an arbitrary point of $[0,L]$, by
\eqref{Cs7-2} and \eqref{Cs8-1} we get

\begin{equation*}
\begin{aligned}
\left|g^{(k)}(\lambda)\right|&=\left|\int_{B_1}
\partial_{\lambda}^kv^T(x,\lambda_0\gamma(x);\lambda)A(x,\lambda_0)u(x,\lambda_0\gamma(x))dx\right|\leq\\&\leq
\int_{B_1}\left|\partial_{\lambda}^kv^T(x,\lambda_0\gamma(x);\lambda)\right||A(x,\lambda_0)||u(x,\lambda_0\gamma(x))|dx\leq\\&\leq
C_1EM_0e^{|\xi|}\rho_0^{-k}k!,
\end{aligned}
\end{equation*}
where $C_1\geq 1$ depends on $M,\rho,c_0, L$ only. 

Summarizing we have

\begin{equation}\label{Cs9-1}
|g(\lambda)|\leq cM_0e^{|\xi|}\varepsilon, \quad\forall
\lambda\in[0,\delta)
\end{equation}
and
\begin{equation}\label{Cs9-2}
\left|g^{(k)}(\lambda)\right|\leq C_1M_0e^{|\xi|}\rho_0^{-k}k!,
\quad\forall \lambda\in[0,L]\mbox{, } k\in \mathbb{N}_0.
\end{equation}
Inequality \eqref{Cs9-2} implies that  $g$ can be extended analytically
in a neighborhood of $[0,L]\times\{0\}\subset \mathbb{C}$. In addition, for any
 $\lambda_{\star}\in[0,L]$, we have that the power series
$$\sum_{k=0}^{\infty}\frac{g^{(k)}(\lambda_{\star})}{k!}\left(z-\lambda_{\star}\right)^k,$$
has the radius of convergence equal to $\rho_0$ and

\begin{equation}\label{Cs9n-2}
\left|\sum_{k=0}^{\infty}\frac{g^{(k)}(\lambda_{\star})}{k!}\left(z-\lambda_{\star}\right)^k\right|\leq
C_1M_0e^{|\xi|}\frac{\rho_0}{\rho_0-\left|z-\lambda_{\star}\right|},
\end{equation}
for $\left|z-\lambda_{\star}\right|<\rho_0$. Therefore the sum of the power series  \eqref{s43-2} is
holomorphic in $B_{\frac{\rho_0}{2}}(\lambda_{\star})$ and the function $g$ can be extended holomorphically in (see Figure 8.1)
\begin{figure}\label{figura-pc10}
	\centering
	\includegraphics[trim={0 0 0 0},clip, width=9cm]{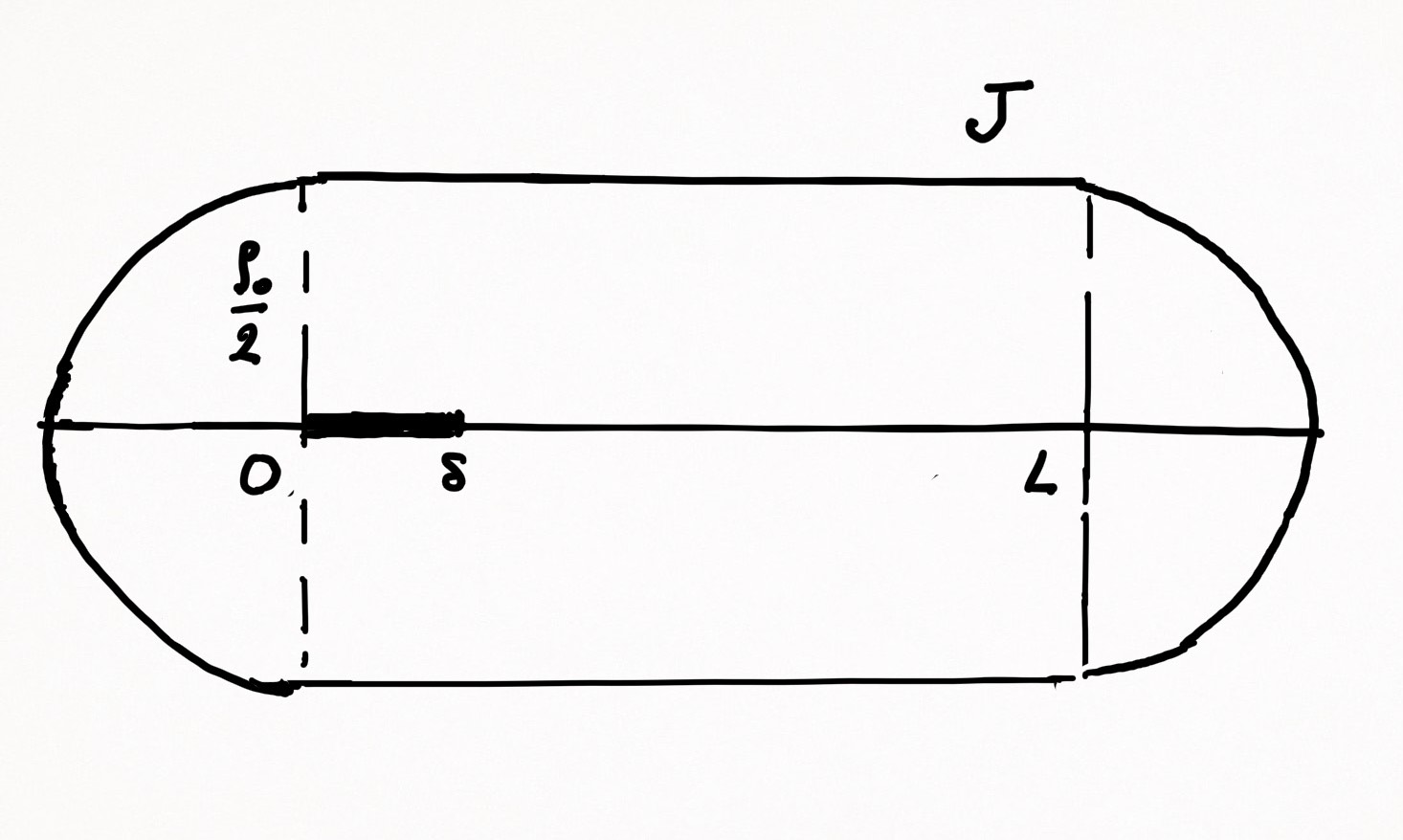}
	\caption{}
\end{figure}

$$J=\left\{z\in \mathbb{C}: \mbox{dist}(z,[0,L]\times\{0\} )<\frac{\rho_0}{2}\right\}.$$
The extension of $g$ to $J$ is formally written as 
$g(\lambda+i\tau)$ and by \eqref{Cs9n-2} we have

\begin{equation*}
 |g(\lambda+i\tau)|\leq 2C_1EM_0e^{|\xi|} , \quad\mbox{for } t+i\tau\in J,
\end{equation*}

\begin{equation*}
 |g(\lambda+i0)|\leq cM_0e^{|\xi|}\varepsilon, \quad \forall
\lambda\in[0,\delta).
\end{equation*}
Now, proceeding in a similar way to what we did to
prove \eqref{s46-1}, we get

\begin{equation*}
\begin{aligned}
 |g(\lambda+i0)|&\leq \left(2C_1EM_0e^{|\xi|}\right)^{1-\vartheta}\left(cEM_0e^{|\xi|}\varepsilon\right)^{\vartheta}\leq\\&\leq
 C_2 e^{|\xi|}E^{1-\vartheta}\varepsilon^{\vartheta},
 \quad\quad\forall \lambda\in[0,L],
\end{aligned}
\end{equation*}
where $\vartheta\in (0,1)$ depends on $L\rho_0^{-1}$ only and
$C_2=2cC_1M_0$. By the definition of $g$ given in \eqref{Cs8-1} and recalling that $w(x)=e^{-|\xi|}e^{-i\xi\cdot x}W$, we have

\begin{equation}\label{Cs10-1}
\left|\int_{B_1} \gamma(x)W^Tu(x,\lambda\gamma(x))e^{-i\xi\cdot
x}dx\right|\leq C_2 e^{|\xi|}E^{1-\vartheta}\varepsilon^{\vartheta},
 \quad\quad\forall \lambda\in[0,L].
\end{equation}

\bigskip

\noindent \textbf{Step 4. Conclusion of the proof.}

Let us fix $\lambda \in [0,L]$. Let $W=e_j$, for $j=1,\cdots,N$ and set

\begin{equation}\label{Cs10-2}
f_j(x)=
\begin{cases}
 \gamma(x)u_j(x,\lambda\gamma(x)),\quad\mbox{for } |x|\leq 1,\\
\\
0, \quad\mbox{for } |x|\geq 1.
\end{cases}
\end{equation}
Let us fix $j=1,\cdots,N$ and, in the sequel, let us omit the index $j$ by
$f_j$. By \eqref{Cs10-1} and \eqref{Cs10-2} we get

\begin{equation}\label{Cs11-1}
\left|\widehat{f}(\xi)\right|=\left|\int_{\mathbb{R}^n}f(x)e^{-i\xi\cdot
x}dx\right|\leq C_2 e^{|\xi|}E^{1-\vartheta}\varepsilon^{\vartheta},
 \quad\forall \xi\in\mathbb{R}^n,
\end{equation}
(where $\widehat{f}$ is the  Fourier transform of $f$).

The proof will be completed as soon as we estimate $|f(x)|$ by
means of \eqref{Cs11-1}. Let us recall that

\begin{equation}\label{Cs11-2}
f(x)=(2\pi)^{-n}\int_{\mathbb{R}^n}\widehat{f}(\xi)e^{-i\xi\cdot
x}d\xi,\quad\forall x\in \mathbb{R}^n.
\end{equation}

Let $s$ be a positive number which we will choose later. By
\eqref{Cs11-1} and \eqref{Cs11-2} we have, for every $x\in
\mathbb{R}^n$

\begin{equation}\label{Cs11-3}
\begin{aligned}
 |f(x)|&\leq
 (2\pi)^{-n}\int_{\mathbb{R}^n}\left|\widehat{f}(\xi)\right|d\xi=\\&=
(2\pi)^{-n}\left(\int_{|\xi|\leq
s}\left|\widehat{f}(\xi)\right|d\xi+\int_{|\xi|>
s}\left|\widehat{f}(\xi)\right|d\xi\right)\leq\\&\leq
(2\pi)^{-n}\left(C_3s^n
e^{s}E^{1-\vartheta}\varepsilon^{\vartheta}+\int_{|\xi|>
s}\left|\widehat{f}(\xi)\right|d\xi\right),
\end{aligned}
\end{equation}
where $C_3=\frac{\omega_n}{n}C_2$ ($\omega_n$ is the measure of unit ball of $\mathbb{R}^n$).

To estimate from above the last integral in \eqref{Cs11-3} we proceed
as follows. First of all we note that by the definition of 
$\gamma$ and $f$ we have $\partial^{\alpha}f(x)=0$ for every
$x\in\partial B_1$ and for every $\alpha\in \mathbb{N}^n_0$ such that $|\alpha|\leq n$ from which we have, for $k=1,\cdots,n$, by using integration by parts

\begin{equation*}
\begin{aligned}
 \left(-i\xi_k\right)^{n+1}\widehat{f}(\xi)&=
 \int_{\mathbb{R}^n}\left(-i\xi_k\right)^{n+1}e^{-i\xi\cdot x}f(x)dx=\\&=
\int_{B_1}\partial^{n+1}_k\left(e^{-i\xi\cdot
x}\right)f(x)dx=\\&=(-1)^{n+1}\int_{B_1}e^{-i\xi\cdot
x}\partial^{n+1}_kf(x)dx.
\end{aligned}
\end{equation*}
Hence \eqref{Cs2-5c} implies

\begin{equation}\label{Cs12-1}
\left|\xi_k\right|^{n+1}\left|\widehat{f}(\xi)\right|\leq
\int_{B_1}\left|\partial^{n+1}_kf(x)\right|dx\leq C_4E,
\end{equation}
where $C_4$ depends on $L$ (and on $n$) only. So that, we have trivially

\begin{equation}\label{Cs12-2}
\left|\xi\right|^{n+1}\left|\widehat{f}(\xi)\right| \leq C_5E,
\end{equation}
where $C_5=n^{\frac{n+1}{2}}C_4$. Now, by
\eqref{Cs12-1} we have

\begin{equation}\label{Cs12-3}
\begin{aligned}
\int_{|\xi|> s}\left|\widehat{f}(\xi)\right|d\xi& =\int_{|\xi|>
s}\left|\xi\right|^{-(n+1)}\left|\xi\right|^{n+1}\left|\widehat{f}(\xi)\right|d\xi\leq\\&\leq
 C_5E\int_{|\xi|>
s}\left|\xi\right|^{-(n+1)}d\xi=\omega_nC_5\frac{E}{s}.
\end{aligned}
\end{equation}

Now we use in \eqref{Cs12-3} what we obtained in \eqref{Cs12-3}  and by the trivial inequality $E<E+2\varepsilon$ we get

\begin{equation}\label{Cs12-4}
 |f(x)|\leq
 C_6(E+2\varepsilon)\left[\left(\frac{\varepsilon}{E+2\varepsilon}\right)^{\vartheta}s^ne^s+\frac{1}{s}\right],\quad\forall x\in
 \mathbb{R}^n, \mbox{ } \forall s>0,
\end{equation}
where $C_6$ depends on $M, L, \rho, c_0$ only.

In order to choose $s$ we proceed as follows. Set
$$\sigma=\left(\frac{\varepsilon}{E+2\varepsilon}\right)^{\vartheta}$$
and rewrite the term on the right--hand side of  \eqref{Cs12-4} as

$$\phi(s):=\exp\left(s+\log s-|\log\sigma|\right)+\frac{1}{s}.$$
Now we choose

\begin{equation}\label{Cs13-1}
 s:s_0=\frac{1}{2}|\log\sigma|,
 \end{equation}
taking into account that $0<\sigma\leq 2^{-\vartheta}<1$, we have

$$\phi(s_0)=\sqrt{\sigma}\left|\log\sqrt{\sigma}\right|+2|\log\sigma|^{-1}\leq
c_{\vartheta}|\log\sigma|^{-1},$$ where $c_{\vartheta}$ depends 
on $\vartheta$. Hence, by \eqref{Cs13-1} and
\eqref{Cs12-4} we have

\begin{equation}\label{Cs13-2-28}
 |f(x)|\leq
 C_7(E+2\varepsilon)\left|\log\left(\frac{\varepsilon}{E+2\varepsilon}\right)\right|^{-1},
 \end{equation}
where $C_7=C_6(c_{\vartheta}+2)\vartheta^{-1}$. Finally,
recalling that $f(x)=\gamma(x)u_j(x,\lambda\gamma(x))$ for $|x|\leq
1$ we have, for any $r\in (0,1)$

\begin{equation}\label{Cs13-2}
 |u(x,\lambda\gamma(x))|\leq
 \frac{C_7N^{1/2}}{(1-r)^{n+1}}(E+2\varepsilon)\left|\log\left(\frac{\varepsilon}{E+2\varepsilon}\right)\right|^{-1},
 \end{equation}
for $|x|\leq 1-r$, from which the thesis follows. $\blacksquare$

\bigskip

\textbf{Remarks.}

We outline some changes that we should make in the case of nonhogeneous system 
\medskip

\noindent \textbf{1. } Let us consider the case in which instead of
\eqref{Cs2-5} we have, for $(x,t)\in
\mathcal{R}_{0,L}$,

\begin{subequations}
\label{Cs14-1}
\begin{equation}
\label{Cs14-1a} u_t(x,t)=\sum_{j=1}^n
A_j(x,t)u_{x_j}(x,t)+A_0(x,t)u(x,t)+F(x,t),
\end{equation}
\begin{equation}
\label{Cs14-1b} \left\Vert
u(\cdot,0)\right\Vert_{L^{\infty}(B_1)}\leq\varepsilon,
\end{equation}
\begin{equation}
\label{Cs14-1c} \left\Vert
u\right\Vert_{C^{n+1}\left(\overline{\mathcal{R}_{0,L}}\right)}\leq
E,
\end{equation}
\begin{equation}
\label{Cs14-1d} \left\Vert
F\right\Vert_{L^{\infty}(\mathcal{R}_{0,L}}\leq\varepsilon_1,
\end{equation}
\end{subequations}
where $F$ \textit{is not necessarily analytic}. Let us continue to denote  (even though $\lambda_1>\lambda_2$) by
$\mathcal{R}_{\lambda_1,\lambda_2}$ the subset of
$\mathbb{R}^{n+1}$ enclosed by $S_{\lambda_1}$ and $S_{\lambda_2}$.
In such a way \eqref{Cs4-1} becomes
\begin{equation}\label{Cs15-1}
	\begin{aligned}
	&\int_{B_1}
	v^T(x,\lambda_2\gamma(x))A(x,\lambda_2)u(x,\lambda_2\gamma(x))dx+\\&
	+\mbox{sgn}\left(\lambda_1-\lambda_2\right)\int_{\mathcal{R}_{\lambda_1,\lambda_2}}v^T(x,t)F(x,t)dxdt=\\&=\int_{B_1}
	v^T(x,\lambda_1\gamma(x))A(x,\lambda_1)u(x,\lambda_1\gamma(x))dx.
	\end{aligned}
\end{equation}
We construct $v(x,t;\lambda)$ likewise the Step 2 of Theorem
\ref{stabilita-john}. Consequently, instead of \eqref{Cs8-1},
 for a fixed $\lambda_0$ in $[0,L]$ and setting

\begin{equation*}
\begin{aligned}
 \widetilde{g}(\lambda)&=\int_{B_1}
w^T(x)u(x,\lambda\gamma(x))dx+\\&
+\mbox{sgn}\left(\lambda_0-\lambda\right)\int_{\mathcal{R}_{\lambda,\lambda_0}}v^T(x,t;\lambda)F(x,t)dxdt,
\end{aligned}
\end{equation*}
we have
\begin{equation}\label{Cs15-2}
\widetilde{g}(\lambda)=\int_{B_1}
v^T(x,\lambda_0\gamma(x);\lambda)A(x,\lambda_0)u(x,\lambda_0\gamma(x))dx,
\end{equation}
for $\left|\lambda-\lambda_0\right|<\delta$ and $\lambda\in
[0,L]$. Exactly like the Step 3 of the proof of Theorem
\ref{stabilita-john} we get the estimate
\begin{equation}\label{Cs15-5}
 \left|\widetilde{g}(\lambda+i0)\right|\leq
 C_2M_0e^{|\xi|}E^{1-\vartheta}\varepsilon^{\vartheta}, \quad
 \forall \lambda\in[0,L].
 \end{equation}
Now by \eqref{Cs7-2} and by \eqref{Cs14-1d} we have

\begin{equation}\label{Cs15-5n}
	\begin{aligned}
		\left|\int_{\mathcal{R}_{\lambda,\lambda_0}}v^T(x,t;\lambda)F(x,t)dxdt\right|&\leq
 \varepsilon_1\int_{\mathcal{R}_{\lambda,\lambda_0}}\left|v^T(x,t;\lambda)\right|dxdt\leq \\&\leq
 cLM_0e^{|\xi|}\varepsilon_1,
 \end{aligned}
 \end{equation}
 where $c$ depends on $n$ only. By \eqref{Cs15-5} and \eqref{Cs15-5n} we have

 \begin{equation*}
 \left|\int_{B_1}\gamma(x) w^T(x)u(x,\lambda\gamma(x))dx\right|\leq
 C_2M_0e^{|\xi|}E^{1-\vartheta}\varepsilon^{\vartheta}+
 cLM_0e^{|\xi|}\varepsilon_1.
 \end{equation*}
From now on we may argue like Step 4 of the proof of Theorem \ref{Cs2-teorema-st} and we find
\begin{equation}\label{Cs16-1}
\left|u(x,t)\right|\leq
\frac{C(E+2\varepsilon_2)}{(1-r)^{n+1}}\left|\log\frac{\varepsilon_2}{E+2\varepsilon_2}\right|^{-1}
\end{equation}
for all $x$ in $\overline{\mathcal{R}_{0,L}}\cap
\left(\overline{B_{1-r}}\times \mathbb{R}\right)$, where $C$ depends on $M, L, \rho, c_0$ and $n$ only and
$$\varepsilon_2=E^{1-\vartheta}\varepsilon^{\vartheta}+
 \varepsilon_1.$$

\medskip

\noindent\textbf{2.} The a priori information  \eqref{Cs2-5c} and
 $\gamma$ (compare \eqref{Cs1-2}) occur in the proof of stability Theorem especially to obtain  \eqref{Cs12-3}, while in the other parts of the proof what is needed
 to know about $u$ is only that
$$\left\Vert
u\right\Vert_{L^{\infty}\left(\mathcal{R}_{0,L}\right)}\leq E.$$ Now we prove
that with some further arrangements we may define

\begin{equation}\label{Cs19-1}
\gamma(x)=1-|x|^2
\end{equation}
and, instead of the a priori information  \eqref{Cs2-5c} we require

\begin{equation}\label{Cs19-2}
\left\Vert
u\right\Vert_{C^{1}\left(\overline{\mathcal{R}_{0,L}}\right)}\leq E.
\end{equation}

First of all, we notice that the functions $f_j$ defined like
\eqref{Cs10-2} (with $\gamma$ given by \eqref{Cs19-1}) satisfy to
(we omit the index $j$)

\begin{equation}\label{Cs19-3}
\left\Vert f\right\Vert_{C^{1}\left(B_1\right)}\leq CE,
\end{equation}
where $C$ depends on $L$ only. Since $f=0$ in
$\mathbb{R}^n\setminus B_1$ we have

\begin{equation*}
-i\xi_k\widehat{f}(\xi)=(-1)\int_{B_1}e^{-i\xi\cdot
x}\partial_kf(x)dx,\quad k=1, \cdots,n.
\end{equation*}
By the latter, taking into account  \eqref{Cs19-3} we have

\begin{equation*}
\int_{\mathbb{R}^n}|\xi|^2\left|\widehat{f}(\xi)\right|^2d\xi=(2\pi)^{-n}\int_{B_1}\left|\nabla
f(x)\right|^2dx\leq \overline{C}E^2,
\end{equation*}
where $\overline{C}$ depends on $L$ (and $n$) only.

Therefore we have (recalling \eqref{Cs11-1}), for every $s>0$

\begin{equation*}
\begin{aligned}
 \int_{\mathbb{R}^n}\left|\widehat{f}(\xi)\right|d\xi&=\int_{|\xi|\leq s}\left|\widehat{f}(\xi)\right|^2d\xi+\int_{|\xi|>
s}\left|\widehat{f}(\xi)\right|^2d\xi\leq\\& \leq c\left(C_2
E^{1-\vartheta}\varepsilon^{\vartheta}\right)^2s^ne^{2s}+\frac{\overline{C}^2E^2}{s^2},
\end{aligned}
\end{equation*}
where $c$ depends on $n$ only. 

All in all, by
$$\int_{B_1}\left|
f(x)\right|^2dx=(2\pi)^{-n}\int_{\mathbb{R}^n}\left|\widehat{f}(\xi)\right|^2d\xi,$$
we have

$$\int_{B_1}\left|
f(x)\right|^2dx\leq \overline{C}_3\left(\left(
E^{1-\vartheta}\varepsilon^{\vartheta}\right)^2s^ne^{2s}+E^2s^{-2}\right),$$
where
$$\overline{C}_3=\max\left\{(2\pi)^{n}cC_2^2,(2\pi)^{n}\overline{C}\right\}.$$
Arguing similarly to the proof of \eqref{Cs13-2}, we get 
\begin{equation}\label{Cs21-1}
\int_{B_1}\left| f(x)\right|^2dx\leq \overline{C}_4
\left(E^2+2\varepsilon^2\right)\left|\log\frac{\varepsilon^2}{E^2+2\varepsilon^2}\right|^{-2},
\end{equation}
where $\overline{C}_4=\overline{c}_{\vartheta}\overline{C}_3$ and
$\overline{c}_{\vartheta}$ depends on $\vartheta$.

By applying  Proposition \ref{Stima2-c18} and by
\eqref{Cs21-1}, we have
\begin{equation}\label{Cs21-2}
\left|u(x,t)\right|\leq
\frac{C(E+\sqrt{2}\varepsilon)}{1-r}\left|\log\frac{\varepsilon}{E+\sqrt{2}\varepsilon}\right|^{-\frac{1}{n+1}},\quad
\end{equation}
for every $x\in \overline{\mathcal{R}_{0,L}}\cap
\left(\overline{B_{1-r}}\times \mathbb{R}\right)$, where $C$ depends on $M, L, \rho, c_0$ and $n$ only.

\medskip

\noindent\textbf{3.} Let us examine the main modifications that we should make in the proof of Theorem
\ref{stabilita-john} to deal with the case where the initial surface in the
Cauchy problem is not a portion of the hyperplane $\{t=0\}$. Let
$\varphi\in C^{2}\left(\overline{B_1}\right)$ satisfy 
$\varphi(0)=|\nabla\varphi(0)|=0$. Let us consider the Cauchy problem
\begin{equation}\label{Cs22-1}
\begin{cases}
u_t=\sum_{j=1}^n A_j(x,t)u_{x_j}+A_0(x,t)u+F(x,t), \\
\\
u(x,\varphi(x))=g(x),\quad x\in B_1,
\end{cases}
\end{equation}
Where $A_j$, $j=0,1,\cdots,n$ are analytic functions. We require that the
surface  $\{t=\varphi(x)|x\in B_1\}$ is noncharacteristic.
This is equivalent to require that the "algebraic" system

\begin{equation}\label{Cs22-1extra}
\begin{cases}
u_t(x,\varphi(x))-\sum_{j=1}^n A_j(x,\varphi(x))u_{x_j}(x,\varphi(x))=\widetilde{f}(x), \\
\\
u_{t}(x,\varphi(x))\varphi_{x_i}(x)+u_{x_i}(x,\varphi(x))=g_{x_i}(x),\quad
i=1,\cdots,n,
\end{cases}
\end{equation}
has a unique solution 
$\left(u_{t}(x,\varphi(x)),u_{x_1}(x,\varphi(x)),\cdots,
u_{x_n}(x,\varphi(x))\right)$ for any $g(x)$, where
$$\widetilde{f}(x)=A_0(x,\varphi(x))g(x)+F(x,\varphi(x)).$$
In turn this is equivalent to require the uniqueness of $u_t, u_{x_i}$, $i=0,1,\cdots, n$ as solution to the system 

\begin{equation*}
\begin{cases}
\left(I+\sum_{j=1}^n\varphi_{x_i}(xx) A_j(x,\varphi(x))\right)u_{t}(x,\varphi(x))=G(x), \\
\\
u_{x_i}(x,\varphi(x))=g_{x_i}(x)-u_{t}(x,\varphi(x))\varphi_{x_i}(x),\quad
i=1,\cdots,n,
\end{cases}
\end{equation*}
where 
$$G(x)=\widetilde{f}(x)+\sum_{j=1}^n A_j(x,\varphi(x))g_{x_i}(x).$$
From which we have that $\{t=\varphi(x)|x\in B_1\}$ is a characteristic surface if and only if 

\begin{equation}\label{Cs22-1extra-24}
\det\left(I+\sum_{j=1}^n\varphi_{x_i}(x) A_j(x,\varphi(x))\right)\neq
0, \quad \forall x\in B_1.
\end{equation}

Let us first consider the case in which $$g\equiv 0.$$ By mean the 
Holmgren transformation, \eqref{3-4-67C}, we may assume that
$\varphi$ is strictly convex. For any $\lambda_1,\lambda_2$
positive numbers, let us denote by $\mathcal{S}_{\lambda_1,\lambda_2}$ the subset of $\mathbb{R}^{n+1}$ enclosed by  hyperplanes 
${t=\lambda_1}$, ${t=\lambda_2}$ and the graph of $\varphi$,
let us suppose that $\lambda_1$ and $\lambda_2$ are small enough in such a way that  $\mathcal{S}_{\lambda_1,\lambda_2}$ has a "lens shape"  and let us apply the Green identity. We get 
\begin{equation}\label{Cs22-2extra}
	\begin{aligned}
	&\int_{B_1}
	v^T(x,\lambda_2)A(x,\lambda_2)u(x,\lambda_2\gamma(x))dx+\\&
	+\mbox{sgn}\left(\lambda_1-\lambda_2\right)\int_{\mathcal{S}_{\lambda_1,\lambda_2}}v^T(x,t)F(x,t)dxdt=\\&
	=\int_{B_1} v^T(x,\lambda_1)A(x,\lambda_1)u(x,\lambda_1)dx,
	\end{aligned}
\end{equation}
Let $v(x,t;\lambda)$ the solution to

\begin{equation}\label{Cs22-3extra}
\begin{cases}
v_t=\sum_{j=1}^n \left(A^T_j(x,t)v\right)_{x_j}-A_0^T(x,t)v, \\
\\
v(x,\lambda)=We^{-i\xi\cdot x},\quad x\in B_1.
\end{cases}
\end{equation}
Now we set
$$g(\lambda):=\int_{B_1}
e^{-i\xi\cdot x}W^Tu(x,\lambda\gamma(x))dx$$ and along the lines
of the proof of Theorem \ref{stabilita-john} we obtain an estimate
like  \eqref{Cs3-1} (reader take care of the details).

In the case in which $g$ does not vanish we may reduce
to the previous case by setting

\begin{equation}\label{Cs22-4extra}
\widetilde{u}(x,t)=u(x,t)-g(x).\end{equation} 
 Let us examine the situation in some detail. First, let us assume, for the sake of brevity,
that in \eqref{Cs22-1} we have $F\equiv 0$. Furthermore, we assume that
the solution $u$ of Cauchy problem \eqref{Cs22-1} there exists in an
open set $D$, we assume that $u\in C^{2}\left(\overline{D}\right)$ and that $u$
satisfies the a priori information

\begin{equation}\label{Cs23-3}
\left\Vert u\right\Vert_{C^{2}\left(\overline{D}\right)}\leq E.
\end{equation}
In addition, let us assume that

\begin{equation}\label{Cs23-2}
\left\Vert g\right\Vert_{L^{\infty}(B_1)}\leq \varepsilon.
\end{equation}
We have that $\widetilde{u}$ satisfies

\begin{equation}\label{Cs23-4}
\begin{cases}
\widetilde{u}_t=\sum_{j=1}^n A_j(x,t)\widetilde{u}_{x_j}+A_0(x,t)\widetilde{u}+\widetilde{F}(x,t), \\
\\
\widetilde{u}(x,\varphi(x))=0,\quad x\in B_1,
\end{cases}
\end{equation}
where

\begin{equation}\label{Cs23-5}
\widetilde{F}(x,t)=\sum_{j=1}^n A_j(x,t)g_{x_j}(x)+A_0(x,t)g(x).
\end{equation}
Now by Proposition \ref{Stima1-c16} we have

\begin{equation}\label{Cs23-6}
\left\Vert \nabla g\right\Vert_{L^{\infty}(B_1)}\leq
c\left(\left\Vert \partial^2
g\right\Vert_{L^{\infty}(B_1)}+\left\Vert
g\right\Vert_{L^{\infty}(B_1)}\right)^{\frac{1}{2}}\left\Vert
g\right\Vert_{L^{\infty}(B_1)}^{\frac{1}{2}},
\end{equation}
where $c$ is a positive constant depending on $n$ only. Inequality
\eqref{Cs23-6} allows us to estimate from above the first derivatives of  $g$  in terms of $\varepsilon$ and the a priori bound \eqref{Cs23-3}. Concerning the latter it suffices to recall that
 $g(x)=u(x,\varphi(x))$ and to calcolate the derivatives of $g$
obtaining

\begin{equation}\label{Cs24-1}
\left\Vert \partial^2 g\right\Vert_{L^{\infty}(B_1)}\leq K_1 E,
\end{equation}
where $K_1$ depends on $\left\Vert
\varphi\right\Vert_{C^{2}\left(\overline{B_1}\right)}$. By
\eqref{Cs23-2}, \eqref{Cs23-5}, \eqref{Cs23-6} and \eqref{Cs24-1} we have
\begin{equation}\label{Cs24-3}
\left\Vert
\widetilde{F}\right\Vert_{C^{2}\left(\overline{D}\right)}\leq K_2
\left(E+\varepsilon\right)^{\frac{1}{2}}\varepsilon^{\frac{1}{2}}.
\end{equation}
Finally, taking into account that by \eqref{Cs22-4extra} and
\eqref{Cs23-3} we have
\begin{equation}\label{Cs24-4}
\left\Vert
\widetilde{u}\right\Vert_{C^{2}\left(\overline{D}\right)}\leq K_3 E.
\end{equation}
By using what is obtained in the case $g\equiv 0$, we get, by \eqref{Cs24-3} and \eqref{Cs24-4}, a stability estimate for $\widetilde{u}$ from which immediately follows a stability estimate for
$u$. We invite the reader to write explicitly a stability estimate for the solution $u$ to problem \eqref{Cs22-1} provided the a priori information \eqref{Cs23-3} is satisfied.  

\medskip

\noindent\textbf{4.} Stability estimate \eqref{Cs3-1} is
a \textit{logarithmic estimate} and, while it is still a stability estimate, it is a rather modest estimate. John, in \cite{Joh60},
called "well-behaved" \index{well--behaved problems}the problems for which a H\"{o}lder conditional stability
holds and  "not well-behaved" the problems for which 
the conditional stability is at best of logarithmic type.
This terminology is still in use today. Of course, in order to be able to say
that a class of problems is "well-behaved" or "not well-behaved"
with respect to certain a priori informations, it needs to be shown that
the estimate in question is optimal in that class of problems with
those certain a priori informations.
Concerning Theorem
\ref{Cs2-teorema-st}, the class of problems is the class of the
Cauchy problems for partial differential equations with analytic
coefficients and the a priori bounds concern a finite numbers of derivatives of the solutions. 
Now, with respect to the class of problems and of the
a priori informations that we have considered above, John himself, in \cite{Joh60},
proved that the Cauchy problem is "not well-behaved."
The example constructed by John  concerns the following Cauchy problem  for the wave equation

\begin{equation}\label{onde-john}
\begin{cases}
u_{xx}+u_{yy}-u_{tt}=0, \quad x^2+y^2<1,\quad t\in \mathbb{R}, \\
\\
u=g_0,\quad x\in \partial B_{\rho}\times \mathbb{R},\\
\\
\frac{\partial u}{\partial \nu}=g_1,\quad x\in \partial B_{\rho}\times \mathbb{R}.\\
\end{cases}
\end{equation}
where $\rho<1$ and $\nu$ is the unit outward normal to $\partial
B_{\rho}\times \mathbb{R}$.

More precisely, set
\begin{equation}\label{onde-1john}
\varepsilon=\left\Vert g_0 \right\Vert_{L^{\infty}\left(\partial
B_{\rho}\times \mathbb{R}\right)}+\left\Vert g_1
\right\Vert_{L^{\infty}\left(\partial B_{\rho}\times
\mathbb{R}\right)},\end{equation} John has proved that for every
$m\in \mathbb{N}$ there exists $u\in C^m\left(\overline{B_{1}}\times
\mathbb{R}\right)$ solution to \eqref{onde-john}, where $g_0,g_1$ satisfy \eqref{onde-1john}, such that
$$ \left\Vert u \right\Vert_{C^{m}\left(
\overline{B_{1}}\times \mathbb{R}\right)}=1$$ and such that
$$\left\Vert u \right\Vert_{L^{\infty}\left(
B_{r}\times \mathbb{R}\right)}\geq
C\left|\log\varepsilon\right|^{-\alpha},$$ where $r\in (\rho,1)$, 
$C>0$ and $\alpha>0$ depend on $r$.$\blacklozenge$

\part{CARLEMAN ESTIMATES AND UNIQUE CONTINUATION PROPERTIES}\label{parte3}

\chapter{PDEs with constant coefficients in the principal part}\label{Nirenberg}

\section{Introduction}\label{Nirenberg-intro}

We begin to study the unique continuation properties for \textbf{operators with non
	analytic coefficients}. First we give some definitions. Let $\Omega$ an open connected set of $\mathbb{R}^n$, we say that the linear differential equation  

\begin{equation}\label{corr:18-4-23-1}
	Lu =0 \quad \mbox{in } \Omega,
\end{equation}
enjoys the \textbf{weak unique continuation property} \index{weak unique continuation property}  if for any open subset $\omega$ of $\Omega$,
$$P(x,\partial)u = 0 \quad \mbox{in } \Omega \mbox{ and }  u =0 \mbox{ in } \omega\quad\Longrightarrow u \equiv 0.$$
We say that equation \eqref{corr:18-4-23-1} enjoys the \textbf{strong unique continuation property} \index{strong unique continuation property} if for any point $x_0\in \Omega$ and for any solution $u$ which satisfies

$$\lim_{r\rightarrow 0}r^{-k}\int_{B_{r}(x_0)}u^2=0,\quad \  \forall k\in \mathbb{N},$$
it follows that
$$u\equiv 0,\quad\mbox{in } \Omega.$$
It is obvious that the strong unique continuation property implies the weak unique continuation property.

As we will see later, the weak unique continuation property is strictly relataded to the uniqueness of the Cauchy problem for equation \eqref{corr:18-4-23-1}.

\bigskip

In the present Chapter we consider the linear differential operators whose \textbf{principal part has constant coefficients}. In other words, we will consider the operators

\begin{equation}\label{1-Nire}
Lu=P(D)u+M(x,D)u,
\end{equation}
where $P(D)$ is a differential operator of order $m$ whose coefficients are  constant (real or complex)
 and $M(x,D)$ is a differential operator of order (less or equal to) $m-1$ whose coefficients belong to $L^{\infty}$ and
$$D_j=\frac{1}{i}\partial_j, \quad j=1,\cdots, n.$$ The latter notation is
very convenient in this Chapter because we will be use
extensively  the Fourier transform.

One of the main purposes of this Chapter is to lay the
ground for the Carleman estimates\index{Carleman estimates}, which will be studied more
systematically in the next chapters. These types of estimates
were introduced by Carleman in \cite{Car1} and \cite{Car2} (in 1933
and 1939 respectively). With these estimates a very important qualitative step
is accomplished in the investigation of the  unique continuation properties
for partial differential equations, particularly for the
Cauchy problem. Indeed, by means of the Carleman estimates, one can prove the
unique continuation properties for differential equations with
nonanalytic coefficients. Actually, the estimates proved in \cite{Car1} and
\cite{Car2} involve partial differential equations of two
variables, but the idea introduced by Carleman has revealed to be very
fruitful leading to the development of a technique that constitutes
certainly the most general and powerful tool, though not unique,
for dealing with unique continuation issues.

The Main Theorem  which we will prove here is due to
\textbf{Nirenberg}, (see Theorem \ref{1.14-corollario-sece}),
\cite{Ni}. Subsequently, we will apply such a Theorem to obtain
the weak unique continuation property for the equation

\begin{equation}\label{2-Nire}
\Delta u-b(x)\cdot\nabla u-c(x)u=0.
\end{equation}
where $b=(b_1,\cdots,b_n)\in L^{\infty}(\mathbb{R}^n, \mathbb{C}^n)$,
$c\in L^{\infty}(\mathbb{R}^n, \mathbb{C})$. Moreover, we will illustrate other applications
and relevant features of Theorem \ref{1.14-corollario-sece}.

\section{The Nirenberg Theorem}\label{Teorema-Nire} Let us introduce
and recall some notations.

Let $P(D)$ be the operator 

\begin{equation}\label{1.8-sece}
P(D)=\sum_{|\alpha|\leq m}a_{\alpha}D^{\alpha},
\end{equation}
where $a_{\alpha}\in \mathbb{C}$, for every $\alpha\in \mathbb{N}^n_0$
satisfying $|\alpha|\leq m$. Let
$$P(\xi)=\sum_{|\alpha|\leq m}a_{\alpha}\xi^{\alpha}, \quad \forall \xi\in \mathbb{R}^n,$$
the symbol of $P(D)$. For each multi-index $\alpha$ we denote 
$$P^{(\alpha)}(\xi)=\partial^{\alpha}_{\xi}P(\xi).$$
We set
$$Q_1=\left\{x\in \mathbb{R}^n:\left|x_j\right|<1,
j=1,\cdots,n\right\}.$$

\medskip

We  prove the following

\medskip

\begin{theo}[\textbf{Nirenberg}]\label{1.14-corollario-sece}
	\index{Theorem:@{Theorem:}!- Nirenberg@{- Nirenberg}}
Let $N\in \mathbb{R}^n$, $|N|=1$. Then there exists a  constant $C$, 
depending on $n$ and $m$ only, such that for every $\alpha\in
\mathbb{N}_0^n$ we have
\begin{equation}\label{1.15-sece}
\int_{Q_1}e^{2\tau N\cdot x}\left|P^{(\alpha)}(D)u\right|^2dx\leq
C\int_{Q_1}e^{2\tau N\cdot x}\left|P(D)u\right|^2dx,
\end{equation}
for every $u \in C^{\infty}_{0}(Q_1,\mathbb{C})$ and for every $\tau\in
\mathbb{R}$.
\end{theo}

\medskip

\textbf{Remark 1.} Estimate \eqref{1.15-sece} is a prototype of
the Carleman estimates. Let us notice that in such an estimate there is a
"weight," $e^{2\tau N\cdot x}$ dependent on a parameter $\tau$, and it is very important that such a parameter can be arbitrarily large.

Let us observe, in particular, that the at right--hand side of \eqref{1.15-sece} it occurs the operator $P(D)$  
applied to an arbitrary $u\in C^{\infty}_{0}(Q_1,\mathbb{C})$,
\textit{not} to a solution of some equation. 
$\blacklozenge$

\bigskip

In order to prove Theorem \ref{1.14-corollario-sece} we need
some preliminary results.

\medskip

First of all, let us recall the following one--dimensinal Poincar\'{e} inequality

\begin{equation}
\label{Poinc-28} \int_{-1}^1 |u|^2dt \leq \frac{4}{\pi^2}\int_{-1}^1
\left| u'\right|^2  dt\mbox{,} \quad \forall u \in
C^{\infty}_{0}((-1,1),\mathbb{C}).
\end{equation}

\medskip

Now let us prove 
\begin{lem}
\label{lem1-sece} There exists $C_0>0$ such that for each $\gamma \in
\mathbb{C}$ we have

\begin{equation}
\label{1.1-sece} \int_{-1}^1 |u|^2dt \leq C_0\int_{-1}^1 \left|
u'-\gamma u\right|^2  dt\mbox{,} \quad \forall u \in
C^{\infty}_{0}((-1,1),\mathbb{C}).
\end{equation}
\end{lem}
\textbf{Proof.} Let $\gamma=\alpha+i\beta$, $\alpha,\beta\in
\mathbb{R}$. We have
\begin{equation*}
\begin{aligned}
\left|u'-\gamma u\right|^2&= \alpha^2\left|u\right|^2
+\left|u'-i\beta u\right|^2-2\alpha\Re\left(\left(u'-i\beta
u\right)\overline{u}\right)=\\&
=\alpha^2\left|u\right|^2+\left|\left(ue^{-i\beta
t}\right)'\right|^2-\alpha\left(|u|^2\right)'.
\end{aligned}
\end{equation*}
Hence, as $u \in C^{\infty}_{0}((-1,1),\mathbb{C})$,
taking into account \eqref{Poinc-28}, we get

\begin{equation*}
\begin{aligned}
\int_{-1}^1 \left| u'-\gamma u\right|^2  dt& =\int_{-1}^1 \left(
\alpha^2\left|u\right|^2+\left|\left(ue^{-i\beta
t}\right)'\right|^2\right) dt \geq\\&
\geq\int_{-1}^1\left|\left(ue^{-i\beta t}\right)'\right|^2 dt\geq
\frac{\pi^2}{4}\int_{-1}^1 \left|ue^{-i\beta t}\right|^2
dt=\\&=\frac{\pi^2}{4}\int_{-1}^1 \left|u\right|^2  dt.
\end{aligned}
\end{equation*}

Therefore inequality \eqref{1.1-sece} is proved with
$C_0=\frac{4}{\pi^2}$. $\blacksquare$

\bigskip

Let $a_1,\cdots, a_k\in \mathbb{C}$, $a_k\neq 0$, and let 
\begin{equation}
\label{1.2-sece} p(\eta)=\sum_{j=0}^k a_j\eta^j, \quad \eta\in
\mathbb{C}.
\end{equation}
Let us consider the differential operator
\begin{equation}
\label{1.3-sece} p(D_t)=\sum_{j=0}^k a_jD_t^j,
\end{equation}
where $D_t=\frac{1}{i}\frac{d}{dt}$. Set
$$p^{\prime}(D_t)=\sum_{j=1}^k ja_jD_t^{j-1}.$$ We have
the following

\begin{lem}
\label{lem2-sece} Let $k\in \mathbb{N}$. Then there exists $C_1>0$ depending on
 $k$ only, such that we have

\begin{equation}
\label{1.4-sece} \int_{-1}^1 \left|p^{\prime}(D_t)u\right|^2dt \leq
C_1\int_{-1}^1 \left|p(D_t)u\right|^2dt\mbox{,} \quad \forall u \in
C^{\infty}_{0}((-1,1),\mathbb{C}).
\end{equation}
\end{lem}

\textbf{Proof.} It is not restrictive to assume that $a_k=1$. Let $\gamma_1,\cdots, \gamma_k\in \mathbb{C}$ 
be the roots of the polynomial $p$, we have
$$p(\eta)=\prod_{1\leq j\leq k}\left(\eta-\gamma_j\right).$$
Set
$$p_l(\eta)=\frac{1}{\left(\eta-\gamma_l\right)}\prod_{1\leq j\leq k}\left(\eta-\gamma_j\right), \mbox{ for } l=1, \cdots, k.$$
We have
\begin{equation}
\label{1.5-sece} p^{\prime}(\eta)=\sum_{l=1}^kp_l(\eta), \ \  \forall \eta\in \mathbb{C}.
\end{equation}
Let  $l\in\{1, \cdots, k\}$ be fixed, $u\in
C^{\infty}_{0}((-1,1),\mathbb{C})$ and let us denote
$$v_l=p_l(D_t)u.$$
Let us observe that
$$\left(D_t-\gamma_l\right) v_l=p(D_t)u.$$
Now we apply Lemma \ref{lem1-sece} to $v_l$ (where
$\gamma=i^{-1}\gamma_l$) and we get

\begin{equation}\label{1.6-sece}
\begin{aligned}
\int_{-1}^1 \left|p_l(D_t)u\right|^2dt& =\int_{-1}^1
\left|v_l\right|^2dt \leq\\&\leq C_0\int_{-1}^1 \left|
\left(D_t-\gamma_l\right) v_l\right|^2  dt=\\&=C_0\int_{-1}^1 \left|
p(D_t)u\right|^2  dt.
\end{aligned}
\end{equation}

By \eqref{1.5-sece} and \eqref{1.6-sece} we have 

\begin{equation}\label{1.7-sece}
\begin{aligned}
\int_{-1}^1 \left|p'(D_t)u\right|^2dt& \leq k\sum_{l=1}^k
\int_{-1}^1 \left|p_l(D_t)u\right|^2dt\leq\\&\leq  C_0k^2\int_{-1}^1 \left|
p(D_t)u\right|^2  dt.
\end{aligned}
\end{equation}

Hence, inequality \eqref{1.4-sece} is proved with
$C_1=k^2C_0$. $\blacksquare$

\bigskip

\begin{theo}[\textbf{H\"{o}rmander}]\label{hormander-sece}
	\index{Theorem:@{Theorem:}!- H\"{o}rmander (for constant differential operators)@{- H\"{o}rmander (for constant differential operators)}}
Let $P(D)$ be a differential operator of order $m$ with constant coefficients. Then there exists a constant $C_2$ which depends on $m$
and on $n$ only, such that we have, for any $\alpha\in \mathbb{N}_0^n$,
\begin{equation}
\label{1.9-sece} \int_{Q_1} \left|P^{(\alpha)}(D)u\right|^2dx \leq
C_2\int_{Q_1} \left|P(D)u\right|^2dx, \quad \forall u \in
C^{\infty}_{0}(Q_1,\mathbb{C}).
\end{equation}
\end{theo}
\textbf{Proof.} First of all we prove \eqref{1.9-sece} for
$\alpha=e_j$, $j=1,\cdots,n$. It is not restrective to assume $j=n$. for any $f\in L^2\left(\mathbb{R}^n\right)$ we set
$$\widehat{f}(\xi',x_n)=\mathcal{F}_{\xi'}(f(\cdot,x_n)=\int_{\mathbb{R}^{n-1}}f(x',x_n)e^{-ix'\dot\xi'}dx',
\quad\forall \xi'\in \mathbb{R}^{n-1}.$$ Let $u \in
C^{\infty}_{0}(Q_1,\mathbb{C})$; we have

\begin{equation}
\label{1.10-sece}
\mathcal{F}_{\xi'}\left(P(D)u\right)=P\left(\xi',D_n\right)\widehat{u}(\xi',x_n)
\end{equation}
and by the Parseval identity, we have
\begin{equation}\label{1.11-sece}
	\begin{aligned}
	&\int_{Q_1} \left|P^{(e_n)}(D)u\right|^2dx=\\&
	=\frac{1}{(2\pi)^{n-1}}\int_{\mathbb{R}^{n-1}}d\xi'\int_{-1}^1\left|P^{(e_n)}\left(\xi',D_n\right)\widehat{u}(\xi',x_n)\right|^2dx_n.
	\end{aligned}
\end{equation}
Now we apply  Lemma \ref{1.6-sece} to the operator
$p(D_n)=P\left(\xi',D_n\right)$, where $\xi'\in \mathbb{R}^{n-1}$ is
fixed. We obtain

\begin{equation}
\label{1.12-sece}
\begin{aligned}
&\int_{-1}^1\left|P^{(e_n)}\left(\xi',D_n\right)\widehat{u}(\xi',x_n)\right|^2dx_n\leq\\&\leq 
C_1\int_{-1}^1\left|P\left(\xi',D_n\right)\widehat{u}(\xi',x_n)\right|^2dx_n.
\end{aligned}
\end{equation}
By \eqref{1.11-sece} and \eqref{1.12-sece} we have

\begin{equation*}
\begin{aligned}
\int_{Q_1} \left|P^{(e_n)}(D)u\right|^2dx&
\leq\frac{C_1}{(2\pi)^{n-1}}\int_{\mathbb{R}^{n-1}}d\xi'\int_{-1}^1\left|P\left(\xi',D_n\right)\widehat{u}(\xi',x_n)\right|^2dx_n=\\&
=C_1\int_{Q_1} \left|P(D)u\right|^2dx.
\end{aligned}
\end{equation*}
Since the previous proof can be repeated for any
indices. We have that for each multi--indices $\alpha$ such that
$|\alpha|=1$ the estimate following holds
\begin{equation}
\label{1.13-sece} \int_{Q_1}\left|P^{(\alpha)}(D)u\right|^2dx\leq
C_1\int_{Q_1}\left|P(D)u\right|^2dx, \quad \forall u \in
C^{\infty}_{0}(Q_1,\mathbb{C}).
\end{equation}
By iteration of \eqref{1.13-sece} we get, for any $\alpha \in \mathbb{N}_0^n$,
\begin{equation*}
\int_{Q_1}\left|P^{(\alpha)}(D)u\right|^2dx\leq
C^{|\alpha|}_1\int_{Q_1}\left|P(D)u\right|^2dx, \quad \forall u \in
C^{\infty}_{0}(Q_1,\mathbb{C}).
\end{equation*}
Hence inequality \eqref{1.9-sece} is now proved with
$C_2=C^{|\alpha|}_1$. $\blacksquare$

\bigskip

\textbf{Proof of Theorem \ref{1.14-corollario-sece}.}

Let $u \in C^{\infty}_{0}(Q_1,\mathbb{C})$. Setting $v=e^{\tau N\cdot
x}u$, we obtain
\begin{equation}\label{per-conclusioni}D_ju=e^{-\tau N\cdot
x}\left(D_j+i\tau N_j \right)v.
\end{equation} For each multi-index $\alpha$, we have
$$D^{\alpha}u=e^{-\tau N\cdot x}\left(D+i\tau N\right)^{\alpha}v.$$
Hence
$$e^{\tau N\cdot x}P(D)u=P(D+i\tau N)v, \quad e^{\tau N\cdot x}P^{(\alpha)}(D)u=P^{(\alpha)}(D+i\tau N)v.$$
Therefore by \eqref{1.9-sece} we get 
\begin{equation*}
\begin{aligned}
\int_{Q_1}e^{2\tau N\cdot x}\left|P^{(\alpha)}(D)u\right|^2dx&
=\int_{Q_1}\left|P^{(\alpha)}(D+i\tau N)v\right|^2dx \leq\\& \leq
C_2\int_{Q_1}\left|P(D+i\tau N)v\right|^2dx=\\&
=C_2\int_{Q_1}e^{2\tau N\cdot x}\left|P(D)u\right|^2dx,
\end{aligned}
\end{equation*}
where $C_2$ is the same constant of \eqref{1.9-sece}.

$\blacksquare$

\bigskip

\textbf{Remark 1.}

Since $P(D)$ is an operator of order $m$, we have that there
exists  $\alpha\in\mathbb{N}_0^n$ such that
$P^{(\alpha)}(\xi)=\alpha !a_{\alpha}\neq 0$. Therefore 
\eqref{1.15-sece} gives, in particular,
\begin{equation}\label{Nire-1}
\int_{Q_1}\left|u\right|^2dx\leq
C_3\int_{Q_1}\left|P(D)u\right|^2dx, \quad \forall u \in
C^{\infty}_{0}(Q_1,\mathbb{C}),
\end{equation}
where $$C_3=C_2(\frac{1}{\alpha !}\min\{|a_{\alpha}|: a_{\alpha}\neq
0, |\alpha|=m \})^2.$$

By the proof of Theorem \ref{1.14-corollario-sece} we
observe that, if $M(\xi)$ is a polynomial for which there exists a constant $C_4>0$ such that

\begin{equation}\label{3-Nire}
\frac{\left|M(\xi+i\tau N)\right|^2}{\sum_{|\alpha|\leq
m}\left|P^{(\alpha)}(\xi+i\tau N)\right|^2}\leq C_4, \quad \forall
\xi\in \mathbb{R}^n,\quad \forall \tau \in \mathbb{R}
\end{equation}
then there exists a constant $C$ such that
\begin{equation}\label{4-Nire}
\int_{Q_1}e^{2\tau N\cdot x}\left|M(D)u\right|^2dx\leq
C\int_{Q_1}e^{2\tau N\cdot x}\left|P(D)u\right|^2dx,
\end{equation}
for every $u \in C^{\infty}_{0}(Q_1,\mathbb{C})$ and for every $\tau\in
\mathbb{R}$.

Likewise, if $M_{\tau}(\xi)$ is a polynomial in the variable
$\xi$ depending by the parameter $\tau$ and if 

\begin{equation}\label{3-Nire-14-8-21}
\sup\left\{\frac{\left|M_{\tau}(\xi+i\tau
N)\right|^2}{\sum_{|\alpha|\leq m}\left|P^{(\alpha)}(\xi+i\tau
N)\right|^2}: (\xi,\tau)\in \mathbb{R}^{n+1}\right\}<+\infty,
\end{equation}
then, for a constant $C$ (independent of $\tau$ and $u$), we have

\begin{equation}\label{4-Nire-14-8-21}
\int_{Q_1}e^{2\tau N\cdot x}\left|M_{\tau}(D)u\right|^2dx\leq
C\int_{Q_1}e^{2\tau N\cdot x}\left|P(D)u\right|^2dx,
\end{equation}
for any $u \in C^{\infty}_{0}(Q_1,\mathbb{C})$ and  any $\tau\in
\mathbb{R}$.

For instance, in the case

\begin{equation*}
P(D)=-\left(D_1^2+\cdots+D_n^2\right)=\Delta,
\end{equation*}
we have

\begin{equation*}
P(\xi)=-\left(\xi_1^2+\cdots+\xi_n^2\right),\quad \sum_{|\alpha|\leq
2}\left|P^{(\alpha)}(\xi)\right|^2=|\xi|^4+4\left(n^2+\left|\xi\right|^2\right),
\end{equation*}
from which we easily obtain that \eqref{3-Nire} is satisfied  for all  $N\in\mathbb{R}^n$, $|N|=1$, provided $M(\xi)=\xi_j$ (as well as, of course, for  $M(\xi)=1$). Hence, we have

\begin{equation}\label{5-Nire-14-8-21-1}
\int_{Q_1}\left(|u|^2+|\nabla u|^2+|D^2 u|^2\right) dx\leq
C\int_{Q_1}\left|\Delta u\right|^2dx,
\end{equation}
for every $u \in C^{\infty}_{0}(Q_1,\mathbb{C})$, where $$|D^2
u|^2=\sum_{j,k=1}^n |D^2_{jk} u|^2=\sum_{j,k=1}^n
|\partial^2_{jk} u|^2.$$ Moreover, we have

$$\sum_{|\alpha|\leq
2}\left|P^{(\alpha)}(\xi+i\tau
N)\right|^2=\left(|\xi|^2-\tau^2\right)^2+4\tau^2(\xi\cdot
N)^2+4\left(|\xi|^2+\tau^2\right)+4n^2$$ and, setting
$$M_{0,\tau}(\xi)=\tau, \quad M_{j,\tau}(\xi)=\xi_j, \quad j=1,\cdots, n,$$
 we have (reader check), for $k=1,\cdots, n$,

\begin{equation*}
\sup\left\{\frac{\left|M_{k,\tau}(\xi+i\tau
N)\right|^2}{\sum_{|\alpha|\leq 2}\left|P^{(\alpha)}(\xi+i\tau
N)\right|^2}: \quad (\xi,\tau)\in \mathbb{R}^{n+1}\right\}<+\infty. 
\end{equation*}
Hence,
\begin{equation}\label{5-Nire-14-8-21-2}
\int_{Q_1}e^{2\tau N\cdot x}\left(\tau^2|u|^2+|\nabla u|^2\right)
dx\leq C\int_{Q_1}e^{2\tau N\cdot x}\left|\Delta u\right|^2dx,
\end{equation}
for every $u \in C^{\infty}_{0}(Q_1,\mathbb{C})$ and for every $\tau\in
\mathbb{R}$. Trivially, also the following estimate holds 

\begin{equation}\label{5-Nire}
\int_{Q_1}e^{2\tau N\cdot x}\left(|u|^2+|\nabla u|^2\right) dx\leq
C\int_{Q_1}e^{2\tau N\cdot x}\left|\Delta u\right|^2dx,
\end{equation}
for every $u \in C^{\infty}_{0}(Q_1,\mathbb{C})$ and for every $\tau\in
\mathbb{R}$.

In what follows, we will exploit estimate \eqref{5-Nire} to prove
some unique continuation property for the equation $$\Delta
U=b(x)\cdot\nabla U+c(x)U,$$ with $b=(b_1,\cdots,b_n)\in
L^{\infty}(\mathbb{R}^n)$, $c\in L^{\infty}(\mathbb{R}^n)$.

As it will be clear later on, the aforesaid
unique continuation property results could be derived with a
slightly less effort by using \eqref{5-Nire-14-8-21-2}
instead of \eqref{5-Nire}. However, part of the arguments that
we will use employing \eqref{5-Nire} can be extended to differential operators which are more general  and this, in a certain
sense, will repay us for the greater effort we will put into using \eqref{5-Nire}. $\blacklozenge$

\bigskip

\section[Application of the Nirenberg Theorem to the Laplace operator]{Application of the Nirenberg Theorem to the Laplace operator}\label{UCP-Delta}
	In this Section we will apply estimate \eqref{5-Nire} to obtain the weak unique continuation property and \textbf{the uniqueness} for the Cauchy problem to equation \eqref{2-Nire}.

The steps to be done are fairly numerous and, in order to highlight the
key points, we proceed gradually. First,
we warn that we should not confuse the solution of equation \eqref{2-Nire} 
with $u$ in inequality \eqref{5-Nire}. Now, let us dwell on inequality
\eqref{5-Nire} and we notice that, since
 $C^{\infty}_{0}(Q_1,\mathbb{C})$ is dense in $H^2_0(Q_1,\mathbb{C})$, estimate \eqref{5-Nire} holds true
also for any $u \in H^2_0(Q_1,\mathbb{C})$.

Hence, we have 
\begin{equation}\label{10-Nire}
\int_{Q_1}e^{2\tau N\cdot x}\left(|u|^2+|\nabla u|^2\right) dx\leq
C\int_{Q_1}e^{2\tau N\cdot x}\left|\Delta u\right|^2dx,
\end{equation}
for every $u \in H^{2}_{0}(Q_1,\mathbb{C})$ and for every $\tau\in
\mathbb{R}$.

Actually, it would not be difficult to prove that 
\eqref{5-Nire} holds for each $u \in H^1_0(Q_1,\mathbb{C})$ which satisfies
 $\Delta u\in L^2(Q_1,\mathbb{C})$. However, at least
for the time being, let us omit further consideration on this
point. 
Let us recall that (Theorem \ref{reg-interno}):
\begin{prop}\label{regolarit-28}
Let $\Omega$ be an open set of $\mathbb{R}^n$, $f\in
L^2(\Omega,\mathbb{C})$ and let $U\in H^1(\Omega,\mathbb{C})$ satisfy
\begin{equation}\label{6-Nire}
\int_{\Omega}\nabla U\cdot \nabla \varphi dx=-\int_{\Omega}f \varphi
dx, \quad \forall\varphi \in  H^1_0(\Omega,\mathbb{C}),
\end{equation}
then $U\in H_{\mbox{loc}}^2(\Omega,\mathbb{C})$.
\end{prop}

\bigskip

Let $b=(b_1,\cdots,b_n)\in L^{\infty}(\mathbb{R}^n;\mathbb{R}^n)$, $c\in
L^{\infty}(\mathbb{R}^n)$ and $f\in L^2(\mathbb{R}^n)$.

Let us start to consider the following Cauchy problem. Let
$$h(x')=1-\sqrt{1-|x'|^2},$$

$$\Lambda=\left\{(x',x_n)\in B_1'\times\mathbb{R}: h(x')<x_n<1
\right\}$$ and

$$\Gamma=\left\{(x',h(x')):x'\in B_1'
\right\}.$$

We say that $U\in H^2(\Lambda)$ is a solution of the
Cauchy problem
\begin{equation}\label{1-CauNire}
\begin{cases}
\Delta U=b(x)\cdot\nabla U+c(x)U+f(x), \quad \mbox{ in } \Lambda, \\
\\
U=0, \quad\mbox{ on } \Gamma,\\
\\
\frac{\partial U}{\partial \nu}=0, \quad \mbox{ on } \Gamma,
\end{cases}
\end{equation} provided $U\in H^2(\Lambda)$ and

\begin{equation}\label{2-CauNire}
\begin{cases}
\Delta U=b(x)\cdot\nabla U+c(x)U+f(x), \quad \mbox{ in } \Lambda \\
\\
U\Psi\in H_0^2(\Lambda) , \quad\forall\Psi\in C^{\infty}(\mathbb{R}^n),\quad\mbox{supp}\Psi\subset \mathbb{R}^{n-1}\times (-\infty,1) .\\
\end{cases}
\end{equation}
Let us observe that in formulation \eqref{2-CauNire}, we express the conditions $U=\frac{\partial U}{\partial \nu}=0$ on $\Gamma$ as "$U\Psi\in H_0^2(\Lambda)$, for every
$\Psi\in C^{\infty}(\mathbb{R}^n)$, such that $\mbox{supp}\Psi\subset
\mathbb{R}^{n-1}\times (-\infty,1)$". Another way to express correctly these initial conditions is through the definition of the traces. In the
case of initial surface $\Gamma$ that we are considering, both the formulations are equivalent (due to the regularity of
$\Gamma$), however formulation \eqref{2-CauNire}
is more elementary because it allows us to dispense with the
notion of the trace.

\medskip

Set

\begin{equation}\label{coeff-CauNire}
K=\left\Vert b\right\Vert_{L^{\infty}(\Lambda)}+\left\Vert
c\right\Vert_{L^{\infty}(\Lambda)}, \quad\quad
\varepsilon=\left\Vert f\right\Vert_{L^{2}(\Lambda)}.
\end{equation}

\medskip

Our goal is to find an estimate that, roughly speaking, tells us that if $\varepsilon$ "is small" then
$\left\Vert U\right\Vert_{L^2(\Lambda)}$ "is small" and tell us that if $\varepsilon=0$ (thus $f\equiv 0$)
then $U\equiv 0$

\medskip
Let us start by considering the simple case in which,  in
\eqref{1-CauNire}, $K$ and $\varepsilon$ are  zero. In such a case
we have

\begin{equation}\label{3-CauNire}
\begin{cases}
\Delta U=0, \quad \mbox{ in } \Lambda, \\
\\
U\Psi\in H_0^2(\Lambda) , \quad\forall\Psi\in C^{\infty}(\mathbb{R}^n),\quad\mbox{supp}\Psi\subset \mathbb{R}^{n-1}\times (-\infty,1) .\\
\end{cases}
\end{equation}
Let $\delta\in (0,\frac{1}{3})$ and let $\zeta\in
C^{\infty}(\mathbb{R})$ satisfy $0\leq \zeta\leq 1$, (Figure 10.1)
\begin{figure}\label{cauchy1}
	\centering
	\includegraphics[trim={0 0 0 0},clip, width=10cm]{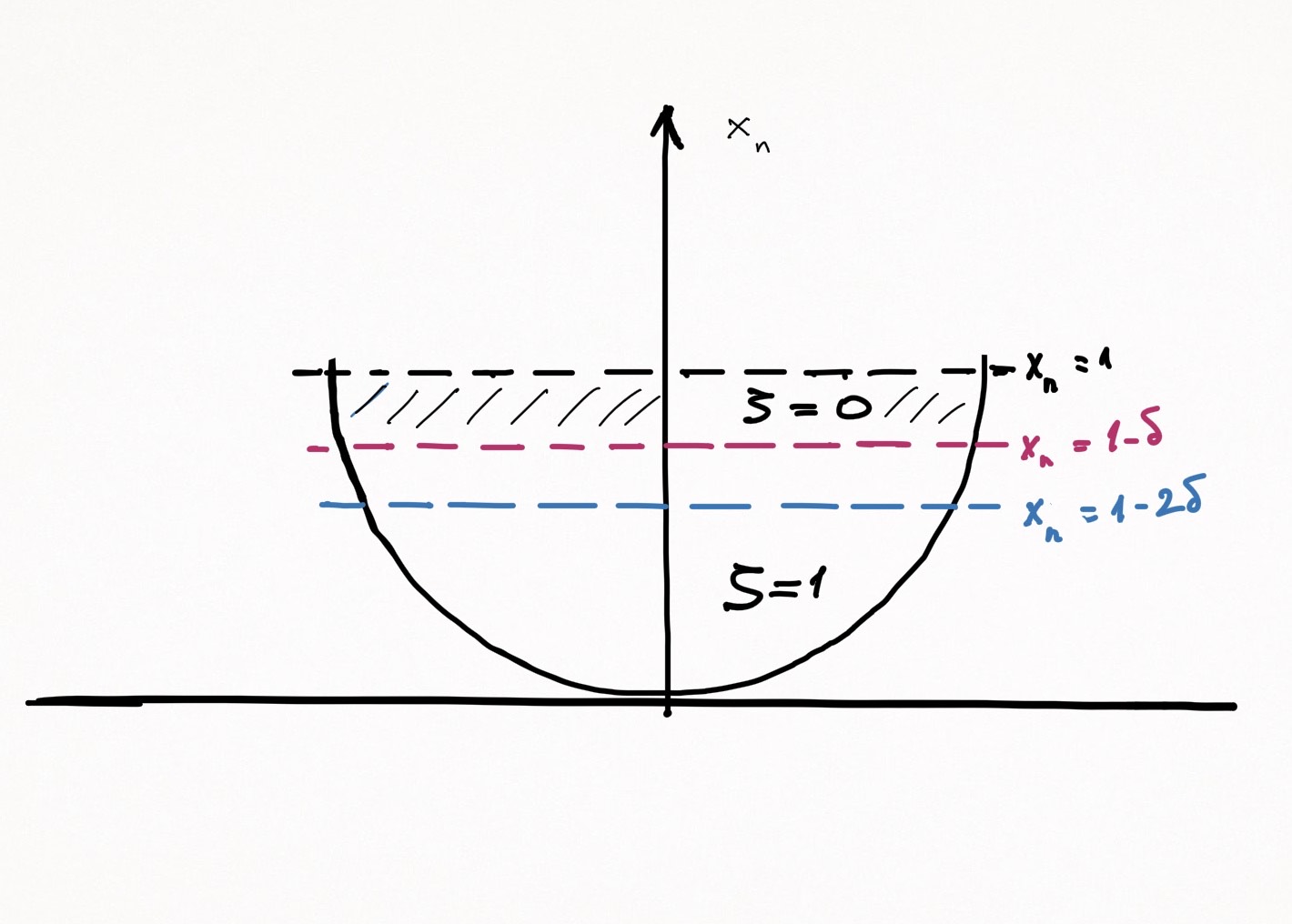}
	\caption{}
\end{figure}

$$\zeta(x_n)=1,\mbox{ for } x_n\leq 1-2\delta,\quad \zeta(x_n)=0\quad\mbox{for }1-\delta\leq
x_n<1$$ and

$$\left|\zeta'(x_n)\right|\leq c\delta^{-1} \mbox{ e }\left|\zeta''(x_n)\right|\leq c\delta^{-2}, \mbox{ for } 1-2\delta\leq
x_n\leq 1-\delta,$$ where $c$ is a constant.

Let us extend $U$ to zero in $\left\{(x',x_n)\in
B_1'\times\mathbb{R}: -1<x_n<h(x') \right\}$, in such a way that this extension
belongs to $H^2(Q_1)$, we continue to denote by $U$ such an extension. We have
$\zeta U\in H_0^2(Q_1)$. Let us apply \eqref{10-Nire} to
$$u=\zeta(x_n) U(x),$$ where $$N=-e_n, \ \ \mbox{and } \tau>0.$$

\noindent We have

$$\Delta (\zeta U)=\zeta(x_n)\Delta U+2\zeta'(x_n)\partial_n U+ \zeta''(x_n) U$$
as $\Delta U=0$, we get

\begin{equation}\label{10n-Nire}
\left|\Delta (\zeta U)\right\vert\leq
\chi_{(1-2\delta,1-\delta)}(x_n)\left(2c\delta^{-1}\left|\partial_n
U\right\vert+ c\delta^{-2}\left|U\right\vert\right).
\end{equation}
Hence by \eqref{10-Nire} and \eqref{10n-Nire} we get

\begin{equation}\label{11-Nire}
\begin{aligned}
&\int_{Q_1}e^{-2 \tau x_n}\left(|\zeta U|^2+|\nabla (\zeta
U)|^2\right) dx \leq C\int_{Q_1}e^{-2 \tau x_n}\left| \Delta (\zeta
U)\right|^2dx \leq\\& \leq
C'\delta^{-4}\int_{Q_1\cap\{1-2\delta<x_n<1-\delta\}}e^{-2 \tau
x_n}\left(|U|^2+|\nabla U|^2\right)dx\leq\\& \leq C'e^{-2 \tau
(1-2\delta)}\delta^{-4}\int_{Q_1\cap\{1-2\delta<x_n<1-\delta\}}\left(|U|^2+|\nabla
U|^2\right)dx\leq\\& \leq C'e^{-2 \tau
(1-2\delta)}\delta^{-4}\left\Vert U\right\Vert_{H^1(Q_1)}^2,
\end{aligned}
\end{equation}
for every $\tau>0$.

Now we have, trivially, for every $\tau>0$

\begin{equation}\label{12-Nire}
\begin{aligned}
&\int_{Q_1}e^{-2 \tau x_n}\left(|\zeta U|^2+|\nabla (\zeta
U)|^2\right) dx \geq \\&\geq  e^{-2 \tau
(1-3\delta)}\int_{\Lambda\cap\{x_n<1-3\delta\}}\left(|U|^2+|\nabla
U|^2\right)dx.
\end{aligned}
\end{equation}

\smallskip

\noindent By \eqref{11-Nire} and \eqref{12-Nire} we obtain

\begin{equation*}
	\begin{aligned}z
&\int_{\Lambda\cap\{x_n<1-3\delta\}}\left(|U|^2+|\nabla
U|^2\right)dx\leq \\&\leq C'e^{-2 \tau \delta}\delta^{-4}\left\Vert
U\right\Vert_{H^1(Q_1)}^2\rightarrow 0, \quad\mbox{as }
\tau\rightarrow+\infty.
\end{aligned}
\end{equation*}
From which we have $U=0$ in $\Lambda\cap\{x_n<1-3\delta\}$ and, as  $\delta$ is arbitrary, we get $U=0$ in $\Lambda$.

\medskip

In the general case, particularly when $b$ or $c$ are not zero,
one would like to argue in a similar manner, but one encounters an obstacle
due to the fact that in expanding $|\Delta(\zeta
U)|$, besides the terms present in
\eqref{10n-Nire}, new terms, depending on $U$ and $\nabla U$, arise. To overcome this difficulty, we first perform a rescaling of estimate \eqref{10-Nire}

For any $r>0$, set

$$Q_r=\left\{x\in \mathbb{R}^n: \quad\left|x_j\right|<r,\quad
j=1,\cdots,n\right\}.$$

\begin{prop}\label{dilat-Nire}
Let $N\in \mathbb{R}^n$ ($|N|=1$). There exists a constant $C>0$ so that
\begin{equation}\label{13-Nire}
\int_{Q_r}e^{2\tau N\cdot x}\left(|u|^2+r^2|\nabla u|^2\right)
dx\leq Cr^4\int_{Q_r}e^{2\tau N\cdot x}\left|\Delta u\right|^2dx,
\end{equation}
for every $r>0$, for every $u \in H^{2}_{0}(Q_r)$ and for every $\tau\in
\mathbb{R}$.
\end{prop}
\textbf{Proof.}

Let $u \in H^{2}_{0}(Q_r)$ and let us denote
$$\widetilde{u}(y)=u(ry),\quad\forall y\in Q_1.$$
We apply \eqref{10-Nire} to $\widetilde{u}$, replacing there $\tau$ by $\tau r$. We get

\begin{equation*}
\int_{Q_1}e^{2\tau r N\cdot y}\left(|u(ry)|^2+r^2|(\nabla
u)(ry)|^2\right) dy\leq Cr^4\int_{Q_1}e^{2\tau r N\cdot
y}\left|(\Delta u)(ry)\right|^2dy.
\end{equation*}
Now, by performing the change of variables $y=r^{-1}x$, we have

\begin{equation*}
\int_{Q_r}e^{2\tau N\cdot x}\left(|u(x)|^2+r^2|\nabla u(x)|^2\right)
r^{-n}dx\leq Cr^4\int_{Q_r}e^{2\tau  N\cdot x}\left|\Delta u(x)\right|^2r^{-n}dy.
\end{equation*}
By last inequality we immediately have \eqref{13-Nire}. $\blacksquare$

\medskip

Now let us go back to Cauchy problem \eqref{1-CauNire}. Let $r\in
(0,1)$ be a number to be choosen later, set
$$\rho=\frac{r^2}{2}, \quad\quad R=\sqrt{2\rho-\rho^2}.$$
It is easy to check that $\rho, R<r$ hence

$$\Lambda_{\rho}:=\left\{(x',x_n)\in B_R'\times\mathbb{R}:
h(x')<x_n<\rho \right\}\subset Q_r.$$ Let $\delta\in
(0,\frac{1}{3})$ and $\zeta\in C^{\infty}(\mathbb{R})$ satisfy  $0\leq
\zeta\leq 1$,

$$\zeta(x_n)=1,\mbox{ for } x_n\leq \rho(1-2\delta),\quad \zeta(x_n)=0\quad\mbox{for }\rho(1-\delta) \leq
x_n<\rho$$ and

$$\left|\zeta'(x_n)\right|\leq c(\delta\rho)^{-1} \mbox{ and }\left|\zeta''(x_n)\right|\leq c(\delta\rho)^{-2}, \mbox{ for } \rho(1-2\delta)\leq
x_n\leq \rho(1-\delta),$$ where $c$ is a constant.

Now we extend $U$ to zero in $\left\{(x',x_n)\in
B_R'\times\mathbb{R}: -\rho<x_n<h(x') \right\}$, this extension belongs to $H^2(Q_r)$, we continue to denote it by $U$, it turns out that
$\zeta U\in H_0^2(Q_r)$. Let us prepare to apply the estimate
\eqref{13-Nire} to $$u=\zeta(x_n) U(x),$$ where $N=-e_n$ and $\tau>0$.

\begin{equation*}
\begin{aligned}
\Delta (\zeta U)& = \zeta(x_n)\Delta U+2\zeta'(x_n)\partial_n U+
\zeta''(x_n) U=\\& =\zeta(x_n)\left(b\cdot\nabla U+cU+f\right)+\\& +
2\zeta'(x_n)\partial_n U+ \zeta''(x_n) U.
\end{aligned}
\end{equation*}
From which we have, taking into account \eqref{coeff-CauNire},

\begin{equation}\label{15-Nire}
\begin{aligned}
\left|\Delta (\zeta U)\right\vert& \leq
K\zeta\left(\left|U\right|+\left|\nabla
U\right|\right)+\zeta\left|f\right|+\\&
+c(\rho\delta)^{-2}\chi_{I}(x_n)(\left|\partial_n
U\right\vert+ \left|U\right|),
\end{aligned}
\end{equation}
where we set $I=(\rho(1-2\delta),\rho(1-\delta))$.

Let us denote

\begin{equation*}
J_1= \int_{Q_r}e^{-2 \tau x_n}
\zeta^2\left(\left|U\right|^2+\left|\nabla U\right|^2\right)dx,
\end{equation*}
and
\begin{equation*}
J_2=\int_{Q_r\cap\{\rho(1-2\delta)<x_n<\rho(1-\delta)\}}e^{-2 \tau
x_n}\left(|U|^2+|\nabla U|^2\right)dx.
\end{equation*}
Hence by \eqref{13-Nire} and \eqref{15-Nire} we have (recall that $C$ denotes a constant that
may change from line to line)

\begin{equation}\label{16-Nire}
\begin{aligned}
&\int_{Q_r}e^{-2 \tau x_n}\left(|\zeta U|^2+r^2|\nabla (\zeta
U)|^2\right) dx\leq \\& \leq Cr^4\int_{Q_r}e^{-2 \tau x_n}\left| \Delta
(\zeta U)\right|^2dx \leq\\& \leq C K^2 r^4 J_1+
Cr^4\left(\rho\delta\right)^{-4} J_2+\\&+ Cr^4\int_{Q_r}e^{-2 \tau
x_n}\zeta^2\left|f\right|^2dx,
\end{aligned}
\end{equation}
for every $\tau>0$.

At this point it is useful to note that compared to 
\eqref{11-Nire}, here we have two new terms: $J_1$ and
$r^4\int_{Q_r}e^{-2 \tau x_n}\zeta^2\left|f\right|^2dx$. Recalling the
second equality in \eqref{coeff-CauNire}, we estimate from above the last term in \eqref{11-Nire} as follows

\begin{equation}\label{17-Nire}
\begin{aligned}
r^4\int_{Q_r}e^{-2 \tau x_n}\zeta^2\left|f\right|^2dx\leq
r^4\int_{Q_r}\left|f\right|^2dx=r^4\varepsilon^2,
\end{aligned}
\end{equation}
for every $\tau>0$.

Before considering the term $J_1$, let us estimate from above $J_2$
basically in the same way as done in \eqref{11-Nire}. We have

\begin{equation}\label{18-Nire}
\begin{aligned}
J_2& =\int_{Q_r\cap\{\rho(1-2\delta)<x_n<\rho(1-\delta)\}}e^{-2 \tau
x_n}\left(|U|^2+|\nabla U|^2\right)dx\leq\\& \leq Ce^{-2 \tau\rho
(1-2\delta)}\int_{Q_r\cap\{\rho(1-2\delta)<x_n<\rho(1-\delta)\}}\left(|U|^2+|\nabla
U|^2\right)dx\leq\\& \leq Ce^{-2 \tau\rho (1-2\delta)}\left\Vert
U\right\Vert_{H^1(Q_r)}^2,
\end{aligned}
\end{equation}
for every $\tau>0$.

Hence, by \eqref{16-Nire}--\eqref{18-Nire}, we get

\begin{equation}\label{19-Nire}
\begin{aligned}
&\int_{Q_r}e^{-2 \tau x_n}\left(|\zeta U|^2+r^2|\nabla (\zeta
U)|^2\right) dx\leq C K^2 r^4 J_1+
\\&+C\left(r\rho^{-1}\delta^{-1}\right)^{4} e^{-2 \tau\rho
(1-2\delta)}\left\Vert U\right\Vert_{H^1(Q_r)}^2+\\&
+Cr^4\varepsilon^2,
\end{aligned}
\end{equation}

\smallskip

\noindent for every $\tau>0$.

By the manner in which $\zeta$ is defined, we have trivially (recall $r<1$)

\begin{equation}\label{20-Nire}
\begin{aligned}
&\int_{Q_r}e^{-2 \tau x_n}\left(|\zeta U|^2+r^2|\nabla (\zeta
U)|^2\right) dx\geq \\&\geq  r^2\int_{Q_r\cap\{x_n<\rho(1-2\delta)\}}e^{-2 \tau
x_n}\left(|U|^2+|\nabla U|^2\right) dx,
\end{aligned}
\end{equation}
for every $\tau>0$.

Concerning $J_1$, let us observe

\begin{equation}\label{21-Nire}
\begin{aligned}
&J_1= \int_{Q_r\cap\{x_n<\rho(1-2\delta)\}}e^{-2 \tau x_n}
\zeta^2\left(\left|U\right|^2+\left|\nabla U\right|^2\right)dx+
\\&
+\int_{Q_r\cap\{\rho(1-2\delta)<x_n<\rho(1-\delta)\}}e^{-2 \tau x_n}
\zeta^2\left(\left|U\right|^2+\left|\nabla U\right|^2\right)dx
\leq\\&\leq \int_{Q_r\cap\{x_n<\rho(1-2\delta)\}}e^{-2 \tau x_n}
\left(\left|U\right|^2+\left|\nabla U\right|^2\right)dx+\\&+e^{-2
\tau\rho (1-2\delta)}\left\Vert U\right\Vert_{H^1(Q_r)}^2,
\end{aligned}
\end{equation}
for every $\tau>0$.

\noindent Now by \eqref{19-Nire}--\eqref{21-Nire}
we have

\begin{equation*}
\begin{aligned}
r^2\left(1-C K^2 r^2\right)\int_{Q_r\cap\{x_n<\rho(1-2\delta)\}}e^{-2
\tau x_n}&\left(|U|^2+|\nabla U|^2\right) dx  \leq
Cr^4\varepsilon^2+\\& +Cr^4\left(\rho^{-1}\delta^{-1}\right)^{4} e^{-2
\tau\rho (1-2\delta)}\left\Vert U\right\Vert_{H^1(Q_r)}^2,
\end{aligned}
\end{equation*}
for every $\tau>0$.

 \noindent Now, let us choose $r=r_0<1$ satisfying $1-C K^2 r_0^2\geq
 \frac{1}{2}$ (here, recall that $C$ does not depends by $r$) and denoting by $\rho_0$ and $R_0$ the values of  $\rho$ and
 $R$ correspondingly to this choice of $r$ we get

\begin{equation}\label{22-Nire}
\begin{aligned}
\int_{Q_{r_0}\cap\{x_n<\rho_0(1-2\delta)\}}e^{-2 \tau
x_n}&\left(|U|^2+|\nabla U|^2\right) dx \leq C\varepsilon^2+\\& +C
e^{-2 \tau\rho_0 (1-2\delta)}\left\Vert
U\right\Vert_{H^1(\Lambda_{\rho_0})}^2,
\end{aligned}
\end{equation}
for every $\tau>0$, where $C$ depends on $K$ and $\delta$ only.

\noindent We have, trivially

\begin{equation*}
\begin{aligned}
&\int_{Q_{r_0}\cap\{x_n<\rho_0(1-2\delta)\}}e^{-2 \tau
x_n}\left(|U|^2+|\nabla U|^2\right) dx \geq \\& \geq
\int_{Q_{r_0}\cap\{x_n<\rho_0(1-3\delta)\}}e^{-2 \tau
x_n}\left(|U|^2+|\nabla U|^2\right) dx\geq \\& 
\\&
\geq e^{-2 \tau\rho_0
(1-3\delta)}\left\Vert
U\right\Vert^2_{H^1(\Lambda_{\rho_0(1-3\delta)})},
\end{aligned}
\end{equation*}
for every $\tau>0$.

\noindent By the last obtained estimate and by \eqref{22-Nire} we get

\begin{equation}\label{23-Nire}
\begin{aligned}
\left\Vert U\right\Vert^2_{H^1(\Lambda_{\rho_0(1-3\delta)})} \leq
C\left[e^{2\tau\rho_0(1-3\delta)}\varepsilon^2+e^{-2\tau\rho_0
\delta}\left\Vert U\right\Vert_{H^1(\Lambda_{\rho_0})}^2\right],
\end{aligned}
\end{equation}
for every $\tau>0$.

Now, let us observe that if $\varepsilon=0$, and if  $\tau$ goes to
$+\infty$, then by \eqref{23-Nire} we get $U=0$ in
$\Lambda_{\rho_0(1-3\delta)}$. Instead, if $\varepsilon>0$ by
\eqref{23-Nire} we can derive a \textbf{stability estimate}
by choosing appropriately $\tau$.

The choice of $\tau$ is driven
by the idea of " balancing" the two right hand addends of \eqref{23-Nire} (or
also, minimize with respect to $\tau$ the right-hand member of
\eqref{23-Nire}). For this purpose it is convenient to rearrange the
inequality and set

$$E=\left\Vert U\right\Vert_{H^1(\Lambda_{\rho_0})}.$$

\noindent By \eqref{23-Nire} we have, trivially 

\begin{equation}\label{24-Nire}
\left\Vert U\right\Vert^2_{H^1(\Lambda_{\rho_0(1-3\delta)})} \leq
C\left[e^{2\tau\rho_0(1-3\delta)}\varepsilon^2+e^{-2\tau\rho_0
\delta}\left(E^2+\varepsilon^2\right)\right],
\end{equation}
for every $\tau>0$.

\noindent Let us choose
$$\tau=\tau_0=\frac{1}{2\rho_0(1-2\delta)}\log\left(\frac{E^2+\varepsilon^2}{\varepsilon^2}\right).$$
We get

$$e^{2\tau_0\rho_0(1-3\delta)}\varepsilon^2=e^{-2\tau_0\rho_0
\delta}\left(E^2+\varepsilon^2\right)=\left(\varepsilon^2\right)^{\mu(\delta)}\left(E^2+\varepsilon^2\right)^{1-\mu(\delta)},$$
where
$$\mu(\delta)=\frac{\delta}{1-2\delta}.$$
Hence, by \eqref{24-Nire} we have the following \textbf{stability estimate}

\begin{equation}\label{24bis-Nire}
\left\Vert U\right\Vert_{H^1(\Lambda_{\rho_0(1-3\delta)})} \leq
2^{3/2}C\varepsilon^{\mu(\delta)}(E+\varepsilon)^{1-\mu(\delta)}.
\end{equation}

\bigskip

It is obvious that \textbf{estimate \eqref{24bis-Nire} implies (putting there $\varepsilon=0$) the uniqueness to Cauchy problem \eqref{1-CauNire}}, however this estimate 
says something more about the Cauchy problem. Precisely
tells us that if we have a bound of the norm $H^1$ of the solution $U$, then we can
estimate the error on the solution $U$ starting from the error on the datum $f$.

\medskip

However, let us leave out the stability issue and we return to the uniqueness question for the Cauchy problem. So far we have obtained a local uniqueness result in the case where the initial surface is a semisphere (the graph of $h$).
Recalling that the operator $\Delta$ is invariant under the rotations,
the result of local uniqueness that we obtained can be
easily extend to $\mathcal{U}\cap \Lambda$, where $\mathcal{U}$ is a neighborhood of the graph of $h$. Actually, by further exploiting the particularity
of the operator $\Delta$ we can easily obtain an intermediate result starting from which it will be easy to reach the global \textbf{uniqueness for the Cauchy problem for quite general initial surfaces}.
  
Such an intermediate result is proved in the following Proposition.

\begin{prop}\label{25nn-Nire}
Let $K, R_2, R_1>0$, $R_1<R_2$. Let us suppose that $U\in
H^2(B_{R_2}(0)\setminus \overline{B_{R_1}})$ satisfy what follows

\begin{equation}\label{26nn-Nire}
\left|\Delta U\right|\leq K\left(\left|\nabla
U\right|+\left|U\right|\right), \quad \mbox{ in }
B_{R_2}\setminus \overline{B_{R_1}}
\end{equation}

\begin{equation}\label{NewC-Nire}
U=0,\quad \frac{\partial U}{\partial \nu}=0, \mbox{ on } \partial
B_{R_2},
\end{equation}
(i.e. $U\Psi\in H_0^2(B_{R_2}\setminus \overline{B_{R_1}})$
for every $\Psi\in H^2(\mathbb{R}^n)$ whose support is contained in $\mathbb{R}^n\setminus\overline{B_{R_1}}$)
then
\begin{equation}\label{28nn-Nire}
 U=0\quad\mbox{in } B_{R_2}\setminus \overline{B_{R_1}}.
\end{equation}
\end{prop}
\textbf{Proof.} It suffices to reduce by dilation to the case
in which $R_2=1$ and observe that to obtain \eqref{24bis-Nire}
(with $\varepsilon=0$) we may just use
\eqref{26nn-Nire} instead of the equation $\Delta U=b\cdot \nabla U+c
U$. Next, starting from \eqref{24bis-Nire} (with
$\varepsilon=0$) and using the invariance of $\Delta$ with respect to the
rotations, we immediately obtain that there exists $r\in (0,R_2)$ such that
$U=0$ in $B_{R_2}\setminus \overline{B_{R_2-r}}$. From which
iterating the obtained result, we get that $U\equiv 0$.
$\blacksquare$

\bigskip

\textbf{Remark 1.} Let us observe that to require that
\begin{equation}\label{26nn-Nire-30-4-23-1}
	\left|\Delta U\right|\leq K\left(\left|\nabla
	U\right|+\left|U\right|\right), \quad \mbox{ in }
	\Omega,
\end{equation}
where $\Omega$ is an open set of $\mathbb{R}^n$ is equivalent to require that there exist \\ $b\in L^{\infty}\left(\Omega; \mathbb{C}^n\right)$ and $c\in L^{\infty}\left(\Omega;\mathbb{C}\right)$ such that $U$ is a a solution of the equation 

\begin{equation}\label{26nn-Nire-30-4-23-2}
	\Delta U=b(x)\cdot \nabla U+c(x)U, \quad \mbox{ in }
	\Omega.
\end{equation}
Indeed, it is obvious that if $U$ satisfies \eqref{26nn-Nire-30-4-23-2} then $U$ satisfies \eqref{26nn-Nire-30-4-23-1} with

$$K=\max\left\{\left\Vert b\right\Vert_{L^{\infty}\left(\Omega; \mathbb{C}^n\right)}, \left\Vert c\right\Vert_{L^{\infty}\left(\Omega;\mathbb{C}\right)} \right\}.$$

Conversely, if $U$ satisfies \eqref{26nn-Nire-30-4-23-1}, then we define

\begin{equation*}
	b(x)=\begin{cases}
		\frac{\overline{\nabla U(x)}(\Delta U(x))}{|\nabla U(x)|^2+|U(x)|^2}, \ \mbox{ for } \ |\nabla U(x)|^2+|U(x)|^2 >0, \\
		\\
		0 , \quad \ \ \  \mbox{for } \ |\nabla U(x)|^2+|U(x)|^2 =0 \\
	\end{cases}
\end{equation*}
and
\begin{equation*}
	c(x)=\begin{cases}
		\frac{\overline{U(x)}(\Delta U(x))}{|\nabla U(x)|^2+|U(x)|^2}, \ \mbox{ for } \ |\nabla U(x)|^2+|U(x)|^2 >0, \\
		\\
		0 , \quad \ \ \  \mbox{for } \ |\nabla U(x)|^2+|U(x)|^2 =0 .\\
	\end{cases}
\end{equation*}
By \eqref{26nn-Nire-30-4-23-1} we have 

$$ \left\Vert b\right\Vert_{L^{\infty}\left(\Omega; \mathbb{C}^n\right)}\leq K, \ \ \mbox{and} \ \  \left\Vert c\right\Vert_{L^{\infty}\left(\Omega;\mathbb{C}\right)}\leq K.$$
In addition, if $|\nabla U(x)|^2+|U(x)|^2 =0$ then \eqref{26nn-Nire-30-4-23-1} implies 
$$\Delta U(x)=0= b(x)\cdot \nabla U(x)+c(x)U(x),$$
and, if $|\nabla U(x)|^2+|U(x)|^2 >0$ then

\begin{equation*}
	\begin{aligned}
	&b(x)\cdot \nabla U(x)+c(x)U(x)=\\&=\frac{\overline{\nabla U(x)}(\Delta U(x))}{|\nabla U(x)|^2+|U(x)|^2}\cdot \nabla U(x)+\frac{\overline{U(x)}(\Delta U(x))}{|\nabla U(x)|^2+|U(x)|^2}U(x)=\\&=\Delta U(x).
\end{aligned}
\end{equation*}
$\blacklozenge$

\bigskip

In the next Proposition we will use \textbf{the Kelvin transform} \index{Kelvin transform} which is defined as follows. Let $u$ be a
sufficiently smooth function, let

\begin{equation}\label{1-Kelvin-Nire}
v(y)=\left|y\right|^{2-n}u\left(y|y|^{-2}\right),
\end{equation}
we have

\begin{equation}\label{2-Kelvin-Nire}
\Delta v(y)=\left|y\right|^{-2-n} (\Delta_xu)\left(y|y|^{-2}\right).
\end{equation}

\bigskip

\begin{theo}[\textbf{weak unique continuation property}]\label{25-Nire}
	\index{Theorem:@{Theorem:}!- weak unique continuation property@{- weak unique continuation property}}
Let $K,R,\rho>0$, $\rho<R$. Let us assume that $U\in H^2(B_{R})$
satisfies the inequality

\begin{equation}\label{26-Nire}
\left|\Delta U\right|\leq K\left(\left|\nabla
U\right|+\left|U\right|\right), \quad \mbox{ in } B_{R}.
\end{equation}
We have that, if 
\begin{equation}\label{27-Nire}
U=0,\quad\mbox{in } B_{\rho},
\end{equation}
then
\begin{equation}\label{28-Nire}
 U=0\quad\mbox{in } B_{R}.
\end{equation}

\end{theo}
\textbf{Proof.} It is not restrictive to assume $R=1$ and, consequently, $\rho<1$. Let us apply the Kelivin transform. Set

\begin{equation}\label{3-Kelvin-Nire}
V(y)=\left|y\right|^{2-n}U\left(y|y|^{-2}\right),\quad \mbox{ in }
\mathbb{R}^n\setminus \overline{B_{1}},
\end{equation}
by \eqref{27-Nire} we have
\begin{equation}\label{29-Nire}
V=0,\quad\mbox{in } \mathbb{R}^n\setminus \overline{B_{1/\rho}}.
\end{equation}
Moreover $V\in H^2\left(B_{1/\rho}\setminus
\overline{B_{1}}\right)$ and

\begin{equation*}
\left|(\nabla_x U)\left(y|y|^{-2}\right)\right|\leq
C\left(|y|^{n-1}|V(y)|+|y|^{n}|\nabla_y V(y)|\right)
\end{equation*}
hence, by this inequality, by \eqref{2-Kelvin-Nire} and by
\eqref{26-Nire} we obtain

\begin{equation}\label{31-Nire}
\left|\Delta V\right|\leq CK\left(\left|\nabla
V\right|+\left|V\right|\right), \quad \mbox{ in }
\mathbb{R}^n\setminus \overline{B_{1}}.
\end{equation}
Now, by applying Proposition \ref{25nn-Nire} we have

\begin{equation}\label{32-Nire}
 V=0\quad\mbox{in } B_{1/\rho}\setminus \overline{B_{1}}.
\end{equation}
From which, taking into account \eqref{1-Kelvin-Nire}, we immediately obtain the
thesis. $\blacksquare$

\bigskip

We now prove the following

\begin{theo}\label{33-Nire}
Let $\Omega$ be a (bounded) connected open set of $\mathbb{R}^n$, let
$\Sigma\subset \partial \Omega$ be a local graph of a function of
class $C^2$. Let us assume that $b\in L^{\infty}(\Omega,
\mathbb{C}^n)$, $c\in L^{\infty}(\Omega, \mathbb{C})$. Moreover, let 
$U\in H^2(\Omega,\mathbb{C})$ satisfy 
\begin{equation}\label{33-CauNire}
\begin{cases}
\Delta U=b(x)\cdot\nabla U+c(x)U, \quad \mbox{ in } \Omega, \\
\\
U=0, \quad\mbox{ on } \Sigma,\\
\\
\frac{\partial U}{\partial \nu}=0, \quad \mbox{ on } \Sigma.
\end{cases}
\end{equation}
Then we have

\begin{equation}\label{34-Nire}
 U=0\quad\mbox{in }\Omega.
\end{equation}
\end{theo}

\medskip

\textbf{Remark 2.} As usual, the conditions $U=\frac{\partial U}{\partial
\nu}=0$ on $\Sigma$ should be understood as: $U\Psi\in
H_0^2(\Omega)$, for every $\Psi\in C^{\infty}(\mathbb{R}^n)$ such that
$\Psi=\frac{\partial \Psi}{\partial \nu}=0$ on $\partial \Omega
\setminus \Sigma$. $\blacklozenge$

\medskip

\textbf{Proof.} For any $M,r>0$ let us denote

$$Q_{r,M}=B'_{r}\times\left(-Mr, Mr\right).$$
Up to rigid transformation of $\mathbb{R}^n$, we may assume
there exist $r_0$ and \\ $h\in C^2\left(\overline{B'_{r_0}}\right)$
such that
$$h(0)=\left|\nabla_{x'}h(0)\right|=0$$ and
$$\mbox{graph}(h)\subset \Sigma,$$
where graph$(h)$ is the graph of $h$. While eventually reducing $r_0$ we may assume that there exists $M_0>0$ such that
$$\left\Vert h\right\Vert_{C^2\left(\overline{B'_{r_0}}\right)}\leq M_0r_0$$
and

$$\Omega\cap Q_{r_0,2M_0}=\left\{x\in B'_{r_0}\times\mathbb{R}:\quad h(x')<x_n<2M_0r_0\right\}.$$
Now, let us denote

$$Q^-_{r_02M_0}=\left\{x\in B'_{r_0}\times\mathbb{R}:\quad -2M_0r_0<x_n\leq
h(x')\right\}$$ and let $$\widetilde{\Omega}=\Omega\cup
Q^-_{r_02M_0}.$$ Moreover, let $\widetilde{U}$ be
the extension of $U$ to $0$ in $Q^-_{r_0,2M_0}$. We have $\widetilde{U}\in H^2\left(\widetilde{\Omega}\right)$ and

$$\Delta \widetilde{U}=\widetilde{b}(x)\cdot\nabla
\widetilde{U}+\widetilde{c}(x)\widetilde{U},$$ where $\widetilde{b}$
e $\widetilde{c}$ are the extensions of $b, c$ to $0$ in
$Q^-_{r_0,2M_0}$. We have trivially $\widetilde{b}\in
L^{\infty}\left(\widetilde{\Omega}, \mathbb{C}^n\right)$,
$\widetilde{c}\in L^{\infty}\left(\widetilde{\Omega},
\mathbb{C}\right)$.

From now on, we argue as we did in the proof of Theorem \ref{UCP-ell-anal}. For the convenience of the reader here we repeat the main steps of this proof.

First of all, we note that $\widetilde{\Omega}$ is connected.

Then, set

\begin{equation}\label{1n-84C-Nire}
A=\left\{x\in \widetilde{\Omega}: \exists \rho_x>0 \mbox{ such that }
\widetilde{U}=0, \mbox{ in }B_{\rho_x}(x)\right\}.
\end{equation}
By the definition of $\widetilde{U}$ we have
$\mbox{Int}\left(Q^-_{r_0,2M_0}\right)\subset A$ hence
$A\neq \emptyset$ and it turns out, trivially, that $A$ is an open set in
$\widetilde{\Omega}$. To prove that A is also closed in
$\widetilde{\Omega}$ (from which, since $\widetilde{\Omega}$ is
connected, we get $A=\widetilde{\Omega}$ and, consequently,
$\widetilde{U}\equiv 0$) we can follow exactly the same
argument followed in the proof of Theorem
\ref{UCP-ell-anal}, using, in this case, Proposition
\ref{25-Nire} instead of Theorem \ref{prop-UC-82C}.
$\blacksquare$

\bigskip

\textbf{Remark 3.} Let us observe that by the weak unique continuation property (Theorem \ref{prop-UC-82C}) we have derived the uniqueness for Cauchy problem \eqref{33-CauNire}. Conversely, if we dispose of the uniqueness for Cauchy problem \eqref{33-CauNire} we may derive the weak unique continuation property for the equation
$$\Delta U=b(x)\cdot\nabla U+c(x)U, \quad \mbox{ in } \Omega, $$
where $\Omega$ is a connected open set of $\mathbb{R}^n$. Indeed, let $\omega$ a subset of $\Omega$ and let us assume that
\begin{equation*}
U=0, \ \ \mbox{in } \omega.
\end{equation*}
Let $B_r(x_0)\Subset \omega$. 

We have that $U$ is a solution to the Cauchy problem 
\begin{equation*}\label{correct:18-4-23-3}
	\begin{cases}
		\Delta U=b(x)\cdot\nabla U+c(x)U, \quad \mbox{ in } \Omega, \\
		\\
		U=0, \quad\mbox{ on } \partial B_r(x_0),\\
		\\
		\frac{\partial U}{\partial \nu}=0, \quad \mbox{ on } \partial B_r(x_0).
	\end{cases}
\end{equation*}
Hence, by the uniqueness for the Cauchy problem, we have
$$U\equiv 0, \ \ \mbox{in } \Omega.$$ $\blacklozenge$
\section{Necessary conditions}\label{Condiz-necessarie}

In the present Section we return to estimate \eqref{1.15-sece} or, more generally,
to estimate \eqref{4-Nire}. We have already observed (Remark 1 of Section \ref{Teorema-Nire}) that if \eqref{3-Nire} holds, then \eqref{4-Nire} holds.
As we will see, the converse is also true. More precisely, we have
the following

\begin{theo}\label{cond-nec-Nire}
Let $M(D)$ and $P(D)$ be two linear differential  operators of order
$r$ and $m$ respectively. The following conditions are equivalent:

\noindent(i) There exists $C>0$ such that

\begin{equation}\label{43-Nire}
\int_{Q_1}e^{2\tau N\cdot x}\left|M(D)u\right|^2dx\leq
C\int_{Q_1}e^{2\tau N\cdot x}\left|P(D)u\right|^2dx,
\end{equation}
for every $u \in C^{\infty}_{0}(Q_1,\mathbb{C})$ and for every $\tau\in
\mathbb{R}$.

\noindent(ii) The following  is true

\begin{equation}\label{44-Nire}
\sup_{(\xi,\tau)\in \mathbb{R}^{n+1}}\left\{\frac{\left|M(\xi+i\tau
N)\right|^2}{\sum_{|\alpha|\leq m}\left|P^{(\alpha)}(\xi+i\tau
N)\right|^2} \right\}<+\infty.
\end{equation}
In particular, if \eqref{43-Nire} holds true, then $r\leq m$.
\end{theo}

\bigskip

In order to prove Theorem \ref{cond-nec-Nire} we need a Lemma 
(Lemma \ref{Simb-Nire}) to which we premise the following \index{extension of the Leibniz formula}

\begin{prop}[\textbf{extension of the Leibniz formula}]\label{L-Niren} Let $P(D)$ be a linear differential  operator of order $m$, we have
\begin{equation}\label{Leibniz-Nire}
P(D)[fu]=\sum_{|\alpha|\leq
m}\frac{D^{\alpha}f}{\alpha!}P^{(\alpha)}(D)u, \quad \forall f,u\in
C^{m}(\mathbb{R}^n,\mathbb{C}).
\end{equation}
\end{prop}

\textbf{Proof.} By the Leibniz formula, we get

\begin{equation}\label{1000}
P(D)[fu]=\sum_{|\alpha|\leq m}(D^{\alpha}f) R_{\alpha}(D)u, \quad
\forall f,u\in C^{m}(\mathbb{R}^n,\mathbb{C}),
\end{equation}
where $R_{\alpha}(D)$ are linear differential operators of order (less or equal) to \\ $m-|\alpha|$. Let $\xi,\eta\in \mathbb{R}^n$ be arbitrary
and let $f(x)=e^{ix\cdot \xi}$, $u(x)=e^{ix\cdot \eta}$. By
\eqref{1000} we have

\begin{equation*}\label{1-Leibniz-Nire}
e^{ix\cdot (\xi+\eta)}P(\xi+\eta)=P(D)[e^{ix\cdot \xi}e^{ix\cdot
\eta}]=e^{ix\cdot (\xi+\eta)}\sum_{|\alpha|\leq m}\xi^{\alpha}
R_{\alpha}(\eta).
\end{equation*}
Hence
$$P(\xi+\eta)=\sum_{|\alpha|\leq m}\xi^{\alpha}
R_{\alpha}(\eta), \quad \forall \xi,\eta \in \mathbb{R}^n.$$ On the other hand, by the Taylor formula we have 
$$P(\xi+\eta)=\sum_{|\alpha|\leq
m}\frac{\xi^{\alpha}}{\alpha!}P^{(\alpha)}(\eta),$$ from which we obtain
$$R_{\alpha}(\eta)=\frac{1}{\alpha!}P^{(\alpha)}(\eta)$$ that gives
\eqref{Leibniz-Nire}. $\blacksquare$

\bigskip

\begin{lem}\label{Simb-Nire}
For any $m \in \mathbb{N}_0$ and for any $\phi\in
C^{\infty}_{0}(Q_1,\mathbb{C})$, which does not vanish identically, there exists a
constant $C_{\phi}\geq 1$ such that for every linear differential operator
$R(D)$ of order $m$ we have

\begin{equation}\label{45-Nire}
C_{\phi}^{-1}\widetilde{R}(\xi)\leq \left(\int_{Q_1}\left\vert
R(D)\left[\phi(x)e^{ix\cdot
\xi}\right]\right\vert^2dx\right)^{1/2}\leq
C_{\phi}\widetilde{R}(\xi),
\end{equation}
for every $\xi \in
\mathbb{R}^n$, where

\begin{equation}\label{47-Nire}
\widetilde{R}(\xi)=\left(\sum_{|\alpha|\leq
m}\left|R^{(\alpha)}(\xi)\right|^2\right)^{1/2}.
\end{equation}

\end{lem}

\textbf{Proof.} Let $\phi\in C^{\infty}_{0}(Q_1,\mathbb{C})
$ be a function  not identically zero

Formula \eqref{Leibniz-Nire} gives

\begin{equation}\label{48-Nire}
R(D)\left[\phi(x)e^{ix\cdot \xi}\right]=e^{ix\cdot
\xi}\sum_{|\alpha|\leq
m}\frac{1}{\alpha!}R^{(\alpha)}(\xi)D^{\alpha}\phi(x).
\end{equation}
Set

\begin{equation}\label{49-Nire}
I_{\alpha\beta}\left(\phi\right)=\int_{\mathbb{R}^n}\frac{1}{\alpha!}D^{\alpha}\phi(x)\frac{1}{\beta!}\overline{D^{\beta}\phi(x)}dx,
\end{equation}
by \eqref{48-Nire}, we have

\begin{equation}\label{50-Nire}
\begin{aligned}
&\int_{\mathbb{R}^n}\left\vert R(D)\left[\phi(x)e^{ix\cdot
\xi}\right]\right\vert^2dx =
\int_{\mathbb{R}^n}\left\vert\sum_{|\alpha|\leq
m}\frac{1}{\alpha!}R^{(\alpha)}(\xi)D^{\alpha}\phi(x)\right\vert^2dx=\\&
= \int_{\mathbb{R}^n}\sum_{|\alpha|,|\beta|\leq
m}R^{(\alpha)}(\xi)\overline{R^{(\beta)}(\xi)}
\frac{1}{\alpha!}D^{\alpha}\phi(x)\frac{1}{\beta!}\overline{D^{\beta}\phi(x)}dx=\\&
=\sum_{|\alpha|,|\beta|\leq
m}R^{(\alpha)}(\xi)\overline{R^{(\beta)}(\xi)}I_{\alpha\beta}\left(\phi\right).
\end{aligned}
\end{equation}
Let us consider the quadratic form

\begin{equation}\label{51-Nire}
H(z)=\sum_{|\alpha|,|\beta|\leq
m}I_{\alpha\beta}\left(\phi\right)z_{\alpha}\overline{z_{\beta}},
\end{equation}
where $z\in \mathbb{C}^{N(m,n)}$ ($N(m,n)$ the number of
multi--indexes $\alpha$ such that $|\alpha|\leq m$). Let us notice that 

\begin{equation}\label{52-Nire}H(z)=\int_{\mathbb{R}^n}\left\vert \sum_{|\alpha|\leq
m}\frac{1}{\alpha!}z_{\alpha}D^{\alpha}\phi(x)\right\vert^2dx.\end{equation}
From which, recalling \eqref{50-Nire} and choosing
\begin{equation}\label{501-Nire}
z^{(0)}_{\alpha}=R^{(\alpha)}(\xi), \quad \mbox{for }\quad
|\alpha|\leq m,
\end{equation}
we have
\begin{equation}\label{5001-Nire}H\left(z^{(0)}\right)=\int_{\mathbb{R}^n}\left\vert
R(D)\left[\phi(x)e^{ix\cdot \xi}\right]\right\vert^2dx.\end{equation}
Of course,  \eqref{52-Nire} gives $H(z)\geq 0$ for every $z\in
\mathbb{C}^{N(m,n)}$. Actually, $H(z)$ is a \textbf{positive--definite form} as we are going to prove. Arguing by contradiction, let us suppose that there exists
$\widetilde{z}\neq 0$ such that $ H\left(\widetilde{z}\right)=0$
then,  applying inequality \eqref{Nire-1} to the operator

$$\widetilde{M}(D)=\sum_{|\alpha|\leq
m}\frac{1}{\alpha!}\widetilde{z}_{\alpha}D^{\alpha},$$ we would have that there exists
a constant $\widetilde{C}$ such that
$$0=\widetilde{C}H\left(\widetilde{z}\right)=\widetilde{C}\int_{Q_1}\left|\widetilde{M}(D)\phi\right|^2dx\geq \int_{Q_1}\left|\phi\right|^2 dx$$
and this would imply $\phi\equiv 0$, which  contradicts the assumption that $\phi$ does not vanish identically. Since $H$ is \textbf{positive--definite}, there exists a constant
$C_{\phi}\geq 1$ (depending by $\phi$) such that

\begin{equation}\label{53-Nire}
C_{\phi}^{-1}\sum_{|\alpha|\leq m}\left|z_{\alpha}\right|^2 \leq
H(z)\leq  C_{\phi}\sum_{|\alpha|\leq m}\left|z_{\alpha}\right|^2,
\quad \forall z \in \mathbb{C}^{N(m,n)}.
\end{equation}
By the latter, recalling \eqref{501-Nire} and \eqref{5001-Nire}
we get \eqref{45-Nire}. $\blacksquare$

\bigskip

\textbf{Proof of Theorem \ref{cond-nec-Nire}.}

We have already seen in Remark 1 of Section \ref{Teorema-Nire} that
\eqref{44-Nire} implies \eqref{43-Nire}. Now we prove the converse. Let $C_{\star}$ be a positive constant, let us fix $\tau
\in\mathbb{R}$ and let us assume that 

\begin{equation}\label{60-Nire}
\int_{Q_1}e^{2\tau N\cdot x}\left|M(D)u\right|^2dx\leq
C_{\star}\int_{Q_1}e^{2\tau N\cdot x}\left|P(D)u\right|^2dx,
\end{equation}
for every $u \in C^{\infty}_{0}(Q_1,\mathbb{C})$. Then, setting
$v=e^{2\tau N\cdot x}u$ and arguing similarly to the proof of
Theorem \ref{1.14-corollario-sece}, we have 
\begin{equation*}
\begin{aligned}
\int_{Q_1}\left|M(D+i\tau N)v\right|^2dx& =\int_{Q_1}e^{2\tau N\cdot
x}\left|M(D)u\right|^2dx\leq\\& \leq C_{\star}\int_{Q_1}e^{2\tau
N\cdot x}\left|P(D)u\right|^2dx=
 \\& =C_{\star}\int_{Q_1}\left|P(D+i\tau N)v\right|^2dx.
\end{aligned}
\end{equation*}

Therefore

\begin{equation}\label{61-Nire}
\int_{Q_1}\left|M(D+i\tau N)v\right|^2dx\leq
C_{\star}\int_{Q_1}\left|P(D+i\tau N)v\right|^2dx,
\end{equation}
for every $v \in C^{\infty}_{0}(Q_1,\mathbb{C})$.

Now let $\phi \in C^{\infty}_{0}(Q_1,\mathbb{C})$ be not identically zero. For instance, let
$$\phi(x)=\begin{cases}
\exp\left(-\frac{1}{1-|x|^2}\right), \quad \mbox{ for } |x|<1, \\
\\
0, \quad \mbox{ for } |x|\geq 1.
\end{cases}
$$
Let $v(x)=\phi(x)e^{ix\cdot \xi}$. By \eqref{45-Nire} and
\eqref{61-Nire} we have, for every $\xi\in\mathbb{R}^n$,

\begin{equation}\label{62-Nire}
 \begin{aligned}
&C^{-1}\sum_{|\alpha|\leq r}\left|M^{(\alpha)}(\xi+i\tau
N)\right|^2\leq \\&\leq\int_{Q_1}\left\vert M(D+i\tau
N)\left[\phi(x)e^{ix\cdot \xi}\right]\right\vert^2dx\leq
\\&\leq C_{\star}\int_{Q_1}\left|P(D+i\tau N)\left[\phi(x)e^{ix\cdot
\xi}\right]\right|^2dx\leq \\& \leq C C_{\star}\sum_{|\alpha|\leq
m}\left|P^{(\alpha)}(\xi+i\tau N)\right|^2.
\end{aligned}
\end{equation}
where $C\geq 1$ \emph{does not} depend on $\tau$. By \eqref{62-Nire}
we obtain

\begin{equation}\label{63-Nire}
\left|M(\xi+i\tau N)\right|^2 \leq C^2C_{\star}\sum_{|\alpha|\leq
m}\left|P^{(\alpha)}(\xi+i\tau N)\right|^2, \quad \forall
(\xi,\tau)\in\mathbb{R}^{n+1}.
\end{equation}
Therefore
\begin{equation*}
\sup_{(\xi,\tau)\in \mathbb{R}^{n+1}}\left\{\frac{\left|M(\xi+i\tau
N)\right|^2}{\sum_{|\alpha|\leq m}\left|P^{(\alpha)}(\xi+i\tau
N)\right|^2} \right\}\leq C^2C_{\star}< +\infty,
\end{equation*}
that concludes the proof.
$\blacksquare$

\bigskip

\textbf{Remark 1.} By reviewing the proof of Theorem
\ref{cond-nec-Nire} it is easily seen that if $\tau_0$ is a
fixed real number the following conditions are equivalent:

\medskip

\noindent(i') There exists $C>0$ such that

\begin{equation}\label{64-Nire}
\int_{Q_1}e^{2\tau_0 N\cdot x}\left|M(D)u\right|^2dx\leq
C\int_{Q_1}e^{2\tau_0 N\cdot x}\left|P(D)u\right|^2dx,
\end{equation}
for every $u \in C^{\infty}_{0}(Q_1,\mathbb{C})$.

\noindent(ii') The following holds true

\begin{equation}\label{65-Nire}
\sup_{\xi\in \mathbb{R}^{n}}\left\{\frac{\left|M(\xi+i\tau_0
N)\right|^2}{\sum_{|\alpha|\leq m}\left|P^{(\alpha)}(\xi+i\tau_0
N)\right|^2} \right\}<+\infty.
\end{equation}

\bigskip

In particular the following conditions are
equivalent:

\medskip

\noindent(i'') There exists $C>0$ such that for every $u \in
C^{\infty}_{0}(Q_1,\mathbb{C})$ we have
\begin{equation*}
\int_{Q_1}\left|M(D)u\right|^2dx\leq
C\int_{Q_1}\left|P(D)u\right|^2dx,
\end{equation*}

\medskip

\noindent(ii'') \begin{equation*} \sup_{\xi\in
\mathbb{R}^{n}}\frac{\left|M(\xi)\right|^2}{\sum_{|\alpha|\leq
m}\left|P^{(\alpha)}(\xi)\right|^2} <+\infty.
\end{equation*}
$\blacklozenge$

\bigskip

\textbf{Remark 2.} In applying Theorem
\ref{1.14-corollario-sece} for proving the uniqueness of
Cauchy problem \eqref{1-CauNire} it was sufficient
to use estimate \eqref{1.15-sece} for $\tau$ sufficiently
large. More precisely, the estimate we actually
used is

\begin{equation}\label{tau-zero}
\int_{Q_1}e^{2\tau N\cdot x}\left|M(D)u\right|^2dx\leq
C\int_{Q_1}e^{2\tau N\cdot x}\left|P(D)u\right|^2dx,
\end{equation}
for every $u \in C^{\infty}_{0}(Q_1,\mathbb{C})$ and for every $\tau\geq
\tau_0$, where $\tau_0$ is a positive number and, as well as $C$,
does not depends on $u$. If $M(\xi)$ and $P(\xi)$ are \textbf{homogeneous polynomials}, it is simple to check that estimate  \eqref{tau-zero}
is equivalent to the estimate

\begin{equation}\label{tau-zero-1}
\int_{Q_1}e^{-2\tau N\cdot x}\left|M(D)u\right|^2dx\leq
C\int_{Q_1}e^{-2\tau N\cdot x}\left|P(D)u\right|^2dx,
\end{equation}
for every $u \in C^{\infty}_{0}(Q_1,\mathbb{C})$ and for every $\tau\geq
\tau_0$. It suffices to consider the simple change of variables
$x\rightarrow -x$ (the reader takes care of the details).
Taking into account what we said in \textbf{Remark 1} of this Section, estimate
\eqref{tau-zero} (hence, estimate \eqref{tau-zero-1}) is
equivalent to

\begin{equation}\label{404-Nire}
\sup\left\{\frac{\left|M(\xi+i\tau N)\right|^2}{\sum_{|\alpha|\leq
m}\left|P^{(\alpha)}(\xi+i\tau N)\right|^2}: (\xi,\tau)\in
\mathbb{R}^{n+1}, |\tau|\geq \tau_0 \right\}<+\infty.
\end{equation}
$\blacklozenge$

\bigskip

\textbf{Remark 3.} If $M_j(\xi)$ are polynomials of degree $r_j$,
for $j=1,\cdots, J$ and $S_k(\xi)$ are polynomials of degree $s_k$, for
$k=1,\cdots, K$, the necessary and sufficient conditions for the
validity of estimate

\begin{equation}\label{tau-zero-2}
\sum_{j=1}^J\int_{Q_1}e^{2\tau N\cdot x}\left|M_j(D)u\right|^2dx\leq
C\sum_{k=1}^K\int_{Q_1}e^{2\tau N\cdot x}\left|S_k(D)u\right|^2dx,
\end{equation}
for every $u \in C^{\infty}_{0}(Q_1,\mathbb{C})$ and for every $\tau\geq
\tau_0$ (or for a fixed $\tau$) can be obtained easily arguing as in the proof of Theorem \ref{cond-nec-Nire}. For instance, the estimate

\begin{equation}\label{tau-zero-3}
\sum_{j=1}^J\int_{Q_1}\left|M_j(D)u\right|^2dx\leq
C\sum_{k=1}^K\int_{Q_1}\left|S_k(D)u\right|^2dx,
\end{equation}
for every $u \in
C^{\infty}_{0}(Q_1,\mathbb{C})$, is equivalent to
\begin{equation}\label{tau-zero-4}
\sup_{\xi\in
\mathbb{R}^{n}}\left\{\frac{\sum_{j=1}^J\left|M_j(\xi)\right|^2}{\sum_{k=1}^K\sum_{|\alpha|\leq
s_k}\left|S_k^{(\alpha)}(\xi)\right|^2} \right\}<+\infty.
\end{equation}
From which, for $r\leq s$, where $r,s$ are nonnegative integer numbers, we easily get the estimate

\begin{equation}\label{tau-zero-5-28}
\int_{Q_1}\left|D^ru\right|^2dx\leq C\int_{Q_1}\left|D^s
u\right|^2dx,\quad \forall u \in C^{\infty}_{0}(Q_1,\mathbb{C}),
\end{equation}
where, for a nonnegative integer $p$, we set
$$\left|D^{p}u\right|^2=\sum_{|\alpha|=p}\left|D^{\alpha}u\right|^2.$$
Similarly it can be proved that
\begin{equation}\label{tau-zero-5}
\int_{Q_1}e^{2\tau N\cdot x}\left|D^ru\right|^2dx\leq
C\int_{Q_1}e^{2\tau N\cdot x}\left|D^s u\right|^2dx,
\end{equation}
for every $u \in C^{\infty}_{0}(Q_1,\mathbb{C})$ and for every $\tau\in
\mathbb{R}$.

\noindent More generally we have
\begin{equation}\label{tau-zero-5-15-8-21}
\tau^{2(s-r)}\int_{Q_1}e^{2\tau N\cdot x}\left|D^ru\right|^2dx\leq
C\int_{Q_1}e^{2\tau N\cdot x}\left|D^s u\right|^2dx,
\end{equation}
for every $u \in C^{\infty}_{0}(Q_1,\mathbb{C})$ and for every $\tau\in
\mathbb{R}$.

\medskip

More attention is required to study the following two estimates. Let $P_m(\xi)$ be a homogeneous polynomial of degree $m\geq 1$ and let us consider the estimates

\begin{equation}\label{tau-zero-6}
\int_{Q_1}\left|D^mu\right|^2dx\leq
C\int_{Q_1}\left|P_m(D)u\right|^2dx,\quad \forall u \in
C^{\infty}_{0}(Q_1,\mathbb{C})
\end{equation}
and

\begin{equation}\label{tau-zero-7}
\int_{Q_1}e^{2\tau N\cdot x}\left|D^mu\right|^2dx\leq
C\int_{Q_1}e^{2\tau N\cdot x}\left|P_m(D)u\right|^2dx,
\end{equation}
for every $u \in C^{\infty}_{0}(Q_1,\mathbb{C})$ and for every
$\tau\geq \tau_0$, where $\tau_0$ is a nonnegative integer number. 
\medskip

Let us begin by \eqref{tau-zero-6}. We distinguish two cases:

\medskip

\noindent \textbf{(a)} $P_m(D)$ is \textbf{elliptic}, that is

\begin{equation}\label{tau-zero-8}
\xi\in \mathbb{R}^n\mbox{, }P_m(\xi)=0\Rightarrow
\xi=0;\end{equation}

\smallskip

\noindent \textbf{(b)} $P_m(D)$ \textbf{is not elliptic},
that is there exists $\xi_0\in \mathbb{R}^n\setminus\{0\}$ such that

\begin{equation}\label{tau-zero-9}
P_m(\xi_0)=0.
\end{equation}

\medskip

Let us check that \textbf{estimate \eqref{tau-zero-6} holds if and only if}
$P_m(D)$ \textbf{is elliptic}.

\noindent Let us denote by

\begin{equation}\label{tau-zero-10}
q(\xi)=\frac{\left|\xi\right|^{2m}}{\sum_{|\alpha|\leq
m}\left|P_m^{(\alpha)}(\xi)\right|^2}.
\end{equation}
Let us assume that (a) holds true. Since $P_m(\xi)$ is a homogeneous polynomial
of degree $m$,  there exists $\lambda\geq 1$
such that

\begin{equation}\label{tau-zero-11}
\lambda^{-1}\left|\xi\right|^m\leq \left|P_m(\xi)\right|\leq
\lambda\left|\xi\right|^m,\quad \forall \xi\in \mathbb{R}^n.
\end{equation}
Hence

\begin{equation*}
q(\xi)\leq \lambda^2,\quad \forall \xi\in \mathbb{R}^n
\end{equation*}
and, by what noted in \textbf{Remark 1} of this Section, we have that estimate
\eqref{tau-zero-6} holds true.

In case (b), let $\xi_0\in \mathbb{R}^n\setminus\{0\}$ satisfy
$P_m(\xi_0)=0$. Let $\mu \in \mathbb{R}$. Then the numerator
of $q(\mu\xi_0)$ is equal to $\mu^{2m}|\xi_0|^{2m}$ and the denominator has degree w.r.t. $\mu$ less or equal to $2m-2$, as 
$P_m(\mu\xi_0)=0$, for every $\mu \in \mathbb{R}$.
Therefore

\begin{equation*}
\lim_{\mu\rightarrow +\infty}q(\mu\xi_0)=+\infty
\end{equation*}
from which we have that in case (b), estimate \eqref{tau-zero-6} does not hold.

\medskip

Let now consider estimate \eqref{tau-zero-7}. We prove that \textbf{it does not hold in any case}.

\noindent For any $\xi\in \mathbb{R}^n$ and $\tau\geq \tau_0$, set

\begin{equation}\label{tau-zero-12}
q(\xi,\tau)=\frac{\left(|\xi|^2+\tau^2\right)^{m}}{\sum_{|\alpha|\leq
m}\left|P_m^{(\alpha)}(\xi+i\tau N)\right|^2}.
\end{equation}

Let us begin by case (a). Let $\xi\in \mathbb{R}^n$ be such that $\xi$ and
$N$ are linearly independent. Let us consider the equation

\begin{equation}\label{tau-zero-13}
P_m(\xi+z N)=0, \quad z\in \mathbb{C}.
\end{equation}
Let $a+ib$ be a solution of \eqref{tau-zero-13}. Then $b\neq 0$,
otherwise we would have $P_m(\xi+a N)=0$, but since $P_m(D)$ is
elliptic, we would have $\xi+a N=0$ which contradicts the assumption of linear independence between
$\xi$ and $N$. Hence, either $b>0$
or $b<0$. Setting $\eta=\xi+a N$, we have $\eta\neq 0$ and

\begin{equation}\label{tau-zero-14}
P_m(\eta+i b N)=0.
\end{equation}
Now, if $b>0$, let $\mu\geq \frac{\tau_0}{b}$ and we get 

\begin{equation}\label{tau-zero-15} P_m(\mu\eta+i\mu b
N)=\mu^mP_m(\eta+i b N)=0.
\end{equation}
Hence the numerator of $q(\mu\eta,\mu\tau)$ is equal to
$\mu^{2m}\left(|\eta|^2+b^2\right)^m$, whereas, by
\eqref{tau-zero-15}, the denominator of $q(\mu\eta,\mu\tau)$ has degree
w.r.t. $\mu$ less or equal to $2m-2$. Therefore

\begin{equation}\label{tau-zero-16}
\lim_{\mu\rightarrow +\infty}q(\mu\eta,\mu b)=+\infty.
\end{equation}
If $b<0$, it suffices to notice that by \eqref{tau-zero-14} it follows
$P_m(-\eta+i (-b) N)=0$ and similarly to \eqref{tau-zero-16} we have
\begin{equation*}
\lim_{\mu\rightarrow +\infty}q(\mu(-\eta),\mu (-b))=+\infty.
\end{equation*}
Hence
\begin{equation} \sup\left\{q(\xi,\tau): (\xi,\tau)\in
\mathbb{R}^{n+1}, |\tau|\geq \tau_0 \right\}=+\infty.
\end{equation}
To conclude, in the elliptic case estimate \eqref{tau-zero-7} does not
hold.

\medskip

Let us consider case (b). Let $\xi_0\in \mathbb{R}^n\setminus\{0\}$ satisfy

\begin{equation*}
P_m(\xi_0)=0.
\end{equation*}
Let $\tau\geq \tau_0$ be fixed. We have

\begin{equation*}
	\begin{aligned}
	\lim_{\mu\rightarrow +\infty}q(\mu\xi_0,\mu \tau)=
	\lim_{\mu\rightarrow +\infty}
	\frac{\left(|\xi_0|^2+\left(\tau\mu^{-1}\right)^2\right)^{m}}{h(\mu)}=+\infty,
	\end{aligned}
\end{equation*}
where
$$h(\mu)=\sum_{|\alpha|\leq
	m-1}\mu^{-2(m-|\alpha|)}\left|P_m^{(\alpha)}(\xi_0+i\left(\tau\mu^{-1}\right)
N)\right|^2+\left|P_m(\xi_0+i\left(\tau\mu^{-1}\right)
N)\right|^2.$$
From which we have that if $P_m(D)$ is not elliptic, estimate \eqref{tau-zero-7} does not hold. $\blacklozenge$

\bigskip

\section{Examples and further considerations. } \label{altre consid}
Remark 1 of Section \ref{Teorema-Nire} implies that,
if $P_m(\xi)$, is a homogeneous polynomial of degree  $m$, then we have
\begin{equation}\label{Nire-64}
\int_{Q_1}e^{2\tau N\cdot x}\left|u\right|^2dx\leq
C\int_{Q_1}e^{2\tau N\cdot x}\left|P_m(D)u\right|^2dx,
\end{equation}
for every $u \in C^{\infty}_{0}(Q_1,\mathbb{C})$ and for every $\tau\in
\mathbb{R}$. 

Let us assume $m\geq 1$. Let $N=-e_n$, let $l$ be a positive number and
$$h:B'_1\rightarrow \mathbb{R}$$ be a \textbf{strictly convex function} which satisfies 

\begin{equation}\label{100-CauNire}
h(0)=0,\quad\mbox{ and }\quad l\leq \inf_{\partial B'_1} h.
\end{equation}
Let

$$\Lambda=\left\{(x',x_n)\in B_1'\times\mathbb{R}: h(x')<x_n<l
\right\}$$ and

$$\Gamma=\left\{(x',h(x')):x'\in B_1'(0)
\right\}.$$ Moreover, let  $a_0\in L^{\infty}(\Lambda)$. By proceeding in
similar manner as we did to prove the uniqueness for
Cauchy problem \eqref{1-CauNire}, w prove the uniqueness for the problem 

\begin{equation}\label{65-CauNire}
\begin{cases}
 P_m(D)U+a_0(x)U=0, \quad \mbox{ in } \Lambda, \\
\\
U\Psi\in H_0^m(\Lambda) , \quad\forall\Psi\in C^{\infty}(\mathbb{R}^n),\quad\mbox{supp}\Psi\subset \mathbb{R}^{n-1}\times (-\infty,1). \\
\end{cases}
\end{equation}
We invite the reader to develop the details (remember to use
the homothetic transformation $x\rightarrow rx$), however, on this issue we
refer to \cite[Theorem 1]{Ni}.

Of course, it is meaningful and interesting to ask what happens regarding
the uniqueness for the Cauchy problem if one perturbs the operator
with operators of order $r$ with $1\leq r\leq m-1$. Keep in mind,
however, that here we are basically considering the case
of uniqueness for the Cauchy problem whose initial surface is a
strictly convex function.

\medskip

Let us consider  the case $r=m-1$. If the estimate holds true 

\begin{equation}\label{101-CauNire}
\int_{Q_1}e^{2\tau N\cdot x}\left|D^{m-1}u\right|^2dx\leq
C\int_{Q_1}e^{2\tau N\cdot x}\left|P_m(D)u\right|^2dx,
\end{equation}
for every $u \in C^{\infty}_{0}(Q_1,\mathbb{C})$ and for every
$\tau\geq \tau_0$ ($\tau_0\geq 0$), then we have the uniqueness of solutions to the Cauchy problem

\begin{equation}\label{102-CauNire}
\begin{cases}
 P_m(D)U+\sum_{|\alpha|\leq
m-1}b_{\alpha}(x)U=0, \quad \mbox{ in } \Lambda, \\
\\
U\Psi\in H_0^m(\Lambda) , \quad\forall\Psi\in C^{\infty}(\mathbb{R}^n),\quad\mbox{supp}\Psi\subset \mathbb{R}^{n-1}\times (-\infty,1),\\
\end{cases}
\end{equation}
where $b_{\alpha}\in L^{\infty}(\Lambda)$, for $|\alpha|\leq m-1$.
As a matter of fact if \eqref{101-CauNire} holds, then by \eqref{tau-zero-5}
we have

\begin{equation}\label{103-CauNire}
\sum_{j=0}^{m-1}\int_{Q_1}e^{2\tau N\cdot
x}\left|D^{j}u\right|^2dx\leq C\int_{Q_1}e^{2\tau N\cdot
x}\left|P_m(D)u\right|^2dx,
\end{equation}
for every $u \in C^{\infty}_{0}(Q_1,\mathbb{C})$ and for every
$\tau\geq \tau_0$. In particular, let us observe that \eqref{103-CauNire} and
\eqref{101-CauNire} are equivalent and, in addition, 
\eqref{103-CauNire} allows us to treat (when $N=-e_n$) problem
\eqref{102-CauNire}  in a manner similar to 
\eqref{1-CauNire}.

\medskip

In \textbf{Remark 3} of Section \ref{Condiz-necessarie}, we have seen, that
necessary and sufficient condition to be hold \eqref{101-CauNire}
is there exists $C>0$ such that, for every $\xi\in \mathbb{R}^n$ and for
every $\tau\geq \tau_0$, we have

\begin{equation}\label{104-CauNire}
q_{m-1}(\xi,\tau)=\frac{\left(|\xi|^2+\tau^2\right)^{m-1}}{\sum_{|\alpha|\leq
m}\left|P_m^{(\alpha)}(\xi+i\tau N)\right|^2}\leq C.
\end{equation}

Now, in the next Proposition we give a simpler formulation of condition
\eqref{104-CauNire}

\begin{prop}\label{105-CauNire}
The following conditions are equivalent:

\smallskip

\noindent(a) Estimate \eqref{103-CauNire} holds true

\noindent(b)  If $(\xi,\tau) \in \mathbb{R}^{n+1}\setminus
\{(0,0)\}$, then
$$P_m(\xi+i\tau N)=0\Rightarrow\sum_{j=1}^n \left|P^{(j)}_m(\xi+i\tau
N)\right|^2> 0,$$ where $P^{(j)}_m(\xi+i\tau N)=P^{(e_j)}_m(\xi+i\tau
N)$.
\end{prop}

\medskip

In order to prove Proposition \ref{105-CauNire} we will use  two 
lemmas.

\begin{lem}\label{106-CauNire}
Let $d$ be a positive integer  number. Let $K$ be a compact set of
$\mathbb{R}^d$ and let $f:K\rightarrow \mathbb{R}$ and $g:K\rightarrow
\mathbb{R}$ be two continuous functions. The following conditions are equivalent:

\smallskip

\noindent(i) $X\in K$, $f(X)=0$ $\Rightarrow$ $g(X)> 0$;

\noindent(ii) There exists $C>0$ such that $C\left(f(X)\right)^2+g(X)>0$
for every $X\in K$.

\end{lem}

\textbf{Proof.} If $K=\emptyset$, the equivalence between (i) and (ii) is trivial. Let us suppose, accordingly, that $K\neq\emptyset$ and that (i) apply. Set
$$K_0=\left\{X\in K: f(X)=0 \right\}.$$ By (i), by the continuity of $g$,
and since $K$ is compact, there exists an open set $V_0$,
of $\mathbb{R}^d,$ which satisfies $K_0\subset V_0$ and $$g(X)>0, \ \ \forall X\in V_0\cap
K.$$ 

If $K\setminus V_0=\emptyset$, then 
(ii) is trivially satisfied. If $K\setminus V_0\neq\emptyset$,
we set
$$M_1=\min_{K\setminus
V_0}f^2>0, \quad M_2=\min_{K\setminus V_0}g$$ and let $C$ be a positive number
such that $CM_1+M_2>0$. We have

$$ C\left(f(X)\right)^2+g(X)\geq g(X)>0, \quad \forall X\in V_0\cap K$$
and
$$ C\left(f(X)\right)^2+g(X)\geq CM_1+M_2>0, \quad \forall X\in K\setminus
V_0,$$ from which (ii) follows.

Let us suppose that (ii) holds, we have trivially that if $f(X)=0$ then

$$g(X)=C\left(f(X)\right)^2+g(X)>0.$$
$\blacksquare$

\bigskip

\begin{lem}\label{1006-CauNire}
Let $\tau_0\in \mathbb{R}$ and let $M_j(\xi)$ be polynomials of degree
$r_j$, for $j=1,\cdots, J$. The following conditions are equivalent:

\smallskip

\noindent (a) There exists $C>0$ such that
\begin{equation}\label{1003-CauNire-a}
\sum_{j=1}^J\int_{Q_1}e^{2\tau_0 N\cdot
x}\left|M_j(D)u\right|^2dx\leq C\int_{Q_1}e^{2\tau_0 N\cdot
x}\left|P_m(D)u\right|^2dx,
\end{equation}
for every $u \in C^{\infty}_{0}(Q_1,\mathbb{C})$;

\noindent (b) there exists $C>0$ such that
\begin{equation}\label{1003-CauNire-b}
\sum_{j=1}^J\int_{Q_1}\left|M_j(D)u\right|^2dx\leq
C\int_{Q_1}\left|P_m(D)u\right|^2dx,
\end{equation}
for every $u \in C^{\infty}_{0}(Q_1,\mathbb{C})$;

\noindent (c) \begin{equation}\label{prop-Nec-005}
\sup_{\xi\in\mathbb{R}^n}\frac{\sum_{j=1}^J\left|M_j(\xi)\right|^2}{\sum_{|\alpha|\leq
m}\left|P_m^{(\alpha)}(\xi)\right|^2}<+\infty.
\end{equation}
\end{lem}
\textbf{Proof.} Let us assume that \eqref{1003-CauNire-a} holds.
For every $u \in C^{\infty}_{0}(Q_1,\mathbb{C})$ we have

\begin{equation}\label{1003-CauNire-1}
 \begin{aligned}
\sum_{j=1}^J\int_{Q_1}\left|M_j(D)u\right|^2dx&=\sum_{j=1}^J\int_{Q_1}e^{-2\tau_0
N\cdot x}e^{2\tau_0 N\cdot x}\left|M_j(D)u\right|^2dx \leq
\\&\leq e^{2|\tau_0|\sqrt{n}}\sum_{j=1}^J\int_{Q_1}e^{2\tau_0 N\cdot x}\left|M_j(D)u\right|^2dx\leq \\&
\leq Ce^{2|\tau_0|\sqrt{n}}\int_{Q_1}e^{2\tau_0 N\cdot
x}\left|P_m(D)u\right|^2dx\leq\\& \leq
Ce^{4|\tau_0|\sqrt{n}}\int_{Q_1}\left|P_m(D)u\right|^2dx.
\end{aligned}
\end{equation}
Hence, (b) follows. Similarly, we can prove that (b) implies (a).

The equivalence between (b) and (c) was proved in
\textbf{Remark 3} of Section \ref{Condiz-necessarie}, see \eqref{tau-zero-4}. $\blacksquare$

\bigskip

\textbf{Remark 1.} Taking into account \textbf{Remark 2} of Section \ref{Condiz-necessarie},
we have that if \eqref{103-CauNire} holds then  
\begin{equation}\label{10031-CauNire}
\sum_{j=0}^{m-1}\int_{Q_1}e^{2\tau N\cdot
x}\left|D^{j}u\right|^2dx\leq C\int_{Q_1}e^{2\tau N\cdot
x}\left|P_m(D)u\right|^2dx,
\end{equation}
for every $u \in C^{\infty}_{0}(Q_1,\mathbb{C})$ and for every $\tau\in
\mathbb{R}$. Hence,  Theorem  \ref{cond-nec-Nire} implies that 
estimate \eqref{103-CauNire} holds if and only if there exists $C>0$ such that for every
 $(\xi,\tau)\in \mathbb{R}^{n+1}$ we have

\begin{equation}\label{condizione-proposiz}
q_{m-1}(\xi,\tau)=\frac{\left(|\xi|^2+\tau^2\right)^{m-1}}{\sum_{|\alpha|\leq
m}\left|P_m^{(\alpha)}(\xi+i\tau N)\right|^2}\leq C.
\end{equation}$\blacklozenge$

\bigskip

\textbf{Proof of Proposition \ref{105-CauNire}.}

In order to prove that (a) implies (b) we argue by contradiction. We assume that (a) holds and that (b) does not hold.
Hence, we assume that there exists $(\xi_{\star}, \tau_{\star})\in
\mathbb{R}^{n+1}\setminus\{(0,0)\}$ satisfying:

$$P_m(\xi_{\star}+i\tau_{\star} N)=0\quad\mbox{and}\quad\sum_{j=1}^n \left|P^{(j)}_m(\xi_{\star}+i\tau_{\star}
N)\right|^2=0.$$ From this, we easily obtain 

$$\lim_{\mu\rightarrow +\infty}q(\mu\xi_{\star},\mu
\tau_{\star})=+\infty$$ which contradicts 
\eqref{condizione-proposiz}. Hence, (a) implies (b).

Now, let us suppose that (b) holds true. Let $$\mathbb{S}^n=\left\{(\xi,\tau)\in
\mathbb{R}^{n+1}:|\xi|^2+\tau^2=1\right\}.$$ By Lemma
\ref{106-CauNire}, there exists $C$, which we may assume larger than $1$,
which satisfies

\begin{equation}\label{condizione-proposiz-1}
C|P_m(\xi+i\tau N)|^2+\sum_{j=1}^n \left|P^{(j)}_m(\xi+i\tau
N)\right|^2>0, \quad \forall (\xi,\tau)\in \mathbb{S}^n.
\end{equation}
Therefore, since the polynomials $P_m$ and $P_m^{(j)}$ are homogeneous of degree
$m$ and $m-1$ respectively, there exists  $\lambda>0$ such that,
for each $(\xi,\tau)\in \mathbb{R}^{n+1}$ we have

\begin{equation}\label{condizione-proposiz-2}
	\begin{aligned}
\gamma(\xi,\tau)&:=\frac{|P_m(\xi+i\tau
N)|^2}{\left(|\xi|^2+\tau^2\right)}+\frac{1}{C}\sum_{j=1}^n
\left|P^{(j)}_m(\xi+i\tau N)\right|^2\geq\\&\geq
\lambda\left(|\xi|^2+\tau^2\right)^{m-1}.
\end{aligned}
\end{equation}
On the other hand, we have trivially that, for some constant
$\widetilde{C}$, we get

\begin{equation*}\label{condizione-proposiz-3}
q_{m-1}(\xi,\tau)\leq \widetilde{C}, \quad |\xi|^2+\tau^2\leq 1
\end{equation*}
and by \eqref{condizione-proposiz-2}, we have

\begin{equation}\label{condizione-proposiz-4}
	\begin{aligned}
q_{m-1}(\xi,\tau)&\leq\frac{\left(|\xi|^2+\tau^2\right)^{m-1}}{\sum_{|\alpha|\leq
m-2}\left|P_m^{(\alpha)}(\xi+i\tau N)\right|^2+\gamma(\xi,\tau)}\leq \\& \leq
\lambda^{-1}, \qquad\qquad\qquad \mbox{for }|\xi|^2+\tau^2\geq 1.
\end{aligned}
\end{equation}
Hence
$$q_{m-1}(\xi,\tau)\leq\max\left\{\widetilde{C},\lambda^{-1}
\right\}, \quad \forall (\xi,\tau)\in \mathbb{R}^{n+1}.$$ Now, taking into account \textbf{Remark 1} of this Section,
(a) follows. $\blacksquare$

\bigskip

\textbf{Examples}

\noindent\textbf{1.} Let $P_2(\xi)=\sum_{j=1}^n\xi_j^2$. We know that the corresponding operator is
$P_2(D)=\sum_{j=1}^nD_j^2=-\Delta$. We already know that the following estimate holds

\begin{equation}\label{010-Nire}
\int_{Q_1}e^{2\tau N\cdot x}\left(|u|^2+|\nabla u|^2\right) dx\leq
C\int_{Q_1}e^{2\tau N\cdot x}\left|\Delta u\right|^2dx,
\end{equation}
for every $u \in C^{\infty}_{0}(Q_1,\mathbb{C})$ and for every $\tau\in
\mathbb{R}$. On the other hand, it is immediate to see that (b)
of Proposition \ref{105-CauNire} is satisfied because

$$\sum_{j=1}^n\left|P^{(j)}_2(\xi+i\tau
N)\right|^2=4\left(|\xi|^2+\tau^2\right)>0, \quad \forall
(\xi,\tau)\in \mathbb{R}^{n+1}\setminus \{(0,0)\}.$$ $\spadesuit$

\medskip

\noindent\textbf{2.} A similar argument applies to the
wave operator $\Box=D^2_0-\Delta$ (here $x_0$ represents the time variable). Since the symbol of the operator $\Box$ is $P_2(\xi)=-\xi_0^2+\sum_{j=1}^{n}\xi_j^2$, we have

$$\sum_{j=0}^n\left|P^{(j)}_2(\xi+i\tau
N)\right|^2=4\left(|\xi|^2+\tau^2\right)>0, \quad \forall
(\xi,\tau)\in \mathbb{R}^{n+1}\setminus \{(0,0)\}.$$ Hence (b)
of Proposition \ref{105-CauNire} is satisfied. $\spadesuit$

\medskip

\noindent\textbf{3.} Whereas, the fourth order operator 
$\Delta^2$, which is called \textit{bilaplacian} or, in
dimension $2$, \textit{plate operator}, does not satisfy
(b) of Proposition \ref{105-CauNire}. In this case we have 
$$P_4(\xi)=\left(\sum_{j=1}^n\xi_j^2\right)^2.$$ Hence
$$P_4(\xi+i\tau N)=\left(|\xi|^2+2i\tau \xi\cdot N-\tau^2\right)^2.$$ Which implies that
$P_4(\xi+i\tau N)=0$ holds if and only if

\begin{equation}\label{esempi-Nire}
\begin{cases}
|\xi|^2-\tau^2=0, \\
\\
\tau \xi\cdot N=0.\\
\end{cases}
\end{equation}
The system above, clearly has non-zero $(\xi,\tau)$ solutions. On the other hand for these values we have

$$\sum_{j=1}^n\left|P^{(j)}_4(\xi+i\tau
N)\right|^2=16\left(|\xi|^2+\tau^2\right)\left\vert|\xi|^2+2i\tau
\xi\cdot N-\tau^2\right\vert^2=0.$$ Therefore $P_4$ does not satisfy (b) of Proposition
\ref{105-CauNire}. $\spadesuit$

\bigskip

\noindent\textbf{4a.} We will now consider more carefully the
case in which $P_2(D)$ is a \textbf{second order elliptic operator with complex coefficients}.
\begin{equation}\label{ell-Nire-1}
P_2(D)=\sum_{j,k=1}^na_{jk}D_{jk}^2,
\end{equation}
where $a_{jk}\in \mathbb{C}$, $a_{jk}=a_{kj}$, for $j,k=1,\cdots, n$.
Set 

\begin{equation}\label{ell-Nire-2}
a(\xi,\eta)=\sum_{j,k=1}^na_{jk}\xi_{j}\eta_{k}, \quad \xi,\eta\in
\mathbb{R}^n.
\end{equation}
Let $N\in \mathbb{R}^n$ a versor and let us suppose that

\begin{equation}\label{ell-Nire-3}
a(N,N)=1.
\end{equation}
Let us notice that  we may always reduce to the situation
\eqref{ell-Nire-3} since by the ellipticity of
$P_2(D)$, we have $a(N,N)=P_2(N)\neq 0$ and, consequently we may divide
from the beginning all the coefficients of $P_2(D)$ by $a(N,N)$
leading us back to the assumption  \eqref{ell-Nire-3}. Hence we have

\begin{equation}\label{ell-Nire-4}
P_2(\xi+i\tau N)=a(\xi,\xi)+2i\tau a(\xi,N)-\tau^2.
\end{equation}

Let us observe that, as  $P_2$ is elliptic, we have that, if 
$(\xi,\tau)\in \mathbb{R}^{n+1}\setminus\{(0,0)\}$ satisfies
\begin{equation}\label{ell-Nire-4-401}
P_2(\xi+i\tau N)=0,
\end{equation} then \textbf{both $\xi$ and $\tau$ need to be
different from zero}. As a matter of fact, if  $\tau=0$ then, as 
 $P_2(D)$ is elliptic, the unique solution of equation 
\begin{equation}\label{ell-Nire-5}
P_2(\xi+i0 N)=0,
\end{equation}
is $\xi=0$. On the other hand, if we had $\xi=0$ then by
the homogeneity of $P_2$ we would have

$$0=P_2(i\tau N)=-\tau^2P_2(N),$$
and by the ellipticity of $P_2$ we would have $\tau=0$.

Moreover, by the homogeneity and the ellipticity of $P_2$, we have that if
$(\xi,\tau)\neq (0,0)$ satisfies \eqref{ell-Nire-4-401}, then $\xi$
and $N$ must be linearly independent. Let us prove the last sentence arguing by contradiction. If $\xi$ 
and $N$ were not linearly dependent then there would exist $a,b\in \mathbb{R}$
not both zero, such that $a\xi+bN=0$. Let us suppose, for instance, $a\neq 0$, then $\xi=-\frac{b}{a}N$. Hence
$$0=P_2(\xi+i\tau N)=P_2\left(\left(-\frac{b}{a}+i\tau\right)
N\right)=\left(-\frac{b}{a}+i\tau\right)^2P_2(N),$$ from which
$P_2(N)=0$, on the other hand this cannot occur because $P_2$ is elliptic and $N\neq 0$. Similarly, we proceed assuming $b\neq
0$.

In the \textbf{case of real coefficients} it is easily seen that
(b) of Proposition \ref{105-CauNire} is satisfied.
Indeed, having to consider only the solutions $(\xi,\tau)\neq
(0,0)$ of equation \eqref{ell-Nire-4-401}, we would get $\xi\neq 0$ and
$\tau\neq 0$ and then to establish (b) it suffices to check that equation \eqref{ell-Nire-4-401}, (considered in the unknown $\tau$)
has no nonzero real double roots. Let, therefore, $\xi_0\neq
0$ and $\tau_0\neq 0$ such that

\begin{equation}\label{ell-Nire-4-401-1}
P_2(\xi_0+i\tau_0 N)=0
\end{equation}
and let us assume that

\begin{equation}\label{ell-Nire-6}
\frac{d}{d\tau}P_2(\xi_0+i\tau
N)_{|\tau=\tau_0}=\sum_{j=1}^nP^{(j)}_2(\xi_0+i\tau N)N_j=0.
\end{equation}
Since
$$\frac{d}{d\tau}P_2(\xi_0+i\tau
N)_{|\tau=\tau_0}=2ia(\xi_0,N)-2\tau_0,$$ by \eqref{ell-Nire-6}
we get $\tau_0=ia(\xi_0,N)$ and, taking into account
\eqref{ell-Nire-4}, we have

$$P_2(\xi_0+i\tau_0 N)=a(\xi_0,\xi_0)-(a(\xi_0,N))^2.$$
On the other hand, since $\xi_0$ e $N$ since $\xi_0$ and $N$ are linearly
independent, by the Cauchy--Schwarz inequality we have
(recall $a(N,N)=1$)
$$a(\xi_0,\xi_0)-(a(\xi_0,N))^2>0,$$ which contradicts
\eqref{ell-Nire-4-401-1}. Therefore, if \eqref{ell-Nire-4-401-1} holds true, then \eqref{ell-Nire-6} cannot be true, consequently
 
$$\sum_{j=1}^n \left|P^{(j)}_2(\xi_0+i\tau_0
N)\right|^2>0.$$

\bigskip

Let us consider now the \textbf{case of complex coefficients}.

\medskip

If $n=2$ and we consider
$$P_2(D)=-D^2_1+2iD^2_{12}+D^2_2=\left(iD_1+D_2\right)^2,$$ then $P_2(D)$ is elliptic, but it is easy to check that (b) of Proposition
\ref{105-CauNire} does not hold. As a matter of fact, we have
$$|P_2(\xi+i\tau N)|^2=\left((\xi_1+\tau N_2)^2+(\xi_2-\tau N_1)^2\right)^2$$
and
$$\sum_{j=1}^2|P^{(j)}_2(\xi+i\tau N)|^2=8\left((\xi_1+\tau N_2)^2+(\xi_2-\tau N_1)^2\right).$$ Hence, if $\xi_0=(-N_2,N_1)$ and $\tau=1$, then we have

\begin{equation*}
|P_2(\xi_0+i N)|^2=0
\end{equation*}
and
$$\sum_{j=1}^2|P^{(j)}_2(\xi_0+i N)|^2=0.$$

Now, let us prove that if $n\geq 3$ then (b) of Proposition
\ref{105-CauNire} is satisfied. We prove, like in the case of the
real coefficients, that if  $(\xi_0,\tau_0)\neq (0,0)$ satisfies
\begin{equation}\label{ell-Nire-7}
P_2(\xi_0+i\tau_0 N)=0,
\end{equation}
then
\begin{equation}\label{ell-Nire-7-401}
\sum_{j=1}^nP^{(j)}_2(\xi_0+i\tau_0
N)N_j=\frac{d}{d\tau}P_2(\xi_0+i\tau N)_{|\tau=\tau_0}\neq 0.
\end{equation}
We argue by contradiction. Let us assume $\xi_0\in \mathbb{R}^n$ and
$\tau_0\in \mathbb{R}$, where $(\xi_0,\tau_0)\neq (0,0)$,  satisfy
\eqref{ell-Nire-7} and let us assume that \eqref{ell-Nire-7-401} does not hold, i.e.  let us assume $\tau_0$ is a double root of equation in $\tau$
\begin{equation}
\frac{d}{d\tau}P_2(\xi_0+i\tau N)_{|\tau=\tau_0}= 0.
\end{equation}
We already noticed that we must have $\xi_0\neq 0$ and
$\tau_0\neq 0$. Since $\tau_0$ is a double solution of the equation
\eqref{ell-Nire-4}, the discriminant of that equation (when
$\xi=\xi_0$) is null. Hence

\begin{equation}\label{ell-Nire-8}
a(\xi_0,\xi_0)=\left(a(\xi_0,N)\right)^2
\end{equation}
and

\begin{equation}\label{ell-Nire-9}
\tau_0=ia(\xi_0,N).
\end{equation}
By \eqref{ell-Nire-8} and \eqref{ell-Nire-9} we have

\begin{equation}\label{ell-Nire-10}
a(\xi_0,\xi_0)=-\tau_0^2,\quad a(\xi_0,N)=-i\tau_0\quad \mbox{ and
(recall) }\quad a(N,N)=1.
\end{equation}
Now, let $\eta$ be  a vector of $\mathbb{R}^n$ such that $\xi_0$, $N$ and
$\eta$ be linearly independent (recall $n\geq 3$).
Let $\mathcal{B}$ a basis of $\mathbb{R}^n$ which complete
$\{\xi_0,N,\eta\}$ and let us write the matrix of bilinear form $a$ w.r.t. 
$\mathcal{B}$. Set
$$\widetilde{a}_{11}=a(\xi_0,\xi_0)=-\tau_0^2,\quad
\widetilde{a}_{12}=\widetilde{a}_{21}=a(\xi_0,N), \quad
\widetilde{a}_{22}=a(N,N),$$
$$\widetilde{a}_{33}=a(\eta,\eta)=\alpha_{33}+i\beta_{33},$$
and
$$\widetilde{a}_{13}=\widetilde{a}_{31}=a(\xi_0,\eta):=\alpha_{13}+i\beta_{13},\quad
\widetilde{a}_{23}=\widetilde{a}_{32}=a(N,\eta):=\alpha_{23}+i\beta_{23},$$
 where $\alpha_{13}, \beta_{13}, \alpha_{23},
\beta_{23},\alpha_{33},\beta_{33},$ are real numbers.

Let us consider the vector $v$ of $\mathbb{R}^n\setminus{0}$ whose
components with respect to the base $\mathcal{B}$ have coordinates
represented by the vector $(x,y,z,0,\cdots, 0)$. Thus, let us note that
$$v\neq 0\Leftrightarrow (x,y,z)\neq (0,0,0).$$
Since $P_2(D)$ is elliptic and $v\neq 0$ we have
\begin{equation}\label{ell-Nire-010}
a(v,v)\neq 0,
\end{equation}
in turn, by \eqref{ell-Nire-10}, this is equivalent to the fact
that for $(x,y,z)\neq (0,0,0)$ we have

\begin{equation*}
	\begin{aligned}
	&x^2-2i\tau_0xy-\tau_0^2y^2+2(\alpha_{13}+i\beta_{13})xz+\\&+2(\alpha_{23}+i\beta_{23})yz+
	(\alpha_{33}+i\beta_{33})z^2\neq 0,	
\end{aligned}
\end{equation*}
for every $(x,y,z)\in \mathbb{R}^3\setminus\{(0,0,0)\}$.  Now, the above condition is equivalent to the fact that $(0,0,0)$ is the unique solution to the algebraic system 

\begin{equation}\label{ell-Nire-11-28}
\begin{cases}
x^2-\tau_0^2y^2+2\alpha_{13}xz+2\alpha_{23}yz+
\alpha_{33} z^2=0,\\
\\
-2\tau_0xy+2\beta_{13}xz+2\beta_{23}yz+
+\beta_{33}z^2=0.\\
\end{cases}
\end{equation}
But this cannot occur.
Let us see why.

First of all, let us recall that $\tau_0\neq 0$. Moreover, let us suppose that
$z\neq 0$. Then, if we set $$X=\frac{x}{z}, \quad Y=\frac{y}{z},$$ 
system \eqref{ell-Nire-11} become

\begin{equation}\label{ell-Nire-11}
\begin{cases}
X^2-\tau_0^2Y^2+2\alpha_{13}X+2\alpha_{23}Y+
\alpha_{33}=0\\
\\
-2\tau_0XY+2\beta_{13}X+2\beta_{23}Y+
\beta_{33}=0\\
\end{cases}
\end{equation}
and it is simple to check that system \eqref{ell-Nire-11}
admits \textit{always} solutions. To convince yourself of this, it suffices to notice
that the asymptotes of the first hyperbole (possibly degenerate) in
\eqref{ell-Nire-11} are parallel to the straight lines $X=\pm \tau_0 Y$ which
must necessarily meet the asymptotes of the second
hyperbola (possibly degenerate) that are parallel to the coordinate axes.
From above it follows that there exists
$v\in\mathbb{R}^n\setminus{0}$ such that $a(v,v)=0$ and this
contradicts \eqref{ell-Nire-010}. In summary if $n\geq 3$ then \eqref{ell-Nire-7-401} must hold, and this implies that (b) of the
Proposition \ref{105-CauNire} is satisfied. $\spadesuit$

\medskip

\noindent\textbf{4b.}  We conclude by considering the case where $P_2(D)$ is a
\textbf{non-elliptic} operator of second order, with
\textbf{real coefficients}

 In such a case there is at least one
 \emph{characteristic direction}. We recall that
 $N\in\mathbb{R}^n\setminus{0}$ is a characteristic direction
 with respect to the operator $P_2(D)$ provided we have
$$P_2(N)=0.$$ In this case, the planes
$${N\cdot x=c},$$ where $c\in \mathbb{R}$ are characteristic surfaces. Let $N$ be a characteristic direction, if we have
$$\sum_{j=1}^n \left|P^{(j)}_2(N)\right|^2>0,$$ we  say
that ${N\cdot x=c}$, $c\in \mathbb{R}$, is a
\textbf{simple characteristic} \index{characteristic:@{characteristic:}!- simple@{- simple}}
If we have $$P_2(N)=0$$ and
$$\sum_{j=1}^n \left|P^{(j)}_2(N)\right|^2=0,$$
we say that $\{N\cdot x=c\},$ is a \textbf{double characteristic} \index{characteristic:@{characteristic:}!- double@{- double}}

For instance, the wave operator has only
simple characteristics (reader check) while
\textbf{the heat operator} and \textbf{the Schr\"{o}dinger operator} \index{Schr\"{o}dinger operator}have double characteristics. The heat operator is given by (for $n>1$)
  
$$P(D)=-\sum_{j=1}^{n-1}D^2_j-iD_n=\Delta_{x'}-\partial_n.$$
So the principal part of the heat operator is
$$P_2(D)=-\sum_{j=1}^{n-1}D^2_j,$$
whose symbol is $-\sum_{j=1}^{n-1}\xi^2_j$.
 It is evident that
the unique characteristic directions are those generated by the versor
$e_n$ which is a double characteristic direction. In the case
of the Schr\"{o}dinger operator we have.

$$P(D)=-\sum_{j=1}^{n-1}D^2_j-D_n=\Delta_{x'}-\frac{1}{i}\partial_n$$
and it is once again clear that the only characteristic directions are those generated by the $e_n$ versor and, as in the case of the heat operator, they are double.

\medskip

Let us check that the operators with simple characteristics
satisfy (b) of Proposition \ref{105-CauNire} for all the versors
$N$.

\noindent Let $N$ be a versor of $\mathbb{R}^n$ such that

\begin{equation}\label{gen-Nire-1}
P_2(N)\neq 0.
\end{equation}
Let $(\xi,\tau)\in \mathbb{R}^{n+1}\setminus\{0\}$ satisfy

\begin{equation}\label{gen-Nire-2}
P_2(\xi+i\tau N)=0.
\end{equation}
Let us denote by $A$ the symmetric matrix  $ \left\{a_{jk}\right\}^n_{
j,k=1}$. Since the coefficients of $P_2(D)$ are real numbers, we have

\begin{equation}\label{gen-Nire3}\sum_{j=1}^n \left|P^{(j)}_2(\xi+i\tau
N)\right|^2=4\left(|A\xi|^2+\tau^2|AN|^2\right).
\end{equation}
Now, if $\tau=0$, then \eqref{gen-Nire-2} implies $P_2(\xi)=0$,
however $P_2(D)$ has only simple characteristics, hence
\begin{equation*}\sum_{j=1}^n \left|P^{(j)}_2(\xi+i 0
N)\right|^2>0.
\end{equation*}
If $\tau\neq0$, then $AN\neq 0$, otherwise, if it were $AN=0$
we would have $P_2(N)=AN\cdot N=0$ which would contradict the
\eqref{gen-Nire-1}. Therefore from \eqref{gen-Nire3} we have

\begin{equation*}\label{gen-Nire-3}\sum_{j=1}^n \left|P^{(j)}_2(\xi+i\tau
N)\right|^2\geq 4\tau^2|AN|^2>0.
\end{equation*}
Hence, if \eqref{gen-Nire-1} holds true, then  (b) of Proposition
\ref{105-CauNire} holds true.

If $N$ is a simple characteristic, then

\begin{equation*}
P_2(N)= 0.
\end{equation*}
and, as $N$ is a simple characteristic, we have

\begin{equation}\label{gen-Nire5}
AN\neq 0.
\end{equation}
Now, let $(\xi,\tau)\in \mathbb{R}^{n+1}\setminus\{0\}$ be a solution to 
\eqref{gen-Nire-2}. If $\tau=0$, we have $P_2(\xi)=0$, hence
\begin{equation*}\sum_{j=1}^n \left|P^{(j)}_2(\xi)\right|^2>0.
\end{equation*}
On the other hand, if $\tau\neq0$ from \eqref{gen-Nire5}, as already seen
above, we have
\begin{equation*}\sum_{j=1}^n \left|P^{(j)}_2(\xi+i \tau
N)\right|^2>0.
\end{equation*}
Therefore, even when $N$ is a simple characteristic direction, (b)
of Proposition \ref{105-CauNire} is satisfied.

Finally, let us consider the case in which $P_2(D)$ has a 
double characteristic direction; be it $\eta$, then any way one chooses the versor
$N$, we have

\begin{equation}\label{gen-Nire-2-16-8}
P_2(\eta+i0 N)=0
\end{equation}
and
\begin{equation}\label{gen-Nire3-28}
\sum_{j=1}^n \left|P^{(j)}_2(\eta+i
0 N)\right|^2=0.
\end{equation}
Therefore, in this case, (b) of Proposition \ref{105-CauNire}
is not satisfied. This implies that the following estimate \textbf{does not hold}

\begin{equation*}
\int_{Q_1}e^{2\tau N\cdot x}\left|u\right|^2dx+ \int_{Q_1}e^{2\tau
N\cdot x}\left|Du\right|^2dx\leq C\int_{Q_1}e^{2\tau N\cdot
x}\left|P_2(D)u\right|^2dx,
\end{equation*}
for every $u \in C^{\infty}_{0}(Q_1,\mathbb{C})$ and for every $\tau\geq
\tau_0$. 

Let us notice that in the case of the heat operator and of
Schr\"{o}dinger operator, the principal part $P_2(D)$ excludes the term
$iD_nu$. Actually, if we employ directly estimate
\eqref{1.14-corollario-sece} we have (reader check)

\begin{equation}\label{gen-Nire5-16-8}
	\begin{aligned}
		&\int_{Q_1}e^{2\tau N\cdot x}\left|u\right|^2dx+ \int_{Q_1}e^{2\tau
N\cdot x}\left|\nabla u\right|^2dx\leq \\&\leq  C\int_{Q_1}e^{2\tau N\cdot
x}\left|\Delta_{x'} u-\partial_n\right|^2dx,
\end{aligned}
\end{equation}
for every $u \in C^{\infty}_{0}(Q_1,\mathbb{C})$ and for every $\tau\in
\mathbb{R}$.

As in previous situations, \eqref{gen-Nire5-16-8} implies the uniqueness for the Cauchy problem with \textbf{strictly convex initial surfaces} for the differential inequalities
\begin{equation}\label{gen-Nire6}
\left|\Delta_{x'} U-\partial_nU\right|\leq
M\left(\left|\nabla_{x'}U\right|+\left| U\right|\right),
\end{equation}
where $M$ is a positive number and $U$ is enough regular.

When $n=2$, one can exploit the particularity of the
dimension two to prove the following unique continuation property, we refer to \cite[Theorem 9]{Ni} for the proof:

Let $\omega$ be an open of $\mathbb{R}^2=\mathbb{R}_x\times
\mathbb{R}_t$ contained in a rectangle $\mathcal{R}$. For $t_0\in
\mathbb{R}$, we denote by $s_{t_0}$ straight line of equation $t=t_0$
and set (Figure 12.2)
\begin{figure}\label{calore1}
	\centering
	\includegraphics[trim={0 0 0 0},clip, width=10cm]{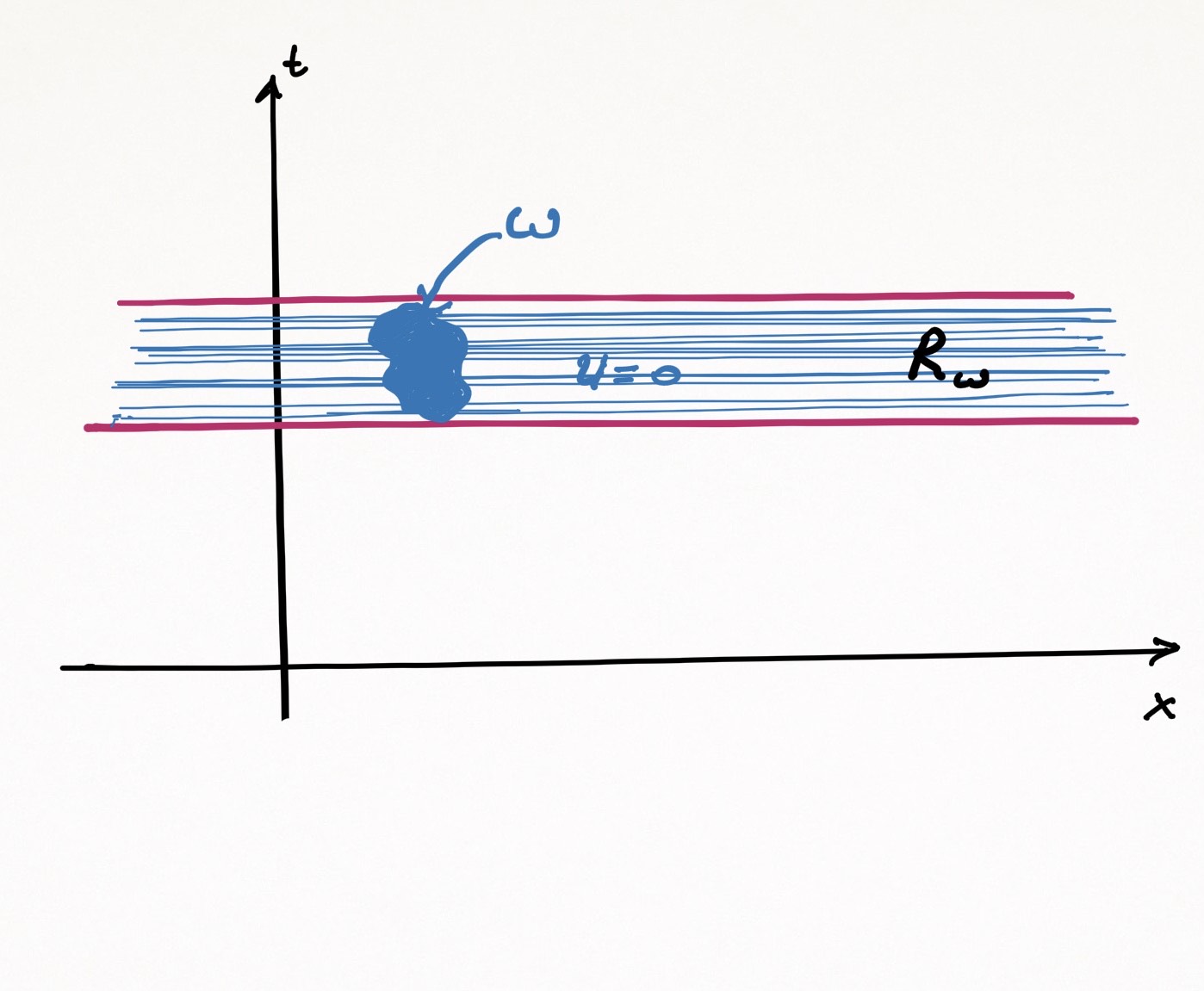}
	\caption{}
\end{figure}

$$\mathcal{R}_{\omega}=\left\{(x,t_0)\in \mathcal{R}: s_{t_0}\cap \omega \neq\emptyset\right\}.$$
Let us assume

\begin{equation}\label{gen-Nire7}
\left| \partial^2_{x}U-\partial_tU\right|\leq
M\left(\left|\partial_{x}U\right|+\left| U \right|\right), \quad
\mbox{ in } \mathcal{R},
\end{equation}
then
\begin{equation}\label{gen-Nire8}
U(x,t)=0, \quad \forall (x,t)\in \mathcal{R}_{\omega}.
\end{equation}

From the previous result it follows, in particular, that if $U$ is
regular enough (it is sufficient that $U,
\partial_xU,\partial_t U, \partial_x^2U$ are continuous in
$(0,1)\times(0,1)$) and if $U$ is a solution of the Cauchy problem

\begin{equation}\label{gen-Nire9}
\begin{cases}
\partial_x^2U-\partial_tU=a(x,t)\partial_xU+b(x,t) U, & \mbox{in}\quad (0,1)\times(0,1), \\
\\
U(0,t)=0,
 &  \mbox{ for } t\in (\alpha,\beta),\\
 \\
 \partial_xU(0,t)=0, &  \mbox{ for } t\in (\alpha,\beta),\\
\end{cases}
\end{equation}
where $\alpha,\beta$ are given numbers which satisfy $0<\alpha<\beta<1$ and $$a,b
\in L^{\infty}((0,1)\times(0,1)),$$ then $$U=0,\quad\mbox{in }
(0,1)\times (\alpha,\beta).$$

We do not enter into the details and refer the interested reader
directly to \cite[Theorem 9]{Ni}. $\spadesuit$

\bigskip

\section{Chapter summary and conclusions} \label{conclusioni}
In this chapter we have proved estimate \eqref{1.15-sece}, in a relatively simple manner. The most relevant peculiarity of such an estimate is that in it there is a "weight" which depends
on a parameter $\tau$ that may be arbitrarily large.

By applying estimate \eqref{1.15-sece} to the Laplace operator and
by exploiting some important invariance properties of this
operator, we have proved, in Theorem \ref{33-Nire},
the global uniqueness for the Cauchy problem for the equation

$$\Delta U=b(x)\cdot\nabla U+c(x)U,$$
where $b,c\in L^{\infty}$.

We proved that the estimate \eqref{1.15-sece}
allows us to prove the global uniqueness for the Cauchy problem with strictly convex initial surface, for the operators $P_m(D)+a_0(x)$, where $P_m(D)$ is a homogeneous operator with
constant coefficients and $a_0\in L^{\infty}$.

We have shown with several examples and remarks related to Theorem
\ref{cond-nec-Nire}, the strict connections that exists between an
estimate of the type \eqref{43-Nire} and some properties of the symbols of operators $M(D)$ and
$PD)$. These connections  makes it possible to  transfer into the
algebraic field the estimates under investigation in this Chapter.

\bigskip

The estimates considered in this Chapter have two remarkable
\textbf{weaknesses} that we now briefly discuss.

\medskip

\textbf{1.} The first weakness lies in the character
of the weight exponent. This is because such an exponent is linear and,
as we have seen, this greatly limits the geometry in which to
apply our estimates. Let us consider, for instance the following Cauchy problem
$$P(D)U=0,\quad \mbox{in } \mathbb{R}^n$$
and $$U=0,\quad \mbox{for } x_n\leq 0,$$ to prove the uniqueness we would be most helped by
a weight whose level surfaces are "curved" with respect to the
$x_n=c$ planes.
 More precisely, if instead of the weight $e^{-2\tau x_n}$ we dispose of estimates with weight $e^{2\tau
	\left(-x_n+\frac{\delta}{2}|x|^2\right)}$, with $\delta>0$ ( even
small), it could be shown that $U$ vanishes in regions of the
type $\left\{-x_n+\frac{\delta}{2} |x|^2<r, x_n>0\right\}$ with
$r>0$ (Figure 12.3).

\begin{figure}\label{D9}
	\centering
	\includegraphics[trim={0 0 0 0},clip, width=10cm]{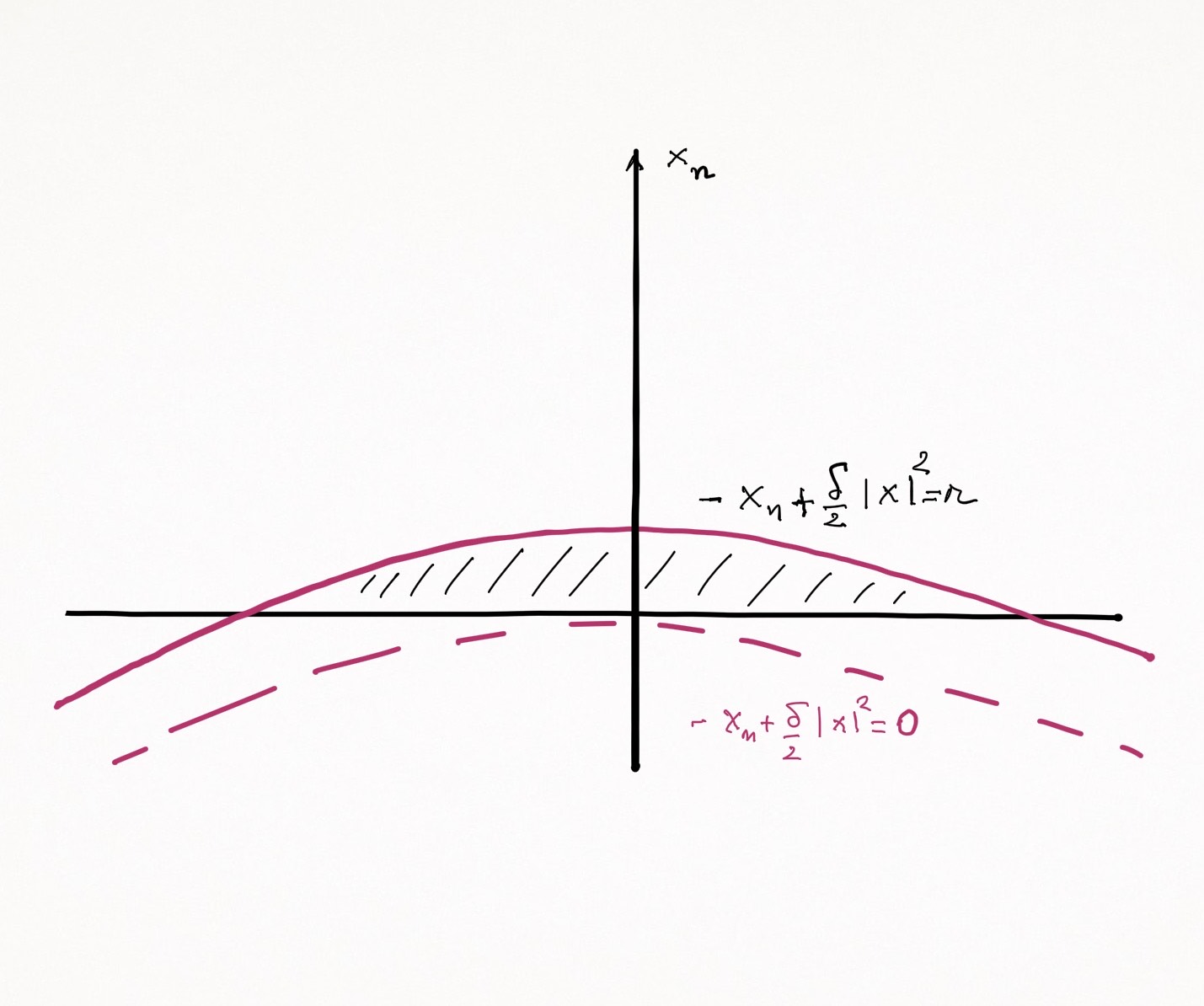}
	\caption{}
\end{figure}

\noindent However, the proof of an estimate that corresponds to
this nonlinear weight would not allow repeat, in a
simple and immediate way, the proof of Theorem
\ref{1.14-corollario-sece}. To realize this, let us observe
that, setting $v=e^{\tau \left(N\cdot
	x+\frac{\delta}{2}|x|^2\right)}u$, instead of the
\eqref{per-conclusioni} we would have
$$D_ju=e^{-\tau N\cdot x}\left(D_j+i\tau \left(N_j +\delta
x_j\right)\right)v$$ and we cannot use in \textit{immediate manner}
the Fourier transform. For the time being, we refer to \cite{Tr1} for
further discussion. 

\medskip

\textbf{2.} The other weakness of estimate \eqref{1.15-sece} consists in the fact that these
estimates hold for operators with constant coefficients in the principal part.

Let $P_m(x,D)$ be the principal part of the operator. One might
be tempted to consider $P_m(x,D)$, in a neighborhood
$\mathcal{U}$ of a point $x_0$, as the pertubation of the operator
with constant coefficients $P_m(x_0,D)$ i.e., we could write
$$P_m(x,D)=P_m(x_0,D)+\left(P_m(x_0,D)-P_m(x,D)\right)$$ and,
by exploiting the regularity of the coefficients of $P_m(x,D)$,
we may consider that, for a suitable neighborhood $\mathcal{U}$, we have

\begin{equation}\label{conclusion}
\left|P_m(x_0,D)U-P_m(x,D)U\right|<\varepsilon|D^m U|, \quad
\forall x\in \mathcal{U}.
\end{equation}
To clarify the idea further let us show, in broad terms, that if
$P_m(x,D)$ is elliptic with continuous coefficients then we
we have

\begin{equation}\label{garding}\int_{Q_1}\left|D^mu\right|^2dx\leq
C\int_{Q_1}\left|P_m(x,D)u\right|^2dx,\quad \forall u \in
C^{\infty}_{0}(Q_1,\mathbb{C}).
\end{equation}
We have already seen that, in the case of elliptic operators with
constant coefficients, \eqref{tau-zero-6} is valid. Now let
$x_0=0$, $\varepsilon>0$ and let $r_{\varepsilon}$ be such that the
\eqref{conclusion} is satisfied for $\mathcal{U}=Q_{r_{\varepsilon}}$.

\noindent Then by estimate \eqref{tau-zero-6} (which we know is
true in the elliptic case) we have 
\begin{equation*}\int_{Q_1}\left|D^mu\right|^2dx\leq
2C\int_{Q_1}\left|P_m(x,D)u\right|^2dx+
2C\varepsilon^2\int_{Q_1}\left|D^mu\right|^2dx,
\end{equation*}
for every $u \in
C^{\infty}_{0}(Q_{r_{\varepsilon}},\mathbb{C})$. It is then evident that by choosing
$\varepsilon=\frac{1}{2\sqrt{C}}$ we get
\begin{equation*}
\int_{Q_{r_{\varepsilon}}}\left|D^mu\right|^2dx\leq
4C\int_{Q_{r_{\varepsilon}}}\left|P_m(x,D)u\right|^2dx,\quad \forall
u \in C^{\infty}_{0}(Q_{r_{\varepsilon}},\mathbb{C}).
\end{equation*}
A similar argument can be made in the neighborhood of the other
points of $Q_1$ and, using a partition of the unity, one can
obtain \eqref{garding}.

It is easily understood that a similar argument for the estimates of the type
\eqref{103-CauNire} -- that is, considering $P_m(x,D)$ as a
pertubation of order $m$ in a neighborhood of $x_0$ of the operator with
constant coefficients $P_m(x_0,D)$ -- is difficult to realize even in the elliptic case. Actually, estimate
\eqref{tau-zero-7} does not hold even for elliptic operators and thus
the following "error term" that would follow from \eqref{conclusion},

 $$C\varepsilon^2\int_{Q_1}e^{2\tau N\cdot x}\left|D^mu\right|^2dx,$$
 cannot be "absorbed" by the terms to the left of the sign of
 inequality. The development of the theory will show that the case of the
 variable coefficients (in the principal part) is more
 tricky than the one with constant coefficients, even in the case of the weight 
 $e^{2{\tau N\cdot x}}$.
 	
 \newpage

 \chapter{Carleman estimates and the Cauchy problem I -- Elliptic operators}\markboth{Chapter 13. Carleman estimates and the Cauchy problem I}{}\label{Carleman}

\section{Introduction}\label{Carleman-intro}

In the previous Chapter we gave a
first insight into the Carleman estimates by showing how they are used in the investigation of the uniqueness of the Cauchy problem and
of the unique continuation property \index{unique continuation property}for operators whose
principal part has constant coefficients. In Section \ref{conclusioni}, we
pointed out some important weaknesses of  Theorem \ref{1.14-corollario-sece}. Such weak points,
 briefly, consist of:

\medskip

\noindent(a) the linear character of the weight exponent;

\medskip

\noindent(b) estimate \eqref{1.15-sece} of Theorem
\ref{1.14-corollario-sece} cannot be easily extended for
operators with variable coefficients.

\medskip

For such reasons, here we begin a more systematic study of the
Carleman estimates \index{Carleman estimates}in order to extend somewhat the uniqueness results we have seen in Chapter \ref{Nirenberg}. Although we will
focus mainly on elliptic operators (Section
\ref{Carleman-ellittico}), the introductory examples (Sect.
\ref{Carlm-esempio-semplice}) and the framework apply
to other types of operators as well. In Chapter
\ref{operatori-2ord} we will consider the \textbf{second order operators whose
principal part has real coefficients and that are not necessarily
elliptic}.

Let us consider the operator

\begin{equation}  \label{stime-Carlm-29-0}
P(x,D)=\sum_{|\alpha|\leq m}a_{\alpha}(x)D^{\alpha}, \quad\mbox{ }
x\in \Omega,
\end{equation}
where $D_j=\frac{1}{i}\partial_j$, $j=1,\cdots,n$, $\Omega$ is an open set
of $\mathbb{R}^n$ and $a_{\alpha}\in
L^{\infty}(\Omega;\mathbb{C})$, for every $\alpha\in \mathbb{N}^n_0$,
$|\alpha|\leq m$.

We will follow the classical approach developed in \cite{HO63} by L.
H\"{o}rmander. This approach is not only more elementary
than the one based on the pseudodifferential operators \index{pseudodifferential operators}
(\cite[vol. IV]{HOII}, \cite{Lern}), but it allows more easily
to reduce on the regularity assumptions of the coefficients
of the principal parts of the operators.

Let us recall that the symbol of the operator $P(x,D)$ is
$$P(x,\xi)=\sum_{|\alpha|\leq m}a_{\alpha}(x)\xi^{\alpha} \ \  \forall \xi\in \mathbb{R}^n.$$
We further denote by $P_m(x,D)$ the
principal part of $P(x,D)$, i.e.

\begin{equation}  \label{stime-Carlm-29-00}
P_m(x,D)=\sum_{|\alpha|= m}a_{\alpha}(x)D^{\alpha}
\end{equation}
(of course we assume that $|\alpha|=m$, $a_{\alpha}$ is not identically zero for at least one $\alpha$ such that
).

Let $\varphi$ be a sufficiently regular real-valued function,
say $\varphi\in C^{\infty}\left(\overline{\Omega}\right)$, however in
many cases a less regularity will suffice.

Let $\mu\geq 0$, we are interested in the Carleman estimates such as

\begin{equation}  \label{stime-Carlm-29-1}
\tau^{\mu}\sum_{|\alpha|\leq m-1}\int_{\Omega}\left\vert
D^{\alpha}u\right\vert^2 e^{2\tau \varphi}dx\leq C
\int_{\Omega}\left\vert P_m(x,D)u\right\vert^2e^{2\tau \varphi}dx,
\end{equation}
for every $u\in C_0^{\infty}(\Omega)$ and for every $\tau\geq \tau_0$,
where $C$ and $\tau_0$ are constants \textit{independent of} $u$
\textit{and} $\tau$.

In this Chapter we will be interested in the case where $$\mu>0.$$ In
this case, it is simple to check that the estimate
\eqref{stime-Carlm-29-1} is equivalent to a similar estimate where
$P(x,D)$ is replaced by $P_m(x,D)$. Indeed, let us suppose that estimate
\eqref{stime-Carlm-29-1} holds and let us denote by

$$R(x,D)=P(x,D)-P_m(x,D),$$
we have

\begin{equation}  \label{stime-Carlm-29-1-00}
\left\vert R(x,D)u\right\vert=\left\vert \sum_{|\alpha|\leq
m-1}a_{\alpha}(x)D^{\alpha}u\right\vert\leq M\sum_{|\alpha|\leq
m-1}\left\vert D^{\alpha}u\right\vert,
\end{equation}
where
$$M=\max_{|\alpha|\leq
m-1}\left\{\left\Vert
a_{\alpha}\right\Vert_{L^{\infty}(\Omega)}\right\}.$$ Hence, by
\eqref{stime-Carlm-29-1} we get
\begin{equation*}
\begin{aligned}
&\tau^{\mu}\sum_{|\alpha|\leq m-1}\int_{\Omega}\left\vert
D^{\alpha}u\right\vert^2 e^{2\tau \varphi}dx\leq C
\int_{\Omega}\left\vert P_m(x,D)u\right\vert^2e^{2\tau
\varphi}dx\leq\\&\leq 2C \int_{\Omega}\left\vert
P(x,D)u\right\vert^2e^{2\tau \varphi}dx+2C\int_{\Omega}\left\vert
R(x,D)u\right\vert^2e^{2\tau \varphi}dx\leq\\&\leq 2C
\int_{\Omega}\left\vert P(x,D)u\right\vert^2e^{2\tau
\varphi}dx+\widetilde{C}M^2\sum_{|\alpha|\leq
m-1}\int_{\Omega}\left\vert D^{\alpha}u\right\vert^2 e^{2\tau
\varphi}dx,
\end{aligned}
\end{equation*}
for every $u\in C_0^{\infty}(\Omega)$ and for every $\tau\geq \tau_0$.
Moving the last sum to the left hand side, we have 
\begin{equation}  \label{stime-Carlm-30-0}
\left(\tau^{\mu}-\widetilde{C}M^2\right)\sum_{|\alpha|\leq
m-1}\int_{\Omega}\left\vert D^{\alpha}u\right\vert^2 e^{2\tau
\varphi}dx\leq C \int_{\Omega}\left\vert
P(x,D)u\right\vert^2e^{2\tau \varphi}dx,
\end{equation}
for every $u\in C_0^{\infty}(\Omega)$ and for every $\tau\geq \tau_0$.
Now let $\tau_1\geq \tau_0$ be a number  such that for every $\tau\geq \tau_1$ we have

$$\tau^{\mu}-\widetilde{C}M^2\geq \frac{\tau^{\mu}}{2},$$
by \eqref{stime-Carlm-30-0} we obtain
\begin{equation}  \label{stime-Carlm-30-00}
\tau^{\mu}\sum_{|\alpha|\leq m-1}\int_{\Omega}\left\vert
D^{\alpha}u\right\vert^2 e^{2\tau \varphi}dx\leq 4C
\int_{\Omega}\left\vert P(x,D)u\right\vert^2e^{2\tau \varphi}dx,
\end{equation}
for every $u\in C_0^{\infty}(\Omega)$ and for every $\tau\geq \tau_1$.
Hence, if \eqref{stime-Carlm-29-1} holds then 
\eqref{stime-Carlm-30-00} holds. The converse (of course with different
values of $C$ and $\tau_0$) can be similarly proved. 

\bigskip

It should be observed at once that estimates of type
\eqref{stime-Carlm-29-1} (or \eqref{stime-Carlm-30-00}) for $\mu>0$
have a local character in the sense specified in the following

\begin{lem}[\textbf{local character of the Carleman estimates}]\label{stime-Carlm-1.3.3-39}
	\index{Lemma:@{Lemma:}!- local character of the Carleman estimates@{- local character of Carleman estimates}}  
Let \\ $\mu>0$. Let $\Omega$ be a bounded open set of  $\mathbb{R}^n$ and let
$P(x,D)$ be a differential operator whose coefficients belong to 
$L^{\infty}(\Omega)$. Let us assume that for each $y\in
\overline{\Omega}$ there exist $\delta_y>0$, $C_y>0$ and $\tau_y\in
\mathbb{R}$ such that

\begin{equation}\label{stime-Carlm-39-1}
	\begin{aligned}
		&\tau^{\mu}\sum_{|\alpha|\leq m-1}\int_{\Omega\cap
B_{\delta_y}(y)}\left\vert D^{\alpha}u\right\vert^2 e^{2\tau
\varphi}dx\leq\\&\leq  C_y \int_{\Omega\cap B_{\delta_y}(y)}\left\vert
P(x,D)u\right\vert^2e^{2\tau \varphi}dx,
\end{aligned}
\end{equation}
for every $u\in C_0^{\infty}\left(\Omega\cap B_{\delta_y}(y)\right)$
and for every $\tau\geq \tau_y$. Then there exist $C>0$ and $\tau_0\in
\mathbb{R}$ such that
\begin{equation}  \label{stime-Carlm-39-2}
\tau^{\mu}\sum_{|\alpha|\leq m-1}\int_{\Omega}\left\vert
D^{\alpha}u\right\vert^2 e^{2\tau \varphi}dx\leq C
\int_{\Omega}\left\vert P(x,D)u\right\vert^2e^{2\tau \varphi}dx,
\end{equation}
for every $u\in C_0^{\infty}(\Omega)$ and for every $\tau\geq \tau_0$.
\end{lem}

\bigskip

\textbf{Proof.} Since $\overline{\Omega}$ is compact, there exists a finite set of points, $y_1,\cdots,y_N$, such that
$$\overline{\Omega}\subset\bigcup_{j=1}^N B_{\delta_{y_j}}(y_j).$$
Let $\left\{\eta_j\right\}_{1\leq j\leq N}$ be a partition of unity (Lemma \ref{Partiz:23-10-22-1}) such that

$$\eta_j\in C_0^{\infty}\left(B_{\delta_{y_j}}(y_j)\right);\quad
0\leq \eta_j\leq 1;\quad \sum_{j=1}^N \eta_j=1, \quad\mbox{in }
\overline{\Omega}. $$ Let $u\in C_0^{\infty}(\Omega)$ and let us denote by

$$u_j=u\eta_j,\quad\quad\mbox{for } j=1,\cdots N.$$
For every $\alpha\in \mathbb{N}^n_0$, by the Cauchy--Schwarz inequality we have

\begin{equation}  \label{stime-Carlm-40-1}
\left\vert D^{\alpha}u\right\vert^2 =\left\vert
\sum_{j=1}^ND^{\alpha}u_j\right\vert^2\leq N \sum_{j=1}^N\left\vert
D^{\alpha}u_j\right\vert^2.
\end{equation}
Now, we notice that 
$$P(x,D)u_j=\eta_jP(x,D)u+R_{m-1}(x,D)u,$$
where $R_{m-1}(x,D)$ is an operator of order $m-1$ whose
coefficients depend on the coefficients of $P(x,D)$ (but not on
their derivatives) and on the functions $\eta_j$, $j=1,\cdots N$ and
their derivatives of order less or equal to $m$. We have

\begin{equation}  \label{stime-Carlm-40-2}
\left\vert P(x,D)u_j\right\vert\leq \eta_j\left\vert
P(x,D)u\right\vert+M\sum_{|\alpha|\leq m-1}\left\vert
D^{\alpha}u\right\vert,
\end{equation}
where $M$ depends on the $L^{\infty}$ norms of the coefficients of
$P(x,D)$ and on the $L^{\infty}$ norms of the derivatives of $\eta_j$ of order  less or equal to $m$.

By \eqref{stime-Carlm-40-2}, \eqref{stime-Carlm-39-1} and \eqref{stime-Carlm-40-1}
 we get

\begin{equation*}
\begin{aligned}
&\tau^{\mu}\sum_{|\alpha|\leq m-1}\int_{\Omega}\left\vert
D^{\alpha}u\right\vert^2 e^{2\tau \varphi}dx\leq\\&\leq 
N^2\tau^{\mu}\sum_{j=1}^N\sum_{|\alpha|\leq m-1}\int_{\Omega\cap
B_{\delta_{y_j}}(y_j)}\left\vert D^{\alpha}u_j\right\vert^2 e^{2\tau
\varphi}dx\leq\\&\leq N^2\sum_{j=1}^N\int_{\Omega\cap
B_{\delta_{y_j}}(y_j)}\left\vert P(x,D)u_j\right\vert^2 e^{2\tau
\varphi}dx\leq \\&\leq C\int_{\Omega}\left\vert P(x,D)u\right\vert^2
e^{2\tau \varphi}dx+C\sum_{|\alpha|\leq m-1}\int_{\Omega}\left\vert
D^{\alpha}u\right\vert^2 e^{2\tau \varphi}dx,
\end{aligned}
\end{equation*}
for every $\tau\geq \tau_0$, where
$$\tau_0=\max_{1\leq j\leq N}\tau_{y_j}$$ and $C$ is a constant.
Now, we move the last integral to the left hand side and we proceed
exactly as did above to prove
\eqref{stime-Carlm-30-00} and we obtain
\eqref{stime-Carlm-39-2}.
$\blacksquare$

\section[Introductory examples]{Introductory examples -- the first order \\ operators}\label{Carlm-esempio-semplice}

This Section has essentially two purposes: the first one consists of
showing, with simple examples concerning the first order operators
(with constant and real coefficients), that certain conditions are needed
(necessary or sufficient) on $\varphi$ in order that it can
be the exponent of a weight in a Carleman estimate of type
\eqref{stime-Carlm-29-1}. The other purpose is to show,
again in the case of first order operators, how to apply the
Carleman estimates to prove the uniqueness of the
Cauchy problem. 

Let $I=(-1,1)$ and $\varphi\in C^2\left(\overline{I}\right)$;
let us begin by considering the following elementary Carleman estimate

\begin{equation}  \label{stime-Carlm-ext-1-1}
\tau^{\mu}\int_{I}\left\vert u(x)\right\vert^2 e^{2\tau
\varphi(x)}dx\leq C \int_{I}\left\vert u'(x)\right\vert^2e^{2\tau
\varphi(x)}dx,
\end{equation}
for every $u\in C_0^{\infty}(I)$ and for every $\tau\geq \tau_0$.
We will show how to derive an estimate like
\eqref{stime-Carlm-ext-1-1} and what conditions on $\mu$ and $\varphi$ are necessary for them to hold.

\bigskip

Let us start by the following
\begin{prop}\label{stime-Carlm-ext-1-3}
Let

\begin{equation}  \label{stime-Carlm-ext-1-A}
A=\left\{x\in \overline{I}: \mbox{ } \varphi'(x)= 0 \right\}.
\end{equation}
If
\begin{equation}  \label{stime-Carlm-ext-1-2}
x\in A\mbox{ }\Rightarrow \mbox{ } \varphi''(x)>0,
\end{equation}
then there exist $\tau_0$ and $C$ such that
\begin{equation}  \label{stime-Carlm-ext-1-3-bis}
\tau\int_{I}\left\vert u(x)\right\vert^2 e^{2\tau \varphi(x)}dx\leq
C \int_{I}\left\vert u'(x)\right\vert^2e^{2\tau \varphi(x)}dx,
\end{equation}
for every $u\in C_0^{\infty}(I)$ and for every $\tau\geq \tau_0$.

\end{prop}

\medskip

\textbf{Remark.} Of course if $A=\emptyset$ 
\eqref{stime-Carlm-ext-1-A} is trivially satisfied. $\blacklozenge$
\medskip

\textbf{Proof.} Set

\begin{equation}\label{stime-Carlm-ext-000-2}
v=e^{\tau \varphi}u,
\end{equation}
we have 

\begin{equation}\label{stime-Carlm-ext-00n-2}
e^{\tau \varphi}u'=e^{\tau \varphi}\left(e^{-\tau
\varphi}v\right)'=v'-\tau\varphi'v.\end{equation} Hence estimate 
\eqref{stime-Carlm-ext-1-3} is equivalent to (we omit for
brevity the integration set)

\begin{equation*} \frac{\tau}{C}\int\left\vert v\right\vert^2dx \leq
\int\left\vert v'-\tau\varphi'v\right\vert^2dx
\end{equation*}
for every $v\in C_0^{\infty}(I)$ and for every $\tau\geq \tau_0$.

Now estimate from below the integral  on the RHS. Spreading the square and
integrating by parts we have

\begin{equation}\label{stime-Carlm-ext-0-2}
\begin{aligned}
\int\left\vert v'-\tau\varphi'v\right\vert^2dx&=\int\left(\left\vert
v'\right\vert^2-2\tau\varphi'\Re
\left(v'\overline{v}\right)+\tau^2\varphi'^2\left\vert
v\right\vert^2\right)dx=\\&= \int\left(\left\vert
v'\right\vert^2-\tau\varphi' \left(\left\vert
v\right\vert^2\right)'+\tau^2\varphi'^2\left\vert
v\right\vert^2\right)dx=\\&= \int\left(\left\vert
v'\right\vert^2+\tau\left(\varphi''+\tau\varphi'^2\right)\left\vert
v\right\vert^2\right)dx.
\end{aligned}
\end{equation}
Now, as \eqref{stime-Carlm-ext-1-2} holds, by applying Lemma \ref{106-CauNire} (with $f=\varphi'$ and $g=\varphi''$) we have that there exists $\tau_0$ such that

$$\tau_0\varphi'^2(x)+\varphi''(x)>0,\quad \forall x\in \overline{I}$$ and,
setting 

$$C^{-1}=\min_{\overline{I}}\left(\tau_0\varphi'^2+\varphi''\right)>0,$$
we have, by \eqref{stime-Carlm-ext-0-2}

\begin{equation}\label{stime-Carlm-ext-bis-2}
\begin{aligned}
\int\left\vert v'-\tau\varphi'v\right\vert^2dx&\geq
\int\tau\left(\varphi''+\tau\varphi'^2\right)\left\vert
v\right\vert^2dx\geq\\& \geq C^{-1}\tau\int\left\vert
v\right\vert^2dx.
\end{aligned}
\end{equation}
From which, taking into account \eqref{stime-Carlm-ext-000-2} and
\eqref{stime-Carlm-ext-00n-2}, 
\eqref{stime-Carlm-ext-1-3-bis} follows.$\blacksquare$

\bigskip

\textbf{Remarks.}

\medskip

\noindent\textbf{1.} If $A=\emptyset$ then a
stronger version of \eqref{stime-Carlm-ext-1-3} holds true. More precisely, we have 

\begin{equation}  \label{stime-Carlm-ext-3-0}
\tau^2\int_{I}\left\vert u(x)\right\vert^2 e^{2\tau
\varphi(x)}dx\leq C \int_{I}\left\vert u'(x)\right\vert^2e^{2\tau
\varphi(x)}dx,
\end{equation}
for every $u\in C_0^{\infty}(I)$ and for every $\tau\geq \tau_0$. In
this case, indeed, we have

$$m:=\min_{\overline{I}}\left\vert \varphi'\right\vert>0$$
as $A=\emptyset$. On the other hand, setting
$$m_1=\left\Vert \varphi''\right\Vert_{L^{\infty}(I)},$$ and taking into account  \eqref{stime-Carlm-ext-0-2}, we get

\begin{equation}\label{stime-Carlm-ext-3-1}
\begin{aligned}
\int\left\vert v'-\tau\varphi'v\right\vert^2dx&= \int\left(\left\vert
v'\right\vert^2+\tau^2\left(\tau^{-1}\varphi''+\varphi'^2\right)\left\vert
v\right\vert^2\right)dx\geq\\&\geq \int\left(\left\vert
v'\right\vert^2+\tau^2\left(-\tau^{-1}m_1+m^2\right)\left\vert
v\right\vert^2\right)dx\geq \\&\geq\frac{\tau^2m^2}{2}\int \left\vert
v\right\vert^2dx,
\end{aligned}
\end{equation}
for every $v\in C_0^{\infty}(I)$ and for every $\tau\geq \tau_0$, where
$\tau_0=\frac{2m_1^2}{m^2}$ e $C=\frac{2}{m^2}$. From which  
\eqref{stime-Carlm-ext-3-0} follows.

\medskip

\noindent\textbf{2.} On the other hand, it is also evident that estimate
\eqref{stime-Carlm-ext-1-1} cannot be true for $\mu>2$. As a matter of fact,
we have

\begin{equation}\label{stime-Carlm-ext-bis-3-1}
\begin{aligned}
\int\left\vert v'-\tau\varphi'v\right\vert^2dx\leq
2\int\left(\left\vert v'\right\vert^2+\tau^2\left\Vert
\varphi'\right\Vert_{L^{\infty}(I)}^2\left\vert
v\right\vert^2\right)dx
\end{aligned}
\end{equation}
and so, if \eqref{stime-Carlm-ext-1-1}, and $\mu>2$, would imply

\begin{equation*}
	\begin{aligned}
		&\int\left\vert
		v\right\vert^2dx\leq \\&\leq 2\tau^{2-\mu}\int\left(\tau^{-2}\left\vert
		v'\right\vert^2+\left\Vert
		\varphi'\right\Vert_{L^{\infty}(I)}^2\left\vert
		v\right\vert^2\right)dx\rightarrow 0,\quad\mbox{as }
		\tau\rightarrow+\infty	
	\end{aligned}
\end{equation*}

$$$$ which is evidently absurd since
$v$ is an arbitrary function of $C_0^{\infty}(I)$. $\blacklozenge$

\bigskip

We establish some necessary conditions for estimate
\eqref{stime-Carlm-ext-1-1}.

\begin{prop}\label{stime-Carlm-extn-1-prop-4}
Let
\begin{equation}  \label{stime-Carlm-extn-0-5}
A_0=\left\{x\in I: \mbox{ } \varphi'(x)=0\right\}.
\end{equation}
If there exist $C$ and $\tau_0$ such that
\begin{equation}  \label{stime-Carlm-extn-1-4}
\int_{I}\left\vert u(x)\right\vert^2 e^{2\tau \varphi(x)}dx\leq C
\int_{I}\left\vert u'(x)\right\vert^2e^{2\tau \varphi(x)}dx,
\end{equation}
for every $u\in C_0^{\infty}(I)$ and for every $\tau\geq \tau_0$, then
\begin{equation}  \label{stime-Carlm-extn-1-5}
x\in A_0\mbox{ }\Rightarrow \mbox{ } \varphi''(x)\geq0.
\end{equation}
\end{prop}

\textbf{Proof.} We first notice that, by
density, \eqref{stime-Carlm-extn-1-4} is satisfied for each
$u\in H_0^{1}(I)$.

In order to prove the Proposition we argue by contradiction and we assume that \eqref{stime-Carlm-extn-1-5}
does not hold. Let $x_0\in I$ satisfy $\varphi'(x_0)=0$ and
$\varphi''(x_0)<0$. For the purpose of simplifying the notations, since
$x_0$ is an interior point of $I$, we assume that $x_0=0$. Hence
we have

\begin{equation}  \label{stime-Carlm-extn-0-6}
\varphi'(0)=0\quad\mbox{and}\quad \varphi''(0)<0.
\end{equation}
Since \eqref{stime-Carlm-extn-1-4} is trivially
equivalent to 
\begin{equation}  \label{stime-Carlm-extn-2-5}
\int_{I}\left\vert u(x)\right\vert^2 e^{2\tau
(\varphi(x)-\varphi(0))}dx\leq C \int_{I}\left\vert
u'(x)\right\vert^2e^{2\tau (\varphi(x)-\varphi(0))}dx,
\end{equation}
for every $u\in H^1_0(I)$ and for every $\tau\geq \tau_0$, we may assume
$$\varphi(0)=0.$$
Set
$$a=-\varphi''(0)>0$$
and let $\psi\in H^1_0(I)$ be a function that we will choose
later.

If $0<\varepsilon\leq \tau_0^{-1/2}$, we have

$$x\rightarrow\psi\left(\varepsilon \sqrt{\tau}x\right)\in H^1_0(I),\quad \forall \tau\geq
\varepsilon^{-2}.$$ Now, introducing the following notation in \eqref{stime-Carlm-extn-2-5}

$$u(x)=\psi\left(\varepsilon \sqrt{\tau}x\right),$$
we have

\begin{equation*}
\int_{I}\left\vert \psi\left(\varepsilon
\sqrt{\tau}x\right)\right\vert^2 e^{2\tau \varphi(x)}dx\leq
C\varepsilon^2\tau \int_{I}\left\vert \psi'\left(\varepsilon
\sqrt{\tau}x\right)\right\vert^2e^{2\tau \varphi(x)}dx.
\end{equation*}
By performing the change of variables $t=\sqrt{\tau}x$, we get

\begin{equation}  \label{stime-Carlm-extn-1-6}
\int_{-\sqrt{\tau}}^{\sqrt{\tau}}\left\vert \psi\left(\varepsilon
t\right)\right\vert^2 e^{2\tau \varphi(t/\sqrt{\tau})}dt\leq
C\varepsilon^2\tau \int_{-\sqrt{\tau}}^{\sqrt{\tau}}\left\vert
\psi'\left(\varepsilon t\right)\right\vert^2e^{2\tau
\varphi(t/\sqrt{\tau})}dt.
\end{equation}
We notice that the Taylor formula gives
$$\varphi(x)=-\frac{a}{2}x^2+\frac{x^2}{2}\omega(x),\quad\forall x\in
I,$$ where

\begin{equation}  \label{lim-omega}
\lim_{x\rightarrow 0}\omega(x)=0.\end{equation} 
 Now we choose 

$$\tau=\varepsilon^{-2},$$
and by \eqref{stime-Carlm-extn-1-6} we get

\begin{equation}  \label{stime-Carlm-extn-2-6}
\int_{-1/\varepsilon}^{1/\varepsilon}\left\vert
\psi\left(\varepsilon t\right)\right\vert^2
e^{-t^2\left(a-\omega(\varepsilon t)\right)}dt\leq
C\int_{-1/\varepsilon}^{1/\varepsilon}\left\vert
\psi'\left(\varepsilon
t\right)\right\vert^2e^{-t^2\left(a-\omega(\varepsilon t)\right)}dt.
\end{equation}
Let us choose $\psi$ such that

\begin{equation}\label{Carlm-psi}
\psi(x)=
\begin{cases}
2(x+1),\quad\mbox{ for } x\in \left[-1,-\frac{1}{2}\right), \\
\\
1,\quad\quad\quad\mbox{ for } x\in \left[-\frac{1}{2},\frac{1}{2}\right),\\
\\
2(-x+1),\mbox{ for } x\in \left[\frac{1}{2},1\right].
\end{cases}
\end{equation}
By \eqref{stime-Carlm-extn-2-6} and \eqref{Carlm-psi} we have

\begin{equation}\label{stime-Carlm-extn-3-6}
\begin{aligned}
\int_{-1/2\varepsilon}^{1/2\varepsilon}
e^{-t^2\left(a-\omega(\varepsilon t\right)}dt&\leq
\int_{-1/\varepsilon}^{1/\varepsilon}\left\vert
\psi\left(\varepsilon t\right)\right\vert^2
e^{-t^2\left(a-\omega(\varepsilon t\right)}dt\leq\\&\leq
C\int_{-1/\varepsilon}^{1/\varepsilon}\left\vert
\psi'\left(\varepsilon
t\right)\right\vert^2e^{-t^2\left(a-\omega(\varepsilon
t\right)}dt=\\&
=8C\int^{1/\varepsilon}_{1/2\varepsilon}e^{-t^2\left(a-\omega(\varepsilon
t\right)}dt.
\end{aligned}
\end{equation}
Passing to the limit as $\varepsilon\rightarrow 0$ and recalling 
\eqref{lim-omega}, we get (by the Dominated Convergence Theorem)

\begin{equation*}
0<\int_{-\infty}^{+\infty} e^{-at^2}dt=\lim_{\varepsilon\rightarrow
0}\int_{-1/2\varepsilon}^{1/2\varepsilon}
e^{-t^2\left(a-\omega(\varepsilon t\right)}dt\leq
\lim_{\varepsilon\rightarrow
0}8C\int^{1/\varepsilon}_{1/2\varepsilon}e^{-t^2\left(a-\omega(\varepsilon
t\right)}dt=0
\end{equation*}
Which is, evidently, absurd. $\blacksquare$

\bigskip

\begin{prop}\label{Carlm-cn-prop-0}
Let $A_0$ be as in Proposition \ref{stime-Carlm-extn-0-5}. If
there exist  $C_0$ e $\tau_0$ such that
\begin{equation}  \label{Carlm-cn-1}
\tau\int_{I}\left\vert u(x)\right\vert^2 e^{2\tau \varphi(x)}dx\leq
C \int_{I}\left\vert u'(x)\right\vert^2e^{2\tau \varphi(x)}dx,
\end{equation}
for every $u\in C_0^{\infty}(I)$ and for every $\tau\geq \tau_0$, then
\begin{equation}  \label{Carlm-cn-2}
x\in A_0\mbox{ }\Rightarrow \mbox{ } \varphi''(x) \geq
\frac{1}{2C}>0.
\end{equation}
\end{prop}

\textbf{Proof.} As already noticed above (proof
of Proposition \ref{stime-Carlm-ext-1-3}), estimate
\eqref{Carlm-cn-1} is equivalent to
\begin{equation}\label{Carlm-cn-2bis}
 \frac{\tau}{C}\int_{I}\left\vert v\right\vert^2dx \leq
\int_{I}\left\vert v'-\tau\varphi'v\right\vert^2dx,
\end{equation}
for every $v\in H_0^{1}(I)$ and for every $\tau\geq \tau_0$. Let us notice that

\begin{equation}\label{Carlm-cn-3}
\begin{aligned}
\int\left\vert v'-\tau\varphi'v\right\vert^2dx&=\int_{I}\left\vert
v'+\tau\varphi'v\right\vert^2dx-4\tau\int_{I}\varphi'\Re(\overline{v}v')dx=\\&
=\int_{I}\left\vert
v'+\tau\varphi'v\right\vert^2dx+2\tau\int_{I}\varphi''\left\vert
v\right\vert^2dx.
\end{aligned}
\end{equation}

We suppose, as in the proof of Proposition
\ref{stime-Carlm-extn-1-prop-4}, that $\varphi'(0)=0$ and we want to prove
that $\varphi''(0) \geq \frac{1}{2C}$. We may also let us assume here
that $\varphi(0)=0$.

Let

$$v(x)=e^{-\tau\varphi(x)}\psi\left(\varepsilon
\sqrt{\tau}x\right),$$ with $\psi\in H_0^{1}(I)$ to be chosen later and with
$\tau\geq \varepsilon^{-2}$, $0<\varepsilon<\tau_0^{-1/2}$. By
\eqref{Carlm-cn-2} and \eqref{Carlm-cn-2bis} we have

\begin{equation}\label{Carlm-cn-4}
\begin{aligned}
 &\frac{1}{C}\int_{I}\left\vert \psi\left(\varepsilon
\sqrt{\tau}x\right)\right\vert^2e^{-2\tau\varphi}dx \leq\\& \leq
2\int_{I}\varphi''(x)\left\vert \psi\left(\varepsilon
\sqrt{\tau}x\right)\right\vert^2e^{-2\tau\varphi}dx+\varepsilon^2\int_{I}\left\vert
\psi'\left(\varepsilon
\sqrt{\tau}x\right)\right\vert^2e^{-2\tau\varphi}dx.
\end{aligned}
\end{equation}
By Proposition
\ref{stime-Carlm-extn-1-prop-4}, if we set
$$\alpha=\varphi''(0),$$
we have
\begin{equation}\label{Carlm-cn-5}
 \alpha\geq 0.
\end{equation}

Now, by the Taylor formula and by performing the change of
variables $t=\sqrt{\tau}x$ we obtain, (we argue as in the proof of Proposition
\ref{stime-Carlm-extn-1-prop-4}),

\begin{equation}\label{Carlm-cn-6}
	\begin{aligned}
		&\frac{1}{C}\int_{-\sqrt{\tau}}^{\sqrt{\tau}}\left\vert \psi\left(\varepsilon
t\right)\right\vert^2
e^{-t^2\left(\alpha+\omega(t/\sqrt{\tau})\right)}dt\leq \\&
\leq
2\int_{-\sqrt{\tau}}^{\sqrt{\tau}}\left(\alpha+\omega_1(t/\sqrt{\tau})\right)\left\vert
\psi\left(\varepsilon
t\right)\right\vert^2e^{-t^2\left(\alpha+\omega(t/\sqrt{\tau})\right)}dt+\\&+\varepsilon^2\int_{-\sqrt{\tau}}^{\sqrt{\tau}}\left\vert
\psi'\left(\varepsilon
t\right)\right\vert^2e^{-t^2\left(\alpha+\omega(t/\sqrt{\tau})\right)}dt,
\end{aligned}
\end{equation}
where $\omega(x)$ and $\omega_1(x)$ go to $0$ as $x$ goes to $0$. Passing to the limit in \eqref{Carlm-cn-6} as
$\tau\rightarrow +\infty$, we have

\begin{equation}\label{Carlm-cn-7}
\begin{aligned}
 &\frac{1}{C}\int_{-\infty}^{+\infty}\left\vert \psi\left(\varepsilon
t\right)\right\vert^2 e^{-\alpha t^2}dt \leq\\&\leq
2\int_{-\infty}^{+\infty}\alpha\left\vert\psi\left(\varepsilon
t\right)\right\vert^2e^{-\alpha
t^2}dt+\varepsilon^2\int_{-\infty}^{\infty}\left\vert
\psi'\left(\varepsilon t\right)\right\vert^2e^{-\alpha t^2}dt.
\end{aligned}
\end{equation}
If it were $\alpha=0$, then \eqref{Carlm-cn-7} would be written

\begin{equation*}
\frac{1}{C}\int_{-\infty}^{+\infty}\left\vert \psi\left(\varepsilon
t\right)\right\vert^2dt\leq
\varepsilon^2\int_{-\infty}^{\infty}\left\vert
\psi'\left(\varepsilon t\right)\right\vert^2dt.
\end{equation*}
The latter, by the change of variable $s=\varepsilon t$,  implies
\begin{equation*}
\frac{1}{C}\int_{I}\left\vert \psi\left( s\right)\right\vert^2ds\leq
\varepsilon^2\int_{I}\left\vert \psi'\left(s\right)\right\vert^2ds
\end{equation*}
which, passing to the limit as
$\varepsilon\rightarrow 0$, leads to an absurd (just choose $\psi$ not identically
null). Therefore necessarily we have
 
$$\alpha>0.$$
At this point, passing to the limit as $\varepsilon\rightarrow 0$
in \eqref{Carlm-cn-7}, we obtain (by the Dominated Convergence Theorem)

\begin{equation*}
 \frac{1}{C}\left\vert \psi\left(
0\right)\right\vert^2\int_{-\infty}^{+\infty} e^{-\alpha t^2}dt\leq
2\left\vert\psi\left(
0\right)\right\vert^2\alpha\int_{-\infty}^{+\infty} e^{-\alpha
t^2}dt.
\end{equation*}
From which we have trivially 
$$2\alpha\geq \frac{1}{C}.$$ $\blacksquare$

\bigskip

\textbf{Remark.} If $1<\mu\leq 2$, in \eqref{stime-Carlm-ext-1-1}, then 
$A_0=\emptyset$  i.e. $\varphi'(x)\neq 0$ for every $x\in (-1,1)$.
As a matter of fact, if $\mu>1$, we would have, for every 
$K>0$
\begin{equation}  \label{stime-Carlm-ext-1-1-bis}
\tau K^{\mu-1}\int_{I}\left\vert u(x)\right\vert^2 e^{2\tau
\varphi(x)}dx\leq C \int_{I}\left\vert u'(x)\right\vert^2e^{2\tau
\varphi(x)}dx,
\end{equation}
for every $u\in C_0^{\infty}(I)$ an for every $\tau\geq \max\{\tau_0,
K\}$. Now, if $A_0\neq\emptyset$, then there exists $x_0\in I$ such that
$\varphi'(x_0)=0$, hence, by \eqref{Carlm-cn-2} we have

$$\varphi''(x_0)\geq\frac{K^{\mu-1}}{2C},\quad\forall K>0$$ from which we have

$$\varphi''(x_0)=+\infty.$$
which contradicts $\varphi\in
C^2\left(\overline{I}\right)$. $\blacklozenge$

\bigskip

\bigskip

We now consider the first-order operator

\begin{equation}  \label{qr1-1}
P_1(\partial)=\sum_{j=1}^na_j\partial_j=a\cdot \nabla.
\end{equation}
where $a=(a_1,\cdots,a_n)\in \mathbb{R}^n\setminus\{0\}$ and $a_j$,
$j=1,\cdots,n$ are constants. Let us propose to transfer
to operator \eqref{qr1-1} what we established above for the
the derivative operator. We will reach a Carleman estimate of the type
\begin{equation}  \label{qr1-carlm-1}
\tau\int\left|u\right|^2e^{2\tau \varphi(x)}dx\leq
C\int\left|P_1(\partial)u\right|^2e^{2\tau \varphi(x)}dx,
\end{equation}
for every $u\in C^{\infty}_0\left(\mathbb{R}^n\right)$ and for every $\tau$
large enough, where $\varphi$ is a function which belongs to 
$C^2\left(\mathbb{R}^n\right)$ on which we will make further assumptions later.

Let us suppose, for instance, that
\begin{equation}  \label{qr1-1-0}
a_n\neq 0.
\end{equation}
Let, for $y\in\mathbb{R}^{n-1}$, $x=X(t,y)$, the equations of characteristic lines satisfying

\begin{equation}\label{qr1-2}
\begin{cases}
\partial_tX(t,y)=a\cdot X(t,y),\\
\\
X(0,y)=(y,0).
\end{cases}
\end{equation}
We have

\begin{equation}\label{qr1-3}
\begin{cases}
X_1(t,y)=a_1t+y_1,\\
\\
\cdots,\\
\\
X_{n-1}(t,y)=a_{n-1}t+y_{n-1},\\
\\
X_{n}(t,y)=a_{n}t.
\end{cases}
\end{equation}
$X$ is a linear and bijective transformation from
$\mathbb{R}^{n}$ in itself since the absolute value of the determinant of the matrix
associated to $X$ is equal to $|a_n|$ and by \eqref{qr1-1-0} we have
$a_n\neq 0$. Moreover, see Section \ref{metodo-caratt}, setting
\begin{equation}\label{qr1-4}
z(t,y)= u((X(t,y)),
\end{equation}

\begin{equation}\label{qr1-5}
\left(a\cdot \nabla u\right)((X(t,y))=\partial_t z(t,y)
\end{equation}
and

\begin{equation*}
\widetilde{\varphi}(t,y)=\varphi(X(t,y)),
\end{equation*}
estimate  \eqref{qr1-carlm-1} is equivalent to the estimate

\begin{equation}  \label{qr1-carlm-2}
\tau\int\left|z\right|^2e^{2\tau \widetilde{\varphi}(t,y)}dtdy\leq
C\int\left|\partial_t z\right|^2e^{2\tau
\widetilde{\varphi}(t,y)}dtdy,
\end{equation}
for every $z\in C^{\infty}_0\left(\mathbb{R}^n\right)$ and for every $\tau$
large enough. Of course, if we are interested in
estimate \eqref{qr1-carlm-1} for $u$ supported in a bounded open
$\Omega$ then estimate \eqref{qr1-carlm-2} will be established for
$z$ supported in a bounded open set. More precisely, set $\widetilde{\Omega}=X^{-1}(\Omega)$, Proposition \ref{stime-Carlm-ext-1-3}  yields what follows:

If

\begin{equation}\label{qr1-6}
\partial_t^2\widetilde{\varphi}(t,y)>0,\quad\mbox{for every } (t,y) \mbox{ such that
} \partial_t\widetilde{\varphi}(t,y)=0,
\end{equation}
then estimate \eqref{qr1-carlm-2} and (consequently) estimate
\eqref{qr1-carlm-1} holds true.

Now we have,

$$\partial_t\widetilde{\varphi}(t,y)=\left(\nabla
\varphi\right)\left(X(t,y)\right)\cdot a
\left(X(t,y)\right)=\sum_{j=1}^n \left(\partial_j \varphi a_j\right)
\left(X(t,y)\right),$$
and

\begin{equation}\label{qr1-6-0}
\begin{aligned}
 \partial_t^2\widetilde{\varphi}(t,y)&=\partial_t\left(\sum_{j=1}^n \left(\partial_j
\varphi\right)\left(X(t,y)\right)a_j\left(X(t,y)\right)\right)=
\\&=\sum_{j,k=1}^n\left(\partial^2_{jk}
\varphi\right)\left(X(t,y)\right)\left(\partial_t
X(t,y)\right)a_j\left(X(t,y)\right)+\\&+
\sum_{j,k=1}^n\left(\partial_j
\varphi\right)\left(X(t,y)\right)\left(\partial_{x_k}a_i\right)(X(t,y))\partial_tX(t,y)=\\&=
\sum_{j,k=1}^n\left(\partial^2_{jk} \varphi a_ja_k\right)(X(t,y)),
\end{aligned}
\end{equation}
(in the second to last step we used that $\partial_{x_k}a_i=0$, as $a$ is a constant vector). Therefore, with respect to the variables $x_1,\cdots, x_n$,
condition \eqref{qr1-6} can be written

\begin{equation}\label{qr1-7}
a\cdot\nabla\varphi(x)=0\quad\Rightarrow
\sum_{j,k=1}^n\partial^2_{jk} \varphi (x) a_ja_k>0.
\end{equation}
So if \eqref{qr1-7} holds, then
for every bounded open set $\Omega$, the following Carleman estimate holds
\begin{equation}  \label{qr1-carlm-3}
\tau\int\left|u\right|^2e^{2\tau \varphi(x)}dx\leq
C\int\left|P_1(\partial)u\right|^2e^{2\tau \varphi(x)}dx,
\end{equation}
for every $u\in C^{\infty}_0\left(\Omega\right)$ and for every $\tau$
large enough.

Before applying this estimate to study of the uniqueness of the
Cauchy problem we provide a geometric interpretation of
condition \eqref{qr1-7}. Let $x_0\in \mathbb{R}^n$ and let us suppose that

\begin{equation}\label{qr1-8}
\nabla\varphi(x_0)\neq 0.
\end{equation}
Now, the assertion  

$$a\cdot\nabla\varphi(x_0)=0,$$ is equivalento to the assertion that the surface
$\{\varphi(x)=\varphi(x_0)\}$ is a characteristic surface for the operator
$P_1(\partial)$ in $x_0$. Regarding the interpretation of the
term
$$\sum_{jk=1}^n\partial^2_{jk} \varphi (x_0) a_ja_k,$$ we are helped by the calculations performed in \eqref{qr1-6-0}. Actually, denoting  by
$x=\gamma(t)$ the parametric equation of the characteristic line
passing through $x_0$, for instance, set $\gamma(0)=x_0$, then we have

$$\frac{d^2\varphi(\gamma(t))}{dt^2}_{|t=0}=\sum_{jk=1}^n\partial^2_{jk} \varphi (x_0)
a_ja_k.$$ Therefore, condition \eqref{qr1-7} states that if
$x_0$ is a characteristic point of the level surface
with respect to the operator $P_1$, then there exists a neighborhood of $0$, $J$,
such that
$$\varphi (\gamma (t))>\varphi (x_0),\quad \forall t\in J\setminus\{0\}$$
that is, the characteristic line $x=\gamma(t)$ remains locally
confined to the region $\{\varphi(x)>\varphi(x_0)\}$ or, in other
words, does not cross the level surface
$$\{\varphi(x)=\varphi(x_0)\}$$ in $x_0$.

\bigskip

Given $\Omega$, an open set of $\mathbb{R}^n$, $x_0\in \Omega$ and $\psi\in
C^1\left(\overline{\Omega}\right)$ a real--valued function such that

\begin{equation}\label{qr1-9}
\nabla\psi(x)\neq 0,\quad\forall x\in \Gamma,
\end{equation}
where
\begin{equation}\label{qr1-10}
\Gamma=\left\{x\in \Omega: \quad \psi(x)=\psi(x_0)  \right\}.
\end{equation} We will say that a real--valued function $\varphi\in
C^0\left(\overline{\Omega}\right)$ with , enjoys the
\textbf{property of convexification with respect to $\Gamma$ in
$x_0$} \index{property of convexification}if $\varphi(x_0)=\psi(x_0)$ and there exists $r>0$ such that (Figure 13.1)

\begin{equation}\label{qr1-11}
	\begin{aligned}
		&\left\{x\in B_r(x_0): \varphi(x)\geq\varphi(x_0)
\right\}\setminus\{x_0\}\subset \\&\subset \left\{x\in B_r(x_0):
\psi(x)>\psi(x_0) \right\}.
\end{aligned}
\end{equation}

\smallskip

\begin{figure}\label{conv}
	\centering
	\includegraphics[trim={0 0 0 0},clip, width=11cm]{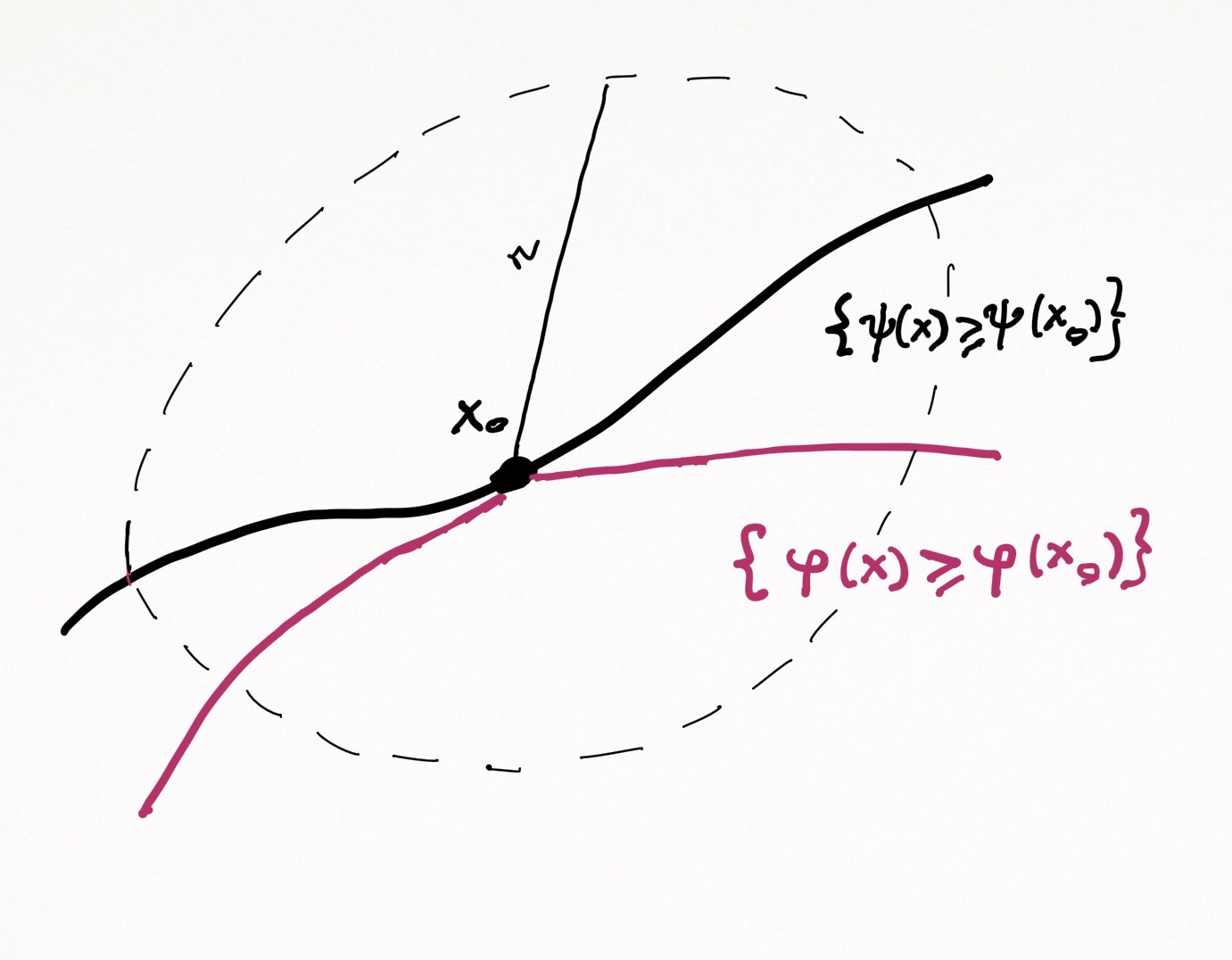}
	\caption{}
\end{figure}

In the next Proposition, we set for a function $f\in
C^2(\Omega)$

$$Q_{f}(x)=\sum_{j,k=1}^n\partial^2_{jk} f(x) a_ja_k$$
The following holds true
\begin{prop}\label{qr-unic-primo ordine}
Let $\Omega$ be an open set of $\mathbb{R}^n$ and let $a\in
\mathbb{R}^n\setminus\{0\}$, $c\in L^{\infty}(\Omega)$ (with values in
$\mathbb{C}$), $x_0\in \Omega$ and $\psi\in
C^2\left(\overline{\Omega}\right)$ a real--valued function satisfying
the following conditions

\begin{equation}\label{qr1-12}
\nabla\psi(x_0)\neq 0
\end{equation}
and let us suppose that
\begin{equation}\label{qr1-13}
a\cdot\nabla\psi(x_0)=0\quad\Rightarrow Q_{\psi}(x_0) >0.
\end{equation}
Let $U\in H^1(\Omega)$ satisfy
\begin{equation}\label{qr1-14}
\begin{cases}
a\cdot\nabla U+c(x)U=0,\quad\quad\mbox{in } \Omega,\\
\\
U(x)=0 \quad\mbox{in } \left\{x\in \Omega: \quad \psi(x)>\psi(x_0)
\right\}.
\end{cases}
\end{equation}
Then there exists a neighborhood $\mathcal{U}_{x_0}$ of $x_0$ such that

\begin{equation}\label{qr1-15}
U=0\quad\mbox{in }\quad \mathcal{U}_{x_0}.
\end{equation}
\end{prop}

\bigskip

\textbf{Proof.} It is not restrictive to assume $x_0=0$,
$$\psi(0)=0$$ and
$$|a|=1.$$
Let

\begin{equation}\label{qr1-16}
\varphi_{\varepsilon}(x)=\psi(x)-\frac{\varepsilon|x|^2}{2},
\end{equation}
where $\varepsilon$ is a positive number to be chosen later.

Now we check that $\varphi_{\varepsilon}$ satisfies \eqref{qr1-7}
in a neighborhood of $0$. We note, to this purpose, that
\eqref{qr1-13} implies that there exists a constant $C_0>0$ such that

\begin{equation*}
M_0:=C_0\left(a\cdot\nabla\psi(0)\right)^2+
 Q_{\psi}(0)>0.
\end{equation*}
We have easily that, for any $\varepsilon\leq
\varepsilon_0=\frac{M_0}{2}$,

\begin{equation}\label{qr1-17}
C_0\left(a\cdot\nabla\varphi_{\varepsilon}(0)\right)^2+
Q_{\varphi_{\varepsilon}}(0) =M_0-\varepsilon\geq \frac{M_0}{2}.
\end{equation}
Let us choose $\varepsilon=\varepsilon_0$ and we omit from now on the
subscript of $\varphi$. By \eqref{qr1-17}, since $\varphi\in
C^2(\Omega)$, there exists $R>0$ such that

\begin{equation}\label{qr1-18}
C_0\left(a\cdot\nabla\varphi(x)\right)^2+ Q_{\varphi}(x)\geq
\frac{M_0}{4}>0,\quad\quad\forall x\in B_{2R}.
\end{equation}
Therefore, \eqref{qr1-7} applies and, consequently, setting

$$Pu=a\cdot\nabla u+c(x)u,$$
we get
\begin{equation}  \label{qr1-carlm-4}
\tau\int\left|u\right|^2e^{2\tau \varphi(x)}dx\leq
C\int\left|Pu\right|^2e^{2\tau \varphi(x)}dx,
\end{equation}
for every $u\in C^{\infty}_0\left(B_{2R}\right)$ and for every $\tau\geq
\tau_0$. By density \eqref{qr1-carlm-4} holds for every $u\in
H^1_0(B_{2R})$.

Let $\eta\in C^{\infty}_0\left(B_{2R}\right)$ be a function such that

$$0\leq \eta(x) \leq 1,\quad\forall x\in B_{R}; \quad \eta(x)=1,\quad \forall x\in B_{R/2}$$
and
$$\mbox{supp } \eta=\overline{B_R}.$$
Let us denote

$$C_1=\left\Vert \nabla
\eta\right\Vert_{L^{\infty}\left(B_{R}\right)}.$$ Now we apply
\eqref{qr1-carlm-4} to $$u=\eta U,$$ since

$$P(\eta U)=\eta\left(a\cdot\nabla U+c(x)U\right)+(a\cdot\nabla
\eta)U,$$ we have
$$|P(\eta U)|\leq C_1\chi_{B_{R}\setminus B_{R/2}} |U|,$$
where $\chi_{B_{R}\setminus B_{R/2}}$ is the characteristic function of $B_{R}\setminus B_{R/2}$. Hence by
\eqref{qr1-carlm-4} we obtain

\begin{equation}  \label{qr1-carlm-5}
\tau\int_{B_R}\left|U\eta\right|^2e^{2\tau \varphi(x)}dx\leq
CC_1^2\int_{B_{R}\setminus B_{R/2}}\left|U\right|^2e^{2\tau
\varphi(x)}dx,
\end{equation}
for every $\tau\geq \tau_0$.

Now, let
$$G=\left(\overline{B_{R}}\setminus B_{R/2}\right)\cap\mbox{supp } U,$$

$$M_1=\max_{G}\varphi$$
and let us prove that
\begin{equation}\label{qr1-19}
M_1<0.
\end{equation}
For this purpose we first observe that $\varphi$ is a
convexification of

$$\Gamma:=\left\{x\in B_R: \quad \psi(x)=0  \right\}.$$ As a matter of fact, let us note that
if

\begin{equation}\label{qr1-20}
x\in\left\{y\in B_R: \varphi(y)\geq 0 \right\}\setminus\{0\}
\end{equation}
then $$\psi(x)\geq
\frac{\varepsilon_0|x|^2}{2}>0,\quad\mbox{for}\quad x\neq 0.$$

Now, arguing by contradiction, let us suppose that \eqref{qr1-19} is false, i.e.
let us suppose that
$$M_1\geq 0.$$ Let $\overline{x}\in G$ satisfy
$\varphi\left(\overline{x}\right)=M_1$. Since $G\subset
\overline{B_R}\setminus B_{R/2}$ we have $\overline{x}\neq
0$. Hence, \eqref{qr1-20} implies
$\psi\left(\overline{x}\right)>0$ from which we get that there exists $\delta>0$
such that
$$\psi(x)>0, \quad \forall x\in B_{\delta}\left(\overline{x}\right)$$
and, recalling that

$$\left\{x\in \Omega: \quad \psi(x)>0  \right\}\subset \left\{x\in \Omega: \quad U(x)=0
\right\},$$ we have $B_{\delta}\left(\overline{x}\right)\subset
\left\{x\in \Omega: \quad U(x)=0 \right\}$. Therefore
$$\overline{x}\notin \mbox{supp } U,$$
that contradicts $\overline{x}\in G\subset$ supp
$U$. Hence \eqref{qr1-19} holds true.

\smallskip

Now, by $\varphi(0)=0$ and by \eqref{qr1-19} we have trivially that
$0$ is an interior point of 

$$\left\{x\in \overline{B_R}:\quad
\varphi(x)>M_1\right\},$$ therefore there exists $r$, $0<r\leq
\frac{R}{2}$, such that
\begin{equation}\label{qr1-21ex}\overline{B_r}\subset
 \left\{x\in \overline{B_R}:\quad
\varphi(x)>M_1\right\}.
\end{equation} Let now $$M_2=\min_{\overline{B_r}}\varphi,$$ by
\eqref{qr1-21ex} we have
$$M_2>M_1.$$

Now, let us come back to  \eqref{qr1-carlm-5}. We have, trivially,

\begin{equation}\label{qr1-21}
\int_{B_{R}\setminus B_{R/2}}\left|U\right|^2e^{2\tau
\varphi(x)}dx\leq e^{2\tau M_1}\int_{G}\left|U\right|^2dx
\end{equation}
and
\begin{equation}\label{qr1-22}
\begin{aligned}
 \int_{B_R}\left|U\eta\right|^2e^{2\tau \varphi(x)}dx&\geq \int_{B_R\cap \left\{
\varphi>M_1\right\}}\left|U\eta\right|^2e^{2\tau
\varphi(x)}dx\geq\\& \geq\int_{B_r} \left|U\right|^2e^{2\tau
\varphi(x)}dx\geq e^{2\tau M_2}\int_{B_r}\left|U\right|^2dx.
\end{aligned}
\end{equation}
By \eqref{qr1-carlm-5}, \eqref{qr1-21} and \eqref{qr1-22} we have

\begin{equation*}
\int_{B_r}\left|U\right|^2dx\leq CC_1^2 e^{-2\tau
(M_2-M_1)}\int_{G}\left|U\right|^2dx
\end{equation*}
for every $\tau\geq \tau_0$, from which, passing to the limit as $\tau$ goes
to $+\infty$, we obtain $U=0$ in $B_r$. $\blacksquare$

\bigskip

\textbf{Remarks.}

\noindent\textbf{1.} The geometric part of the proof of
Proposition \eqref{qr-unic-primo ordine} is to be considered
standard and will occur again even in the case of more general operators than the ones considered so far.  On the contrary, the path that
we followed to arrive to estimate
\eqref{qr1-carlm-3} is not extendable (or, at least, is not
easily extendable) to more general situations. However, the reader
can easily repeat, for the variable coefficient operator

\begin{equation}\label{qr1-23}
P_1(x,\partial)=a(x)\cdot \nabla,
\end{equation}
the calculations we did in the case where $a$ is a constant
vector. Of course, assuming, for instance,
that 

\begin{equation*}
a_n(0)\neq 0,
\end{equation*}
instead of \eqref{qr1-2}, the reader may  consider

\begin{equation*}
\begin{cases}
\partial_tX(t,y)=a(X(t,y)),\\
\\
X(0,y)=(y,0).
\end{cases}
\end{equation*}
obtaining, unlike the case in which $a$ is constant, a local change of coordinates and, instead of \eqref{qr1-7}, it will be found
(compare with \eqref{qr1-6-0}) the following condition in
$\overline{\mathcal{U}_0}$, where $\mathcal{U}_0$ is a neighborhood of
$0$.

\begin{equation}\label{qr1-24}
\begin{aligned}
&a(x)\cdot\nabla \varphi(x)=0\quad\Rightarrow\\& \Rightarrow
\sum_{j,k=1}^n\partial^2_{jk} \varphi (x)
a_j(x)a_k(x)+\sum_{j,k=1}^n\partial_{k} a_j(x)a_k(x)
\partial_{j}\varphi(x)>0.
\end{aligned}
\end{equation}

\bigskip

By the procedure that we have outlined above, we have that if \eqref{qr1-24} holds then

\begin{equation}  \label{qr1-carlm-6}
\tau\int\left|u\right|^2e^{2\tau \varphi(x)}dx\leq
C\int\left|P_1(x,\partial)u\right|^2e^{2\tau \varphi(x)}dx,
\end{equation}
for every $u\in C^{\infty}_0\left(\mathcal{U}_0\right)$ and for every $\tau$
large enough.

\medskip

\noindent \textbf{2.}  We now derive estimate \eqref{qr1-carlm-6} by means of
a procedure based on integrations by parts. This
procedure will be extended in the next Sections to more general operators.

Let

$$P_{1,\tau}(x,\partial)v:=e^{\tau\varphi}P_{1}(x,\partial)\left(e^{-\tau\varphi}v\right)=a\cdot\nabla v-\tau
(a\cdot\nabla \varphi)v,$$ estimate \eqref{qr1-carlm-6} is
equivalent to
\begin{equation}  \label{qr1-carlm-7}
\tau\int\left|v\right|^2dx\leq
C\int\left|P_{1,\tau}(x,\partial)v\right|^2dx,
\end{equation}
for every $v\in C^{\infty}_0\left(\mathcal{U}_0\right)$ and for every $\tau$
large enough. Since the coefficients of
$P_{1}(x,\partial)$ are real-valued, to prove the
\eqref{qr1-carlm-7} it suffices to consider $v$ real-valued.  We have

\begin{equation}\label{qr1-24-bis}
\begin{aligned}
 \int\left|P_{1,\tau}(x,\partial)v\right|^2dx&= \int\left|a\cdot\nabla v-\tau
(a\cdot\nabla \varphi)v\right|^2dx=\\& =\int\left|a\cdot\nabla
v\right|^2dx+ \tau^2\int\left|(a\cdot\nabla
\varphi)v\right|^2dx-\\&- 2\tau\int \left(a\cdot\nabla
\varphi\right)\left(a\cdot\nabla v\right)vdx.
\end{aligned}
\end{equation}

Now, let us consider the third integral on the right hand side in \eqref{qr1-24};
Integrating by parts we have

\begin{equation*}
\begin{aligned}
 - 2\tau\int \left(a\cdot\nabla
\varphi\right)\left(a\cdot\nabla v\right)vdx&=-\tau\int
a\cdot\nabla\left(v^2\right)\left(a\cdot\nabla \varphi\right)dx=
\\&=\tau\int\mbox{div }\left(\left(a\cdot\nabla
\varphi\right)a\right)v^2dx.
\end{aligned}
\end{equation*}
By the the just obtained equality and by \eqref{qr1-24} we have.

\begin{equation}\label{qr1-25}
\begin{aligned}
 \int\left|P_{1,\tau}(x,\partial)v\right|^2dx\geq\int\left(\tau^2\left(a\cdot\nabla \varphi\right)+\tau
 \widetilde{Q}_{\varphi}\right)v^2dx,
\end{aligned}
\end{equation}
where
\begin{equation*}
\begin{aligned}
 \widetilde{Q}_{\varphi}&=\mbox{div }\left(\left(a\cdot\nabla
\varphi\right)a\right)=\\&=\sum_{j,k=1}^n\partial^2_{jk} \varphi
a_ja_k+\sum_{j,k=1}^n\partial_{k} a_ja_k
\partial_{j}\varphi+(a\cdot\nabla \varphi) (\mbox{div }a).
\end{aligned}
\end{equation*}
Now, proceeding as in the proof of
Proposition \ref{stime-Carlm-ext-1-3} we have that if 
\begin{equation}\label{qr1-26}
a(x)\cdot\nabla\varphi(x)=0\quad\Rightarrow
\widetilde{Q}_{\varphi}>0,\quad\mbox{in} \ \  \overline{\mathcal{U}_0}
\end{equation}
then, taking into account \eqref{qr1-25},

\begin{equation*}
\begin{aligned}
 \int\left|P_{1,\tau}(x,\partial)v\right|^2dx\geq\frac{\tau}{C}\int|v|^2dx,
\end{aligned}
\end{equation*}
for every $v\in C^{\infty}_0\left(\mathcal{U}_0\right)$ and for
$\tau\geq \tau_0$ ($\tau_0$ independent of $v$). From the latter
\eqref{qr1-carlm-6} follows. Let us note that \eqref{qr1-24} and
\eqref{qr1-26} are equivalent.

\section{Quadratic differential form and their integration by parts}\label{forme quadratiche}
As we have already seen in the simple examples of the
previous Section, the first steps one makes to prove a
Carleman estimate consists in
setting, for an arbitrary function $u\in C^{\infty}_0(\Omega)$,
$$v=e^{\tau\varphi}u.$$
In this way, denoting, for the sake of brevity, by $P$ the differential operator
 $P_m(x,D)$ (principal part of the operator $P(x,D)$)
we introduce \textbf{the conjugate operator of $P$}\index{conjugate operator} which is defined by

\begin{equation}\label{coniugato}
P_{\tau}v=e^{\tau\varphi}P\left(e^{-\tau\varphi}v\right)=e^{\tau\varphi}P_m(x,D)\left(e^{-\tau\varphi}v\right);
\end{equation}
after that, since
$$\int_{\Omega}\left\vert P_m(x,D)u\right\vert^2e^{2\tau \varphi(x)}dx=\int_{\Omega}\left\vert
P_{\tau}v\right\vert^2dx,$$ we wish to prove a suitable 
 estimate from below of $\int_{\Omega}\left\vert P_{\tau}v\right\vert^2dx$, if $\tau$ is
 \textit{large enough}. Of course, it is precisely
 this estimate from below the most tricky part of the
 proof of a Carleman estimate. In the examples we
 encountered in the previous Section we first spread the square $\left\vert P_{\tau}v\right\vert^2$ and then we integrate it by parts, but it is evident that unless appropriate arrangements are made, this procedure leads to great difficulty for the operators just a little slightly more general than those encountered in the previous Section. These difficulties also arise in the case of the Laplace operator. As a matter of fact, set $$P=\Delta=-\sum_{j=1}^nD^2_j,$$ we have

\begin{equation}\label{coniugato-laplaciano}
P_{\tau}v=e^{\tau\varphi}\Delta\left(e^{-\tau\varphi}v\right)=\Delta
v-\tau\nabla\varphi\cdot\nabla
v-\tau\Delta\varphi+\tau^2\left|\nabla\varphi\right|^2v,
\end{equation}
and, spreading the square we obtain an expression like
$$\sum_{|\alpha|,|\beta|\leq
2}\tau^{\gamma_{\alpha,\beta}}a_{\alpha\beta}(x)\partial^{\alpha}v\partial^{\beta}v.$$
To handle this kind of expressions, in the present Section
we will study the quadratic forms
\begin{equation}\label{fqd-1}
\sum_{\alpha,\beta}a_{\alpha\beta}(x)D^{\alpha}u\overline{D^{\beta}u},
\end{equation}
where the sum is finite, $a_{\alpha\beta}$ are complex--valued functions. We are particularly interested in
the integration by parts of forms \eqref{fqd-1}.

\bigskip

We recall that a \textbf{sesquilinear form}  \index{form:@{form:}!- sesquilinear@{- sesquilinear}}on a complex vector space $V$ is a function $$\Phi:V\times
V\rightarrow V$$ such that $f(\cdot,v)$ is linear for every $v\in V$
and $\Phi(u,\cdot)$ is antilinear for every $u\in V$. We say that a
sesquilinear form on $V$ is  \textbf{hermitian}\index{form:@{form:}!- hermitian@{- hermitian}}, if

\begin{equation}\label{fqd-herm}
\Phi(u,v)=\overline{\Phi(v,u)},\quad\quad\forall u,v\in V.
\end{equation}

In the sequel of this Section we will denote by $\mathcal{V}$ the
space of sesquilinear forms on
$C^{\infty}\left(\mathbb{R}^n,\mathbb{C}\right).$

Let us consider the  sesquilinear forms
\begin{equation}\label{fqd-2}
\Phi_{\alpha\beta}:
C^{\infty}\left(\mathbb{R}^n,\mathbb{C}\right)\times
C^{\infty}\left(\mathbb{R}^n,\mathbb{C}\right)\rightarrow
C^{\infty}\left(\mathbb{R}^n,\mathbb{C}\right),
\end{equation}

\begin{equation}\label{fqd-3}
\Phi_{\alpha\beta}(u,v)=D^{\alpha}u\overline{D^{\beta}v},\quad
\ \ \forall u,v\in C^{\infty}\left(\mathbb{R}^n,\mathbb{C}\right).
\end{equation}

In what follows, for any $\zeta\in \mathbb{C}^n$,
$\zeta=(\zeta_1,\cdots,\zeta_n)$ and $x\in \mathbb{R}^n$,
\\ $x=(x_1,\cdots,x_n)$, we will denote by

$$\zeta\cdot x=\sum_{j=1}^n\zeta_jx_j.$$

\bigskip

The following Proposition holds true

\begin{prop}\label{fqd-4}
The family of sesquilinear forms
$\left\{\Phi_{\alpha\beta}\right\}_{\alpha\beta\in \mathbb{N}^n_0}$
is linearly independent in $\mathcal{V}$.
\end{prop}

\textbf{Proof.} Let us consider a finite linear combination of the forms $\Phi_{\alpha\beta}$, $\alpha,\beta\in \Lambda$, where
$\Lambda$ is a finite subset of 
$\mathbb{N}^n_0\times\mathbb{N}^n_0$ and let us assume that it vanishes identically. We have, for some $c_{\alpha\beta}\in \mathbb{C}$,

\begin{equation}\label{fqdr-1}
\sum_{\alpha,\beta\in \Lambda}c_{\alpha\beta}\Phi_{\alpha\beta}=0.
\end{equation}
Now we prove that

\begin{equation}\label{fqdr-1-0}
c_{\alpha\beta}=0,\quad\quad \forall \alpha, \beta\in \Lambda.
\end{equation}
Let us notice that \eqref{fqdr-1} is equivalent to
\begin{equation}\label{fqdr-5}
\sum_{\alpha,\beta\in
\Lambda}c_{\alpha\beta}D^{\alpha}u\overline{D^{\beta}v}=0,\quad\forall
u,v\in C^{\infty}\left(\mathbb{R}^n,\mathbb{C}\right).
\end{equation}
Now, let $\zeta,\eta\in \mathbb{C}^n$ be arbitrary and put 
$u=e^{i\zeta\cdot x}$, $v=e^{i\eta\cdot x}$ in \eqref{fqdr-5}.
We get

\begin{equation*}
0=\sum_{\alpha,\beta\in
\Lambda}c_{\alpha\beta}D^{\alpha}\left(e^{i\zeta\cdot
x}\right)\overline{D^{\beta}\left(e^{i\eta\cdot x}\right)}=
e^{i\zeta\cdot x}\overline{e^{i\eta\cdot x}}\sum_{\alpha,\beta\in
\Lambda}c_{\alpha\beta}\zeta^{\alpha}\overline{\eta^{\beta}}.
\end{equation*}
Therefore
\begin{equation*}
\sum_{\alpha,\beta\in
\Lambda}c_{\alpha\beta}\zeta^{\alpha}\overline{\eta^{\beta}}=0,\quad\forall
\zeta,\eta\in \mathbb{C}^n.
\end{equation*}
From which we obtain \eqref{fqdr-1-0}. $\blacksquare$

\bigskip

Let us denote by $\mathcal{W}$ the subspace of $\mathcal{V}$
generated by $\left\{\Phi_{\alpha\beta}\right\}_{\alpha\beta\in
\mathbb{N}^n_0}$. We will call \textbf{sesquilinear differential form}\index{form:@{form:}!- sesquilinear differential@{- sesquilinear differential}} any element of the space $\mathcal{W}$. Thus an
arbitrary element of $\mathcal{W}$ is 
\begin{equation}\label{fqd-6}
\Phi(u,v)=\sum_{\alpha,\beta}a_{\alpha\beta}D^{\alpha}u\overline{D^{\beta}v},
\end{equation}
where the sum is finite, and  $a_{\alpha\beta}\in \mathbb{C}$.

\bigskip

The following Proposition holds true.
\begin{prop}\label{fqd-equivalenze}
Let $\Phi\in \mathcal{W}$ be given by \eqref{fqd-6} the following conditions are equivalent.
\begin{subequations}
\label{fqd-equivalenze-1}
\begin{equation}
\label{fqd-equivalenze-1a} \Phi\mbox{ is an hermitian form},
\end{equation}
\begin{equation}
\label{fqd-equivalenze-1b} \Phi(u,u)\in \mathbb{R},\quad\quad\forall
u\in C^{\infty}\left(\mathbb{R}^n,\mathbb{C}\right),
\end{equation}
\begin{equation}
\label{fqd-equivalenze-1c}
a_{\alpha\beta}=\overline{a_{\beta\alpha}},\quad\quad \forall
\alpha,\beta\in \mathbb{N}^n_0,
\end{equation}
\begin{equation}
\label{fqd-equivalenze-1d}
\sum_{\alpha,\beta}a_{\alpha\beta}\zeta^{\alpha}\overline{\zeta^{\beta}}\in
\mathbb{R},\quad\forall \zeta\in \mathbb{C}^n.
\end{equation}
\end{subequations}
\end{prop}

\textbf{Proof}. We follow the pattern
$$\eqref{fqd-equivalenze-1a}\Longleftrightarrow\eqref{fqd-equivalenze-1b}\Longleftrightarrow\eqref{fqd-equivalenze-1c}\Longleftrightarrow\eqref{fqd-equivalenze-1d}.$$

\medskip

\noindent implication \eqref{fqd-equivalenze-1a}$
\Longrightarrow$ \eqref{fqd-equivalenze-1b} is trivial.

\medskip

\noindent Let us prove that  \eqref{fqd-equivalenze-1b}$ \Longrightarrow$
\eqref{fqd-equivalenze-1a}.

\medskip

Let $u,v\in C^{\infty}\left(\mathbb{R}^n,\mathbb{C}\right)$
be arbitraries. By \eqref{fqd-equivalenze-1b} we have
\begin{equation}\label{fqd-7}
\Phi(u+v,u+v)\in \mathbb{R},\quad \Phi(u+iv,u+iv)\in \mathbb{R}.
\end{equation} Setting
$$z=\Phi(u,v),\quad\quad w=\Phi(v,u),$$
we get, by \eqref{fqd-7},
$$z+w=\Phi(u+v,u+v)-\Phi(u,u)-\Phi(v,v)\in
\mathbb{R}$$ and

$$-iz+iw=\Phi(u+iv,u+iv)-\Phi(u,u)-\Phi(v,v)\in
\mathbb{R}.$$ From which we have $\Im (z+w)=0$ and $\Re (z-w)=0$; that is
\begin{equation*}
\begin{cases}
z+w-\overline{(z+w)}=0,\\
\\
z-w+\overline{(z-w)}=0.
\end{cases}
\end{equation*}
and adding member to member we have $$2z-2\overline{w}=0.$$ Hence

$$\Phi(u,v)=\overline{\Phi(v,u)}.$$

\medskip

\noindent Let us prove that \eqref{fqd-equivalenze-1b}$ \Longrightarrow$
\eqref{fqd-equivalenze-1c}.

\medskip

Let us assume that \eqref{fqd-equivalenze-1b} holds true. By the equivalence
proved previously we have
$$\Phi(u,v)=\overline{\Phi(v,u)}, \quad\quad\forall u,v\in
C^{\infty}\left(\mathbb{R}^n,\mathbb{C}\right).$$ Hence, for any $u,v\in C^{\infty}\left(\mathbb{R}^n,\mathbb{C}\right)$
we get

\begin{equation*}
\begin{aligned}
 \sum_{\alpha,\beta}a_{\alpha\beta}D^{\alpha}u\overline{D^{\beta}v}&=\overline{\sum_{\alpha,\beta}a_{\alpha\beta}D^{\alpha}v\overline{D^{\beta}u}}=\\&=
 \sum_{\alpha,\beta}\overline{a}_{\alpha\beta}\overline{D^{\alpha}v}D^{\beta}u=\sum_{\alpha,\beta}\overline{a}_{\beta\alpha}\overline{D^{\beta}v}D^{\alpha}u
\end{aligned}
\end{equation*}
(the last step is a mere change of indices). From what
just obtained and setting

$$c_{\alpha,\beta}=a_{\alpha\beta}-\overline{a}_{\beta\alpha},$$
we have
\begin{equation}\label{fqd-8}
\sum_{\alpha,\beta}c_{\alpha\beta}D^{\alpha}u\overline{D^{\beta}v}=0
\end{equation}
and Proposition \ref{fqd-4} gives
$$a_{\alpha\beta}-\overline{a}_{\beta\alpha}=c_{\alpha,\beta}=0$$ for any $\alpha,\beta\in \mathbb{N}_0^n$. Therefore
\eqref{fqd-equivalenze-1c} holds.

\medskip

\noindent Let us prove that \eqref{fqd-equivalenze-1c} $\Longrightarrow$
\eqref{fqd-equivalenze-1b}.

\medskip

Let us suppose 

$$a_{\alpha\beta}=\overline{a}_{\beta\alpha}.$$

Let $u\in C^{\infty}\left(\mathbb{R}^n,\mathbb{C}\right)$
be arbitrary. We get

\begin{equation*}
\begin{aligned}
 \sum_{\alpha,\beta}a_{\alpha\beta}D^{\alpha}u\overline{D^{\beta}u}&=\sum_{\alpha,\beta}\overline{a}_{\beta\alpha}D^{\alpha}u\overline{D^{\beta}u}=\\&
 =\overline{\sum_{\alpha,\beta}a_{\beta\alpha}\overline{D^{\alpha}u}D^{\beta}u}=
 \overline{\sum_{\alpha,\beta}a_{\alpha\beta}\overline{D^{\beta}u}D^{\alpha}u}=\overline{\Phi(u,u)}.
\end{aligned}
\end{equation*}

\medskip

\noindent the implication \eqref{fqd-equivalenze-1c}$
 \Longrightarrow$ \eqref{fqd-equivalenze-1d} can be proved in a similar way to the previous one (just replace $\zeta$ to $D$).

\medskip

\noindent Let us prove that
\eqref{fqd-equivalenze-1d}$\Longrightarrow$\eqref{fqd-equivalenze-1c}.

Let us assume that \eqref{fqd-equivalenze-1d} holds. Then

\begin{equation*}
\begin{aligned}
\sum_{\alpha,\beta}a_{\alpha\beta}\zeta^{\alpha}\overline{\zeta^{\beta}}&=\overline{\sum_{\alpha,\beta}a_{\alpha\beta}\zeta^{\alpha}\overline{\zeta^{\beta}}}=\\&=
\sum_{\alpha,\beta}\overline{a}_{\alpha\beta}\overline{\zeta^{\alpha}}\zeta^{\beta}=\sum_{\alpha,\beta}\overline{a}_{\beta\alpha}\overline{\zeta^{\beta}}\zeta^{\alpha}.
\end{aligned}
\end{equation*}
From what obtained above and setting

$$c_{\alpha,\beta}=a_{\alpha\beta}-\overline{a}_{\beta\alpha},$$
we get
\begin{equation*}
\sum_{\alpha,\beta}c_{\alpha\beta}\zeta^{\alpha}\overline{\zeta^{\beta}}=0,\quad\quad\forall
\zeta\in \mathbb{C}^n.
\end{equation*}
From which we have, for any $\gamma,\delta\in \mathbb{N}_0^n$,

\begin{equation}\label{fqd-8-bis}
\gamma!\delta!c_{\gamma\delta}=\partial_{\zeta}^{\gamma}\partial_{\overline{\zeta}}^{\delta}\left(
\sum_{\alpha,\beta}c_{\alpha\beta}\zeta^{\alpha}\overline{\zeta^{\beta}}
\right)_{|\zeta=0}=0,
\end{equation} where, $\zeta=(\zeta_1, \cdots, z_n)$,
$z_j=\xi_j+i\eta_j$,

\begin{equation*}
\partial_{\zeta_j}=\frac{1}{2}\left(\partial_{\xi_j}-i\partial_{\eta_j}\right),
\quad\quad
\partial_{\overline{\zeta}_j}=\frac{1}{2}\left(\partial_{\xi_j}+i\partial_{\eta_j}\right),\quad j=1,\cdots,n.
\end{equation*}
Finally, \eqref{fqd-8-bis} gives
\eqref{fqd-equivalenze-1c}. $\blacksquare$

\bigskip

From here on, it is convenient to use the following notations: to denote
a sesquilinear form with constant coefficients we will write

$$F(D,\overline{D})[u,\overline{v}]:=\Phi(u,v)=\sum_{\alpha,\beta}a_{\alpha\beta}D^{\alpha}u\overline{D^{\beta}v},$$
\index{$F(D,\overline{D})[u,\overline{v}]$}where $a_{\alpha\beta}\in \mathbb{C}$ are null except for a
finite set of multi--indices.

We will call \textbf{differential quadratic form} \index{form:@{form:}!- differential quadratic@{- differential quadratic}}with constant coefficients, the following form

\begin{equation}\label{fqd-9-28}
F(D,\overline{D})[u,\overline{u}]:=\Phi(u,u)=\sum_{\alpha,\beta}a_{\alpha\beta}D^{\alpha}u\overline{D^{\beta}u}.
\end{equation}
In what follows we will do the convention of the "regrouped terms" according to
to which the terms with the same indices $\alpha$ and $\beta$
occur only one time. With this convention, the following polynomial in $\zeta$ and $\overline{\zeta}$

\begin{equation}\label{fqd-9}
F(\zeta,\overline{\zeta})=\sum_{\alpha,\beta}a_{\alpha\beta}\zeta^{\alpha}\overline{\zeta^{\beta}},
\end{equation}
is uniquely associated to the form $F$. As a matter of fact, it turns out

$$F(\zeta,\overline{\zeta})=e^{-2(\Im\zeta)\cdot
x}F(D,\overline{D})\left[e^{i\zeta\cdot x},\overline{e^{i\zeta\cdot
x}}\right].$$ Hence, if
\begin{equation}\label{fqd-9-401}
F(D,\overline{D})[u,\overline{u}]=0,\quad \forall u\in
C^{\infty}\left(\mathbb{R}^n\right),
\end{equation}
then
\begin{equation}\label{fqd-9-401-1}
F(\zeta,\overline{\zeta})=0,\quad \forall\zeta\in
\mathbb{C}^n.
\end{equation}
Conversely, if \eqref{fqd-9-401-1} holds, then we have

$$\gamma!\delta!a_{\gamma\delta}=\partial_{\zeta}^{\gamma}\partial_{\overline{\zeta}}^{\delta}\left(
\sum_{\alpha,\beta}a_{\alpha\beta}\zeta^{\alpha}\overline{\zeta^{\beta}}
\right)_{|\zeta=0}=\partial_{\zeta}^{\gamma}\partial_{\overline{\zeta}}^{\delta}F(\zeta,\overline{\zeta})_{|\zeta=0}=0,$$
from which we get \eqref{fqd-9-401}. All in all \eqref{fqd-9-401} and
\eqref{fqd-9-401-1} are equivalent. 

We call the polynomial $F(\zeta,\overline{\zeta})$ the
\textbf{symbol of the  differential quadratic form} $F(D,\overline{D})$.\index{symbol of differential quadratic form}

\bigskip

The following Proposition will be useful later on.
\begin{prop}\label{parseval-fq}
Let $F(D,\overline{D})[u,\overline{u}]$ be a differential quadratic form. We have

\begin{equation}\label{note-42}
	\begin{aligned}
		&\int_{\mathbb{R}^n}
F(D,\overline{D})[u,\overline{u}]dx=\\&=\frac{1}{(2\pi)^n}\int_{\mathbb{R}^n}F(\xi,\xi)\left\vert
\widehat{u}(\xi)\right\vert^2 d\xi,\quad\forall u\in
C^{\infty}_0\left(\mathbb{R}^n\right).
\end{aligned}
\end{equation}
\end{prop}

\textbf{Proof.} Let $u\in C^{\infty}_0\left(\mathbb{R}^n\right)$ and
$$F(D,\overline{D})[u,\overline{u}]=\sum_{\alpha,\beta}a_{\alpha\beta}D^{\alpha}u\overline{D^{\beta}u}.$$
From the Parseval identity we have

\begin{equation*}
\begin{aligned}
\int_{\mathbb{R}^n}
F(D,\overline{D})[u,\overline{u}]dx&=\sum_{\alpha,\beta}a_{\alpha\beta}\int_{\mathbb{R}^n}D^{\alpha}u\overline{D^{\beta}u}dx=\\&=
\frac{1}{(2\pi)^n}\sum_{\alpha,\beta}a_{\alpha\beta}\int_{\mathbb{R}^n}\widehat{D^{\alpha}u}\overline{\widehat{D^{\beta}u}}d\xi=\\&=
\frac{1}{(2\pi)^n}\sum_{\alpha,\beta}a_{\alpha\beta}\int_{\mathbb{R}^n}\xi^{\alpha+\beta}\left\vert
\widehat{u}(\xi)\right\vert^2
d\xi=\\&=\frac{1}{(2\pi)^n}\int_{\mathbb{R}^n}F(\xi,\xi)\left\vert
\widehat{u}(\xi)\right\vert^2 d\xi.
\end{aligned}
\end{equation*}
$\blacksquare$

\bigskip

We are interested to examine \textbf{under what conditions $F(D,\overline{D})[u,\overline{u}]$ can be "written as a divergence" of some vector field}.

More precisely, we are interested in examining under
which conditions there exist some differential quadratic forms
$G_k(D,\overline{D})$, $k=1,\cdots,n$ with constant coefficients, such
that

\begin{equation}\label{fqd-10}
F(D,\overline{D})[u,\overline{u}]=\sum_{k=1}^n\partial_k\left(G_k(D,\overline{D})[u,\overline{u}]\right).
\end{equation}
For this purpose we consider the  differential quadratic form with constant coefficients

\begin{equation*}
G(D,\overline{D})[u,\overline{u}]=\sum_{\alpha,\beta}c_{\alpha\beta}D^{\alpha}u\overline{D^{\beta}u}
\end{equation*}
and we wish to express the symbol of
$\partial_k \left(G(D,\overline{D})[u,\overline{u}]\right)$ through the symbol of
$G(D,\overline{D})[u,\overline{u}]$. We have

\begin{equation*}
\begin{aligned}
 \partial_k \left(G(D,\overline{D})[u,\overline{u}]\right)&=\sum_{\alpha,\beta}c_{\alpha\beta}
 \partial_k\left(D^{\alpha}u\overline{D^{\beta}u}\right)=\\&=
\sum_{\alpha,\beta}c_{\alpha\beta}\left(\partial_kD^{\alpha}u\overline{D^{\beta}u}+D^{\alpha}u\partial_k\overline{D^{\beta}u}\right)=\\&=
\sum_{\alpha,\beta}c_{\alpha\beta}\left(iD_kD^{\alpha}u\overline{D^{\beta}u}-iD^{\alpha}u\overline{D_kD^{\beta}u}\right).
 \end{aligned}
\end{equation*}
Hence, the symbol associated to
$\partial_k \left(G(D,\overline{D})[u,\overline{u}]\right)$ is

\begin{equation*}
\begin{aligned}
 \sum_{\alpha,\beta}c_{\alpha\beta}\left(i\zeta_k\zeta^{\alpha}\overline{\zeta^{\beta}}-i\zeta^{\alpha}\overline{\zeta_k\zeta^{\beta}}\right)&=
 i\left(\zeta_k-\overline{\zeta}_k\right)\sum_{\alpha,\beta}c_{\alpha\beta}\zeta^{\alpha}\overline{\zeta^{\beta}}=\\&=
 i\left(\zeta_k-\overline{\zeta}_k\right)G(\zeta,\overline{\zeta}).
 \end{aligned}
\end{equation*}
Therefore, in order to have \eqref{fqd-10} it is necessary that

\begin{equation}\label{fqd-10-0}F(\zeta,\overline{\zeta})=i\sum_{k=1}^n\left(\zeta_k-\overline{\zeta}_k\right)G_k(\zeta,\overline{\zeta}).
\end{equation}
Now, setting $\zeta=\xi+i\eta$, where $\xi,\eta\in\mathbb{R}^n$,
we have

\begin{equation}\label{fqd-11}
F(\xi+i\eta,\xi-i\eta)=-2\sum_{k=1}^n\eta_kG_k(\xi+i\eta,\xi-i\eta).
\end{equation}
In particular, we have
\begin{equation}\label{fqd-12}
F(\xi,\xi)=0,\quad\quad \forall \xi\in \mathbb{R}^n
\end{equation}
and
\begin{equation}\label{fqd-13}
G_k(\xi,\xi)=-\frac{1}{2}\frac{\partial}{\partial\eta_k}F(\xi+i\eta,\xi-i\eta)_{|\eta=0},
\quad\quad\forall\xi\in\mathbb{R}^n.
\end{equation}
Therefore a \textbf{necessary condition to be true
\eqref{fqd-10}} is that \eqref{fqd-12} be true. Below we will see that
this condition is also sufficient, but first we give
the definition of the double and total order of a differential quadratic form.

\begin{definition}\label{ordine-fq}
	\index{Definition:@{Definition:}!- orders of a differential quadratic form@{- orders of a differential quadratic form}}
Let
\begin{equation}\label{fqd-14}
F(D,\overline{D})[u,\overline{u}]=\sum_{\alpha,\beta}a_{\alpha\beta}D^{\alpha}u\overline{D^{\beta}u},
\end{equation}
be a differential quadratic form with constant coefficients. We say that $F$ has \textbf{double order} $(\mu;m)$ provided
\begin{equation}\label{fqd-15}
a_{\alpha\beta}\neq 0\quad \Longrightarrow \quad
|\alpha|+|\beta|\leq \mu;\quad |\alpha|,|\beta|\leq m.
\end{equation}
$\mu$ is called the \textbf{total order} and $m$ is called
the \textbf{separated order} of the differential quadratic form $F$.
\end{definition}
It is evident that
$$\mu\leq 2m.$$
Moreover, when adopting the convention of the grouped terms, the
previous definition uniquely determines the order of the differential 
quadratic form. Here and in the sequel we will naturally extend the notions of
the double order and the total order to the symbol of a differential quadratic form.

\bigskip

\begin{lem}\label{integr-parti-cost}
Let $F(D,\overline{D})[u,\overline{u}]$ be a differential quadratic form with constant copefficients. Let us suppose that
\begin{equation}\label{fqd-16}
F(\xi,\xi)=0,\quad\quad \forall \xi\in \mathbb{R}^n,
\end{equation}
then there exist $n$ differential quadratic forms
$G_k(D,\overline{D})[u,\overline{u}]$, $k=1,\cdots,n$, such that
\begin{equation}\label{fqd-17}
F(D,\overline{D})[u,\overline{u}]=\sum_{k=1}^n\partial_k\left(G_k(D,\overline{D})[u,\overline{u}]\right)
\end{equation} and we have

\begin{equation}\label{fqd-13-ext}
G_k(\xi,\xi)=-\frac{1}{2}\frac{\partial}{\partial\eta_k}F(\xi+i\eta,\xi-i\eta)_{|\eta=0},
\quad\quad\forall\xi\in\mathbb{R}^n.
\end{equation}

In addition, let us assume that $F(D,\overline{D})[u,\overline{u}]$ has a double order $(\mu;m)$, $m>0$ then

\medskip

\noindent (a) if $\mu<2m$, the forms $G_k$, $k=1,\cdots,n$, can be choosen of double order $(\mu-1;m-1)$;

\medskip

\noindent (b) if $\mu=2m$, the forms $G_k$, $k=1,\cdots,n$, can be choosen of double order $(\mu-1;m)$.
\end{lem}

\bigskip

\textbf{Proof.} Let us assume that \eqref{fqd-16} holds.
Let us consider the polynomial

$$\eta\rightarrow F(\xi+i\eta,\xi-i\eta)$$ and let us apply the Taylor formula
at $\eta=0$. We have, for suitable polynomials $f_k$,
$k=1,\cdots,n$,

$$F(\xi+i\eta,\xi-i\eta)=\sum_{k=1}^n\eta_kf_k\left(\xi,\eta\right)=-\frac{i}{2}\sum_{k=1}^n \left(\zeta_k-\overline{\zeta}_k\right)f_k
\left(\frac{\zeta+\overline{\zeta}}{2},\frac{\zeta-\overline{\zeta}}{2i}\right).$$
Set

$$G_k\left(\zeta, \overline{\zeta}\right)=-\frac{1}{2}
f_k\left(\frac{\zeta+\overline{\zeta}}{2},\frac{\zeta-\overline{\zeta}}{2i}\right),$$
we get

\begin{equation}\label{fqd-17-0}
\begin{aligned}
 F(\zeta,\overline{\zeta})=i\sum_{k=1}^n\left(\zeta_k-\overline{\zeta}_k\right)G_k(\zeta,\overline{\zeta})=-2\sum_{k=1}^n\eta_kG_k(\xi+i\eta,\xi-i\eta),
 \end{aligned}
\end{equation}
from which we have

\begin{equation}\label{fqd-18}
\frac{\partial}{\partial\eta_k}F(\xi+i\eta,\xi-i\eta)=-2G_k(\xi+i\eta,\xi-i\eta),
\quad\quad\forall\xi\in\mathbb{R}^n.
\end{equation}
Hence
\begin{equation*}\label{40122}
G_k(\xi,\xi)=-\frac{1}{2}\frac{\partial}{\partial\eta_k}F(\xi+i\eta,\xi-i\eta)_{|\eta=0},
\quad\quad\forall\xi\in\mathbb{R}^n.
\end{equation*}
Moreover, from the first equality in \eqref{fqd-17-0} (retracing to
backward the calculations that led to \eqref{fqd-10-0}) we have

\begin{equation*}
 F(D,\overline{D})[u,\overline{u}]=i\sum_{k=1}^n\left(D_k-\overline{D}_k\right)\left(G_k(D,\overline{D})[u,\overline{u}]\right)=
 \sum_{k=1}^n\partial_k\left(G_k(D,\overline{D})[u,\overline{u}]\right).
\end{equation*}

\medskip

\textbf{Proof of (a) and (b)}.

\noindent \textbf{Case (a)}, $\mu<2m$.
Let us show that if $\alpha',\alpha'',\beta',\beta''$ are multi--indices
such that

\begin{equation}\label{fqd-20}
\begin{cases}
\alpha'+\beta'=\alpha''+\beta''\\
\\
\left|\alpha'\right|+\left|\beta'\right|=\left|\alpha''\right|+\left|\beta''\right|\leq\mu\leq
2m-1\\
\\
\left|\alpha'\right|,\mbox{ }\left|\beta'\right|, \mbox{
}\left|\alpha''\right|,\mbox{ }\left|\beta''\right|\leq m,
\end{cases}
\end{equation}
then

\begin{equation}\label{fqd-21}
\zeta^{\alpha'}\overline{\zeta^{\beta'}}=
\zeta^{\alpha''}\overline{\zeta^{\beta''}}+i\sum_{j=1}^n\left(\zeta_j-\overline{\zeta}_j\right)h_{j}(\zeta,\overline{\zeta}),
\end{equation}
where $h_{j}(\zeta,\overline{\zeta})$, $j=1,\cdots,n$, have the total order
less or equal than $\mu-1$ and the separated order less or equal than  $m-1$.

Notice that by \eqref{fqd-20} we have either $\left|\alpha'\right|<m$
or $\left|\beta'\right|<m$ (likewise for
$\left|\alpha''\right|$ and $\left|\beta''\right|$).

The proof consists of repeatedly applying both
simple identities.

\begin{subequations}
\label{fqd-id-1}
\begin{equation}
\label{fqd-id-1a}
\overline{\zeta}_j=\zeta_j-\left(\zeta_j-\overline{\zeta}_j\right),
\end{equation}
\begin{equation}
\label{fqd-id-1b}
\zeta_j=\overline{\zeta}_j+\left(\zeta_j-\overline{\zeta}_j\right).
\end{equation}
\end{subequations}
Let us consider $\zeta^{\alpha'}\overline{\zeta^{\beta'}}$ and let us suppose
$\left|\alpha'\right|<m$. Identity \eqref{fqd-id-1a} allows us to move the factors from $\overline{\zeta^{\beta'}}$ to
$\zeta^{\alpha'}$ as long as the exponent of $\zeta$ does not have modulus
$m$, when this occurs identity \eqref{fqd-id-1b} is used. Let us see more in detail. If
$\left|\alpha'\right|<m$ e $\beta'\neq 0$, for instance let
$\beta'_j>0$, then

\begin{equation*}
\begin{aligned}
\zeta^{\alpha'}\overline{\zeta^{\beta'}}&=\zeta^{\alpha'}\overline{\zeta^{\beta'-e_j}}\overline{\zeta}_j
=\zeta^{\alpha'}\overline{\zeta^{\beta'-e_j}}\left(\zeta_j-\left(\zeta_j-\overline{\zeta}_j\right)\right)=\\&=
\zeta^{\alpha'+e_j}\overline{\zeta^{\beta'-e_j}}-\left(\zeta_j-\overline{\zeta}_j\right)\zeta^{\alpha'}\overline{\zeta^{\beta'-e_j}}.
\end{aligned}
\end{equation*}
Let us notice that $\zeta^{\alpha'}\overline{\zeta^{\beta'-e_j}}$ has total order
$\left|\alpha'\right|+\left|\beta'-e_j\right|\leq
\mu-1$ and separated order less or equal than $m-1$ (recall
$\left|\alpha'\right|<m$ and $\left|\beta'\right|\leq m$). If,
on the other hand $\left|\beta'\right|<m$
 (and this includes the case
 $\beta'=0$ which was neglected previously) we use identity
\eqref{fqd-id-1b} and proceeding as above we reach (assuming
$\alpha'_j>0$, for some $j$) to

\begin{equation}\label{fqd-21-ext}\zeta^{\alpha'}\overline{\zeta^{\beta'}}=\zeta^{\alpha'-e_j}\overline{\zeta^{\beta'-e_j}}+
\left(\zeta_j-\overline{\zeta}_j\right)\zeta^{\alpha'-e_j}\overline{\zeta^{\beta'}}.
\end{equation}
Similarly to the case $\left|\alpha'\right|<m$, we have that 
$\zeta^{\alpha'-e_j}\overline{\zeta^{\beta'}}$ has total order 
less or equal than $\mu-1$ and separated order less or equal than
$m-1$. Repeatedly applying the procedure used above we arrive to \eqref{fqd-21}, which in turn implies that there exists $\widetilde{F}$ such that

\begin{equation}\label{fqd-XXII}
F(\zeta,\overline{\zeta})=\widetilde{F}(\zeta,\overline{\zeta})
+i\sum_{j=1}^n\left(\zeta_j-\overline{\zeta}_j\right)G_j(\zeta,\overline{\zeta}),
\end{equation}
where $G_{j}(\zeta,\overline{\zeta})$, $j=1,\cdots,n$, has double order
$(\mu-1;m-1)$ and
\begin{equation}\label{fqd-23}
\widetilde{F}(\zeta,\overline{\zeta})=\sum_{(\alpha,\gamma)\in\Lambda}c_{\gamma}\zeta^{\alpha}\overline{\zeta^{\gamma-\alpha}},
\end{equation}
where $c_{\gamma}\in \mathbb{C}$
$$\Lambda=\left\{(\alpha,\gamma)\in \mathbb{N}_0^n\times \mathbb{N}_0^n: \mbox{ } |\gamma|\leq m,\mbox{
}\alpha\leq\gamma\right\}$$ and also $\widetilde{F}$ has double order $(\mu;m)$. Notice that, thanks to \eqref{fqd-21} and \eqref{fqd-21-ext}, the summation in \eqref{fqd-23}
is written in such a way that for a given sum of the multi-indices
only one addend occurs.

By \eqref{fqd-16} and \eqref{fqd-XXII} we have
\begin{equation}\label{fqd-30}
0=\widetilde{F}(\xi,\xi)=\sum_{|\gamma|\leq
m}N_{\gamma}c_{\gamma}\xi^{\gamma},\quad\forall \xi\in \mathbb{R}^n,
\end{equation}
where $N_{\gamma}$ is the cardinality of the set
$\left\{\alpha\in \mathbb{N}_0^n:\mbox{ } \alpha\leq\gamma\right\}$
, hence, $c_{\gamma}=0$  from which we have  $\widetilde{F}\equiv 0$.
Therefore

\begin{equation*}\label{fqd-31}
F(\zeta,\overline{\zeta})=i\sum_{j=1}^n\left(\zeta_j-\overline{\zeta}_j\right)G_j(\zeta,\overline{\zeta})
\end{equation*}
and \eqref{fqd-17} is proved in case (a).

\medskip

\noindent \textbf{Case (b)}, $\mu=2m$.

Obviously, we can handle the terms of
\begin{equation*}
F(\zeta,\overline{\zeta})=\sum_{\alpha,\beta}a_{\alpha\beta}\zeta^{\alpha}\overline{\zeta^{\beta}}.
\end{equation*}
satisfying $|\alpha|+|\beta|<2m$ (and $|\alpha|,|\beta|\leq m$) in the same way of case (a). Let us examine in which a way we can handle the terms such that $|\alpha|+|\beta|=2m$. Since
$|\alpha|,|\beta|\leq m$ we have $|\alpha|=|\beta|=m$.

By identities \eqref{fqd-id-1a} and \eqref{fqd-id-1b} we have

$$\zeta_j\overline{\zeta}_k=\zeta_k\overline{\zeta}_j+\left(\zeta_j-\overline{\zeta}_j\right)\overline{\zeta}_k-\left(\zeta_k-\overline{\zeta}_k\right)\overline{\zeta}_j,$$
for $j,k=1,\cdots,n$.

Now, let us suppose that $\alpha',\alpha'',\beta',\beta''$
satisfy
\begin{equation*}
\begin{cases}
\alpha'+\beta'=\alpha''+\beta'',\\
\\
\left|\alpha'\right|=\left|\beta'\right|=\left|\alpha''\right|=\left|\beta''\right|=m,
\end{cases}
\end{equation*}
then we have (if $\alpha'_j>0$ and $\beta'_k>0$)
\begin{equation*}
\begin{aligned}
\zeta^{\alpha'}\overline{\zeta^{\beta'}}&=\zeta^{\alpha'-e_j}\overline{\zeta^{\beta'-e_k}}\zeta_j\overline{\zeta}_k=\zeta^{\alpha'-e_j+e_k}\overline{\zeta^{\beta'-e_k+e_j}}+\\&
+\left(\zeta_j-\overline{\zeta}_j\right)h_1(\zeta,\overline{\zeta})+
\left(\zeta_k-\overline{\zeta}_k\right)h_2(\zeta,\overline{\zeta}),
\end{aligned}
\end{equation*}
where
$$h_1(\zeta,\overline{\zeta})=\zeta^{\alpha'-e_j}\overline{\zeta^{\beta'}},\quad\quad
h_2(\zeta,\overline{\zeta})=-\zeta^{\alpha'-e_j}\overline{\zeta^{\beta'-e_k+e_j}},$$
have double order $(\mu-1,m)$. From now on, one may
proceed as in case (a) and we reach the conclusion.
$\blacksquare$

\bigskip

In the case of a differential quadratic forms with variable coefficients we have the
following

\begin{lem}\label{integr-parti-var}
Let \index{$F(x,D,\overline{D})[u,\overline{u}]$}
\begin{equation}\label{fqd-1-bis}
F(x,D,\overline{D})[u,\overline{u}]=\sum_{\alpha,\beta}a_{\alpha\beta}(x)D^{\alpha}u\overline{D^{\beta}u}
\end{equation}
be a differential quadratic form with variable coefficients $a_{\alpha\beta}\in
C^{s}(\Omega)$, where $s\in \mathbb{N}$. Let us suppose that $F$ has
double order $(\mu;m)$, $m>0$, and that
\begin{equation}\label{fqd-16-var}
F(x,\xi,\xi)=0,\quad\quad \forall x\in \Omega,\mbox{ } \forall\xi\in
\mathbb{R}^n.
\end{equation}
Then there exists a differential quadratic form 
$G(x,D,\overline{D})[u,\overline{u}]$ whose coefficients belong to 
$C^{s-1}(\Omega)$ such that
\begin{equation}\label{fqd-17-var}
\int_{\Omega}F(x,D,\overline{D})[u,\overline{u}]dx=\int_{\Omega}G(x,D,\overline{D})[u,\overline{u}]dx,\quad\forall
u\in C_0^{\infty}(\Omega)
\end{equation}
and such that:

\medskip

\noindent (a) if $\mu<2m$, then $G(x,D,\overline{D})$ can be
chosen of double order $(\mu-1;m-1)$;

\medskip

\noindent (b) if $\mu=2m$, then $G(x,D,\overline{D})$ can be
chosen of double order $(\mu-1;m)$. Moreover

\begin{equation}\label{fqd-18-var}
G(x,\xi,\xi)=\frac{1}{2}\sum_{k=1}^n\partial_{x_k\eta_k}^2F(x,\xi+i\eta,
\xi-i\eta)_{|\eta=0}.
\end{equation}
\end{lem}
\textbf{Proof.} Let $F_1,\cdots, F_N$ be a basis of the
vector space (of finite dimension) of all quadratic forms $H$
of double order $(\mu;m)$ with constant coefficients and satisfying
$$H(\xi,\xi)=0,\quad\quad\forall \xi \in \mathbb{R}^n.$$ By Lemma \ref{integr-parti-cost} there exist differential quadratic forms with
constant coefficients $G^k_j$, $j=1,\cdots, N$, $k=1,\cdots, n$, of
double order $(\mu-1;m')$, with $m'=m-1$, provided $\mu<2m$, and $m'=m$
provided $\mu=2m$, such that

$$F_j(D,\overline{D})[u,\overline{u}]=\sum_{k=1}^n\partial_{x_k}\left(G^k_j(D,\overline{D})[u,\overline{u}]\right),\quad
j=1,\cdots, N.$$ Now, \eqref{fqd-16-var} implies that there exist
$c_j\in C^{s}(\Omega)$, $j=1,\cdots, N$, such that

$$F(x,D,\overline{D})[u,\overline{u}]=\sum_{j=1}^N c_j(x)
F_j(D,\overline{D})[u,\overline{u}].$$ Hence, if $u\in
C_0^{\infty}(\Omega)$, then integration by parts yields

\begin{equation*}
\begin{aligned}
\int_{\Omega}F(x,D,\overline{D})[u,\overline{u}]dx&=
\sum_{j=1}^N\sum_{k=1}^n\int_{\Omega}c_j(x)\partial_{x_k}\left(G^k_j(D,\overline{D})[u,\overline{u}]\right)dx=\\&=
-\sum_{j=1}^N\sum_{k=1}^n\int_{\Omega}\partial_{x_k}c_j(x)G^k_j(D,\overline{D})[u,\overline{u}]dx.
\end{aligned}
\end{equation*}
Thus, we can choose
$$G(x,D,\overline{D})[u,\overline{u}]=\sum_{j=1}^N\sum_{k=1}^n-\partial_{x_k}c_j(x)G^k_j(D,\overline{D})[u,\overline{u}].$$
From which we have 
\begin{equation}\label{fqd-19-var}
G(x,\xi,\xi)=\sum_{j=1}^N\sum_{k=1}^n-\partial_{x_k}c_j(x)G^k_j(\xi,\xi).
\end{equation}
On the other hand by \eqref{fqd-13-ext} we have
\begin{equation*}
G_j^k(\xi,\xi)=-\frac{1}{2}\frac{\partial}{\partial\eta_k}F_j(\xi+i\eta,\xi-i\eta)_{|\eta=0},
\quad\quad\forall\xi\in\mathbb{R}^n.
\end{equation*}
Now, by the last equality we get

\begin{equation*}
\begin{aligned}
G(x,\xi,\xi)&=\sum_{k=1}^n\sum_{j=1}^N-\partial_{x_k}c_j(x)G^k_j(\xi,\xi)=\\&=
\frac{1}{2}\sum_{k=1}^n\left(\partial^2_{x_k\eta_k}\sum_{j=1}^Nc_j(x)F_j(\xi+i\eta,\xi-i\eta)
\right)_{|\eta=0}=\\&=
\frac{1}{2}\sum_{k=1}^n\partial^2_{x_k\eta_k}F(\xi+i\eta,\xi-i\eta)_{|\eta=0}.
\end{aligned}
\end{equation*}
$\blacksquare$
\section[The conjugate of $P_m(x,D)$]{The conjugate of $P_m(x,D)$ -- Set up of a Carleman estimate}\label{sec-coniugato} 
In this Section we will consider the
conjugate of the operator $P_m(x,D)$ which, we recall, is defined
by

\begin{equation}\label{sec-coniugato-1-28}
P_{\tau}v=e^{\tau\varphi}P_m(x,D)\left(e^{-\tau\varphi}v\right).
\end{equation}
We first observe that
\begin{equation}\label{sec-coniugato-1}
P_{\tau}v=P_m(x,D+i\tau\nabla\varphi(x))v.
\end{equation}
As a matter of fact we have
$$e^{\tau\varphi}D_j\left(e^{-\tau\varphi}v\right)=D_jv-\tau
(D_j\varphi)v=D_jv+i\tau (\partial_j\varphi)v,$$
$$e^{\tau\varphi}D_kD_j\left(e^{-\tau\varphi}v\right)=e^{\tau\varphi}D_k\left(e^{-\tau\varphi}(D_jv+i\tau
(\partial_j\varphi)v)\right)=(D_k+i\tau
(\partial_k\varphi))(D_j+i\tau (\partial_j\varphi))v,$$
$$\vdots$$
$$e^{\tau\varphi}D^{\alpha}\left(e^{-\tau\varphi}v\right)=\left(D+i\tau\nabla \varphi(x)\right)^{\alpha}v,$$
for every  $\alpha\in \mathbb{N}_0^n$. Now, since
\begin{equation}\label{sec-coniugato-1bis}
P_m(x,D)=\sum_{|\alpha|=m}a_{\alpha}(x)D^{\alpha},
\end{equation} 
we have
\begin{equation*}
\begin{aligned}
e^{\tau\varphi}P_m(x,D)\left(e^{-\tau\varphi}v\right)&=\sum_{|\alpha|=m}a_{\alpha}(x)e^{\tau\varphi}D^{\alpha}\left(e^{-\tau\varphi}v\right)=\\&=
\sum_{|\alpha|=m}a_{\alpha}(x)\left(D+i\tau\nabla
\varphi(x)\right)^{\alpha}v=\\&=P_m(x,D+i\tau\nabla \varphi(x))v.
\end{aligned}
\end{equation*}

Now, let us consider the  polynomial $P_m(x,\xi+i\tau\nabla \varphi(x))$ in the variable $\xi$ and let us denote by
\begin{equation}\label{note-1-49} \index{$p_m(x,D,\tau)$}
p_m(x,D,\tau)\mbox{ \textbf{the operator whose symbol is} }
P_m(x,\xi+i\tau\nabla \varphi).\end{equation} Let us note that, in general, operator \eqref{note-1-49} does not equal to the operator $P_m(x,D+i\tau\nabla \varphi)$.
 For instance, if
$$P_2(x,D)=D^2_{1}$$
we have
$$P_2(x,D+i\tau\nabla \varphi(x))v= D^2_{1}v+2i\tau\partial_1\varphi
D_1v-\tau^2(\partial_1\varphi)^2v+\tau\left(\partial^2_1\varphi\right)v$$
hence

$$p_2(x,D,\tau)=P_2(x,D+i\tau\nabla \varphi(x))-\tau\partial^2_1\varphi(x).$$
In general we have
\begin{equation}\label{note-2-49}
p_m(x,D,\tau)=\sum_{|\alpha|+j=m}\tau^jb_{\alpha j}(x) D^{\alpha},
\end{equation}
where the coefficients $b_{\alpha j}(x)$ depend on $\nabla \varphi$,
on coefficients of $P_m(x,D)$ (but not on their derivatives) and \textbf{do not
depend} on $\tau$. In addition we have

\begin{equation}\label{note-3-49}
P_m(x,D+i\tau\nabla \varphi(x))=p_m(x,D,\tau)+R_{m-1,\tau}(x,D,\tau),
\end{equation}
where
$$R_{m-1,\tau}(x,D)=\sum_{|\alpha|+j\leq m-1}\tau^j\widetilde{b}_{\alpha
j}(x) D^{\alpha},$$ the coefficients $\widetilde{b}_{\alpha j}$
depend on the coefficients of $P_m(x,D)$ (but not on their
derivatives), on $\nabla \varphi$ and on the higher-order derivatives
of $\varphi$ and \textbf{do not depend on} $\tau$. The second term on the right--hand side in
\eqref{note-3-49}, as we will realize soon, may
be regarded as a harmless perturbation of the operator
$p_m(x,D,\tau)$.

Now let us deal with the square in the integral
\begin{equation}
\int\left\vert P_m(x,D+i\tau\nabla\varphi(x))v\right\vert^2dx,
\end{equation}
From here on, since $v$ is supported in $\Omega$,
we omit the set of integration. First we notice that by
\eqref{note-3-49} we have

\begin{equation}\label{note-4-49}
\begin{aligned}
\int\left\vert P_m(x,D+i\tau\nabla\varphi(x))v\right\vert^2dx&\geq
\frac{1}{2}\int\left\vert
p_m(x,D,\tau)v\right\vert^2dx-\\&-\int\left\vert
R_{m-1,\tau}(x,D)v\right\vert^2dx\geq\\&\geq
\frac{1}{2}\int\left\vert
p_m(x,D,\tau)v\right\vert^2dx-\\&-C\sum_{|\alpha|\leq
m-1}\tau^{2(m-|\alpha|)-2}\int\left\vert D^{\alpha}v\right\vert^2dx,
\end{aligned}
\end{equation}
where $C$ depends by the $L^{\infty}$ norms of the coefficients of
$P_m(x,D)$.
 
 \medskip

 We need some additional notation. Let $M(x,\xi)$ be
a polynomial with respect to the variable $\xi$, let us suppose that the
coefficients of $M(x,\xi)$ are differentiable. Let us set 
$$M^{(j)}(x,\xi)=\partial_{\xi_j}M(x,\xi),\quad\quad
M_{,j}(x,\xi)=\partial_{x_j}M(x,\xi),\quad j=1,\cdots,n.$$
Let us denote by $$\overline{M}(x,\xi)$$ the polynomial in $\xi$ whose 
coefficients are the complex conjugate of the coefficients of
$M(x,\xi)$. Keep in mind that if $\zeta\in \mathbb{C}^n$, then
$$\overline{M(x,\zeta)}=\overline{M}\left(x,\overline{\zeta}\right).$$
If $L(x,\xi)$ and $M(x,\xi)$ are two polynomials in the variable $\xi$ with 
 differentiable coefficients, we define their \textbf{Poisson brackets} \index{Poisson brackets}

\begin{equation}\label{note-parentesi-p}
	\begin{aligned}
&\left\{L(x,\xi),M(x,\xi)\right\}=\\&=\sum_{j=1}^n\left(L^{(j)}(x,\xi)M_{,j}(x,\xi)-L_{,j}(x,\xi)M^{(j)}(x,\xi)\right).
\end{aligned}
\end{equation}

Now we anticipate that in points 3 and 4 of the Remarks of the present Section, we will observe that $\overline{p}_m(x,D,\tau)$ is a suitable approximation of the adjoint of the operator $p_m(x,D,\tau)$. 

Let
\begin{subequations}
\label{note-5-49}
\begin{equation}
\label{note-5-49a}
S(x,D,\tau)=\frac{1}{2}\left(p_m(x,D,\tau)+\overline{p}_m(x,D,\tau)\right),
\end{equation}
\begin{equation}\label{note-5-49b}
A(x,D,\tau)=\frac{1}{2}\left(p_m(x,D,\tau)-\overline{p}_m(x,D,\tau)\right),
\end{equation}
\end{subequations}

  We have trivially
\begin{equation} \label{note-5-49a-2911}
p_m(x,D,\tau)=S(x,D,\tau)+A(x,D,\tau).
\end{equation}
 Hence
\begin{equation}
\begin{aligned}\label{note-1-50}
&\int\left\vert p_m(x,D,\tau)v\right\vert^2dx=\\&=\int\left\vert
S(x,D,\tau)v\right\vert^2dx+\int\left\vert
A(x,D,\tau)v\right\vert^2dx+\\&+2\int\Re\left(S(x,D,\tau)v\overline{A(x,D,\tau)v}\right)dx.
\end{aligned}
\end{equation}

\medskip

A \textbf{crucial point in the proof of a Carleman estimate} consists in finding an appropriate estimate from below of the third integral on the right hand side of \eqref{note-1-50}.

Now, let us consider the differential quadratic form

\begin{equation}
\label{note-1ext-50}
F(x,D,\overline{D},\tau)\left[v,\overline{v}\right]=2\Re\left(S(x,D,\tau)v\overline{A(x,D,\tau)v}\right),
\end{equation}
whose symbol is
$$F(x,\zeta,\overline{\zeta},\tau)=2\Re\left(S(x,\zeta,\tau) \overline{A(x,\zeta,\tau)}\right),\quad \forall \zeta\in \mathbb{C}^n.$$
By the definition of $p_m(x,D,\tau)$ we have that the symbol of 
$\overline{p}_m(x,D,\tau)$ is given by $\overline{P}_m(x,\xi-i\tau\nabla\varphi)$, from which we have, for
$\zeta=\xi+i\eta\in\mathbb{C}^n$

\begin{equation*}
\begin{aligned}
&F(x,\zeta,\overline{\zeta},\tau)=2\Re\left(S(x,\zeta,\tau)\overline{A(x,\zeta,\tau)}\right)=\\&=
\frac{1}{2}\Re\left(\left(P_m(x,\zeta+i\tau\nabla\varphi(x))+\overline{P}_m(x,\zeta-i\tau\nabla\varphi(x))\right)\times\right.\\&\left.\times\left(\overline{P_m(x,\zeta+i\tau\nabla\varphi(x))-
	\overline{P}_m(x,\zeta-i\tau\nabla\varphi(x))}\right)\right)=\\&=
\frac{1}{2}\left(\left\vert
P_m(x,\zeta+i\tau\nabla\varphi(x))\right\vert^2-\left\vert
\overline{P}_m(x,\zeta-i\tau\nabla\varphi(x))\right\vert^2\right).
\end{aligned}
\end{equation*}
Hence
\begin{equation}\label{note-2-50}
	\begin{aligned}
		&F(x,\zeta,\overline{\zeta},\tau)=\\&=
		\frac{1}{2}\left(\left\vert
P_m(x,\zeta+i\tau\nabla\varphi(x))\right\vert^2-\left\vert
\overline{P}_m(x,\zeta-i\tau\nabla\varphi(x))\right\vert^2\right).
\end{aligned}
\end{equation}
By this equality we get

\begin{equation*}
\begin{aligned}
 F(x,\xi,\xi,\tau)&
 =\frac{1}{2}\left(\left\vert
P_m(x,\xi+i\tau\nabla\varphi(x))\right\vert^2-\left\vert
\overline{P}_m(x,\xi-i\tau\nabla\varphi(x))\right\vert^2\right)=\\&=
\frac{1}{2}\left(\left\vert
P_m(x,\xi+i\tau\nabla\varphi(x))\right\vert^2-\left\vert
\overline{P_m(x,\xi+i\tau\nabla\varphi(x))}\right\vert^2\right)=0.
\end{aligned}
\end{equation*}
Hence, \textbf{assuming that the coefficients of $P_m(x,D)$ belong to
$C^1\left(\overline{\Omega}\right)$} we can apply Lemma
\ref{integr-parti-var}. Using formula \eqref{fqd-18-var} and
denoting by

\begin{equation}\label{note-3-50}
G(x,\xi,\xi,\tau):=\frac{1}{2}\sum_{k=1}^n\partial_{x_k\eta_k}^2F(x,\xi+i\eta,
\xi-i\eta,\tau)_{|\eta=0},
\end{equation}
we have

\begin{equation}\label{note-4-50}
\begin{aligned}
2\int\Re\left(S(x,D,\tau)v\overline{A(x,D,\tau)v}\right)dx&=\int
F(x,D,\overline{D},\tau)[v,\overline{v}]dx=\\&=\int
G(x,D,\overline{D},\tau)[v,\overline{v}]dx,
\end{aligned}
\end{equation}
for every $v\in C_0^{\infty}(\Omega)$.

Now we calculate the expression on the right hand side in \eqref{note-3-50}. Although
the calculation is elementary, let us perform it in detail. We have

\begin{equation}\label{note-5-50}
\begin{aligned}
&\partial_{\eta_k}\frac{1}{2}\left(\left\vert
P_m(x,\xi+i\eta+i\tau\nabla\varphi(x))\right\vert^2-\right.\\&\left.-\left\vert
\overline{P}_m(x,\xi+i\eta-i\tau\nabla\varphi(x))\right\vert^2\right)_{|\eta=0}=\\&=
\partial_{\eta_k}\frac{1}{2}\left(\left\vert
P_m(x,\xi+i\eta+i\tau\nabla\varphi(x))\right\vert^2-\right.\\&\left.-\left\vert
P_m(x,\xi-i\eta+i\tau\nabla\varphi(x))\right\vert^2\right)_{|\eta=0}=\\&=
2\Re\left(-iP_m(x,\xi+i\tau\nabla\varphi(x))\overline{P^{(k)}_m(x,\xi+i\tau\nabla\varphi(x))}\right).
\end{aligned}
\end{equation}
Set for short
$$\zeta=\xi+i\tau\nabla \varphi (x)$$
and let us differentiate  what obtained in \eqref{note-5-50} w.r.t. $x_k$.
By \eqref{note-3-50} we get 

\begin{equation}\label{note-6-50}
\begin{aligned}
&G(x,\xi,\xi,\tau)=\\&=\tau\sum_{j,k=1}^n\partial^2_{x_jx_k}\varphi(x)
P_m^{(j)}(x,\zeta)\overline{P_m^{(k)}(x,\zeta)}
+\\&+\Im\left(\sum_{k=1}^nP_{m,k}(x,\zeta)\overline{P_m^{(k)}(x,\zeta)}\right)+\\&+\Im\left[
P_m(x,\zeta)\left(\sum_{k=1}^n\overline{P_{m,k}^{(k)}(x,\zeta)}
-i\tau\sum_{j,k=1}^n\overline{P_m^{(k,j)}(x,\zeta)}\partial^2_{x_jx_k}\varphi (x)\right)\right].
\end{aligned}
\end{equation}
Let us observe that we have

\begin{equation}\label{note-6-50-28}
\begin{aligned}
&P_m(x,\xi+i\tau\nabla\varphi)=0\mbox{ }
\Rightarrow\\&\Rightarrow\mbox{ }
G(x,\xi,\xi,\tau)=\frac{i}{2}\left\{P_m(x,\xi+i\tau\nabla\varphi),\overline{P_m(x,\xi+i\tau\nabla\varphi)}\right\},
\end{aligned}
\end{equation}
where $\{\cdot,\cdot\}$ is the Poisson bracket defined in
\eqref{note-parentesi-p}.  In order to check  \eqref{note-6-50-28}
it suffices to develop the Poisson bracket in \eqref{note-6-50-28}, and
to notice that the third term on the right hand side in \eqref{note-6-50} vanishes when $P_m(x,\xi+i\tau\nabla\varphi(x))=0$.

\medskip

Let us notice that $G(x,\xi,\tau)$ is
a homogeneous polynomial of degree $2m-1$ in the variables $(\xi,\tau)$.
However, we will be interested in more precise information about the
differential quadratic form $G(x,D,\overline{D},\tau)$ or, equivalently on
its symbol $G(x,\zeta,\overline{\zeta},\tau)$, to this end
we prove

\begin{prop}\label{note-52-prop}
Let $P_m(x,D)$ the differential operator
\begin{equation*}
P_m(x,D)=\sum_{|\alpha|=m}a_{\alpha}(x)D^{\alpha},
\end{equation*}
where $a_{\alpha}\in C^1\left(\overline{\Omega},\mathbb{C}\right)$,
for $|\alpha|=m$. Let $F(x,\zeta,\overline{\zeta},\tau)$ be defined by
\eqref{note-2-50}. Then there exists a differential quadratic form
$G\left(x,D,\overline{D},\tau \right)$ such that

\begin{equation}\label{note-52-prop-1}
\begin{aligned}
\int F\left(x,D,\overline{D},\tau\right)[v,\overline{v}]dx=\int
G\left(x,D,\overline{D},\tau \right)[v,\overline{v}]dx.
\end{aligned}
\end{equation}
Moreover
\begin{equation}\label{note-52-prop-2}
G\left(x,D,\overline{D},\tau \right)=\sum_{h=0}^{2m-1}\tau^{h}
G^{(h)}\left(x,D,\overline{D}\right),
\end{equation}
where $G^{(h)}\left(x,D,\overline{D}\right)$ is a differential quadratic form
which has double order  $\left(2m-h-1;m\right)$, for $h=0,1 \cdots,
2m-1$.

 If the \textbf{coefficients $a_{\alpha}$, for $|\alpha|=m$, are real valued functions} then  \eqref{note-52-prop-1} continues to hold true,
 but instead of \eqref{note-52-prop-2} we have
\begin{equation}\label{note-52-prop-3}
G\left(x,D,\overline{D},\tau \right)=\tau\sum_{h=0}^{2m-2}\tau^{h}
G^{(h)}\left(x,D,\overline{D}\right),
\end{equation}
where $G^{(h)}\left(x,D,\overline{D}\right)$ is a differential quadratic form
which has double order  $\left(2m-h-2;m\right)$, $h=0, \cdots,
2m-2$.

In  any  case $G(x,\xi,\xi,\tau)$ is given by \eqref
{note-6-50}.
\end{prop}

\textbf{Proof.} By the Taylor formula we get

\begin{equation}\label{note-6quater-50}
P_m(x,\zeta+i\tau\nabla\varphi(x))=\sum_{k=0}^m\tau^kq_{m-k}(x,\zeta),\quad\forall\zeta
\in \mathbb{C}^n,
\end{equation}
where, for $k=0,1,\cdots,m$, $q_{m-k}(x,\zeta)$ are polynomials in the variable
$\zeta$ of degree  $m-k$. Moreover the  coefficients of $q_{m-k}(x,\zeta)$
are of class $C^1\left(\overline{\Omega}\right)$. We have

\begin{equation*}
\left\vert
P_m(x,\zeta+i\tau\nabla\varphi(x))\right\vert^2=\sum_{k,j=0}^m\tau^{k+j}q_{m-k}(x,\zeta)\overline{q}_{m-j}(x,\overline{\zeta}).
\end{equation*}
Hence, by \eqref{note-2-50} we get

\begin{equation*}
F\left(x,\zeta,\overline{\zeta},\tau\right)=\sum_{k,j=0}^m\tau^{k+j}F_{kj}\left(x,\zeta,\overline{\zeta}\right),
\end{equation*}
where, for $j,k=1,\cdots,m$
\begin{equation}\label{note-star-52}
F_{kj}\left(x,\zeta,\overline{\zeta}\right)=q_{m-k}(x,\zeta)\overline{q}_{m-j}(x,\overline{\zeta})-
q_{m-k}(x,\overline{\zeta})\overline{q}_{m-j}(x,\zeta).\end{equation}
Each of the forms $F_{kj}$ has double order
$\left(2m-(k+j);m\right)$, furthermore by \eqref{note-star-52},
since $q_{00}$ has degree $0$, we have
\begin{equation}\label{note-starbis-52}
F_{mm}\left(x,\zeta,\overline{\zeta}\right)=0,\quad\forall\zeta \in
\mathbb{C}^n.
\end{equation}
Moreover

$$F_{kj}\left(x,\xi,\xi\right)=0,\quad\quad \forall \xi\in
\mathbb{R}^n$$ and by Lemma \eqref{integr-parti-var} -- case (b)
-- we have that, for $j,k=1,\cdots,m$, where either $j$ or $k$ are different from $m$,
there exist a differential quadratic form $G_{kj}$ which have  double order
$\left(2m-(k+j)-1;m\right)$ and satisfying

$$\int F_{kj}\left(x,D,\overline{D}\right)[v,\overline{v}]dx=\int
G_{kj}\left(x,D,\overline{D}\right)[v,\overline{v}]dx,\quad\forall
v\in C_0^{\infty}(\Omega),$$ of course, since
\eqref{note-starbis-52} holds, we may choose

\begin{equation*}
G_{mm}\equiv 0.
\end{equation*}
Therefore, by the last obtained equality and setting
$$G^{(h)}\left(x,D,\overline{D}\right)=\sum_{k+j=h}G_{kj}\left(x,D,\overline{D}\right), \quad h=1,\cdots, 2m-1,$$

$$G\left(x,D,\overline{D},\tau\right)[v,\overline{v}]=\sum_{h=0}^{2m-1}\tau^h
G^{(h)}\left(x,D,\overline{D}\right)[v,\overline{v}],$$ we have that
$G^{(h)}$ is a differential quadratic form which has double order
$\left(2m-h-1;m\right)$ and
\begin{equation}\label{note-6extra-50}
\begin{aligned}
\int
F\left(x,D,\overline{D},\tau\right)[v,\overline{v}]dx&=\sum_{k,j=0}^m\tau^{k+j}\int
F_{kj}\left(x,D,\overline{D}\right)[v,\overline{v}]dx=\\&=
\sum_{k,j=0}^m\tau^{k+j}\int
G_{jk}\left(x,D,\overline{D}\right)[v,\overline{v}]dx=\\&=
\sum_{h=0}^{2m-1}\tau^{h}\int
G^{(h)}\left(x,D,\overline{D}\right)[v,\overline{v}]dx=\\&= \int
G\left(x,D,\overline{D},\tau\right)[v,\overline{v}]dx.
\end{aligned}
\end{equation}

If \textbf{the coefficients of $P_m(x,D)$ are real--valued} then
also the coefficients of the polynomials $q_j(x,\zeta)$ in  \eqref{note-6quater-50} are real--valued and \eqref{note-star-52} can be written as

\begin{equation}\label{note-1-51}
F_{kj}\left(x,\zeta,\overline{\zeta}\right)=q_{m-k}(x,\zeta)q_{m-j}(x,\overline{\zeta})-
q_{m-k}(x,\overline{\zeta})q_{m-j}(x,\zeta).\end{equation} Hence,
besides \eqref{note-starbis-52},  we have 
\begin{equation}\label{note-3a}
F_{00}\left(x,\zeta,\overline{\zeta}\right)=0,\quad\forall \zeta\in
\mathbb{C}^n.
\end{equation}
Therefore we have
\begin{equation*}
F\left(x,\zeta,\overline{\zeta},\tau\right)=\tau\sum_{h=0}^{2m-2}\tau^h
F^{(h)}\left(x,\zeta,\overline{\zeta}\right) ,
\end{equation*}
where
\begin{equation*}
F^{(h)}\left(x,\zeta,\overline{\zeta}\right)=\sum_{k+j=h+1}
F_{kj}\left(x,\zeta,\overline{\zeta}\right).
\end{equation*}
Hence $F^{(h)}$, for $h=0,\cdots,2m-2$, has double order
$(2m-h-1;m)$. Therefore by applying Lemma \eqref{integr-parti-var}
-- case (a) -- there exist
$G^{(h)}\left(x,D,\overline{D}\right)$, differential quadratic forms which have double
order $(2m-h-2;m-1)$, such that

$$\int F^{(h)}\left(x,D,\overline{D}\right)[v,\overline{v}]dx=\int
G^{(h)}\left(x,D,\overline{D}\right)[v,\overline{v}]dx,\quad\forall
v\in C_0^{\infty}(\Omega),$$ and setting
$$G\left(x,D,\overline{D},\tau\right)[v,\overline{v}]=\tau\sum_{h=0}^{2m-2}\tau^h
G^{(h)}\left(x,D,\overline{D}\right)[v,\overline{v}]$$ we get

$$\int F\left(x,D,\overline{D},\tau\right)[v,\overline{v}]dx=\int
G\left(x,D,\overline{D},\tau\right)[v,\overline{v}]dx.$$
$\blacksquare$

\bigskip

\bigskip

Now, we  \textbf{broadly outline} the main ideas that
are involved in proving a Carleman estimate.  We come back,
then, to the third integral in \eqref{note-1-50}. Let $x_0\in
\overline{\Omega}$ we have from \eqref{note-4-50} and \eqref{note-42},

\begin{equation}\label{note-6bis-50}
\begin{aligned}
2\int\Re&\left(S(x,D,\tau)v\overline{A(x,D,\tau)v}\right)dx=\int
G(x,D,\overline{D},\tau)[v,\overline{v}]dx=\\&=\int
G(x_0,D,\overline{D},\tau)[v,\overline{v}]dx+\\&+\int
\left(G(x,D,\overline{D},\tau)-G(x_0,D,\overline{D},\tau)\right)[v,\overline{v}]dx=\\&=
\frac{1}{(2\pi)^n}\int
G(x_0,\xi,\xi,\tau)\left|\widehat{v}(\xi)\right|^2d\xi+\\&+
\underset{\mathcal{R}}{\underbrace{\int
\left(G(x,D,\overline{D},\tau)-G(x_0,D,\overline{D},\tau)\right)[v,\overline{v}]dx}}.
\end{aligned}
\end{equation}

The main idea that we will follow consits essentially in what follows:

\medskip

\noindent\textbf{(a) Choice of $\varphi$. } The choice of $\varphi$
will be made so that we have
\begin{equation}\label{note-6ter-50}
	\begin{aligned}
		&P_m(x,\xi+i\tau\nabla\varphi(x))=0\mbox{ }\Rightarrow\\&
\Rightarrow G(x,\xi,\xi,\tau)>0, \quad (\xi,\tau)\in \mathbb{R}^n\times
\mathbb{R}\setminus(0,0), \tau>0
\end{aligned}
\end{equation}

\medskip

\noindent hence, by the homogeneity
of $G$ w.r.t. $(\xi,\tau)$ we get
\begin{equation}\label{note-6ter-50-28}
	\begin{aligned}
&P_m(x,\xi+i\tau\nabla\varphi(x))=0\mbox{ }\Rightarrow\\&
\Rightarrow G(x,\xi,\xi, \tau)\geq C\left(|\xi|^2+\tau^2\right)^{m-\frac{1}{2}},
\quad \forall x\in \overline{\Omega}
\end{aligned}
\end{equation}
for every $(\xi,\tau)\in \mathbb{R}^{n+1}$, $\tau>0$.

\medskip

\noindent\textbf{(b) Next steps.} Keeping in mind Lemma
\ref{stime-Carlm-1.3.3-39}, we will exploit the local character of a
Carleman estimate to focus on the case where
\textbf{the support of $v$ (and hence of $u$) is sufficiently small}.
With this expedient, the term $\mathcal{R}$ on the right-hand side in \eqref{note-6bis-50} can be treated as a kind of rest and can be efficiently estimated from below by simultaneously exploiting the continuity of the coefficients of the quadratic form $G(x,D,\overline{D},\tau)$ and Proposition \ref{note-52-prop}.
 \textbf{In coarse words}, where $P_m(x,\xi+i\tau\nabla\varphi(x))=0$ \eqref{note-6ter-50} is used and where
$P_m(x,\xi+i\tau\nabla\varphi(x))\neq 0$ (so where one will not be able to
exploit the property \eqref{note-6ter-50}) we will exploit the
specific character of the operator $P_m$

\bigskip

\textbf{Remarks.}

\smallskip

\noindent\textbf{1.} Let us notice that if
$P_m(x_0,\xi+i\tau\nabla\varphi(x_0))$ has some zero of multiplicity larger than $1$ in $(\xi_0,\tau_0)\neq (0,0)$, then \eqref{note-6ter-50} cannot be true. As a matter of fact in this case we have 

$$P^{(j)}_m(x_0,\xi_0+i\tau_0\nabla\varphi(x_0))=0, \quad\quad j=1,,\cdots,
n,$$ hence by \eqref{note-6-50} we have $G(x_0,\xi_0,\xi_0,
\tau_0)=0$.

\medskip

\noindent\textbf{2.} Let us notice that if the coefficients of $P_m(x,D)$
are constant, then we have

\begin{equation}\label{note-65}
\begin{aligned}
&\frac{i}{2\tau}\left\{P_m(x,\xi+i\tau\nabla\varphi(x)),\overline{P_m(x,\xi+i\tau\nabla\varphi(x))}\right\}=\\&=\frac{1}{\tau}
G(x,\xi,\xi,\tau)=\\&= \sum_{j,k=1}^n\partial^2_{x_jx_k}\varphi(x)
P_m^{(j)}(x,\xi+i\tau\nabla\varphi(x))\overline{P_m^{(k)}(x,\xi+i\tau\nabla\varphi(x))}.
\end{aligned}
\end{equation}

\noindent\textbf{3.} Operators \eqref{note-5-49a},
\eqref{note-5-49b} are the "first-order approximations"
respectively, of the symmetric and the antisymmetric parts
of the operator $p_m(x,D,\tau)$.  Let us examine this issue in more
detail. Let us suppose that the operator $P_m(x,D)$ has very regular coefficients, say $C^{\infty}$, and let us write $p_m(x,D,\tau)$ as follows

$$p_m(x,D,\tau)=\sum_{|\alpha|+j=
m}c_{\alpha,j}(x)\tau^jD^{\alpha}=\tau^m\sum_{|\alpha|+j=
m}c_{\alpha,j}(x)\left(\tau^{-1} D\right)^{\alpha}.$$ Let us consider
the formal adjoint of $p_m(x,D,\tau)$ i.e. the operator
$p^{\star}_m(x,D,\tau)$ such that

$$\int (p_m(x,D,\tau)v \overline{w}dx=\int v \left(\overline{p^{\star}_m(x,D,\tau)w}\right)dx,\quad\quad\forall u,w\in C^{\infty}_0(\Omega).$$
We have, integrating by parts,

$$p^{\star}_m(x,D,\tau)v=\tau^m\sum_{|\alpha|+j= m}\left(\tau^{-1}
D\right)^{\alpha}\left(\overline{c}_{\alpha,j}v\right)=\sum_{|\alpha|+j=
m}\tau^jD^{\alpha}\left(\overline{c}_{\alpha,j}v\right). $$
 On the other hand
$$\overline{p}_m(x,D,\tau)=\sum_{|\alpha|+j=
m}\overline{c}_{\alpha,j}\tau^jD^{\alpha}.$$ Hence we have

\begin{equation}\label{note-7-50}
\begin{aligned}
&\left(p^{\star}_m(x,D,\tau)-\overline{p}_m(x,D,\tau)\right)v=\\&=\frac{1}{i}\sum_{|\alpha|+j=
m-1}\tau^j\sum_{k=1}^n\binom{{\alpha}}{{e_k}}\partial_{x_k}c_{\alpha,j}(x)
D^{\alpha-e_k}v+\\&+r_{m-2}(x,D,\tau)v,
\end{aligned}
\end{equation}

\medskip

\noindent where
$$r_{m-2}(x,D,\tau)v=\sum_{|\alpha|+j\leq
m-2}\tau^j\widetilde{c}_{\alpha j}(x)D^{\alpha}v$$ and
$\widetilde{c}_{\alpha j}$ are suitable coefficients. By expressing
\eqref{note-7-50} by means of the symbols of the operator, we have 
\begin{equation}\label{note-8-50}
p^{\star}_m(x,\xi,\tau)=\overline{p_m(x,\xi,\tau)}+\frac{1}{i}\sum_{k=1}^np^{(j)}_{m,j}(x,\xi,\tau)+r_{m-2}(x,\xi,\tau).
\end{equation}
Relationship \eqref{note-7-50} and \eqref{note-8-50}  expresses in a precise manner that $\overline{p}_m(x,D,\tau)$ approximates
$p^{\star}_m(x,D,\tau)$ to the first order.

\medskip

\noindent\textbf{4.} It can be noticed that by spreading the square  in \eqref{note-1-50} in a standard way one would
leads to conclusions not unlike those seen above, in
particular, with regard to \eqref{note-6-50}. Here we give
a brief mention referring the interested reader to \cite[Ch.
4]{Lern}. We warn, however, that this approach requires
generally, assumptions of greater regularity on the coefficients of
$P_m(x,D)$ than we will make in this Chapter.

Set

\begin{equation*}
s(x,D,\tau)=\frac{1}{2}\left(p_m(x,D,\tau)+p^{\star}_m(x,D,\tau)\right),
\end{equation*}
\begin{equation*}
a(x,D,\tau)=\frac{1}{2}\left(p_m(x,D,\tau)-p^{\star}_m(x,D,\tau)\right).
\end{equation*}
We have trivially
\begin{equation}\label{note-50-2911}
p_m(x,D,\tau)=s(x,D,\tau)+a(x,D,\tau).
\end{equation}
Denoting by $p$, $s$ and $a$, respectively, $p_m(x,D,\tau)$,
$s(x,D,\tau)$, $a(x,D,\tau)$ and denoting by $\langle \cdot,\cdot\rangle$ the scalar
product in $L^2(\Omega, \mathbb{C})$, we have:

\begin{equation}\label{note-9-50}
\begin{aligned}
\left\Vert p(v)\right\Vert^2_{L^2(\Omega)}&=\langle
p(v),p(v)\rangle=\\&=\langle (s+a)(v),(s+a)(v)\rangle=\\&=\left\Vert
s(v)\right\Vert^2_{L^2(\Omega)}+\left\Vert
a(v)\right\Vert^2_{L^2(\Omega)}+2\Re\langle s(v),a(v)\rangle.
\end{aligned}
\end{equation}
Let us note that, denoting by
$$[s,a]=sa-as,$$ the \textbf{commutator of $a$ and $s$} \index{commutator of two operators}and taking into account that
$$s^{\star}=s,\quad\quad a^{\star}=-a,$$ we have
\begin{equation}\label{note-10-50}
\begin{aligned}
2\Re\langle s(v),a(v)\rangle&=\langle a(v),s(v)\rangle+\langle
s(v),a(v)\rangle=\\&=\langle s^{\star}a(v),v\rangle+\langle
a^{\star} s(v),v\rangle=\\&=\langle sa(v),v\rangle-\langle
as(v),v\rangle=\\&=\langle [s,a](v) ,v\rangle.
\end{aligned}
\end{equation}
Now, by \eqref{note-9-50} we have

\begin{equation}\label{note-11-50}
\begin{aligned}
\int\left\vert p_m(x,D,\tau)v\right\vert^2dx&=\int\left\vert
s(x,D,\tau)v\right\vert^2dx+\int\left\vert
a(x,D,\tau)v\right\vert^2dx+\\&+2\int\Re\left(s(x,D,\tau)v\overline{a(x,D,\tau)v}\right)dx.
\end{aligned}
\end{equation}
this, by \eqref{note-10-50}, can be written as
\begin{equation}\label{note-12-50}
\begin{aligned}
\int\left\vert p_m(x,D,\tau)v\right\vert^2dx&=\int\left\vert
s(x,D,\tau)v\right\vert^2dx+\int\left\vert
a(x,D,\tau)v\right\vert^2dx+\\&+2\int\left(\left[s(x,D,\tau),a(x,D,\tau)\right]v\right)\overline{v}dx.
\end{aligned}
\end{equation}
Now let us compare the integrals

$$2\int\Re\left(S(x,D,\tau)v\overline{A(x,D,\tau)v}\right)dx,\quad\quad
2\int\Re\left(s(x,D,\tau)v\overline{a(x,D,\tau)v}\right)dx$$ which 
occur, respectively, as the third term on the right hand side in
\eqref{note-1-50} and the third term on the right hand side in \eqref{note-11-50}.

 Set

$$R(x,D,\tau)=\frac{1}{2}\left(\frac{1}{i}\sum_{k=1}^np^{(j)}_{m,j}(x,D,\tau)+r_{m-2}(x,D,\tau)\right),$$
we have
$$s(x,D,\tau)=S(x,D,\tau)+R(x,D,\tau)\quad\mbox{ and }\quad
a(x,D,\tau)=S(x,D,\tau)-R(x,D,\tau)$$ then

\begin{equation*}
\begin{aligned}
2\Re\left(s(x,D,\tau)v\overline{a(x,D,\tau)v}\right)&=2\Re\left(S(x,D,\tau)v\overline{A(x,D,\tau)v}\right)-\\&-q(x,D,\overline{D},\tau)[v,\overline{v}]
-\left|R(x,D,\tau)v\right|^2,
\end{aligned}
\end{equation*}
where
\begin{equation*}
q(x,D,\overline{D},\tau)=
2\Re\left(S(x,D,\tau)v\overline{R(x,D,\tau)v}-R(x,D,\tau)v\overline{A(x,D,\tau)v}\right).
\end{equation*}
Recalling \eqref{note-5-49a} and \eqref{note-5-49b} we get

\begin{equation*}
\begin{aligned}
q(x,D,\overline{D},\tau)[v,\overline{v}]=q_1(x,D,\overline{D},\tau)[v,\overline{v}]+q_2(x,D,\overline{D},\tau)[v,\overline{v}],
\end{aligned}
\end{equation*}
where
\begin{equation*}
q_1(x,D,\overline{D},\tau)[v,\overline{v}]=
2\Re\left(p_m(x,D,\tau)v\overline{R(x,D,\tau)v}\right)
\end{equation*}
and
\begin{equation*}
	\begin{aligned}
&q_2(x,D,\overline{D},\tau)[v,\overline{v}]=\\&=
\Re\left(\overline{p}_m(x,D,\tau)v\overline{R(x,D,\tau)v}-R(x,D,\tau)v\left(\overline{\overline{p}_m(x,D,\tau)v}\right)\right)=0.
\end{aligned}
\end{equation*}
Therefore
\begin{equation*}
\begin{aligned}
&2\Re\left(s(x,D,\tau)v\overline{a(x,D,\tau)v}\right)=\\&=2\Re\left(S(x,D,\tau)v\overline{A(x,D,\tau)v}\right)-\\&-2\Re\left(p_m(x,D,\tau)v\overline{R(x,D,\tau)v}\right)
-\left|R(x,D,\tau)v\right|^2
\end{aligned}
\end{equation*}
and by \eqref{note-4-50} we have
\begin{equation}\label{note-13-50}
\begin{aligned}
&2\int\Re\left(s(x,D,\tau)v\overline{a(x,D,\tau)v}\right)dx=\\&=\int
G(x,D,\overline{D},\tau)[v,\overline{v}]dx-\\&-2\Re\int\left(p_m(x,D,\tau)v\overline{R(x,D,\tau)v}\right)dx-\\&-\int
\left|R(x,D,\tau)v\right|^2dx.
\end{aligned}
\end{equation}

\medskip

\noindent This relationship allows (see \textbf{Exercise} subsequent to the
proof of Theorem \ref{Carlm-ell-teor}) to consider
equivalent the approach we are following with the one outlined in
this Remark (of course, when the coefficients
of the operator are sufficiently regular). $\blacklozenge$

\bigskip

We conclude this Section with some lemma that will be useful later on. 
\begin{lem}\label{lemma-carlm-ell-I}
Let $\varphi\in C^{\infty}\left(\overline{\Omega}\right)$. Then for every
$m\in \mathbb{N}_0$ there exists a constant $C>1$ such that
\begin{equation}\label{lemma-carlm-ell-I-1}
\begin{aligned}
C^{-1}\sum_{|\alpha|\leq m} \tau^{2(m-|\alpha|)}\left\vert
D^{\alpha}\left(e^{\tau\varphi}u\right)\right\vert^2 &\leq
\sum_{|\alpha|\leq m}\tau^{2(m-|\alpha|)} \left\vert
D^{\alpha}u\right\vert^2e^{2\tau\varphi}\leq\\&\leq
C\sum_{|\alpha|\leq m} \tau^{2(m-|\alpha|)}\left\vert
D^{\alpha}\left(e^{\tau\varphi}u\right)\right\vert^2,
\end{aligned}
\end{equation}
for every $u\in C^{\infty}\left(\overline{\Omega}\right)$ and for every
$\tau\geq 1$.
\end{lem}
\textbf{Proof.} Both the inequalities are proved
easily by means of Leibniz formula. Here we limit ourselves to
prove 

\begin{equation}\label{lemma-carlm-ell-I-1-23-4-23}
	\begin{aligned}
		 \sum_{|\alpha|\leq m}\tau^{2(m-|\alpha|)} \left\vert
		D^{\alpha}u\right\vert^2e^{2\tau\varphi}\leq
		C\sum_{|\alpha|\leq m} \tau^{2(m-|\alpha|)}\left\vert
		D^{\alpha}\left(e^{\tau\varphi}u\right)\right\vert^2.
	\end{aligned}
\end{equation}

We use the induction principle. If $m=0$, then \eqref{lemma-carlm-ell-I-1-23-4-23} is trivial. Let us suppose that

\begin{equation*}\label{lemma-carlm-ell-I-2}
\begin{aligned}
\sum_{|\alpha|\leq m}\tau^{2(m-|\alpha|)} \left\vert
D^{\alpha}u\right\vert^2e^{2\tau\varphi}\leq C_m\sum_{|\alpha|\leq
m} \tau^{2(m-|\alpha|)}\left\vert
D^{\alpha}\left(e^{\tau\varphi}u\right)\right\vert^2,
\end{aligned}
\end{equation*}
where $C_m\geq 1$ and we have

\begin{equation}\label{lemma-carlm-ell-I-3}
\begin{aligned}
\sum_{|\alpha|\leq m+1} \tau^{2(m+1-|\alpha|)}\left\vert
D^{\alpha}\left(e^{\tau\varphi}u\right)\right\vert^2&=\sum_{|\alpha|=
m+1} \left\vert
D^{\alpha}\left(e^{\tau\varphi}u\right)\right\vert^2+\\&+
\tau^2\sum_{|\alpha|\leq m} \tau^{2(m-|\alpha|)}\left\vert
D^{\alpha}\left(e^{\tau\varphi}u\right)\right\vert^2\geq\\&\geq
\sum_{|\alpha|= m+1} \left\vert
D^{\alpha}\left(e^{\tau\varphi}u\right)\right\vert^2+\\&+C_m^{-1}\sum_{|\alpha|\leq
m}\tau^{2(m+1-|\alpha|)} \left\vert
D^{\alpha}u\right\vert^2e^{2\tau\varphi}.
\end{aligned}
\end{equation}
Let $\delta\in (0,1)$ be to choose. Using the Leibniz formula
 we have, for $\tau\geq 1$,

\begin{equation}\label{lemma-carlm-ell-I-4}
\begin{aligned}
\sum_{|\alpha|= m+1} \left\vert
D^{\alpha}\left(e^{\tau\varphi}u\right)\right\vert^2&\geq \delta
\sum_{|\alpha|= m+1} \left\vert
D^{\alpha}\left(e^{\tau\varphi}u\right)\right\vert^2\geq\\&\geq
\delta \sum_{|\alpha|= m+1} \left\vert
D^{\alpha}u\right\vert^2e^{2\tau\varphi}-\\&-\delta
\widetilde{C}_m\sum_{|\alpha|\leq m}\tau^{2(m+1-|\alpha|)}
\left\vert D^{\alpha}u\right\vert^2e^{2\tau\varphi},
\end{aligned}
\end{equation}
where $\widetilde{C}_m\geq 1$ is a suitable constant depending on $m$. By
\eqref{lemma-carlm-ell-I-3} and \eqref{lemma-carlm-ell-I-4} we get

\begin{equation*}\label{lemma-carlm-ell-I-5}
\begin{aligned}
\sum_{|\alpha|\leq m+1}\tau^{2(m+1-|\alpha|)} \left\vert
D^{\alpha}\left(e^{\tau\varphi}u\right)\right\vert^2&\geq \delta
\sum_{|\alpha|=m+1} \left\vert D^{\alpha}u\right\vert^2e^{2\tau\varphi}+\\&+
\left(C_m^{-1}-\delta \widetilde{C}_m\right) \sum_{|\alpha|\leq
m}\tau^{2(m+1-|\alpha|)} \left\vert
D^{\alpha}u\right\vert^2e^{2\tau\varphi}.
\end{aligned}
\end{equation*}
Now, we choose $\delta=\frac{1}{2C_m\widetilde{C}_m}$ and we get
\begin{equation*}
\begin{aligned}
\sum_{|\alpha|\leq m+1}\tau^{2(m+1-|\alpha|)} \left\vert
D^{\alpha}u\right\vert^2e^{2\tau\varphi}\leq 2C_m\sum_{|\alpha|\leq
m+1}  \tau^{2(m+1-|\alpha|)}\left\vert
D^{\alpha}\left(e^{\tau\varphi}u\right)\right\vert^2,\end{aligned}
\end{equation*}
which concludes the proof.$\blacksquare$

\bigskip

\begin{lem}\label{lemma-carlm-ell-II}
Let $\varphi\in C^{\infty}\left(\overline{\Omega}\right)$. Then for each
$m\in \mathbb{N}_0$ there exists a constant $C>1$ such that for every
$v\in C_0^{\infty}\left(\Omega\right)$ and for every $\tau\geq 1$ we have
\begin{equation}\label{lemma-carlm-ell-II-1}
\begin{aligned}
C^{-1}\sum_{|\alpha|\leq m}\tau^{2(m-|\alpha|)}\int \left\vert
D^{\alpha}v\right\vert^2dx & \leq \int
\left(|\xi|^2+\tau^2\right)^{m}\left|\widehat{v}(\xi)\right|^2d\xi\leq\\&\leq
C\sum_{|\alpha|\leq m}\tau^{2(m-|\alpha|)} \int\left\vert
D^{\alpha}v\right\vert^2dx.
\end{aligned}
\end{equation}
\end{lem}

\textbf{Proof.} We start with the first inequality in
\eqref{lemma-carlm-ell-II-1}. By Lemma \ref{lemma-carlm-ell-I} and the
Parseval identity we have

\begin{equation*}
\begin{aligned}
\sum_{|\alpha|\leq m} \int\tau^{2(m-|\alpha|)}\left\vert
D^{\alpha}v\right\vert^2dx&=\frac{1}{(2\pi)^n}\int
\sum_{|\alpha|\leq m}\tau^{2(m-|\alpha|)}\left\vert
\xi^{\alpha}\right\vert^2\left|\widehat{v}(\xi)\right|^2d\xi\leq\\&\leq
C\int
\left(|\xi|^2+\tau^2\right)^{m}\left|\widehat{v}(\xi)\right|^2d\xi.
\end{aligned}
\end{equation*}
Concerning the  second inequality in
\eqref{lemma-carlm-ell-I-1}, we have similarly
\begin{equation*}
\begin{aligned}
\int
\left(|\xi|^2+\tau^2\right)^{m}\left|\widehat{v}(\xi)\right|^2d\xi
&=\sum_{k=0}^m\binom{{m}}{{k}}\tau^{2(m-k)}\int|\xi|^{2k}\left|\widehat{v}(\xi)\right|^2d\xi\leq\\&\leq
2^m\sum_{k=0}^m\tau^{2(m-k)}\int\sum_{|\alpha|=k}\left|\xi^{\alpha}\right|^2\left|\widehat{v}(\xi)\right|^2d\xi=\\&=2^m\int\sum_{|\alpha|\leq
m}\tau^{2(m-|\alpha|)}\left\vert
\xi^{\alpha}\right\vert^2\left|\widehat{v}(\xi)\right|^2d\xi=\\&=
2^m(2\pi)^n\sum_{|\alpha|\leq m}\tau^{2(m-|\alpha|)} \int\left\vert
D^{\alpha}v\right\vert^2dx.
\end{aligned}
\end{equation*}
$\blacksquare$

\bigskip

\begin{lem}\label{lemma-carlm-ell-III}
Let us assume that the coefficients of operator
\eqref{sec-coniugato-1bis} belong to
$C^{0}\left(\overline{\Omega}\right)$. Let $\varphi\in
C^{\infty}\left(\overline{\Omega}\right)$ and let $p_m(x,D,\tau)$
be the operator defined by \eqref{note-1-49}. Then for every 
$\varepsilon>0$ there exists $\delta>0$ such that
\begin{equation}\label{lemma-carlm-ell-III-1}
\begin{aligned}
\left\vert\frac{1}{(2\pi)^n}\int
\left|p_m(x_0,\xi,\tau)\right|^2\left|\widehat{v}(\xi)\right|^2d\xi-\int
\left|p_m(x,D,\tau)v\right|^2dx\right\vert\leq &\\ \leq\varepsilon
\sum_{|\alpha|\leq m}\tau^{2(m-|\alpha|)}\int \left\vert
D^{\alpha}v\right\vert^2dx,
\end{aligned}
\end{equation}
for every $v\in C_0^{\infty}(B_{\delta}(x_0)\cap \Omega)$, for every
$\tau\in \mathbb{R}$ and for every $x_0\in \overline{\Omega}$.
\end{lem}

\textbf{Proof. } Let $x_0\in \overline{\Omega}$. Let recall that
by \eqref{note-2-49} we have

\begin{equation}\label{lemma-carlm-ell-III-2}
p_m(x,D,\tau)=\sum_{|\alpha|+j=m}\tau^jb_{\alpha j}(x) D^{\alpha},
\end{equation}
by the assumptions on $\varphi$ and on the coefficients of $P_m(x,D)$, we have \\
$b_{\alpha j}\in C^0\left(\overline{\Omega}\right)$, for any $\alpha$ and
$j$ such that $|\alpha|+j=m$.

Let $\varepsilon>0$ and $\delta>0$ be such that for any $\alpha$ and $j$
satisfying $|\alpha|+j=m$ we have

$$\left|b_{\alpha j}(x)-b_{\alpha j}(x_0)\right|<\varepsilon, \quad
\forall x\in B_{\delta}(x_0)\cap\overline{\Omega}$$ ($\delta$
indipendent of $x_0$). We obtain 

\begin{equation}\label{lemma-carlm-ell-III-3}
\begin{aligned}
&\left\vert p_m(x,D,\tau)v-p_m(x_0,D,\tau)v\right\vert\leq\\&\leq
\sum_{|\alpha|+j=m}|\tau|^j\left\vert b_{\alpha j}(x)-b_{\alpha
j}(x_0) \right\vert\left| D^{\alpha}v\right|\leq \\&\leq
C\varepsilon\sum_{|\alpha|+j=m}|\tau|^j\left| D^{\alpha}v\right|, \ \ \forall x\in B_{\delta}(x_0)\cap\overline{\Omega}.
\end{aligned}
\end{equation}
On the other hand
\begin{equation}\label{lemma-carlm-ell-III-4}
\left\vert p_m(x,D,\tau)v\right\vert\leq C
\sum_{|\alpha|+j=m}|\tau|^j\left| D^{\alpha}v\right|,\quad\forall x\in
\overline{\Omega}
\end{equation}
Now, taking into account the elementary inequality $$\left||z|^2-|w|^2\right|\leq
\left(|z|+|w|\right)|z-w|,\quad\forall z,w\in \mathbb{C},$$ we have
by \eqref{lemma-carlm-ell-III-3} and \eqref{lemma-carlm-ell-III-4}, for every $x\in B_{\delta}(x_0)\cap\overline{\Omega}$,

\begin{equation}\label{lemma-carlm-ell-III-5}
\begin{aligned}
\left\vert\left\vert p_m(x,D,\tau)v\right\vert^2-\left\vert
p_m(x_0,D,\tau)v\right\vert^2\right\vert\leq C\varepsilon
\sum_{|\alpha|\leq m}\tau^{2(m-|\alpha|)}\left|
D^{\alpha}v\right|^2.
\end{aligned}
\end{equation}
Therefore, for every $v\in C_0^{\infty}(B_{\delta}(x_0)\cap \Omega)$

\begin{equation}\label{lemma-carlm-ell-III-6}
\begin{aligned}
&\frac{1}{(2\pi)^n}\int
\left|p_m(x_0,\xi,\tau)\right|^2\left|\widehat{v}(\xi)\right|^2d\xi-\int
\left|p_m(x,D,\tau)v\right|^2dx =\\&=\int
\left(\left|p_m(x_0,D,\tau)v\right|^2-\left|p_m(x,D,\tau)v\right|^2\right)dx\leq
\\&
\leq C\varepsilon \sum_{|\alpha|\leq
m}\tau^{2(m-|\alpha|)}\int\left| D^{\alpha}v\right|^2dx,
\end{aligned}
\end{equation}
and similarly, for every $v\in C_0^{\infty}(B_{\delta}(x_0)\cap \Omega)$,

\begin{equation}\label{lemma-carlm-ell-III-7}
\begin{aligned}
\frac{1}{(2\pi)^n}\int
\left|p_m(x_0,\xi,\tau)\right|^2&\left|\widehat{v}(\xi)\right|^2d\xi-\int
\left|p_m(x,D,\tau)v\right|^2dx \geq\\&\geq -C\varepsilon
\sum_{|\alpha|\leq m}\tau^{2(m-|\alpha|)}\int\left|
D^{\alpha}v\right|^2dx.
\end{aligned}
\end{equation}
Finally
\eqref{lemma-carlm-ell-III-6} and \eqref{lemma-carlm-ell-III-7} implies
\eqref{lemma-carlm-ell-III-1}. $\blacksquare$

\section{Carleman estimates for the elliptic operators}\label{Carleman-ellittico} Let $m\in \mathbb{N}$,
and let $\Omega$ be a bounded open set of $\mathbb{R}^n$. Let $a_{\alpha}$
complex--valued functions. We recall that the operator
\begin{equation}\label{Carleman-ellittico-1}
P(x,D)= \sum_{|\alpha|\leq m}a_{\alpha}(x)D^{\alpha},
\end{equation}
is elliptic in a point $x_0$ if
\begin{equation}\label{Carleman-ellittico-2}
P_m(x_0,\xi)= \sum_{|\alpha|= m}a_{\alpha}(x_0)\xi^{\alpha}\neq
0,\quad\quad \forall \xi\in \mathbb{R}^n\setminus\{0\}.
\end{equation}
We also say that $P(x,D)$ is elliptic in $\overline{\Omega}$ if
\eqref{Carleman-ellittico-2} holds for every $x_0\in \Omega$. Let us note that if $a_{\alpha}\in
C^0\left(\overline{\Omega},\mathbb{C}\right)$, for $|\alpha|=m$, the ellipticity condition  for the operator $P(x,D)$ is equivalent to the existence of a constant 
$\lambda>0$ such that
\begin{equation}\label{Carleman-ellittico-2bis}
\left\vert P_m(x,\xi)\right\vert\geq \lambda \left\vert
\xi\right\vert^{m},\quad\quad \forall \xi\in \mathbb{R}^n, \quad
\forall x\in \overline{\Omega}.
\end{equation}

\medskip

For the sake of brevity, in the proof of Theorem below, for an
open $\omega\subset \Omega$ we will identify $C_0^{\infty}(\omega)$
with the function space

$$\left\{u\in C_0^{\infty}(\Omega):\quad \mbox{supp }u\subset \omega\right\}.$$

\medskip

\begin{theo}[\textbf{Carleman--H\"{o}rmander}]\label{Carlm-ell-teor}
	\index{Theorem:@{Theorem:}!- Carlemann--H\"{o}rmander@{- Carleman--H\"{o}rmander}}
Let $\varphi\in C^{\infty}\left(\overline{\Omega}\right)$ be a 
real--valued function which satisfies
\begin{equation}\label{Carleman-ellittico-3}
\nabla \varphi(x)\neq 0,\quad\quad\forall x\in \overline{\Omega}.
\end{equation}
Let $P(x,D)$ be an operator of order $m$ whose coefficients belong to
$L^{\infty}(\Omega,\mathbb{C})$. Let us assume that the coefficients of the principal part  $P_m(x,D)$ belong to
$C^1\left(\overline{\Omega},\mathbb{C}\right)$. Let us suppose that
$P(x,D)$ satisfies the ellipticity condition \eqref{Carleman-ellittico-2bis} and that the following condition is satisfied:

\medskip

\noindent ($\bigstar$) If

\begin{equation}\label{Carleman-ellittico-4}
\begin{aligned}
\begin{cases}
P_m(x,\xi+i\sigma\nabla\varphi(x))=0,\\
\\
x\in \overline{\Omega},\\
\\
(\xi,\sigma)\in \mathbb{R}^{n+1}\setminus\{(0,0)\},
\end{cases}
\end{aligned}
\end{equation}
then
\begin{equation}\label{Carleman-ellittico-5}\frac{i}{2\sigma}\left\{P_m(x,\xi+i\sigma\nabla\varphi(x)),\overline{P_m(x,\xi+i\sigma\nabla\varphi(x))}\right\}>0,
\end{equation}
where $\{\cdot,\cdot\}$ is the Poisson bracket defined in
\eqref{note-parentesi-p}.

\medskip

Then there exist constants  $C$ and $\tau_0$  such that

\begin{equation}\label{Carleman-ellittico-6}
\sum_{|\alpha|\leq m}\tau^{2(m-|\alpha|)-1}\int \left\vert
D^{\alpha}u\right\vert^2e^{2\tau\varphi}dx\leq C\int \left\vert
P(x,D)u\right\vert^2e^{2\tau\varphi}dx,
\end{equation}
for every $u\in C_0^{\infty}(\Omega)$ and for every $\tau\geq \tau_0$.

Moreover $C$ and $\tau_0$ depend on $\lambda$, on 
the $L^{\infty}(\Omega,\mathbb{C})$ norms of $a_{\alpha}$, $|\alpha|\leq
m$, on the $L^{\infty}(\Omega,\mathbb{C})$ norms of $\nabla a_{\alpha}$, $|\alpha|= m$, and on the moduli of
continuity  of $\nabla a_{\alpha}$, for $|\alpha|= m$.
\end{theo}

\medskip
\textbf{Remark 1.} Let us notice that requiring that
$(\xi,\sigma)\neq (0,0)$ in \eqref{Carleman-ellittico-4} is equivalent
to require that both $\xi$ and $\sigma$ are different from zero.
As a matter of fact, if $\xi=0$ then, by $P_m(x,\xi+i\sigma\nabla\varphi(x))=0$,
we have $(i\sigma)^mP_m(x,\nabla\varphi(x))=0$ in addition, since
$P_m(x,\nabla\varphi(x))\neq 0$ and since $P_m(x,D)$ is elliptic and
$\nabla\varphi(x)\neq 0$, we have $\sigma=0$. Similarly, if $\sigma=0$ by
the ellipticity of $P_m(x,D)$ we have $\xi=0$. $\blacklozenge$

\medskip

\textbf{Remark 2.} Taking into account Remark 2 of
Section \ref{sec-coniugato}, if the coefficients of $P_m(x,D)$
are constants (let us rename it $P_m(D)$),  condition ($\bigstar$)
become:

\medskip

\noindent(($\bigstar$)\textbf{ -- constant coefficients})

If

\begin{equation}\label{Carleman-ellittico-4bis}
\begin{aligned}
\begin{cases}
P_m(\xi+i\sigma\nabla\varphi(x))=0,\\
\\
(\xi,\sigma)\in \mathbb{R}^{n+1}\setminus\{(0,0)\},
\end{cases}
\end{aligned}
\end{equation}
then
\begin{equation*}
\sum_{j,k=1}^n\partial^2_{x_jx_k}\varphi(x)
P_m^{(j)}(\xi+i\tau\nabla\varphi(x))\overline{P_m^{(k)}(\xi+i\tau\nabla\varphi(x))}>0.
\end{equation*}
In particular, if the Hessian matrix of $\varphi$ is positive definite then  ($\bigstar$)\textbf{ -- constant coefficients} is satisfied. $\blacklozenge$

\bigskip

\textbf{Proof of Theorem \ref{Carlm-ell-teor}.}

Let $u\in C_0^{\infty}(\Omega)$, set
$$v=e^{-\tau\varphi}u.$$
As observed in the previous Section, we have

\begin{equation}\label{dim-1}
e^{\tau\varphi}P_m(x,D)u=e^{\tau\varphi}P_m(x,D)\left(e^{-\tau\varphi}v\right)=P_m(x,D+i\tau\nabla\varphi(x))v
\end{equation}
and, denoting by $p_m(x,D,\tau)$ the operator whose symbol is
$P_m(x,\xi+i\tau\nabla \varphi(x))$, by \eqref{note-4-49} we get

\begin{equation}\label{dim-2}
\begin{aligned}
&\int\left\vert P_m(x,D+i\tau\nabla\varphi(x))v\right\vert^2dx\geq\\&\geq 
\frac{1}{2}\int\left\vert
p_m(x,D,\tau)v\right\vert^2dx-\\&-C_1\sum_{|\alpha|\leq
m-1}\tau^{2(m-|\alpha|)-2}\int\left\vert D^{\alpha}v\right\vert^2dx,
\end{aligned}
\end{equation}
where $C_1$ depends by the $L^{\infty}$ norms of the coeifficients of
$P_m(x,D)$.

We now derive an appropriate estimate from below of the first term on the right--hand side
in \eqref{dim-2}. 

By \eqref{note-1-50} and \eqref{note-4-50} we get
(by multiplying both equalities by $\tau$)

\begin{equation}
\begin{aligned}\label{dim-3}
\tau\int\left\vert p_m(x,D,\tau)v\right\vert^2dx&\geq
2\tau\int\Re\left(S(x,D,\tau)v\overline{A(x,D,\tau)v}\right)dx=\\&=\tau\int
G(x,D,\overline{D},\tau)[v,\overline{v}]dx,
\end{aligned}
\end{equation}
where $G(x,D,\overline{D},\tau)$ has been  defined in Proposition
\ref{note-52-prop}.

Let now $x_0\in\overline{\Omega}$ be a fixed point. We may assume $x_0=0\in
\overline{\Omega}$. We get

\begin{equation}\label{dim-4}
\begin{aligned}
&\tau\int G(x,D,\overline{D},\tau)[v,\overline{v}]dx=\tau\int
G(0,D,\overline{D},\tau)[v,\overline{v}]dx+\\&+\tau\int
\left(G(x,D,\overline{D},\tau)-G(0,D,\overline{D},\tau)\right)[v,\overline{v}]dx=\\&=
\frac{\tau}{(2\pi)^n}\int
G(0,\xi,\xi,\tau)\left|\widehat{v}(\xi)\right|^2d\xi+\tau
\mathcal{R},
\end{aligned}
\end{equation}
where

\begin{equation*}
\mathcal{R}= \int
\left(G(x,D,\overline{D},\tau)-G(0,D,\overline{D},\tau)\right)[v,\overline{v}]dx.
\end{equation*}
By \eqref{note-52-prop-2} we have
\begin{equation}\label{dim-5}
G\left(x,D,\overline{D},\tau \right)=\sum_{h=0}^{2m-1}\tau^{h}
G^{(h)}\left(x,D,\overline{D}\right),
\end{equation}
where
$$G^{(h)}\left(x,D,\overline{D}\right)[v,\overline{v}]=\sum_{(\alpha,\beta)\in\Lambda_h}c^{(h)}_{\alpha\beta}(x)D^{\alpha}v\overline{D^{\beta}v}$$
and $$\Lambda_h=\left\{(\alpha,\beta)\in \mathbb{N}_0^n:\quad
|\alpha|\leq m, \mbox{ } |\beta|\leq m, \mbox{ }|\alpha|+|\beta|\leq
2m-h-1\right\},$$ for $h=0,1 \cdots, 2m-1$ and, further,
$c^{(h)}_{\alpha\beta}\in
C^0\left(\overline{\Omega},\mathbb{C}\right)$ for
$(\alpha,\beta)\in\Lambda_h$. Let $\varepsilon$ be a positive number
that we will choose later and let $\rho_1>0$ be such that

$$\left\vert
c^{(h)}_{\alpha\beta}(x)-c^{(h)}_{\alpha\beta}(0)\right\vert<\varepsilon,\quad
\forall x\in B_{\rho_1}\cap \overline{\Omega},\mbox{ }
(\alpha,\beta)\in\Lambda_h,\mbox{ } h=0,1 \cdots, 2m-1.$$ We have,
for every $\tau\geq 1$ and for every $x\in B_{\rho_1}\cap \overline{\Omega}$,

\begin{equation*}
	\begin{aligned}
&\left\vert\tau^{h+1}\left(G^{(h)}\left(x,D,\overline{D}\right)-G^{(h)}\left(0,D,\overline{D}\right)\right)[v,\overline{v}]\right\vert\leq
\\&
\leq \sum_{(\alpha,\beta)\in\Lambda_h}\tau^{h+1}\left\vert
c^{(h)}_{\alpha\beta}(x)-c^{(h)}_{\alpha\beta}(0)\right\vert
 \left\vert D^{\alpha}v\right\vert
\left\vert\overline{D^{\beta}v}\right\vert\leq\\& \leq \varepsilon
\sum_{(\alpha,\beta)\in\Lambda_h}\tau^{2m-(|\alpha|+|\beta|)}\left\vert
D^{\alpha}v\right\vert \left\vert D^{\beta}v\right\vert=\\&
= \varepsilon\sum_{(\alpha,\beta)\in\Lambda_h}
\left(\tau^{m-|\alpha|}\left\vert
D^{\alpha}v\right\vert\right)\left(\tau^{m-|\beta|}\left\vert
D^{\beta}v\right\vert\right)\leq \\&\leq C\varepsilon \sum_{|\alpha|\leq
m}\tau^{2(m-|\alpha|)}\left\vert D^{\alpha}v\right\vert^2,
	\end{aligned}
\end{equation*}
where C depends on $m$ only. Let us notice that in the second inequality we have
exploited that, for $\tau\geq 1$,
$$(\alpha,\beta)\in\Lambda_h\Rightarrow h+1\leq
2m-(|\alpha|+|\beta|)\Rightarrow \tau^{h+1}\leq
\tau^{2m-(|\alpha|+|\beta|)}.$$ Therefore, for every  $x\in B_{\rho_1}\cap \overline{\Omega}$,

\begin{equation*}
\left\vert\tau\left( G\left(x,D,\overline{D},\tau
\right)-G\left(0,D,\overline{D},\tau
\right)\right)[v,\overline{v}]\right\vert\leq C\varepsilon
\sum_{|\alpha|\leq m}\tau^{2(m-|\alpha|)}\left\vert
D^{\alpha}v\right\vert^2.
\end{equation*}
Hence, by Lemma \ref{lemma-carlm-ell-II} we have for any $\tau\geq 1$,
\begin{equation}\label{dim-6}
	\begin{aligned}
\left|\tau\mathcal{R}\right|&\leq C\varepsilon\sum_{|\alpha|\leq
m}\tau^{2(m-|\alpha|)}\int \left\vert D^{\alpha}v\right\vert^2dx\leq \\&\leq
C\varepsilon\int \left(|\xi|^2+\tau^2
\right)^{m}\left|\widehat{v}(\xi)\right|^2d\xi,
\end{aligned}
\end{equation}
for every $v\in C_0^{\infty}(B_{\rho_1}\cap \Omega)$,

Now, by \eqref{dim-3}, \eqref{dim-4} and \eqref{dim-6} we get

\begin{equation}
\begin{aligned}\label{dim-7}
\tau\int\left\vert p_m(x,D,\tau)v\right\vert^2dx&\geq
\frac{\tau}{(2\pi)^n}\int
G(0,\xi,\xi,\tau)\left|\widehat{v}(\xi)\right|^2d\xi-\\&-C\varepsilon\int
\left(|\xi|^2+\tau^2 \right)^{m}\left|\widehat{v}(\xi)\right|^2d\xi,
\end{aligned}
\end{equation}
for every $v\in C_0^{\infty}(B_{\rho_1}\cap\Omega)$ and for every
$\tau\geq 1$.

\bigskip

Now we prove the following

\medskip

\noindent \textbf{Claim.}

Set
$$N=\nabla\varphi(0),$$
there exist two positive constants $C_1$ and $C_2$ such that

\begin{equation}\label{dim-8}
	\begin{aligned}
		&C_1\left|\xi+i\sigma N\right|^{2m}\leq \\&\leq \sigma
G(0,\xi,\xi,\sigma)+C_2\left|P_m\left(0,\xi+i\sigma
N\right)\right|^2,\mbox{ } \forall (\xi,\sigma)\in \mathbb{R}^{n+1}.
\end{aligned}
\end{equation}

\textbf{Proof of the Claim.} Let us denote 
$$\mathbb{S}^n=\left\{(\xi,\sigma)\in\mathbb{R}^{n+1}:\quad \left|\xi+i\sigma N\right|=1
\right\}$$
and
$$\eta=\frac{\xi}{\left|\xi+i\sigma N\right|},\quad\quad \mu=\frac{\sigma}{\left|\xi+i\sigma
N\right|},$$ 
let us note that \eqref{note-6-50} implies, by
homogeneity, that condition $(\bigstar)$ is equivalent to the following one

\medskip

\noindent($\bigstar'$) If

\begin{equation}\label{Carleman-ellittico-4bis2}
	\begin{aligned}
		\begin{cases}
			P_m(0,\eta+i\mu N)=0,\\
			\\
			(\eta,\mu)\in \mathbb{S}^n,
		\end{cases}
	\end{aligned}
	\end{equation}
then (recall Remark 1)
\begin{equation*}
\mu G(0,\eta,\eta,\mu)>0.
\end{equation*} Furthermore, \eqref{dim-8}
is equivalent to
\begin{equation}\label{dim-9}
C_1\leq \mu G(0,\eta,\eta,\mu)+C_2\left|P_m\left(0,\eta+i\mu
N\right)\right|^2,\mbox{ } \forall (\eta,\mu)\in \mathbb{S}^n.
\end{equation}

\medskip

Now, by \eqref{Carleman-ellittico-2bis}, we have
$$\left(\mu G(0,\eta,\eta,\mu)+\left|P_m\left(0,\eta+i\mu
N\right)\right|^2\right)_{|\mu=0}\geq \lambda,\quad \mbox{for } |\eta|=1.$$ By the compactness of
$\mathbb{S}^n$ there exists $\mu_0>0$ such that 

\begin{equation}\label{dim-10}
\mu G(0,\eta,\eta,\mu)+\left|P_m\left(0,\eta+i\mu
N\right)\right|^2\geq\frac{\lambda}{2},
\end{equation}
for every $(\eta,\mu)\in \mathbb{S}^n\cap\{|\mu|\leq \mu_0\}$.
\medskip

\noindent Let us denote by $K$ the compact set
$$K=\mathbb{S}^n\cap\{|\mu|\geq \mu_0\}$$ (of course, if $K=\emptyset$
the proof would be concluded). Since  ($\bigstar'$) gives
trivially

$$\left|P_m\left(0,\eta+i\mu
N\right)\right|=0,\mbox{ } (\eta,\mu)\in K\Longrightarrow \mu
G(0,\eta,\eta,\mu)>0,$$ by Lemma \ref{106-CauNire} we have that
there exists $C>0$ such that 
\begin{equation}\label{dim-11}
\mu G(0,\eta,\eta,\mu)+C\left|P_m\left(0,\eta+i\mu
N\right)\right|^2>0,\quad \forall
(\eta,\mu)\in K.
\end{equation}
By \eqref{dim-10} and \eqref{dim-11} we have

\begin{equation}\label{dim-12}
\mu G(0,\eta,\eta,\mu)+(C+1)\left|P_m\left(0,\eta+i\mu
N\right)\right|^2>0,\quad \forall (\eta,\mu)\in \mathbb{S}^n
\end{equation}
and \eqref{dim-9} follows with
$$C_1=\min_{(\eta,\mu)\in
\mathbb{S}^n}\left(G(0,\eta,\eta,\mu)+(C+1)\left|P_m\left(0,\eta+i\mu
N\right)\right|^2\right)$$ and $$C_2=C+1.$$ The proof
of the Claim is concluded.

\medskip

Now, we set
$$\gamma=\min\left\{1, \min_{\overline{\Omega}}|\nabla
\varphi|\right\},$$using \eqref{dim-8} in \eqref{dim-7}
we get

\begin{equation}
\begin{aligned}\label{dim-13}
\tau\int\left\vert p_m(x,D,\tau)v\right\vert^2dx&\geq
(2\pi)^{-n}C_1\gamma^2\int \left(|\xi|^2+\tau^2
\right)^{m}\left|\widehat{v}(\xi)\right|^2d\xi-\\&-(2\pi)^{-n}C_2\int\left|p_m\left(0,\xi,\tau
\right)\right|^2\left|\widehat{v}(\xi)\right|^2d\xi-\\&-C\varepsilon\int
\left(|\xi|^2+\tau^2 \right)^{m}\left|\widehat{v}(\xi)\right|^2d\xi,
\end{aligned}
\end{equation}
for every $v\in C_0^{\infty}(B_{\rho_1}\cap\Omega)$ and for every
$\tau\geq 1$.  By Lemma \ref{lemma-carlm-ell-III} there
exists $\rho_2\leq \rho_1$ such that for every $v\in
C_0^{\infty}(B_{\rho_2}\cap \Omega)$ and for every $\tau\geq 1$ we have

\begin{equation}\label{dim-14}
\begin{aligned}
(2\pi)^{-n}\int
\left|p_m(0,\xi,\tau)\right|^2\left|\widehat{v}(\xi)\right|^2d\xi\leq
&\int \left|p_m(x,D,\tau)v\right|^2dx+\\& +C\varepsilon \int
\left(|\xi|^2+\tau^2 \right)^{m}\left|\widehat{v}(\xi)\right|^2d\xi.
\end{aligned}
\end{equation}
By \eqref{dim-13} and \eqref{dim-14} we have

\begin{equation}
\begin{aligned}\label{dim-15}
&\tau\int\left\vert p_m(x,D,\tau)v\right\vert^2dx\geq\\&\geq
\left((2\pi)^{-n}C_1\gamma^2-C\varepsilon\right)\int
\left(|\xi|^2+\tau^2
\right)^{m}\left|\widehat{v}(\xi)\right|^2d\xi-\\&-(2\pi)^{-n}C_2\int\left|p_m\left(x,D,\tau
\right)v\right|^2dx,
\end{aligned}
\end{equation}
for every $v\in
C_0^{\infty}(B_{\rho_2}\cap \Omega)$ . Now, let us choose
$$\varepsilon=\varepsilon_0:=\frac{(2\pi)^{-n}C_1\gamma^2}{2C}$$ and let us denote by $\overline{\rho}$ the value of $\rho_2$ when
$\varepsilon=\varepsilon_0$. Moving the last integral of \eqref{dim-15} to the left--hand side and  recalling Lemma \ref{lemma-carlm-ell-II}, we have

\begin{equation}
\begin{aligned}\label{dim-16}
C_3\tau\int\left\vert p_m(x,D,\tau)v\right\vert^2dx&\geq
\varepsilon_0C^{-1}\sum_{|\alpha|\leq m}\tau^{2(m-|\alpha|)}\int
\left\vert D^{\alpha}v\right\vert^2dx,
\end{aligned}
\end{equation}
($C_3=1+(2\pi)^{-n}C_2$) for every $v\in
C_0^{\infty}(B_{\overline{\rho}}\cap \Omega)$, for every $\tau\geq
1$.

At this point we use \eqref{note-4-49} and we have

\begin{equation}
\begin{aligned}\label{dim-17}
&C\sum_{|\alpha|\leq m-1}\tau^{2(m-|\alpha|)-1}\int\left\vert
D^{\alpha}v\right\vert^2dx+\\&+2C_3\tau\int\left\vert
P_m(x,D+i\tau\nabla\varphi(x))v\right\vert^2dx\geq\\&\geq
\varepsilon_0C^{-1}\sum_{|\alpha|\leq m}\tau^{2(m-|\alpha|)}\int
\left\vert D^{\alpha}v\right\vert^2dx,
\end{aligned}
\end{equation}
for every $v\in C_0^{\infty}(B_{\overline{\rho}}\cap \Omega)$ and for every
$\tau\geq 1$. Now, in \eqref{dim-17} we move on the right--hand side  the first term which is on the left--hand side and we get

\begin{equation*}
\begin{aligned}
2C_3\tau\int\left\vert
P_m(x,D+i\tau\nabla\varphi(x))v\right\vert^2&dx\geq
\varepsilon_0C^{-1}\sum_{|\alpha|= m}\tau^{2(m-|\alpha|)}\int
\left\vert D^{\alpha}v\right\vert^2dx+\\& +\sum_{|\alpha|\leq
m-1}\tau^{2(m-|\alpha|)}\left(C^{-1}\varepsilon_0-C\tau^{-1}\right)\int\left\vert
D^{\alpha}v\right\vert^2dx,
\end{aligned}
\end{equation*}
for every $v\in C_0^{\infty}(B_{\overline{\rho}}\cap \Omega)$ and for every $\tau\geq 1$. Hence, if $\tau\geq \tau_0$, where
$\tau_0=\max\{2C^2\varepsilon_0^{-1},1\}$, we have

\begin{equation*}
\begin{aligned}
2C_3\tau\int\left\vert
P_m(x,D+i\tau\nabla\varphi(x))v\right\vert^2&dx\geq
\frac{\varepsilon_0C^{-1}}{2}\sum_{|\alpha|\leq
m}\tau^{2(m-|\alpha|)}\int \left\vert D^{\alpha}v\right\vert^2dx,
\end{aligned}
\end{equation*}
for every $v\in C_0^{\infty}(B_{\overline{\rho}}\cap \Omega)$ and for every $\tau\geq \tau_0$. By using Lemma \ref{lemma-carlm-ell-I}
and by recalling (compare with \eqref{sec-coniugato-1})
\begin{equation*}
P_m(x,D+i\tau\nabla\varphi(x))v=e^{\tau\varphi(x)}P_m(x,D)u,
\end{equation*}
we have

\begin{equation}\label{dim-18}
\sum_{|\alpha|\leq m}\tau^{2(m-|\alpha|)-1}\int \left\vert
D^{\alpha}u\right\vert^2e^{2\tau\varphi}dx\leq C\int \left\vert
P_m(x,D)u\right\vert^2e^{2\tau\varphi}dx.
\end{equation}

Estimate \eqref{Carleman-ellittico-6} follows by Lemma
\ref{stime-Carlm-1.3.3-39} and by the comments made at the beginning
of the introduction to this Chapter. $\blacksquare$

\bigskip

\underline{\textbf{Exercise.}} Prove Theorem \ref{Carlm-ell-teor} (assuming $C^{\infty}$ coefficients in the principal part)
by employing  decomposition \eqref{note-50-2911} instead of 
decomposition \eqref{note-5-49a-2911}. [Hint: recall
\eqref{note-13-50} and use

$$-2\Re\left(z\overline{w}\right)\geq -|z|^2-|w|^2,$$ for $z,w\in \mathbb{C}$]. $\clubsuit$

\subsection{Elliptic operators
with Lipschitz continuous coefficients and the Cauchy problem}\label{Carleman-ellittico-28}
In Theorem \ref{Carlm-ell-teor} we have assumed that the coefficients of the
principal part are of class $C^{1}(\overline{\Omega})$ and it turns out that
the constants, $C$ and $\tau_0$, in the estimate \eqref{Carleman-ellittico-6}
depend on the \textbf{modulus of continuity} of the gradients of these
coefficients. We will now see that with a relatively modest effort
we can prove a Carleman estimate for elliptic operators
with Lipschitz  continuous coefficients in principal part. In this regard, it is useful to point out that this assumption cannot be substantially reduced as has been shown in the
counterexamples of Mandache's Mandache \cite{Mand} and of Pl\u{\i}s \cite{Pl}.

 \medskip

For any $x\in \mathbb{R}^n$ and $R>0$ let us denote by $$Q_R(x)=\left\{y\in
\mathbb{R}^n:\quad |y_j-x_j|<R, \quad j=1,\cdots, n\right\}.$$
Let us introduce a special \textbf{partition of unity}.

\medskip

Let $\vartheta_0\in C^\infty_0(\mathbb{R})$ satisfy
\begin{equation*}
\vartheta_0(t)=
\begin{cases}
1,\quad \mbox{for }|t|\leq 1,\\
\\
0, \quad \mbox{for }|t|\geq 3/2.
\end{cases}
\end{equation*}
Let, further, $0\leq\vartheta\leq 1$ such that
$$\vartheta(x)=\vartheta_0(x_1)\cdots\vartheta_0(x_{n}),$$ we have
\begin{equation*}
\vartheta(x)=
\begin{cases}
1,\quad \mbox{for } x\in Q_1(0),\\
\\
0, \quad \mbox{for } x\in\mathbb{R}^{n}\setminus \overline{Q_{3/2}(0)}.
\end{cases}
\end{equation*}
For any $\mu\geq 1$ and $g\in \mathbb{Z}^{n}$, let us denote
$$x_g=g/\mu$$
and
$$\vartheta_{g,\mu}(x)=\vartheta(\mu(x-x_g)).$$
Hence, we have 
\begin{equation*}
{\rm supp}\,\vartheta_{g,\mu}\subset
\overline{Q_{3/2\mu}(x_g)}\subset Q_{2/\mu}(x_g)
\end{equation*}
and
\begin{equation}\label{6.4}
|D^k\vartheta_{g,\mu}|\leq C_1\mu^k(\chi_{
Q_{3/2\mu}(x_g)}-\chi_{Q_{1/\mu}(x_g)}),\quad k=0,1,\cdots,m,
\end{equation}
where $C_1\geq 1$ depends on $n$ only.

\medskip

For any $g\in \mathbb{Z}^{n}$, set $$A_g=\{g'\in \mathbb{Z}^{n}\,|\,{\rm
supp}\,\vartheta_{g',\mu}\cap{\rm
supp}\,\vartheta_{g,\mu}\neq\emptyset\},$$ then
\begin{equation}\label{cardAg}
card(A_g)\mbox{  depends only on } n.\end{equation}
 Therefore we can define
\begin{equation}\label{6.5}
\widetilde{\vartheta}_{\mu}(x):=\sum_{g\in
\mathbb{Z}^n}\vartheta_{g,\mu}(x)\geq 1,\quad \forall x\in \mathbb{R}^{n}.
\end{equation}
By \eqref{6.4}, we get
\begin{equation}\label{6.6}
|D^k\widetilde{\vartheta}_{\mu}|\leq C_2\mu^k,
\end{equation}
where $C_2\geq 1$ depends on $n$ only. Define
$$\eta_{g,\mu}(x)=\vartheta_{g,\mu}(x)/\widetilde{\vartheta}_{\mu}(x),\quad \forall x\in \mathbb{R}^{n},$$
we have thus
\begin{equation}\label{6.7}
\begin{cases}
\eta_{g,\mu}\geq 0,\\
\\
\sum_{g\in \mathbb{Z}^{n}}\eta_{g,\mu}= 1,\quad \mbox{in } \mathbb{R}^{n},\\
\\
{\rm supp}\,\eta_{g,\mu}\subset \overline{Q_{3/2\mu}(x_g)}\subset Q_{2/\mu}(x_g),\\
\\
|D^{\alpha}\eta_{g,\mu}|\leq
C_3\mu^{|\alpha|}\chi_{\overline{Q_{3/2\mu}(x_g)}},\quad \forall
\alpha \in \mathbb{N}^n, \mbox{ }1\leq|\alpha|\leq m,
\end{cases}
\end{equation}
where $C_3\geq 1$ depends on $n$ only.

\bigskip

Let $m\in\mathbb{N}$ and let 
\begin{equation}\label{qv-oper}
P(x,D)= \sum_{|\alpha|\leq m}a_{\alpha}(x)D^{\alpha},
\end{equation}
be an elliptic operator.
Let $M_0, M_1, \lambda$ be positive constants and let us suppose
that

\begin{subequations}
\label{qv-1-1w}
\begin{equation}
\label{qv-1-1wa} \left\Vert
a_{\alpha}\right\Vert_{L^{\infty}(Q_1)}\leq M_0, \quad \mbox{for }
|\alpha|\leq m,
\end{equation}
\begin{equation}
\label{qv-1-1wb} \left\vert
a_{\alpha}(x)-a_{\alpha}(y)\right\vert\leq M_1|x-y|,\quad \forall x,y\in Q_1, \quad\mbox{for }
|\alpha|=m,
\end{equation}
\begin{equation}
\label{qv-1-1wc} \left\vert P_m(x,\xi)\right\vert\geq \lambda
\left\vert \xi\right\vert^{m},\quad\quad \forall \xi\in
\mathbb{R}^n, \quad \forall x\in \overline{Q}_1.
\end{equation}
\end{subequations}

\medskip

Let $\varphi\in C^{\infty}\left(\overline{Q}_1\right)$ and for $x,y\in
\overline{Q}_1$, $(\xi,\sigma)\in\mathbb{R}^{n+1}$ set

\begin{equation}\label{qv-1-6}
	\begin{aligned}
		&\mathcal{G}(x,y;\xi,\sigma)=\\&=\sum_{j,k=1}^n\partial^2_{x_jx_k}\varphi(x)
P_m^{(j)}(y,\xi+i\sigma\nabla\varphi(x))\overline{P_m^{(j)}(y,\xi+i\sigma\nabla\varphi(x))}.
\end{aligned}
\end{equation}

\begin{theo}\label{qv-prop-2w}
Let us assume that operator \eqref{qv-oper} satisfies ellipticity condition
\eqref{qv-1-1wc} and that its
coefficients satisfy conditions \eqref{qv-1-1wa} and
\eqref{qv-1-1wb}. Moreover, let us assume that

\begin{equation}\label{qv-1-2w}
\begin{aligned}
\begin{cases}
P_m(0,\xi+i\sigma\nabla\varphi(0))=0,\\
\\
(\xi,\sigma)\in \mathbb{R}^{n+1}\setminus\{(0,0)\},
\end{cases}\mbox{ }\Longrightarrow \mbox{ }\mathcal{G}(0,0;\xi,\sigma)>0.
\end{aligned}
\end{equation}
Then there exist $\overline{R}\in (0,1]$, $\delta_0\in(0,1]$ $C_0\geq
1$ and $\tau_0$ such that
\begin{equation}\label{qv-2-2w}
\sum_{|\alpha|\leq m}\tau^{2(m-|\alpha|)-1}\int \left\vert
D^{\alpha}u\right\vert^2e^{2\tau\varphi}dx\leq C_0\int \left\vert
P_m(\delta x,D)u\right\vert^2e^{2\tau\varphi}dx,
\end{equation}
for every $\delta\in (0,\delta_0]$, for every $u\in
C_0^{\infty}\left(Q_{\overline{R}}\right)$ and for every $\tau\geq
\tau_0$.
\end{theo}

\textbf{Proof.} Since
$|P_m(0,\xi+i\sigma\nabla\varphi(0))|^2$ and
$\left(|\xi|^2+\tau^2\right)\mathcal{G}(0,0;\xi,\sigma)$ are homogeneous
polynomials of degree $2m$, \eqref{qv-1-2w} is
equivalent to the following property: there exist positive constants $C_1$ and
$C_2$ such that (Lemma \ref{106-CauNire})
\begin{equation}\label{qv-3-2w-28}
\begin{aligned}
C_2\left|P_m(0,\xi+i\sigma\nabla\varphi(0))\right|^2+\left(|\xi|^2+\sigma^2\right)&\mathcal{G}(0,0;\xi,\sigma)\geq\\&\geq
C_1 \left(|\xi|^2+\sigma^2\right)^m. \end{aligned}
\end{equation}
for every $(\xi,\sigma)\in \mathbb{R}^{n+1})$.

 Let us denote
\begin{equation}\label{qv-3bis-2w}
H(x,y;\xi,\sigma)=C_2\left|P_m(y,\xi+i\sigma\nabla\varphi(x))\right|^2+\left(|\xi|^2+\sigma^2\right)\mathcal{G}(x,y;\xi,\sigma)
\end{equation}
and let us notice that $H$ is a continuous function. Moreover, by \eqref{qv-3-2w-28} we have trivially

$$H(0,0;\xi,\sigma)\geq C_1, \quad \mbox{for all } (\xi,\sigma) \mbox{ such that } |\xi|^2+\sigma^2=1.$$ Hence, the continuity of $H$ implies that
there exists $\overline{R}_1\in(0,1]$ such that

$$H(x,y;\xi,\sigma)\geq \frac{C_1}{2}, \quad \mbox{for all } (\xi,\sigma) \mbox{ such that } |\xi|^2+\sigma^2=1,$$ for every $x,y\in
\overline{Q}_{\overline{R}_1}$. 

Therefore

\begin{equation}\label{qv-1-3w}
\begin{aligned}
C_2\left|P_m(y,\xi+i\sigma\nabla\varphi(x))\right|^2+\left(|\xi|^2+\sigma^2\right)&\mathcal{G}(x,y;\xi,\sigma)\geq\\&\geq
\frac{C_1}{2} \left(|\xi|^2+\sigma^2\right)^m,
\end{aligned}
\end{equation}
for every $x,y\in \overline{Q}_{\overline{R}_1}$. By the previuos inequality we have that for every $\widetilde{y}\in Q_{\overline{R}_1}$ it occurs

\begin{equation}\label{qv-1-3wbis}
\begin{aligned}
\begin{cases}
P_m(\widetilde{y},\xi+i\sigma\nabla\varphi(x))=0,\\
\\
(\xi,\sigma)\in \mathbb{R}^{n+1}\setminus\{(0,0)\}%
\end{cases}\mbox{ }\Longrightarrow \mbox{
}\mathcal{G}(x,\widetilde{y};\xi,\sigma)>0.
\end{aligned}
\end{equation}
Now, for $y\in \overline{Q}_{\overline{R}_1}$ \textit{fixed} and
$\delta\in(0,1]$ to be chosen, let us consider the operator with
constant coefficients w.r.t. the variable $x$
\begin{equation}\label{qv-2-3w}
P_m(\delta y, D_x)=\sum_{|\alpha|=m}a_{\alpha}(\delta
y)D^{\alpha}_x.
\end{equation}
Of course $\delta y\in \overline{Q}_{\overline{R}_1}$ (as
$\delta\in(0,1]$). Now, \eqref{qv-1-3wbis} (considered for
$\widetilde{y}=\delta y$) is nothing but (compare
\textbf{Remark 2} after Theorem \ref{Carlm-ell-teor})
 condition ($\bigstar$) of Theorem \ref{Carlm-ell-teor}. Hence
there exist $C_3>0$ e $\tau_1$ such that

\begin{equation}\label{qv-3-2w}
	\begin{aligned}
	&\sum_{|\alpha|\leq m}\tau^{2(m-|\alpha|)-1}\int \left\vert
D^{\alpha}u\right\vert^2e^{2\tau\varphi(x)}dx\leq \\&\leq C_3\int \left\vert
P_m(\delta y,D)u\right\vert^2e^{2\tau\varphi(x)}dx,
\end{aligned}
\end{equation}
for every $u\in C_0^{\infty}\left(Q_{\overline{R}_1}\right)$ and for every $\tau\geq \tau_1$. Moreover $C_3>0$ and $\tau_1$ \textbf{do not 
depend neither on $y\in \overline{Q}_{\overline{R}_1}$ nor on
$\delta\in(0,1]$}.

\medskip
We now use the partition of unity introduced above with

\begin{equation}\label{mu-w}
\mu=\sqrt{\varepsilon \tau},
\end{equation} for $\tau\geq
\tau^{(\varepsilon)}:=\max\left\{\varepsilon^{-1},\tau_1\right\}$
where $\varepsilon\in(0,1]$ is to be chosen.

Let $u\in C_0^{\infty}\left(Q_{\overline{R}_1}\right)$. By the first relation of  \eqref{6.7} we have
\begin{equation}\label{qv-2-4w}
u=\sum_{g\in \mathbb{Z}^n}u\eta_{g,\mu}.
\end{equation}
Now we apply \eqref{qv-3-2w} (for $y=x_{g}\in
Q_{\overline{R}_1}$). We have, for every $\tau\geq
\tau^{(\varepsilon)}$

\begin{equation}\label{qv-3-4w}
\begin{aligned}
&\sum_{|\alpha|\leq m}\tau^{2(m-|\alpha|)-1}\int \left\vert
D^{\alpha}u\right\vert^2e^{2\tau\varphi}dx\leq \\&\leq c\sum_{g\in
\mathbb{Z}^n}\sum_{|\alpha|\leq
m}\tau^{2(m-|\alpha|)-1}\int\left|D^{\alpha}\left(u\eta_{g,\mu}\right)\right|^2e^{2\tau\varphi}dx\leq
\\& \leq cC_3 \sum_{g\in \mathbb{Z}^n}\int\left\vert P_m\left(\delta
x_g,D\right)\left(u\eta_{g,\mu}\right)\right\vert^2e^{2\tau\varphi}dx,
\end{aligned}
\end{equation}

\medskip

\noindent where the constant $c$ appearing in the second
inequality, by \eqref{cardAg}, depends on $n$ only.

Now, let us estimate form above the last term on the right--hand side of \eqref{qv-3-4w}.
We have

\begin{equation}\label{qv-4-4w}
\begin{aligned}
\left\vert P_m\left(\delta
x_g,D\right)\left(u\eta_{g,\mu}\right)\right\vert^2&\leq 2\left\vert
P_m\left(\delta x,D\right)\left(u\eta_{g,\mu}\right)\right\vert^2
+\\&+2\left\vert \left(P_m\left(\delta x_g,D\right)-P_m\left(\delta
x,D\right)\right)\left(u\eta_{g,\mu}\right)\right\vert^2.
\end{aligned}
\end{equation}
In order to estimate the \textbf{first term on the right--hand side in \eqref{qv-4-4w}}
we notice that

\begin{equation}\label{qv-4-4w-28}
\begin{aligned}
P_m\left(\delta
x,D\right)\left(u\eta_{g,\mu}\right)&=\eta_{g,\mu}\sum_{|\alpha|\leq
m}a_{\alpha}(\delta x)D^{\alpha}u+\\&+\sum_{|\alpha|\leq
m}a_{\alpha}(\delta
x)\sum_{\beta<\alpha}\binom{{\alpha}}{{\beta}}D^{\beta}u
D^{\alpha-\beta}\eta_{g,\mu}=\\&=\eta_{g,\mu}P_m\left(\delta
x,D\right)u+\widetilde{P}(x,D,\mu)u
\end{aligned}
\end{equation}
where we set

$$\widetilde{P}(x,D,\mu)u=\sum_{|\alpha|\leq
m}a_{\alpha}(\delta
x)\sum_{\beta<\alpha}\binom{{\alpha}}{{\beta}}D^{\beta}u
D^{\alpha-\beta}\eta_{g,\mu}.$$ Let us note that this operator has order
$m-1$. Moreover by \eqref{6.7} and \eqref{qv-1-1wa}
we get

\begin{equation}\label{qv-4bis-4w}
\left|\widetilde{P}(x,D,\mu)u\right|\leq
CM_0\chi_{Q_{2/\mu}(x_g)}\sum_{|\beta|\leq
m-1}\left|D^{\beta}u\right|\mu^{m-|\beta|}.
\end{equation}

From \eqref{qv-4-4w}, \eqref{qv-4bis-4w} and recalling \eqref{mu-w}
we have (for the first term on the right we use the trivial inequality
$\eta_{g,\mu}^2\leq \eta_{g,\mu}$)

\begin{equation}\label{qv-1-5w}
\begin{aligned}
\left\vert P_m\left(\delta
x,D\right)\left(u\eta_{g,\mu}\right)\right\vert^2&\leq
\eta_{g,\mu}\left\vert P_m\left(\delta
x,D\right)u\right\vert^2+\\&+CM^2_0\chi_{Q_{2/\mu}(x_g)}\sum_{|\alpha|\leq
m-1}(\varepsilon\tau)^{m-|\alpha|}\left|D^{\alpha}u\right|^2.
\end{aligned}
\end{equation}
We now estimate the second term on the right--hand side in \eqref{qv-4-4w}.
Proceeding in a similar way as above, we have
\begin{equation*}
\begin{aligned}
&\left\vert\left(P_m\left(\delta x_g,D\right)-P_m\left(\delta
x,D\right)\right)\left(u\eta_{g,\mu}\right)\right\vert\leq\\&\leq
\sum_{|\alpha|= m}\left| \left(a_{\alpha}(\delta
x)-a_{\alpha}(\delta
x_g)\right)\right|\left|D^{\alpha}\left(u\eta_{g,\mu}\right)\right|=\\&=
\eta_{g,\mu}\sum_{|\alpha|= m}\left| \left(a_{\alpha}(\delta
x)-a_{\alpha}(\delta
x_g)\right)\right|\left|D^{\alpha}u\right|+\\&+CM_0\chi_{Q_{2/\mu}(x_g)}\sum_{|\alpha|\leq
m-1}\left|D^{\alpha}u\right|\mu^{m-|\alpha|}.
\end{aligned}
\end{equation*}
In order to estimate the second-to-last term, it must be taken into account that the estimate has to be
done in the support of $\eta_{g,\mu}$. By \eqref{qv-1-1wb} we get,
therefore,

\begin{equation}\label{qv-2-5w}
\begin{aligned}
&\left\vert\left(P_m\left(\delta x_g,D\right)-P_m\left(\delta
x,D\right)\right)\left(u\eta_{g,\mu}\right)\right\vert^2\leq
\\&\leq CM^2_1\eta_{g,\mu}\frac{\delta^2}{\varepsilon\tau}\sum_{|\alpha|=m}
\left|D^{\alpha}u\right|^2+\\&+CM^2_0\chi_{Q_{2/\mu}(x_g)}\sum_{|\alpha|\leq
m-1}(\varepsilon\tau)^{m-|\alpha|}\left|D^{\alpha}u\right|^2.
\end{aligned}
\end{equation}

\smallskip

\noindent Now, we insert \eqref{qv-1-5w} and \eqref{qv-2-5w} into
\eqref{qv-4-4w} and we get

\begin{equation*}
\begin{aligned}
\left\vert P_m\left(\delta
x_g,D\right)\left(u\eta_{g,\mu}\right)\right\vert^2&\leq
2\eta_{g,\mu}\left\vert P_m\left(\delta
x,D\right)u\right\vert^2+\\&+
CM^2_1\eta_{g,\mu}\frac{\delta^2}{\varepsilon\tau}\sum_{|\alpha|=m}
\left|D^{\alpha}u\right|^2+\\&+CM^2_0\chi_{Q_{2/\mu}(x_g)}\sum_{|\alpha|\leq
m-1}(\varepsilon\tau)^{m-|\alpha|}\left|D^{\alpha}u\right|^2.
\end{aligned}
\end{equation*}
Inserting the latter into \eqref{qv-3-4w} we obtain

\begin{equation*}
\begin{aligned}
\sum_{|\alpha|\leq m-1}\tau^{2(m-|\alpha|)-1}\int \left\vert
D^{\alpha}u\right\vert^2e^{2\tau\varphi}dx&+\frac{1}{\tau}\sum_{|\alpha|=
m}\int \left\vert D^{\alpha}u\right\vert^2e^{2\tau\varphi}dx
\leq\\&\leq cC_3\int \left\vert P_m\left(\delta
x,D\right)u\right\vert^2e^{2\tau\varphi}dx+\\&+
C_4\frac{\delta^2}{\varepsilon\tau}\sum_{|\alpha|= m}\int \left\vert
D^{\alpha}u\right\vert^2e^{2\tau\varphi}dx+\\&+C_5\sum_{|\alpha|\leq
m-1}\int(\varepsilon\tau)^{m-|\alpha|}\left|D^{\alpha}u\right|^2e^{2\tau\varphi}dx,
\end{aligned}
\end{equation*}
where $C_4$ depends on $M_1$ only and $C_5$ depends by $M_0$ only. From which we have

\begin{equation}\label{qv-2bis-4w}
\begin{aligned}
\sum_{|\alpha|\leq
m-1}\tau^{(m-|\alpha|)}&\left(\tau^{m-|\alpha|-1}-C_5\varepsilon^{m-|\alpha|}\right)\int
\left\vert
D^{\alpha}u\right\vert^2e^{2\tau\varphi}dx+\\&+\frac{1}{\tau}\left(1-C_4\frac{\delta^2}{\varepsilon}\right)\sum_{|\alpha|=
m}\int \left\vert D^{\alpha}u\right\vert^2e^{2\tau\varphi}dx
\leq\\&\leq cC_3\int \left\vert P_m\left(\delta
x,D\right)u\right\vert^2 e^{2\tau\varphi}dx.
\end{aligned}
\end{equation}
Let us choose
$$\varepsilon=\varepsilon_0:=\frac{1}{2C_5},$$
$$\delta\leq\delta_0:=\sqrt{\frac{\varepsilon_0}{2C_4}}$$
and by \eqref{qv-2bis-4w} we get
\begin{equation}\label{qv-1-6w}
\sum_{|\alpha|\leq m}\tau^{2(m-|\alpha|)-1}\int \left\vert
D^{\alpha}u\right\vert^2e^{2\tau\varphi}dx\leq 2cC_3\int \left\vert
P_m(\delta x,D)u\right\vert^2e^{2\tau\varphi}dx,
\end{equation}
for every $u\in C_0^{\infty}\left(Q_{\overline{R}_1}\right)$ and for every
$\tau\geq \tau^{(\varepsilon_0)}$. Estimate \eqref{qv-2-2w} is
proved. $\blacksquare$

\medskip

\textbf{Remark.} Let us notice that (reader check) by the change of variables $X=\delta x$,  \eqref{qv-2-2w} si become

\begin{equation}\label{qv-2-2w-cambiov}
\begin{aligned}
\sum_{|\alpha|\leq m}\tau^{2(m-|\alpha|)-1}&\delta^{m-|\alpha|}\int
\left\vert
D^{\alpha}u\right\vert^2e^{2\tau\varphi\left(\delta^{-1}X\right)}dX\leq\\&\leq
C_0\int \left\vert
P_m(X,D)u\right\vert^2e^{2\tau\varphi\left(\delta^{-1}X\right)}dX,
\end{aligned}
\end{equation}
for every $u\in C_0^{\infty}\left(Q_{\delta\overline{R}}\right)$ and
for every $\tau\geq \tau_0$. $\blacklozenge$

\bigskip

In the following Theorem we will apply estimate \eqref{qv-2-2w-cambiov}
to prove a uniqueness result for the Cauchy  problem.

\begin{theo}\label{Cauchy-ell-teor}
Let $\psi\in C^1\left(\overline{Q}_1\right)$ be real--valued function
such that
\begin{equation}\label{Cauchy-ell-1}
\nabla\psi(0)\neq 0.
\end{equation}
Let $P(x,D)$ be operator \eqref{qv-oper} and let us suppose that
\eqref {qv-1-1w} holds true. Let $U\in H^m\left(Q_{1}\right)$ satisfy
\begin{equation}\label{Cauchy-ell-2}
\begin{cases}
P(x,D)U=0,\quad\mbox{in }\quad Q_{1},\\
\\
U(x)=0 \quad\mbox{in } \left\{x\in \overline{Q}_1: \quad
\psi(x)>\psi(0) \right\}.
\end{cases}
\end{equation}
Let us suppose that for every $\xi\in \mathbb{R}^n\setminus\{0\}$ we have
\begin{equation}\label{Cauchy-ell-4}
\sigma \rightarrow P_m(0,\xi+i\sigma \nabla \psi(0))\mbox{ }\mbox{ has no real multiple roots
	 }.
\end{equation}
Then there exist a neighborhood $\mathcal{U}_{0}$ of $0$ such that

\begin{equation}\label{Cauchy-ell-5}
U=0\quad\mbox{in }\quad \mathcal{U}_{0}.
\end{equation}
\end{theo}
\textbf{Remark.} As it is easily checked, condition
\eqref{Cauchy-ell-4} can be expressed equivalently as
follows

\begin{equation}\label{oper-2ord-343-n}
	\begin{aligned}
		&\begin{cases}
P_m(0,\xi+i\tau \nabla \psi(0))=0,\\
\\
(\xi,\tau)\neq (0,0),
\end{cases} \Longrightarrow \\& \Longrightarrow \sum_{j=1}^nP_m^{(j)}(0,\xi+i\tau \nabla
\psi(0))\partial_j\psi(0)\neq 0.
\end{aligned}
\end{equation}

We further observe that if $m=2$ and the coefficients of $P_2(x,D)$ are
real then the \eqref{Cauchy-ell-4} is satisfied (see Example
4a, Section \ref{altre consid}). $\blacklozenge$

\bigskip

\textbf{Proof of Theorem \ref{Cauchy-ell-teor}.}

 It is not restrictive to assume $\psi(0)=0$ and, since $\nabla\psi(0)\neq 0$,  we may 
 reduce to consider, up to isometries, the case where,
 for an appropriate $r_0>0$, we have

\begin{equation}\label{Cauchy-ell-6}
\left\{x\in Q_{r_0}:\quad \psi(x)=0\right\}=\left\{(x',f(x')):\quad
x'\in Q'_{r_0}(0)\right\},
\end{equation}
($Q'_{r_0}=(-r_0,r_0)^{n-1}$) where $f\in C^1(Q'_{r_0})$, $f(0)=|\nabla
f(0)|=0$ and

$$\left\{x\in Q_1: \quad \psi(x)>0  \right\}\cap Q_{r_0}=\left\{(x',x_n)\in Q_{r_0}:\quad
x_n<f(x')\right\}.$$ Let us notice that in this way condition
\eqref{Cauchy-ell-4} becomes.

\begin{equation}\label{Cauchy-ell-7}
\sigma\rightarrow P_m(0,\xi-i\sigma e_n)\mbox{ has no real multiple roots }.
\end{equation}

Now we use \textbf{Holmgren transformation} introduced in
\eqref{3-4-67C}, that is we consider the transformation
\begin{equation}\label{Cauchy-ell-8}
\Lambda:\mathbb{R}^n_x\rightarrow \mathbb{R}^n_y,\quad x\rightarrow
y=\Lambda(x',x_n)=\left(x',x_n+\frac{A}{2}|x'|^2\right),
\end{equation}
(recall that $\Lambda$ is a diffeomorphism) where $A>0$ satisfies

\begin{equation}\label{correct:19-3-23-1}
A>\left\Vert\partial^2f\right\Vert_{L^{\infty}(B'_{r_0})}
\end{equation}
and where $\partial^2f$ is the Hessian matrix of $f$. Let us fix $A$ that satisfies \eqref{correct:19-3-23-1} and we recall that, with this choice, the function

\begin{equation}\label{Cauchy-ell-9}
g(x')=f(x')+\frac{A}{2}|x'|^2,
\end{equation}
is strictly convex and
\begin{equation}\label{2611-1}
g(0)=|\nabla g(0)|=0.
\end{equation}
Let us denote by $\widetilde{P}(y,D_y)$ the transformed operator
of $P(x,D_x)$ by mean of $\Lambda$. Since
$$P_m(x,\xi)=i^m\sum_{|\alpha|=
m}a_{\alpha}(x)\xi^{\alpha}$$ we have (compare with \eqref{4-11C})

\begin{equation}\label{Cauchy-ell-10}
\widetilde{P}_m(y,\eta)=i^m\sum_{|\alpha|=
m}a_{\alpha}(\Lambda^{-1}(y))\left(\left(\partial_x\Lambda(x)\right)^t\eta\right)_{|x=\Lambda^{-1}(y)}^{\alpha}.
\end{equation}
So the condition \eqref{Cauchy-ell-7} is written (reader
check)
\begin{equation}\label{Cauchy-ell-11}
	\begin{aligned}
		&\mbox{for fixed } \xi\in\mathbb{R}^n\setminus \{0\}, \\& \widetilde{P}_m(0,\xi-i\sigma
e_n)\mbox{ }\mbox{ has no real multiple roots. }
\end{aligned}
\end{equation}

\smallskip

Moreover, \eqref{Cauchy-ell-2} implies

\begin{equation}\label{Cauchy-ell-12}
\begin{cases}
\widetilde{P}(y,D_y)\widetilde{U}=0,\quad\mbox{in }\quad Q_{r_0},\\
\\
\widetilde{U}(y)=0 \quad\mbox{in } \left\{(x',x_n)\in Q_{r_0}:\quad
x_n<g(x')\right\}.
\end{cases}
\end{equation}
where $\widetilde{U}(y)=U\left(\Lambda^{-1}(y)\right)$. It turns out
$\widetilde{U}\in H^m\left(Q_{r_0}\right)$.

We agree from here on to omit " $\widetilde{  \mbox{
  }}$ " from $P$ and  $U$, and  to rename "$x$" the variable "$y$". Let
$$h(x_n)=-x_n+\frac{x_n^2}{2}$$
(let us notice that $h$ is strictly decreasing in $0\leq x_n\leq
1$)

\noindent and
$$\varphi(x)=h(\delta_0x),$$
where $\delta_0$ is defined in Theorem \ref{qv-prop-2w}.
We have  $$\nabla\varphi(0)=-\delta_0e_n.$$ We have
$$\nabla\varphi(0)=-\delta_0e_n$$ and also, 
\eqref{Cauchy-ell-11} implies that if $\xi\in
\mathbb{R}^n\setminus\{0\}$ and

$$P_m(0,\xi+i\sigma\nabla\varphi(0))=P_m(0,\xi-i\sigma\delta_0e_n)=0,$$
then

\begin{equation*}
\begin{aligned}
&\mathcal{G}(0,0;\xi,\sigma)=\\&=\sum_{j,k=1}^n\partial^2_{x_jx_k}\varphi(0)
P_m^{(j)}(0,\xi+i\sigma\nabla\varphi(0))\overline{P_m^{(j)}(0,\xi+i\sigma\nabla\varphi(0))}=\\&=
\delta_0^2\left|P_m^{(n)}(0,\xi+i\sigma\nabla\varphi(0))\right|^2=\delta_0^2\left|P_m^{(n)}(0,\xi-i\sigma\delta_0e_n)\right|^2>0.
\end{aligned}
\end{equation*}
Therefore the assumptions of Theorem \ref{qv-prop-2w} are satisfied. Then
Carleman estimate \eqref{qv-2-2w-cambiov} holds. We may write such a Carleman estimate as (setting $\delta=\delta_0$)

\begin{equation}\label{Cauchy-ell-13}
\begin{aligned}
\sum_{|\alpha|\leq m}\tau^{2(m-|\alpha|)-1}\int \left\vert
D^{\alpha}u\right\vert^2&e^{2\tau h\left(x_n\right)}dx\leq\\&\leq
C\int \left\vert P(x,D)u\right\vert^2e^{2\tau h\left(x_n\right)}dx,
\end{aligned}
\end{equation}
for every $u\in C_0^{\infty}\left(Q_{\delta_0\overline{R}}\right)$ and
for every $\tau\geq \overline{\tau}_0$ for a certain
$\overline{\tau}_0\geq\tau_0$. Set
$$r_1=\min\left\{r_0,\delta_0\overline{R}\right\}$$
and, for $\rho>0$,
$$E_{\rho}=\left\{(x',x_n)\in Q_{r_1}:\quad
g(x')<x_n<\rho\right\}.$$ For the strict convexity of $g$ and by
\eqref{2611-1}, we have that there exists $\rho_1>0$ such that

$$ \overline{E}_{\rho_1}\subset Q_{r_1}.$$

Let $\rho_2\in \left(0,\rho_1\right)$ Let $\eta\in C^{\infty}\left(\mathbb{R}\right)$ be a function such that

$$0\leq \eta(x_n) \leq 1,\quad\forall x\in \mathbb{R}; \quad \eta(x_n)=1,\quad \forall x_n\leq \rho_2; \quad \eta(x)=0,\quad \forall x_n\geq \rho_1.$$

Let us assume that

$$\eta^{(k)}(x_n)\leq C\left(\rho_1-\rho_2\right)^{-k}.$$
By density,  \eqref{qv-2-2w-cambiov} holds for every
$u\in H^m_0(Q_{r_1})$, hence, in particular,
\eqref{qv-2-2w-cambiov} holds  for $u(x)=U(x)\eta(x_n)$. From now on, the
proof is quite standard, we present it for
completeness. Since $P(x,D)U=0$ in $Q_{r_1}$ we have.

$$|P(x,D)(U\eta)|\leq CM_0\chi_{\mathbb{R}\setminus(\rho_2,\rho_1)}\sum_{|\alpha|\leq
m-1}(\rho_1-\rho_2)^{-|\alpha|}|D^{\alpha}U|$$ and, for any
$0<\rho<\rho_2$, by \eqref{qv-2-2w-cambiov}, we get  

\begin{equation*}
\begin{aligned}
&e^{2\tau h\left(\rho\right)}\sum_{|\alpha|\leq
m}\tau^{2(m-|\alpha|)-1}\int_{E_{\rho}} \left\vert
D^{\alpha}U\right\vert^2dx\leq\\&\leq \sum_{|\alpha|\leq
m}\tau^{2(m-|\alpha|)-1}\int_{E_{\rho_1}} \left\vert
D^{\alpha}(U\eta)\right\vert^2e^{2\tau
h\left(x_n\right)}dx\leq\\&\leq C\int_{E_{\rho_1}} \left\vert
P(x,D)(U\eta)\right\vert^2e^{2\tau h\left(x_n\right)}dx\leq\\&\leq
CM_0^2e^{2\tau h\left(\rho_2\right)}\sum_{|\alpha|\leq
m-1}\int_{E_{\rho_1}\setminus
E_{\rho_2}}(\rho_1-\rho_2)^{-|\alpha|}|D^{\alpha}U|dx,
\end{aligned}
\end{equation*}
for every $\tau\geq \overline{\tau}_0$. Hence

\begin{equation*}
\begin{aligned}
&\sum_{|\alpha|\leq m}\tau^{2(m-|\alpha|)-1}\int_{E_{\rho}}
\left\vert D^{\alpha}U\right\vert^2dx\leq
\\&\leq CM_0^2e^{-2\tau\left(
h\left(\rho\right)-h\left(\rho_2\right)\right)}\sum_{|\alpha| \leq
m-1}\int_{E_{\rho_1}\setminus
E_{\rho_2}}(\rho_1-\rho_2)^{-|\alpha|}|D^{\alpha}U|dx
\end{aligned}
\end{equation*}
for every $\tau\geq \overline{\tau}_0$. Passing to the limit as $\tau\rightarrow +\infty$
which goes to infinity, and and taking into account that
$h\left(\rho\right)-h\left(\rho_2\right)>0$ we have $U=0$ in $Q_{r_1}\cap\{x_n\leq \rho\}$. Theorem is proved. $\blacksquare$

\chapter{Carleman estimates and the Cauchy problems II -- Second order operators} \markboth{Chapter 14. Carleman estimates and the Cauchy problem II}{}\label{operatori-2ord}

\section{Introduction}\label{corr:26-4-23} In this Chapter we will consider the second-order operators
whose principal part (not necessarily elliptic) is given by

\begin{equation}\label{oper-2ord-1}
P_2(x,\partial)=\sum_{j,k=1}^ng^{jk}(x)\partial^{2}_{x_jx_k},
\end{equation}
where the matrix of coefficients
$\left\{g^{jk}(x)\right\}_{j,k=1}^n$ is a symmetric and invertible matrix, whose entries are the real--valued functions $g^{jk}$  defined on a bounded open set $\Omega\subset\mathbb{R}^n$, on which we will make appropriate regularity assumptions. When
we will refer to the symbol of the operator \eqref{oper-2ord-1}, here we will
always refer to the polynomial in the variable $\xi$
\begin{equation}\label{oper-2ord-2}
P_2(x,\xi)=\sum_{j,k=1}^ng^{jk}(x)\xi_j\xi_k.
\end{equation}
We note that, with the notation for the derivatives used in the previous Sections, operator
 \eqref{oper-2ord-1} can be written 

\begin{equation}\label{oper-2ord-3}
P_2(x,\partial)=-\sum_{j,k=1}^ng^{jk}(x)D^{2}_{x_jx_k}
\end{equation}
and so polynomial \eqref{oper-2ord-2} is simply the
symbol of operator \eqref{oper-2ord-3} with the sign changed.
This abuse of notation will not create major problems, in particular, it will not create problems when we compare procedures
and results found in this Section with those in the previous sections.

The purpose of this Section is to derive the Carleman estimates
in a more direct fashion than the last two
sections. This will allow us, in particular, to write down
explicitly the quadratic form $G(x,D,\overline{D},\tau)$
given by Proposition \eqref{note-52-prop} making it more
easy to obtain estimates for the operators with $C^{0,1}$ coefficients. 

\section{The case of the Laplace operator}\label{24-4-23:laplace-operator}

The case of the Laplace operator will serve us somewhat as a model
for more general operators of type \eqref{oper-2ord-1}.

\medskip

We begin by the following

\begin{lem}[\textbf{The Rellich identity}]\label{rellich}
	\index{Lemma:@{Lemma:}!- Rellich identity@{- Rellich identity}}
Let $\beta \in C^{0,1} (\Omega, \mathbb{R}^n)$,
$\beta=\left(\beta^1,\cdots, \beta^n \right)$ and $v \in C^2
(\Omega)$, then

\begin{equation}\label{oper-2ord-4}
\begin{aligned}
2(\beta \cdot \nabla v) \Delta v &= \mbox{div} \left(2(\beta \cdot
\nabla v ) \nabla v - \beta |\nabla v|^2\right) + \\&+(\mbox{div }
\beta) |\nabla v|^2 - 2  \partial_k \beta^j
\partial_j v
\partial_k v, \ \ \mbox{a.e. } x\in \Omega,
\end{aligned}
\end{equation}
(in \eqref{oper-2ord-4} we have used the Einstein notation
of repeated indices).
\end{lem}

\textbf{Proof.}  We have
\begin{equation*}
\begin{aligned}
&2(\beta \cdot \nabla v) \Delta v  = 2 \left( \beta^j \partial_j v
\right) \Delta v= 2 \left( \partial_k \left(\beta^j \partial_j v
\partial_k v\right) -
\partial_k \left(\beta^j \partial_j v\right) \partial_k v\right) = \\ & =
2\left(
\partial_k \left(\beta^j \partial_j v \partial_k v -
\left(\partial_k \beta^j\right) \partial_j v \partial_k v-\beta^j\partial^2_{jk} v \partial_k v\right)\right)= \\
& = 2 \mbox{div } \left[(\beta \cdot \nabla v) \nabla v\right] -
2\left(\partial_k \beta^j\right)
\partial_j v
\partial_k v - \beta^j \partial_j\left(|\nabla v|^2\right)= \\ & =
2 \mbox{div } \left[(\beta \cdot \nabla v) \nabla v\right] -
2\left(\partial_k \beta^j\right)
\partial_j v
\partial_k v-\partial_j (\beta^j |\nabla v|^2)+(\mbox{div }\beta) |\nabla v|^2=\\& =\mbox{div } \left[2(\beta \cdot \nabla v ) \nabla v - \beta
|\nabla v|^2\right]  - 2 (\partial_k \beta^j)
\partial_j v
\partial_k v+(\mbox{div }\beta) |\nabla v|^2.
\end{aligned}
\end{equation*}
aalmost everywhere in $\Omega$. $\blacksquare$

\bigskip

\textbf{Remark.} By \eqref{oper-2ord-4} we have immediately 

\begin{equation}\label{oper-2ord-5}
\begin{aligned}
\int_{\Omega}2(\beta \cdot \nabla v) \Delta v
dx&=\int_{\Omega}\left((\mbox{div } \beta) |\nabla v|^2 - 2
\partial_k \beta^j
\partial_j v
\partial_k v\right)dx,
\end{aligned}
\end{equation}
for every $v\in C^{\infty}_0(\Omega)$. On the other hand, as can be
easily checked, if $v$ is a real--valued function we have

$$2(\beta \cdot \nabla v) \Delta
v=F\left(x,D,D\right)\left[v,\overline{v}\right],$$ where

$$F\left(x,D,D\right)\left[v,\overline{v}\right]=i\sum_{j,k=1}^n\beta^j\left(D^2_kv\overline{D_jv}-D_jv\overline{D^2_kv}\right)$$
which satisfies condition \eqref{fqd-16-var} of Lemma
\ref{integr-parti-var}. In our case \eqref{fqd-17-var} takes
the form \eqref{oper-2ord-5}. Let us notice that by
\eqref{fqd-18-var} we obtain
$$G(x,\xi,\xi)=(\mbox{div
}\beta)|\xi|^2-2\partial_k\beta^j\xi_j\xi_k.$$  $\blacklozenge$

\bigskip

\textbf{Let us review some key steps of the proof of
	Theorem \ref{Carlm-ell-teor} using Rellich identity
	\eqref{rellich} to perform the integrations by parts}.

We begin by rewriting the statement of Theorem \ref{Carlm-ell-teor}
in the case of the Laplace operator

\begin{theo}\label{teo2bis} Let $\Omega$ be a bounded open set of $\mathbb{R}^n$,
and let $\varphi \in C^{\infty}(\bar{\Omega})$ be a real--valued function such that
$\nabla \varphi \neq 0$ on $\bar{\Omega}$. Let us assume that the following implication holds true

\begin{equation}\label{pseudconv-Delta}
\begin{aligned}
\begin{cases}
|\xi|^2=\tau^2| \nabla \varphi(x)|^2, \\
\\
\xi \cdot \nabla \varphi (x)= 0,\\
\\
\tau\neq 0,
\end{cases}\mbox{ }
\Longrightarrow \mbox{ }&Q(x,\xi,\tau)=\sum_{j,k=1}^{n}
\partial^2_{jk} \varphi (x) \xi_j\xi_k+\\&+\tau^2\sum_{j,k=1}^{n}
\partial^2_{jk}(x) \varphi
\partial_j\varphi(x)\partial_k\varphi(x)>0.
\end{aligned}
\end{equation}

\medskip

Then there exist  constants $C$ and $\tau_0$  such that
\begin{equation*}
\begin{aligned}
& \tau^3 \int_{\Omega} |u|^2 e^{2 \tau \varphi}dx + \tau
\int_{\Omega}
 |\nabla u|^2 e^{2 \tau \varphi}dx +
\tau^{-1}\int_{\Omega} |\partial^2 u|^2 e^{2 \tau \varphi}dx  \leq C
\int_{\Omega} |\Delta u|^2 e^{2 \tau \varphi}dx,
\end{aligned}
\end{equation*}
for every $u \in C^{\infty}_0(\Omega)$ and for every $\tau\geq \tau_0$.
\end{theo}

\bigskip

\textbf{Proof.} First, we observe that \eqref{pseudconv-Delta} is simply  the rewriting of condition $(\bigstar)$ of Theorem
\ref{Carlm-ell-teor} in the case of the Laplace operator.

Let us denote by $L$ the operator $\Delta$. Let $u\in
C^{\infty}_0(\Omega)$. Set $v= e^{\tau \varphi} u$. Let us calculate
\begin{equation*}
\begin{aligned}
&\partial_j u= e^{-\tau \varphi} \left(\partial_j v - \tau \partial_j \varphi v\right),\\
&\partial^2_j u = e^{-\tau \varphi} \left(\partial^2_j v -2\tau
\partial_j \varphi \partial_j v - \tau \partial^2_j \varphi v + \tau^2 v
(\partial_j \varphi)^2\right).
\end{aligned}
\end{equation*}
We obtain
\begin{equation}\label{oper-2ord-6}
L_\tau v=e^{\tau \varphi} L(e^{-\tau \varphi} v)= \Delta v - \tau
(\Delta \varphi) v - 2 \tau \nabla \varphi \cdot \nabla v + \tau^2
|\nabla \varphi|^2 v.
\end{equation}
Now in the setting provided in Section \ref{sec-coniugato}, in the
first line, (see in particular \eqref{note-3-49}) we have
neglected the term $- \tau (\Delta \varphi) v$ and, in the middle part of the proof of Theorem
\ref{Carlm-ell-teor}, we have focused on the operator $p_m(x,D,\tau)$, which in this special 
case is given by

\begin{equation}\label{oper-2ord-7}
\widetilde{L}_\tau v= \Delta v - 2 \tau \nabla \varphi \cdot \nabla
v + \tau^2 |\nabla \varphi|^2 v.
\end{equation}
Consequently, operators \eqref{note-5-49a} and \eqref{note-5-49b}
in the case of the Laplace operator, are, respectively, given by

\begin{equation*}
 S_\tau v= \Delta v + \tau^2 |\nabla
\varphi |^2 v,
\end{equation*}
\begin{equation*}
A_\tau v = -2 \tau \nabla \varphi \cdot \nabla v.
\end{equation*}
Of course, $$\widetilde{L}_\tau  v = S_\tau v + A_\tau v,$$ which implies
\begin{equation}\label{oper-2ord-8}
\int \left|\widetilde{L}_\tau v\right|^2dx= \int |S_\tau v|^2 dx +
\int |A_\tau v|^2dx + 2 \int S_\tau v A_\tau vdx,
\end{equation}
(for brevity, we omit the domain of integration). As
we saw, a crucial point in the proof of Theorem
\ref{Carlm-ell-teor} consists to handle the third integral to the
right--hand side of \eqref{oper-2ord-8}, which we will pursue here using
the identity \eqref{oper-2ord-5}, where

$$\beta=\nabla \varphi.$$
We have
$$2 \int (\nabla \varphi \cdot \nabla v) \Delta v \; dx = \int [\Delta \varphi |\nabla v|^2 - 2 \partial^2_{jk} \varphi \partial_j v \partial_k v] \; dx.$$
Integrating by parts, we get
\begin{equation*}
\begin{aligned}
&2 \int S_\tau v A_\tau vdx = -2 \tau \int \left(2(\nabla \varphi
\cdot \nabla v) \Delta v + 2 \tau^2 |\nabla \varphi|^2 (\nabla
\varphi \cdot \nabla v) v\right)dx =\\&= -2 \tau \int \left(\Delta
\varphi |\nabla v|^2 - 2
\partial^2_{jk} \varphi
\partial_j v \partial_k v +\tau^2 |\nabla \varphi|^2 \nabla \varphi \cdot
\nabla (v^2)\right)dx =\\ & = -2 \tau \int \left(\Delta \varphi |\nabla
v|^2 - 2 \partial^2_{jk} \varphi \partial_j v
\partial_k v - 2
\tau^2\partial^2_{jk} \varphi \partial_j \varphi \partial_k \varphi v^2 -\tau^2 |\nabla \varphi|^2 \Delta \varphi v^2 \right) dx= \\
& = -2 \tau \int \left( \Delta \varphi (|\nabla v|^2-\tau^2 |\nabla
\varphi|^2 v^2) - 2
\partial^2_{jk} \varphi
\partial_j v \partial_k v - 2 \tau^2 (\partial^2_{jk} \varphi
\partial_j \varphi \partial_k \varphi) v^2 \right)
dx.
\end{aligned}
\end{equation*}
Hence
\begin{equation}\label{oper-2ord-9}
\begin{aligned}
2 \int S_\tau v A_\tau vdx&=4\tau\int \left(\partial^2_{jk} \varphi
\partial_j v \partial_k v+\tau^2 (\partial^2_{jk} \varphi
\partial_j \varphi \partial_k \varphi) v^2\right)dx-\\&-2\tau\int \Delta \varphi (|\nabla v|^2-\tau^2
|\nabla \varphi|^2 v^2)dx:=I.
\end{aligned}
\end{equation}
At this point one could continue without involving the
Fourier transform, but for this approach we refer to
\cite{Bal}.

Now, for any $x_0\in \overline{\Omega}$, the integral on the right--hand side
\eqref{oper-2ord-9}, which we denoted by $I$, can be written as

\begin{equation*}
\begin{aligned}
I=I_{x_0}+ R_{x_0},
\end{aligned}
\end{equation*}
where
\begin{equation*}
\begin{aligned}
I_{x_0}&=4\tau\int \left(\partial^2_{jk} \varphi(x_0)
\partial_j v \partial_k v+\tau^2 \partial^2_{jk} \varphi(x_0)
\partial_j \varphi(x_0) \partial_k(x_0)
v^2\right)dx-\\&-2\tau\int \Delta \varphi(x_0) (|\nabla v|^2-\tau^2
|\nabla \varphi(x_0)|^2 v^2)dx,
\end{aligned}
\end{equation*}
and, of course, $$R_{x_0}=I-I_{x_0}.$$ Now, for $\varepsilon>0$ to be chosen, there exists $\rho_1>0$ such that

\begin{equation*}
\begin{aligned}
\left|R_{x_0}\right|\leq \varepsilon \int\left(\tau|\nabla
v|^2+\tau^3|v|^2\right)dx, \ \ \forall v\in C_0^{\infty}(B_{\rho_1}(x_0)\cap \Omega).
\end{aligned}
\end{equation*}
Now by the Parseval identity we have

\begin{equation*}
\begin{aligned}
(2\pi)^n I_{x_0}=\int
q(\xi,\tau)\left|\widehat{v}(\xi)\right|^2d\xi,
\end{aligned}
\end{equation*}
where
\begin{equation*}
\begin{aligned}
q(\xi,\tau)&=4\tau\left(\partial^2_{jk} \varphi(x_0) \xi_j
\xi_k+\tau^2
\partial^2_{jk} \varphi(x_0)
\partial_j \varphi(x_0) \partial_k(x_0)\right)-
\\&-2\tau\Delta
\varphi(x_0) (|\xi|^2-\tau^2 |\nabla \varphi(x_0)|^2).
\end{aligned}
\end{equation*}
Now \eqref{teo2bis} implies there exist positive constants
$C_1$ and $C_2$, such that
\begin{equation*}
C_1\left|\xi+i\tau \nabla\varphi(x_0)\right|^{4}\leq \tau
q(\xi,\tau)+C_2\left|\sum_{j=1}^n\left(\xi+i\tau\partial_j\varphi(x_0)\right)^2\right|^2,\mbox{
} \forall (\xi,\tau)\in \mathbb{R}^{n+1}.
\end{equation*}
(as in \eqref{dim-8}). At this point the most challenging part
of the proof is done and it is not difficult to put the
together the various "pieces" as in the proof of Theorem
\ref{Carlm-ell-teor}, we invite the reader to do so.
 $\blacksquare$

\bigskip

\underline{\textbf{Examples.}} Let us examine some example of function satisfying
\eqref{pseudconv-Delta}.

\medskip

\textbf{Example 1.} Let $\varphi \in C^{\infty}\left(\bar{\Omega}\right)$ satisfy $\nabla \varphi \neq 0$ in $\bar{\Omega}$ and let $\Omega$ be a bounded open set of $\mathbb{R}^2$. We wish to prove that condition
\eqref{pseudconv-Delta} holds for each and only the functions $\varphi$
which satisfy

\begin{equation}\label{oper-2ord-9-28}
 \Delta \varphi >
0,\quad\forall x\in \overline{\Omega}.
\end{equation}
Let us denote
$$N(x)=\frac{\nabla\varphi(x)}{|\nabla\varphi(x)|}$$
and 

$$Q(x,\xi,\tau)=\sum_{j,k=1}^{2} \partial^2_{jk} \varphi (\xi_j \xi_k + N_j N_k).$$
It is clear that \eqref{pseudconv-Delta}  is equivalent to the following condition

\medskip

\noindent ($\bigstar''$) if

\begin{equation}\label{(sub2)}
\begin{aligned} \begin{cases}
|\xi|^2 = 1, \\
\\
N(x) \cdot \xi = 0,
\end{cases}
\end{aligned}
\end{equation}
then  $$Q(x,\xi,\tau)>0.$$

\medskip

Let us suppose, then, that  \eqref{(sub2)} is true.

Since $\nabla \varphi \neq 0$, we may assume $N_1 \neq 0$.
Recalling that  $|N|=1$, by \eqref{(sub2)} we obtain
\begin{equation*}
\begin{aligned} \begin{cases}
\xi_1^2 + \xi_2^2 = 1, \\
\\
\xi_1 N_1 + \xi_2 N_2 = 0,
\end{cases} \Longleftrightarrow
\begin{cases}
\xi_2^2 \left( \frac{N_2^2}{N_1^2} + 1 \right) = 1, \\
\\
\xi_1= - \frac{\xi_2 N_2}{N_1},
\end{cases}  \Longleftrightarrow
\begin{cases}
|\xi_2| = |N_1|, \\
\\
\xi_1= - \frac{\xi_2 N_2}{N_1}.
\end{cases}
\end{aligned}
\end{equation*}
Let us suppose that $\xi_2=|N_1|$, which implies $\xi_1 =
-\frac{N_2}{N_1}|N_1|$, we have
\begin{equation*}
\begin{aligned}
Q & = \partial^2_1 \varphi \xi^2_1 + 2 \partial^2_{12} \varphi \xi_1
\xi_2 + \partial^2_2 \varphi \xi^2_2 + \partial^2_1 \varphi N^2_1 +
2
\partial^2_{12} \varphi N_1 N_2 + \partial^2_2 \varphi N^2_2 = \\ & =
\partial^2_1 \varphi N^2_2 + 2 \partial^2_{12} \varphi \left(
-\frac{N_2}{N_1}|N_1| \right) |N_1| + \partial^2_2 \varphi N^2_1 +
\partial^2_1 \varphi N^2_1 + 2 \partial^2_{12} \varphi N_1 N_2 +
\partial^2_2 \varphi N^2_2 = \\ & = \partial^2_1 \varphi - 2
\partial^2_{12} \varphi N_1 N_2 + \partial^2_2 \varphi + 2 \partial^2_{12}
\varphi N_1 N_2 = \Delta \varphi.
\end{aligned}
\end{equation*}
If $\xi_2=-|N_1|$, we get a similar result. Therefore
($\bigstar''$) is equivalent to \eqref{oper-2ord-9-28}. $\spadesuit$

\bigskip

\textbf{Example 2.} Let $\varphi = e^{\lambda \psi}$, where $\lambda
\in \mathbb{R}$ and $$|\nabla \psi (x)| \neq 0,\quad \forall x \in
\bar{\Omega}.$$ Let us look at whether there are any $\lambda$ for which \eqref{pseudconv-Delta} applies.

Let us calculate
$$\partial_j \varphi= \lambda e^{\lambda \psi} \partial_j \psi, $$
$$ \partial^2_{jk} \varphi = \lambda e^{\lambda \psi} \partial^2_{jk} \psi + \lambda^2 e^{\lambda \psi} \partial_j \psi \partial_k \psi. $$
 Hence \eqref{pseudconv-Delta} becomes
\begin{equation}
\label{(3-26.2)}
\begin{aligned} \begin{cases}
|\xi|^2 = \tau^2 \lambda^2 e^{2 \lambda \psi} |\nabla \psi|^2, \\
\\
\nabla \psi \cdot \xi = 0, \\
\\
\tau \neq 0,
\end{cases} \quad\Longrightarrow \quad Q_\lambda (x,\xi,\tau)>0,
\end{aligned}
\end{equation}
where

\begin{equation}
\label{Bal-1}
\begin{aligned}
Q_\lambda (x,\xi,\tau)&= \sum_{j,k=1}^{n}\lambda e^{\lambda \psi}
\partial^2_{jk} \psi\xi_j\xi_k + \tau^2 \sum_{j,k=1}^{n}\lambda^3 e^{3 \lambda
\psi}\partial^2_{jk} \psi \partial_j \psi
\partial_k \psi + \\ & +\lambda^2 e^{\lambda \psi}  (\nabla \psi\cdot\xi)^2 + \tau^2\lambda^4 e^{3 \lambda
\psi} |\nabla \psi|^2.
\end{aligned}
\end{equation}

\medskip

In order to  examine \eqref{(3-26.2)} let us suppose 
\begin{equation}\label{Bal-2}
\begin{cases}
|\xi|^2 = \tau^2 \lambda^2 e^{2 \lambda \psi} |\nabla \psi|^2, \\
\\
\nabla \psi \cdot \xi = 0, \\
\\
\tau \neq 0,
\end{cases}
\end{equation}
from which we have  $\nabla\psi\cdot \xi=0$. Hence,
if \eqref{Bal-2} holds true, we have
$$ Q_\lambda = \sum_{j,k=1}^{n} \lambda e^{\lambda \psi} \left(\partial^2_{jk} \psi \xi_j \xi_k +
\lambda^2 \tau^2 e^{2\lambda \psi}\partial^2_{jk} \psi \partial_j
\psi \partial_k \psi \right)+ \tau^2 \lambda^4 |\nabla \psi |^2 e^{3
\lambda \psi}.$$ Set
$$\eta = \frac{\xi}{\tau \lambda e^{\lambda \psi} |\nabla \psi|}.$$
In this way \eqref{Bal-2} it is rewritten as
\begin{equation*}
\begin{aligned} \begin{cases}
|\eta|^2 = 1 \\
\\
\nabla \psi \cdot \eta = 0 \\
\end{cases} \Longrightarrow \tilde{Q}_\lambda =  \sum_{j,k=1}^{n} \lambda \left(\partial^2_{jk} \psi \eta_j \eta_k +
\partial^2_{jk} \psi \frac{\partial_j \psi}{|\nabla \psi|} \frac{\partial_k \psi}{|\nabla \psi|} \right)+\lambda^2 >0.
\end{aligned}
\end{equation*}

It is clear now that if $|\lambda|$ is sufficiently large, then \eqref{(3-26.2)} is satisfied. More precisely, set
$$M= \max_{|\eta| =1} \left\vert \sum_{j,k=1}^{n}\partial^2_{jk} \psi \left( \eta_j \eta_k +
\frac{\partial_j \psi}{|\nabla \psi|} \frac{\partial_k \psi}{|\nabla
\psi|} \right) \right \vert,$$  we have that, if $$|\lambda|>M$$ then
condition \eqref{pseudconv-Delta} of Theorem \ref{teo2bis} is
satisfied. $\spadesuit$

\bigskip

 \textbf{Example 3.} Let us consider a radial function
 \begin{equation}\label{radial}
 \varphi(x)= f(|x|),\end{equation}
 in $\Omega= B_1 \setminus \overline{B_r}$, $r\in (0,1)$.
 
 We would like to find some functions $f$ such that $\varphi$
 satisfies condition \eqref{pseudconv-Delta}. We will see that in
 some cases this is not possible. Let us proceed in a similar manner
 to the previous example. Let us calculate
  $$\partial_j \varphi(x) =
\frac{x_j}{|x|} f'(|x|),$$
$$\partial^2_{jk} \varphi(x) = \left( \frac{\delta_{jk}}{|x|} - \frac{x_j x_k}{|x|^3} \right) f'(|x|) + \frac{x_j x_k}{|x|^2} f''(|x|),$$
$$\partial^2_{jk} \varphi(x) \xi_j \xi_k = \left( \frac{|\xi|^2}{|x|} - \frac{(\xi \cdot x)^2}{|x|^3} \right) f'(|x|) + \frac{(\xi \cdot x)^2}{|x|^2} f''(|x|),$$
\begin{equation*}
	\begin{aligned}
	\partial^2_{jk} \varphi(x) \partial_j \varphi(x) \partial_k \varphi(x) &= \left( \frac{|\nabla \varphi(x)|^2}{|x|} - \frac{(\nabla \varphi(x) \cdot x)^2}{|x|^3} \right) f'(|x|)
	+ \\&+\frac{(\nabla \varphi(x) \cdot x)^2}{|x|^2} f''(|x|).	
\end{aligned}
\end{equation*}

\medskip

 Since
$$|\nabla \varphi(x)|^2 = f'^2(|x|)$$ and  $$\nabla \varphi(x) \cdot
x= |x| f'(|x|)$$ we have
$$\partial^2_{jk} \varphi(x) \partial_j \varphi(x) \partial_k \varphi (x)  = (f'(|x|))^2 f''(|x|).$$
Then by \eqref{pseudconv-Delta} we can write  $Q$ as follows
$$Q = \left( \frac{|\xi|^2}{|x|} - \frac{(\xi \cdot x)^2}{|x|^3} \right) f'(|x|) + \frac{(\xi \cdot x)^2}{|x|^2} f''(|x|) + \tau^2 (f'(|x|))^2 f''(|x|).$$
Let us suppose that the antecedent of condition \eqref{pseudconv-Delta} holds, i.e., let us suppose that
$$\begin{cases}
|\xi|^2 = \tau^2 |\nabla \varphi(x)|^2=\tau^2 (f'(|x|))^2, \\
\\
\tau \nabla \varphi(x) \cdot \xi = \tau \frac{(\xi \cdot x)^2}{|x|^2} f'(|x|) = 0, \\
\\
 \tau \neq 0,
\end{cases}$$
namely
$$\begin{cases}
|\xi|^2 =\tau^2 (f'(|x|))^2, \\
\\
\xi \cdot x = 0,\\
\\
\tau\neq 0,
\end{cases}$$
Hence $Q$ can be written as follows
$$Q = \frac{\tau^2 (f'(|x|))^2}{|x|} f'(|x|) + \tau^2 (f'(|x|))^2 f''(|x|).$$
We get
$$\tau^{-2} Q = \frac{(f'(|x|))^2}{|x|} f'(|x|) + (f'(|x|))^2 f''(|x|).$$
We characterize the functions $$f:(0,1)\rightarrow (0,+\infty)$$ which satisfy
\eqref{pseudconv-Delta} and for which we have

$$\lim_{t\rightarrow 0} f(t)=+\infty$$ and $$f'<0, \ \ \mbox{in } (0,1).$$ The condition $Q>0$, for $t \in (0,1)$, becomes
$$Q = \frac{f'}{t} + f''>0.$$
To solve this differential inequality we set
$$f(t)= \psi (\log t),$$
from which we have
$$f'(t)= \psi'(\log t) \frac{1}{t},$$
$$f''(t)= \psi''(\log t) \frac{1}{t^2} - \psi'(\log t) \frac{1}{t^2}.$$
Hence the differential inequality can be written
$$\frac{f'}{t} + f'' = \frac{1}{t^2} \psi''(\log t)>0 \Longleftrightarrow \psi''(   \log t)>0.$$
Let $s=\log t$. Then, as $t \in (0,1)$, $s \in (-\infty ,
0)$, condition $Q>0$ becomes:
$$ \psi''(s)>0 \Longleftrightarrow \frac{d^2}{ds^2}(f(e^s)) >0 \; \; \; \forall s \in (-\infty,0).$$

Let us observe that there are functions that do not satisfy this condition.
For $\alpha >0$, we consider functions $f (t)$ of the
type

$$f(t)= \left( \log \frac{1}{t} \right)^{\alpha}.$$
Let us calculate
$$f'(t)= \alpha \left( -\frac{1}{t} \right) \left( \log \frac{1}{t} \right)^{\alpha-1},$$
$$f'(t)= \frac{\alpha }{t^2} \left( \log \frac{1}{t} \right)^{\alpha-1} + \frac{\alpha (\alpha -1)}{t^2} \left( \log \frac{1}{t} \right)^{\alpha-2},$$
from which it follows that the condition on $Q$ can be written as
$$Q = \frac{f'}{t} + f''= \frac{\alpha (\alpha -1)}{t^2} \left( \log \frac{1}{t} \right)^{\alpha-2}>0.$$
Consequently if $\alpha <1$, we have
$$ \frac{f'}{t} + f'' < 0.$$
If $\alpha=1$, we have
$$\varphi= \log \frac{1}{|x|},$$
Therefore $$Q=0.$$ $\spadesuit$

\bigskip

\underline{\textbf{Exercise.}} Prove that if $U\in H^2(B_1)$ is a 
solution to the equation
$$\Delta U=b(x)\cdot \nabla U+c(x) U=0,\quad \mbox{in } B_1,$$
where $b\in L^{\infty}(B_1;\mathbb{R}^n)$ and $c\in
L^{\infty}(B_1)$ then the following three sphere inequalityT holds true

\begin{equation}\label{esercizio-tre-sfere}
\left\Vert U \right\Vert_{L^2(B_{\varrho}(0))}\leq C \left\Vert U
\right\Vert_{L^2(B_{r})}^{\vartheta}\left\Vert U
\right\Vert_{L^2(B_{1})}^{1-\vartheta},
\end{equation}
for $0<r\leq \varrho \leq C^{-1}$, where $C\leq 1$ e $\vartheta\in
(0,1)$ are constants  depending on $\left\Vert b
\right\Vert_{L^{\infty}(B_{1};\mathbb{R}^n)}$, $\left\Vert c
\right\Vert_{L^{\infty}(B_{1})}$, on  $\varrho$ and on $r$.

[Hint: apply  Theorem \ref{teo2bis} where $\varphi$
is a suitable radial function.] $\clubsuit$

\section[Second order operators I]{Second order operators I -- constant coefficients in the principal part}\label{24-4-23:second-order}
We now consider a more general case related to the operator
\eqref{oper-2ord-1}. Let $M_0$ and let us $M_1$ be positive numbers given and
let us assume that

\begin{subequations}
\label{oper-2ord-101}
\begin{equation}
\label{oper-2ord-101a} \left\Vert
g^{jk}\right\Vert_{L^{\infty}(\Omega)}\leq M_0, \quad \mbox{for }
j,k=1,\cdots, n,
\end{equation}
\begin{equation}
\label{oper-2ord-101b} \left\vert g^{jk}(x)-g^{jk}(y)\right\vert\leq
M_1|x-y|,\quad \mbox{for } j,k=1,\cdots, n,\mbox{ } \forall x,y\in
\Omega.
\end{equation}
\end{subequations}
Let us note that with conditions \eqref{oper-2ord-101a} and
\eqref{oper-2ord-101b}, operator \eqref{oper-2ord-1} is not
necessarily elliptic. As we noted in Section
\ref{Carleman-intro}, to establish a Carleman estimate for
operator \eqref{oper-2ord-1}, under assumption
\eqref{oper-2ord-101b}, is equivalent to establish a Carleman estimate
for the operator

\begin{equation}
\label{oper-2ord-102}
L_gu=\partial_j\left(g^{jk}(x)\partial_ku\right).
\end{equation}
As a matter of fact we have
$$(L_g-P_2)u=\partial_j\left(g^{jk}(x)\right)\partial_ku,$$
which is a first order operator with bounded coefficients.

We begin by establishing an identity analogous to \eqref{rellich}. To this purpose we introduce some notations. We set
$$\xi^{(g)}=\left\{g^{jk}(x)\xi_k \right\}_{j=1}^n$$
and for any function $v$, sufficiently regular, we set
$$\nabla^{(g)}v=\left\{g^{jk}(x)\partial_k v\right\}_{j=1}^n.$$
We further set,

$$\mathbf{g}(\xi,\eta)=g^{jk}(x)\xi_j\eta_k.$$
Using these notations we have

\begin{equation}
\label{oper-2ord-103} L_gu=\mbox{div}\left(\nabla^{(g)}u\right).
\end{equation}
Let us notice that if  $g^{jk}(x)=\delta^{jk}$ then
$\mathbf{g}(\xi,\eta)=\xi\cdot\eta$ and $L_g=\Delta$ and, if
$$\left\{g^{jk}(x)\right\}_{j,k=1}^n=\mbox{diag }(1,\cdots,1,-1),$$
then $$\mathbf{g}(\xi,\eta)=\xi'\cdot\eta'-\xi_n\eta_n$$ and
$$L_gu=\Delta_{x'}u-\partial^2_nu=\square u.$$

\bigskip

\begin{lem}[\textbf{generalized Rellich identity}]\label{re91221}
	\index{Lemma:@{Lemma:}!- generalized Rellich identity@{- generalized Rellich identity}}
Let $\beta \in C^{0,1} (\Omega, \mathbb{R}^n)$,
$\beta=\left(\beta^1,\cdots, \beta^n \right)$ and $v \in C^2
(\Omega)$, then we have

\begin{equation}\label{oper-2ord-104}
\begin{aligned}
2\mathbf{g}\left(\beta,\nabla v\right)L_gv&= \mbox{div }
\left(2\mathbf{g}\left(\beta,\nabla v\right)\nabla
v-\mathbf{g}\left(\nabla v,\nabla v\right)\beta \right)+\\&+
(\mbox{div } \beta) \mathbf{g}\left(\nabla v,\nabla
v\right)-2\partial_l\beta^kg^{lj}\partial_jv\partial_kv+\\&+\beta^k\left(\partial_k
g^{lj}\right)\partial_lv\partial_jv,\quad \mbox{a.e. } x \in
\Omega
\end{aligned}
\end{equation}
and

\begin{equation}\label{oper-2ord-104-24-4-23}
	\begin{aligned}
		2\int_{\Omega}\mathbf{g}\left(\beta,\nabla v\right)L_gvdx &= \int_{\Omega}
		(\mbox{div } \beta) \mathbf{g}\left(\nabla v,\nabla
		v\right)-2\partial_l\beta^kg^{lj}\partial_jv\partial_kv+\\&+\int_{\Omega}\beta^k\left(\partial_k
		g^{lj}\right)\partial_lv\partial_jv dx,
	\end{aligned}
\end{equation}
\end{lem}

\textbf{Proof.} First, we observe that identity \eqref{oper-2ord-104-24-4-23} is am immediate consequence of \eqref{oper-2ord-104} after its integration over $\Omega$. Hence, it suffices to prove \eqref{oper-2ord-104}.

We have

\begin{equation}\label{oper-2ord-105}
\begin{aligned}
2\mathbf{g}\left(\beta,\nabla v\right)L_gv&=2(\beta^k\partial_k
v)\partial_l\left(g^{lj}\partial_jv\right)=\\&= 2\partial_l
\left(\beta^k\partial_kvg^{lj}\partial_jv\right)-2\partial_l\left(\beta^k\partial_kv\right)g^{lj}\partial_jv=\\&=
\mbox{div } \left(2\mathbf{g}\left(\beta,\nabla v\right)\nabla^{(g)}
v\right)-2\partial_l\left(\beta^k\partial_kv\right)g^{lj}\partial_jv=\\&=
\mbox{div } \left(2\mathbf{g}\left(\beta,\nabla
v\right)\nabla^{(g)}v\right)-2\left(\partial_l\beta^kv\right)g^{lj}\partial_kv\partial_jv-\\&-
2\beta^kg^{lj}\partial_{lk}^2v\partial_jv.
\end{aligned}
\end{equation}
Now, we notice that
\begin{equation*}
	\begin{aligned}
		\partial_k\left(g^{lj}\partial_{l}v\partial_jv\right)&=g^{lj}\partial^2_{lk}v\partial_jv+g^{lj}\partial_{l}v\partial^2_{jk}v+\partial_k\left(g^{lj}\right)\partial_{l}v\partial_jv=\\&=2g^{lj}\partial_{lk}^2v\partial_jv+\partial_k\left(g^{lj}\right)\partial_{l}v\partial_jv,
	\end{aligned}
\end{equation*}
from which we have
$$2g^{lj}\partial_{lk}^2v\partial_jv=
\partial_k\left(g^{lj}\partial_{l}v\partial_jv\right)-\partial_k\left(g^{lj}\right)\partial_{l}v\partial_jv.$$
Using this identity, we transform the last term on the
right--hand side of \eqref{oper-2ord-105}. We have

\begin{equation*}
\begin{aligned}
- 2\beta^kg^{lj}\partial_{lk}^2v\partial_jv&=-\beta^k
\partial_k\left(g^{lj}\partial_{l}v\partial_jv\right)+\beta^k\partial_kg^{lj}\partial_{l}v\partial_jv=\\&=
-\partial_k\left(\beta^kg^{lj}\partial_{l}v\partial_jv\right)+(\mbox{div
}\beta)g^{lj}\partial_{l}v\partial_jv+\beta^k\partial_kg^{lj}\partial_{l}v\partial_jv.
\end{aligned}
\end{equation*}
Hence

\begin{equation*}
\begin{aligned}
- 2\beta^kg^{lj}\partial_{lk}^2v\partial_jv&=-\mbox{div
}\left(\mathbf{g}\left(\nabla v,\nabla
v\right)\right)+\\&+(\mbox{div }\beta)\mathbf{g}\left(\nabla
v,\nabla v\right)+\beta^k\partial_kg^{lj}\partial_{l}v\partial_jv
\end{aligned}
\end{equation*}
and using the just obtain equality in \eqref{oper-2ord-105} we obtain
\eqref{oper-2ord-104}. $\blacksquare$

\bigskip

In the present Subsection we consider the \textbf{case of constant coefficients}. In this case 
the operator is given by
\begin{equation}\label{oper-2ord-80122}
P_2(\partial)=g^{jk}\partial_{jk}^2,
\end{equation}
where $g^{jk}=g^{kj}$ are \textbf{real constants}. We begin by giving
some definitions whose geometric meaning will be explained
later (see Section \ref{cond-nec-Carlm}).

\begin{definition}[\textbf{pseudo--convex functions}]\label{psedoconvess}
	\index{Definition:@{Definition:}!- pseudo--convex functions@{- pseudo--convex functions}} Let $\Omega$ be a bounded open set of $\mathbb{R}^n$
and let $\phi\in C^2\left(\overline{\Omega}\right)$ satisfy
\begin{equation}\label{oper-2ord-106-28}
\nabla \phi(x)\neq 0,\quad\quad \forall x\in \overline{\Omega}.
\end{equation}
We say that $\phi$ is \textbf{pseudo--convex w.r.t. the operator} \eqref{oper-2ord-80122} in the point $x\in
\overline{\Omega}$, if we have

\begin{equation}\label{oper-2ord-106}
\begin{aligned}
\begin{cases}
P_2(\xi)=0, \\
\\
P_2^{(j)}(\xi)\partial_j\phi(x)= 0,\\
\\
\xi\neq 0,
\end{cases}\mbox{ }
\Longrightarrow \mbox{ }&
\partial^2_{jk}\phi(x)P_2^{(j)}(\xi)P_2^{(k)}(\xi)>0.
\end{aligned}
\end{equation}
We say that $\phi$ is \textbf{pseudo--convex w.r.t. operator} $P_2(\partial)$ if \eqref{oper-2ord-106} holds true
for every $x\in \overline{\Omega}$.
\end{definition}

\bigskip

\textbf{Remark.} Using the notations introduced above, \eqref{oper-2ord-106} can be written

\begin{equation}\label{oper-2ord-106-bis}
\begin{aligned}
\begin{cases}
\mathbf{g}(\xi,\xi)=0,\\
\\
\mathbf{g}(\xi,\nabla\phi(x))= 0,\\
\\
\xi\neq 0,
\end{cases}\mbox{ }
\Longrightarrow \mbox{ }&
\partial^2\phi(x)\xi^{(g)}\cdot\xi^{(g)}>0.
\end{aligned}
\end{equation}
Let us notice that \textbf{if the matrix
$\left\{g^{jk}\right\}_{j,k=1}^n$ is singular do not exist any
pseudo--convex functions} because there exists $\xi \in \mathbb{R}^n\setminus\left\{0\right\}$ such that $\xi^g=g\cdot \xi=0$. It should also be noticed that if the operator
$P_2(\partial)$ is elliptic, then condition \eqref{oper-2ord-106} is
trivially satisfied since the antecedent is false.
$\blacklozenge$

\bigskip

\begin{definition}[\textbf{strong pseudo--convex functions}]\label{forte-psedoconvess}
	\index{Definition:@{Definition:}!- strong pseudo--convex functions@{- strong pseudo--convex functions}} Let $\Omega$, $\phi$ and $P_2(\partial)$ be as in
Definition \ref{psedoconvess}. We say that $\phi$ is
\textbf{strongly pseudo--convex w.r.t. the operator $P_2$}
in the point $x\in \overline{\Omega}$ if $\phi$ is pseudo--convex w.r.t.  $P_2$ and, further, we have

\begin{equation}\label{oper-2ord-107}
\begin{aligned}
&\begin{cases}
P_2(\xi+i\tau\nabla \phi (x))=0, \\
\\
P_2^{(j)}(\xi+i\tau \nabla \phi (x))\partial_j\phi(x)= 0,\\
\\
\tau\neq 0,
\end{cases}\mbox{ }\Longrightarrow\\&
\quad\quad\quad\quad\Longrightarrow \mbox{ }
\partial^2_{jk}\phi(x)P^{(j)}_2(\xi+i\tau\nabla \phi (x))\overline{P^{(k)}_2(\xi+i\tau\nabla \phi (x))}>0.
\end{aligned}
\end{equation}
We say that $\phi$ is \textbf{strongly pseudo--convex w.r.t. the operator} $P_2(\partial)$, if it is strongly pseudo--convex w.r.t. the operator $P_2(\partial)$ in each point
$x\in \overline{\Omega}$.

\end{definition}

\bigskip

\textbf{Remarks.}

\medskip

\noindent \textbf{1.}  With the notations introduced above,  \eqref{oper-2ord-107} can be written

\begin{equation}\label{oper-2ord-108}
\begin{aligned}
&\begin{cases}
\mathbf{g}(\xi,\xi)=\tau^2\mathbf{g}(\nabla\phi(x),\nabla\phi(x)), \\
\\
\mathbf{g}(\xi,\nabla\phi(x))= 0,\\
\\
\mathbf{g}(\nabla\phi(x),\nabla\phi(x))=0,\\
\\
\tau\neq 0
\end{cases}\mbox{ }
\Longrightarrow \mbox{ }\\& \quad\quad\quad\Longrightarrow \mbox{ }
\partial^2\phi(x)\xi^{(g)}\cdot\xi^{(g)}+\tau^2\partial^2\phi(x)\nabla^{(g)}\phi(x)\cdot\nabla^{(g)}\phi(x)>0.
\end{aligned}
\end{equation}

\medskip

\noindent \textbf{2.} As we will easily check, \textbf{in the case of the real coefficients we are considering, the definitions of pseudo--convexity and strong pseudo--convexity are equivalent}. 
As a matter of fact, if $\phi$ is strongly pseudo-convex
with respect to $P_2(\partial)$, it is trivially pseudo--convex. We prove
that if $\phi$  is pseudo--convex then it is strongly
pseudo--convex.
 Let us suppose, hence, that \eqref{oper-2ord-106}
is satisfied in $x_0\in \overline{\Omega}$ and let us prove that
\eqref{oper-2ord-108} is satisfied in $x_0$. If
$\mathbf{g}(\nabla\phi(x_0),\nabla\phi(x_0))\neq 0$, then 
\eqref{oper-2ord-108} is trivially satisfied as
the antecedent of the implication \eqref{oper-2ord-108} is false.
If, on the other hand, we have

$$\mathbf{g}(\nabla\phi(x_0),\nabla\phi(x_0))= 0,$$ then the antecedent of condition
\eqref{oper-2ord-108} becomes (in $x_0$)

\begin{equation}\label{oper-2ord-109}
\begin{aligned}
\begin{cases}
\mathbf{g}(\xi,\xi)=0, \\
\\
\mathbf{g}(\xi,\nabla\phi(x_0))= 0,\\
\\
\mathbf{g}(\nabla\phi(x_0),\nabla\phi(x_0))=0,\\
\\
\tau\neq 0.
\end{cases}
\end{aligned}
\end{equation}

Now, by the first two conditions in \eqref{oper-2ord-109} and the
pseudo--convexity of $\phi$ we have

\begin{equation}\label{oper-2ord-110}
\partial^2\phi(x_0)\xi^{(g)}\cdot\xi^{(g)}>0.
\end{equation}
Moreover, since $\mathbf{g}(\nabla\phi(x_0),\nabla\phi(x_0))=0$
and $\nabla\phi(x_0)\neq 0$, setting $\xi_0=\nabla\phi(x_0)$ we get by \eqref{oper-2ord-106} trivially

$$\begin{cases}
\mathbf{g}(\xi_0,\xi_0)=0, \\
\\
\mathbf{g}(\xi_0,\nabla\phi(x_0))= 0,\\
\\
\xi_0\neq 0,
\end{cases}$$
hence by \eqref{oper-2ord-106-bis} we have 
$$\partial^2\phi(x_0)\nabla^{(g)}\phi(x_0)\cdot\nabla^{(g)}\phi(x_0)>0$$ and taking into account that (by
\eqref{oper-2ord-109}) $\tau\neq 0$, we have

\begin{equation}\label{oper-2ord-111}
\tau^2\partial^2\phi(x_0)\nabla^{(g)}\phi(x_0)\cdot\nabla^{(g)}\phi(x_0)>0.
\end{equation}
Now, by \eqref{oper-2ord-110} and \eqref{oper-2ord-111} we have
$$\partial^2\phi(x_0)\xi^{(g)}\cdot\xi^{(g)}+\tau^2\partial^2\phi(x_0)\nabla^{(g)}\phi(x_0)\cdot\nabla^{(g)}\phi(x_0)>0.$$
Hence, we have proved that \eqref{oper-2ord-108} is
satisfied and, therefore, we have the equivalence of the
definitions \ref{psedoconvess} and \ref{forte-psedoconvess}

Let us note that \textbf{in the elliptic case}, each $\phi\in
C^2\left(\overline{\Omega}\right)$ such that
\begin{equation*}
\nabla \phi(x)\neq 0,\quad\quad \forall x\in \overline{\Omega}
\end{equation*}
is trivially pseudo--convex (hence, it is strongly
pseudo--convex). $\blacklozenge$

\bigskip

\noindent \textbf{\underline{Warning about definitions \ref{psedoconvess} and
\ref{forte-psedoconvess}.}} It is important to point out that generally the
definitions of pseudo--convexity and strong pseudo--convexity
are referred to the \textbf{level surfaces}
$\left\{\phi(x)=\phi(x_0)\right\}$, where $x_0\in \Omega$. And using the
term "surface" we also want to emphasize the
invariant character of the definition (see \cite[\S 8.6]{HO63}).
Our modification of the terminology is only due to the purpose of
to lighten the exposition a little. $\blacktriangle$

\bigskip

To obtain the Carleman estimate of Theorem
\ref{Carl-2ord-costanti} (see below) we need a condition more
stringent than the strong pseudo-convexity. This condition, has the same form
of condition $(\bigstar)$ of Theorem \ref{Carlm-ell-teor}.

\begin{definition}\label{condizione-S} 
	\index{Definition:@{Definition:}!- condition \textbf{S}@{- condition \textbf{S}}}
	Let $\Omega$, $\phi$ and $P_2$ as in Definition \ref{forte-psedoconvess}
We say that $\phi$ satisfies condition $(\mathbf{S})$ w.r.t. the operator $P_2(\partial)$, if $\phi$ is pseudo--convex
w.r.t. $P_2(\partial)$ and we have

\begin{equation}\label{oper-2ord-112}
\begin{aligned}
&\begin{cases}
P_2(\xi+i\tau\nabla \phi (x))=0, \\
\tau\neq 0,
\end{cases}\mbox{ }\Longrightarrow\\&
\quad\quad\quad\quad\Longrightarrow \mbox{ }
\partial^2_{jk}\phi(x)P^{(j)}_2(\xi+i\tau\nabla \phi (x))\overline{P^{(k)}_2(\xi+i\tau\nabla \phi (x))}>0.
\end{aligned}
\end{equation}
\end{definition}

\bigskip

\textbf{Remark.} With the notations introduced above
\eqref{oper-2ord-107},  we write  condition \eqref{oper-2ord-112} as follows

\begin{equation}\label{oper-2ord-113}
\begin{aligned}
&\begin{cases}
\mathbf{g}(\xi,\xi)=\tau^2\mathbf{g}(\nabla\phi(x),\nabla\phi(x)), \\
\\
\mathbf{g}(\xi,\nabla\phi(x))= 0,\\
\\
\tau\neq 0,
\end{cases}
\Longrightarrow \\& \Longrightarrow\mbox{ }
Q_{\phi}:=\partial^2\phi(x)\xi^{(g)}\cdot\xi^{(g)}+\tau^2\partial^2\phi(x)\nabla^{(g)}\phi(x)\cdot\nabla^{(g)}\phi(x)>0.
\end{aligned}
\end{equation} $\blacklozenge$

\bigskip

It is evident that if $\phi$ satisfies condition $(\mathbf{S})$
then it is strongly pseudo-convex. However, the
converse does not hold. Let us consider, for instance, the function
$\phi(x)=\log\frac{1}{|x|}$; this function is strongly
pseudo--convex, but, as we saw in \textbf{Example 3} of
this Section, it does not satisfy condition $(\mathbf{S})$.

\medskip

However the following Proposition holds

\begin{prop}\label{oper-2ord-114}
 Let $\Omega$ be a bounded open set of $\mathbb{R}^n$ and let $\phi\in
C^2\left(\overline{\Omega}\right)$ strongly pseudo--convex
w.r.t. operator \eqref{oper-2ord-80122}, then
 $\varphi=e^{\lambda \phi}$ satisfies condition
$(\mathbf{S})$ if $\lambda$ is large enough.
\end{prop}
\textbf{Proof.} If  $P_2(\partial)$ is
elliptic we can readily reduce to what we have done in Example 2. Hence let us suppose that $P_2(\partial)$ is \textbf{not elliptic}. We first prove that $\varphi=e^{\lambda \phi}$
is pseudo--convex. Let us calculate
$$\nabla \varphi= \lambda e^{\lambda \phi} \nabla \phi, $$
$$ \partial^2_{jk} \varphi = \lambda e^{\lambda \phi} \partial^2_{jk} \phi + \lambda^2 e^{\lambda \phi} \partial_j \phi \partial_k \phi,\quad j,k=1,\cdots, n $$
and
\begin{equation}\label{oper-2ord-334}
\begin{aligned}
&\partial^2_{jk}\varphi(x)P_2^{(j)}(\xi)P_2^{(k)}(\xi)=\\&=\lambda e^{\lambda\phi}
\left[\partial^2_{jk}\varphi(x)P_2^{(j)}(\xi)P_2^{(k)}(\xi)+\lambda\left|P_2^{(j)}(\xi)\partial_j\phi\right|^2\right].
\end{aligned}
\end{equation}
Set

$$X=\overline{\Omega}\times\left\{\xi\in \mathbb{R}^n: \quad P_2(\xi)=0, \quad
|\xi|=1\right\}.$$ Since $P_2(\partial)$ is not ellptic we have
that $X$ is a nonempty compact subset of $\mathbb{R}^{2n}$. Now,
by the definition of pseudo--convexity in $\overline{\Omega}$
we have

\begin{equation*}
\begin{aligned}
(x,\xi)\in X,\quad P_2^{(j)}(\xi)\partial_j\phi(x)=0\mbox{ }
\Longrightarrow\mbox{
}\partial^2_{jk}\phi(x)P_2^{(j)}(\xi)P_2^{(k)}(\xi)>0
\end{aligned}
\end{equation*}
and this, by Lemma \ref{106-CauNire}, implies that there exists
$\lambda_0>0$ such that
\begin{equation}\label{oper-2ord-80122-1}
\lambda_0
\left|P_2^{(j)}(\xi)\partial_j\phi(x)\right|^2+\partial^2_{jk}\phi(x)P_2^{(j)}(\xi)P_2^{(k)}(\xi)>0,
\end{equation}
for every $(x,\xi)\in X$. 

Inequality \eqref{oper-2ord-80122-1}, in turn, implies (taking into account that the polinomial (in $\xi$) on the left--hand side is homogeneous of degree $2$)

\begin{equation}\label{oper-2ord-335}
\lambda
\left|P_2^{(j)}(\xi)\partial_j\phi(x)\right|^2+\partial^2_{jk}\phi(x)P_2^{(j)}(\xi)P_2^{(k)}(\xi)>0
\end{equation}
for every $x\in \overline{\Omega}$, for every $\xi\in \mathbb{R}^n\setminus\{0\}$ and for every $\lambda\geq \lambda_0$. By \eqref{oper-2ord-334} and
\eqref{oper-2ord-335} we have that 
\eqref{oper-2ord-106} is satisfied by the function $\varphi$.

\medskip

Now we prove \eqref{oper-2ord-112}. Let us introduce the
following notation. For every $\xi\in \mathbb{R}^n$, for every $t\in
\mathbb{R}$ and for every function $f\in
C^2\left(\overline{\Omega}\right)$ such that $\nabla f\neq 0$, in $
\overline{\Omega}$, set

$$\zeta_{t,f}=\xi+it\nabla f.$$
Similarly to \eqref{oper-2ord-334} we have

\begin{equation}\label{oper-2ord-336}
\begin{aligned}
&\partial^2_{jk}\varphi
P_2^{(j)}(\zeta_{\tau,\varphi})\overline{P_2^{(k)}(\zeta_{\tau,\varphi})}=\\&=\lambda
e^{\lambda\phi}
\left[\partial^2_{jk}\phi(x)P_2^{(j)}(\zeta_{\tau\lambda,\phi})\overline{P_2^{(k)}(\zeta_{\tau\lambda,\phi})}+
\lambda\left|P_2^{(j)}(\zeta_{\tau\lambda,\phi})\partial_j\phi\right|^2\right].
\end{aligned}
\end{equation}
Set

\begin{equation}\label{oper-2ord-337}
X_1=\left\{(x,\xi,t)\in \overline{\Omega}\times\mathbb{R}^n: \quad
P_2(\zeta_{t,\phi(x)})=0, \quad |\zeta_{t,\phi(x)}|=1\right\}.
\end{equation}
It turns out that $X_1\neq \emptyset$ and that $X_1$ is a compact of
$\mathbb{R}^{2n+1}$ (since $\nabla\phi(x)\neq 0$ in
$\overline{\Omega}$). We now check that we have
\begin{equation}\label{oper-2ord-338}
\begin{aligned}
(x,\xi,t)\in X_1,\quad&
P_2^{(j)}(\zeta_{t,\phi(x)})\partial_j\phi(x)=0\mbox{ }
\Longrightarrow\\&\mbox{ }\Longrightarrow \mbox{ }
\partial^2_{jk}\phi(x)P_2^{(j)}(\zeta_{t,\phi(x)})\overline{P_2^{(k)}(\zeta_{t,\phi(x)})}>0.
\end{aligned}
\end{equation}
As a matter of fact, for $t=0$, \eqref{oper-2ord-338} is nothing but
\eqref{oper-2ord-106} and so it is satisfied. Now, if $t\neq 0$ then
\eqref{oper-2ord-338} is satisfied because $\phi$ is strongly
pseudo--convex.  Hence by Lemma \ref{106-CauNire},
there exists $\lambda_1\geq\lambda_0$ such that

\begin{equation*}
\begin{aligned}
\lambda_1
\left|P_2^{(j)}(\zeta_{t,\phi(x)})\partial_j\phi(x)\right|^2+\partial^2_{jk}\phi(x)P_2^{(j)}(\zeta_{t,\phi(x)})\overline{P_2^{(k)}(\zeta_{t,\phi(x)}})>0,
\end{aligned}
\end{equation*}
for every $(x,\xi,t)\in X_1$, which in turn implies
\begin{equation}\label{oper-2ord-339}
\begin{aligned}
\lambda
\left|P_2^{(j)}(\zeta_{t,\phi(x)})\partial_j\phi(x)\right|^2+\partial^2_{jk}\phi(x)P_2^{(j)}(\zeta_{t,\phi(x)})\overline{P_2^{(k)}(\zeta_{t,\phi(x)}})>0,
\end{aligned}
\end{equation}
for every $\lambda\geq \lambda_1$, for every $x\in \overline{\Omega}$
and for every $(\xi,t)\in \mathbb{R}^n\setminus\{(0,0)\}$. Trivially 
\eqref{oper-2ord-339} holds true for $t=\tau\lambda e^{\lambda \phi(x)}$, $\tau\neq 0$.
Hence, recalling \eqref{oper-2ord-336}, we have that $\varphi$
satisfies the condition

\begin{equation}\label{oper-2ord-340}
\begin{aligned}
&\begin{cases}
P_2(\xi+i\tau\nabla \varphi (x))=0, \\
\\
\tau\neq 0,
\end{cases}\mbox{ }\Longrightarrow\\&
\quad\quad\quad\quad\Longrightarrow \mbox{ }
\partial^2_{jk}\varphi(x)P^{(j)}_2(\xi+i\tau\nabla \varphi (x))\overline{P^{(k)}_2(\xi+i\tau\nabla \varphi (x))}>0.
\end{aligned}
\end{equation}
The proof is concluded. $\blacksquare$

\bigskip

In the sequel we will need some notations and  a Lemma.

\medskip

Let $k\in \mathbb{R}$ and let $N\in \mathbb{R}^n\setminus \{0\}$.
For any $f\in C^{\infty}_0(\Omega)$, let us denote by 

\begin{equation}\label{oper-2ord-201}
\left\Vert
f\right\Vert^2_{k,\tau}=\frac{1}{(2\pi)^n}\int_{\mathbb{R}^n}\left\vert
\xi+i\tau N\right\vert^{2k}\left\vert
\widehat{f}(\xi)\right\vert^2d\xi.
\end{equation}

\medskip

\textbf{Remark.} By \eqref{oper-2ord-201} and by the 
Parseval identity we have

\begin{equation}\label{oper-2ord-203}
\left\Vert f\right\Vert_{-1,\tau}\leq |\tau N|^{-1}\left\Vert
f\right\Vert_{L^2(\Omega)},
\end{equation}
for every $f\in C^{\infty}_0(\Omega)$. $\blacklozenge$

\bigskip

The following Lemma holds
\begin{lem}\label{Lemma-lipsch}
Let $\rho>0$ and let $h$ be a Lipschitz continuous function defined in $B_{\rho}(x_0)$. Let $$A=[h]_{0,1,\overline{B_{\rho}(x_0)}},$$
the Lipschitz constant of $h$. Let us suppose that
$$h(x_0)=0.$$
We have, for every $w\in C_0^{\infty}(B_{\rho}(x_0))$

\begin{equation}\label{oper-2ord-202}
\left\Vert h(\partial_jw-\tau N_jw)\right\Vert_{-1,\tau}\leq
A\left(\rho+|\tau N|^{-1}\right) \left\Vert
w\right\Vert_{L^2(B_{\rho}(x_0))}.
\end{equation}
\end{lem}
\textbf{Proof.} We have

$$h(x)(\partial_jw-\tau N_jw)=(\partial_j-\tau
N_j)(hw)-w\partial_jh.$$ Hence

\begin{equation}\label{oper-2ord-300}
\left\Vert h(\partial_j-\tau N_jw)\right\Vert_{-1,\tau}\leq
\left\Vert (\partial_j-\tau N_j)(hw)\right\Vert_{-1,\tau}+\left\Vert
w\partial_jh\right\Vert_{-1,\tau}.
\end{equation}
Let us notice that
\begin{equation}\label{oper-2ord-299}
\left\vert h(x)\right\vert= \left\vert h(x)-h(x_0)\right\vert\leq
A\rho,\quad \forall x\in B_{\rho}(x_0).
\end{equation}
Now, the triangle inequality and \eqref{oper-2ord-299}
yield
\begin{equation}\label{oper-2ord-301}
\begin{aligned}
&\left\Vert (\partial_jw-\tau
N_j)(hw)\right\Vert^2_{-1,\tau}=\\&=\frac{1}{(2\pi)^n}\int_{\mathbb{R}^n}\frac{\left\vert
\xi_j+i\tau N_j\right\vert^2\left\vert
\widehat{hw}(\xi)\right\vert^2}{\left\vert \xi+i\tau
N\right\vert^2}d\xi\leq \\&\leq
\frac{1}{(2\pi)^n}\int_{\mathbb{R}^n}\left\vert
\widehat{hw}(\xi)\right\vert^2d\xi=\\&=\int_{\mathbb{R}^n}
\left\vert h(x)w(x)\right\vert^2 dx\leq A^2\rho^2\left\Vert
w\right\Vert^2_{L^2(B_{\rho}(x_0))}.
\end{aligned}
\end{equation}
On the other hand
\begin{equation}\label{oper-2ord-302}
\begin{aligned}
\left\Vert
w\partial_jh\right\Vert^2_{-1,\tau}&=\frac{1}{(2\pi)^n}\int_{\mathbb{R}^n}\frac{\left\vert\widehat{w\partial_j
h}(\xi)\right\vert^2}{\left\vert \xi+i\tau N\right\vert^2}d\xi\leq \\&\leq \frac{|\tau
N|^{-2}}{(2\pi)^n}\left\Vert
\widehat{w\partial_jh}\right\Vert^2_{L^2(\mathbb{R}^n)}=\\&= |\tau
N|^{-2}\left\Vert
w\partial_jh\right\Vert^2_{L^2(B_{\rho}(x_0))}\leq\\&\leq A^2|\tau
N|^{-2}\left\Vert w\right\Vert^2_{L^2(B_{\rho}(x_0))}.
\end{aligned}
\end{equation}
Therefore, by  \eqref{oper-2ord-300}, \eqref{oper-2ord-301} and
\eqref{oper-2ord-302} we get \eqref{oper-2ord-202}.
$\blacksquare$

\bigskip

We now state and prove the following

\begin{theo}\label{Carl-2ord-costanti} Let $\Omega$ be a bounded open set of
$\mathbb{R}^n$ and let $$P_2(\partial)=g^{jk}\partial_{jk}^2,$$ where
the matrix  $\left\{g^{jk}\right\}_{j,k=1}^n$ is
real, constant and symmetric. Let us suppose that  $\varphi \in
C^{2}(\bar{\Omega})$ satisfies condition $(\mathbf{S})$ w.r.t.  $P_2(\partial)$.

Then there exist constants $C$ and
$\tau_0$ such that
\begin{equation}\label{oper-2ord-115}
\begin{aligned}
& \tau^3 \int_{\Omega} |u|^2 e^{2 \tau \varphi}dx + \tau
\int_{\Omega}
 |\nabla u|^2 e^{2 \tau \varphi}dx \leq C
\int_{\Omega} |P_2 (\partial)u|^2 e^{2 \tau \varphi}dx,
\end{aligned}
\end{equation}
for every $u \in C^{\infty}_0(\Omega)$ and for every $\tau\geq \tau_0$.
\end{theo}

\textbf{Proof.} Let us denote
$$L=P_2(\partial).$$ 
Let $u\in C^{\infty}_0(\Omega)$ and set $v=
e^{\tau \varphi} u$. We have
\begin{equation*}
\begin{aligned}
&P_{2, \tau}v:=e^{\tau \varphi} L(e^{-\tau
\varphi} v)= \\&=L(v)
+\tau^2\mathbf{g}(\nabla\varphi,\nabla\varphi) v - 2 \tau
\mathbf{g}(\nabla\varphi,\nabla v)-\tau L(\varphi) v.
\end{aligned}
\end{equation*}
let us define
$$S_\tau v=L(v)  +
\tau^2\mathbf{g}(\nabla\varphi,\nabla\varphi) v$$ and
$$A_{\tau}v=- 2 \tau \mathbf{g}(\nabla\varphi,\nabla v)$$
and set
\begin{equation}\label{oper-2ord-303}
p_2(x,\partial,\tau) v=S_\tau v+ A_{\tau}v,
\end{equation}
Let us note that
\begin{equation}\label{oper-2ord-116}
P_{2,\tau} v-p_2(x,\partial,\tau) v=- \tau
L(\varphi) v.
\end{equation}

In what follows, for the sake of brevity, we will omit the domain and the element of
integration.

Let $\mu$ be a positive number  which we will choose later, 
we have trivially

\begin{equation}\label{oper-2ord-117}
\begin{aligned}
&\int \left|p_2(x,\partial,\tau) v\right|^2e^{2\mu\varphi}=\\&= \int
|S_\tau v|^2 e^{2\mu\varphi}+ \int |A_\tau v|^2e^{2\mu\varphi} +\\&+
2 \int \left(S_\tau v A_\tau v \right)e^{2\mu\varphi}\geq\\&\geq 2
\int \left(S_\tau v A_\tau v \right) e^{2\mu\varphi}=\\&=
-2\tau\int2\left(L(v)+\tau^2\mathbf{g}(\nabla\varphi,\nabla\varphi)
v\right)\mathbf{g}(\nabla\varphi,\nabla v)e^{2\mu\varphi}:=\\&:=
-2\tau J_1-2\tau^3 J_2,
\end{aligned}
\end{equation}
where

\begin{equation}\label{oper-2ord-118}
\begin{aligned}
J_1=\int2 L(v)\mathbf{g}(\nabla\varphi,\nabla v)e^{2\mu\varphi}
\end{aligned}
\end{equation}
and
\begin{equation}\label{oper-2ord-119}
\begin{aligned}
J_2=\int2 \mathbf{g}(\nabla\varphi,\nabla\varphi)
\mathbf{g}(\nabla\varphi,\nabla v)ve^{2\mu\varphi}.
\end{aligned}
\end{equation}

\bigskip

\noindent \textbf{Let us examine} $J_1$.

\medskip

Let us apply \eqref{oper-2ord-104} with

$$\beta=e^{2\mu\varphi}\nabla^{(g)} \varphi.$$
We get
\begin{equation}\label{oper-2ord-120}
\begin{aligned}
&J_1= \int\left(\mbox{div } (e^{2\mu\varphi}\nabla^{(g)}
\varphi)\right) \mathbf{g}\left(\nabla v,\nabla
v\right)-\\&-2\partial_l\left(e^{2\mu\varphi}g^{kj}\partial_j\varphi\right)
g^{ls}\partial_sv\partial_kv=\\&= \int \left(
e^{2\mu\varphi}(L(\varphi))\mathbf{g}\left(\nabla v,\nabla
v\right)+2\mu e^{2\mu\varphi}\mathbf{g}\left(\nabla \varphi,\nabla
\varphi\right)\mathbf{g}\left(\nabla v,\nabla
v\right)-\right.\\&\quad\left.-
2e^{2\mu\varphi}g^{kj}\partial^2_{lj}\varphi
g^{ls}\partial_sv\partial_kv-4\mu e^{2\mu\varphi}\partial_l\varphi
g^{kj}\partial_j\varphi g^{ls}\partial_s v\partial_kv\right)=\\&=
\int e^{2\mu\varphi}\left((L(\varphi))\mathbf{g}\left(\nabla
v,\nabla
v\right)-2\partial^2\varphi\nabla^{(g)}v\cdot\nabla^{(g)}v+\right.\\&\quad\left.+2\mu
\mathbf{g}\left(\nabla \varphi,\nabla
\varphi\right)\mathbf{g}\left(\nabla v,\nabla v\right)-4\mu
\left(\mathbf{g}\left(\nabla \varphi,\nabla
v\right)\right)^2\right).
\end{aligned}
\end{equation}

\bigskip

\noindent \textbf{Let us examine} $J_2$.

\medskip

\begin{equation*}
\begin{aligned}
&J_2= \int 2e^{2\mu\varphi} \mathbf{g}(\nabla\varphi,\nabla\varphi)
\mathbf{g}(\nabla\varphi,\nabla v)v=\\&= \int e^{2\mu\varphi}
\mathbf{g}(\nabla\varphi,\nabla\varphi)g^{ls}\partial_l\varphi\partial_s(v^2)=\\&=
-\int
\partial_s\left(e^{2\mu\varphi}\mathbf{g}(\nabla\varphi,\nabla\varphi)g^{ls}\partial_l\varphi\right)v^2=\\&=
-\int
e^{2\mu\varphi}\left(2\partial^2\varphi\nabla^{(g)}\varphi\cdot\nabla^{(g)}\varphi+\mathbf{g}(\nabla\varphi,\nabla\varphi)L(\varphi)+2\mu
\left(\mathbf{g}(\nabla \varphi,\nabla \varphi)\right)^2\right)v^2.
\end{aligned}
\end{equation*}

\medskip

\noindent By the above obtained equality and by \eqref{oper-2ord-120} we have

\begin{equation}\label{oper-2ord-122}
\begin{aligned} 2 \int \left(S_\tau v A_\tau v \right) e^{2\mu\varphi}&= -2\tau J_1-2\tau^3 J_2=\\&= 2\tau\int
e^{2\mu\varphi}q_{\mu}(x,v,\nabla v,\tau),
\end{aligned}
\end{equation}
where

\begin{equation}\label{oper-2ord-123-28}
\begin{aligned}
q_{\mu}(x,v,\nabla v,\tau)&=q_{0}(x,v,\nabla v,\tau)+\mu
q_{1}(x,v,\nabla v,\tau),
\end{aligned}
\end{equation}

\begin{equation}\label{oper-2ord-137}
\begin{aligned}
&q_{0}(x,v,\nabla
v,\tau)=\\&=2\left[\partial^2\varphi(x)\nabla^{(g)}v\cdot\nabla^{(g)}v+\tau^2\partial^2\varphi(x)\nabla^{(g)}\varphi(x)\cdot\nabla^{(g)}\varphi(x)
v^2\right] + \\&
+L(\varphi)\left[\tau^2\mathbf{g}(\nabla\varphi(x),\nabla\varphi(x))v^2-\mathbf{g}(\nabla
v,\nabla v)\right],
\end{aligned}
\end{equation}
and

\begin{equation}\label{oper-2ord-138}
\begin{aligned}
&q_{1}(x,v,\nabla v,\tau)=\\&=-2 \mathbf{g}(\nabla \varphi(x),\nabla
\varphi(x))\mathbf{g}(\nabla v,\nabla v)+4\left(\mathbf{g}(\nabla
\varphi (x),\nabla v)\right)^2+\\&+2\tau^2\left(\mathbf{g}(\nabla
\varphi(x),\nabla \varphi(x))\right)^2v^2.
\end{aligned}
\end{equation}

\medskip

\noindent Trivially, for any $x_0\in \overline{\Omega}$, the right--hand side
of \eqref{oper-2ord-122} can be written 

\begin{equation*}
\begin{aligned}
2\tau\int e^{2\mu\varphi}q_{\mu}(x,v,\nabla v,\tau)=I_{x_0}+
R_{x_0},
\end{aligned}
\end{equation*}
where
\begin{equation}\label{oper-2ord-140}
\begin{aligned} I_{x_0}&= 2\tau e^{2\mu\varphi(x_0)}\int
q_{\mu}(x_0,v,\nabla v,\tau),
\end{aligned}
\end{equation}
and \begin{equation}\label{oper-2ord-141}
 R_{x_0}=2\tau \int \left(
e^{2\mu\varphi(x)}q_{\mu}(x,v,\nabla v,\tau)-
e^{2\mu\varphi(x_0)}q_{\mu}(x_0,v,\nabla
v,\tau)\right).\end{equation} Now, let $\rho\in (0,1]$, to be chosen.
We get easily

\begin{equation}\label{oper-2ord-123}
\begin{aligned}
\left|R_{x_0}\right|\leq C\rho (\mu+1)^2 e^{2\mu
\Phi_0}\int\left(\tau|\nabla v|^2+\tau^3|v|^2\right)dx,
\end{aligned}
\end{equation}
for every $v\in C_0^{\infty}(B_{\rho}(x_0)\cap \Omega)$, where
$$\Phi_0=\max_{x\in \overline{\Omega}}|\varphi|,$$
$C$ does not depend by $\mu$, but depends on the $C^2\left(\overline{\Omega}\right)$ norm of $\varphi$ 
 and by
$$M_0=\max_{1\leq j,k\leq n}\left|g^{jk}\right|.$$
 Now, the Parseval identity gives 
\begin{equation}\label{oper-2ord-124}
\begin{aligned}
 I_{x_0}= \frac{2\tau e^{2\mu\varphi(x_0)}}{(2\pi)^n}\int
Q_{\mu}(\xi,\tau)\left|\widehat{v}(\xi)\right|^2d\xi,
\end{aligned}
\end{equation}
where
\begin{equation}\label{oper-2ord-136}
\begin{aligned}
Q_{\mu}(\xi,\tau)&=Q_0(\xi,\tau) +\mu Q_1(\xi,\tau),
\end{aligned}
\end{equation}

\smallskip

\begin{equation*}
\begin{aligned}
Q_0(\xi,\tau)=&
2\left[\partial^2\varphi(x_0)\xi^{(g)}\cdot\xi^{(g)}+\tau^2\partial^2\varphi(x_0)\nabla^{(g)}\varphi(x_0)\cdot\nabla^{(g)}
\varphi(x_0)\right]+ \\&
+L(\varphi)(x_0)\left[\tau^2\mathbf{g}(\nabla\varphi(x_0),\nabla\varphi(x_0))-\mathbf{g}(\xi,\xi)\right]
\end{aligned}
\end{equation*}
and
\begin{equation*}
\begin{aligned}
Q_1(\xi,\tau)=&-2 \mathbf{g}(\nabla \varphi(x_0),\nabla
\varphi(x_0))\left(\mathbf{g}(\xi ,\xi)-\tau^2\mathbf{g}(\nabla
\varphi(x_0),\nabla
\varphi(x_0))\right)+\\&+4\left(\mathbf{g}(\nabla
\varphi(x_0),\xi)\right)^2.
\end{aligned}
\end{equation*}

\medskip

\textbf{Claim.}

 There exist $C_1$, $C_2$ and $\mu$, positive number, such that

\begin{equation}\label{oper-2ord-125}
\begin{aligned}
\left|\xi+i\tau \nabla\varphi(x_0) \right|^2\leq
C_1Q_{\mu}(\xi,\tau)+C_2\frac{\left|P_2(\xi+i\tau\nabla\varphi(x_0))\right|^2}{\left|\xi+i\tau
\nabla\varphi(x_0) \right|^2},
\end{aligned}
\end{equation}
for every $(\xi,\tau)\in \mathbb{R}^{n+1}$.

\medskip

\textbf{Proof of the  Claim.}

\medskip

\textbf{If $P_2$ is elliptic}, we can proceed in a manner
similar to what we did in the proof of Theorem \ref{Carlm-ell-teor}. For completeness we provide the proof. First, we notice that if $\xi=0$ then \eqref{oper-2ord-112}--\eqref{oper-2ord-113} are
trivially satisfied in $x_0$, because $P_2$ is elliptic and if we had $\xi=0$ we would have in the first condition of  \eqref{oper-2ord-112}--\eqref{oper-2ord-113},
$P_2(i\tau\nabla\varphi(x_0))=0$ from which we would have $\tau=0$
arriving to a contradiction.

  Then for $\xi=0$ the antecedent of
 \eqref{oper-2ord-112}--\eqref{oper-2ord-113} is false. Now, we choose $\mu=0$ and set
 
$$\widetilde{\Sigma}=\left\{(\xi,\tau)\in \mathbb{R}^{n+1}: \quad
\left|\xi+i\tau\nabla\varphi(x_0)\right|=1\right\}.$$ By
\eqref{oper-2ord-112} we have

\begin{equation*}
(\xi,\tau)\in \widetilde{\Sigma},\quad
P_2(\xi+i\tau\nabla\varphi(x_0))=0\quad \Longrightarrow\quad
Q_0(\xi,\tau)>0.
\end{equation*}
Hence, by Lemma \ref{106-CauNire} we have that there exists $B>0$ such that
\begin{equation}\label{oper-2ord-131-712}
\frac{B\left|P_2(\xi+i\tau\nabla\varphi(x_0))\right|^2}{\left|\xi+i\tau\nabla\varphi(x_0)\right|^2}+
Q_0(\xi,\tau)>0,\quad\forall (\xi,\tau)\in \widetilde{\Sigma}
\end{equation}
and, by homogeneity,  \eqref{oper-2ord-125} follows, with $\mu=0$.

\medskip

\textbf{If $P_2$ is not elliptic}, the set

$$K_1=\left\{\xi\in \mathbb{R}^n: \quad \mathbf{g}(\xi,\xi)=0,\quad
|\xi|=1\right\},$$ is nonempty. Let us denote
$$K_2=\left\{\xi\in K_1: \quad \mathbf{g}(\xi,\nabla\varphi(x_0))=
0\right\}.$$ If $K_2=\emptyset$, then we have trivially 
$$\left|\mathbf{g}(\xi,\nabla\varphi(x_0))\right|>0, \quad\forall
\xi\in K_1$$ hence, by the compactness of $K_1$, there exists $m_1>0$
such that
$$\left|\mathbf{g}(\xi,\nabla\varphi(x_0))\right|\geq m_1, \quad\forall
\xi\in K_1.$$ Now we set

$$C_0=1+\max_{|\xi|=1}\left\vert\partial^2\varphi(x_0)\xi^{(g)}\cdot\xi^{(g)}\right\vert$$
abd we have, for any $\xi\in K_1$ and $$\mu=\mu_0:=\frac{C_0}{m_1^2},$$

\begin{equation}\label{oper-2ord-200}
\begin{aligned}
Q_{\mu}(\xi,0)&=2\partial^2\varphi(x_0)\xi^{(g)}\cdot\xi^{(g)}+4\mu\left(\mathbf{g}(\nabla
\varphi(x_0),\xi)\right)^2\geq\\&\geq -2C_0+4\mu m^2_1\geq 2C_0.
\end{aligned}
\end{equation}
Let us fix $\mu=\mu_0$ and let $(\xi,\tau)\in \widetilde{\Sigma}$ satisfy
$$P_2(\xi+i\tau\nabla\varphi(x_0))=0;$$ we have what follows:

\smallskip

\noindent (a) if $\tau=0$ then we have, by \eqref{oper-2ord-200},
$Q_{\mu}(\xi,0)>0$;

\noindent (b) if $\tau\neq0$, then we have, by \eqref{oper-2ord-112},
$Q_{\mu}(\xi,\tau)>0$.

\smallskip

Hence

\begin{equation*}
(\xi,\tau)\in \widetilde{\Sigma},\quad P_2(\xi+i\tau\nabla\varphi(x_0))=0\quad
\Longrightarrow\quad Q_{\mu}(\xi,\tau)>0
\end{equation*}
and by Lemma \ref{106-CauNire} we derive that there exists  $B>0$ such that

\begin{equation}\label{oper-2ord-131}
\frac{B\left|P_2(\xi+i\tau\nabla\varphi(x_0))\right|^2}{\left|\xi+i\tau\nabla\varphi(x_0)\right|^2}+
Q_{\mu}(\xi,\tau)>0,\quad\forall (\xi,\tau)\in \widetilde{\Sigma}.
\end{equation}
Hence, if $K_2=\emptyset$, then \eqref{oper-2ord-125} is
satisfied for $\mu=\mu_0$.

\medskip

\noindent Now, let us suppose that $K_2\neq\emptyset$. By
\eqref{oper-2ord-106-bis} we have

\begin{equation}\label{oper-2ord-126}
\xi\in K_2\quad \Longrightarrow\quad Q_0(\xi,0)>0.
\end{equation}
Since $Q_0$ is continuous and $K_2$ is compact, 
there exists  $\delta_1>0$ such that
\begin{equation}\label{oper-2ord-127}
\xi\in K_2\quad \Longrightarrow\quad Q_0(\xi,0)\geq \delta_1.
\end{equation} By compactness of $K_2$, by continuity of
$Q_0$ and by \eqref{oper-2ord-127} it follows that there exists $d_0>0$ such that

\begin{equation}\label{oper-2ord-128}
\xi\in K_1,\mbox{ }
\left|\mathbf{g}(\xi,\nabla\varphi(x_0))\right|\leq d_0
\Longrightarrow\quad Q_0(\xi,0)\geq \frac{\delta_1}{2}.
\end{equation}
 Let us denote, like before,

$$C_0=1+\max_{|\xi|=1}\left\vert\partial^2\varphi(x_0)\xi^{(g)}\cdot\xi^{(g)}\right\vert.$$ We notice that if
$$\xi\in K_1\quad\mbox{ e }\quad \left|\mathbf{g}(\xi,\nabla\varphi(x_0))\right|\geq
d_0,$$ then, for any $\mu\geq \frac{C_0}{d^2_0}$ we have

\begin{equation}\label{oper-2ord-129}
\begin{aligned}
Q_{\mu}(\xi,0)&=2\partial^2\varphi(x_0)\xi^{(g)}\cdot\xi^{(g)}+4\mu\left(\mathbf{g}(\nabla
\varphi(x_0),\xi)\right)^2\geq\\&\geq -2C_0+4\mu d^2_0\geq 2C_0>0.
\end{aligned}
\end{equation}
Now, let us fix
$$\mu=\mu_1:=\frac{C_0}{d^2_0}$$ and by \eqref{oper-2ord-128} and
\eqref{oper-2ord-129} we obtain

\begin{equation}\label{oper-2ord-130}
\xi\in K_1\quad \Longrightarrow\quad Q_{\mu}(\xi,0)\geq \delta_2>0,
\end{equation} where
$$\delta_2=\min\left\{2C_0,
\frac{\delta_1}{2}\right\}.$$  From now on, we proceed
as we already did to prove the \eqref{oper-2ord-131}.
\textbf{The proof of the Claim is concluded}.

\medskip
From now on, we \textbf{fix a value} $\mu$ for which
\eqref{oper-2ord-125} is satisfied (recalling that, however,
$\mu$ depends on $x_0$). Set $N=\nabla\varphi(x_0)$. Using the
notations introduced in \eqref{oper-2ord-201}, we have by
\eqref{oper-2ord-125}

\begin{equation}\label{oper-2ord-132}
\begin{aligned}
\left\Vert v\right\Vert^2_{1,\tau}&\leq \frac{C_1}{(2\pi)^n}\int
Q_{\mu}(\xi,\tau)\left|\widehat{v}(\xi)\right|^2d\xi+ C_2\left\Vert
p_2(x_0,\partial,\tau)v\right\Vert^2_{-1,\tau},
\end{aligned}
\end{equation}
where $p_2(x,\partial,\tau)$ is defined in \eqref{oper-2ord-303}.
We have easily

\begin{equation}\label{oper-2ord-133}
\begin{aligned}
\left\vert
p_2(x,\partial,\tau)v-p_2(x_0,\partial,\tau)v\right\vert&\leq
C\rho\left(\tau |\nabla v|+\tau^2|v|\right),
\end{aligned}
\end{equation}
for every $v\in C_0^{\infty}(B_{\rho}(x_0)\cap \Omega)$. Now, by
\eqref{oper-2ord-203}, \eqref{oper-2ord-133} and by the triangle inequality, we get

\begin{equation*}
\begin{aligned}
&\left\Vert p_2(x_0,\partial,\tau)v\right\Vert_{-1,\tau}\leq
\frac{1}{\tau|\nabla\varphi(x_0)|}\left\Vert
p_2(x_0,\partial,\tau)v\right\Vert_{L^2(\Omega)}\leq\\&\leq
\frac{1}{\tau|\nabla\varphi(x_0)|}\left(\left\Vert
p_2(\cdot,\partial,\tau)v-
p_2(x_0,\partial,\tau)v\right\Vert_{L^2(\Omega)}+\left\Vert
p_2(\cdot,\partial,\tau)v\right\Vert_{L^2(\Omega)}\right)\leq
\\&\leq
\frac{1}{\tau|\nabla\varphi(x_0)|}\left[C\rho\tau\left\Vert |\nabla
v|+\tau|v|\right\Vert_{L^2(\Omega)}+\left\Vert
p_2(\cdot,\partial,\tau)v\right\Vert_{L^2(\Omega)}\right]\leq\\&\leq
C\rho\left\Vert
v\right\Vert_{1,\tau}+\frac{1}{\tau|\nabla\varphi(x_0)|}\left\Vert
p_2(\cdot,\partial,\tau)v\right\Vert_{L^2(\Omega)}.
\end{aligned}
\end{equation*}
In sum, we have that
\begin{equation}\label{oper-2ord-134}
\begin{aligned}
\left\Vert p_2(x_0,\partial,\tau)v\right\Vert_{-1,\tau}\leq
C\left(\rho\left\Vert v\right\Vert_{1,\tau}+\frac{1}{\tau}\left\Vert
p_2(\cdot,\partial,\tau)v\right\Vert_{L^2(\Omega)}\right),
\end{aligned}
\end{equation}
for every $v\in C_0^{\infty}(B_{\rho}(x_0)\cap \Omega)$, where $C$
depends on $C^1\left(\overline{\Omega}\right)$ norm of $\varphi$,
on $M_0$ and on $$m_2=\min_{\overline{\Omega}}|\nabla\varphi|.$$ By
\eqref{oper-2ord-134} and by \eqref{oper-2ord-132} we obtain

\begin{equation}\label{oper-2ord-135}
\begin{aligned}
2\tau\left\Vert v\right\Vert^2_{1,\tau}\leq &\frac{ 2\tau
C_1}{(2\pi)^n}\int
Q_{\mu}(\xi,\tau)\left|\widehat{v}(\xi)\right|^2d\xi+\\&+ 2\tau
C\rho^2\left\Vert v\right\Vert^2_{1,\tau}+C\tau^{-1}\left\Vert
p_2(\cdot,\partial,\tau)v\right\Vert^2_{L^2(\Omega)},
\end{aligned}
\end{equation}
for every $v\in C_0^{\infty}(B_{\rho}(x_0)\cap \Omega)$.

Now we need to estimate the first term on the right--hand side in
\eqref{oper-2ord-135}. From \eqref{oper-2ord-140},
\eqref{oper-2ord-141} and \eqref{oper-2ord-136} we have

\begin{equation}\label{oper-2ord-142}
\begin{aligned}
&\frac{2\tau}{(2\pi)^n}\int
Q_{\mu}(\xi,\tau)\left|\widehat{v}(\xi)\right|^2d\xi=2\tau\int
q_{\mu}(x_0,v,\nabla v,\tau)dx=\\&= 2\tau e^{-2\mu\varphi(x_0)}\int
e^{2\mu\varphi(x_0)}q_{\mu}(x_0,v,\nabla v,\tau)dx=\\&=2\tau
e^{-2\mu\varphi(x_0)}\int
\left(e^{2\mu\varphi(x_0)}q_{\mu}(x_0,v,\nabla
v,\tau)-\right. \\& \left. -e^{2\mu\varphi(x)}q_{\mu}(x,v,\nabla
v,\tau)\right)dx+\\&+2\tau e^{-2\mu\varphi(x_0)}\int
e^{2\mu\varphi(x)}q_{\mu}(x,v,\nabla v,\tau)dx=\\&= R_{x_0}+2\tau
e^{-2\mu\varphi(x_0)}\int e^{2\mu\varphi(x)}q_{\mu}(x,v,\nabla
v,\tau)dx.
\end{aligned}
\end{equation}

\smallskip

From what we have just obtained, from \eqref{oper-2ord-122} and from
\eqref{oper-2ord-123} we have (recall that we have fixed $\mu$)

\begin{equation}\label{oper-2ord-304}
\begin{aligned}
&\frac{2\tau}{(2\pi)^n}\int
Q_{\mu}(\xi,\tau)\left|\widehat{v}(\xi)\right|^2d\xi\leq \\& \leq C\rho
(\mu+1)^2 e^{2\mu \Phi_0}\int\left(\tau|\nabla
v|^2+\tau^3|v|^2\right)dx+\\&+C \int \left(S_\tau v A_\tau v\right)
e^{2\mu\varphi}dx\leq\\&\leq C\rho\tau \left\Vert
v\right\Vert^2_{1,\tau}+C\int \left|p_2(x,\partial,\tau)
v\right|^2e^{2\mu\varphi}dx,
\end{aligned}
\end{equation}

\smallskip

\noindent for every $v\in C_0^{\infty}(B_{\rho}(x_0)\cap \Omega)$. In the last estimate from above we used that (see \eqref{oper-2ord-117})
$$2 \int \left(S_\tau v A_\tau v\right)
e^{2\mu\varphi}dx\leq \int \left|p_2(x,\partial,\tau)
v\right|^2e^{2\mu\varphi}dx.$$ By \eqref{oper-2ord-135} and
\eqref{oper-2ord-304} (recalling that $\tau\geq 1$) we get

\begin{equation*}
\begin{aligned}
2\tau\left\Vert v\right\Vert^2_{1,\tau}\leq C\rho\tau \left\Vert
v\right\Vert^2_{1,\tau}+C\left(1+\tau^{-1}\right)\int
\left|p_2(x,\partial,\tau) v\right|^2dx,
\end{aligned}
\end{equation*}
for every $v\in C_0^{\infty}(B_{\rho}(x_0)\cap \Omega)$. Now, let us choose

$$\rho=\rho_0:=\min\left\{\frac{1}{C},1\right\}$$ and we get

\begin{equation*}
\begin{aligned}
\tau\left\Vert v\right\Vert^2_{1,\tau}\leq C\int
\left|p_2(x,\partial,\tau) v\right|^2dx,
\end{aligned}
\end{equation*}
for every $v\in C_0^{\infty}(B_{\rho_0}(x_0)\cap \Omega)$ and for every
$\tau\geq 1$.

Now, recalling \eqref{oper-2ord-116} and applying  Lemmas
\ref{lemma-carlm-ell-I} and \ref{lemma-carlm-ell-II} we obtain
\begin{equation}\label{oper-2ord-305}
\begin{aligned}
& \tau^3 \int_{\Omega} |u|^2 e^{2 \tau \varphi}dx + \tau
\int_{\Omega}
 |\nabla u|^2 e^{2 \tau \varphi}dx \leq C
\int_{\Omega} |P_2 (\partial)u|^2 e^{2 \tau \varphi}dx,
\end{aligned}
\end{equation}
for every $u\in C_0^{\infty}(B_{\rho_0}(x_0)\cap \Omega)$ and for every
$\tau\geq \tau^*$, where $C$ and $\tau^*$ are suitable positive numbers depending by $x_0$. At this point we have only to
apply  Lemma \ref{stime-Carlm-1.3.3-39} to conclude the proof. $\blacksquare$

\section[Second order operators II]{Second order operators II -- Lipschitz coefficients in the principal part}\label{25-4-23:second-order}

In the following Theorem we will consider the operator
\begin{equation}\label{oper-2ord-332-1}
P_2(x,\partial)=g^{jk}(x)\partial_{jk}^2,\quad \mbox{in } B_1,
\end{equation}
where $\left\{g^{jk}(x)\right\}_{j,k=1}^n$ is a real symmetric matrix--valued function. Recall that $\phi\in
C^2\left(\overline{\Omega}\right)$ satisfies
\begin{equation*}
\nabla \phi(x)\neq 0,\quad \forall x\in \overline{B_1}.
\end{equation*}
Moreover, $\phi$ satisfies condition $(\mathbf{S})$ in $0$ w.r.t. the operator
$P_2(x,\partial)$ if we have

\begin{equation}\label{oper-2ord-332-2}
\begin{aligned}
&\begin{cases}
P_2(0,\xi)=0, \\
\\
P_2^{(j)}(\xi)\partial_j\phi(0)= 0,\\
\\
 \xi\neq 0,
\end{cases}\mbox{ }\Longrightarrow\\&
\quad\quad\quad\quad\Longrightarrow \mbox{ }
\partial^2_{jk}\phi(0)P^{(j)}_2(0,\xi)P^{(k)}_2(0,\xi)>0
\end{aligned}
\end{equation}
and
\begin{equation}\label{oper-2ord-332-3}
\begin{aligned}
&\begin{cases}
P_2(0,\xi+i\tau\nabla \phi (0))=0, \\
\\
\tau\neq 0,
\end{cases}\mbox{ }\Longrightarrow\\&
\Longrightarrow \mbox{ }
\partial^2_{jk}\phi(0)P^{(j)}_2(0,\xi+i\tau\nabla \phi (0))\overline{P^{(k)}_2(0,\xi+i\tau\nabla \phi (0))}>0.
\end{aligned}
\end{equation}

\medskip

\noindent Concerning \eqref{oper-2ord-332-2} and 
\eqref{oper-2ord-332-3}, keep in mind, respectively, the
Remarks that follow definitions
\ref{psedoconvess} and \ref{condizione-S}.

\bigskip

\begin{theo}\label{oper-2ord-306} Let
$P_2(x,\partial)$ be operator \eqref{oper-2ord-332-1}. Let us assume that \eqref{oper-2ord-101a} and \eqref{oper-2ord-101b} are satisfied with 
$\Omega=B_1$. Let us suppose that $\varphi \in
C^{2}\left(\overline{B_1}\right)$ satisfies  condition
$(\mathbf{S})$ w.r.t. $P_2(x,\partial)$ in $x=0$.

Then there exist $\rho_0\in (0,1]$, $\delta_0\in(0,1]$, $C\geq 1$ and
$\tau_0\geq 1$ such that

\begin{equation}\label{oper-2ord-307}
\begin{aligned}
& \tau^3 \int_{B_1} |u|^2 e^{2 \tau \varphi}dx + \tau
\int_{B_1}
 |\nabla u|^2 e^{2 \tau \varphi}dx \leq \\&\leq C
\int_{B_1} |P_2 (\delta x,\partial)u|^2 e^{2 \tau \varphi}dx,
\end{aligned}
\end{equation}
for every $\delta\in(0,\delta_0]$, for every $u \in
C^{\infty}_0(B_{\rho_0})$ and for every $\tau\geq \tau_0$.
\end{theo}
\textbf{Proof.} For the most part of the proof
we repeat the steps we did in the proof of
Theorem \ref{Carl-2ord-costanti} by paying special attention to the
additional terms that come up because now the coefficients
of the operator are variables.

\smallskip

Let us denote 

$$g_{\delta}^{jk}(x)=g^{jk}(\delta x),\quad j,k=1,\cdots,n.$$
Let $u\in C^{\infty}_0(B_1)$ and set $v= e^{\tau \varphi}u$.
We have
\begin{equation*}
\begin{aligned}
&P_{2}(\delta x,\partial,\tau) v:=e^{\tau \varphi} P_2(\delta
x,\partial)(e^{-\tau \varphi} v)= \\&= P_2(\delta x,\partial)v -
\tau (P_2(\delta x,\partial)\varphi) v
+\tau^2\mathbf{g}_{\delta}(\nabla\varphi,\nabla\varphi) v - 2 \tau
\mathbf{g}_{\delta}(\nabla\varphi,\nabla v)=\\&=
\partial_j\left(g_{\delta}^{jk}(x)\partial_k\right)- 2 \tau
\mathbf{g}_{\delta}(\nabla\varphi,\nabla v)-\\& - \tau (P_2(\delta
x,\partial)\varphi) v-\partial_j\left(g^{jk}_{\delta}(
x)\right)\partial_kv.
\end{aligned}
\end{equation*}
Let us denote
\begin{equation}\label{oper-2ord-308}
L_{\delta}(v)=\partial_j\left(g_{\delta}^{jk}(x)\partial_k v\right),
\end{equation}

\begin{equation}\label{oper-2ord-309}
p_2(x,\partial,\tau) v=S_\tau v+ A_{\tau}v,
\end{equation}
where
$$S_\tau v=L(v)  +
\tau^2\mathbf{g}_{\delta}(\nabla\varphi,\nabla\varphi) v$$ and
$$A_{\tau}v=- 2 \tau \mathbf{g}_{\delta}(\nabla\varphi,\nabla v).$$
Let us notice that
\begin{equation}\label{oper-2ord-310}
P_2(\delta x,\partial,\tau) v-p_2(x,\partial,\tau) v=- \tau
(P_2(\delta x,\partial)\varphi) v-\partial_j\left(g_{\delta}^{jk}(
x)\right)\partial_kv.
\end{equation}

\smallskip

For the sake of brevity, we omit the domain and the element of
integration.

Let $\mu$ be a positive number that we will choose later.
We have trivially

\begin{equation}\label{oper-2ord-311}
\begin{aligned}
\int \left|p_2(x,\partial,\tau) v\right|^2e^{2\mu\varphi}&= \int
|S_\tau v|^2 e^{2\mu\varphi}+ \int |A_\tau v|^2e^{2\mu\varphi} +\\&+
2 \int \left(S_\tau v A_\tau v\right) e^{2\mu\varphi}\geq\\&\geq 2
\int \left(S_\tau v A_\tau v \right)e^{2\mu\varphi}.
\end{aligned}
\end{equation}

Set

\begin{equation}\label{oper-2ord-312}
\begin{aligned}
J_1=\int2\left(
L_{\delta}v\right)\mathbf{g}_{\delta}(\nabla\varphi,\nabla
v)e^{2\mu\varphi}
\end{aligned}
\end{equation}
and
\begin{equation}\label{oper-2ord-313}
\begin{aligned}
J_2=\int2 \mathbf{g}_{\delta}(\nabla\varphi,\nabla\varphi)
\mathbf{g}_{\delta}(\nabla\varphi,\nabla v)ve^{2\mu\varphi},
\end{aligned}
\end{equation}
we have

\begin{equation}\label{oper-2ord-314}
\begin{aligned}
2 \int \left(S_\tau v A_\tau v\right) e^{2\mu\varphi}=-2\tau
J_1-2\tau^3 J_2.
\end{aligned}
\end{equation}
Now, to handle $J_1$, we apply \eqref{oper-2ord-104} with

$$\beta_{\delta}=e^{2\mu\varphi}\nabla^{(g_{\delta})} \varphi.$$
It should be kept in mind that now, with respect to proof of the
Theorem \ref{Carl-2ord-costanti}, the coefficients of the operator $P_2(\delta
x,\partial)$ depend on $x$, more precisely they are of the type
$f(\delta x)$ where $f$ is a Lipschitz continuous  function. Hence we have

\begin{equation*}
\begin{aligned}
&J_1\leq\int
e^{2\mu\varphi}\left((L_{\delta}(\varphi))\mathbf{g}_{\delta}\left(\nabla
v,\nabla
v\right)-2\partial^2\varphi\nabla^{(g_{\delta})}v\cdot\nabla^{(g_{\delta})}v+\right.\\&\quad\left.+2\mu
\mathbf{g}_{\delta}\left(\nabla \varphi,\nabla
\varphi\right)\mathbf{g}_{\delta}\left(\nabla v,\nabla v\right)-4\mu
\left(\mathbf{g}_{\delta}\left(\nabla \varphi,\nabla
v\right)\right)^2\right)+ C\delta\int e^{2\mu\varphi}|\nabla v|^2
\end{aligned}
\end{equation*}
and, similarly,
\begin{equation*}
\begin{aligned}
J_2\leq&-\int
e^{2\mu\varphi}\left(2\partial^2\varphi\nabla^{(g_{\delta})}\varphi\cdot\nabla^{(g_{\delta})}\varphi+\mathbf{g}_{\delta}(\nabla\varphi,\nabla\varphi)L_{\delta}(\varphi)+
\right.\\&\left.+2\mu \left(\mathbf{g}_{\delta}(\nabla
\varphi,\nabla \varphi)\right)^2\right)v^2+C\delta\int e^{2\mu\varphi} v^2,
\end{aligned}
\end{equation*}
in the last two inequalities, $C$ depends on $M_0$, $M_1$ and
$\left\Vert\varphi\right\Vert_{C^2\left(\overline{B_1}\right)}$,
but does not depend on $\mu$. By these inequalities and by
\eqref{oper-2ord-314} we get

\begin{equation}\label{oper-2ord-315}
\begin{aligned} 2\int  \left(S_\tau v
A_\tau v \right) e^{2\mu\varphi}&= -2\tau J_1-2\tau^3 J_2\geq\\&\geq
2\tau\int e^{2\mu\varphi}q^{(\delta)}_{\mu}(x,v,\nabla
v,\tau)-\\&-C\delta\tau\int e^{2\mu\varphi}\left(\tau^2v^2+|\nabla v|^2\right),
\end{aligned}
\end{equation}
where

\begin{equation}\label{oper-2ord-316}
\begin{aligned}
q^{(\delta)}_{\mu}(x,v,\nabla v,\tau)&=q^{(\delta)}_{0}(x,v,\nabla
v,\tau)+\mu q^{(\delta)}_{1}(x,v,\nabla v,\tau),
\end{aligned}
\end{equation}

\begin{equation}\label{oper-2ord-317}
\begin{aligned}
&q^{(\delta)}_{0}(x,v,\nabla
v,\tau)=\\&=2\left[\partial^2\varphi(x)\nabla^{(g_{\delta})}v\cdot\nabla^{(g_{\delta})}v+\right. \\& \left.\tau^2(\partial^2\varphi(x)\nabla^{(g_{\delta})}\varphi(x)\cdot\nabla^{(g_{\delta})}
\varphi(x)) v^2\right] + \\&
+L_{\delta}(\varphi)\left[\tau^2\mathbf{g}_{\delta}(\nabla\varphi(x),\nabla\varphi(x))v^2-\mathbf{g}_{\delta}(\nabla
v,\nabla v)\right]
\end{aligned}
\end{equation}
and
\begin{equation}\label{oper-2ord-318}
\begin{aligned}
&q^{(\delta)}_{1}(x,v,\nabla v,\tau)=\\&=-2 \mathbf{g}_{\delta}(\nabla
\varphi(x),\nabla \varphi(x))\mathbf{g}_{\delta}(\nabla v,\nabla
v)+4\left(\mathbf{g}_{\delta}(\nabla \varphi (x),\nabla
v)\right)^2+\\&+2\tau^2\left(\mathbf{g}_{\delta}(\nabla
\varphi(x),\nabla \varphi(x))\right)^2v^2.
\end{aligned}
\end{equation}

\medskip

Now we \textbf{examine  the first addend on the right--hand side of
\eqref{oper-2ord-315}}, namely
$$2\tau\int e^{2\mu\varphi}q^{(\delta)}_{\mu}(x,v,\nabla
v,\tau).$$ Likewise to the proof of
Theorem \ref{Carl-2ord-costanti} we write

\begin{equation*}
\begin{aligned}
2\tau\int e^{2\mu\varphi}q^{(\delta)}_{\mu}(x,v,\nabla
v,\tau)=I_{0}+ R_{0},
\end{aligned}
\end{equation*}
where
\begin{equation}\label{oper-2ord-319}
\begin{aligned} I_{0}&= 2\tau e^{2\mu\varphi(0)}\int
q^{(\delta)}_{\mu}(0,v,\nabla v,\tau)
\end{aligned}
\end{equation}
and \begin{equation}\label{oper-2ord-320}
 R_{0}=2\tau \int \left(
e^{2\mu\varphi(x)}q^{(\delta)}_{\mu}(x,v,\nabla v,\tau)-
e^{2\mu\varphi(0)}q^{(\delta)}_{\mu}(0,v,\nabla
v,\tau)\right).\end{equation} Let $\rho\in (0,1]$, to be chosen;
we get easily (recall $\delta\leq 1$)

\begin{equation}\label{oper-2ord-321}
\begin{aligned}
\left|R_{0}\right|\leq C\rho (\mu+1)^2 e^{2\mu
\Phi_0}\int\left(\tau|\nabla v|^2+\tau^3|v|^2\right)dx,
\end{aligned}
\end{equation}
for every $v\in C_0^{\infty}(B_{\rho})$, where $$\Phi_0=\max_{x\in
\overline{B_1}}|\varphi|,$$ $C$ does not depend on $\mu$, but depends
on the $C^2\left(\overline{B_1}\right)$ norm  of $\varphi$ and on
$M_0$.

Let us denote  $$g_0^{jk}=g^{jk}(0).$$By the Parseval identity we have
\begin{equation}\label{oper-2ord-322}
\begin{aligned}
 I_{0}= \frac{2\tau e^{2\mu\varphi(0)}}{(2\pi)^n}\int
Q_{\mu}(\xi,\tau)\left|\widehat{v}(\xi)\right|^2d\xi,
\end{aligned}
\end{equation}
where
\begin{equation}\label{oper-2ord-323-28}
\begin{aligned}
Q_{\mu}(\xi,\tau)&=Q_0(\xi,\tau) +\mu Q_1(\xi,\tau),
\end{aligned}
\end{equation}

\smallskip

\begin{equation*}
\begin{aligned}
Q_0(\xi,\tau)=&
2\left[\partial^2\varphi(0)\xi^{(g_0)}\cdot\xi^{(g_0)}+\tau^2\partial^2\varphi(0)\nabla^{(g_0)}\varphi(0)\cdot\nabla^{(g_0)}
\varphi(0)\right]+ \\&
+L_{0}(\varphi)(0)\left[\tau^2\mathbf{g}_0(\nabla\varphi(0),\nabla\varphi(0))-\mathbf{g}_0(\xi,\xi)\right]
\end{aligned}
\end{equation*}
and
\begin{equation*}
\begin{aligned}
Q_1(\xi,\tau)=&-2 \mathbf{g}_0(\nabla \varphi(0),\nabla
\varphi(0))\left(\mathbf{g}_0(\xi ,\xi)-\tau^2\mathbf{g}_0(\nabla
\varphi(0),\nabla \varphi(0))\right)+\\&+4\left(\mathbf{g}_0(\nabla
\varphi(0),\xi)\right)^2.
\end{aligned}
\end{equation*}

\medskip

\noindent Similarly to what was done in the proof of Theorem \ref{Carl-2ord-costanti}  can be proved
the existence of constants $C_1$, $C_2$ and $\mu$ such that

\begin{equation}\label{oper-2ord-323}
\begin{aligned}
\left|\xi+i\tau \nabla\varphi(0) \right|^2\leq
C_1Q_{\mu}(\xi,\tau)+C_2\frac{\left|P_2(0,\xi+i\tau\nabla\varphi(0))\right|^2}{\left|\xi+i\tau
\nabla\varphi(0) \right|^2},
\end{aligned}
\end{equation}
for every $(\xi,\tau)\in \mathbb{R}^{n+1}$.

\textbf{From now on we fix a value $\mu$ for which
\eqref{oper-2ord-323} is satisfied}. Adopting the notations
introduced in \eqref{oper-2ord-201} with $N=\nabla\varphi(0)$, by
\eqref{oper-2ord-323} we have

\begin{equation}\label{oper-2ord-324}
\begin{aligned}
&2\tau\left\Vert v\right\Vert^2_{1,\tau}\leq\\&\leq  \frac{2\tau
C_1}{(2\pi)^n}\int
Q_{\mu}(\xi,\tau)\left|\widehat{v}(\xi)\right|^2d\xi+ 2\tau
C_2\left\Vert p_2(0,\partial,\tau)v\right\Vert^2_{-1,\tau}.
\end{aligned}
\end{equation}

\medskip

Now, \textbf{we estimate from above the last term on the right--hand side in
\eqref{oper-2ord-324}}. For this purpose, we will repeatedly use
the triangle inequality and Lemma  \ref{Lemma-lipsch}.

\smallskip

First of all, set

\begin{equation*}
\begin{aligned}
\widetilde{p}_2(x,\partial,\tau)v=g_{\delta}^{jk}(x)\partial^2_{jk}v-2\tau
g_{\delta}^{jk}(x)\partial_j\varphi(0)\partial_kv+\tau^2g_{\delta}^{jk}(x)\partial_j\varphi(0)\partial_k\varphi(0)v
\end{aligned}
\end{equation*}
and let us note that

\begin{equation*}
\begin{aligned}
\widetilde{p}_2(0,\partial,\tau)=p_2(0,\partial,\tau).
\end{aligned}
\end{equation*}
By this equality and by the triangle inequality we have

\begin{equation}\label{oper-2ord-325}
\begin{aligned}
\left\Vert p_2(0,\partial,\tau)v\right\Vert^2_{-1,\tau}=&\left\Vert
\widetilde{p}_2(0,\partial,\tau)v\right\Vert^2_{-1,\tau}\leq
2\left\Vert
\widetilde{p}_2(\cdot,\partial,\tau)v\right\Vert^2_{-1,\tau}+\\&+
2\left\Vert
\widetilde{p}_2(\cdot,\partial,\tau)v-\widetilde{p}_2(0,\partial,\tau)v\right\Vert^2_{-1,\tau}.
\end{aligned}
\end{equation}
Let us estimate from above the first term on the right--hand side in \eqref{oper-2ord-325}.
By the triangle inequality and by the definition of
$\left\Vert\cdot \right\Vert_{-1,\tau}$ we obtain

\begin{equation}\label{oper-2ord-328}
\begin{aligned}
&\left\Vert
\widetilde{p}_2(\cdot,\partial,\tau)v\right\Vert^2_{-1,\tau}\leq
C\tau^{-2}\left\Vert
\widetilde{p}_2(\cdot,\partial,\tau)v\right\Vert^2_{L^2(B_1)}\leq\\&\leq
C\tau^{-2}\left\Vert
\widetilde{p}_2(\cdot,\partial,\tau)v-p_2(\cdot,\partial,\tau)v\right\Vert^2_{L^2(B_1)}+\\&+C\tau^{-2}\left\Vert
p_2(\cdot,\partial,\tau)v\right\Vert^2_{L^2(B_1)},
\end{aligned}
\end{equation}
on the other hand,

\begin{equation*}
\begin{aligned}
&\widetilde{p}_2(x,\partial,\tau)v-p_2(x,\partial,\tau)v=\\&=-\partial_j\left(g_{\delta}^{jk}(x)\right)\partial_kv+2\tau
g_{\delta}^{jk}(x)\left(\partial_j\varphi(x)-\partial_j\varphi(0)\right)\partial_kv+\\&+
\tau^2 g_{\delta}^{jk}(x)\left(\partial_j \varphi(x)\partial_k
\varphi(x)-\partial_j \varphi(0)\partial_k \varphi(0)\right)v.
\end{aligned}
\end{equation*}
From the latter and the triangle inequality, we deduce easily that

\begin{equation}\label{oper-2ord-329}
\begin{aligned}
&\tau^{-2}\left\Vert
\widetilde{p}_2(\cdot,\partial,\tau)v-p_2(\cdot,\partial,\tau)v\right\Vert^2_{L^2(B_1)}\leq\\&\leq 
C\delta^2\tau^{-2}\left\Vert \nabla
v\right\Vert^2_{L^2(B_1)}+C\rho^2\left\Vert
v\right\Vert^2_{1,\tau}
\end{aligned}
\end{equation}
for every $v\in C_0^{\infty}(B_{\rho})$.

\medskip

 By \eqref{oper-2ord-328} and \eqref{oper-2ord-329} we get
\begin{equation}\label{oper-2ord-330}
\begin{aligned}
\left\Vert
\widetilde{p}_2(\cdot,\partial,\tau)v\right\Vert^2_{-1,\tau}&\leq
C\tau^{-2}\left\Vert \nabla
v\right\Vert^2_{L^2(B_1)}+C\rho^2\left\Vert
v\right\Vert^2_{1,\tau}+\\&+C\tau^{-2}\left\Vert
p_2(\cdot,\partial,\tau)v\right\Vert^2_{L^2(B_1)}.
\end{aligned}
\end{equation}

\medskip

Now, let us estimate from above the second term on the right--hand side in \eqref{oper-2ord-325}.

\noindent We notice
\begin{equation*}
\begin{aligned}
\widetilde{p}_2(x,\partial,\tau)v&-\widetilde{p}_2(0,\partial,\tau)v=\\&=\left(g_{\delta}^{jk}(x)-g_{\delta}^{jk}(0)\right)
\left(\partial_j-\tau\partial_j\varphi(0)\right)\left(\partial_k-\tau\partial_k\varphi(0)\right)v.
\end{aligned}
\end{equation*}
and, for a fixed  $k=1,\cdots, n$, we set
$$w_k=\left(\partial_k-\tau\partial_k\varphi(0)\right)v.$$ By applying Lemma \ref{Lemma-lipsch} we have, for every $v\in
C_0^{\infty}(B_{\rho})$,

\begin{equation}\label{oper-2ord-326}
\begin{aligned}
&\left\Vert \left(g_{\delta}^{jk}(x)-g_{\delta}^{jk}(0)\right)
\left(\partial_j-\tau\partial_j\varphi(0)\right)w_k\right\Vert^2_{-1,\tau}\leq\\&
\quad\quad\quad\quad\quad\quad\quad\quad \leq
C\left(M_1^2\delta^2\rho^2+|\tau
\nabla\varphi(0)|^{-2}\right)\left\Vert
w_k\right\Vert^2_{L^2(B_1)}.
\end{aligned}
\end{equation}
On the other hand, for $k=1,\cdots, n$,
\begin{equation*}
\begin{aligned}
\left\Vert w_k\right\Vert^2_{L^2(B_1)}&= \left\Vert
\left(\partial_k-\tau\partial_k\varphi(0)\right)v\right\Vert^2_{L^2(B_1)}\leq
C\left\Vert v\right\Vert^2_{1,\tau}.
\end{aligned}
\end{equation*}
By the latter and by \eqref{oper-2ord-326} we get
\begin{equation}\label{oper-2ord-327}
\begin{aligned}
\left\Vert
\widetilde{p}_2(\cdot,\partial,\tau)v-\widetilde{p}_2(0,\partial,\tau)v\right\Vert^2_{-1,\tau}\leq
C\left(M_1^2\rho^2+\tau^{-2}\right)\left\Vert
v\right\Vert^2_{1,\tau}.
\end{aligned}
\end{equation}
By \eqref{oper-2ord-325}, \eqref{oper-2ord-330} and
\eqref{oper-2ord-327} we have

\begin{equation*}
\begin{aligned}
\left\Vert p_2(0,\partial,\tau)v\right\Vert^2_{-1,\tau}\leq
C\left(\rho^2+\tau^{-2}\right)\left\Vert
v\right\Vert^2_{1,\tau}+C\tau^{-2}\left\Vert
p_2(\cdot,\partial,\tau)v\right\Vert^2_{L^2(B_1)},
\end{aligned}
\end{equation*}
for every $v\in C_0^{\infty}(B_{\rho})$, where $C$ depends on $M_0$ and $M_1$.

By the latter and by \eqref{oper-2ord-324} we have

\begin{equation}\label{oper-2ord-331}
\begin{aligned}
2\tau\left\Vert v\right\Vert^2_{1,\tau}&\leq \frac{2\tau
C_1}{(2\pi)^n}\int
Q_{\mu}(\xi,\tau)\left|\widehat{v}(\xi)\right|^2d\xi+
\\&+C\left(\rho^2\tau+\tau^{-1}\right)\left\Vert
v\right\Vert^2_{1,\tau}+C\tau^{-1}\left\Vert
p_2(\cdot,\partial,\tau)v\right\Vert^2_{L^2(B_1)},
\end{aligned}
\end{equation}
for every $v\in C_0^{\infty}(B_{\rho})$.

\medskip

In order to estimate from above the first term on the right--hand side in \eqref{oper-2ord-331} we proceed as in the proof of Theorem \ref{Carl-2ord-costanti} (see \eqref{oper-2ord-142}) and
by \eqref{oper-2ord-321}, \eqref{oper-2ord-322}  we get

\begin{equation*}
\begin{aligned}
&\frac{2\tau}{(2\pi)^n}\int
Q_{\mu}(\xi,\tau)\left|\widehat{v}(\xi)\right|^2d\xi= R_{0}+2\tau
e^{-2\mu\varphi(0)}\int
e^{2\mu\varphi(x)}q^{(\delta)}_{\mu}(x,v,\nabla v,\tau)dx\leq\\&\leq
C\rho\tau\left\Vert v\right\Vert^2_{1,\tau}+2\tau
e^{-2\mu\varphi(0)}\int
e^{2\mu\varphi(x)}q^{(\delta)}_{\mu}(x,v,\nabla v,\tau)dx,
\end{aligned}
\end{equation*}
on the other hand, by \eqref{oper-2ord-315} we have

\begin{equation*}
\begin{aligned} 2\tau\int e^{2\mu\varphi}q^{(\delta)}_{\mu}(x,v,\nabla
v,\tau)dx\leq C\delta\tau \left\Vert v\right\Vert^2_{1,\tau}+2\int
\left(S_\tau v A_\tau v \right) e^{2\mu\varphi},
\end{aligned}
\end{equation*}
hence, recalling \eqref{oper-2ord-311},

\begin{equation*}\label{oper-2ord-332}
\begin{aligned}
\frac{2\tau}{(2\pi)^n}\int
Q_{\mu}(\xi,\tau)&\left|\widehat{v}(\xi)\right|^2d\xi&\leq
C(\delta+\rho)\tau\left\Vert v\right\Vert^2_{1,\tau}+C\int
\left|p_2(x,\partial,\tau) v\right|^2e^{2\mu\varphi},
\end{aligned}
\end{equation*}
for every $v\in C_0^{\infty}(B_{\rho})$. By inserting the latter in
\eqref{oper-2ord-331} we have

\begin{equation*}
\begin{aligned}
2\tau\left\Vert v\right\Vert^2_{1,\tau}&\leq
C\left((\rho+\delta)\tau+\tau^{-1}\right)\left\Vert
v\right\Vert^2_{1,\tau}+C\left(1+\tau^{-1}\right)\left\Vert
p_2(\cdot,\partial,\tau)v\right\Vert^2_{L^2(B_1)},
\end{aligned}
\end{equation*}
for every $v\in C_0^{\infty}(B_{\rho})$. Let
$$\rho_0=\delta_0=\frac{1}{4C},$$ then

\begin{equation}\label{oper-2ord-333}
\begin{aligned}
\tau\left\Vert v\right\Vert^2_{1,\tau}\leq 2C\left\Vert
p_2(\cdot,\partial,\tau)v\right\Vert^2_{L^2(B_1)},
\end{aligned}
\end{equation}
for every $v\in C_0^{\infty}(B_{\rho_0})$, for every $\tau\geq 2$
and for every $\delta\leq \delta_0$. Taking into account  \eqref{oper-2ord-310},
by \eqref{oper-2ord-333} we easily deduce that there exists  $\tau*$
such that

\begin{equation*}
\begin{aligned}
\tau\left\Vert v\right\Vert^2_{1,\tau}\leq 2C\left\Vert
P_2(x,\partial,\tau)v\right\Vert^2_{L^2(B_1)},
\end{aligned}
\end{equation*}
for every $v\in C_0^{\infty}\left(B_{\rho_0}\right)$ and for every $\tau\geq \tau*$. Finally, applying Lemmas \ref{lemma-carlm-ell-I} and
\ref{lemma-carlm-ell-II} we obtain \eqref{oper-2ord-307}.
$\blacksquare$

\subsection{Application to the Cauchy problem}\label{applic-Cauchy}

By Carleman estimate
\eqref{oper-2ord-307}, we can obtain an uniqueness result for the Cauchy problem for the operator

\begin{equation}\label{oper-2ord-342}
P(x,\partial)u=g^{jk}(x)\partial^2_{jk}u+b_k(x)\partial_ku+c(x)u,
\end{equation}
where matrix $\left\{g^{jk}\right\}$ has real entries and
satisfy conditions \eqref{oper-2ord-101a} and
\eqref{oper-2ord-101b}, $b_k\in L^{\infty}(B_1, \mathbb{C})$,
$k=1,\cdots, n$, $c\in L^{\infty}(B_1, \mathbb{C})$. 
Below we state the theorem; the proof is only briefly mentioned
as it is carried out in an analogous way to that of Theorem
\ref{Cauchy-ell-teor}. Theorem \ref{oper-2ord-341} has been
proved by \textbf{Calder\'{o}n} in 1957 for
more general (but with coefficients $C^{\infty}$) than those
considered here (for further discussion, we refer to \cite[Ch.
3]{Lern}).

\begin{theo}[\textbf{Calder\'{o}n}]\label{oper-2ord-341}
	\index{Theorem:@{Theorem:}!- Calder\'{o}n@{- Calder\'{o}n}} 
Let $\psi\in C^1\left(\overline{B_1}\right)$ be a real--valued function such that
\begin{equation*}
\nabla\psi(0)\neq 0.
\end{equation*}
Let $P(x,\partial)$ be operator \eqref{oper-2ord-342}. Let $U\in
H^2\left(B_1\right)$  such that
\begin{equation*}
\begin{cases}
P(x,\partial)U=0,\quad\mbox{in }\quad B_{1},\\
\\
U(x)=0, \quad\mbox{in } \left\{x\in B_{1}: \quad \psi(x)>\psi(0)
\right\}.
\end{cases}
\end{equation*}
Let us suppose that
\begin{equation}\label{oper-2ord-343}
	\begin{aligned}
		&\begin{cases}
P_2(0,\xi+i\tau \nabla \psi(0))=0,\\
\\
(\xi,\tau)\neq (0,0),
\end{cases}\Longrightarrow \\& \Longrightarrow P^{(j)}_2(0,\xi+i\tau \nabla
\psi(0))\partial_j\psi(0)\neq 0,
\end{aligned}
\end{equation}
($P_2(x,\xi)=g^{jk}(x)\xi_j\xi_k$).

Then there exists a neighborhood $\mathcal{U}_{0}$ of $0$ such that

\begin{equation*}
U=0\quad\mbox{in }\quad \mathcal{U}_{0}.
\end{equation*}
\end{theo}

\textbf{Proof.} The proof is carried out in a manner similar to
that of Theorem \ref{Cauchy-ell-teor}. Therefore, first of all we
write Carleman estimate \eqref{oper-2ord-307} in the form

\begin{equation*}
\begin{aligned}
\tau^{3}\delta^{4}&\int \left\vert
u\right\vert^2e^{2\tau\varphi\left(\delta^{-1}X\right)}dX+\tau\delta^{2}\int
\left\vert \nabla
u\right\vert^2e^{2\tau\varphi\left(\delta^{-1}X\right)}dX
\leq\\&\leq C_0\int \left\vert
P_2(X,\partial)u\right\vert^2e^{2\tau\varphi\left(\delta^{-1}X\right)}dX,
\end{aligned}
\end{equation*}
for every $\delta\in(0,\delta_0]$, for every $u \in
C^{\infty}_0(B_{\widetilde{\rho}_0}(0))$ and for every $\tau\geq
\widetilde{\tau}_0$.

Next, by means of a diffeomorphism we reduce to the case
in which $\left\{x\in B_{r_0}: \quad \psi(x)\leq 0 \right\}$ is
the epigraph of a function $f$ strictly convex and such that
$f(0)=|\nabla f(0)|=0$. Like the proof of Theorem
\ref{Cauchy-ell-teor} we introduce the function

$$\varphi(x)=h(\delta_0x),$$
where

$$h(x_n)=-x_n+\frac{x_n^2}{2}$$
and we check that $\varphi$ satisfies condition ($\mathbf{S}$). To
check this, first we check that $\varphi$ is
pseudo--convex. For this purpose, it suffices to notice that by
\eqref{oper-2ord-343} we have that if $\tau=0$ then $\xi\neq 0$.  
Therefore

\begin{equation*}
\begin{cases}
P_2(0,\xi)=0,\\
\\
\xi\neq 0,
\end{cases}\Longrightarrow P^{(n)}_2(0,\xi)\neq 0,
\end{equation*}
the latter, in particular, implies that the antecedent of implication
\eqref{oper-2ord-332-2} is not satisfied by $\varphi$ in $0$
which, in turn, implies that the condition \eqref{oper-2ord-106}
holds. Regarding \eqref{oper-2ord-332-3},
we have that if (recall $\nabla\varphi(0)=-\delta_0e_n$)
\begin{equation*}
\begin{aligned}
&\begin{cases}
P_2(0,\xi-i\tau\delta_0e_n)=0, \\
\\
\tau\neq 0,
\end{cases}
\end{aligned}
\end{equation*}
then, \eqref{oper-2ord-343} implies

$$P^{(n)}_2(0,\xi-i\tau\delta_0e_n)\neq 0,$$ hence
$$\partial^2_{jk}\varphi(0)P^{(j)}_2(0,
\xi-i\tau\delta_0e_n)\overline{P^{(k)}_2(\xi-i\tau\delta_0e_n)}=\delta_0^2\left|P^{(n)}_2(0,\xi-i\tau\delta_0e_n)\right|^2>0.$$
The remaining part of the proof is identical to that of
Theorem \ref{Cauchy-ell-teor}. $\blacksquare$

\bigskip

\underline{\textbf{Let us examine condition \eqref{oper-2ord-343}}}.

Let us first notice that if $N:=\nabla\psi(0)$ satisfies
\eqref{oper-2ord-343}, then $N$ cannot be a characteristic direction. As a matter of fact, let us suppose the opposite, that is let us suppose that

\begin{equation}\label{oper-2ord-345}
P_2(0,N)=0.
\end{equation}
Let $\xi=0$, then for every $\tau\neq 0$ we have
$$P_2(0,0+i\tau N)=0,$$
hence the antecedent of \eqref{oper-2ord-343} holds true, but
(by Euler Theorem on homogeneous function)
$$P^{(j)}_2(0,0+i\tau N)N_j=2i\tau P_2(0,N)=0.$$
Therefore, \textbf{a necessary condition in order that \eqref{oper-2ord-343} holds true is that the surface
$\left\{\psi(x)=\psi(0)\right\}$ is noncharacteristic in $0$}.
Nevertheless, as we are going to see, the converse is not true, i.e. the condition
$P_2(0,N)\neq 0$ is not sufficient for the validity of
\eqref{oper-2ord-343}.

Let us start with the following
\begin{prop}\label{oper-2ord-346}
If	$N:=\nabla\psi(0)$ satisfies
\begin{equation}\label{oper-2ord-345-ex}
P_2(0,N)\neq 0.
\end{equation}
then condition \eqref{oper-2ord-343} is equivalent to
 (we set $g^{jk}=g^{jk}(0)$, $j,k=1,\cdots, n$)
\begin{equation}\label{oper-2ord-347}
\begin{cases}
\mathbf{g}(\xi,\xi)=0,\\
\\
\xi\nparallel N,
\end{cases}\Longrightarrow \mathbf{g}(\xi,N)\neq 0,
\end{equation}
$\xi\nparallel N$ mean "$\xi$ and $N$ linearly independent" \index{$\nparallel$}.
\end{prop}
\textbf{Proof.} We begin by noticing that the
\eqref{oper-2ord-343} is equivalent to
\begin{equation}\label{oper-2ord-348}
\begin{cases}
P_2(0,\xi+i\tau N)=0,\\
\\
\xi\nparallel N,\\
\\
(\xi,\tau)\neq (0,0),
\end{cases}\Longrightarrow P^{(j)}_2(0,\xi+i\tau N)N_j\neq 0.
\end{equation}
Indeed, if \eqref{oper-2ord-343} holds then trivially 
\eqref{oper-2ord-348} holds. Let us suppose now that 
\eqref{oper-2ord-348} is true, we have again that \eqref{oper-2ord-343} is trivially satisfied, for any $\xi\nparallel N$. Instead, if
$\xi=\lambda N$, $\lambda \in \mathbb{R}$, we have that
$$P_2(0,\xi+i\tau N)=0,$$
implies
$$0=P_2(0,\lambda N+i\tau N)=(\lambda+i\tau)^2P_2(0,N)$$
from which, taking into account that $P_2(0,N)\neq 0$, we have $\tau=\lambda=0$,
that is $\xi=0$ and $\tau=0$, consequently the antecedent of
\eqref{oper-2ord-343} is false and thus the condition
\eqref{oper-2ord-343} is satisfied.

\medskip

Now we prove that \eqref{oper-2ord-347} and \eqref{oper-2ord-348}
are equivalent. Let us begin assuming that
\eqref{oper-2ord-348} holds true. To prove \eqref{oper-2ord-347}
let us suppose that $\xi\nparallel N$ and that
$$\mathbf{g}(\xi,\xi)=0.$$
Now, if it were

\begin{equation}\label{oper-2ord-349}
\mathbf{g}(\xi,N)=0,
\end{equation}
we would have at the same time $$\xi\nparallel N,$$ $$P_2(0,\xi+i0
N)=\mathbf{g}(\xi,\xi)=0$$ and $$P^{j}_2(0,\xi+i0
N)N_j=2\mathbf{g}(\xi,N)=0.$$ So there would be a contradiction with
\eqref{oper-2ord-348}, therefore \eqref{oper-2ord-349} does not hold. Thus, if  \eqref{oper-2ord-348} holds true then \eqref{oper-2ord-347} holds true.

Now let us suppose that \eqref{oper-2ord-347} holds and let us suppose that
\begin{equation}\label{oper-2ord-350}
\begin{cases}
P_2(0,\xi+i\tau N)=0,\\
\\
\xi\nparallel N,\\
\\
(\xi,\tau)\neq (0,0).
\end{cases}
\end{equation}
If $\tau=0$, then, \eqref{oper-2ord-350} implies $\xi\nparallel
N$ and $\mathbf{g}(\xi,\xi)=P_2(0,\xi+i0 N)=0$, moreover 
\eqref{oper-2ord-347} implies $\mathbf{g}(\xi,N)\neq 0$. Hence
$$P^{(j)}_2(0,\xi+i0 N)=2\mathbf{g}(\xi,N)\neq 0.$$
If $\tau\neq 0$, recalling that $N$ is a noncharacteristic direction -- so that $\mathbf{g}(N,N)=P_2(0, N)\neq 0$ -- we have
$$P^{(j)}_2(0,\xi+i\tau N)N_j=2\mathbf{g}(\xi,N)+i\tau
\mathbf{g}(N,N)\neq 0.$$ All in all, if 
\eqref{oper-2ord-347} holds then \eqref{oper-2ord-348} holds. The proof is complete. $\blacksquare$

\bigskip

\underline{\textbf{Remarks and Examples.}}

\medskip

\noindent\textbf{1.} If $P_2(x,\partial)$ is ellipic (with real coefficients), then \eqref{oper-2ord-343} is
satisfied, as already proved in Example 4a of Section \ref{altre consid}.

\medskip

\noindent\textbf{2.} Let us consider the wave operator

\begin{equation}\label{oper-2ord-351}
P_2(\partial)=\Delta_{x'}-\partial^2_{x_n}.
\end{equation}
Let us check for which  $N\neq 0$ condition \eqref{oper-2ord-343} 
is satisfied. It is not restrictive to assume 
\begin{equation}\label{oper-2ord-351-bis}
|N|=1.
\end{equation}
 We first need to assume that $N$ is a noncharacteristic direction for
$P_2(\partial)$, i.e.

\begin{equation}\label{oper-2ord-352}
P_2(N)=|N'|^2-N_n^2\neq 0.
\end{equation}
Now, let us see when \eqref{oper-2ord-347} is satisfied (which,
we recall, is equivalent to \eqref{oper-2ord-343}). This
condition can be written

\begin{equation}\label{oper-2ord-353}
\begin{cases}
|\xi'|^2-\xi_n^2=0,\\
\\
\xi\nparallel N,
\end{cases}\Longrightarrow \xi'\cdot N'-\xi_n N_n\neq 0.
\end{equation}
Let us first examine the \textbf{case $n=2$}. In this case
condition \eqref{oper-2ord-352} can be written
\begin{equation}\label{oper-2ord-354}
N_1^2-N_2^2\neq 0.
\end{equation}
the first condition of the antecedent of \eqref{oper-2ord-353} can be written
$\xi_1^2-\xi_2^2=0$  and it is equivalent to
$$\xi_2=\pm \xi_1.$$
From which,  by \eqref{oper-2ord-352} and $\xi\neq
0$, we have
$$\xi_1 N_1-\xi_2 N_2=\xi_1\left(N_1\mp N_2\right)\neq 0.$$

Therefore, if $n=2$, \eqref{oper-2ord-343} is satisfied for
all $N\in \mathbb{R}^2$ that satisfy \eqref{oper-2ord-354}, i.e. that are not a characteristic direction.

\medskip

Let us consider now the \textbf{case $n\geq 3$}. Since the scalar product in $\mathbb{R}^{n-1}$ is invariant w.r.t. the rotations, we may assume that

$$N=N_1e_1+N_ne_n.$$
In this way, conditions \eqref{oper-2ord-351-bis} and
\eqref{oper-2ord-352}, can be written, respectively;
$$N_1^2+N_n^2=1$$
and

\begin{equation*}
N_1^2-N_n^2\neq 0.
\end{equation*}

Let us distinguish two cases

\medskip

\noindent \textbf{(a)} $N_1^2-N_n^2< 0$, i.e. $|N'|<|N_n|$;

\medskip

\noindent \textbf{(b)} $N_1^2-N_n^2>0$, i.e. $|N'|>|N_n|$;

\medskip

\noindent \textbf{Case (a)}. Let $\xi\nparallel N$ such that
$$|\xi'|^2-\xi_n^2=0.$$
In particular, we have  $\xi_n\neq 0$ and $\xi'\neq 0$ (as a matter of fact,
if one of them is zero the other is also zero) moreover
$$|\xi'|=|\xi_n|.$$
Hence

$$|\xi'\cdot N'|\leq
|\xi'||N'|=|\xi'||N_1|<|\xi'||N_n|=|\xi_n||N_n|,$$ from which we have

$$\xi'\cdot N'-\xi_n N_n\neq 0.$$
Therefore, in case (a), \eqref{oper-2ord-347} is satisfied.

In case (b) it is simple to check that  \eqref{oper-2ord-347}
is not satisfied. To check this, Let

$$\xi_0=e_1N_n+e_2\sqrt{N_1^2-N_n^2}+e_n N_1.$$
We have

\begin{equation*}
\begin{cases}
|\xi'_0|^2-\xi_{0,n}^2=0,\\
\\
\xi_0\nparallel N,
\end{cases}
\end{equation*}
but
$$\xi'\cdot N'-\xi_n N_n=N_nN_1-N_1N_n=0.$$

\medskip

Now it is interesting to point out (we refer to \cite[Ch.
6]{Lern}) that it has been proved that there exist $u, q\in
C^{\infty}(\mathbb{R}^3, \mathbb{C})$ such that

\begin{equation}\label{oper-2ord-354-29}
\begin{cases}
\partial^2_{t}u-\partial^2_{x_1}u-\partial^2_{x_2}u+q(x,t)u=0,\\
\\
\mbox{supp } u=\left\{x_2\geq 0\right\},
\end{cases}
\end{equation}
or, in other words, \textbf{does not hold uniqueness} for the Cauchy problem with initial surface  $\Gamma:=\left\{x_2 =0\right\}$, for the equation
$$\partial^2_{t}u-\partial^2_{x_1}u-\partial^2_{x_2}u+q(x,t)u=0$$ where
$$q\in C^{\infty}(\mathbb{R}^3, \mathbb{C}).$$ Keep in mind that, when $q$ is
\textbf{analytic}, the Holmgren Theorem provides  uniqueness
for the Cauchy problem for the equation
$$\partial^2_{t}u-\partial^2_{x_1}u-\partial^2_{x_2}u+q(x,t) u=0,$$
with initial surface $\Gamma$ (as it is a noncharacteristic surface).

\bigskip

\noindent\textbf{3.} We check that if $\theta\in (0,1)$, then
\begin{equation}\label{oper-2ord-355}
\phi(x)=\frac{1}{2}\left(|x'|^2-\theta^2x_n^2\right),
\end{equation}
is a  pseudo--convex function w.r.t. wave operator
\eqref{oper-2ord-351}  in the open set $$\Omega=\left(B_R\setminus
\overline{B_r}\right)\times (-T,T),$$ for every $0<r<R$ and $T>0$.

We only need to check \eqref{oper-2ord-106-bis}. Calculate
$$\nabla\phi(x)=\left(x',-\theta^2x_n\right)\neq 0,\quad\mbox{in } \overline{\Omega},$$
$$\partial^2\phi(x)=\mbox{diag
}\left(1,\cdots,1,-\theta^2x_n\right).$$ Now let us suppose
\begin{equation}\label{oper-2ord-356}
\begin{aligned}
\begin{cases}
|\xi'|^2-\xi^2_n=0, \\
\\
\xi'\cdot x'+\theta^2\xi_nx_n= 0,\\
\\
\xi\neq 0
\end{cases}
\end{aligned}
\end{equation}
and let us check that
\begin{equation}\label{oper-2ord-357}
\partial^2\phi(x)\widetilde{\xi}\cdot \widetilde{\xi}>0,
\end{equation}
where
$$\widetilde{\xi}=\left(\xi',-\theta^2\xi_n\right).$$ We have
\begin{equation*}
\partial^2\phi(x)\widetilde{\xi}\cdot
\widetilde{\xi}=|\xi'|^2-\theta^2\xi^2_n.
\end{equation*}
On the other hand, by the first condition of \eqref{oper-2ord-356}
we have $|\xi'|^2=\xi^2_n$ hence, taking into account that if $\xi$
satisfies at the same time $\xi\neq 0$ and $|\xi'|^2-\xi^2_n=0$, then $\xi'\neq 0$, we have

\begin{equation*}
\partial^2\phi(x)\widetilde{\xi}\cdot
\widetilde{\xi}=|\xi'|^2\left(1-\theta^2\right)>0.
\end{equation*}
Therefore, we have proved that \eqref{oper-2ord-356} implies
\eqref{oper-2ord-357}. Hence $\varphi$ is a pseudo--convex  function in $\Omega$. $\blacklozenge$

\bigskip

In the following Theorem we prove the uniqueness for a Cauchy problem under the assumption of pseudo-convexity of the initial surface.

\begin{theo}\label{241221-1}
Let us suppose that the coefficients of the principal part 
$P_2(\partial)$, of operator \eqref{oper-2ord-342} are
costants. Let $U\in H^2\left(B_1(x_0)\right)$ satisfy
\begin{equation*}
\begin{cases}
P(\partial)U=0,\quad\mbox{in }\quad B_{1}(x_0)\\
\\
U(x)=0 \quad\mbox{in } \left\{x\in B_{1}(x_0): \quad
\psi(x)>\psi(x_0) \right\}.
\end{cases}
\end{equation*}
Let $\psi\in C^2\left(\overline{B_1(x_0)}\right)$ be a real--valued function and pseudo--convex w.r.t. $P_2(\partial)$ in $x_0$. Then there exists a neighborhood $\mathcal{U}_{x_0}$ of $0$ such that
\begin{equation*}
U=0\quad\mbox{in }\quad \mathcal{U}_{x_0}.
\end{equation*}
\end{theo}

\textbf{Proof.} The proof is very similar to that
of Proposition \ref{qr1-12}. Therefore, here we merely
point out the most important differences inviting the reader to
care the details.

It is not restrictive to assume that $x_0=0$ and
$$\psi(0)=0.$$ Since $\psi$ is pseudo--convex in $0$, by Proposition \ref{oper-2ord-114} we have that the function
$$\varphi(x)=e^{\lambda \psi(x)}-1,$$
satisfies condition $(\mathbf{S})$ in $0$ for any $\lambda$
large enough. Let

\begin{equation}\label{241221-2}
\varphi_{\varepsilon}(x)=\varphi(x)-\frac{\varepsilon|x|^2}{2},
\end{equation}
where $\varepsilon$ is a positive number that can be
chosen in such a way that $\varphi_{\varepsilon}$ satisfies
condition $(\mathbf{S})$ (the reader cure the details). We fix this $\varepsilon$ and from Theorem \ref{Carl-2ord-costanti} we have
that there exists $R\in (0,1/2)$ and there exist two constants $C$ and $\tau_0$
such that

\begin{equation}\label{241221-3}
\begin{aligned}
& \tau^3 \int_{B_1} |u|^2 e^{2 \tau \varphi_{\varepsilon}}dx +
\tau \int_{B_1}
 |\nabla u|^2 e^{2 \tau \varphi_{\varepsilon}}dx \leq C
\int_{B_1} |P_2 (\partial)u|^2 e^{2 \tau \varphi_{\varepsilon}}dx,
\end{aligned}
\end{equation}

\smallskip

\noindent for every $u \in C^{\infty}_0(B_{2R}(0))$ and for every $\tau\geq
\tau_0$. Starting from this point one repeats, with obvious
modifications, what we did in the proof of Proposition 
\ref{qr-unic-primo ordine}. $\blacksquare$

\section{Stability estimate for the wave equation in a cylinder}\label{corr:26-4-23}

In this Section we adopt the traditional notations: the "spatial coordinates" are denoted by
$x_1,\cdots,x_n$, the time coordinate is denoted by $t$ and,
therefore, the wave operator is

\begin{equation}\label{oper-2ord-358}
\square=\partial_t^2-\Delta_x=\partial_t^2-\left(\partial_{x_1}^2+\cdots+\partial_{x_1}^2\right).
\end{equation}
Let us denote by
$$\nabla_x=\left(\partial_{x_1},\cdots,\partial_{x_1}\right),\quad\quad
\nabla_{x,t}=(\nabla_x,\partial_t).$$ Let $T>1$, set
$$S_T=\mathbb{R}^n\times (-T,T).$$

The following Theorem holds true (see also \cite{RA})

\bigskip

\begin{theo}[\textbf{stability estimate for the wave equation}]\label{oper-2ord-stab}
		\index{Theorem:@{Theorem:}!- stability estimate for the wave equation@{- stability estimate for the wave equation}}
	
Let $a\in L^{\infty}\left(S_T,\mathbb{R}^n\right)$, $b\in
L^{\infty}\left(S_T\right)$ and $c\in L^{\infty}\left(S_T\right)$. Let
$M\geq 1$. Let us assume that
\begin{equation}\label{oper-2ord-359}
\left\Vert
a\right\Vert_{L^{\infty}\left(S_T,\mathbb{R}^n\right)}+\left\Vert
b\right\Vert_{L^{\infty}\left(S_T\right)}+\left\Vert c
\right\Vert_{L^{\infty}\left(S_T\right)}\leq M.
\end{equation}
Let $F\in L^{2}\left(S_T\right)$ and  $U\in
C^{\infty}\left(S_T\right)$ satisfy
\begin{equation}\label{oper-2ord-360}
\begin{cases}
\Box U+a\cdot \nabla_{x}U+
b\partial_tU+cU=F,\quad \mbox{in }
S_T,\\
\\
U(x,t)=0,\quad \mbox{for }|x|>1,\mbox{ } t\in (-T,T).
\end{cases}
\end{equation}

\medskip 

Then
\begin{equation}\label{oper-2ord-361}
\begin{aligned}
\int^{T}_{-T}\int_{B_1} \left(|U|^2+|\nabla_{x,t} U|^2\right)
dxdt\leq C \left\Vert F\right\Vert^2_{L^2\left(B_1\times
(-T,T)\right)},
\end{aligned}
\end{equation}
where $C$ depends on $M$ and $T$.
\end{theo}

\bigskip

\noindent The proof of Theorem \ref{oper-2ord-stab} is
based on what follows:

\smallskip

\noindent(a) an \textbf{energy estimate} for the equation in
\eqref{oper-2ord-360}

\noindent(b) Carleman estimate \eqref{oper-2ord-115}.

\smallskip

\begin{lem}[\textbf{energy estimate}]\label{stima-energia}
		\index{Lemma:@{Lemma:}!- energy estimate for the wave equation@{- energy estimate for the wave equation}}
Let $U\in C^2\left(\overline{B_1}\times (-T,T)\right)$ satisfy

\begin{equation}\label{oper-2ord-360-29}
\begin{cases}
\Box U+a\cdot \nabla_{x}U+
b\partial_tU+cU=F,\quad \mbox{in }
B_1\times
(-T,T),\\
\\
U(x,t)=0,\quad \mbox{for }(x,t)\in\partial B_1\times (-T,T),
\end{cases}
\end{equation}

\medskip

\noindent where $a$, $b$, $c$ and $F$ satisfy the same assumption
of Theorem \ref{oper-2ord-stab}, then the following inequality holds true

\begin{equation}\label{oper-2ord-362}
\begin{aligned}
&\int_{-T}^T\int_{B_1}\left(U^2(x,t)+U_t^2(x,t)+\left| \nabla_x
U(x,t)\right|^2\right)dxdt\leq \\&\leq
CT\delta^{-1}\int^{\delta}_{-\delta}
\int_{B_1}\left(U^2(x,t)+U_t^2(x,t)+\left| \nabla_x
U(x,t)\right|^2\right)dxdt+\\&+ C\left\Vert
F\right\Vert^2_{L^2\left(B_1\times (-T,T)\right)},
\end{aligned}
\end{equation}
where $C$ depends on $M$ and $T$ only.
\end{lem}
\textbf{Proof of Lemma \ref{stima-energia}. } By applying 
Lemma \ref{re91221} with
$$\beta=(0,\cdots,0,1)$$
and
$$g=\mbox{diag }(-1,\cdots,-1,1),$$
we get

\begin{equation}\label{oper-2ord-363}
(\Box U)U_t=\frac{1}{2}\partial_t\left(U^2_t+\left| \nabla_x
U\right|^2\right)-\mbox{div }_x\left(U_t\nabla_x U\right).
\end{equation}
Moreover, by \eqref{oper-2ord-360-29} and taking into account 
\eqref{oper-2ord-359} we have easily

\begin{equation}\label{oper-2ord-364}
\left|(\Box U)U_t\right|\leq C(M+1)\left(U^2+U^2_t +\left| \nabla_x
U\right|^2+F^2\right).
\end{equation}
Let $s,\sigma\in (-T,T)$,  $s\leq\sigma$. By integrating both the sides of \eqref{oper-2ord-363} over $B_1\times [s, \sigma]$,
we have from the divergence Theorem and from \eqref{oper-2ord-364}

\begin{equation}\label{oper-2ord-365}
\begin{aligned}
&\int_{B_1}\left(U^2_t(x,\sigma)+\left| \nabla_x
U(x,\sigma)\right|^2\right)dx-\\&-\int_{B_1}\left(U^2_t(x,s)+\left|
\nabla_x U(x,s)\right|^2\right)dx=\\&=
\int_s^{\sigma}\int_{B_1}\partial_t\left(U^2_t+\left| \nabla_x
U\right|^2\right)dxdt=\\&=2\int_s^{\sigma}\int_{B_1}(\Box
U)U_t dxdt\leq\\&\leq C(M+1)\int_s^{\sigma}\int_{B_1}\left(U^2+U^2_t
+\left| \nabla_x U\right|^2+F^2\right)dxdt.
\end{aligned}
\end{equation}
Let us note that by the first Poincar\'{e} inequality (Theorem \ref{Poincar}) we have

\begin{equation}\label{oper-2ord-366}
\begin{aligned}
C_*^{-1}\int_{B_1}\left| \nabla_x U(x,t)\right|^2dx&\leq
\int_{B_1}\left(U^2(x,t)+\left| \nabla_x
U(x,t)\right|^2\right)dx\leq
\\&\leq C_* \int_{B_1}\left| \nabla_x
U(x,t)\right|^2dx,
\end{aligned}
\end{equation}
 for every $t\in (-T,T)$, where $C\geq 1$ is a constant.
 Therefore, setting

 $$E(t)=\int_{B_1}\left(U_t^2(x,t)+\left| \nabla_x
U(x,t)\right|^2\right)dx,$$ by \eqref{oper-2ord-365} e
\eqref{oper-2ord-366} we have

\begin{equation*}
\begin{aligned}
E(\sigma)\leq E(s)+C_1\left\Vert
F\right\Vert^2_{L^2\left(B_1\times
(-T,T)\right)}+C_1\int^{\sigma}_s E(t)dt
\end{aligned}
\end{equation*}
($C_1=C_*C(M+1)$) and by the Gronwall inequality we obtain

\begin{equation}
\begin{aligned}
E(\sigma)\leq \left(E(s)+C_1\left\Vert
F\right\Vert^2_{L^2\left(B_1\times
(-T,T)\right)}\right)e^{2C_1T}.
\end{aligned}
\end{equation}
We notice that if $\sigma\leq s$ we obtain similarly the previous estimate
 if that we integrate both the sides of \eqref{oper-2ord-363}
over $B_1\times [\sigma, s]$ and we interchange
$\sigma$ and $s$. Therefore we have, for each $s,\sigma\in (-T,T)$,

\begin{equation}
\begin{aligned}\label{oper-2ord-367}
E(\sigma)\leq \left(E(s)+C_1\left\Vert
F\right\Vert^2_{L^2\left(B_1\times
(-T,T)\right)}\right)e^{2C_1T}.
\end{aligned}
\end{equation}
Now, by integrating with respect to $s$ both members of
\eqref{oper-2ord-367} over $(-\delta, \delta)$, where $\delta\in (0,T)$,
we have, for each $\sigma\in (-T,T)$,

\begin{equation}
\begin{aligned}\label{oper-2ord-368}
2\delta E(\sigma)\leq 2C_2\int^{\delta}_{-\delta} E(s)ds +2\delta
C_2\left\Vert F\right\Vert^2_{L^2\left(B_1\times (-T,T)\right)}
\end{aligned}
\end{equation}
($C_2=e^{2C_1T}$). Finally, by integrating both the sides of
\eqref{oper-2ord-368} with respect to $\sigma$ over $(-T,T)$ and taking into account \eqref{oper-2ord-366}, we obtain \eqref{oper-2ord-362}. $\blacksquare$

\bigskip

\textbf{Remark 4.} As can be seen immediately from the
proof,  it is not necessary for Lemma \ref{oper-2ord-362}  that $T$
be greater than $1$. $\blacklozenge$

\bigskip

\textbf{Proof of Theorem \ref{oper-2ord-stab}}.
In Remark 3 of the previous Section we have proved that

$$\phi(x,t)=-\theta^2t^2+|x|^2,$$
is a pseudo--convex funtion if $\theta\in (0,1)$ in
$\left(B_R\setminus \overline{B_r}\right)\times (-T,T)$ for every $0<r<R$. Let us fix $\theta$ in such a way that 

\begin{equation}\label{oper-2ord-370}
\frac{1}{T}<\theta<1.
\end{equation}
(recall that $T>1$) and let $\rho\in (0,1/10)$ satisfy
\begin{equation}\label{oper-2ord-372}
\rho<\frac{\theta T-1}{8}.
\end{equation} 
Proposition
\ref{oper-2ord-114} implies that if $\lambda$ is sufficiently
large, then the functions

\begin{equation}\label{oper-2ord-371}
\varphi_0(x,t)=e^{\lambda \phi(x,t)},\quad\quad
\varphi_{1}(x,t)=e^{\lambda \phi(x-8\rho e_1,t)},
\end{equation}
satisfy condition $(\mathbf{S})$ in $\left(B_R\setminus
\overline{B_r}\right)\times (-T,T)$ for every $0<r<1<R$.

Let us fix $\lambda>0$ in such a way that $\varphi_0$ satisfy condition
$(\mathbf{S})$ in $\left(B_2\setminus
\overline{B_{\rho}}\right)\times (-T,T)$ (and, consequently
$\varphi_1$ satisfies condition $(\mathbf{S})$ in
$\left(B_2(8\rho e_1)\setminus \overline{B_{\rho}(8\rho
e_1)}\right)\times (-T,T)$).

For any $s>0$ set

$$Z_{0, s}=\left\{(x,t)\in \overline{B_1}\times
(-T,T):\quad \varphi_0(x,t)\geq e^{\lambda s}\right\}$$ and

$$Z_{1, s}=\left\{(x,t)\in \overline{B_1}\times
(-T,T):\quad \varphi_{1}(x,t)\geq e^{\lambda s}\right\}.$$

Let us check that

\begin{subequations}\label{oper-2ord-374}
\begin{equation}\label{oper-2ord-374a}
\varphi_0(x,\pm T)\leq 1,\quad \varphi_{1}(x,\pm T)\leq 1,\quad
\forall x\in B_1,
\end{equation}
\begin{equation*}
\left.\right.
\end{equation*}
\begin{equation}
\label{oper-2ord-374b} \overline{B_1}\times [-2\rho,2\rho]\subset
Z_{0, s} \cup Z_{1, s},\quad \forall s\in(0, 12\rho^2].
\end{equation}
\end{subequations}

\medskip

The first of \eqref{oper-2ord-374a} is an immediate consequence
of \eqref{oper-2ord-370}. Concerning the second
of \eqref{oper-2ord-374a}, we observe that from
\eqref{oper-2ord-372} we have, for each $x\in B_1$

$$|x-8\rho e_1|\leq 1+8\rho <\theta T$$ which implies

$$\varphi_{1}(x,\pm T)=e^{\lambda\left(-\theta^2T^2+|x-8\rho
e_1|^2\right)}\leq 1,$$ for every $x\in B_1$.

Now let us check \eqref{oper-2ord-374b}. Set
$$\psi(x)=\max\left\{|x|^2,|x-8\rho e_1|^2\right\},$$ we obtain easily

\begin{equation}\label{oper-2ord-375}
Z_{0,s} \cup Z_{1, s}=\left\{(x,t)\in \overline{B_1}\times
(-T,T):\quad e^{\lambda(-\theta^2t^2+\psi(x))}\geq e^{\lambda
s}\right\}.
\end{equation}
Let us note now that $\psi$ can be written as 

\begin{equation*}
\psi(x)=\begin{cases} |x|^2,\quad \mbox{for }x_1\geq 4\rho,\\
\\
(x_1-8\rho)^2+ x_2^2+\cdots+ x_n^2,\quad\mbox{for } x_1<4\rho,
\end{cases}
\end{equation*}
from which we have
\begin{equation*}
\psi(x)\geq 16\rho^2,\quad\forall x\in \mathbb{R}^n.
\end{equation*}
Now, if $(x,t)\in \overline{B_1}\times [-2\rho,2\rho]$, then

\begin{equation*}
-\theta^2t^2+\psi(x)\geq
-\theta^2t^2+16\rho^2>-4\rho^2+16\rho^2=12\rho^2\geq s, \quad
\forall s\in(0,12\rho^2],
\end{equation*}
by the latter and by \eqref{oper-2ord-375}, we get $(x,t)\in
Z_{0,s} \cup Z_{1, s}$ for $0\leq s\leq 12\rho^2$. Hence
\eqref{oper-2ord-374b} is proved.

\bigskip

Let us apply Carleman estimate \eqref{oper-2ord-115}
to the operator $\Box$ where $\varphi=\varphi_0$. Set
$$Q_T=\overline{B_1}\times (-T,T),$$ we get

\begin{equation}\label{oper-2order-376}
\begin{aligned}
 \tau^3 \int_{Q_T} |u|^2 e^{2 \tau \varphi_0}dxdt &+ \tau \int_{Q_T}
 |\nabla_{x,t} u|^2 e^{2 \tau \varphi_0}dxdt \leq \\& \leq C
\int_{Q_T} |\Box u|^2 e^{2 \tau \varphi_0}dxdt,
\end{aligned}
\end{equation}

\medskip

\noindent for every $u \in C^{\infty}_0\left(\mathbb{R}^{n+1}\right)$, such that supp $u$ $\subset
Q_T=\overline{B_1}\times (-T,T) $ and for every $\tau\geq \tau_0$.
Let $\widetilde{\eta}\in C^{\infty}(\mathbb{R})$ satisfy

$$\widetilde{\eta}(r)=0,\mbox{ } r\leq 9\rho^2;\quad 0\leq \widetilde{\eta}(r)\leq
1,\mbox{ } 9\rho^2<r<10\rho^2;\quad \widetilde{\eta}(r)=1, \mbox{ }
r\geq 10\rho^2;$$

\begin{equation}\label{oper-2order-377}
\left|\frac{d\widetilde{\eta}}{dr}\right|\leq C\rho^{-2},\quad
\left|\frac{d^2\widetilde{\eta}}{dr^2}\right|\leq C\rho^{-4},
\end{equation}
where $C$ is a constant (independent by $\rho$). Set

$$\eta(x,t)=\widetilde{\eta}(-\theta^2t^2+|x|^2)$$ and let us apply estimate \eqref{oper-2order-376} to $U\eta$. By \eqref{oper-2ord-359} we have

\begin{equation*}
\begin{aligned}
\left|\Box(\eta U)\right|&=\left|\eta\Box
U+2\left(\partial_t\eta\partial_tU-\nabla_x\eta\cdot\nabla_xU\right)+U(\Box
\eta)\right|\leq\\&\leq \eta
M(|\nabla_{x,t}U|+|U|)+\eta|F|+C\rho^{-2}\chi_{Z_{0,9\rho^2}\setminus
Z_{0,10\rho^2}}|\nabla_{x,t}U|+\\&+C\rho^{-4}\chi_{Z_{0,9\rho^2}\setminus
Z_{0,10\rho^2}}|U|,
\end{aligned}
\end{equation*}
where $C$ depends on $T$. Let us observe that (by inserting what was obtained
in \eqref{oper-2order-376}), we have

\begin{equation}\label{oper-2order-378}
\begin{aligned}
& \int_{Q_T} \left(\tau^3|U\eta|^2+\tau|\nabla_{x,t} (\eta
U)|^2\right) e^{2 \tau \varphi_0}dxdt \leq\\&\leq CM^2 \int_{Q_T}
\eta^2 \left(|\nabla_{x,t}U|^2+|U|^2\right) e^{2 \tau
\varphi_0}dxdt+\\&+C\int_{Q_T}|F|^2e^{2 \tau \varphi_0}dxdt+\\&
+C\rho^{-8}\int_{Z_{0,9\rho^2}\setminus Z_{0,10\rho^2}}|U|^2e^{2
\tau \varphi_0}dxdt+\\&+C\rho^{-4}\int_{Z_{0,9\rho^2}\setminus
Z_{0,10\rho^2}}|\nabla_{x,t}U|^2e^{2 \tau \varphi_0}dxdt,
\end{aligned}
\end{equation}
for every $\tau\geq \tau_0$.

Now, let us estimate from below the left--hand side of \eqref{oper-2order-378}

\begin{equation}\label{oper-2order-379}
\begin{aligned}
 &\int_{Q_T} \left(\tau^3|U\eta|^2+\tau|\nabla_{x,t} (\eta
U)|^2\right) e^{2 \tau \varphi_0}dxdt \geq\\&\geq
 \tau\int_{Z_{0,10\rho^2}} \left(|U|^2+|\nabla_{x,t}
U|^2\right) e^{2 \tau \varphi_0}dxdt
\end{aligned}
\end{equation}
and let us estimate from above first integral on the right--hand side as follows

\begin{equation}\label{oper-2order-380}
\begin{aligned}
&\int_{Q_T} \eta^2 \left(|\nabla_{x,t}U|^2+|U|^2\right) e^{2 \tau
\varphi_0}dxdt=\\&= \int_{Z_{0,10\rho^2}}
\left(|\nabla_{x,t}U|^2+|U|^2\right) e^{2 \tau
\varphi_0}dxdt+\\&+\int_{Z_{0,9\rho^2}\setminus Z_{0,10\rho^2}}
\eta^2 \left(|\nabla_{x,t}U|^2+|U|^2\right) e^{2 \tau
\varphi_0}dxdt\leq \\&\leq \int_{Z_{0,10\rho^2}}
\left(|\nabla_{x,t}U|^2+|U|^2\right) e^{2 \tau \varphi_0}dxdt+\\&+
e^{2 \tau e^{10\lambda \rho^2}}\int_{Z_{0,9\rho^2}\setminus
Z_{0,10\rho^2}} \left(|\nabla_{x,t}U|^2+|U|^2\right) dxdt.
\end{aligned}
\end{equation}
Using \eqref{oper-2order-379} and \eqref{oper-2order-380} in
\eqref{oper-2order-378}, we obtain, by simple calculations (recall
$\rho<1$)

\begin{equation}\label{oper-2order-381}
\begin{aligned}
&(\tau-CM^2)\int_{Z_{0,10\rho^2}} \left(|U|^2+|\nabla_{x,t}
U|^2\right) e^{2 \tau \varphi_0}dxdt \leq \\&\leq C\int_{Q_T}|F|^2e^{2 \tau
\varphi_0}dxdt+\\& +C\rho^{-8}e^{2 \tau e^{10\lambda
\rho^2}}\int_{Z_{0,9\rho^2}\setminus
Z_{0,10\rho^2}}(\left(|U|^2+|\nabla_{x,t}U|^2\right)dxdt,
\end{aligned}
\end{equation}
for every $\tau\geq \tau_0$. We set $\tau_1=\max\left\{\tau_0, (2C
M^2)^{-1} \right\}$ and by \eqref{oper-2order-381} we obtain 

\begin{equation}\label{oper-2order-382}
\begin{aligned}
&\int_{Z_{0,10\rho^2}} \left(|U|^2+|\nabla_{x,t} U|^2\right) e^{2
\tau \varphi_0}dxdt \leq C\int_{Q_T}|F|^2e^{2 \tau
\varphi_0}dxdt+\\& +C\rho^{-8}e^{2 \tau e^{10\lambda
\rho^2}}\int_{Z_{0,9\rho^2}\setminus
Z_{0,10\rho^2}}(\left(|U|^2+|\nabla_{x,t}U|^2\right)dxdt,
\end{aligned}
\end{equation}
for every $\tau\geq \tau_1$. Now, in \eqref{oper-2order-382}, we estimate trivially from below the integral on the left--hand side and we estimate trivially the integrals on the right--hand side. We get
\begin{equation*}
\begin{aligned}
&e^{2 \tau e^{12\lambda \rho^2}}\int_{Z_{0,12\rho^2}}
\left(|U|^2+|\nabla_{x,t} U|^2\right)
dxdt\leq\\&\leq\int_{Z_{0,10\rho^2}} \left(|U|^2+|\nabla_{x,t}
U|^2\right) e^{2 \tau \varphi_0}dxdt \leq \\&\leq Ce^{2 \tau e^{\lambda
}}\int_{Q_T}|F|^2e^{2 \tau \varphi_0}dxdt+\\& +C\rho^{-8}e^{2 \tau
e^{10\lambda \rho^2}}\int_{Z_{0,9\rho^2}\setminus
Z_{0,10\rho^2}}(\left(|U|^2+|\nabla_{x,t}U|^2\right)dxdt,
\end{aligned}
\end{equation*}
which implies

\begin{equation}\label{oper-2order-383-29}
\begin{aligned}
&\int_{Z_{0,12\rho^2}} \left(|U|^2+|\nabla_{x,t} U|^2\right)
dxdt\leq\\&\leq\int_{Z_{0,10\rho^2}} \left(|U|^2+|\nabla_{x,t}
U|^2\right) e^{2 \tau \varphi_0}dxdt \leq \\&\leq  Ce^{2 \tau
\left(e^{\lambda }-e^{12\lambda
\rho^2}\right)}\int_{Q_T}|F|^2dxdt+\\& +C\rho^{-8}e^{2 \tau
\left(e^{10\lambda \rho^2}-e^{12\lambda
\rho^2}\right)}\int_{Z_{0,10\rho^2}\setminus
Z_{0,9\rho^2}}(\left(|U|^2+|\nabla_{x,t}U|^2\right)dxdt\leq\\&\leq
Ce^{2 \tau \left(e^{\lambda }-e^{12\lambda
\rho^2}\right)}\int_{Q_T}|F|^2dxdt+\\& +C\rho^{-8}e^{2 \tau
\left(e^{10\lambda \rho^2}-e^{12\lambda
\rho^2}\right)}\int_{Q_T}\left(|U|^2+|\nabla_{x,t}U|^2\right)dxdt,
\end{aligned}
\end{equation}

\medskip

for every $\tau\geq \tau_1$. By Lemma \ref{stima-energia} and by
\eqref{oper-2order-383-29} we have

\begin{equation}\label{oper-2order-384}
\begin{aligned}
&\int_{Z_{0,12\rho^2}} \left(|U|^2+|\nabla_{x,t} U|^2\right)
dxdt\leq Ce^{2 \tau \left(e^{\lambda }-e^{12\lambda \rho^2}\right)}
\left\Vert F\right\Vert^2_{L^2\left(B_1\times (-T,T)\right)}+\\&+
C\rho^{-9}e^{2 \tau \left(e^{10\lambda \rho^2}-e^{12\lambda
\rho^2}\right)}\int^{2\rho}_{-2\rho}\int_{B_1}\left(|U|^2+|\nabla_{x,t}U|^2\right)dxdt,
\end{aligned}
\end{equation}

\smallskip

\noindent for every $\tau\geq \tau_1$, where $C$ depends on $M$ and $T$. 
At this point we note that, by using \eqref{oper-2ord-115} for
operator $\Box$ with $\varphi=\varphi_1$, we obtain an estimate similar to
\eqref{oper-2order-384} .
More precisely we have 

\begin{equation}\label{oper-2order-383}
\begin{aligned}
&\int_{Z_{1,12\rho^2}} \left(|U|^2+|\nabla_{x,t} U|^2\right)
dxdt\leq \\&\leq Ce^{2 \tau \left(e^{\lambda }-e^{12\lambda \rho^2}\right)}
\left\Vert F\right\Vert^2_{L^2\left(B_1\times (-T,T)\right)}+\\&+
C\rho^{-9}e^{2 \tau \left(e^{10\lambda \rho^2}-e^{12\lambda
\rho^2}\right)}\int^{2\rho}_{-2\rho}\int_{B_1}\left(|U|^2+|\nabla_{x,t}U|^2\right)dxdt,
\end{aligned}
\end{equation}
for every $\tau\geq \tau_1$. By \eqref{oper-2ord-374b}, \eqref{oper-2order-384} e
\eqref{oper-2order-383} we have

\begin{equation*}
\begin{aligned}
&\int^{2\rho}_{-2\rho}\int_{B_1} \left(|U|^2+|\nabla_{x,t}
U|^2\right) dxdt \leq \\&\leq  Ce^{2 \tau \left(e^{\lambda }-e^{12\lambda
\rho^2}\right)} \left\Vert F\right\Vert^2_{L^2\left(B_1\times
(-T,T)\right)}+\\&+ C\rho^{-9}e^{2 \tau \left(e^{10\lambda
\rho^2}-e^{12\lambda
\rho^2}\right)}\int^{2\rho}_{-2\rho}\int_{B_1}\left(|U|^2+|\nabla_{x,t}U|^2\right)dxdt,
\end{aligned}
\end{equation*}
For every $\tau\geq \tau_1$.

Let us choose $\tau=\tau_2\geq\tau_1$, where $\tau_2$  satisfies

$$C\rho^{-9}e^{2 \tau_2 \left(e^{10\lambda
\rho^2}-e^{12\lambda \rho^2}\right)}\leq\frac{1}{2}.$$ 
By this choice of $\tau$ the second term on the right--hand side of
\eqref{oper-2order-384} is absorbed on the left--hand side, and we get
 \eqref{oper-2ord-361}. $\blacksquare$

\bigskip

\textbf{Remark.} Theorem \ref{oper-2ord-stab} can
be proved under less restrictive assumptions on $U$, but here we do not go into this question. Instead, we want to illustrate a simple and direct application to the proof \textbf{of the uniqueness of the following inverse problem} -- which we discuss here only at the formal level -- let $F$ a function depending only on the variable $x$, let $U$ be the solution to the following direct  problem

\begin{equation}\label{oper-2order-385}
\begin{cases}
\partial^2_t U-\Delta U=F(x),\quad
\mbox{in }(x,t)\in B_1\times
(-T,T),\\
\\
U(x,t)=h_0,\quad \mbox{for } (x,t)\in\partial B_1\times
(-T,T),\\
\\
U(x,0)=U_{0}(x),\quad U_t(x,0)=U_1(x), \quad \mbox{for } x\in B_1,
\end{cases}
\end{equation}
where $h_0, U_0, U_1$ and $F$ given functions. The proof
of the uniqueness and the existence of solution to the direct problem \eqref{oper-2order-385}, for a given $F\in L^2(B_1)$, can be found, for instance, in \cite[Ch. 7]{EV}.
Let us assume now to know also 

$$\frac{\partial U}{\partial\nu}_{|\partial B_1\times
(-T,T)}=h_1,$$ 
we wish to determine $F$. Here we give a sketch of the
proof of the uniqueness for the problem of determining $F$
(inverse problem) by mean of $h_0,h_1,U_0,U_1$. Since the problem is 
linear, it is enough to prove that if
$$U_0=U_1=0,\quad h_0=h_1=0,$$ then $F\equiv 0$.

\medskip

Let $U$ satisfy
\begin{equation}\label{oper-2order-386}
\begin{cases}
\partial^2_t U-\Delta U=F(x),\quad
\mbox{in }(x,t)\in B_1\times
(-T,T),\\
\\
U(x,t)=0,\quad \frac{\partial U}{\partial\nu}=0 \quad \mbox{for }
(x,t)\in\partial B_1\times
(-T,T),\\
\\
U(x,0)=0,\quad U_t(x,0)=0, \quad \mbox{for } x\in B_1
\end{cases}
\end{equation}
Set
$$w=\partial_tU$$
and differentiate both the sides of the equazion with respect to $t$. Since
$F$ does not depend on  $t$ we have

\begin{equation}\label{oper-2order-387}
\begin{cases}
\partial^2_t w-\Delta w=0,\quad
\mbox{in }(x,t)\in B_1\times
(-T,T),\\
\\
w(x,t)=0,\quad \frac{\partial w}{\partial\nu}=0 \quad \mbox{for }
(x,t)\in\partial B_1\times (-T,T).
\end{cases}
\end{equation}
Applying to $w$  estimate \eqref{oper-2ord-361} we get

$$w\equiv 0, \quad \mbox{in } B_1\times (-T,T).$$
Hence
$$U(x,t)=U(x,0)+\int^t_0w(x,s)ds=0,\quad \mbox{in } B_1\times (-T,T).$$
Therefore by \eqref{oper-2order-386} we have
$$F(x)=\partial_t^2U-\Delta U=0.$$
This proves the uniqueness for the inverse problem. The idea that
we have outlined is only a miniature of a more general method
of proving uniqueness and stability results for
inverse problems related to evolution equations in which it is required the determination of the time iondepent coefficients of some equation. For further study we refer to  \cite{KL-TI}, \cite{Is2}.
$\blacklozenge$

\section[Geometric meaning of the pseudo --
convexity]{Geometric meaning of the pseudo --
	convexity. Some remarks on the necessary conditions.}\label{cond-nec-Carlm}

In this Section we briefly discuss
some necessary conditions on the weight exponent
for a Carleman estimate to be valid. We state without
proof a H\"{o}rmander Theorem  \cite[Theorem
8.1.1]{HO63}. In some sense, Proposition \ref{stime-Carlm-extn-1-prop-4} is a "miniature" of such a Theorem
\begin{theo}[\textbf{necessary condition for the Carleman estimate}]
	\index{Theorem:@{Theorem:}!- necessary condition for the Carleman estimate@{- necessary condition for the Carleman estimate}}
	
\label{cond-nec-Carlm-1} Let $\Omega$ be a bounded open set of
$\mathbb{R}^n$. Let $\varphi \in
C^{\infty}\left(\bar{\Omega}\right)$ satisfy $\nabla \varphi \neq
0$ in $\bar{\Omega}$. Let
\begin{equation}
\label{cond-nec-Carlm-2} P(x,D)= \sum_{|\alpha| \leq m}
a_{\alpha}(x) D^{\alpha},
\end{equation}
be a linear differential operator of order $m$, such that $a_{\alpha}\in L^{\infty}(\Omega)$, $|\alpha|\leq m $ and $a_{\alpha}\in C^1\left(\bar{\Omega}\right)$, $|\alpha|=m$.

Let us suppose that there exist $K_1>0$ and $\tau_0$ such that for every  $u \in
C^{\infty}_{0}(\Omega, \mathbb{C})$ and for everey $\tau\geq\tau_0$ we have

\begin{equation}
\label{cond-nec-Carlm-3} \tau \sum_{|\alpha| \leq
m-1}\binom{m-1}{\alpha} \int_{\Omega} |D^{\alpha} u |^2 e^{2 \tau
\varphi} dx  \leq K_1 \int_{\Omega} |P(x, D) u |^2 e^{2 \tau
\varphi}dx.
\end{equation}
Then, setting $\zeta=\xi+i \sigma \nabla \varphi(x)$ where $x
\in \Omega$, $\xi \in \mathbb{R}^n$, $\sigma \in \mathbb{R}$,
we have

\begin{equation}
\label{cond-nec-Carlm-4}
\begin{aligned} &\begin{cases}
P_m(x, \zeta) = 0 \\
\\
\sigma \neq 0
\end{cases}  \Longrightarrow \\& \Longrightarrow|\zeta |^2 \leq 2K_1 \left[ \sum_{j,k=1}^{n} \partial_{jk}^2 \varphi(x) P^{(j)}_m(x, \zeta) \overline{P^{(k)}_m(x, \zeta)} + \right.
\\ & \left. + \sigma^{-1} \Im \sum_{k=1}^{n} P_{m,k}(x, \zeta) \overline{P^{(k)}_m (x, \zeta)} \right].
\end{aligned}
\end{equation}
\end{theo}

\medskip

As already noted in \eqref{note-6-50} the expression on the right--hand side
of the implication \eqref{cond-nec-Carlm-4} can be written
by Poisson brackets and it is equal to

\begin{equation*}
\frac{i}{2\sigma}\left\{P_m(x,\xi+i\sigma\nabla\varphi),\overline{P_m(x,\xi+i\sigma\nabla\varphi)}\right\}.
\end{equation*}

\bigskip

\textbf{Remark.} Let us observe that if $P(x,D)$ is an elliptic operator and estimate
 \eqref{cond-nec-Carlm-3} holds true, then Theorem
\ref{cond-nec-Carlm-1} implies that necessary condition \eqref{cond-nec-Carlm-4} holds true in any open set $\widetilde{\Omega}$ compactly contained in $\Omega$
 this, in turn, implies that condition
($\bigstar$) of Theorem \ref{Carlm-ell-teor} is satisfied in
$\overline{\widetilde{\Omega}}$. Therefore such a Theorem
\ref{Carlm-ell-teor} gives the estimate

\begin{equation*}
\sum_{|\alpha|\leq m}\tau^{2(m-|\alpha|)-1}\int \left\vert
D^{\alpha}u\right\vert^2e^{2\tau\varphi}dx\leq C\int \left\vert
P(x,D)u\right\vert^2e^{2\tau\varphi}dx,
\end{equation*}
for every $u\in C_0^{\infty}(\widetilde{\Omega})$ and for every $\tau$
large enough. $\blacklozenge$

\bigskip

In Theorem \ref{cond-nec-Carlm-1} we provided necessary conditions for the validity
of Carleman estimates. To
establish necessary conditions for the unique continuation property we would need deeper the exploration of the notion of pseudo--convexity which 
we mentioned in the previous Section for second-order operators. For
this kind of further investigation we refer to Chapters 4, 5,
6 of the book \cite{Lern}.

\medskip

 \underline{\textbf{Geometric interpretation of pseudo--convexity condition.}}

Let

\begin{equation}\label{cond-nec-Carlm-8}
P_2(x,\xi)=\sum_{j,k=1}^ng^{jk}(x)\xi_j\xi_k.
\end{equation}
where  $\left\{g^{jk}\right\}_{j,k=1}^n$
is a nonsingular symmetric matrix, whose entries are \textbf{real constants}. Let $\phi \in C^{2}(\bar{\Omega})$
where $\Omega$ is a bounded open set of $\mathbb{R}^n$, and
let us suppose that $\nabla \phi \neq 0$ in $\bar{\Omega}$. Let $x_0\in
\Omega$ and set
\begin{equation}\label{cond-nec-Carlm-7}\xi_0=\nabla \phi(x_0).
\end{equation}
Let us consider the Hamiltonian system  associated to $P_2(x,\xi)$ which,
as operator $P_2$ has constant coefficients, is 
\begin{equation}\label{cond-nec-Carlm-9}
\begin{cases}
\overset{\cdot }{x}^j=P^{(j)}_2(\xi(t)), \quad j=1,\cdots, n,\\
\\
\overset{\cdot }{\xi}_j=-P_{2,j}(\xi(t))\mbox{ }(=0), \quad j=1,\cdots, n.\\
\end{cases}
\end{equation}
A solution to \eqref{cond-nec-Carlm-9} on an interval $J$ is
the line $(x(t),\xi(t))$ of $\mathbb{R}^n_x\times
\mathbb{R}^n_{\xi}$ which is called  \textbf{bicharacteristic line}
of the operator $P_2(D)$. Let us note that by \eqref{cond-nec-Carlm-9}
we have

\begin{equation}\label{cond-nec-Carlm-10}
\frac{d}{dt}P_2(\xi(t))=\overset{\cdot
}{\xi}_jP^{(j)}_2(x(t),\xi(t))=0,\quad \forall t\in J.
\end{equation}
If
$$P_2(\xi(t))=0, \quad \forall t\in J,$$
we say that $(x(t),\xi(t))$ is a \textbf{null bicharacteristic line}. In particular, we have that if $(x(t),\xi(t))$ is 
a null bicharacteristic line in a point $t_0$ then it is
null bicharacteristic line in the whole interval $J$. As a matter of fact, if
$$P_2(\xi(t_0))=0,$$
then \eqref{cond-nec-Carlm-10} implies $P_2(\xi(t))=0$, for every
$t\in J$. Let us notice that in each point $(x(t),\xi(t))$ of a
null bicharacteristic line, $\xi(t)$ is a characteristic direction for operator $P_2$ in the point $x(t)$. The projection $x(\cdot)$, on
$\mathbb{R}^n_x$ of a null bicharacteristic line is called
\textbf{ray} of $P_2$ (see Section \ref{NonlinEq}). Now, if $x(t)$ is a ray
of $P_2$, then

\begin{equation}\label{cond-nec-Carlm-11}
\frac{d}{dt}\phi(x(t))=\partial_j\phi(x(t))\frac{dx^j}{dt}=P^{(j)}_2(\xi(t))\partial_j\phi(x(t))
\end{equation}
and
\begin{equation}\label{cond-nec-Carlm-12}
\begin{aligned}
&\frac{d^2}{dt^2}\phi(x(t))=\frac{d}{dt}\left(P^{(j)}_2(\xi(t))\partial_j\phi(x(t))\right)=\\&=
\partial^2_{jk}\phi(x(t))\frac{dx^{k}}{dt}P^{(j)}_2(\xi(t))+P^{(jk)}_2(\xi(t))\frac{d\xi_k}{dt}\partial_j\phi(x(t))=\\&=
\partial^2_{jk}\phi(x(t))P^{(j)}_2(\xi(t))P^{(k)}_2(\xi(t)).
\end{aligned}
\end{equation}
Let us suppose that in $\xi_0\in \mathbb{R}^n\setminus\{0\}$ we have
\begin{equation}\label{cond-nec-Carlm-12n}
P_2(\xi_0)=0,\end{equation}
 let $x_0\in \Omega$ and let us consider the
solution $(x(t),\xi(t))$ of system \eqref{cond-nec-Carlm-9} satisfying the initial condition in $t_0\in J$
\begin{equation}\label{cond-nec-Carlm-12nn}
x(t_0)=x_0,\quad \xi(t_0)=\xi_0.
\end{equation}
 We quickly notice that
\eqref{cond-nec-Carlm-11} allows us to write  condition
\eqref{oper-2ord-343} of Theorem \ref{oper-2ord-341} (see also
Proposition \ref{oper-2ord-346}) assuming there
\eqref{cond-nec-Carlm-12n} and \eqref{cond-nec-Carlm-12nn} as follows

\begin{equation}\label{cond-nec-Carlm-14}
\frac{d}{dt}\phi(x(t))_{|t=t_0}\neq 0.
\end{equation}
This is equivalent to the fact that the ray $x(t)$ passing through $x_0$ is
transverse to the level surface
$$\Gamma_{\phi}=\left\{x\in
\Omega :\quad \phi(x)=\phi(x_0)\right\}.$$

Furthermore, \eqref{cond-nec-Carlm-11} and \eqref{cond-nec-Carlm-12}
allow us to write the \eqref{oper-2ord-106} (in $x_0$) in the
form

\begin{equation}\label{cond-nec-Carlm-13}
\begin{aligned}
\frac{d}{dt}\phi(x(t))_{|t=t_0}= 0\mbox{ } \Longrightarrow \mbox{ }&
\frac{d^2}{dt^2}\phi(x(t))_{|t=t_0}>0,
\end{aligned}
\end{equation}
this implies that if $t_0$ is a critical point of $\phi(x(t))$
then it is a proper minimum point of the function $\phi(x(t))$. Condition
\eqref{cond-nec-Carlm-13} can also be formulated in the following
manner: let us consider the above defined level surface
$\Gamma_{\phi}$ and the level set
$$\Omega^+_{\phi}=\left\{x\in
\Omega :\quad \phi(x)>\phi(x_0)\right\},$$ then \eqref{cond-nec-Carlm-13}
says that if the ray $x(t)$ is tangent in $x_0$ to the level
surface $\Gamma_{\phi}$ (this is expressed
by the antecedent of implication \eqref{cond-nec-Carlm-13}), then
for every $t$ in a neighborhood of $t_0$, we have, for $t\neq t_0$, $x(t)\in
\Omega^+_{\phi}$. In other words, the ray $x(t)$ cannot
"cross" the level surface $\Gamma_{\phi}$ at points where $x(t)$ is tangent to $\Gamma_{\phi}$.

\chapter{Optimal three sphere and doubling inequality for second order elliptic equations}\markboth{Chapter 15. Optimal three sphere inequality and doubling inequality}{}\label{tre sfere-ellittiche}
\section{Introduction}\label{tre sfere-ellittiche-intro}

In this Chapter we will prove the \textbf{strong unique  continuation property} \index{strong unique  continuation property} for the second order elliptic equations  with real coefficients (in the principal part). Now we  recall
briefly this property and provide an introduction to the
Chapter.

Let $\left\{a^{ij}(x)\right\}^n_{i,j=1}$ be a \textbf{symmetric matrix of
real--valued functions}. We assume that 
the following uniform ellipticity condition is satisfied 

\begin{equation*}
\lambda^{-1}\left|\xi\right|^{2}\leq
\sum_{i,j=1}^na^{ij}(x)\xi_i\xi_j\leq\lambda\left|\xi\right|^{2},\quad
\forall\xi\in\mathbb{R}^{n} \mbox{, }\forall x\in B_1,
\end{equation*}
where $\lambda\geq 1$. Let us assume that the function $a^{ij}$, $i,j=1,\cdots,n$ are Lipschitz continuous

\begin{equation}\label{intro14-4An}
\left|a^{ij}(x)-a^{ij}(y)\right|\leq\Lambda |x-y|, \quad\mbox{for
}i,j\in\left\{1,\cdots, n\right\}, \quad\forall x,y\in B_1.
\end{equation}

\medskip

Let $b^i\in L^{\infty}\left(B_1\right)$, $i=1,\cdots, n$ and
$c\in L^{\infty}\left(B_1\right)$ (these coefficients can
also be complex--valued) satisfy

\begin{equation*}
\left\Vert b^i\right\Vert_{L^{\infty}\left(B_1\right)}\leq
M,\quad\mbox{for } i=1,\cdots,n
\end{equation*}
and
\begin{equation}\label{intro14-4Ann}
\left\Vert c\right\Vert_{L^{\infty}\left(B_1\right)}\leq M.
\end{equation}

We recall that the equation

\begin{equation}  \label{intro14-equ}
Lu=\sum_{i,j=1}^na^{ij}(x)\partial^2_{x^ix^j}u+\sum_{i=1}^nb^i(x)\partial_{x^i}u+c(x)u=0,\quad\mbox{in
} B_1,
\end{equation}
enjoys the \textbf{strong unique continuation property} provided that any solution $u$ to \eqref{intro14-equ} satisfying the conditions 

\begin{equation}  \label{intro14-equ-1}
\int_{B_r}|u|^2dx=\mathcal{O}\left(r^m\right),\quad\mbox{as }
r\rightarrow 0,\quad \forall m\in \mathbb{N},
\end{equation}
identically vanishes.
We will prove the strong unique continuation property as a
consequence of an \textbf{optimal three sphere inequality}. A
prototype of such an inequality is the Hadamard three circle inequality for the holomorphic functions that we first encountered in
Section  \ref{esempi-stime-errore}.

Generally speaking, a three sphere inequality for solutions to the equation
\eqref{intro14-equ} is an inequality of the type

\begin{equation}  \label{intro14-equ-2}
\int_{B_{\rho}}|u|^2dx\leq
C\left(\int_{B_{R}}|u|^2dx\right)^{1-\theta}\left(\int_{B_{r}}|u|^2dx\right)^{\theta},
\end{equation}
where $0<r<\rho<R\leq 1$, $C$ and $\theta\in (0,1)$ depend by
$\lambda$, $\Lambda$, $M$ and $R$, $\rho$ ($C$ and $\theta$ do not
depend on $u$).

We say that \eqref{intro14-equ-2} is an \textbf{optimal} three spheres inequality provided $C$ \textbf{does not depend} on $r$ and, for fixed
$R,\rho$, we have

\begin{equation}  \label{intro14-equ-3}
\theta\sim\left|\log r\right|^{-1} ,\quad\mbox{as } r\rightarrow 0.
\end{equation}
We recall that by $f(r)\sim g(r)$, as $r\rightarrow
0$ we means

$$0<\lim_{r\rightarrow 0}\frac{f(r)}{g(r)}<+\infty.$$

First author that proved \eqref{intro14-equ-2} was
Landis in \cite{LAN}.

Arguing like in Remark 5 of Section \ref{esempi-stime-errore}, it is easy to check
that if a function $u$
satisfies an optimal three sphere inequality, then whenever $u$ satisfies  \eqref{intro14-equ-1} it
vanishes identically. Therefore, if an optimal three sphere inequality holds true for equation
\eqref{intro14-equ-1} then such an equation satisfies the strong unique continuation property.

Another type of inequality that implies the
strong unique continuation property is the so-called \textbf{doubling inequality}\index{doubling inequality}. Such an inequality occurs in the form

\begin{equation}  \label{intro14-equ-4}
\int_{B_{2r}}|u|^2dx\leq K\int_{B_{r}}|u|^2dx, \quad\forall r\in \left(0,\frac{1}{2}\right),
\end{equation}
where  $K$ depends on $u$ but \textbf{does not depend on} $r$. We will prove later on in which a way \eqref{intro14-equ-4} implies the strong unique continuation property. The main idea may be expressed as
follows. Iterating inequality \eqref{intro14-equ-4} we have, for
every $j\in \mathbb{N}$

\begin{equation*}
\int_{B_{1/2}}|u|^2dx\leq K\int_{B_{1/4}}|u|^2dx\leq \cdots
\leq K^{j-1}\int_{B_{1/2^j}}|u|^2dx
\end{equation*}
which, together with \eqref{intro14-equ-1}, provides, for each
 $j,m\in
\mathbb{N}$ ($C_m$ depends on $m$ only),

\begin{equation}\label{intro14-equ-5}
\int_{B_{1/2}}|u|^2dx\leq
K^{j-1}C_m\left(\frac{1}{2^j}\right)^m=C_mK^{-1}\left(\frac{K}{2^m}\right)^j.
\end{equation}
Let now $m$ satisfy $$2^m>K$$ and passing to the limit as
$j\rightarrow \infty$, we get by  \eqref{intro14-equ-5} 

\begin{equation*}
\int_{B_{1/2}(0)}|u|^2dx=0.
\end{equation*}

Both the optimal three sphere inequality and the doubling inequality
will be obtained from an appropriate Carleman estimate for the elliptic operator $L$ (or, equivalently, for the principal part of that operator). We will approach this question in
two phases: first we will study the case of the Laplace operator $\Delta$ and
then we will study the case of equations with variable coefficients. In both the cases, the
proofs of the Carleman estimates will start from
rewriting the elliptic operators in polar coordinates:
in the case of Laplace operator, in Euclidean polar coordinates; in the case of 
variable coefficients, in geodesic polar coordinates w.r.t. the
Riemannian structure induced by a metric conforming to

$$a_{ij}dx^{i}\otimes dx^{j},$$ ($\left\{a_{ij}(x)\right\}^n_{i,j=1}$
is the inverse of the matrix $\left\{a^{ij}(x)\right\}^n_{i,j=1}$).
We warn that throughout this Chapter we will use the Einstein convention
of repeated indices. We will, in addition, adhere more
scrupulously to the notation on indices of the components of a
tensor. Actually, these notations are mostly needed in the
Sections \ref{geodesic-polar} and \ref{aronszajn-general}, but
we will adopt it in the preceding sections as well. The proof in the
case of the Laplace operator presents most of the main
difficulties that we will encounter in the case of variable coefficient,
which, of course, presents additional technical difficulties.
In addition to the proofs that we give here,
There exist other proofs in the literature, e.g. \cite{E-V}, \cite{HO}.

\section{Formulas for the change of variables of second order operators}\label{cambiamento-variabili} We begin by deriving
a formula to the change of variables of the operator

\begin{equation}\label{ias-0-1}
\mbox{div}\left(A(x)\nabla u(x)\right),
\end{equation}
where $A(x)=\left\{a^{ij}(x)\right\}^n_{i,j=1}$ is a symmetric matrix, $a^{ij}$ are sufficiently smooth functions for $i,j=1,\cdots, n$ (will suffice
$a^{ij}\in C^{0,1}\left(\mathbb{R}^n\right)$).

Let us consider the case where $A(x)$ is the identity matrix. Let
$\Lambda$ be an open set of $\mathbb{R}^n$ and let $\Phi\in
C^1\left(\overline{\Lambda}, \mathbb{R}^n\right)$ be a injective
map such that

\begin{equation}\label{ias-1n-1}
\det(J\Phi(x))\neq 0,\quad \forall x\in \overline{\Lambda},
\end{equation}
where $J\Phi(x)$ is the jacobian matrix of $\Phi$ in $x$
 $$J\Phi(x)=\left(
\begin{array}{ccc}
\partial_{x^1}\Phi^1((x) & \cdots  & \partial_{x^n}\Phi^1((x) \\
\vdots & \ddots  & \vdots \\
\partial_{x^1}\Phi^n((x) & \cdots  & \partial_{x^n}\Phi^n((x)%
\end{array}
\right). $$ Set $$\Omega=\Phi(\Lambda)$$ and
$$\Psi=\Phi^{-1},$$ consequently $\Psi \in C^1\left(\overline{\Omega}, \mathbb{R}^n\right)$. Let us prove
that, if $u\in C^2\left(\overline{\Lambda}\right)$, then

\begin{equation}\label{ias-3-1}
(\Delta u)(\Psi(y))=\frac{1}{\left\vert
\det(J\Psi(y))\right\vert}\mbox{div}\left(B(y)v(y)\right),
\end{equation}
where
\begin{equation}\label{ias-4-1}
v(y)=u(\Psi(y))
\end{equation}
and
\begin{equation}\label{ias-5-1}
B(y)=\left|\det(J\Psi(y))\right|(J\Psi(y))^{-1}\left((J\Psi(y))^{-1}\right)^T.
\end{equation}

More  generally, for operator \eqref{ias-0-1} we have

\begin{equation}\label{ias-1-3}
\mbox{div}\left(A(x)\nabla
u(x)\right)_{|x=\Psi(y)}=\frac{1}{\left\vert
\det(J\Psi(y))\right\vert}\mbox{div}\left(\widetilde{A}(y)v(y)\right),
\end{equation}
where $v$ is given by \eqref{ias-4-1} and

\begin{equation}\label{ias-1-3n}
\widetilde{A}(y)=\left|\det(J\Psi(y))\right|(J\Psi(y))^{-1}A(\Psi(y))\left((J\Psi(y))^{-1}\right)^T.
\end{equation}
Since \eqref{ias-1-3} can be proved similarly, we limit ourselves to prove \eqref{ias-3-1}.

\bigskip

\textbf{Proof of \eqref{ias-3-1}.} Let $w\in
C^{\infty}_0\left(\Omega\right)$ and let

$$\widetilde{w}(x)=w(\Phi(x)), \quad\forall x\in \Lambda.$$

Integrating by parts, we have

\begin{equation}\label{ias-1-3nn}-\int_{\Lambda}\Delta u(x)\widetilde{w}(x)dx=\int_{\Lambda}\nabla
u(x)\cdot\nabla\widetilde{w}(x)dx\end{equation} and by the formula of
change of variables for multiple integrals, we have

\begin{equation}\label{ias-0-2}
	\begin{aligned}
		&\int_{\Lambda}\nabla_x
u(x)\cdot \nabla_x\widetilde{w}(x)dx=\\&=\int_{\Omega}\left(\nabla_x
u\right)(\Psi(y))\cdot\left(\nabla_x\widetilde{w}\right)(\Psi(y))\left|\det(J\Psi(y))\right|
dy.
\end{aligned}
\end{equation} Now, \eqref{ias-4-1} gives

$$\nabla_yv(y)=(J\Psi(y))^T\left(\nabla_x
u\right)(\Psi(y)),$$ ($\nabla_x u$ and $\nabla_yv$ are 
column vectors). Hence

\begin{equation}\label{ias-1-2}
\left(\nabla_x
u\right)(\Psi(y))=\left((J\Psi(y))^{-1}\right)^T\nabla_yv(y),
\end{equation}
similarly

\begin{equation}\label{ias-2-2}
\left(\nabla_x
\widetilde{w}\right)(\Psi(y))=\left((J\Psi(y))^{-1}\right)^T\nabla_y
w(y).
\end{equation}
Substituting \eqref{ias-1-2} and \eqref{ias-2-2} into the integral on the
right--hand side of \eqref{ias-0-2} and recalling \eqref{ias-1-3nn}, we have

\begin{equation}\label{ias-3-2}
 \begin{aligned}
-\int_{\Lambda}\Delta
u(x)\widetilde{w}(x)dx&=\int_{\Omega}B(y)\nabla_y v(y)\cdot\nabla_y
w(y)dy=\\&= -\int_{\Omega}\mbox{div}_y \left(B(y)\nabla_y
v(y)\right) w(y)dy.
\end{aligned}
\end{equation}
Now, we change the variables in the integral on the
left--hand side of \eqref{ias-3-2} and we find

\begin{equation*}
-\int_{\Omega}(\Delta u)(\Psi(y))w(y)\left|\det(J\Psi(y))\right|dy=
-\int_{\Omega}\mbox{div}_y \left(B(y)\nabla_y v(y)\right) w(y)dy.
\end{equation*}
Since $w$ is arbitrary in
$C^{\infty}_0\left(\Omega\right)$, we get \eqref{ias-3-1}.

\bigskip

\underline{\textbf{Exercise 1.}} Apply formula \eqref{ias-3-1} to write
the Laplace operator, in $\mathbb{R}^2$, in polar coordinates showing that, setting
$$(x^1,x^2)=\Psi(\varrho,\vartheta)=(\varrho \cos\vartheta, \varrho\sin \vartheta), \quad
(\varrho,\vartheta)\in (0,+\infty)\times[0,2\pi)$$ and
$$v(\varrho,\vartheta)=u(\Psi(\varrho,\vartheta))=u(\varrho \cos\vartheta, \varrho\sin
\vartheta),$$ we have

\begin{equation}\label{ias-0-6}
(\Delta u)(\varrho \cos\vartheta, \varrho\sin
\vartheta)=\partial^2_{\varrho}v+\frac{1}{\varrho}\partial_{\varrho}v+\frac{1}{\varrho^2}\partial^2_{\vartheta}v.
\end{equation} $\clubsuit$

\bigskip

\underline{\textbf{Exercise 2.}} Apply  formula \eqref{ias-3-1} to
write the Laplace operator, in $\mathbb{R}^3$, in polar coordinates 
showing that, setting
\begin{equation}\label{ias-1-6}(x^1,x^2,x^3)=\Psi(\varrho,\vartheta,\phi)=(\varrho
\sin\vartheta\cos\phi, \varrho\sin\vartheta\sin\phi,
\varrho\cos\vartheta), 
 \end{equation} where
for $(\varrho,\vartheta,\phi)\in (0,+\infty)\times) (0,\pi)\times (0,2\pi)$ and
$$v(\varrho,\vartheta,\phi)=u(\Psi(\varrho,\vartheta\phi))=u(\varrho \sin\vartheta\cos\phi, \varrho\sin\vartheta\sin\phi, \varrho\cos\vartheta),$$ we have
\begin{equation}\label{ias-2-6}
	\begin{aligned}
	&(\Delta u)(\varrho \sin\vartheta\cos\phi,
	\varrho\sin\vartheta\sin\phi, \varrho\cos\vartheta)=\\&=
	\partial^2_{\varrho}v+\frac{2}{\varrho}\partial_{\varrho}v+\frac{1}{\varrho^2\sin
		\vartheta}\partial_{\vartheta}\left(\sin \vartheta
	\partial_{\vartheta}v\right)+\frac{1}{\varrho^2\sin^2
		\vartheta}\partial^2_{\phi}v.
	\end{aligned}
\end{equation}
 $\clubsuit$

\section[Polar coordinates in $\mathbb{R}^n$]{Polar coordinates in $\mathbb{R}^n$. The Laplace--Beltrami operator on the sphere}\label{n-euclid-coordinate polari} In
the $n$--dimensional case, it is possible to express the Laplace operator in polar coordinates  by means some formulas similar to \eqref{ias-0-6} and to \eqref{ias-2-6}.
Neverthless, in this way we come to write quite
cumbersome formulas . For this reason we present here a procedure,  which consists, first of all, in defining
the Laplace operator on the sphere
$$\Sigma=\left\{x\in \mathbb{R}^n| \mbox{ } |x|=1\right\}.$$
Concerning the polar coordinates in $\mathbb{R}^n$ if one wishes to
follow the path suggested in Exercise 1 and 2 of the previous Section, it would be convenient to write a transformation similar to \eqref{ias-1-6}. For this purpose it would suffice to use the following
recursive formula

\begin{equation}\label{ias-1-7}
\begin{cases}
r_n=|x|\\
\\
x^n=r_n\cos \vartheta^n, \quad \vartheta^n\in [0,2\pi),\\
\\
r_{n-1}=r_n\sin \vartheta^n, \quad \vartheta^n\in [0,2\pi),\\
\\
x^{n-1}=r_{n-1}\cos \vartheta^{n-1}, \quad \vartheta^{n-1}\in [0,\pi)\
\\
\cdots\\
\\
x^1=r_1 \cos \vartheta^{2}, \quad \vartheta^{2}\in [0,\pi)
.\end{cases}
\end{equation}
In this way, the "angular coordinates" are
$\left(\vartheta^{2},\vartheta^{3},\cdots,\vartheta^{n}\right)$. 
By using \eqref{ias-1-7} one reaches
to the following formula of change of variables in polar coordinates ($f$ integrable function)\index{integral in Euclidean polar coordinates}

\begin{equation}\label{ias-2-7}
	\begin{aligned}
		&\int_{B_R}f(x)dx=\int^R_0\left(\int_{\Sigma}f(\rho\omega)\rho^{n-1}dS\right)d\rho=\\&=\int^R_0\left(\int_{\partial
B_{\rho}}f(y)dS_y\right)d\rho,
\end{aligned}
\end{equation}
where $R>0$.

In order to write the Laplace operator in polar coordinates we proceed in two steps:

\medskip

\noindent \textbf{Step I.} We will write the Laplace--Beltrami operator on the sphere  $\Sigma$;

\medskip

\noindent \textbf{Step II.} We will complete the transformation begun in Step
I and we will prove Theorem \ref{ias-1p-18}.

\bigskip

\noindent \underline{\textbf{Step I.}} Let  $u(\omega)$ be a real--valued function on
$\Sigma$. We associate to $u$ the homogeneous function of degree $0$ defined by
 $$\widetilde{u}:\mathbb{R}^n\setminus\{0\}\rightarrow
\mathbb{R}.$$ 

\begin{equation}\label{ias-1-8}
\widetilde{u}(x)=u\left(\frac{x}{|x|}\right), \quad x\in
\mathbb{R}^n\setminus\{0\}.
\end{equation}

If $\widetilde{u}\in C^k\left(\mathbb{R}^n\setminus\{0\}\right)$, we say  that $u$ belongs to $C^k(\Sigma)$. We define the scalar product of  two functions $u,v\in C^0(\Sigma)$ as

\begin{equation}\label{ias-2-8}
(u,v)_{\Sigma}=\int_{\Sigma}u(\omega)v(\omega) d\omega
\end{equation}
($d\omega:=dS$) and we set

\begin{equation}\label{ias-3-8}
\left\Vert u\right\Vert_{L^2(\Sigma)}=\int_{\Sigma}\left\vert
u(\omega)\right\vert^2 d\omega;
\end{equation}
Hence $C^0(\Sigma)$, equipped with the scalar scalar product \eqref{ias-2-8}, is a prehilbertian space and $L^2(\Sigma)$
is the completion of that space.

\bigskip

For any $1\leq i\leq n$ we define the operator $d_i$ as follows.
Let $u\in C^1(\Sigma)$, let us denote

\begin{equation}\label{ias-4-8}
\left(d_i u\right)(\omega)=\left(\partial_{x^i}
\widetilde{u}(x)\right)_{|\mbox{ }x=\omega}, \quad \omega\in \Sigma.
\end{equation}
Other symbols which are used in literature to denote  $d_iu$ are: $\Omega_iu$,
$\partial_{\omega^i}u$.

\bigskip

\textbf{Remarks.}

\noindent\textbf{1.} By \eqref{ias-4-8} we have for $ j,k=1,\cdots, n$

\begin{equation}\label{ias-0-9}
	\begin{aligned}
	&\left(d_i
	u\right)(\omega)=\left(\left(\partial_{x^i}\frac{x^k}{|x|}\right)\left(\partial_{x^k}\widetilde{u}\right)\left(\frac{x}{|x|}\right)\right)_{|\mbox{ }x=\omega}=\\&
	=\frac{1}{|x|}\left(\left(\partial_{x^i}\widetilde{u}\right)\left(\frac{x}{|x|}\right)-
	\frac{x^i}{|x|}\frac{x^k}{|x|}\left(\partial_{x^k}\widetilde{u}\right)\left(\frac{x}{|x|}\right)\right)_{|\mbox{ }x=\omega}.
	\end{aligned}
\end{equation}
Let us note that the function written in the brackets in the second line
is the $i$--th component of

\begin{equation}\label{ias-00-9}
\left(\nabla
\widetilde{u}\right)\left(\frac{x}{|x|}\right)-\frac{x}{|x|}\left(\frac{x}{|x|}\cdot
\left(\nabla \widetilde{u}\right)\left(\frac{x}{|x|}\right)\right),
\end{equation}
which, in turn, is the $i$--th  tangential component on
$\{|x|=1\}$ of the gradient of $\widetilde{u}$.

\medskip

\noindent\textbf{2.} The operators $d_i$, $1\leq i\leq n$,
\textbf{are not independent}. As a matter of fact, by
\eqref{ias-4-8} and by the fact that $\widetilde{u}$ is a
homogeneous function of degree $0$, we have
\begin{equation}\label{ias-1-9}
\sum_{i=1}^n\omega^i\left(d_iu\right)(\omega)=\sum_{i=1}^n
\left(\frac{x^i}{|x|}\left(\partial_{x^i}\widetilde{u}\right)\left(\frac{x}{|x|}\right)\right))_{|\mbox{
}x=\omega}=0.
\end{equation}

Moreover, by \eqref{ias-4-8} we have, for every $u,v \in
C^1(\Sigma)$,
\begin{equation}\label{ias-2-9}
d_i(uv)=\left(d_iu\right)v+u\left(d_iv\right), \quad 1\leq i\leq n.
\end{equation}

Also let us note

\begin{equation}\label{ias-1-10}
\partial_{x^i}\widetilde{u}(x)=\frac{1}{|x|}\left(d_iu\right)\left(\frac{x}{|x|}\right) , \quad x\in
\mathbb{R}^n\setminus\{0\}\mbox{, } 1\leq i\leq n.
\end{equation}
As a matter of fact, by homogeneity  we have

$$\widetilde{u}(x)=\widetilde{u}(\lambda x)\quad \forall \lambda\in \mathbb{R}\setminus\{0\}.$$
Hence

$$\partial_{x^i}\widetilde{u}(x)=\lambda\left(\partial_{x^i}\widetilde{u}\right)(\lambda
x),\quad 1\leq i\leq n$$ from which, choosing $\lambda=\frac{1}{|x|}$, we get

$$\partial_{x^i}\widetilde{u}(x)=\frac{1}{|x|}\left(\partial_{x^i}\widetilde{u}\right)\left(\frac{x}{|x|}\right)=\frac{1}{|x|}\left(d_{i}u\right)\left(\frac{x}{|x|}\right)
 ,\quad 1\leq i\leq n,$$
from which  \eqref{ias-1-10} follows. $\blacklozenge$

\bigskip

\begin{prop}\label{ias-prop1-10}
Let $n>1$. We have

\begin{equation}\label{ias-2-10}
\left(d_iu,v\right)_{\Sigma}=\left(u,d^{\ast}_iv\right)_{\Sigma},\quad
\forall u,v \in C^1(\Sigma) \mbox{, } 1\leq i\leq n,
\end{equation}
where

\begin{equation}\label{ias-3-10}
 d^{\ast}_i=(n-1)\omega^i-d_i, \quad 1\leq i\leq n.
\end{equation}
$d^{\ast}_i$ is the formal adjoint of $d_i$ with respect to the scalar product
 \eqref{ias-2-8}.
\end{prop}

\medskip

\textbf{Proof.} First, let us notice that, for any $1\leq i\leq n$,
$\partial_{x_i}\widetilde{u}(x)$ is  integrable over $B_1$ because

$$\left|\partial_{x_i}\widetilde{u}(x)\right|\leq
\frac{1}{|x|}\max_{\Sigma}\left|\nabla\widetilde{u}\right|, \quad
1\leq i\leq n$$ and $n>1$.

Now, let us check that, for any $1\leq i\leq n$, we have

\begin{equation}\label{ias-1-11}
 \int_{\Sigma}d_iu(\omega)v(\omega)d\omega=(n-1)\int_{B_1}\partial_{x^i}\widetilde{u}(x) \widetilde{v}(x)dx.
\end{equation}
We check \eqref{ias-1-11}. If $\varepsilon\in (0,1)$ and $1\leq i\leq n$, we have by 
\eqref{ias-2-7} and by \eqref{ias-1-10},
\begin{equation*}
 \begin{aligned}
\int_{B_1\setminus
B_{\varepsilon}}\partial_{x^i}\widetilde{u}(x)
\widetilde{v}(x)dx&=\int_{B_1\setminus
B_{\varepsilon}}\frac{1}{|x|}\left(d_iu\right)\left(\frac{x}{|x|}\right)v\left(\frac{x}{|x|}\right)dx=\\&=
\int_{\varepsilon}^1\rho^{n-2}\int_{\Sigma}d_iu(\omega)v(\omega)d\omega=\\&=
\frac{1-\varepsilon^{n-1}}{n-1}\int_{\Sigma}d_iu(\omega)v(\omega)d\omega.
\end{aligned}
\end{equation*}
Passing to the limit as $\varepsilon\rightarrow 0$ we obtain
\eqref{ias-1-11}.

Now, by \eqref{ias-1-11} and integrating by parts we have

\begin{equation}\label{ias-1n-11}
 \begin{aligned}
&\frac{1}{n-1}\left(d_iu,v\right)_{\Sigma}=\int_{B_1}\partial_{x^i}\widetilde{u}(x)
\widetilde{v}(x)dx=\\&=\lim_{\varepsilon\rightarrow
0}\int_{B_1\setminus
B_{\varepsilon}}\partial_{x^i}\widetilde{u}(x)
\widetilde{v}(x)dx=\\&=\lim_{\varepsilon\rightarrow
0}\int_{B_1\setminus
B_{\varepsilon}}\left(\partial_{x^i}\left(\widetilde{u}\widetilde{v}\right)-\widetilde{u}\partial_{x^i}\widetilde{v}\right)dx=\\&=
\int_{\Sigma}u(\omega)v(\omega)\omega^id\omega-\\&-\lim_{\varepsilon\rightarrow
0}\left(\int_{\partial
B_{\varepsilon}}\widetilde{u}(x)\widetilde{v}(x)\frac{x^i}{|x|}dS+\int_{B_1\setminus
B_{\varepsilon}}\widetilde{u}(x)\partial_{x^i}\widetilde{v}(x)dx\right).
\end{aligned}
\end{equation}
On the other hand

$$\int_{\partial
B_{\varepsilon}}\widetilde{u}(x)\widetilde{v}(x)\frac{x^i}{|x|}dS=\mathcal{O}\left(\varepsilon^{n-1}\right),
\quad\mbox{as } \varepsilon\rightarrow 0,$$ hence we have

$$\lim_{\varepsilon\rightarrow
0}\int_{\partial
B_{\varepsilon}}\widetilde{u}(x)\widetilde{v}(x)\frac{x^i}{|x|}dS=0.$$

Therefore \eqref{ias-1n-11} gives

$$\frac{1}{n-1}\left(d_iu,v\right)_{\Sigma}= \int_{\Sigma}u(\omega)v(\omega)\omega^id\omega-\int_{B_1}\widetilde{u}(x)\partial_{x^i}\widetilde{v}(x)dx$$
and employing \eqref{ias-1-11} (interchange $u$ with $v$ in
the latter) we get

\begin{equation*}
 \begin{aligned}
\frac{1}{n-1}\left(d_iu,v\right)_{\Sigma}&
=\int_{\Sigma}u(\omega)\left(v(\omega)\omega^i-\frac{1}{n-1}d_iv(\omega)\right)d\omega=\\&=
\frac{1}{n-1}\left(u,(n-1)\omega^iv-d_iv\right)_{\Sigma}.
\end{aligned}
\end{equation*}
By the just obtained equality, \eqref{ias-2-10} immediately follows .
$\blacksquare$

\bigskip

\begin{definition}[\textbf{Laplace -- Beltrami operator on the sphere}]\label{ias-def-12}
	\index{Definition:@{Definition:}!- Laplace -- Beltrami operator on the sphere@{- Laplace -- Beltrami operator on the sphere}} The Laplace--Beltrami operator on the sphere is defined
as
\begin{equation}\label{ias-lap-bel-12}
 \Delta_{\Sigma}=\sum_{i=1}^nd^2_i.
\end{equation}
\end{definition}

\bigskip

\begin{prop}\label{ias-prop2-13}
We have

\begin{equation}\label{ias-1-13}
\left(\Delta_{\Sigma}u,v\right)_{\Sigma}=-\sum_{i=1}^n\left(d_iu,d_iv\right)_{\Sigma},\quad
\forall u,v \in C^2(\Sigma).
\end{equation}
\end{prop}

\bigskip

\textbf{Proof.} By \eqref{ias-1-9} and \eqref{ias-2-10} we get
\begin{equation*}
 \begin{aligned}
\left(\Delta_{\Sigma}u,v\right)_{\Sigma}&=\sum_{i=1}^n\left(d^2_iu,v\right)_{\Sigma}=\sum_{i=1}^n\left(d_iu,d_i^{\ast}v\right)_{\Sigma}=\\&
=\sum_{i=1}^n\left(d_iu,(n-1)\omega^iv-d_iv\right)_{\Sigma}=\\&=(n-1)\left(\sum_{i=1}^n\omega^id_iu,v\right)_{\Sigma}-\sum_{i=1}^n\left(d_iu,d_iv\right)_{\Sigma}
=\\&=-\sum_{i=1}^n\left(d_iu,d_iv\right)_{\Sigma}.
\end{aligned}
\end{equation*}$\blacksquare$

\bigskip

\textbf{Remark.}  Proposition \eqref{ias-prop2-13}
implies

\begin{equation}\label{ias-2-13}
\left(\Delta_{\Sigma}u,u\right)_{\Sigma}=-\sum_{i=1}^n\left\Vert
d_iu\right\Vert^2_{L^2(\Sigma)}\leq 0,\quad \forall u \in
C^2(\Sigma)
\end{equation}
and
\begin{equation}\label{ias-3-13}
\int_{\Sigma}\Delta_{\Sigma}u(\omega)d\omega=\left(\Delta_{\Sigma}u,1\right)_{\Sigma}=0.
\end{equation}
$\blacklozenge$

\bigskip

\begin{prop}\label{ias-es-14}
Si ha, per $i,j=1, \cdots, n$,
\begin{equation}\label{ias-1-14}
d_jd_iu-d_id_ju=\omega_jd_iu-\omega_id_ju,\quad \forall u \in
C^2(\Sigma).
\end{equation}
\end{prop}

\bigskip

\textbf{Proof.} By \eqref{ias-1-10} we have

\begin{equation*}
\left(d_iu\right)\left(\frac{x}{|x|}\right)=|x|\partial_{x^i}\widetilde{u}(x),
\quad 1\leq i\leq
n, \mbox{ } \forall x\in \mathbb{R}^n\setminus\{0\}.
\end{equation*}

Hence

\begin{equation}\label{ias-2-14}
 \begin{aligned}
d_jd_iu(\omega)&=\partial_{x^j}\left(|x|\partial_{x^i}\widetilde{u}(x)\right)_{|\mbox{
}x=\omega}=\\&
=\left(\left(\partial_{x^j}|x|\right)\partial_{x^i}\widetilde{u}(x)+|x|\partial^2_{x^ix^j}\widetilde{u}(x)\right)_{|\mbox{
}x=\omega}=\\&=\left(\frac{x^j}{|x|}\partial_{x^i}\widetilde{u}(x)\right)_{|\mbox{
}x=\omega}+\left(\partial^2_{x^ix^j}\widetilde{u}\right)(\omega)=\\&=
\omega^jd_iu(\omega)+\left(\partial^2_{x^ix^j}\widetilde{u}\right)(\omega)
\end{aligned}
\end{equation}
and, similarly,

\begin{equation}\label{ias-3-14}
d_id_ju(\omega)=\omega^id_ju(\omega)+\left(\partial^2_{x^jx^i}\widetilde{u}\right)(\omega).
\end{equation}
Subtracting \eqref{ias-3-14} by \eqref{ias-2-14} we obtain
\eqref{ias-1-14}.$\blacksquare$

\bigskip

\bigskip

\noindent \textbf{Step II.} We associate to $u\in
C^2\left(B_R\setminus\{0\}\right)$, where $R>0$, the function $U\in
C^2\left((0,R)\times \left(\mathbb{R}^n\setminus\{0\}\right)\right)$

\begin{equation}\label{ias-0-15}
U(\rho,y)=u\left(\rho\frac{y}{|y|}\right), \quad \rho\in(0,R),
\mbox{ } y\in \mathbb{R}^n\setminus\{0\}.
\end{equation}
It is obvious that if $\rho=|x|$ and $\omega=\frac{x}{|x|}$, then we have

\begin{equation}\label{ias-00-15}
U(\rho,\omega)=u\left(|x|\frac{x}{|x|}\right)=u(x).
\end{equation}

\begin{prop}\label{ias-1prop-16}
Let $U$ be defined by \eqref{ias-0-15}. Then we have, for any $i=1,
\cdots, n$,
\begin{equation}\label{ias-1-16}
\partial_{x^i}u(x)=\omega^i\partial_{\rho}U(\rho,\omega)+\frac{1}{\rho}d_i
U(\rho,\omega).
\end{equation}
where $\rho=|x|$ and $\omega=\frac{x}{|x|}$
\end{prop}

\bigskip

\textbf{Proof.} By \eqref{ias-00-15} we have

\begin{equation}\label{ias-1n-16}
 \begin{aligned}
&\partial_{x^i}u(x)=\left(\partial_{x^i}|x|\right)\partial_{\rho}U\left(|x|,\frac{x}{|x|}\right)+\\&+
\sum_{k=1}^n\partial_{x^i}\left(\frac{x_k}{|x|}\right)\partial_{y_k}U\left(|x|,\frac{x}{|x|}\right)=\\&=\frac{x^i}{|x|}\partial_{\rho}U\left(|x|,\frac{x}{|x|}\right)+\\&+
\sum_{k=1}^n\left[\frac{\delta^{ik}}{|x|}-\frac{x^ix^k}{|x|^3}\right]\partial_{y^k}U\left(|x|,\frac{x}{|x|}\right).
\end{aligned}
\end{equation}
Now, \eqref{ias-0-9} implies

$$\sum_{k=1}^n\left[\frac{\delta^{ik}}{|x|}-\frac{x^ix^k}{|x|^3}\right]\partial_{y^k}U\left(|x|,\frac{x}{|x|}\right)=\frac{1}{|x|}d_iU(\rho,\omega).$$
By the latter and by \eqref{ias-1n-16} we get
\eqref{ias-1-16}.$\blacksquare$

\bigskip

\begin{theo}\label{ias-1p-18}
Let $u\in C^2\left(B_R\setminus\{0\}\right)$ and let $U$ be defined by
\eqref{ias-0-15}. Set
\begin{equation}\label{ias-0-18}
U(\rho,\omega)=u(\rho\omega),\quad \quad \forall (\rho,\omega)\in (0,R)\times\Sigma.
\end{equation}
We have $v\in
C^2\left((0,R)\times \Sigma\right)$ and
\begin{equation}\label{ias-1-18}
(\Delta
u)(\rho\omega)=\partial^2_{\rho}v(\rho,\omega)+\frac{n-1}{\rho}\partial_{\rho}v(\rho,\omega)+\frac{1}{\rho^2}\Delta_{\Sigma}v(\rho,\omega)
\end{equation}
\end{theo}

\bigskip

\textbf{Proof.} We denote, for the sake of brevity,
$$\rho(x)=|x|,\quad \omega(x)=\frac{x}{|x|}.$$
Let $i\leq i\leq n$; \eqref{ias-1-16} implies

\begin{equation}\label{ias-2-16}
\partial^2_{x^i}u(x)=\underset{J_1
}{\underbrace{\partial_{x^i}\left(\omega^i(x)\partial_{\rho}U(\rho(x),\omega(x))\right)}}+\underset{J_2
}{\underbrace{\partial_{x^i}\left(\frac{1}{\rho(x)}d_i
U(\rho(x),\omega(x))\right)}}.
\end{equation}

\medskip

\textbf{Computation of $J_1$.}

We have

\begin{equation}\label{ias-3-16}
 \begin{aligned}
&J_1=\partial_{x^i}\left(\omega^i(x)\right)\partial_{\rho}U(\rho(x),\omega(x))+\omega^i(x)\partial_{x^i}\left(\partial_{\rho}U(\rho(x),\omega(x))\right)=\\&
=\left(\frac{1}{|x|}-\frac{(x^i)^2}{|x|^3}\right)\partial_{\rho}U(\rho(x),\omega(x))+\\&+\omega^i(x)\left[\omega^i(x)\partial^2_{\rho}U(\rho(x),\omega(x))+
\frac{1}{\rho(x)}\left(d_i\partial_{\rho}U\right)(\rho(x),\omega(x))\right]=\\&=\left(\frac{1}{\rho(x)}-\frac{\omega^2_i(x)}{\rho(x)}\right)\partial_{\rho}U(\rho(x),\omega(x))
+\\&
+(\omega^i(x))^2\partial^2_{\rho}U(\rho(x),\omega(x))+\frac{\omega^i(x)}{\rho(x)}\left(d_i\partial_{\rho}U\right)(\rho(x),\omega(x)).
\end{aligned}
\end{equation}

\medskip

\textbf{Computation of $J_2$.}

We have

\begin{equation}\label{ias-1-17}
 \begin{aligned}
&J_2=\partial_{x^i}\left(\frac{1}{\rho(x)}\right)d_iU(\rho(x),\omega(x))+\frac{1}{\rho(x)}\partial_{x^i}\left((d_iU)(\rho(x),\omega(x))\right)=\\&=
-\frac{\omega^i(x)}{\rho^2(x)}d_iU(\rho(x),\omega(x))+\\&+\frac{\omega^i(x)}{\rho(x)}\partial_{\rho}\left(d_i
U\right)(\rho(x),\omega(x))+\frac{1}{\rho^2(x)}\left(d^2_iU\right)(\rho(x),\omega(x)),
\end{aligned}
\end{equation}
in the last equality we have applied  \eqref{ias-1-16}
to $(d_iU)(\rho(x),\omega(x))$.

\medskip

\textbf{Computation of $\Delta$.}

Now, adding up \eqref{ias-3-16} and \eqref{ias-1-17} we obtain
(we omit, the variables)

\begin{equation}\label{ias-1n-17}
 \begin{aligned}
&\Delta u=\sum_{i=1}^n\partial^2_{x^i}u=\\&=
\frac{n-1}{\rho}\partial_{\rho}U+\partial^2_{\rho}U+\frac{1}{\rho}\sum_{i=1}^n\omega^id_i\left(\partial_{\rho}U\right)-\\&
-\frac{1}{\rho^2}\sum_{i=1}^n\omega^id_iU+\frac{1}{\rho}\sum_{i=1}^n\omega^i\partial_{\rho}\left(d_i
U\right)+\frac{1}{\rho^2}\Delta_{\Sigma}U.
\end{aligned}
\end{equation}
Now,  \eqref{ias-1-9} implies

$$\sum_{i=1}^n\omega_id_iU=0,\quad
\sum_{i=1}^n\omega^id_i\left(\partial_{\rho}U\right)=0;$$ by the first equality we get
$$\sum_{i=1}^n\omega_i\partial_{\rho}\left(d_i
U\right)=\partial_{\rho}\left(\sum_{i=1}^n\omega^id_i U\right)=0.$$
Therefore,  \eqref{ias-1n-17} gives

\begin{equation*}
\Delta
u(x)=\partial^2_{\rho}U(\rho(x),\omega(x))+\frac{n-1}{\rho(x)}\partial_{\rho}U(\rho(x),\omega(x))+\frac{1}{\rho^2(x)}\left(\Delta_{\Sigma}
U\right)(\rho(x),\omega(x))
\end{equation*}
and by \eqref{ias-0-18} we get  \eqref{ias-1-18}.
$\blacksquare$

\bigskip

\textbf{Remarks.}

\noindent\textbf{1.}  Generally speaking, even  if $u\in C^2\left(B_R\right)$, $U(\rho,\omega)$ is not differentiable w.r.t.  $\rho$ in
$0$. However \eqref{ias-1-9} and \eqref{ias-1-16} give

\begin{equation}\label{ias-2-18}
\partial_{\rho}U(\rho(x),\omega(x))=\nabla u(x)\cdot \omega(x),
\end{equation}
from which we have

 \begin{equation}\label{ias-2n-18}
\left|\partial_{\rho}U(\rho,\omega)\right|\leq \max_{\overline{B_r}}|\nabla
u|,\quad\mbox{for } \rho\leq r\leq R.
\end{equation}

\medskip

\noindent\textbf{2.} Let us observe that 
\begin{equation}\label{ias-2nn-18}
\begin{aligned}
\partial^2_{\rho}U(\rho(x),\omega(x))&=\sum_{ij=1}^n\partial^2_{x^ix^j}u(x)\frac{x^ix^j}{|x|^2}=
\sum_{|\alpha=2|}\frac{2!}{\alpha!}\partial^{\alpha}u(x)\left(\frac{x}{|x|}\right)^{\alpha}=\\&=
\underset{\mbox{Hessian matrix}}{\underbrace{\partial^2u(x)}}\frac{x}{|x|}\cdot
\frac{x}{|x|}.
\end{aligned}
\end{equation}
Indeed we have
\begin{equation*}
\begin{aligned}
\partial^2_{\rho}U(\rho(x),\omega(x))&=\sum_{k=1}^n\omega^k(x)\partial_{x^k}\left(\sum_{j=1}^n\omega^j(x)\partial_{x^j}u\right)=\\&
=\sum_{k=1}^n\frac{x^k}{|x|}\sum_{j=1}^n\left[\left(\frac{\delta^{jk}}{|x|}-\frac{x^jx^k}{|x|^3}\right)\partial_{x^j}u+\frac{x^j}{|x|}\partial^2_{x^jx^k}u\right]=\\&=
\sum_{k=1}^n\frac{x^k}{|x|}\left[\frac{1}{|x|}\partial_{x^k}u-\frac{x^k}{|x|^2}\left(\frac{x}{|x|}\cdot
\nabla
u\right)+\sum_{j=1}^n\frac{x^j}{|x|}\partial^2_{x^jx^k}u\right]=\\&=\sum_{ij=1}^n\partial^2_{x^ix^j}u(x)\frac{x^ix^j}{|x|^2}.
\end{aligned}
\end{equation*}
Similarly can be checked that

\begin{equation}\label{ias-1-19}
\partial^m_{\rho}U(\rho(x),\omega(x))=\frac{1}{|x|^m}
\sum_{|\alpha|=m}\frac{m!}{\alpha!}\partial^{\alpha}u(x)x^{\alpha}.
\end{equation}

\medskip

\noindent\textbf{3.} Of course, by \eqref{ias-2-18} we have (for
$r\leq R$)

\begin{equation}\label{ias-2-19}
\frac{\partial u}{\partial
\nu}(x)=\partial_{\rho}U(r,\omega),\quad\mbox{for } x\in \partial B_r,
\end{equation}
where $\nu$ is the unit outward normal to $\partial B_r$. Similarly we have by \eqref{ias-1-19},

\begin{equation*}
\frac{\partial^m u}{\partial
\nu^m}(x):=\frac{d^m}{dt^m}u(x+t\nu)_{|t=0}=\partial^m_{\rho}U(r,\omega).
\end{equation*} $\blacklozenge$

\bigskip

\underline{\textbf{Exercise 1.}} Let $u\in C^2\left(B_R\right)$ and
$U(\rho,\omega)=u(\rho\omega)$. Apply formula
\eqref{ias-1-18} to prove that for any $r<R$ we have
\begin{equation*}
\int_{\Sigma}\partial_{\rho}U(\rho,\omega)d\omega=\frac{1}{r^{n-1}}\int_{B_r}\Delta
u(x) dx.
\end{equation*} $\clubsuit$

\bigskip

\underline{\textbf{Exercise 2.}} Let $u$ be a harmonic function in $B_R$.
Prove that for any $r<R$ we have
\begin{equation*}
\int_{\partial B_r} \frac{\partial^m u}{\partial \nu^m}(x)dS=0.
\end{equation*} $\clubsuit$

\section[The case of Laplace leading operator]{The case of Laplace leading operator} \label{aronszajn-constant}

In this Section and in the
sequel we will use the following \textbf{Hardy inequality}, \cite{l:HLP34}, \cite{l:T67}.
\begin{lem}[\textbf{the Hardy inequality}]
    \label{prop:Hardy}
    \index{Lemma:@{Lemma:}!- Hardy inequality@{- Hardy inequality}}
        If $f\in C^{\infty}_0(0,+\infty)$, then
\begin{equation}
    \label{eq:24.1}
        \int_0^{+\infty} \frac{f^2(s)}{s^2}ds\leq 4 \int_0^{+\infty} (f^{\prime}(s))^2ds.
\end{equation}
\end{lem}

\medskip

\textbf{Note.} \textit{the number $4$ on the right--hand side of \eqref{eq:24.1}
is the best constant.}

\medskip

\textbf{Proof.} We set
\begin{equation}
    \label{g}
        g(s)=s^{-\frac{1}{2}}f(s).
\end{equation}
and we have $$f(s)=s^{\frac{1}{2}}g(s),\quad
f^{\prime}(s)=\frac{1}{2}s^{-\frac{1}{2}}g(s)+s^{\frac{1}{2}}g^{\prime}(s).$$
Hence \eqref{eq:24.1} is equivalent to

\begin{equation}
    \label{eq:24.1-aeq}
        \int_0^{+\infty} s^{-1}g^2(s)ds\leq 4 \int_0^{+\infty} \left(\frac{1}{2}s^{-\frac{1}{2}}g(s)+s^{\frac{1}{2}}g^{\prime}(s)\right)^2ds.
\end{equation}
Now, we have
\begin{equation*}
\begin{aligned}
&4 \int_0^{+\infty}
\left(\frac{1}{2}s^{-\frac{1}{2}}g(s)+s^{\frac{1}{2}}g^{\prime}(s)\right)^2ds=\\&=
4 \int_0^{+\infty}\left(\frac{1}{4}s^{-1}g^2(s)+g(s)g^{\prime}(s)+s
\left(g^{\prime}(s)\right)^2\right)ds=\\& =4
\int_0^{+\infty}\left(\frac{1}{4}s^{-1}g^2(s)+\frac{1}{2}\left(g^2\right)^{\prime}(s)+s
\left(g^{\prime}(s)\right)^2\right)ds=\\&
=\int_0^{+\infty}\left(s^{-1}g^2(s)+4s
\left(g^{\prime}(s)\right)^2\right)ds\geq\\&\geq \int_0^{+\infty}
s^{-1}g^2(s)ds.
\end{aligned}
\end{equation*}
The proof is complete. $\blacksquare$


\bigskip

\begin{theo}[{\bf Carleman estimate for  $\Delta$}]\label{prop:Carlm-delta}
Let $\epsilon\in(0,1]$. Let us define
\begin{equation}
    \label{peso}
        \rho(x) = \phi_{\epsilon}\left(|x|\right), \ \ \mbox{for } x\in B_1,
\end{equation}
where
\begin{equation}
    \label{eq:24.3-delta}
        \phi_{\epsilon}(s) = \frac{s}{\left(1+s^{\epsilon}\right)^{1/\varepsilon}}.
\end{equation}
Then there exist $\tau_1>1$ and $C>1$, depending on $\epsilon$ only, such that

\begin{equation}\label{Carlm-delta-0}
\begin{aligned}
\tau^{3}\int\rho^{\epsilon-2\tau}|u|^2dx+&
\tau\int\rho^{2+\epsilon-2\tau}|\nabla u|^2dx\leq\\&\leq
C\int\rho^{4-2\tau}|\Delta u|^2dx,
\end{aligned}
\end{equation}
for every $u\in
C^\infty_0\left(B_1\setminus \{0\}\right)$ and for every $\tau\geq \tau_1$ .

Moreover, there exist $\tau_2>1$, $C>1$, depending on $\epsilon$ only,
such that
\begin{equation}\label{Carlm-delta}
\begin{aligned}
\tau^{3}\int\rho^{\epsilon-2\tau}|u|^2dx+&
\tau\int\rho^{2+\epsilon-2\tau}|\nabla u|^2dx+\\&
+\tau^2r\int\rho^{-1-2\tau}u^2dx\leq C\int\rho^{4-2\tau}|\Delta
u|^2dx,
\end{aligned}
\end{equation}
for every $r\in (0,1)$, for every $u\in
C^\infty_0\left(B_1\setminus \overline{B}_{r/4}\right)$ and for every $\tau\geq \tau_2$.
\end{theo}

\bigskip
\textbf{Remark.} We have
\begin{equation}
    \label{peso-osserv}
        \frac{|x|}{2^{1/\epsilon}}\leq \rho(x) \leq |x|, \quad\forall x\in B_1.
\end{equation}
$\blacklozenge$

\bigskip

\textbf{Proof.} It is not restrictive to assume  that $u$ is a real--valued function.
 First we prove \eqref{Carlm-delta-0}, afterwards, with a few modifications we will prove \eqref{Carlm-delta}.

Let $u$ be an arbitrary function of $C^\infty_0\left(B_1\setminus
\{0\}\right)$ and let us consider the \\ $n$-dimensional Laplace operator  in the 
polar coordinates  $(\varrho,\omega)$, that is (recalling that
$\Sigma=\partial B_1$)
\begin{equation}\label{laplace-polar-1}
    \Delta u=u_{\varrho\varrho}+\frac{n-1}{\varrho}u_{\varrho}+\frac{1}{\varrho^2}\Delta_{\Sigma} u, \ \ \forall  (\varrho,\omega)\in (0,\infty)\times \Sigma.
\end{equation}
Let us perform the following change of variables
$$\varrho=e^t,\quad\quad\widetilde{u}(t,\omega)=u\left(e^t,\omega\right),\quad\forall (t,\omega)\in (-\infty,0)\times \Sigma.$$
We have, for every $(t,\omega)\in (-\infty,0)\times \Sigma$,
\begin{equation}\label{laplace-polar-2}
    e^{2t}(\Delta u)(e^t,\omega)=\mathcal{L}\widetilde{u}:=\widetilde{u}_{tt}+(n-2)\widetilde{u}_{t}+
    \Delta_{\Sigma}\widetilde{u}.
\end{equation}
For the sake of brevity, for any function $h\in
C^\infty_0\left(B_1\setminus \{0\}\right)$ we will write  $h^{\prime}$,
$h^{\prime\prime}$, ... instead of $h_t$, $h_{tt}$, ... .

By \eqref{peso} we have (we omit the subscript $\epsilon$ from now on)
\begin{equation}\label{peso-1}
\varphi(t):=\log(\phi(e^t))=t-\epsilon^{-1}\log\left(1+e^{\epsilon
t}\right), \quad \forall t\in(-\infty,0).
\end{equation}
We get
\begin{equation}\label{peso-derivate}
\varphi^{\prime}(t)=\frac{1}{1+e^{\epsilon t}}, \quad
\varphi^{\prime\prime}(t)=-\frac{\epsilon e^{\epsilon
t}}{(1+e^{\epsilon t})^2},\quad \forall  t\in(-\infty,0).
\end{equation}
Let
\begin{equation*}
f(t,\omega)=e^{-\tau\varphi}\widetilde{u}(t,\omega) , \quad\forall  (t,\omega)\in(-\infty,0)\times \Sigma.
\end{equation*}
We have
\begin{equation}\label{L-tau}
\mathcal{L}_{\tau}f:=e^{-\tau\varphi}\mathcal{L}(e^{\tau\varphi}f)=\underset{\mathcal{A}_{\tau}f}{\underbrace{b_0f+b_1f^{\prime}}}+
\underset{\mathcal{S}_{\tau}f}{\underbrace{a_0f+f^{\prime\prime}+\Delta_{\Sigma}f}},
\end{equation}
where
\begin{equation}\label{coefficienti}
a_0=\tau^2\varphi^{\prime^2}+\tau(n-2), \quad b_0=\tau
\varphi^{\prime\prime},\quad b_1=2\tau\varphi^{\prime}+(n-2).
\end{equation}
Let us denote by $\int(\cdot)$ the integral
$\int^{0}_{-\infty}\int_{\Sigma}(\cdot) d\omega dt$ and set
\begin{equation}\label{gamma}
\gamma:=\frac{1}{\varphi'}=1+e^{\epsilon t}.
\end{equation}
We obtain
\begin{equation}\label{square}
\int\gamma\left\vert \mathcal{L}_{\tau}f\right\vert^2\geq
2\int\gamma\mathcal{A}_{\tau}f\mathcal{S}_{\tau}f+\int\gamma\left\vert\mathcal{A}_{\tau}f\right\vert^2
\end{equation}
and
\begin{equation}\label{commutator}
\begin{aligned}
2\int\gamma\mathcal{A}_{\tau}f\mathcal{S}_{\tau}f&=\underset{I_1}{\underbrace{2\int\gamma\left(b_0f+b_1f^{\prime}\right)\Delta_{\Sigma}f}}+\\&
+\underset{I_2}{\underbrace{2\int\gamma\left(b_0f+b_1f^{\prime}\right)
\left(a_0f+f^{\prime\prime}\right) }}.
\end{aligned}
\end{equation}

\textbf{We examine $I_1$}.

For any function $f$, $g$ on $\in
C^\infty_0\left(B_1\setminus \{0\}\right)$, let us denote by

$$\langle \nabla _{\Sigma }f,\nabla _{\Sigma }g\rangle
=\sum_{i=1}^nd_i f d_ig,\quad\mbox{ }\left\vert
\nabla_{\Sigma}f\right\vert^2=\langle \nabla _{\Sigma }f,\nabla
_{\Sigma }f\rangle.$$

Integraing by parts, using Proposition
\eqref{ias-prop2-13} and taking into account \eqref{gamma}, we have
\begin{equation*}
\begin{aligned}
I_1&=2\int\left(\gamma b_0f\Delta_{\Sigma}f+\gamma
b_1f^{\prime}\Delta_{\Sigma}f\right)=\\& =2\int\left(-\gamma
b_0\left\vert \nabla_{\Sigma}f\right\vert^2 - \gamma b_1\left\langle
\nabla _{\Sigma }f,\nabla _{\Sigma
}f^{\prime}\right\rangle\right)=\\& =2\int\left(-\gamma b_0
\left\vert \nabla_{\Sigma}f\right\vert^2- \frac{1}{2}\gamma
b_1\left(\left\vert
\nabla_{\Sigma}f\right\vert^2\right)^{\prime}\right)=\\&=2\int\left(-\gamma
b_0 + \frac{1}{2}\left(\gamma b_1\right)^{\prime}\right)\left\vert
\nabla_{\Sigma}f\right\vert^2.
\end{aligned}
\end{equation*}
By \eqref{peso-derivate}, \eqref{coefficienti} and \eqref{gamma} 
we obtain
\begin{equation}\label{calcolo-1}
\begin{aligned}
-\gamma b_0 + \frac{1}{2}\left(\gamma
b_1\right)^{\prime}&=-\gamma\tau\varphi^{\prime\prime}+\frac{1}{2}\gamma
b_1^{\prime}+\frac{1}{2}\gamma^{\prime}b_1=\\&
=-\gamma\tau\varphi^{\prime\prime}+\gamma\tau\varphi^{\prime\prime}+\frac{1}{2}\left(2\tau\varphi^{\prime}+n-2\right)\gamma^{\prime}=\\&
=\frac{\epsilon e^{\epsilon t}}{2}\left(\frac{2\tau}{1+ e^{\epsilon
t}}+n-2\right)\geq \frac{\tau}{2}\epsilon e^{\epsilon t}, \ \  \forall  \tau>0.
\end{aligned}
\end{equation}
Hence, we have
\begin{equation}\label{I-1}
I_1\geq\int\tau\epsilon e^{\epsilon t}\left\vert
\nabla_{\Sigma}f\right\vert^2, \quad\forall  \tau>0.
\end{equation}

\textbf{Now we examine $I_2$}.

Integration by parts gives
\begin{equation}\label{I-2}
\begin{aligned}
I_2&=2\int\gamma\left(b_0f+b_1f^{\prime}\right)\left(a_0f+f^{\prime\prime}\right)=\\&
=2\int\gamma\left(a_0b_0f^2+b_0ff^{\prime\prime}+b_1a_0f^{\prime}f+b_1f^{\prime}f^{\prime\prime}\right)=\\&
=2\int\gamma a_0b_0f^2-\left(\gamma
b_0f\right)^{\prime}f^{\prime}+\frac{1}{2}\gamma
b_1a_0\left(f^2\right)^{\prime}+\frac{1}{2}\gamma b_1
\left(f^{\prime^2}\right)^{\prime}=\\&=2\int\left[\gamma
a_0b_0-\frac{1}{2}\left(\gamma b_1
a_0\right)^{\prime}\right]f^2-\left(\gamma
b_0\right)^{\prime}ff^{\prime}-\\&-\int \gamma b_0f^{\prime^2}
+\frac{1}{2}\left(\gamma
b_1\right)^{\prime}f^{\prime^2}=\\&=2\int\underset{H_1}{\underbrace{\left[\gamma
a_0b_0-\frac{1}{2}\left(\gamma b_1
a_0\right)^{\prime}+\frac{1}{2}\left(\gamma
b_0\right)^{\prime\prime}\right]}}f^2-\\&-\int\underset{H_2}{\underbrace{\left[\gamma
b_0+\frac{1}{2}\left(\gamma b_1\right)^{\prime}\right]}}f^{\prime^2}.
\end{aligned}
\end{equation}
Let now  examine $H_1$.

Since by \eqref{coefficienti},  $H_1$ is a polynomial of third degree
w.r.t. $\tau$, we begin by evaluating the coefficient of
$\tau^3$.

Let us notice that the terms of $H_1$ have the following behavior, as $\tau \rightarrow +\infty$

\begin{equation*}
\begin{aligned}
&\gamma a_0b_0=\mathcal{O}\left(\tau^3\right),  \\&
-\frac{1}{2}\left(\gamma b_1 a_0\right)^{\prime}=\mathcal{O}\left(\tau^3\right),  \\&
\frac{1}{2}\left(\gamma b_0\right)^{\prime\prime}
=\mathcal{O}\left(\tau\right).
\end{aligned}
\end{equation*}
Hence, let us first examine the term
$$\widetilde{H}_1:=\gamma a_0b_0-\frac{1}{2}\left(\gamma b_1 a_0\right)^{\prime}.$$
By \eqref{peso-derivate}, \eqref{coefficienti} and \eqref{gamma}
we have
\begin{equation}\label{H-tilde-1-0}
\begin{aligned}
\widetilde{H}_1&=\gamma a_0b_0-\frac{1}{2}\left(\gamma b_1
a_0\right)^{\prime}=\\&
=\gamma\left(\tau^2\varphi^{\prime^2}+\tau(n-2)\varphi^{\prime}\right)\tau\varphi^{\prime\prime}-\frac{1}{2}b_1^{\prime}\left(\gamma
a_0\right)-\frac{1}{2}b_1\left(\gamma
a_0\right)^{\prime}=\\&=\gamma\left(\tau^2\varphi^{\prime^2}+\tau(n-2)\varphi^{\prime}\right)\tau\varphi^{\prime\prime}
-\\&-\frac{1}{2}\left(2\tau\varphi^{\prime}+(n-2)\right)^{\prime}\gamma\left(\tau^2\varphi^{\prime^2}+\tau(n-2)\varphi^{\prime}\right)-\frac{1}{2}b_1\left(\gamma
a_0\right)^{\prime}= \\& =-\frac{1}{2}b_1\left(\gamma
a_0\right)^{\prime} =\\&=-\frac{1}{2}\left(\frac{2\tau}{1+e^{\varepsilon
t}}+n-2\right)\left[\gamma\left(\tau^2\varphi^{\prime^2}+\tau(n-2)\varphi^{\prime}\right)\right]^{\prime}=\\&
=-\left(\frac{\tau}{1+e^{\epsilon
t}}+\frac{n-2}{2}\right)\left(\tau^2\frac{1}{1+e^{\epsilon
t}}+\tau(n-2)\right)^{\prime}=\\&=
\left(\tau^3\frac{1}{1+e^{\epsilon
t}}+\tau^2\frac{n-2}{2}\right)\frac{\epsilon e^{\epsilon
t}}{(1+e^{\epsilon t})^2}.
\end{aligned}
\end{equation}
Then, using the trivial inequality
$$\frac{1}{1+e^{\epsilon t}}\geq \frac{1}{2},\quad\quad \forall t\in (-\infty,0),$$ we have (for $t\in (-\infty,0)$)

\begin{equation}\label{H-tilde-1}
\widetilde{H}_1\geq \frac{\tau^3}{8} \epsilon e^{\epsilon t}, \quad
\forall \tau>0.
\end{equation}
and
\begin{equation}\label{altrotermine}
\frac{1}{2}\left(\gamma
b_0\right)^{\prime\prime}=\frac{\tau}{2}\left(\frac{-\epsilon
e^{\epsilon t}}{1+e^{\epsilon
t}}\right)^{\prime\prime}=-\frac{\tau}{2}\frac{\epsilon^3
e^{\epsilon t}\left(1-e^{\epsilon t}\right)}{(1+e^{\epsilon
t})^3}\geq-\frac{\tau}{2}\epsilon^3e^{\epsilon t}, \ \forall
\tau>0.
\end{equation}
Inequalities \eqref{H-tilde-1} and \eqref{altrotermine} give
\begin{equation}\label{H-1}
\begin{aligned}
H_1=\widetilde{H}_1+\frac{1}{2}\left(\gamma
b_0\right)^{\prime\prime}\geq \frac{\tau^3}{8} \epsilon e^{\epsilon
t}-\frac{\tau}{2}\epsilon^3e^{\epsilon t}\geq \frac{\tau^3}{16}
\epsilon e^{\epsilon t}, \quad \forall \tau>\sqrt{8}\epsilon.
\end{aligned}
\end{equation}
Now, let us consider $H_2$.
\begin{equation}\label{H-2-0}
\begin{aligned}
&H_2=\gamma b_0+\frac{1}{2}\left(\gamma b_1\right)^{\prime}=\\&
=\tau(1+e^{\epsilon
t})\varphi^{\prime\prime}+\frac{1}{2}\left[2\tau(1+e^{\epsilon
t})\varphi^{\prime}+(n-2)(1+e^{\epsilon t})\right]^{\prime}=\\&
=-\tau\frac{\epsilon e^{\epsilon t}}{1+e^{\epsilon
t}}+\frac{1}{2}\left[2\tau+(n-2)(1+e^{\epsilon
t})\right]^{\prime}=\\& =-\tau\frac{\epsilon e^{\epsilon
t}}{1+e^{\epsilon t}}+\frac{(n-2)\epsilon e^{\epsilon t}}{2}.
\end{aligned}
\end{equation}
Hence
\begin{equation}\label{H-2}
-2H_2\geq\frac{\tau}{2} \epsilon e^{\epsilon t},\quad \forall
\tau\geq 2(n-2).
\end{equation}

By \eqref{I-2}, \eqref{H-1} and \eqref{H-2} we have
\begin{equation}\label{I-2-1}
I_2\geq \int \frac{\tau^3}{8}\epsilon e^{\epsilon t}f^2+\frac{\tau}{2}\epsilon
e^{\epsilon t}f^{\prime^2}, \quad \forall \tau\geq \tau_1,
\end{equation}
where $\tau_1=\max\{\sqrt{8}\varepsilon,2(n-2) \}$.

Now, \eqref{commutator}, \eqref{I-1} and \eqref{I-2-1} give

\begin{equation}\label{commutator-1}
2\int\gamma\mathcal{A}_{\tau}f\mathcal{S}_{\tau}f\geq
\frac{\epsilon}{8}\int \left(\tau^3f^2+\tau
\left(f^{\prime^2}+\left\vert
\nabla_{\Sigma}f\right\vert^2\right)\right)e^{\epsilon t},\quad
\forall \tau\geq \tau_1.
\end{equation}
By \eqref{square} and \eqref{commutator-1} we have
\begin{equation}\label{commutator-final}
\int\gamma\left\vert \mathcal{L}_{\tau}f\right\vert^2\geq
\frac{\epsilon \tau^3}{8}\int f^2+\frac{\epsilon \tau}{8}
\int\left(f^{\prime^2}+\left\vert
\nabla_{\Sigma}f\right\vert^2\right)e^{\epsilon t},
\end{equation}
for every $\tau\geq \tau_1$ and for every $f\in
C_0^{\infty}((-\infty,0)\times\Sigma)$.

Now, in order to obtain \eqref{Carlm-delta-0} we come back to $u$ and to the
original variables. Let us recall that 
$f(t,\omega)=e^{-\tau\varphi}u(e^t,\omega)$. By using
\eqref{peso}, \eqref{laplace-polar-2} and \eqref{L-tau}, we get

\begin{equation}\label{final-1}
\begin{aligned}
\int^{0}_{-\infty}\int_{\Sigma}\left\vert
\mathcal{L}_{\tau}f\right\vert^2d\omega
dt=&\int^{0}_{-\infty}\int_{\Sigma}e^{-2\tau\varphi(t)}e^{4t}|(\Delta
u)(e^t,\omega)|^2d\omega dt= \\& =
\int^{1}_{0}\int_{\Sigma}e^{-2\tau\varphi(\log\varrho)}\varrho^{3}|(\Delta
u)(\varrho,\omega)|^2d\omega
d\varrho=\\&=\int_{B_1}\rho^{-2\tau}|x|^{4-n}|\Delta u|^2dx.
\end{aligned}
\end{equation}

Similarly, we get
\begin{equation}\label{final-3}
\int^{0}_{-\infty}\int_{\Sigma} f^2e^{\epsilon t}d\omega
dt=\int_{B_1}\rho^{-2\tau}|x|^{\epsilon-n} u^2dx.
\end{equation}
Concerning the second integral on the right--hand side of
\eqref{commutator-final}, let $\delta\in (0,1)$ be a number that we will choose later, we have
\begin{equation}\label{final-4}
\begin{aligned}
&\int^{0}_{-\infty}\int_{\Sigma} e^{\epsilon
t}\left(f^{\prime^2}+\left\vert
\nabla_{\Sigma}f\right\vert^2\right)d\omega dt\geq
\\&\geq \delta\int^{0}_{-\infty}\int_{\Sigma} e^{\epsilon
t}\left(f^{\prime^2}+\left\vert
\nabla_{\Sigma}f\right\vert^2\right)d\omega dt\geq \\& \geq
\frac{\delta}{2}\int^{0}_{-\infty}\int_{\Sigma}e^{\epsilon
t}e^{-2\tau\varphi(t)}\left(|u_{\varrho}(e^t,\omega)|^2e^{2t}+\right.
\\& \left.+ \left\vert
\nabla_{\Sigma}u(e^t,\omega)\right\vert^2-2\tau^2|u(e^t,\omega)|^2\right)
d\omega dt=\\&
=\frac{\delta}{2}\int_{B_1}\rho^{-2\tau}|x|^{\epsilon-n}
\left(|x|^{2}|\nabla u|^2-2\tau^{2}|u|^2\right)dx.
\end{aligned}
\end{equation}

\medskip

Now, let us choose $\delta=\frac{1}{2}$ so that, by \eqref{commutator-final} and
\eqref{final-1}--\eqref{final-4}, we have

\begin{equation}\label{final-5}
\begin{aligned}
\int_{B_1}\rho^{-2\tau}|x|^{4-n}|\Delta
u|^2dx\geq&\frac{\varepsilon\tau}{16}\int_{B_1}\rho^{-2\tau}|x|^{\epsilon+2-n}
|\nabla u|^2dx+ \\&
+\frac{\epsilon\tau^3}{16}\int_{B_1}\rho^{-2\tau}|x|^{\epsilon-n}
u^2dx,
\end{aligned}
\end{equation}
for every $u\in C^\infty_0(B_1\setminus
\{0\})$ and for every $\tau\geq \tau_1$. Finally, by \eqref{peso-osserv},  we substitute in 
\eqref{final-5} $\tau$ by $(\tau-\frac{n}{2})$ and we find inequality \eqref{Carlm-delta-0}.

\bigskip

Let us now look at \eqref{Carlm-delta}. Let $u\in
C^\infty_0(B_1\setminus \overline{B}_{r/4})$.

By \eqref{square} and \eqref{commutator-final} we have

\begin{equation}\label{commutator-final-1}
\int\left\vert \mathcal{L}_{\tau}f\right\vert^2\geq
\frac{\epsilon}{8}\int \left(\tau^3f^2+\tau
\left(f^{\prime^2}+\left\vert
\nabla_{\Sigma}f\right\vert^2\right)\right)e^{\epsilon t}+
\int\gamma\left\vert\mathcal{A}_{\tau}f\right\vert^2.
\end{equation}
To obtain the first term on the left--hand side of \eqref{Carlm-delta}
we estimate from below the last term on the right--hand side of \eqref{commutator-final-1}.

Let us note that by the trivial inequality
$(a+b)^2\geq\frac{1}{2}a^2-b^2$ and by \eqref{peso-derivate},
\eqref{gamma} we have

\begin{equation}\label{antisimm-1}
\begin{aligned}
\int\gamma\left\vert\mathcal{A}_{\tau}f\right\vert^2&\geq\frac{1}{2}\int\gamma
\left(2\tau\varphi^{\prime}+n-2\right)^2f^{\prime^2}-\int\gamma\tau^2
\varphi^{\prime\prime^2}f^2\geq\\& \geq\tau^2\int f^{\prime^2}-
\epsilon^2\tau^2\int e^{2\epsilon t}f^2, \quad \forall \tau>0
\end{aligned}
\end{equation}
Using inequality \eqref{antisimm-1} in
\eqref{commutator-final-1}, we have

\begin{equation}\label{antisimm-2}
\begin{aligned}
\int\left\vert \mathcal{L}_{\tau}f\right\vert^2\geq &\tau^2\int
f^{\prime^2} +\epsilon\tau^3\int
\left(\frac{1}{8}-\epsilon\tau^{-1}e^{\epsilon t}\right) e^{\epsilon
t}f^2+\\&+\frac{\epsilon}{8}\tau\int \left(f^{\prime^2}+\left\vert
\nabla_{\Sigma}f\right\vert^2\right)e^{\epsilon t}, \ \  \forall \tau\geq\tau_1.
\end{aligned}
\end{equation}
Now, since $\left(\frac{1}{8}-\epsilon\tau^{-1}e^{\epsilon
t}\right)\geq \frac{1}{16}$ for every $\tau\geq 4\epsilon$, by \eqref{antisimm-2},
we get 
\begin{equation}\label{antisimm-3}
\begin{aligned}
\int\left\vert \mathcal{L}_{\tau}f\right\vert^2\geq &\tau^2\int
f^{\prime^2}+ \frac{\epsilon\tau^3}{16}\int e^{\epsilon t}f^2
+\frac{\epsilon}{8}\tau\int \left(f^{\prime^2}+\left\vert
\nabla_{\Sigma}f\right\vert^2\right)e^{\epsilon t},
\end{aligned}
\end{equation}
for every  $\tau\geq\tau_2$, where $\tau_2=\max\{4\epsilon, \tau_1\}$.

Proposition \ref{prop:Hardy} implies
\begin{equation}\label{new-Hardy}
\begin{aligned}
\int^0_{-\infty} f^2(t,\omega)e^{-t}dt=&\int^1_{0} s^{-2}f^2(\log
s,\omega)ds\leq\\&\leq  4\int^1_{0} \left\vert\frac{\partial}{\partial
s}f(\log s,\omega)\right\vert^2 ds= \\& = 4\int^0_{-\infty}
f^{\prime^2}(t,\omega)e^{-t}dt, \quad \forall \omega\in
\Sigma.
\end{aligned}
\end{equation}

\medskip

On the other hand, since $u\in C^\infty_0(B_1\setminus
\overline{B}_{r/4})$, we have $f(t,\omega)=0$ for every $t\leq \log
(r/4)$ and for every $\omega\in \Sigma$, \eqref{new-Hardy} gives
$$\int^0_{-\infty} f^2(t,\omega)e^{-t}dt
\leq 4\int^{\log\frac{r}{4}}_{-\infty}
f^{\prime^2}(t,\omega)e^{-t}dt\leq\frac{16}{r}\int^{0}_{-\infty}
f^{\prime^2}(t,\omega)dt, \quad \forall \omega\in
\Sigma.$$ Let us integrate over $\Sigma$ both the sides of the just obtained inequality
 and let us use \eqref{antisimm-3} to obtain
\begin{equation}\label{antisimm-4}
\int f^2e^{-t}\leq\frac{16}{r}\int f^{\prime2}\leq
\frac{16}{\tau^2r}\int\left\vert \mathcal{L}_{\tau}f\right\vert^2, \ \ 
\forall \tau\geq \tau_2.
\end{equation}

By \eqref{antisimm-4} and \eqref{antisimm-3} we have
\begin{equation}\label{antisimm-5}
\begin{aligned}
C\int\left\vert \mathcal{L}_{\tau}f\right\vert^2 &\geq 
\epsilon\tau^3\int e^{\epsilon t}f^2+\epsilon\tau\int
\left(f^{\prime^2}+\left\vert
\nabla_{\Sigma}f\right\vert^2\right)e^{\epsilon t}+\\& +\tau^2r\int
f^2e^{-t}, \quad \ \ \forall \tau\geq \tau_2,
\end{aligned}
\end{equation}
where $C$ is a constant.

Finally, by \eqref{final-5}, \eqref{peso-osserv} and by 
\begin{equation}\label{final-2}
\int^{0}_{-\infty}\int_{\Sigma} f^2e^{-t}d\omega
dt=\int_{B_1}\rho^{-2\tau}|x|^{-1-n} u^2dx,
\end{equation}
we obtain \eqref{Carlm-delta}. $\blacksquare$

\section{Proof of the optimal three sphere and the doubling inequality}\label{Dim-tre-sfere-doubling}

In the next Theorem we will consider the equation  
\begin{equation}
    \label{2-1-appdoub}
    \Delta U=b(x)\cdot\nabla U+c(x)U,\quad\mbox{in } B_1,
\end{equation}
where $b\in L^{\infty}\left(B_1;\mathbb{R}^n\right)$ and $c\in
L^{\infty}\left(B_1\right)$. Moreover, set
\begin{equation}
    \label{limitaz-coeff}
M=\max\left\{\left\Vert
b\right\Vert_{L^{\infty}\left(B_1;\mathbb{R}^n\right)},\left\Vert
c\right\Vert_{L^{\infty}\left(B_1\right)}\right\}.
\end{equation}

\bigskip

\begin{theo}[\textbf{optimal three sphere and doubling inequality}]\label{theo:40.teo}
	\index{Theorem:@{Theorem:}!- optimal three sphere and doubling inequality@{- optimal three sphere and doubling inequality}}
	Let us assume that $U\in H^2\left(B_{1}\right)$
    is a solution to equation \eqref{2-1-appdoub}.
Let $x_0\in B_1$ and $0<R_0\leq 1-|x_0|$. Then there exists $C>1$
depending on $M$ only, such that, if $0<2r<R<\frac{R_0}{2}$ then

\begin{equation}
    \label{three-sphere-enunciato}
    \int_{B_{R}(x_0)}U^2
    \leq
    C\left(\frac{R_0}{R}\right)^{C}\left(\int_{B_{r}(x_0)}U^2\right)^{\theta}\left(\int_{B_{R_0}(x_0)}U^2\right)^{1-\theta},
\end{equation}
where
\begin{equation}
    \label{three-sphere-5}
\theta=\frac{\log\frac{R_0}{2R}}{\log\frac{2R_0}{r}}.
\end{equation}

Moreover, if $U$ does not vanish identically in $B_{R_0/4}(x_0)$
then the following doubling inequality\index{doubling inequality} holds 
\begin{equation}\label{eq:10.6.1102-cube}
    \int_{B_{2r}(x_0)}U^2\leq C N_{x_0,R_0}^3\int_{B_{r}(x_0)}U^2,
\end{equation}
where
\begin{equation}
    \label{eq:10.6.1108-cube}
    N_{x_0,R_0}=\frac{\int_{B_{R_0}(x_0)}U^2}{\int_{B_{R_0/4}(x_0)}U^2}.
    \end{equation}
\end{theo}

\bigskip

In order to prove Theorem \ref{theo:40.teo} we need the following

\begin{lem}\label{newlemma1-cube} Under the same assumption of Theorem
\ref{theo:40.teo} we have, for every $x_0\in B_1$, $R$ and for every $R$, $r$ such that $0<2r<R<\frac{R_0}{2}$, where $R_0\leq 1-|x_0|$,
\begin{equation} \label{Lemma-doub-cube}
	\begin{aligned}
	 &R(2r)^{-2\tau}\int_{B_{2r}(x_0)}U^2+R^{1-2\tau}\int_{B_{R}(x_0)}
	U^2\leq\\&
	\leq C \overline{M}^2\left[\left(\frac{r}{4}\right)^{-2\tau}\int_{B_r(x_0)}U^2+
	\left(\frac{R_0}{2}\right)^{-2\tau}\int_{B_{R_0}(x_0)}U^2
	\right],	
\end{aligned}
\end{equation}
for every $\tau\geq \widetilde{\tau}_2$, where $\widetilde{\tau}_2$ and
$C\geq 1$ depend on $M$ only.
\end{lem}

\bigskip

\textbf{Proof of Lemma.} By a translation we may assume that $x_0=0$. Let $r,R$ satisfy
\begin{equation}
    \label{eq:25.1-cube}
0<2r<R<\frac{R_0}{2}.
\end{equation}
Let $\eta\in C^\infty_0((0,R_0))$ satisfy
\begin{equation}
    \label{eq:25.2-cube}
    0\leq \eta\leq 1,
\end{equation}
\begin{equation}
    \label{eq:25.4-cube}
\eta=0, \quad \hbox{ in }\left(0,\frac{r}{4}\right)\cup
\left(\frac{2R_0}{3}, R_0\right); \quad \eta=1, \quad \mbox{in
}\left[\frac{r}{2}, \frac{R_0}{2}\right],
\end{equation}
\begin{equation}
    \label{eq:25.6-cube}
\left|\frac{d^k\eta}{dt^k}(t)\right|\leq C r^{-k}, \quad \hbox{ in
}\left(\frac{r}{4}, \frac{r}{2}\right),\quad\mbox{for } 0\leq k\leq
2,
\end{equation}
\begin{equation}
    \label{eq:25.7-cube}
\left|\frac{d^k\eta}{dt^k}(t)\right|\leq CR_0^{-k}, \quad \mbox{in
}\left(\frac{R_0}{2}, \frac{2R_0}{3}\right),\quad\hbox{for } 0\leq
k\leq 2.
\end{equation}
We define
\begin{equation}
    \label{eq:25.5-cube}
\xi(x)=\eta(|x|).
\end{equation}
Exploiting Carleman estimate \eqref{Carlm-delta} and fixing there $\epsilon=1$, we get

\begin{equation}\label{Lemma-doub:1-1}
\begin{aligned}
\tau^{3}\int_{B_{R_0}}\rho^{1-2\tau}u^2dx&+
\tau\int_{B_{R_0}}\rho^{3-2\tau}|\nabla u|^2+\\&
+\tau^2r\int_{B_{R_0}}\rho^{-1-2\tau}u^2\leq
C\int_{B_{R_0}}\rho^{4-2\tau}|\Delta u|^2,
\end{aligned}
\end{equation}
for every $u\in
C^\infty_0\left(B_{R_0}\setminus \overline{B}_{r/4}\right)$ and for every $\tau\geq \tau_2$  
(we recall that $\tau_2$ and $C$ depend neither on $r$ nor on
$R_0$ and that the value $C$ may change from line to line).

\medskip

Since $\xi U\in H^2_0(B_{R_0})$, by density we can
apply Carleman estimate \eqref{Lemma-doub:1-1} to $u=\xi U$. Hence we find
\begin{equation}\label{eq:26.1-cube}
\begin{aligned}
&\tau^{3}\int_{B_{R_0}}\rho^{1-2\tau}\xi^2U^2+
\tau\int_{B_{R_0}}\rho^{3-2\tau}|\nabla (\xi U)|^2+\\&
+\tau^2r\int_{B_{R_0}}\rho^{-1-2\tau}\xi^2U^2\leq
C\int_{B_{R_0}}\rho^{4-2\tau}|\Delta (\xi U)|^2,
\end{aligned}
\end{equation}
 for every $\tau\geq \tau_2$. 
 
 Since we have

\begin{equation}\label{Lemma-doub:2-2}
\left|\Delta (\xi U)\right|^2\leq 2\xi^2\left|\Delta
U\right|^2+C\left( \left|\partial^{2} \xi\right|^2U^2+ \left|\nabla
\xi\right|^2\left|\nabla U\right|^2\right),
\end{equation}
setting

\begin{equation}
    \label{eq:27.1-cube}
    J_0 =\int_{B_{r/2}\setminus B_{r/4}}\rho^{4-2\tau}
    \left(r^{-4} U^2+r^{-2}|\nabla U|^2\right),
\end{equation}
\begin{equation}
    \label{eq:27.2-cube}
    J_1 =\int_{B_{2R_0/3}\setminus B_{R_0/2}}\rho^{4-2\tau}
   \left(U^2+|\nabla U|^2\right),
\end{equation}
we get
\begin{equation}
\begin{aligned}\label{eq:27.3-cube}
\tau^{3}\int_{B_{R_0}}\rho^{1-2\tau}\xi^2U^2dx&+
\tau\int_{B_{R_0}}\rho^{3-2\tau}|\nabla (\xi
U)|^2dx+\tau^2r\int_{B_{R_0}}\rho^{-1-2\tau}\xi^2U^2dx\leq \\& \leq
C\int_{B_{R_0}}\rho^{4-2\tau}|\Delta U|^2+C\left(J_0+J_1\right),
\end{aligned}
\end{equation}
for every $\tau\geq \tau_2$.

Now we perform what follows: we use \eqref{eq:25.1-cube}--\eqref{eq:25.5-cube} and
\eqref{Lemma-doub:2-2}, we estimate trivially from below the left--hand side of
\eqref{eq:27.3-cube} and we estimate trivially from above the right--hand side of
\eqref{eq:27.3-cube}, obtaining
\begin{equation}\label{eq:33.1-cube}
	\begin{aligned}
	&\tau^{3}\int_{B_{R_0/2}
		\setminus B_{r/2}}\rho^{1-2\tau}U^2+ \tau\int_{B_{R_0/2} \setminus
		B_{r/2}}\rho^{3-2\tau}|\nabla U|^2+\\&
	+\tau^2r\int_{B_{R_0}}\rho^{-1-2\tau}\xi^2U^2dx\leq \\&\leq 
	CM^2\int_{B_{R_0/2} \setminus B_{r/2}}
	\rho^{4-2\tau}\left(U^2+|\nabla U|^2\right)+\\&
	+C\overline{M}^2(J_0+J_1),	\end{aligned}
\end{equation}
for every $\tau\geq \tau_2$, where $\overline{M}=\sqrt{M^2+1}$.

\medskip

By \eqref{eq:33.1-cube}, we get
\begin{equation}\label{Lemma-doub:6}
	\begin{aligned}
	&\int_{B_{R_0/2}
		\setminus B_{r/2}}\left(\tau^{3}-CM^2\rho^3\right)\rho^{1-2\tau}U^2+\\&+
	\int_{B_{R_0/2} \setminus B_{r/2}}\left(\tau-CM^2\rho\right)\rho^{3-2\tau}|\nabla U|^2+\\
	& +\tau^2r\int_{B_{R_0}}\rho^{-1-2\tau}\xi^2U^2\leq
	C\overline{M}^2(J_0+J_1).	\end{aligned}
\end{equation}
By the latter, taking into account that $\rho\leq 1$ in $B_{R_0}$,
we have

\begin{equation} \label{Lemma-doub:1-7}
	\begin{aligned}
	 &\frac{\tau^3}{2}\int_{B_{R_0/2}\setminus B_{r/2}}\rho^{1-2\tau}U^2+
	\frac{\tau}{2}\int_{B_{R_0/2} \setminus B_{r/2}}\rho^{3-2\tau}|\nabla U|^2+\\
	& +\tau^2r\int_{B_{R_0}}\rho^{-1-2\tau}\xi^2U^2dx\leq
	C\overline{M}^2(J_0+J_1),	
\end{aligned}
\end{equation}
for every $\tau\geq \widetilde{\tau}_2$, where (recall
$\overline{\tau}\geq 1$)
$$\widetilde{\tau}_2=\min\left\{2CM^2, \tau_2\right\}.$$

Now, we estimate from above $J_0$ and $J_1$. By the Caccioppoli inequality (Theorem \ref{dis-Caccioppoli})
and recalling \eqref{peso-osserv}, we have

\begin{equation}
 \begin{aligned}\label{0-9}
J_0&= \int_{B_{r/2}\setminus B_{r/4}}\rho^{4-2\tau}
    \left(r^{-4} U^2+r^{-2}|\nabla U|^2\right)\leq\\&\leq
    C\left(\frac{r}{4}\right)^{-2\tau}\int_{B_{r/2}} \left(U^2+r^{2}|\nabla U|^2\right)\leq C\left(\frac{r}{4}\right)^{-2\tau}\int_{B_{r}}U^2,
\end{aligned}
\end{equation}
where $C$ depends on $M$ only.

Similarly we get

\begin{equation}
    \label{eq:37.1-cube}
     J_1
     \leq C\left(\frac{R_0}{2}\right)^{-2\tau}\int_{B_{R_0}}U^2.
\end{equation}
By \eqref{Lemma-doub:1-7} -- \eqref{eq:37.1-cube}
we have
\begin{equation}\label{Lemma-doub:2-9}
	\begin{aligned}
		&\tau^2r\int_{B_{R_0}}\rho^{-1-2\tau}\xi^2U^2+\tau^3\int_{B_{R_0/2}\setminus B_{r/2}}\rho^{1-2\tau}U^2dx \leq\\& \leq
		C\overline{M}^2\left(\left(\frac{r}{4}\right)^{-2\tau}\int_{B_{r}}U^2+\left(\frac{R_0}{2}\right)^{-2\tau}\int_{B_{R_0}}U^2\right),	
	\end{aligned}
\end{equation}
for every $\tau\geq \widetilde{\tau}_2$.

\medskip

Now, recalling that $2r<R<\frac{R_0}{2}$, by \eqref{eq:25.4-cube} we have

\begin{equation}
\label{Lemma-doub:2n-9}
    \tau^2r\int_{B_{R_0}}\rho^{-1-2\tau}\xi^2U^2\geq (2r)^{-2\tau}\int_{B_{2r}\setminus B_{r/2}} U^2,
\end{equation}
and
\begin{equation}
\label{Lemma-doub:2nn-9}
    \tau^3\int_{B_{R_0/2}\setminus B_{r/2}}\rho^{1-2\tau}U^2\geq
R^{1-2\tau}\int_{B_{R}\setminus B_{r/2}} U^2.
\end{equation}
By \eqref{Lemma-doub:2-9}, \eqref{Lemma-doub:2n-9} and
\eqref{Lemma-doub:2-9} we have
\begin{equation}\label{10.6.9-cube}
	\begin{aligned}
	&(2r)^{-2\tau}\int_{B_{2r}\setminus B_{r/2} }U^2+R^{1-2\tau}\int_{B_{R} \setminus B_{r/2}}
	U^2\leq\\ &
	\leq C \overline{M}^2\left[\left(\frac{r}{4}\right)^{-2\tau}\int_{B_r}U^2+
	\left(\frac{R_0}{2}\right)^{-2\tau}\int_{B_{R_0}}U^2
	\right],	
\end{aligned}
\end{equation}
for every $\tau\geq \widetilde{\tau}$. Now, we add to both the sides of \eqref{10.6.9-cube} the quantity
$$R(2r)^{-2\tau}\int_{B_{r/2} }U^2+(R)^{1-2\tau}\int_{B_{r/2} }U^2,$$ 
 and we find  \eqref{Lemma-doub-cube} for $r<R/2$
and $R<R_0/2$. $\blacksquare$

\bigskip

\textbf{Proof of Theorem \ref{theo:40.teo}.}

Let us suppose $x_0=0$, \eqref{Lemma-doub-cube} gives, for
$0<2r<R<\frac{R_0}{2}$,

\begin{equation}
    \label{three-sphere}
    R^{1-2\tau}\int_{B_{R}}U^2
    \leq C \overline{M}^2\left[\left(\frac{r}{4}\right)^{-2\tau}\int_{B_r}U^2+
    \left(\frac{R_0}{2}\right)^{-2\tau}\int_{B_{R_0}}U^2
    \right],
\end{equation}
for every $\tau\geq \widetilde{\tau}_2$.

Set

$$A(s)=\int_{B_{s}}U^2.$$ By \eqref{three-sphere} we have
\begin{equation}
    \label{three-sphere-1}
    A(R)
    \leq C R^{-1}\overline{M}^2\left[\left(\frac{4R}{r}\right)^{2\tau}A(r)+
    \left(\frac{2R}{R_0}\right)^{2\tau}A(R_0)
    \right],
\end{equation}
for every $\tau\geq \widetilde{\tau}_2$. Let

\begin{equation}
    \label{three-sphere-2}
    \widehat{\tau}=\frac{\log
    \frac{A(R_0)}{A(r)}}{2\log\frac{2R_0}{r}}.
\end{equation}
If
\begin{equation}
    \label{three-sphere-3}
\widehat{\tau}\geq \widetilde{\tau}_2,\end{equation} we choose
$\tau=\widehat{\tau}$ in \eqref{three-sphere-1} and since

$$\left(\frac{4R}{r}\right)^{2\widehat{\tau}}A(r)=
    \left(\frac{2R}{R_0}\right)^{2\widehat{\tau}}A(R_0),$$
we have

\begin{equation}
    \label{three-sphere-4}
    \begin{aligned}
    &A(R)
    \leq C R^{-1}\overline{M}^2\left[\left(\frac{4R}{r}\right)^{2\widehat{\tau}}A(r)+
    \left(\frac{2R}{R_0}\right)^{2\widehat{\tau}}A(R_0)\right]=\\&=
    2CR^{-1}\overline{M}^2\left(\frac{4R}{r}\right)^{2\widehat{\tau}}A(r)=R^{-1}\overline{M}^2\left(A(r)\right)^{\theta}\left(A(R_0)\right)^{1-\theta},
\end{aligned}
\end{equation}
where $\theta$ is given by \eqref{three-sphere-5}. Whereas, if 
\eqref{three-sphere-3} does not hold, then

\begin{equation*}\label{three-sphere-6}
    \frac{\log
    \frac{A(R_0)}{A(r)}}{\log\frac{2R_0}{r}}\leq 2\widetilde{\tau}_2
\end{equation*}
and multiplying both the sides of the last inequality by
$\log\frac{R_0}{2R}$, we have

\begin{equation}\label{three-sphere-7}
\left(A(R_0)\right)^{\theta}\leq
\left(\frac{R_0}{2R}\right)^{2\widetilde{\tau}_2}\left(A(r)\right)^{\theta},
\end{equation}
from which we have trivially

\begin{equation}\label{three-sphere-8}
\begin{aligned}
A(R)&\leq
A(R_0)=\left(A(R_0)\right)^{\theta}\left(A(R_0)\right)^{1-\theta}\leq\\&\leq
\left(\frac{R_0}{2R}\right)^{2\widetilde{\tau}_2}\left(A(r)\right)^{\theta}\left(A(R_0)\right)^{1-\theta},
\end{aligned}
\end{equation}
which, togheter with \eqref{three-sphere-4}, gives
\eqref{three-sphere-enunciato}.

\medskip

Now we prove \eqref{eq:10.6.1102-cube}.

Let us fix $R=\frac{R_0}{4}$ in \eqref{Lemma-doub-cube}. We have
\begin{equation} \label{eq:10.6.952-cube}
	\begin{aligned}
	&\frac{(2r)^{-2\tau}}{4} \int_{B_{2r}}U^2+\left(\frac{R_0}{4}\right)^{1-2\tau}\int_{B_{R_0/4}}
	U^2\leq\\&
	\leq C \overline{M}^2\left[\left(\frac{r}{4}\right)^{-2\tau}\int_{B_r}U^2+
	\left(\frac{R_0}{2}\right)^{-2\tau}\int_{B_{R_0}}U^2
	\right],
	\end{aligned}
\end{equation}
for every $\tau\geq \widetilde{\tau}_2$.

Now, by choosing $\tau=\tau_0$, where

\begin{equation*}
   \label{eq:10.6.1005-cube}
    \tau_0=\widetilde{\tau}+\log_4\left(4C \overline{M}^2N\right)
\end{equation*}
and

\begin{equation}
    \label{eq:10.6.1008-cube}
    N=\frac{\int_{B_{R_0}}U^2}{\int_{B_{R_0/4}}
    U^2},
\end{equation}
we have
$$\left(\frac{R_0}{4}\right)^{1-2\tau_0}\int_{B_{R_0/4}}
    U^2\geq
C\overline{M}^2\left(\frac{R_0}{2}\right)^{-2\tau_0}\int_{B_{R_0}}U^2.$$
Hence, by \eqref{eq:10.6.952-cube}, we obtain
\begin{equation}
    \label{eq:10.6.1014-cube}
    \frac{(2r)^{-2\tau_0}}{4}\int_{B_{2r}}U^2
    \leq C
    \overline{M}^2\left(\frac{r}{4}\right)^{-2\tau_0}\int_{B_r}U^2.
\end{equation}
 By using \eqref{eq:10.6.1008-cube} and
\eqref{eq:10.6.1014-cube}, we have

\begin{equation}
    \label{eq:10.6.1024-cube}
    \int_{B_{2r}}U^2
    \leq CN^3\int_{B_r}U^2,
\end{equation}
where $C$ depends on $M$ only.

The proof is complete. $\blacksquare$

\begin{cor}[\textbf{strong unique continuation for Laplace operator}]\label{SUCP}
	\index{Corollary:@{Corollary:}!- strong unique continuation for Laplace operator@{- strong unique continuation for Laplace operator}}
	
Let $U\in H^2\left(B_{1}\right)$ be a solution to equation
\eqref{2-1-appdoub}. Let $x_0\in B_1$ and $0<R_0\leq 1-|x_0|$.

If $U$ does not vanish identically in $B_{R_0/4}(x_0)$, then we have, for every $r<s\leq\frac{R_0}{16}$,

\begin{equation}\label{SUCP-1}
    \int_{B_{s}(x_0)}U^2\leq
    CN_{x_0,R_0}^3\left(\frac{s}{r}\right)^{\log_2(CN_{x_0,R_0}^3)}\int_{B_r(x_0)}U^2,
\end{equation}
where $N_{x_0,R_0}$ is defined by \eqref{eq:10.6.1108-cube}.

Moreover, if

\begin{equation}  \label{SUCP-2}
\int_{B_r(x_0)}U^2=\mathcal{O}\left(r^m\right),\quad\mbox{as }
r\rightarrow 0,\quad \forall m\in \mathbb{N},
\end{equation}
then

\begin{equation}  \label{SUCP-tesi}
U\equiv 0, \quad\mbox{in } B_1.
\end{equation}
\end{cor}

\textbf{Proof.} We prove \eqref{SUCP-1}. Let us suppose that $x_0=0$ and let
$r<s\leq\frac{R_0}{16}$. Set
$j=\left[\log_2\left(sr^{-1}\right)\right]$ (we recall that
$[a]$ is the integer part of $a$). We have
$$2^jr\leq s <2^{j+1}r$$ and applying repeatedly
\eqref{eq:10.6.1102-cube} we obtain

\begin{equation*}
    \int_{B_{s}}U^2\leq \int_{B_{2^{j+1}r}}U^2
    \leq \left(C N^3\right)^{j+1}\int_{B_r}U^2\leq CN^3\left(\frac{s}{r}\right)^{\log_2(CN^3)}\int_{B_r}U^2.
\end{equation*}
From which we get \eqref{SUCP-1}.

Now, let us suppose that \eqref{SUCP-2} holds true. Hence let us suppose that there
exists a sequence $C_m$ such that

\begin{equation}  \label{SUCP-3}
\int_{B_r}U^2\leq C_mr^m,\quad\mbox{for } r<1,\quad \forall m\in
\mathbb{N}.
\end{equation}

Set

\begin{equation}  \label{SUCP-4}
r_0=\sup\left\{r\in [0,1]:\mbox{ }\int_{B_r}U^2=0\right\}
\end{equation}
(let us note that in \eqref{SUCP-4} the "sup" is, actually, the
maximum).

We distinguish two cases

\medskip

\noindent (i) $r_0=0$,

\medskip

\noindent (ii) $r_0\in (0,1]$.

\bigskip

In case (i) we have

\begin{equation}  \label{SUCP-5}
\int_{B_r}U^2>0,\quad\forall r\in (0,1].
\end{equation}
Hence, setting

$$K=\log_2(CN_{0,1}^3)$$
by \eqref{SUCP-1} and \eqref{SUCP-3}, we get, for $r,s$ such that
$r<s\leq \frac{1}{4}$
\begin{equation}\label{SUCP-6}
    \int_{B_{s}}U^2\leq
    CN_{0,1}^3\left(\frac{s}{r}\right)^{K}\int_{B_r}U^2\leq
    CC_mN_{0,1}^3s^{K}r^{m-K}.
\end{equation}
Let $m>K$. Passing to the limit in \eqref{SUCP-6} as $r$ goes to
 $0$. We obtain
\begin{equation*}
    \int_{B_{s}}U^2=0, \quad \mbox{for } s\leq \frac{1}{4}
\end{equation*}
which contradicts \eqref{SUCP-5}.

Let us consider case (ii).

If $r_0=1$ there is nothing to prove. Let, therefore, $r_0\in
(0,1)$. By the definition of $r_0$ we have

\begin{equation}\label{SUCP-7}
    \int_{B_{r_0}}U^2=0.
\end{equation}
Let $$\delta<\min\left\{r_0,\frac{1-r_0}{15}\right\}$$ and let
$\overline{x}$ be a point such that $|\overline{x}|=r_0-\delta$. Setting
$\overline{R}=1-|\overline{x}|$ we have
$$r_1:=r_0-\delta+\frac{\overline{R}}{16}=r_0+\frac{1-r_0-15\delta}{16}>r_0.$$ Now,
since \eqref{SUCP-7} trivially implies 

\begin{equation*}
    \int_{B_{r}\left(\overline{x}\right)}U^2=\mathcal{O}\left(r^m\right),\quad\mbox{as
    } r\rightarrow 0, \mbox{ } \forall m\in \mathbb{N},
\end{equation*}
repeating the argument of case (i) in the ball
$B_{\overline{R}}\left(\overline{x}\right)$, we reach

\begin{equation*}
    \int_{B_{\overline{R}/16}\left(\overline{x}\right)}U^2=0.
\end{equation*}
Finally, since this equality holds for each $\overline{x}$ such that
$|\overline{x}|=r_0-\delta$, taking into account \eqref{SUCP-7},
we have
\begin{equation*}
    \int_{B_{r_1}}U^2=0,
\end{equation*}
which, as $r_1>r_0$, contradicts the definition of $r_0$ given in \eqref{SUCP-4}.$\blacksquare$

\bigskip

\textbf{Remarks and comments.}

\noindent\textbf{1.} The optimal three sphere inequality can be obtained by the doubling inequality, \eqref{SUCP-1}, in a simple and direct way.
As a matter of fact, by \eqref{SUCP-1}, using the elementary properties of the logarithmic function, we have, for $2r\leq s\leq
\frac{R_0}{16}$,

\begin{equation}\label{SUCP-8}
    \int_{B_{s}(x_0)}U^2\leq
    \left(CN_{x_0,R_0}^3\right)^{2\log_2\frac{s}{r}}\int_{B_r(x_0)}U^2.
\end{equation}
Now, by \eqref{eq:10.6.1108-cube} we have trivially

\begin{equation*}
    N_{x_0,R_0}=\frac{\int_{B_{R_0}(x_0)}U^2}{\int_{B_{R_0/4}(x_0)}U^2}\leq
    \frac{\int_{B_{R_0}(x_0)}U^2}{\int_{B_{s}(x_0)}U^2}.
\end{equation*}
By the latter and by \eqref{SUCP-8} we get

\begin{equation*}
    \left(\int_{B_{s}(x_0)}U^2\right)^{1+6\log\frac{s}{r}}\leq
    \left(C\int_{B_{R_0}(x_0)}U^2\right)^{6\log_2\frac{s}{r}}\int_{B_r(x_0)}U^2,
\end{equation*}
which gives
\begin{equation*}
    \int_{B_{s}(x_0)}U^2\leq
    \left(C\int_{B_{R_0}(x_0)}U^2\right)^{1-\widetilde{\theta}}\left(\int_{B_r(x_0)}U^2\right)^{\widetilde{\theta}}
\end{equation*}
where

\begin{equation*}
   \widetilde{\theta}=\frac{1}{1+6\log_2\frac{s}{r}}.
\end{equation*}
Let us notice that also $\widetilde{\theta}$ is an optimal exponent
in the sense that, for fixed $s$, \eqref{intro14-equ-3} holds.

\medskip

\noindent\textbf{2.} Let us observe that three sphere inequality \eqref{three-sphere-enunciato}  has been proved using Carleman estimate \eqref{Carlm-delta-0}, while to prove the doubling inequality we have used Carleman estimate \eqref{Carlm-delta}, which differs from the estimate \eqref{Carlm-delta-0}
for the occurrence of the term

$$\tau^2r\int\rho^{-1-2\tau}u^2dx.$$

The idea of including this term is indebted to Bakri
\cite{l:Bakri1}, \cite{l:Bakri2} and this idea simplifies the
proof of the doubling inequality with respect to the
previous proofs based on the Carleman estimates. It should also be
pointed out that in the literature there are other methods to prove 
the doubling inequality, see for instance \cite{Ga_Li} ,
\cite{KUK}. We will briefly discuss the main ideas underlying such methods in \ref{Misc:3-11-22-2-0}. $\blacklozenge$

\section{The geodesic polar coordinates}\label{geodesic-polar} In
this Section we give the definitions and the main properties of
geodesic polar coordinates\index{geodesic polar coordinates} introduced by Aronszajn, Krzywicki and
Szarski in \cite{AKS}.

Let $n$ be an integer number,  $n\geq 2$. For any $r>0$ we denote by $\widetilde{B}_{r}$ the set
$B_{r}\setminus \{0\}$.

For any $A=\left\{ a_{ij}\right\} _{i,j=1}^{n}$ real matrix,
we denote by $\left| A\right| $ the norm of $A$, i.e. 
$$\left|
A\right| ^{2}=\sum_{i,j=1}^{n}a_{ij}^{2}.$$

In what follows we will use the Einstein convention of the repeated indices.

By $I_{n}$ we denote the $n\times n$ identity matrix. Given two vectors $x,y\in\mathbb{R}^{n}$, $x=\left(x^1,\cdots,x^n\right)$,
$y=\left(y^1,\cdots,y^n\right)$ we denote by
$$x\cdot y=\delta^{ij}x^iy^j=x^iy^i,$$ their Euclidean scalar product and by $|x|=\sqrt{x\cdot x}$ the Euclidean norm.

Let $\lambda $, $\lambda \geq 1$, $\Lambda $ be positive numbers. Let
$G\left(x\right)=\left\{g_{ij}\left(x\right) \right\}_{i,j=1}^{n}$ be
a nonsingular symmetric real matrix whose entries are functions that belong to
$C^{\infty}\left(\overline{B_{2}}\right)$. Let us denote by
$G^{-1}\left( x\right) =\left\{ g^{ij}\left( x\right)
\right\}_{i,j=1}^{n}$ the inverse of $G(x)$. Let us suppose that

\begin{equation}  \label{2A}
\lambda^{-1}\left|\xi\right|^{2}\leq
G\left(x\right)\xi\cdot\xi\leq\lambda\left|\xi\right|^{2},\quad
\forall\xi\in\mathbb{R}^{n} \mbox{, }\forall x\in B_{2},
\end{equation}
\begin{equation}\label{4A}
\left|\partial_{x^k}G\right|\leq\Lambda, \quad\mbox{for
}k\in\left\{1,\cdots n\right\}, \quad\mbox{in } B_{2}.
\end{equation}

Set $$r(x)=|x|$$ and let us denote by
\begin{equation}\label{8A}
\mu _{0}(x)=G^{-1}(x)\nabla r(x)\cdot\nabla r(x), \mbox{ for }
x\in\widetilde{B}_{2},
\end{equation}
and
\begin{equation}\label{9A}
\widetilde{G}(x)=\mu_{0}(x)G(x), \mbox{ for }x\in\widetilde{B}_{2}.
\end{equation}
Let us denote by $\widetilde{g}_{ij}(x)$, $i$, $j\in\left\{ 1,\cdots,n\right\}$ the entries of the matrix $\widetilde{G}\left( x\right) $.

\bigskip

We wish to introduce the geodesic polar coordinates with respect to the metric tensor  $$\widetilde{g}_{ij}( x) dx^{i}\otimes
dx^{j}.$$ Then we will express in such geodesic polar coordinates the
Laplace--Beltrami operator \index{$\Delta_{g}$}
\begin{equation} \label{3B}
\Delta_{g}(\cdot)=\frac{1}{\sqrt{g(x)}}\partial_{x^{i}}\left(\sqrt{g\left(x\right)}g^{ij}(x)
\partial_{x^{j}}\cdot \right), \mbox{ in }B_{2},
\end{equation}
where $$g(x) =\det G(x).$$

\noindent To perform the above mentioned transformation we assume also
\begin{equation}\label{3A}
G(0) =I_{n},
\end{equation}

\medskip

For any $\overline{x}\in \overline{B}_{1}\setminus\{0\} $
let us denote by $\Gamma(\sigma ;\overline{x})$ the global solution of
the following Cauchy problem
\begin{equation}\label{10Aa}
\begin{cases}
\overset{\cdot }{\Gamma }(\sigma ;\overline{x})=\widetilde{G}^{-1}(\Gamma(\sigma;\overline{x}))\nabla r(\Gamma(\sigma ;\overline{x})),\\
\\
\Gamma(r(\overline{x}),\overline{x})=\overline{x},\end{cases}
\end{equation}
where $\overset{\cdot }{\Gamma}$ is the derivative of $\Gamma(\sigma
;\overline{x})$ w.r.t. $\sigma$.

\medskip

\textbf{Remark.} Let us observe that for any
$\overline{x}\in\overline{B_{1}}\setminus\left\{ 0\right\} $,
$\Gamma$ is a geodesic line w.r.t. the Riemannian metric
$\widetilde{g}_{ij}\left(x\right)dx^{i}\otimes dx^{j}$. Indeed, by
\eqref{8A} and \eqref{9A} we have trivially that function $r(\cdot)$ is a solution to eikonal equation 
$$\widetilde{g}^{ij}\left(
x\right)\partial_{x^{i}}r\partial_{x^{j}}r=1.$$ Hence, denoting by $p(\sigma)=\nabla r(\Gamma(\sigma;\overline{x}))$ we have that, see Section \ref{NonlinEq}, $(\Gamma, p)$ is a solution to Hamilton-Jacobi equations 

\begin{equation*}
\left\{\begin{array}{ll}
\overset{\cdot}{\Gamma}=\nabla_{p}H\left(\Gamma,p\right),\\
\\
\overset{\cdot}{p}=-\nabla_{x} H\left(\Gamma,
p\right),\end{array}\right.
\end{equation*}
where $$H\left(
x,p\right)=\dfrac{1}{2}\widetilde{g}^{ij}\left(x\right)p_{i}p_{j}$$
is the Hamiltonian. Therefore  $\Gamma $ solves the Euler equation, see Section \ref{appendice},
$$\frac{d}{d\sigma } \nabla_{q}L\left(\Gamma,\overset{\cdot}{\Gamma}
\right)=L_{x}\left(\Gamma,\overset{\cdot}{\Gamma}\right),$$ where
$$L\left( x,q\right) =\dfrac{1}{2}\widetilde{g}_{ij}\left( x\right)
q^{i}q^{j}.$$ Hence $\Gamma(\cdot,\overline{x})$ is a geodesic line
w.r.t. the metric
$$\widetilde{g}_{ij}\left(x\right)dx^{i}\otimes dx^{j}.$$ $\blacklozenge$

\bigskip

The following Proposition holds true.

\medskip

\begin{prop}\label{Propos-2}
Let $\Gamma(\cdot ;\overline{x})$ be the global solution to Cauchy problem
\eqref{10Aa}. Then
\begin{equation}\label{30519-1}
\Gamma\left( \cdot ;\overline{x}\right) \mbox{ is defined
in the interval }( 0,2) ,
\end{equation}
and
\begin{equation}\label{10sAb}
r\left(\Gamma\left(\sigma ;\overline{x}\right)\right)=\sigma, \mbox{
for every }\sigma \in \left(0,2\right) .
\end{equation}
\end{prop}

\textbf{Proof.} The proof is the same of that given in the Remark after Theorem \ref{prop5-geod}, we repeat with different notation for the reader's convenience. Let us denote by $J$ the interval on which
$\Gamma\left( \cdot ;\overline{x}\right) $ is defined.

\medskip

\textbf{Claim.}

We have
\begin{equation}\label{11na}
r\left(\Gamma\left( \sigma ;\overline{x}\right) \right) =\sigma,
\quad \forall\sigma \in J.
\end{equation}

\medskip

\textbf{Proof of the Claim}

To prove \eqref{11na} we note that equation (\ref{10Aa})
gives (we omit $\overline{x}$ in $\Gamma$)
\begin{equation*}
\frac{d}{d\sigma
}r\left(\Gamma(\sigma)\right)=\frac{d\Gamma^{i}(\sigma)
}{d\sigma}\partial_{x^{i}}r(\Gamma(\sigma))=\widetilde{g}^{ij}(
\Gamma(\sigma))
\partial_{x^{i}}r(\Gamma(\sigma))r\partial_{x^{j}}r(\Gamma(\sigma))=1,
\end{equation*}
for every $\sigma\in J$. Therefore, there exists a constant $c$, such that
\begin{equation*}
r\left(\Gamma\left(\sigma \right)\right)=\sigma+c,\quad\forall \sigma\in J.
\end{equation*}
By initial condition we have
$$\Gamma(r(\overline{x}),\overline{x})=\overline{x},$$ hence
$$r\left(\overline{x}\right)
=r(\Gamma(r(\overline{x}),\overline{x}))=r\left(\overline{x}\right)
+c,$$ consequently $c=0$ which implies \eqref{11na}. Claim is proved.

\medskip 

By \eqref{11na} and by standard results of general theory of ordinary differential equations we have that $\Gamma$ can be defined in the whole interval $(0,2)$
hence \eqref{30519-1} is proved.  $\blacksquare$

\bigskip

 In order to introduce the geodesic polar coordinates we need some additional notations. Set
$\Sigma =\partial B_{1}$. Let $\left\{ U_{\alpha },\varphi _{\alpha
}\right\} _{\alpha \in \mathcal{J} }$ be a finite family of local maps
which define an oriented $C^{\infty }$ differentiable structure on
$\Sigma $. For any $\alpha \in \mathcal{J} $, set $V_{\alpha
}=\varphi _{\alpha }\left(U_{\alpha }\right) $. Let us denote by
$\Phi$ the map
\begin{equation}\label{11s1A}
\Phi:B_{1}\setminus\left\{ 0\right\} \rightarrow \left( 0,1\right)
\times \Sigma,
\end{equation}
such that
\begin{equation}\label{11s2A}
\Phi\left(x\right) =\left(\left|x\right|
,\Gamma\left(1;x\right)\right), \quad\forall x\in
B_{1}\setminus\left\{0\right\} .
\end{equation}
By \eqref{10Aa} we have easily that $\Phi$ is bijective,
moreover
\begin{equation}\label{11s3A}
\Phi^{-1}\left(
\varrho,p\right)=\Gamma\left(\varrho;p\right),\quad\forall\left( \varrho,p\right) \in\left( 0,1\right)\times\Sigma.
\end{equation}
For any
$\alpha\in\mathcal{J}$, let us consider the following geodesic sector
\begin{equation}
\mathcal{I}\left( U_{\alpha }\right)
=\left\{\Gamma\left(\varrho;p\right) : \
\varrho \in \left( 0,1\right) \text{, }p\in U_{\alpha }\right\} ,
\label{12A}
\end{equation}
let us denote by $\Phi_{\alpha }$ the map $\Phi$ expressed in the local coordinates, i.e.
\begin{equation}\label{13A}
\Phi_{\alpha } :\mathcal{I}\left(U_{\alpha }\right)\rightarrow
\left( 0,1\right) \times V_{\alpha},
\end{equation}
\begin{equation}\label{14A}
\Phi _{\alpha } \left( x\right) =\left( \left| x\right| ,\varphi
_{\alpha }\left(\Gamma\left( 1;x\right) \right) \right),\quad
\forall x\in \mathcal{I}\left(U_{\alpha}\right).
\end{equation}
For any $\alpha \in \mathcal{J}$ let us denote by

\begin{equation}\label{16A}
\Gamma_{\alpha }\left(\varrho,\theta\right) =\Gamma\left(
\varrho,\varphi _{\alpha }^{-1}\left(\theta\right)\right),
\quad\forall\left( \varrho,\theta\right) \in \left(
0,1\right) \times V_{\alpha } .
\end{equation}
We have
\begin{equation}\label{15A}
\Phi_{\alpha }^{-1}\left( \varrho,\theta \right) =\Gamma_{\alpha
}\left(\varrho,\theta\right), \quad\forall\left(\varrho,\theta \right) \in (0,1) \times V_{\alpha }.
\end{equation}
Let us note that, for any $\alpha ,\alpha ^{\prime }\in \mathcal{J}$,
we have
\begin{equation*}
 (\Phi _{\alpha }\circ \Phi
_{\alpha ^{\prime }}^{-1})
\left(\varrho,\theta\right)=\left(\varrho,\left( \varphi _{\alpha
}\circ \varphi _{\alpha ^{\prime }}^{-1}\right)
(\theta)\right),\quad \forall
\left(\varrho,\theta\right)\in\Phi_{\alpha}\left(\mathcal{I}\left(U_{\alpha}\cap
U_{\alpha^{\prime}}\right)\right).
\end{equation*}
Let us note that $ \mathcal{I}( U_{\alpha })\cap \mathcal{I}(U_{\alpha
^{\prime }})=\mathcal{I}\left( U_{\alpha }\cap U_{\alpha ^{\prime
}}\right)$, hence $\left\{\mathcal{I}\left(U_{\alpha }\right)
,\Phi_{\alpha }\right\}_{\alpha \in \mathcal{J}}$ defines an oriented $C^{\infty }$ differentiable structure on $B_{1}\setminus\{0\}$.

\medskip

\noindent Let us observe that by \eqref{10Aa} and \eqref{16A} we have  
\begin{equation}\label{16sAa}
	\begin{aligned}
		\partial_{\varrho}\Gamma_{\alpha}(\varrho,\theta)&=\widetilde{G}^{-1}(\Gamma_{\alpha}(\varrho,\theta))\nabla r(\Gamma_{\alpha}(\varrho,\theta))=\\&=\frac{1}{\varrho}\widetilde{G}^{-1}(\Gamma_{\alpha}(\varrho,\theta))\Gamma_{\alpha}(\varrho,\theta)
		\end{aligned}
\end{equation}
and
\begin{equation}\label{16sAb}
\Gamma_{\alpha }\left( 1,\theta\right)=\varphi_{\alpha }^{-1}(
\theta).
\end{equation}
Moreover, by \eqref{10sAb} we have
\begin{equation}\label{16sAc}
\left| \Gamma_{\alpha }(\varrho,\theta)\right| =\varrho.
\end{equation}

\bigskip

To save the sum index convention, 
in the next Proposition and in the sequel, we denote by $\theta
^{j+1}$ the $j-$th component of $\theta$, so $\theta
=\left(\theta ^{2},\cdots,\theta ^{n}\right)$.

Let $\alpha \in
\mathcal{J}$ be fixed. Let $\eta =\left\{\eta ^{h}\right\}
_{h=2}^{n}$ a vector, set
\begin{equation}\label{39A}
y_{\eta }\left( \varrho,\theta\right)
=\frac{1}{\varrho}\partial_{\theta}\Gamma_{\alpha}(\varrho,\theta)
\eta,
\end{equation}
where $\partial_{\theta}(\cdot)$ denotes the  jacobian matrix w.r.t. the variables $\theta_2,\cdots,\theta_n $. Let us check that by \eqref{16sAa} and
\eqref{16sAb} we obtain, respectively,
\begin{equation}\label{39sAa}
\partial_{\varrho }y_{\eta }=\frac{1}{\varrho}\left(\widetilde{G}
^{-1}\left(\Gamma_{\alpha}\right)-I_{n}\right) y_{\eta
}+\frac{1}{\varrho} (\partial_{x^{k}}\widetilde{G}^{-1}(
\Gamma_{\alpha })) y_{\eta }^{k}\Gamma_{\alpha},
\end{equation}
and
\begin{equation}\label{39sAb}
y_{\eta }(1,\theta)=\partial_{\theta} \varphi_{\alpha
}^{-1}(\theta)\eta.
\end{equation}
Equality \eqref{39sAb} is an immediate consequence of \eqref{16sAb}. Concerning \eqref{39sAa}, first we set
$$\widetilde{y}_h=\frac{1}{\varrho}\partial_{\theta^h}\Gamma_{\alpha}(\varrho,\theta),\quad h=2,\cdots, n,$$ and we have, for $h=2,\cdots, n$,
\begin{equation*}
	\begin{aligned}
		\partial_{\varrho}\widetilde{y}_h&=\partial_{\varrho}\left(\frac{1}{\varrho}\partial_{\theta^h}\Gamma_{\alpha}(\varrho,\theta)
		\right)=\\&=-\frac{1}{\varrho^2}\partial_{\theta^h}\Gamma_{\alpha}(\varrho,\theta)
		+\frac{1}{\varrho}\partial^2_{\varrho\theta^h}\Gamma_{\alpha}(\varrho,\theta)
		=\\&=-\frac{1}{\varrho}\widetilde{y}_h+\frac{1}{\varrho}\partial_{\theta^h}\left(\partial_{\varrho}\Gamma_{\alpha}(\varrho,\theta)\right)=\\&=-\frac{1}{\varrho}\widetilde{y}_h+\frac{1}{\varrho^2}\partial_{\theta^h}\left(\widetilde{G}^{-1}(\Gamma_{\alpha}(\varrho,\theta))\Gamma_{\alpha}(\varrho,\theta)\right)=\\&
		=-\frac{1}{\varrho}\widetilde{y}_h+\frac{1}{\varrho^2}\partial_{\theta^h}\left(\widetilde{G}^{-1}(\Gamma_{\alpha}(\varrho,\theta))\right)\Gamma_{\alpha}(\varrho,\theta)+\\&+\frac{1}{\varrho^2}\widetilde{G}^{-1}(\Gamma_{\alpha}(\varrho,\theta))\partial_{\theta^h}\Gamma_{\alpha}(\varrho,\theta)=\\&=-\frac{1}{\varrho}\widetilde{y}_h+\frac{1}{\varrho^2}\left(\partial_{x^k}\widetilde{G}^{-1}(\Gamma_{\alpha}(\varrho,\theta))\partial_{\theta^h}\Gamma^k_{\alpha}(\varrho,\theta)\right)\Gamma_{\alpha}(\varrho,\theta)+\\&+\frac{1}{\varrho}\widetilde{G}^{-1}(\Gamma_{\alpha}(\varrho,\theta))\widetilde{y}_h.
		\end{aligned}.
\end{equation*}
Hence, for $h=2,\cdots, n$, we have  

\begin{equation*}
	\begin{aligned}
		\partial_{\varrho}\widetilde{y}_h=&\frac{1}{\varrho}\left(\widetilde{G}^{-1}(\Gamma_{\alpha}(\varrho,\theta))-I_n\right)\widetilde{y}_h+\frac{1}{\varrho^2}\left(\partial_{x^k}\widetilde{G}^{-1}(\Gamma_{\alpha}(\varrho,\theta))\partial_{\theta^h}\Gamma^k_{\alpha}(\varrho,\theta)\right)\Gamma_{\alpha}(\varrho,\theta).
	\end{aligned}
\end{equation*}
By multlying both the sides of the last equality by $\eta^h$ and adding up the index $h$, we obtain \eqref{39sAa}. $\blacksquare$

\bigskip

\begin{lem}\label{ar-lemma}
Let us assume that \eqref{2A}, \eqref{4A} and \eqref{3A} hold true. Let 
$\Gamma_{\alpha}$, $\alpha\in\mathcal{J}$ be defined by \eqref{16A}.
We have
\begin{equation}\label{40A}
C^{-1}\left| \partial_{\theta}\varphi _{\alpha }^{-1}(\theta)\eta
\right|\leq\left|y_{\eta}(\varrho,\theta)\right|\leq C\left|
\partial_{\theta}\varphi_{\alpha }^{-1}(\theta)\eta\right|,
\end{equation}
for every $\left(\varrho,\theta\right)\in(0,1) \times V_{\alpha }$ and
for every $\eta \in \mathbb{R}^{n-1}$, where $C$ and $C\geq 1,$ depend on
$\lambda $ and $\Lambda $ only (here and in the sequel we omit the dependence on
$n$).
\end{lem}

\textbf{Proof.} Let us omit the index $\alpha $. By \eqref{2A}, \eqref{9A}, \eqref{3A},
and \eqref{16sAc} we have
\begin{equation*}
\left| \partial_{x^{k}}\widetilde{G}^{-1}(\Gamma)y_{\eta
}^{k}\right| \leq C\left|y_{\eta}\right| ,
\end{equation*}
\begin{equation*}
\left| \widetilde{G}^{-1}(\Gamma)-I_{n}\right| \leq C\varrho,
\end{equation*}
where $C$ depends on $\lambda$ and $\Lambda$ only. Therefore, by
\eqref{39sAa} we have
\begin{equation}\label{X1}
\left|\partial_{\varrho}y_{\eta}\right|\leq C\left|y_{\eta }\right|,
\end{equation}
where $C$ depends on $\lambda$ and $\Lambda$ only.

\noindent By \eqref{39sAb} and \eqref{X1} we have
\begin{equation}\label{44A}
\left|y_{\eta}\left( \varrho,\theta\right) \right|
\leq\left|\partial_{\theta}\varphi^{-1}( \eta \right|
+C\int\nolimits_{\varrho}^{1}\left| y_{\eta }\left( s,\theta
\right) \right| ds,\quad\forall\varrho\in(0,1],
\end{equation}
By \eqref{44A} and by the Gronwall inequality we get the
second inequality of \eqref{40A}.

Now we prove the first inequality of \eqref{40A}. Let
$\overline{\varrho}$  be fixed in $(0,1)$. We have
\begin{equation*}
\left|y_{\eta}(\varrho,\theta)\right|\leq\left|
y_{\eta}(\overline{\varrho},\theta)\right|+\int_{\overline{\varrho}}^{\varrho}\left|
\partial_s y_{\eta}(s,\theta)\right| ds,\quad \forall\varrho\in
[\overline{\varrho},1] ,
\end{equation*}
hence by \eqref{X1} we get
\begin{equation*}
\left|y_{\eta}(\varrho,\theta)\right|\leq\left|y_{\eta}\left(\overline{\varrho},\theta\right)\right|+C\int_{\overline{\varrho
}}^{\varrho}\left|y_{\eta }(s,\theta)\right|ds,\quad\forall\varrho\in\left[\overline{\varrho},1\right].
\end{equation*}
By the Gronwall inequality we have 
\begin{equation}\label{45A}
\left|y_{\eta}\left(\varrho,\theta\right)\right|\leq\left|y_{\eta}\left(
\overline{\varrho},\theta\right)\right|e^C,\quad\forall\varrho\in \left[\overline{\varrho},1\right] ,
\end{equation}
where $C$ depends on $\lambda $ and $\Lambda $ only. By \eqref{39sAb} and
\eqref{45A} we have

\begin{equation*}
\left|\partial_{\theta}\varphi^{-1}(\theta)\eta\right|e^{-C}\leq
\left|y_{\eta }\left( \overline{\varrho},\theta\right) \right|,
\end{equation*}
that gives the first inequality of \eqref{40A}. $\blacksquare$

\bigskip

The following Proposition holds true.

\medskip

\begin{prop}\label{Propos-4}
For any $\alpha\in\mathcal{J}$ let us denote by $\widetilde{b}_{\alpha
,hk}$, $h,$ $k\in \left\{1,\cdots, n\right\} $, the components of the metric 
tensor $\widetilde{g}_{ij}\left( x\right) dx^{i}\otimes
dx^{j}$ with respect to the local coordinates $\left(\mathcal{I}\left(
U_{\alpha }\right) ,\Phi_{\alpha } \right).$ We have
\begin{equation}\label{17A}
\widetilde{b}_{\alpha ,hk}\left( \varrho,\theta\right)
=\widetilde{g}_{ij}\left(\Gamma_{\alpha }\left(
\varrho,\theta\right)\right) \partial_{\theta ^{h}} \Gamma_{\alpha
}^{i} \partial_{\theta ^{k}}\Gamma_{\alpha }^{j}, \quad\mbox{for }
h,k\in \left\{2,\cdots, n\right\},
\end{equation}
\begin{equation}\label{18A}
\widetilde{b}_{\alpha ,h1}\left( \varrho,\theta\right)
=\widetilde{b} _{\alpha ,1h}\left(
\varrho,\theta\right)=0,\quad\mbox{for } h\in \left\{2,\cdots,
n\right\}
\end{equation}
and 
\begin{equation}\label{19A}
\widetilde{b}_{\alpha ,11}\left( \varrho,\theta\right)=1.
\end{equation}
\end{prop}

\textbf{Proof.} To simplify the notazion, in what follows we omit the index
$\alpha$. Equality \eqref{17A} are nothing but the rule of 
transformation of the components of the metric tensor.

Let us prove \eqref{18A}.

By \eqref{16sAa} and \eqref{16sAc} we get, for every $h\in
\left\{2,\cdots n\right\}$,

\begin{equation*}
\begin{aligned}
\widetilde{b}_{h1}&=\widetilde{g}_{ij}(\Gamma)\partial_{\theta^{h}}\Gamma^{i}\partial_{\varrho}\Gamma^{j}=\widetilde{g}_{ij}
(\Gamma)\partial_{\theta^{h}}\Gamma^{i}
\widetilde{g}^{jk}(\Gamma)\partial_{x^{k}}r(\Gamma)=\\&
=\delta_{i}^{k}\partial_{\theta^{h}}\Gamma^{i}\partial_{x^{k}}r(\Gamma)=\partial_{\theta^{k}}\left(r(\Gamma)\right)=\partial_{\theta^{k}}\varrho=0,
\end{aligned}
\end{equation*}
since $\widetilde{b}_{ij}$  is symmetric, we obtain
\eqref{18A}.

Now, let us prove \eqref{19A}. We have
\begin{equation*}
\begin{aligned}
\widetilde{b}_{11}&=\widetilde{g}_{ij}(\Gamma)\partial_{\varrho}\Gamma^{i}\partial_{\varrho}\Gamma^{j}=\widetilde{g}
_{ij}(\Gamma)\widetilde{g}^{ik}(\Gamma) \partial_{x^{k}}r(\Gamma)
\widetilde{g}^{jh}(\Gamma)\partial_{x^{h}}r(\Gamma)= \\&
=\delta_{j}^{k}\widetilde{g}^{jh}(\Gamma) \partial_{x^{k}}r(\Gamma)
\partial_{x^{h}}r(\Gamma) =\widetilde{g}^{kh}(\Gamma)
\partial_{x^{k}}r(\Gamma)\partial_{x^{h}}r(\Gamma)=1.
\end{aligned}
\end{equation*}
$\blacksquare$

\bigskip

In formulas \eqref{24A}--\eqref{31A} below, we introduce some notations.

Set, for any $\alpha\in\mathcal{J}$,
\begin{equation}\label{24A}
\mu_{\alpha}=\mu_{0}\circ \Phi_{\alpha }^{-1},
\end{equation}
\begin{equation}\label{27A}
\widetilde{b}_{\alpha }=\det\left\{\widetilde{b}_{\alpha
,ij}\right\}_{i,j=1}^{n},
\end{equation}
\begin{equation}\label{28A}
\beta_{\alpha ,hk}=\frac{1}{\varrho^{2}}\widetilde{b}_{\alpha
,hk},\mbox{ for }h,k\in \left\{ 2,\cdots, n\right\} ,
\end{equation}
\begin{equation}\label{28Ab}
\left\{\beta_{\alpha
}^{hk}\right\}_{h,k=2}^{n}=\left(\left\{\beta_{\alpha
,hk}\right\}_{h,k=2}^{n}\right) ^{-1},
\end{equation}
\begin{equation}\label{29A}
\beta_{\alpha }=\det\left\{\beta_{\alpha,hk}\right\}_{h,k=2}^{n}.
\end{equation}
In addition, let
\begin{equation}\label{31A}
\beta _{\alpha }=\varrho^{-2\left( n-1\right) }\widetilde{b}_{\alpha
}.
\end{equation}

\bigskip

The following Proposition holds true.

\medskip

\begin{prop}\label{Propos-5}
For every $\overline{\varrho}\in (0,1)$,
$\beta_{\alpha,hk}\left(\overline{\varrho},\theta\right) ,$ $h,k\in
\left\{ 2,\cdots,n\right\}$, are the components of a metric tensor
on $\Sigma $ with respect to the local maps $\left(U_{\alpha
},\varphi_{\alpha}\right)$.
\end{prop}
\textbf{Proof.} Let $\overline{p}\in \Sigma $ and let $\left( U_{\alpha },\varphi
_{\alpha }\right) $, $\left( U_{\alpha ^{\prime }},\varphi_{\alpha
^{\prime }}\right) $ be two coordinate neighborhoods such that $\overline{p}\in
U_{\alpha }\cap U_{\alpha ^{\prime }}$. Let $p$ be an arbitrary point of
 $U_{\alpha }\cap U_{\alpha^{\prime }}$. Set
$$\theta^{(p)}=\varphi_{\alpha}(p),\quad\quad\widehat{\theta}^{(p)}=\varphi_{\alpha^{\prime}}( p).$$ We have trivially,
$$p=\varphi_{\alpha}^{-1}\left(\theta^{\left(p\right)}\right)=\varphi_{\alpha^{\prime}}^{-1}\left(
\widehat{\theta}^{(p) }\right).$$ Moreover, since
$\Gamma_{\alpha}\left(\cdot,\theta^{(p)}\right)$ and
$\Gamma_{\alpha^{\prime}}\left(\cdot,\widehat{\theta}^{\left(p\right)}\right)$
are solutions to equation \eqref{16sAa} and
$$\Gamma_{\alpha}\left(1,\theta^{(p) }\right)
=\Gamma_{\alpha^{\prime}}\left(1,\widehat{\theta}^{(p)}\right)=p,$$
we have
$$\Gamma_{\alpha}\left(\cdot,\theta^{(p)}\right)=\Gamma_{\alpha^{\prime}}\left(\cdot,\widehat{\theta}^{(p)}\right).$$
Therefore
$$\Gamma_{\alpha^{\prime}}\left(\overline{\varrho},\widehat{\theta
}^{(p)}\right)
=\Gamma_{\alpha}\left(\overline{\varrho},\left(\varphi_{\alpha}\circ\varphi_{\alpha^{\prime}}^{-1}\right)
\left(\widehat{\theta}^{(p)}\right) \right).$$ Hence, if
$$\theta\in\varphi_{\alpha^{\prime}}\left(U_{\alpha }\cap
U_{\alpha^{\prime}}\right),$$ then
$$\Gamma_{\alpha^{\prime}}\left(\overline{\varrho},\theta\right)
=\Gamma_{\alpha}\left(\overline{\varrho}
,\left(\varphi_{\alpha}\circ\varphi_{\alpha^{\prime}}^{-1}\right)
\left(\theta\right)\right).$$ Differentiating w.r.t. $\theta^{l}$ both the sides of the last equality, we obtain
\begin{equation}\label{35A}
\partial_{\theta^{l}}\Gamma_{\alpha^{\prime}}\left(\overline{\varrho},\theta\right)=\left(\partial_{\theta^k}\Gamma_{\alpha}\right)
\left(\overline{\varrho},
\left(\varphi_{\alpha}\circ\varphi_{\alpha^{\prime}}^{-1}\right)\left(\theta\right)\right)\partial_{\theta
^{l}}\left(\varphi_{\alpha}\circ\varphi_{\alpha
^{\prime}}^{-1}\right)^{k}(\theta),
\end{equation}
for every $l\in \left\{ 2,\cdots, n\right\} $, $\theta \in
\varphi_{\alpha ^{\prime}}\left(U_{\alpha}\cap
U_{\alpha^{\prime}}\right)$.

Now
\begin{equation}\label{35Ab}
\beta_{\alpha^{\prime},lm}\left(\overline{\varrho},\theta\right) =
\overline{\varrho}^{-2}\widetilde{g}_{ij}\left(\Gamma_{\alpha^{\prime}}\left(\overline{\varrho},\theta\right)\right)
\partial_{\theta^{l}}\Gamma_{\alpha^{\prime}}^i \partial_{\theta^{m}}\Gamma_{\alpha^{\prime}}^j,
\end{equation}
for $l,m\in \left\{ 2,\cdots,n\right\} $, $\theta \in \varphi
_{\alpha ^{\prime }}\left( U_{\alpha }\cap U_{\alpha ^{\prime
}}\right)$.

 Therefore by \eqref{35A} and \eqref{35Ab} we obtain
\begin{equation}\label{ar-extra1}
	\begin{aligned}
		&\beta_{\alpha^{\prime},lm}\left(\overline{\varrho},\theta\right)=\\&=\beta_{\alpha^{\prime},ks}\left(\overline{\varrho},
\left(\varphi_{\alpha}\circ\varphi_{\alpha
^{\prime}}^{-1}\right)(\theta)\right)\partial_{\theta^{l}}\left(\varphi_{\alpha}\circ\varphi_{\alpha^{\prime}}^{-1}\right)^{k}(\theta)
\partial_{\theta^{m}}\left(\varphi_{\alpha}\circ\varphi_{\alpha^{\prime}}^{-1}\right)^{s}(\theta),
\end{aligned}
\end{equation}
for every $l,m\in \left\{ 2,\cdots,n\right\} $ and for every
$\theta\in\varphi_{\alpha^{\prime}}\left(U_{\alpha}\cap
U_{\alpha^{\prime}}\right)$.

Equality \eqref{ar-extra1} proves that
$\beta_{\alpha,hk}\left(\overline{\varrho},\theta\right) ,$ $h,k\in
\left\{ 2,\cdots,n\right\}$, are the components of a tensor which is
a metric tensor because the matrix
$\left\{\beta_{\alpha,hk}\left(\overline{\varrho},\theta
\right)\right\}_{h,k=2}^{n}$ is symmetric and positive. The proof is complete. $\blacksquare$

\bigskip

Now we begin to derive the expression of operator \eqref{3B}
in the polar coordinates introduced above.

Let $u\in $C$^{\infty }\left( B_{2}\right) $ and let us denote by $w$ the function
\begin{equation}\label{def-w}
(0,1) \times\Sigma\ni\left(\varrho,p\right) \rightarrow
w\left(\varrho,p\right)=u\left(\Gamma\left(\varrho,p\right)\right).
\end{equation}
Set
\begin{equation}\label{4B}
w_{\alpha }=u\circ \Phi _{\alpha }^{-1},
\end{equation}
by \eqref{15A}  we have
\begin{equation}\label{5B}
w\left( \varrho,\varphi _{\alpha }^{-1}\left( \theta \right)\right)
=w_{\alpha }\left( \varrho,\theta\right) , \quad\forall
\left(\varrho,\theta\right) \in \left(0,1\right)\times V_{\alpha},
\end{equation}
hence, for any fixed $\varrho\in \left( 0,1\right)$, $w_{\alpha
}\left(\varrho,\cdot\right)$ is the expression of
$w\left(\varrho,\cdot\right) $ in the local coordinates $\left(
U_{\alpha },\varphi_{\alpha }\right)$.

Now, for any fixed  $\overline{\varrho}\in(0,1)$, Proposition
\ref{Propos-5} allows us to define on $\Sigma $ the Riemannian structure, \cite{BOO}, induced by metric tensor whose components with respect coordinate neighborhood $\left(U_{\alpha },\varphi
_{\alpha}\right)$ are equals to $\beta _{\alpha
,hk}\left(\overline{\varrho},\theta\right) $ for $h,k\in
\left\{2,\cdots,n\right\}$.

Let us denote by $\left\langle\cdot,\cdot\right\rangle$ and
$\left|_{\cdot }\right|_{\Sigma }$, respectively, the inner product and the associated norm on the above defined Riemannian structure.
Let us denote by $\nabla_{\Sigma}$ and div$_{\Sigma }$ the gradient and the divergence operators on $\Sigma $ respectively.

We have
\begin{equation}\label{6B}
\left( \nabla_{\Sigma}w\left( \overline{\varrho},\cdot\right)
\right)\circ \varphi_{\alpha }^{-1}\left(\theta
\right)=\left\{\beta_{\alpha }^{hk}\left(\overline{\varrho},\theta
\right)
\partial_{\theta ^{k}} w_{\alpha }\left( \overline{\varrho},\theta\right) \right\} _{h=2}^{n},\quad\forall\theta \in V_{\alpha
}.
\end{equation}
Set
\begin{equation*}
\mu
\left(\varrho,p\right)=\mu_{0}\left(\Gamma\left(\varrho,p\right)\right),\quad
\forall\left(\varrho,\theta\right) \in (0,1)\times\Sigma,
\end{equation*}
we have
\begin{equation}\label{7B}
	\begin{aligned}
	&\mbox{div}_{\Sigma }\left( \mu ^{1-\frac{n}{2}}\left(
	\overline{\varrho},\cdot , \right)\nabla_{\Sigma
	}w(\overline{\varrho},\cdot)\right)\circ\varphi_{\alpha
	}^{-1}\left(\theta \right) =\\&
	=\frac{1}{\sqrt{\beta_{\alpha }\left(\overline{\varrho},\theta
			\right) }}\partial_{\theta^{h}}\left( \mu _{\alpha }^{1-
		\frac{n}{2}}\left( \overline{\varrho},\theta\right) \sqrt{\beta
		_{\alpha }\left( \overline{\varrho},\theta\right)}\beta_{\alpha
	}^{hk}\left(\overline{\varrho},\theta\right) \partial_{\theta
		^{k}}w_{\alpha}\left(\overline{\varrho},\theta\right)\right),
	\end{aligned}
\end{equation}
for every $\theta\in V_{\alpha }.$

Let us note that Proposition \eqref{Propos-4} implies that the derivatives 
$$\partial_{\varrho}\log\sqrt{\beta_{\alpha}(\overline{\varrho},\theta)}, \ \ \mbox{for } \alpha\in \mathcal{J},$$ are the expression in the coordinate neighborhoods of
a  $C^{\infty}$ function on $\Sigma$.

Let us observe that by the equality  $$\Gamma_{\alpha}\left(
\varrho,\theta\right) =\Gamma_{\alpha^{\prime
}}\left(\varrho,\left(\varphi_{\alpha }\circ\varphi_{\alpha ^{\prime
}}^{-1}\right) \left(\theta \right)\right) \ \ \forall\theta\in
\varphi_{\alpha^{\prime}}\left(U_{\alpha }\cap U_{\alpha ^{\prime
}}\right),$$ we get
\begin{equation*}
\partial_{\theta^{h}}\Gamma_{\alpha}\left( \varrho,\theta\right)=\partial_{\theta^{k}}\Gamma_{\alpha^{\prime}}\left(\varrho,
\left(\varphi_{\alpha}\circ\varphi_{\alpha^{\prime}}^{-1}\right)(\theta)\right)\partial_{\theta^{h}}
\left(\varphi_{\alpha}\circ\varphi_{\alpha^{\prime}}^{-1}\right)^k(\theta).
\end{equation*}

Let us denote, respectively, by $\mathcal{L}_{1}$,
$\mathcal{L}_{2}$ and $\mathcal{L}$ the operators
\begin{equation}\label{81B}
\mathcal{L}_{1}w=\mu \left( \frac{\partial ^{2}w}{\partial
\varrho^{2}}+\frac{ n-1}{\varrho}\frac{\partial w}{\partial
\varrho}+\frac{1}{\varrho^{2}\mu ^{1-n/2}} div_{\Sigma }\left(\mu
^{1-n/2}\nabla _{\Sigma }w\right) \right),
\end{equation}
\begin{equation}\label{91B}
\mathcal{L}_{2}w=\mu \frac{\partial }{\partial \varrho}\left( \log \left(\mu ^{1-n/2}%
\sqrt{\beta }\right)\right)\partial_{\varrho}w
\end{equation}
and
\begin{equation}\label{92B}
\mathcal{L}=\mathcal{L}_{1}+\mathcal{L}_{2}.
\end{equation}

We have the following

\begin{prop}\label{Propos-6}
\textbf{(Tranformation of the operator }$\Delta_{g}$\textbf{)} The following equality holds true
\begin{equation}\label{101B}
\left(\left(\Delta_{g}u\right)
\circ\Phi_{\alpha}^{-1}\right)\left(\varrho,\theta\right)=\left(\mathcal{L}w\right)\left(\varrho,\varphi_{\alpha}^{-1}\left(\theta\right)\right),
\quad\forall(\varrho,\theta) \in (0,1)\times V_{\alpha}.
\end{equation}
\end{prop}

\textbf{Proof.} Let us denote by $\left\{ b_{\alpha }^{ij}\right\} _{i,j=1}^{n}$ the inverse matrix
of $\left\{\mu_{\alpha }^{-1}\widetilde{b}_{\alpha
,ij}\right\}_{i,j=1}^{n}$ (recall that
$\left\{\widetilde{b}_{\alpha ,ij}\right\} _{i,j=1}^{n}$ is
defined in Proposition \eqref{Propos-4}). Let us recall that
$w_{\alpha }=u\circ \Phi _{\alpha }^{-1}$. By \eqref{9A} and
\eqref{18A} we have

\begin{equation}\label{102B}
	\begin{aligned}
	&\left(\Delta_g
u\right)\left(\Phi_{\alpha}^{-1}(\varrho,\theta)\right)=\\&=\frac{1}{\sqrt{b_{\alpha
}}}\left(\partial_{\varrho}\left( \sqrt{b_{\alpha }}\mu _{\alpha
}\partial_{\varrho} w_{\alpha }\right) + \partial_{\theta
^{h}}\left(\sqrt{b_{\alpha }} b_{\alpha
}^{hk}\partial_{\theta^{k}}w_{\alpha }\right)\right),
	\end{aligned}
\end{equation}
where $b_{\alpha }=\det \left\{\mu_{\alpha
}^{-1}\widetilde{b}_{\alpha ,ij}\right\}_{i,j=1}^{n}$.

We have
\begin{equation}\label{103B}
b_{\alpha }=\mu _{\alpha }^{-n}\varrho^{2\left( n-1\right)
}\beta_{\alpha},
\end{equation}
\begin{equation}\label{104B}
b_{\alpha }^{hk}=\varrho^{-2}\mu _{\alpha }\beta _{\alpha }^{hk},\quad\mbox{for }%
h\text{, }k\in \left\{ 2,\ldots, n\right\} .
\end{equation}
By \eqref{7B}, \eqref{102B}, \eqref{103B} and \eqref {104B} we get
\begin{equation*}
	\begin{aligned}
	&\left( \Delta_g u\right) \circ \Phi _{\alpha
	}^{-1}=\\&=\left(\mu\left(\partial_{\varrho}^{2}w
	+\frac{n-1}{\varrho}\partial_{\varrho}w +\frac{1}{\varrho^{2}\mu
		^{1-n/2}}\mbox{div}_{\Sigma }\mu
	^{1-n/2}\nabla_{\Sigma}w\right)\right)\circ \varphi _{\alpha }^{-1}
	\\&
	+\left(\mu \partial_{\varrho}\left(\log
	\mu^{1-n/2}\sqrt{\beta}\right)\partial_{\varrho}
	w\right)\circ\varphi_{\alpha }^{-1}.
	\end{aligned}
\end{equation*}
$\blacksquare$

\bigskip

In the next propositions we will estimate the tensors which occur in the transformed operator $\mathcal{L}$.

\begin{prop}\label{Propos-7}
Let $\mu _{\alpha }$, $\left\{ \beta _{\alpha ,hk}\right\}
_{h,k=2}^{n}$, $\beta_{\alpha}$ be defined by \eqref{24A},
\eqref{28A}, \eqref {31A}. For any $\alpha \in \Gamma $ and
$\left(\varrho,\theta\right)\in(0,1)\times V_{\alpha}$, we have
\begin{equation}\label{mu}
\lambda^{-1}\leq \mu_{\alpha}(\varrho,\theta)\leq \lambda,
\end{equation}
\begin{equation}\label{65A}
\left|\partial_{\varrho} \beta _{\alpha ,hk}\eta ^{h}\eta
^{k}\right| \leq  C
\beta_{\alpha,hk}(\varrho,\theta)\eta^{h}\eta^{k},\quad\mbox{for }
\eta \in\mathbb{R}^{n-1},
\end{equation}
\begin{equation}\label{67A}
\left|\partial_{\varrho}\log \sqrt{\beta _{\alpha }}\right| \leq C,
\end{equation}
\begin{equation}\label{71sAa}
\left|\partial_{\varrho}\mu _{\alpha }(\varrho,\theta) \right| \leq
C,
\end{equation}
where $C$ depends on $\lambda $ and $\Lambda $ only.
\end{prop}

\textbf{Proof.} Let us omit the index $\alpha$. By \eqref{4A} and \eqref{8A} we have
$$\lambda^{-1}\leq \mu_0(x)\leq \lambda,$$ by these inequalities and by
\eqref{24A} we obtain \eqref{mu}.

 Now, we  prove \eqref{65A}. For any vector $\left\{\eta ^{h}\right\}
_{h=2}^{n}$, let $y_{\eta}$ be defined by \eqref{39A}.

 By \eqref{39A}, \eqref{17A} and \eqref{28A} we have
\begin{equation}\label{72A}
\beta
_{hk}\left(\varrho,\theta\right)\eta^{h}\eta^{k}=\widetilde{G}\left(\Gamma\left(\varrho,\theta\right)\right)y_{\eta}\cdot
y_{\eta}.
\end{equation}
Therefore
\begin{equation}\label{73A}
\partial_{\varrho }\beta_{hk}(\varrho,\theta)\eta^{h}\eta^{k}=
2\widetilde{G}(\Gamma) \partial_{\varrho} y_{\eta }\cdot y_{\eta}+
\left(\partial_{x^{k}}\widetilde{G}(\Gamma) \partial_{\varrho
}\Gamma^{k}y_{\eta}\right)\cdot y_{\eta }.
\end{equation}
Let us recall that $$\left|\partial_{\varrho}y_{\eta}\right| \leq
C\left|y_{\eta}\right|.$$ By \eqref{2A}), \eqref{4A} \eqref{9A},
\eqref{16sAc}, \eqref{73A} and by the last inequality we get
\begin{equation}\label{73Ab}
	\begin{aligned}
		\left|\partial_{\varrho}
\beta_{hk}(\varrho,\theta)\eta^{h}\eta^{k}\right|&\leq
C\left|y_{\eta}\right|^{2}\leq \\&\leq C\lambda^{2}\widetilde{G}(\Gamma)
y_{\eta}\cdot y_{\eta
}=\\&=C\lambda^{2}\beta_{hk}\left(\varrho,\theta\right)\eta^{h}\eta^{k},
\end{aligned}
\end{equation}
where $C$ depends on $\lambda $ and $\Lambda $ only. Therefore
\eqref{65A} is proved.

In order to prove \eqref{67A}, recall that if $A(s)$ is a matrix--valued function of class $C^{1}$ such that $\det A(s)\neq 0$ then we have the following equality (we denote by tr$(\cdot)$ the trace of the
matrix in the brackets)
\begin{equation*}
\frac{d}{ds}\log \left| \det A\left(s\right) \right|
=\text{tr}\left( \frac{dA\left( s\right) }{ds}A^{-1}\left( s\right)
\right).
\end{equation*}
This equality and \eqref{65A} give
\begin{equation*}
\begin{aligned}
\left|\partial_{\varrho}\log
\sqrt{\beta}\right|&=\frac{1}{2}\left|(\partial_{\varrho}\beta_{ij})\beta^{ji}\right|
=\frac{1}{2}\left|
(\partial_{\varrho}\beta_{hk})\delta_{i}^{h}\delta _{j}^{k}\beta
^{lm}\delta_{l}^{i}\delta_{m}^{j}\right|\leq\\&\leq C\beta
_{hk}\delta _{i}^{h}\delta _{j}^{k}\beta ^{lm}\delta_{l}^{i}\delta
_{m}^{j}=C\beta_{ij}\beta ^{ji}=C\left(n-1\right),
\end{aligned}
\end{equation*}
where $C$ depends on $\lambda $ and $\Lambda $ only. Therefore 
\eqref{67A} follows. $\blacksquare$

\bigskip

In order to prove Propositions \ref{Propos-8} and \ref{Propos-9}
stated below we need a partition of unity $\left\{ \zeta _{\alpha }\right\}_{\alpha\in \mathcal{J}}$
subordinate to the (finite) covering  $\left\{\left( U_{\alpha
},\varphi_{\alpha}\right)\right\}_{\alpha\in \mathcal{J}}$. By this, we mean that for each $\alpha\in\mathcal{J}$,
$\zeta_{\alpha }\in C^{\infty }\left( \Sigma \right)$,
$\zeta_{\alpha }\geq 0$, supp $\zeta_{\alpha}\subset U_{\alpha}$ and
$$\sum_{\alpha \in \mathcal{J}}\zeta_{\alpha}(p) =1,\quad\forall  p\in
\Sigma.$$ Let us denote by $\widehat{\zeta}_{\alpha}$ the function
$\zeta_{\alpha }\circ \varphi_{\alpha}^{-1}$ and set
$\widetilde{\zeta}_{\alpha}(x)=\zeta_{\alpha}\left(\Gamma (1;x,)
\right)$. We have $$\sum_{\alpha \in \mathcal{J}}\widetilde{\zeta
}_{\alpha }(x) =1,\quad\forall x\in \widetilde{B}_{1}.$$

\begin{prop}\label{Propos-8}
For every $\varrho\in (0,1) $, let us denote by $d\Omega_{\varrho}$
the element of volume on $\Sigma$.

Let $f$ be a function belonging to $C^{0}(\overline{B_1})$. We have
\begin{equation}
\int_{B_{1}}f(x) \sqrt{\widetilde{g}(x)
}dx=\int_{0}^{R_{0}}d\varrho\int_{\Sigma }f(\Gamma(\varrho,p))
\varrho^{n-1}d\Omega_{\varrho}.
\end{equation}
\end{prop}

\textbf{Proof.} For any  $\sigma \in (0,1) $ we have
\begin{equation}\label{112B}
\int_{B_{1}\setminus B_{\sigma}}f\left(x\right)
\sqrt{\widetilde{g}(x)}dx=\sum_{\alpha\in\mathcal{J}}
\int_{\mathcal{I}\left(U_{\alpha }\right)\setminus B_{\sigma }}f(x)
\sqrt{\widetilde{g}(x)}\widetilde{\zeta}_{\alpha}(x) dx.
\end{equation}
Now, in the integral on the right--hand side of \eqref{112B}, we perform the following change of variables: $$x=\Phi_{\alpha}^{-1}(\varrho,\theta).$$
By such a change of variables and by \eqref{15A}, \eqref{27A}, \eqref{31A} we get
\begin{equation}\label{113B}
	\begin{aligned}
	&\int\limits_{\mathcal{I}\left( U_{\alpha }\right)\setminus B_{\sigma
	}}f(x) \sqrt{\widetilde{g}(x)
	}\,\widetilde{\zeta}_{\alpha}\left(x\right) dx=\\&
	=\int_{\sigma }^{1}\int\limits_{V_{\alpha}}f(\Gamma_{\alpha}
	(\varrho;\varphi_{\alpha}^{-1}(\theta )) \varrho^{n-1} \sqrt{\beta
		_{\alpha }(\varrho,\theta) }\,\widehat{\zeta}_{\alpha }(\theta)
	d\theta d\varrho.
	\end{aligned}
\end{equation}
By \eqref{112B} and \eqref{113B} we get
\begin{equation*}\label{114B}
	\begin{aligned}
	&\int_{B_{1}\setminus B_{\sigma}}f(x)\sqrt{\widetilde{g}(x)}dx=
	\\&
	=\int_{\sigma }^{1}d\varrho\sum\limits_{\alpha\in\mathcal{J}
	}\int\limits_{V_{\alpha}}f\left(\Gamma\left(\varrho;\varphi_{\alpha
	}^{-1}(\theta)\right) \right)\varrho^{n-1}\sqrt{\beta _{\alpha
		}\left( \varrho,\theta\right)}\,\widehat{\zeta }_{\alpha }\left(
	\theta \right) d\theta=
	\\&
	=\int_{\sigma }^{1}d\varrho\int\limits_{\Sigma }f\left( \Gamma
	\left(\varrho,p\right)\right) \varrho^{n-1}d\Omega_{\varrho}.
	\end{aligned}
\end{equation*}
Finally, by the latter we have
\begin{equation*}
	\begin{aligned}
		\int_{B_{1}}f(x)\sqrt{\widetilde{g}(x)}dx&=\lim_{\sigma \rightarrow
0}\int_{B_{1}\setminus B_{\sigma}}f( x)\sqrt{\widetilde{g}(x)}dx
=\\&=\int_{0}^{1}d\varrho\int_{\Sigma}f\left(\Gamma
(\varrho,p)\right)\varrho^{n-1}d\Omega_{\varrho}.
\end{aligned}
\end{equation*}
$\blacksquare$

\bigskip

For any $\varrho\in(0,1)$ denotes by $\Xi ^{(\varrho)}$ the covariant 
tensor of order $2$ whose components with respect to 
coordinate neighborhoods $\left(U_{\alpha},\varphi_{\alpha}\right)$ are equal to
$$ \partial_{\varrho}\beta_{\alpha }^{lm}(\varrho,\cdot) \beta_{\alpha,lh}(\varrho,\cdot)\beta_{\alpha,mk}(\varrho,\cdot),\quad\mbox{for } h,k=2,\cdots,n.$$ Let us denote by $\ell(\varrho,\cdot)$  the  $C^{\infty }(\Sigma)$ function whose 
expressions with respect to the 
coordinate neighborhoods  $\left( U_{\alpha
},\varphi _{\alpha }\right)$ are equal to
$$\partial_{\varrho}\log\sqrt{\beta _{\alpha}(\varrho,\cdot)}.$$

\medskip

In the following Proposition we will denote by $\int \left(\cdot\right)
$ the integral $\int_{0}^{\varrho_{0}}d\varrho\int_{\Sigma }(\cdot)
d\Omega_{\varrho}$.

\medskip

\begin{prop}\label{Propos-9}
Let $v_{1},v_{2}\in C^{\infty }((0,\varrho_{0}) \times \Sigma)$.
Let us suppose either $v_{1}$ or $v_{2}$ of compact support.
Let $h\in C^{\infty}((0 ,\varrho_{0}))$.

Then we have
\begin{equation} \label{15B}
\int v_{1}\partial_{\varrho} v_{2}=-\int \left(v_{1}\ell
+\partial_{\varrho}v_{1}\right)v_2.
\end{equation}
 If $v_{1}$  has compact support then
\begin{equation}\label{17B}
\begin{aligned}
&\int h\left( \varrho\right)
\partial_{\varrho}v_{1}\mbox{div}_{\Sigma}\left(
\mu^{1-\frac{n}{2}}\nabla_{\Sigma}v_{1}\right) =\\&= \frac{1}{2} \int
h\left(\varrho\right)\mu^{1-\frac{n}{2}}\Xi^{\left(\varrho\right)
}\left(\nabla_{\Sigma }v_{1},\nabla_{\Sigma }v_{1}\right)+\\&
+\frac{1}{2}\int \left(\partial_{\varrho}\left( h(\varrho)
\mu^{1-\frac{n}{2}}\right)+h(\varrho)\mu^{1-\frac{n}{2}}\ell\right)
\left| \nabla _{\Sigma }v_{1}\right|_{\Sigma }^{2}.
\end{aligned}
\end{equation}
\end{prop}
\textbf{Proof.} For every $\alpha \in \mathcal{J}$ we denote, respectively, by
$v_{1,\alpha }( \varrho,\theta)$ and
$v_{2,\alpha}\left(\varrho,\theta\right)$, the funzctions
$v_{1}\left(\varrho,\varphi_{\alpha }^{-1}( \theta)\right) $ and
$v_{2}\left(\varrho,\varphi_{\alpha }^{-1}( \theta)\right)$.

Let us prove \eqref{15B}.

We have
\begin{equation*}
	\begin{aligned}
	&\int_{0}^{\varrho_{0}}d\varrho\int\limits_{\Sigma
	}v_{1,\alpha}\partial_{\varrho}v_{2,\alpha}d\Omega_{\varrho}=\sum\limits_{\alpha
		\in\mathcal{J}}\int_{0}^{\varrho_{0}}d\varrho\int\limits_{V_{\alpha}}
	v_{1,\alpha}(\partial_{\varrho}v_{2,\alpha})
	\sqrt{\beta_{\alpha}}\widehat{\zeta}_{\alpha}d\theta= \\&
	=-\sum\limits_{\alpha\in\mathcal{J}}\int_{0
	}^{\varrho_{0}}d\varrho\int\limits_{V_{\alpha}}\left(v_{1,\alpha
	}\partial_{\varrho}\log\sqrt{\beta_{\alpha}}+\partial_{\varrho}v_{1,\alpha
	}
	\right)v_{2,\alpha }\sqrt{%
		\beta_{\alpha }}\widehat{\zeta}_{\alpha }d\theta= \\&
	=-\int_{0}^{\varrho_{0}}d\varrho\int\limits_{\Sigma }\left(
	v_{1}\ell+\partial_{\varrho}v_{1}\right) v_{2}d\Omega _{\varrho}.
	\end{aligned}
\end{equation*}

\medskip

Let us prove \eqref{17B}. Let us suppose
$v_{1}$ with compact support. By the divergence Theorem on the
Riemannian manifold $\Sigma $ we have
\begin{equation}\label{19B}
	\begin{aligned}
	&\int_{0}^{\varrho_{0}}d\varrho\int\limits_{\Sigma }f( \varrho)
	\partial_{\varrho}v_{1} div_{\Sigma
	}\left(\mu^{1-\frac{n}{2}}\nabla_{\Sigma }v_{1}\right) d\Omega
	_{\varrho}=
	\\&
	=-\int_{0}^{\varrho_{0}}d\varrho\int\limits_{\Sigma }f\left(
	\varrho\right)\mu^{1-\frac{n}{2}}\left\langle \nabla _{\Sigma
	}v_{1},\nabla _{\Sigma }\partial_{\varrho}v_{1}\right\rangle d\Omega
	_{\varrho}.
	\end{aligned}
\end{equation}
Now
\begin{equation*}
\left\langle \nabla _{\Sigma }v_{1},\nabla_{\Sigma
}\partial_{\varrho}v_{1}\right\rangle=\frac{1}{2}\partial_{\varrho}
\left\langle \nabla _{\Sigma }v_{1},\nabla _{\Sigma
}v_{1}\right\rangle - \frac{1}{2}\Xi
^{(\varrho)}\left(\nabla_{\Sigma }v_{1},\nabla _{\Sigma
}v_{1}\right) ,
\end{equation*}
that, with \eqref{15B} gives \eqref{17B}. $\blacksquare$
\section{The  case of variable coefficients} \label{aronszajn-general}

In this Section we will prove the Carleman estimate of Aronszajn--Krzywicki--Szarski, \cite{AKS}. Basically we will proceed in a similar way to Section \ref{aronszajn-constant}, however, instead of the
Euclidean polar coordinates we will use the geodesic polar coordinates introduced
in Section \ref{geodesic-polar}. Compared with the original proof
of \cite{AKS}, the one we will prove here has some
simplification.

\medskip

Precisely we prove

\begin{theo}[{\bf Carleman estimate for $\Delta_g$}]\label{prop:Carlm-delta-g}
	\index{Theorem:@{Theorem:}!- Carleman estimate for $\Delta_g$@{- Carleman estimate for $\Delta_g$}}
Let us suppose that the matrix $G=\left\{g_{ij}\left( x\right)
\right\}_{i,j=1}^{n}$ satisfies to \eqref{2A}, \eqref{4A}, \eqref{3A}
and $g_{ij}\in C^{\infty}\left(B_2\right)$, for $i,j=1,\cdots, n$.

 Let $\epsilon\in(0,1]$.We define
\begin{equation}
    \label{peso-g}
        \rho(x) = \phi_{\epsilon}\left(|x|\right),\quad \forall x\in B_1\setminus \{0\},
\end{equation}
where
\begin{equation}
    \label{eq:24.3-delta-g}
        \phi_{\epsilon}(s) = \frac{s}{\left(1+s^{\epsilon}\right)^{1/\epsilon}}.
\end{equation}
Then there exist  $r_0\in (0,1)$, $\overline{\tau}>1$ and $C>1$,
which depend on $\epsilon$, $\lambda$ and $\Lambda$ only, such that 

\begin{equation}\label{Carlm-delta-g}
\begin{aligned}
\tau^{3}\int\rho^{\epsilon-2\tau}|u|^2dx&+
\tau\int\rho^{2+\epsilon-2\tau}|\nabla u|^2dx+\\&
+\tau^2r\int\rho^{-1-2\tau}u^2dx\leq C\int\rho^{4-2\tau}|\Delta_g
u|^2dx,
\end{aligned}
\end{equation}
for every $r\in (0,r_0)$, for every $\tau\geq \overline{\tau}$ and for every $u\in C^\infty_0\left(B_{r_0}\setminus
\overline{B}_{r/4}\right)$.
\end{theo}

\bigskip

We start by the following simple estimation of the first order operator
 $\mathcal{L}_{2}$ defined in \eqref{91B}.

\begin{prop}\label{estimate-L-2}
The following estimate holds true
\begin{equation}\label{91B-estimate}
\left\vert\mathcal{L}_{2}w\right\vert\leq
C\left\vert\partial_{\varrho}w\right\vert\quad \forall w\in
C^{\infty}\left((0,1) \times\Sigma\right),
\end{equation}
where $C$ depends on $\lambda$ and $\Lambda$ only.
\end{prop}
\textbf{Proof. } By \eqref{67A} and \eqref{71sAa} we have
\begin{equation*}
\begin{aligned}
&\left\vert\left(\mathcal{L}_2w\right)\left(\varrho,\varphi_{\alpha}^{-1}\left(\theta\right)\right)\right\vert=\\&
=
\left\vert\mu_{\alpha}\left(\left(1-\frac{n}{2}\right)\partial_{\varrho}
\log \mu_{\alpha}+\partial_{\varrho} \log
\sqrt{\beta_{\alpha}}\right)\partial_{\varrho}w_{\alpha}\right\vert
\leq\\& \leq
\left(\left\vert\mu_{\alpha}\left(1-\frac{n}{2}\right)\partial_{\varrho}
\log \mu_{\alpha}\right\vert+\left\vert\partial_{\varrho} \log
\sqrt{\beta_{\alpha}}\right\vert\right)\left\vert\partial_{\varrho}w_{\alpha}\right\vert\leq\\&
\leq C\left\vert\partial_{\varrho}w_{\alpha}\right\vert,
\end{aligned}
\end{equation*}
where $C$ depends on $\lambda$ and $\Lambda$ only. $\blacksquare$

\bigskip

Let $u$ be an arbitrary function that belongs to
$C^{\infty}_0\left(B_1\setminus\{0\}\right)$ and let us denote by $w$ the function (recall \eqref{def-w})
\begin{equation}\label{def-w-1}
(0,1) \times\Sigma\ni\left(\varrho,p\right) \rightarrow
w\left(\varrho,p\right)=u\left(\Gamma\left(\varrho,p\right)\right) .
\end{equation}
where $\Gamma$ is defined in Proposition \ref{Propos-2}.

Now, we carry out the following change of variables
$$\varrho=e^t,\quad\quad\widetilde{w}(t,p)=w\left(e^t,p\right),\quad\forall(t,p)\in (-\infty,0)\times \Sigma$$
and we adopt the following conventions: for any function
$h(\varrho,\cdot)$ (or for every tensor) in which the variable $\varrho$ occurs , we denote by $\widetilde{h}(t,\cdot)$ the function (or tensor) $h(e^t,\cdot)$. We will continue to denote by
$\left\langle\cdot,\cdot\right\rangle$ and $\left|_{\cdot
}\right|_{\Sigma }$, respectively, the inner product and the norm
associated with it in the Riemannian structure induced by the metric tensor
$$\left\{\widetilde{\beta}_{\alpha,hk}\right\}_{h,k=2}^n.$$
We will still denote by $\nabla_{\Sigma}$, div$_{\Sigma }$ and
$d\Omega_t$, respectively, the gradient, the divergence operators and the element of volume on $\Sigma$ in the above mentioned structure. In particular, the local expression of $d\Omega_t$
is equal to $\sqrt{\widetilde{\beta}_{\alpha}}d\theta$. Moreover we 
set
\begin{equation}\label{operatore-M}
\mathcal{M}(\cdot)=\mbox{div}_{\Sigma}\left(\widetilde{\mu}^{1-\frac{n}{2}}\nabla_{\Sigma}\cdot\right).
\end{equation}
We have
\begin{equation}\label{laplace-g-polar-2}
    e^{2t}\widetilde{\mu}^{-1}(\mathcal{L}_1 w)(e^t,p)=\mathcal{P}\widetilde{w}(t,p),\quad\forall (t,p)\in (-\infty,0)\times \Sigma,
\end{equation}
where
\begin{equation}\label{laplace-g-polar-3}
    \mathcal{P}\widetilde{w}=\widetilde{w}_{tt}+(n-2)\widetilde{w}_{t}+\frac{1}{\widetilde{\mu}^{1-n/2}}
    \mathcal{M}\widetilde{w}.
\end{equation}

For the reader's convenience, we reformulate Propositions \ref{Propos-7}
and \ref{Propos-9} in the Riemannian structure induced by the metric tensor
$\left\{\widetilde{\beta}_{\alpha,hk}\right\}_{h,k=2}^n$.
Recall that $\left\{ U_{\alpha },\varphi _{\alpha
}\right\}_{\alpha \in \mathcal{J}}$ is a (finite) family  of
coordinate neighborhoods  defining on $\Sigma $ a structure of
$C^{\infty }$ oriented differentiable manifold, where, for any
$\alpha \in\mathcal{J}$, we set $V_{\alpha }=\varphi _{\alpha
}\left( U_{\alpha }\right)$.

\begin{prop}\label{Propos-7-bis}
For any  $\alpha \in\mathcal{J}$ and for any
$\left(t,\theta\right)\in(-\infty,0)\times V_{\alpha}$, we have
\begin{equation}\label{65A-bis}
\left|\partial_t \widetilde{\beta}_{\alpha ,hk}\eta ^{h}\eta
^{k}\right| \leq  C e^t
\widetilde{\beta}_{\alpha,hk}(\varrho,\theta)\eta^{h}\eta^{k},\quad \forall\eta \in\mathbb{R}^{n-1},
\end{equation}
\begin{equation}\label{67A-bis}
\left|\partial_t\log \sqrt{\widetilde{\beta} _{\alpha }}\right|\leq
Ce^t, 
\end{equation}
\begin{equation}\label{71sAa-bis}
\left|\partial_t\widetilde{\mu} _{\alpha }(t,\theta) \right| \leq
Ce^t,
\end{equation}
where $C$ depends on $\lambda $ and $\Lambda$ only.
\end{prop}

By the convention introduced above, for any $t\in(-\infty ,0)$,
let us denote by $\widetilde{\Xi}^{(t)}$ the covariant tensor satisfying
$$\widetilde{\Xi}^{(t)}=\Xi^{(e^t)}$$ and let us denote by
$$\widetilde{\ell}(t,\cdot)=\ell(e^t,\cdot).$$

\medskip

In the next Proposition let us denote by $\int (\cdot)$
the integral $\int_{-\infty}^{t_{0}}dt\int_{\Sigma }(\cdot)
d\Omega_{t}$.

\begin{prop}\label{Propos-9-bis}
Let $v_{1},$ $v_{2}\in C^{\infty }((-\infty,t_{0}) \times
\Sigma)$. Let us suppose either $v_{1}$ or $v_{2}$ of compact support. Let  $h\in
C^{\infty}((-\infty ,t_{0}))$. Then we have
\begin{equation} \label{15B-bis}
\int v_{1}\partial_t v_{2}=-\int \left(
e^tv_{1}\widetilde{\ell}+\partial_t v_{1}\right)v_2.
\end{equation}
If $v_{1}$ has compact support then 
\begin{equation}\label{17B-bis}
\begin{aligned}
&\int h(t) \partial_t v_{1}\mbox{div}_{\Sigma}\left(
\widetilde{\mu}^{1-\frac{n}{2}}\nabla_{\Sigma}v_{1}\right)=\\&=
\frac{1}{2} \int e^t
h(t)\widetilde{\mu}^{1-\frac{n}{2}}\widetilde{\Xi}^{(t)
}\left(\nabla_{\Sigma }v_{1},\nabla_{\Sigma }v_{1}\right)+\\& +
\frac{1}{2}\int \left[\partial_t\left( h(t)
\widetilde{\mu}^{1-\frac{n}{2}}\right)+e^t
h(t)\widetilde{\mu}^{1-\frac{n}{2}}\widetilde{\ell}\right]
\left\vert \nabla_{\Sigma }v_{1}\right\vert_{\Sigma }^{2}.
\end{aligned}
\end{equation}
\end{prop}

\bigskip

\begin{prop}\label{Propos-10}
We have
\begin{equation} \label{ell}
\left\vert\widetilde{\ell}(t,p)\right\vert \leq C, \quad\forall (t,p)\in (-\infty,0) \times \Sigma,
\end{equation}
\begin{equation}\label{Chi}
\left\vert\widetilde{\Xi}^{(t)}(\nabla_{\Sigma}v,\nabla_{\Sigma}v)\right\vert\leq C
\left\vert\nabla_{\Sigma}v\right\vert^2_{\Sigma}, \ \ \forall
v\in C^{\infty }((-\infty,0)\times \Sigma),
\end{equation}
where $C$ depends on $\lambda$ and $\Lambda$ only.
\end{prop}
\textbf{Proof.} Inequalities \eqref{ell} and \eqref{Chi} are an immediate
consequences of Proposition \ref{Propos-7}.$\blacksquare$

\bigskip

\textbf{Proof of Theorem \ref{prop:Carlm-delta-g}.}

For any smooth function $v$ we write $v^{\prime}$,
$v^{\prime\prime}$, ... instead of $\partial_tv$, $\partial_{tt}v$,
... .

By \eqref{peso-1} we have, similarly to the proof
of Theorem \ref{prop:Carlm-delta-g} (here and in the sequel we omit
the subscript $\epsilon$)
\begin{equation}\label{peso-1-g}
\varphi(t):=\log(\phi(e^t))=t-\epsilon^{-1}\log\left(1+e^{\epsilon
t}\right), \quad\forall t\in(-\infty,0)
\end{equation}
and
\begin{equation}\label{peso-derivate-g}
\varphi^{\prime}(t)=\frac{1}{1+e^{\varepsilon t}}, \quad
\varphi^{\prime\prime}(t)=-\frac{\epsilon e^{\varepsilon
t}}{(1+e^{\epsilon t})^2},   \quad \forall t\in(-\infty,0).
\end{equation}
Let
\begin{equation}\label{def-f-g}
f(t,p)=e^{-\tau\varphi}w(t,p) , \quad\forall
(t,p)\in(-\infty,0)\times \Sigma,
\end{equation}
where $w$ is defined in \eqref{def-w-1}.

 We have
\begin{equation}\label{P-tau}
\mathcal{P}_{\tau}f:=e^{-\tau\varphi}\mathcal{P}(e^{\tau\varphi}f)=\underset{\mathcal{A}_{\tau}f}{\underbrace{b_0f+b_1f^{\prime}}}+
\underset{\mathcal{S}_{\tau}f}{\underbrace{a_0f+f^{\prime\prime}+\frac{1}{m}\mathcal{M}f}},
\end{equation}
where

\begin{equation*}
m=\widetilde{\mu}^{1-\frac{n}{2}}
\end{equation*}
and
\begin{equation}\label{coefficienti-g}
a_0=\tau^2\varphi^{\prime^2}+\tau(n-2), \quad b_0=\tau
\varphi^{\prime\prime},\quad b_1=2\tau\varphi^{\prime}+(n-2).
\end{equation}
Let us note that \eqref{mu} gives
\begin{equation}\label{stima-m}
\lambda^{\frac{n}{2}-1}\leq m\leq\lambda^{1-\frac{n}{2}}.
\end{equation}
Set
\begin{equation}\label{gamma-g}
\gamma:=\frac{1}{\varphi'}=1+e^{\epsilon t}.
\end{equation}
We have
\begin{equation}\label{square-g}
\int m\gamma\left\vert \mathcal{P}_{\tau}f\right\vert^2\geq 2\int
m\gamma \mathcal{A}_{\tau}f\mathcal{S}_{\tau}f+\int
m\gamma\left\vert\mathcal{A}_{\tau}f\right\vert^2,
\end{equation}

\begin{equation}\label{commutator-g}
\begin{aligned}
2\int
m\gamma\mathcal{A}_{\tau}f\mathcal{S}_{\tau}f&=\underset{I_1}{\underbrace{2\int
\gamma\left(b_0f+b_1f^{\prime}\right)\mathcal{M}f}}+\\&
+\underset{I_2}{\underbrace{2\int
m\gamma\left(b_0f+b_1f^{\prime}\right)
\left(a_0f+f^{\prime\prime}\right) }}.
\end{aligned}
\end{equation}

\textbf{We examine $I_1$}.

We have
\begin{equation}\label{split-I-1}
I_1=2\int\left(\gamma b_0f\mathcal{M}f+\gamma
b_1f^{\prime}\mathcal{M}f\right)
=\underset{I_{11}}{\underbrace{2\int\gamma b_0f\mathcal{M}f}}
+\underset{I_{12}}{\underbrace{2\int\gamma
b_1f^{\prime}\mathcal{M}f}}.
\end{equation}
By the divergence Theorem we obtain

\begin{equation}\label{I-11}
I_{11}=2\int\gamma
b_0f\mbox{div}_{\Sigma}\left(m\nabla_{\Sigma}f\right) =-2\int
m\gamma b_0\left\vert \nabla_{\Sigma}f\right\vert^2_{\Sigma}.
\end{equation}
By \eqref{17B-bis} we get
\begin{equation}\label{I-12}
\begin{aligned}
I_{12}&=2\int\gamma b_1f^{\prime}\mathcal{M}f =
 \int \left[\left( m \gamma b_1\right)'+e^t m \gamma b_1\widetilde{\ell}\right] \left\vert \nabla_{\Sigma }f\right\vert_{\Sigma}^{2}+\\&+\int e^t m \gamma b_1\widetilde{\Xi}^{(t)
}\left(\nabla_{\Sigma }f,\nabla_{\Sigma }f\right).
\end{aligned}
\end{equation}
By \eqref{split-I-1}, \eqref{I-11} and \eqref{I-12} we have

\begin{equation}\label{I-1-somma}
\begin{aligned}
I_1=I_{11}+I_{12}&= \int \left\{-2m\gamma b_0+\left[\left(m\gamma
b_1\right)^{\prime}+e^t m \gamma b_1\widetilde{\ell}\right]
\right\}\left\vert \nabla_{\Sigma }f\right\vert_{\Sigma}^{2}+\\&
+\int e^t m \gamma b_1\widetilde{\Xi}^{(t) }\left(\nabla_{\Sigma
}f,\nabla_{\Sigma }f\right).
\end{aligned}
\end{equation}
We get (compare with \eqref{calcolo-1})

\begin{equation}\label{I-1-primo-termine}
\begin{aligned}
&-2m\gamma b_0+\left[\left( m \gamma b_1\right)^{\prime}+e^t m \gamma
b_1\widetilde{\ell}\right]=\\&=2m\left(-\gamma b_0 +
\frac{1}{2}\left(\gamma b_1\right)^{\prime}\right)+m^{\prime}\gamma
b_1+e^t m \gamma b_1\widetilde{\ell}\geq \\&\geq \tau\epsilon
e^{\epsilon t}++m^{\prime}\gamma b_1+e^t m \gamma
b_1\widetilde{\ell}.
\end{aligned}
\end{equation}

Now, by \eqref{71sAa-bis} we have
\begin{equation}\label{m-derivata}
\left\vert m^{\prime}\right\vert\leq Ce^t,
\end{equation}
where $C$ depends on  $\lambda$ and $\Lambda$ only. By this
inequality, by \eqref{Chi} and by \eqref{I-1-somma}, taking into account \eqref{ell}, we have
\begin{equation}\label{I-1-somma-29}
I_1\geq \tau\int e^{\epsilon t}\left(\epsilon
-C_{\star}e^{(1-\epsilon )t}\right)\left\vert \nabla_{\Sigma
}f\right\vert_{\Sigma}^{2}, \quad\forall \tau\geq 1,
\end{equation}
where $C_{\star}$ depends on $\lambda$ and $\Lambda$ only.

Now, for every $t\geq
T_1:=\frac{1}{1-\varepsilon}\log\frac{\varepsilon}{2 C_{\star}}$
we have 
\begin{equation*}
\epsilon -C_{\star}e^{(1-\epsilon )t}\geq \frac{\epsilon}{2}
\end{equation*}
By this inequality and by \eqref{I-1-somma-29} we have, for every
$f\in C_0^{\infty}((-\infty,T_1)\times \Sigma)$,
\begin{equation}\label{I-1-somma-1}
I_1\geq \frac{\epsilon \tau }{2}\int e^{\epsilon t}\left\vert
\nabla_{\Sigma }f\right\vert_{\Sigma}^{2}, \quad\forall
\tau\geq 1.
\end{equation}

\textbf{We examine $I_2$.}

By \eqref{15B-bis} we have
\begin{equation*}
	\begin{aligned}
&I_2=2\int
m\gamma\left(b_0f+b_1f^{\prime}\right)\left(a_0f+f^{\prime\prime}\right)=\\&=2\int
m\gamma\left(a_0b_0f^2+b_0ff^{\prime\prime}+b_1a_0f^{\prime}f+b_1f^{\prime}f^{\prime\prime}\right)=\\&
=2\int m\gamma a_0b_0f^2-\left(m\gamma
b_0f\right)^{\prime}f^{\prime}+\frac{1}{2}m\gamma
b_1a_0\left(f^2\right)^{\prime}+\frac{1}{2}m\gamma b_1
\left(f^{\prime^2}\right)^{\prime}-\\&-2\int e^tm\gamma
b_0\widetilde{\ell}ff^{\prime}=\\&=2\int m\gamma
a_0b_0f^2-\left(m\gamma
b_0f\right)^{\prime}f^{\prime}-\frac{1}{2}\left(m\gamma
b_1a_0\right)^{\prime}f^2-\frac{1}{2}\left(m\gamma
b_1\right)^{\prime} f^{\prime^2}-\\& -2\int e^tm\gamma
\widetilde{\ell} \left(b_0ff^{\prime}+\frac{1}{2}
b_1a_0f^2+\frac{1}{2}b_1 f^{\prime^2}\right)==2\int
m\left[\left(\gamma a_0b_0 +\frac{1}{2}\left(\gamma
b_1a_0\right)^{\prime}\right)f^2 +\right.\\& \left.+\left(\gamma
b_0+\frac{1}{2}\left(\gamma b_1\right)^{\prime}\right)
f^{\prime^2}-\left(\gamma b_0\right)^{\prime}ff^{\prime}\right]-\\&
-2\int e^tm\gamma \widetilde{\ell} \left(b_0ff^{\prime}-\frac{1}{2}
b_1a_0f^2-\frac{1}{2}b_1 f^{\prime^2}\right)-\\&-2\int m^{\prime}\gamma
\left(b_0 ff^{\prime}+\frac{1}{2}b_1a_0f^2+\frac{1}{2}b_1 f^{\prime^2}\right)=\\&
\end{aligned}
\end{equation*}
\begin{equation*}
\begin{aligned}
&=2\int m\left[\left(\gamma a_0b_0 +\frac{1}{2}\left(\gamma b_1a_0\right)^{\prime}\right)f^2
+\right. \\&\left.+\left(\gamma b_0+\frac{1}{2}\left(\gamma b_1\right)^{\prime}\right)
f^{\prime^2}-\left(\gamma b_0\right)^{\prime}ff^{\prime}\right]-\\&
-2\int \gamma\left(me^t\widetilde{\ell}+m^{\prime}\right)
 \left(b_0 ff^{\prime}+\frac{1}{2}b_1a_0f^2+\frac{1}{2}b_1 f^{\prime^2}\right)=I_{21}+I_{22}.
\end{aligned}
\end{equation*}
Where we set
\begin{equation*}
	\begin{aligned}
		&I_{21}=2\int m\left[\left(\gamma a_0b_0 +\frac{1}{2}\left(\gamma b_1a_0\right)^{\prime}\right)f^2
		+\right. \\&\left.+\left(\gamma b_0+\frac{1}{2}\left(\gamma b_1\right)^{\prime}\right)
		f^{\prime^2}-\left(\gamma b_0\right)^{\prime}ff^{\prime}\right]
	\end{aligned}
\end{equation*}
and
\begin{equation*}
	\begin{aligned}
		&I_{22}=-2\int \gamma\left(me^t\widetilde{\ell}+m^{\prime}\right)
		\left(b_0 ff^{\prime}+\frac{1}{2}b_1a_0f^2+\frac{1}{2}b_1 f^{\prime^2}\right)
	\end{aligned}
\end{equation*}

Now, let us estimate $I_{21}$ e $I_{22}$.

We have (compare with \eqref{H-tilde-1-0} and \eqref{H-tilde-1}), for every $\tau>0$
\begin{equation}\label{H-tilde-1-0-g}
\gamma a_0b_0-\frac{1}{2}\left(\gamma b_1 a_0\right)^{\prime}=
\left(\tau^3\frac{1}{1+e^{\epsilon
t}}+\tau^2\frac{n-2}{2}\right)\frac{\epsilon e^{\epsilon
t}}{(1+e^{\epsilon t})^2}\geq \frac{\tau^3}{8} \epsilon e^{\epsilon
t}, 
\end{equation}
and (compare with \eqref{H-2-0} and \eqref{H-2}) fpr every $\tau>2(n-2)$
\begin{equation}\label{H-2-0-g}
\begin{aligned}
-\left(\gamma b_0+\frac{1}{2}\left(\gamma b_1\right)^{\prime}\right)
=\tau\frac{\epsilon e^{\epsilon t}}{1+e^{\epsilon
t}}-\frac{(n-2)\epsilon e^{\epsilon t}}{2}\geq\frac{\epsilon\tau
}{4} e^{\epsilon t}.
\end{aligned}
\end{equation}
Moreover

\begin{equation*}
\left\vert(\gamma b_0)'ff^{\prime}\right\vert =\tau\frac{\epsilon
e^{\epsilon t}}{\left(1+e^{\epsilon t}\right)^2}\left\vert
ff^{\prime}\right\vert\leq \epsilon e^{\epsilon
t}\left(\tau^2f^2+f^{\prime^2}\right).
\end{equation*}
By this inequality and by \eqref{H-tilde-1-0-g},  \eqref{H-2-0-g} we have

\begin{equation}\label{I-21}
I_{21}\geq \epsilon\int m e^{\epsilon
t}\left(\frac{\tau^3}{8}f^2+\frac{\tau}{4}f^{\prime^2}\right), \quad
\forall \tau\geq\overline{\tau}_1,
\end{equation}
where $\overline{\tau}_1=\max\{16,2(n-2)\}$.

Now we estimate  $I_{22}$ e $I_{23}$.

By \eqref{ell}, \eqref{coefficienti-g}, \eqref{stima-m},  and
\eqref{m-derivata} we have

\begin{equation}\label{I-22}
\left\vert I_{22}\right\vert\leq C_{\ast}\int
e^t\left(\tau^3f^2+\tau f^{\prime^2}\right), \quad  \forall \tau\geq
\overline{\tau}_1,
\end{equation}
where $C_{\ast}$ depends on $\lambda$ and $\Lambda$ only.

By \eqref{I-21} and \eqref{I-22} we have, for every $\tau\geq
\overline{\tau}_1$

\begin{equation}\label{I-2-2}
I_2\geq I_{21}-\left\vert I_{22}\right\vert\geq \int e^{\epsilon
t}\left(\frac{\epsilon
\lambda^{\frac{n}{2}-1}}{8}-C_{\ast}e^{(1-\epsilon)t}\right)\left(\tau^3
f^2+\tau f^{\prime^2}\right).
\end{equation}
Notice that for every $$t\geq
T_2:=\frac{1}{1-\varepsilon}\log\left(\frac{\varepsilon
\lambda^{\frac{n}{2}-1}}{16 C_{\ast}}\right)$$ we have

\begin{equation*}
\frac{\epsilon
\lambda^{\frac{n}{2}-1}}{8}-C_{\ast}e^{(1-\epsilon)t}\geq
\frac{\epsilon \lambda^{\frac{n}{2}-1}}{16},
\end{equation*}
this inequality and \eqref{I-2-2} imply that
\begin{equation}\label{I-2-3}
I_2\geq \frac{\varepsilon \lambda^{\frac{n}{2}-1}}{16}\int
e^{\varepsilon t}\left(\tau^3 f^2+\tau f^{\prime^2}\right),
\end{equation}
for every $f\in C_0^{\infty}((-\infty,T_2)\times \Sigma)$ and for every $\tau\geq\overline{\tau}_1$. By \eqref{commutator-g}, \eqref{I-1-somma-1} and \eqref{I-2-3}
we have,

\begin{equation}\label{commutator-final-g-1}
2\int m\gamma\mathcal{A}_{\tau}f\mathcal{S}_{\tau}f\geq
\frac{\varepsilon \tau }{2}\int e^{\epsilon t}\left\vert
\nabla_{\Sigma }f\right\vert_{\Sigma}^{2}+\frac{\epsilon
\lambda^{\frac{n}{2}-1}}{16}\int e^{\epsilon t}\left(\tau^3 f^2+\tau
f^{\prime^2}\right) ,
\end{equation}
for every $\tau\geq \overline{\tau}_1$ and for every $f\in
C_0^{\infty}((-\infty,T_3)\times \Sigma)$, where
$T_3=\min\{T_1,T_2\}$.

\medskip

Set
$$\epsilon_0=\epsilon\min\{\frac{1}{2},\frac{\lambda^{\frac{n}{2}-1}}{16}\},$$
by \eqref{square-g} and \eqref{commutator-final-g-1} we get
\begin{equation}\label{commutator-final-g-2}
	\begin{aligned}
&\int m\gamma\left\vert \mathcal{P}_{\tau}f\right\vert^2\geq\\&\geq 
\epsilon_0 \int \left(\tau^3f^2+\tau \left(f^{\prime^2}+\left\vert
\nabla_{\Sigma}f\right\vert_{\Sigma}^2\right)\right)e^{\epsilon t}
+\int m\gamma\left\vert\mathcal{A}_{\tau}f\right\vert^2,
\end{aligned}
\end{equation}
for every $\tau\geq \tau_1$ and for every $f\in
C_0^{\infty}((-\infty,T_3)\times \Sigma)$.

\medskip

In order to obtain the third term on the left--hand side of
\eqref{Carlm-delta-g} we argue as in the proof of 
Theorem \ref{prop:Carlm-delta}. For the sake of clarity, we repeat the
most important steps.

By the trivial inequality $(a+b)^2\geq\frac{1}{2}a^2-b^2$ we have
\begin{equation*}\label{antisimm-1-g}
\begin{aligned}
\int
m\gamma\left\vert\mathcal{A}_{\tau}f\right\vert^2&\geq\frac{\lambda^{\frac{n}{2}-1}}{2}\int
\gamma
\left(2\tau\varphi^{\prime}+n-2\right)^2f^{\prime^2}-\lambda^{\frac{n}{2}-1}\int
\gamma\tau^2 \varphi^{\prime\prime^2}f^2\geq\\&
\geq\lambda^{\frac{n}{2}-1}\tau^2\int f^{\prime^2}-
\epsilon^2\lambda^{\frac{n}{2}-1}\tau^2\int e^{2\epsilon t}f^2,
\quad \forall \tau>0.
\end{aligned}
\end{equation*}
Plugging this inequality into \eqref{commutator-final-g-2}
we have

\begin{equation}\label{antisimm-2-g}
\begin{aligned}
\int m\gamma\left\vert \mathcal{P}_{\tau}f\right\vert^2&\geq
\lambda^{\frac{n}{2}-1}\tau^2\int f^{\prime^2} +\tau^3\int
\left(\epsilon_0-\epsilon\lambda^{\frac{n}{2}-1}\tau^{-1}e^{\epsilon
t}\right) e^{\epsilon t}f^2+\\&+\epsilon_0\tau\int
\left(f^{\prime^2}+\left\vert
\nabla_{\Sigma}f\right\vert^2\right)e^{\epsilon t},
\end{aligned}
\end{equation}
for every $\tau\geq \tau_1$ and for every $f\in
C_0^{\infty}((-\infty,T_3)\times \Sigma)$.

By \eqref{antisimm-2-g} we have

\begin{equation}\label{antisimm-3-g}
\begin{aligned}
&\int m\gamma\left\vert \mathcal{P}_{\tau}f\right\vert^2\geq
\\&\geq \lambda^{\frac{n}{2}-1}\tau^2\int f^{\prime^2}+
\frac{\varepsilon_0\tau^3}{2}\int e^{\varepsilon t}f^2
+\varepsilon_0\tau\int \left(f^{\prime^2}+\left\vert
\nabla_{\Sigma}f\right\vert_{\Sigma}^2\right)e^{\varepsilon t},
\end{aligned}
\end{equation}

\medskip

\noindent for every  $\tau\geq\overline{\tau}_2$, where
$\overline{\tau}_2=\max\{2\epsilon\epsilon_0^{-1}\lambda^{\frac{n}{2}-1},
\overline{\tau}_1\}$ and for every $f\in
C_0^{\infty}((-\infty,T_3)\times \Sigma)$.

Let $r\leq 4 e^{T_3}$. Since $u=0$ in $B_{r/4}$ we have 
$f(t,p)=0$ for every $t\leq \log (r/4)$ and for every $p\in \Sigma$.
Proceeding as in the proof of \eqref{antisimm-5} we obtain

\begin{equation}\label{antisimm-5-g}
\begin{aligned}
C\int\left\vert \mathcal{P}_{\tau}f\right\vert^2 \geq &  \tau^3\int
e^{\varepsilon t}f^2+\tau\int \left(f^{\prime^2}+\left\vert
\nabla_{\Sigma}f\right\vert_{\Sigma}^2\right)e^{\varepsilon t}+\\&
+\tau^2r\int f^2e^{-t}, \quad\forall \tau\geq
\overline{\tau}_2
\end{aligned}
\end{equation}
where $C$ depends on $\epsilon$, $\lambda$ and $\Lambda$ only.

\bigskip

Now we come back to the original coordinates. By \eqref{def-w-1},
\eqref{laplace-g-polar-2}, \eqref{def-f-g} and \eqref{P-tau} we have

\begin{equation}\label{final-1-g}
\begin{aligned}
&\int^{0}_{-\infty}dt\int_{\Sigma}\left\vert
\mathcal{P}_{\tau}f\right\vert^2 d\Omega_{t}
=\\&=\int^{0}_{-\infty}dt\int_{\Sigma}e^{-2\tau\varphi(t)}e^{4t}|\widetilde{\mu}^{-1}(\mathcal{L}_1
w)(e^t,p)|^2d\Omega_{t}= \\& =
\int^{1}_{0}d\varrho\int_{\Sigma}e^{-2\tau\varphi(\log\varrho)}\varrho^{3}|\mu^{-1}\mathcal{L}_1
w(\varrho,p)|^2d\Omega_{\varrho} .
\end{aligned}
\end{equation}
By proposizions \ref{Propos-6} and \ref{estimate-L-2} and by
\eqref{final-1-g} we have

\begin{equation}\label{final-2-g}
\begin{aligned}
\int^{0}_{-\infty}dt\int_{\Sigma}\left\vert
\mathcal{P}_{\tau}f\right\vert^2 d\Omega_{t}&\leq
C\int^{1}_{0}d\varrho\int_{\Sigma}e^{-2\tau\varphi(\log\varrho)}\varrho^{3}|\mathcal{L}
w|^2d\Omega_{\varrho}+\\&
+C\int^{1}_{0}d\varrho\int_{\Sigma}e^{-2\tau\varphi(\log\varrho)}\varrho^{3}|\partial_{\varrho}w|^2d\Omega_{\varrho},
\end{aligned}
\end{equation}
where $C$ depends on $\lambda$ and $\Lambda$ only.

Moreover we have
\begin{equation}\label{final-3-g}
\int^{0}_{-\infty}dt\int_{\Sigma} f^2e^{-t}d\Omega_{t}
=\int^{1}_{0}d\varrho\int_{\Sigma}e^{-2\tau\varphi(\log\varrho)}\varrho^{-2}
|w|^2d\Omega_{\varrho}
\end{equation}
and
\begin{equation}\label{final-4-g}
\int^{0}_{-\infty}dt\int_{\Sigma} f^2e^{\epsilon t}d\Omega_{t}
=\int^{1}_{0}d\varrho\int_{\Sigma}e^{-2\tau\varphi(\log\varrho)}\varrho^{\epsilon-1}
|w|^2d\Omega_{\varrho}.
\end{equation}
Concerning the second integral on the right--hand side of
\eqref{final-1-g}, let $\delta\in (0,1)$ be to choose later, we get
\begin{equation}\label{final-5-g}
	\begin{aligned}
&\int^{0}_{-\infty}dt\int_{\Sigma} e^{\epsilon
	t}\left(f^{\prime^2}+\left\vert
\nabla_{\Sigma}f\right\vert_{\Sigma}^2\right)d\Omega_{t}\geq\\&\geq
\delta\int^{0}_{-\infty}dt\int_{\Sigma} e^{\epsilon
	t}\left(f^{\prime^2}+\left\vert
\nabla_{\Sigma}f\right\vert_{\Sigma}^2\right)d\Omega_{t} \geq \\& \geq
\frac{\delta}{2}\int^{0}_{-\infty}dt\int_{\Sigma}e^{\epsilon
	t}e^{-2\tau\varphi(t)}\left(|w_{\varrho}(e^t,p)|^2e^{2t}+\right.  \\&\left. +\left\vert
\nabla_{\Sigma}w(e^t,p)\right\vert_{\Sigma}^2-2\tau^2|w(e^t,p)|^2\right)
d\Omega_{t}=\\&
=\frac{\delta}{2}\int^{1}_{0}d\varrho\int_{\Sigma}e^{-2\tau\varphi(\log\varrho)}\left(|w_{\varrho}(\varrho,p)|^2+\right.\\&\left.+\varrho^{-2}\left\vert
\nabla_{\Sigma}w(\varrho,p)\right\vert_{\Sigma}^2\right)\varrho^{\epsilon+1}d\Omega_{\varrho}-\\&
-\frac{\delta\tau^2}{2}\int^{1}_{0}d\varrho\int_{\Sigma}e^{-2\tau\varphi(\log\varrho)}|w(\varrho,p)|^2
\varrho^{\epsilon-1}d\Omega_{\varrho}.
	\end{aligned}
\end{equation}
Now, plugging \eqref{final-2-g}, \eqref{final-3-g},
\eqref{final-4-g} and \eqref{final-5-g} into \eqref{antisimm-5-g}, we have 
\begin{equation}\label{final-6-g}
	\begin{aligned}
&\tau^3\left(1-\frac{\delta}{2}\right)\int^{1}_{0}d\varrho\int_{\Sigma}e^{-2\tau\varphi(\log\varrho)}|w|^2
\varrho^{\epsilon-1}d\Omega_{\varrho}+\\&
 +\tau^2r\int^{1}_{0}d\varrho\int_{\Sigma}e^{-2\tau\varphi(\log\varrho)}|w|^2 \varrho^{-2}d\Omega_{\varrho}+\\&
 +\frac{\delta}{2}\tau\int^{1}_{0}d\varrho\int_{\Sigma}e^{-2\tau\varphi(\log\varrho)}\left(|w_{\varrho}|^2+
 \varrho^{-2}\left\vert \nabla_{\Sigma}w\right\vert_{\Sigma}^2\right)\varrho^{\epsilon+1}d\Omega_{\varrho}\leq\\&
 \leq C\int^{1}_{0}d\varrho\int_{\Sigma}e^{-2\tau\varphi(\log\varrho)}\varrho^{3}|\mathcal{L} w|^2d\Omega_{\varrho}
 +\\&+C\int^{1}_{0}d\varrho\int_{\Sigma}e^{-2\tau\varphi(\log\varrho)}\varrho^{3}|\partial_{\varrho}w|^2d\Omega_{\varrho},
 \end{aligned}
\end{equation}
for every $\tau\geq \tau_2$ and for every $f\in
C_0^{\infty}((0,r_0)\times \Sigma)$, where $r_0=e^{T_3}$. Now,
let us choose $\delta=\frac{1}{2}$. It turns out that the last integral on the right--hand side of \eqref{final-6-g} can be absorbed by the third integral 
on the left--hand side. Finallly, applying Proposition \ref{Propos-8} and,
taking into account that $2^{-1/\epsilon}|x|\leq \rho(x)\leq |x|$ and
replacing $\tau$ by $(\tau-\frac{n}{2})$, we obtain
inequality \eqref{Carlm-delta-g}.$\blacksquare$

\bigskip

\begin{cor}\label{prop:Carlm-delta-Lip}
Let us assume that the entries of the matrix $G=\left\{g_{ij}\left(
x\right) \right\}_{i,j=1}^{n}$ are of class
$C^{0,1}\left(B_2\right)$ and that $G$ satisfies \eqref{2A},
\eqref{3A} and (instead of \eqref{4A}) satisfies 
\begin{equation}\label{4A-Lip}
\left|G(x)-G(y)\right|\leq\Lambda\left|x-y\right|, \quad\forall x,y\in B_{2},
\end{equation}
then Carleman estimate \eqref{Carlm-delta-g} continue to hold.

More precisely, there exist $C$ and
$\widehat{\tau}_{\ast}$, depending on $\lambda$ and $\Lambda$ only,
such that

\begin{equation}\label{Carlm-delta-g-Lip-5}
\begin{aligned}
\tau^{3}\int\rho^{\epsilon-2\tau}|u|^2dx&+
\tau\int\rho^{2+\epsilon-2\tau}|\nabla u|^2dx+\\&
+\tau^2r\int\rho^{-1-2\tau}u^2dx\leq C\int\rho^{4-2\tau}|\Delta_g
u|^2dx,
\end{aligned}
\end{equation}
for every $r\in (0,r_0)$, for every $\tau\geq \widehat{\tau}_{\ast}$ and for every $u\in C^\infty_0\left(B_{r_0}\setminus
\overline{B}_{r/4}\right)$.
\end{cor}

\textbf{Proof.} Let $\psi\in C_0^{\infty}\left(\mathbb{R}^n\right)$ satisfy supp
$\psi\subset B_1$, $\psi\geq 0$ and $\int_{\mathbb{R}^n}\psi dx =1$.
Set
$$\psi_{\nu}(x)=\nu^n\psi(\nu x),\quad\nu\in \mathbb{N}$$ and

$$G_{\nu}(x)=\left(G\star
\psi_{\nu}\right)(x)=\int_{\mathbb{R}^n}G(x-y)\psi_{\nu}(y)dy,\quad\nu\in
\mathbb{N}.$$ We have that $G_{\nu}$ satisfies to \eqref{2A},
\eqref{4A}, \eqref{3A}. Moreover $$G_{\nu}\in
C^{\infty}\left(B_2\right)$$ and

\begin{equation}\label{Carlm-delta-g-Lip-0}
\left\Vert G_{\nu}-G\right\Vert_{L^{\infty}(B_1)}\rightarrow
0,\quad\mbox{as } \nu\rightarrow \infty. \end{equation}

Let $r\in(0,r_0)$ and let $u$ be an arbitrary function belonging to
$C^\infty_0\left(B_{r_0}\setminus \overline{B}_{r/4}\right)$. By
\eqref{Carlm-delta-g} we have, for every $\nu\in \mathbb{N}$ and for every
 $\tau\geq \overline{\tau}$,
\begin{equation}\label{Carlm-delta-g-Lip}
\begin{aligned}
&\tau^{3}\int\rho^{\epsilon-2\tau}|u|^2dx+
\tau\int\rho^{2+\epsilon-2\tau}|\nabla u|^2dx+\\&
+\tau^2r\int\rho^{-1-2\tau}u^2dx\leq
C\int\rho^{4-2\tau}\left|\Delta_{g_{\nu}} u\right|^2dx.
\end{aligned}
\end{equation}
On the other hand, by \eqref{4A} we have (by using the convention of repeated index)

\begin{equation*}
\int\rho^{4-2\tau}\left|\Delta_{g_{\nu}} u\right|^2dx\leq
2\int\rho^{4-2\tau}\left|g^{ij}_{\nu} \partial^2_{x^i
x^j}u\right|^2dx+C\int\rho^{4-2\tau}\left|\nabla u\right|^2dx,
\end{equation*}
where $C$ \textbf{depends only on} $\lambda$ and $\Lambda$. By the just obtained inequality and by 
\eqref{Carlm-delta-g-Lip} we have, for every $\nu\in \mathbb{N}$,

\begin{equation}\label{Carlm-delta-g-Lip-2}
\begin{aligned}
&\tau^3\int\rho^{\epsilon-2\tau}|u|^2dx+
\int\rho^{2+\epsilon-2\tau}\left(\tau-C\rho^{2-\epsilon}\right)|\nabla
u|^2dx+\\& +\tau^2r\int\rho^{-1-2\tau}u^2dx\leq
C\int\rho^{4-2\tau}\left|g^{ij}_{\nu} \partial^2_{x^i
x^j}u\right|^2dx,
\end{aligned}
\end{equation}
where $C$ depends on $\lambda$ and $\Lambda$ only. Now let
$\overline{\tau}_{\ast}\geq\overline{\tau}$ satisfy  (recall that
$\rho\leq 1$ in $B_1$) for every $\tau\geq \overline{\tau}_{\ast}$,
$$\tau-C\rho^{2-\epsilon}\geq \frac{\tau}{2}.$$ By
\eqref{Carlm-delta-g-Lip-2} we have, for every $\nu\in \mathbb{N}$,
\begin{equation}\label{Carlm-delta-g-Lip-3}
\begin{aligned}
&\tau^3\int\rho^{\epsilon-2\tau}|u|^2dx+
\frac{\tau}{2}\int\rho^{2+\epsilon-2\tau}|\nabla u|^2dx+\\&
+\tau^2r\int\rho^{-1-2\tau}u^2dx\leq
C\int\rho^{4-2\tau}\left|g^{ij}_{\nu} \partial^2_{x^i
x^j}u\right|^2dx.
\end{aligned}
\end{equation}
Passing to the limit  as
$\nu\rightarrow\infty$ in \eqref{Carlm-delta-g-Lip-3} we obtain

\begin{equation*}
\begin{aligned}
&\tau^3\int\rho^{\epsilon-2\tau}|u|^2dx+
\frac{\tau}{2}\int\rho^{2+\epsilon-2\tau}|\nabla u|^2dx+\\&
+\tau^2r\int\rho^{-1-2\tau}u^2dx\leq C\int\rho^{4-2\tau}\left|g^{ij}
\partial^2_{x^i x^j}u\right|^2dx.
\end{aligned}
\end{equation*}
By the just obtained inequality, employing

\begin{equation*}
\int\rho^{4-2\tau}\left|g^{ij} \partial^2_{x_i x_j}u\right|^2dx\leq
2\int\rho^{4-2\tau}\left|\Delta_{g}
u\right|^2dx+C\int\rho^{4-2\tau}\left|\nabla u\right|^2dx
\end{equation*}
and repeating the arguments already used above, we have that there exists
$\widehat{\tau}_{\ast}\geq \overline{\tau}_{\ast}$, where
$\widehat{\tau}_{\ast}$ depends on $\lambda$ and $\Lambda$ only,
such that

\begin{equation*}
\begin{aligned}
&\tau^{3}\int\rho^{\epsilon-2\tau}|u|^2dx+
\tau\int\rho^{2+\epsilon-2\tau}|\nabla u|^2dx+\\&
+\tau^2r\int\rho^{-1-2\tau}u^2dx\leq C\int\rho^{4-2\tau}|\Delta_g
u|^2dx,
\end{aligned}
\end{equation*}
for every $r\in (0,r_0)$, for every $\tau\geq \widehat{\tau}_{\ast}$ and for every $u\in C^\infty_0\left(B_{r_0}\setminus
\overline{B}_{r/4}\right)$. $\blacksquare$

\bigskip

Now we can state the analog of Theorem \eqref{theo:40.teo} and
of Corollary \ref{SUCP} for the solutions $U\in
H^2\left(B_1\right)$ to the equation

\begin{equation}  \label{14-equ-fin}
LU=\sum_{i,j=1}^na^{ij}(x)\partial^2_{x^ix^j}U+\sum_{i=1}^nb^i(x)\partial_{x^i}U+c(x)U=0,\quad\mbox{in
} B_1,
\end{equation}
where $A(x)=\left\{a^{ij}(x)\right\}^n_{i,j=1}$ is a symmetric matrix whose entries are
real--valued functions, such that 
\begin{equation}  \label{4-2A-fin}
\lambda^{-1}\left|\xi\right|^{2}\leq
\sum_{i,j=1}^na^{ij}(x)\xi_i\xi_j\leq\lambda\left|\xi\right|^{2},\quad
\forall\xi\in\mathbb{R}^{n} \mbox{, }\forall x\in B_1,
\end{equation}
where $\lambda\geq 1$. Also we assume 
\begin{equation}\label{intro14-4An-fin}
\left|a^{ij}(x)-a^{ij}(y)\right|\leq\Lambda |x-y|, \quad\mbox{for
}i,j\in\left\{1,\cdots, n\right\}, \quad\forall x,y\in B_1.
\end{equation}
Moreover, $b^i\in L^{\infty}\left(B_1\right)$, $i=1,\cdots, n$ and
$c\in L^{\infty}\left(B_1\right)$ (these may be complex--valued coefficients) and set

\begin{equation}
    \label{limitaz-coeff-fin}
M=\max\left\{\left\Vert
b\right\Vert_{L^{\infty}\left(B_1,\mathbb{R}^n\right)},\left\Vert
c\right\Vert_{L^{\infty}\left(B_1\right)}\right\},
\end{equation}
where $b=\left(b^1,\cdots,b^n\right)$.

\begin{theo}
    \label{theo:40.teo-fin}
    Let us assume that $U\in H^2\left(B_{1}\right)$ is a solution to the equation \eqref{4-2A-fin}.
Let $x_0\in B_1$ and $0<R_0\leq 1-|x_0|$. Then there exist $C_1\geq
C\geq 1$ depending on $\lambda$, $\Lambda$ and $M$ only, such that,
if $0<r<\frac{R}{C}<\frac{R_0}{C_1}$ then

\begin{equation}
    \label{three-sphere-enunciato-fin}
    \int_{B_{R}(x_0)}U^2
    \leq
    C\left(\frac{R_0}{R}\right)^{C}\left(\int_{B_{r}(x_0)}U^2\right)^{\theta}\left(\int_{B_{R_0}(x_0)}U^2\right)^{1-\theta},
\end{equation}
where
\begin{equation}
    \label{three-sphere-5-fin}
\theta=\frac{\log\frac{R_0}{C_1R}}{\log\frac{C_1R_0}{r}}.
\end{equation}

Moreover, if $U$ does not vanish identically in $B_{R_0/C}(x_0)$
the following doubling inequality holds true 
\begin{equation}\label{eq:10.6.1102-cube-fin}
    \int_{B_{2r}(x_0)}U^2\leq C N_{x_0,R_0}^k\int_{B_{r}(x_0)}U^2,
\end{equation}
where
\begin{equation}
    \label{eq:10.6.1108-cube-fin}
    N_{x_0,R_0}=\frac{\int_{B_{R_0}(x_0)}U^2}{\int_{B_{R_0/C}(x_0)}U^2}
    \end{equation}
and $k$ is a positive number ($k\geq 3$).
\end{theo}
\textbf{Proof.} For fixed $x_0\in B_1$, since
$A(x_0)$ is a symmetric matrix, there exists a linear map 
$S$ such that
$$SA(x_0)S^T=I_n.$$ Hence, we may perform the change of variables  $y=S^{-1}(x-x_0)$ in the equation \eqref{14-equ-fin} that
allows to  apply, after some simple modifications, the Carleman estimate
\eqref{Carlm-delta-g-Lip-5} in a manner quite similar to what was done
in the proof of Theorem \eqref{theo:40.teo}. To obtain the analogon of Lemma \ref{newlemma1-cube}, we may fix, for instance $\epsilon=\frac{1}{2}$. We leave the details to the reader. $\blacksquare$

\medskip

Also, the following Corollary can be proved similarly to 
Corollary \ref{SUCP}

\begin{cor}[\textbf{strong unique continuation for elliptic equations}]\label{SUCP-fin}
	\index{Corollary:@{Corollary:}!- strong unique continuation for elliptic equations@{- strong unique continuation for elliptic equations}}
	
Let $U\in H^2\left(B_{1}\right)$ be a solution to equation
\eqref{4-2A-fin}. Let $x_0\in B_1$ and $0<R_0\leq 1-|x_0|$. There exists
$C$ depending on $\lambda$, $\Lambda$ and $M$ only such that we have what follows. If $U$ does not vanish identically in
$B_{R_0/C}(x_0)$ then we have, for every $r<s\leq\frac{R_0}{C}$,

\begin{equation}\label{SUCP-1-fin}
    \int_{B_{s}(x_0)}U^2\leq
    CN_{x_0,R_0}^k\left(\frac{s}{r}\right)^{\log_2(CN_{x_0,R_0}^k)}\int_{B_r(x_0)}U^2,
\end{equation}
where $N_{x_0,R_0}$ is defined by \eqref{eq:10.6.1108-cube-fin} and $k$ is the same number that occurs in \eqref{eq:10.6.1102-cube-fin}.

Moreover, if

\begin{equation}  \label{SUCP-2-fin}
\int_{B_r(x_0)}U^2=\mathcal{O}\left(r^m\right),\quad\mbox{as }
r\rightarrow 0,\quad \forall m\in \mathbb{N},
\end{equation}
then

\begin{equation}  \label{SUCP-tesi-fin}
U\equiv 0, \quad\mbox{in } B_1.
\end{equation}
\end{cor}

\chapter{Miscellanea}\label{Misc:27-11-22}
\section{Introduction}\label{Misc:27-11-22-1}
In this final Chapter we will first give (Section \ref{Misc:27-11-22-2}) a brief outline of two methods, alternative to the Carleman estimates, for dealing within the unique continuation issue. These methods are generally called the log -- convexity and the frequency function method. We will see that they are intimately related. Next, in Section \ref{Ap:7-12-22} , we will give a little mention of $A_p$ weights, pointing out some applications of them to inverse problems. In Section \ref{Rung:9-12-22} we will consider the Runge property for the Laplace operator. 

\section{The backward problem for the heat equation}\label{Misc:27-11-22-2-0}

Let us consider a rod of heat conducting material. Let $\pi$ be the length of the rod, let us assume that its temperature is zero at its extremes and that the heat flows only in the direction of the axis of the rod.
Let $u(x,t)$ be the temperature of the rod at the point $x$ and the time $t$. If the initial temperature is $f(x)$, then $u$ is a solution of the following + initial--boundary value  problem for the heat equation

\begin{equation}\label{Misc:27-11-22-2-01}
	\begin{cases}
		u_t-u_{xx}=0,\quad \mbox{for } (x,t)\in (0,\pi)\times (0,+\infty),\\
		\\
		u(0,t)=u(\pi,t)=0, \quad \mbox{for } t\in [0,+\infty),\\
		\\
		u(x,0)=f(x), \quad \mbox{for } x\in [0,\pi].
	\end{cases}
\end{equation}

Also, we assume 
\begin{equation}\label{Misc:28-11-22-2-1}
	f\in C^1([0,\pi]),\quad \mbox{and}\quad f(0)=f(\pi)=0. 
\end{equation} 
Let us recall that, \cite{We}, there exists an unique solution to problem \eqref{Misc:27-11-22-2-01} in the class $C^0([0,\pi]\times [0,+\infty))\cap C^2((0,\pi)\times (0,+\infty))$ and it is given by
\begin{equation}\label{Misc:28-11-22-2-2}
	u(x,t)=\sum_{k=1}^{\infty}f_k\sin kx e^{-k^2t},
\end{equation} 
where 
\begin{equation}\label{Misc:28-11-22-2-3}
	f_k=\frac{2}{\pi}\int^{\pi}_0f(x)\sin kx dx.
\end{equation} 

\medskip

Moreover we have
\begin{equation*}
	\int^{\pi}_0u^2(x,t)dx\leq \frac{\pi}{2}\int^{\pi}_0f^2(x)dx, \quad \forall t\geq 0
\end{equation*} 
and, more generally
\begin{equation*}
	\int^{\pi}_0|\partial^m_xu(x,t)|dx\leq C_{m,t}\int^{\pi}_0f^2(x)dx, \quad \forall t> 0.
\end{equation*} 
These inequalities imply a continuous dependence of the solution of problem \eqref{Misc:27-11-22-2-01} by the initial datum $f$.

In the \textbf{backward problem} \index{backward problem}we are interested in determining the temperature $u$, if we know, instead of initial temperature, the temperature at an instant $t>0$,  say $t=1$. Set
\begin{equation}\label{Misc:30-11-22-2-0}
	g(x)=u(x,1), \quad \mbox{in } [0,\pi].
\end{equation} 
It is evident that, for $t>1$ by the translation $t'=t-1$ we reduce to problem \eqref{Misc:27-11-22-2-01}. when $t<1$ we will examine what happens for what concerns  the uniqueness and continuous dependence of $u$ by $g$.

\medskip

\underline{\textbf{Uniqueness.}} By the linearity of the problem, it suffices to check that if $g\equiv 0$ then $u\equiv 0$. Now, since $u$ is given by \eqref{Misc:28-11-22-2-2}, we may consider the equation (of the unkwnown $f$)

\begin{equation}\label{Misc:30-11-22-2-1}
	\sum_{k=1}^{\infty}f_k\sin kx e^{-k^2}=0, \quad \forall x\in [0,\pi],
\end{equation} 
from which, multiplying both the sides by $\sin mx$, for $m\in \mathbb{N}$ and integrating over $[0,\pi]$, we obtain
$$\frac{\pi}{2}f_me^{-m^2}=0, \quad \forall m\in \mathbb{N}.$$
Therefore $f_m=0$ for every $m\in \mathbb{N}$. Hence 
$$u\equiv 0.$$ 

\medskip
 
\underline{\textbf{Continuous dependence and conditional stability.}}
Let us consider the sequence of functions
\begin{equation*}
g_{\nu}(x)=e^{\nu}\sin \pi\nu x, \quad \nu \in \mathbb{N}.
\end{equation*} 
It is easily checked that $$u_{\nu}(x,t)=e^{\nu^2(1-t)}e^{\nu}\sin \pi\nu x, \quad \nu \in \mathbb{N}, \quad \nu \in \mathbb{N}.$$
Hence
\begin{equation*}
	\left\Vert g^{(n)}_{\nu} \right\Vert_{L^2(0,\pi)}\rightarrow 0, \quad \mbox{as }\nu \rightarrow \infty, \quad \forall n\in \mathbb{N},
\end{equation*} 
($g^{(n)}_{\nu}$ is the $n$--th derivative of $g$), but

\begin{equation*}
	\left\Vert u_{\nu}(\cdot, t) \right\Vert_{L^2(0,\pi)}\rightarrow +\infty, \quad \mbox{as } \nu \rightarrow \infty, \quad \forall t\in (0,1).
\end{equation*} 
In plain words, even if we had the estimates of the error of all derivatives of the datum $g$ we could not control the error on $u(\cdot,t)$ when $t<1$.

\medskip

Let us denote by 

\begin{equation}\label{Misc:30-11-22-2-2}
	\varepsilon:=\left\Vert g \right\Vert_{L^2(0,\pi)}
\end{equation}
and let us suppose the temperature at the initial time is bounded (in the $L^2(0,\pi)$ norm) by a known constant. More precisely, we suppose that
\begin{equation}\label{Misc:30-11-22-2-3}
	\left\Vert f \right\Vert_{L^2(0,\pi)}=\left\Vert u(\cdot, 0) \right\Vert_{L^2(0,\pi)}\leq E,
\end{equation}
where $E>0$ is known. By \eqref{Misc:28-11-22-2-2} and \eqref{Misc:30-11-22-2-0} we have
 
\begin{equation*}
	\sum_{k=1}^{\infty}f_k\sin kx e^{-k^2}=g(x), \quad \forall x\in [0,\pi],
\end{equation*} 
from which we have

\begin{equation*}
	f_m=\frac{2}{\pi}g_me^{-m^2}, \quad \forall m\in \mathbb{N},
\end{equation*} 
where
\begin{equation*}
	g_m=\frac{2}{\pi}\int^{\pi}_0f(x)\sin mx dx,\quad \forall m\in \mathbb{N}.
\end{equation*} 
Therefore  \eqref{Misc:30-11-22-2-2} and  condition \eqref{Misc:30-11-22-2-3} are expressed, respectively, by
\begin{equation*}
	\frac{\pi}{2}\sum_{k=1}^{\infty}f_k^2e^{-2k^2}=\varepsilon^2
\end{equation*} 
and
\begin{equation*}
	\frac{\pi}{2}\sum_{k=1}^{\infty}f_k^2\leq E^2.
\end{equation*} 
On the other hand, we are interested in estimating  

\begin{equation*}
\left\Vert u(\cdot, t) \right\Vert^2_{L^2(0,\pi)}=	
\frac{\pi}{2}\sum_{k=1}^{\infty}f_k^2e^{-2k^2t}
\end{equation*} 
for $t\in (0,1)$.

Applying the H\"{o}lder inequality we get 

\begin{equation*}
\begin{aligned}
		\frac{\pi}{2}\sum_{k=1}^{\infty}f_k^2e^{-2k^2t}&= \frac{\pi}{2}\sum_{k=1}^{\infty}\left|f_k\right|^{2(1-t)}\left(\left|f_k\right|^2e^{-2k^2}\right)^t\leq \\&\leq \frac{\pi}{2}\left(\sum_{k=1}^{\infty}\left|f_k\right|^2\right)^{1-t}\left(\sum_{k=1}^{\infty}\left|f_k\right|^2e^{-2k^2}\right)^t\leq\\&\leq E^{2(1-t)}\varepsilon^{2t}.
\end{aligned}
\end{equation*}
Hence we have proved the following conditional stability estimate

\begin{equation}\label{Misc:30-11-22-2-4}
	\left\Vert u(\cdot, t) \right\Vert_{L^2(0,\pi)}\leq E^{1-t}\varepsilon^{t}, \quad \forall t\in (0,1).
\end{equation}

\bigskip

\textbf{Remark.} It is easily checked that estimate  \eqref{Misc:30-11-22-2-4} cannot be improved. Furthermore, \eqref{Misc:30-11-22-2-4} implies the log--convexity of the function

$$[0,\pi]\ni t\rightarrow \left\Vert u(\cdot, t) \right\Vert_{L^2(0,\pi)}$$ By this we mean that the function
$$F(t)=\log  \left\Vert u(\cdot, t) \right\Vert_{L^2(0,\pi)},$$ 
is convex. $\blacklozenge$

\section{The log--convexity method and the frequency function method}\label{Misc:27-11-22-2}
\subsection{The log--convexity method}\label{Misc:3-11-22-2-0}
At the base of the \textbf{log--convexity method} \index{log--convexity method}there are the following simple considerations.

Let us suppose that $F\in C^2([0,1])$ is a nonnegative function, $F''\geq 0$ in $[0,1]$ and let us suppose
\begin{equation}\label{Misc:27-11-22-3}
	F''(t)F(t)-F'^2(t)\geq 0, \quad \forall t\in [0,1].
\end{equation}
It is immediately checked that this inequality is equivalent to

$$F''(t)(F(t)+\gamma)-F'^2(t)\geq 0, \quad \forall t\in [0,1], \ \ \forall\gamma>0$$ 
which, in turn,  is equivalent to the log--convexity of $F+\gamma$ in $[0,1]$. As a matter of fact we have, 
$$\frac{d^2}{dt^2}\log (F(t) +\gamma)=\frac{F''(t)(F(t) +\gamma)-F'^2(t)}{F^2(t)}\geq 0.$$ Now, the log--convexity of $F+\gamma$ is equivalent to inequality
\begin{equation*}
	F(t)+\gamma\leq (F(0)+\gamma)^{1-t}(F(1)+\gamma)^{t}, \quad \forall t\in [0,1], \ \ \forall \gamma>0
\end{equation*}
hence 
\begin{equation}\label{Misc:27-11-22-4}
	F(t)\leq (F(0))^{1-t}(F(1))^{t}, \quad \forall t\in [0,1].
\end{equation}
By the latter inequality we derive that if  one of values $F(0)$, $F(1)$ is zero then $F$ vanishes identically. 

Let us at once see an application of the aforementioned idea for proving the uniqueness and a conditional stability estimate for the following backward problem

\begin{equation}\label{Misc:27-11-22-5}
	\begin{cases}
u_t-\left(a(x)u_x\right)_{x}=0,\quad \mbox{for } (x,t)\in (0,1)\times (0,T),\\
\\
u(0,t)=u(1,t)=0, \quad \mbox{for } t\in  [0,T],\\
\\
u(x,T)=g(x),  \quad \mbox{for } x\in [0,1]
	\end{cases}
\end{equation}
where $T>0$. Let us suppose that $a\in C^1([0,1])$  and that there exists $u\in C^2([0,1]\times [0,T])$ solution to \eqref{Misc:27-11-22-5}. 

We define

\begin{equation}\label{Misc:27-11-22-6}
	F(t)=\int^{1}_0 u^2(x,t)dx, \quad t\in [0,T]
\end{equation}
and we have
\begin{equation}\label{Misc:27-11-22-7}
	F'(t)=2\int^{1}_0 u(x,t)u_t(x,t)dx.
\end{equation}
Now, by the equation $u_t-\left(a(x)u_x\right)_{x}=0$, taking into account that  $u(0,t)=u(1,t)=0$ and integrating by parts, we get

\begin{equation*}
	2\int^{1}_0 u(x,t)u_t(x,t)dx=2\int^{1}_0 u(x,t)\left(a(x)u_x\right)_{x}dx=-2\int^{1}_0a(x)u^2_{x}(x,t)dx.
\end{equation*}
Hence
\begin{equation*}
	F'(t)=-2\int^{1}_0a(x)u^2_{x}(x,t)dx.
\end{equation*}
By using the just obtained equality, we calculate the second derivative of $F$  
\begin{equation*}
	F''(t)=-4\int^{1}_0a(x)u_{x}(x,t)u_{xt}(x,t)dx.
\end{equation*}
Now we integrate by parts, and we recall that $u_t(0,t)=u_t(1,t)=0$ (which is obtained by differentiating $u(0,t)=u(1,t)=0$ with respect to $t$), by using again the equation, we get
\begin{equation*}
	\begin{aligned}
	-4\int^{1}_0a(x)u_{x}(x,t)u_{xt}(x,t)dx&=4\int^{1}_0 \left(a(x)u_{x}(x,t)\right)_xu_{t}(x,t)dx=\\&=4\int^{1}_0u^2_{t}(x,t)dx.
		\end{aligned}
\end{equation*}
Hence 
\begin{equation}\label{Misc:27-11-22-8}
	F''(t)=4\int^{1}_0u^2_{t}(x,t)dx\geq 0.
\end{equation}
Now, by \eqref{Misc:27-11-22-6}, \eqref{Misc:27-11-22-7} and \eqref{Misc:27-11-22-8} we have 
\begin{equation*}
\begin{aligned}
F''(t)F(t)-F'^2(t)&=4\int^{1}_0u^2_{t}(x,t)dx\int^{1}_0u^2(x,t)dx-\\&-4\left(\int^{1}_0 u(x,t)u_t(x,t)dx\right)^2\geq 0,
\end{aligned}	
\end{equation*}
where the last inequality follows by the Cauchy--Schwarz inequality. Hence, by \eqref{Misc:27-11-22-4} we have

\begin{equation*}
	\int^{1}_0u^2(x,t)dx\leq \left(\int^{1}_0u^2(x,0)dx\right)^{1-\frac{t}{T}}\left(\int^{1}_0u^2(x,T)dx\right)^{\frac{t}{T}}, \quad \forall t\in [0,T]
\end{equation*}
from which, recalling $u(x,T)=g(x)$ in $[0,1]$,  
\begin{equation}\label{Misc:27-11-22-9}
	\int^{1}_0u^2(x,t)dx\leq \left(\int^{1}_0u^2(x,0)dx\right)^{1-\frac{t}{T}}\left(\int^{1}_0g^2(x)dx\right)^{\frac{t}{T}},
\end{equation}
for every $t\in [0,T]$. By the last estimate we obtain the uniqueness for the backward problem \eqref{Misc:27-11-22-5}. As a matter of fact, if $g\equiv 0$ then \eqref{Misc:27-11-22-9} implies $u\equiv 0$. Moreover, if we have the information

$$\int^{1}_0g^2(x)dx\leq \varepsilon^2, \quad \mbox{(error)}$$ and
$$\int^{1}_0u^2(x,0)dx\leq E^2, \quad \mbox{(a priori information)}$$ then we get the following conditional stability estimate

\begin{equation}\label{Misc:2-12-22-1-0}
	\left(\int^{1}_0u^2(x,t)dx\right)^{1/2}\leq E^{1-\frac{t}{T}}\varepsilon^{\frac{t}{T}},\quad \forall t\in [0,T].
\end{equation}

\bigskip

\textbf{Remark.} Unlike the method based on the Carleman estimates, in the log--convexity method, the equation and the initial and boundary data are used directly. We refer to \cite{Ag} and \cite{Pa} for more details on this topic. The elegance of the method and the simplicity of the proof that we have just given should not lead us to believe that the procedure does not have its asperities. To make a rough comparison with the method based on Carleman estimates, one could say that as, in the latter, the choice of weight is crucial (and non  trivial), in the log--convexity method, the choice of the function $F$ is crucial (and non  trivial). $\blacklozenge$

\bigskip 

We reconsider the backward problem with $a$ depending on $x$ and $t$. That is, we consider

  \begin{equation}\label{Misc:2-12-22-1}
  	\begin{cases}
  		u_t-\left(a(x,t)u_x\right)_{x}=0,\quad \mbox{for } (x,t)\in (0,1)\times (0,T),\\
  		\\
  		u(0,t)=u(1,t)=0, \quad \mbox{for } t\in [0,T] ,\\
  		\\
  		u(x,T)=g(x), \quad \mbox{for } x\in [0,1]
  	\end{cases}
  \end{equation}
where $a\in C^1([0,1]\times [0,+\infty))$. We denote

\begin{equation}\label{Misc:2-12-22-3-0}
\lambda=\min_{(x,t)\in [0,1]\times[0,T]}a(x,t)>0
\end{equation} 
and

\begin{equation}\label{Misc:2-12-22-3-01}
	M=\max_{(x,t)\in [0,1]\times[0,T]}\left|a_t(x,t)\right|.
\end{equation} 

We are searching for a function $\mu$, such that
\begin{equation}\label{Misc:2-12-22-2}
\mu:[0,1]\rightarrow [0,T],
\end{equation}       
bijective, increasing, which satisfies $\mu \in C^2([0,1])$ and such that 
$$[0,1]\ni s\rightarrow \Phi(s):=F(\mu(s)),$$ is log--convex, where $F$ is given by  

\begin{equation*}
	F(t)=\int^{1}_0 u^2(x,t)dx, \quad t\in [0,T].
\end{equation*}
We get
\begin{equation*}
	\overset{\cdot }{\Phi }(s)=F'(\mu(s))\overset{\cdot }{\mu }(s), \quad \forall s\in [0,1],
\end{equation*}
where $\overset{\cdot }{\Phi}$ denotes the derivative of $\Phi$ w.r.t. $s$. In addition
\begin{equation}\label{Misc:2-12-22-3}
	\overset{\cdot \cdot}{\Phi }(s)=F''(\mu(s))\overset{\cdot }{\mu }^2(s)+F'(\mu(s))\overset{\cdot\cdot }{\mu }(s), \quad \forall s\in [0,1].
\end{equation}
Hence
\begin{equation}\label{Misc:2-12-22-4}
	\begin{aligned}
	\overset{\cdot \cdot}{\Phi }(s)\Phi(s)-	\overset{\cdot }{\Phi}^2(s)&=\left(	F''(\mu(s))F(\mu(s))-F'^2(\mu(s))\right)\overset{\cdot}{\mu}^2(s)(s)+\\&+
	F'(\mu(s))F(\mu(s))\overset{\cdot\cdot }{\mu }(s), \quad \quad \forall s\in [0,1].
	\end{aligned}
\end{equation}       
Now, by the equation $u_t=\left(a(x,t)u_x\right)_{x}$, we get
\begin{equation*}
	\begin{aligned}
	F'(t)&=2\int^{1}_0 u(x,t)u_t(x,t)dx=\\&=2\int^{1}_0 u(x,t)\left(a(x,t)u_x(x,t)\right)_{x}dx=\\&=-2\int^{1}_0 a(x,t)u_x^2(x,t)dx
	\end{aligned}
\end{equation*}
and
\begin{equation*}
	\begin{aligned}
		F''(t)&=-4\int^{1}_0 a(x,t)u_x(x,t)u_{xt}(x,t)dx-2\int^{1}_0 a_t(x,t)u_x^2(x,t)dx=\\&=
		4\int^{1}_0 u_t(x,t)\left(a(x,t)u_x(x,t)\right)_{x}dx-2\int^{1}_0 a_t(x,t)u_x^2(x,t)dx=\\&=
		4\int^{1}_0 u^2_t(x,t)dx-2\int^{1}_0 a_t(x,t)u_x^2(x,t)dx.
	\end{aligned}
\end{equation*}
Recalling \eqref{Misc:2-12-22-3-01}, we have

\begin{equation}\label{Misc:2-12-22-5}
	F''(t)\geq 4\int^{1}_0 u^2_t(x,t)dx-2M\int^{1}_0 u_x^2(x,t)dx, \quad \forall t\in [0,T].
\end{equation}      
 By the last equality and by \eqref{Misc:2-12-22-4} we get (for the sake of brevity, we omit the variables)
 
 \begin{equation*}
 	\begin{aligned}
 		\overset{\cdot \cdot}{\Phi }(s)\Phi(s)-	\overset{\cdot }{\Phi}^2(s)&\geq \left[\left(4\int^{1}_0 u^2_tdx-2M\int^{1}_0 u_x^2dx\right)\left(\int^{1}_0 u^2dx\right)-\right.\\&\left.-4\left(\int^{1}_0 u_tudx\right)^2\right]\overset{\cdot}{\mu}^2(s)+\\&+
 		2\left(\int^{1}_0 u_tudx\right)\left(\int^{1}_0 u^2dx\right)\overset{\cdot\cdot}{\mu}(s)=\\&=
 		4\left[\left(\int^{1}_0 u^2_tdx\right)\left(\int^{1}_0 u^2dx\right)-\left(\int^{1}_0 u_tudx\right)^2\right]\overset{\cdot}{\mu}^2(s)+\\&+
 		\left[-4M\overset{\cdot}{\mu}^2(s)\int^{1}_0 u_x^2dx+2\left(\int^{1}_0 u_tudx\right)\overset{\cdot\cdot}{\mu}(s)\right]\left(\int^{1}_0 u^2dx\right).
 		\end{aligned}
 \end{equation*}
By applying the the Cauchy--Schwarz inequality to the expression in the first square bracket, we obtain

\begin{equation}\label{Misc:2-12-22-6}
	\begin{aligned}
		\overset{\cdot \cdot}{\Phi }(s)\Phi(s)-	\overset{\cdot }{\Phi}^2(s)&\geq 
		\left[-4M\overset{\cdot}{\mu}^2(s)\int^{1}_0 u_x^2dx+\right. \\&
		\left.	+2\left(\int^{1}_0 u_tudx\right)\overset{\cdot\cdot}{\mu}(s)\right]\left(\int^{1}_0 u^2dx\right).
	\end{aligned}
\end{equation}
On the other hand
\begin{equation*}
	\int^{1}_0 u_tudx=\int^{1}_0\left(a(x,t)u_x\right)_{x}udx=-\int^{1}_0a(x,t)u^2_xdx\geq -\lambda \int^{1}_0u^2_xdx.
\end{equation*}
Now, proposing to find $\mu$ \textbf{concave}, by the last inequality and by \eqref{Misc:2-12-22-6} we have  
 
\begin{equation*}
		\overset{\cdot \cdot}{\Phi }(s)\Phi(s)-	\overset{\cdot }{\Phi}^2(s)\geq 
		-2\left[2M\overset{\cdot}{\mu}^2(s)+\lambda\overset{\cdot\cdot}{\mu}(s)\right]\left(\int^{1}_0 u^2dx\right)\left(\int^{1}_0 u_x^2dx\right).
\end{equation*}
Hence, in order to  
\begin{equation}\label{Misc:2-12-22-7}
\overset{\cdot \cdot}{\Phi }(s)\Phi(s)-\overset{\cdot }{\Phi}^2(s)\geq 0
\end{equation}
it suffices that $\mu$ satisfies the following conditions
\begin{equation}\label{Misc:2-12-22-8}
	\begin{cases}
	2M\overset{\cdot}{\mu}^2(s)+\lambda\overset{\cdot\cdot}{\mu}(s)\leq 0, \quad \forall s\in [0,T],\\
	\\
	\overset{\cdot}{\mu}(s)\geq 0, \quad \quad \ \ \ \ \ \ \ \ \ \ \ \ \ \  \forall s\in [0,T],\\
	\\
	\overset{\cdot\cdot}{\mu}(s)\leq 0, \quad \quad \ \ \ \ \ \ \ \ \ \ \ \ \ \  \forall s\in [0,T],\\
	\\
	\mu(0)=0,\quad \mu(1)=T.
	\end{cases}
\end{equation}
 Let us notice that the condition  $\overset{\cdot\cdot}{\mu}\leq 0$ implies (see \eqref{Misc:2-12-22-3}) $\overset{\cdot\cdot}{\Phi}\geq 0$. 
 
 Setting
 $$\alpha=\frac{2M}{\lambda},$$
 it is simple to check that 
 \begin{equation}\label{Misc:2-12-22-9}
 	\mu(s)=T+\frac{1}{\alpha}\log\left[e^{-\alpha T} +\left(1-e^{-\alpha T}\right)\right], \quad s\in [0,1],
 \end{equation}
 satyisfies all conditions \eqref{Misc:2-12-22-8}. In particular we get
 $$2M\overset{\cdot}{\mu}^2(s)+\lambda\overset{\cdot\cdot}{\mu}(s)= 0, \quad \forall s\in [0,T].$$ 
All in all, if we have

$$\int^{1}_0g^2(x)dx\leq \varepsilon^2$$ and
$$\int^{1}_0u^2(x,0)dx\leq E^2, $$ we obtain the conditional stability estimate 

\begin{equation}\label{Misc:2-12-22-10}
	\left(\int^{1}_0u^2(x,t)dx\right)^{1/2}\leq E^{1-\mu^{-1}(t)}\varepsilon^{\mu^{-1}(t)},\quad \forall t\in [0,T].
\end{equation}
where 
\begin{equation*}
\mu^{-1}(t)=\frac{e^{-\alpha (T-t)}-e^{-\alpha t}}{1-e^{-\alpha T}}.
\end{equation*}
Let us note that as $\alpha$ goes to $0$ (corresponding to the case in which $a$ does not depend on $t$) $\mu^{-1}(t)$ goes to $\frac{t}{T}$, i.e. the  exponent of the estimate \eqref{Misc:2-12-22-1-0}.

\bigskip

The log--convexity method can also be applied to prove the uniqueness and conditional stability for the Cauchy problem. Perhaps, the first author to use it was M. M. Lavrent'ev in 1956, \cite{Lav}. He applied the method to the Cauchy problem for the Laplace equation in a convex region. With some minor simplification, the situation considered is as follows (we consider only the uniqueness) 

\begin{equation}\label{Misc:2-12-22-11}
	\begin{cases}
	u_{yy}(x,y)+\Delta_{x}u(x,y)=0, \quad \forall(x,y)\in B_1\times (0,1),\\
	\\
	u(x,y)=0	\quad \forall(x,y)\in \partial B_1\times [0,1],\\
	\\
	u(x,0)=u_y(x,0)=0, \quad \forall x\in \overline{B_1},
	\end{cases}
\end{equation}    
 where 
 $$\Delta_{x}u(x,y)=\sum_{j=1}^nu_{x_jx_j}(x,y)$$ and we suppose that $u\in C^2\left(\overline{B_1}\times [0,1]\right)$.
 
 Set
 \begin{equation}\label{Misc:2-12-22-13}
 	F(y)=\int_{B_1}u^2(x,y)dx.
 \end{equation}
 We have
 \begin{equation}\label{Misc:2-12-22-14}
 	F'(y)=2\int_{B_1}u(x,y)u_y(x,y)dx
 \end{equation}
 and
\begin{equation}\label{Misc:2-12-22-15}
	F''(y)=2\int_{B_1}\left(u_y^2(x,y)+u(x,y)u_{yy}(x,y)\right)dx.
\end{equation} 

\medskip

Now, let us prove 
\begin{equation}\label{Misc:2-12-22-16}
	\int_{B_1}u(x,y)u_{yy}(x,y)dx=\int_{B_1}u_y^2(x,y)dx.
\end{equation}

\medskip

First, we note that, by the equation $u_{yy}+\Delta_{x}u=0$ and by the condition $u(x,y)=0$ on $\partial B_1\times [0,1]$ we have 
  
\begin{equation}\label{Misc:2-12-22-17}
	\begin{aligned}
		\int_{B_1}u(x,y)u_{yy}(x,y)dx&=-\int_{B_1}\Delta_{x}u(x,y)u(x,y)dx=\\&=\int_{B_1}|\nabla_x u(x,y)|^2dx.	\end{aligned}
	\end{equation}
 On the other hand 
 
 \begin{equation*}
 	\begin{aligned}
 		\frac{d}{dy}\int_{B_1}u_y^2dx&=2\int_{B_1}u_{yy}u_ydx=\\&=
 		-2\int_{B_1}(\Delta_{x}u)u_ydx=\\&=
 		-2\int_{B_1}\left[\mbox{div}_x\left(\nabla_xu u_y\right)-\nabla_x u\cdot \nabla_x u_y\right]dx=\\&=
 		2\int_{B_1}\nabla_x u\cdot \nabla_x u_ydx=
 		\frac{d}{dy}\int_{B_1}|\nabla_x u(x,y)|^2dx.
 	\end{aligned}
 \end{equation*}
 By the just obtained equality and by  \eqref{Misc:2-12-22-17} we get
 
 \begin{equation}\label{correct:30-4-23}
 		\frac{d}{dy}\left(\int_{B_1}u_y^2(x,y)dx-\int_{B_1}|\nabla_x u(x,y)|^2dx\right)=0.
 \end{equation}
 Now, since we have $$u(x,0)=u_y(x,0)=0$$, we get 
 \begin{equation*}
 	\int_{B_1}u_y^2(x,0)dx-\int_{B_1}|\nabla_x u(x,0)|^2dx=0
 \end{equation*}
and by \eqref{correct:30-4-23} we have
 
 \begin{equation*}
 \int_{B_1}u_y^2(x,y)dx=\int_{B_1}|\nabla_x u(x,y)|^2dx, \ \ \forall y\in (0,1).
 \end{equation*}
By the latter and by \eqref{Misc:2-12-22-17} we have \eqref{Misc:2-12-22-16} which, in turn (recall 
\eqref{Misc:2-12-22-15}), gives

\begin{equation*}
	F''(y)=4\int_{B_1}u_y^2(x,y)dx.
\end{equation*} 
By the just obtained equality, by \eqref{Misc:2-12-22-13},  \eqref{Misc:2-12-22-14} and by the Cauchy--Schwarz inequality we have 
\begin{equation*}
	F''(y)F(y)-F'^2(y)=4\left(\int_{B_1}u_y^2dx\right)\left(\int_{B_1}u^2dx\right)-4\left(\int_{B_1}uu_ydx\right)^2\geq 0.
\end{equation*} 
Hence $F(y)$ is a log--convex function, consequently we have
\begin{equation*}
	\int_{B_1}u^2(x,y)dx\leq \left(\int_{B_1}u^2(x,0)dx\right)^y\left(\int_{B_1}u^2(x,1)dx\right)^{1-y}=0
\end{equation*} 
for every $y\in [0,1]$. 

Therefore
$$u\equiv 0.$$

\subsection{The frequency function method}\label{Misc:3-12-22-1}

In Section \ref{esempi-stime-errore} we saw various versions of the Hadamard three circle inequality for the holomorphic functions and for the harmonic functions in two variables (for the latter, see \eqref{3-12-22-s36n-1-0} and \eqref{3-12-22-s36n-1}).
In Chapter \ref{tre sfere-ellittiche}, using the Carleman estimates, we have extended this inequality to the solutions of the second-order elliptic equations and we proved the doubling inequality. The frequency function method was used for the first time in \cite{Ga_Li} for the second-order elliptic equations with variable coefficients. Here, in order to present the main ideas of this method, we consider the case of the Laplace equation. 

\medskip

Let us propose to prove inequality \eqref{3-12-22-s36n-1-0} for the harmonic functions of $n$ variables.  Then, let $u$ be a solution of the Laplace equation.

\begin{equation}\label{laplace}
	\Delta u=0, \quad  \mbox{ in } B_{R_0}\subset \mathbb{R}^n.
\end{equation}
The inequality we are interested in can be written    
\begin{equation}\label{3-12-22-1}
	H\left(r_2\right)\leq
	\left(H\left(r_1\right)\right)^{\theta_0}\left(H\left(r_3\right)\right)^{1-\theta_0},
\end{equation}
for $0<r_1<r_2< r_3<R_0$, where
\begin{equation}\label{Misc:I-H}
H(r)=\int_{\partial B_r}u^2dS,  
\end{equation}
and
\begin{equation}\label{Misc:theta}
	\vartheta=\frac{\log\left(\frac{r_3}{r_2}\right)}{\log\left(\frac{r_3}{r_1}\right)}.
\end{equation}
Now it is important to observe that \eqref{3-12-22-1} is equivalent to the log--convexity of 
 
\begin{equation}\label{Misc:3-12-22-2}
	F(t)=H\left(e^t\right), \quad t\in \left(-\infty,\log R_0\right).
\end{equation}
 Therefore the derivative of $\log F(t)$ needs to be increasing. Let us calculate such a derivative
 \begin{equation*}
 	\frac{d}{dt}\log F(t)=\frac{F'(t)}{F(t)}=\frac{e^tH'\left(e^t\right)}{H\left(e^t\right)}, \quad t\in \left(-\infty,\log R_0\right).
 \end{equation*}
This equality implies that \eqref{3-12-22-1} is equivalent to the fact that the function
\begin{equation}\label{Misc:3-12-22-3-0}
	\frac{rH'\left(r\right)}{H\left(r\right)}, \quad r\in \left(0,R_0\right),
\end{equation}
is increasing. Now we calculate $H'\left(r\right)$. First we notice that

\begin{equation}\label{Misc:3-12-22-3}
	H(r)=\frac{1}{r}\int_{B_r} \mbox{div}(x u^2)dx,
\end{equation}
as a matter of fact, by the divergence Theorem we obtain
\begin{equation*}
	H(r)=\frac{1}{r}\int_{\partial B_r} \left(x\cdot\frac{x}{|x|}\right) u^2dS=\frac{1}{r}\int_{\partial B_r} (x\cdot\nu) u^2dS=\frac{1}{r}\int_{B_r} \mbox{div}(x u^2)dx.
\end{equation*}
Now, set
\begin{equation}\label{Misc:3-12-22-4}
	I(r)=\int_{B_r}|\nabla u|^2dx, \quad r\in (0,R_0) 
\end{equation}
and let us notice that
\begin{equation}\label{Misc:4-12-22-1}
	I(r)=\int_{B_r}\nabla u\cdot \nabla u dx=\int_{B_r}\mbox{div} \left(u\nabla u\right) dx= \int_{\partial B_r}u\frac{\partial u}{\partial\nu}dS.
\end{equation}
 We have

\begin{equation*}
	\begin{aligned}
		&H'(r)=\frac{1}{r}\int_{\partial B_r} \mbox{div}(x u^2)dS -\frac{1}{r^2}\int_{B_r} \mbox{div}(x u^2)dx\\
		&=\frac{1}{r}\left\{n\int_{\partial B_r} u^2dS+2\int_{\partial B_r} (x\cdot \nabla u) udS\right\}-\frac{1}{r}H(r)\\
		&=\frac{n-1}{r}\int_{\partial B_r} u^2dS+2\int_{\partial B_r}u\frac{\partial u}{\partial\nu}dS=\frac{n-1}{r}H(r)+2I(r).
	\end{aligned}
\end{equation*}
Therefore

\begin{equation}\label{Misc:3-12-22-5}
	H'(r)=\frac{n-1}{r}H(r)+2I(r), \quad r\in (0,R_0).
\end{equation}
Hence, function \eqref{Misc:3-12-22-3-0} can be written
\begin{equation*}
\frac{rH'\left(r\right)}{H\left(r\right)}=n-1+2\frac{r I(r)}{H(r)}, \quad r\in (0,R_0)
\end{equation*} 
In sum, the log--convexity  of $F$ is equivalent to the fact that the function\begin{equation}\label{Misc:3-12-22-6}
	N(r)=\frac{r I(r)}{H(r)}, \quad r\in (0,R_0)
\end{equation}
is increasing. $N(\cdot)$ is called the \textbf{frequency function} of $u$. 
 The frequency function was introduced by Almgren (1977), \cite{ALM} and has yaken this name because for a homogeneous harmonic polynomial of degree $m$, we have \\ $N(r)=m$ for all $r$. For instance, in dimension $2$, the homogeneous harmonic polynomials of degree $m$ are (in polar coordinates) of the type $$p_m(\varrho,\theta)=A\varrho^m \cos m\phi+B\varrho^m \sin m\phi$$ and it is easy to verify what has been asserted.
 
 In the sequel, we will assume that $H(r)>0$ for every $r\in (0,R_0)$, otherwise we can employ the device shown at the beginning of Section \ref{Misc:3-11-22-2-0}. We can also notice that if there exists $\overline{r}\in (0,R_0)$ such that $H\left(\overline{r}\right)=0$ then $u=0$ on $\partial B_r$ from which, $u=0$ in $B_r$ and by the unique continuation property we have $u\equiv 0$, which make the \eqref{3-12-22-1} trivial. The second device (i.e., assuming a unique continuation property) is certainly legitimate, but it is somewhat reductive because by the method we are illustrating we can obtain independent proof of unique continuation property.
 
 \bigskip 
   
At this point, we state and prove the following.
\begin{prop}\label{N-cresce}
	If $u$ is a non identically zero solution to \eqref{laplace} then $N(r)$ is an increasing function.
\end{prop}

\textbf{Proof.} 
First we prove that 

\begin{equation}\label{I'}
	I'(r)=2\int_{\partial B_r}\left(\frac{\partial u}{\partial\nu}\right)^2 dS+\frac{n-2}{r}I(r).
\end{equation} 
For this purpose we use the Rellich identity (Lemma \ref{rellich}) for $v \in C^2 \left(B_{R_0}\right)$ and $\beta(x)=x$,
\begin{equation}\label{Misc: 4-12-22-rellich}
	2 (x\cdot \nabla v) \Delta v = \mbox{div} [2 (x\cdot \nabla v) \nabla v - x |\nabla v|^2 ] +(n-2)|\nabla v|^2.
\end{equation}

We get by \eqref{Misc:3-12-22-4}

\begin{equation}\label{I'-1}
	I'(r)=\int_{\partial B_r}|\nabla u|^2dS.
\end{equation}

Now we apply identity \eqref{Misc: 4-12-22-rellich} to the function $u$, we integrate both the sides of \eqref{Misc: 4-12-22-rellich} over $B_r$ and we recall that $\Delta u=0$, obtaining
\begin{equation*}
	\begin{aligned}
		&0=\int_{\partial B_r}\left\{2 (x\cdot \nabla u) \nabla u\cdot\frac{x}{|x|} - x\cdot \frac{x}{|x|} |\nabla u|^2\right\}dS+ (n-2)\int_{B_r} |\nabla u|^2dx\\
		&=\int_{\partial B_r}\left\{2 r \left(\frac{\partial u}{\partial \nu}\right)^2 - r |\nabla u|^2\right\}dS+ (n-2)\int_{B_r} |\nabla u|^2dx,
	\end{aligned}
\end{equation*}
from which, taking into account \eqref{I'-1}, we have
\begin{equation}\label{I'-2}
	\begin{aligned}
		&I'(r)=\int_{\partial B_r}|\nabla u|^2dS=2\int_{\partial B_r} \left(\frac{\partial u}{\partial \nu}\right)^2dS+ \frac{n-2}{r}\int_{B_r} |\nabla u|^2dx\\
		&=2\int_{\partial B_r}\left(\frac{\partial u}{\partial\nu}\right)^2 dS+\frac{n-2}{r}I(r).
	\end{aligned}
\end{equation}
hence \eqref{I'} is proved.

\bigskip

Differentiating both the sides of \eqref{Misc:3-12-22-6}, we get
\begin{equation}\label{N'}
	N'(r)=\frac{I(r)}{H(r)}+r\frac{I'(r)}{H(r)}-r\frac{I(r)H'(r)}{H^2(r)}.
\end{equation}
Unless $u$ is constant (and, therefore \eqref{3-12-22-1} becomes trivial) we have $$N(r)>0$$ and we divide both the sides of \eqref{N'} by $N(r)$, obtaining

\begin{equation*}
	\frac{N'(r)}{N(r)}=\frac{1}{r}+\frac{I'(r)}{I(r)}-\frac{H'(r)}{H(r)},
\end{equation*}
from which, taking into account \eqref{Misc:3-12-22-5} and \eqref{I'}, we have

\begin{equation}\label{N'-2}
	\begin{aligned}
		&\frac{N'(r)}{N(r)}=\frac{1}{r}+\frac{I'(r)}{I(r)}-\frac{\frac{n-1}{r}H(r)+2I(r)}{H(r)}=\\
		&=\frac{2-n}{r}+\frac{I'(r)}{I(r)}-2\frac{I(r)}{H(r)}=\\
		&=2\left\{\frac{\int_{\partial B_r}\left(\frac{\partial u}{\partial\nu}\right)^2 dS}{I(r)}-\frac{I(r)}{H(r)}\right\}.
	\end{aligned}
\end{equation}
Now we express $I(r)$ by means of \eqref{Misc:4-12-22-1} and, by \eqref{N'-2}, we get

\begin{equation*}
	\begin{aligned}
		&\frac{N'(r)}{N(r)}=2\left\{\frac{\int_{\partial B_r}\left(\frac{\partial u}{\partial\nu}\right)^2 dS}{\int_{\partial B_r}\frac{\partial u}{\partial\nu}u dS}-\frac{\int_{\partial B_r}\frac{\partial u}{\partial\nu}u dS}{\int_{\partial B_r}u^2 dS}\right\}=\\
		&=\frac{2}{I(r)H(r)}\left\{\left(\int_{\partial B_r}\left(\frac{\partial u}{\partial\nu}\right)^2 dS\right)\left(\int_{\partial B_r}u^2 dS\right)-\left(\int_{\partial B_r}\frac{\partial u}{\partial\nu}u dS\right)^2\right\}.
	\end{aligned}
\end{equation*}
Now, by the Cauchy--Schwarz inequality we have 

\begin{equation*}
	\left(\int_{\partial B_r}\left(\frac{\partial u}{\partial\nu}\right)^2 dS\right)\left(\int_{\partial B_r}u^2 dS\right)-\left(\int_{\partial B_r}\frac{\partial u}{\partial\nu}u dS\right)^2\geq 0
\end{equation*}
from which the thesis follows.
$\blacksquare$

\bigskip

Let us summarize what has been obtained so far. By Proposition \ref{N-cresce} we get the log--convexity of function \eqref{Misc:3-12-22-2}, from which inequality \eqref{3-12-22-1} follows. By the latter, proceeding as done to prove \eqref{s36n-1} and using by \eqref{3-12-22-1}, we obtain, 

\begin{equation}\label{Misc:6-12-22-1}
	\int_{B_{r_2}}u^2dx\leq
	\left(\int_{B_{r_1}}u^2dx\right)^{\theta_0}\left(\int_{B_{r_3}}u^2dx\right)^{1-\theta_0},
\end{equation}
for $0<r_1<r_2< r_3<R_0$, where $\theta$ is given by \eqref{Misc:theta}. 
    
\bigskip

By the properties of the frequency function we derive a doubling inequality. Indeed the following holds true.

\begin{prop}\label{Misc:6-12-22-2}
Let $u$ be a non identically zero solution to \eqref{laplace} then
\begin{equation}\label{Misc:6-12-22-2-0}
\int_{B_{2r}}u^2dx\leq 2^{2^n N(R_0)+1}\int_{B_{r}}u^2dx, \quad \forall r\in \left(0,\frac{R_0}{2}\right]
\end{equation}
and
\begin{equation}\label{Misc:6-12-22-2-1}
	\int_{B_{2r}}u^2dx\leq \frac{2\int_{B_{R_0}}u^2dx}{\int_{B_{R_0/4}}u^2dx} \int_{B_{r}}u^2dx,\quad \forall r\in \left[0,\frac{R_0}{4}\right],
\end{equation}
\end{prop}

\bigskip

In the proof of Proposition \ref{Misc:6-12-22-2} we use

\begin{lem}\label{lemma:6-12-22-1}
Let $u$ be a solution to \eqref{laplace} then we have

\begin{equation}\label{lemma:6-12-22-2}
\frac{1}{r}\int_{B_r}u^2dx\leq \int_{\partial  B_r}u^2dS\leq \frac{C}{r}\int_{B_{2r}}u^2dx, \quad r\in \left(0,\frac{R_0}{2}\right],
\end{equation}
where $C$ depends on $n$ only.
\end{lem}
\textbf{Proof of Lemma.} Let us prove the first inequality of \eqref{lemma:6-12-22-2}. By \eqref{Misc:3-12-22-5} we have $H'(r)\geq 0$. Hence 

\begin{equation}
	\begin{aligned}
	\int_{B_r}u^2dx&=\int^{r}_0\left(\int_{\partial B_t}u^2dS\right)dt\leq \\&\leq 
	\int^{r}_0\left(\int_{\partial B_r}u^2dS\right)dt=\\&=r\int_{\partial B_r}u^2dS.
	\end{aligned}
	\end{equation}

Concerning the second inequality of \eqref{lemma:6-12-22-2}, by \eqref{Misc:3-12-22-3} and by the  Caccioppoli inequality, \eqref{dis-Caccioppoli}, we have
\begin{equation*}
	\begin{aligned}
	H(r)&=\frac{1}{r}\int_{B_r} \mbox{div}(x u^2)dx=\\&=
\frac{n}{r}\int_{B_r} u^2dx+\frac{1}{r}\int_{B_r}u\left(x\cdot\nabla u\right)dx\leq \\&\leq
\frac{n}{r}\int_{B_r} u^2dx+\left(\int_{B_r} u^2dx\right)^{1/2}\left(\int_{B_r} |\nabla u|^2dx\right)^{1/2}\leq \\&\leq \frac{n}{r}\int_{B_r} u^2dx+\frac{\widetilde{C}}{r}\left(\int_{B_r} u^2dx\right)^{1/2}\left(\int_{B_{2r}} u^2dx\right)^{1/2}\leq \\&\leq \frac{\widetilde{C}+n}{r}\int_{B_{2r}} u^2dx,
\end{aligned}
\end{equation*}
where $\widetilde{C}$ depends on $n$ only. The second inequality in \eqref{lemma:6-12-22-2} is proved with $C=\widetilde{C}+n$. $\blacksquare$

\bigskip

\textbf{Proof of Proposition \ref{Misc:6-12-22-2}.} By \eqref{Misc:3-12-22-5} and \eqref{Misc:3-12-22-6} we get
\begin{equation}\label{Misc:6-12-22-3}
\frac{d}{dr}\left(\log\frac{H(r)}{r^{n-1}}\right)=\frac{2N(r)}{r}.
\end{equation}
Let $\rho\in \left(0,\frac{R_0}{2}\right]$. Integrating both the sides of \eqref{Misc:6-12-22-3} over $[\rho,2\rho]$ and recalling that $N$ is increasing, we have

\begin{equation*}
	\log\left(\frac{H(2\rho)}{2^{n-1}H(\rho)}\right)=2\int^{2\rho}_{\rho}\frac{2N(r)}{r}\leq (2\log2)N(R_0).
\end{equation*}
Hence

\begin{equation*}
H(2\rho)\leq 2^{2^n N(R_0)}H(\rho), \quad\forall \rho\in \left(0,\frac{R_0}{2}\right], 
\end{equation*}
Integrating over $[0,r]$, $r\in \left(0,\frac{R_0}{2}\right]$, we get \eqref{Misc:6-12-22-2-0}.

\medskip 

Let now prove \eqref{Misc:6-12-22-2-1}. Let $\rho\in \left(0,\frac{R_0}{4}\right]$ and $R\in\left[\frac{R_0}{4},\frac{R_0}{2}\right]$. Integrating both the sides of \eqref{Misc:6-12-22-3} over $[\rho,2\rho]$ and recalling that $N$ is an increasing function, we have
\begin{equation*}
	\log\left(\frac{H(2\rho)}{2^{n-1}H(\rho)}\right)\leq (2\log2)N(R).
\end{equation*}
By the just obtained inequality, taking into account \eqref{Misc:6-12-22-3}, we have
\begin{equation*}
\frac{1}{R}	\log\left(\frac{H(2\rho)}{2^{n-1}H(\rho)}\right)\leq (2\log2)\frac{N(R)}{R}=(\log2)\frac{d}{dR}\left(\log\frac{H(R)}{R^{n-1}}\right).
\end{equation*}
Again we integrate both the sides of the last inequality w.r.t. $R$ over $\left[\frac{R_0}{4},\frac{R_0}{2}\right]$ so that we have 
\begin{equation*}
	\log\left(\frac{H(2\rho)}{2^{n-1}H(\rho)}\right)\leq \log\left(\frac{H\left(\frac{R_0}{2}\right)}{2^{n-1}H\left(\frac{R_0}{4}\right)}\right).
\end{equation*}
Now, Lemma \ref{lemma:6-12-22-1} gives 
\begin{equation*}
H(2\rho)\leq \frac{\int_{B_{R_0}}u^2dx}{\int_{B_{R_0/4}}u^2dx} H(\rho),\quad \forall \rho\in \left(0,\frac{R_0}{4}\right].
\end{equation*}
From which integration of both the sides over $[0,r]$, $r\in \left[0,\frac{R_0}{4}\right]$, of last inequality gives 
\begin{equation*}
\int_{B_{2r}}u^2dx\leq \frac{2\int_{B_{R_0}}u^2dx}{\int_{B_{R_0/4}}u^2dx} \int_{B_{r}}u^2dx,\quad \forall r\in \left[0,\frac{R_0}{4}\right].
\end{equation*}
Therefore, we get \eqref{Misc:6-12-22-2-1}. $\blacksquare$

\bigskip

\textbf{Final Remark.} The proof of Proposition \ref{N-cresce}, based on the Rellich identity, differs from the proof given in \cite{Ga_Li} which is based on the transformation of the elliptic operator in polar coordinates. The proof based on the Rellich identity was given by \cite{KUK} for the second--order elliptic operators. $\blacklozenge$ 

\section{A brief review about the $A_p$ weights.}\label{Ap:7-12-22} 
In this Section we will provide a brief summary of the $A_p$ weights by referring to \cite{Co-Fe} and \cite[Chapter 4]{GC-RdF} for further reading.
In Section \ref{Diff-Leb:15-10-22-1} we have introduced the maximal function related to a function $f\in L^1\left(\mathbb{R}^n\right)$. Basically equivalent to it is the following definition. In what follows we will denote by $Q$ a closed cube whose sides are parallel to the axes. Let $f\in L^1\left(\mathbb{R}^n\right)$. We define the maximal Hardy-Littlewood function as
\begin{equation}\label{Ap:7-12-22-1}
\mathcal{M}(f)(x)=\sup\left\{\dashint_{Q}|f(y)|dy:\mbox{ } Q\ni x\right\},
\end{equation} 
where, we recall, 
$$\dashint_{Q}|f(y)|dx=\frac{1}{|Q|}\int_Q|f(y)|dy.$$
It is simple to check that 

\begin{equation}\label{Ap:7-12-22-2}
\frac{\omega_n}{2^{n}n}M(f)(x)\leq \mathcal{M}(f)(x)\leq  \omega_nn^{\frac{1}{2}n-1} M(f)(x), \quad \forall x\in \mathbb{R}^n,\end{equation} 
 where $M(f)$ is defined in \eqref{Diff-Leb:15-10-22-1}. From inequalities \eqref{Ap:7-12-22-2} one can prove for $\mathcal{M}(f)$ properties similar to those of $M(f)$ some of which have been proved or presented in Section \ref{Diff-Leb}. In particular, the following apply
 (compare with Lemma \ref{Leb:16-10-22-1-0} and \eqref{7-12-1922}, respectively) 
 
\begin{equation*}
	\left|E_t\right|\leq \frac{C_n}{t}\int_{\mathbb{R}^n}|f(x)|dx, \quad \forall f\in L^1\left(\mathbb{R}^n\right), \forall t>0
\end{equation*}
($C_n$ depends on $n$ only), where 
$$E_t=\left\{x\in \mathbb{R}^n: \mbox{ } M(f)(x)>t \right\}$$ 
and
\begin{equation*}
	\left\Vert M(f)\right\Vert_{L^p\left(\mathbb{R}^n\right)}\leq C_{n,p} \left\Vert f\right\Vert_{L^p\left(\mathbb{R}^n\right)},\quad\forall f\in L^p\left(\mathbb{R}^n\right).
\end{equation*}
($C_{n,p}$ depends on $n$ and $p$ only). 

\medskip

The $A_p$ weight were introduced by Mucknhoupt, \cite{Muck}, and by Coifman and Fefferman, \cite{Co-Fe} to answer to the following question:
 
Let $p\in (1,+\infty)$, determine all measurable and nonnegative functions $w$ such that 

\begin{equation}\label{Ap:7-12-22-3}
\int_{\mathbb{R}^n}|\mathcal{M}(f)(x)|^pw(x)dx\leq C\int_{\mathbb{R}^n}|f(x)|^pw(x)dx, \quad \forall f\in L^p\left(\mathbb{R}^n\right),
\end{equation}
where $C$ depends on $w$ only. The functions $w$ which enjos property \eqref{Ap:7-12-22-3} are called \textbf{$A_p$ weight} \index{weight $A_p$}. If $p=1$, then we say that $w$ is an $A_1$ weight if 
\begin{equation*}
	\int_{\widetilde{E}_t}w(x)dx\leq C\int_{\mathbb{R}^n}|f(x)|w(x)dx, \quad \forall f\in L^1\left(\mathbb{R}^n\right),
\end{equation*} 
where 
$$\widetilde{E}_t=\left\{x\in \mathbb{R}^n: \mbox{ } \mathcal{M}(f)(x)>t \right\}.$$

The following Theorem can be proved  (\cite[Chap. 4, Sect 2]{GC-RdF}).

\begin{theo}\label{Ap:7-12-22-4}
	Let $w\in L^1_{loc}\left(\mathbb{R}^n\right)$, $w\geq 0$. The following conditions are equivalent
	
	\smallskip
	
	(a) there exists $p\in [1,+\infty)$ such that $w\in A_p$;
	
	 \smallskip
	 
	 (b) if $p\in (1,+\infty)$, then there exists $C>0$ such that
	 \begin{equation}\label{Ap:7-12-22-5}
	 	\left(\dashint_{Q}w(x)dx\right)\left(\dashint_{Q}w^{-\frac{1}{p-1}}(x)dx\right)^{p-1}\leq C,\quad \mbox{for every cube } Q
 		 \end{equation} 
	  
	  if $p=1$, then we have
	\begin{equation}\label{Ap:7-12-22-6}
		\left(\dashint_{Q}w(x)dx\right)\mbox{ess}\sup \left(w^{-1}\right)\leq C, \quad \mbox{for every cube } Q;
	\end{equation} 

\smallskip

(c) there exist $\delta>0$ and $C>0$  such that

\begin{equation}\label{Ap:7-12-22-7}
	\left(\dashint_{Q}w^{1+\delta}(x)dx\right)^{\frac{1}{1+\delta}}\leq C\dashint_{Q}w(x)dx\quad \mbox{for every cube } Q;
\end{equation} 
\smallskip

(d) there exist $s>0$ and $C>0$  such that, for every cube $Q$ and for every $E\subset Q$, $E$ Lebesgue measurable set, we have

\begin{equation}\label{Ap:7-12-22-8}
\frac{|E|}{|Q|}\leq C\left(\frac{\int_Ew(x)dx}{\int_Qw(x)dx}\right)^s.
\end{equation} 
\end{theo}
 
 \bigskip
 
 \textbf{Comments and Remarks.} An inequality like \eqref{Ap:7-12-22-7} is called \textbf{"reverse H\"{o}lder inequality"}. It can be proved that Theorem \ref{Ap:7-12-22-4} can be reformulated for weight functions $w$ defined on an open set $\Omega$, as long as we replace, in (b), (c) and (d) "for every cube $Q$" by "for every cube $Q\subset \Omega$.
The following remark turns out to be very useful in proving quantitative estimates in inverse problems. Below we illustrate this idea in a simplified way. 
 In Theorem \ref{theo:40.teo-fin} we proved that a nonzero if $u$ is a solution of the elliptic equation \index{equation:@{equation:}!- elliptic@{- elliptic}}
 
 \begin{equation}\label{Misc:8-12-22-2}  
 	\sum_{i,j=1}^na^{ij}(x)\partial^2_{x^ix^j}u+\sum_{i=1}^nb^i(x)\partial_{x^i}u+c(x)u=0,
 \end{equation}
 whose coefficients satisfy the hypotheses \eqref{4-2A-fin}, \eqref{intro14-4An-fin} and \eqref{limitaz-coeff-fin} then $u$ satisfies doubling inequality \eqref{eq:10.6.1102-cube-fin}. Let us suppose that $u$ satisfies \eqref{Misc:8-12-22-2} in $B_{\overline{R}}$, $\overline{R}>1$ and let us write inequality  \eqref{eq:10.6.1102-cube-fin} in the form
    
 \begin{equation}\label{Misc:8-12-22-3}
 	\int_{B_{2r}(x_0)}u^2\leq C N_{x_0}^k\int_{B_{r}(x_0)}u^2,\quad x_0\in B_1.
 \end{equation}
 
 \begin{equation}
 	\label{Misc:8-12-22-4}
 	N_{x_0}=\frac{\int_{B_{1}(x_0)}u^2dx}{\int_{B_{r_0}(x_0)}u^2dx}.
 \end{equation}
 where $r_0$ is a suitable point of $(0,1)$ (here, $R_0=1$ and $r_0=1/C$, in \eqref{eq:10.6.1102-cube-fin}). Let us assume for simplicity
 $$\overline{R}=\max\{1+2\sqrt{n}, 16\}.$$ 
 Let us denote by 
 \begin{equation}\label{Misc:8-12-1922-5}
 	F=\frac{\int_{B_{\overline{R}}}u^2dx}{\int_{B_{r_0/2}}u^2dx}.
 \end{equation}
By the Caccioppoli inequality and the Sobolev Embedding Theorem, proceeding similarly to what was done in Lemma \ref{imm-Sobol}, we have

 \begin{equation}\label{Misc:8-12-22-5}
 	\left(\dashint_{Q_r(x_0)}|u|^{q}dx\right)^{1/q}\leq \left(\dashint_{Q_{2r}(x_0)}u^{2}dx\right)^{1/2}, \quad \forall x_0\in B_{r_0/2},
 \end{equation}
  where $q$ is an arbitrary number of $(1,+\infty)$ when $n=2$, and it is equal to $\frac{2n}{n-2}$ when $n\geq 3$ . On the other hand (see Corollary \ref{SUCP-fin}), 
  \begin{equation*}
  	\begin{aligned}
  \left(\dashint_{Q_{2r}(x_0)}u^{2}dx\right)^{1/2}\leq& \left(\dashint_{B_{2\sqrt{n}r}(x_0)}u^{2}dx\right)^{1/2}\leq\\&\leq \left(CN_{x_0}\right)^{\widetilde{k}} \left(\dashint_{B_{r}(x_0)}u^{2}dx\right)^{1/2}\leq \\&\leq \left(CN_{x_0}\right)^{\widetilde{k}} \left(\dashint_{Q_{r}(x_0)}u^{2}dx\right)^{1/2}, \quad \forall x_0\in B_{r_0/2},
  	\end{aligned}
  \end{equation*}
where $\widetilde{k}=\frac{k}{2}\left(1+\log_2(2\sqrt{2n})\right)$. By the just obtained inequality and by \eqref{Misc:8-12-22-5}, we have
\begin{equation*}
	\left(\dashint_{Q_r(x_0)}|u|^{q}dx\right)^{1/q}\leq \left(CN_{x_0}\right)^{\widetilde{k}} \left(\dashint_{Q_{r}(x_0)}u^{2}dx\right)^{1/2}, \quad \forall  x_0\in B_{r_0/2}. 
\end{equation*}
 Moreover, by \eqref{Misc:8-12-1922-5} we have trivially 
 \begin{equation*}
 	N_{x_0}=\frac{\int_{B_{1}(x_0)}u^2dx}{\int_{B_{r_0}(x_0)}u^2dx}\leq \frac{\int_{B_{\overline{R}}}u^2dx}{\int_{B_{r_0/2}}u^2dx}.  
 \end{equation*}
 Therefore we have, in particular, for any cube  $Q\subset B_{r_0/2}$
 \begin{equation}\label{Misc:8-12-22-6}
 	\left(\dashint_{Q}|u|^{q}dx\right)^{1/q}\leq F^{\widetilde{k}} \left(\dashint_{Q}u^{2}dx\right)^{1/2}.
 \end{equation}
 Recalling that $q>2$, we get by  \eqref{Misc:8-12-22-6} that $u^2$ satisfies a reverse H\"{o}lder inequality, consequently  $u^2$ is an $A_p$ weight. In particular, Theorem  \ref{Ap:7-12-22-4} yields that  for every $Q\subset B_{r_0/2}$ and for every $E\subset Q$, $E$ Lebesgue measurabile set, we have
 
 \begin{equation}\label{Misc:8-12-22-7}
	\frac{|E|}{|Q|}\leq C\left(\frac{\int_Eu^2dx}{\int_Qu^2dx}\right)^s, 
\end{equation} 
where $C$ depends by $F$.
Let us suppose, now that the set $E$ has positive measure, then, if we have some bounds on
 $F$ (generally obtainable from values at the boundary of $u$), estimate \eqref{Misc:8-12-22-7} can be trivially rewritten

 \begin{equation}\label{Misc:8-12-22-8}
	\int_Qu^2dx\leq \left(\frac{C|Q|}{|E|}\right)^{1/s}\int_Eu^2dx
\end{equation} 

This estimate implies, in particular, that if $u$ vanishes on a set $E$ of positive measure then $u$ vanishes identically in $B_{\overline{R}}$. Actually, \eqref{Misc:8-12-22-8} also allows us to control, in terms of the measure $E$ only, the propagation of the error
 $$\int_Eu^2dx\leq \varepsilon^2$$
 On a "small" cube $Q$ and from there on the whole
  $B_{\overline{R}}$. 
 
 Another trivial translation of \eqref{Misc:8-12-22-7} is 
 \begin{equation*}
 	|E|\leq C|Q|\left(\frac{\int_Eu^2dx}{\int_Qu^2dx}\right)^s, 
 \end{equation*} 
which allows us to estimate the Lebesgue measure of $E$ by the integral of $u^2$ on $E$ itself. The latter observation is useful for finding \textbf{size estimates} \index{size estimates} of unknown inclusions in problems of the type considered in Section \ref{size-est} (for details, see \cite{A-R-M}).
$\blacklozenge$

\section{The Runge property}\label{Rung:9-12-22}
As an introduction to the main topic of this Section, we show by an example that it is not always possible to extend a harmonic function $u$ from $B_1$ to an open set, $\Omega$,  such that $B_1\Subset \Omega$
The example we present here is due to Hadamard, \cite{Ha1}.  
Let us consider the function whose expression in polar coordinates is given by
\begin{equation}\label{Rung:9-12-22-1}
	u(\rho,\phi)=\sum_{n=1}^{\infty}2^{-n}\rho^{4^n}\sin (4^n\phi).
\end{equation}
It is simple to check that $u\in C^0\left(\overline{B}_1\right)\cap C^2\left(B_1\right)$ and that $u$ is harmonic in $B_1$. Now we check that $$u\notin C^1\left(\overline{B}_1\right).$$ For this purpose we show that 
\begin{equation}\label{Rung:10-12-22-1}
	\lim_{r\rightarrow 1^-}\int_{B_r}|\nabla u|^2dxdy=\lim_{r\rightarrow 1^-}\int^r_0d\rho\int^{2\pi}_0\left(u^2_{\rho}+\rho^{-2}u^2_{\phi}\right)\rho d\phi=+\infty.
\end{equation}
We get

\begin{equation*}
	 \int^{2\pi}_0\left(u^2_{\rho}+\rho^{-2}u^2_{\phi}\right)\rho d\phi=2\pi \sum_{n=1}^{\infty}4^n\rho^{2\cdot 4^n-1}
\end{equation*}
Hence 
\begin{equation*}
	\int^r_0d\rho\int^{2\pi}_0\left(u^2_{\rho}+\rho^{-2}u^2_{\phi}\right)\rho d\phi=\pi \sum_{n=1}^{\infty} r^{2\cdot 4^n},
\end{equation*}
from which \eqref{Rung:10-12-22-1} follows. It is therefore evident that $u$ cannot be extended to a harmonic function in an open set containing $\overline{B}_1$. 

However, it is of interest to know whether $u$ can be approximated by functions that are harmonic in an open set containing $\overline{B}_1$. A property of this kind is called \textbf{Runge property} \index{Runge property}for the operator $\Delta$. This issue has been studied for operators which are more general than the Laplace operator, but here we limit ourselves to the Laplace operator operator only, referring to the final comments for hints on further consideration. 

\bigskip

We have

\begin{theo}\label{Runge:12-12-22-1}
Let $\Omega_1\Subset \Omega_2$ be two open sets of $\mathbb{R}^n$, where $\partial \Omega_2$ is of class $C^{1,1}$. Let us assume that $\Omega_2\setminus\overline{\Omega}_1$ connected. Then for every $u$ such that
\begin{equation*}
	\Delta u=0, \quad \mbox{in}\quad \Omega_1
\end{equation*}   
and for every $\varepsilon>0$ there exists $v\in H^1(\Omega_2)$ such that
\begin{equation*}
	\Delta v=0, \quad \mbox{in}\quad \Omega_2
\end{equation*}   
and
\begin{equation*}
	\left\Vert u-v\right\Vert_{L^2\left(\Omega_1\right)}<\varepsilon.
\end{equation*}     
\end{theo}    
\textbf{Proof.} Set
\begin{equation*}
	\mathcal{S}_1=\left\{u\in H^1\left(\Omega_1\right):\mbox{ } \Delta u=0 \mbox{ in } \Omega_1\right\}
\end{equation*}
and
\begin{equation}\label{Runge:13-12-22-1}
\mathcal{S}_2=\left\{v_{|_{\Omega_1}}:\mbox{ } \Delta v=0 \mbox{ in } \Omega_2\right\}.
\end{equation}
The property that we wish to prove is equivalent to the fact that $\mathcal{S}_2$ is dense in $\mathcal{S}_1$, with respect to the topology induced by $L^2\left(\Omega_1\right)$. We now prove this density property. To this aim, it suffices to prove that if $u\in \mathcal{S}_1$ and
\begin{equation}\label{Runge:14-12-22-1}
	\int_{\Omega_1}uvdx=0,\quad \forall v\in \mathcal{S}_2,
\end{equation}
then
$$u\equiv 0, \quad \mbox{in}\quad \Omega_1.$$ Let $\widetilde{u}$ 

\begin{equation*}
\widetilde{u}(x)=	\begin{cases}
		u,\quad \mbox{in}\quad \Omega_1,\\
		\\
		0, \ \ \mbox{in } \Omega_2\setminus \Omega_1.
	\end{cases}
\end{equation*} 
and let $w$ satisfy
\begin{equation}\label{Runge:14-12-22-2}
	\begin{cases}
		\Delta w= \widetilde{u},\quad \mbox{in}\quad \Omega_2,\\
		\\
		w\in H^1_0\left(\Omega_2\right).
	\end{cases}
\end{equation} 
Since $\widetilde{u}\in L^2\left(\Omega_2\right)$ and $\partial\Omega_2$ is of class $C^{1,1}$, by Theorem \ref{reg-bordo} we have $$w\in H^2\left(\Omega_2\right).$$ Now for an arbitrary $v\in \mathcal{S}_2$ we have
\begin{equation}\label{Runge:14-12-22-3-0}
	\begin{aligned}
		&0=\int_{\Omega_1}uvdx=\int_{\Omega_2}\widetilde{u} vdx=\int_{\Omega_2}\Delta w vdx=\\&=
	\int_{\partial \Omega_2}\frac{\partial w}{\partial \nu}vdS+\int_{\Omega_2}w\Delta vdx=\\&=\int_{\partial \Omega_2}\frac{\partial w}{\partial \nu}vdS.
	\end{aligned}
\end{equation}
Hence
\begin{equation*}
	\int_{\partial \Omega_2}\frac{\partial w}{\partial \nu}vdS=0,\quad \forall v\in \mathcal{S}_2.
\end{equation*}

By the latter and by Theorem \ref{teorema-es-unic} we have

\begin{equation}\label{Runge:14-12-22-3-1}
	\int_{\partial \Omega_2}\frac{\partial w}{\partial \nu}\varphi dS=0,\quad \forall \varphi\in H^{1/2}\left(\partial\Omega_2\right),
\end{equation}
hence for $\varphi = \frac{\partial w}{\partial \nu}$ we have
\begin{equation*}
	\frac{\partial w}{\partial \nu}=0,\quad \mbox{on}\quad \partial\Omega_2.
\end{equation*}
Therefore we have

\begin{equation}\label{Runge:14-12-22-3}
	\begin{cases}
		\Delta w= \widetilde{u},\quad \mbox{in}\quad \Omega_2\setminus\overline{\Omega}_1,\\
		\\
		w=0, \quad \mbox{on}\quad \partial\Omega_2,\\
		\\
		\frac{\partial w}{\partial \nu}=0, \quad \mbox{on}\quad \partial\Omega_2.
	\end{cases}
\end{equation} 
Now, since $\Omega_2\setminus\overline{\Omega}_1$ is connected, \eqref{Runge:14-12-22-3} implies, by the unique continuation property, 

\begin{equation*}
	w\equiv 0, \quad \mbox{in}\quad \Omega_2\setminus\overline{\Omega}_1.
\end{equation*}
By this relationship, recalling that $w\in H^2\left(\Omega_2\right)$ we have that $w\in H^2_0\left(\Omega_1\right)$ from which we have
\begin{equation*}
\int_{\Omega_1}u^2dx=\int_{\Omega_1}u\Delta wdx=\int_{\Omega_1}\Delta u wdx=0
\end{equation*}
Hence
$$u\equiv 0,\quad \mbox{in}\quad \Omega_1.$$
As we desired to prove. $\blacksquare$

\bigskip

\textbf{Remarks.} 

\smallskip

\textbf{1.} Theorem \ref{Runge:12-12-22-1} continues to be true also if, instead of $\partial\Omega_2\in C^{1,1}$, we assume only $\partial\Omega_2\in C^{0,1}$
To prove this, one must first keep in mind that in \eqref{Runge:14-12-22-3-0} occurs $$\langle \frac{\partial w}{\partial \nu}, v \rangle_{H^{-1/2}(\partial \Omega_2), H^{1/2}(\partial \Omega_2)}$$ instead of $$\int_{\partial \Omega_2}\frac{\partial w}{\partial \nu}vdS$$ and, consequently, instead of \eqref{Runge:14-12-22-3-1}, we have

\begin{equation*}
\langle \frac{\partial w}{\partial \nu}, \varphi \rangle_{H^{-1/2}(\partial \Omega_2), H^{1/2}(\partial \Omega_2)}=0,\quad\forall \varphi \in H^{1/2}\left(\partial\Omega_2\right).
\end{equation*}
Therefore, we likewise have that $\frac{\partial w}{\partial \nu}=0$. However, since $\frac{\partial w}{\partial \nu} \in H^{-1/2}(\partial \Omega_2)$, it will be necessary to first reformulate Cauchy problem \eqref{Runge:14-12-22-3} in a weak form and then to prove the uniqueness for such a Cauchy problem, for both of which we refer the reader to the paper \cite{A-R-R-V}.

\smallskip

\textbf{2.} The assumption that $\Omega_2\setminus\overline{\Omega}_1$ is connected cannot be dropped as the following simple \textbf{counterexample} shows.
 Let $$\Omega_1=\left(B_7\setminus \overline{B}_5\right)\cup \left(B_3\setminus \overline{B}_1\right),\quad\quad \Omega_2=B_8,$$ let us notice that $\Omega_2\setminus\overline{\Omega}_1$ is not  connected. 

Let $u$ be the following function

\begin{equation}\label{Runge:15-12-22-1}
	u(x)=\begin{cases}
	1,\quad \mbox{for } x\in B_3\setminus \overline{B}_1,\\
	\\
	0, \quad \mbox{for } x\in B_7\setminus \overline{B}_5,
\end{cases}
\end{equation}
 of course $u$ is harmonic in $\Omega_1$. Let $\varepsilon$ be a given positive number which we will choose later and let $v\in H^1\left(\Omega_2\right)$ be a harmonic function satisfying
 
 \begin{equation}\label{Runge:15-12-22-2}
 	\left\Vert u-v\right\Vert_{L^2\left(\Omega_1\right)}<\varepsilon.
 \end{equation} 
Let $x_0\in \partial B_2$, taking into account that $B_1(x_0)\subset \Omega_1$, by \eqref{Runge:15-12-22-2} we have

\begin{equation*}
	\begin{aligned}
	|v(x_0)-1|&=\left|\frac{1}{\left|B_1(x_0)\right|}\int_{B_1(x_0)}(v-1)dx\right|\leq\\&\leq  \frac{1}{\left|B_1(x_0)\right|^{1/2}} \left(\int_{B_1(x_0)}(v-1)^2dx\right)^{1/2}\leq\\&\leq	\frac{\varepsilon}{c_n},
\end{aligned}
\end{equation*}
where $c_n=|B_1|^{1/2}$. Hence we have
\begin{equation}\label{Runge:16-12-22-1}
	1-\frac{\varepsilon}{c_n}<v(x)<1+\frac{\varepsilon}{c_n}, \quad \forall x\in \partial B_2 .
\end{equation}
Similarly we get
\begin{equation}\label{Runge:16-12-22-2}
	-\frac{\varepsilon}{c_n}<v(x)<\frac{\varepsilon}{c_n}, \quad \forall x\in \partial B_6 .
\end{equation}
Now, by \eqref{Runge:16-12-22-1}, \eqref{Runge:16-12-22-2} and by using maximum principle we have

$$1-\frac{\varepsilon}{c_n}<v(x)\leq \max_{x\in \partial B_6} v< \frac{\varepsilon}{c_n}, \quad \forall x\in \partial B_6,$$ that leads to a contradiction provided $\varepsilon<\frac{c_n}{2}$. 

\smallskip

\textbf{3.} Theorem \ref{Runge:12-12-22-1} can be extended to a large class of operators; in particular, it can be extended to the second-order elliptic operators with real coefficients whose formal adjoint enjoys the unique continuation property.
 Let us specify this a little. Let us consider, for instance, the operator 
\begin{equation}
	Lu=\mbox{div}\left(A(x)\nabla u\right),
\end{equation}  
where $A\in L^{\infty}\left(\mathbb{R}^n;\mathbb{M}(n)\right)$ is a not necessarily symmetric matrix and such that ($\lambda\geq 1$)
 \begin{equation*}
 	\lambda^{-1}|\xi|^2\leq A(x)\xi\cdot \xi,\quad \forall \xi\in \mathbb{R}^n, \mbox{ } \forall x\in \mathbb{R}^n.
 \end{equation*}  
 Then the formal adjoint of $L$ is the operator
 
 \begin{equation}
 	L^{\star}u=\mbox{div}\left(A^T(x)\nabla u\right).
 \end{equation}  
 We say that $L$ enjoys  the \textbf{Runge property} provided it occurs what follows. 
 
 Let $\Omega_1\Subset \Omega_2$ be two open sets of $\mathbb{R}^n$, like Theorem \ref{Runge:12-12-22-1}, then for every $u$ such that
  \begin{equation*}
 	L u=0, \quad \mbox{in}\quad \Omega_1
 \end{equation*}   
 and for every $\varepsilon>0$ there exists  $v\in H^1(\Omega_2)$ such that
 \begin{equation*}
 	L v=0, \quad \mbox{in}\quad \Omega_2
 \end{equation*}   
 and
 \begin{equation*}
 	\left\Vert u-v\right\Vert_{L^2\left(\Omega_1\right)}<\varepsilon.
 \end{equation*}      
 Lax in \cite{Lax} has proved that the following conditions are equivalent
 
 \smallskip
 
 (a) $L$ enjoys Runge property 
 
 \smallskip
 
 (b) $L^{\star}$ enjoys the unique continuation property.
 
 \smallskip
 
The proof of (b) $\Longrightarrow$ (a) is analogous to the proof of Theorem 
\ref{Runge:12-12-22-1} and it s left to the reader as an exercise. Concerning the implication (a) $\Longrightarrow$ (b), we refer to  \cite{Lax}.

\smallskip

\textbf{4.} The quantitative versions of the Runge property are also of interest (especially in the stability issue of inverse problems). That is, it is of interest to estimate appropriately from above in terms of $\varepsilon$, the quantity 
$$\left\Vert v\right\Vert_{H^{1/2}\left(\partial\Omega_2\right)},$$
(by Hadamard example, illustrated at the beginning of this Section, we should expect that, in general, as $\varepsilon$ goes to $0$ we should have
 $\left\Vert v\right \Vert_{H^{1/2}\left(\partial\Omega_2\right)}$ goes to infinity). A result in this regard is proven in \cite{Ru-Sa}. $\blacklozenge$

\printindex

\bigskip


\begin{thebibliography}{99}

\bibitem{Ag} S. Agmon,  Unicit\'{e} et Convexit\'{e} dans les Probl\`{e}mes Diff\'{e}rentiels (S\'{e}minaire de Math\'{e}matiques Sup\'{e}rieures, No. 13 (\'{E}t\'{e}, 1965). Les Presses de l'Universit\'{e} de Montr\'{e}al, Montreal, Que), 1966.

\bibitem{A-R-R-V} G. Alessandrini, L. Rondi, E. Rosset,
S. Vessella, The stability for the Cauchy problem for elliptic
equations, Inverse Problems 25 (2009), 1-47.

\bibitem{A-R-M} G. Alessandrini, A. Morassi, E. Rosset, Size estimates, In "Inverse problems: theory and
applications", Contemp. Math. 333 (2003), 1--33.

\bibitem{A-R-M-08} G. Alessandrini, A. Morassi, E. Rosset, The linear constraint in the Poincar\'{e} and Korn type inequalities, Forum Math., 20 (2008), pp. 557 -- 569.

\bibitem{ALM} F. I. Almgren Jr., Dirichlet's problem for multiple valued functions and the regularity of mass minimizing integral currents, in "Minimal Submanifolds and Gedesic", (M. Obata Ed) pp. 1--6, North--Holland, Amsterdam, (1979).

\bibitem{Ar}  N. Aronszajn, A unique continuation theorem for solutions of elliptic partial differential equations or inequalities of second order, J. Math. Pures Appl. (9) 36 (1957), 235–249.

\bibitem{AKS}  N. Aronszajn, A. Krzywicki and J. Szarski, A Unique
Continuation Theorem for Exterior Differential Forms on Riemannian
manifolds, Ark. for Matematik, 4, (34), 1962, 417-453.

\bibitem{l:Bakri1} L. Bakri,  Quantitative uniqueness for Schr\"odinger operator, Indiana Univ. Math. J. 61 (2012), no. 4, 1565--1580.

\bibitem{l:Bakri2} L. Bakri, Carleman estimates for the Schr\"odinger operator. Applications to quantitative uniqueness, Comm. Partial Differential Equations 38 (2013), no. 1, 69--91.

\bibitem{Bal} L. Baldassari, Stime di Carleman per l'operatore di
Laplace e applicazioni allo studio della continuazione unica e della
stabilit\`{a}, Tesi di Laurea Magistrale in Matematica , Marzo 2017.

\bibitem{BOO}  W. M. Boothby, An Introduction to Differentiable Manifolds
and Riemannian Geometry. Academic Press, New York, 1975.

\bibitem{Br} H. Brezis, Analisi funzionale, Ed. Liguori, Napoli,
1986.


\bibitem{Car0} T. Carleman, Les Fonctions quasi-analitiques, Gauthier-Villars, Paris (1926), pp. 3--5.

\bibitem{Car1} T. Carleman, Sur les syst\`{e}mes lin\'{e}aires aux d\`{e}riv\`{e}es
partielles du primier ordre \`{a} deux variables, C. R. Acad. Sci.
Paris, 19, (1933), 471--474.

\bibitem{Car2} T. Carleman, Sur un probl\`{e}me d'unicit\'{e} pour
les syst\`{e}mes d'\`{e}quations aux d\`{e}riv\`{e}es partielles
\`{a} deux variables ind\'{e}pendentes, Ark. Mat. Astr. Fys., 26B
(1939), 1--9.

\bibitem{CTV} V. Checcucci, A. Tognoli, E. Vesentini, Lezioni di topologia generale, Feltrinelli, Milano, 1968.

\bibitem{Co-Fe} R. Coifman and C. Fefferman. Weghted norm inequalities for maximal functions and singular
integrals. Studia Math. 5L (1974),24I 250.

\bibitem{Co-Hi} R. Courant and D. Hilbert, Methods of Mathematical Phisics, Vol. II, John Wiley \& Sons Inc. New York 1989.

\bibitem{DeFi-Go} D. G. de Figueredo, J. P. Gossez, Strict
monotonicity of eigenvalues and unique continuation, Comm. Part.
Diff. Equat. 17 (1, 2), 339--346 (1992).

\bibitem{Du} P. L. Duren, Theory of $H^p$ spaces, Academic Press,
New York, 1970.


\bibitem{Fo} G. Folland, Introduction to Partial Differential Equations, Princeton University Press, 1995.

\bibitem{E-V} L. Escauriaza and S. Vessella,  Optimal Three Cylinder Inequalities for Solutions to Parabolic Equations with Lipschitz Leading Coefficients (Inverse Problems Theory and Applications Contemporary Mathematics vol 333) ed G. Alessandrini and G. Uhlmann (Providence, RI: American  Mathematical Society) pp 79--87, (2003).

\bibitem{EV}  L. C. Evans, Partial Differential Equations, AMS, Providence,
R. I., 1998.


\bibitem{EV-GA} L. C. Evans and R. E. Gariepy, Measure Theory and Fine Properties of Functions, CRC
Press, Boca Raton, FL, 1992.

\bibitem{GC-RdF} J. Garc\`{i}a-Cuerva and J. L. Rubio De Francia. Weighted norm inequalities and related topics, Amsterdam: North-Holland, 1985.

\bibitem{Garding} L. Garding, Some Points in Analysis ant Their Hystory, University Lectures Series, vol. 11, American Mathematical Society, 1998.


\bibitem{Ga_Li} N. Garofalo and F. H. Lin,  Monotonicity properties of variational
integrals, Ap weights and unique continuation Indiana Univ. Math. J.
35, 245 -- 268, (1986).



\bibitem{G-T} D. Gilbarg and N. Trudinger, Elliptic Partial Differential Equations of Second Order, (Grundlehren der
Mathematischen Wissenschaften vol 224) 2nd edn (Berlin: Springer),
1983.

\bibitem{Giusti} E. Giusti, Equazioni ellittiche del secondo ordine. Quaderni dell'Unione Matematica Italiana, Ed. Pitagora, Bologna 1978.  


\bibitem{Gorny} A. Gorny, Contribution a l'\`{e}tude des fonctions d\'{e}rivables d'une variable r\'{e}elle. Acta Math., 71 (1939), pp. 257--268.

\bibitem{Ha1} J. Hadamard, Sur le principe de Dirichlet,
Bulletin de la S. M. F., tome 34 (1906), p. 135-138.

\bibitem{Ha} J. Hadamard, Le probl\`{e}me de Cauchy et les \'{e}quations aux d\'{e}riv\'{e}es partielles lin\'{e}aires hyperboliques, Herman Paris, 1932


\bibitem{l:HLP34}
G.H. Hardy, J.E. Littlewood and G. P\'olya, Inequalities.
Cambridge University Press, Cambridge, 1934.



\bibitem{HO63}  L. H\"{o}rmander, Linear Partial Differential Operators,
Springer, New York, 1963.

\bibitem{HO}  L. H\"{o}rmander, Uniqueness Theorem for Second Order Elliptic
Differential Equations, Comm. Part. Diff. Equations, 8 (1),1983,
21-64.

\bibitem{HOII} L. H\"{o}rmander, The Analysis of Linear Partial Differential Operators, vol. I-IV, Springer, New York, 1985.

\bibitem{KL-TI} M. V. Klibanov and A. Timonov, Carleman Estimates for
Coefficient Inverse Problems and Numerical Applications, VSP,
Utrecht--Boston, 2004.

\bibitem{Is1} V. Isakov, Inverse Source Problems (Mathematical Surveys and Monographs vol 34) (Providence, RI: American Mathematical Society), 1990.

\bibitem{Is2} V. Isakov, Inverse Problems for Partial Differential Equations (Applied Mathematical Sciences vol 127) 2nd edn (New York: Springer), 2006.

\bibitem{Joh60} F. John, Continuous dependence on data for solutions of partial differential equations with a prescribed bound Commun. Pure Appl. Math., VOL. XIII, 551-586 (1960).

\bibitem{Joh}  F. John, Partial Differential Equations. Springer-Verlag, New York, 1982.

\bibitem{KUK} I. Kukavica, Quantitative uniqueness for second order
elliptic operators, Duke Math. J. 91, 1998, 225--240.

\bibitem{K-J-F} A. Kufner, O. John, S. Fu\v{c}ik, Function spaces. Noordhoff International Publishing, Leyden; Academia, Prague,
1977.

\bibitem{Kur} C. Kuratowski, Topologie. Polska Akademia Nauk, Warzawa 1958.

\bibitem{LAN}  E. M. Landis, A three sphere theorem, Dokl. Akad. Nauk SSSR
148 (1963), 277-279, Engl. trans. Soviet Math. Dokl. 4 (1963),
76--78.

\bibitem{Lang} S. Lang, Real Analysis, 2nd edition, Addison--Wesley, London, 1983.


\bibitem{Lav} M. M. Lavrent'ev, On the Cauchy problem for Laplace equation Izv. Akad. Nauk SSSR. Ser. Mat.
20 819--42, 1956.

\bibitem{L-R-S} M. M. Lavrent'ev, V. G. Romanov V G and S. P. \v{S}i\v{s}atski\u{\i},  Nekorrektnye Zadachi Matematicheskoi Fiziki i
Analiza, (Moscow: Nauka), 1980.

M. M. Lavrent'ev, V. G. Romanov V G and S. P.
\v{S}i\v{s}atski\u{\i}, Problemi non ben posti in fisica matematica
e analisi Pubblicazioni dell'Istituto di Analisi Globale e
Applicazioni 12 (Firenze: IAGA) (Italian Translation), 1983

M. M. Lavrent'ev, V. G. Romanov V G and S. P.
\v{S}i\v{s}atski\u{\i},  Ill-Posed Problems of Mathematical Physics
and Analysis (Translations of Mathematical Monographs vol 64)
(Providence, RI: American Mathematical Society) (Engl. Transl.),
1986.


\bibitem{Lax} P. D. Lax, A stability theorem for solutions of abstract differential equations and its application to the study of local behavior of solutions of elliptic equations, Comm. Pure Appl. Math. 9, (1956) 747-766.
 
\bibitem{Lern} N. Lerner, Carleman inequalities -- An introduction
and more, Springer, Berlin 2019.



\bibitem{Ma} F. Maggi, Sets of finite perimeter and geometric
variational problems -- an introduction to geometric measure theory,
Cambridge University Press, 2012.

\bibitem{Magn1} R. Magnanini, Dispense del Corso di Analisi III, Corso di Laurea Triennale in Matematica,
 
 https://web.math.unifi.it/users/magnanin/Istit/a3gsm15.pdf, 2015.

\bibitem{Magn2} R. Magnanini, Dispense del Corso di
Istituzioni di Analisi Superiore,
Laurea Magistrale in Matematica,  

https://web.math.unifi.it/users/magnanin/Istit/DispenseIAS2016.pdf, 2016.


\bibitem{Mand} N. Mandache N, On a counterexample concerning unique continuation
for elliptic equations in divergence form Mat. Fiz. Anal. Geom. 3,
(1996) 308-31.

\bibitem{Met} G. M\'{e}tevier, Counterexample to Holmgren's uniqueness theorem for analytic non--linear Cauchy problems, Inv. Math. 112 (1993), 217-222.

\bibitem{Mizo} S. Mizohata, Lectures on Cauchy Problem, Tata Institute of Fundamental Research, Bombay, 1965.

\bibitem{Mo-Ni} C. B. Morrey, L. Nirenberg, On the analyticity of the solutions of linear elliptic systems
of partial differential equations. Comm. Pure Appl. Math. 10 (1957),
271--290.

\bibitem{Muck} Muckenhoupt, Weighted norm inequalities for the Hardy maximal function, Trans. Amer. Math. Soc. 165 (1972) 207-226.

\bibitem{Ne} J. Ne\c{c}as,  Direct methods in the theory of elliptic equations, Springer-Verlag Berlin Heidelberg, 2012.

\bibitem{Ni} L. Nirenberg, Uniqueness in Cauchy Problems for
Differential Equations with Constant Leading Coefficients, Comm.
Pure Appl. Math., 10, (1957) 89-105.

\bibitem{Pa} L. E. Payne  Improperly Posed Problems in Partial Differential Equations (Philadelphia, PA: Society for Industrial and Applied Mathematics), 1975.

\bibitem{Petr} I. G. Petrovsky, Lectures on Partial Differential Equations. Interscience Publishers, New York, 1954.

\bibitem{Pi-St-Vi} L. C. Piccinini, G. Stampacchia, G. Vidossich,
Equazioni differenziali ordinarie in $\mathbb{R}^n$. Liguori
Editore, Napoli, 1978.

L. C. Piccinini, G. Stampacchia and G. Vidossich,
Ordinary differential equations in  $\mathbb{R}^n$. Springer-Verlag New York, 1984.


\bibitem{Pl} A. Pl\u{\i}s, On non-uniqueness in Cauchy problem for an elliptic
second order differential equation Bull. Acad. Pol. Sci. S\'{e}r.
Sci. Math. Astron. Phys. 11, (1963) 95 -- 100.


\bibitem{Pu} C. Pucci, Istituzioni di Analisi Superiore, Unione Matematica Italiana, Bologna, 2013


\bibitem{RA} Rakesh, Carleman estimates for second order PDEs, http://www.bcamath.org/documentos\_public/courses/3.\_Carleman\_estimates.\\ RAKESH.pdf


\bibitem{Rellich} F. Rellich, Darstellung der eigenwerte von
$\Delta u+ \lambda u=0$ durch ein Randitegral, Math. Zeit., 46,635
(1940) -- 636.

\bibitem{Roy} H. L. Royden, Real Analysis, Macmillan P. C., New York 1988. 


\bibitem{Ru} W. Rudin, Functional Analysis, McGraw--Hill Publishing Company, New York 1975.

\bibitem{Ru-Sa} A. R\"{u}land and M. Salo, Quantitative Runge Approximation and Inverse Problems, International Mathematics Research Notices, Vol. 2019, No. 20, pp. 6216--6234 Advance Access Publication January 19, 2018, doi:10.1093/imrn/rnx301.


\bibitem{ST} E. M. Stein, Singular integrals and differentiability properties of functions, Princeton Univ. Press,
Princeton, NJ, 1970.

\bibitem{l:T67} G. Talenti, Sopra una disuguaglianza integrale, Annali Scuola Normale Superiore Pisa, classe di Scienze, 21(2) (1967) 167--188.

\bibitem{Tik} A. N. Tikhonov and V. A. Arsenin, Solutions of Ill--Posed Problems, John Wiley and Sons, Toronto 1977.

\bibitem{Tik1} A. Tychonoff, Th\'{e}oremes d'unicit\'{e} pour l'\'{e}quation de la chaleur, Mat. Sb., 42, 199--216 (1935).


\bibitem{Tr1} F. Treves, Partial Differential Equations with Constants Coefficients, Gordon and Breach, New York, 1966.

\bibitem{Tr} F. Treves, Basic Linear Partial Differential Equations, Academic Press, 1975.

\bibitem{We} H. F. Weinberger, A First Course in Partial Differential Equations, Dover, 1995.

\bigskip



\end{thebibliography}
\end{document}